\documentclass[10pt]{amsbook} 
\usepackage[T1]{fontenc}
\usepackage{amsmath}
\usepackage{mathrsfs}
\usepackage{amssymb}
\usepackage{latexsym}
\usepackage{framed}
\usepackage{emptypage}
\usepackage{graphicx}
\usepackage{imakeidx}
\usepackage{framed}
\usepackage{tikz}
\usepackage[all]{xy}
\makeindex
\numberwithin{equation}{chapter} 

\newtheorem{theorem}{Theorem}[section]
\newtheorem{lemma}[theorem]{Lemma}
\theoremstyle{definition}
\newtheorem{definition}[theorem]{Definition}
\newtheorem{example}[theorem]{Example}

\newtheorem{observation}[theorem]{Observation}
\theoremstyle{remark}
\newtheorem{remark}[theorem]{Remark}
\theoremstyle{corollary}
\newtheorem{corollary}[theorem]{Corollary}
\theoremstyle{proposition}
\newtheorem{proposition}[theorem]{Proposition}

\newcommand{\im}{\operatorname{Im}}
\newcommand{\re}{\operatorname{Re}}

\newcommand{\Rn}{\mathbb{R}^n}
\newcommand{\R}{\mathbb{R}}
\newcommand{\Cn}{\mathbb{C}^n}
\newcommand{\C}{\mathbb{C}}

\newcommand{\Z}{\mathbb{Z}}

\newcommand{\N}{\mathbb{N}}

\newcommand{\abs}[1]{\left|#1\right|}
\newcommand{\norm}[1]{\left\|#1\right\|}

\usepackage{mathtools}

\begin{document}

\title{Theory of Polyanalytic functions}
\author{Abtin Daghighi} 
%
%
%

\frontmatter

\maketitle

\cleardoublepage 
\vspace*{\stretch{1}}
\begin{flushright}
\itshape
{\Large Dedicated to my parents and my brother}
\end{flushright}
\vspace{\stretch{3}}


\tableofcontents

\newpage

\mainmatter

\chapter*{Introduction}
There are different equivalent ways to define complex analytic functions. From the perspective of partial differential equations they can be identified as the
set of functions annihilated by the Cauchy-Riemann operator $\frac{\partial}{\partial \bar{z}}$.
Powers of the Cauchy-Riemann operator where present in relation to elasticity problems studied by Kolosov \cite{kolosov} (1909), whose studies involved so called bianalytic functions which are functions of the form $a(z)+\bar{z}b(z),$ where $a$ and $b$ are holomorphic. Before that Goursat \cite{goursat} had studied so called biharmonic functions (constituting the kernel of the square of the Laplace operator $\Delta^2$) in particular proving that any biharmonic function $u$ can be identified as the real part of what was later called bianalytic functions. Burgatti \cite{burgatti} and Theodorescu \cite{theodorescu} where among the first to initiate a general study of the kernels (by kernel here we mean nullspace) of the operators $\left(\frac{\partial}{\partial \bar{z}}\right)^q$ for positive integers $q$ (i.e.\
powers of the Cauchy-Riemann operator).
As it turns out, many properties of complex analytic functions (the case $q=1$) 
have similar counterparts in the higher order situation ($q>1$)
for example integral and power series representations, automatic real-analyticity of solutions to the defining equations and the structure of the sets of uniqueness. 
One way to characterize polyanalyticity in finite dimension is precisely by replacing the system
$\partial_{\bar{z}_j} f=0,$ $j=1,\ldots,n,$ for functions defined on open subsets of $\Cn$, by the system
$\partial_{\bar{z}_j}^{\alpha_j} f=0,$ $j=1,\ldots,n,$ for $\alpha\in \Z_+^n.$ Although, as we shall see
(see Section \ref{ahernbrunasec} of Chapter \ref{basicpropsec})
there is in finite higher dimensional complex analysis an alternative notion, due to Ahern \& Bruna \cite{ahernbruna},
which may be more natural and which conforms to 
the possibility of representing complex analytic functions in terms of homogeneous series (which in particular renders the property that the restriction to each complex slice is again complex analytic).
It has been about thirty years since a textbook (a comprehensive survey) appeared on the subject of polyanalytic functions and even then some developments that where already present (e.g.\ the generalization to several complex variables) where excluded in that survey, by Mark Balk \cite{ca1}. This book will present most of the main developments since then and also attempt to lay a groundwork for a fundamental understanding of the concept of polyanalyticity by approaching from different perspectives.
It may be appropriate to subcategorize the study of polyanalytic functions into two major categories:
one approaching the subject from classical function theory and the other from operator theory. Both will of course intersect the theory
of partial differential equations.
Although recently, the study of the higher power case has been mentioned also in the context of other more applied fields of research for instance time-frequency analysis and
wavelets
and eigenspaces of the so-called magnetic Schr\"odinger operator (see e.g.\ Askour \cite{askour} and Hedenmalm \cite{hedenhaimi}), such applications are intentionally not covered here. 
The focus of this book is not on any applied mathematics but strictly on pure mathematics.
One way the zeros of $q$-analytic functions arise naturally in complex analysis is 
in terms of {\em analytic arcs}\index{Analytic arc}.
An analytic arc, $\gamma$, is the image of an interval $[a,b]$ under a one-to-one
conformal map $A$ ($A$ is sometimes called the {\em Schwarz function of $\gamma$}).
If a point $z=A(t_0)$ 
lies in a neighborhood of the analytic arc, $\gamma$, then the reflection $z^*$ of $z$, in $\gamma$, 
is defined as $\overline{A(t_0)}.$ A point in the image of $\phi$, thus belongs to
$\gamma$ if and only if it equals its reflection with respect to $\gamma,$ i.e.\
if it satisfies
$z=\overline{A(t_0)}$ which is equivalent to
\begin{equation}
f(z):=B(z)+\bar{z}=0
\end{equation}
where $B(z)=-A(z)$ is a conformal one-to-one map.
Hence, an analytic arc is precisely the zero set of a $2$-analytic function (also called bianalytic function).
So for each $q\in \Z_+$ consider the set of (holomorphic) functions of the form
$F(z)=\sum_{j=0}^{q-1}a_j(z)(A(z))^j,$  for holomorphic $a_j(z),$ and a one-to-one conformal map $A$.
The restriction of the holomorphic function $F(z)$,  
to an analytic arc, is unique, but at each point of $\gamma=\{A(z)=\bar{z}\}$ the value of $F(z)$ coincides with that of the $q$-analytic
function
$\sum_{j=0}^{q-1}a_j(z)\bar{z}^j.$ Schwarz functions have been heavily utilized for proving results on the properties of boundary values of
$q$-analytic functions e.g.\ when considering the Dirichlet problem for bounded domains with rational boundary.
$\alpha$-analytic functions on domains in $\Cn$ can always be identified as the restriction to a complex $n$-dimensional submanifold, $M$ of $\C^{n+1},$ of so-called
{\em pseudopolynomials} on $\C^{n+1}$ (see Definition \ref{pseudopolynomdef}),
and this line of thinking opens up the door to hypoanalytic theory and 
$CR$ geometry, see Chapter \ref{hypoanalsec}.
An arbitrary linear higher order complex partial differential equation has the form
\begin{equation}\label{equation1}
\left(\frac{\partial}{\partial z}\right)^m\left(\frac{\partial}{\partial \bar{z}}\right)^n f+
\sum_{\abs{(\mu,\nu)}<m+n}\left(\frac{\partial}{\partial z}\right)^\mu\left(\frac{\partial}{\partial \bar{z}}\right)^\nu f =0
\end{equation}
Below we shall for convenience of notation
denote $\partial_z:=\frac{\partial}{\partial z}$,
 $\partial_{\bar{z}}:=\frac{\partial}{\partial \bar{z}}$.
Recall that a fundamental solution for the Cauchy-Riemann operator is $\frac{-1}{\pi z}$ and by iteration a fundamental solution for 
$\partial_{\bar{z}}^n$ is given by $-\frac{1}{\pi}\frac{\bar{z}^{n-1}}{(n-1)!z}.$
A fundamental solution for the Laplacian is $\log\abs{z}^2$ and
 a fundamental solution for $\partial_{z}^m\partial_{\bar{z}}^n$ is
\begin{equation}
-\frac{1}{\pi}\frac{z^{m-1}}{(m-1)!}\frac{\bar{z}^{n-1}}{(m-1)!\pi}
\left(\log\abs{z}^2-\sum_{\mu=1}^{m-1}\frac{1}{\mu}-
\sum_{\nu=1}^{n-1}\frac{1}{\nu}\right)\end{equation}
Generalizing a technique used by Vekua \cite{vekua}, 
Begehr \& Vaitekhovic \cite{begehrvaitekovic}, point out that Eqn.(\ref{equation1}) can on a plane domain $\Omega\subset\C$ be
rewritten as
\begin{multline}\label{equation2}
\partial_{z}^m\partial_{\bar{z}}^n f+\sum_{\stackrel{\mu+\nu=m+n}{(\mu,\nu)\neq(m,n)}}\left(q_{\mu\nu}\partial_z^\mu\partial_{\bar{z}}^\nu f +\hat{q}_{\mu\nu}\overline{\partial_z^\mu\partial_{\bar{z}}^\nu f}\right) +\\
\sum_{\stackrel{\mu+\nu<m+n}{(\mu,\nu)\neq(m,n)}}
\left(a_{\mu\nu}\partial_z^\mu\partial_{\bar{z}}^\nu f +b_{\mu\nu}\overline{\partial_z^\mu\partial_{\bar{z}}^\nu f}\right) +c=0
\end{multline}
where $a_{\mu\nu},b_{\mu\nu},c\in L^p(\Omega,\C), p>1,$
and $\sum_{\mu+\nu=m+n\\ (\mu,\nu)\neq(m,n)}(\abs{q_{\mu\nu}(z)}+\abs{\hat{q}_{\mu\nu}(z)})<1,$
and Eqn.(\ref{equation2})
can in turn via the representation 
$w:=\varphi+T_{m,n}\rho,$
$\partial_{z}^m\partial_{\bar{z}}^n \varphi=0,$
$\partial_{z}^m\partial_{\bar{z}}^n w=\rho,$ 
be transformed into the singular integral equation
\begin{multline}\label{equation3}
\rho+\sum_{\stackrel{\mu+\nu=m+n}{(\mu,\nu)\neq(m,n)}}\left(q_{\mu\nu}T_{m-\mu,n-\nu}\rho +\hat{q}_{\mu\nu}\overline{T_{m-\mu,n-\nu}\rho}\right) +\\
\sum_{\mu+\nu<m+n}\left(a_{\mu\nu}T_{m-\mu,n-\nu}\rho +b_{\mu\nu}\overline{T_{m-\mu,n-\nu}\rho}\right)+\\
\sum_{\stackrel{\mu+\nu=m+n}{(\mu,\nu)\neq(m,n)}}\left( q_{\mu\nu}\partial_z^\mu\partial_{\bar{z}}^\nu \varphi
+\hat{q}_{\mu\nu}\overline{\partial_z^\mu\partial_{\bar{z}}^\nu \varphi}\right)+\\
\sum_{\mu+\nu<m+n}\left(a_{\mu\nu}\partial_z^\mu\partial_{\bar{z}}^\nu \varphi+
b_{\mu\nu}\overline{\partial_z^\mu\partial_{\bar{z}}^\nu \varphi}
\right) +c=0
\end{multline}
where
\begin{equation}\label{pomieueq1}
T_{m,n}f(z):=\int_{\Omega}K_{m,n}(z-\zeta)f(\zeta)d\xi d\eta,\quad f\in L^1(\Omega,\C)
\end{equation}
$T_{0,0}f:=f,$ and
\begin{equation}\label{pomieueq2}
K_{m,n}(z)\colon=\left\{
\begin{array}{l}
\frac{(-1)^m (-m)!}{(n-1)!\pi}z^{m-1}\bar{z}^{n-1},\quad m\leq 0,\\
\frac{(-1)^m(-n)!}{(m-1)!\pi}z^{m-1}\bar{z}^{n-1},\quad n\leq 0,\\
\frac{z^{m-1}}{(m-1)!}\frac{\bar{z}^{n-1}}{(m-1)!\pi}
\left(\log\abs{z}^2-\sum_{\mu=1}^{m-1}\frac{1}{\mu}-
\sum_{\nu=1}^{n-1}\frac{1}{\nu}\right),\quad 0<m,n,
\end{array}
\right.
\end{equation}
for $0\leq m+n,0<m^2+n^2.$
The $K_{m,n}$ are called Cauchy-Poisson kernels and the $T_{m,n}$ are higher order generalizations of the so-called Pompieu operator.
By construction (the choice of $\hat{q}_{\mu\nu},q_{\mu\nu}$) the first sum in Eqn.(\ref{equation3}) determines a contraction in $L^p(\Omega,\C), p>2,$ for sufficiently small $(p-2)$, and the second is a compact operator on $L^p(\Omega,\C).$
Under these conditions, {\em Fredholm theory}\index{Fredholm theory} applies, see Section \ref{fredholmsec} of Chapter \ref{reproducingsec}, and Begehr \cite{begehrbok}, p.227, regarding the application of Fredholm theory to such singular integral equations.
The equation for $\rho$, stemming from the term $\left(\frac{\partial}{\partial z}\right)^m\left(\frac{\partial}{\partial \bar{z}}\right)^n f$, in the 
general equation Eqn.(\ref{equation2}),
be recognized, for $m\leq n$ as the inhomogeneous system of a polyanalytic and a polyharmonic equation, namely,
\begin{equation}
\left(\frac{\partial}{\partial \bar{z}}\right)^{m-n} f=w, \left(\frac{\partial}{\partial z}\frac{\partial}{\partial \bar{z}}\right)^m w=f
\end{equation}
This illustrates how the study of the kernel of the operators $\left(\frac{\partial}{\partial \bar{z}}\right)^q$ 
for positive integers $q$, arise naturally as a cornerstone in the analysis of an arbitrary higher order linear complex partial differential equation. Studying the kernel of each
$\left(\frac{\partial}{\partial \bar{z}}\right)^q$ (and regarding boundary value problems, also the inhomogeneous counterpart) is thus an important step toward 
the study of 
the general system. 
Balk \& Zuev \cite{balkzuev} (1970), and later Balk \cite{ca1} (1991) collected the main results, up to that period, on polyanalytic functions ($q$-analytic functions) of one complex variable in survey form. Since then the theory of $\alpha$-analytic functions has developed in various directions simultaneously. 
The aim of the current text is to provide an introduction to the modern theory of $\alpha$-analytic functions, bringing up to date past surveys regarding the theory in one variable by presenting the new settings to which $\alpha$-analytic functions have more recently been generalized. Emphasis has been put on pure mathematics, and only selected topics from the periphery of hard analysis have been included.
The bulk of the classical results on the topic are those in the vein of the classical Russian school of $q$-analytic functions, whose fore-figure seems to have been Mark Benevich Balk (the results of the author can be regarded as belonging to this category). It seems that Vazgain Avanissian was the first to seriously consider the generalization to several complex variables (i.e.\ replacing the integer $q$ by a multi-integer $\alpha$), although it should be mentioned that Balk \& Zuev \cite{balkzuev}, Parag. 8, gave an alternative (equivalent) definition of the $\alpha$-analytic functions, not based upon partial differential equations, and before that 
Balk \cite{balk69} introduced such functions for the case of $\C^2$ and proved a generalized uniqueness theorem for such functions. 
However, there exist also fairly developed branches of the current theory of $\alpha$-analytic functions, the main one being
boundary value problems for polyanalytic functions 
(where pioneering work is due to Heinrich Begehr).
Another developed branch is approximation, where where many important results seem to be inspired by methods used by Anatoli Georgievich Vitushkin used in the theory of holomorphic functions (here important current contributions are due to Joan Josep Carmona). 
Other existing branches are
operator theoretical analysis of polyanalytic functions (mainly generalizations of the theory of Bergman spaces) and hypercomplex analysis (e.g.\ the theory of $k$-monogenic functions).
In a remarkable paper, Ahern \& Bruna \cite{ahernbruna} present a notion of polyanalytic type in several complex variables, that is in a sense isotropic. This notion turn out to be highly natural from the perspective of extending polyanalyticity to both
hypoanalytic theory and infinite-dimensional complex analysis. Hence for future developments of the theory these authors may have played a fundamental role.
It is important to point out that from the perspective partial differential equations, the theory of $q$-analytic functions could be thought of as a small subcategory of the theory of elliptic equations, thus many deep and powerful results are rendered automatically (e.g.\ the Malgrange-Lax theorem) due to known results which one could for instance look for in the excellent book of Tarkhanov \cite{tarkhanov}.
However, it will become clear to the reader that
the potential of future developments of the theory of $\alpha$-analytic functions is vast.
Namely, any time there exists an analogue of the Cauchy-Riemann operator there usually exists a natural
way to introduce the concept of $q$-analyticity, which is exemplified for the case of infinite dimensional complex
analysis and analysis on the Gaussian integers. 
Another way the theory can be generalized is of course to replace the monomial defining partial differential 
operators with polynomial differential operators (with renders the theory of meta-analytic functions), such a theory is deserving of a monograph in its own right and will not be discussed much in this book.

\section*{Some notation}
Given an open subset, $\Omega$, of a manifold $M$, we denote the set of real-analytic functions on $\Omega$ by $C^\omega (\Omega).$ The set of continuous functions on $\Omega$ is denoted by $C^0(\Omega), $ and the set of $k$-times continuously differentiable functions by $C^k(\Omega)$. 
 The set of holomorphic functions on an open subset $\Omega\subset \Cn$ is denoted $\mathscr{O}(\Omega).$ If not stated otherwise we denote by $z=x+iy$ the complex coordinate on $\C^n,$ where $z=(z_1,...,z_n),$ $z_j=x_j+iy_j,$ and $(x,y)$ are Euclidean coordinates for $\R^{2n}.$ We also denote by $\re z$ and $\im z$ the real part and the imaginary part respectively, of the complex vector $z$. 
 Throughout the book we shall, depending upon the context, use the notation
 $d\mu(z)$ for the area measure on $\Cn$, 
 for example, in the complex plane, depending upon the context, this shall stand for Lebesgue measure $dx\wedge dy$.
 The Cauchy-Riemann operator\index{Cauchy-Riemann operator} on $\C$ is denoted $\overline{\partial}:=\frac{\partial}{\partial \bar{z}}=\frac{1}{2}\left(\frac{\partial}{\partial x} +i\frac{\partial}{\partial y}\right).$ 
 Alternatively $\partial_{\bar{z}}:=\frac{\partial}{\partial \bar{z}},$ in order to conform to some cited literature.
 For a positive integer $n$, a Euclidean coordinate $t$ on $K^n$ (where $K$ is the field $\R$ or $\C$) and a multi-index $\alpha=(\alpha_1,...,\alpha_n)\in \N^n,$
we define the monomial $t^\alpha:=t_1^{\alpha_1}\cdots t_n^{\alpha_n}.$ Denote $\abs{\alpha}=\sum_j \alpha_j$ and by $\abs{t}$ the Euclidean norm. Let $X,Y,Z$ be differential manifolds and let $f:X\to Y,$ $g:Y\to Z$ be maps. The {\em pullback}\index{Pullback} of $g$ by $f$
is the map $f^*g=g\circ f.$ Let $p\in X.$ The {\em pushforward}\index{Pushforward} of a vector $v\in T_p X$
is the vector $f_*v\in T_{f(p)}Y$ defined by $f_*v(h)=v(h\circ f)$, for all smooth $g:Y\to \R.$
Recall that a tangent vector can be identified via an equivalence class of curves $[\gamma]$, $\gamma:(-1,1)\to Y,$ 
where $\frac{\partial(h\circ\gamma)}{\partial t}(0)=v.$
The pushforward of an equivalence class $[\gamma]$ of curves is given by
$f_*v=f_*[\gamma]=[f\circ \gamma].$
All manifolds are assumed to be Hausdorff and paracompact. We always assume the ZFC system of axioms (in particular the axiom of choice is assumed).
We shall assume the reader is familiar with the definition of distributions and distribution solutions (see e.g. Friedlander \& Joshi \cite{friedlander}).
Recall that the {\em dual}, $X^*$ of a real or complex vector space $X$ is the set of linear functionals $X\to K$ (where $K$ is the field $\R$ or $\C$). 
Given a real differentiable manifold $M$ the dual, $T^*_p M,$ of the tangent space $T_p M$, at $p\in M,$ is called the {\em cotangent space}\index{Cotangent} at $p$ and its elements are called cotangent vectors. A real differential form (of degree one) on $M,$ is a section of the cotangent bundle $\bigcup_{p\in M} T_p^* M.$ In this text \index{Section} sections of a vector bundle, $B$, over a domain $\Omega$ will be denoted $\Gamma(B,\Omega)$ and when it is clear from the context what $B$ is we simply write $\Gamma(\Omega).$ The $k$:th {\em exterior power}\index{Exterior power} $\Lambda^k T_p^* M,$ of $T^*_p M$ is defined as follows. $\Lambda^0 T_p^* M :=\{0\},$ and for $k>0,$ $\Lambda^k T_p^* M$ is defined as the set of {\em alternating} $k$-multilinear maps $\eta\colon \overbrace{T_p^* M \times\cdots \times T_p^* M}^{k\mbox{ times}}\to \R$,
where alternating means that 
$\eta(v_{\pi_1},\ldots,v_{\pi_k})=\mbox{sign} \pi \eta(v_1,\ldots,v_k),$ for all
$v_j\in T^*_p M$ and all permutations $\pi$ of $\{1,\ldots,k\}.$ The {\em wedge product}\index{Wedge product} of an $r$-form $\xi$ and an $s$-form $\eta$ is an $(r+s)$-form, denoted $\xi\wedge \eta$, given by
$\frac{1}{r!s!} \sum_{\pi} \mbox{sign} \pi \xi(v_{\pi_1},\ldots,v_{\pi_r})\eta(v_{\pi_{r+1}},\ldots,v_{\pi_{r+s}}),$ where summation is taken over all permutations $\pi$ of $\{1,\ldots,r+s\}.$
For the Euclidean space $\R^{m}$ with coordinates $(x_1,\ldots,x_m)$
 we can define the complex differential forms of degree $k$ (i.e.\ each section of the cotangent bundle) as $\phi=\sum_{\beta}\phi_\beta dx_{\beta_1}\wedge \cdots \wedge dx_{\beta_{k}},$ (here the $\phi_{\beta}$ are complex valued and skew-symmetric in the indices $\beta_1,\ldots,\beta_{k}$,
 and the $\binom{m}{k}$ elements $dx_{\beta_1}\wedge\cdots dx_{\beta_k}$, $1\leq \beta_1<\cdots <\beta_k\leq m$ form a basis for the set of differential forms of degree $k$). 
 Exterior differentiation\index{Exterior differentiation}, $d$, is a linear operator acting on the set of differential forms of order $k$, sending these to the set of differential forms of order $k+1$ according to:
 $d\colon f\mapsto df,$ where $df_p:v\mapsto \mbox{Jac}(f)(p)\cdot v,$ for a real valued $0$-form $f$ and where $\mbox{Jac}$ denotes the Jacobian. If $\xi$ and $\eta$ are forms of order $l$ and $m$ respectively then we define $d(xi\wedge\eta)=(d\xi)\wedge \eta+(-1)^l\xi\wedge d\eta .$ It has the well-known property $d\circ d=0.$
 To see this, let $(x_1,\ldots,x_n),$ denote Euclidean coordinates on the domain of a $0$-form $f$ and write $d(df)=d\sum_j \partial_j f dx_j=
 \sum_{l<j} (\partial_l\partial_j f-\partial_j\partial_l f) dx_l\wedge dx_l =0.$ For a basis element $\eta=fdx_\iota$ (where $\iota$ is an $n$-multi-integer) of the set of $k$-forms we have $d(d\eta)=d(df)\wedge dx_{\iota} +df\wedge  d(dx_\iota),$ and $d(dx_\iota)=d(1)\wedge dx_{\iota}.$
 A $k$-form $\xi$ is called {\em closed} if $d\xi=0,$ and it is called {\em exact} if there exists a $k-1$ form $\nu$ such that $d\nu=\xi.$
  When $m=2n$ and we have Euclidean coordinates $(x_1,y_1,\ldots,x_n,y_n)$ then we identify $dz=dx_j+idy_j$ and $d\bar{z}=dx_j-idy_j$ and thus obtain a representation of any complex differential form of order $k=p+q$ as
 $\phi=\sum_{(\alpha,\beta)}\phi_{(\alpha,\beta)} dz_{\beta_1}\wedge \cdots \wedge dz_{\beta_{q}}\wedge d\bar{z}_{\alpha_1}\wedge \cdots \wedge d\bar{z}_{\alpha_{p}},$ (here the $\phi_{(\alpha,\beta)}$ are complex valued and skew-symmetric in the indices $(\alpha_1,\ldots,\alpha_p,\beta_1,\ldots,\beta_{q}).$ 
Such $\phi$ is said to have {\em bidegree $(p,q)$}.
Denote, for a given open subset $\Omega\subset\Cn$, the set of smooth complex $(p,q)$-forms by $\mathcal{A}^{p,q}(\Omega).$
When $\phi$ is a $0$-form, we can express $d\phi$ as a linear combination fo the differentials $dz_j,d\bar{z}_j,$
$d\phi=\partial \phi +\overline{\partial}\phi,$
where
$\partial\phi=\sum_{j=1}^n \frac{\partial \phi}{\partial z_j} dz_j,$
$\overline{\partial}\phi=\sum_{j=1}^n \frac{\partial \phi}{\partial \bar{z}_j} d\bar{z}_j.$
This decomposes $d$ into a $\C$-linear and an anti-$\C$-linear part respectively.
When $\phi=\sum\phi_{\alpha,\beta} dz^\alpha\wedge d\bar{z}^\beta,$ is a $(p,q)$-form $d\phi$ decomposes into
$\partial\phi :=\sum_{\alpha,\beta} \partial \phi_{\alpha,\beta} dz^\alpha\wedge d\bar{z}^\beta$
and $\overline{\partial}\phi :=\sum_{\alpha,\beta} \partial \phi_{\alpha,\beta} dz^\alpha\wedge d\bar{z}^\beta$ respectively.
The property $d=\overline{\partial}+\partial$ together with $d^2=0$ implies $\partial\overline{\partial} =-\overline{\partial}\partial.$
Next note that by definition we have the following property for $d$ and two forms $\xi$ and $\eta$ of orders $l$ and $m$ respectively
$d(\xi\wedge \eta)=d \xi \wedge \eta +(-1)^{l}\xi \wedge d \eta,$ and a direct consequence is that
for $\alpha\in \mathcal{A}^{p,q}$
\begin{equation}
\partial (\alpha\wedge \beta)=\partial \alpha \wedge \beta +(-1)^{p+q}\alpha \wedge \partial \beta
\end{equation}
\begin{equation}
\overline{\partial} (\alpha\wedge \beta)=\overline{\partial} \alpha \wedge \beta +(-1)^{p+q}\alpha \wedge 
\overline{\partial} \beta
\end{equation}
Let us recall a classical result on closed differential forms.
\begin{proposition}\label{primitihorm}
Let $v$ be an $n+1$-form with $C^k$-smooth coefficients, $k\geq 1$ on an open convex subset $U$ of a finite dimensional vector space $V.$ If $dv=0$ then there exists an $n$-form $u$ on $U$ with $C^k$-smooth coefficients such that $du=v.$
\end{proposition}
\begin{proof}
We may assume $0\in U$ and set $U':=\{(x,t)\in V\times \R:tx\in U\}$ which is an open neighborhood of $U\times [0,1].$ Set $f(x,t)=f_t(x)=tx$
and $f^*v=f^*_t v+dt\wedge _t$ where $f^*_t$ is defined regarding $t$ as a parameter so that $f^*_t v$ and $w_t$ are differential forms that do not contain a factor $dt$ but only differentials of the coordinates in $U.$ Since $df^*v=f^*dv=0$ we have $0=dt\wedge(\partial_t(f_t^* v)-d_x w_t)+R$ where $R$ is a form which does not contain $dt$ and $d_x w_t$ denotes the differential of $w$ then $t$ is regarded as a parameter. This yields $\partial_t(f_t^* v)=d_x w_t$ and integration from $t=0$ to $t=1$ yields
$f_1^* v-f_0^*v=du$ where $u=\int_0^1 w_t dt.$ Since $f_1$ is the identity and $f_0$ maps $U$ to $0$ this implies $v=du$. This completes the proof.
\end{proof}

Let $V_j,1\leq j\leq k,$ be vector spaces over $\R,$ 
Then the tensor product\index{Tensor product}, $V_1 \otimes \cdots \otimes V_k,$ 
is defined as the quotient,
$P/G,$ where $P$ is the space of formal linear combinations over $\R,$ of elements of $V_1 \otimes \cdots \otimes V_k,$ 
(i.e.\ elements of the form $\sum_{j=1}^n c_j(v_{1,j},\ldots,v_{k,j})$, for some $n\in \N$, which is a real vector space when we define $(0,\ldots,0)$ as the additive identity, $1\cdot (v_1,\ldots,v_p)=(v_1,\ldots,v_p)$
and $0\cdot (v_1,\ldots,v_p)=0$) and where $G$ is the vector 
subspace generated by elements of the form,
\begin{eqnarray}
(i)\quad (v_1 ,\cdots ,v_{i-1} ,x,v_{i+1} ,\ldots ,v_k) +
(v_1 ,\cdots ,v_{i-1} ,y,v_{i+1} ,\ldots ,v_k)- \nonumber\\
\nonumber
(v_1 ,\cdots ,v_{i-1} ,x+y,v_{i+1} ,\ldots ,v_k),\\ 
\nonumber
(ii)\quad c(v_1 ,\cdots ,v_k)-
(v_1 ,\cdots ,v_{i-1} ,c,v_{i+1} ,\ldots ,v_k) 
\end{eqnarray}
In the quotient space we have,
\begin{multline}
(v_1 \otimes\cdots \otimes v_{i-1} \otimes x\otimes v_{i+1} \otimes \ldots \otimes v_k) +
(v_1 \otimes\cdots \otimes v_{i-1} \otimes y\otimes v_{i+1} \otimes \ldots \otimes v_k)= \\
(v_1 \otimes\cdots \otimes v_{i-1} \otimes (x+y)\otimes v_{i+1} \otimes \ldots \otimes v_k),\mbox{  and}
\\
c(v_1 \otimes \ldots \otimes v_k) =(v_1 \otimes\cdots \otimes v_{i-1} \otimes cv_i\otimes v_{i+1} 
\otimes \ldots \otimes v_k),\quad i=1,\ldots k.
\end{multline}
The element in $P/G$ corresponding to the equivalence class of $(v_1,\ldots,v_k)\in V_1\times\cdots\times V_k$
in $P/G$
is denoted $v_1\otimes\cdots\otimes v_k$ (this is called a simple tensor).
In the tensor product $V_1\otimes V_2$ we have for $(u,v),(u',v')\in V_1\times V_2$ and a scalar $c,$
that $u\otimes v=v\otimes u,$ $u\otimes v+u'\otimes v=(u+u')\otimes v$, 
$u\otimes v+u\otimes v'=u\otimes (v+v')$ and $c u\otimes v=u\otimes cv.$
Here we could view $V_1\otimes V_2\simeq  P/G$
where $P$ is interpreted as a linear map on $V_1\times V_2$
and $G$ is the set of all formal sums of the form 
$(cu,v)-c(u,v),$ $(u,cv)-c(u,v),$, $(u+u',v)-(u,v)-(u',v),$ $(u,v+v')-(u,v)-(u,v'),$
where we denote $(u,v):=1(u,v)$. Clearly, $P/G$ is an $\R$-module since both $P$ and $G$ are and $G$ is a submodule of $P.$
Denote the quotient map $q:V_1\times V_2 \to V_1\otimes V_2$, $q(u,v)=u\otimes v$
(by which we mean $q(u,v)=u\otimes v +G$ so $u\otimes v$ is the residue in the quotient target space).
A map
$\phi:M_1\to M_2$ between $\R$-modules is called {\em bilinear}
if for each $u\in M_1$ the map $v\mapsto \phi(u,v)$ is a linear map $M_2\to Y,$ and for each
$v\in M_2$ the map $u\mapsto \phi(u,v)$ is a linear map $M_1\to Y.$
The choice of $G$ makes $q$ a bilinear map. 
If $V$ is a real vector space then as a real vector space, the tensor product $\C\otimes V$ is generated
over $\R$ by $1\otimes v$ and $i\otimes v$ for $v\in V.$ Then $\C\otimes V$ can be made into a complex vector space by defining $a(b\otimes v)=
ab\otimes v$ for all $a,b\in C$ and $v\in V.$ As a complex vector space $\C\otimes V$ is generated by $1\otimes V.$ The natural 
conjugation for $\C\otimes V$ is, if we identify $1v=1\otimes v$, and write $av$ for $a\otimes v,$ given by $\overline{av}=\bar{a}v=\bar{a}\otimes v.$
Tensor products can also be defined in terms of universal maps.
This is particularly simple in the case of tensor product of vector spaces.
In particular, if $V_1$ and $V_2$ respectively are finite dimensional vector spaces with basis
$\{e_1,\ldots,e_N\}$ and $\{f_1,\ldots,f_M\}$ for positive integers $N,M$ then
the tensors $\{e_i\otimes f_j\}_{i,j}$ form a basis for $V_1\otimes V_2,$
i.e.\ any tensor $u\otimes v\in V_1\otimes V_2$ takes the form
$u\otimes v=\sum_{i,j} c_{ij} (e_i\otimes f_j)$ for constant scalars $c_{ij}.$
Recall that a linear map is uniquely determined by its action on the basis vectors and it 
sends bases to bases i.e.\ the elements $e_i\otimes f_j$, $i=1,\ldots,N,$ $j=1,\ldots,M$ form a basis
for $V_1\otimes V_2$ which thus is an $NM$ real vector space.
Denote by $B$ the map which sends for arbitrary constants $a_i,b_j$ the elements
$(\sum_{i=1}^N a_ie_i,\sum_{j=1}^M b_jf_j)$ to $\sum_{i,j} a_ib_j (e_i\otimes f_j).$
One verifies 
immediately that $B$ is a bilinear map.
Suppose $T$ is a bilinear map from $P$ (the formal finite linear combination of pairs in 
$V_1\times V_2$) to $V_1\otimes V_2$. Define the linear map $A:V_1\otimes V_2\to Y$
according to 
$A(\sum_{i,j} c_{ij} (e_i\otimes f_j))=\sum_{i,j} c_{ij} T(e_i\otimes f_j)$.
Then we have $T=A\circ B.$ Hence any arbitrary bilinear $T:P\to V_1\otimes V_2$ factors
through the {\em universal} map $B$ for some linear map $A$ (this result also holds for the more general case of
tensor product of $R$-modules, $M_1,M_2$ over a unital ring $R$).
If $M$ is a finite dimensional manifold then for each point 
$p\in M$ the cotangent space $T^*_p M$ is a vector space with same dimension as $M$
so we can
form the tensor product of $T^*_p M\otimes T^*_p M$ and thereby define
so called tensor fields\index{Tensor field} as smooth sections of the bundles $T^*M\otimes T^*M$ where we
require smoothness with respect to the base point, of the assignment of a tensor to each $T^*_p M\otimes T^*_p M$. 

\chapter{Basic properties}\label{basicpropsec}

\section{Defining equations and immediate properties} 

Most of the current literature on $\alpha$-analytic functions regard the case of one complex dimension and from an introductory viewpoint it makes sense to describe this case separately before we introduce the more general case.
In this text we shall predominantly use the following terminology proposed by Burgatti \cite{burgatti} in 1922.
\begin{definition}[$q$-analytic functions]\label{nolkoldef}
Let $\Omega\subset \C$ be an open subset and let $q\in \Z_+.$ 
A distribution solutions, $f,$ to
$\partial_{\bar{z}}^q f=0$ on $\Omega$ 
is called a $q$-analytic function on $\Omega.$
\end{definition}
As pointed out by Balk \cite{ca1} there are many different equivalent definitions of
and terminology used simultaneously for what we call $q$-analytic functions and 
Balk \cite{ca1} uses interchangeably $n$-analytic and polyanalytic of order $n$ 
with a preference for the latter. 
\begin{definition}[Polyanalytic functions of order $q$]\label{endim}
Let $\Omega\subseteq\C$ be a domain and let $q\in \Z_+$. A function $f$ is called {\em polyanalytic of order $q$ at $p_0$} if it can, near $p_0$, be represented in the form
\begin{equation}\label{eq1}
f(z)=\sum_{j=0}^{q-1}a_j(z)\bar{z}^j
\end{equation} 
where $a_j, j=0,\ldots , q-1,$ are holomorphic functions near $p_0.$ The case $a_{q-1}\equiv 0$ is not excluded.
When $a_{q-1}\not\equiv 0$ the number $q$ is called the {\em exact} order of polyanalyticity of $f$.
A function is called polyanalytic of order $q$ on $\Omega$ if it is polyanalytic of order $q$ at each point of $\Omega.$
The space of polyanalytic functions on $\Omega$ is denoted by $\mbox{PA}_q(\Omega)$. The functions $a_j$ are called the analytic (or holomorphic)\index{Analytic component} components of $f$.
\end{definition}
Note that since $\Omega$ is connected, each of the local analytic components of $f(z)$
extend to a global (on $\Omega$) analytic component.
We shall prove below that Definition \ref{nolkoldef} is equivalent to
Definition \ref{endim}.
\begin{definition}
A function $f$ is called {\em countably analytic} on an open subset $U\subset \C,$ if for every point $p\in U$ there is an open
neighborhood $U_p$ of $p$ such that on $U_p$, $f$ can be represented by a uniformly convergent series $f(z)=\sum_{j=0}^{\infty} a_{p,j}(z)\bar{z}^j$ for 
holomorphic $a_{p,j}(z)$ on $U_p.$
\end{definition}

\begin{definition}
	Let $\Omega\subset\C$ be a domain and let $q\in \Z_+$.
	A polyanalytic function $g(z)=\sum_{j=0}^{q-1}a_j(z)\bar{z}^j$, for holomorphic $a_j(z)$ on $\Omega$, is called {\em reduced}\index{Reduced polyanalytic function} if
	there exists holomorphic functions $a'_j(z)$ such that $a_j(z)=z^ja'_j(z),$ $j=0,\ldots,q-1.$
\end{definition}

\begin{definition}[Elliptic operator]
Let $\Omega\subset \R^n$ be an open subset, with Euclidean coordinate $z,$ and let $L$ be a partial differential operator of order $m$ with complex valued coefficients. Assume $L$ can be represented as a differential polynomial $P(\partial):=\sum_{\abs{\beta}\leq m} c_{\beta}(z) \partial^\beta,$ where the $c_\beta(z)$ are complex valued coefficients (for convenience we have here used the notation $\partial$ not to be confused with the $\C$-linear $\partial$ operator acting on differential $(p,q)$-forms to yield $(p+1,q)$-forms).
$L$ is called {\em elliptic} if the {\em symbol}\index{Symbol of differential operator}, $\sigma_{L}(z,\xi):=\sum_{\abs{\beta}= m} c_{\beta} \xi^\beta,$
satisfies $\sigma_{L}(z,\xi)\neq 0$ for all $(z,\xi)\in \Omega\times (\Cn\setminus \{0\}).$ Here $\xi$ denotes a coordinate in a fiber of the cotangent bundle.
\end{definition}
For Euclidean coordinates $(x_1,x_2)$ in $\R^2$, we can identify the Cauchy-Riemann operator in $\C$ as 
$\overline{\partial}(\frac{\partial}{\partial x_1},\frac{\partial}{\partial x_2})=\frac{1}{2}\left(\frac{\partial}{\partial x_1} +i\frac{\partial}{\partial x_2}\right),$
where $z=x_1+ix_2$ is a Euclidean coordinate for $\C.$
The symbol $\sigma_{\overline{\partial}}(z,\xi)=\frac{1}{2}\left(\xi_1 +i\xi_2\right),$ vanishes precisely for $\xi=0,$ thus $\partial_{\bar{z}}$ is elliptic.
\begin{definition}[Hypoelliptic operator]\index{Hypoelliptic operator}
Let $\Omega\subset \R^n$ be an open subset and let $L$ be a partial differential operator with complex coefficients. $L$ is called {\em hypoelliptic at $p_0\in \Omega$} if there is an open $U\ni p_0$ such that any distribution solution $u$ to $Lu=0$ can be represented in terms of a $C^\infty$-smooth function on $U.$ $L$ is called {\em hypoelliptic} if it is hypoellitic at each point.
\end{definition}
The well-known Elliptic Regularity Theorem\index{Elliptic Regularity Theorem} (see Theorem \ref{elliptictheorem}) states that any distribution solution $u$ to $Lu=0,$ where $L$ is elliptic, can be represented by a $C^\infty$-smooth function
(i.e.\ elliptic operators are hypoelliptic). 
Let $\Omega\subseteq\R^n$ be a domain with Euclidean coordinate $x$, let $q\in \Z_+,$ and let $L$ be an elliptic operator of order $m$, with complex coefficients, on $\Omega.$ Then
$L^q$ is elliptic: Indeed,
we can write $L^q = \left(\sum_{\abs{\alpha}=qm} a_\alpha(x) \partial^\alpha\right) +\mbox{lower order terms}.$
We can thus identify
$\sigma_{L^q}(x,\xi)=\left(\sigma_{L}(x,\xi)\right)^q,$ where each factor on the right hand side is a polynomial in $\xi$
that is nonzero for all nonzero $\xi.$
\\
In particular, the operators $(\partial_{\bar{z}})^q$ are elliptic for each $q\in \Z_+.$
\begin{proposition}\label{prop2}
Let $\Omega\subseteq\C$ be a domain and let $q\in \Z_+$. Then $\mbox{PA}_q(\Omega)$ coincides with the set of 
$q$-analytic functions on $\Omega$ and
the representation in Eqn.(\ref{eq1}) is unique.
\end{proposition}
\begin{proof}
The members of $\mbox{PA}_q(\Omega)$ are real-analytic.
Let us for brevity (in this proof) denote by $\mbox{Ker}(\partial_{\bar{z}}^q)$
the set of $q$-analytic functions on $\Omega$.
$\mbox{PA}_q(\Omega)\subset \mbox{Ker}(\partial_{\bar{z}}^q)$ follows from the fact that smooth solutions are also distribution solutions.
We prove $\mbox{Ker}(\partial_{\bar{z}}^q)\subset \mbox{PA}_q(\Omega)$ by induction in the positive integer $q.$
The case $q=1$ follows from ellipticity together with the identity theorem for holomorphic functions. Let $m\geq 2$ and
assume the statement holds true for $q=1,\ldots, m-1.$
Let $q=m.$ Note that 
$\partial_{\bar{z}}^{q-1}(\partial_{\bar{z}}f)\equiv 0,$
so by the induction hypothesis we have the following representation near $p_0,$ $\partial_{\bar{z}}f=\sum_{j=0}^{q-2}b_j(z)\bar{z}^j,$ for holomorphic $b_j(z)$ on $\Omega.$
Thus $\partial_{\bar{z}}(f(z)-\bar{z}h(z))\equiv 0,$ 
for $h(z):=\sum_{j=0}^{q-2}b_j(z)\frac{\bar{z}^j}{j+1}$.
This implies 
$f(z)=\bar{z}h(z)+a_0(z),$ for a holomorphic $a_0(z)$ on $\Omega$. This completes the induction regarding the existence of the given local representation. 
Finally, given a representation of the form in Eqn.(\ref{eq1}), and assuming that the given solution $f$ vanishes near a point $p_0\in \Omega,$ the equations
$\partial_{\bar{z}}^j f\equiv 0,$ $1\leq j\leq q-1,$ imply that 
$a_j\equiv 0.$ By the identity theorem for holomorphic functions each component thus vanishes on all of $\Omega,$ implying uniqueness of the components in the representation of $f$ (and the same argument applies if we replace $f$ with the difference of two solutions whose representations agree near some point $p_0\in \Omega$). This completes the proof.
\end{proof}

In particular, Proposition \ref{prop2} implies $ \mbox{PA}_q(\Omega)\subset C^\omega(\Omega)$ for all $q\in \Z_+$ (where $C^\omega(\Omega)$ denotes the set of real-analytic functions on $\Omega$). Note that the proof of Proposition \ref{prop2} establishes (since the representation is unique) the following identity theorem: 
\begin{proposition}
Let $\Omega\subseteq\C$ be a domain and let $q\in \Z_+$. If $f\in \mbox{PA}_q(\Omega)$ and $f(z)=0$ for all $z\in V$ for some open subset $V\subset \Omega.$ Then $f\equiv 0$ on $\Omega.$
\end{proposition}
Note also that for a given domain, $\mbox{PA}_q$ is a module over the set of holomorphic functions $\mathscr{O},$ i.e.\ by abuse of notation, assuming that a global holomorphic coordinate, $z$ is given we can write
$\mbox{PA}_q=\oplus^{q}_{j=1} \bar{z}^{j-1}\mathscr{O}$.
The $q$-analytic counterpart to the identity theorem that establishes vanishing of a complex analytic function by vanishing of all complex derivatives at a single point is the following (a modification is of course required since for example the function $f(z)=z\bar{z}$ satisfies $\partial_z^j f(0)=0$ for all $j\in \Z_{\geq 0}$ but $f\not\equiv 0$).
\begin{proposition}\label{uniquepointprop}
Let $q\in \Z_+,$ let $\Omega\subset\C$ be a domain and let $f\in \mbox{PA}_q(\Omega).$ If there exists a point $p_0\in \Omega$ such that $\partial_z^j\partial_{\bar{z}}^k f(p_0)=0$
for all $j,k\in \Z_{\geq 0},$ then $f\equiv 0.$
\end{proposition}
\begin{proof}
Without loss of generality suppose $p_0=0.$ By definition $f$ has a representation of the form
$f(z)=\sum_{j=0}^{q-1} a_j(z)\bar{z}^j,$ for holomorphic $a_j.$ Then $(q-1)!a_j(z)=\partial_{\bar{z}}^{q-1} f(z).$
Thus $\partial_{z}^{j} a_j(0)=0$ for all $j\in \Z_{\geq 0}$ which implies $a_j\equiv 0.$ Thus 
$f(z)=\sum_{j=0}^{q-2} a_j(z)\bar{z}^j,$ and the exact same process can be repeated to show that $a_{q-2}\equiv 0$ and so on. This completes the proof.
\end{proof}
The generalization of $q$-analytic functions to $\Cn$ was introduced properly, in terms of partial differential equation, by Avanissian \& Traor\'e \cite{avan1}, see also Traor\'e \cite{traore94}, \cite{traore97} \cite{traore02}, although it should be mentioned that Balk \& Zuev \cite{balkzuev}, Parag. 8, gave an alternative (equivalent) definition of the $\alpha$-analytic functions not based upon partial differential equations, and before that 
Balk \cite{balk69} introduced such functions for the case of $\C^2$ and proved a generalized a uniqueness theorem for such functions. We present here both definitions and we shall later prove their equivalence.
\begin{definition}[Balk \& Zuev \cite{balkzuev}]
	Let $\omega\subset\Cn$ and $\alpha\in \Z_+^n.$ 
	A function $f$ on $\Omega$ is called {\em polyanalytic of vectorial order $\alpha$}
		if it is a polynomial with respect to $\bar{z}$ of degree $\alpha_j-1$ with respect to $\bar{z}_j$,
		$j=1,\ldots,n$ and with holomorphic coefficients with respect to $z$ on $\omega.$
		A function $f$ on a domain $\Omega\subset \Cn$ is called {\em areolar} at $p\in\omega$
		if it is expressible in some polycylinder $\delta\subset \Omega$ with center $p$ as
		\begin{equation}
		\sum_{\abs{\beta}=0}^\infty (\bar{z}-\bar{p})^\beta \phi_{\beta,p}(z)
		\end{equation} 
		for holomorphic $\phi_{\beta,p}$ on $\delta.$
\end{definition} 
\begin{definition}[$\alpha$-analytic functions]\label{avandefhel}
Let $\Omega\subset \Cn$ be a domain and let $\alpha\in \Z^n_+.$ 
A function $f$ on $\Omega$ is called {\em $\alpha$-analytic} on 
$\Omega,$
if it is a
distribution solutions to
$\partial_{\bar{z}_1}^{\alpha_1}f=\cdots= \partial_{\bar{z}_n}^{\alpha_n} f=0$ on $\Omega.$
The space of $\alpha$-analytic functions (or {\em polyanalytic functions of order $\alpha$}) on $\Omega$ is denoted by $\mbox{PA}_\alpha(\Omega)$. 
\end{definition}
In Definition \ref{avandefhel}, the order of an $\alpha$-analytic function is said to be {\em exact} if $f$ has, near each $p_0\in \Omega,$ a representation (see Proposition \ref{prop1} for the motivation of this) of the form
\begin{equation}\label{eq1a}
f(z)=\sum_{0\leq \beta_j<\alpha_j}a_\beta(z)\bar{z}^\beta
\end{equation} 
where the $a_\beta$ are holomorphic functions near $p_0,$ such that
$a_{\alpha_j-1}\not\equiv 0, j=1,\ldots,n$.

\begin{proposition}\label{systemab}
	Let $\alpha,\beta\in \Z_+^n,$
	$\beta_j\leq \alpha_j,$ $j=1,\ldots,n.$ Let $\Delta(p,r)$ be a polydisc in $\Cn,$ with center $p\in \C,$ where $r=(r_1,\ldots,r_n).$
	Let $f_1,\ldots,f_n$ be polyanalytic of order $\leq\alpha$ (by which we mean that the separate order with respect to $z_j$ is $\leq \alpha_j,$ $j=1,\ldots,n$) on $\Delta(p,r)$
	such that
	\begin{equation}\label{homoab} \partial_{z_j}^{\beta_j}f_k=\partial_{z_k}^{\beta_k}f_j,\quad j,k=1,\ldots,n\end{equation}
	Then the system
	\begin{equation}\label{lossytsab00} \partial_{z_j}^{\beta_j}f=f_j,\quad j=1,\ldots,n\end{equation}
	has a solution that is polyanalytic of order $\leq \alpha$ on $\Omega$, unique up to
	addition by a polyanalytic polynomial of order $\leq \alpha$, that is of degree $<\beta_j$ with respect to $z_j,$
	$j=1,\ldots,n.$ 
\end{proposition}
\begin{proof}
	Suppose $f_j$ is $\kappa_j$-analytic, $\kappa\in \Z_+^n$. 
	We have, for each $j=1,\ldots,n,$ a representation $f_j=\sum_{\kappa_j\leq \gamma_j} a_{\kappa,j}(z)\bar{z}^\kappa,$ for holomorphic $a_{\kappa,j}.$ The system defined by Eqn.(\ref{homoab})
	is equivalent to
	\begin{equation} \partial_{z_i}^{\beta_i}a_{\kappa,j}=\partial_{z_k}^{\beta_k}a_{\kappa,i},\quad i, j=1,\ldots,n
	\end{equation}
	To find a solution $f=\sum_{\kappa} a_{\kappa}\bar{z}^\kappa$ we must solve
	\begin{equation}\label{lossytsab} \partial_{z_i}^{\beta_i}a_{\kappa}=a_{\kappa,i},\quad i=1,\ldots,n, \kappa_j\leq \alpha_j-1, j=1,\ldots,n\end{equation}
	
	\begin{lemma}\label{inforlossyst}
		Let $\alpha\in \Z_+^n.$ and let $\Delta(p,r)$ be a polydisc in $\Cn,$ with center $p\in \C,$ where $r=(r_1,\ldots,r_n).$
		Given holomorphic functions $f_1,\ldots,f_n,$ satisfying, on the polydisc the system
		\begin{equation}\label{condabbecr} \partial_{z_j}^{\alpha_j}f_k=\partial_{z_k}^{\alpha_k}f_j,\quad j,k=1,\ldots,n\end{equation} 
		there exist a unique (up to a holomorphic polynomial of degree $\leq \alpha_j$ with respect to $z_j$, $j=1,\ldots,n$) holomorphic solution $f$
		to the system
		\begin{equation}\label{lossytsab001} \partial_{z_j}^{\alpha_j}f=f_j,\quad j=1,\ldots,n\end{equation}
	\end{lemma}
	\begin{proof}
		We use induction in $n.$ For $n=1$ we have $f_1=\sum_j a_j (z-p_1)^j$ so we may take $f=\sum_j \frac{a_j}{(j+k)\cdots(j+1)}(z-p_1)^{j+k}.$
		So suppose $n>1$ and that the result holds true for all dimensions $<n.$ If
		$f_n=\sum_j a_j(z_1,\ldots,z_{n-1})(z_n-p_n)^j$ we set $\phi_n:=\sum_j \frac{a_j(z_1,\ldots,z_{n-1})}{(j+\alpha_n)\cdots (j+1)}
		(z_n-p_n)^{j+\alpha_n},$ which gives
		$\partial_{z_n}^{\alpha_n} \phi_n =f_n.$ Setting $\psi_i=f_i-\partial_{z_i}^{\alpha_i}\phi_n$, $i=1,\ldots,n$ we have $\psi_n=0$
		and by the conditions of Eqn.(\ref{condabbecr}) we have $\partial_{z_i}^{\alpha_i}\psi_j=\partial_{z_j}^{\alpha_i}\psi_i,\quad j=1,\ldots,n.$
		This implies $\partial_{z_i}^{\alpha_i}\psi_n=0$ for all $i$, which yields
		$\partial_{z_n}^{\alpha_n}\psi_n=0$, $i=1,\ldots,n-1.$ Hence
		$\psi_i=\sum_{j=0}^{\alpha_n-1} \psi_{ij}(z_1,\ldots,z_{n-1})z_n^j,$ $i=1,\ldots,n-1,$ for holomorphic $\psi_{ij}$ on $\abs{z_i-p_i}<r_i,$
		satisfying the conditions of Eqn.(\ref{condabbecr}) for each fixed $j.$ By the induction hypothesis we have that, for each fixed $j$,
		the system
		$\partial_{z_i}^{\alpha_i}f=\psi_{ij},$ $i=1,\ldots,n-1,$ has a solution $g_j( z_1,\ldots,z_{n-1})$. Hence $g:=\sum_j g_j z_n^j$ is a solution to the system
		$\partial_{z_i}^{\alpha_i}g=\psi_{i},$ $i=1,\ldots,n-1.$ This means that the function $f=g+\phi_n$ is a solution to the system defined by Eqn.(\ref{lossytsab001})
		This proves Lemma \ref{inforlossyst}.
	\end{proof}
	By Lemma \ref{inforlossyst} we have, for each fixed $\kappa,$ a solution $a_{\kappa}$, to the system defined by Eqn.(\ref{lossytsab}) i.e.\
	$\partial_{z_i}^{\beta_i}a_{\kappa}=a_{\kappa,i},\quad i,j=1,\ldots,n,$ 
	where $a_\kappa$ is a holomorphic function, uniquely determined up to a polynomial of order $\leq \beta_i$ with respect to $z_i.$
	This completes the proof.
\end{proof}
Recall that even for holomorphic functions there exists functions
that are real-analytic almost everywhere and satisfies the Cauchy-Riemann equation
everywhere but are not holomorphic.
\begin{example}
Consider 
\begin{equation}
f(z)=\left\{
\begin{array}{ll}
\exp(-z^{-4}) & ,z\neq 0\\
0 & ,z=0
\end{array}
\right.
\end{equation}
Looman \cite{looman} pointed out that although $\partial_{\bar{z}} f\equiv 0,$
it is not analytic at $z=0$ since $f(z)/z\to \infty$ as $z\to 0.$ along arg$z=\pi/4.$
In particular, $f$ has an essential singularity at $0$.
\end{example}
Note that $\partial_{\bar{z}_1}^{\alpha_1}\cdots \partial_{\bar{z}_n}^{\alpha_n}$ is in general not an elliptic operator.
\begin{example}
Let $n=2,$ $\alpha=(1,1)$ and the operator $L=\partial_{\bar{z}_1}^{\alpha_1} \partial_{\bar{z}_2}^{\alpha_2}=$
$\partial_{\bar{z}_1}\partial_{\bar{z}_2}$ on $\Cn.$
Then
$\sigma_L(z,\xi)=\frac{1}{4}\left(\xi_1 +i\xi_2\right)\left(\xi_3 +i\xi_4\right).$
Thus $\sigma_L(z,(\xi_1,\xi_2,0,0))=\sigma_L(z,(0,0,\xi_3,\xi_4))=0,$
so $L$ is not elliptic.
\end{example}
However, the system $\partial_{\bar{z}_1}^{\alpha_1}f=\cdots =\partial_{\bar{z}_n}^{\alpha_n}f=0$ is in a sense hypoelliptic as was proved by Avanissian \& Traor\'e \cite{avan2}.
\begin{theorem}\label{cool}
Assume $\Omega\subset\Cn,$ is a domain and $u$ a distribution $u\in \mathcal{D}(\Omega)$
 satisfying in distribution sense the equations
$\partial_{\bar{z}_1}^{\alpha_1}u=\cdots= \partial_{\bar{z}_n}^{\alpha_n}u=0$
for a multi-integer $\alpha\in \Z_+^n$.
Then there exists an $\alpha$-analytic $f$ on $\Omega$ such that $u=f$ almost everywhere.
\end{theorem}
\begin{proof}
Let $\theta_r(z)=\theta_r(\abs{z})\geq 0$ be a $C^\infty$-smooth function with compact support 
(denoted $\theta_r\in C^\infty_c(\Cn)$), 
$\theta_r\in \mathcal{D}(\Omega)$ such that
$\theta_r =0$ for $\abs{z}\geq r$ and $\int_{\Cn} \theta_r =1.$ 
Recall that the so called regularization by the test function $\theta_r$, is given by
$f_r *\theta_r:=\int f_r(z-\zeta)\theta_r(\zeta)d\mu(\zeta)$ and it 
is well-known that this convolution is $C^\infty$-smooth (see Theorem \ref{regularizationthm}) 
and satisfies $\partial^\beta (u*\theta_r)=u*\partial^\beta \theta_r$ for multi-integers $\beta$
(see e.g.\ Friedlander \& Joshi \cite{friedlander}).
Hence $f_r *\theta_r$ satisfies the system
$\partial_{\bar{z}_1}^{\alpha_1}f=\cdots =\partial_{\bar{z}_n}^{\alpha_n}f=0$
whenever $u$ does and therefore
$f_r$ is $\alpha$-analytic in the usual sense, so there exists holomorphic
$a_{r,\beta}(z),$ $\beta\in \N^n,$ $\beta_j=0,\ldots,\alpha_j-1$, such that
\begin{equation}
f_r(z)=\sum_{\stackrel{0\leq \beta_j\leq \alpha_j-1}{1\leq j\leq n}} a_{r,\beta}(z) \bar{z}^{\beta}
\end{equation}
Furthermore, in the sense of distributions $f_r \to u$ and $\partial_{\bar{z}_1}^{\beta_1}\cdots \partial_{\bar{z}_n}^{\beta_n} f_r \to \partial^\beta u$ (and
$\partial_{z_1}^{\beta_1}\cdots \partial_{z_n}^{\beta_n} f_r\to \partial_{z_1}^{\beta_1}\cdots \partial_{z_n}^{\beta_n} u$)
for all multi-integers $\beta\in N^n.$
We thus obtain that each $a_{r,\beta}$ converges in distribution sense, to say $a_{\beta}$ and furthermore
\begin{equation}
0=\partial_{\bar{z}_j}a_{r,\beta}(z)\to \partial_{\bar{z}_j}a_{\beta}(z) 
\end{equation}
for $j=1,\ldots,n,$
thus $a_{\beta}$ is almost everywhere equal to a holomorphic function. This implies that
$f:=\sum_{\stackrel{0\leq \beta_j\leq \alpha_j-1}{1\leq j\leq n}} a_{\beta}(z)\bar{z}^\beta$
has the wanted properties. 
Here we use the fact that a distribution which satisfies the Cauchy-Riemann equation is a.e.\ equal to a holomorphic function.
This completes the proof.
\end{proof}
Furthermore, the following generalization of Hartogs' theorem on separate analyticity is known. 
\begin{theorem}[Avanissian \& Traor\'e~{\cite[Theorem 1.3, p.\,264]{avan2}}]\label{hartog1}
  Let $\Omega\subset\Cn$ be a domain and let $z=(z_1,\ldots,z_n),$
  denote holomorphic coordinates in $\Cn$ with $\re z=:x, \im z=y$. Let
  $f$ be a function which, for each $j$, is smooth in $x_j,y_j$ and polyanalytic of order
  $\alpha_j$ in the variable $z_j=x_j+iy_j$ (in such case we shall
  simply say that $f$ is separately polyanalytic of order
  $\alpha$). Then $f$ is jointly smooth with respect to
  $(x,y)$ on $\Omega$ and furthermore is polyanalytic of order
  $\alpha=(\alpha_1,\ldots,\alpha_n).$ 
\end{theorem}
\begin{proof}
Suppose $f$ satisfying the conditions of theorem is also bounded on $\Omega,$
in particular a  
By Theorem \ref{cool}
there exists a $alpha$-analytic $\tilde{f}$ on $\Omega$ such that $\tilde{f}-f=0$ almost everywhere.
Let $\Omega_1,\Omega_2$ be open subsets of $\Rn$ and let $\Phi(x,y)$ be a 
real-valued function on $\Omega_1\times \Omega_2$, such that for each $y\in \Omega_2$
$\Phi$ is continuous in the variable $x.$ 
It follows from Fubinis theorem that if $\Phi=0$
almost everywhere then $\Phi\equiv 0.$
Replacing $\Phi$ with $\re(\tilde{f}-f)$ and $\im(\tilde{f}-f)$ respectively
yields the conclusion of the theorem for the case that $f$ is assumed to be bounded.

\begin{lemma}\label{lemk2}
Let $z=(z',z'')\in \C^{k}\times\C^{n-k}$ and let $\Omega=\Omega'\times \Omega''\subset\Cn$ be an open subset. 
Assume a function $f(z)$ is $q'$-analytic with respect to the variables $z'=(z_1,\ldots,z_k)$
for any fixed $z''=(z_{k+1},\ldots,zn),$ such that $(z',z'')\in \Omega$
and $q''$-analytic with respect to the variables $z''=(z_{k+1},\ldots,zn),$
for any fixed $z'=(z_1,\ldots,z_k)$ such that $(z',z'')\in \Omega$.
If $B$ is a polydisc with with nonempty interior and $\overline{B}\subset \Omega$, then there exists 
polydisc $B'\subset B$ with nonempty interior such that $f$ is bounded on $B'.$
\end{lemma}
\begin{proof}
Let $A_C:=\{ \abs{f(z',z'')}\leq C, \mbox{ for all }z''\in \Omega''\}.$
For fixed $z''\in \Omega''$ the function $z'\mapsto f(z',z'')$ is continuous so each $A_C$ is closed. Also
For fixed $z'\in \Omega'$ the function $z''\mapsto f(z',z'')$ is continuous thus bounded. 
Furthermore
$\bigcup_{C=1}^\infty A_C =\Omega'.$ By the Baire category theorem some $A_{C}$ must have nonempty interior, thus must contain a nonempty
closed polydisc, say $B_1\subset \Omega'$ such that $f$ is bounded on $B'=B_1\times\Omega''.$
This completes the proof of the lemma.
\end{proof}
\begin{lemma}\label{lem3k}
Let $z=(z',z_n)\in \C^{n-1}\times\C^{n}$ and let $\Omega\subset\Cn$ be the polydisc 
$\{\abs{z_j}<R,1\leq j\leq n\}$. 
If a function $f(z)$ is $q'$-analytic with respect to the variables $z'=(z_1,\ldots,z_{n-1})$
for any fixed $z_n,$ such that $(z',z_n)\in \Omega$ and if there exists $r>0$ such that
$f$ is $(\alpha',1)$-analytic and bounded on the polydisc
\begin{equation}
B'=\{\abs{z_j}<r,1\leq n-1;\abs{z_n}<R\}
\end{equation}
then $f$ is $(\alpha',1)$-analytic on $B.$
\end{lemma}
\begin{proof}
The function $f$ can be represented as
\begin{equation}
f(z',z_n)=\sum_\beta a_\beta(z',z_n)\bar{z'}^{\beta}
\end{equation}
with $\alpha_\beta$ holomorphic in $z'$.
It is known (as in the case of holomomorphic functions of several variables) 
that there is a constant $c<\infty$ not depending on $f$ such that (see Avanissian \& Traour\'e \cite{avan1})
\begin{equation}
\sup_{\abs{z'}\leq r_1 <r} \abs{a_\beta (z',z_n)}\leq c\sup_{B'}\abs{f}
\end{equation}
Indeed, on $B'$ the Cauchy integral with respect to $z_n$ can be applied to each $a_\beta$ 
(and the holomorphic dependence on the additional variables will be maintained as can be seen from Theorem \ref{dolbeaultlemma}). A power series expansion can be obtained
using for $\abs{z_n}<\abs{\zeta_n},$ $\frac{1}{\zeta_n -z_n}=\sum_{j\geq 0}\frac{z_n^j}{\zeta_n^{j+1}},$
which can be substituted in the Cauchy integral, and by Fubinis theorem
summation and integration can be interchanged. By Hartogs' lemma 
in the theory of holomorphic functions applied componentwise
to the $a_\beta$, this implies that $f$ is bounded and separately $(\alpha',1)$-analytic on the
polydisc $B$ thus 
$f$ is $(\alpha',1)$-analytic on $B.$
\end{proof}

Hence starting from a function $f$ satisfying the conditions of the theorem with $z''=z_n$, we can first conclude by Lemma \ref{lemk2} that
near any point, $Q,$ there is a small polydisc with center $Q$, on which $f$ is bounded and by the introductory remark
$f$ is $\alpha$-analytic on that disc. Furthermore, by Lemma \ref{lem3k} $f$ extends to
the largest polydisc $U_Q$, with center $Q$, whose closure is contained in $\Omega.$
But then we can use induction in the number of elements of $z''$. When that number is $1$, we have $z''=z_n$
which has already been handled so assume the statement holds true when $z'=(1,\ldots,k),$ $k< n-1.$
For fixed $z_n=p_n,$ $f(z_1,\ldots,z_{n-1},p_n)$ can, by what we have already done
be realized as a $(\alpha_1,\ldots,\alpha_{n-1})$ analytic function in $(z_1,\ldots,z_{n-1}).$
Thus replacing the role of $z'$ with $(z_1,\ldots,z_{n-1})$ and that of $z''$ with $z_n$
we obtain that $f$ is $\alpha$-analytic on $U_Q$. This completes the induction. This completes the proof.
\end{proof}
We have the following generalization of Proposition \ref{prop2} to several complex variables.
\begin{theorem}[Avanissian \& Traour\'e \cite{avan1}, Prop.2.1, p.744]\label{prop1} Let $\alpha\in \Z_+^{n}.$ The condition that $f$ satisfies $\partial_{z_1}^{\alpha_1}f=\cdots =\partial_{z_n}^{\alpha_n} f=0$ on a domain $\Omega\subset\Cn ,$ in distribution sense, is equivalent to $f$ having on the given domain, near each point $p_0$, the unique representation $\sum_{0\leq \beta_j<\alpha_j}a_{\beta}(z)\bar{z}^\beta$, where $a_{\beta}(z)$ are holomorphic functions near $p_0.$
\end{theorem}
\begin{proof}
If $f$ has a representation of the form $\sum_{0\leq \beta_j<\alpha_j}a_{\beta}\bar{z}^\beta$, then it is obviously separately
$\alpha_j$-analytic with respect to $z_j$, $j=1,\ldots,n.$ By Theorem \ref{hartog1} it is thus jointly
$\alpha$-analytic. Conversely, suppose $f$ is $\alpha$-analytic.
We use induction in $q=\sum \alpha_j-n$. $q=0$ 
implies that $f$ is holomorphic o the sought representation is obvious. Assume $m>1$ and that the result holds true for $q=m-1$. Let $q=m$. W.l.o.g.\ assume $\alpha_1> 1$. By hypoellipticity we can assume $f$ is $C^\infty$-smooth. 
Then $\partial_{\bar{z}_1}f(z)$ satisfies the induction hypothesis, since it is $\alpha'$-analytic
with $\sum \alpha'_j-n=(\alpha_1-1)+\sum_{j=2}^n \alpha_j-n=q-1.$
Thus it has, by the induction hypothesis, a representation of the form 
$\sum_{0\leq \beta_j<\alpha'_j}a'_{\beta}\bar{z}^\beta$, where $a'_{\beta}(z)$ are holomorphic functions near $p_0,$
where $\alpha_j'=\alpha_j,$ $j\geq 2,$ $\alpha_1'=\alpha_1-1.$
Set $g:=\sum_{0\leq \beta_j<\alpha'_j}\frac{1}{\beta_1+1} a'_{\beta}(z)\bar{z}^\beta$. Then
$\partial_{\bar{z}_1}(f-\bar{z}_1 g)=0,$ and by the representation of $g$, $\partial_{\bar{z}_j}^{\alpha_j}(f-\bar{z}_1 g)=0$, $j=2,\ldots,n.$
Hence $f=(f-\bar{z}_1g)+\bar{z}_1g,$ is the sum of two functions both of which have the sought representation.
Also the identity theorem holds true for holomorphic functions in several variables and if $f$ vanishes on an open subset then the equations $\partial_{\bar{z}_1}^{\beta_1}f=\cdots =\partial_{\bar{z}_n}^{\beta_n} f=0, \beta_k \leq \alpha_k$,
$k=1,\ldots,n,$ imply that the components $a_\alpha(z)$ all vanish on an open neighborhood, thus on all of $\Omega.$ This verifies uniqueness of the representation. This completes the proof.
\end{proof}
Proposition \ref{prop1} can be proved more easily using an integral (Cauchy type) representation result, see Theorem \ref{theodorescuthm}.
Again note that we actually obtain a real-analytic representation and not merely $C^\infty$-smooth.

\section{$q$-analyticity in the sense of Ahern-Bruna}\label{ahernbrunasec}
Ahern \& Bruna \cite{ahernbruna}, p.132, give an alternative notion to that of polyanalyticity in several complex variables.
\begin{definition}[$q$-analytic functions in the sense of Ahern-Bruna]\label{ahernbrunadef}
	Let $\Omega\subset\Cn$ be an open subset and let $q\in \Z_+$. Denote by $\mbox{AB}_q(\Omega)$\index{$\mbox{AB}_q$, $q$-analytic functions in the sense of Aher-Bruna} the set of
	functions $f\in C^q(\Omega)$ that satisfy
	\begin{equation}
	\partial_{\bar{z}_1}^{\beta_1}\cdots \partial_{\bar{z}_n}^{\beta_n} f=0 \quad \forall \beta\in \N^n \mbox{ such that }\abs{\beta}=q
	\end{equation}
	We shall call such functions {\em $q$-analytic in the sense of Ahern-Bruna}\index{$q$-analyticity in the sense of Ahern-Bruna}.
	Ahern \& Bruna \cite{ahernbruna} point out that these are precisely the functions $f$
	that are annihilated by $\partial_{\bar{\zeta}}^q$ on the intersection of $\Omega$ with any complex line, parametrized by 
	$\zeta.$ Specifically for each $z_0\in \Omega$ and each $w_0\in \Cn$ we have
	\begin{equation}
	\partial_{\bar{\zeta}}^q f(z_0+\zeta w_0)=0, z_0+\zeta w_0 \in \Omega
	\end{equation}
\end{definition}
This (in a sense isotropic, i.e.\ independent of the individual coordinate components) notion will become interesting when studying polyanalyticity from the perspective of hypoanalytic theory (and also in the case of infinite dimensional
polyanalyticity).
Note that $\alpha$-analytic functions are defined via a system of elliptic equations whereas the defining equations of
$q$-analytic functions form a more complicated compound system. 
\begin{remark}
	Let $\Omega\subset\Cn$ be a domain subset and let $q\in \Z_+$. It is obvious that
	for $\alpha\in \Z_+^n,$ a function $f$ on $\Omega$ is $\alpha$-analytic on $\Omega$ if
	$f$ is $q$-analytic in the sense of Ahern-Bruna for $q=\sum_j \alpha_j.$
	But note that for each $n\in \Z_+$ and each $q\in \Z_+$ the function
	$f(z)=\sum_{j=1}^n \bar{z}_j^{q-1}$ is $(q,\ldots,q)$-analytic in $\Cn$ with $\abs{\alpha}=nq,$ whereas it is merely
	$q$-analytic in the sense of
	Ahern-Bruna.
	It is clear that any function $f(z)$ is $\alpha$-analytic on $\Omega$ if
	$f$ is $q$-analytic in the sense of Ahern-Bruna for $q=\max_j \alpha_j.$
	Not only is this condition not sharp but it will not imply that it is $q$-analytic on each complex line. Indeed, take the function $f(z)=\bar{z}_1^2 +\bar{z}_2^2
	+\bar{z}_1^2\bar{z}_2^2$ on $\C^2.$
	Clearly, $f$ is $(3,3)$-analytic. But on $\{z_2=z_1\}$ the term $\bar{z}_1^2\bar{z}_2^2$
	renders the restriction of $f$ to be $5$-analytic.
\end{remark}
We can give a sharp representation as follows.
\begin{proposition}\label{aherbrunacharact1prop}
	Let $\Omega\subset\Cn$ be a domain subset and let $q\in \Z_+$.
	A function $f$ on $\Omega$ is $q$-analytic in the sense of Ahern-Bruna
	if and only if $f$ has a representation of the form
	\begin{equation}\label{ekccharacahern}
	f(z)=\sum_{\beta\in \N^n,\abs{\beta}<q} a_\beta(z)\bar{z}^\beta, \quad \abs{\beta}:=\sum_{k=1}^n \beta_k
	\end{equation}
\end{proposition}
\begin{proof}
	Suppose $f$ has the representation in Eqn.(\ref{ekccharacahern}).
	Note that any function of the form $\bar{z}^\beta$ with $\sum_j\beta_j<q$
	is annihilated by
	$\partial_{\bar{z}_1}^{\gamma_1}\cdots \partial_{\bar{z}_n}^{\gamma_n}$
	for any $\gamma\in \Z_+^n$ satisfying $\sum_j \gamma_j = q.$
	Indeed, if $\partial_{\bar{z}_1}^{\gamma_1}\cdots \partial_{\bar{z}_n}^{\gamma_n} \bar{z}^\beta\not\equiv 0$
	then $\gamma_j\leq \beta_j$ for all $j=1,\ldots,n,$
	which would imply $q = \sum_j \gamma\leq \sum_j\beta_j<q.$
	This implies that each term in the representation of $f$ is annihilated by
	$\partial_{\bar{z}_1}^{\gamma_1}\cdots \partial_{\bar{z}_n}^{\gamma_n}$
	for any $\gamma\in \Z_+^n$ satisfying $\sum_j \gamma_j = q,$ hence $f$ is $q$-analytic in the sense of Ahern-Bruna.
	Conversely, suppose $f$ is $q$-analytic in the sense of Ahern-Bruna. Clearly, $f$ is, in particular separately
	$(q,\ldots,q)$-analytic, thus jointly $\alpha$-analytic with 
	$\alpha=(q,\ldots,q).$
	On a domain, $\Omega$, $f$ thus has a representation of the form 
	$f=\sum_{j=1}^n\sum_{\beta_j<q} a_\beta(z)\bar{z}^\beta.$ 
	Now any operator
	$\partial_{\bar{z}_1}^{\gamma_1}\cdots \partial_{\bar{z}_n}^{\gamma_n}$ for any $\gamma\in \Z_+^n$ satisfying $\sum_j \gamma_j = q,$
	must annihilate $f$, and since we have already seen that such operators annihilate
	all terms with $\abs{\beta}<q$ this gives
	$0\equiv \sum_{\abs{\beta}\geq q} a_\beta(z)\bar{z}^\beta.$
	Hence
	$f=\sum_{\abs{\beta}< q} a_\beta(z)\bar{z}^\beta.$ 
	This completes the proof.
\end{proof}

The following follows immediately from Proposition \ref{aherbrunacharact1prop}
\begin{proposition}
	Let $\Omega\subset\Cn$ be a domain subset and let $q\in \Z_+$.
	If $f$ is $\alpha$-analytic for some $\alpha\in \Z_+^n$ then it is $q$-analytic in the sense of
	Ahern-Bruna with $q=(\sum_j \alpha_j)-(n-1)$.
	Conversely if $f$ is $q$-analytic in the sense of
	Ahern-Bruna the it is $\alpha$-analytic for some $\alpha\in \Z_+^n$ satisfying 
	$q=(\sum_j \alpha_j)-(n-1).$
\end{proposition}

\section{The areolar derivative}\ref{areolarsec}
Polyanalytic functions can also be described in terms of so called areolar derivatives, see e.g.\ Pascali \cite{pascali}.
Pompieu \cite{pompieu} introduced in 1912 the {\em areolar derivative at $z=z_0$} of a continuous function $f(z)=u(x,y)+iv(x,y),$
$z=x+iy,$ according to
\begin{equation}
\frac{\partial f}{\partial \bar {z}}|_{z=z_0}=\lim_{\gamma\to z_0} \frac{\frac{1}{2\pi i}\int_{\gamma} f(z)dz}{\frac{1}{\pi}
\int_{\sigma} dx\wedge dy}
\end{equation}
where the limit is taken on all closed simple rectifiable curves $\gamma$ enclosing a domain $\sigma$, and continuously contractible to the point $z_0.$
If $f$ is $C^1$-smooth (for short $f\in C^1$) then
\begin{equation}
\frac{\partial f}{\partial \bar {z}}=\frac{1}{2}\left(\frac{\partial f}{\partial x}+i\frac{\partial f}{\partial y}
\right)=\frac{1}{2}\left(\left(\frac{\partial u}{\partial x}-\frac{\partial v}{\partial y}\right)+i\left(\frac{\partial u}{\partial y}-\frac{\partial v}{\partial x}\right)
\right)
\end{equation}
In the sense of distributions (i.e.\ independent of the existence of the partial derivatives of $u$ and $v$) the inhomogeneous equation
$\frac{\partial f}{\partial \bar {z}}=g(z)$ with a continuous function $g$, in a plane domain $\Omega$ takes the form
\begin{equation}
\frac{1}{2i}\int_{\gamma} f(z)dz-
\int_{\sigma} g(z)dx\wedge dy=0
\end{equation}
which holds for all pairs $(\gamma,\sigma)$ in $\Omega.$
Theodorescu \cite{theodorescu} extended the following representation to functions $f$ with merely continuous areolar derivative,
\begin{equation}
\frac{1}{2\pi i}\int_{\Omega} \frac{f(t)}{t-z}dt-
\frac{1}{\pi}\int_{\Gamma} \frac{\partial f}{\partial \bar{t}}\frac{d\eta d\nu}{t-z}=f(z),\quad z\in \Omega 
\end{equation}
where $t=\eta+ i\nu$ and $\Gamma$ consists of a finite number of simple closed rectifiable curves.

\subsection{The Pompieu operator}\label{areolarsec}
Theodorescu \cite{theodorescu} considered repeated application of the so called Pompieu integral operator\index{Pompieu integral operator} (see Pompieu \cite{pompieu}) defined for a bounded
measurable function $f$ on a domain $\Omega\subset \C$ by
the formula
\begin{equation}\label{pompdef0}
Tf(z)=-\frac{1}{\pi}\int_{\Omega} \frac{f(\zeta)}{\zeta-z}d\mu(\zeta)
\end{equation}
If $\Omega\subset\C$ is a domain bounded by a rectifiable Jordan curve $\gamma$
the Green formula yields for a function $f\in C^0(\Omega)\cap C^1(\Omega)$
\begin{equation}
\int_\Omega \partial_{\bar{z}} f dx\wedge dy=
\frac{1}{2i}\int_\gamma f(z)dz
\end{equation} 
\begin{equation}
\int_\Omega \partial_{z} f dx\wedge dy=%
-\frac{1}{2i}\int_\gamma f(z)d\bar{z}
\end{equation} 
thus for a fixed $z\in \Omega$ and sufficiently small $\epsilon>0$ we have
\begin{equation}
\int_\Omega \frac{\partial_{\bar{z}}f(\zeta)d\mu(\zeta)}{\zeta-z}-\frac{1}{2i}\int_{\abs{\zeta-z}=\epsilon}\frac{f(\zeta)d\zeta}{\zeta-z}=\int_{\Omega\cap \{\abs{\zeta-z}>\epsilon\}} \partial_{\bar{z}} f\frac{d\mu(\zeta)}{\zeta-z}
\end{equation} 
Letting $\epsilon\to 0$
we obtain
\begin{equation}\label{cauchypompformulaendim}
f(z)=\frac{1}{2\pi i}\int_\gamma \frac{f(\zeta)d\zeta}{\zeta-z}-\frac{1}{\pi}\int_\Omega \partial_{\bar{z}} f\frac{d\mu(\zeta)}{\zeta-z}
\end{equation} 
So if $f$ satisfies $\partial_{\bar{z}} f=g$ for $g\in C^1(\Omega)$ we have
\begin{equation}\label{tgreffen}
f(z)=\frac{1}{2\pi i}\int_\gamma \frac{f(\zeta)d\zeta}{\zeta-z}-\frac{1}{\pi }\int_\Omega g\frac{d\mu(\zeta)}{\zeta-z}
\end{equation} 
where we recognize $Tg$ in the last term.
\begin{lemma}\label{lemmanvibalktus0}
	Let $\Omega\subset\C$ be a domain. If $f\in L^1(\overline{\Omega})$ then 
	\begin{equation}
	\int_\Omega Tf \partial_{\bar{\zeta}}\phi 
	=-\int_\Omega f \phi , 
	\quad \forall \phi\in C^1_c(\Omega)
	\end{equation}
\end{lemma}
\begin{proof}
	By Eqn.(\ref{tgreffen}) we have (since $\phi$ vanishes on the boundary)
	\begin{equation}
	\phi(z)=-\frac{1}{\pi}\int_\Omega \frac{\partial_{\bar{\zeta}}\phi}{\zeta -z}d\mu(\zeta) =T(\partial_{\bar{z}}\phi)
	\end{equation}
	This implies 
	\begin{equation}
	\int_\Omega Tf \partial_{\bar{\zeta}}\phi =\frac{1}{\pi}\int_\Omega f(\zeta)d\mu(\zeta) \int_\Omega\frac{ \partial_{\bar{\zeta}}\phi}{\zeta -z} d\mu(\zeta)
	=-\int_\Omega f\phi 
	\end{equation}
	This completes the proof.
\end{proof}
\begin{lemma}\label{lemmanvibalktus}
	Let $\Omega\subset\C$ be a domain.
	If $f=\partial_{\bar{z}} g\in L^1(\overline{\Omega})$ then
	\begin{equation}
	g(z)=h(z)-\frac{1}{\pi}\int_\Omega \frac{f(\zeta)d\mu(\zeta)}{\zeta -z} (=h(z)+Tf(z))
	\end{equation}
	for a function $h$ holomorphic on $\Omega.$
	Conversely, if $h\in \mathscr{O}(\Omega)$ and $f\in L^1(\overline{\Omega})$ then the function
	$g=h+Tf$ has satisfies $\partial_{\bar{z}} g\in L^1(\overline{\Omega})$
	and $\partial_{\bar{z}} g =f.$
\end{lemma}
\begin{proof}
	The second assertion is obvious.
	For the first assertion we have 
	by Lemma \ref{lemmanvibalktus0}
	$\partial_{\bar{z}}(g-Tf)=f-f=0.$
	This completes the proof.
\end{proof}
\begin{theorem}\label{theodorescuthm}
For a $q$-analytic function $f$ at each point of $\overline{\Omega}$ (for $q\in \Z_+$ and a domain $\Omega$ bounded by a rectifiable closed contour $\Gamma$)
\begin{equation}\label{theodorescueq}
f(z)=\frac{1}{2\pi i}\sum_{j=0}^{q-1}\int_{\Gamma} \frac{(\bar{t}-\bar{z})^j}{j!(t-z)} \frac{\partial^j f}{\partial \bar{t}^j} dt, \quad z\in \Omega 
\end{equation}
\end{theorem}
The result is a fundamental proof-tool. We will give a proof in Section \ref{bvpsec} that includes this result, as it shows that polyanalytic functions are determined by the boundary values of their
$\overline{\partial}$ derivatives.
Here we give also the original method of proof of Theodorescu.
\begin{proof}
By definition there exists holomorphic functions $a_j(z)$ on a neighborhood of $\overline{\Omega}$ such that
$f(z)=\sum_{j=0}^{q-1} a_j(z)\bar{z}^j$ at each point of $\overline{\Omega}.$
Define for $j=0,\ldots,q-1$
\begin{equation}
\phi_j[f]:=f(z)-\frac{\bar{z}^1}{1!}\partial_z f(z)+\frac{\bar{z}^2}{2!}\partial_z^2 f(z)\cdots +(-1)^{q-1}\frac{\bar{z}^{q-1-j}}{(q-1-j)!}\partial_z^{q-1} f(z)
\end{equation}
This implies that
\begin{equation}
\phi_j[f]=\phi_0[\partial_z^j f]
\end{equation}
and thus
\begin{equation}
\frac{1}{j!}\phi_j[f]=a_j(z)
\end{equation}
which yields
\begin{equation}\label{theodoreqqen}
f(z)=\phi_0[f]+\frac{\bar{z}^1}{1!}\phi_1[f]+\cdots +\frac{\bar{z}^{q-1}}{(q-1)!}\phi_{q-1}[f]
\end{equation}
Applying the Cauchy formula to each $\phi_j[f]$ in Eqn.(\ref{theodoreqqen}) and changing the order of the terms appropriately this renders
\begin{multline}
f(z)=\frac{1}{2\pi i}\int_\Gamma \frac{f(\zeta)}{\zeta-z}d\zeta + \frac{1}{2\pi i}\int_\Gamma \frac{(\bar{z}-\bar{\zeta})^1}{1!}
\frac{f(\zeta)}{\zeta-z}d\zeta+\\
\cdots +\frac{1}{2\pi i}\int_\Gamma \frac{(\bar{z}-\bar{\zeta})^{q-1}}{(q-1)!}\frac{\partial_z^{q-1} f(\zeta)}{\zeta-z}d\zeta
\end{multline}
This completes the proof.
\end{proof}
Theodorescu \cite{theodorescu}, p.20, also points out that if conversely $f_0,\ldots,f_{q-1}$ are continuous functions on a domain $\Omega$ such that $\int_\gamma \phi_j[f(z)]dz=0$,
$\phi_j[f]:=f_j-\bar{z}f_{j-1}+\cdots+(-1)^{q-1-j}\frac{z^{q-1-j}}{(q-1-j)!}f_{q-1}(z),$ for each closed rectifiable contour $\gamma$ in $\Omega$, then $f_0$ is $q$-analytic and the functions $f_1,\ldots,f_{q-1}$ are the successive $\partial_{\bar{z}}$ derivatives of $f_0$. Indeed, since $\phi_{q-1}[f]=f_{q-1}(z)$ we have that $f_{q-1}$ is holomoprhic, thus by $\phi_{q-2}[f]=f_{q-2}-\bar{z}f_{q-1}$ we obtain that $\phi_{q-1}[f]$ is holomorphic, which implies that $f_{q-2}(z)$ is $2$-analytic so on. 
This is a Morera type of theorem for $q$-analytic functions.
\\
Now a slightly refined version of Theorem \ref{theodorescuthm} holds true.
\begin{theorem}\label{theodorescuthmrefin}
Let $\Omega$ be a domain  bounded by a rectifiable closed contour $\Gamma$, let $q\in \Z_+$.
Suppose $f$ is a $q$-analytic function on $\Omega$ such that each
$\partial_{\bar{z}}^j f$, $j=0,\ldots,q-1$ has continuous extension to
 $\overline{\Omega}$ (naturally the condition $\partial_{\bar{z}}^q f=0$ on $\Omega$
 implies that the latter has trivial continuous extension to $\overline{\Omega}$). Then
\begin{equation}
f(z)=\frac{1}{2\pi i}\sum_{j=0}^{q-1}\int_{\Gamma} \frac{(\bar{t}-\bar{z})^j}{j!(t-z)} \frac{\partial^j f}{\partial \bar{t}^j} dt, \quad z\in \Omega 
\end{equation}
\end{theorem}
\begin{proof}
First consider the case $q=2$. Let $f,g$ be 
$2$-analytic functions
on a bounded domain $\Omega\subset\C$ with boundary given by a rectifiable Jordan curve $\gamma,$
such that $\partial_{\bar{z}}^j f$, $\partial_{\bar{z}}^j g,$ $j=0,1,2$, have continuous extension to
$\overline{\Omega}$. 
Since we have $g\partial_{\bar{z}}^2 f-f\partial_{\bar{z}}^2 g =0$ on $\Omega.$
\begin{equation}
\int_\gamma (g\partial_{\bar{z}}^2 f-f\partial_{\bar{z}}^2 g)d\zeta =0
\end{equation}
Let $z\in \Omega.$ Let 
\begin{equation}
g(z,\zeta)=\frac{\bar{\zeta}-\bar{z}}{2(\zeta -z)}
\end{equation}
which is obviously a $2$-analytic function of $\zeta$ for $\zeta\neq z$ satisfying
\begin{equation}
\partial_{\bar{\zeta}} g(\zeta,z)=\frac{1}{\zeta -z}
\end{equation}
Note also that $g$ is bounded also when $z=\zeta.$
We replace $\Omega$ by $\Omega_\epsilon\setminus \{\abs{\zeta-z}<\epsilon\},$
so the new boundary is $\gamma\cup \{ \abs{\zeta-z}=\epsilon\},$ and the restrictions of $f,g$ to $\overline{\Omega}_\epsilon$ are now used (where $g$ is now $2$-analytic on $\Omega_\epsilon$). 
Thus for sufficiently small $\epsilon>0$ we have
\begin{multline}
\int_\gamma \left(\frac{\bar{\zeta}-\bar{z}}{2(\zeta -z)}  \partial_{\bar{\zeta}} f(\zeta) -f(\zeta) \frac{1}{\zeta -z}
\right)d\zeta=\\
\int_{\{ \abs{\zeta-z}=\epsilon\}} \frac{\bar{\zeta}-\bar{z}}{2(\zeta -z)}  \partial_{\bar{\zeta}} f(\zeta) d\zeta
-\int_{\{ \abs{\zeta-z}=\epsilon\}} \frac{f(\zeta)d\zeta}{\zeta -z}d\zeta
\end{multline}
where the integrand $\frac{\bar{\zeta}-\bar{z}}{2(\zeta -z)}  \partial_{\bar{\zeta}} f(\zeta)$
can be written as 
$\exp(2i\mbox{Arg}(\zeta-z)) \partial_{\bar{\zeta}} f(\zeta)$ thus, letting $\epsilon \to 0$ shows that the right hand side reduces to the Cauchy
integral reproducing $2\pi i f(z)$.
This proves that for $z\in \Omega$ we have
\begin{equation}
f(z)=\frac{1}{2\pi i} \int_\gamma \left(
\frac{\bar{\zeta}-\bar{z}}{2(\zeta -z)}  \partial_{\bar{\zeta}} f(\zeta) -f(\zeta) \frac{1}{\zeta -z}
\right)d\zeta
\end{equation}
Now consider the case of functions $f,g$, $q$-analytic in $\Omega$ such that
$\partial_{\bar{z}}^j f$, $\partial_{\bar{z}}^j g,$ $j=0,\ldots,q-1$, have continuous extension to
$\overline{\Omega}$.
 We have
\begin{multline}
g\partial_{\bar{z}}^q f=
\partial_{\bar{z}}(g\partial_{\bar{z}}^{q-1} f -
(\partial_{\bar{z}}g) \partial_{\bar{z}}^{q-1} f=\\
\partial_{\bar{z}}\left( g\partial_{\bar{z}}^{q-1} f- \partial_{\bar{z}}g\partial_{\bar{z}}^{q-2} f+\cdots+
(-1)^{q-1}(\partial_{\bar{z}}^{q-1} g)f\right)
+(-1)^q (\partial_{\bar{z}}^{q} g)f
\end{multline}
and using the notation $\partial_{\bar{z}}^{0}=1,$
\begin{equation}
g\partial_{\bar{z}}^q f +(-1)^{q-1} (\partial_{\bar{z}}^{q}g) f =
\partial_{\bar{z}}\left(\sum_{j=0}^{q-1}(-1)^{j}(\partial_{\bar{z}}^{j} g)
\partial_{\bar{z}}^{q-1-j} f
\right)
\end{equation}
We define the auxiliary functions
\begin{equation}
g_j(z,\zeta):=\frac{(\bar{\zeta}-\bar{z})^{j-1}}{2^{j-1}(j-1)!(\zeta -z)}
\end{equation}
which is obviously bounded on bounded sets, $j$-analytic when $z\neq \zeta$ and satisfies
\begin{equation}
\partial_{\bar{\zeta}} g_j =g_{j-1}
\end{equation}
Hence, completely analogous to the case $q=2$ we obtain
\begin{equation}
f(z)=\frac{(-1)^{q-1}}{2\pi i}\int_\gamma \left(\sum_{j=0}^{q-1} (-1)^j g_{q-j}(z,\zeta) \partial_{\bar{\zeta}}^{q-1-j} f(\zeta)\right) d\zeta
\end{equation}
This completes the proof.
\end{proof}

It was proved by Plemelj-Privalov (\cite{plemelj1}, \cite{plemelj2})
that for $\Gamma=\{\abs{z}=1\},$ the
Cauchy principal value integral
\begin{equation}
f(t)\mapsto \frac{1}{\pi i}\int_{\Gamma}\frac{f(\zeta)}{\zeta -t}, \quad t\in \Gamma 
\end{equation}
behaves invariantly with respect to the H\"older class $C^{0,\alpha}(\Gamma),$ $\alpha\in (0,1).$
This was generalized to piecewise smooth curves by Muskhelishvili \cite{muskesh} and by Salaev \cite{salaev} to 
closed rectifiable Jordan curves satisfying that there is a constant $c >0$ such that for all $p_0\in \Gamma$ and
$r>0$, the arc-length measure of $\Gamma\cap \{\abs{z-p_0}<r\}$ is at most $cr$. 
Vekua \cite{vekua} proved for $p>2,$ $T$ is a continuous linear operator $L^p(\Omega)\to C^{\kappa}(\overline{\Omega})$
the H\"older space with exponent $\kappa=\frac{p-2}{p}.$ This can be derived from the following theorem.
\begin{theorem}\label{vekuaineqthm}
Let $\Omega\subset\C$ be a bounded domain and $f\in L^p(\overline{\Omega}),$ for some integer $p>2.$
Then the Pompieu operator acting on $f$ according to Eqn.(\ref{pompdef0}) 
satisfies that there exists constants $M_1,M_2>0$, $M_1$ depending only on $p,\Omega$ and $M_2$ depending only on $p$, such that
\begin{equation}\label{vekuaforsta0}
\abs{Tf(z)}\leq M_1 \norm{f}_{L^p(\overline{\Omega})}
\end{equation}
\begin{equation}\label{vekuaandra0}
\abs{Tf(z_1)-Tf(z_2)}\leq M_2 \norm{f}_{L^p(\overline{\Omega})}\abs{z_1-z_2}^\alpha,\quad \alpha=\frac{p-2}{p}
\end{equation}
for arbitrary points $z_1,z_2$.
\end{theorem}
\begin{proof}
	By the H\"older inequality we have for $q$ satisfying $1/p+1/q =1$
	\begin{equation}
	\abs{Tf(z)}\leq M_1\frac{1}{\pi}\left(\int_\Omega \abs{f}^p\right)^{\frac{1}{p}}\left(\int_\Omega \abs{\zeta-z}^{-q}\right)^{\frac{1}{q}}
	\end{equation}
Since $q<2$ we have
	\begin{equation}
\frac{1}{\pi}\left(\int_\Omega \abs{\zeta-z}^{-q}\right)^{\frac{1}{q}}\leq 
\frac{1}{\pi}\left(\frac{2\pi}{\alpha q}\right)^{\frac{1}{q}} (\mbox{diam}(\Omega))^\alpha=:M_1
\end{equation}
where $\mbox{diam}(\Omega)$ denotes the diameter and $\alpha=(p-2)/p.$
Hence Eqn.(\ref{vekuaforsta0}) follows from Eqn.(\ref{vekuaandra0}).
Applying the H\"older inequality to 
\begin{equation}
Tf(z_1)-Tf(z_2)=\frac{z_1-z_2}{\pi}\int\frac{f(\zeta)d\mu(\zeta)}{(\zeta-z_1)(\zeta -z_2)},\quad z_1\neq z_2
\end{equation}
gives
\begin{equation}\label{veksexfyra}
\abs{Tf(z_1)-Tf(z_2)}\leq \norm{f}_{L^p(\overline{\Omega})}\frac{\abs{z_1-z_2}}{\pi}
\left(\int_\Omega \abs{\zeta-z_1}^{-q}\abs{\zeta-z_2}^{-q}d\mu(\zeta)\right)^{\frac{1}{q}}
\end{equation}
Let $C_1=\{\abs{\zeta-z_1}= 2\abs{z_1-z_2}$
and $r_0>0$ such that $\overline{\Omega}\subset C_0:=\{\abs{\zeta-z_1}= r_0\}.$
For
$\zeta$ outside of $C_1$ we have $2\abs{\zeta-z_2}\geq \abs{\zeta-z_1},$ which implies
that for $\alpha<2,\beta<2,$
\begin{multline}
A_1:=\int_{C_0\setminus C_1}\Omega \abs{\zeta-z_1}^{-\alpha}\abs{\zeta-z_2}^{-\beta}d\mu(\zeta)\leq \pi 2^{1+\beta}\int_r^{2r_0} r^{1-\alpha-\beta} dr<\\
\left\{
\begin{array}{ll}
\frac{8\pi\abs{z_1-z_2}^{2-\alpha-\beta}}{\alpha+\beta-2} & , \alpha+\beta >2\\
8\pi \log\frac{r_0}{\abs{z_1-z_2}} & , \alpha+\beta =2\\
\frac{32}{2-\alpha-\beta}r_0^{2-\alpha-\beta} & , \alpha+\beta <2\\
\end{array}
\right.
\end{multline}
Also we have for a constant $M_{\alpha,\beta}>0$ depending on $\alpha,\beta$
\begin{multline}
A_2:=\int_{C_1} \abs{\zeta-z_1}^{-\alpha}\abs{\zeta-z_2}^{-\beta}d\mu(\zeta)=\\
\frac{1}{\abs{z_1-z_2}^{\alpha+\beta-2}}\int_{\abs{\zeta}<2} \frac{d\mu(\zeta)}{\abs{\zeta}^\alpha \abs{\zeta -\exp(i\theta)}^\beta}\leq
	\frac{M_{\alpha,\beta}}{\abs{z_1-z_2}^{\alpha+\beta-2}}
		\end{multline}
		Now
		\begin{equation}
		J_{\alpha,\beta}:=\int_{\Omega} \abs{\zeta-z_1}^{-\alpha}\abs{\zeta-z_2}^{-\beta}d\mu(\zeta)\leq A_1+A_2
		\end{equation}
		implies that there exist constants $M_{\alpha,\beta}'>0,$ $M_{\alpha,\beta,\Omega}''>0$ such that
		\begin{equation}\label{casevek123}
		{J_{\alpha,\beta}\leq}
		\left\{
		\begin{array}{ll}
		M_{\alpha,\beta}'\abs{z_1-z_2}^{2-\alpha-\beta} & , \alpha+\beta >2\\
	M_{\alpha,\beta,\Omega}''+8\pi\abs{ \log\abs{z_1-z_2}} & , \alpha+\beta =2\\
	M_{\alpha,\beta,\Omega}'' & , \alpha+\beta <2\\
	\end{array}
	\right.
	\end{equation}
	Since $1<q<2$ the first case in Eqn.(\ref{casevek123}) together with Eqn.(\ref{veksexfyra}) yields
	for a constant $M_p>0$ depending upon $p$
	\begin{equation}
	\frac{\abs{z_1-z_2}}{\pi}
	\left(\int_\Omega \abs{\zeta-z_1}^{-q}\abs{\zeta-z_2}^{-q}d\mu(\zeta)\right)^{\frac{1}{q}}\leq M_p\abs{z_1-z_2}^{\frac{p-2}{p}}
	\end{equation}
	This completes the proof.
\end{proof}
\begin{definition}[See De la Cruz et al.\ \cite{delacruz}]
Let $E\subset \R^2$ be a closed subset and let $k\in \N$ and $\alpha\in (0,1].$
A real-valued function $f$ on $E$ is said to belong to $\mbox{Lip}(E,k+\alpha)$
if there exists real-valued bounded functions $f^{(j)}$, $0<\abs{j}\leq k,$ on $E$, with $f^{(0)}=f,$ such that
\begin{equation}
R_j(x,y)=f^{(j)}-\sum_{\abs{j+l}\leq k} \frac{f^{(j+1)}(y)}{l!}(x-y)^l,\quad x,y\in E
\end{equation}
\begin{multline}\label{normekv}
\abs{f^{(j)}(x)}\leq c,\quad \abs{R_j(x,y)}\leq c\abs{x-y}^{k+\alpha-\abs{j}},\quad x,y\in E,\\
\abs{j}\leq k-\sum_{\abs{j+l}\leq k} \frac{f^{(j+1)}(y)}{l!}(x-y)^l,\quad x,y\in E
\end{multline}
for a positive constant $c$. 
$\mbox{Lip}(E,k+\alpha)$ is equipped with the norm $\norm{\cdot}_{k,\alpha}$ defined as the smallest $c$ 
in Eqn.(\ref{normekv}), and turns $\mbox{Lip}(E,k+\alpha)$ into a Banach space.
\end{definition}
Based upon Theorem \ref{theodorescuthm}, De la Cruz et al.\ \cite{delacruz} 
define the operator
\begin{equation}
\mathcal{S}^{(k)} f(t):=\sum_{j=0}^k \frac{1}{\pi i} \int_\Gamma \frac{(\overline{t-\zeta})^j}{j!(\zeta-t)} 
f^{(0,j)}(\zeta)d\zeta,\quad t\in \Gamma
\end{equation}
and prove the following extension of a theorem of Plemelj-Privalov to the polyanalytic case.
\begin{theorem}
$\mathcal{S}^{(k)}(\mbox{Lip}(E,k+\alpha))\subset \mbox{Lip}(E,k+\alpha),\quad \alpha \in (0,1)$
\end{theorem}
See Chapter \ref{bvpsec} and Chapter \ref{reproducingsec} for more on the use of versions of the Pompieu operators and integral kernels in the theory of $q$-analytic functions.

\section{$\alpha$-analytic functions on complex manifolds}

\begin{definition}
Let $X$ be a connected, second countable, Hausdorff space. 
A {\em local chart} on $X$ at $p\in X$ is a pair $(U, \phi)$ where $p\in U,$ $U \subset X$ is an open subset and $\phi : U \to \phi(U)\subseteq K^n$ (where $K$ is the field $\C$ or $\R$ respectively) is a homeomorphism for some positive integer $n$. the components of $\phi= (\phi_1,\ldots ,\phi_n)$ are called local coordinates on $U$. When $K=\R,$ $\phi$ is a diffeomorphism and when $K=\C$, $\phi$ is a biholomorphism.
A $C^k$-smooth (or $C^{\omega}$ or holomorphic) {\em atlas}, $\mathcal{A},$ 
on $X$ is given by a collection of charts $\{(\phi^j,U^j)\}_{j\in I}$ for an index set $I$ such that each chart has the same dimension, say $n$, in the target space, $X$ is
covered by the $U_j$ and the transition maps 
$\phi^k (U_j\cap U_k)\circ (\phi^j)^{-1} (U_j\cap   U_k),$ are $C^k$-smooth (or $C^{\omega}$ or biholomorphic) for all $j,k\in I.$ In that case we say that the transition maps are are {\em compatible}. Here $C^{\omega}$ denotes the set of real-analytic functions.
Of course it is required that $K=\C$ (i.e.\ the transition maps are biholomorphisms) for the atlas to be called holomorphic/complex analytic. 
A manifold of dimension $n$ is a doublet $(M,\mathcal{A})$ where $M$ is a connected, second countable, Hausdorff space and $\mathcal{A}$ a given atlas on $M$. When $K=\C$, $(M,\mathcal{A})$ is called a complex manifold of dimension $n.$,
A $C^{\omega}$ or holomorphic chart $(\phi,U)$ on $M$ is called 
{\em compatible} with the atlas, $\mathcal{A},$ of $M$ if given any chart
$(\psi,V)\in \mathcal{A},$ such that $U\cap V\neq \emptyset$, 
$\phi (U\cap V)\circ\psi^{-1} (U\cap V)$ is a homeomorphisms.
\end{definition}
In this book we shall always assume any manifold is Hausdorff and paracompact.

\begin{definition}[Submanifolds]
An $m$ dimensional {\em submanifold} of an $n$ dimensional manifold $(M,\mathcal{A}_M)$ is a 
closed connected subset $N\subseteq M$ 
with atlas $\mathcal{A}_N$, each chart of which
can be extended to a chart on $M$ compatible with 
$\mathcal{A}_M$. For the subset $N\subset M$ there exists a natural atlas obtained by the restrictions of the charts from 
$\mathcal{A}_M,$ this atlas obviously turns $N$ into a submanifold. 
Thus whenever a submanifold is presented it is important that the associated atlas be given. If we do not specify the atlas then we shall always assume that the submanifold has the {\em natural} atlas, by which we mean that each of its charts is induced from the supermanifold as a restriction (see above). 
When this natural choice of atlas for $N$ is made and when it is clear from the context what the atlas of $M$ is, then we abbreviate $N\subset M$ and call $N$ a {\em natural} submanifold of $M$. If we do not specify an atlas for a submanifold $N$ of $M$ then we shall always assume the natural (induced by restriction) atlas and simply say that $N\subset M$ is a {\em submanifold}.
When it is clear from the context we shall not always write out the atlas of standard manifolds such as $\Cn$.
As an example, if we write that
$M\subset \Cn$ is a submanifold without specifying any atlases then we assume that $M$ is a natural submanifold and as such
any local chart of $M$ is given by the restriction of the ambient coordinate function for $\Cn.$
\end{definition} 
\index{Manifold}\index{Local chart}\index{Transition function}
\index{Atlas}\index{$C^k$ functions on manifolds}

\begin{definition}[Holomorphic function on manifold]
Let $\Omega$ be a complex manifold and let $f\in C^1(\Omega).$
Then $f$ is called holomorphic at $p_0\in U$ if
for some (and hence any) local chart $(U,\phi)$ in the atlas of $\Omega$, such that 
$p_0\in U,$ the function $f\circ \phi^{-1}$ is holomorphic in the usual sense, near $\phi(p_0)$.
\end{definition}
\begin{remark}
It is clear that this is independent of local chart, i.e.\
if $(V,\psi)$ is another local chart in the atlas of $\Omega$, such that 
$p_0\in U,$ then on $U\cap V$ we have after a coordinate change given by the diffeomorphism $\psi\circ\phi^{-1}$
that the function $f\circ \psi^{-1}$ has a local convergent power series expansion at $p_0$, in the new variable $(\psi\circ\phi^{-1})(\phi (p))=\psi(p),$ for $p$ near $p_0$ in $\Omega.$
\end{remark}
For a complex manifold $M$ it is well known (see e.g.\ Gunning \cite{gunning}, p.25) that the atlas of $M$ induces a complex atlas for the {\em (tangent) vector bundle,} $TM:=\bigcup_{p\in M} T_p M$. The projection $\pi\colon TM\to M$ is a smooth map, and a {\em section} of $TM$ is
is a smooth map $X\colon M \to TM$ satisfying
$\pi\circ X =\mbox{id}_M.$
A section of a vector bundle is also called a {\em vector field}.
When $\Omega$ is a complex manifold and there is a {\em global frame}, i.e.\ global sections $Z_1,\ldots, Z_n$, of the tangent bundle, s.t.\ $dZ_1\wedge\ldots\wedge dZ_n (p)\neq 0$ at each point $p\in\Omega$, then we can globally define the operators $Z_j^{\alpha_j}$ (which plays the role of $\partial_{\bar{z}_j}^{\alpha_j}$), and in a natural way call functions $\alpha$-analytic whenever they are annihilated by the conjugate counterpart of these differential operators $j=1,\ldots,n$ (the point is that the order is fixed, the $\alpha_j$ are uniquely associated to a fixed basis vector field $Z_j$). For instance, if in one complex dimension we have $s=\partial x +i\partial y$ as a global $C^\infty$-smooth section of $\C\otimes T\Omega$, then by definition this varies $C^\infty$-smoothly with respect to the base point and it is a linear map $C^\infty(\Omega) \to C^\infty(\Omega)$, and the same holds for its multi-powers in $(\Gamma(\Omega,\C\otimes T\Omega))^\alpha$. 
Obstruction to global sections can be described in terms of de Rham cohomology groups (or alternatively Chern classes), see e.g.\ Wells \cite{wells} p.84.
\\
Let $M$ be any complex manifold. Recall that
the complex tangent space $\C\otimes TM$ has a splitting into a $\C$-linear ($T^{1,0}M$) and anti-$\C$-linear ($T^{0,1}M$) part respectively, i.e.\
$\C\otimes M=T^{1,0}M + T^{0,1}M.$
Holomorphic functions are precisely those that are annihilated by all smooth local sections of $T^{0,1}M.$ Note that this is independent of local choice of charts. As we have pointed out, for any such local section, say on $V\subset M$, its 
multi-power in $(\Gamma(V,T^{0,1}M))^\alpha$ is a
linear map $C^\infty(V) \to C^\infty(V)$, and determines locally a differential operator that generalizes the standard operator $\overline{\partial}^{\alpha}.$
\begin{definition}[$q$-analytic function on a complex one dimensional manifold]\index{$q$-analytic function on a complex one dimensional manifold} 
Let $M$ be a complex one dimensional manifold and let $q\in \Z_+.$ 
Recall that
the complex tangent space $\C\otimes TM$ has a splitting into a $\C$-linear ($T^{1,0}M$) and anti-$\C$-linear ($T^{0,1}M$) part respectively, i.e.\
$\C\otimes M=T^{1,0}M + T^{0,1}M.$
A $C^\infty$-smooth function $f$ is called $q$-analytic at $p\in M$ if $L^q f=0$ near $p$ in $M$, for any smooth local section $L$ of $T^{0,1}M$ near $p.$ $f$ is called $q$-analytic on an open subset $U\subseteq M$ if it is $q$-analytic at each point of $U.$
\end{definition} 
In contrast to the case of holomorphy, being
$\alpha$-analytic requires a specific ordered basis for
the local sections of the anti-$\C$-linear part of the tangent space, in particular the condition required will be atlas dependent.
\begin{definition}[$\alpha$-analytic function on complex finite dimensional manifolds]\index{$\alpha$-analytic function on a complex submanifold of $\Cn$}  
Let $\Omega\subset \Cn$ be a complex submanifold 
and let $\alpha\in \Z_+^{n}.$ 
A $C^\infty$-smooth function $f$ is called $\alpha$-analytic at $p$ if 
there exists a 
chart at $p$, in the atlas (with induced coordinate $z$) of $\Omega$, with respect to which 
$f$ is a solution near $p$ to the system $\partial_{\bar{z}_1}^{\alpha_1}f=\cdots =\partial_{\bar{z}_n}^{\alpha_n}f=0$.
\end{definition} 
\begin{remark}[Notations and assumptions regarding atlases of ambient space and submanifolds]
In the context of $\alpha$-analytic functions it is important to be aware of the atlas of a given manifold. We shall try to specify the atlas as often as possible but for convenience we shall make the following assumptions:
Whenever we write $\Cn$ without specifying an atlas, then we shall assume that this is $\Cn$ equipped with the standard complex structure (in particular equipped with an atlas consisting of a single {\em global} complex analytic chart).
Whenever $(M,\mathcal{A})$ is a given manifold and we call
$N\subset M$ a (not necessarily complex) {\em submanifold} without specifying an atlas for $N$, then we assume that $N$ is a natural submanifold (by which we mean that it is equipped with the atlas induced from $M$, in which each chart of the atlas of $N$ can be identified as the restriction of a chart in the atlas of $M$).
For example, if we write that $M\subset \Cn$ is a complex submanifold of dimension $n$, 
then the restriction, to $M$, of any $\alpha$-analytic function on $\Cn$, is again an $\alpha$-analytic function.
\end{remark}
Clearly, on each local chart Proposition \ref{prop2} yields unique representation, and the local uniqueness in the case of Euclidean manifolds immediately yields the identity principle.
\begin{corollary}[to Proposition \ref{prop2}]
Let $\Omega$ be a connected complex $n$-dimensional manifold. Any 
$\alpha$-analytic function on $\Omega$ has a real-analytic representation on $\Omega.$ Furthermore, any $\alpha$-analytic function on $\Omega$ that vanishes on an open subset vanishes identically. 
\end{corollary}

\section{Classical extension}\label{extensionsec}
As a consequence of the unique decomposition of $\alpha$-analytic functions, as described in Proposition \ref{prop2}, we have the following extension theorem.
\begin{theorem}\label{extensionthm1}
Given a domain $\Omega\subset \Cn,$ $n\geq 2,$ and a compact subset $K \subset \Omega$ such that $\Omega \setminus K$ is connected, Any bounded member of $\mbox{PA}_\alpha(\Omega \setminus K)$, extends uniquely to a member of $\mbox{PA}_\alpha(\Omega).$ 
\end{theorem}
\begin{proof}
The result is well-known in the case of $\abs{\alpha_j}=1,$ for all $1\leq j\leq n.$ For general $\alpha\in\Z_+^n$, assume without loss of generality $p_0=0\in \Omega$ and 
let $f\in\mbox{PA}_\alpha(\Omega \setminus K)$ have representation $f=\sum_{1\leq \beta_j <\alpha_j} a_\beta(z) \bar{z}^\beta.$ By the known result for holomorphic functions, each $a_\beta(z)$ extends uniquely to a holomorphic function $\tilde{a}_\alpha$ on $\Omega$. Hence the function $\tilde{f}:=\sum \tilde{a}_\beta(z) \bar{z}^\beta$ defines an $\alpha$-analytic extension of $f$. This completes the proof.
\end{proof}
Obviously, obstruction to holomorphic extension implies obstruction to $\alpha$-analytic extension for any $\alpha\in\Z_+^n$ where a non-holomorphic counterexample function can be found for each $\alpha$ satisfying $\max_j \alpha_j >1.$
\begin{definition}[Pseudoconvex domain]\label{pseudoconvexdef}\index{Pseudoconvex domain}
Let $\Omega\subset \Cn$ be a domain with smooth boundary and smooth defining function
$\rho$, in the sense that
$\Omega :=\{ \rho(z)<0\},$ $\abs{d\rho}\neq 0$ on $\partial\Omega.$
Recall that $h\colon V\times V\to \C$
on a complex vector space $V$ is called a {\em Hermitian form}\index{Hermitian form} if it is linear in the first coordinate and such that
$h(X,Y)=\overline{h(Y,X)}.$ Every Hermitian form
has an associated Hermitian matrix, $A,$ such that $h(X,Y)=X A \overline{Y}^T$.
For each $p_0\in \partial \Omega$, the matrix
\begin{equation}
\begin{bmatrix}
\frac{\partial^2\rho(p_0)}{\partial z_j\partial \bar{z}_k}
\end{bmatrix}_{1\leq j,k\leq n}
\end{equation}
defines a Hermitian form on $T_{p_0}^{(1,0)}(\Cn),$ and its restriction to $\partial \Omega$ defines
a Hermitian form, $\mathcal{L}_{p_0},$ on $T_{p_0}^{(1,0)}(\partial\Omega),$ defined as follows:
For $X=\sum_{j=1}^n a_j \partial_{z_j}$ and $Y=\sum_{j=1}^n b_j \partial_{z_j}$ are two complex tangent vectors, at $p_0,$
$\mathcal{L}_{p_0}(X,Y):=\sum_{j,k=1}^n a_j\bar{b}_k\frac{\partial^2\rho(p_0)}{\partial z_j\partial \bar{z}_k}$.
This Hermitian form is real-valued\footnote{In fact can be identified as $\mathcal{L}_{p_0}(X,Y):=\frac{1}{2i}
[\overline{\tilde{X}},\tilde{Y}]_{p_0}$mod$\C\otimes T^c_{p_0} \partial\Omega,$ where $T^c_{p_0}\partial \Omega=T_{p_0}\partial \Omega \cap J(T^c_{p_0}\partial \Omega),$ for the, fiberwise defined, complex structure map $J.$}
If $\mathcal{L}_{p_0}$ is non-negative definite (positive definite) at $p_0\in \partial\Omega,$ then 
$\Omega$
 is called {\em pseudoconvex} ({\em strongly pseudoconvex}) at $p_0$\index{Pseudoconvexity}.
$\Omega$ is called pseudoconvex (strongly pseudoconvex) is $\mathcal{L}_{p_0}$ is non-negative definite (positive definite)
at each point $p\in\partial\Omega.$
\end{definition}
\begin{definition}[Domain of holomorphy]
A domain $\Omega \subset\Cn$ is called a {\em domain of holomorphy} if there exists a holomorphic function $f:\Omega\to \C$ such that for every
$z_0\in \Omega,$ the function $f|_{B(z_0,r_0)}$, where $r_0=\max\{ r>0: B(z_0,r)\subset\Omega\},$ $B(z_0,r):=\{\abs{z-z_0}<r\},$
cannot be extended holomorphically to $\Omega\cup \{\abs{z-z_0}<r'\},$ for some $r'>r_0.$
\end{definition}
\begin{theorem} A domain $\Omega\subset\Cn$ is (Levi) pseudoconvex iff $\Omega$ is domain of holomorphy.
\index{Domain of holomorphy}
\end{theorem}

If $\Omega\subset \C$ is a non-pseudoconvex domain, 
then it is a is a well-known result in complex analysis of several variables that it is not a domain of holomorphy, thus there exists a point of the boundary such that all holomorphic functions extend to an open neighborhood of that point. 
As pointed out by Gunning \cite{gunning}, it is not necessarily true that this implies that the functions are extended to an open set whose interior lies in the complement of $\Omega$ since e.g.\ $\Omega \cap U$ can be dense in $\Omega$. Furthermore, holomorphic extension and in particular maximal holomorphic extension renders in general many-valued functions and in order to continue working with single-valued functions one can define Riemann domains as the natural structures for envelopes of holomorphic.\index{Envelope of holomorphy} 
We can point out that for 
certain special types of domains in $\Cn$, their envelopes of holomorphy has a global chart and therefore
the proof of Theorem \ref{extensionthm1} can be mimicked to obtain extension to the envelope of holomorphy. 
Two such examples are \index{Reinhart domain}Reinhart domains and tube domains\index{Tube domain} (the latter being unbounded) respectively, see e.g.\ Gunning \cite{gunning}, Vol 1, p.71.
More can be said regarding the case of {\em bounded} domains.
First recall the definition of Schlicht Riemann domains.
\begin{definition}
A {\em Riemann domain}\index{Riemann domain} of dimension $n$ is a connected Hausdorff topological space $M$
together with a 
local homeomorphism 
$\pi\colon M\to \Cn,$ called the {\em projection}\index{Projection} of $M$ over $\Cn.$ 
A Riemann domain is called {\em schlicht}\index{Schlicht Riemann domain} if it is single-sheeted.
\end{definition}
\index{Schlicht Riemann domain}
\begin{theorem}\label{extensiontm2} Let $\Omega\subset\Cn$ be a bounded domain and assume that the envelope of holomorphy $\widehat {\Omega}$ of $\Omega$ is schlicht over $\Omega$. Then $\alpha$-analytic extension to the envelope of holomorphy holds true.
\end{theorem} 
\begin{proof}
Jupiter \cite{jupiter} has shown that if $\widehat{\Omega}$ is schlicht over $\Omega$ then it is a schlicht Riemann domain.
As pointed out in Shabat \cite{shabat}, p.213, any function $g\in \mathscr{O}(\Omega)$ satisfies
$g(\widehat{\Omega})\subseteq g(\Omega),$ and this applied to the coordinate functions $g=z_j,$ $j=1,\ldots ,n,$ implies (since $\Omega$ is bounded) that $\widehat{\Omega}$ is a bounded domain in $\Cn.$ In particular, there exists a {\em global} holomorphic (coordinate) chart, say $z$, for $\widehat{\Omega}$. As in the proof of Theorem \ref{extensionthm1}, each $a_\beta(z)$ extends uniquely to a holomorphic function $\tilde{a}_\alpha$ on $\Omega$. Hence the function $\tilde{f}:=\sum \tilde{a}_\beta(z) \bar{z}^\beta$ defines an $\alpha$-analytic extension of $f$. This completes the proof. 
\end{proof}
Proposition \ref{uniquepointprop} obviously has an immediate extension to the several variable case with analogous proof.
As a consequence we obtain the following.
\begin{proposition}
Let $\Omega\subset\Cn$ be a domain, let $\alpha\in \Z_+^n$ and let $f$ and $g$ be two $\alpha$-analytic functions on $\Omega.$ The set $A$ of points, $p$ such that $f=g$ on a neighborhood of $p$, is open and closed.
\end{proposition}
\begin{proof}
Obviously $A$ is open.
If $\{z_j\}_{j\in \Z_+}$ is a sequence $A$ such that $z_j\to p\in \Omega$ then the limit of each $\partial_z^\beta\partial_{\bar{z}}^\gamma f(z_j)$
and $\partial_z^\beta\partial_{\bar{z}}^\gamma g(z_j)$ respectively
tends to $\partial_z^\beta\partial_{\bar{z}}^\gamma f(p)$
and $\partial_z^\beta\partial_{\bar{z}}^\gamma g(p)$
as $z_j\to p.$ Hence 
$\partial_z^\beta\partial_{\bar{z}}^\gamma (f-g)(p)=
0$
for all multi-integers
$\beta,\gamma\in \Z_{\geq 0}^n,$ which implies that $(f-g)\equiv 0$ near $p$ thus $p\in A.$ This completes the proof. 
\end{proof}

\begin{proposition}
Let $\Omega\subset\Cn$ be a domain, let $\Omega'\subset\Omega$ be a subdomain and let $\alpha\in \Z_+^n$. If $f$ is an $\alpha$-analytic function of exact order $\alpha$ (i.e.\ for each $j$ there is an analytic component of $f$ corresponding to some $\bar{z}^\beta$ such that
$\beta_j=\alpha_j-1$ which does not vanish identically) on $\Omega$
then $f|_{\Omega'}$ is an $\alpha$-analytic function of exact order $\alpha$. Conversely if $g$ is an $\alpha$-analytic function of exact order $\alpha$ on $\Omega'$ which has a (necessarily unique) polyanalytic extension, $G$, to
$\Omega$, then $G$ is $\alpha$-analytic.
\end{proposition}
\begin{proof}
We know that $f$ has a representation of the form 
$f(z)=\sum_{\beta_j<\alpha_j} a_\beta(z) \bar{z}^\beta,$
where the $a_\beta$ are the holomorphic components. 
If the restriction to $\Omega'$ was $\alpha'$-analytic where
for some $j_0$ we had $\alpha'_{j_0}<\alpha_{j_0},$
then the restriction of the analytic components corresponding to $\bar{z}^\beta$ such that
$\beta_{j_0}=\alpha_{j_0}-1$ would vanish on an open subset thus identically, contradicting that the order of $f$ is exact. Similarly, if $G$ where to increase the order of analyticity of $g$ then it would require the introduction of a new analytic component that vanished on $\Omega'.$ This completes the proof. 
\end{proof}
Tarkhanov \cite{tarkhanov}, Sec 1.1.1, gives a clear description on how to naturally equip the sheaves of germs of solutions to a partial differential operator $L$ with a topology that will be Hausdorff iff and only if each homogeneous (distribution) solution $u$, to $Lu=0$
satisfies that if $u$ vanishes on an open subset then $u\equiv 0$ (an identity principle which we have established for the case of $\alpha$-analytic functions). 
Let $\Omega\subset\Rn$ be an open subset, let $L$ be a differential operator on $\Omega$ and denote by $P_0$ the set of distribution solutions $u$ to $Lu=0.$ The correspondence $\Omega \to P_0$ determines a sheaf and the covering space, $P,$ over $\Omega$ consists of germs of distribution solutions to $Lu=0.$ The topology on $P_0$ can be defined as follows. Let $p\in \Omega$ and $[u]_p$ be a germ in $P$
and let $(u,U)$ be a representative, i.e.\ $U\subset\Omega$ is an open neighborhood of $p$. For each $x\in U$ denote by $u_x$ the germ of the solution $u$ at $x.$ Then $u_p$ has a neighborhood in $P$ of the form $\bigcup_{x\in U} u_x$. The projection $\pi:P\to \Omega$ sends the germ $u_x$ in $P$ to $x$ and this projection will be a local homeomorphism. The operation of addition and multiplication of germs by scalars are readily induced and a continuous section of $P$ on an open set $U$ can be identified as a distribution solution to $Lu=0$ on $U.$
Let us run through the necessary steps for the particular case of $\alpha$-analytic functions.
\begin{definition}\label{germsdefforsta}
For $\alpha\in \Z_+^n$ and $p\in \Cn$ denote by $O_p^\alpha$ the set of germs of $\alpha$-analytic functions at $p$.
 Define $O^\alpha:=\bigcup_{p} O_p^\alpha,$ $p\in \Cn.$ An equivalence class, $[f]_p,$ for a germ at $p$ has a representative in terms of a doublet $(f,U)$ where $U$ is an open neighborhood of $p$ and
 $f$ an $\alpha$-analytic function on $U.$ 
By definition, the set $\{ f_z :z\in U\}$ defines an open neighborhood of $[f]_p$, which thus renders a fundamental system of open neighborhoods of $f_a$ as $(f,U)$ varies over all doublets in the definition of $[f]_p.$ By what we have described in the preceding paragraph this induces a topology
on $O^\alpha$. We also introduce the space $O_p$ of germs of polyanalytic functions (of arbitrary order) at $p$ and set
$O:=\bigcup_p O_p.$ 
It is clear that $O=\bigcup_{\alpha\in \Z_+^n}O^\alpha.$ The space $O$ can completely analogously be equipped with a natrular topology by the same procedure as for $O^\alpha.$
\end{definition}
\begin{proposition}
The map $\pi:O \to \Cn$, $\pi([f]_p)=p$ is a local homeomorphism.
\end{proposition}
\begin{proof}
Let $U$ be an open connected neighborhood of $p$ in the definition of $[f]_p,$ and let $V$ be another open connected neighborhood of $p$. Then the set $\{ [f]_z:z\in U\cap V\}$ 
defines an open neighborhood of $[f]_p$ with image in $V$ under $\pi.$ Thus $\pi$ is continuous. 
Now let $(f,U)$ be a representative of $[f]_p$. The $U\ni p$ contains a connected open neighborhood
of $p$ on which $\pi$ is bijective with inverse given by $z\mapsto [f]_z,$ and this inverse is continuous. 
This completes the proof.
\end{proof}
\begin{proposition}
The topological space $O$ is Hausdorff.
\end{proposition}
\begin{proof}
Let $[f]_p,[g]_r$ be two elements in $O$, $[f]_p\neq [g]_r$. If $p\neq r$ there are two disjoint domains $U_p \ni p,$
$U_r \ni r$ such that $\pi^{-1}(U_p)\cap \pi^{-1}(U_r)=\emptyset$ and 
$\pi^{-1}(U_p)$ is an open neighborhood of $[f]_p$ and $\pi^{-1}(U_r)$ is an open neighborhood of $[g]_r.$
If $p=r$ let $(U,f),$ $(U',g)$ be members of $[f]_p, [g]_p$ respectively. Then $U\cap U'$ contains a domain $V\ni p.$
The sets $W_p=\{[f]_z :z\in V\},$ $W'_p=\{[g]_z :z\in V\}$ respectively, are two open neighborhoods of $[f]_p$,$[g_p]$ respectively. If there exists $[h]_z\in W_p\cap W'$ then $[f]_z=[g]_z$ and since $V$ is connected 
this would imply $f\equiv g$ on $V$, i.e.\ the contradiction $[f]_p=[g]_p$. Hence we have $W_p\cap W'=\emptyset.$
This completes the proof.
\end{proof}
\begin{proposition}
	The topological space $O$ is paracompact (i.e.\ each open cover has a locally finite refinement in the sense that each member of the new cover is a subset of a member in the original cover and each point has an open neighborhood that intersects only finitely many members of the refinement).
\end{proposition}
\begin{proof}
	We shall need the following which we state without proof.
	\begin{theorem}[Poincar\'e-Volterra, see e.g.\ Narasimhan \cite{narasimhannieve}]\label{poincarevolterra}
		Let $X$ be a Hausdorff topological space with a countable basis of open sets and let $Y$ be a connected variety. If
		$f:X\to Y,$ is a continuous function such that for each $x\in X,$
		$f^{-1}(y)$ is discrete then $Y$ has a countable basis of open sets.
	\end{theorem}
Let $Y$ be a connected component of $O.$ For each $x\in \Cn$, $\pi^{-1}(x)$ equals $O_x$ which is discrete
so the restriction $\pi|_Y$ satisfies the conditions of Theorem \ref{poincarevolterra}, hence
each connected component of $O$ has a countable basis of open sets.
By a known result of Dieudonn\'e \cite{dieudonne}, this implies that $O$ is paracompact. This completes the proof.	
\end{proof}
\begin{remark}
The topology on $O^\alpha$ is induced by the topology on $O.$ To see this note that each open subset of $O^\alpha$ is an open subset of $O$. Conversely suppose $U$ is an open subset of $O$. Consider $U\cap O^\alpha$.
If $[f]_p\in U\cap O^\alpha$ such that $f$ is polyanalytic of order $\leq \alpha$ in a neighborhood $W$ of $p$ then
$\{[f]_z:z\in W\}$ is an open neighborhood of $[f]_p$ in $U\cap O^\alpha$, so $U\cap O^\alpha$ is open in $O^\alpha.$
This also shows that $O^\alpha$ is open in $O.$
 We claim that also $O^\alpha$ is closed in $O.$
 Let $[g]_x\in \overline{O\alpha},$ where $g$ is a polyanalytic on a domain $U$ containing $x.$ 
 Then there exists $y\in U$ and a function $f$ that is polyanalytic of order $\leq \alpha$ near $y$ such that $[f]_y=[g]_y$ hence $g$ is of order $\leq \alpha$ near $y$ thus also near $x.$ This implies $[g]_x\in O^\alpha$ so
  $O^\alpha$ is closed in $O.$
\end{remark}
Note that the map $q:O\to \Z_+^n$, $[f]_p\mapsto \alpha$ whenever $f$ is $\alpha$-analytic, is a locally constant map. Hence it
is constant on each connected component of $O$.
Let $\pi_i:\R^n\to \R$ be the natural projection. Then $q_i\circ q :O\to \R$,
$[f]_p \mapsto \alpha_i$ is continuous for each $i=1,\ldots,n.$ If $K\subset O$ is compactthen $q_i\circ q$
is bounded, thus $K$ belongs to some $O^\alpha$. Hence
a subset $K\subset O$ is compact iff there exists $\alpha\in \Z_+^n$ such that
$K$ belongs to some $O^\alpha$ and is compact in $O^\alpha.$
\begin{remark}
The sheaf $(O,\pi)$ can be identified as the inductive limit of the $(O^\alpha,\pi)$ with respect to the canonical injections
$j_\alpha :O^\alpha \to O.$ If $\alpha\leq \beta$ the injection $j_{\alpha,\beta}:O^\alpha \to O^\beta$ is a 
sheaf homeomorphism, so the topological space $O$ is the inductive limit of the $O^\alpha.$ A subset $U\subset O$ is open iff $U\cap O^\alpha$
is open in $O^\alpha$ for some $\alpha\in \Z_+^n.$
\end{remark}
\begin{definition}
	Let $(M,\pi)$ be a Riemann domain of dimension $n$, i.e.\ $\pi:M\to \Cn$ is a local homeomorphism and	$M$ is a connected Hausdorff topological space. A function $f:M\to \C$ is called
	{\em $\alpha$-analytic} on $(M,\pi)$ if for all $p\in M$ there exists an open neighborhood $U$ such that (i) $\pi|_U$ is a homeomorphism of $U$ onto an open subset $V\subset\Cn$ and (ii) $f\circ\pi^{-1}$ is $\alpha$-analytic. 
\end{definition}
Note that by the definition if $f$ is an $\alpha$-analytic function on a Riemann domain $(M,\pi)$ over $\Cn$
then if $g$ is a such that locally near some point $f=g\circ \pi$ then $g$ is $\alpha$-analytic. We define the derivatives
$\partial^\alpha f(x):=\partial^\alpha g(z),$ $\overline{\partial}^\alpha f(x):=\overline{\partial}^\alpha g(z),$
$z=\pi(x).$ 
\begin{remark}\label{prevriemrem}
If $(M,\pi)$ is a Riemann domain over $\Cn$ and $f$
is an $\alpha$-analytic function $\Omega$ then the zero set $f^{-1}(0)$ is obviously open. Let $p\in \overline{f^{-1}(0)}$ (by which we mean the closure).
Then $\pi(p)=q\in \overline{\pi(f^{-1}(0))}.$ But if locally $f=g\circ \pi$ then the set of points on which
$g=0$ is closed, thus $q\in \pi(f^{-1}(0)).$ Hence $g=0$ near $q$ and therefore $f=0$ near $p,$ thus $p\in f^{-1}(0).$
Since $f^{-1}(0)$ is closed and open we have that $f$ vanishes identically on a connected component whenever $f$ vanishes on an open subset
of the component. In particular, if two $\alpha$-analytic functions on a Riemann domain coincide on an open subset then they coincide everywhere
since $M$ is connected.
\end{remark}

\begin{definition}
	Let $(X,\pi)$ be a Riemann domain over $\Cn$, let $\Omega\subset \Cn$ be a domain and let 
	$f$ be an $\alpha$-analytic function $\Omega$.
	We say that $f$ has $\alpha$-analytic extension\index{$\alpha$-analytic extension on Riemann domains} in $(X,\pi)$ if there exists a homeomorphism $j:\Omega\to X$ onto an open subset of $X$
	and an $\alpha$-analytic function $g$ (necessarily unique by Remark \ref{prevriemrem}) on $(X,\pi)$ such that (i) $\pi\circ j=\mbox{Id}$ and (ii) $g\circ j=f.$
\end{definition}

\begin{theorem}
	Let $\Omega\subset \Cn$ be a domain and let 
	$f$ be an $\alpha$-analytic function $\Omega$. Then there exists a Riemann domain $(X,\pi)$ over $\Cn$ and an $\alpha$-analytic extension
	of $f$ on $(X,\pi)$ such that for any other Riemann domain $(X',\pi')$  and any $\alpha$-analytic extension $g'$ of $f$ on $(X',\pi')$
	there exists a local homeomorphism
	$\phi:X'\to X$ such that
	\begin{equation}
	\pi\circ \phi=\pi',\quad g\circ \phi=g',\quad \phi\circ j'=j
	\end{equation}
	\end{theorem}
	\begin{proof}
	Consider the sheaf $(O,\pi)$ of germs of polyanalytic functions. The set $\Omega'=\{[f]_z:z\in \Omega\}$ is open in $O.$
	If $X$ is the connected component of $\Omega'$ in $O$ then the restriction $\pi'|_{\Omega'}$ is a homeomorphism with inverse $j:z\mapsto [f]_z.$ The restriction $\pi|_X$ defines a domain in $(X,\pi)$
	and $\pi\circ j=\mbox{Id},$ where $j(\Omega)=\Omega'$ is open and connected in $X.$
	If $g:X\to \C,$ $[u]_z\mapsto u(z)$ for a germ $[u]_z$ of a polyanalytic function near $z$
	in $\Cn$ then for all $z\in \Omega$ we have $g\circ j(z)=g([f]_z)=f(z).$ Hence $g\circ j=f$ and $g$ is an extension of $f$ on $(X,\pi).$
	Also for each $[u]_z$ we have $g\circ\pi^{-1}(z)=g([u]_z)=u(z)$ where $u$ is $\alpha$-analytic near $z$ since by our previous remarks
	$\alpha$ is locally constant and $\alpha$ is the order of $g([f]_z)=f(z)$ so the extension $g$
	is $\alpha$-analytic.
	Now suppose 
	$X',\pi',g',j'$ define another $\alpha$-analytic extension.
	For $x'\in X'$ and $p'=\pi'(x')$ and an open neighborhood $U'$ of $p'$
	we know that $\pi'(U')$ is an open neigbhorhood (the homeomorphic image) of $p'$. Set $h=g'\circ (\pi')^{-1}$ in $U.$
	With $\phi(x'):=[h]_{x'}$ we have a map $X'\to O$. If $x'=j'(z),$ $z\in \omega$ then $p'=z$ and $\phi(x')=
	[g'\circ j']_z=[f]_z\in \Omega'$. So we have $\phi(j'(\Omega))\subset X$. Since $\phi$
	is continuous and $X'$ is connected we have $\phi(X')\subset X,$ and $\phi$ can be verified to satisfied the
	sought conditions. This completes the proof. 
	\end{proof}

\section{The real and imaginary parts of $\alpha$-analytic functions}\label{realandimagsec}

We begin with a well-known result.
\begin{proposition}
	Suppose $u(x,y)$ is a real valued harmonic function on an open subset $\Omega\subset\C.$ Then $u$ is locally 
	the real part of a holomorphic function which is determined up to an
	additive 
	complex constant. 
\end{proposition}
\begin{proof}
	Since $\partial_{\bar{z}}\partial_z u=0$ we know that
	$\partial_z u$ is holomorphic thus has a complex primitive, $f(z)$ (up to an additive constant),
	on $\Omega,$
	satisfying $df=2\partial_z g(z) dz$. This shows that $\overline{\partial} f=0$ so $f$ is holomorphic. Furthermore, taking the conjugate yields
	(since $g$ is real so the conjugate og
	$\partial_z g$ is $\partial_{\bar{z}} g$)
	$d\bar{f}=2\partial_{\bar{z}}gd\bar{z}$. We thus
	obtain
	\begin{equation}
	\frac{1}{2}d(f+\bar{f})=dg
	\end{equation}
	thus $g$ is the real part of $f$ up to an additive constant.
	Next suppose $f_1,f_2$ are two holomorphic functions whose real part coincide.
	Let $h=f_1-f_2$. Then $d(h+\bar{h})=0$ i.e.\
	\begin{equation}
	\partial_z hdz+\partial_{\bar{z}}\bar{h} d\bar{z} =0
		\end{equation}
		which implies that $\partial_z h=0$ and $\partial_{\bar{z}}\bar{h}$ thus $h$ is a constant.
			This completes the proof.
			\end{proof}
			Note that the result is local. For example, if $\Omega=\C\setminus \{0\}$ then the function $\log \abs{z}$
			is not the real part of a holomorphic function on $\Omega$ since the logarithm of $z$ has no single-valued branch
			on $\Omega.$ However, because for a harmonic function $g$ on a {\em simply connected} domain, 
			the differential form $2\partial_z g dz$ does have a primitive in $\Omega$ we have the following.
			\begin{proposition}
			If $\Omega\subset\C$ is a domain then any real valued harmonic $g$ is the real part of a holomorphic function
			on $\Omega.$
			\end{proposition}
			Suppose $f(z)$ is holomorphic with power expansion on a disc $\{\abs{z}<r\}$ for some $r>0$
			given by
			$f(z)=\sum_{j=0}^\infty c_j z^j,$ in particular $\bar{f}(\bar{z})=\sum_{j=0}^\infty \bar{c}_j \bar{z}^j$ 
	converges also.
	Consider
	\begin{equation}
	g(x,y):=f(x+iy)+\bar{f}(x-iy)
	\end{equation}
	We have
	\begin{equation}
	2g\left(\frac{z}{2},\frac{z}{2i}\right)=f(z)+\bar{f}(0)
	\end{equation}
	In particular,
	\begin{equation}
	2g(0,0)=f(0)+\bar{f}(0)
	\end{equation}
	thus
	\begin{equation}
	2g\left(\frac{z}{2},\frac{z}{2i}\right)-g(0,0)=f(z)+\frac{1}{2}(\bar{f}(0)-f(0))
	\end{equation}
	which shows that $f(z)$ is up to an additive pure imaginary constant
	equal to
	$2g\left(\frac{z}{2},\frac{z}{2i}\right)-g(0,0).$
	Now let
	$g(x,y)$ be a real valued harmonic function and suppose $g$ is the real part of a holomorphic function 
	$f(z)=\sum_{j=0}^\infty c_jz^j$ on a disc $\{\abs{z}<\rho\}$, some $\rho>0$, and we suppose that $\im c_0=0$
	(which as we have shown above completely determines $f(z)$, given $g$).
	For $r<\rho$ 
	\begin{equation}
	g(r\cos\theta,r\sin\theta)=c_0+\frac{1}{2}\sum_{j=1}^\infty r^j (a_j\exp(ij\theta)+\bar{c}_j\exp(-ij\theta))
	\end{equation}
	with normal convergence with respect to $\theta\in [0,2\pi)$ and the Fourier coefficients
	are given by
	\begin{equation}
	c_0=\frac{1}{2\pi} \int_0^{2\pi} g(r\cos\theta,r\sin\theta)d\theta
	\end{equation}
	\begin{equation}
	c_j=\frac{1}{2\pi} \int_0^{2\pi} \frac{g(r\cos\theta,r\sin\theta)}{(r\exp(i\theta))^j} d\theta ,\quad j\in \Z_+
	\end{equation}
	This yields for $\abs{z}<r$
	\begin{equation}
	f(z)=\frac{1}{2\pi} \int_0^{2\pi}g(r\cos\theta,r\sin\theta)
	\left(1+2\sum_{j=1}^\infty \frac{z^j}{(r\exp(i\theta))^j}\right) d\theta
	\end{equation}
	Due to normal convergence 
	we may change the order of summation and integration. Also
	\begin{equation}
	\left(1+2\sum_{j=1}^\infty \frac{z^j}{(r\exp(i\theta))^j}\right)=\frac{r\exp(i\theta)+z}{r\exp(i\theta)-z}
	\end{equation}
	thus
	\begin{equation}\label{forstabvpekvationenpoisson}
	f(z)=\frac{1}{2\pi} \int_0^{2\pi}g(r\cos\theta,r\sin\theta)
	\frac{r\exp(i\theta)+z}{r\exp(i\theta)-z} d\theta ,\quad r<\rho
	\end{equation}
	This expresses $f(z)$ in terms of its real part on the boundary of the disc $\abs{z}<r$.
	Equating the real parts in Eqn.(\ref{forstabvpekvationenpoisson}) gives
	\begin{equation}\label{forstabvpekvationenpoisson1}
	g(x,y)=\frac{1}{2\pi} \int_0^{2\pi} g(r\cos\theta,r\sin\theta)
	\frac{r^2 -\abs{z}^2}{\abs{r\exp(i\theta)-z}^2} d\theta ,\quad r<\rho
	\end{equation}
	We have thus derived the {\em Poisson kernel}\index{Poisson kernel}, $\frac{r^2 -\abs{z}^2}{\abs{r\exp(i\theta)-z}^2}$ 
	for the disc,
	as the real part of the holomorphic function
	$\frac{r\exp(i\theta)+z}{r\exp(i\theta)-z}$. In particular, the Poisson kernel 
	is harmonic and is zero at all points of $\abs{z}=r$ except $z=r\exp(i\theta).$
	Furthermore, it is positive for $\abs{z}<r$ and for fixed $r$ and $z$ with $\abs{z}<r$ it is a periodic function of $\theta$
	defining a positive mass distribution of the unit circle
	and using the harmonic function $g\equiv 1$ in Eqn.(\ref{forstabvpekvationenpoisson1}) gives
	\begin{equation}
	\frac{1}{2\pi} \int_0^{2\pi}
	\frac{r^2 -\abs{z}^2}{\abs{r\exp(i\theta)-z}^2} d\theta =1
	\end{equation} 

Recall that in complex analysis of one variable there exists a unique harmonic conjugate to
a harmonic function $u\in C^0(\overline{\Omega})$ where $\Omega\subset \C$ is a simply connected bounded domain. 
Let 
\begin{equation}\label{greeneq}
g(\zeta,z_0)=-\log \abs{w(\zeta,z_0)}
\end{equation} 
where $w$ is the Riemann mapping (biholomorphism) from a domain $\Omega\subset\C$
to $\Omega$,
$z_0\in \Omega$ satisfying $w(\zeta,z_0)=0,$ $w'(z_0,z_0)\neq 0$. 
The so called Green function is usually defined by the continuous extension to $\overline{\Omega}$ (Caratheodory theorem), of $g$, obtained via the
continuous extension of the Riemann mapping.
This extension will satisfy $\abs{w(\zeta,z)}=1$ for $\zeta\in \gamma.$
Specifically, a Green function of $\Omega$ is required to be a continuous function $g(\zeta,z)$, defined for 
$z\neq \zeta$,
$z\in \overline{\Omega}$ 
$\zeta\in\overline{\Omega}$, 
satisfying that for each $z\in \Omega$
\begin{equation}\label{greeneqdisken}
g(\zeta,z)=\ln \abs{\zeta -z} +h_z(\zeta)
\end{equation}
where $h_z(\zeta)$ is harmonic in $\Omega$ and continuous on $\overline{\Omega}$ such that
if $\zeta\in \gamma=\partial\Omega$, $z\in \gamma$, $\zeta\in \overline{\Omega}$ 
then
\begin{equation}
g(\zeta,z)=0
\end{equation}
i.e.\ for fixed $z_0\in \Omega$, we have that, as a function of $\zeta\in \gamma$, $h(\zeta,z_0)=-\ln \abs{\zeta -z}$ (clearly, the right hand side will be
continuous if the boundary is continuous, thus the solution to the harmonic version of the Dirichlet problem will render an appropriate $h$).
Existence of the Green function for simply connected domains with continuous boundary is in a sense equivalent the existence of the Riemann mapping as 
described above. 
For the unit disc 
\begin{equation}
\frac{\partial}{\partial n}=\zeta\frac{\partial}{\partial\zeta}+
\bar{\zeta}\frac{\partial}{\partial\bar{\zeta}}=2\re(\zeta\frac{\partial}{\partial\zeta})
\end{equation} and
\begin{equation}
g(\zeta,z)=-\ln\abs{\frac{\zeta-z}{1-\bar{z}\zeta}}
\end{equation}
More generally, for a disc $\{\abs{\zeta}<R\}$ we have
\begin{equation}\label{goldbergnevekv4}
w(\zeta,z)=\frac{R(\zeta-z)}{R^2-\bar{z}\zeta}
\end{equation} and
\begin{equation}
g(\zeta,z)=\ln\abs{\frac{R^2-\zeta\bar{z}}{R(\zeta -z)}}
\end{equation}
Note that
\begin{equation}
\frac{\partial_\zeta w(\zeta,z)}{w(\zeta,z)}d\zeta =i\frac{\partial g}{\partial n}ds
\end{equation}
and (for the case $\{\abs{\zeta}<R\}$ and $\zeta=R\exp(i\theta)$
\begin{equation}\label{goldbergnevekv8}
\frac{\partial g(\zeta,z)}{\partial n}ds =-i\partial_\zeta \left(\ln \frac{R(\zeta -z)}{R^2 -\zeta \bar{z}}\right)d\zeta =
\frac{R^2-\abs{z}^2}{\abs{\zeta -z}^2}d\theta
\end{equation}
where for $z=r\exp(i\varphi)$
\begin{equation}\label{goldbergnevekv9}
\frac{R^2-\abs{z}^2}{\abs{\zeta -z}^2} =\frac{R^2-r^2}{R^2+r^2-2Rr\cos(\varphi-\theta)}=\re\frac{\zeta +z}{\zeta -z}
\end{equation}
If $\Omega$ is simply connected let $\gamma_1,\ldots,\gamma_n$ are analytic arcs forming $\gamma=\partial\Omega$ and let $a_j,$ 
$j=1,\ldots,n$ be the common end point of $\gamma_j$ and $\gamma_{j+1}$ ($\gamma_{n+1}=\gamma_1$ for $j=n$)
Let $\pi\alpha_j$, $0\leq \alpha_l\leq 2$ be the angles between $\gamma_j$ and $\gamma_{j+1}.$
By the Riemann-Schwarz symmetry principle
the Riemann mapping $w(\zeta,z)$ has for fixed $z\in \Omega$
analytic extension to a domain containing $\overline{\Omega}\setminus\{a_1,\ldots,a_n\}.$
We shall later need to consider $\Omega$ satisfying the following property (which is obvious when the boundary consists of a single analytic curve): 
for a fixed $z$ we have for sufficiently small $\epsilon>0$ that the Riemann mapping $w(\zeta,z)$
has on $\{\abs{\zeta-a_j}<\epsilon\}\cap \Omega$ the representation
\begin{equation}\label{goldberganvsen}
w(z,\zeta)=(\zeta-a_j)^{\frac{1}{\alpha_j}}\varphi_j(\zeta)+w_j,\quad j=1,\ldots,n
\end{equation}
for analytic $\varphi_j(\zeta)$ on $\{\abs{\zeta-a_j}<\epsilon\}$ with $\varphi_j(a_j)\neq 0$ and $\abs{w_j}=1,$
$j=1,\ldots,n.$
\\
The Green function is obviously unique since if $g_1,g_2$ are two Green functions then
$g_1-g_2$ is, for fixed $z\in \Omega$, a difference of two harmonic functions in $\zeta$ on $\Omega$
and continuous on $\overline{\Omega},$ so the second condition on Green functions implies that
$g_1-g_2$ vanishes on the boundary thus by the maximum principle
it vanishes identically.
Note that the first condition implies that for $z\in \Omega$ there exists an $\epsilon>0$ such that
$g(\zeta,z)>0$ for $\zeta$ in $\{\abs{\zeta-z}<\epsilon\}$. By the minimum principle for $\Omega\setminus  \{\abs{\zeta-z}<\epsilon\}$
we have $g(\zeta,z)>0$ for all $\zeta\in \Omega,z\in \Omega.$
\\
Recall that the solution for the Dirichlet problem for harmonic functions on $\Omega$
is given by
\begin{equation}
u(z)=-\frac{1}{2\pi}\int_{\partial \Omega} a(\zeta)\frac{\partial g(\zeta,z)}{\partial n}d\zeta, \quad z\in  \Omega
\end{equation}
for a given boundary function $a\in C^0(\partial \Omega,\R),$
where $n$ denotes the outward pointing unit normal.
Locally near a point $p_0\in \Omega,$ a harmonic conjugate, $h(z,p_0)$, to $g$
is given by 
\begin{equation}
h(z,p_0)=\int_{p_0}^z \frac{\partial g(\zeta,p_0)}{\partial n} +\mbox{const.}
\end{equation}
Then a harmonic conjugate to $u$ in $D$ which is continuous on $\partial D$ is given by
\begin{equation}
v(z)=-\frac{1}{2\pi}\int_{\partial D} u(\zeta)\frac{\partial h(\zeta,z)}{\partial n}d\zeta, \quad z\in D
\end{equation} 
Obviously, the pair $u,v$ will still satisfy the Cauchy-Riemann equations if a constant is added to either. Thus for fixed real part
there is a family of harmonic conjugates where the members differ by imaginary constants.
Consider the equation for $q$-harmonic functions $u$ for $q\geq 1,$
\begin{equation}\label{polyharmeq}
(\partial_z\partial_{\bar{z}})^q u=0
\end{equation}
Let $f$ be a $q$-analytic function. Then obviously 
$\overline{\partial}^q (\re f)=-i\overline{\partial}^q (\im f)$ so that
\begin{equation}
\Delta^q \re f=-i\Delta^q \im f
\end{equation}
whence both $\re f$ and $\im f$ are $q$-harmonic.
Conversely (following e.g.\ Gakhov \cite{gakhov}, p.250)
Eqn.(\ref{polyharmeq}) is equivalent to
\begin{equation}
(\partial_z\partial_{\bar{z}})^{q-1} u=\varphi(z)+\psi(\bar{z}),
\end{equation}
where $\varphi,\psi$ are holomorphic functions of their arguments. Since the left hand side is real-valued
we have $\varphi(z)+\psi(\bar{z})=\overline{\varphi(z)+\psi(\bar{z})}$ thus $\varphi(z)+\overline{\psi(\bar{z})}=
2\re \varphi(z),$ which determines $\varphi(z):=\sum_{j=0}^\infty z^j,$ up to an imaginary term. 
Integrating both sides with respect to $z$ and then with respect to $\bar{z}$
gives
\begin{equation}
(\partial_z\partial_{\bar{z}})^{q-2} u=z\bar{z}(\varphi_1(z)+\overline{\varphi_1(\bar{z})})+\psi_1(z)+\overline{\psi_1(z)}=2\re(z\bar{z}\varphi_1(z)+\psi_1(z))
\end{equation}
where $\varphi_1(z)=\sum_{j=0}^\infty \frac{a_j}{j+1}z^j,$ and is determined by $u$ up to an imaginary term.
Iterating this $q-2$ times renders the
representation
\begin{equation}
u(z)=\re \left( \sum_{j=0}^{q-1}(z\bar{z})^j\varphi_j(z)\right)
\end{equation}
for holomorphic $\varphi_j$ determined by $u$ up to arbitrary imaginary terms.
Clearly, $u(z)$ is then identified as the real part of a {\em reduced} $q$-analytic function
and we call
\begin{equation}
v(z):=\im \left( \sum_{j=0}^{q-1}(z\bar{z})^j\varphi_j(z)\right)
\end{equation}
a {\em polyharmonic conjugate to $u$}. Note however
that the function $f(z)=u(z)+iv(z)$ is not necessarily reduced.
Suppose $f$ is an entire $q$-analytic function with 
the representation $\sum_{j=0}^q a_j(z)\bar{z}^j$ for entire $a_j.$
This implies that $\re f$ and $\im f$ can be represented as real-analytic functions of 
$x=(z+\bar{z})/2,$ $y=(z-\bar{z})/2i$, which have holomorphic extensions to $\C^2$,
in the sense that there exists a holomorphic functions $U(z,w),V(z,w)$ in two variables
such that $U(z,w)|_{w=\bar{z}}=\re f(z,\bar{z})$, $V(z,w)|_{w=\bar{z}}=\im f(z,\bar{z})$.
Let $O_p^\alpha$ be the set of germs of polyanalytic functions of order $\leq \alpha$ (see Definition \ref{germsdefforsta})
and set $O^\alpha=\bigcup_{p} O_p^\alpha$
and assume the topologies defined in Section \ref{extensionsec}, in particular $[f]_p\mapsto p$ is a map $\pi:O^\alpha \to \Cn$ that will be a local homeomorphic if we use as fundamental system of open neighborhoods of $[f]_p$ the sets
$V(U,f):=\{[f]_p,p\in U\}$ for each $(U,f)$, $U$ an open neighborhood of $p.$
\begin{proposition}
	Let $\alpha,\beta\in \Z_+^n,$
	$\beta_j\leq \alpha_j,$ $j=1,\ldots,n.$ Let $\Omega\subset\Cn$ be a simply connected domain and let 
	$f_1,\ldots,f_n$ be polyanalytic of order $\leq\alpha$ (by which we mean that the separate order with respect to $z_j$ is $\leq \alpha_j,$ $j=1,\ldots,n$) on $\Omega$
	such that
	\begin{equation} \partial_{z_j}^{\beta_j}f_k=\partial_{z_k}^{\beta_k}f_j,\quad j,k=1,\ldots,n\end{equation}
	Then the system
	\begin{equation} \partial_{z_j}^{\beta_j}f=f_j,\quad j=1,\ldots,n\end{equation}
	has a solution that is polyanalytic of order $\leq \alpha$ on $\Omega$, unique up to
	addition by a polyanalytic polynomial of order $\leq \alpha$, that is of degree $<\beta_j$ with respect to $z_j,$
	$j=1,\ldots,n.$ 
\end{proposition}
\begin{proof}
	Denote by $O^{\alpha,n}$ the cartesian product $O^\alpha\times \cdots\times O^\alpha$ of $n$ factors, whose members are of the form $(\phi_1,\ldots,\phi_n),$
	$\pi(\phi_1)=\cdots=\pi(\phi_n).$ Then the extension of $\pi$ defines a local homeomorphism $\pi: \bigcup_\alpha O^{\alpha,n} \to \Cn,$ where $\overbrace{O\times\cdots\times O}^{n-\mbox{times}}=\bigcup_\alpha O^{\alpha,n}.$
	Set
	\begin{equation}
	W_0:=\{ (\phi_1,\ldots,\phi_n)\in O^{\alpha,n},\partial_{z_i}^{\beta_i} \phi_j=\partial_{z_j}^{\beta_j} \phi_i,i,j=1,\ldots,n\}
	\end{equation}	
	where the partial derivatives of a germ are equal to the germs of the partial derivatives. Then $W_0$ is open and closed in $O^{\alpha,n}$.
	Then $\pi_0=\pi|_{W_0}$ is again a local homeomorphis. Set $\psi:\Omega\to W_0$, $\mapsto ([f_{1}]_z,\ldots, [f_{n}]_z),$ where $[f_i]_z$ is the germ of $f_i$ at $z$ in $\Omega.$
	If $W=\psi(\Omega)$ then $\psi$ is a homeomorphism of $\Omega$ onto $E$ and the inverse is $\pi_0|_W,$ hence $W$ is also simply connected.
	Set
	\begin{equation}
	\eta:\overbrace{O\times\cdots\times O}^{n-\mbox{times}}\to W_0, [f]_z\mapsto \left(\partial_{z_1}^{\beta_1}[f]_z,\ldots,\partial_{z_n}^{\beta_n}[f]_z \right)
	\end{equation} 
	If $([\phi_1]_p,\ldots,[\phi_n]_p))\in W_0$ and $\Delta(p,r)$ is a polydisc in $\Cn$ on which each $\phi_i$ is polyanalytic then by Proposition \ref{systemab} there is a polyanalytic function $f$ on the polydisc such that
	$\partial_{z_i}^{\beta_i} f=\phi_i,$ $i=1,\ldots,n.$
	Set $V:=\{[\phi_{1}]_z,\ldots, [\phi_{n}]_z:z\in \Delta(p,r)\}.$ Then each connected component of $\eta^{-1}8V)$ is of the form $U=\{[f+Q]_z:z\in \Delta(p,r)\},$ for a polyanalytic polynomial $Q$ (see Proposition \ref{systemab}) and $\eta|_U$ is a homeomorphism. If $C$ is a connected component of $\eta^{-1}(W)$ then, since $W$ is simply connected, $\eta|_C$ i a covering map thus a homeomorphism. Setting $F(z):=\eta^{-1}(\psi(z))(z)$ to be the value of the germ $[f+Q]_z$, i.e.\ $F(z)=f(z)+Q(z),$ we have
	$\partial_{z_i}^{\beta_i} F=f_i,$ $i=1,\ldots,n.$ This completes the proof.
\end{proof}

\begin{definition}
	Let $\Omega\subset\Cn$ be a simply connected domain and let $\alpha\in \Z_+^n.$. A smooth function $f:\Omega\to \C$ is called {\em pluriharmonic of order $\alpha$}\index{Pluriharmonic function}
	if $\alpha_1,\ldots,\alpha_n$ 
	satisfy 
	\begin{equation} \partial_{z_j}^{\alpha_j}\partial_{\bar{z}_k}^{\alpha_k}f=0,\quad j,k=1,\ldots,n\end{equation}  
\end{definition}
\begin{proposition}
	A smooth real function $u:\Omega\to \R$ on a simply connected domain 
	$\Omega\subset\Cn$ is pluriharmonic of order $\alpha$ iff $u$ is the (not necessarily unique) real part of an $\alpha$-analytic function.	
\end{proposition}
\begin{proof}
	If $f$ is $\alpha$-analytic then clearly it is pluriharmonic of order $\alpha$ and so is its real part
	$2\re f=f+\overline{f}.$	
	Conversely, suppose $u:\Omega\to \C$ satisfies
	\begin{equation} \partial_{z_j}^{\alpha_j}\partial_{\bar{z}_k}^{\alpha_k}u=0,\quad j,k=1,\ldots,n\end{equation}  
	For fixed $k$ we consider the equations
	\begin{equation} \partial_{z_j}^{\alpha_j}f=f_k,\quad j=1,\ldots,n\end{equation} 
	where $f$ denotes a general solution. 
	We have that 
	\begin{equation} \partial_{z_k}^{\alpha_k}f_k=0, \quad  k=1,\ldots,n\end{equation}  
	so by Proposition \ref{systemab} there is an $\alpha$-analytic function $g$ on $\Omega$ such that
	$\partial_{z_k}^{\alpha_k}g=f_k.$ This implies that 
	\begin{equation} \partial_{z_j}^{\alpha_j}(f-g)=0, \quad  j=1,\ldots,n\end{equation}
	thus $f=g+\overline{h}$ for some 
		for an $\alpha$-analytic function $h.$
	Since $u$ is a real-valued we have
	\begin{equation} u=g+\overline{h}=\overline{g}+h=\frac{1}{2}\re(g+\overline{h})\end{equation}
	and since $u$ is pluriharmonic of order $\alpha$ we have that
	$g+\overline{h}$ is $\alpha$-analytic. This completes the proof.
\end{proof}
\begin{remark}
	If $f,g$ are two $\alpha$-analytic functions on a simply connected domain then so is $f-g$ and $\bar{f}+\bar{g}$ so the condition with $u=\re f=\re g$ implies that
	$f-g$ is a pure imaginary polynomial, thus any $\alpha$-analytic function $h$ having the same real part as $f$ has the form
	$f+iP$ for some real polynomial $P$ of (not necessarily strict) order $\alpha$.
\end{remark}
\section{Local properties of classical kind}
We use here the definitions and notations on germs of polyanalytic functions given in
Section \ref{extensionsec}. 
In particular, given a function, $f$, polyanalytic on a domain $\Omega\subset\Cn$, $0\in \Omega$,
of a given arbitrary order, near the origin $z=0$
we shall here denote by $[f]_0$ the germ of $f$ at $0$ in the set, 
$\mbox{PA}(\Omega)$ of polyanalytic functions of any arbitrary order, near $0$. 
We denote by $\mbox{PA}(0)$ the ring of germs of polyanalytic functions at $0$.
Note that an element $[f]_0\in \mbox{PA}(0)$ is invertible iff $f$ is analytic at $0$ and $f(0)\neq 0.$
Classically one can study the ring structure of power series algebras using the theory of polynomial algebras, 
(even if the former are not in general principal ideals), 
if the ideals are finitely generated. 
We shall need to recall some basic definitions and we point out that pseudopolynomials will be recurring as an important proof-tool in upcoming chapters.
\begin{definition}[Fritzche \& Grauert \cite{frigrau}, p.124]\label{pseudopolynomdef}
	If $\Omega\subset\Cn$ is a domain then a {\em pseudopolynomial},\index{Pseudopolynomials} $F$, of order (or degree) $q$, is defined as
	a function $F:\Omega\times \C \to \C$, $(z,w)\mapsto \sum_{j=0}^{q} a_j(z)w^j,$
	for holomorphic functions $a_j(z).$
\end{definition}

Pseudopolynomials are certainly not new to complex variables as they appear already
in the works of Weierstrass (see Section \ref{plemeljsec} of the appendix).
As usual, an analytic function $f(\zeta)$ is called divisible by another analytic function $g(\zeta)$
at a point $\zeta=p,$ $\zeta\in \Cn$, if there exists an analytic function $h(\zeta)$ (a divisor of $f$ at $p$) near $p$
such that near $p$ we have $f=gh.$ 
Division by the zero function is excluded.
An analytic $f(\zeta)$ such that $f(0)=0$ is called irreducible at $p\in \Cn$
if $f=f_1 f_2$ near $p$, for analytic $f_1,f_2$ such that $f_j(0)=0$
and $f_j\not\equiv 0,$ $j=1,2.$ An irreducible factor is called a prime factor.
If $f$ is analytic such that $f(p)=0$ but $f\not\equiv 0$ then
$f$ has (up to equivalent factors) a unique decomposition into prime factors at $p$.
By equivalent factors we mean the following: two analytic functions
$g,h$ with $g(p)=h(p)=0$, are called equivalent at $p\in \Cn$ if each is divisible by the other near $p$.
For this it is necessary and sufficient that there exist an analytic $Q$, $Q(p)\neq 0$
such that near $p$, $g=Q h.$ If two analytic functions vanishing at $p$ but not identically have no
common divisor at $p$ they are called relatively prime. These definitions are invariant under nonsingular analytic
transformations.
A pseudopolynomial $f(z,w)$ is called {\em divisible}, at a point $p$ of its domain, by a
pseudopolynomial $g(z,w)$ if there exists a pseudopolynomial $h(z,w)$ near $p$ such that
$f(z,w)=g(z,w)h(z,w)$ near $p.$ Here division by the zero function is excluded.

\begin{proposition}
	A pseudopolynomial $f(z,w)=\sum_{j=0}^m a_j(z_1,\ldots,z_n)w^j$ (where the $a_j(z)$ are analytic functions
	and $(z,w)\in \C^n\times\C$) is divisible at the origin
	by a function $\phi(z)$ then each coefficient function $a_j(z)$ is divisble by
	$\phi(z).$
\end{proposition}
\begin{proof}
	We have by hypothesis an analytic function $g(z,w)$ near the origin such that $f(z,w)=g(z,w)\phi(z).$
	We expand $g(z,w)$ with respect to $w$ near the origin according to
	$g(z,w)=\sum_{j=0}^\infty c_j w^j.$ Substituting this into the expression for 
	$f(z,w)$ and comparing coefficients proves the proposition.
\end{proof}
A pseudopolynomial is called {\em primitive} if it is not divisible by any function $g(z)$
at the origin with $g(0)=0.$
If $f(z,w)\not\equiv 0$ is a pseudopolynomial in $w$ (with analytic coefficients in $z$) then it can be written as the product of two
psedudopolynomials one of which is of degree zero and the second is primitive, namely if $f(z,w)$
is itself not primitive
then we can factor out the greatest common divisor of the
analytic coefficients in $z.$ Also the same proof as in the case of holomorphic functions
shows that if the
product of two pseudopolynomials in $w$ is divisible (at the origin) by a function $\phi(z)$, irreducible at the origin, then one of the
two factors is divisible by $\phi(z)$ at the origin.
Important examples are of course pseudomonomials at a point $p\in \C^m$,
namely
$f(z,w)=a(z)(w-p)^N$ for an integer $N\in \N$ and analytic $a(z).$ These are obviously primitive.

\begin{proposition}
If $\Omega\subset\C$ is a domain then the ring of holomorphic funcitons
$\mathscr{O}(\Omega)$ is an integral domain.	
\end{proposition}	
\begin{proof}
	Assume $f_1\not\equiv 0,f_2\not\equiv 0$ are two
	holomorphic functions on $\Omega.$ Since their zero sets have only isolated points
	and since there exists $p_0\in \Omega$ such that $f_1(p_0)f_2(p_0)\neq 0$
	we have $f_1\cdot f_2\not\equiv 0.$ 
\end{proof}	
\begin{remark}
	A direct consequence is that the set $\mathscr{O}(\Omega)[w]$ of pseudopolynomials
	also has no zero divisor. Let $Q$ be the quotient field of $\mathscr{O}(\Omega)$
	and let $Q^0[w]$ denote the group of units in the 
	integral domain that consists of the nonzero polynomials of degree $0$. Then 
	$Q^0[w]\cap \mathscr{O}(\Omega)$ is the multiplicative subgroup of non-identically vanishing holomorphic functions
	on $\Omega.$
\end{remark}

\begin{definition}
A ring $R$ is called noetherian\index{Noetherian ring} if every ideal in $R$ is finitely generated.
An $R$-module is called noetherian if every $R$-submodule is finitely generated. 
A {\em monoid}\index{Monoid}, $M$, is a commutative semi-group with $1$ such that $ab=ac\Leftrightarrow b=c$ for all
$a,b,c\in M.$ If $a,b$ are non-units in $M$ we call $a$ a {\em divisor}\index{Divisor in a monoid}
of $b$ (denoted $a|b$) if $b=ac$ for some $c\in M.$ The element $a$ is called a {\em proper divisor} of $b$ if
$b=ac$ for some non-unit $c\in M.$ An element $a\in M$ is called {\em irreducible}\index{Irreducible element in a monoid}
if $a$ has no proper divisors. An element $a\in M$ is called {\em prime}\index{Prime element in a monoid}
$a|bc\Rightarrow a|b$ or $a|c.$
\end{definition}
Obviously prime elements are irreducible since if $a=bc$ and $a|b$ then $b=ad$ for some $d$, hence $a=adc$ and $dc=1$
so $c$ is a unit. In a principal ideal domain the unique factorization theorem holds true.
\begin{definition}
A monoid $M$  is called factorial if each non-unit of $M$ is the product of finitely many prime elements. A ring is called factorial if $R\setminus \{0\}$ is a factorial monoid.
\end{definition}
In particular, factorial rings have no non-trivial zero divisors.
\begin{proposition}
A monoid $M$ is factorial iff it satisfies the property (C) every non-unit of $M$ is a product of finitely many irreducible elements, and the
factorization is unique up to order and up to units.
\end{proposition}
\begin{proof}
	If $M$ is factorial then in particular it has the property 
	(A) each nonzero non-unit is a finite product of irreducible elements. Also if 
	$u$ is irreducible the factorial property implies that $u=v_1,\cdots v_n$ for primes $v_j$, $j=1,\ldots,n$, thus we must
	have $n=1$, so we have the property (B) any irreducible element is prime. 
	Conversely it is clear that (A)$\wedge$(B) together imply the factorial property.
	So ($M$ is factorial)$\Leftrightarrow$(A)$\wedge$(B).
	Now suppose (C) holds true. Clearly, (C) implies (A).
	If $u\in M$ is irreducible there
	exists $a,b,c\in M$ such that $cu=ab$. Now the elements $a,b,c$ have by (A)
	finite decomposition into irreducible factors, which inserted into $ab=cu$
	shows that $u$ divides both $a$ and $b.$ We conclude from this property (B).
	Hence (C)$\Rightarrow$(A)$\wedge$(B).
	Now suppose (A)$\wedge$(B). 
	Let $a$ be a nonzero non-unit in $M.$ Each decomposition of $a$ into irreducible factors is by (B) a decomposition into prime factors.
	It remains to verify that this is unique up to order and up to units.
	For a prime $u\in M$ denote by $B_u(a)$ the number $m$ such that $a\in M u^m$, $a\notin Mu^{m+1}.$
	If $a_1,\ldots,a_n$ are nonzero in $R$ then we claim
	\begin{equation}
	B_u(a_1\cdots a_n)=\sum_{j=1}^n B_u(a_j)
	\end{equation}
	It suffices to verify this for the case $n=2.$ Let $r=B_u(a_1)$ and $s=B_u(a_2).$ Then there exist $f,g\in M$ such that $f,g\notin Mu,$
	so $a_1=f u^r,$ $a_2=g u^s.$ It follows that $a_1a_2=fg u^{r+s}$ thus $fg$ does not belong to the prime ideal $Mu.$
	This implies that $a_1a_2$ do not belong to $Mu^{r+s+1}$. 
	Indeed, if
	there exist a $z\in M$ such that $a_1a_2=zu^{r+s+1}$ then 
	$(zu)u^{r+s}=(fg) u^{r+s},$ contradicting the fact that $fg\notin uM.$
	This in turn implies $B_u(a_1a_2)=B_u(a_1)+B_u(a_2).$
	By induction we may verify the claim.
	Let $v$ be an irreducible element.
	If $v\notin Mu$ then $B_u(v)=0.$ If $v\in Mu$ then $v=fu$ where $f$ is necessarily a unit so that $B_u(f)=1.$
	Now suppose $a=v_1\cdots v_n$ where the $v_1,\ldots,v_n$ are irreducible factors.
	Then $B_u(a)=\sum_{j=1}^n B_u(v_j)$, which shows that the decomposition of $a$ is unique up to reordering and up to units.
	This shows that (A)$\wedge$(B)$\Rightarrow$(C).
	This completes the proof.
	\end{proof}
The proof follows that of Reiffen, Scheja \& Vetter \cite{reiffen}, Thm. 160.
\begin{definition}
	Let $R$ be a factorial ring and $K$ its field of quotients. Let $a\in K,$
	$a\neq 0.$ If $p\in R$ is prime then $a=p^r b$ for $b\in K,$ 
	$r$ an integer
	and $p$ does not divide the numerator or denominator of $b.$ Using the unique factorization in $R$ we have that $r$ is uniquely determined by $a$.
	The integer $r$ is called the order of $a$ at $p$ and denoted $r=\mbox{ord}_p a.$
	If $a=0$ we define the order to be $\infty.$
	If $f(x)\in K[x]$ is a polynomial $f(x)=a_0+\cdots a_nx^n$ we define the order as follows: if $f=0$
	the order is $\infty,$ if $f\not\equiv 0$ then the order of $f$ at $p$, denoted $\mbox{ord}_p f$, is defined as $\min_i \mbox{ord}_p a_i$
	where the minimum is taken over all $i$ such that $a_i\not\equiv 0.$
	If $r=\mbox{ord}_p f$ we call $up^r$ a $p$-content for $f$, if $u$ is any unit of $R.$
	The {\em content} of $f$, denoted $\mbox{cont}(f)$ is defined as $\Pi p^{\mbox{ord}_p} f$, where the product is taken over all $p$ such that $\mbox{ord}_p f\neq 0.$, or any multiple of this product by a unit of $R.$
		Then the content is well-defined up to multiplication by a unit of $R.$ 
	\end{definition}
	\begin{remark}\label{remglem}
		Note that if $b\in K,$
		$b\neq 0$ then $\mbox{cont}(bf)=b\mbox{cont}(f),$ hence $f(x)=c\cdot f_1(x)$
		where $c=\mbox{cont}(f)$ and $f_1(x)$ has content $1$. In particular, all coefficients of $f_1$ lie in $R$
		and their g.c.d\ is $1.$
		\end{remark}
\begin{proposition}\label{notherfactorprop}
The ring $\mbox{PA}(0)$ is noetherian and factorial.
\end{proposition}
\begin{proof}
Denote by $\mathbf{j}$ the map which sends each polyanalytic function near the origin, $f(z)=
\sum_{k} a_k(z)\bar{z}^k$ (for holomorphic $a_k$ near $z=0$) to $\mathbf{j}(f):=\sum_k [a_k]_0 \omega^k,$
in the ring $\mathscr{O}_0$ of germs of analytic functions at $0.$ Then $\mathbf{j}$ is 
ring isomorphism. 
\begin{lemma}[Hilbert's basis theorem]\label{hilbertbasis}
	A polynomial ring $R[x_1,\ldots,x_N]$ 
	over a noetherian ring $R$, is noetherian.
\end{lemma}
\begin{proof}
We only prove the result for one indeterminant $x_1=x$, the general case follows by induction.
Let $I$ be an ideal in $R[x]$. 
Let $\{f_j\}_{j\in N}$ (for an index set $N\subseteq\Z_+$) be a sequence in $I$ such that $f_1$ is nonzero and of least possible degree
in $R$ and for $j\geq 1$, if $(f_1,\ldots,f_j)\neq I$ we choose $f_{j+1}$
to be an element in $I$ of least possible degree such that $f_{j+1}\notin (f_1,\ldots,f_j)$.
If $N$ is finite we are done. So assume $N=\Z_+.$
We can relabel and reorder the sequence such that $f_j=\sum_{k=0}^j a_k x^k$ for coefficients $a_k\in R.$
Since $R$ is noetherian the ideal $J=(a_1,a_2,\ldots)$ is finitely generated by say $\{a_1,\ldots,a_m\}$ for a fixed $m\in \Z_+.$
Suppose (in order to reach a contradiction) that $I\neq (f_1,\ldots,f_m).$
Then $f_{m+1}a_{m+1}\in J,$ so that
$a_{m+1}=\sum_{k=1}^m b_k a_k$ for $b_k\in R,$ $k=1,\ldots,m.$ Set
$g=\sum_{k=1}^m b_k f_j x^{\mbox{deg}f_{m+1}-\mbox{deg}f_k}$.
	Then $g\in (f_1,\ldots,f_m)$ and $g$ has the same leading term as $f_{m+1}.$ Thus $f_{m+1}-g\in I$ and $f_{m+1}-g\notin (f_1,\ldots,f_m)$
	but has degree strictly less than $\mbox{deg}f_{m+1},$ which is a contradiction since $f_{m+1}$ 
	had minimal degree. This completes the proof. 
	\end{proof}
\begin{theorem}[Gauss]\label{gaussfactorial}
	A ring $R$ is factorial iff $R[x_1,\ldots,x_n]$ is factorial for each $n\in \N.$
\end{theorem}
\begin{proof}
We prove the result for $x=x_1$. The general case can be proved by induction.
\begin{lemma}\label{glemma}
	Let $R$ be factorial and let $K$ be the quotient field of $R$. Let $f,g\in K[x]$
	be polynomials in one variable.
	Then the $\mbox{cont}(fg)=\mbox{cont}(f)\mbox{cont}(g).$
\end{lemma}
\begin{proof}
	By Remark \ref{remglem} we can write $f=c f_1,$ $g=dg_1$ where
	$c=\mbox{cont}(f),$
	$d=\mbox{cont}(g)$, thus it suffices to prove that if $f,g$ have content $1$ then also $fg$ has content $1$.
	For this it suffices to prove that for each prime $p$,
	$\mbox{ord}_p(fg)=0$. Let $f(x)=\sum_{j=0}^n a_j x^j,$ $a_n\neq 0,$
	$g(x)=\sum_{j=0}^n b_j x^j,$ $b_n\neq 0,$
	be polynomials of content $1$ and let $p\in R$ be a prime element. It suffices to prove that $p$ does not divide all coefficients of $fg$.
	Let $r$ be the largest integer such that $0\leq r\leq n,$ $a_r\neq 0$ and $p$ does not divide $a_r$.
	Let $b_s$ be the coefficient of $g$ with least $s$ such that $b_s\neq 0$ such that $p$ does not divide
	$b_s.$
	The coefficient of $x^{r+s}$ in $f(x)g(x)$ is $c=a_rb_s+a_{r+1}b_{s-1}+\cdots+a_{r-1}b_{s+1}+\cdots$
	and $p$ does not divide $a_r b_s.$ On the other hand $p$ divides every nonzero term in the latter sum since
	each term is some coefficient $a_i$ or $b_j$ with $i$ or $j$ less than $r$ or $s$ respectively.
	Hence $p$ does not divide $c$. This completes the proof. 
\end{proof}

Let $f\in R[x]$, $f\neq 0$. By the unique factorization in $\C[x]$ we have a factorization
\begin{equation}
f(x)=c\cdot p_1(x)\cdots p_r(x),\quad c\in R
\end{equation} 
where $p_j\in R[x], j=1,\ldots,r$ are polynomials which are irreducible in $K[x].$ W.l.o.g.\ 
we may (after extracting their contents) assume that $\mbox{cont}(p_i)=1$ for each $i$. Then $c=\mbox{cont}(f)$ by Lemma \ref{glemma}.
This proves the existence of factorization. Suppose we have another factorization
\begin{equation}
f(x)=d\cdot q_1(x)\cdots q_s(x),\quad d\in R
\end{equation} 
By the unique factorization in $K[x]$ we have $r=s$ and after permutation of factors we have 
$p_i=a_i q_i,$ $i=1,\ldots,r,$ for $a_i\in K.$ Since both $p_i,q_i$ are assumed to have content $1$,
we have $a_i\in R$ and $a_i$ is a unit for $i=1,\ldots,r.$ This completes the proof.
\end{proof}
\begin{lemma}
$\mathscr{O}_0$ is noetherian and factorial.
\end{lemma}
\begin{proof}
We use induction in the complex dimension $n$ of the ambient space $\Cn$.
	For $n=0$, i.e.\ $\mathscr{O}_0=\C$ we have a field and trivially a noetherian ring that is a unique factorization domain. 
	Assume as induction hypothesis that for dimensions $0,\ldots,n-1$, $\mathscr{O}_0$ is noetherian and let $I$ be an ideal.
	Let $g\in I$ be nonzero. Then (if necessary after a nonsingular linear transformation and if necessary multiplication of $g$ by a unit) we can assume that $g\in I\cap \mathscr{O}_0^{(n-1)}[z_n]$ is a Weierstrass polynomial, where $\mathscr{O}_0^{(n-1)}$ denotes the ring of germs of holomorphic functions at $0$ with respect to the first $n-1$ complex variables. By the Hilbert basis theorem $\mathscr{O}_0^{(n-1)}[z_n]$ is noetherian, hence the ideal $I\cap \mathscr{O}_0^{(n-1)}[z_n]$ is finitely generated by say $g_1,\ldots,g_m.$
	For any $f\in I$ we have by the Weierstrass division theorem that $f=gh +r$ for some $h\in \mathscr{O}_0$ and $r\in \mathscr{O}_0^{(n-1)}[z_n].$ But $r\in I$ so we also have $r=\sum_j g_j h_j$ for some $h_j\in \mathscr{O}_0^{(n-1)}[z_n],$ i.e.\ $f=gh+\sum_j g_jh_j.$ Hence $g,g_1,\ldots,g_m$ generate $I$. This proves that $\mathscr{O}_0$ is noetherian. Now Assume as induction hypothesis that for dimensions $0,\ldots,n-1$, $\mathscr{O}_0$ is factorial. Let $f\in \mathscr{O}_0$.
	Then (if necessary after a nonsingular linear transformation and if necessary multiplication of $g$ by a unit) we can assume that $f=uh$ for a unit $u\in \mathscr{O}_0$ and a Weierstrass polynomial $h\in \mathscr{O}_0^{(n-1)}[z_n]$. By Theorem \ref{gaussfactorial} $\mathscr{O}_0^{(n-1)}[z_n]$ is factorial hence $h$ can be written (uniquely up to order of factors and up to units) as a product of irreducible polynomials.
	This proves that $\mathscr{O}_0$ is factorial. This completes the proof.
	\end{proof}
Since $\mathscr{O}_0$ is noetherian it follows from Lemma \ref{hilbertbasis} that
$\mbox{PA}(0)$ is noetherian. 
Also by Theorem \ref{gaussfactorial} , $R[x_1,\ldots,x_N]$ is factorial if $R$ is factorial hence
$\mbox{PA}(0)$ is factorial.	
This completes the proof.
\end{proof}
Hence any polyanalytic function $f$ near the origin satisfying
$f(0)=0$ can be factorized into a finite product of irreducible factors, uniquely up to order and invertible factors.
However, some of the factors may be nonzero at the origin.
It is clear that the zeros of a polyanalytic function $f(z)=\sum_{j=0}^q a_j \bar{z}^j$ on a domain $\Omega\subset\Cn$ is the union of two
analytic sets, namely one complex analytic set given  by $a(z):=(a_0,\ldots,a_{q-1})=0$ and one real analytic set
given by the equation
\begin{equation}
\norm{a\wedge Z}(z)=(\norm{a}\norm{Z})(z)
\end{equation}
where
\begin{equation}
Z=(1,z^1,\dots,z^{n-1}),\quad \norm{a\wedge Z}^2:=\sum_{j,k}\abs{a_z^k -a_kz^j}^2
\end{equation}
In particular, on any given compact there can only exists finitely many isolated zeros.
\begin{definition}
Given two subsets $A_1,A_2\subset\Cn$, we say that $A_1$ and $A_2$ are equivalent at a point $p\in A_1\cap A_2$ if there exists an open neighborhood $U$ of $p$ in $\Cn$
such that $A_1\cap U=A_2\cap U.$ This renders an equivalence relation on subsets and we call the equivalence class $[A_1]_p$ at $p$ the germ of $A_1$ at $p$.
Let $\Omega\subseteq\Cn$ be a domain. A subset $A\subset \Omega$ is called a {\em polyanalytic subset}\index{Polyanalytic subset}
if it is closed in $\Omega$ and coincides with common set of zeros of finitely many polyanalytic functions
on $\Omega$. Note that the germ of a polyanalytic subset can be identified with the common zeros of germs of the defining polyanalytic functions.
\end{definition}

\begin{proposition}
If $[A]_0$ is a germ of a polyanalytic subset $A$ at $0$ then $[A]_0$ is the union of finitely many irreducible 
germs of polyanalytic subsets $[A_1]_0,\ldots,[A_N]_0,$ (for a fixed $N\in \N$) and the decomposition is unique under the condition that
$[A_i]_0\not\subset [A_j]_0$ for $i\neq j.$
\end{proposition}
\begin{proof}
This is immediate from Proposition \ref{notherfactorprop} since the ring of germs of polyanalytic functions at $0$	
is noetherian and factorial.
\end{proof}
\begin{remark}
If $[A]_0$ is a germ of a polyanalytic subset $A$ at $0$ we denote by $I([A]_0)$ the ideal in $\mbox{PA}(0)$ consisting of the set of
 germs of polyanalytic functions near $0$ that vanish on $A.$
It is clear that if $[A_1]_0$ and $[A_2]_0$ are germs of polyanalytic subsets $A_1,A_2$ at $0$ then
$[A_1]_0\subseteq [A_2]_0\Rightarrow I([A_2]_0)\subseteq I([A_1]_0)$ and furthermore
$[A_1]_0$ is prime iff $I([A_1]_0)$ is prime.
Note that in particular we have 
$[A_1]_0 = [A_2]_0\Rightarrow I([A_2]_0) = I([A_1]_0)$.
It is important to note that the germ of a polyanalytic subset is not uniquely induced, namely if $f$ is a polyanalytic function
with polynomial analytic components then $\bar{f}$ is also polyanalytic and with the same zero set, however the ideal generated by $[f]_0$ does not
necessarily contain $[\bar{f}]_0.$
\end{remark}
\begin{theorem}
Suppose $f$ is a polyanalytic function near $0$ in $\C$ such that $f(0)=0$, that $0$ is a non-isolated zero of $f$ and that $[f]_0$ is irreducible.
Denote by $Z([f]_0)$ the germ induced by the zero set of $f.$ Then $[f]_0$ divides any germ $[F]_0$ such that $[F]_0$ that vanishes on $Z([f]_0).$ 
\end{theorem}
\begin{proof}
Let $f$ be $q$-analytic near $0$ and
there have representation
$f(z)=\sum_{j=0}^{q-1} a_j(z)\bar{z}^j$, for holomorphic $a_j$, $j=0,\ldots,q-1.$
Let $[F]_0$ be zero on $Z([f]_0),$ $F=\sum_{j=0}^{m-1} b_j(z)\bar{z}^j$ for holomorphic $b_j,$ 
$j=0,\ldots,m-1.$ Suppose (in order to reach a contradiction) that $[f]_0$
does not divide $[F]_0.$ Then $[f]_0$ and $[F]_0$ have g.c.d.\ $1$. Consider the polynomials
$\phi=\sum_{j=0}^{q-1} [a_j]_0 \omega^j$ associated to $f$, and
$\psi=\sum_{j=0}^{m-1} [b_j]_0 \omega^j$ associated to $F$, both
obtained by application of the map $\mathbf{j}.$ These also must have g.c.d.\ $1$. Let $A$
be the ring of germs of holomorphic functions near $0$ and let $\alpha$
be the field of quotients of $A.$ Then $\phi$ and $\psi$ have g.c.d.\ $1$ also in 
$\alpha[\omega].$ By the Euclidean algorithm there exists $\lambda,\mu\in \alpha[\omega]$ such that $\lambda \phi+\mu\psi=1.$ After 
rewriting the terms to having the same denominator we see that there exists a germ $[g]_0\neq 0$
of a holomorphic function in $z$ together with 
$\lambda',\mu'\in \alpha[\omega]$ such that
$\lambda'\phi +\mu'\psi=[g]_0.$ Set $\phi_1(z,\omega):=\sum_{j=0}^{q-1} a_j(z)\omega^j$
and $\psi_1(z,\omega):=\sum_{j=0}^{m-1} b_j(z)\omega^j$.
This shows that there exists a neighborhood, $U$ of $0$ together with a holomorphic function $g_1$
on $U$, $g_1\not\equiv 0$
and polynomials with holomorphic coefficients, $\lambda_1,\mu_1$ such that on $U$
we have $g_1=\lambda_1\phi_1+\mu_1\psi_1.$ This implies (replacing the roles of $\omega$ and $\bar{z}$)
that there exists polyanalytic $\lambda_1',\mu_1'$ such that 
\begin{equation}
\lambda_1'(z)f(z)+\mu'_1(z)F(z)=g_1(z)
\end{equation}
It follows that the zeros of $g_1(z)$ have $z=0$ as accumulation point which means that $g_1\equiv 0$
which is a contradiction. This completes the proof.
\end{proof}
\begin{corollary}
If $f$ is a polyanalytic function in one variable, near $0$, such that $[f]_0$ is irreducible and if $z=0$ is a non-isolated zero of $f$ then $Z([f]_0)$ is irreducible. 
\end{corollary}

\begin{corollary}
	Suppose $f$ is a polyanalytic function in one variable, near $0$, such that $[f]_0=\Pi_{i=1}^m [p_i]^{n_i}$ 
	where each factor is irreducible. If for each $i$, $0$ is a non-isolated zero of $p_i$, then
	$Z([f]_0)=\bigcup_{i=1}^m Z([p_i]_0)$ is a decomposition with irreducible components. 
\end{corollary}
Let $p\in \C[x,y]$ be irreducible such that $0$ is a non-isolated zero of $p.$ in $\R^2.$ Let $Z$ denote the set of zeros of $p$
in $\R^2.$ If $q\in \C[x,y]$ such that there exists an open neighborhood $U$ of $0$ satisfying $q|_{U\cap Z}=0$, then there exists
$\lambda\in \C[x,y]$ such that $q=\lambda p.$ 
Setting $z=x+iy,$ $\bar{z}=x-iy$ we see that $p(z,\bar{z})$ becomes a polyanalytic polynomial in one complex variable. The ring of germs of polyanalytic polynomials is a subring of $\mbox{PA}(0)$ thus by the preceeding results apply.

\chapter{Zero sets}
\section{Some preliminary proof-tools}
In this section we describe some important tools appearing in the study of values of $q$-analytic functions and in particular in the study of their zero sets, namely Schwarz functions.
This can be used to identify associated functions to a given polyanalytic function, that agree with the function on a given analytic arc. As Schwarz functions shall also be used
as proof tool in upcoming chapters we shall present also several notable properties of these functions that are not necessarily used in the results on zero sets.
We shall also need to be familiar with some notions from Nevanlinna theory (these will reappear in Chapter \ref{polyentiresec}).

\subsection{Schwarz functions}\label{schwarzsec}
A closed {\em analytic(Jordan) curve}\index{Analytic curve} is a 
curve which can be parametrized by a function $z= x + iy = z(w)$,
analytic in a domain which includes the unit circle $\abs{w}= 1.$
Let $\gamma$ be a closed Jordan analytic curve that bounds a domain $\Omega\subset\C.$
Let $w=\phi(z)$ be a one-sheeted conformal mapping of $\Omega$ onto the unit disc and let $\psi$ be its inverse mapping.
Then the equation of $\gamma$ can be written $A_\gamma(z) =\bar{z},$ where $A_\gamma(z):=
\overline{\psi(1/\overline{\phi(z)})}$ is called the {\em Schwarz function}
of $\gamma.$ $A_\gamma(z)$ is holomorphic and one-sheeted in an open neighborhood of $\gamma.$ 
See e.g\ Davis \cite{davis} for the theory of Schwarz functions.
A way to define an {\em analytic arc}, $\gamma$ in $\R^2\simeq \C,$ is as follows.
$\gamma$ is such an arc if and only if it is parametrized by $z=\phi(t),$ $t\in [0,1],$ 
for a function $f$, analytic on an open neigbhorhood of $[0,1]$ that is bijective
on $[0,1]$ and satisfies $\phi'(t)\neq 0,$ $t\in [0,1].$ Equivalently,
$\gamma$ is such an arc if and only if for each point, $p_0=x_0+iy_0\in \gamma$ 
there exists an open neighborhood 
$U\subset\C^2$ of the point  $(x_0+i0,y_0+i0)\in \C^2$
and a function $\phi\in \mathscr{O}(U)$ such that 
$\phi(v,w)$ is real for real arguments $v,w$, ($v=x+ix',w=y+iy'$) and such that
\begin{equation}
\gamma\cap U=\{ (v,w)\in U\cap \R^2 :\phi(v,w)=0
\end{equation}
and the gradient of $\phi$ is nonzero at $(0,0).$

\begin{proposition}\label{schwarzfuncprop}
	Let $\gamma$ be a nonsingular analytic Jordan arc in $\C.$ Then there exists an open neighborhood $U$
	of $\gamma$ and a uniquely determined holomorphic function
	$S(z)$ on $U$ such that $\gamma=\{z:S(z)=\bar{z}\}.$
\end{proposition}
\begin{proof}
	It suffices to construct $S$ near a single point $p_0\in \gamma$ since the relation
	$S(z)=\bar{z}$ will allow us to glue together these local determinations of $S.$
	By definition there exists a holomorphic function $\phi(s,t)$, $s=x+ix', t=y+iy',$ on an open neighborhood, $V$, of $p_0$,
	with nonvanishing gradient at $p_0$ such that $V\cap \gamma=\{(x,y)\in \R^2: \phi(x,y)=0\}.$
	Set
	\begin{equation}
	f(z,w):=\phi\left(\frac{z+w}{2},\frac{z-w}{2i}\right)
	\end{equation}
	where $(z,w)$ varies in a neighborhood of $(p_0,\overline{p_0})$ in $\C^2.$
	Then $f(p_0,\overline{p_0})=0$ and 
	\begin{equation}
	\partial_w f =\frac{1}{2} \partial_s\phi\left(\frac{z+w}{2},\frac{z-w}{2i}\right)
	-\frac{1}{2i} \partial_t\phi\left(\frac{z+w}{2},\frac{z-w}{2i}\right)
	\end{equation}
	Hence
	\begin{equation}
	\partial_w f(p_0,\overline{p_0})=\frac{1}{2}(\partial_x\phi(p_0,\overline{p_0})+i\partial_y\phi(p_0,\overline{p_0})\neq 0
	\end{equation}
	By the complex version of the implicit function theorem there exists a unique complex-valued $w=S(z)$
	defined in an open neighborhood of $z=p_0$ such that
	$S(p_0)=\overline{p_0}$ and $f(z,S(z))=0,$ and moreover $S$ is analytic. Now for $z\in \gamma$ and near $p_0$
	$f(z,\bar{z})=\phi(z)=0$ on $\gamma,$ i.e.\
	the equation $f(z,w)=0$ is solved by $w=\bar{z}.$ This implies that $S(z)=\bar{z}$ for $z\in \gamma$, $z$ sufficiently near $p_0.$ 
	By the beginning remark, this proves existence. Uniqueness follows from the fact that
	holomorphic functions do not vanish on sets with an accumulation point unless they vanish identically. This completes the proof.
\end{proof}

In the theory of Schwarz functions it is well-known that if $\gamma$ is a closed analytic arc bounding a 
simply connected bounded
domain then the Schwarz function $S(z)$ of $\gamma$ has at least a singular point in $\Omega$ 
and if $S(z)$ is meromorphic on $\Omega$ then $\gamma$ is algebraic (see e.g.\ Davis \cite{davis} and Mathurin \cite{mathurin}). Let give some
of the details.

\begin{theorem}
	The Schwarz function of an analytic arc $\gamma$ is rational in $z$ if and only if 
	$\gamma$ is an arc of a circle or a straight line.
\end{theorem}
\begin{proof}
	The Schwarz function of a circle $\abs{z-p_0}=\rho$ for some $\rho>0$ and fixed $p_0\in \C$
	is given by
	$S(z)=\frac{\rho^2}{z-p_0}+\bar{p}_0$, thus rational. Conversely, suppose $S(z)=P(z)/Q(z)$
	where
	$P(z)=a_0z^m+\cdots,$
	$Q(z)=b_0z^n+\cdots ,$ $a_0b_0\neq 0,$ where greatest common divisor of $P$ and $Q$ is $1.$
	We have $\bar{R}(R(z))=z$ which implies
	\begin{equation}
	\frac{\bar{a}_0R^m+\bar{a}_1 R^{m-1}+\cdots}{\bar{b}_0 R^n+\bar{b}_1 R^{n-1}+\cdots}=z
	\end{equation}
	or
	\begin{equation}
	(\bar{a}_0R^m+\bar{a}_1 R^{m-1}+\cdots)=z(\bar{b}_0 R^n+\bar{b}_1 R^{n-1}+\cdots)
	\end{equation}
	Suppose $m>n$.
	Then 
	\begin{equation}
	\bar{a}_0P^m+\bar{a}_1 P^{m-1}+\cdots=Q^m(\bar{b}_0 R^n+\bar{b}_1 R^{n-1}+\cdots)
	\end{equation}
	so $Q$ divides $\bar{a}_0 P^m$. Since $P$ and $Q$ have greatest common divisor $1$
	we obtain that $Q\equiv$constant. Thus $R=\alpha_0 z^n+\alpha_1 z^{m-1} +\cdots$. Since $\bar{R}(R(z))=z$
	\begin{equation}
	z=\bar{\alpha}_0(\alpha_0 z^m+\alpha_1 z^{m-1} +\cdots )^m+\bar{\alpha}_1(\alpha_0 z^m +\cdots )^{m-1}+\cdots 
	\end{equation}
	Equating coefficients gives
	$\alpha_0\bar{\alpha}_0z^{m^2}=z$ thus $m=1$ and $\alpha_0\bar{\alpha}_0=1.$ Thus
	$R(z)=\alpha z+\beta,$ $\alpha\bar{\alpha}=1,$ $\re (\bar{\alpha}\beta )=0.$ This takes care of the case $m>n.$
	\\
	Suppose instead $m<n.$ Then 
	$Q$ divides $\bar{b}_0 z P^m$ so since $P$ and $Q$ have greatest common divisor $1$
	we obtain that $Q\equiv z\cdot\mbox{const.}$
	Since $P$ has smaller degree than $Q$ this implies $R(z)=c/z$. 
	Since $\bar{R}(R(z))=z$ this implies $c=\bar{c}.$ This takes care of the case $m<n.$
	\\
	\\
	Finally, suppose $m=n.$
	\begin{equation}
	\bar{a}_0P^m+\bar{a}_1 P^{m-1}+\cdots=z(\bar{b}_0 P^m+\bar{b}_1 P^{m-1}+\cdots)
	\end{equation}
	This implies $\bar{a}_0-\bar{b}_0 z=0$mod$Q$, hence
	$Q=cz+d$ which gives $P=az+b$ for constants $a,b,c,d\in \C.$
	Set
	Then $\bar{R}(R(z))=z$ implies
	\begin{equation}
	\begin{bmatrix}
	\bar{a} & \bar{b}\\
	\bar{c} & \bar{d}
	\end{bmatrix}\begin{bmatrix}
	a & b\\
	c & d
	\end{bmatrix}=\lambda I
	\end{equation}
	We consider matrices 
	\begin{equation}
	M=\begin{bmatrix}
	a & b\\
	c & d
	\end{bmatrix},\quad \bar{M}M=I
	\end{equation}
	and we
	may assume det$(M)=1.$ In particular, we must have
	\begin{equation}
	\begin{bmatrix}
	\bar{a} & \bar{b}\\
	\bar{c} & \bar{d}
	\end{bmatrix}=
	\begin{bmatrix}
	d & -b\
	-c & a
	\end{bmatrix}
	\end{equation}
	Thus
	\begin{equation}
	M=\begin{bmatrix}
	a & ib\
	ic & \bar{a}
	\end{bmatrix}
	\end{equation}
	where $a\bar{a}+bc=1.$ Hence $R(z)$ takes the form
	\begin{equation}
	R(z)=\frac{\rho^2}{z-p_0}+\bar{p}_0
	\end{equation}
	with $\rho=\frac{1}{c},$ $p_0=\frac{i\bar{a}}{c}.$
	This case care of the case $m=n.$
	This completes the proof.
\end{proof}

\begin{theorem}\label{schwarzapplic1}
	Let $\Omega$ be a bounded, simply connected domain whose boundary is an analytic curve, $\gamma,$
	with Schwarz function $S(z).$ Then $S(z)$ is meromorphic on $\Omega$ if and only if there exists a linear functional $L$
	defined on the algebra, $\mathscr{O}(\Omega)\cap C^0(\overline{\Omega})$, of functions, $f$, analytic in $\Omega$ and continuous on $\overline{\Omega},$ such that
	$L$ can be represented as
	\begin{equation}
	L(f)=\sum_{j=1}^N\sum_{k=0}^{j_k} a_{jk} f^{(k)}(z_j)
	\end{equation}
	where $z_1,\ldots,z_N$ are distinct points in $\Omega$ and $a_{jk}$ are constants independent of $f$;
	and such that for all $f\in \mathscr{O}(\Omega)\cap C^0(\overline{\Omega})$
	\begin{equation}\label{tqrepppp}
	\int_\Omega f(z)d\mu(z)=L(f)
	\end{equation}
\end{theorem}
\begin{proof}
	Let $g_1,f$ belong to $\mathscr{O}(\Omega)\cap C^0(\overline{\Omega})$.
	In the Green formula
	\begin{equation}
	\int_\gamma Pdy-Qdx=\int\Omega (P_x+Q_y)dx\wedge dy
	\end{equation}
	set $P=\frac{1}{2}\bar{g_1}f,$ $Q=\frac{i}{2}\bar{g}_1f$.
	Then $P_x=\frac{1}{2}(\partial_z+\partial_{\bar{z}})\bar{g}_1f$ and
	$Q_y=-\frac{1}{2}(\partial_z -\partial_{\bar{z}})\bar{g}_1f$ so $P_x+Q_y=\partial_{\bar{z}} (\bar{g}_1f)=
	\partial_{\bar{z}} (\bar{g}_1(\bar{z})f(z))=\overline{g_1'(z)}g_2(z).$ This gives
	\begin{multline}
	\int_\Omega \overline{g_1'(z)}f(z)  =\frac{1}{2i}\int_\gamma \overline{g_1(z)}f(z)dz=
	\frac{1}{2i}\int_\gamma \bar{g}_1(\bar{z}z)f(z)dz
	\end{multline}
	Since on the boundary we have $\bar{z}=S(z)$ this gives
	\begin{equation}
	\int_\Omega \overline{g_1'(z)}f(z)=\frac{1}{2i}\int_\gamma \bar{g}_1(S(z))f(z)dz
	\end{equation}
	Setting $g_1\equiv z$ yields
	\begin{equation}
	\int_\Omega f(z)=\frac{1}{2i} \int_\gamma S(z)f(z)dz
	\end{equation}
	Suppose $S(z)$ is meromorphic. Then $S(z)$ has only finitely 
	many poles $z_1,\ldots,z_N$ in $\Omega.$
	Hence for sufficiently small circles, $C_j,$ in $\Omega$, centered at $z_j,$ $j=1,\ldots,N$ we have
	\begin{equation}
	\frac{1}{2i}\int_\gamma S(z)f(z)dz=\sum_{j=1}^N \frac{\pi}{2\pi i} \int_{C_j} S(z)f(z)dz
	\end{equation}
	Inside $C_j$, $S(z)$ has an expansion 
	\begin{equation}
	S(z)=\zeta_j(z)+\frac{B_{1j}}{z-z_j}+\frac{B_{2j}}{(z-z_j)^2}+\cdots +\frac{B_{p_j j}}{(z-z_j)^{p_j}}
	\end{equation}
	for regular $\zeta_j(z)$ on the closed disc defined by $C_j.$ This gives
	\begin{equation}
	\frac{\pi}{2\pi i} \int_{C_j} S(z)f(z)dz    =\pi\sum_{k=1}^{p_j}\frac{B_{k j}f^{(k-1)}(z_j)}{(k-1)!}
	\end{equation}
	This proves Eqn.(\ref{tqrepppp}).
	\\
	Conversely suppose Eqn.(\ref{tqrepppp}) holds true for all $f\in \mathscr{O}(\Omega)\cap C^0(\overline{\Omega})$ 
	and for some $L$ as in the theorem.
	Setting
	\begin{equation}
	R(z)=\frac{1}{\pi} \sum_{j=1}^N\sum_{j=1}^{j_k} \frac{k!c_{jk}}{(z-z_j)^{k+1}}
	\end{equation}
	we have for all $f\in \mathscr{O}(\Omega)\cap C^0(\overline{\Omega})$ 
	\begin{equation}
	\int_\gamma (S(z)-R(z))f(z)dz =0
	\end{equation}
	By a known result on conditions for continuous functions to be boundary values of holomorphic functions, Corollary \ref{schwarzapplic1lem1}
	(in this book), we know that
	$S(z)-R(z)|_\gamma$ determines the boundary value of a function
	$\zeta(z)$ in $\mathscr{O}(\Omega)\cap C^0(\overline{\Omega})$.
	But this implies by the uniqueness of holomorphic functions
	that $S(z)-R(z)$ share the same region of analyticity, thus 
	$S(z)=R(z)+\zeta(z)$ on $\Omega$ thus is meromorphic on $\Omega.$
	This completes the proof.
\end{proof}

\begin{theorem}\label{schwarzapplic2}
	Let $\Omega$ be a bounded, simply connected domain whose boundary is an analytic curve, $\gamma,$
	with Schwarz function $S(z).$ Suppose $z=0$ belongs to $\Omega$ and $z=\phi(w), \phi(0)=0$ maps the unit disc conformally one-to-one
	onto $\Omega$. Then $S(z)$ is meromorphic on $\Omega$ if and only if 
	$\phi(w)$ is rational.
\end{theorem}
\begin{proof}
	Suppose $\phi$ is rational and write
	\begin{equation}
	\phi(w)=\frac{aw(1-\beta_1w)\cdots(1-\beta_p w)}{(1-\alpha_1w)\cdots(1-\alpha_n w)}=aw+\cdots
	\end{equation}
	where the complex constants $\alpha_j$ are different from the $\beta_k$.
	We have
	\begin{multline}
	\int_\Omega f(z) dx\wedge dy =\frac{1}{2\pi i} \int_\gamma \bar{z} f(z)dz=\\
	\frac{1}{2\pi i} \int_{\abs{w}=1} \overline{\phi(w)}f(\phi(w))\phi'(w)dw=
	\frac{1}{2\pi i} \int_{\abs{w}=1} \bar{\phi}\left(\frac{1}{w}\right)f(\phi(w))\phi'(w)dw
	\end{multline}
	where 
	\begin{equation}
	\bar{\phi}\left(\frac{1}{w}\right)=\frac{\bar{a}w^{n-p-1}(1-\bar{\beta}_1w)\cdots(1-\bar{\beta}_p w)}
	{(1-\bar{\alpha}_1w)\cdots(1-\bar{\alpha}_n w)}
	\end{equation}
	If $n\geq p+1$ then $\bar{\phi}(1/w)$ has poles $\bar{\alpha}_1,\ldots,\bar{\alpha}_n$
	and no other singularities. 
	If the $\alpha_j$ are distinct the by the residue theorem
	\begin{equation}
	\int_\Omega f(z)dx\wedge dy =\pi\bar{a}\sum_{k=1}^m f(\phi(\bar{\alpha}_k))\phi'(\bar{\alpha}_k)
			\frac{\bar{\alpha}_k^{n-p-1}(\bar{\alpha}_k-\bar{\beta}_1)\cdots(\bar{\alpha}_k-\bar{\beta}_p)}{\bar{P}'(\alpha_k)}
			\end{equation}
			where
			\begin{equation}
			P(w)=(w-\alpha_1)\cdots (w-\alpha_n)
			\end{equation}
			Thus
			$\int_\Omega f(z) dx\wedge dy$
			takes the form
			$\sum_{k=1}^m c_k f(z_k)$ where the $c_k$ and the abscissas $z_k=\phi(\bar{\alpha}_k)$
			are independent of $f(z).$ If the $\alpha_j$ are not distinct then each point of higher order multiplicity, say $\tau_k$,
			contributes a differential operator of order $\tau_k-1$ evaluated at $\bar{\alpha}_k$ and we again obtain
			the wanted representation for $f(z).$
			This proves the result under the assumption that 
			$n\geq p+1.$
			If instead $n<p+1$ then
			it also has a pole of order $p+1-n$ at the origin, $\alpha=0$ and since $\phi(0)=0$, $f(0)$,
			we have that $f(z)$ in the case $n-p-1=-1$, is represented by itself, and for
			$n-p-1<-1$ is represented by its higher order derivatives.
			In both cases $\int_\Omega f(z) dx\wedge dy=L(f)$ for a linear functional
			satisfying the conditions of Theorem \ref{schwarzapplic1}. Hence $S(z)$ must be meromorphic
			on $\Omega.$ We have thus proved that if $\phi$ is rational then $S(z)$ is meromorphic on $\Omega.$
			\\
			Conversely suppose that $S(z)$ is meromorphic on $\Omega.$ Then by 
			Theorem \ref{schwarzapplic1} we have a linear functional satisfying the conditions of Theorem \ref{schwarzapplic1} such that
			$\int_\Omega fdx\wedge dy=L(f)$ for all
			$f\in \mathscr{O}(\Omega)\cap C^0(\overline{\Omega})$.			
			Denote by $M$ the map of $\Omega$ conformally one-to-one onto the unit disc
			$\abs{w}\leq 1$ with $M(0)=0.$
			By the Riesz representation theorem a linear functional $L$ on $L^2$
			has a representation for $f\in L^2$ according to
			$L(f)=\langle f,r\rangle=\int_\Omega f(z)\overline{r(z)} dx\wedge dy$.
			If furthermore $L(f)=\int_\Omega f(z) dx\wedge dy$
			then $r(z)\equiv 1.$
			Now the Bergman reproducing kernel for the unit disc
			takes the form (see Chapter \ref{reproducingsec}, which in Section
			\ref{realandimagsec} is exemplified by the Poisson kernel and in
			Chapter \ref{bvpsec} by e.g.\ the Schwarz kernel) 
			\begin{equation}
			K(\bar{z},w)=\frac{1}{\pi}\frac{\overline{M'(z)}M'(w)}{(1-\overline{M(z)}M(w))^2}
			\end{equation}
			which implies
			\begin{equation}
			\overline{L_\eta K(\bar{z},\eta)}\equiv 1
			\end{equation}
			This is a functional equation for a linear functional whose Riesz representation is the 
			constant function $1.$ This implies
			\begin{equation}
			1\equiv \frac{1}{\pi}\overline{L_\eta \left(\frac{\overline{M'(z)}M'(\eta)}{(1-\overline{M(z)}M(\eta))^2}\right)}
				\end{equation}
				or in terms of a series, an absolute uniformly convergent on compacts of $\Omega$ given by
				\begin{equation}
				1\equiv \frac{1}{\pi}\sum_{j=0}^\infty (j+1)(M(z))^j M'(z)\overline{L((M(\eta))^j M'(\eta))}
				\end{equation}
				Thus
				\begin{equation}
				1\equiv \frac{1}{\pi} \frac{d}{dz} \sum_{j=0}^\infty (M(z))^{j+1} M'(z)\overline{L((M(\eta))^j M'(\eta))}
				\end{equation}
				This yields
				\begin{equation}
				\pi z=M(z) \sum_{j=0}^\infty (M(z))^{j} M'(z)\overline{L((M(\eta))^j M'(\eta))}
				\end{equation}
				thus 
				\begin{equation}
				\bar{z}=\overline{M(z)} L_\eta \left(\frac{M'(\eta)}{1-\overline{M(\eta)}M(z)}\right)
			\end{equation}
			Since $z=\phi(w),\phi(0)=0$ is the inverse to $w=M(z)$ on $\Omega$
			this yields
			\begin{equation}
			\pi \overline{\phi(w)}=\bar{w} L_\eta \left(\frac{M'(\eta)}{1-\bar{w}M(z)}\right)
		\end{equation}
		Since $L$ is a point differential functional this implies
		$L_\eta((M'(\eta))/(1-\bar{w}M(\eta)))$, and therefore also $\overline{\phi(w)}$, is a rational function of $\bar{w}$,
		thus $\phi(w)$ is rational.
		This completes the proof.
	\end{proof}
\subsection{Some lemmas from Nevanlinna theory}\label{nevanlinnasec}
We need to be familiar with some notions from Nevanlinna theory (these will reappear in Chapter \ref{polyentiresec}).
See e.g.\ Goldeberg \& Otrovskii \cite{goldostrov}. Here we only 
point out, for future reference (see e.g.\ Chapter \ref{polyentiresec}), the following results from Nevanlinna theory
Let $\Omega\subset\C$ be a bounded domain whose boundary, $\gamma$, consists of finitely many
piecewise analytic curves.
Denote by $\frac{\partial}{\partial n}$ the operator of differentiation in the direction of the
inward directed normal to $\gamma.$
Let $\sigma$ be the Lebesgue measure in the plane and let $s$ denote the one dimensional
Lebesgue measure on $\gamma.$ Let $f(z)=f(x,y)=u(x,y)+iv(x,y)$ be a $C^2$ function
on $\overline{\Omega}.$
The Green formula in the plane
\begin{equation}
\int_\Omega (\partial_x Q -\partial_y P)d\sigma =\int_\gamma Pdx+Qdy
\end{equation}
yields with $P=-u\partial_y v,$ $Q=u\partial_x v$
\begin{multline}
\int_\Omega (u\Delta v+\partial_x u\partial_x v+\partial_y u\partial_x v)d\sigma =\\
\int_\gamma 
-u\partial_y v dx +u\partial_x v dy =-\int_\gamma u \frac{\partial v}{\partial n}ds
\end{multline}
Subtracting from this the formula which is obtained after interchanging $u$ and $v$ yields
\begin{equation}\label{greenekv2}
\int_\Omega (u\Delta v-v\Delta u)d\sigma =-\int_\gamma \left(u\frac{\partial v}{\partial n}
-v\frac{\partial u}{\partial n}\right)ds
\end{equation}
\begin{theorem}\label{goldbergnevanthm1}
	Let $\Omega$ be a simply connected domain with piecewise analytic boundary $\gamma$, let 
	$\gamma_1,\ldots,\gamma_n$ be analytic arcs forming $\gamma=\partial\Omega$ and let $a_j,$ 
	$j=1,\ldots,n$ be the common end point of $\gamma_j$ and $\gamma_{j+1}$ ($\gamma_{n+1}=\gamma_1$ for $j=n$)-
	Let $\pi\alpha_j$, $0+\leq \alpha_l\leq 2$ be the angles between $\gamma_j$ and $\gamma_{j+1},$ and suppose
	$\Omega$ has
	the property defined by Eqn.(\ref{goldberganvsen}).
	Let $u(z)$ be a $C^2$-smooth 
	function on a neighborhood of $\overline{\Omega}\setminus \{c_1,\ldots,c_q\}$ for complex points $c_j\in \overline{\Omega}$.
	Suppose that in a neighborhood of each $c_j$ the function $u(z)$ has the form
	\begin{equation}\label{goldnevanekv12}
	u(z)=d_j\ln\abs{z-c_j}+u_j(z),\quad j=1,\ldots,q
	\end{equation}
	for constants $d_j$ and $u_j(z)$ $C^2$-smooth in a neighborhood of $c_j.$
	Then
	\begin{equation}\label{goldnevanekv13}
	u(z)+\frac{1}{2\pi}\int_\Omega g(\zeta,z)\Delta u(\zeta)d\sigma = \frac{1}{2\pi}\int_\gamma u(\zeta)
	\frac{\partial g}{\partial n}ds -\sum_{j=1}^q d_j g(c_j,z)
	\end{equation}
	for $z\in \Omega\setminus \{c_1,\ldots,c_q\}.$
\end{theorem}
\begin{proof}
	Let $\epsilon$ be sufficiently small such that the sets $\{\abs{\zeta-c_j}<\epsilon\},$ $j=1,\ldots,q,$
	$\{\abs{\zeta-z}<\epsilon\}$, $\{\abs{\zeta-a_j}<\epsilon\}$, $j=1,\ldots,n$ are disjoint
	and if any of the centers belongs to $\Omega$ then the disc with that center of radius $\epsilon$
	also belongs to $\Omega.$ Denote by $\Omega_\epsilon$ be the domain
	obtained by removing from $\Omega$ the union of the above discs. Denote by $\gamma_\epsilon$ the part
	of $\gamma$ which does not belong to any of the discs, and denote by $C(\epsilon,a)=\{z:\abs{z-a}=\epsilon\}\cap \Omega.$
	By Eqn.(\ref{greenekv2}) with $u=u(\zeta)$ and $v=g(\zeta,z)$ (in particular $\Delta v=0$) 
	we have for the domain $\Omega_\epsilon$ (where if some of the points $c_j,a_k$ appear more than once we only write the associated integral once)
	\begin{multline}\label{goldnevanekv14}
	\int_{\Omega_\epsilon} v\Delta d d\sigma=
	\\
	\left(\int_{\Omega_\epsilon} +\int_{C(\epsilon,z)}+\sum_{j=1}^q \int_{C(\epsilon,c_j)} +\sum_{j=1}^n \int_{C(\epsilon,a_j)}\right)
	\left(u\frac{\partial v}{\partial n}-v\frac{\partial u}{\partial n}\right)ds
	\end{multline}
	Since $v=0$ for $\zeta\in \gamma$ we have
	\begin{equation}
	\lim_{\epsilon\to 0}\int_{\gamma_\epsilon} \left(u\frac{\partial v}{\partial n}-v\frac{\partial u}{\partial n}\right)ds
	=\int_{\gamma_\epsilon} u(\zeta)\frac{\partial g}{\partial n}ds
	\end{equation}
	By the mean value theorem we have for each $a\in \overline{\Omega}$ with $(2\pi\epsilon)$ being realized as the length of $C(\epsilon,a)$
	\begin{equation}
	\int_{C(\epsilon,a)} \left(u\frac{\partial v}{\partial n}-v\frac{\partial u}{\partial n}\right)ds=(2\pi\epsilon)
	\left(u\frac{\partial v}{\partial n}-v\frac{\partial u}{\partial n}\right)|_{\zeta^*}
	\end{equation}
	where $\zeta^*$ is a point in $C(\epsilon,a).$
	We have
	$v=\ln \frac{1}{\epsilon} +O(1),$
	$\frac{\partial v}{\partial n}=\frac{\partial v}{\partial \epsilon}=-\frac{1}{\epsilon} +O(1),$
	$u=u(z)+o(1)$ and $\frac{\partial u}{\partial n}=O(1)$.
	This implies
	\begin{equation}
	\lim_{\epsilon\to 0}\int_{C(\epsilon,z)} \left(u\frac{\partial v}{\partial n}-v\frac{\partial u}{\partial n}\right)ds
	=-2\pi u(z)
	\end{equation}
	For $c_j\in \Omega$, $j=1,\ldots,q,$
	\begin{equation}
	\lim_{\epsilon\to 0}\int_{C(\epsilon,c_j)} \left(u\frac{\partial v}{\partial n}-v\frac{\partial u}{\partial n}\right)ds
	=-2\pi d_j g(c_j,z)
	\end{equation}
	If $a$ is any one of the points $c_j$ belonging to $\gamma$ or one of the points
	$a_j$, $j=1,\ldots,n$ then the length of $C(\epsilon,a)$ is $O(\epsilon)$ thus on 
	$C(\epsilon,a)$ we have by Eqn.(\ref{goldberganvsen})
	as $\epsilon\to 0$ that
	$u=O(\abs{\ln \epsilon}),$
	$\frac{\partial u}{\partial n}=O\left(\frac{1}{\epsilon}\right)$, $v=o(1)$ and
	$\frac{\partial v}{\partial n}=O\left(\frac{1}{\sqrt{\epsilon}}\right)$.
	\begin{equation}
	\lim_{\epsilon\to 0}\int_{C(\epsilon,a)} \left(u\frac{\partial v}{\partial n}-v\frac{\partial u}{\partial n}\right)ds=0
	\end{equation}
	Letting $\epsilon\to 0$ in Eqn.(\ref{goldnevanekv14}) yields the wanted result.
	\end{proof}
	
	\begin{theorem}\label{goldbergnevanthm21}
		Let $\Omega\subset\C$ be a simply connected domain with piecewise analytic boundary $\gamma$ and $f\not\equiv 0$ a meromorphic function
		on $\overline{\Omega}.$ Let $a_1,\ldots,a_m$ denote the zeros 
		of $f$ in $\Omega$
		and let $b_1,\ldots,b_n$ denote the poles of $f$ in $\Omega$.
		Then
		\begin{equation}\label{goldnevanekv21}
		\ln \abs{f(z)}=\frac{1}{2\pi} \int_\gamma \ln\abs{f(\zeta)}\frac{\partial g}{\partial n}ds-
		\sum_{j=1}^m g(a_j,z) +\sum_{j=1}^n g(b_j,z)
		\end{equation}
	\end{theorem}
	\begin{proof}
		By Theorem \ref{goldbergnevanthm1} with $u(z):=\ln\abs{f(z)}$ we have
		that in Eqn.(\ref{goldnevanekv13}) the left hand side vanishes, because if 
		$z$ is neither a zero or a pole of $f(z)$ we have $\Delta u(z)=0.$ If $c_j$ is a zero (pole) of $f$ of order $p_j$ then
		Eqn.(\ref{goldnevanekv12}) with $d_j=p_j$ ($d_j=-p_j$) holds true in a neighborhood of $c_j.$
		Thus the right hand side of Eqn.(\ref{goldnevanekv13}) is equal to the right hand side of
		Eqn.(\ref{goldnevanekv21}).
		This completes the proof.
		\end{proof}
		Note that applying Eqn.(\ref{goldbergnevekv4}), Eqn.(\ref{goldbergnevekv8}) and Eqn.(\ref{goldbergnevekv9})
		to the case of $\{\abs{z}<R\}$ in Theorem \ref{goldbergnevanthm21} yields the following (this is sometimes called the Poisson-Jensen formula\index{Poisson-Jensen formula}).
		\begin{corollary}
		Let $f(z)\not\equiv 0$ be a meromorphic function on $\{\abs{z}\leq R\}$. Then
		\begin{multline}\label{goldnevanekv22}
		\ln\abs{f(z)}=\frac{1}{2\pi} \int_0^{2\pi} \ln \abs{f(R\exp(i\theta)}\re \frac{R\exp(i\theta)+z}{R\exp(i\theta)-z}d\theta\\
		-\sum_{j=1}^m \ln \abs{\frac{R^2-\bar{a}_j z}{R(z-a_j)}} +\sum_{j=1}^n \ln \abs{\frac{R^2-b_j z}{R(z-b_j)}}
		\end{multline}
		where the $a_j$, $j=1,\ldots,m$ are the zeros in inside the disc and the $b_j$, $j=1,\ldots,n$, are the poles inside the disc.
		\end{corollary}
		
		\begin{corollary}
			Let $f(z)\not\equiv 0$ be a meromorphic function on $\{\abs{z}\leq R\}$. Then for a real constant $c$ we have
			\begin{multline}\label{goldnevanekv24}
			\ln f(z)=\frac{1}{2\pi} \int_0^{2\pi} \ln \abs{f(R\exp(i\theta)}\re \frac{R\exp(i\theta)+z}{R\exp(i\theta)-z}d\theta\\
			-\sum_{j=1}^m \ln \abs{\frac{R^2-\bar{a}_j z}{R(z-a_j)}} +\sum_{j=1}^n \ln \abs{\frac{R^2-b_j z}{R(z-b_j)}} +ic
			\end{multline}
			where the $a_j$, $j=1,\ldots,m$ are the zeros in inside the disc and the $b_j$, $j=1,\ldots,n$, are the poles inside the disc.
		\end{corollary}
		\begin{proof}
			Note that both sides of Eqn.(\ref{goldnevanekv24}) are holomorphic with respect to $z.$
			By Theorem \ref{goldbergnevanthm21} the real parts of both sides of Eqn.(\ref{goldnevanekv24}) coincide.
			By the Cauchy-Riemann equations they are equal up to a pure imaginary additive constant.
			This completes the proof.
		\end{proof}

		\begin{theorem}\label{goldbergnevanthm21}
			Let $f(z)\not\equiv 0$ be a meromorphic function on $\{\abs{z}\leq R\}$ with Laurent expansion at $z=0$ given by
			\begin{equation}\label{goldnevanekv21}
			f(z)=c_\lambda z^\lambda +c_{\lambda +1} z^{\lambda +1}+\cdots,\quad c_\lambda\neq 0
			\end{equation}
			Then we have ({\em Jensen formula}\index{Jensen formula}) 
			\begin{multline}
			\ln \abs{c_\lambda}=\frac{1}{2\pi} \int_0^{2\pi} \ln \abs{f(R\exp(i\theta)}d\theta\\
			-\sum_{0<\abs{a_j}<R} \ln \frac{R}{\abs{a_j}} +\sum_{\abs{b_j}<R} \ln \frac{R}{\abs{b_j}} -\lambda \ln R
					\end{multline}
					where the $a_j$ are the zeros and the $b_j$ are the poles of $f$ in
					$\{0\abs{z}<R\}$.
				\end{theorem}
				\begin{proof}
					Letting $z\to 0$ in the Poisson-Jensen formula we have for the sums in the right hand side
					\begin{multline}
					-\sum_{\abs{a_j}<R} \ln \abs{\frac{R^2-\bar{a}_j z}{R(z-a_j)}} +\sum_{\abs{b_j}<R} \ln \abs{\frac{R^2-b_j z}{R(z-b_j)}}
					=\\
					-\sum_{0<\abs{a_j}<R} \ln \abs{\frac{R^2-\bar{a}_j z}{R(z-a_j)}} +\sum_{0<\abs{b_j}<R} \ln \abs{\frac{R^2-b_j z}{R(z-b_j)}} +\lambda\ln\frac{\abs{z}}{R}\\
					=-\sum_{0<\abs{a_j}<R} \ln\frac{R}{\abs{a_j}} +\sum_{0<\abs{b_j}<R} \ln\frac{R}{\abs{b_j}} +o(1)+\lambda\ln \frac{\abs{z}}{R}
					\end{multline}
					From the left hand side in the Poisson-Jensen formula we get $\ln \abs{f(z)}=\lambda \ln \abs{z} +\abs{c_\lambda} +o(1).$
					This completes the proof.
				\end{proof}
\begin{definition}
	Let $f(z)$ be a meromorphic function on $\{\abs{z}<r\}$
	for some $r>0$ and let $n(r,f)$ denote the number of poles 
	of $f$ counting multiplicity, in $\{\abs{z}<r\}.$
	The Nevanlinna function\index{Nevanlinna characteristic} $N(r,0,f)$ of $f$ at $0$
	is defined as
	\begin{equation}
	N(r,f)=\int_0^r (n(t,f)-n(0,f))\frac{dt}{t} +n(0,f)\ln r
	\end{equation}
	Setting $m(r,0,f):=\frac{1}{2\pi}\int_0^{2\pi} \ln^+\abs{f(r\exp(i\theta))}d\theta$
	the Nevanlinna characteristic at $z=0$
	is defined as 
	\begin{equation}\label{charactnevek}
	T(r,f)=m(r,0,f)+N(r,0,f)
	\end{equation}
	For a point $a\neq 0$ we let $n(r,a)$ be the number of solutions to $f(z)=a$ in the disc $\{\abs{z}<r\}$ counting multiplicity. 
	If $f(0)\neq a$ and $f(0)<\infty,$ we set $N(r,a)=\int_0^r \frac{n(t,a)}{t}dt,$
	otherwise we use
	$N(r,a)=\int_0^r (n(t,a)-n(0,a))t^{-1} dt +n(0,a)\ln a.$
\end{definition} 
Let $f(z)$ be a meromorphic function on $\{\abs{z}<r\}$ for some
$r>0.$ By the argument principle and the Cauchy-Riemann equations we have
for $a\neq \infty$
\begin{multline}
n(r,a)-n(r,\infty)=\frac{1}{2\pi i}\int_{\abs{z}=r} \frac{f'(z)}{f-a}dz=\\
\frac{r}{2\pi}\frac{d}{dr}\int_{-\pi}^\pi \ln \abs{f(r\exp(i\theta))-a}d\theta
\end{multline}
Dividing by $r$ and integrating with respect to $r$
we have if $f(0)\neq a$ and $\abs{f(0)}<\infty,$
\begin{equation}\label{jensekk}
\frac{1}{2\pi}\int_{-\pi}^\pi \ln \abs{f(r\exp(i\theta))-a}d\theta =\ln \abs{f(0)-a}+N(r,a,f)-N(r,\infty,f)
\end{equation}
If $f(0)=a$ and $f(z)-a=cz^m+\cdots,$ 
$c\neq 0$ then $\ln\abs{f(0)-a}$ has to be replaced by $\ln \abs{c}$
and we must use $N(r,a,f)=\int_0^r (n(t,a)-n(0,a))t^{-1} dt +n(0,a)\ln a.$
If we set $m(r,\infty,f)=\frac{1}{2\pi}\int_{-\pi}^\pi \ln^+ \abs{f(r\exp(i\theta))}d\theta$
	and \begin{equation}m(r,a,f)=m(r,\infty,(f-a)^{-1})=\frac{1}{2\pi}\int_{-\pi}^\pi \ln^+ \frac{1}{\abs{f(r\exp(i\theta))-a}}d\theta\end{equation}
		then the relation
		\begin{equation}
		\log^+\abs{x\pm y}\leq \log^+\abs{x} +\log^+ \abs{y} +\log 2
		\end{equation}
		implies that
		\begin{equation}
		m(r,\infty,f)=m(r,\infty,f-a)+O(1)
		\end{equation}
		which implies that Eqn.(\ref{jensekk}) can be written as
		\begin{equation}
		m(r,\infty,f)+N(r,\infty,f)=m(r,a,f)+N(r,a,f)+O(1), \quad r\to \infty
		\end{equation}
		confer Eqn.(\ref{charactnevek}), i.e.\ the definition of
		the Nevanlinna characteristic $T(r,f)$ is in some sense independent
		of $a.$ The term $N(r,a,f)$ counts the $a$-points 
		and the term $m(r,a,f)$ is a measure of proximity of $f(z)$ to $a$
		on the circle $\abs{z}=r.$

\section{Polyanalytic version of the fundamental theorem of algebra}
According to the fundamental theorem of algebra, a nonconstant complex polynomial must have a zero.
This is not the case for polyanalytic functions with polynomial polyanalytic components, i.e.\
functions of the form $f(z)=\sum_{j=0}^{q-1} p_j(z)\bar{z}^j$, for polynomials $p_j(z)$.

Take e.g.\ $f(z)=1+(z+\bar{z})=1+2i\im z.$
\begin{theorem}[Hurwitz]\label{hurwitsthmzero}
	Let $\{f_j\}_{j\in \N}$ be a sequence of analytic functions on a domain $\Omega\subset\C$
	such that the sequence converges uniformly on compacts of $\Omega$ to a function $f(z)$ such that
	$f(z)$ has a zero of order $n$ at $p_0\in\Omega.$  Then there exists $\epsilon>0$, $N\in \N$,
	such that for $j\geq N$ each $f_j(z)$ has exactly $n$ zeros in
	$\{ \abs{z-p_0}<\epsilon\}$, counting multiplicity,and these zeros converge to $p_0$ as $j\to \infty.$
\end{theorem}
\begin{proof}
	Let $\epsilon>0$ be sufficiently small such that $\{\abs{z-p_0}\leq r\}$ is contained in $\Omega$ so that $f(z)\neq 0$ for $0<\abs{z-p_0}\leq r.$
	Let $\delta>0$ be such that $\abs{f(z)}\geq \delta$ on $\{\abs{z-p_0}=r\}.$
	Since $f_j\to f$ uniformly on $\{ \abs{z-p_0}\leq \epsilon\}$ we have
	for sufficiently large $N$ that for $j\geq N$,
	$\abs{f_j(z)}>\epsilon/2$ for $\abs{z-p_0}= \epsilon$. This implies that $f_j'(z)/f_j(z)$
	is well defined for all $z\in \{ \abs{z-p_0}= \epsilon\}$. By Morera's theorem
	$f_j'(z)/f_j(z)$ converges uniformly to $f'(z)/f(z)$ for $\abs{z-p_0}= \epsilon$.
	Hence we have 
	\begin{equation}\label{hurwzeros}
	\frac{1}{2\pi i}\int_{\abs{z-p_0}=\epsilon} \frac{f'_j(z)}{f_j(z)}dz\to \frac{1}{2\pi i}\int_{\abs{z-p_0}=\epsilon} \frac{f'(z)}{f(z)}dz
	\end{equation}
	By the argument principle we have that the number $n_j$($n$ respectively) of the zeros of $f_j$ ($f$ respectively) inside $\{ \abs{z-p_0}<\epsilon\}$
	is given by
	\begin{equation}
	n_j=\frac{1}{2\pi i}\int_{\abs{z-p_0}=\epsilon} \frac{f'_j(z)}{f_j(z)}dz,\quad n=\frac{1}{2\pi i}\int_{\abs{z-p_0}=\epsilon} \frac{f'(z)}{f(z)}dz
	\end{equation}
	By Eqn.(\ref{hurwzeros}) we have $n_j\to n$ as $j\to \infty$, hence for sufficiently large $N$ we have for $j\geq N$ 
	that $f_j(z)$ has $n$ zeros on $\{ \abs{z-p_0}< \epsilon\}.$ Since $\epsilon>0$ can be chosen arbitrarily small
	this implies that the zeros of $f_j(z)$ accumulate to $p_0.$ This completes the proof.
\end{proof}

Balk \cite{balkfundament} proved the following theorem on existence of zeros for a subclass of polyanalytic polynomials.
\begin{theorem}\label{fundamentalthmbalk}
	Suppose $P(z,\bar{z})$ is a polynomial
	such that the exact degree of $P(v,w)$ as a joint polynomial in complex variables $v,w$ is $n\in \N$
	such that $n> 2\max\{n_1,n_2\}$ where $n_1$ is the exact degree of $P$ with respect to $v$ and 
	$n_2$ is the exact degree of $P$ with respect to $w$. Then $P(z,\bar{z})$ has at least one complex root.
\end{theorem}
\begin{proof}
	We may write $P(z,\bar{z})=\sum_{j=0}^{n_2} \bar{z}^j P_j(z),$ where $P_j(z)=\sum_{k=0}^{n-j} a_{j,k} z^k$ such that not all
	$a_{n-j,j}\in \C$ are zero. 
	Suppose (in order to reach a contradiction) that $P(z,\bar{z})$
	has no zeros. Then 
	the polynomial $\phi(z,t):=\sum_{j=0}^{n_2} t^{2j}z^{n_2 -j} P_j(z)$ has
	precisely $n_2$ zeros in $\{\abs{z}<t\}$ for any $t>0.$ For 
	$\abs{z}<1$ and $t\to \infty$ we have
	\begin{equation}
	\phi(z,t)=(a_{n,0}z^{n_2+n}+a_{n-1,1}z^{n+n_2-2}+\cdots a_{n-n_2,n_2} z^{n+n_2})+o(1)
	\end{equation}
	By Theorem \ref{hurwitsthmzero} $\phi(z,t)$ has at least $n-n_2$ zeros in $\{\abs{z}<1\}$ for large $t.$
	This implies that $n_2\geq n-n_2,$ $n\leq 2n_2.$ Now $\overline{P(z,\bar{z})}$ has degree $n_1$ with respect to
	$\bar{z}$ which implies $n\leq 2n_1$. This completes the proof.
\end{proof}

\section{Inner uniqueness}
Let us begin with the an observation on sets of uniqueness.
\begin{definition}
	Let $\Omega\subset\Cn$ be a domain and let $\alpha\in \Z_+^n.$
	We say that a set $E\subset \Omega$
	is a {\em set of uniqueness} for $\alpha$-analytic functions if
	any $\alpha$-analytic function $f$ on $\Omega$ satisfying $f|_E= 0$ vanishes identically on $\Omega$. 
\end{definition}

\begin{lemma}\label{203lemmabalk}
	Let $\Omega\subset\mathbb{C}$ be a domain, $0\in \omega$, let $q\in \mathbb{Z}_+$ and let $z=x+iy$ denote the complex coordinate in $\mathbb{C}$.
	Let $c_0,\ldots,c_{q-1}$ be distinct positive real numbers such that 
	$\{|z| =c_k\} \cap\Omega\neq \emptyset,$ $k=0,\ldots,q-1$.
	Let $f(z)=\sum_{j=0}^{q-1} a_j(z)|z|^{2j}$
	for holomorphic functions $a_j(z)$ on $\Omega,$ $j=0,\ldots,q-1$ (i.e.\ $f$ is, up to translation, an arbitrary 
	reduced $q$-analytic function
	on $\Omega$). Let \begin{equation}
	E:=\bigcup_{k=0}^{q-1} \Omega\cap \{ |z|=c_k\}
	\end{equation}
	If $f|_E=0$ then $f\equiv 0$ on $\Omega$ (in other words $E$ is a set of uniqueness
	for reduced $q$-analytic functions).
\end{lemma}
\begin{proof}
	In polar form we shall use the notation $z=|z|\exp(i\theta).$ 
	The condition $f|_E=0$
	can be written for $\tilde{a}_{j,k}(\theta):=a_j(c_k\exp(i\theta)),$ and all $\theta$ such that $|z|\exp(i\theta)\in \Omega$
	\begin{equation}\label{ekv1}
	0=\sum_{j=0}^{q-1} \tilde{a}_{j,k}(\theta)c_k^{2j}=:\phi_k(\theta),\quad k=0,\ldots,q-1
	\end{equation}
	where each $\tilde{a}_k(\theta)$ is the restriction to a nonempty intersection of a circle with $\Omega$, of a unique holomorphic function $a_j(z)$,
	so that each $\phi_k(\theta)$ is the restriction to a nonempty intersection of a circle with $\Omega$, of a unique holomorphic function 
	\begin{equation}\label{zsbalkekv1}
	\Phi_k(z):= \sum_{j=0}^{q-1} a_{j}(z)c_k^{2j}, \quad \Phi_k|_{ \Omega\cap \{ |z|=c_k\}}=\phi_k,
	\quad k=0,\ldots,q-1
	\end{equation}
	Since $\phi_k\equiv 0\Rightarrow \Phi_k\equiv 0$ we have by Eqn.(\ref{zsbalkekv1}) (since $\Omega$ is connected)
	the following system for each 
	\begin{equation}\label{zsbalkekv2}
	0=\sum_{j=0}^{q-1} a_{j}(z)c_k^{2j},\quad k=0,\ldots,q-1, \quad z\in \Omega
	\end{equation}
	This is a system which can be written in matrix form
	\begin{equation}
	[a_1(z),\ldots,a_{q-1}(z)]M=0
	\end{equation}
	where the coefficient matrix is clearly a Vandermonde matrix
	\begin{equation}
	M:=\begin{bmatrix}
	1 & 1 & \cdots & 1\\
	(c_0^2)^1 & (c_1^2)^1 & \cdots & (c_{q-1}^2)^1\\
	\vdots & \vdots & \vdots & \vdots\\
	(c_0^2)^{(q-1)} & (c_1^2)^{(q-1)} & \cdots & (c_{q-1}^2)^{(q-1)}\\
	\end{bmatrix}
	\end{equation}
	and $M$ has determinant given by
	\begin{equation}
	\mbox{det}(M)=\Pi_{i<j} ((c_i^2)-(c_j^2))
	\end{equation}
	Since the $c_k$, $k=0,\ldots,q-1$ are distinct the determinant is nonzero.
	Hence $M$ is invertible and thus  
	$a_j\equiv 0$, $j=0,\ldots,q-1.$ Hence $f\equiv 0.$
	This proves Lemma \ref{203lemmabalk}.
	\end{proof}

\begin{proposition}\label{uniqueness203}
	Let $\Omega\subset\mathbb{C}$ be a domain and let $z=x+iy$ denote the complex coordinate in $\mathbb{C}.$ 
	Let $p_0\in \Omega.$
	Let $c_0,\ldots,c_{q-1}$ be distinct nonzero real numbers.
	Set for $j=0,\ldots,q-1$
	\begin{equation}
	E^1_j:=\{ y=c_j\},\quad E^2_j:=\{ x=c_j\},\quad E^3_j:=\{ |z-p_0|^2=c_j^2\}
	\end{equation}
	\begin{equation}
	E^k:=\bigcup_{j=0}^{q-1}\Omega\cap E_j^k,\quad k=1,2,3
	\end{equation}
	If $E_j^k\cap \Omega \neq\emptyset,$ $j=0,\ldots,q-1,$ 
	then $E^k$ is a set of uniqueness for $q$-analytic functions, $k=1,2,3$.
\end{proposition}
\begin{proof}
	It suffices (by translation, since translation preserved $q$-analyticity) to prove the result for the case that $p_0=0\in \Omega$
	and verify uniqueness with respect to any $q$-analytic function of the form
	$f(z)=\sum_{j=0}^{q-1} a_j(z)\bar{z}^j$
	for holomorphic functions $a_j(z),$ $j=0,\ldots,q-1.$
	Note that $f|_{E^3}=0$ iff $(z^{q-1} f)|_{E^3} =0.$
	Since $z^{q-1}f$ is reduced we have by Lemma \ref{203lemmabalk}	that if
	$E_j^3\cap \Omega \neq\emptyset,$ $j=0,\ldots,q-1,$ then $E^3$ is 
	a set of uniqueness for $f$ which was (up to translation) an arbitrary $q$-analytic functions.
	For the case of $E^1$ we rewrite
	\begin{equation}
	f(z)=\sum_{j=0}^{q-1} a_j(z)\bar{z}^j=\sum_{j=0}^{q-1} a_j(z)(z-2iy)^j
	\end{equation}
	which can be expanded into
	\begin{equation}\label{zsbalkekv3}
	f(z)=\sum_{j=0}^{q-1} \sum_{k=0}^j y^k \binom{j}{k} (-2i)^k z^{j-k}a_j(z)=\sum_{j=0}^{q-1} b_j(z)y^j
	\end{equation}
	for holomorphic $b_j(z),$ $j=0,\ldots,q-1$.
	Similarly, for the case $E^2$ we rewrite
	\begin{equation}
	f(z)=\sum_{j=0}^{q-1} a_j(z)\bar{z}^j=\sum_{j=0}^{q-1} a_j(z)(-z+2x)^j=\sum_{j=0}^{q-1} d_j(z)x^j
	\end{equation}
	for holomorphic $d_j(z),$ $j=0,\ldots,q-1$.
	In the case of $E^1$, the proof of Lemma \ref{203lemmabalk} can be repeated
	after in appropriate places replacing $a_j(z),|z|^2,\theta$ with $b_j(z),y,x$, and we then obtain a system whose coefficient matrix is given by
	\begin{equation}
	M:=\begin{bmatrix}
	1 & 1 & \cdots & 1\\
	c_0 & c_1 & \cdots & c_{q-1}\\
	\vdots & \vdots & \vdots & \vdots\\
	c_0^{(q-1)} & c_1^{(q-1)} & \cdots & c_{q-1}^{(q-1)}\\
	\end{bmatrix}
	\end{equation}
	and has determinant given by
	\begin{equation}
	\mbox{det}(M)=\Pi_{i<j} (c_i-c_j)
	\end{equation} 
	Since the $c_k$, $k=0,\ldots,q-1$ are distinct and nonzero the determinant is nonzero.	
	This implies (since $\Omega$ is connected) that $b_j\equiv 0$, $j=0,\ldots,q-1,$ which in turn implies by Eqn.(\ref{zsbalkekv3}) that $f\equiv 0.$
	Analogously, in the case of $E^2$, the proof of Lemma \ref{203lemmabalk} can be repeated
	after in appropriate places replacing $a_j(z),|z|^2,\theta$ with $d_j(z),x,y$, and repeating the analogous arguments given for the case of $E^1$,
	we get $d_j\equiv 0$, $j=0,\ldots,q-1,$ which in turn implies $f\equiv 0.$
	We conclude that $E^1$ and $E^2$ are both set of uniqueness for $q$-analytic functions.
	This completes the proof.
	\end{proof}
Let $E\subset \C$ be an arbitrary set of points and $a\in \C$ any fixed point (not necessarily in $E$).
Let $\ell$ be a ray emanating from $a$ and let the equation of this ray 
be $\mbox{arg}(z-a)=\alpha.$
The set $E$ is said to be {\em condensed} to the point $a$ along the ray $\ell$
if there is a sequence of points $z_j$ in $E$ such that for
$j\to \infty$ $\lim_j z_j = a$ and $\lim_j \mbox{arg}(z_j-a) =\alpha.$
we see that $E$ is condensed to the point $a$ along $\ell$
if every angle with a vertex at the point $a$ containing the ray $\ell$
also contains a part of the set $E$ for which $a$ is a limit point. 
A point a is called a condensation point of order $k$ for the set $E$, if $E$ condenses to the point $a$  
along $k$ distinct lines. 
We can phrase this as follows.
	\begin{definition}\label{deff0}(See Balk~\cite[p.~4]{balk1965}.)
	Let $U\subset\C$ be a domain and let $p\in E\subset U$. 
	We say that the line $\ell :=\{z\in \C\colon z=p+te^{i\theta},\abs{t}<\infty, t\in \R\}$, $p$ and $\theta$ constants, 
	is a \emph{limiting direction of the set $E$ at $p$} if $E$ contains a sequence of points $z_j =p+t_j e^{i\theta_j}$,
	$t_j\to 0,\theta_j\to \theta,t_j\neq 0$. The point $p$ is called a \emph{condensation point of order $k$ of $E$} if there are $k$ different lines through $p$ which are limiting directions of $E$. 
	If in this case there are not $k+1$ distinct straight lines along which the set $E$ is condensed to $a$, 
	then we say that the condensation point $a$ has {\em exact}
	\index{Exact order of condensation} order $k$.
\end{definition}
\begin{example}
	The origin is a condensation point of order $1$ for 
	$y=\sin(x),$ whereas it is a condensation point of order $2$ for the lemniscate
	$(x^2+y^2)=2(x^2-y^2)$ and 
	a condensation point of order $3$ for $y=x\sin(1/x),$ $x\neq 0.$
\end{example}
\begin{proposition}\label{balkuniquelemmat}
	If $f$ is a polyanalytic function of order $q$ in some disk $\{\abs{z -a}<r\}$ for a fixed $r>0,$
	and vanishes on
	a set $E$ having a condensation point of order $q$ at $a$, then $f\equiv 0.$
\end{proposition}
\begin{proof}
	W.l.o.g.\ we can assume $a=0.$ We know that on the disc $\{\abs{z}<r\}$, $f$ has a representation of the form
	$f(z)=\sum_{j=0}^{q-1} a_j(z)\bar{z}^j$ for holomorphic $a_j(z),$ $j=0,\ldots,q-1,$
	which in polar coordinates $z=r\exp(i\theta)$ and after Taylor expansion near $0$ yields
	\begin{multline}
	f(z)=\sum_{j=0}^{q-1} \bar{z}^j \sum_{k=0}^\infty b_{j,k}z^k =\sum_{j=0}^{q-1}\sum_{k=0}^\infty
	r^{j+k} \exp(i(k-j)\theta)=\\
	\sum_{m=0}^\infty
	r^{m}  \sum_{j=0}^j b_{j,m} \exp(i(m-2j)\theta)
	\end{multline}
	for complex constants $b_{j,k},$ with
	\begin{equation}\label{uniqsju}
	b_{j,k}=0,\quad k<0\mbox{ or }j\geq q
	\end{equation}
	Setting 
	\begin{equation}
	B_m(\nu):=\sum_{j=0}^m b_{j,m-j}\nu^{m-j}, \quad \nu:=\exp(2i\theta)
	\end{equation}
	gives
	\begin{equation}\label{uniqnio}
	f(z)=\sum_{m=0}^\infty r^m B_m(\nu)\exp(-im\theta) 
	\end{equation}
	Assume (in order to reach a contradiction) that $f\not\equiv 0.$ Then there exists $m_0$
	such that $B_{m}(\nu)\equiv 0$ for $m<m_0$ and $B_{m_0}(\nu)\equiv 0.$
	By Eqn.(\ref{uniqnio}) we obtain
	\begin{equation}\label{uniqelva}
	r^{-m_0}f(z)=\exp(-im_0\theta) B_{m_0}(\nu)+r\sum_{m=m_0+1}^\infty r^{m-p_0-1} B_m(\nu)\exp(-im\theta) 
	\end{equation}
	By the conditions of the proposition there exists $q$ distinct rays (emanating from $z=0$)
	say given by $\mbox{arg} z=\theta_1,\ldots,\mbox{arg}z=\theta_q$ such that $E$ condenses along each of the rays to $z=0$.
	In particular, there exists $\{z_j\}_{j\in \N}$, $z_j\neq 0$, $\lim_{j\to \infty} z_j =0,$
	$\lim_{j\to \infty}\mbox{Arg}z_j=\theta_1$.  Clearly, $f(z_j)=0$. Setting in Eqn.(\ref{uniqelva}) $z=z_j$ and letting
	$j\to \infty$ yields
	\begin{equation}\label{uniqtolv}
	B_{m_0}(\exp(-2i\theta_1)=0 
	\end{equation}
	Hence the polynomial $B_{m_0}(\nu)$ has a root $\exp(2i\theta_1)$. Repeating this for $\theta_j$
	$j=2,\ldots,q$
	we obtain that $B_{m_0}(\nu)$ has the $q$ distinct roots
	$\exp(2i\theta_j),$ $j=1,\ldots,q.$ By Eqn.(\ref{uniqsju}) we can write
	\begin{equation}\label{uniqtretton}
	B_{m_0}(\nu)=\nu^{m_0-q+1}A_{m_0}(\nu),\quad A_{m_0}(\nu)=\sum_{j=0}^{q-1} b_{j,q-j}\nu^{q-1-j} 
	\end{equation}
	Hence, the polynomial $A_{m_0}(\nu)\not\equiv 0$ and has $q$ distinct roots 
	$\exp(2i\theta_j),$ $j=1,\ldots,q.$ But this is impossible, 
	since deg$A_{m_0}(\nu)\leq q-1$. Thus by contradiction the assumption $f\not \equiv 0$ must be false.
	This completes the proof.
\end{proof}
Obviously, this includes the well-known result that open sets are sets of uniqueness for
polyanalytic functions, namely
if $f(z)$ is a polyanalytic of order $q$ in some disk $D$ with center $a$, and 
vanishes on some arbitrarily small disc $E$ such that $a\in E\cup\partial E$
then clearly $E$ condenses to $a$ along an infinite set of rays emanating from $a$
thus $f\equiv 0.$ 
As a direct consequence of Proposition \ref{balkuniquelemmat}, Balk \cite{balk1965} deduces the following theorem.
\begin{theorem}\label{balkthmforstauniq}
	Let $\Omega\subset\C$ be a domain and $q\in \Z_+.$
	Let $E\subset\Omega$ be a set which has in $\Omega$ a condensation
	point of order $q$. Then $E$ is a set of uniqueness for $q$-analytic functions on 
	$\Omega$
	(in the sense that if $f,g$ are $q$-analytic on $\Omega$
	such that $f|_E=g|_E$ then $f\equiv g$ on $\Omega$).
\end{theorem}
\begin{proof}
Let $a\in E$ be a condensation point of order $q$.
It suffices to consider, for each given $b\in \Omega$, 
chains of finitely many discs $D_1,\ldots,D_m$ relatively compact
in $\Omega$ such that $D_1$ has center $a$ and $D_m$ contains $b$ and such that
the center of $D_j$ belongs to the closure, $\overline{D}_{j-1}$, of $D_{j-1}.$
Let $\phi=f-g$. By Proposition \ref{balkuniquelemmat} we have $\phi\equiv 0$ on $D_1$. This implies that
$\phi=0$ near the center of $D_2$, thus $\phi\equiv 0$ on $D_2.$
By repeated application of Lemma \ref{balkuniquelemmat} we obtain $\phi\equiv 0$ on $D_m$,
thus on an open neighborhood of the arbitrarily chosen $b\in \omega.$
This completes the proof.	
\end{proof}
Balk \cite{balk1965} points out that the proofs of Proposition \ref{balkuniquelemmat} and Theorem \ref{balkthmforstauniq}
carry over with minimal adaptation to the case of countably analytic functions, and as we have pointed out in the
first chapter, the order of polyanalyticity is constant for a function that is polyanalytic of finite order on a domain.
\begin{theorem}
	Let $\Omega\subset\C$ be a simply connected domain, let $f$ and $g$ be analytic functions
	on $\Omega$ and let $P(x,y)$ and $Q(x,y)$ be real polynomials of degree at most $q$ satisfying that there exists a point $x_1+iy_1\in \C$
	such that  
	\begin{equation}\label{balkuniqintesant}
	P(x_1,y_1)=Q(x_1,y_1)\neq 0
	\end{equation}
	Suppose $E\subset \Omega$ has a condensation point of order $q+1.$
	If
	\begin{equation}\label{balkuniqintesant}
	P(x,y)f(z)=Q(x,y)g(z) \mbox{ on } E
	\end{equation}
	Then $f\equiv g$ on $\Omega.$
\end{theorem}
\begin{proof}
	If both $f\equiv 0$ and $g\equiv 0$ we are done so assume w.l.o.g.\ $g\not\equiv 0.$
	Rewriting $2x=z+\bar{z}$ and $2iy=\bar{z}-z$ we see that any polynomial in $x,y$ of order $q$
	is an entire $q+1$-analytic function, thus $P(x,y)f(z)$ and $Q(x,y)g(z)$ are both $q+1$-analytic.
	By assumption there is a point $(x_1,y_1)\in \Omega$ satisfying Eqn.(\ref{balkuniqintesant})
	Then there exists a disc $D_1\subset\Omega$, on which $g(z)$ has no zeros.
	Since by assumption the polyanalytic polynomial $P(x,y)\not\equiv 0$
	there exists a disc $D_2\subset D_1$ on which $P\neq 0,$ thus we have
	\begin{equation}
	\frac{f(z)}{g(z)}=-\frac{Q(x,y)}{P(x,y)},\quad \quad z=x+iy\in D_2
	\end{equation}
	where the left hand side is holomorphic whereas the right hand side is real-valued. Hence
	$f(z)/g(z)=c,$ for some real constant $c$, for all $z\in D_2,$ and thus $Q(x,y)\equiv cP(x,y)$ on $D_2.$
	But $Q$ and $P$ are the two polyanalytic functions that coincide on the open set $D_2$
	so (because $\Omega$ is connected) they coincide on $\Omega.$ 
	Since $P,Q$ are polynomials with respect to their arguments we must have $c=1,$  and in particular
	\begin{equation}
	c=\frac{Q(x_1,y_1)}{P(x,y)}=1,\quad \quad z=x+iy\in D_2
	\end{equation}
	Hence $f|_{D_2}\equiv g|_{D_2}$ and by connectedness $f\equiv g$ on $\Omega.$
	This completes the proof.
	\end{proof}
	
	\begin{theorem}
		If $f$ and $g$ are two polyanalytic functions that are equal on a subset $E$ of an open analytic 
		arc $\gamma$, having a point of accumulation, $p_0$, on $\gamma$, then they are equal on the whole arc. 
	\end{theorem}
	\begin{proof}
		By the definition it is possible to express $\gamma$ as $z=\lambda(t)=\lambda_1(t)+t\lambda_2(t),$
		$t\in (a,b)$ for a real interval $(a,b)$ and real-analytic $\lambda_1,\lambda_2,$
		$\lambda'(t)\neq 0$ on $(a,b).$
		Let $p_0=\lambda(t_0),$ for some $t_0\in (a,b)$ and choose $(\alpha,\beta)$ such that
		$a<\alpha<t_0<\beta<b.$ Then $\lambda([\alpha,\beta])$ is a closed arc $\gamma'\subset\gamma.$
		Now replacing the argument, $t,$ of $\lambda$ by a complex variabel, $\tau$ (i.e.\ extending the domain of
		$\lambda$) on an open neighborhood $\Delta$ of $[\alpha,\beta]$ we can find a choice of $\Delta$
		such that the extensions $\lambda_1(\tau),\lambda_2(\tau)$ and $\lambda(\tau)$ are holomorphic on $\Delta.$
		Let $t\in [\alpha,\beta]$ be a fixed point and $p=\lambda(t)$ the associated point on $\gamma.$
		In some neighborhood $D_1(p)$ of $p$ the functions $f(z)$ and $g(z)$ 
		can be represented in the form of uniformly convergent series
		\begin{equation}
		f(z)=\sum_{j=0}^\infty (\bar{z}-\bar{p})^j \phi_j(p;z),\quad g(z)=\sum_{j=0}^\infty (\bar{z}-\bar{p})^j \psi_j(p;z)
		\end{equation}
		for holomorphic $\phi_j\psi_j,$ $j\in \N.$ Since $\lambda'(t)\neq 0$ there exists a neighborhood $\Delta(t)$ of $t$ sufficiently small so that
		$z=\lambda(\tau)$ maps $\Delta(t)$ one-to-one onto some subdomain $\delta(p)$ of $D_1(p).$
		Now the infinite family of choices of $\Delta(t)$ covering $[\alpha,\beta]$, has (by compactness of the interval)
		a finite subfamily that covers $[\alpha,\beta]$ and we can always assume the disc $\Delta_0:=\Delta(t_0)$ is included.
		Denote the subfamily (ordered by increasing abscissa of their centers) $\Delta_j:=\Delta(t_j),$
		$j=-m,-m+1,\ldots,-1,0,1,\ldots,q,$ and denote the corresponding domains $\delta_j:=\delta(p_j),$ $p_j:=\lambda(t_j).$
		The art of $(a,b)$ belonging to $\Delta_j$ is an interval, which we denote $i_j$, so that
		$\{i_j\}_{j}$ cover $[\alpha,\beta].$
		The map $z=\lambda(\tau)$ maps $i_j$ to some arc $\iota_j$ belonging to $\gamma$, and they cover $\gamma'.$
		Further note that for each $j,$ $-m,\ldots,q-1$ we have that the intervals $i_j$ ans $i_{j+1}$ 
		have a common interval and $\iota_j$ and $\iota_{j+1}$ have a common arc.
		Now on $\delta_0$ the functions $f$ and $g$ can be represented as
		\begin{equation}
		f(z)=\sum_{j=0}^\infty (\bar{z}-\bar{p_0})^j \phi_j(p_0;z),\quad g(z)=\sum_{j=0}^\infty (\bar{z}-\bar{p_0})^j \psi_j(p_0;z)
		\end{equation}
		in particular they are, on $\delta_0,$ real-analytic with respect to $(x,y),$ say $f=\phi(x,y),$ $g=\psi(x,y).$
		For $x=\lambda_1(\tau),$ $y=\lambda_2(\tau)$ consider on $\Delta_0$ the functions
		$\Phi(\tau):=\phi(\lambda_1(\tau),\lambda_2(\tau)),$
		$\Psi(\tau):=\phi(\lambda_1(t),\lambda_2(t)).$ If $t\in [\alpha,\beta]$ then 
		\begin{equation}\label{balkuniq25}
		\Phi(\tau)=f(\lambda(\tau)),\quad \Psi(\tau)=g(\lambda(\tau))
		\end{equation}
		(however for complex $\tau$, $\Phi(\tau)$ may differ from $f(\lambda(\tau))$,
		and likewise for $\Psi(\tau)$ may differ from $g(\lambda(\tau))$).  
		By assumption $E$ contains a sequence of points $\{z_j\}_{j\in \Z_+},$ $z_j\to p_0$ as $j\to \infty,$
		such that staring for a sufficiently large $\nu\in \Z_+$, say $\nu\geq \nu_0$, each $z_\nu$, belongs to $\delta_0.$
		Since the map $z=\lambda(\tau)$
		is one-to-one on $\Delta_0$ there exists a point $\tau_\nu$ such that $\lambda(\tau_\nu)=z_\nu$ and we have
		$\tau_\nu\to t_0$ as $\nu\to \infty$ (and $z_\nu\to p_0$) and furthermore $\Phi(\tau_\nu)=\Psi(\tau_\nu),$ $\nu\geq \nu_0.$
		By the uniqueness property of holomorphic functions $\Phi(\tau)\equiv\Psi(\tau)$ on $\Delta_0$ in particular on $i_0.$
		By Eqn.(\ref{balkuniq25}) we have
		\begin{equation}
		f(\lambda(\tau))\equiv g(\lambda(\tau)),\mbox{ on the interval }i_0
		\end{equation} 
		which yields
		\begin{equation}
		f(z)\equiv g(z),\mbox{ on the arc }\gamma_0
		\end{equation}
		Repeating the arguments for the discs $\Delta_1,\Delta_{-1}$ yields 
		$f(z)\equiv g(z),$ on the arcs $\gamma_1,\gamma_{-1}$, and by iteration we obtain $f(z)\equiv g(z)$ onall the arcs
		$\gamma_j,$ $j=-m,-m+1,\ldots,-1,0,1,\ldots,q.$ Since we can always choose a segment $[\alpha,\beta]$ such that $\gamma'$ 
		contains a chosen point of $\gamma$ we obtain that $f(z)=g(z)$ on $\gamma.$
		This completes the proof.
	\end{proof}

Balk \cite{ca1} states on p.203 (although we have been able to locate the proof) that
Proposition \ref{uniqueness203} generalizes as follows.
\begin{definition}
A set $E_1$ is called {\em polyanalytically separable}\index{Polyanalytically separable} for a set $E_2$ if there exists a polyanalytic function which vanishes at all points on $E_1$ and not more than finitely many points of $E_2.$
\end{definition}

\begin{proposition}
	Let $\Omega\subset\C$ be a domain.
If $E=\bigcup_{j=1}^m E_j$ where the $E_j\subset\Omega$ are pairwise polyanalytically separable, and if each $E_j$ has a condensation point of order $n_j$
then $E$ is a set of uniqueness for $n$-analytic functions, $n:=\sum_{j=1}^m n_j.$
\end{proposition}

We shall however use a method of proof based upon
the notion of Schwarz functions in order to prove a similar result on sufficient conditions (not involving
condensation points of higher order) for a subset of a domain to b a set of uniqueness
for $q$-analytic functions.

\begin{definition}\label{oscardef}
	We denote by an {\em OSCAR}\index{OSCAR}, an {\em open} simple curve, $\gamma,$ that is analytic and regular. By this we means that $\gamma$ is the homeomorphic image in $\C$, of some real interval $(a,b)$, such that $\gamma$ is
	defined in $\C$ by the equation $z=w(t),$ where $w(t)$ is holomorphic on a neighborhood of $(a,b)$, with $w'(t)\neq 0,$
	for $t\in (a,b).$ By a {\em closed} regular analytic arc\index{Closed regular analytic arc} we mean a homeomorphic image $\Gamma$ of some segment $[\alpha,\beta]$, 
	defined in $\C$ by the equation $z=w(t),$ where $w(t)$ is holomorphic on a neighborhood of $[a,b]$. 
	Two OSCAR's 
	are called {\em non-contiguous}\index{Non-contiguous OSCAR:s} if their Schwarz functions
	are not analytic continuations of each other.
\end{definition}
\begin{remark}\label{noncontigrem}
	Let $q\in \Z_+,$ let $\Omega\subset \mathbb{C}$ be a domain, $\omega\subset \Omega$ a subdomain, let $S(z)$ be a holomorphic function on $\omega$
	and let $a_j(z),$ $j=0,\ldots,q-1$ be holomorphic functions on $\omega$.
	Suppose the function $\sum_{j=0}^{q-1} a_j(z) S^j(z)\equiv 0$ on $\omega.$ Then there exists a discrete set $D\subset\Omega$ (i.e.\ $D$ has no accumulation points in $\Omega$)
	such that $S(z)$ has analytic continuation to any simply connected domain $\omega'\subset \Omega\setminus D,$ such that $\omega'\cap \omega\neq \emptyset.$
	In particular, if $S(z)$ is the Schwarz function of an OSCAR, $E$, in $\Omega$ then, if there exists a polyanalytic function that vanishes on
	$E$, we must have that there exists a $q\in \Z_+$ together with holomorphic $a_j(z)$, $j=0,\ldots,q-1$, satisfying that
	$\sum_{j=0}^{q-1} a_j(z) S^j(z)=0$ for all $z\in E.$ Since $E$ has an accumulation point this implies 
	$\sum_{j=0}^{q-1} a_j(z) S^j(z)\equiv 0$ and thus the domain of the Schwarz function can be extended to any simply connected domain in the connected component
	in $\Omega\setminus D$, of a given point $p\in E,$ for some discrete set $D.$
\end{remark}

\begin{theorem}
	Let $\Omega\subset \mathbb{C}$, be a domain and $q\in \mathbb{Z}_+$.
	If
	$E_1,\ldots,E_{q}$ is a collection of pairwise non-contiguous OSCAR:s in $\Omega$,
	then the set $E=\bigcup_{j=1}^q E_j$ is a set of uniqueness for $q$-analytic functions on $\Omega.$ 
\end{theorem}
\begin{proof}
	By Remark \ref{noncontigrem} there exists points
	$p_1,\ldots,p_q$, $p_j\in E_j,$ such that the  
	the domains, $U_j$, of the Schwarz functions, $S_j(z)$, of $E_j$, near $p_j,$ $j=1,\ldots,q$, 
	(i.e.\ $S_j(z)$ is holomorphic on $U_j$ such that $U_j\cap \{S_j(z)=\bar{z}\}=E_j\cap U_j$)
	satisfy
	\begin{equation}
	U:=\bigcap_{j=1}^q U_j \neq \emptyset 
	\end{equation}
	where obviously $U\subset \Omega$ is an open subset.
	Let $f(z)=\sum_{j=0}^{q-1} a_j(z)\bar{z}^j$
	be a $q$-analytic function on $\Omega$ such that $f|_E=0.$
	and set
	\begin{equation}\label{triggen0}
	f_k(z):=\sum_{j=0}^{q-1} a_j(z)(S_k(z))^j, \quad z\in U,\quad k=1,\ldots,q
	\end{equation} 
	Clearly, each $f_k(z)$ is a holomorphic function on $U$ and we have
	\begin{equation}\label{triggen1}
	0=f(z)=f_k(z)=\sum_{j=0}^{q-1} a_j(z)(S_k(z))^j,\quad z\in E_k\cap U_k,\quad k=1,\ldots,q
	\end{equation}
	Now in the system defined by Eqn.(\ref{triggen1}), the right hand sides are holomorphic, thus
	vanishing on any set containing an accumulation point (in particular on $U_j\cap E_j$) implies vanishing on $U_j$ thus on $U.$
	We thus obtain the system 
	\begin{equation}\label{triggen2}
	0\equiv \sum_{j=0}^{q-1} a_j(z)(S_k(z))^j,\quad z\in U, k=1,\ldots,q
	\end{equation}
	\begin{lemma}\label{triggenlemma0}
		There exists a point $p_0\in U$ and $\epsilon>0$ such that
		the complex numbers 
		\begin{equation}
		S_1(p),\ldots,S_q(p)
		\end{equation}
		are nonzero and distinct for all $p\in U$ such that $|p-p_0|<\epsilon$
		(where $\epsilon$ is sufficiently small such that $\{|z-p_0|<\epsilon\}\subset U$).
	\end{lemma}
	\begin{proof}
	Set $Z=U\cap \bigcup_{j=1}^q S_j^{-1}(0)$ and $W:=U\setminus Z.$
	Note that, since each $S_j(z)$ is holomorphic (and we can assume it is $\not\equiv 0$ since it defines an OSCAR)
	we have that $W$ is a nonempty domain.
	Define the holomorphic function
	\begin{equation}
	\Theta(z):=\Pi_{1\leq j<k\leq q}(S_j(z)-S_k(z))
	\end{equation}
	Note that on the domain of $\Theta$ we have for each $z\in W$ that 
	the collection $S_j(z)$, $j=1,\ldots,q$, is a set of nonzero complex numbers.
	Since $\Theta$ is holomorphic it must have isolated zeros or vanish identically on $W.$
	If $p_0\in W$ such that $\Theta(p_0)\neq 0$ then by continuity
	there exists $\epsilon>0$ such that
	$\Theta(z)\neq 0$ for all $z\in \{\abs{z-p_0}<\epsilon\}\subset W$,
	which implies that
	\begin{equation}
	S_j(z)-S_k(z)\neq 0,\quad  1\leq j<k\leq q,\quad  z\in \{|z-p_0|<\epsilon\}
	\end{equation}
	and we are done.
	On the other hand, each $\phi_{jk}(z):= S_j(z)-S_k(z)$ is holomorphic thus has isolated zeros or vanishes identically.
	If $\Theta\equiv 0$ on $U$ then there must exists
	$1\leq j_0<k_0\leq q$ such that $S_{j_0}(z)\equiv S_{k_0}(z).$
	But this would contradict that $E_{j_0}$ and $E_{k_0}$ are non-continugous OSCAR:s.
	This proves Lemma \ref{triggenlemma0}.
	\end{proof}
Let $p_0$ and $\epsilon$ be as in Lemma \ref{triggenlemma0} and set $V:=\{\abs{z-p_0}<\epsilon\}$.
Then the system defined by Eqn.(\ref{triggen2}) restricted to $V\subset U,$
\begin{equation}\label{triggen2}
0\equiv \sum_{j=0}^{q-1} a_j(z)(S_k(z))^j,\quad z\in V, k=1,\ldots,q
\end{equation}
can we written in matrix form as $[a_0(z),\ldots,a_{q-1}(z)]M$ where
\begin{equation}
M:=\begin{bmatrix}
1 & 1 & \cdots & 1\\
S_1(z) & S_2(z) & \cdots & S_{q}(z)\\
\vdots & \vdots & \vdots & \vdots\\
S_1^{(q-1)}(z) & S_2^{(q-1)}(z) & \cdots & S_{q}^{(q-1)}(z)\\
\end{bmatrix}
\end{equation}
and has determinant given by
\begin{equation}
\mbox{det}(M)=\Pi_{i<j} ((S_i(z))-(S_j(z)))
\end{equation} 
Since the $S_k(z)$, $k=0,\ldots,q-1$ are distinct for all $z\in V,$ the determinant is nonzero.
This proves that $a_j\equiv 0$, $j=0,\ldots,q-1$ on $V$ thus on all of $\Omega.$
This completes the proof.	
\end{proof}

\section{Rad\'o's theorem}
Rad\'o's theorem states that a continuous function on an open subset of $\Cn$ that is holomorphic off its zero set extends to a holomorphic function on the given open set. For the one-dimensional result see Rad\'o \cite{t1}, and for a generalization to several variables, see e.g.~Cartan \cite{cartan}.
\begin{definition}
	Let $\Omega\subset\Cn$ be an open subset and let $(z_1,\ldots,z_n)$ denote holomorphic coordinates for $\Cn.$
	A function $f$, on $\Omega,$ is said to be {\em separately $C^{k}$-smooth with respect to the $z_j$-variable}, if
	for any fixed $(c_1,\ldots,c_{n-1})\in \C^{n-1},$ chosen such that
	the function 
	\[ z_j \mapsto f(c_1,\ldots,c_{j-1},z_j,c_j,\ldots,c_{n-1}), \]
	is well-defined (i.e.\ such that $(c_1,\ldots,c_{j-1},z_j,c_j,\ldots,c_{n-1})$ belongs to the domain of $f$) is $C^{k}$-smooth with respect to $\re z_j, \im z_j$.
	For $\alpha\in \Z_+^n$ we say that $f$ is separately $\alpha$-smooth if $f$ is separately
	$C^{\alpha_j}$-smooth with respect to $z_j$ for each $1\leq j\leq n$.
\end{definition}

\begin{definition}
	Let $\Omega\subseteq\Rn$ be an open subset. For a fixed $p>1,$
	 the {\em $p$-Laplacian} of a real-valued function $u$ on $\Omega$ is defined as\index{$p$-Laplacian}
\begin{equation}
\Delta_p:=\mbox{div}(\abs{\nabla u}^{p-2}\nabla u)
\end{equation}
The operator can also be defined for $p=1$ (it is then the negative of the so-called mean curvature operator) and $p=\infty$ but we shall not concern ourselves with such cases.
\end{definition}

\begin{remark}\label{regrem}
	Note the subtle similarity between the notation for the $p$-Laplacian  
	\begin{equation}
	\Delta_p=\mbox{div}(\abs{\nabla u}^{p-2}\nabla u)
	\end{equation}
	and that of the $p$:th power of the Laplacian $\Delta^p$.
	We have that $\Delta_{2} =\Delta.$ 
	At least they both share the property of being elliptic operators. In the case of $\Delta^p$ this is a direct consequence of the fact that $\Delta$ is an elliptic operator and therefore any finite power is also, in particular the elliptic regularity theorem
	applies to $\Delta^p$ and to $\Delta_p$, and implies that any real-valued 
	distribution solution $u$ to $\Delta^p u=0$ (or to $\Delta_p$) on a domain $\Omega\subset\R^n$ is Lebesgue a.e.\ equal to a $C^\infty$-smooth solution
	$\tilde{u}$ to $\Delta^p \tilde{u}=0$ (or to $\Delta_p\tilde{u}=0$) on $\Omega.$
\end{remark}
Kilpel\"ainen \cite{kilpelainen} proved the following.
\begin{theorem}\label{kipelthm}
	If $\Omega\subset\R^2$ is a domain and if $u\in C^1(\Omega)$ satisfies the $p$-Laplace equation
	$\mbox{div}(\abs{\nabla}^{p-2}\nabla u)=0$ on $\Omega\setminus u^{-1}(0)$ then $u$ is a solution to the $p$-Laplacian on $\Omega.$
\end{theorem}
We mention that, more recently, Tarkhanov \& Ly \cite{tarkhanovly} proved the following related result in higher dimension.
\begin{theorem}\label{tarkhanovthm}
	Let $\Omega\subseteq\Rn$ be an open subset. If $u\in C^{1,\frac{1}{p-1}}(\Omega)$ such that
	$\mbox{div}(\abs{\nabla}^{p-2}\nabla u)=0$ on $\Omega\setminus u^{-1}(0)$ then this holds true on all of $\Omega.$
\end{theorem}

The following is a version of a result of Daghighi \& Wikstr\"om \cite{frankabberado}.
\begin{proposition}\label{halvmainlemma}
	Let $\Omega\subset\C$ be a simply connected domain, let $q\in \Z_+$ and let $f\in C^q(\Omega)$
	be a $q$-analytic function on $\Omega\setminus f^{-1}(0).$
	Then $f$ is $q$-analytic on $\Omega.$
\end{proposition}
\begin{proof}
	If $f\equiv 0$ then we are done, so assume $f\not\equiv 0.$
	Since $f$ is $C^q$-smooth the function
	$\partial_{\bar{z}}^q f$ is continuous.
	By assumption $\partial_{\bar{z}}^q f=0$ on $\Omega\setminus f^{-1}(0).$
	Set $Z:=(f^{-1}(0))^\circ$ and $X:=\{f\neq 0\}\cup Z.$ Now $f|_Z$ clearly satisfies $\partial_{\bar{z}}^q f=0$.
	Let $p\in \partial X.$ If $p$ is an isolated zero of $f$, then
	by continuity we have $\partial_{\bar{z}}^q f(p)=0.$ Suppose $p$ is a non-isolated zero.
	We have for each sufficiently large $j\in \Z_+$
	that $\{\abs{z-p}<1/j\}\cap X\neq \emptyset$.
	This implies that for there exists a sequence $\{z_j\}_{j\in \Z_+}$
	of points $z_j\in X$ such that $z_j\to p$ as $j\to \infty.$
	By continuity we have 
	\begin{equation}\partial_{\bar{z}}^q f(p)=\lim_{j\to \infty} \partial_{\bar{z}}^q f(z_j) =0\end{equation}
	This completes the proof.
\end{proof}
\begin{theorem}\label{halvmain}%
	Let $\Omega\subset \C$ be a bounded, simply connected domain, let $q\in \Z_+$
	and let $f$ be a function $q$-analytic on $\Omega\setminus f^{-1}(0)$. 
	Suppose at least one of the following conditions holds true:\\
	(i) $f\in C^\kappa(\Omega),$ for $\kappa=\min\{1,q-1\},$ and $\re f$ ($\im f$ respectively) is a solutions to the $p'$-Laplace equation ($p''$-Laplace equation respectively) on $\Omega\setminus f^{-1}(0)$, for some $p',p''>1$.\\
	(ii) $f\in C^{q}(\Omega).$
	\\
	Then $f$ agrees (Lebesgue) a.e.\ with a function that is $q$-analytic on $\Omega.$ 
\end{theorem}
\begin{proof}
	The case (ii) follows from Proposition \ref{halvmainlemma}. So suppose (i) holds true.
		If $q=1$ the theorem is well-known and due to Rad\'o \cite{t1}, so assume $q\geq 2.$ 
	Let $f=u+iv$ where $u=\re f,$ $v=\im f.$ 
	Now $f^{-1}(0)= u^{-1}(0) \cap v^{-1}(0)$, whence $u$ (and $v$ respectively) is 
	a solution to the $p'$-Laplace equation ($p''$-Laplace equation respectively)
	on $\Omega\setminus u^{-1}(0)$ ($\Omega\setminus v^{-1}(0)$ respectively). 
	If $f\in C^{\kappa}(\Omega)$ and $q\geq 2$ then $u$ and $v$ respectively are at least $C^1$-smooth thus satisfy the conditions of 
	Theorem \ref{kipelthm}. Hence
	it follows that $u$ ($v$ respectively) are 
	solutions to the $p'$-Laplace equation ($p''$-Laplace equation respectively) on all of $\Omega$.
	By Remark \ref{regrem} (in particular Elliptic regularity) it follows that $u$ and $v$ respectively agree (Lebesgue) a.e.\ on $\Omega$ with
	$C^\infty$-smooth functions $\tilde{u}$ and $\tilde{v}$ respectively.
	This implies that the function $\tilde{f}:=\tilde{u}+i\tilde{v}$ is $C^\infty$-smooth on $\Omega$
	and agrees (Lebesgue) a.e.\ on $\Omega$ with $f.$
Suppose there exists a point $p_0\in \Omega$ 
such that $\partial_{\bar{z}}^q \tilde{f}(p_0)\neq 0.$ Set $Z:=(f^{-1}(0))^\circ$ and $X:=\{f\neq 0\}\cup Z.$ By continuity there exists an open neighborhood $U_{p_0}$ of $p_0$ in $\Omega$
such that $\partial_{\bar{z}}^q \tilde{f}\neq 0$ on the open subset $U_{p_0}\cap X.$ By the definition of $\tilde{f}$ there exists a set $E$ of zero measure such that
on $V_{p_0}:=(X\cap U_{p_0})\setminus E$ we have that $\partial_{\bar{z}}^q f$ exists (since $X$ contains no point of $f^{-1}(0)\setminus Z$) and satisfies
$0=\partial_{\bar{z}}^q f=\partial_{\bar{z}}^q \tilde{f}$ on $V_{p_0},$ which could only happen if
$V_{p_0}$ is empty which is impossible since $E$ cannot possess interior points.
We conclude that
$\partial_{\bar{z}}^q \tilde{f}=0$ on $\Omega.$ This completes the proof.
\end{proof}

\begin{theorem}[Rad\'o's theorem for polyanalytic functions in several complex variables]\label{frabsmooth}
	Let $\Omega\subset \Cn$ be a bounded $\mathbb{C}$-convex domain. Let $\alpha\in \Z_+^n$. 
	Suppose $f$ is $\alpha$-analytic on $\Omega\setminus f^{-1}(0)$ such that one of the following conditions hold true:\\
	(i) For each $j=1,\ldots,n$ $f$ is separately 
	$C^{\kappa_j}$-smooth with respect to $z_j$ (i.e.\ for each fixed value of the remaining variables $z_k,$ $k\neq j$,
	$f$ becomes a $C^{\kappa_j}$-smooth function of $z_j$), $\kappa_j=\min\{1,\alpha_j-1\}$ and $\re f$ ($\im f$ respectively) are
	solutions to the $p'$-Laplace equation ($p''$-Laplace equation respectively) for some $p',p''>1.$
	\\
	(ii) For each $j=1,\ldots,n$ $f$ is separately 
	$C^{\alpha_j}$-smooth with respect to $z_j$.\\
	Then $f$ agrees (Lebesgue) a.e.\ with a function that is $\alpha$-analytic on~$\Omega$. 
\end{theorem}
\begin{proof}
	Denote for a fixed $c\in \C^{n-1}$, $\Omega_{c,k}:=\{ z\in \Omega
	:z_j=c_j,j<k, z_j=c_{j-1},j>k \}$. Since $\Omega$ is $\C$-convex, $\Omega_{c,k}$ is simply connected. Consider the function
	$f_c(z_k):=$ $f(c_1,\ldots,c_{k-1},z_k,c_{k},\ldots,c_{n-1})$.
	Clearly, $f_c$ is $\alpha_k$-analytic on $\Omega_{c,k}\setminus
	f^{-1}(0)$ for any $c\in \C^{n-1}.$ Since $f^{-1}_c(0)\subseteq
	f^{-1}(0)$, Theorem~\ref{halvmain} applies to $f_c$ meaning that
	$f$ agrees a.e.\ with a function $\tilde{f}$ that is {\em separately} polyanalytic of order $\alpha_j$ in the
	variable $z_j, 1\leq j\leq n$.  
	By Theorem~\ref{hartog1} the function $\tilde{f}$ must be polyanalytic of order $\alpha$ 
	on $\Omega$.
	This completes the proof.
\end{proof}

\begin{corollary}
	Let $\Omega\subset \C$ be a bounded $\C$-convex domain and let $\alpha\in \Z_+^n$. 
	Suppose $f$ is separately $C^{\alpha_j}$-smooth with respect to $z_j,$ $j=1,\ldots,n.$
	If $f$ is $\alpha$-analytic on $\Omega\setminus f^{-1}(0),$
	then $f$ agrees (Lebesgue) a.e.\ with a function that is $\alpha$-analytic on~$\Omega$. 
\end{corollary}

Note that $q-1$ cannot be replaced by $q-2$.
\begin{example}
	Define
	\[
	f(z) := \begin{cases} \abs{z}^2 - 1, & \abs{z} \ge 1 \\ 1-\abs{z}^2, & \abs{z} < 1. \end{cases}
	\]
	Then $f$ is continuous and $2$-analytic off its zero set $\{ \abs{z} = 1 \}$, but not $2$-analytic.
\end{example}

It is possible to give a proof for a similar result using other types of conditions as follows (see Daghighi \& Wikstr\"om \cite{daghighiwikstromjpde}).
\begin{definition}\label{bmp}
	Let $\Omega\subset \Cn$ be a submanifold and let
	$\mathfrak{M}\subset C(\Omega,\C)$ be a family of functions.  Let
	$D$ denote the unit disc, $D:=\{\zeta\in \C:\abs{\zeta}<1\}$.  The
	family $\mathfrak{M}$ is said to have the \emph{one dimensional
		boundary maximum modulus property} (\textsc{1d-bmmp}) if given
	$f\in \mathfrak{M}$, then for any $\psi\in \mathscr{O}(D,\Omega)\cap
	C(\overline{D},\Omega)$,
	\begin{equation}
	\max_{\zeta\in \overline{D}}\abs{f\circ \psi(\zeta)}\leq
	\max_{\zeta\in \partial D}\abs{f\circ \psi(\zeta)}.
	\end{equation}
\end{definition}
\begin{example}
	Let $\Omega\subset \Cn$ be a complex submanifold. Then
	$\mathscr{O}(\Omega)\subset C(\Omega,\C)$ clearly has the one
	dimensional boundary maximum modulus property.
\end{example}
We shall use the following notation:
Let $U\subset\Cn$ be a submanifold. Denote by
$\mathscr{M}_{\alpha}(U)$ the set of restrictions to $U$ of
$\alpha$-analytic functions on $\Cn$ that obey the
one dimensional boundary maximum modulus property of
Definition~\ref{bmp}.
An immediate consequence of the fact that the set of holomorphic
functions obey a strong maximum principle, is that, if $U\subset\Cn$
is a complex submanifold (not necessarily of dimension $n$), then
$\mathscr{O}(U)\subseteq\mathscr{M}_{\alpha}(U)$. The set $\mathscr{M}_{\alpha}$
also contains non-holomorphic functions. For example, harmonic ($\alpha$-analytic)
functions  belong to $\mathscr{M}_\alpha$, so if
$U \subset \mathbb{C}$, then $\mathscr{M}_{2}(U)$ contains functions such as
$z + \bar z$. There are also non-harmonic examples, such as $z\bar z \in \mathscr{M}_2$.
Unfortunately, we are not aware of any simple characterization of $\mathscr{M}_\alpha$, but
let us note that $\mathscr{M}_\alpha$ has no obvious algebraic structure. In particular, the sum
of two functions in $\mathscr{M}_\alpha$ is not necessarily an $\mathscr{M}_\alpha$-function. (For example
$1$ and $-z\bar z$ are in $\mathscr{M}_2$, but $1-z\bar z$ is not.) 
\\
\\
Kaufman \cite{kauf} makes use of a maximum principle
and combines this with an approximation property on the boundary to prove
Rad\'o's theorem and we shall use a similar method of proof for
the following generalization for the case of bounded
domains.
We will need a fairly weak assumption on our domains. In fact, we need to
approximate continuous functions on $\partial \Omega$ uniformly by global
harmonic functions. There is a harmonic version of Mergelyan's theorem that takes
care of this. In order to state the theorem, let us begin with some preliminary
notation. If $K \subset \mathbb{C}$ is compact, we let $\hat K$ denote the union of
$K$ and all bounded components of $\mathbb{C} \setminus K$.

\begin{definition}
	Let $\Omega \subset \mathbb{C}$ be a bounded domain.
	We say that $\Omega$ is \emph{harmonically fat}
	if $\partial (\widehat{\partial\Omega}) = \partial \Omega$.
	If $\Omega \subset \Cn$ ($n > 1$), we say that $\Omega$ is harmonically fat
	if the intersection of $\Omega$ with any complex line is (empty or)
	harmonically fat.
\end{definition}

In particular, note that simply connected domains in $\mathbb{C}$ (and hence
$\mathbb{C}$-convex domains in $\Cn$) are harmonically fat. One harmonic version
of Mergelyan's theorem~\cite[Corollary 1.16]{gardiner} can now be stated as:

\begin{theorem}\label{harmonic-mergelyan}
	Let $\Omega \subset \mathbb{C}$ be harmonically fat and let $u \in C(\partial \Omega)$.
	Then, for every $\varepsilon > 0$, there is an entire harmonic function $v$ such that
	$|v-u| < \varepsilon$ on $\partial \Omega$.
\end{theorem}
We begin with the single variable case.
\begin{theorem}\label{radoprop}
	Let $\Omega\subset\C$ be a bounded harmonically fat domain. Let $q$ be a positive
	integer and let $f \in C^{q-1}(\overline{\Omega})$ such that $f\in
	\mathscr{M}_q(\Omega\setminus f^{-1}(0))$. Then $f\in
	\mathscr{M}_q(\Omega)$.
\end{theorem}

\begin{proof}
	Note that $q=1$ corresponds to the well-known Rad\'o's theorem in one
	complex variable.
	
	\begin{remark}\label{remarkrado}
		If $U\subset \Cn$ a submanifold and $G$ is a continuous function on
		$\overline{U}$ which satisfies the boundary maximum modulus principle on the
		(necessarily open) subset $U\setminus G^{-1}(0)$ then $G$ satisfies
		the boundary maximum modulus principle on $U$. Indeed, let $V\Subset
		U,$ be a domain. If $V\cap G^{-1}(0)=V$ then, by continuity,
		$\max_{z\in \overline{V}}\abs{G(z)}=0$, so we are done. If instead
		the (necessarily open) set $V\cap \{G\neq 0\}$ is nonempty, then
		$\max_{z\in \overline{V}}\abs{G(z)}=$ $\max_{z\in \overline{V\cap
				\{G\neq 0\}}}\abs{G(z)}=$ $\max_{z\in \partial V\cap \{G\neq
			0\}}\abs{G(z)}$ $=\max_{z\in \partial V}\abs{G(z)}$.
	\end{remark}
	
	By Remark \ref{remarkrado} it is sufficient to show that $f$ extends
	to a $q$-analytic function on $\Omega$. For the sake of clarity we
	first prove the case $q=2$, i.e.\ assume $f\in
	\mathscr{M}_q(\Omega\setminus f^{-1}(0))$ (only small modifications
	are then required to prove the cases $2<q<\infty$).  
	By a well-known result (we shall prove this later, see Corollary \ref {ahernbrunakonsekvens} or Theorem \ref{balkuniformconv}) 
	we have that $f$ is $2$-analytic on $\Omega\setminus
	f^{-1}(0)$, and by definition we know that $f$ and $\frac{1}{f}$ both
	satisfy the boundary maximum modulus principle on $\Omega\setminus
	f^{-1}(0)$. Set
	\begin{equation}
	u(z):=\frac{\partial}{\partial\overline{z}} f(z)=\partial_{\bar{z}} f(z).
	\end{equation}
	Since $f\in C^{2-1}(\Omega)$, the function $u$ is well-defined on
	$\Omega$. Let $U = \Omega \setminus f^{-1}(0)$ for convenience of
	notation. By definition, $u$ is holomorphic on $U$. Furthermore, $u =
	0$ on the interior of $f^{-1}(0)$ so,
	\begin{equation}
\partial_{\bar{z}} u(z) = 0, \qquad \forall z \in \Omega \setminus \partial U.
	\end{equation}
	
	We shall need the following lemma.
	
	\begin{lemma}\label{flem}
		Let $g$ be a function continuous on $\overline{U}$ and holomorphic
		on $U$.  Then for all $z\in U$,
		\begin{equation}
		|g(z)| \leq \sup_{\zeta \in \partial U \cap \partial \Omega} |g(\zeta)|.
		\end{equation}
	\end{lemma}
	
	\begin{proof}
		First of all, $\sup_{\partial U\cap \Omega^{\circ} }\left| f\right|
		=0$ since $f=0$ on the given set.  Secondly, note that $f$ and $g^j$
		each satisfies the boundary maximum modulus principle applied to
		$U$.  Thus for any $z\in U$,
		\begin{equation}
		\left|f(z)\right| \left|g(z)\right|^j  \leq \left(\sup_{\partial U\cap \overline{\Omega} } \left| f\right|\right)\cdot \left(
		\sup_{\partial U\cap \overline{\Omega}} \left| g\right|
		\right)^j\leq
		\left(\sup_{\partial U\cap \partial \Omega } \left| f\right|\right)\cdot \left(
		\sup_{\partial U\cap \partial \Omega} \left| g\right|
		\right)^j,
		\end{equation}
		which in turn implies, after taking $1/j$\,th power and passing to the limit
		as $j\to \infty$,
		\begin{equation}
		\left| g(z)\right|\leq \sup_{w\in \partial
			U\cap \partial \Omega } \left|g(w) \right|,\qquad \text{for $z\in
			\Omega\setminus \partial U$.}\qedhere
		\end{equation}
	\end{proof}
	By Lemma \ref{flem},
	\begin{equation}\label{ek00}
	\left| u(z)\right|\leq \sup_{w\in \partial U\cap \partial \Omega } \left|u (w) \right|, z\in U.
	\end{equation}
	Also, we know that $\abs{u(z)}=0$ for all $z$ in the interior of
	$f^{-1}(0)$. In fact, given Lemma \ref{flem}, a verbatim repetition of
	an argument which can be found in e.g.\ Kaufmann \cite{kauf}, proves
	that $U$ must be a dense subset of $\Omega$, in particular $f^{-1}(0)$
	must have \emph{empty} interior.  This together with
	Eqn.\eqref{ek00} gives,
	\begin{equation}\label{ek01}
	\left| u(z)\right|\leq \sup_{w\in \partial U\cap \partial \Omega } \left|u(w) \right|,
	\quad z\in \Omega\setminus \partial U.
	\end{equation}
	Next we show that $\partial u/\partial\bar{z}$ is harmonic, and since
	it is zero on a dense open subset of $\Omega$ it vanishes identically.
	To see that $\partial u/\partial\bar{z}$ is harmonic, set $u = w+iv$,
	and show that $w$, $v$ are harmonic as follows:
	Since we assume that $\Omega$ is harmonically fat, we can
	find a sequence $\{P_j\}_{j\in \N}$,
	of entire functions such that $\re (P_j-u)\to 0$
	on $\partial \Omega$. We have $\abs{e^{P_j-u}}=e^{\re(P_j-u)}$, and
	also that $e^{P_j-u}$ and $e^{u-P_j}$ are holomorphic on $U$. Thus the
	maximum principle of Eqn.\eqref{ek01} applies to both
	$e^{P_j-u}$ and $e^{u-P_j}$. We can choose $P_j$ such that $\abs{P_j
		-u} < \frac{1}{j}$. Consequently,
	\begin{equation}\label{ineq0}
	\abs{e^{(P_j(z) -u(z))}} < e^{\frac{1}{j} } \text{ and }
	\abs{e^{(u(z)- P_j(z))}} < e^{\frac{1}{j} },
	\quad z\in \partial \Omega ,
	\end{equation}
	and by the maximum principle of Eqn.\eqref{ek01}, the
	inequalities Eqn.\eqref{ineq0} continue to hold in $\Omega$. This in turn
	implies that
	\begin{equation}\label{ineqw00}
	e^{\frac{1}{j}}>\abs{e^{P_j(z)-u(z)}}=e^{\re(P_j(z)-u(z))}, \quad\forall z\in \Omega.
	\end{equation}
	After taking logarithms,
	\begin{equation}\label{ineqw0}
	\frac{1}{j}>\re(P_j(z)-u(z)),\quad\forall z\in \Omega.
	\end{equation}
	
	Since the real part of each $P_j$ is harmonic, this uniform
	convergence implies that $\re u(z)=w(z)$ is harmonic on
	$\Omega$. Analogously one shows that $v$ is harmonic.  It then follows
	that $\partial u/\partial\bar{z}$ has harmonic real and imaginary
	parts since, $\Delta \frac{\partial u}{\partial\bar{z}} =
	\frac{\partial }{\partial x}(\Delta w)- \frac{\partial }{\partial
		y}(\Delta v) +i\left( \frac{\partial }{\partial y}(\Delta w)-
	\frac{\partial }{\partial x}(\Delta v) \right) =0$. A harmonic
	function (in particular the real and imaginary part respectively of
	$\frac{\partial u}{\partial\bar{z}}$) vanishing on a dense open subset
	vanishes identically thus $\frac{\partial u}{\partial\bar{z}}$
	vanishes identically on $\Omega$ i.e.\ $u$ is holomorphic on $\Omega$
	meaning that,
	\begin{equation}
	0=\partial_{\bar{z}} u(z)=\partial_{\bar{z}}^2 f(z),\quad\forall z\in \Omega,
	\end{equation}
	i.e.\ $f$ is bianalytic on $\Omega$. This proves
	Theorem~\ref{radoprop} for $q=2$.
	\\
	\\
	Now we can easily adapt the proof to the cases $q>2$. Assume $f\in
	C^{q-1}(\Omega)$ is $q$-analytic on $U$ and as before let
	$u(z):=\partial_{\bar{z}}^{q-1} f(z)$.  If $f^{-1}(0) \cap \Omega$ has nonempty
	interior, then also $\partial_{\bar{z}}^q f(z)=0$ on $\left( \Omega\cap
	f^{-1}(0)\right)^{\circ}$, thus $\partial_{\bar{z}} u(z) =0$ on
	$\Omega\setminus\left(\partial f^{-1}(0)\right)$, which is a dense
	open subset of $\Omega.$ Applying the same arguments to $u$ as for the
	case $q=2$ we obtain that $\partial_{\bar{z}} u(z)$ vanishes identically on
	$\Omega,$ in particular $u$ is differentiable\footnote{In fact, $u$
		must be $C^{\infty}$-smooth a.e.\ (due to 
		ellipticity).} on $\Omega$
	thus $\partial_{\bar{z}}^{q-1} f(z)$ is well-defined and differentiable on
	$\Omega$, meaning that we can write $\partial_{\bar{z}}^q f=
	\frac{\partial}{\partial\overline{z}}\left(\partial_{\bar{z}}^{q-1}
	f(z)\right)=0$, $\forall z\in \Omega$. This completes the proof.
\end{proof}

Note the importance of the starting function $f$ to be $C^{q-1}(\Omega)$
instead of merely continuous, namely, we need in the proof for $q>2$
that $u$ be continuous on $\Omega$.
\begin{example}
	Let $\Omega:=\{\abs{z}<2\}\subset\C,$ and set,
	\begin{equation}
	f(z):=
	\left\{
	\begin{array}{ll}
	1-z\bar{z}, & \abs{z}\leq 1, \\
	z\bar{z}-1, & 1< \abs{z}< 2.
	\end{array}
	\right.
	\end{equation}
	The function $f$ is clearly $2$-analytic on the open subset
	$\Omega\setminus f^{-1}(0)$ (which consists of two disjoint domains),
	and we have $f\in C^0(\Omega)$, $0=q-2$. However $f$ is \emph{not}
	$2$-analytic at any point of $\{\abs{z}=1\}$.
\end{example}	
	Note that this example breaks \emph{both} the regularity assumption
	\emph{and} the boundary maximum modulus principle required in
	Theorem~\ref{radoprop}. 
We immediately obtain the following consequence of Theorem~\ref{radoprop}.

\begin{theorem}\label{radoprop1}
	Let $\Omega\subset\Cn$ be a bounded harmonically fat domain and let $(z_1,\ldots
	,z_n)\in \Cn$ denote holomorphic variables. Let
	$\alpha=(\alpha_1,\ldots,\alpha_n)\in \Z_+^n$. Then any
	function $f$ which is $C^{\alpha_j-1}$-smooth in the $z_j$ variable,
	for each $j$, up to the boundary
	and which is a member of $\mathscr{M}_\alpha(\Omega\setminus
	f^{-1}(0))$ is automatically a member of $\mathscr{M}_\alpha(\Omega).$
\end{theorem}

\begin{proof}
	Recall that when $\Omega$ is a complex manifold of same dimension as
	the ambient space, $\mathscr{M}_\alpha$ reduces to a subspace of
	$\alpha$-analytic functions which satisfy the boundary maximum
	modulus principle.  Denote for a fixed $c\in \C^{n-1}$,
	$\Omega_{c,k}:=\{ z\in \Omega :z_j=c_j,j<k, z_j=c_{j-1},j>k
	\}$. Consider the function $f_c(z_k):=$
	$f(c_1,\ldots,c_{k-1},z_k,c_{k},\ldots,c_{n-1})$. Note that the
	restriction of a function in $n$ complex variables which satisfies
	the boundary maximum modulus principle, to a complex one-dimensional
	submanifold, must also satisfy the boundary maximum modulus
	principle.  Clearly, $f_c$ is $\alpha_k$-analytic on
	$\Omega_{c,k}\setminus f^{-1}(0)$ for any $c\in \C^{n-1}.$ Since
	$f^{-1}_c(0)\subseteq f^{-1}(0)$, Theorem~\ref{radoprop} applies
	to $f_c$ meaning that $f$ is \emph{separately} polyanalytic of order
	$\alpha_j$ in the variable $z_j$, $1\leq j\leq n$. By
	Theorem~\ref{hartog1} the function $f$ must be polyanalytic of order
	$\alpha$ 
	on $\Omega.$ By
	Remark \ref{remarkrado}, $f$ satisfies the boundary maximum modulus
	principle. This completes the proof.
\end{proof}

\section{Non-isolated zeros in the entire case }
The following result can rather quickly be proved using the theory of semi-analytic sets (see
e.g.\ the proof of Theorem \ref{mpmainresult}, 
however we give here the original proof
of Balk \& Zuev \cite{balkzuev1972}. 
\begin{proposition}\label{numbervalueslem1}
	If $f(z)$ is an $n$-analytic function on a domain $\Omega\subseteq \C$
	with a non-isolated zero at a point $p_0\in \Omega$ then 
	there exists analytic arc 
	$\Gamma$, near $p_0$, in $f^{-1}(0)$ 
	with end point $p_0$.
\end{proposition}
\begin{proof}
	W.l.o.g.\ assume $p_0=0$ and 
	suppose $f(z)$ has an nonisolated zero at $z=0.$
	Since $\Omega$ is a domain
	we may write $f(z)=\sum_{j=0}^{n-1} A_j(z)\bar{z}^j$
	for holomorphic $A_j$, $j=0,\ldots,n-1.$
	and by assumption this has a 
	nonisolated zero. Define 
	\begin{equation}
	g(z,w):=\sum_{j=0}^{n-1} A_j(z)w^j
	\end{equation}
	By the Weierstrass preparation theorem it is possible to write
	\begin{equation}
	g(z,w)=z^\mu F(z,w)\Pi_{j=1}^s P_j(z,w)
	\end{equation}
	for some holomorphic $F(z,w)$ that is nonzero 
	in a neighborhood of $(0,0)$, $\mu\geq 0$ and irreducible
	pseudopolynomials $P_j,$ $j=1,\ldots,s.$
	Since $f(z)$ has a nonisolated zero then at least one of the functions $P_j(z,\bar{z})$
	has a nonisolated zero. Hence it suffices to prove the lemma for
	the restriction to $z=\bar{z}$ of an irreducible
	pseudopolynomial of the form
	\begin{equation}
	P(z,w)=w^m+\sum_{j=0}^{m-1} a_j(z)w^j
	\end{equation}
	for some holomorphic $a_j(z)$ near $z=0.$
	It is well-known (see e.g.\ Fuks \cite{fuks1962}, p.98) 
	that the analytic set $P(z,w)=0$
	can, in some neighborhood $\Delta$ of $(0,0)$ be described as
	\begin{equation}
	w-A(z^{\frac{1}{m}})=0
	\end{equation}
	for $A(\zeta)$ holomorphic in $\zeta$ near $\zeta=0.$ 
	For each point $(z,w)\in \Delta$ belonging to $w=\bar{z}$
	and belonging to the zero set of $P(z,w)$ we have that
	\begin{equation}
	\bar{z}-A(z^{\frac{1}{m}})=0
	\end{equation}
	has at least one root $z^{\frac{1}{m}}.$
	With the new variable $\zeta$ defined by $z=\zeta^m$
	the latter condition becomes
	$\bar{\zeta}^m -A(\zeta)=0.$
	Since by assumption $z=0$ is a nonisolated zero, 
	there exists a sequence $\{z_j\}_{j\in \N}$ in 
	the zero set of $P(z,\bar{z})$ going to $z=0$, and associated
	sequence $\{\zeta_j\}_{j\in \N}$ going to zero such that
	\begin{equation}
	\zeta_j\to 0, j\to \infty,\quad \bar{\zeta}_j^m-A(\zeta_j) =0
	\end{equation}
	From this it is possible to deduce 
	$A(\zeta)\equiv \zeta^m a(\zeta)$ for an analytic function $a(\zeta)$ near $\zeta=0$
	with $a(0)\neq 0,$
	since if we Taylor expand $A(\zeta)=\sum_{j\geq k_0} A_{j}z^j$ about the origin,
	we have that 
	$\zeta^m (\exp(-2mi\mbox{arg}\zeta)+  \sum_{j\geq k_0} A_{j}z^{j-m})$
	is bounded at $0$, by continuity.  
	Hence it is possible to decompose $\bar{\zeta}^m-A(\zeta)$ according to
	\begin{equation}\label{alskdjfhg}
	\bar{\zeta}^m-A(\zeta) =\Pi_{j=1}^m(\bar{\zeta}-\lambda_j\theta(\zeta))
	\end{equation}
	Where $\lambda_j$, $j=1,\ldots,m$ are roots of unity and $\theta(\zeta)$
	is one of the branches of $(a(\zeta))^{\frac{1}{m}}.$
	Since the point $\zeta=0$ is a nonisolated zero of
	Eqn.(\ref{alskdjfhg}) it is a nonisolated zero of at least one of the bianalytic functions
	$\bar{\zeta}-\lambda_j\theta(\zeta),$ $j=1,\ldots,m.$
	As for each $j$, this is the equation for a Schwarz function,
	there exists
	an analytic arc $\gamma:=\{\zeta=\lambda(t),\alpha\leq t\leq \beta\}$
	each point of which is a zero of
	Eqn.(\ref{alskdjfhg}) and we may assume $\lambda(\alpha)=0$, i.e.\ one end of the arc is an
	accumulation point of the zero set (see e.g.\ Balk \cite{balkbianal64}).
	Then the set $\gamma:=\{z=\zeta^m\}$ for $\zeta$ in the domain of $A$,
	maps to an analytic arc $\Gamma$ with end point $p_0=0.$
	This proves Proposition \ref{numbervalueslem1}.
\end{proof}

\begin{definition}
	A function $f(z)$ on a domain $\Omega\subset \C$ is said to take the value $a\in \C$
	in a non-isolated way if the set of roots of the equation $f(z)=a$ has
	an accumulation point.
\end{definition}

\begin{corollary}\label{numbervalueslem01}
	If $f(z)$ is an $n$-analytic function 
	takes a values $a$ in a non-isolated way then 
	for each nonisolated preimage $p_0$
	there exists analytic arc
	contained in the preimage $f^{-1}(a)$
	with end point being $p_0.$
\end{corollary}

The following can be found in Borel \cite{borel1897}.
\begin{proposition}\label{borel1897lemma}
	Let $Q(z)$ be a nonconstant entire function and let $\lambda(r):=\max_{\abs{z}=r} \abs{Q(z)}$ be the real valued increasing positive function of $r:=\abs{z}.$
	Let $G_1(z),\ldots,G_n(z)$ be entire functions such that
	each $\abs{G_j(z)}$ grows slower than $\exp(\lambda(r))$ as $r\to \infty$, $j=1,\ldots,n,$ (in particular this is the case when
	each $G_j$ is a polynomial).
	Let $H_1(z),\ldots,H_n(z)$ be entire functions such that 
	each $\abs{H_j(z)}$ grows slower than $\exp(\lambda(r))$ or there exists an increasing function $\nu(r)$ dominating $\lambda(r)$
	such that $\abs{H_i-H_j}$ lies between $(\nu(r))^2$ and $\exp(\nu(r))$. 
Suppose that each 
$\abs{H_i-H_j}$ grows
	faster than $\left(\lambda(r)\right)^2.$
	Then the equations
	\begin{equation}\label{borelekv1}
	\sum_{j=1}^n G_j(z)\exp H_j(z)\equiv 0
	\end{equation}
	imply
	\begin{equation}
	G_j(z)\equiv 0,\quad j=1,\ldots,n
	\end{equation}
\end{proposition}
\begin{proof}
	We use induction in the number of members of the collection $n\geq 1.$
	The result is obvious for the case $n=1$. Suppose the result holds true for collections with $1,\ldots,n-1$ members.
	dividing by 
	$G_n\exp(H_n)$ the equation $\sum_{j=1}^n G_j(z)\exp H_j(z)$
	$=0$ and differentiating with respect to $z$ yields
	\begin{equation}\label{borelekv2}
	(G_1'G_n-G_1 G_n'+G_1G_n H_1'-G_1G_nH_n')\exp(H_1-H_n)+\cdots =0
	\end{equation}
	The coefficients of the exponentials in Eqn.(\ref{borelekv2}) then satisfy the growth conditions of the lemma
	thus they all must vanish by the induction hypothesis since there
	are only $n-1$ exponentials in Eqn.(\ref{borelekv2}).
	We show that this implies that all $G_j$ vanish identically.
	Indeed, if the first coefficient in Eqn.(\ref{borelekv2}) vanishes then
	\begin{equation}\label{borelekv00}
	\frac{G_1}{G_n}\exp(H_1-H_n)\equiv\mbox{const.}
	\end{equation}
	If neither $G_1$ nor $G_n$ is identically zero
	then Eqn.(\ref{borelekv00}) is impossible, since 
	then we may find an infinite set of circles on which
	$G_1/G_n$ is greater in modulus than $\exp(-2\lambda(r)),$ and where
	$\abs{\exp(H_1-H_2)}=\exp(\re(H_1-H_2))$ exceeds the function $(\exp(\lambda(r))^2)$ 
		thus their product can be made arbitrary large, hence the product cannot be constant.
	However if it where the case that at least $G_1$ nor $G_n$ is identically zero then
	Eqn.(\ref{borelekv2}) would reduce to one with $n-1$ terms thus the induction hypothesis would apply.
	This completes the proof for the case when each $H_j(z)$ grows slower than $\exp(\lambda(r))$, $j=1,\ldots,n.$
	Now consider the case when we do not require that each $H_j(z)$ grows slower than $\exp(\lambda(r))$, $j=1,\ldots,n.$
	If there exists an increasing function $\nu(r)$ dominating $\lambda(r)$
	such that $\abs{H_i-H_j}$ lies between $(\nu(r))^2$ and $\exp(\nu(r))$ 
	then the result follows by replacing $\lambda$ by $\nu.$
	This completes the proof.
\end{proof}
In fact, the result holds true even if one drops the conditions
that 
each $H_j(z)$ grows slower than $\exp(\lambda(r))$ or that there exists an increasing function $\nu(r)$ dominating $\lambda(r)$
such that $\abs{H_i-H_j}$ lies between $(\nu(r))^2$ and $\exp(\nu(r))$, 
Borel \cite{borel1897}.
We shall only need a special case.
\begin{lemma}\label{borellemcor}
	if $H_j(z)$, $j=1,\ldots,n$ are entire functions such that
\begin{equation}\label{jockoolinik}
H_j(z)-H_l(z)\not\equiv \mbox{const.},\quad j\neq l
\end{equation}
and if $p_j(z)$ $j=1,\ldots,n$  are polynomials then 
\begin{equation}
\sum_{j=1}^n p_j(z)\exp H_j(z)\equiv 0 \Rightarrow p_j(z)\equiv 0,\quad j=1,\ldots,n
\end{equation}
\end{lemma}
\begin{proof}
	The conditions of Eqn.(\ref{jockoolinik}) implies (since there are finitely many) that
	 there exists $i_0,j_0$ such that with
	$\lambda(r)^2=\max_{\abs{z}=r} \abs{H_{i_0}-H_{j_0}}$
	each
$\abs{H_i-H_j}$ grows
faster than $\left(\lambda(r)\right)^2$ as $r\to \infty.$ Furthermore, since $H_{i_0}-H_{j_0}$ is entire and nonconstant 
any polynomial $p_j(z)$ will satisfy that $\abs{p_j(z)}$ grows slower than $\exp(\lambda(r))$ as $r\to\infty.$
By Proposition \ref{borel1897lemma} this proves the lemma.
\end{proof}
Let us call a function of the form $B(z)=b_0(z)+\bar{z}b_1(z)$ a meromorphic bianalytic function if $b_0,b_1$ are meromorphic on $\C.$
Balk \cite{balkbianal64} proved for the case $q=2$ that a meromorphic bianalytic function, $B(z)$, on $\C$ 
which has a nonisolated zero can be represented as $B(z)=M(z)V(z)$ for a meromorphic $M(z)$ and a degenerate (in the sense that its image does not contain interior points) bianalytic 
$V(z).$ Utilizing the known possible representations of entire bianalytic functions in Balk \cite{balkbianal64b} 
it is possible to deduce that $V(z)$ is a real-valued polynomial. Furthermore, it is clear that
if $F(z)$ is a product of entire bianalytic functions then the same type of decomposition will remain valid in the sense that if
each bianalytic factor has a nonisolated zero then $F(z)$ can be decomposed into $M_1(z)V_1(z)$
for meromorphic $M_1(z)$ and $V_1(z)$ a real-valued polynomial.
Balk \& Zuev \cite{balkzuev1972}
generalize this as follows.

\begin{theorem}\label{bztheorem1}
	Every function $F(z)=\sum_{j=0}^{n-1} A_j(z)\bar{z}^j$ where the $A_j$ are meromorphic on $\C$, can be represented as a product 
	$M(z)V(z)$ where $M(z)$ is meromorphic with only isolated zeros and $V(z)$ is a real-valued polyanalytic polynomial.
\end{theorem}
\begin{proof}
	Let for entire functions $e_0,\ldots,e_{n-1}$
	\begin{equation}
	\Pi(z)=\sum_{j=0}^{n-1} e_j(z)\bar{z}^{n-1} 
	\end{equation}
	be an $n$-analytic function, with exact order $n\geq 2$ which we assume to be irreducible in the sense that
	$\Pi(z)$ cannot be written as a product of two polyanalytic functions that are not analytic.
	Suppose that $\Pi(z)$ has a nonisolated zero. We may write
	\begin{equation}\label{bzekv3} 
	\Pi(z)=e(z)\pi(z),\quad \pi(z):=\sum_{j=0}^{n-1} a_j(z)\bar{z}^{n-1} 
	\end{equation}
	where the analytic components $a_j(z)$ of $\pi(z)$ have no common root and where $\pi(z)$ is again irreducible in 
	the sense that
	$\Pi(z)$ cannot be written as a product of two polyanalytic functions that are not analytic.
	Define the associated pseudopolynomials
	\begin{equation}\label{bzekv4}
	A(z,w)=\sum_{j=0}^{n-1} w^j a_j(z)
	\end{equation}
	\begin{equation}
	B(z,w)=\overline{A(\bar{z},\bar{w})}=\sum_{j=0}^{n-1} w^j b_j(w)
	\end{equation}
	where $b_j(w):=\overline{a_j(\bar{w})}.$
	\begin{lemma}\label{visaclaimbz6}
		There exists an entire function $h(z,w)$ such that
		\begin{equation}\label{visaclaimbzekv6}
		A(z,w)=B(z,w)\exp(h(z,w))
		\end{equation}
	\end{lemma}
	\begin{proof}
		If $G(z,w)$ is an entire function we say that $D(z,w)$ is an entire divisor of $G(z,w)$ if $D(z,w)$ is entire and 
		there exists an entire $H(z,w)$ such that $G=D\cdot H.$
		We shall use the fact that (see Zuev \cite{zuev1970pseuopolyn}) that each 
		entire divisor of an entire pseudopolynomial is, up to an irreducible factor which has no roots, an entire pseudopolynomial.
		It follows (see Osgood \cite{osgood1917} and \cite{osgood1929}, p.219) that there exists entire
		functions $a(z,w),b(z,w), D(z,w)$ such that
			\begin{equation}
		A(z,w)=D(z,w)a(z,w),\quad B(z,w)=D(z,w)b(z,w)
		\end{equation}
		Consider the set of points $(z_0,w_0)$ where $a(z,w)$ and $b(z,w)$ vanish simultaneously.
				If $D(z,w)$ does not vanish anywhere then for each root, $z_0$, of $\pi(z)$, the
		point $(z_0,\bar{z}_0)$ would be a common root of $a(z,w)$ and $b(z,w)$, where we know that 
		the function
		$A(z,w)/B(z,w)=a(z,w)/b(z,w)$ is meromorphic.
		By Proposition \ref{numbervalueslem1}
		we have that (since the set of all roots of $\pi(z)$ will then contain a terminal condensation point) the zero set of $\pi$ contains an analytic arc, which is impossible since
		the set of points $(z_0,\bar{z}_0)$ which in $\C^2$ are indeterminant points for a meromorphic
		function ($A(z,w)/B(z,w)$) is discrete (see Zuev \cite{zuev1970pseuopolyn}, p.157).
		Hence there exists a point $(z,w)$ at which $D(z,w)$ vanishes.
		It follows that there exists a pseudopolynomial $d(z,w)$ (finite degree
		with respect to $w$ with analytic coefficients with respect to $z$)
		and an entire $H(z,w)$ function without zeros such that
		\begin{equation}
		D(z,w)=H(z,w)d(z,w)
		\end{equation}
		which in turn implies that there exists an entire function $\gamma(z,w)$
		such that
		\begin{equation}
		D(z,w)=d(z,w)\exp(\gamma(z,w))
		\end{equation}
		Now if $a(z,w)$ would have a zero in $\C^2$
		then there would exist an entire $\delta(z,w)$ and a pseudopolynomial $\alpha(z,w)$ 
		(with analytic coefficients with respect to $z$) vanishing at least at one point,
		such that
		$a(z,w)$ could be represented as
		\begin{equation}
		a(z,w)=\alpha(z,w)\exp(\delta(z,w))
		\end{equation}
		But that would imply
		\begin{equation}\label{bzvaluesekv9}
		A(z,w)=d(z,w)\alpha(z,w)\exp(E(z,w)),\quad E=\gamma+\delta
		\end{equation}
		Since for each $z$ the analytic component $a_j(z)$ of the pseudopolynomial $A(z,w)$ does not vanish by assumption
		the same is true for $d/z,w)$ and $\alpha(z,w)$ and therefore their exact degree
		with respect to $w$ is greater than $0$. The function $E(z,w)$ can be decomposed into
		\begin{equation}
		E(z,w)=\sum_{j=0}^\infty E_j(z)w^k
		\end{equation}
		for entire $E_j(z).$
		It is not possible that $E_j\not\equiv 0$ for some $k\geq 1$ because then by choosing $z_0$
		such that at least one of $E_j(z_0)$ is nonzero we would get according to Eqn.(\ref{bzvaluesekv9})
		\begin{equation}\label{bzvaluesekv9}
		\exp(E(z_0,w)=\frac{A(z_0,w)}{d(z_0,w)\alpha(z_0,w)}
		\end{equation}
		which is impossible since a transcendental function can not equal a rational function. Thus
		we must have
		\begin{equation}
		E(z,w)=E_0(z),\quad A(z,w)=d(z,w)\alpha(z,w)\exp(E_0(z))
		\end{equation}
		This shows that the pseudopolynomial $A(z,w)$ is in the ring of entire pseudopolynomials thus the
		polyanalytic function $\Pi(z)$ is reducible
		in the class entire polyanalytic functions which is a contradiction to our assumption in the beggining of the proof.
		Thus we conclude, by this contradiction (which was obtained by assuming that $a(z,w)$ had a zero)
		that $a(z,w)$ has no zeros in $\C^2$.
		Thus there exists an entire $h_1(z,w)$ such that 
		\begin{equation}
		a(z,w)=\exp(h_1(z,w))
		\end{equation}
		Applying the analogous arguments to $b(z,w)$ taking into account that $B(z,w)$ is a pseudopolynomial in $z$ we obtain that there exists
		an entire $h_2(z,w)$ such that $b(z,w)=\exp(h_2(z,w)$. We thus obtain with $h(z):=h_1(z,w)-h_2(z,w)$
		\begin{equation}
		A(z,w)=B(z,w)\exp(h(z,w))
		\end{equation}
		for an entire $h(z,w),$ which is Eqn.(\ref{visaclaimbzekv6}).
			This proves Lemma \ref{visaclaimbz6}.
	\end{proof}
	
	\begin{lemma}\label{visaclaimbzny}
		It is possible to choose numbers $w_1,\ldots,w_n$ such that
		\begin{equation}
		h(z,w_j)-h(z,w_1),\quad j=2,3,\ldots,n
		\end{equation}
		are constant, i.e.\ independent of $z.$
	\end{lemma}
	\begin{proof}
		Assume not (in order to reach a contradiction). Then we have that for any choice of $w_{2n}$
		there must exists
		$2n-1$ numbers $w_1,\ldots,w_{2n-1}$ such that
		\begin{equation}
		h(z,w_j)-h(z,w_l)\not\equiv\mbox{const.},\quad k,l=1,\ldots,2n, \quad k\neq l
		\end{equation}
		By Eqn.(\ref{bzekv4}) and Eqn.(\ref{visaclaimbzekv6}) we get two systems
		\begin{equation}\label{bzekv11}
		\sum_{j=0}^{n-1} a_j(z) w_k^j =B(z,w_k) \exp( h(z,w_k)),\quad k=1,\ldots,n
		\end{equation}
		\begin{equation}\label{bzekv12}
		\sum_{j=0}^{n-1} a_j(z) w_k^j =B(z,w_k) \exp( h(z,w_k)),\quad k=n+1,\ldots,2n
		\end{equation}
		and we solve for $a_{n-1}(z)$ from each of these to get
		\begin{equation}\label{bzekv13}
		a_{n-1}(z)=\sum_{j=0}^{n} \lambda_j B(z,w_j) \exp h(z,w_j)=,
		\sum_{j=0}^{n} \lambda_{n+j} B(z,w_{n+j}) \exp h(z,w_{n+j})
		\end{equation}
		where $\lambda_j$ are constants and $B(z,w_j)$ polynomials with respect to $z.$
		By Lemma \ref{borellemcor} 
		we have that
		if $H_j(z)$, $j=1,\ldots,m$ are entire functions such that
		\begin{equation}
		H_j(z)-H_l(z)\not\equiv \mbox{const.},\quad j\neq l
		\end{equation}
		and if $p_j(z)$ are polynomials then 
		\begin{equation}
		\sum_{j=1}^m p_j(z)\exp H_j(z)\equiv 0 \Rightarrow p_j(z)\equiv 0,\quad j=1,\ldots,m
		\end{equation}
		Applying this to Eqn.(\ref{bzekv13}) we obtain
		\begin{equation}
		\lambda_{2n} B(z,w_{2n})\equiv 0
		\end{equation}
		By Eqn.(\ref{bzekv12}) and Eqn.(\ref{bzekv13}) for $\lambda_{2n}\neq 0$, and
		denoting by \\$V(w_{n+1},w_{n+2},\ldots,w_{2n-1})$, $V(w_{n+1},w_{n+2},\ldots,w_{2n})$ the standard Vandermonde determinants
		we have
		\begin{equation}
		\lambda_{2n} =\frac{V(w_{n+1},w_{n+2},\ldots,w_{2n-1})}{V(w_{n+1},w_{n+2},\ldots,w_{2n})}
		\end{equation}
		Thus $B(z,w_{2n})\equiv 0$. Since this holds true for 
		for any $w_{2n}$ and any $z$ we have $B(\bar{z},z)\equiv 0$that is $\pi(z)\equiv 0$. Thus a contradiction to the assumption
		that $\pi(z)$ is a polyanalytic function of order $n\geq 2$ (see Eqn.(\ref{bzekv3})). This implies that there exists constants
		$c_1,\ldots,c_n$ such that
		\begin{equation}
		h(z,w_j)-h(z,w_1)\equiv c_j,\quad k=0,\ldots,n-1
		\end{equation}
			This proves Lemma \ref{visaclaimbzny}.
	\end{proof}
	By Lemma \ref{visaclaimbzny} we have that the function $h(z,w)$ in Eqn.(\ref{visaclaimbzekv6}) does not depend on $w.$
	From the system defined by Eqn.(\ref{bzekv11}) we have constants $\lambda_j^{(k)}$, $k=0,\ldots,n-1,$ such that
	\begin{equation}
	a_k(z)=\sum_{j=0}^{n} \lambda_j^{(k)} B(z,w_k) \exp( h(z,w_k)) =q_k(z)\exp\theta(z)
	\end{equation}
	where $\theta(z)\equiv h(z,w_1)$ and $q_k(z)$ polynomials in $z.$ By Eqn.(\ref{bzekv3})
	we obtain
	\begin{equation}\label{bzekv14}
	\Pi(z)=E(z)q(z,\bar{z})
	\end{equation}
	where 
	\begin{equation}
	E(z)=e(z)\exp(\theta(z))
	\end{equation}
	are entire functions and
	\begin{equation}
	q(z,\bar{z})=\sum_{j=0}^{n-1} q_j(z)\bar{z}^j
	\end{equation}
	an entire polyanalytic polynomial.
	Now $q(z,\bar{z})$ can be represented in the form
	\begin{equation}
	\varphi(x,y)+i\psi(x,y)
	\end{equation}
	for real valued polynomials $\varphi,\psi.$ Since by assumption $\Pi(z)$ has a nonisolated zero
	the equations (curves)
	$\varphi(x,y)=0$ and $\psi(x,y)=0$ have infinitely many common zeros.
	It is known (see Borel \cite{borel1934}, parag.75) that this is impossible if the
	polynomials $\varphi,\psi$ are relatively prime. 
	Hence there exists real valued polynomials 
	$\Phi(x,y),\Psi(x,y)$ and
	$s(x,y)\equiv V(z,\bar{z})$ for relatively prime $\Phi(x,y),\Psi(x,y)$ with
	\begin{equation}
	\phi(x,y)=s(x,y)\Phi(x,y)),\quad \Psi(x,y)=s(x,y)\Psi(x,y)
	\end{equation}
	But then
	\begin{equation}
	q(z,\bar{z})=(\Phi(x,y)+i\Psi(x,y))V(z,\bar{z})
	\end{equation}
	From Eqn.(\ref{bzekv14}) we obtain
	\begin{equation}
	\Pi(z)=G(z)V(z,\bar{z})
	\end{equation}
	where
	\begin{equation}
	G(z)=E(z)(\Phi(x,y)+i\Psi(x,y))
	\end{equation}
	is an entire polyanalytic function with only isolated zeros and $V(z,\bar{z})$ a real valued polynomial.
	Since $\Pi(z)$ is an irreducible polyanalytic function which has a nonisolated zero
	$G(z)$ is an entire (analytic) function.
	Thus the theorem is proved for the case of entire irreducible polyanalytic functions that have a nonisolated zero.
	This in turn immediately implies the theorem for arbitrary entire polyanalytic functions. Finally, a meromorphic entire polyanalytic function
	follows by factoring out an entire polyanalytic factor.
		This completes the proof.
\end{proof}

\begin{theorem}\label{bztheorem2}
	Let
	$f(z)=\sum_{j=0}^{n-1} a_j(z)\bar{z}^j$ for entire $a_j(z)$ and let $p_j(z)$,
	$j=1,\ldots,n$ be polynomials such that $p_i(z)\not\equiv p_k(z),$ $i\neq k,$
	and let $E_j$ be the set $\{z: f(z)-p_j(z)=0\}$.
	If each $E_j$, $j=1,\ldots,n$ has at least one accumulation point, then $f(z)$ is a polyanalytic polynomial.
\end{theorem}
\begin{proof}
	Define the polyentire functions
	\begin{equation}
	\varphi_j(z)=f(z)-p_j(z),\quad j=1,\ldots,n
	\end{equation}
	This can be represented as a finite product of polyentire functions where at least one of them, 
	which we denote $\theta_j(z)$, must have a nonisolated zero (since $\varphi_j(z)$ does).
	By Theorem \ref{bztheorem1} we have
	\begin{equation}
	\theta_j(z):= G_j(z)V_j(z,\bar{z}) 
	\end{equation}
	where $G_j$ is an entire function and $V_j$ a real valued polynomial.
	Clearly, $V_j(z,\bar{z})$ has a nonisolated zero and we denote its zeros
	set by $e_j(\subset E_j).$
	We may write
	\begin{equation}
	V_j(z,\bar{z}) =\sum_{k=0}^{m} q_k(z)\bar{z}^k
	\end{equation}
	for polynomials $q_k(z).$
	By Corollary \ref{numbervalueslem1} $e_j$ contains
	an analytic arc
	$\gamma_j:=\{\bar{z}=A_j(z)\}$ where $A_j$ is the Schwartz function of $\gamma_j$ defined and analytic 
	on a domain $\Delta_j.$
	On $\gamma_j$ (and therefore on $\Delta_j$)
	we have
	\begin{equation}
	V_j(z,A_j(z)) \equiv 0
	\end{equation}
	This implies that $A_j(z)$ can be analytically continued to each point of the $z$-plane except
	the finite set $\delta_j$ consisting of the zeros of the polynomial $q_m(z)$ and the discriminant set.
	We keep the same notation for the analytic continuation, $w=A_j(z)$, and this satisfies
	\begin{equation}\label{bzekv16}
	V_j(z,w) = 0
	\end{equation}
	Consider the associated function to $f(z)$
	\begin{equation}
	F(z,w)=\sum_{j=0}^{n-1} a_j(z)w^j
	\end{equation}
	On $\gamma_j$ we have
	\begin{equation}
	F(z,A_j(z)) =F(z,\bar{z})=f(z)=p_j(z)
	\end{equation}
	which implies
	\begin{equation}\label{bzekv17}
	F(z,A_j(z)) \equiv p_j(z),\quad j=1,\ldots,n
	\end{equation}
	If for some $j,k$, $A_j(z)\equiv A_k(z)$ for $j\neq k$ then by Eqn.(\ref{bzekv17})
	$p_j(z)\equiv p_k(z)$ which contradicts the conditions of the theorem.
	This ensures that all the
	$A_j(z),$ $j=1,\ldots,n$ are different.
	Let $D$ be a simply connected domain in
	$\C-\sum_{j=1}^n \delta_j$.
	We choose for each $A_j(z)$, one unique branch, $A_j^0(z)$, of $A_j(z)$, in $D.$
	The system of Eqn.(\ref{bzekv17}) 
	\begin{equation}
	\sum_{k=0}^{n-1} a_k(z)(A^0_j(z))^k =p_j(z),\quad j=,\ldots,n
	\end{equation}
	shows that each $a_j(z)$ is expressed rationally in terms of the $A_j^0(z)p_j(z).$
	Since the st of algebraic functions form a field each $a_j(z)$, $j=0,\ldots,n-1$,
	is algebraic. But by assumption they are entire, hence they must be polynomials.
	This completes the proof of Theorem \ref{bztheorem2}.
\end{proof}

Using a theorem on factorization of entire functions with bounded zero set, Theorem \ref{entirepathm} 
we are able to obtain the following result.

\begin{proposition}\label{bzlemma2}
	If a transcendental polyanalytic function $F(z)$ is divisble (in the sense
	of divisibility in the class of polyentire functions) by a polyanalytic polynomial of exact order $m>1$ then
	the set of points where $F(z)$ coincides with a nonzero polyanalytic polynomial of order smaller than $m$, is unbounded.
\end{proposition}
\begin{proof}
	Let 
	\begin{equation}\label{bzekv18}
	F(z)=p(z,\bar{z})\varphi(z,\bar{z})
	\end{equation} 
	for a polyanalytic polynomial $p(z,\bar{z})$ of order $m>0$
	and a polyentire $\varphi(z,\bar{z}).$ Let $q(z,\bar{z})$ be a polyanalytic polynomial of exact degree $s<m$ and suppose
	(in order to reach a contradiction) that
	$\{F=q\}$ is bounded. By Theorem \ref{entirepathm} (proved in Chapter \ref{polyentiresec}) there exists an entire $A(z)\not\equiv\mbox{const.}$
	and a polyanalytic polynomial such that
	\begin{equation}\label{bzekv19}
	F(z)-q(z,\bar{z})=h(z,\bar{z})\exp(A(z))
	\end{equation} 
	By Eqn.(\ref{bzekv18}) and Eqn.(\ref{bzekv19}) we have
	\begin{equation}\label{bzekv20}
	p(z,\bar{z})\varphi(z,\bar{z})-h(z,\bar{z})\exp(A(z))-q(z,\bar{z})\equiv 0
	\end{equation} 
	Equating the analytic components yields by uniqueness of analytic components
	\begin{equation}\label{bzekv20}
	\varphi(z,\bar{z})\equiv R(z,\bar{z})\exp(A(z))
	\end{equation} 
	For some polyanalytic polynomial $R(z,\bar{z})$. By Eqn.(\ref{bzekv20}) $q(z,\bar{z})$ is a transcendental polyanalytic 
	function which is a contradiction. This completes the proof.
\end{proof}

\begin{corollary}\label{bztheorem3}
	If the set $\{z:F(z)-P(z)=0\}$ for a polyentire $F(z)$ and an analytic polynomial $P(z)$ is bounded
	and if $\{z:F(z)-Q(z)=0\}$ for another analytic polynomial $Q(z)$, has a condensation point, then $F$ is a polyanalytic polynomial.
\end{corollary}
\begin{proof}
	By Theorem \ref{bztheorem1} $F(z):=f(z)-Q(z)$ is divisible by a polyanalytic polynomial of order $m>1.$
	By Proposition \ref{bzlemma2} applied to the function $F(z)$ have with $P(z)\equiv a$, for a complex constant $a$
	and $Q(z)\equiv b$
	for a complex constant $b$
	shows that $f(z)$ is a polyanalytic polynomial.
	This completes the proof.
\end{proof}
The theorem does not hold true in general if we merely require that $P,Q$ be polyanalytic polynomials, take e.g.\
$P(z)=\bar{z}z-1,$ $Q(z)=zP(z)$ and the polyanlaytic $f(z)=(\exp(z)+1)(\bar{z}z-1)$ which is clearly not a polyanalytic polynomial
but $f=Q$ contains the unit disc and $\{f=P\}$ belongs to $\{\abs{z}<2\}.$ 

\begin{corollary}\label{bztheorem4}
	Let $f(z)$ be a transcendental analytic function and let $H(z,\bar{z})$ be a polyanalytic non-analytic polynomial.
	The set $\{f(z)-H(z,\bar{z})=0\}$ is unbounded and consists only of isolated point.
\end{corollary}
\begin{proof}
	Under the assumption that the function $f(z)-H(z,\bar{z})$ has a bounded zero set or at least one 
	nonisolated zero at least one of the following equations hold true
	\begin{equation}\label{bzekv21}
	f(z)-H(z,\bar{z})\equiv P(z,\bar{z})\exp(g(z))
	\end{equation}
	\begin{equation}\label{bzekv22}
	f(z)-H(z,\bar{z})\equiv Q(z,\bar{z})\psi(z,\bar{z})
	\end{equation}
	for an entire analytic $g(z)$, polyanalytic polynomials $P(z),$ $Q(z)$ where the exact order of $Q$ is $>1$ and $\psi(z,\bar{z})$ a polyentire function.
	By uniqueness of analytic components we may equate the analytic components in 
	each of Eqn.(\ref{bzekv21}) and Eqn.(\ref{bzekv22}) the proof of
	Proposition \ref{bzlemma2} shows that in each case the equation is only possible if $f(z)$ is a polyanalytic polynomial.
	This completes the proof.
\end{proof}

\subsection{Petrov's growth condition}
Petrov \cite{petrov1966} proved the following.
\begin{theorem}\label{petrovthm1}
	Let $\Omega$ be a Jordan domain whose boundary $\partial\Omega$ contains a segment $L$ of a straight line or a circle
	given by the equation $\bar{z}=G(z),$ where $G(z)$ is the Schwarz function of $L.$ Let $f(z)$ be an $n$-analytic function
	on $\Omega$ for some $n\in \Z_+.$
	If the function
	\begin{equation}\label{petrovekv1}
	\frac{f(z)}{(\bar{z}-G(z))^{n-1}}
	\end{equation}
	has boundary values that vanish on $L$ then $f\equiv 0.$
\end{theorem}
\begin{proof}
	First consider the case where $L$ is a segment of $\{x=0\}$, where $z=x+iy.$ Then the function in Eqn.(\ref{petrovekv1}) takes the form
	\begin{equation}\label{holoulu}
	\frac{f(z)}{x^{n-1}}
	\end{equation}
	We use induction in the order of analyticity $n.$ The case $n=1$ is a direct corollary to the Luzin-Privalov theorem for holomorphic functions.
	Suppose $n>1$ and that the theorem holds true for $(n-1)$-analytic functions.
	Now on $\Omega$ there are analytic components $a_j(z),$ $j=0,\ldots,q-1$ such that
	\begin{equation}
	f(z)=\bar{z}^{n-1} a_{n-1}(z)+\sum_{j=0}^{n-2} a_j(z)\bar{z}^j=(z-(z+\bar{z}))^{n-1}a_{n-1}(z)+\sum_{j=0}^{n-2} a_j(z)\bar{z}^j 
	\end{equation}
	thus, since $(z-(z+\bar{z}))^{n-1}=(z+\bar{z})^{n-1}+\phi(z,\bar{z})$ for a function $\phi$ polyanalytic of order $n-1$, we have that 
	$f(z)$ can be represented on $\Omega$ as
	\begin{equation}\label{petrovekv2}
	f(z)=u(z)+x^{n-1}h(z)
	\end{equation}
	where $h(z)$ is holomorphic, $u(z)$ polyanalytic of order $q-1$. Set
	\begin{equation}\label{petrovekv3}
	u(z)=u_0(z)+x^{n-1}u_1(z),\quad h(z)=h_0(z)+x^{n-1}h_1(z)
	\end{equation}
	for polyharmonic $u_0,u_1,$ and harmonic $h_0,h_1$ so that
	\begin{equation}\label{petrovekv3}
	\re f(z)=u_0(z)+x^{n-1}h_0(z),\quad \im f(z)=u_1(z)+x^{n-1}h_1(z)
	\end{equation}
	By the conditions of the theorem the functions $\frac{\re f(z)}{x^{n-1}}$ and $\frac{\im f(z)}{x^{n-1}}$ have zero boundary values
	on $L$.
	This in turn implies (by a known result of Huber \cite{huber})
	that the functions $\frac{u_0(z)}{x^{n-1}}$ and $\frac{\im u_1(z)}{x^{n-1}}$ have zero boundary values
	on $L$.
	Consequently the same holds true for the $n-1$-analytic function $u(z)$. By the induction hypothesis $u\equiv 0$ on $\Omega$.
		Thus on $\Omega$
	\begin{equation}
	f(z)=x^{n-1}h(z)
	\end{equation}
	for a holomorphic $h(z)$. But this implies that by the condition that the expression in Eqn.(\ref{holoulu}) take zero boundary values on $L$
	that also $h\equiv 0.$ This shows that $f\equiv 0.$
	This proves the theorem for the case when $L$ is a segment of the imaginary axis. 
	\\
	Now let $L$ be an arbitrary segment of the curve $px+ql+l=0.$ We may assume that $L$ is not parallel to the real axis (otherwise the proof is simplyfied).
	Then the expression of Eqn.(\ref{petrovekv1}) becomes
	\begin{equation}
	\frac{f(z)}{(px+qy+l)^{n-1}}
	\end{equation}
	We may express $L$ in complex form with $z=a+\exp(i\alpha) t,$ $a=-\frac{l}{p},$ $\alpha=\mbox{atan}\left(-\frac{p}{q}\right),$ $t\in (-\infty,\infty)$
	and the transformation of the plane given by
	\begin{equation}
	w=(z-a)\exp\left(t\left(\frac{\pi}{2}-\alpha\right)\right)
	\end{equation}
	we may assume that $L$ is transformed to a segment of $t\exp\left(i\frac{\pi}{2}\right),$ $\re w=0,$ and
	$\Omega$ is transformed to some domain $\Omega'$ whose boundary contains a segment $L'$ of the imaginary axis and
	the function
	$\varphi(w)=f(w^{-1})$ is $n$-analytic on $\Omega'.$ Since
	\begin{equation}
	\re w=\frac{px+qy+l}{-\sqrt{p^2+q^2}}
	\end{equation}
	the function 
	\begin{equation}
	\frac{\varphi(w)}{(\re w)^{n-1}}
	\end{equation}
	has zero boundary values on $L'$ and consequently $\varphi(w)\equiv 0$ on $\Omega'$, by the case we have already proved, 
	which in turn implies $f(z)\equiv 0$ on $\Omega.$ This proves the theorem for the case when $L$ is a line segment.
	\\
	Now suppose $\partial\Omega$ contains an arc in some circle $\{\abs{z-a}=R\},$ $R>0$ and that $L$ is a segment of that arc.
	The condition of Eqn.(\ref{petrovekv1}) becomes that
	\begin{equation}\label{petrovprimekv}
	\frac{f(z)}{(R-\abs{z-a})^{n-1}}
	\end{equation} 
	has zero angular boundary values on $L$, and we may assume w.l.o.g.\ that $L$ belongs to
	the circle
	$\abs{z-1/2}=1/2.$ Set $z'=\frac{1}{\bar{z}}$ and
	$\rho:=\abs{z-a},$ $z'=r'\exp(i\alpha')=x'+iy$,
	$z=r\exp(i\alpha).$
	This transformation send $L$ into a segment of the line $\re z'=1$, and $\Omega$ to a domain $\omega'$ symmetric with respect to the unit circle
	such that the boundary $\partial\Omega'$ contains a segment of $x'-1=0.$ Define the $n$-analytic function on $\Omega'$
	\begin{equation}
	F(z'):=(\bar{z}')^{n-1}\cdot f\left(\frac{1}{z'}\right)
	\end{equation}
	We can calculate that
	\begin{equation}
	x'-1=r'^2\left(\frac{1}{4}-\rho^2\right)
	\end{equation}
	\begin{equation}
	\frac{F(z')}{(x'-1)^{n-1}}=\left(\frac{r^2}{z}\right)^{n-1}\overline{\left(\frac{f(z)}{\left(\frac{1}{4}-\rho^2\right)^{n-1}}\right)}
	\end{equation}
	By the condition given by Eqn.(\ref{petrovprimekv}) it follows that $\frac{F(z')}{(x'-1)^{n-1}}$
	vanishes on $L'$, thus by the previous case of the theorem we have that $F(z')\equiv 0$ on $\Omega'$. This renders $f(z)\equiv 0$.
	This completes the proof of Theorem \ref{petrovthm1}.
	\end{proof}
\section{A condition for unbounded zero sets}
We shall need the following.
\begin{theorem}\label{narasimzerooflimithm00}
	Let $\Omega\subset \Cn$ be a domain. Let $\{\phi_j\}_{j\in \N}$ be a sequence of continuous open maps  
	on $\Omega$
	converging uniformly on compacts to a function $\phi(z).$ Suppose that for some $p_0\in \Omega$ $p_0$ is an isolated point
	of $\phi^{-1}\phi(p_0)$. Then for any neighborhood $U$ of $p_0$ there exists $j_0$ such that $\phi(p_0)\in \phi_j(U)$ for $j\geq j_0.$
\end{theorem}
\begin{proof}
	Assume not. If necessary after passing to a subsequence we can suppose $U\Subset\Omega$, $\phi(p_0)\notin \phi(\partial U)$.
	Hence there is a neighborhood $V$ of $\phi(\partial U)$ and a polydisc $P$ centered at $\phi(p_0)$ such that $P\cap V=\emptyset.$
	If $j_0$ is sufficiently large $\phi_j(\partial U)\subset V$ for $j\geq j_0.$ Since $\phi_j$ is open we have 
	$\partial \phi_j(U)\subset \phi_j(\partial U)$ (here $\partial \phi_j(U)$ denotes the boundary of the set $\phi_j(U)$). Hence $\phi_j(U)\Subset\Cn$
	with $\partial\phi_j (U)\subset V$. We show that for sufficiently large $j$ we have
	$(\partial\phi_j(U))\cap P\neq\emptyset$, which will render a contradiction. To see this note that since $\phi(p_0)\in P$ and $\phi_j(p_0)\to \phi(p_0)$ we have that if $j$ is sufficiently large $\phi_j(p_0)\in P$. However $\phi(p_0)\notin \phi_j(U)$ by assumption. If $(\partial \phi_j(U))\cap P=\emptyset$
	then we would have
	$P=(\phi_j(U)\cap P)\cup ((\Cn \setminus \overline{\phi_j(U)})\cap P)$ and each of these open sets
		is nonempty since $\phi_j(p_0)$ belongs to the first and $\phi(p_0)$ belongs to the second.
		This would contradict the fact that $P$ is connected. Hence
		$(\partial\phi_j(U))\cap P\neq \emptyset.$ This completes the proof.
	\end{proof}
	
	\begin{theorem}[Hurwitz]\label{narasimzerooflimithm0}
		Let $\Omega\subset \Cn$ be a domain and let $\{f_j\}_{j\in \N}$ be a sequence of holomorphic functions on $\Omega$
		converging uniformly on compacts to a function $f(z).$ If $f_j(z)\neq 0$ for all $j\in \N$ and all $z\in \Omega$ then
		$f(z)\neq 0$ for all $z\in \Omega$
		\end{theorem}
		\begin{proof}
			Suppose $f(a)=0$ for some $p_0\in \Omega$ and let $\delta>0$ such that $V:=\{\abs{z-p_0}<\delta\}\subset \Omega.$
			We can assume $f\not\equiv 0$ (we already know $f$ is holomorphic).
			Let $p_1\in V$ such that $f(p_1)\neq 0.$
			Set $U:=\{\lambda\in \C:p_0+\lambda(p_1-p_0)\in V\}.$ Then $U$ is a convex, hence connected open subset. Let $\phi_j(z):=f_j(p_0+\lambda(p_1-p_0)),$
			$\phi(\lambda)=f(p_0+\lambda (p_1-p_0)).$
			Then $\phi(0)=0,$ $\phi(1)=f(b)\neq 0$ hence $\phi$ is nonconstant on $U$ thus for sufficiently large $j$ also each $\phi_j(z)$ is nonconstant, 
			thus an open map of $U$ into $\C.$
			This implies by Theorem \ref{narasimzerooflimithm00} that for sufficiently large $j$ we have
			$\{0\}\subset \phi(\Omega)\subset f(\Omega)$ which is a contradiction.
			This completes the proof.
		\end{proof}
		
		\begin{corollary}\label{narasimzerooflimithm}
			Let $\Omega\subset \Cn$ be a domain and let $\{f_j\}_{j\in \N}$ be a sequence of holomorphic functions on $\Omega$
			converging uniformly on compacts to a function $f(z).$ Suppose $f(z)$ has a zero in $\Omega.$ Then there exists
			a subsequence
			$\{f_{j_l}\}_{l\in \N}$ of $\{f_j\}_{j\in \N}$ such that $f_{j_l}\to f$ uniformly on compacts and such that
			each $f_{j_l}$ has a zero in $\Omega.$
		\end{corollary}
		\begin{proof}
			The contrapositive to Theorem \ref{narasimzerooflimithm0} is that if $f(z)$ has a zero in $\Omega$ then 
			there exists at least one $f_{j}(z)$ which also has a zero in $\Omega.$ Let $S:=\{f_{j_l}\}_l$ be the set of
				members of $\{f_j\}_{j\in \N}$ that each have a zero in $\Omega.$ Then $S$ must contain infinitely many members
				for otherwise the set $\{f_j:f_j\notin S\}$ would contain a subsequence converging uniformly on compacts of $\Omega$ to $f(z)$ such that each member of the subsequence 
				would be zero free on $\Omega$ (a contradiction to the contrapositive of Theorem \ref{narasimzerooflimithm0}).
				But if $S$ is an infinite sequence (which was obtained from a sequence converging uniformly on compacts of $\Omega$ to $f(z)$) then it contains a
				subsequence converging uniformly on compacts to a function $f(z).$ This completes the proof.
			\end{proof}

Balk \cite{balkbianalpicard65}, Lemma 1, proved
the following for $2$-analytic functions
and quite possible the same proof can be tweeked to hold for arbitrary order. We shall later present a different proof for polyentire functions of arbitrary finite order for the case
of the zero value.
\begin{proposition}\label{balkunboundedzerolem0}
	Suppose $f(z)=a_0(z) +\bar{z}a_1(z)$ for entire $a_0(z),a_1(z)$. If $a_1^{-1}(0):=\{z:a_1(z)=0\}$ is unbounded
	then for all $\alpha\in \C$ the set $f^{-1}(\alpha):=\{z:f(z)=\alpha\}$ is unbounded.
\end{proposition}
\begin{proof}
	Note that (since the zeros of an entire $\not\equiv 0$) function are isolated an entire ($\not\equiv 0$) function has unbounded zero set iff it has an infinite set of zeros.
	We show that  if $a_1(z)\not\equiv 0$ has infinitely many zeros then
	for all $\alpha\in\C$ the set $f^{-1}(\alpha)$ is unbounded.
	Let
	\begin{equation}
	g(z)=z(f(z)-\alpha)=:\Phi +\abs{z}^2a_1(z),\quad \Phi(z)=z(a_0(z)-\alpha)
	\end{equation}
	We shall show that $g^{-1}(0)$ is unbounded.
	Assume (in order to reach a contradiction) that this is false. Then there exists $d>0$ such that all the zeros of $g$ lie
	inside the oriented contour $\gamma:=\{\abs{z}=d\}$ and
	\begin{equation}
	\frac{1}{2\pi}\Delta_\gamma \mbox{arg} g(z)=m<\infty
	\end{equation}
	where $\frac{1}{2\pi}\Delta_\gamma \mbox{arg} g(z)$ denotes the change in argument of $g(z)$ as $z$ traverses $\gamma$ in the positive direction once,
	and the value is finite since $g(z)$ is continuous on the bounded disc defined by $\gamma$ and for appropriate choice of $d$ is nonzero on $\gamma.$
	By the condition of the proposition there exist $c>d$ such that for the oriented $\Gamma:=\{\abs{z}=c\}$, $a_1(z)$ has $s>m$ zeros inside $\Gamma$ and $a_1\neq 0$ on $\Gamma.$ Since the continuous function $f(z)$ has no zeros
	on $d\leq\abs{z}\leq c$ we have (see e.g.\ Szilard \cite{szilard})
	\begin{equation}
	\frac{1}{2\pi}\Delta_\Gamma \mbox{arg} g(z)=m
	\end{equation}
	Set for $n\in \Z_+$, $F_n(z):=\Phi(z)+n^2 a_1(z)$ and
	$\Gamma_n:=\{\abs{z}=n\}.$ Then $F_n|_{\Gamma_n}=g|_{\Gamma_n}$ and for $n>c$ we have
	\begin{equation}\label{balkhot005}
	\frac{1}{2\pi}\Delta_\Gamma \mbox{arg} F_n(z)=m
	\end{equation}
	Let $D'$ denote the set obtained by removing the zeros of $a_1(z)$ from $\{\abs{z}< c\}$.
	Then $\{F_n\}_{n}$ converges uniformly to $\infty$ on $D'$ thus $\{F_n\}_{n}$ is quasinormal 
	and the irregular points of the sequence will be the set of zeros of $a_1(z).$
	We surround these zeros by pairwise nonintersecting circles that also do not intersect $\Gamma.$
	By a known theorem of Montel (see Theorem \ref{a4prop}) 
	there exists $N$ such that for $n\geq N$, each $F_n(z)$ will take the value $0$ in each of the discs enclosed by a circle as above, thus $F_n$ has at least $s$ zeros inside $\Gamma$. This implies that
	\begin{equation}\label{balkhot006}
	\frac{1}{2\pi}\Delta_\Gamma \mbox{arg} F_n(z)\geq s>n
	\end{equation}
	which for $n>\max\{c,N\}$ contradicts Eqn.(\ref{balkhot005}).
	This proves that the zeros of $g$ must be unbounded. This in turn implies that
	$f^{-1}(\alpha)$ is also unbounded. This completes the proof.
	\end{proof}
For the case of the value $\alpha=0$ we generalize this to arbitrary polyentire functions.
\begin{proposition}\label{balkunboundedzerolem}
	Suppose $f(z)=\sum_{j=0}^{q-1} a_j(z)\bar{z}^j$ for entire $a_j(z),$ $j=0,\ldots,q-1$. 
	If $a_{q-1}^{-1}(0):=\{z:a_{q-1}(z)=0\}$ is unbounded
	then $f^{-1}(0):=\{z:f(z)=0\}$ is unbounded.
\end{proposition}
\begin{proof}
	Let $\{z_j\}_{j\in \Z_+}\subset a_{q-1}^{-1}(0)$ such that
	$c_j:=\abs{z_j}\to \infty$ as $j\to \infty$ and let each fixed $\epsilon>0,$
	$\{r_{j,\epsilon}\}_{j\in \Z_+}$ be a sequence of positive real numbers such that
	the sets $U_{j,\epsilon}:=\{ \abs{z-a_j}<r_{j,\epsilon}\}$ are pairwise disjoint and define 
	the functions
	\begin{equation}
	g_j(z) :=\frac{1}{c_{j,\epsilon}^{2k}}\left(\sum_{k=0}^{q-1} \frac{a_k(z)}{z^k} c^{2k}_{j,\epsilon}\right)=:\frac{1}{c_{j,\epsilon}^{2k}}\hat{f}_{j,\epsilon}(z)
	\end{equation}
	By Theorem \ref{narasimzerooflimithm} (recall that $a_{q-1}(z)$ has a zero, $z_{j_0}$, in $U_{j_0,\epsilon}$) there exists, for each $j_0\in \Z_+$ such that $z_{j_0}\neq 0$
	a subsequence $\{g_{j_l}\}_{l\in \Z_+}$ of $\{g_j\}_{j\in \Z}$ such that
	\begin{equation}
	g_{j_l}\to \frac{a_{q-1}(z)}{z^{q-1}},\mbox{ uniformly on }U_{j_0,\epsilon}
	\end{equation}
	satisfying that each $g_{j_l}(z)$ has a zero in $U_{j_0,\epsilon}.$
	Note that if we replace $\epsilon$ by a smaller value $\epsilon'$ then the same functions $g_j$ and $g_{j_l}$ can be used in our arguments.
	In particular, we can consider limits with respect to $\epsilon\to 0$ because
	the $g_{j_l}$ will be independent of $\epsilon.$
	Now there exists an $\epsilon_0>0$ such that if $g_{j_l}$ has a zero in $U_{j_0,\epsilon}$ then obviously
	we have for $\epsilon<\epsilon_0$ (and assuming $z_{j_0}\neq 0$) that
	$\hat{f}_{j,\epsilon}|_{U_{j_0,\epsilon}}$ also has a zero in  $U_{j_0,\epsilon}$ for all $\epsilon<\epsilon_0.$
	Furthermore, for $k=0,\dots,q-1$ we have
	\begin{equation}
	\lim_{\epsilon\to 0} \frac{c^{2k}_{j,\epsilon}}{z^k}=\lim_{\epsilon\to 0}\frac{\bar{z}^k z^k}{z^k}=\bar{z}_{j_0}^k
	\end{equation}
	hence
	\begin{equation}
	\lim_{\epsilon\to 0} \hat{f}_\epsilon(z)=\sum_{k=0}^{q-1} a_k(z_{j_0})\bar{z}_{j_0}^k=f(z_{j_0})
	\end{equation}
	This implies by continuity that
	$z_{j_0}\in f^{-1}(0).$ Since $z_{j_0}$ was arbitrary in the unbounded set $\{z_j\}_{j\in \Z_+}\setminus \{0\}$ this completes the proof.
\end{proof}

\section{Polyanalytic functions on the unit disc}
\begin{definition}
	Let $z$ denote the complex coordinate in $\C$ with polar representation $z=r\exp(i\theta).$
	Denote by $\Pi(\{\abs{z}<1\})$ the set of holomorphic functions $f(z)$ on $\{\abs{z}<1\}$ that 
	satisfy
	\begin{equation}
	T(r,\pi):=O\left(\ln \frac{1}{1-r}\right)\mbox{ as }r\to \infty
	\end{equation}
	Suppose 
	$\Omega(z)=(a_0(z),\ldots,a_q(z))$ is a holomorphic curve on $\{\abs{z}<1\}$ with $a_q(z)\not\equiv 0$.
	We define
	\begin{equation}
	\norm{\Omega(z)}:=\sum_{j=0}^q \abs{a_j(z)},\quad T(r,\Omega)=\frac{1}{2\pi}\int_{0}^{2\pi} \ln \norm{\Omega(r\exp(i\theta)}d\theta
	\end{equation}
\end{definition}

By the fundamental theorem of algebra any polynomial $P(z)$ can be written
\begin{equation}
P(z)=cz^p\Pi_{r=p+1}^n \left(1-\frac{z}{z_r}\right)
\end{equation}
where the $z_r$ are the zeros of $P(z)$ outside the origin.
There exists a generalization of this representation to certain meromorphic functions.
\begin{definition}
	For $q\in \Z_+$ define the {\em primary Weierstrass multiplier}\index{Primary Weierstrass multiplier} as
	\begin{equation}
	E(z,q):=(1-z)\exp(z+\frac{1}{2}z^2+\cdots +\frac{1}{q}z^q),\quad E(z,0)=1-z
	\end{equation}
\end{definition}
\begin{theorem}[See Hayman \cite{hayman}, Thm. 1.9]
	If $f(z)$ is meromorphic in $\C$ and has zeros $a_\mu,$ and poles $b_\nu$ such that it satisfies one of the conditions
	\begin{equation}
	(a)\qquad \limsup_{R\to \infty} \frac{T(R,f)}{R^q}=0
	\end{equation}
	\begin{equation}
	(b)\qquad \liminf_{R\to \infty} \frac{T(R,f)}{R^q}=0
	\end{equation}
	then
	\begin{equation}
	f(z)=z^p\exp(P_{q-1}(z))\left(
	\frac{\Pi_{\abs{a_\mu}<R} E\left(\frac{z}{a_\mu},q-1\right)}{\Pi_{\abs{b_\mu}<R} E\left(\frac{z}{b_\nu},q-1\right)}
	\right)
	\end{equation}
	where the limit as $R\to \infty$ is taken through all values in case (a) and through a suitable sequence of values in case (b). Here $p\in \Z_+$ and
	$P_{q-1}(z)$ is a polynomial of degree at most $q-1.$
\end{theorem}
\begin{proof}
Suppose first that $f(0)\neq 0$ and $f(0)\neq \infty.$ By the Poisson-Jensen formula (Theorem \ref{goldbergnevanthm21}) we have with $z=r\exp(i\theta)$
\begin{multline}
\ln \abs{f(z)} =
\frac{1}{2\pi}\int_0^{2\pi} \ln\abs{f(R\exp(i\theta))} \frac{(R^2-r^2)}{R^2-2Rr\cos(\phi-\theta)+r^2}d\phi
+\\
\sum_{\abs{b_\nu}<R}\ln \abs{\frac{R^2-\bar{b}_{\nu}z}{R(z-b_{\nu})}}+
\sum_{\abs{a_\mu}<R}\ln \abs{ \frac{R(z-a_{\mu})}{R^2-\bar{a}_{\mu}z)}}
\end{multline}
Both sides are equal to a harmonic function, say $v(z$, near any $z=r\exp(i\theta)$ where $f(r\exp(i\theta))\notin \{0,\infty\}.$
We may assume that there are no poles of $f$ on $\abs{z}=R.$
Applying $\partial_x-i\partial_y$ to both sides, differentiating under the integral sign and using
\begin{equation}
\re\left(\frac{R\exp(i\phi)+z}{R\exp(i\phi)-z}\right)=\frac{(R^2-r^2)}{R^2-2Rr\cos(\phi-\theta)+r^2}
\end{equation}
yields 
\begin{multline}
\frac{f'(z)}{f(z)} =
\frac{1}{2\pi}\int_0^{2\pi} \ln\abs{f(R\exp(i\theta))} \frac{2R\exp(i\phi)}{(R\exp(i\phi)-z)^2}d\phi
+\\
\sum_{\abs{b_\nu}<R}\left(\frac{1}{b_{\nu}-z}-\frac{\bar{b}_\nu}{R^2-\bar{b}_\nu z}\right)-
\sum_{\abs{a_\mu}<R}\left(\frac{1}{a_{\mu}-z}-\frac{\bar{a}_\mu}{R^2-\bar{a}_\mu z}\right)
\end{multline}
Differentiating $q-1$ times we obtain\begin{multline}\label{hayman191ekv117}
\frac{d^{q-1}}{dz^{q-1}}\left(\frac{f'(z)}{f(z)}\right)=\frac{q!}{\pi}
\int_0^{2\pi}\frac{\ln \abs{f(R\exp(i\theta))} R\exp(i\theta))d\theta}{(R\exp(i\theta)-z)^{q+1}}+\\
(q-1)!\sum_{\abs{b_\nu}<R}\left(\frac{1}{(b_\nu -z)^q} -\frac{\bar{b}_\nu^q}{(R^2-\bar{b}_\nu z)^q}\right)-\\
(q-1)!
\sum_{\abs{b_\nu}<R}\left(\frac{1}{(a_\mu -z)^q} -\frac{\bar{a}_\mu^q}{(R^2-\bar{a}_\mu z)^q}\right)
\end{multline}
Suppose $T(2r)/r^q\to 0$ as $r\to \infty$ either through all values or through a suitable sequence of values, say $R_k$ (which tends to $\infty$ with $k$), such a sequence
exists by our assumptions. We may choose $\{R_k\}_{k\in \N}$ such that $f(z)\neq 0$ on $\abs{z}=R_k.$ Set
$R=R_k$ in Eqn.(\ref{hayman191ekv117}) we have, using $m(r,f)\geq 0$
\begin{equation}
T(2R_k,f)\geq N(2R_k,f)\geq \int_{R_k}^{2R_k} \frac{n(t,f)}{t}dt\geq n(R_k,f)\ln 2
\end{equation}
Note that with the notation 
\begin{equation}
N(R,f)=\sum \ln\abs{\frac{R}{b_\nu}}=\int_0^R n(t,f)\frac{dt}{t}
\end{equation}
where $n(t,f)$ is the number of poles of $f$ in $\abs{z}<t$
the (Jensen) formula
\begin{equation}
\ln \abs{f(0)}=\frac{1}{2\pi}\int_0^{2\pi}\ln \abs{f(R\exp(i\theta))}d\theta +\sum_{\mu=1}^M\ln \frac{\abs{a_\mu}}{R}
-\sum_{\nu=1}^N\ln \frac{\abs{b_\nu}}{R}
\end{equation}
where $N$ are the number of zeros and $M$ the number of poles of $f$ in $\abs{z}<R$
becomes
\begin{equation}
\ln \abs{f(0)}=m(R,f)-m(R,\frac{1}{f})+N(R,f)-N(R,\frac{1}{f})
\end{equation}
so that
\begin{equation}
m(R,f)+N(R,f)=\ln \abs{f(0)}+m(R,\frac{1}{f})+N(R,\frac{1}{f})
\end{equation}
So that the characteristic function $T(R,f)=m(R,f)+N(R,f)$ satisfies equals $T(R,\frac{1}{f})+\ln\abs{f(0)}$ thus
$T\left(R_k,\frac{1}{f}\right)=T(R_k,f)+O(1)$. This implies that
\begin{equation}
\frac{n(R_k,f)}{R^q_k}\to 0,\quad \frac{n\left(R_k,\frac{1}{f}\right)}{R^q_k}\to 0,\mbox{ as }k\to \infty
\end{equation}
Suppose $\abs{z}<\frac{R_k}{2}$. Then $\abs{\bar{b}_\nu z}<\frac{R_k^2}{2}$ for 
$\abs{b_\nu}<R_k$
and $\abs{R^2_k-\bar{b}z}\geq R_k^2-\abs{\bar{b}z}>\frac{R^2_k}{2}.$ This implies
\begin{equation}
\abs{\frac{\bar{b}_\nu^q}{(R^2_k -\bar{b}_\nu z)^q}}<\frac{R^q_k}{\left(\frac{R^2_k}{2}\right)^2}=\frac{2^q}{R^q_k}
\end{equation}
for all poles $b_\nu$ such that $\abs{b_\nu}<R_k.$ Summation over such poles yields
\begin{equation}
\abs{\sum_{\abs{b_\nu}<R_k} \frac{\bar{b}_\nu^q}{(R^2_k -\bar{b}_\nu z)^q}}\leq
\sum_{\abs{b_\nu}<R_k} 
\abs{\frac{\bar{b}_\nu^q}{(R^2_k -\bar{b}_\nu z)^q}}<\frac{2^q n(R_k,f)}{R_k^q}
\end{equation}
where the right hand side goes to zero uniformly on bounded subsets as $k\to\infty.$
The analogous arguments render that the summation of the terms corresponding to the
zeros (i.e.\ corresponding to the $a_\mu$) of the function go to zero on bounded subsets as $k\to \infty.$
Now
\begin{equation}
\end{equation}
The modulus of the right hand side of Eqn.(\ref{hayman191ekv117}) is bounded by
\begin{multline}
q!\frac{R_k}{\pi}\frac{2^{q+1}}{R^{q+1}_k}\int_0^{2\pi}\abs{\ln \abs{f(R_k\exp(i\theta)}}d\theta=\\
\frac{q!}{\pi}\frac{2^{q+1}}{R_k^q}\left(\int_0^{2\pi}\ln^+\abs{f(R_k\exp(i\theta))}d\theta +
\int_0^{2\pi}\ln^+\abs{\frac{1}{f(R_k\exp(i\theta))}}d\theta\right)=\\
\frac{q!}{\pi}\frac{2^{q+1}}{R_k^q}(m(R_k,f)+m(R_k1/f))=O\left(\frac{T(R,f)}{R_k^q}\right)
=\\
O\left(\frac{T(2R,f)}{R_k^q}\right)\to 0,\mbox{ as }k\to \infty
\end{multline}
Thus Eqn.(\ref{hayman191ekv117}) can be written
\begin{equation}
\frac{d^{q-1}}{dz^{q-1}} \frac{f'(z)}{f(z)}=\lim_{k\to \infty} S_k(z)
\end{equation}
where
\begin{equation}
S_k(z)=(q-1)!\left(\sum_{\abs{b_\nu}<R_k}\frac{1}{(b_\nu-z)^q}- \sum_{\abs{a_\nu}<R_k}\frac{1}{(a_\nu-z)^q}\right)
\end{equation}
and where the convergence is uniform for any bounded set of values of $z$ that does not contain
zeros or poles of $f(z).$
By uniform convergence we may integrate both sides $(q-1)$ times along an appropriate path from $0$ to $z$ in order to obtain
\begin{multline}
\frac{f'(z)}{f(z)}=
\lim_{k\to \infty}\left(\sum_{\abs{b_\nu}<R_k}\left(\frac{1}{b_\nu-z}-\frac{1}{b_\nu}-\frac{z}{b_\nu^2}-\cdots\right.\right.\\
\left.\left. -\frac{z^{q-2}}{b^{q-1}_\nu}\right)-
\left(\frac{1}{a_\nu-z}-\frac{1}{a_\nu}-\frac{z}{a_\nu^2}-\cdots-\frac{z^{q-2}}{a^{q-1}_\nu}\right)\right)
+P_{q-2}(z)
\end{multline}
where $P_{q-2}(z)$ is a polynomial of degree at most $q-2$.
Integrating both sides from $0$ to $z$ and taking exponentials proves the theorem for the case when $f(0)\neq 0$, $f(0)\neq \infty.$ If $0$ is a pole of $f$ of order $p$ we may apply this to the function $f(z)/z^p$ to once again obtain the wanted result. This completes the proof. 	
\end{proof}
We shall need the following result due to Goldberg \cite{goldbergconvex}.
\begin{proposition}
	\label{goldberghay}
	Let $k=1,\ldots,q$, let $y_k(x)$ be even, bounded, with period $2\pi$, and monotone increasing for $0\leq x\leq \pi,$ and let $\Phi(y)$
	be logarithmically convex ($0<y<\infty$). Then for arbitrary real numbers $x_k$, $j=1,\ldots,q,$ we have
	\begin{equation}
	\int_{-\pi}^\pi \Phi\left(\Pi_{k=1}^q y_k(x+x_k)\right)dx\leq \int_{-\pi}^\pi \Phi\left(\Pi_{k=1}^q y_k(x)\right)dx
	\end{equation}
\end{proposition}
\begin{proof}
First of all note that if $\psi(s)$ is a convex function of $s\in \R$
 since for 
$s_1<s_2<s_3$ we have
\begin{equation}
(t_3-t_1)\psi(t_2)\leq (t_2-t_1)\psi(t_3)+(t_3-t_2)\psi(t_1)
\end{equation}	
setting $t_1=t-h_2,$ $t_3=t+h_2$ and $t_2=t+ h_1$ ($t_2=t-h_1$ respectively) where $0\leq h_1<h_2$
gives
\begin{equation}\label{qqqqqq777}
2h_2\psi(t+h_1)\leq  (h_1+h_2)\psi(t+h_2) +(h_2-h_1)\psi(t-h_2) \end{equation}
and for the case $t_2=t-h_1$
\begin{equation}\label{qqqqqq777a}
2h_2\psi(t-h_1)\leq (h_2-h_1)\psi(t+h_2) +(h_2+h_1)\psi(t-h_2) \end{equation}
Adding Eqn.(\ref{qqqqqq777}) and Eqn.(\ref{qqqqqq777a}) yields
\begin{equation}
2h_2(\psi(t+h_1)+\psi(t-h_1))\leq 2h_2(\psi(t+h_2)+\psi(t-h_2) \end{equation}
thus $\psi(s+h)+\psi(s-h)$ is increasing in $h$ for $h\geq 0$
Since $\Phi(y)$ is convex in $\log y$ we have
setting $a=\exp(t_1-h_1), A=\exp(t_1+h_1),$
$b=\exp(t_2-h_2), B=\exp(t_2+h_2)$, and setting $\psi(t)=\Phi(\exp(t))$, the above 
deduction 
yields that
\begin{multline}
\psi(t_1+t_2+h_2-h_1)+\psi(t_1+t_2+h_1-h_2)\leq\\
\psi(t_1+t_2+h_2+h_1)+\psi(t_1+t_2-h_1-h_2)
\end{multline}
We conclude that
for $0<a<A<\infty,$ $0<b\leq B<\infty$
\begin{equation}\label{haymanlem49ekv}
\Phi(Ab)+\Phi(aB)\leq \Phi(ab)+\Phi(AB)
\end{equation}
Now since the $y_k$ are periodic with period $2\pi$ we can w.l.o.g.\ assume $0<x_k<2\pi$, $k=1,\ldots,q.$
Suppose $x_\mu>0$. Then the periodic property of $y_j(x)$ yields
\begin{multline}
\int_{-\pi}^\pi \Phi\left(\Pi_{j=1}^q y_j(x+x_j)\right)dx =
\int_{-\pi-\frac{x_\mu}{2}}^{\pi-\frac{x_\mu}{2}} \Phi\left(\Pi_{j=1}^q y_j(x+x_j)\right)dx=\\
\int_{-\pi}^\pi \Phi\left(\Pi_{j=1}^q y_j(x+k_j)\right)dx
\end{multline}
where $k_j=x_j-\frac{1}{2} x_\mu$(mod $2$).
Now 
\begin{equation}
\cos(x-k_j)-\cos(x+k_j)=2\sin(x)\sin(k_j)
\end{equation}
so $\cos(x-k_j)<\cos(x+k_j)$ for $\pi<k_j$ and
$\cos(x-k_j)>\cos(x+k_j)$ for $\pi>k_j$.
Since the $y_j(x)$ are decreasing with respect to $\cos x$
we have that if $0\leq x\leq \pi$ then $y_j(x-k_j)\leq y_j(x+k_j)$ if $0\leq k_j<\pi$ and
$y_j(x+k_j)\leq y_j(x-k_j)$ if $\pi\leq k_j<2\pi$.
Set
\begin{equation}
f_1(x):=\Pi_{0\leq k_j<\pi} y_j(x+k_j)
\end{equation}
\begin{equation}
f_2(x):=\Pi_{\pi\leq k_j<2\pi} y_j(-x+k_j)=\Pi_{\pi\leq k_j<2\pi} y_j(x-k_j)
\end{equation}
For $0\leq x\leq \pi$ we have $f_1(-x)\leq f_1(x)$ and $f_2(-x)\leq f_2(x)$.
By Eqn.(\ref{haymanlem49ekv}) we thus obtain
\begin{multline}
\int_{-\pi}^\pi \Phi\left(\Pi_{j=1}^q y_j(x+x_j)\right)dx=
\int_{-\pi}^\pi \Phi\left(f_1(x)f_2(-x)\right)dx=\\
\int_{0}^\pi \Phi\left(f_1(x)f_2(-x)\right)+\Phi\left(f_1(-x)f_2(x)\right)dx\leq\\
\int_{0}^\pi \Phi\left(f_1(x)f_2(x)\right)+\Phi\left(f_1(-x)f_2(-x)\right)dx=
\int_{-\pi}^\pi \Phi\left(f_1(x)f_2(x)\right)dx=\\
\int_{-\pi}^\pi \Phi\left(\Pi_{0\leq k_j<\pi} y_j(x+k_j)\Pi_{\pi\leq k_j<2\pi} y_j(x-k_j)\right)dx=
\int_{-\pi}^\pi \Phi\left(\Pi_{j=1}^q  y_j(x+x'_j)\right)dx
\end{multline}
where
$x'_j=k_j-\frac{1}{2}x_\mu$ (mod$2\pi$), $0\leq x_j'<2\pi$, if $0\leq k_j<\pi$;
and $x'_j=-k_j-\frac{1}{2}x_\mu$ (mod$2\pi$), $0\leq x_j'<2\pi$, if $\pi\leq k_j<2\pi$;
Furthermore, $k_\mu=\frac{1}{2}x_\mu<\pi,$ so that
$x'_\mu=0.$ If for some $j_0$ we have that $x'_{j_0}$ is nonzero then we may repeat the procedure in order to obtain a new set $x''_j$, $j=1,\ldots,q$,
which are zero for each $j$ such that $x'_j$ is zero and for one extra value of $j$.
After at most $q$ steps we obtain that all the $x_j$ are zero.
Since
\begin{equation}
\int_{-\pi}^\pi \Phi\left(\Pi_{j=1}^q y_j(x+k_j)\right)dx
\end{equation}
is non-decreasing from stage to stage, the proof is complete.
\end{proof}

Balk \cite{balkdisc1977} proved the following.
\begin{theorem}\label{balkdiscthm}
	Let $q\in \Z_+$ and $f(z)=\sum_{j=0}^{q}a_j(z),$ for holomorphic $a_j(z),$ on $\{\abs{z}<1\},$ $j=0,\ldots,q.$
	Suppose the zero set of $f(z)$ is compact in $\{\abs{z}<1\}$. Then $f(z)$ can be represented on $\{\abs{z}<1\}$ as
	\begin{equation}
	f(z)=G(z)\sum_{j=0}^{q-1}\pi_j(z)\bar{z}^j
	\end{equation}
	for a holomorphic function $G(z)$ and for $\pi_j(z)\in \Pi(\{\abs{z}<1\}).$\\
	Furthermore, if $\Omega(z)=(a_0(z),\ldots,a_q(z))$ with $a_q(z)\not\equiv 0$ then
	\begin{equation}
	\limsup_{r\to 1}\left(T(r,\Omega)\ln \frac{1}{1-r}\right)\leq q
\end{equation}
\end{theorem}
\begin{proof}
The proof consists of a sequence of lemmas.
\begin{lemma}\label{balkdisclemma1}
	In $\C^{n+1}\setminus\{0\}$ the following function is continuous
	\begin{equation}
	J(w)=J(w_0,\ldots,w_n):=\int_{1/2}^1 \abs{\sum_{k=0}^n r^{2k}w_k}^{-\alpha} dr, \quad \alpha\in (0,1/n)
	\end{equation} 
\end{lemma}
\begin{proof}
	Let $(a_0,\ldots,a_n)\in \C^{n+1},$ $a\neq 0.$If $Q(r,a)=\sum_{k=0}^n r^{2k}a_k$ has no zeros on $[1/2,1]$ then then we have $\lim_{w\to a} J(w)=J(a).$
	So we can assume
	$Q(r,a)=P(r,a)(r-\lambda_1)\cdots (r-\lambda_s),$ for
	$1/2 \leq \lambda_1\leq \cdots\leq \lambda_s \leq 1$ for a polynomial $P(r,a)$ which has no zeros
	on $[1/2,1].$
	Since $Q(r,a)$ is an even function of $r$ we have $s\leq n.$ There exists $p>0$ such that $P(r,a)>2p$ for 
	$r\in [1/2,1]$. 
	By a theorem of Hurwitz 
	(see, e.g.\ Marden \cite{marden}, p.4) 
	the polynomial $Q(r,w)$ can for $w$ near $a$ be written
	$Q(r,w)=P(r,w)(r-\lambda_1(w))\cdots (r-\lambda_s(w))$
	where
	$\lambda_k(w)\to \lambda_k(a)=\lambda_k$ as $w\to a,$ $1\leq k\leq s.$ and $P(r,w)\to P(r,a)$
	uniformly with respect to $r$, $r\in [1/2,1]$ as $w\to a.$
	We assume $w$ is sufficiently close to $a$ such that $\abs{P(r,w)}>p$ for $r\in [1/2,1].$
	Let $\eta\in (0,1/2).$ Pick neighborhoods $A_k$ in $\C$ of radius $2\eta$ of $\lambda_k,$
	$k=1,\ldots_s$, with $\eta$ sufficiently small such that the $A_k$ are pairwise disjoint if they have different centers and 
	$\delta$ such that $\abs{\lambda_k(w)-\lambda_k}<\delta<\eta$, $k=1,\ldots,s.$
	Denote
	$\Delta_k:=(\lambda_k-\delta,\lambda_k +\delta)\cap [1/2,1]$
	and $A:=\bigcup_{k=1}^s \Delta_k,$
	$B:=[1/2,1]\setminus A.$ Let
	$\lambda_1=\cdots \lambda_q<\lambda_{q+1}.$ 
	Then setting $x_k=\lambda_l-\re \lambda_k(w),$ $\abs{x_k}<\delta,$ $k=1,\ldots,q$ we have
	\begin{multline}\label{balkdisc8}
	\int_{\Delta_1} \abs{Q(r,w)}^{-\alpha} dr \leq p^{-\alpha} \int_{\lambda_1 -\delta}^{\lambda_1 +\delta} \frac{dr}{\abs{r-\lambda_1(w)}^\alpha \cdots \abs{r-\lambda_s(w)}^\alpha}\leq \\
p^{-\alpha} (2\eta)^{-\alpha(s-q)} \int_{\lambda_1 -\delta}^{\lambda_1 +\delta}
\frac{dr}{\abs{r-\re \lambda_1(w)}^\alpha \cdots \abs{r-\re \lambda_s(w)}^\alpha}= \\
p^{-\alpha} (2\eta)^{-\alpha(s-q)} \int_{-\delta}^{\delta}
\frac{dr}{\abs{x-x_1}^\alpha \cdots \abs{x+x_q}^\alpha}
\end{multline}
Let $\psi(x)$ be a periodic function with period $2\delta$ such that
$\psi(x)=\abs{x}^{-\alpha}$ on $[-\delta,\delta.$ By the inequality $\abs{x+x_k}^{-\alpha}\leq\psi(x+x_k),$ $k=1,\ldots,q$
we have
\begin{equation}\label{balkdisc9}
\int_{-\delta}^{\delta} \frac{dr}{\abs{x-x_1}^\alpha \cdots \abs{x+x_q}^\alpha}\leq
\int_{-\delta}^{\delta} \psi(x+x_1)\cdots \psi(x+x_q)dx
\end{equation}
By Proposition \ref{goldberghay}
\begin{equation}
\int_{-\delta}^{\delta} \left(\Pi_{k=1}^q y_k(x+x_k)\right)dx\leq
\int_{-\delta}^{\delta} \abs{y(x)}^q dx
\end{equation}
Setting $y(x)=\psi(x)$, $\alpha\in (0,1/n),$ $\alpha q<1$ we obtain
\begin{equation}\label{balkdisc10}
\int_{-\delta}^{\delta} \psi(x+x_1)\cdots\psi(x+x_q)dx \leq \int_{-\delta}^{\delta} \abs{\psi(x)}^q dx
\leq \int_{-\delta}^{\delta} \abs{x}^{-\alpha q} dx=\frac{2}{1-\alpha q} \delta^{1-\alpha q}
\end{equation}
Similarly, we obtain bounds for the integrals over $\Delta_k$, $k=2,\ldots,s.$
By Eqn.(\ref{balkdisc8}), Eqn.(\ref{balkdisc9}) and Eqn.(\ref{balkdisc10})
we obtain
\begin{equation}\label{balkdisc11}
\int_A \abs{Q(r,w)}^{-\alpha} dr\leq sp^{-\alpha} (2\eta)^{-\alpha(s-1)} \frac{2}{1-\alpha s} \delta^{1-\alpha s}
\end{equation}
Let $\epsilon>0$ and choose $\delta$ sufficiently small such that
\begin{equation}
sp^{-\alpha} (2\eta)^{-\alpha(s-1)} \frac{2}{1-\alpha s} \delta^{1-\alpha s}\leq \frac{\epsilon}{3}
\end{equation}
For $r\in B$ we have uniform convergence
\begin{equation}
\lim_{w\to a} \abs{Q(r,w)}^{-\alpha} =\abs{Q(r,a)}^{-\alpha}
\end{equation}
Let $0<\delta_1<\delta$ be sufficiently small such that for $w\in \{\abs{w-a}<\delta_1\}$ we have
\begin{equation}
\abs{\int_B \abs{Q(r,w)}^{-\alpha} dr-\int_B \abs{Q(r,a)}^{-\alpha} dr}<\frac{\epsilon}{3}
\end{equation}
This implies
\begin{multline}
\abs{J(w)-J(a)}\leq \int_A \abs{Q(r,w)}^{-\alpha} dr
+\int_A \abs{Q(r,a)}^{-\alpha} dr+\\
\abs{\int_B \abs{Q(r,w)}^{-\alpha} dr-\int_B \abs{Q(r,w)}^{-\alpha} dr}<\epsilon
\end{multline}
This proves Lemma \ref{balkdisclemma1}.
\end{proof}

\begin{lemma}\label{balkdisclem2}
	If $f_j(z),$ $j=0,\ldots,n,$ be holomorphic functions on $\{\abs{z}<1\}$ and assume the reduced $n+1$-analytic function 
	$f(z)=\sum_{j=0}^n f_j(z)\bar{z}^{2j}$ have zeros contained in $D_\rho:=\{\abs{z}<\rho\}$ for some $\rho\in (0,1)$
	Then the number of zeros, $n(c,0,\phi_c)$, in $D_c$, $\rho<c<1,$ of the function
	\begin{equation}
	\phi_c(z)=\sum_{j=0}^n f_j(z)c^{2j}
	\end{equation}
	is independent of the choice of $c$, i.e.\ $h=n(c,0,\phi_c)$ is constant.
\end{lemma}
\begin{proof}
	Define for each curve $\gamma$, $V(f,\gamma):=\frac{1}{2\pi}\int_\gamma d\mbox{arg} f(z).$	
Since $\phi_c(z)=f(z)$ on $\gamma_c:=\{\abs{z}=c\}$ and $f(z)$ is continuous and without zeros in $\{r<\abs{z}<1\}$
we have that $V(f,\gamma_c)$ is independent of $c$ for $\rho\leq c<1.$ Hence
$n(c,0,\phi_c)=V(\phi_c,\gamma_c)=V(f,\gamma_c)=V(f,\gamma_\rho)=h$
is constant.
This completes the proof of Lemma \ref{balkdisclem2}.	
\end{proof}

We shall need the following classical result, see e.g.\ Goldberg \& Ostrovskii \cite{goldostrov}.
\begin{lemma}\label{goldnevanlem}
\begin{equation}
N\left(r,\frac{1}{f}\right)-N(r,f)=\frac{1}{2\pi}\int_{0}^{2\pi} \ln \abs{f(r\exp(i\theta))}d\theta -\ln \abs{c_\lambda}
\end{equation}
\end{lemma}
\begin{proof}
Let $f(z)$, $f\not\equiv 0$ on the disc $\{\abs{z}\leq R\}$ for some $R>0$ and let
have Laurent expansion near $z=0$ given by
\begin{equation}
f(z)=c_\lambda z^\lambda +c_{\lambda +1}z^{\lambda +1} +\cdots ,\quad c_\lambda\neq 0
\end{equation}
By the Poisson-Jensen formula 
we have
\begin{multline}
\ln \abs{f(z)} =
\frac{1}{2\pi}\int_0^{2\pi} \ln\abs{f(R\exp(i\theta))}\re \frac{R\exp(i\theta))+z}{R\exp(i\theta))-z}d\theta
-\\
\sum_{\abs{a_m}<R}\ln \abs{ \frac{R^2-\bar{a}_{m}z}{R(z-a_m)}} +
\sum_{\abs{b_n}<R}\ln \abs{\frac{R^2-b_{n}z}{R(z-b_n)}}
\end{multline}
where the $a_m$ are zeros of $f(z)$ and where the $b_n$ are the poles of $f(z)$
Letting $z\to 0$ we obtain in the right hand side
\begin{multline}
-\sum_{\abs{a_m}<R} +\sum_{\abs{b_n}<R}=
-\sum_{0<\abs{a_m}<R} +\sum_{0<\abs{b_n}<R} +\lambda\ln\frac{\abs{z}}{R}=\\
-\sum_{0<\abs{a_m}<R}\ln\frac{R}{\abs{a_m}} +\sum_{0<\abs{b_n}<R}\ln\frac{R}{\abs{b_n}} +o(1)+\lambda\ln\frac{\abs{z}}{R}=\\
\end{multline}
and in the left hand side \begin{equation}\ln \abs{f(z)}=\lambda \ln\abs{z} +\ln \abs{c_\lambda} +o(1)\end{equation}
This implies
\begin{multline}
\ln \abs{c_\lambda}=\frac{1}{2\pi}\int_0^{2\pi} \ln\abs{f(R\exp(i\theta))}d\theta -\\
-\sum_{0<\abs{a_m}<R}\ln\frac{R}{\abs{a_m}} +\sum_{0<\abs{b_n}<R}\ln\frac{R}{\abs{b_n}}-\lambda\ln R
\end{multline}
which thus determines the $c_\lambda.$
Now by the properties of the Stieltjes integral we have
\begin{multline}
\sum_{0<\abs{b_n}<r}\ln \frac{r}{\abs{b_n}}=\int_{+0}^r \ln\frac{r}{t} dn(t,f)=
\int_0^r \ln \frac{r}{t} d(n(t,f)-n(0,f))=\\
\ln\frac{r}{t}(n(t,f)-n(0,f))|_{+0}^r +\int_0^r \frac{n(t,f)-n(0,f)}{t}dt=\\
\int_0^r \frac{n(t,f)-n(0,f)}{t}dt
\end{multline}
The lemma now follows from the equalities
\begin{equation}
\sum_{0<\abs{b_n}<r}\ln\frac{r}{\abs{b_n}}=N(r,f)-n(0,f)\ln r
\end{equation}
\begin{equation}
\sum_{0<\abs{a_m}<r}\ln\frac{r}{\abs{a_m}}=N\left(r,\frac{1}{f}\right)-n\left(0,\frac{1}{f}\right)\ln r
\end{equation}
\begin{equation}
-\lambda=n\left(0,f\right)-n\left(0,\frac{1}{f}\right)
\end{equation}
where the $a_m$ are zeros of $f(z)$ and where the $b_n$ are the poles of $f(z)$ and
$\lambda$ determined as above.
This proves Lemma \ref{goldnevanlem}.
\end{proof}

\begin{lemma}\label{balkdisclem3}
If $\Omega(z)=(\phi_0,\ldots,\phi_n(z))$ with $\phi_n(z)\not\equiv 0,$
and all $\phi_j(z),$ $j=0,\ldots,n$ holomorphic such that
\begin{equation}
T(r,\Omega)=O\left(\ln \frac{1}{1-r}\right),\quad r\to 1
\end{equation}
then there exists $H(z)\in \Pi(\{\abs{z}<1\})$ such that 
$H(z)$ has the same zeros, with the same multiplicity, as $\phi_n(z).$
\end{lemma}
\begin{proof}
If $\phi_n(z)$ has finitely many zeros we are done, so suppose $\phi_n(z)$ has
infinitely many zeros. By Lemma \ref{goldnevanlem} (application of the Jensen formula)
we have for $r\in (0,1),$ $r\to 1$ that there exists $c'>0$ such that
\begin{multline}
N(r,0,\phi_n)=\frac{1}{2\pi} \int_0^{2\pi} \ln\abs{\phi_n (e\exp(i\theta))}d\theta +c'\leq\\
T(r,\Omega)+O(1)=O\left(\ln \frac{1}{1-r}\right)
\end{multline}
W.l.o.g.\ we we assume $\phi_n(0)\neq 0$ and let $\{a_k\}_{k\in \N}$ be the zeros of $\phi_n(z).$ We have
\begin{equation}
N\left(\frac{r+1}{2},0,\phi_n\right)\geq =\int_r^{\frac{r+1}{2}}\frac{n(r,0,\phi_n)}{t}dt  \geq n(r,0,\phi_n) \frac{1-r}{1+r}
\end{equation}
\begin{equation}
n(r,0,\phi_n)\leq \frac{1+r}{1-r} N\left(\frac{r+1}{2},0,\phi_n\right)=O\left(
\frac{1}{1-r}\ln \frac{1}{1-r}\right)
\end{equation}
This implies
\begin{equation}
\sum_{k=1}^\infty (1-\abs{a_k})^2 =2\int_0^1 n(t,0,\phi_n)(1-t)dt <\infty
\end{equation}
Set
\begin{equation}
H(z):=\Pi_{k=1}^\infty E(L_k(z),1),\quad L_k(z)=(1-\abs{a_k}^2)/(1-\bar{a}_k z)
\end{equation}
where $E(\zeta,1)$ is a primary Weierstrass multiplier. 
For $r=\abs{z},$ $1-\rho=(1-r)/4,$ and $1>\abs{a_k}\geq \rho$ we have
\begin{equation}
\abs{L_k(z)}<\frac{1-\rho^2}{1-r}<2\frac{1-\rho^2}{1-r}=\frac{1}{2}
\end{equation}
Now $\abs{E(\zeta,1)}\leq \abs{\zeta}^2$ for $\abs{\zeta}\leq 1/2$ which yields
\begin{equation}
\ln^+\abs{E(L_k(z),1)}\leq \abs{L_k(z)}^2
\end{equation}
\begin{multline}
m(r,E(L_k(z),1),1)\leq (1-\abs{a_k}^2)^2 \frac{1}{2\pi} \int_0^{2\pi} \frac{d\theta}{\abs{1-\abs{a_k}r\exp(i\theta))}^2} 
=\\
\frac{(1-\abs{a_k}^2)^2}{1-\abs{a_k}^2r^2} <
\frac{4(1-\abs{a_k}^2)}{1-r}
\end{multline}
where we assume that $\re L_k(z)>0$ for $\abs{z}<1,$ $\abs{a_k}<1$ and we have
\begin{equation}
\ln \abs{E(L_k(z),1)}=\ln \abs{1-L_k(z)}+\re L_k(z)\leq \re L_k(z)
\end{equation}
This implies for $\abs{a_k}<\rho$
\begin{equation}
m(r,E(L_k(z),1))\leq  \int_0^{2\pi} \re L_k(r\exp(i\theta))d\theta =1-\abs{a_k}^2\leq 2(1-\abs{a_k})
\end{equation}
Hence
\begin{multline}
m(r,H)\leq 2\sum_{\abs{a_k}<\rho}(1-\abs{a_k})+\frac{4}{1-r}\sum_{\abs{a_k}\geq \rho} (1-\abs{a_k})^2=\\
2\int_0^\rho (1-t)dn(t) +\frac{1}{1-\rho}\int_\rho^1(1-t)^2 dn(t)\leq\\
2\int_0^\rho dn(t)dt+(1-\rho)n(\rho) +\frac{2}{1-\rho}\int_\rho^1 (1-t)n(t)dt
\end{multline}
Since $(1-\rho)n(\rho)=O\left(\ln \frac{1}{1-\rho}\right)$ as $\rho \to 1$ we have
\begin{equation}
\int_0^\rho dn(t)dt \leq N(\rho,0,\phi_n)= O\left(\ln \frac{1}{1-\rho}\right)
\end{equation}
\begin{equation}
\int_\rho^1 dn(t)dt = O\left(\int_\rho^1 \ln  \frac{1}{1-t}dt\right)
= O\left((1-\rho)\ln \frac{1}{1-\rho}\right)
\end{equation}
Hence 
\begin{equation}
m(r,H)=O\left(\ln\frac{1}{1-\rho}\right)=O\left(\ln\frac{1}{1-r}\right),\mbox{ as }r\to 1
\end{equation}
Since $T(r,H)=m(r,H)$, $r\in (0,1)$, the function $H(z)$ satisfies the wanted properties
of Lemma \ref{balkdisclem3}. This completes the proof of Lemma \ref{balkdisclem3}. 	
\end{proof}

Now for the $n+1$-analytic function $\phi(z)=\sum_{j=0}^n \phi_j(z)\bar{z}^j$, set
\begin{equation}
f(z):=\phi(z)z^n=\sum_{j=0}^n f_j(z)\abs{z}^{2j}
\end{equation}
Let $a$ be the first nonzero coefficient of the Taylor expansion of
$f_n(z)$ and let $F(z)=(f_0(z),\ldots,f_n(z))$, in particular
\begin{equation}\label{balkdisc19}
\norm{F(z)}=\sum_{j=0}^n \abs{f_j(z)},\quad T(r):=T(r,F)=\frac{1}{2\pi}\int_0^{2\pi} \ln \norm{F(r\exp(i\theta))}d\theta
\end{equation}
Since each pole of $f_j(z)/f_n(z)$ in $\{\abs{z}<1\}$ is a pole of $1/f_n(z)$ we have
\begin{equation}
N\left(r,\frac{f_j}{f_n}\right)-n\left(0,\frac{f_j}{f_n}\right)\ln r\leq N\left(r,\frac{1}{f_n}\right)-n\left(0,\frac{1}{f_n}\right)\ln r
\end{equation}
Hence for $r\in (0,1)$ we have
\begin{equation}
N\left(r,\frac{f_j}{f_n}\right) \leq N\left(r,\frac{1}{f_n}\right)+n\left(0,\frac{1}{f_n}\right)
\ln \frac{1}{r}
\end{equation}
\begin{multline}
T\left(r,\frac{f_j}{f_n}\right)=\frac{1}{2\pi} \int_0^{2\pi} \ln^+\abs{\frac{f_j(r\exp(i\theta))}{f_n(r\exp(i\theta))}}d\theta +N\left(r,\frac{f_j}{f_n}\right)\leq\\
	\int_0^{2\pi} \ln \norm{F(r\exp(i\theta))}d\theta - \int_0^{2\pi} \ln \abs{f_n(r\exp(i+theta))}d\theta+N\left(r,\frac{1}{f_n}\right)+n\left(0,\frac{1}{f_n}\right)\ln\frac{1}{r}
	\end{multline}
	By the Jensen formula 
	we have 
	\begin{equation}
	-\frac{1}{2\pi}\int_0^{2\pi} \ln \abs{f_n(r\exp(i\theta))}d\theta +N\left(r,\frac{1}{f_n}\right)=-\ln\abs{a}
	\end{equation}
	and for $r\to 1$
	\begin{equation}\label{balkdisc20}
	T\left(r,\frac{f_j}{f_n}\right)\leq T(r)-\ln\abs{a}+n\left(0,\frac{1}{f_n}\right)\ln\frac{1}{r}\leq T(r)+O(1)
	\end{equation}
	Set
	\begin{equation}\label{balkdisc21}
	\varphi_c(z):=\sum_{j=0}^n c^{2j} f_j(z),\quad H_c(z)=\frac{\norm{F(z)}}{\abs{\varphi_c(z)}}
		\end{equation}
		\begin{equation}\label{balkdisc22}
		M_c(t)=\frac{1}{2\pi}\int_0^{2\pi} \ln H_c(t\exp(i\theta))d\theta,\quad t\in (0,1)
		\end{equation}
		By the Jensen formula applied to $\varphi_c(z)$ we have
		\begin{equation}\label{balkdisc23}
		\frac{1}{2\pi}\int_0^{2\pi} \ln \abs{\varphi_c(t\exp(i\theta))}d\theta=N(t,0,\varphi_c)+\ln \abs{a(c)}
			\end{equation}
			where $a(c)$ is the first nonzero Taylor coefficient of the expansion of $\varphi_c(z)$ with respect to $z.$
			Let $1/2\leq \rho \leq r \leq c<1$.
			By Eqn.(\ref{balkdisc19}),Eqn.(\ref{balkdisc21}), Eqn.(\ref{balkdisc22}), Eqn.(\ref{balkdisc23}) together 
			with $H_c(<)\geq 1$ and Lemma
			\ref{balkdisclem2}
			we have
			\begin{equation}
			T(r)=M_c(r)+N(r,0,\varphi_c)+\ln\abs{a(c)}
			\end{equation}
			\begin{equation}
			T(\rho)=M_c(\rho)+N(\rho,0,\varphi_c)+\ln\abs{a(c)}\geq N(\rho,0,\varphi_c)+\ln\abs{a(c)}
			\end{equation}
			\begin{equation}
			N(r,0,\varphi_c)-N(r,0,\varphi_c)=\int_\rho^r\frac{n(t,0,\varphi_c)}{t}dt \leq
			\int_\rho^r\frac{n(c,0,\varphi_c)}{t}dt<h\ln\frac{1}{\rho}
			\end{equation}
			This implies that there is a constant $C>0$ independent of $r$ and $c$ such that
			\begin{equation}
			T(r)<M_c(r)+C
			\end{equation}
			This gives for constants, $\mbox{const.}$, independent of $r$ and $c$ 
			\begin{equation}\label{balkdisc25}
			T(r)<\frac{1}{1-r}\int_r^1 M_c(r)dc +\mbox{const.}=\frac{1}{2\pi}\int_0^{2\pi} d\theta \left(\frac{1}{1-r}\int_r^1 \ln H_c(r\exp(i\theta))dc\right)+\mbox{const.}
			\end{equation}
			Now we have
			$\varphi_c(z)/\norm{F(z)}=\sum_{j=0}^n w_k c^{2j}$ where $\sum_{j=0}^n \abs{w_j}=1.$ Set
			$\sigma:=\{w\in \C^{n+1}:\sum_{j=0}^n \abs{w_j}=1\}.$ By Eqn.(\ref{balkdisc25})
			we have
			\begin{equation}\label{balkdisc26}
			T(r)<\lambda(r)+\mbox{const.}
			\end{equation}
			where $\lambda(r)=\sup_{w\in \sigma} \left(\frac{1}{1-r} \int_r^1 \ln \abs{\sum_{j=0}^n w_j c^{2j}}^{-1} dc\right)$
			By an analogue of the Cauchy inequalites (see Polya \& Szeg\"o \cite{polyaszego}, p.77)
			\begin{equation}\label{balkdisc27}
			\frac{1}{b-a}\int_a^b \ln \psi(t)dt \leq \ln\left(\frac{1}{b-a}\int_a^b \psi(t)dt\right)
			\end{equation}
			for real positive $\psi(t)$ on $[a,b].$ Let $\alpha\in (0,1/n).$ By Eqn.(\ref{balkdisc25}) and Eqn.(\ref{balkdisc27})
			\begin{multline}\label{balkdisc28}
			\frac{1}{1-r}\int_r^1 \ln \abs{\sum_{j=0}^n w_j c^{2j}}^{-1} dc\leq \frac{1}{\alpha} \ln \left(\frac{1}{1-r}\int_r^1 
			\abs{\sum_{j=0}^n w_j c^{2j}}^{-\alpha}dc\right)\leq\\
			\frac{1}{\alpha} \ln \frac{1}{1-r} +\frac{1}{\alpha}\ln J(w)
			\end{multline}
					By Lemma \ref{balkdisclemma1}, $J(w)$ is continuous on $\sigma$ thus
			$J(w)$ is bounded from above on $\sigma$ which gives
			\begin{equation}
			\lambda(r)\leq \frac{1}{\alpha}\ln\frac{1}{1-r} +O(1),\mbox{ as} r\to 1
			\end{equation}
		By Eqn.(\ref{balkdisc26}) and Eqn.(\ref{balkdisc28}) we have
			\begin{equation}
			\limsup_{r\to 1} \frac{T(r,F)}{\ln\frac{1}{1-r}} \leq n
			\end{equation}
			Thus for $r\in (0,1)$ we have
			\begin{equation}
			T(r,F)\leq T(r,\Omega)\leq T(r,F)+n\ln\frac{1}{r} =T(r,F)+o(1)
			\end{equation}
			This proves the first part of the theorem.
			We have also shown, By Eqn.(\ref{balkdisc26}), Eqn.(\ref{balkdisc28}) and Eqn.(\ref{balkdisc20}) that
			\begin{equation}
			T\left(r,\frac{f_j}{f_n}\right)=O\left(\ln \frac{1}{1-r}\right),\mbox{ as }r\to 1
			\end{equation}
			By Lemma \ref{balkdisclem2} there exists $H(z)$ that is holomorphic in $\{\abs{z}<1\}$
			and has the same zeros as $\phi_n(z)$ in the unit disc
			such that $H(z)\in \Pi(\{\abs{z}<1\}).$
			Then the function $G(z):=\phi_n(z)/H(z)$
			is holomorphic on the unit disc and has no zeros there. Setting $\pi_j(z):=\phi_j(z)/G(z)$ we obtain $\pi_j(z)\in \Pi(\{\abs{z}<1\}).$
			This completes the proof of the first part of the theorem.
This proves Theorem \ref{balkdiscthm}.
\end{proof}
Balk \cite{balkdisc1977} also points out that Theorem \ref{balkdiscthm} is no longer true
if the condition $T(r,\pi)=O\left(\ln\frac{1}{1-r}\right)$, as $r\to 1$,
is replaced by $T(r,\pi)=o\left(\ln\frac{1}{1-r}\right)$, as $r\to 1$.

\section{Further notable results}

Some further results on zero sets for polyentire functions are given in Chapter \ref{polyentiresec}.
Balk \cite{balkdisc1977} proved the following.
\begin{theorem}\label{balk1977}
	Suppose $f_j(z)$, $j=0,\ldots,n$ are holomorphic on the unit disc and suppose that
	the zeros in the unit disc, of
	$f(z):=\sum_{j=0}^n f_j(z)\abs{z}^{2j}$, belong compactly to the unit disc. Then
	$f(z)$ has on the unit disc the representation
	$f(z)=G_1(z)\sum_{j=k}^n g_j(z)(1-\abs{z}^2)^j$
	for polynomials $g_j(z)\in \Pi(\{\abs{z}<1\}),$ $g_k(z)\not\equiv 0$,
	and a holomorphic $G_1(z)$ that does not have zeros in the unit disc.
	\end{theorem}
	\begin{definition}
		Let $p\in \Omega,$ for a domain $\Omega\subset\C$ and 
		let $f(z)=\sum_{j=0}^{q-1}a_j(z)(\bar{z}-\bar{p})^j$, for holomorphic $a_j(z)$, be a $q$-analytic function $\Omega.$
		The function $f(z)$ is called {\em pseudoanalytically irreducible}\index{Pseudoanalytically irreducible} at $p$ if the function of two complex variables
		$F(z,w):=\sum_{j=0}^{q-1}a_j(z)(w-\bar{p})^j$ is pseudoanalytically irreducible at $(p,\bar{p}).$
	\end{definition}
	Balk proves the following theorem (see Balk \cite{ca1}, p.207, Theorem 1.7, Balk \cite{balk1974irred} and Schopf \cite{schopf}).
	\begin{theorem}
		Let $f$ be a polyanalytic function that is pseudoanalytically irreducible at a point $a$, with $f\not\equiv 0.$ If $a$ is a non-isolated zero of $f$ then the set of all zeros of $f$ in a 
		sufficiently small neighborhood of $a$ coincides with the set of all points of an open simple (but possibly non-regular) analytic arc; the point $a$ splits this arc into two
		OSCAR:s.
	\end{theorem}
Balk \cite{ca1}, Theorem 1.7 (i), states without proof that a corollary is the following information on zero sets of countably analytic functions.
\begin{theorem}[See Balk \cite{ca1}, p.208, (i)]\label{balk208bra}
	Let $\Omega\subset\C$ be a domain and let $f$
	be a countably analytic function on $\Omega.$ Denote $f^{-1}(0):=\{z\in \Omega:f(z)=0\},$
	and denote by $Z$ the subset of all non-isolated points of $f^{-1}(0).$
	For each $p\in Z$ there exists an open neighborhood $U_p$ of $p$ in $\Omega$
	such that $Z\cap U_p$ is a union of at most $2q-2$ OSCARs (see Definition \ref{oscardef}) emanating from $p$.
\end{theorem}

The following result is due to Balk \& Zuev \cite{balkzuev1972}.
\begin{theorem}\label{balkzuev1972thm}
	Let $q\in \Z_+$ and $f(z)=\sum_{j=0}^{q-1} a_j(z)\bar{z}^j$ for entire $a_j(z),$ $j=0,\ldots, q-1.$
	If $f$ is a transcendental $q$-analytic function then $f^{-1}(a)$
	contains non-isolated points for at most $(q-1)$-values $a.$
\end{theorem}
\begin{example}
	The function $f(z)=(z+\bar{z})+\Pi_{j=1}^{q-1}(z+\bar{z}-j)$ takes the value $j$ on $\re z=j,$
	$j=1,\ldots,q-1$. This show that the bound in Theorem \ref{balkzuev1972thm} is sharp.
\end{example}

\begin{theorem}[Balk \cite{balk1965}, Theorem 3]\label{balkuniformconv}
Let $\{f_j(z)\}_{j\in\N}$ be a sequence of $q$-analytic functions $q\in \Z_+,$
on a domain $\Omega\subset\C.$ Suppose $\{f_j(z)\}_{j\in\N}$ is uniformly bounded on $\Omega$
and converges on a subset $E\subset\Omega$ having a condensation point of order $q.$
Then the sequence converges uniformly on $\Omega$ to a $q$-analytic function.
\end{theorem}

\chapter{Quasinormal families and algebras of $q$-analytic functions}
Let $X\subset\C$ be a subset and let $q\in \Z_+$. We denote by $A_q(X)$ $(P_q(X))$ the set of uniform limits, on $X$, of continuous functions on $X$ that are $q$-analytic functions ($q$-analytic functions with polynomial analytic components) on the interior $\mbox{int}(X)$. For a subset $X\subset\Cn$ we denote by $A_1(X)$ the set of uniform limits, on $X$, of continuous functions on $X$ that are holomorphic on the interior $\mbox{int}(X).$ We denote by $C^0(X)$ the set of continuous functions on $X$ equipped with the sup-norm when $X$ is bounded.
We shall describe some results on quasinormal families which are a generalization of {\em normal} families so 
for the readers convenience we recall here first some well-known theory on normal 
families originally due to Montel \cite{montel1907} (but the term {\em normal} seems to appear first in Montel \cite{montel1911}). 
\begin{definition}\index{Normal family}
	A family $F$ of analytic functions on a domain $\Omega\subseteq\C$ is called {\em normal} if every sequence $\{f_j\}_{j\in \N}$ in $F$ contains either a subsequence which converges to a limit function $f\not\equiv \infty$ uniformly on compacts of $\Omega$,
	or a subsequence which converges uniformly to $\infty$ on compacts of $\Omega.$
	$F$ is called normal {\em at a point} $z_0\in \Omega$ if it is normal in a neighborhood of $z_0$ in $\Omega.$
	The family $F$ is called normal {\em at} $\infty$ if the corresponding family $\{g:f(z)=f(1/z)\}$ is normal at $z=0$.
\end{definition}
Note that normality at $\infty$ means in some sense normality in some neighborhood of $\infty$, $\Delta(\infty,R):=\{z:\abs{z}>R\}\cup\{\infty\},$ $R>0.$
We define $F$ to be normal in a domain in $\C\cup\{\infty\}$ if $F$ is normal in $\Omega\setminus\{\infty\}$ in the usual sense, as well as normal at $\infty$ in the above sense. If $F$ is not normal in a domain $\Omega$ as above then there exists at least one point in $\Omega$ at which $F$ is not normal we call such a point an {\em irregular point}\index{Irregular point}.
\begin{theorem}
	A family $F$ of analytic functions on a domain $\Omega\subseteq\C$ is called {\em normal} if and only if it is 
	normal at each point of $\Omega.$
\end{theorem}
\begin{proof}
	Necessity is obvious. For sufficiency suppose that $F$ is normal at each point of $z\in\Omega.$
	Let $\{z_n\}_{n\in \Z_+}$ be a dense subset of $\Omega$ say $z_n=x_n+iy_n$ for rational $x_n,y_n.$ Denote by $r_n>0$ the largest
	radius such that $F$ is normal on $\{\abs{z-z_n}<r_n\}.$ Note that if $z_{n_k}\to \zeta\in \Omega$ then
	$\zeta\in \{\abs{z-z_{n_k}}<r_{n_k}\}$ for sufficiently large $k,$ since $r_{n_k}\to 0$ only if $\zeta\in \partial\Omega.$
	Hence $\bigcup_{n=1}^\infty \{\abs{z-z_n}<r_n\}$ covers $\Omega.$ For each sequence $\{f_n\}_n$ in $F$ we can pick 
	a convergent subsequence $\{f^{(1)}_{n_k}\}_k$ which converges uniformly in $\{\abs{z-z_1}<r_1/2\}$
		to an analytic function or to $\infty.$
		The sequence $\{f^{(1)}_{n_k}\}_k$ in turn has a subsequence $\{f^{(2)}_{n_k}\}_k$
				which converges uniformly in $\{\abs{z-z_2}<r_2/2\}$
				to an analytic function or to $\infty.$ By iteration we obtain a diagonal sequence
				$\{f^{(k)}_{n_k}\}_k$
					which converges uniformly in $\{\abs{z-z_n}<r_n/2\}$, $n\in \Z_+$
					to either an analytic function or to $\infty.$ This yields a division of the points $z\in \Omega$ into two classes, which we denote 
					$\Omega_0$ and $\Omega_\infty$ respectively, which are open such that $\Omega_0\cap \Omega_\infty=\emptyset,$
					$\Omega_0\cup\Omega_\infty=\Omega.$ If $K\subset\Omega$ is a compact then the family
					$\{ \{\abs{z-z_n}<r_n/2\}\}$ forms an open cover of $K$ and since $K$ can be covered by a finite subcover this completes the proof
				\end{proof}
				
				Using the notion of normality we can give the following result on uniform convergence (cf. Theorem \ref{hormkonvthm}).
				\begin{theorem}[Vitali]\label{vitaliporterthm}
					Let $\Omega\subset \Cn$ be a domain and let $\{f_j\}_{j\in \N}$ be a sequence of locally bounded holomorphic functions on 
					a domain $\Omega\subset\C$. Suppose that $\lim_{n\to \infty} f_n(z)$ exists for each $z$ in a set $E\subset\omega$ which has an accumulation point
					in $\Omega.$ Then $\{f_j\}_{j\in \N}$ converges uniformly on compacts of $\Omega$ to an analytic function.
					converging uniformly on compacts to an analytic function $f(z).$
				\end{theorem}
				\begin{proof}
					The sequence $\{f_j\}_{j\in \N}$ is normal thus there exists a subsequence $\{f_{j_k}\}_{k\in \N}$ that converges normally to an analytic function $f.$ Then 
					$\lim_{k\to \infty} f_{j_k}(z)=f(z)$ for each $z\in E.$ Suppose (in order to reach a contradiction) 
					that $\{f_j\}_{j\in \N}$ does not converge uniformly on compacts of $\omega.$ to $f.$
					Then there exists $\epsilon>0$ together with a compact $K\subset \Omega$ and a subsequence $\{f_{m_l}\}_{l\in \N}$ and points
						$z_l\in K$ such that
						\begin{equation}\label{suchthatererer}
						\abs{f_{m_l}(z_l)-f(z_l)}\geq\epsilon,\quad j\in N
						\end{equation} 
						But then $\{f_{m_l}\}_{l\in \N}$ has a subsequence which converges uniformly on compacts of $\Omega$ to an analytic function $g$ and $g\not\equiv f$ by
							Eqn.(\ref{suchthatererer}). Now $f|_E=g|_E$ implies by the identity theorem that $f\equiv g$ on $\Omega$ which is a contradiction.
							This completes the proof.
						\end{proof}

\section{Quasinormal families}

\begin{definition}
	Let $\Omega\subset\C$ be a domain.
	We say that a holomorphic function $f$ admits the exceptional value\index{Exceptional value}
	$a$ on $\Omega$ if $f(z)-a$ is nowhere zero on $\Omega.$
	Note that this is equivalent to
	the condition that
	the linear combination $\lambda_0 +\lambda_1f(z)$ is nowhere zero on $\Omega$ as soon as
	$\lambda_1=-a\lambda_0.$
	We say that $f$ admits the exceptional combination \index{Exceptional combination}
	\begin{equation}
	\lambda_0+\lambda_1 f(z),\quad \lambda_11\neq 0
	\end{equation}
	if for an exceptional value $a$ we have $\lambda_0=-a\lambda_1.$
	\end{definition}
	Montel \cite{montel1926}
		considered on a domain $\Omega\subset\C$ linear combinations 
	\begin{equation}
	F=\lambda_0 +\lambda_1 f_1+\cdots +\lambda_m f_m
	\end{equation}
	of nonvanishing holomorphic functions $f_j$, and constant coefficients $\lambda_j$, and calls these
	{\em exceptional combinations}. He points out that
	such combinations play for systems of functions $\mathbf{f}=(f_1,\ldots,f_m)$, the role that
	the linear combinations $\lambda_0+\lambda_1 f$ plays for finding the exceptional values of the function $f.$
	We say that the system $\mathbf{f}$ {\em admits the exceptional combination $F$}
	if $F$ does not vanish on $\Omega.$
	We say that the system $\mathbf{f}$ admits $m$ distinct exceptional combinations
	\begin{equation}
	\begin{array}{l}
	\lambda_0^1+\lambda_1^1 f_1+\cdots\lambda_m^1 f_m,\\
	\qquad \vdots\\
	\lambda_0^m+\lambda_1^m f_1+\cdots\lambda_m^m f_m
	\end{array}
	\end{equation}
	if the determinant 
	\begin{equation}
	\abs{
	\begin{array}{lll}
\lambda_1^1 & \cdots &\lambda_m^1\\
\vdots & \vdots & \vdots\\
\lambda_1^m & \cdots &\lambda_m^m
\end{array}
}
\end{equation}
is nonzero.
We say that $m+1$ exceptional combinations 
\begin{equation}
\begin{array}{l}
\lambda_0^1+\lambda_1^1 f_1+\cdots\lambda_m^1 f_m,\\
\qquad \vdots\\
\lambda_0^m+\lambda_1^{m+1} f_1+\cdots\lambda_m^{m+1} f_m
\end{array}
\end{equation}
are {\em distinct} for some system $\mathbf{f}$
if the determinant 
\begin{equation}
\abs{
\begin{array}{lll}
\lambda_0^1 & \cdots &\lambda_m^1\\
\vdots & \vdots & \vdots\\
\lambda_0^{m+1} & \cdots &\lambda_m^{m+1}
\end{array}
}
\end{equation}
is nonzero.
Suppose the functions $f_1,\ldots,f_m$ are entire. If the expression of a combination has finitely many zeros
then if is exceptional outside a sufficiently large disc in $\C$
thus reduces to
a polynomial $P(z)$ or to the product
$P(z)\exp(Q(z)$ for an entire function $Q.$
In the case that
an exceptional combination is of {\em first type}\index{First type}
if it is a polynomial and of the {\em second type}\index{Second type}
otherwise (i.e.\ when it is transcendent).
\begin{theorem}
When the functions $f_1,\ldots,f_m$ are entire such that at least one is not a polynomial 
then the number of exceptional combinations
is $\leq 2m-1.$ There cannot exist more than $m-1$ combinations of first type and there cannot exist more than $m$ combinations of second type.
\end{theorem}
\begin{proof}
Suppose that at least one is not a polynomial.
If there existed $m$ distinct combinations of first type then
there are polynomials $P_j$ such that
\begin{equation}
\begin{array}{l}
\lambda_0^1+\lambda_1^1 f_1+\cdots\lambda_m^1 f_m =P_1,\\
\qquad \vdots\\
\lambda_0^m+\lambda_1^{m} f_1+\cdots\lambda_m^{m} f_m=P_m
\end{array}
\end{equation}
such that the determinant $\delta$ of the system defined by the coefficients $\lambda_i^k$ is nonzero.
Solving the system with respect to $f_1,\ldots,f_m$ enders that these are all polynomials which is a contradiction to the assumption that
at least one of the $f_j$ is not a polynomial. This proves that there cannot exists more than $m-1$ combinations of first type.
Now suppose there exists $m+1$ combinations of second type.
Then there are polynomials $P_1,\ldots,P_{m+1}$ and nonconstant entire functions $Q_1,\ldots,Q_{m+1}$ such that
\begin{equation}
\begin{array}{l}
F_1:=\lambda_0^1+\lambda_1^1 f_1+\cdots\lambda_m^1 f_m =P_1\exp(Q_1),\\
\qquad \vdots\\
F_{m+1}:=\lambda_0^m+\lambda_1^{m+1} f_1+\cdots\lambda_m^{m+1} f_m=P_{m+1}\exp(Q_{m+1})
\end{array}
\end{equation}
Since the combinations are supposed to be distinct, the determinant $\Delta$ 
\begin{equation}
\Delta=\abs{
	\begin{array}{ll}
	\delta & A\\
	\lambda & P_{m+1}\exp(Q_{m+1})
	\end{array}
}
\end{equation}
is nonzero, where  
\begin{equation}
A=
\begin{bmatrix}
P_1\exp(Q_1)\\
\vdots\\
P_m\exp(Q_{m})
\end{bmatrix},\quad \lambda=(\lambda_1^{m+1},\ldots,\lambda_m^{m+1})
\end{equation}
and thus we have constants $C_1,\ldots,C_{m+1}$ such that
\begin{equation}
\Delta=\sum_{j=1}^{m+1}C_j P_j(z)\exp(Q_j(z))
\end{equation}
Since $\Delta\neq 0$ this is impossible by a known result of Borel \cite{borel1896}.
This proves that there cannot exists more than $m$ combinations of second type.
Finally, the total number of possible combinations must be less then the sum of combinations of first type and combinations of second type i.e\
$m+m-1=2m-1.$ This completes the proof.
	\end{proof}
The notion of quasinormality was introduced by Montel \cite{montelquasi}.
\begin{definition}[Quasinormal family]
Let $F$ be a family of analytic functions on a domain $\Omega\subset\C$. $F$
is called {\em quasinormal} on $\Omega$ if each sequence $\{f_j\}_{j\in \N}$ in $F$
contains a subsequence which converges uniformly on compacts of $\Omega\setminus E$,
where $E$ is a (possibly empty) finite set of points in $\Omega.$ 
The choice of $E$ may vary with the choice of subsequence. The points of $E$ are called {\em irregular points}.
\index{Irregular points of a quasi-normal family} Note that the convergence is either to an analytic function
or to the function $f\equiv \infty$ on $\Omega\setminus E.$ $F$ is said to have order of quasi-normality $q$ if the set of irregular points
always has $\leq q$ elements and there exists one subsequence for which the set of irregular points has $q$ points.
\index{Order of quasi-normality}
\end{definition}

\begin{theorem}\label{a2thm}
	Let $F$ be a quasi-normal family of analytic functions on a domain $\Omega\subset\C.$ Let $\{f_j\}_{j\in \N}$ be a sequence such that there is a subsequence $\{f_{j_k}\}_{k\in \N}$ which converges uniformly on compacts of $\Omega\setminus E,$ but not on any compact subset
		containing points of $E$, where $E$ is a finite nonempty set. Then the limit function, $f(z),$ of the subsequence must be $\equiv \infty$ on $\Omega\setminus E.$
	\end{theorem}
	\begin{proof}
		Suppose $f(z)$ is analytic on $\Omega\setminus E.$ Let $p_0\in E$ and $r>0$ sufficiently small such that
		$\{\abs{z-p_0}=r\}\subset\Omega,$ where the circle $\{\abs{z-p_0}=r\}$ does not contain any other points of $E.$
		Since $f_{j_k}\to f$ uniformly on  $\{\abs{z-p_0}=r\}$ we have by the Weierstrass theorem
		that $f_{j_k}\to f$ uniformly on
		$\{\abs{z-p_0}\leq r/2\}$ which is a contradiction. Hence
		$f$ cannot be analytic, whereas the uniform limit of analytic functions must be analytic,
		thus $f\equiv \infty.$
	\end{proof}
\begin{proposition}\label{a4prop}
	Let $\{f_j\}_{j\in\N}$ be a sequence in a quasinormal family $F$ on a domain $\Omega\subset\C.$
	Assume that for some $p_0\in \Omega$ no subsequence of $\{f_j\}_{j\in\N}$
	converges uniformly on compacts of $\Omega$ containing $p_0.$ Then, if $a\in\C,$ and $U$ is a fixed neighborhood of $p_0$, each
	equation $f_j(z)-a$ has a root in $U,$ for sufficiently large $j.$
\end{proposition}
\begin{proof}
	Assume not. Then there exists $a\in \C$, an $r>0$ such that $U=\{\abs{z-p_0}\leq r\}\subseteq\Omega$ and a subsequence
	$\{f_{j_k}\}_{k\in \N}$ such that $f_{j_k}(z)\neq a$ on $U.$ Noe there is a subsequence
	$\{f_{j_{k_l}}\}_{l\in \N}$ which converges uniformly on compacts of $\Omega\setminus E$ for a finite set $E$ of irregular points.
	Since $p_0\in E$ Theorem \ref{a2thm} implies that $f_{j_{k_l}}\to \infty$ uniformly on compacts of $\Omega\setminus E,$
		in particular on $\{\abs{z-p_0}=r\}.$ Then the function $\phi_l(z)=1/(f_{j_{k_l}}(z)-a)$ is analytic on
		$\{\abs{z-p_0}< r\}$ and goes to $0$ uniformly on $\{\abs{z-p_0}=r\}$ thus $\phi_l(z)\to 0$ as $l\to \infty.$
		This implies that $f_{j_{k_l}}\to \infty$ uniformly in $\{\abs{z-p_0}< r\}$ which is a contradiction. This completes the proof.
	\end{proof}
	The following is a quasinormal version of the so-called fundamental normality test.
	\begin{proposition}\label{montelpropenquasinormal}
		Let $F$ be a family of analytic functions
		on a domain $\Omega\subset\C$ which do not take the value $a\in \C$ more than $n_1$ times and do not take the value $b\in \C$
		more than $n_2$ times. Then $F$ is quasinormal of order $\leq \min\{n_1,n_2\}.$ 
	\end{proposition}
	\begin{proof}
		W.l.o.g.\ we cn assume $a=0,b=1.$
		Let $\{f_j\}_{j\in\N}$ be a seq
		uence in $F$ and $p_1\in \Omega$ an accumulation point of the zeros of the $f_j$ if such a point exists.
		Let $\{f_{j_k}\}_{k\in \N}$ be a subsequence such that each member has a zero in some fixed arbitrarily small disc $\Delta_1$
		centered at $p_1.$
		Suppose that the set of zeroes has another accumulation point $p_2$ in $\Omega$, we may assume (if necessary after removing a finite number of terms
		in the beginning of $\{f_{j_k}\}_{k\in \N}$ we can take the radius of $\Delta_1$ arbitrarily small) 
		$p_2\notin \Delta_1$. Let $\Delta_2$ be a disc centered at $p_2$ and disjoint from $\Delta_1.$
		Let $\{f_{j_{k_l}}\}_{l\in \N}$ be a subsequence such that each member has a zero in $\Delta_2.$ This process can be repeated in order to find
		$n_1'\leq n_1$ points $p_1,\ldots,p_{n_1'}$ and a subsequence $\{f_m\}_{m\in N_1},$ for an infinite set $N_1\subset\N$, such that
		such that each member has a zero in the discs $\Delta_1,\ldots,\Delta_{n_1'}.$ Let $\Omega'\Subset \Omega$
		and $\bigcup_{i=1}^{n_1'}\Delta_i\subseteq \Omega'.$ Set $\Omega''=\Omega'\setminus \bigcup_{i=1}^{n_1'}\Delta_i.$
		Then $\{f_m\}_{m\in N_1}$ has only finitely many zeros in $\Omega''$ so for sufficiently large $m$ $\{f_m\}_{m\in N_1}$ has no zeros
		in $\Omega''.$ Repeating the above procedure with $\{f_m\}_{m\in N_1}$ for the equations $f_m(z)-1=0$ we obtain
		$n_2'\leq _2$ points $q_1,\ldots,q_{n_2'}$ (some of which may coincide with some of $p_1,\ldots,p_{n_1'}$) and a subsequence $\{f_n\}_{n\in N_2},$ for an infinite set $N_2\subset\N$
		such that each member takes the value $1$ in a small neighborhood of the points
		$q_1,\ldots,q_{n_2'}$. Since $\{f_n\}_{n\in N_2}$ may be assumed to have no zeros in $\Omega''$
		and each member takes the value $1$ at most $n_2$ times the standard (Montel) 
		normality test implies that there is a subsequence $\{f_l\}_{l\in N_3},$ for a subset $N_3\subset\N$,
		converging uniformly on compacts of $\Omega''.$ \\
		\\
		\textit{Case (i).} Assume the limit function is analytic in $\Omega''.$ As in the proof of Theorem \ref{a2thm}
		$\{f_l\}_{l\in N_3}$ converges to an analytic function on each $\Delta_i$, thus on compacts of $\Omega'$ and similarly on compacts of
			$\Omega.$ Hence no irregular point arises from Case (i).
			\\
			\\
			\textit{Case (ii).} Assume the limit function is $\equiv \infty$ on $\Omega''$ Then (since the discs $\Delta_i$ can be chosen arbitrarily small for sufficiently large $l$) $\{f_l\}_{l\in N_3}$
			converges to $\infty$ on compacts of $\Omega\setminus\{p_1,\ldots,p_{n_1'}\}.$
			If no subsequence
			$\{f_l\}_{l\in N_3}$ converges uniformly on compacts of $\Omega$ containing $p_j$ (i.e.\ $p_j$ is an irregular point) then (as in the proof of Proposition
			\ref{a4prop}) we have for sufficiently large $l$ that $f_(z)-1$ has a root in a neighborhood of $p_j.$ Since $p_j$ does not belong to
			$\{q_1,\ldots,q_{n_2'}\}$ there cannot exists such root for sufficiently large $l$. Hence 
			Hence if $p_j$ does not belong to $\{q_1,\ldots,q_{n_2'}\}$ then $p_j$ is not an irregular point. This proves that $F$ is quasinormal and any irregular point must be
			one belonging to both $\{q_1,\ldots,q_{n_2'}\}$ and $\{p_1,\ldots,p_{n_1'}\}$ simultaneously, hence the order is at most $\min\{n_1,n_2\}.$
			This completes the proof.
			\end{proof}

The uniqueness of analytic components if a $q$-analytic function makes it tempting to consider for each $q$-analytic function on a domain in $\C$ given by $f(z)=\sum_{j=0}^{q-1} f_j(z)\bar{z}^j$, the associated vector function
$(f_1(z),\ldots,f_{q-1}(z)).$
Dufresnoy \cite{dufresnoy} studied families of $q$-tuples of holomorphic functions, $(f_0,\ldots,f_{q-1})\in \mathscr{O}^q$, and their uniform limits. 
He uses the metric (distance) defined for $f,g\in \mathscr{O}^p$ by
\begin{equation}
\mbox{dist}(f,g)=\sqrt{\frac{\sum\abs{f_i g_j-f_j g_i}^2}{\left(\sum\abs{f_i}^2\right)\left(\sum\abs{g_i}^2\right)}}
\end{equation}
Convergence and the Cauchy convergence criterion is defined in a standard way with respect to this metric, which is complete in the sense that Cauchy sequences
are convergent. Similarly, continuity and uniform continuity is defined on this metric space
as usual, and this also induces convergence and Cauchy sequences of sequences $\{f^\nu\}_\nu$ in
a family $F$ of functions in $\mathscr{O}^q$.
\begin{theorem}[Dufresnoy \cite{dufresnoy}, p.18-19]\label{dufresnoylem}
	Let $F$ be a family of $n$-dimensional vector functions 
	\begin{equation}
	\phi_\nu(z)=(\phi_{\nu,0}(z),\ldots,\phi_{\nu,n-1}(z)),\quad \nu\in \N
	\end{equation}
	that are holomorphic on a domain $\Omega\subset \C.$
	For any $\nu$ consider the set, $\lambda_\nu^n$, of $2n-1$ linear combinations
	from the components of the vector function 
	$\phi_\nu(z)$ with constant coefficients (independent of $z,\nu$)
	\begin{equation}\label{balksystem}
	f_{\nu,k}=\sum_{\mu=0}^{n-1} b_{k,\mu}\phi_{\nu,\mu},\quad k=1,\ldots,2n-1
	\end{equation}
	such that for any fixed $\mu$ and any $n$ combinations out of the $2n-1$ combinations fo the system in Enq.(\ref{balksystem}),
	are linearly independent.
	Suppose that for any domain $\Omega'$ that compactly belongs to $\Omega$, all except possibly finitely many of the functions
	$f_\nu,k(z)$, $k=1,\ldots,2n-1$, $\nu\in \N$, do not have zeros in $\Omega'.$
	Then there exists a number $m,$ $1\leq m\leq n$ such that each of the following families is compact
	\begin{equation}
	\{ f_{\nu,k}/f_{\nu,m}\},\quad k=1,\ldots,n,\nu\in \N
		\end{equation}
	\end{theorem}
\begin{definition}
Let $q\in \Z_+,$ let $\Omega\subset\C$ be a domain and $K\subset\Omega$ a compact subset.
Let $\{f_\nu\}_{\nu\in\Z_{\geq 0}}$ be a family of $q$-analytic functions on $\Omega.$ In particular, there exists holomorphic components 
$a_{\nu,k}(z)$ on $\Omega$ such that each $f_\nu(z)=\sum_{k=0}^{q-1} a_{\nu,k}(z)\bar{z}^k$ on $\Omega.$
For any complex number $c$ there exists holomorphic $b_{\nu,k}(z)$ on $\Omega$ such that we can rewrite 
$f(z)=\sum_{k=0}^{q-1} b_{\nu,k}(z)(\bar{z}-\bar{c})^k$.
Such a sequence $\{f_\nu\}_{\nu\in \N}$ is called {\em regular}\index{Regular sequence} of rank $p$, $p\geq 0$,
on a compact $K$ if the sequence $\{a_{\nu,p}(z)\}_{\nu\in \N}$ is uniformly convergent on $K$ to $\infty$, and for
any fixed $k\in\{0,\ldots,q-1\}$ the sequence $\{ a_{\nu,k}(z)/a_{\nu,q}(z)\}_{0\leq k<q,\nu\in \N}$ is uniformly convergent on $K$
to a finite limit $\theta_k(z),$ and to zero when $k>p.$ The sequence is said to be regular of rank $-1$ if
for any fixed $k\in\{0,\ldots,q-1\}$ the sequence $\{a_{\nu,k}(z)\}_{\nu\in \N}$ is uniformly convergent on $K$ to
a finite limit $\psi_K(z).$ 
$\{f_\nu\}_{\nu\in \N}$ is called regular of rank $p$, $-1\leq p\leq q-1$ on $K$ {\em with respect to the vertex $c$}
if $q$ sequences
$\{ b_{\nu,k}(z,c)\}_{0\leq k<q,\nu\in \N}$ satisfy the conditions given above with $a_{\nu,k}(z)$
replaced by $b_{\nu,k}(z,c)$.
\end{definition}
Now
\begin{equation}
b_{\nu,k}=\sum_{\mu=k}^{q-1} \lambda_{k,\mu}a_{\nu,k},\quad a_{\nu,k}=\sum_{\mu=k}^{q-1} \lambda'_{k,\mu}b_{\nu,k}
\end{equation}
for constants $\lambda_{k,\mu},$ $\lambda'_{k,\mu},$ depending only on $c,k,\mu,q,$ with $\lambda_{k,k}\neq 0,$ 
$\lambda'_{k,k}\neq 0.$
This implies that if $\{f_\nu\}_{\nu\in \N}$ is regular of rank $p$ on a compact $K$ with respect to some vertex $c$ then it is regular
of rank $p$ on $K$ with respect to any other vertex also.
Furthermore, if $\{f_\nu\}_{\nu\in \N}$ is regular of rank $p$ on each of two compacts $K_1,K_2$ with $K_1\cap K_2\neq \emptyset$ then it is regular
on $K_1\cup K_2$ and if it is regular of rank $p(K)$ on every compact $K$ i a domain $\Omega\subset\C,$
and if on one of these the rank is $p_0$ then it is regular of rank $p_0$ on all the compacts.
\begin{definition}
For a domain $\Omega\subset\C,$ a sequence $\{f_\nu\}_{\nu\in \N}$ of $q$-analytic functions on $\Omega$ is called regular of rank $p$ on $\Omega$ if it is regular
of rank $p$ on each compact subset of $\Omega.$
\end{definition}
Now let $D_1,\ldots,D_m$ be a sequence of open discs in $\C$ with finitely many points removed. Suppose 
 $\{f_\nu\}_{\nu\in \N}$ is a sequence of $q$-analytic functions on $\Omega=\cup_{j=1}^m D_j$ that is regular of rank $p$ on each $D_j$. 
For each $j=1,\ldots,m$ let $\{\tilde{D}_{j,l}\}_l$ be a sequence of closed discs where for each removed point $v\in D_j$
we have removed an open disc $d_v$ centered at $v$, such that $d_v\setminus \{v\}\subset D_j,$ in such a way that
the sequence of compacts $\{\tilde{D}_{j,l}\}_l$ exhausts $D_j.$ Then the sequence of compacts
$\{\tilde{D}_{1,l}\cup \cdots \cup\tilde{D}_{m,l}\}_l$ exhausts $\Omega.$ Hence we conclude that
$\{f_\nu\}_{\nu\in \N}$ is regular of rank $p$ on $\Omega.$ 

\begin{definition}
For a domain $\Omega\subset\C,$ a sequence $\{f_\nu\}_{\nu\in \N}$ of $q$-analytic functions on $\Omega$ is called {\em quasiregular}\index{Quasiregular} 
of rank $p$ on $\Omega$ if it is regular
of rank $p$ on $\Omega\setminus E$ for any discrete set $E\subset \Omega$ (i.e.\ $E$ does not have accumulation points in $\Omega$).
If the sequence is not regular in $\Omega\setminus E'$ for any $E'\subset E$ then any point of $E$ is an {\em irregular} point for
$\{f_\nu\}_{\nu\in \N}$ in $\Omega$ and $E$ is called an {\em irregular set}\index{Irregular set} for $\{f_\nu\}_{\nu\in \N}$ in $\Omega.$
\end{definition}

\begin{definition}
For a domain $\Omega\subset\C,$ a family $F:=\{f_\eta\}_{\eta\in I}$ of $q$-analytic functions on $\Omega$, where $I$ is infinite and not necessarily
countable, is called {\em quasinormal}\index{Quasinormal} 
if any sequence of functions in $F$ has a subsequence that is quasiregular on $\Omega.$ $F$ is called quasinormal at a point $c$ if it is quasinormal in a neighborhood of $c.$
\end{definition}
\begin{remark}\label{quasirem1}
Suppose $F$ is quasinormal on $\Omega$, then clearly it is quasinormal at each point. Conversely assume it is quasinormal at each point of $\Omega$.
For each point associate to it an open disc on which $F$ is quasinormal and cover $\Omega$ by such discs. Then by the Heine-Borel lemma $F$ can be deduced to be quasinormal on $\Omega.$
\end{remark}
It is clear that if a sequence $\{f_\nu\}_{\nu\in \N}$ of $q$-analytic functions is quasiregular of rank $p=-1$ on a domain $\Omega$ then 
it is uniformly convergent to a $q$-analytic function (see Corollary \ref{ahernbrunakonsekvens}).
On the other hand if $p\geq 0$ then the representation
\begin{equation}
f_\nu(z)=a_{\nu,p}\left(\frac{a_{\nu,0}}{a_{\nu,p}}+\cdots+\frac{a_{\nu,p-1}}{a_{\nu,p}}\bar{z}^{p-1}+\bar{z}^p+\frac{a_{\nu,p+1}}{a_{\nu,p}}\bar{z}^{p+1}+\cdots +\frac{a_{\nu,q-1}}{a_{\nu,p}}\bar{z}^{q-1} \right)
\end{equation}
shows that the sequence will diverge to $\infty$ on an open set of the form $\Omega\setminus U_p$, where $U_p$ is the set of all zeros of a polyanalytic function of exact order $(p+1)$.
\begin{definition}
For a continuous and nonzero function $f(z)$ on an arc $\Gamma$ we set
\begin{equation}
V(f,\Gamma)=\frac{1}{2\pi}\int_\Gamma 
\end{equation}
The {\em multiplicity of an isolated zero}\index{Multlicity of an isolated zero} $p_0$ of a polyanalytic function $f(z)$ is defined
by the Poincar\'e index
\index{$\mbox{Ind}(p_0,f)$, Poincar\'e index}
\begin{equation}
\mbox{Ind}(p_0,f)=V(f,\gamma)
\end{equation}
where $\gamma:=\{\abs{z-p_0}=r\}$ is a sufficiently small circle centered at $p_0$ such that $f$ has no zeros (except $p_0$) or discontinuities on 
$\{\abs{z-p_0}\leq r\}$.
A polyanalytic function $f$ is said to {\em take the value $a$ at most $m$ times} in a domain $\Omega$ if $(z)-a$ has finitely may zeros $z_1,\ldots,z_l$ in $\Omega$ such that
\begin{equation}
n(f-a,\Omega):=\sum_{k=1}^l \abs{\mbox{Ind}(z_k,f-a)}\leq m
\end{equation}
If a function $f(z)=\sum_{j=0}^{q-1} a_j(z)(\bar{z}-\bar{c})^j$ for fixed $c\in \C$ is polyanalytic on $\gamma:=\{\abs{z-c}\leq r\}$, then for a fixed $r>0$ it coincides with a holomorphic function 
\begin{equation}
f(z,\gamma):=\sum_{j=0}^{q-1} a_j(z)r^{2j} \frac{a_j(z)}{(z-c)^j}
\end{equation}
on $\{\abs{z-c}= r\}$ We say that $f(z,\gamma)$ is {\em matched} with $f$ on $\gamma$.\index{Matched function}
\end{definition}

\begin{theorem}\label{quasithm1}
	Let $F$ be a family of $q$-analytic functions on a domain $\Omega\subset \C.$
	Suppose that there exists constants $a,b\in \C$, $a\neq b,$ and natural numbers $n_1,n_2$ 
	such that each member of $F$ takes the value $a$ at most $n_1$ times and the value $b$ at most $n_2$ times
	on $\Omega.$ Then $F$ is quasicontinuoous.
	\end{theorem}
	\begin{proof}
	W.l.o.g.\ we can assume $a=0$ and $b=1.$ It suffices to show that $F$ is quasicontinuous at each point say $c\in \Omega$ (see Remark \ref{quasirem1}).
	Let $\{f_\nu\}_\nu$ be a sequence in $F$, 
	such that each $f_\nu(z)=\sum_{k=0}^{q-1} a_{\nu,k}(z)\bar{z}^k$
	for holomorphic components 
	$a_{\nu,k}(z)$ on $\Omega$. Let $0<R<R',$ chosen such that $\Delta'=\{ \abs{z-c}<R'\}\subset\Omega.$
	The set of all points $z\in \Omega$ at which at least one $a_{\nu,k}(z)$ (some $\nu$) takes the value $0$ or $1$
	is at most countable. This implies that there exists $3q-2$ distinct circles $\Gamma_\mu=\{\abs{z}=R_\mu\},$
	$\mu=1,\ldots,3q-2,$ $R<R_\mu <R'$, on which none of the $f_\nu$ take the value $0$ or $1$.
	For any fixed $\mu$ we consider the family 
	\begin{equation}\label{quasifive}
	f_{\nu}(z,\Gamma_\mu)=\sum_{j=0}^{q-1} R_\mu^{2j} \frac{a_{\nu,k}(z)}{z^k}\end{equation}
		of matched functions with respect to $\Gamma_\mu.$
Let $n_{\nu,\mu}$ be the number of zeros (counting multiplicity) 
of $f_\nu(z,\gamma_\mu)$ inside $\gamma_\nu$ and
let $n_{\nu,\mu}(R)$ be the number of zeros inside $\{\abs{z-c}<R\}.$
Let $p_{\nu,\mu}$ be the order of the pole of $f_\nu(z,\gamma_\mu)$
at $z=0$ (when $f_\nu(z,\gamma_\mu)$ is holomorphic near $0$ we set
$p_{\nu,\mu}=0$). Then
\begin{equation}
n_{\nu,\mu}(R)-p_{\nu,\mu}\leq \abs{n_{\nu,\mu}-p_{\nu,\mu}}=\abs{V(f_\nu(z,\gamma_\mu),\gamma_\mu)}=
	\abs{V(f_\nu(z),\gamma_\mu)}\leq n_1
	\end{equation}
	(where $n_1$ is the bound for the number of times $f$ attains $a$ in $\Omega$).
	Hence
	\begin{equation}
	n_{\nu,\mu}(R)\leq n_1+q-1=:p_1
	\end{equation}
	For each $\nu\in \N$
	the function $f_\nu(z,\gamma_\mu)$ has at most $p_1$ zeros in $\{\abs{z-c}<R\}\setminus\{0\}.$ Similarly
	($f_\nu(z,\gamma_\mu)-1)$ has at most $p_1q_1:=n_2+q-1$ zeros in $\{\abs{z-c}<R\}\setminus\{0\}$
	(where $n_2$ is the bound for the number of times $f$ attains $b$ in $\Omega$).
	By Montel's criterion on quasinormality (see Montel \cite{montel}, \cite{montelquasi} and Shiff \cite{shiff}) we have that
	for fixed $\mu=1,\ldots,3n-2$, the family $\{f_\nu(z,\gamma_\mu)\}_{\nu\in \N}$ is quasinormal on $\{\abs{z-c}<R\}\setminus \{0\},$
	and the order of quasinormality is $\leq s:=\min(p_1,q_1).$ Hence there exists an infinite sequence of numbers
	$N_1\subset \N$ such that for each $\mu=1,\ldots,3n-2$ the sequence
	$\{f_\nu(z,\gamma_\mu)\}_{\nu\in N_1}$ (where $N_1$ has induced order and $\nu\in N_1$ is assumed to follow increasing order) 
	is uniformly convergent to either $\infty$ or a holomorphic function on a domain of the form
	$\Omega\setminus P$, where $P$
	is a finite set containing $0$ and at most $(3q-2)s$ other points.
	\\
	\\
	\textit{Case (1).} At least $2q-1$ of the $(3q-2)s$ sequences $\{f_\nu(z,\gamma_\mu)\}_{\nu\in\N}$ are uniformly convergent to $\infty$
	on $\Omega\setminus P.$\\
	We denote for the $2q-1$ sequences the indices $\mu_1,\ldots,\mu_{2q-1}.$
	For fixed $\nu$ in Eqn.(\ref{quasifive}) we pick out the system
	obtained for $\nu=1,\ldots,q.$
	This renders a system of $q$ linear equations for the functions
	$a_{\nu,k}(z)/z^k$, and the determinant of the system is the so-called Vandermonde determinant
	of the numbers $R_1^2,\ldots,R_q^2.$
	Recall that a Vandermonde matrix is a matrix of the form
	\begin{equation}
	\begin{bmatrix}
	1 & 1 & \cdots & 1\\
	R_1^2 & R_2^2 & \cdots & R_q^2\\
	\vdots & \vdots & \vdots & \vdots\\
	R_1^{2(q-1)} & R_2^{2(q-1)} & \cdots & R_q^{2(q-1)}\\
	\end{bmatrix}
	\end{equation}
	and has determinant given by
	\begin{equation}
	\mbox{det}\begin{bmatrix}
	1 & 1 & \cdots & 1\\
	R_1^2 & R_2^2 & \cdots & R_q^2\\
	\vdots & \vdots & \vdots & \vdots\\
	R_1^{2(q-1)} & R_2^{2(q-1)} & \cdots & R_q^{2(q-1)}\\
	\end{bmatrix}
	=
	\Pi_{i<j} (R^2_i-R^2_j)
	\end{equation}
	Solving the system yields
	\begin{equation}\label{hatjo}
	\frac{a_{\nu,k}(z)}{z^k}=\sum_{\mu=1}^q c_{k,\mu} f_\nu(z,\gamma_\mu),\quad k=1,\ldots,q-1,\nu\in N_1
	\end{equation}
	where the $c_{k,\mu}$ depend upon $R_1,\ldots,R_q$ but are independent of $\nu$, and where the
	rank of the matrix $[c_{k,\mu}]_{0\leq k\leq q-1,1\leq \mu\leq q}$ is $q.$
	Setting in Theorem \ref{dufresnoylem}, $f_{\nu,\mu}(z)=_\nu(z,\gamma_\mu),$ $\mu=1,\ldots,2q-1,\mu\in \N$ we obtain
		that there exists an $m\in \{1,\ldots,q\}$ such that each of the families
		\begin{equation}
		F_\mu=\{f_\nu(z,\gamma_\mu)/f_\nu(z,\gamma_m)\},\quad \mu=1,\ldots,q,\nu\in N_1
		\end{equation}
		is compact in $\{\abs{z-c}<R\}\setminus P.$ W.l.o.g.\ we can assume $m=1.$ This implies that there 
		are functions $\psi_1,\ldots,\psi_q$
		holomorphic in $\{\abs{z-c}<R\}\setminus P$ and an infinite sequence of integers $N_2\subset N_1\subset\N$ such that
		the sequence $f_\nu(z,\gamma_\mu)/f_\nu(z,\gamma_1)$
		converges uniformly on compacts of $\{\abs{z-c}<R\}\setminus P$, to $\psi_\mu,$
		$\mu=1,\ldots,q,$ $\nu\in N_2.$
		By Eqn.(\ref{hatjo}) we obtain for $k=0,\ldots,q-1$ on $\{\abs{z-c}<R\}\setminus P$ as $\nu\to \infty,$ $\nu\in N_2$
		\begin{equation}\label{quasinine}
		\frac{a_{\nu,k}(z)}{f_\nu(z) }\to b_k(z),\quad b_k(z)=z^k\sum_{\mu=1}^q c_{k,\mu}\psi_\mu(z)
		\end{equation}
		If all $b_k(z)$, $k=1,\ldots,q-1$ would vanish identically
		then $\psi_\mu\equiv 0$ for $\mu=1,\ldots,q$ which is  contradiction since $f_\nu(z,\gamma_\mu)/f_\nu(z,\gamma_1)\to \psi_\mu$ on 
		compacts of $\{\abs{z-c}<R\}\setminus P$, implies that $\psi_1\equiv 1.$
		Hence not all $b_k$ vanish identically and we pick out the largest integer, $n_1$, such that
		$b_{n_1}\not\equiv 0.$ Let $Q$ be the set of zeros of $b_{n_1}$ on $\{\abs{z-c}\}<R\}.$
	
	So by Eqn.(\ref{quasinine}) we have as $N_1\ni \nu\to \infty$ that on $\{\abs{z-c}<R\}\setminus (P\cup Q)$,
	$\{a_{\nu,n_1}(z)$ converges to $\infty$ on compacts and 
	$a_{\nu,k}(z)/a_{\nu,n_1}$ converges uniformly on compacts to $b_k(z)/b_{n_1}(z).$
	Hence the sequence $\{f_\nu\}_{\nu\in \N}$ chosen from $F$ has a subsequence regular of rank $n_1$ in $\Omega.$
	In particular, in Case (1) a sequence from $F$ contains a subsequence that is quasiregular in $\Omega.$
	\\
	\\
	\textit{Case (2).} There exists at least $q$ sequences (among the $3q-2$ sequences) that are uniformly convergent to holomorphic functions on $\Omega\setminus P.$
	\\
	By Eqn.(\ref{hatjo}) we have that for each $k=0,\ldots,q-1,$ the family
	$\{a_{\nu,k}(z)\}_{\nu\in N_1}$ is compact in $\Omega\setminus P.$ By the Weierstrass theorem the family is compact in $\Omega.$
	Hence the sequence $\{f_\nu\}_{\nu\in \N}$ chosen from $F$ has a subsequence regular of rank $-1$ in $\Omega.$
	In particular, also in Case (2) a sequence from $F$ contains a subsequence that is quasiregular in $\Omega.$
	\\
	\\
	Case (1) and Case (2) together imply that $F$ is quasinormal at $c$. Since $c$ was arbitrary this shows that $F$ is quasinormal on $\Omega.$
	This completes the proof.
\end{proof}

This can be used to prove the following.
\begin{theorem}
	If $f(z)$ is an entire polyanalytic function that is 
	{\em transcendental} and 
	that takes a value $a$ only on a bounded set of points in $\C$, then there exists
	a ray
	$L$ originating at $z=0$ such that in any angular neighborhood, $\delta(L)$, of the ray the function $f(z)$ takes any value $b\neq a$, 
	and the set $M_b$ of all $b$-points of $f(z)$ that lie in $\delta(L)$ is unbounded and consists of isolated points.
	\end{theorem}
	\begin{proof}
	It suffices to prove the theorem for the case $a=0$.
	Let
	\begin{equation}\label{quasitio}
	f(z)=\sum_{k=0}^{q-1} A_k(z)\bar{z}^k
	\end{equation} 
	
	\begin{lemma}\label{quasiapplem}
		The family $\{f_\nu\}_{\nu\in \N}$ defined by
		\begin{equation}\label{quasielva}
		f_\nu(z)=f(2^\nu z)=\sum_{k=0}^{q-1} a_{\nu,k}(z)\bar{z}^k,\nu\in \N,\quad a_{\nu,k}(z)=2^{\nu k}A_k(2^{\nu k} z)
		\end{equation}
		is not quasinormal on $\Omega:=\{0<\abs{z}<\infty\}.$
	\end{lemma}
	\begin{proof}
		Suppose (in order to reach a contradiction) that the family of Eqn.(\ref{quasielva}) is quasiregular in $\Omega=\{0\leq\abs{z}<\infty\}.$
		So it has a subsequence
		\begin{equation}\label{quasitretton}
		\{f_\nu\}_{\nu\in N1},\quad N_1\subset\N\mbox{ an infinite subset}
		\end{equation}
		that is quasiregular in $\Omega.$ 
		If it is quasiregular of rank $-1$ then on the circle $\gamma=\{\abs{z}=R>1\}$
		we would have
		$\abs{a_{\nu,k}(z)}<C_1$ for a constant $C_1,$ $k=0,\ldots,q-1,$ $\nu\in N_1$.
		By Eqn.(\ref{quasielva}) this implies that for $k=0,\ldots,q-1$ we have $\abs{A_k(z)}<C_2$ for a constant $C_2$
			on the circles
			$\{\gamma_\nu=\{\abs{z}=2^\nu R\},$ $\nu\in N_1.$ Hence the function given by Eqn.(\ref{quasitio}) would be rational 
			and in fact a constant,
			which contradicts the condition
			of the theorem. If instead the subsequence would be quasiregular of rank $n_1\geq 0$ on $\Omega$ then on a circle $\gamma=\{\abs{z}=R>1\}$
			we would have that $\{a_{\nu,n_1}(z)\}$ goes to $\infty$ uniformly, $\nu\in N_1,$ $\nu\to \infty$.
			For some subsequence $N_2\subset N_1\subset\N$ we thus obtain that on $\gamma$, 
			$\abs{a_{\nu,n_1}(z)}>1,$ i.e.\ $2^{\nu n_1}>1,$ $\nu\in N_2.$
		This implies that on the circles $\gamma_\nu=\{\abs{z}=2^\nu R\},$ $\nu\in N_2$ we have
		\begin{equation}\label{quasifjorton}
		\abs{z^{n_1} A_{n_1}} >R^{n_1}
		\end{equation}
		By a known result of Balk \cite{balk1966} (see Chapter \ref{polyentiresec})
		$f(z)$ takes the form
		\begin{equation}\label{quasifemton}
		f(z)=\exp(g(z))\sum_{j=0}^{q-1} \pi_j(z)\bar{z}^j
		\end{equation}
		for polynomials $\pi_j(z),$ in $z$, $j=0,\ldots,q-1,$ and an entire $g(z).$
		By Eqn.(\ref{quasitio}), Eqn.(\ref{quasifjorton}) and Eqn.(\ref{quasifemton})
		we have
		\begin{equation}
		A_{n_1}(z)=\exp(g(z))\pi_{n_1}(z),\quad \abs{\exp(-g(z))}\leq \frac{\abs{z^{n_1}\pi_{n_1}(z)}}{R^{n_1}},\mbox{ on }\gamma_\nu, \nu\in N_2
			\end{equation}
					This implies that $g\equiv\mbox{const.}$ thus $f(z)$ is a polynomial, which is a contradiction.		
		This proves Lemma \ref{quasiapplem}
	\end{proof}
By Lemma \ref{quasiapplem} 
the family in Eqn.(\ref{quasielva}) cannot be quasinormal in $\Omega$ thus there exists $c\in \omega$ at which it is not quasinormal.
Denote by $L$ a ray originating at $z=0$ and passing $c$. Let $\delta(L)$ be an angular neighborhood of $c$ and $\rho>0$ sufficiently small such that
a circle
$\omega_0$ centered at $c$ with radius $\rho$ belongs to $\delta(L).$
For each $\nu\in \N$ let $\omega_\nu$ be a circle (homothetically) similar to $\omega_0$ center of similarity (homothetic center) $z=0$ with coefficient $2^\nu.$
We will show that for the sequence of circles $\{\omega_\nu\}_{\nu\in \N}$
there exists a subsequence (independent of the choice of $b$) $\{\omega_{\nu\in N'},$ for an infinite subset $N'\subset\N$ such that 
for any $b$, $b\neq a$ the function $f$ takes the value $b$ in all but finitely many circles of the subsequence. The subsequence will be chosen such that
for any $b,$ $b\neq a,$
\begin{equation}\label{quasiarton}
\lim n(f-b,\omega_\nu)=\infty,\quad \nu\to \infty, \nu\in N'
\end{equation}
By a result of Balk \& Zuev \cite{balkzuev1972} (p.323, Thm 3), see Theorem \ref{balkzuev1972thm}, the preimage set of a point $b$ (also called $b$-points),\index{Preimage of a point ($a$-points)} with $b\neq a$, of an entire transcendental polyanalytic function $f(z)$ is discrete in 
$\C.$ This implies that $f$ takes the value $b$ only at an isolated point. 
Suppose that the sequence of circles $\{\omega_\nu\}_{\nu\in \N}$ does not contain a subsequence
$\{\omega_{\nu\in N'},$ for an infinite subset $N'\subset\N$ such that 
Eqn.(\ref{quasiarton}) holds true for any $b,$ $b\neq a.$
Choose in Eqn(\ref{quasielva}) a sequence $\{f_\nu(z)\}_{\nu\in N_1}$, $N_1\subset\N$ and let
$\{\omega_\nu\}_{\nu}$ be the corresponding sequence of circles. By assumption Eqn.(\ref{quasiarton}) fails for 
$\nu\in N_1,$ for all $b$, $b\neq a.$ This implies that
there exists $b$ a constant $C$ and a sequence of circles $\{\omega_\nu\}_{\nu\in N_2}$ for an infinite subset $N_2\subset N_1\subset\N,$
such that $n(f-b,\omega_\nu)\leq C$ which is equivalent to
\begin{equation}
n(f_\nu-b,\omega_0)\leq C
\end{equation}
However, boundedness of the preimage set of $a$ of $f_\nu$ implies that for all but finitely many $\nu$ we have $n(f_\nu-a,\omega_0)=n(f-a,\omega_\nu)=0.$
By Theorem \ref{quasithm1} the family $\{f_\nu\}_{\nu\in N_2}$ is quasinormal in $\omega_0,$ thus contains a subsequence quasiregular in $\omega_0.$
This implies that that any sequence from the family in Eqn.(\ref{quasielva}) contains a subsequence quasiregular in $\omega_0$, thus quasiregularity at $c$, which is a contradiction. 	
This completes the proof.
\end{proof}

\section{Algebras of $q$-analytic functions}
The context of Banach and uniform algebras has many useful application in the theory of holomorphic functions, see e.g.\ Stout \cite{stout}, Gamelin \cite{gamelinbok} and Brouwder \cite{brouwderbok}.
If $q\in \Z_+,$ and $f,g$ are two $q$-analytic functions on a domain $\Omega,$ such that the order $q$ is exact for both $f$ and $g$,
then $\partial_{\bar{z}}^q (fg)\equiv 0$ if and only if $q=1.$ Hence, for the case of $q>1$ we do not
have algebras of $q$-analytic functions. In fact, the subalgebra of continuous functions on $\Omega$ that is generated by the
$2$-analytic functions is all of $C^\omega(\Omega)$ (by which we mean the space of real-analytic functions).
\begin{definition}
A {\em Banach algebra}\index{Banach algebra} is a commutative $\C$-algebra, $A$, with identity (a so called unital algebra)
that at the same time is a complex Banach space with respect to a norm $\norm{\cdot}$ that satisfies $\norm{xy}\leq \norm{x}\norm{y}$
for all $x,y\in A,$ and $\norm{1}=1$ (where $1$ in the left hand side denotes the identity in $A$ and $1$ in the right hand side denotes the 
identity in $\C$).
\end{definition}
\begin{definition}[See e.g.\ Stout \cite{stout}, Def.1.2.6]
A {\em uniform algebra}\index{Uniform algebra} on a compact Hausdorrf space $X$ is a uniformly closed, point-separating subalgebra with identity of the algebra $C(X).$
\end{definition}
\begin{theorem}[Stone-Weierstrass] Let $X$ be a compact Hausdorff space. Let
$A\subset C^0(X,\R)$ be a subalgebra, which contains the constant function $1$ (i.e\ a unital algebra). Assume $A$
separates the points of $X$, i.e. for any $p_1,p_2\in X$, with $p_1\neq p_2$, there exists $f\in A$ such
that $f(p_1)\neq f(p_2).$ Then $A$ is dense in $C^0(X,\R)$, in the norm topology.
\end{theorem}
The complex version where $X\subset\C$ and $A$ is a subalgerba of $C^0(X,\C)$, 
includes the condition that $A$ be a $C^*$-subalgebra
(which implies that $f\in A\Rightarrow \bar{f}\in A$). 
\begin{definition}
A $C^*$-algebra, $A$, is a Banach algebra over the field of complex numbers together with an involution operator
$*$ on $A$ satisfying $\norm{x^*x}=\norm{x}^2,$ for all $x\in A.$ 
\end{definition}
In the given context the involution is identified as conjugation.
Consider the case of the disc algebra, i.e.\ $X=\{\abs{z}\leq 1\},$ and the subalgebra
$A_1(X).$ Denote $\mathbb{T}=\{\abs{z}=1\}.$ Then the maximum principle for holomorphic functions continuous on the disc 
implies that the restriction map $\Phi: A_1(X)\to C^0(\mathbb{T}),$ $f\mapsto f|_{\mathbb{T}}$ is an isometric unital
algebra homomorphism. In particular, the range $\mbox{Ran}(\Phi)$ is a closed subalgebra of $C^0(\mathbb{T}),$
which contains the constant function $1$. Denote $\phi(\zeta):=\zeta,$ $\zeta\in\mathbb{T}.$
One can verify that $\Phi$ is an isomorphism which implies that 
the spectrum $\mbox{Spec}_{\mbox{Ran}(\Phi)}(\phi)=\mbox{Spec}_{A_1(X)}(\tilde{\phi}),$
where $\tilde{\phi}\in A_1(X),$ is the function $\tilde{\phi}(\zeta)=\zeta,$ $\zeta\in  \{\abs{z}\leq 1\}.$
This implies $\mbox{Spec}_{\mbox{Ran}(\Phi)}(\phi)=\{\abs{z}\leq 1\}.$
The standard example given by the uniform limits on $\{\abs{z}\leq 1\}$ of complex polynomials 
(which contains the constant $1$ and separates the points of $\{\abs{z}\leq 1\}$) is not a $C^*$-subalgebra, 
since it does not satisfy $f\in A\Rightarrow \bar{f}\in A$ (indeed, the uniform limit must be holomorphic on the interior of the disc). 
Similarly, the complex version of the Stone-Weierstrass theorem
is inapplicable to any kind of algebra that would contain a $q$-analytic function with a nonpolynomial analytic component but not its conjugate (which would of course not be a $k$-analytic function for any $k\in \Z_+$). 
We shall introduce a multiplication that allows for the construction of a kind of uniform algebra of $q$-analytic functions, for any fixed $q\in \Z_+.$
For $q=1$ this reduces precisely to the case of holomorphic functions. 
To do this we need some preparatory work.

\section{Behaviour of analytic components under uniform limits}
Ahern \& Bruna \cite{ahernbruna} prove the following.
\begin{theorem}\label{ahernbrunaprop}
Let $U\subset \C$ be an open subset and let $q\in \Z_+.$ 
Let $p_0\in U$ and let $r>0$ be any number such that $\{\abs{z-p_0}<r\}\subset U.$
Then for each pair of integers $l,m\in \Z_{\geq 0}$ there exists a constant $C_{lm}$
such that for any $f\in \mbox{PA}_q(U)\cap L^1(U)$ we have
\begin{equation}
\abs{\partial_z^l\partial^m_{\bar{z}}f(p_0)}\leq \frac{C_{lm}}{r^{l+m+2}}
\int_{\{\abs{z-p_0}<r\}} \abs{f(\zeta)}d\mu(\zeta)
\end{equation}
where $d\mu$ is the standard area measure in $\C.$
\end{theorem}
\begin{proof}
It suffices to prove the result for $p_0=0$ and $r=1$ because
the inequality will continue to hold after dilation and translation.
Assume $f\in \mbox{PA}_q(\{\abs{z}<R\}),$ for $R>1$
with representation $f(z)=\sum_{j=0}^{q-1} a_j(z)\bar{z}^j$
for holomorphic $a_j(z).$
Make the ansatz 
\begin{equation}
P_{l,m}(z)=z^m\bar{z}^l \sum_{j=0}^{q-1} c_j\abs{z}^{2\nu}
\end{equation}
where we shall determine the constants $c_j\in C,$ such that
\begin{equation}\label{ahernvill}
\partial_z^l\partial_{\bar{z}}^m f(0)=\int_{\{\abs{z}<1\}} f(z)P_{lm}(z)d\mu(z)
\end{equation}
so we can assume $m\leq q-1$ and $\partial_z^l\partial_{\bar{z}}^m f(0)=q! \partial_z^l a_q(0).$
This will prove the wanted result with $C_{lm}:=\sup_{\abs{z}<1} \abs{P_{l,m}(z)}.$
Denoting by $\hat{a}_j(n)$ the Fourier coefficients, we have the power series expansion
\begin{equation}
a_j(z)=\sum_{n=0}^\infty \hat{a}_j(n)z^n
\end{equation}
which gives with $\abs{z}\exp(iarg(z))=:r\exp(i\theta)$ respect to Eqn.(\ref{ahernvill})
\begin{multline}
\sum_{j,\nu=0}^{q-1} c_\nu \int_0^1 r^{1+j+m+l+2\nu} \int_0^{2\pi}\exp(i(q-j-l)\theta)\sum \hat{a}_j(n)r^n\exp(in\theta)\frac{d\theta}{2\pi}dr\\
=\sum_{j,\nu=0}^{q-1} c_\nu \hat{a}_j(j+l-m)\int_0^1 r^{1+j+m+l+2\nu +j+l-m}dr=\\
\sum_{j=0}^{q-1}\left(\sum_{\nu=0}^{q-1} \frac{c_\nu}{2(1+j+l+\nu)}\right)\hat{a}_j(j+l-m)
\end{multline}
For a fixed $l\geq 0$, invertibility of the following Hilbert matrix
\begin{equation}
(H_{j\nu}):= \left(\frac{1}{1+j+\nu+l}\right)_{j,\nu=0}^{k-1}
\end{equation}
implies that there are numbers $\{c_\nu\}$ such that
\begin{equation}
\sum_{\nu=0}^{q-1} \frac{c_\nu}{1+j+\nu+l} =\left\{
\begin{array}{ll}
2m!l! & , j=q\\
0 & , j\neq q
\end{array}
\right.
\end{equation}
This completes the proof.
\end{proof}
Applying the lemma to a disc with radius $\mbox{dist}(p_0),\partial U)$ gives.
\begin{corollary}
Let $\Omega\subset\C$ be an open subset. Then there exists a constant $C=C_{lm}$ such that for any $f\in \mbox{PA}_q(U)\cap L^1(U)$, 
$p_0\in \Omega,$ 
\begin{equation}
\partial_z^l\partial_{\bar{z}}^m f(p_0) \leq \frac{C}{(\mbox{dist}(p_0),\partial U))^{k+l+2}}\norm{f}_1
\end{equation}
\end{corollary}
Another result that follows immediately is Corollary \ref{balk206b} stating that if
$\Omega\subset\C$ is a domain $q\in \Z_+,$
and we have a sequence $\{f_j\}_{j\in \Z_+}$ of $q$-analytic functions converging uniformly inside $\Omega$
then the limit function is also $q$-analytic in $\Omega.$
We mention that this result extends to $\alpha$-analytic functions in $\Cn$, in light of the 
$\alpha$-analytic version of 
Hartog's separate analyticity theorem (see Theorem \ref{hartog1}).
\begin{corollary}
Let $\Omega\subset\Cn$, $\alpha\in \Z_+^n$ and let $\{P_k\}_{k\in \Z_+}$ be a sequence of $\alpha$-analytic functions converging uniformly to a function
$f$ on $\Omega.$ Then $f$ is $\alpha$-analytic.
\end{corollary} 
\begin{proof}
For each $c\in \C^{n-1},$ set 
\begin{equation}
M_{j,c}:=\{ z: z_j\in \C, z_l=c_l, l<j, z_l=c_{l-1}, l>j\}
\end{equation}
On each $M_j\cap \Omega$ we have that the restrictions of each $P_k$ is
$\alpha_j$-analytic with respect to $z_j$, thus these restrictions converge uniformly on $M_j\cap \Omega$
to an $\alpha_j$-analytic function with respect to $z_j.$ Since $f$ is separately $\alpha$-analytic, applying Theorem \ref{hartog1} we obtain the wanted result.
\end{proof}     
We have the following. 
   \begin{corollary}\label{ahernbrunakonsekvens}
          Let $\Omega\subset \C$ be a bounded simply connected domain bounded by a rectifiable Jordan curve and let $q\in \Z_+.$
          Let $\{f_j\}_{j\in \Z_+}$ be a sequence of continuous functions on the closure $K:=\Omega\cup\partial\Omega$ such that 
          each $f_j$ is of the form $f_j(z)=\sum_{k=0}^{q-1}\bar{z}^k a_{j,k}(z),$
          $a_{j,k}\in A_1(K)$ (here $\bar{z}^j$ denotes the restriction of $\bar{z}^j$ to $K$) and
          satisfying $f_j\to f$ uniformly on $K$. Then $f$ is a
          continuous function on $K$, $q$-analytic on $\Omega,$ and of the form $f(z)=\sum_{k=0}^{q-1}\bar{z}^k g_{k}(z),$
          for $g_{k}\in A_1(K)$, $k=0,\ldots,q-1.$
          \end{corollary}
         \begin{proof}
                   Recall that each member of $A_1(K)$ is automatically a bounded holomorphic function on $\Omega$ and has unique continuous boundary values, when $K$ is the closed unit disc $A_1(K)$ is the standard disc algebra.
                   That $f$ is continuous is immediate from the fact that $C^0(K)$ is closed under the sup-norm. 
                   That $f$ is $q$-analytic on $\Omega,$ follows from Corollary \ref{balk206b} (which as we have stated follows from  Proposition \ref{ahernbrunaprop}).
                   By Proposition \ref{ahernbrunaprop} (alternatively using Corollary \ref{balk206b})
                   we have for each $m=0,\ldots,q-1,$
                   \begin{equation}\label{crocat}
                   \partial_{\bar{z}}^m f_j\to \partial_{\bar{z}}^m f\mbox{ uniformly on }\Omega
                   \end{equation}
                   Indeed, either we apply Corollary \ref{balk206b} directly and Eqn.(\ref{crocat}) is immediate, or we note that, since $f_j-f$ has continuous extension to the boundary it is $L^1$ and
                   we use the fact that for each $z_0\in \Omega$ there is an $\epsilon_0>0$ such that
                   $\{\abs{z-z_0}<\epsilon_0\}\subset \Omega,$ so by Proposition \ref{ahernbrunaprop}
                   \begin{equation}
                   \abs{\partial_{\bar{z}}^m(f_j-f)(z_0)}\leq \frac{C_{m}}{r^{m+2}}
                   \int_{\{\abs{z-p_0}<r\}} \abs{(f_j-f)(\zeta)}d\mu(\zeta)
                   \end{equation}
                   and since $f_j\to f$ uniformly on $K$ Eqn.(\ref{crocat}) follows.                                                      
                   Also note that the fact that each of the holomorphic components $a_{j,k}(z)$ of $\Omega$ have (since they belong to $A_1(K)$ by definition) continuous extension to $K$, implies that the same holds true
                   for $\partial_{\bar{z}}^m f_j=\sum_{k=0} a_{j,k}(z)(\partial_{\bar{z}}^m \bar{z}^k)$ which obviously will have continuous extension to $K$ as soon as the $a_{j,k}(z)$ do.
                   Now $\{a_{j,k}(z)\}_{j\in \Z_+}$ is a sequence in $A_1(K)$, and the latter space is closed with respect to the sup-norm.
                   We claim that if the uniform limit $g_k:=\lim_j a_{j,k}$, on $\Omega$, exists then the limit 
                   also extends continuously to the boundary, i.e.\ $g_k\in A_1(K).$
                   (If we merely had $a_{j,k}\in H^\infty(\Omega)$ this is false, take e.g.\ the $a_{j,k}$
                   to be a family converging to a Riemann map of $\Omega$ onto a bounded domain with discontinuous boundary).
                   For fixed $k$ denote each boundary function $h_{j}(\zeta):=a_{j,k}(\zeta),$ $\zeta\in \partial\Omega.$
                   Let $\epsilon>0$ be given.
                   Choose $N_\epsilon \in \Z_+$ sufficiently large such that for all $j,l>N_\epsilon$ we have
                   $\sup_{\Omega} \abs{(a_{j,k}-a_{l,k})(z)}<\epsilon.$
                   Since $a_{j,k}(z)-a_{l,k}(z)\in A_1(K),$ we have, by the maximum principle, for all $j,l>N_\epsilon,$
                   \begin{equation}
                   \epsilon> \sup_{\Omega} \abs{(a_{j,k}-a_{l,k})(z)}=\max_{\zeta\in \partial\Omega} \abs{(h_j-h_l)(\zeta)}
                   \end{equation}
                   Hence $\{ h_j\}_{j\in \Z_+}$ is a Cauchy sequence in $C^0(\partial\Omega)$ with respect to the sup-norm.
                   Since $C^0(\partial\Omega)$ is complete, we can (if necessary after replacing $\{ a_{j,k}\}_{j\in \Z_+}$
                   by a subsequence, keeping the same notation) conclude that
                   the uniform limit $\lim_j a_{j,k}(\zeta)|_{\partial\Omega}=\lim_j h_j(\zeta)$
                   exists as a member of $C^0(\partial\Omega)$ and can be identified with boundary
                   values of $g_k$ (the uniform limit on $\Omega$ of the $a_{j,k}$ as $j\to \infty$).
                   This proves the claim that $g_k$ can be identified as a member of $A_1(K).$
                   In particular, the uniform limit
                   $\frac{1}{(q-1)!}\partial_{\bar{z}}^{q-1} f_j\to g_{q-1}$, $j\to\infty,$ is a holomorphic function
                   on $\Omega$ coinciding with $\lim_j a_{j,q-1}$ (see Theorem \ref{hormkonvthm}) 
                   and has continuous extension to the boundary, thus we can identify $g_{q-1}$ as a member of $A_1(K)$.
                   This in turn can be iterated for $f_j$ replaced by $(f_j-\bar{z}^{q-1}g_{q-1})$
                   in order to obtain that $g_{q-2}$ is a holomorphic function
                   on $\Omega$ coinciding with $\lim_j a_{j,q-2}$ with continuous extension to $K$, and so on.
                   Hence $f$ can, on $\Omega$ be identified as
                   the function $f(z)=\sum_{k=0}^{q-1} g_k(z)\bar{z}^k$ for holomorphic
                   $g_k(z)$ such that each $g_k$ has continuous extension to $K.$ 
                   This completes the proof.
                   \end{proof}

\section{The algebra $\mathbf{PA}_{q}(K)$}
\begin{definition}
Let $U\subset\C$ be a bounded simply connected domain such that $\partial U$ is a rectifiable Jordan curve and set
$K=U\cup \partial U.$
Let $A_1(K)$ be the uniform algebra of functions holomorphic on the interior of
$K$ and continuous on $K.$ Let $\mbox{PA}_q(K)$\index{$\mbox{PA}_q(X)$} denote the
space of functions $f$ of the form $f(z)=\sum_{j=0}\bar{z}^j a_j(z),$
$a_j\in A_1(X)$ (here $\bar{z}^j$ denotes the restriction of $\bar{z}^j$ to $X$). In other words $\mbox{PA}_q(X)$ is the module over $A_1(X)$
generated by $\{1,\bar{z},\ldots,\bar{z}^{q-1}\}$ (notice that a compact $X$ is not an open subset
so this notation must be distinguished from the case when $X\subset\C$ is an open subset, where in the latter case 
this will coincide with the set of uniform limits of continuous functions that are $q$-analytic on $X$. 
Note specifically that there is a distinct difference in general for the case of open $X$ and compact $X$. For example when $X:=\{\abs{z}\leq 1\}$
the function $f(z):=(\bar{z}^2-1)/\sqrt{1-z}$ is continuous on $X$, polyanalytic on $X$ but its analytic components are not all continuous on $X$).
Denote by $\mathbf{PA}_{q}(K)$\index{$\mathbf{PA}_q(X)$} the space $\mbox{PA}_q(K)$ equipped 
with the operation, $\diamond_q ,$ of multiplication modulo $(\bar{z}^q)$ where $\bar{z}^q$ denotes the
element defined by the restriction
of $\bar{z}^q$ to $K,$ and $(\bar{z}^q)$ denotes the ideal generated by 
$(\bar{z}^q)$
in the space, $\mbox{CA}(K)$\index{$\mbox{CA}(K)$}, of functions countably analytic on int$K$ and continuous on $K.$ 
That is, if $f,g\in \mathbf{PA}_q(X),$ we define $f\diamond_q g =f(z)\cdot g(z)\mbox{ mod }(\bar{z}^q).$
\end{definition}
We can write down this multiplication more explicitly as follows.
Let $\Omega\subset\C$ be a bounded simply connected domain bounded by a rectifiable Jordan curve and let $q\in \Z_+.$
Introduce the 'multiplication' operation
\begin{equation}
\diamond_q\colon  \mbox{PA}_q(\Omega)\times \mbox{PA}_q(\Omega) \to \mbox{PA}_q(\Omega)
\end{equation}
as follows: Let $f,g\in \mbox{PA}_q(\Omega).$ By definition there exists holomorphic functions $a_j,b_j$ on 
$\Omega$ such that
$f=\sum_{j=0}^{q-1}a_j(z)\bar{z}^j,$ $g=\sum_{j=0}^{q-1}b_j(z)\bar{z}^j.$ Define 
\begin{equation}\label{equniformforsik}
(f\diamond_q g)(z):= \sum_{j=0}^{q-1} \left(\sum_{k+l=j} a_k(z) b_l(z)\right)\bar{z}^j
\end{equation}
Let $\mbox{PA}_{q,\diamond_q}(\Omega)$ denote the space consisting of the members of $\mbox{PA}_q(\Omega)$
equipped with the usual addition but with multiplication replaced by $\diamond_q$.
If $K\subset \C$ is the closure of a bounded simply connected domain, $U$, with boundary given by a rectifiable a Jordan curve, 
we have $\mathbf{PA}_{q}(K)=\mbox{PA}_{q,\diamond_q}(\mbox{int}K)\cap C^{q-1}(K)$
(by which we mean the set of functions 
that are $q$-analytic on $U$ such that their holomorphic components have
continuous extension to $K$, in the sense that they belong to $A_1(K)$). 
Obviously, we could instead define $\mbox{PA}_{q,\diamond_q}$ for compacts $K$
and then not have to worry about making sure to include the $C^{q-1}$-smoothness criterion.
In particular, if $f\in \mathbf{PA}_{q}(K)$, then on $U$, $f$ has the representation $f=\sum_{j=0}^{q-1}a_j(z)\bar{z}^j,$
for $a_j\in A_1(K)\subset C^0(K)$ and obviously $f\in C^0(K).$
The operation $\diamond_q$ is extended to $\mathbf{PA}_{q}(K)$ as follows: For $f,g\in \mathbf{PA}_{q}(K)$ define
\begin{equation}
f\diamond g(z):= 
\left\{
\begin{array}{ll}
(f\diamond g)(\zeta) & , z\in U\\
\lim_{U\ni \zeta\to z} (f\diamond g)(\zeta) & , z\in K\setminus U
\end{array}
\right.
\end{equation}
where the limit makes sense because each analytic component of the respective functions $f,g$ have continuous extension to $K$.
It is clear that
\begin{equation}
f\diamond_q g=g\diamond_q f
\end{equation}
Since $\mathbf{PA}_{q}(K)\subseteq C^0(K),$ we can define on $\mathbf{PA}_{q}(K)$ the sup-norm
\begin{equation}
\norm{f}_{\mathbf{PA}_{q}(K)}=\norm{f}_\infty=\sup_{z\in K} \abs{f(z)}
\end{equation}
By Theorem \ref{ahernbrunakonsekvens},
(alternatively Balk \cite{ca1}, Cor 4, p.206).
$\mathbf{PA}_{q}(K)$ is closed with respect to $\norm{\cdot }_{\mathbf{PA}_{q}(K)}$.
It is clear that $\diamond$ is distributive over addition and associative
since
\begin{multline}
f\diamond_q (g\diamond_q h)=f\diamond_q
\left(\sum_{j=0}^{q-1}\left(\sum_{k+l=j}b_kc_l\right)\bar{z}^j
\right)=\\
\sum_{j=0}^{q-1} \left(\sum_{m+n=j} a_n\sum_{k+l=m}b_kc_l\right)\bar{z}^j=
\sum_{j=0}^{q-1} \left(\sum_{k+l+n=j} a_n b_kc_l\right)\bar{z}^j
=(f\diamond_q g)\diamond_q h
\end{multline}
Let $U$ be a simply connected bounded domain as above with boundary given by a rectifiable Jordan curve and denote by $K$ the closure of $U.$
Obviously, for $q=1$, $\diamond_1$ is the usual multiplication $\cdot$ and $\mathbf{PA}_{q}(K)$ reduces to
the uniform algebra of functions holomorphic on $U$ and continuous on $K.$
Let us look at the properties of $\diamond_q$ in the more general case.
\begin{definition}
Given a $q$-analytic function $f=\sum_{j=0}^{q-1}a_j(z)\bar{z}^j,$
we denote by
ord$(f)$ the nonnegative integer $\max \{ j\colon a_j\not\equiv 0\}.$
\end{definition}

Let $f,g\in \mathbf{PA}_{q}(K)$ such that on $U$,
$f=\sum_{j=0}^{q-1}a_j(z)\bar{z}^j,$ $g=\sum_{j=0}^{q-1}b_j(z)\bar{z}^j,$ 
$h=\sum_{j=0}^{q-1}c_j(z)\bar{z}^j$
for
holomorphic functions $a_j,b_j,c_j$ on 
$U$. Let $n_1:=\mbox{ord}(f),$ $n_2=\mbox{ord}(g),$
$n_3=\mbox{ord}(h).$
If $n_1=0,$
we have for $z\in U$
\begin{equation}(f\diamond_q g)(z)=\sum_{j=0}^{q-1} \left(\sum_{k+l=j} a_k(z) b_l(z)\right)\bar{z}^j=
\sum_{j=0}^{q-1} a_0(z) b_j(z)\bar{z}^j=fg
\end{equation}
In particular, the space $\mathbf{PA}_{q}(K)$
has a unit with respect to $\diamond_q$ given by the usual multiplicative unit $g(z)\equiv 1.$
Obviously, we have for a complex constant $k_1,$
\begin{equation}(h\diamond_q k_1 g)(z)=
(k_1 h\diamond_q g)(z)=k_1 (h\diamond_q g)(z)
\end{equation}
Also
\begin{multline}
(f\diamond_q (g+h))(z)=
\sum_{j=0}^{q-1} \left(\sum_{k+l=j} a_k(z) (b_l(z)+c_l(z))\right)\bar{z}^j=\\
\sum_{j=0}^{q-1} \left(\sum_{k+l=j} a_k(z) b_l(z)\right)\bar{z}^j +
\sum_{j=0}^{q-1} \left(\sum_{k+l=j} a_k(z) c_l(z)\right)\bar{z}^j=\\
(f\diamond_q g)(z) +(f\diamond_q h)(z)
\end{multline}
Hence we have an associative commutative algebra over $\C$.
However, we do not have a Banach algebra (note that we include automatically commutativity and the condition of the existence of a unit in the algebra, which is not always the case
in various literature, in such cases it is important to write unital and commutative Banach algebra).
because the condition $\norm{fg}\leq \norm{f}\norm{g}$ will not continue to hold true when we replace the multiplication by $\diamond_q.$
\begin{example}\label{banachcounter}
Let $K:=\{z\colon \abs{z-3/4}\leq 1/4\}.$
 There exists $f,g\in \mbox{PA}_2(K)$ such that $\norm{ f\diamond_2 g}_K>\norm{f}_K\cdot \norm{g}_K.$ Indeed, let $f(z):=z\bar{z},$ $g(z):=1-z\bar{z}.$
 Then $f\diamond_2 g=z\bar{z},$ so that $\norm{f}_K=1,$ $\norm{g}_K=1-1/4=3/4,$ $\norm{f\diamond_2 g}_K=1>1\cdot 3/4=\norm{f}_K\cdot \norm{g}_K$.
\end{example}

\begin{remark}\label{innanvarje}
Due to (the counter-) Example \ref{banachcounter} it makes sense to consider replacing, in the definition of uniform algebra, "Banach algebra" by "algebra with the standard addition, that is a Banach space under the sup-norm". 
\end{remark}

\begin{remark}
Note that the definition of a uniform algebra, $A$, requires that the algebra be a subalgebra of $C(X)$, so formally this would require
that the multiplication in $A$ coincide with that in $C(X).$ We do however have the following.
Let $X\subset \C$ be a compact. A function $P\in \mbox{PA}_q(\mbox{int}(X))$ is a polyanalytic polynomial
(in the sense that its analytic components are complex polynomials) if and only if
$P(z)$ is a "polyanalytic polynomial" with respect to the multiplication in 
$\mathbf{PA}_{q}(\mbox{int}(X))$, in the sense that for some $N,M\in \Z_+,$ and complex constants $c_j$
\begin{equation}\label{parepbf}
P(z)=\sum_{j=0}^N \left(\sum_{k=0}^{M} d_{j,k} \underbrace{(z \diamond_q \cdots \diamond_q z)}_{k-\mbox{times}}  \right) \underbrace{(\bar{z} \diamond_q \cdots \diamond_q \bar{z})}_{j-\mbox{times}} 
\end{equation}
This follows immediately from the fact that 
\begin{equation}
\underbrace{(z \diamond_q \cdots \diamond_q z)}_{k-\mbox{times}} =z^k,\quad \underbrace{(\bar{z} \diamond_q \cdots \diamond_q \bar{z})}_{j-\mbox{times}}=\bar{z}^j\mbox{ mod}\bar{z}^q
\end{equation}
Denote by $\mathbf{P}_q(X)$ the set of uniform limits on $X$, of $C^{q-1}(X)$-functions having the representation in Eqn.(\ref{parepbf}). 
Then
\begin{equation}
P_q(X)\cap \cap C^{q-1}(X)=\mathbf{P}_q(X)
\end{equation}
Similarly, it is easily realized that for a simply connected bounded domain $U\subset\C$ such that $\partial U$ is a rectifiable Jordan curve, we have for $X=U\cup\partial U,$
\begin{equation}
A_q(X)\cap C^{q-1}(X)=\mathbf{PA}_q(X)
\end{equation}
\end{remark}
\begin{definition}[Uniform algebra in the generalized sense]
A {\em uniform algebra in the generalized sense}, $A$, on a compact Hausdorrf space $X$ is a uniformly closed, point-separating (commutative and unital) algebra over $\C$, with identity and with the standard addition, that is a Banach space under the sup-norm
such that $A\subset C^0(X)$ (i.e.\ we do not require that the multiplication in $A$ coincide with that in $C^0(X)$).
\end{definition}

\begin{observation}
For a bounded simply connected domain $U\subset\C$ such that $\partial U$ is a rectifiable Jordan curve, we have for $K=U\cup\partial U,$ that
$\mathbf{PA}_{q}(K)$ is a uniform algebra in the generalized sense.
\end{observation}
\begin{proof}
We have already proved that 
$\mathbf{PA}_{q}(K)$ is a unital commutative function
algebra with the standard addition, contained in $C(X)$, closed with respect to the sup-norm and contains the constants. We
need to verify that it separates points of $K$, but this follows from the fact
that it contains the set of functions holomorphic on $U$ and continuous on $K$ and
that is sufficient to separate points. This proves the observation.
\end{proof}
This opens up the question on whether a spectral analysis can be interesting in the theory of polyanalytic functions.
The spectrum as we shall see can be identified as a subset of the complex dual space equipped with the so-called weak-star topology (see below).
This topic is 
only recently been considered so here we shall merely give some background on what is known for the case of $q=1$ and some remarks on the difficulties in the polyanalytic (nonanalytic cases), see Daghighi \& Gauthier \cite{daghighiuniform}.
\\
\\
Recall that the norm of a bounded linear functional $\phi$ on a Banach space (we denote $\phi\in B'$ where $B'$ is the complex dual), $B$, is given by 
\begin{equation}
\norm{\phi}:=\norm{\phi}_{B'}=\sum_{\{f\neq 0\}} \frac{\abs{\phi(f)}}{\norm{f}}
\end{equation}
\begin{proposition}
The complex dual $B'$ is also a Banach space.
\end{proposition}
\begin{proof} 
Indeed, $B'$ is a linear space
and with the given norm
\begin{equation}
\norm{\phi_1 +\phi_2}\leq \sum_{\{f\neq 0\}} \frac{\abs{\phi_1(f)}}{\norm{f}}+\sum_{\{f\neq 0\}} \frac{\abs{\phi_2(f)}}{\norm{f}}
\leq \norm{\phi_1}+\norm{\phi_2}
\end{equation}
Also if $\{\phi_j\}_{j\in \Z_+}$ is a Cauchy sequence in $B'$ then for each $f\in B$ we have\\
$\abs{\phi_j(f)-\phi_k(f)}\leq \norm{\phi_j-\phi_k}\norm{f}$ so the sequence $\{\phi_j(f)\}_{j\in \Z_+}$ is Cauchy,
so we can define $\phi(f)=\lim_{j\to\infty} \phi_j(f).$ Linearity of $\phi$ follows from that of the $\phi_j.$
If $N\in \Z_+$ is such that $\abs{\phi_j-\phi_k}<1$ for $j,k>N$ then for $f\in B,$
\begin{multline}
\abs{\phi(f)}\leq \abs{\phi(f)-\phi(_N(f)}+\abs{\phi_N(f)}\leq \\
\limsup_{n\to\infty} \norm{\phi_j +\phi_N}\norm{f} +\norm{\phi_N}\norm{f}\leq (1+\norm{\phi_N})\norm{f}
\end{multline}
thus $\phi\in B'$. Also given $\epsilon >0$ we can choose $N$ sufficiently large such that $\norm{\phi_j -\phi_k}<\epsilon$ for $j,k>N$ so that for $f\in B$
\begin{equation}
\abs{(\phi-\phi_j)(f)}\leq \abs{(\phi-\phi_k)(f)}+\abs{(\phi_k-\phi_j)(f)}\leq 
\abs{(\phi +\phi_k)(f)} +\epsilon\norm{f}
\end{equation}
Since $\lim \abs{(\phi-\phi_j)(f)}=0$ we have $\norm{\phi -\phi_k}<\epsilon$ thus 
the arbitrary Cauchy sequence $\{\phi_j(f)\}_{j\in \Z_+}$ converges to a $\phi\in B'$ so $B'$ is complete.
This completes the proof.
\end{proof}

\begin{definition}
Let $A$ be a (commutative unital) Banach algebra. Denote by $A^{-1}$ the multiplicative group of invertible elements
in $A$ (i.e.\ $f\in A^{-1}$ if there exists $g\in A$ such that with the multiplication in $A$ we have $f\cdot g=1$).
We denote by $f^{-1}$ the inverse element of $f$.
A complex number is said to belong to the {\em resolvent}\index{Resolvent} of $f\in A$ if $\lambda\cdot 1 -f$ is invertible. 
The set of complex $\lambda$ for which $\lambda\cdot 1 -f$ is not invertible, is called the {\em spectrum} of $f$ (in $A$)
and is denoted $\sigma_A(f).$ 
\end{definition}

For our purposes (see Remark \ref{innanvarje}) it makes sense to consider the following.
\begin{definition}
Let $A$ be a uniform algebra in the generalized sense. Denote by $A^{-1}$ the multiplicative group of invertible elements
in $A$ (i.e.\ $f\in A^{-1}$ if there exists $g\in A$ such that with the multiplication in $A$ we have $f\cdot g=1$).
We denote by $f^{-1}$ the inverse element of $f$.
A complex number is said to belong to the {\em resolvent}\index{Resolvent} of $f\in A$ if $\lambda\cdot 1 -f$ is invertible. 
The set of complex $\lambda$ for which $\lambda\cdot 1 -f$ is not invertible, is called the {\em spectrum} of $f$ (in $A$)
and is denoted $\sigma_A(f).$ 
\end{definition}

\begin{proposition}
Let $A$ be a (commutative unital) Banach algebra. Then $A^{-1}$ is open in $A.$
\end{proposition}
\begin{proof} If $\norm{a-1}<1$ then
comparison to the Geometric series shows that
$\sum_{j=0}^\infty (-1)^j (a-1)^j$ converges in $A$. So $a=1+(a-1)$ shows that
the series represents $a^{-1}$ hence $a^{-1}\in A^{-1}.$ Also if $b\in A^{-1}$
then for each $x\in A$ we have
$x=(b(b^{-1}x-1)+1).$ Since $b$ is invertible and for $x$ sufficiently near $b$
$b^{-1}x-1$ has norm $<1$ we have $(b^{-1}x-1)+1\in A^{-1}.$ Thus for all $x$ sufficiently near $b$,
$x\in A^{-1}.$ This completes the proof.
\end{proof}
Note that the first part of the proof does not verbatim carry over to
uniform algebras in the generalized sense. Take e.g.\
$q=2$,
$K=\{\abs{z}\leq \frac{1}{2}\}$ and $b=\overline z.$ Note that $\overline z^0=1$ and for $k=1,2\ldots,$
$$
\underbrace{(\overline{z} \diamond_q \cdots \diamond_q \overline{z})}_{2k-\mbox{times}}=0
\quad \mbox{and} \quad
\underbrace{(\overline{z} \diamond_q \cdots \diamond_q \overline{z})}_{2k-1-\mbox{times}}=\overline{z}.
$$
Thus, 
$$
\sum_{j=0}^{2n} (-1)^j\underbrace{(\overline{z} \diamond_q \cdots \diamond_q \overline{z})}_{j-\mbox{times}}=1-n\overline{z}
\quad \mbox{and} \quad
\sum_{j=0}^{2n-1} (-1)^j\underbrace{(\overline{z} \diamond_q \cdots \diamond_q \overline{z})}_{j-\mbox{times}}=1-(n-1)\overline{z}.
$$
Hence, the series
$$
\sum_{j=0}^\infty (-1)^j\underbrace{(\overline{z} \diamond_q \cdots \diamond_q \overline{z})}_{j-\mbox{times}}
$$
diverges, and this is
$$
\sum_{j=0}^\infty (-1)^j 
\underbrace{(b \diamond_q \cdots \diamond_q b)}_{j-\mbox{times}}
$$ 
for a $b$ with $\norm{b}<1.$

\begin{theorem}[Gelfand-Mazur theorem]
If a (commutative unital) Banach algebra, $A$, is a division ring then it is isomorphic to the field $\C.$
\end{theorem}
\begin{proof}
That $A$ is a division ring implies that for each $a\in A$ 
and $\lambda_a \in \sigma_A(a)$ then, because $\lambda_a -a$ is noninvertible, we have $\lambda_a \cdot 1 =a,$ so $\lambda_a$ is unique (an isolated point in $\C$). 
The map $A\to \C,$
$a\mapsto \lambda_a$ is an isomorphism of $\C$-algebras.
This completes the proof.
\end{proof}

\begin{theorem}\label{paunifthm45}
Let $A$ be a (commutative unital) Banach algebra. For each $f\in A,$ $\sigma_A(f)$ is non-empty and compact and 
the inverse element $(\lambda\cdot 1 -f)^{-1}$ depends analytically on $\lambda$ for $\lambda$ belonging to the resolvent of $f.$
If $\lambda\in \sigma(f)$ then $\abs{\lambda}\leq \norm{f}.$ If $\lambda$ belongs to the resolvent
of $f$ then $\mbox{dist}(\lambda,\sigma(f))\geq \norm{(\lambda\cdot 1-f)^{-1}}^{-1}.$
\end{theorem}
\begin{proof}
For $\abs{\lambda}>\norm{f},$ the series $\sum_{j=0}^\infty f^j/\lambda^{j+1}$ converges to a function $g(\lambda)$
analytic at $\infty$ satisfying $g(\lambda)(\lambda \cdot 1-f)=1.$ Since $g(\lambda)=(\lambda\cdot 1-f)^{-1}$ we have 
$\sigma_A(f)\subset \{\abs{z}\leq \norm{f}\},$ thus is compact. Let $\lambda_0$ belong to the resolvent of $f$.
Then $\sum_{j=0}^\infty (\lambda_0\cdot 1 -\lambda\cdot 1)^j(\lambda_0\cdot 1 -f)^{-(j+1)}$ converges to a function $h(\lambda)$ 
analytic in $\{ \abs{\lambda-\lambda_0}< \norm{(\lambda_0\cdot 1 -f)^{-1}}^{-1} \}$ satisfying $h(\lambda)=(\lambda\cdot 1 -f)^{-1}.$ Thus
the resolvent of $f$ is open and $(\lambda\cdot 1 -f)^{-1}$ is analytic on the resolvent of $f.$ For each $\phi$ be a continuous linear functional on $A$ we have that 
$\phi((\lambda\cdot 1 -f)^{-1})$ is analytic with respect to $\lambda$ on the resolvent of $f$ and
vanishes at $\infty.$ By Lioville's theorem $\sigma(f)=\emptyset$ would imply $\phi\equiv 0$ which by the Hahn-Banach theorem would imply 
$(\lambda\cdot 1 -f)^{-1}\equiv 0$ which is impossible. Hence $\sigma(f)\neq \emptyset$.
This completes the proof.
\end{proof}
Let us look at the case when $A$ in Theorem \ref{paunifthm45} is replaced by the uniform algebra in the generalized sense, $\mathbf{PA}_q(K)$, for a fixed $q\in \Z_+,$ and a $K:=\overline{U}$ for a
bounded simply connected domain
$U$ whose boundary is a rectifiable Jordan curve.
The proof breaks down because we have that the sum 
$\sum_{j=0}^\infty f^j/\lambda^{j+1}$ may not be $q$-analytic with the usual multiplication and
may not converge with the multiplication $\diamond_q.$
Let us look at the case when $A$ in Theorem \ref{paunifthm45} is replaced by the uniform algebra in the generalized sense, $\mathbf{PA}_q(K)$, for a fixed $q\in \Z_+,$ and a $K:=\overline{U}$ for a
bounded simply connected domain
$U$ whose boundary is a rectifiable Jordan curve.
The proof breaks down because we have that the sum 
$\sum_{j=0}^\infty f^j/\lambda^{j+1}$ may not be $q$-analytic with the usual multiplication and
may not converge with the multiplication $\diamond_q.$
Let us look at what the inverse of a member $f\in \mathbf{PA}_q(K)$
will look like on $U$.
Let $\lambda\in \C$. 
As we would expect, it is easy to verify elements with point spectrum.
Let $q=2$ and let $f(z):=1+\bar{z}.$
Let $\lambda\neq 1, \lambda\in \C$. Then $g(z)\equiv 1-\bar{z}$ satisfies
$(1+\bar{z})\diamond_2 (1-\bar{z})=1+\bar{z}-\bar{z}=1.$
Furthermore, we may set $h(z):=\frac{1}{\lambda-1}\left(1-\bar{z} c_1(z)\right)$ and solve for $c_1(z)$ in
\begin{equation}
1\equiv \frac{1}{\lambda-1}\left(1-\bar{z} c_1(z)\right)\diamond_2 (\lambda-1-\bar{z})=\frac{1}{\lambda -1}(\lambda -1-\bar{z}+(1-\lambda)\bar{z}c_1(z))
\end{equation}
which gives $c_1(z)=(1-\lambda)^{-1},$ thus
\begin{equation}
h(z):=\frac{1}{\lambda-1}-\frac{\bar{z}}{(\lambda -1)(1-\lambda)}\Rightarrow (\lambda-f)\diamond_2 h(z)\equiv 1
\end{equation}
Also if $\lambda=1$ then $1-f(z)=-\bar{z}$ then the equation $(c_0(z)+\bar{z}c_1(z))\diamond_2 (1-f(z))\equiv 1$ has no solution
(since that would require $-c_1(z)\bar{z}\equiv 1$).
Thus $\{1\}=\sigma_{\mathbf{PA}_q(K)}(f).$
\\
Let us go back to the more general case of $q$.
We have on $U$ a representation of the form
$f(z)=\sum_{j=0}^{q-1}a_j(z)\bar{z}^j$ for holomorphic $a_j(z)$ (which also have continuous extension up to the boundary).
Suppose $g(z)\in \mathbf{PA}_q(K)$ such that
$g(z)\diamond_q f(z)\equiv 1.$ Let $g(z)=\sum_{j=0}^{q-1}b_j(z)\bar{z}^j$ on $U$,
for holomorphic $b_j(z)$ on $U$ (which also have continuous extension up to the boundary). By Eqn.(\ref{equniformforsik}) we have
$1\equiv \sum_{j=0}^{q-1} \left(\sum_{k+l=j} a_k(z) b_l(z)\right)\bar{z}^j.$
Let us suppose $a_0(z)$ is nonzero on $\overline{U}.$
By uniqueness of analytic coefficients this gives
\begin{equation}
b_0(z)\equiv \frac{1}{a_0(z)},\quad b_j(z)\equiv -\frac{1}{a_0(z)}\sum_{k+l=j,k>0} a_k(z) b_l(z),  \quad j=1,\ldots,q-1
\end{equation}
that is, on $U$ we have
\begin{equation}
g(z)\equiv\frac{1}{a_0(z)}\left(1-\sum_{j=1}^q\bar{z}^j\sum_{k+l=j,k>0} a_k(z) b_l(z)\right)
\end{equation}
Similarly, if $\lambda\in \C$ is such that $h:=(\lambda\cdot 1-f)^{-1}$ exists in 
$\mathbf{PA}_q(K)$, with representation on $U$ given by $\sum_{j=0}^{q-1}c_j(z)\bar{z}^j$.
Now $\lambda-f=(\lambda-a_0(z))+\sum_{j=1}^{q-1}a_j(z)$. Suppose 
\begin{equation}
\lambda\notin a_0(U)
\end{equation} 
Then (since $\lambda-f=(\lambda-a_0(z))+\sum_{j=1}^{q-1}a_j(z)$) we may repeat the derivation of the expression for 
the case when $g$ existed in order to obtain
\begin{equation}\label{cjekvannen}
c_0(z)\equiv \frac{1}{\lambda-a_0(z)},\quad c_j(z)\equiv -\frac{1}{\lambda-a_0(z)}\sum_{k+l=j,k>0} a_k(z) c_l(z), \quad j=1,\ldots,q-1
\end{equation}
and
\begin{equation}
h(z)\equiv\frac{1}{\lambda-a_0(z)}\left(1-\sum_{j=1}^q\bar{z}^j\sum_{k+l=j,k>0} a_k(z) c_l(z)\right)
\end{equation}
for holomorphic $c_j(z)$ on $U$ (which also have continuous extension up to the boundary).
Since all involved analytic components of $h$ together with all $a_j,c_j$, $j=0,\ldots,q-1$, 
have continuous extension to $K$ we cannot have that $a_0(z)=\lambda$ for $z\in \partial U$. 
This yields an expression of the inverses $h=(\lambda-f)^{-1}$ for $\lambda$ belonging to the resolvent
such that $\lambda\notin a_0(K).$
We have the following.
\begin{proposition}
	Consider the uniform algebra in the generalized sense, $\mathbf{PA}_q(K)$, for a fixed $q\in \Z_+,$ and $K=\overline{U}$ for a
	bounded simply connected domain
	$U$ whose boundary is a rectifiable Jordan curve. Let $f(z)$ have on $U$ the representation $\sum_{j=0}^{q-1}a_0(z)\bar{z}^j$, where each $a_j(z)$ is holomorphic and has continuous extension to $K.$
	Then $\sigma_{\mathbf{PA}_q(K)}(f)\subseteq a_0(K)$.
	In particular the spectrum is closed and bounded, i.e.\ compact.
	Furthermore, $(\lambda -f)^{-1}$ is an analytic function of $\lambda$
	on $\{\C\setminus a_0(K)\}.$
\end{proposition}
\begin{proof}
Let $f\in \mathbf{PA}_q(K)$ with representation on $U$ given by $f(z)=\sum_{j=0}^{q-1}a_j(z)\bar{z}^j$, $a_0\not\equiv 0,$ and let 
\begin{equation}
\lambda\notin a_0(K)
\end{equation}
Let $h\in \mathbf{PA}_q(K)$, with representation on $U$ given by $\sum_{j=0}^{q-1}c_j(z)\bar{z}^j$. 
Then it is sufficient for $h:=(\lambda\cdot 1-f)^{-1}$ that we are able to solve for $c_0(z),\ldots,c_{q-1}(z)$ in 
the system given by (using uniqueness of analytic coefficients) Eqn.(\ref{equniformforsik}), i.e.\
$1\equiv \sum_{j=0}^{q-1} \left(\sum_{k+l=j} \tilde{a}_k(z) c_l(z)\right)\bar{z}^j,$
with 
\begin{equation}
\tilde{a}_0(z):=\lambda -a_0(z),\quad \tilde{a}_j(z):=a_j(z),\quad j=1,\ldots,q-1
\end{equation}
such that the $c_j(z)$ have continuous extension to $K.$
The system can, for $z\in U$, be written 
\begin{equation}
\underbrace{\begin{bmatrix}
	\tilde{a}_0(z) & 0     & 0 & 0 &\cdots & 0\\
	\tilde{a}_{1}(z) & \tilde{a}_0(z)   & 0 &  0  &\cdots &  0\\
	\tilde{a}_{2}(z) & \tilde{a}_{1}(z)   & \tilde{a}_0(z) & 0   &\cdots &  0\\
	\vdots &  & & & \ddots & \vdots \\
	\tilde{a}_{q-1}(z) & \tilde{a}_{q-2}(z)   & \tilde{a}_{q-3}(z)     &   &\cdots   & \tilde{a}_0(z)\\
	\end{bmatrix}}_{=:A}
\begin{bmatrix}
c_0(z)\\
\vdots\\
c_{q-1}(z)
\end{bmatrix}=
\begin{bmatrix}
1\\
0\\
\vdots\\
0
\end{bmatrix}
\end{equation}
which is a lower triangular system with determinant (for any fixed $z$)
given by 
$(\lambda -a_0(z))^{q},$ which is nonzero since by assumption $\lambda\notin a_0(K).$ Hence, for each $\lambda\notin a_0(K)$
the determinant is nonzero and we may solve for the $c_j(z)$
in terms of the $a_j(z)$, $j=0,\ldots,q-1$ (and $\lambda$).
Next note that the explicit calculation of the
matrix inverse
yields for 
\begin{equation}
A^{-1}=\frac{1}{\mbox{det}(A)}\mbox{adj}(A)=(\lambda -a_0(z))^{-q}\mbox{adj}(A)
\end{equation}
where $\mbox{adj}(A)$ is the transpose of the cofactor matrix of $A$ and 
the elements of the cofactor matrix are finite sums of finite multiples of the
$a_j(z)$, $j=0,\ldots,q-1$ (and $\lambda$).
In particular, the resulting expression for $h$ on $U$ will be an analytic map with respect to $\lambda$ and furthermore each
$c_j(z)$ will have continuous extension with respect to $z$, to $K$ as soon as all the $a_j(z),$
$j=0,\ldots,q-1$ have continuous extension to $K$.
We have also obtained the implication
\begin{equation}
\lambda\notin a_0(K)\Rightarrow \lambda\notin \sigma_{\mbox{PA}_q(K)}(f) 
\end{equation}
which yields
\begin{equation}\label{spectruandraeqw}
\sigma_{\mathbf{PA}_q(K)}(f) \subseteq  a_0(K)
\end{equation}
This completes the proof.	
\end{proof}
Note that in the case $q=1$ we have $f(z)=a_0(z)$ so clearly
$\lambda\in \sigma_{\mathbf{PA}_1(K)}(f)$ implies that
$\abs{\lambda}\leq \max_{z\in K}\abs{f(z)}=\norm{f}.$
Since for a $q$-analytic $f(z)$ on a domain with representation $f(z)=\sum_{j=0}^{q-1}a_j(z)\bar{z}^j$
we have an associated function of two complex variables
\begin{equation}
F(z,w):=\sum_{j=0}^{q-1}a_j(z)w^j 
\end{equation}
This implies
$a_0(z)=F(z,0).$
We immediately have the following
\begin{corollary}
	Let $q\in \Z_+$ and $K=\overline{U}$ for a
	bounded simply connected domain
	$U$ whose boundary is a rectifiable Jordan curve
	and let $f\in \mathbf{PA}_q(K)$. Then
	$\sigma_{\mathbf{PA}_q(K)}(f)\subseteq \{\abs{z}\leq \norm{F(z,0)}\}.$
\end{corollary}
Investigating the Shilov boundary does not seem fruitful since given any
point $p_0\in C$ the entire $2$-analytic function $f(z)=1-\abs{z-p_0}^2$ peaks at $p_0.$
But we could look at character spaces.
\begin{definition}
A character\index{Character of a Banach algebra} $\varphi$ of a Banach algebra $F(+,\, '\cdot')$ (where $'\cdot'$ denotes the multiplication in $F$) 
is a $\C$-linear multiplicative functional on $F$,
i.e.\ $\varphi\in F^*$ (where $F^*$ denotes the complex dual) such that for all $f,g\in F,$
$\varphi(f\, '\cdot' g)=\varphi(f)\varphi(g).$ We denote the set of characters of $F$ by Char$(F).$
\end{definition}
In light of Remark \ref{innanvarje} our version of this definition will be the following.
\begin{definition}
A character $\varphi$ of a uniform algebra in the generalized sense, $F$,
is a $\C$-linear multiplicative functional on $F$,
i.e.\ $\varphi\in F^*$ (where $F^*$ denotes the complex dual) such that for all $f,g\in F,$
$\varphi(f \, '\cdot' g)=\varphi(f)\varphi(g).$ We denote the set of characters of $F$ by Char$(F).$
\end{definition}

A {\em net} is a function $\phi:I\to J$ between a directed ordered set $I$ and a topological space $A.$ By convention
the collection of the elements in the image is written with the induced order as $\{\phi_\tau\}_{\tau\in I},$ where each $\phi_\tau=\phi(\tau).$
A net $\{\phi_\tau\}_{\tau\in I},$ is said to converge to an element $t\in A$ if 
for each neighborhood $U$ of $t$, $\phi_\tau$ is eventually in $U.$
Recall that given a normed space $X$ the map $X\to X\to (X')',$ $x\mapsto \iota,$ $\iota_x(\phi)=\phi(x),$ is an embedding. The {\em weak-star
topology}\index{Weak-star-topology} on $X'$ is the coarsest topology such that the maps $\iota_x :X'\to \C$ are continuous.
A net $\{\phi_\tau\}_{\tau\in I}$ in $X'$ is said to  converge in the weak-star-topology if for all $x\in X,$ $\phi_\tau (x)\to \phi(x).$
Recall also that the kernel of each character is closed and a maximal ideal.  

\begin{lemma}\label{zornlemma}
Let $A$ be a commutative unital ring. Then every proper ideal is is contained in a maximal ideal and
an ideal $J$ is maximal if and only if $A/J$ is a field.
\end{lemma}
\begin{proof}
Let $I\subset A$ be a proper ideal and let $\mathcal{I}$ be the set of all proper ideals that contain $I.$
Then $\mathcal{I}$ can be equipped with the partial order $\subseteq.$ 
Let $\{I_\alpha\}_\alpha$ be a totally ordered subset of $\mathcal{I}$
and let $J=\cup_\alpha I_\alpha.$ Then $J$ is non-empty since $0\in J_\alpha$ for all $\alpha.$.
If $x,y\in J$ such that $x\in J_{\alpha_1}, y\in J_{\alpha_2}$ with $J_{\alpha_1}\subseteq J_{\alpha_2}$
then $x,y\in J_{\alpha_2},$ thus $x+y,x-y\in J_{\alpha_2}\subseteq J$ so $J$ forms an additive group. For $z\in J$, $z\in J_{\alpha}$
for some $\alpha$ and $a\in A$ then $az,za\in J_{\alpha}\subseteq J$ so $J$ is an ideal. If $J=A$ then $1\in J$ which implies $1\in J_\alpha$ for some $\alpha$
so that $J_\alpha =A$ which contradicts that $J_\alpha$ is proper. Also since $I\subset J_\alpha$
for all $\alpha$ we have $I\subset J.$ Hence by Zorn's lemma $\mathcal{I}$ has a maximal element.
Now let $J$ be a maximal ideal. 
Let $a+J$ be a nonzero element in $A/J.$ Then $a+J\neq J$ so $a\notin J.$ Thus $J$ is properly contained in $aA+J.$
Since $J$ is maximal this implies $aA +J=A.$ Thus there exists $b\in A$ and $c\in J$ such that
$ba+c=1$ thus $ba\equiv 1$mod$J$ so $b+J$ is the multiplicative inverse of $a+J,$ showing that $A/J$ is a field. 
Conversely suppose $A/J$ is a field. Let $a\in A\setminus J.$ Then $a+J\neq 0+J$ in $A/J.$
Thus $a+J$ has a multiplicative inverse, $b+J$ i.e.\ $(a+J)(b+J)=(1+J).$ This implies $ba+J=1+J$ or $1\in ba+J,$
thus $aA+J=A.$ Since $aA+J$ is the smallest ideal containing $J$ and $a$ it is maximal.
This completes the proof.
\end{proof}

\begin{theorem}\label{mazurclosed}
Let $A$ be a (commutative unital) Banach algebra. Every maximal ideal, $J$, of $A$ is closed and
$A/J$ is isometrically isomorphic to the field $\C.$
\end{theorem}
\begin{proof}
For each $f\in A$ such that $\norm{1-f}<1$ we have 
that $1$ belongs to the resolvent of $1-f$ and $f=1-(1-f) \in A^{-1}.$
If $I$ is a proper ideal and $f\in I$ then $f\notin A^{-1}$ thus $\norm{1-f}\geq 1.$ This inequality holds true for $f$ belonging
to the closure of $I$. Thus the closure of any proper ideal is a proper ideal thus a maximal ideal must be closed.
Now let $J$ be a maximal ideal. Since $J$ is closed by Lemma \ref{zornlemma} $A/J$ 
becomes a Banach space with norm given by $\norm{f+J}=\inf\{\norm{f+g}:g\in J\}.$
Since $\norm{fg +J}\leq \norm{f+J}\norm{g+J}$ for $f,g\in A,$ 
$A/J$ is commutative.
Since $\norm{1+g}\geq 1$ for all $g\in J$, $\norm{1+J}=1$ so that
$1+J$ is an identity for $A/J.$
By the Gelfand-Mazur theorem $A/J$ is isometrically isomorphic to the field $\C.$
This completes the proof.
\end{proof}

\begin{proposition}
Let $A$ be a uniform algebra (on a compact Hausdorff space $X$) 
and $a\in A.$ Then $a\in A^{-1}$ if and only if no character
of $A$ annihilates $a$.
\end{proposition}
\begin{proof}
Suppose that no character
of $A$ annihilates $a$. Assume (in order to reach a contradiction) that $a\in A^{-1}.$
The principal ideal $(a)$ generated by $a$ is contained in a maximal ideal, $m_a$, by Lemma \ref{zornlemma}.
$m_a$ is closed by Theorem \ref{mazurclosed}. $A/m_a$ is a field 
with the quotient norm
\begin{equation}\norm{h+m_a}_{A/m_a} =\inf\{\norm{h+f}_X \colon f\in m_a\}\end{equation}
By Theorem \ref{mazurclosed} the field is isomorphic to $\C.$ Hence the quotient map $A\mapsto A/m_a$ is a character
which annihilates $a.$ A contradiction.\\
Conversely let $a\in A^{-1}$ and assume there exists a character $\phi$ such that $\phi(a)=0.$
Then $\phi(1)=\phi(a)\phi(a^{-1})=0$ a contradiction.
This completes the proof.
\end{proof}
\begin{proposition}
For a compact Hausdorff space $X$ and $\phi\in$Char$C^0(X)$ there exists a unique point $x\in X$ such that $\phi(f)=f(x)$ for each $f\in C^0(X).$
\end{proposition}
\begin{proof}
If $x$ exists then uniqueness is evident.
Assume (in order to reach a contradiction) that existence fails. Then compactness of $X$ implies that there exists 
$f_1,\ldots,f_m\in C^0(X)$ such that $\phi(f_j)=0$ and for each $x\in X,$ $f_j(x)\neq 0$ for some
$j\in \{1,\ldots,m\}.$ Set $g_j=f_j / \left(\sum_{j=1}^m f_j\bar{f}_j\right).$ 
Then
$1=\sum_{j=1}^n f_j g_j,$ so
$1=\phi(1)=\phi(\sum_{j=1}^n f_j g_j)=\sum_{j=1}^n \phi(f_j)\phi( g_j)=0,$ a contradiction. This proves, given $\phi,$ existence of $x$ 
such that $\phi(f)=f(x)$ for each $f\in C^0(X).$ This completes the proof.
\end{proof}
This implies that for a compact Hausdorff space $X$ and $\phi\in$Char$C^0(X)$
we have
\begin{equation}
\abs{\phi(f)}\leq \norm{f}_X
\end{equation}
in particular each character has norm $1$ since $\phi(1)=1.$

\begin{proposition}\label{normchar}
For a (unital and commutative) Banach algebra $A$ any character is continuous and has norm $1.$
\end{proposition}
\begin{proof}
Suppose $A$ is a (unital) Banach algebra. Let $\phi\in$Char$(A).$
Then $\norm{1}=1$ since $(\phi(1))^2=\phi(1),$ and $\phi(1)\neq 0$ for $\phi\not\equiv 0.$ 
If $\abs{\lambda}>\norm{f}$ then $\lambda\cdot 1-f$ is invertible so
$\phi(\lambda\cdot 1-f)\phi((\lambda\cdot 1-f)^{-1})=\phi(1)=1$ implies $\phi(\lambda\cdot 1-f)\neq 0$ or $\phi(f)\neq\lambda.$
Hence $\abs{\phi(f)}\leq \norm{f}.$ Since this holds for all $f\in A,$ $\phi$ is continuous and $\norm{\phi}\leq 1.$ Since $\phi(1)=1,$ $\norm{\phi}=1.$
This completes the proof.
\end{proof}

For the proof above to extend to the cases $q>1$ we would need that 
$\lambda\cdot 1-f$ is invertible whenever $\abs{\lambda}>\norm{f}$.
The difficulty for a uniform algebra in the generalized sense, $A$, on a compact Hausdorff space $X,$ is as before the fact that power series may not converge. Indeed, 
let $f\in A$ such that $\norm{f}_X<1$.
If $\phi\in$Char$(A)$ such that $\phi(f)=c$ for some $c$ with $\abs{c}\geq 1$
then $\norm{f}_X<1$ implies that there exists an $\epsilon>0$ such that $\abs{f(x)}<1-\epsilon$
for all $x\in X.$ If we wish to replicate the known techniques for the case $q=1$ we would like to consider the sum $\frac{1}{c}\sum_{j=0}^\infty (f/c)^j$ (i.e.\ with the usual multiplication), which in that case converges uniformly on $X$ to some $g\in A$ and since $1=(c-f)g$ this would imply $1=\phi(1)=(c-\phi(f)\phi(g))=0,$ a contradiction, and we would obtain that
for each $\phi\in$Char$(A)$ we have $\abs{\phi(f)}<1.$ 
However we cannot control the convergence of 
$\sum_{j=0}^\infty \overbrace{(f/c)\diamond_q\cdots\diamond_q (f/c)}^{j-\mbox{times}}$.
Hence a different type of analysis is necessary.

\begin{proposition}
Let $A$ be a unital Banach algebra and equip Char$(A)$ with the weak-star-topology.
Then Char$(A)$ is compact.
\end{proposition}
\begin{proof}
Since the unit ball in $A'$ (the complex dual) is compact (the Banach-Alaoglu theorem), in the weak-star-
topology it suffices to show that Char$(A)$ is weak-star-closed
in the unit ball in $A^*$. Take a net $\{\lambda_\iota\}_{\iota\in I}\subset$Char$(A)$, (for an index set $I$) 
and $\phi\in$Char$(A)$ with $\phi$ belonging to the unit ball in $A'$ 
such that $\phi:=\lim^*_{\iota\in I} \lambda_\iota$ (where $\lim^*$ denotes the weak-star-limit). By the definition of the weak-star topology 
$\phi(a)=\lim_{\iota\in I} \lambda_\iota(a),$ for all $a\in A.$ This implies that $\phi$ is linear and multiplicative and 
since $\phi(1)=\lim_{\iota\in I} \lambda_\iota(1)=1,$ we have $\phi\not\equiv 0,$ so that $\phi\in$Char$(A).$ Since the sequence was arbitrary with limit in the unit 
ball in $A'$, this completes the proof.
\end{proof}
We are interested in the case when a character $\phi$ is a multiplicative linear functional on $\mathbf{PA}_{q}(K)$ where 
the multiplication 
is given by $\diamond_q.$
For any $f\in C^{q}(K)$ the restriction of $f$ to $U$ satisfies that each $f_j:=\partial_{\bar{z}}^j f$ has continuous extension to $K.$
Recall that by Theorem \ref{theodorescuthm}, $f\in \mathbf{PA}_{q}(K)$ is given, on any subdomain $U'\subset U$, bounded by a rectifiable Jordan curve, by the formula 
\begin{equation}
f(z)=\frac{1}{2\pi i}\sum_{j=0}^{q-1} \int_{\partial U'} \frac{1}{j!(\zeta-z)}f_j(\zeta)(\bar{\zeta}-\abs{z})^j d\zeta
\end{equation}
and we have seen that each $\partial_{\bar{z}}^j f$ has continuous extension to the boundary iff
each analytic component does, so by the refinement following that theorem, if the
$f_j$ have continuous extension to $K$, the formula 
extends to the case where $U'$ is replaced by $U.$
We only need however the holomorphic case of this.
Denote by $(z,w)$ the Euclidean coordinate in $\C^2$ and for a compact $S\subset\C^2,$ denote 
by $A_1(S)$ the set of functions, 
holomorphic for $(z,w)\in \mbox{int}S$ and continuous on $S.$
For a subset $K\subset \C$ denote by $K^*:=\{\bar{z} \colon z\in K\}\subset \C.$   
We shall denote by $\widetilde{A_1(S)/(w^q)}$ the set of Weierstrass polynomials
$\sum_{j=0}^{q-1} a_j(z)w^j.$
\begin{proposition}\label{forstasaken}
Let $K$ be the closure of a bounded simply connected domain $U\ni 0$ with continuous boundary (given by a rectifiable Jordan curve).
Let $q\in \Z_+.$
Set $S:=K\times K^*\subset\C^2.$
Then
\begin{equation}
\mbox{Char}(\mathbf{PA}_{q}(K))=\{\varphi\in \mbox{Char}(A_1(S))\colon \varphi(w^q)=0\}=\mbox{Char}(A_1(S))|_{\widetilde{A_1(S)/(w^q)}}
\end{equation}
where $\mbox{Char}(A_1(S))|_{\widetilde{A_1(S)/(w^q)}}$ denotes, as above, the restriction of the characters of $A_1(S)$ 
to the representatives of the elements in the quotient by the ideal $(w^q)$ generated by the monomial $w^q.$
\end{proposition} 
\begin{proof}
Let $g(z,w)\in A_1(S).$  
Since $g$ is bounded and holomorphic on int$S$ and with continuous boundary values, 
we have that the generalized version of the Cauchy integral formula, (see Theorem \ref{theodorescuthmrefin}, alternatively H\"ormander \cite{hormander}, Thm 2.2.1, for the case of the unit disc,
the case of bounded simply connected domain with boundary given by a rectifiable Jordan curve follows from applying
the following generalization of the Riemann mapping theorem given by the Carath\'eodoy theorem) holds with respect to
the rectifiable Jordan curve $\partial (U^*),$ 
for $(z,w)\in K\times U^*$ 
\begin{equation}
g(z,w)=\frac{1}{2\pi i}\int_{\partial (U^*)}g(z,\zeta) \frac{1}{\zeta-w}d\zeta
\end{equation}
Since $w\neq \zeta$ for all $\zeta\in \partial (U^*)$ we can
Taylor expand the expression $(\zeta-w)^{-1}$ up to order $q$ in order to obtain functions
$a_j(z)\in A_1(K),$ $j=0,\ldots,q-1$, and $R(z,w)\in A_1(S),$ such that
\begin{equation}\label{weierguru}
g(z,w)=\frac{1}{2\pi i}\int_{\partial (U^*)}g(z,\zeta) \frac{1}{\zeta-w}d\zeta
=\sum_{j=0}^{q-1} a_j(z)w^j + w^q R(z,w)
\end{equation} 
Here $w^j$ denotes the restriction of the function $w^j$ to $K^*\subset \C,$ and can thus be identified as 
a function (independent of $z$) in $A_1(S)$,
in particular this implies that $w^q R(z,w)\in (w^q)$, the ideal generated by $w^q\in A_1(S)$.
(Note that locally in int$S$, Eqn.(\ref{weierguru}) is precisely the
Weierstrass division theorem in two complex variables, for the Weierstrass polynomial $w^q$ of order $q$).
If we introduce on $A_1(S),$ the equivalence relation $F\sim G \Leftrightarrow F-G\in (w^q),$
and denote the equivalence class of $F\in A_1(S)$ by $[F]$ then we know that each $[F]$ can be identified as an element of $\mathbf{PA}_{q}(K)$
which is then the natural quotient space under the given equivalence relation 
and this induces the wanted identification of character spaces. 
Let us give the details.
First of all we can
define the map $\Psi: (A_1(S))'\to (\mathbf{PA}_{q}(K))'$ as follows. 
For $F(z,w)\in A_1(S)$
and $\varphi\in (A_1(S))',$ 
we have 
$F=F_1+F_2,$ $F_1=\sum_{k=0}^{q-1} a_k(z)w^k,$ for $a_k\in A_1(K)$ and
$F_2\in (w^k)\subset A_1(S),$ so that
$\varphi=\varphi_1+\varphi_2,$ with
$\varphi_1(F)=\varphi(F_1),$ $\varphi_2(F)=\varphi(F_2).$
Define for each $\varphi\in (A_1(S))',$ 
\begin{equation}
(\mathbf{PA}_{q}(K))'\ni \psi(z):=\varphi_1(z,\bar{z}),\quad \Psi(\varphi):=\psi
\end{equation} 
where $(\mathbf{PA}_{q}(K))'$ denotes the complex dual space.
Then $\Psi$ is a linear mapping satisfying Ker$(\Psi)=(w^k)$ and the image of $\Psi$ is Im$(\Psi)=(\mathbf{PA}_{q}(K))'.$
$\Psi$ is a 
a bounded linear operator so we have a vector space isomorphism $(A_1(S))'/\mbox{Ker}\Psi = \mbox{Im}\Psi$
which shows $(A_1(S))'|_{\widetilde{A_1(S)/(w^q)}}=(\mathbf{PA}_{q}(K))'.$
We are interested in the subset 
$\mbox{Char}(\mathbf{PA}_{q}(K))\subset (\mathbf{PA}_{q}(K))'$
of multiplicative functionals. However, as we have seen,
each $\psi\in (\mathbf{PA}_{q}(K))'$ corresponds uniquely to
some $\varphi_1(z,\bar{z})$, where $\varphi_1(z,w)$ is the restriction of an element $\varphi=\varphi_1+\varphi_2\in (A_1(S))'$ to $\widetilde{A_1(S)/(w^q)}.$
So if $\psi$ is multiplicative, i.e.\
$\psi(f)\psi(g)=\psi(f\diamond_q g)=\psi(fg\mbox{ mod } \bar{z}^q),$
we have  
that $\varphi_1$ is multiplicative (with respect to the usual $\cdot$)
when restricted to $\widetilde{A_1(S)/(w^q)},$ hence has an extension to a multiplicative $\varphi$ in $\mbox{Char}(A_1(S)).$
This gives an injection $\mbox{Char}(\mathbf{PA}_{q}(K))\to 
\mbox{Char}(A_1(S))|_{A_1(S)/(w^q)}.$
Conversely, note that if we start from a multiplicative $\varphi\in \mbox{Char}(A_1(S))$ then
its restriction, $\varphi_1$, to $\widetilde{A_1(S)/(w^q)}$, must also be multiplicative, whence
the associated functional $\varphi_1(z,\bar{z})=\psi\in (\mathbf{PA}_{q}(K))'$
will multiplicative  and thereby, by definition, belong to $\mbox{Char}(\mathbf{PA}_{q}(K)).$
This completes the proof.
\end{proof}
Character spaces of spaces of type $A_1(S)$ have previously been examined. Indeed
since $\C^2\setminus S$ is connected and $S$ is compact we have
$P_1(S)=A_1(S)$ so it suffices to describe Char$P_1(S).$
\begin{theorem}[See e.g.\ Stout \cite{stout}, Thm 1.2.9]
If $S\subset\Cn$ is compact then every character of $P_1(S)$ is of the form $f\mapsto \hat{f}(z)$ for 
a unique $z\in \hat{S}.$
\end{theorem}
This completely characterizes the space Char$A_1(S)$ as the evaluations at points of the polynomially convex hull of $S.$ Theorem \ref{forstasaken} thus characterizes the
Char($\mathbf{PA}_q(K)$) in terms of the subset of such characters that annihilate $w^q.$
Now character spaces can be used in the study of the spectrum.
The Gelfand-Mazur theorem renders the following characterization of Char$(\mathbf{PA}_q(K)).$
\begin{proposition}\label{abtinssatsomuniform}
Let $K=U\cup\partial U$ for a bounded simply connected domain $U$ such that $\partial U$ is a rectifiable Jordan curve.
Let $S=K\times K^*\subset\C^2,$ (where $K^*:=\{w \in \C \colon \bar{w}\in U\}$) with Euclidean coordinates $(z,w)$ for $\C^2.$
Then 
\begin{equation}\label{rhsofeqnref}
\mbox{Char}(\mathbf{PA}_q(K))=\{A_1(S)/J\colon J\mbox{ a maximal ideal with }w^q\in J\}
\end{equation}
\end{proposition}
\begin{proof}
We have proved that $\mathbf{PA}_q(K)=A_1(S)/(w^q).$ By
the Gelfand-Mazur theorem any character of $A_1(S)$ takes the form $A_1(S)/J$ for a maximal ideal $J$ in $A_1(S).$
Thus the set in the right hand side of Eqn.(\ref{rhsofeqnref}) is precisely
the set of characters of $A_1(S)$ that annihilate $(w^q)$. By Proposition \ref{forstasaken}, this completes the proof.
\end{proof}
In a more general setting, if $Y$ is a set and $T$ a topological space and $F$
a family of functions $Y\to T.$ The {\em weak topology}\index{Weak topology}
on $Y$ induced by $F$ is the smallest (coarsest) topology, $\tau$, on $Y$
for which the functions of $F$ are continuous. Then $\tau$ is generated by (i.e.\ has neighborhood basis
at each point that belongs to) the sets
\begin{equation}
\{f^{-1}(U)\colon f\in F,U\mbox{ open in }T\}, \quad f^{-1}(U):=\{y\in Y\colon f(y)\in U\}
\end{equation}
Recall that convergence of a net $\{x_\alpha\}_\alpha$ is defined by $\lim_{\alpha} x_\alpha =x$ if and only if $\lim_\alpha f(x_\alpha)=f(x)$ for each $f\in F.$
If $T$ is Hausdorff and $F$ separates the points of $Y$ then the weak topology is Hausdorff.
\begin{definition}
Let $Y$ be a Banach space. For each $f\in Y$ denote by $\hat{f}$ the function on the complex dual $Y'$
defined by \begin{equation}\label{gelftrans}
\hat{f}(\phi)=\phi(f) 
\end{equation}
The weak-star topology\index{Weak-star topology} on $Y'$ is the weak topology on $Y'$ induced by the family $\{ \hat{f}\colon f\in Y\}.$
\end{definition}

\begin{proposition}\label{characthausdorf}
The weak-star topology of the complex dual $A'$ of a Banach space $A$, is Hausdorff.
\end{proposition}
\begin{proof}
Let $\phi_1,\phi_2\in A'$.
If $\phi_1\neq \phi_2$ there exists $f\in A$ such that $\phi_1(f)\neq \phi_2(A).$
Hence $\hat{f}(\phi_1)\neq \hat{f}(\phi_2)$ which implies that the functions $\{\hat{f}\colon f\in A\}$ separate the points of $A'$.
Since convergence of a net $\{x_\alpha\}_\alpha$ is defined by $\lim_{\alpha} x_\alpha =x$ if and only if $\lim_\alpha \phi(x_\alpha)=\phi(x)$ for each $\phi\in A',$
the weak topology is the topology of pointwise convergence and since $\C$ is Hausdorff the fact that
$A'$ separates the points of $A$ is sufficient for $A'$ to be Hausdorff.
\end{proof}

\begin{definition}
The spectrum $\Sigma(A)$\index{$\Sigma(A)$} of a uniform algebra in the generalized sense, $A$, on a compact Hausdorff
space $X$ is the set of all continuous characters of $A$ equipped with the weak-star topology
(so that a net $\{\phi_\iota\}_\iota$ in $\Sigma(A)$ converges to $\phi\in \Sigma(A)$ if for each $f\in A$
the net
$\{\phi_\iota (f)\}_\iota$ converges to $\phi(f)$).
\end{definition}
A neighborhood basis of $\phi\in \Sigma(A)$ is given by sets of the form
\begin{equation}
\{ \phi\in \Sigma(A)\colon \abs{\phi(f_k)-\psi(f_k)}<\epsilon, k=1,\ldots,r\}
\end{equation}
where $\epsilon$ ranges over $(0,\infty),$ $r$ ranges over $\Z_+$ and $f_j$ ranges over $A$.

\begin{proposition}
The spectrum $\Sigma(A)$ of a uniform algebra in the generalized sense $A$ on a compact Hausdorff
space $X$, is a compact Hausdorff space.
\end{proposition}
\begin{proof}
By Proposition \ref{characthausdorf} the space of characters with the weak-star topology is Hausdorff. We must prove compactness.
Denote for each $f\in A$, by $\C^f$, a copy of $\C$ and denote by $\C^A$ the cartesian product $\pi_{f\in A}\C^f.$ 
Set $\Phi:\Sigma(A)\to \C^A,$ $\Phi(\phi)=\{\phi(f)\colon f\in A\}.$ By the definition of the topologies
$\Phi$ is injective and continuous. Furthermore
$\Phi(\Sigma(A))$ is a subset of $\Pi_{f\in A} \{\zeta\colon \abs{\zeta}\leq \norm{f}_X\}$ which is compact since
by Tychonoff's theorem the product of compact sets is compact
in the product topology. 
Let $\{z_f:f\in A\}$ be a limit point of $\Phi(\Sigma(A))$ i.e.\ there exists a net $\{\phi_\alpha\}_\alpha$ in $\Sigma(A)$,
$\phi_\alpha(f)\to z_f$ for each $f\in A.$ Define $\varphi:A\to \C,$ $\varphi(f):=z_f.$ Then $\varphi(1)=1$ and $\varphi$ is $\C$-linear and multiplicative
so $\varphi\in \Sigma(f)$. Thus Ran$(\Phi)$ is closed in a compact so it is itself a compact space.
If $\{\Phi(\phi_\alpha)\}_\alpha$ converges to $\Phi(\varphi),$ then $\phi_\alpha(f)\to \varphi(f)$ thus $\phi_\alpha\to \varphi$ in $\Sigma(A).$ 
Hence $\Phi^{-1}$ is continuous proving that $\Phi$ is a homeomorphism which in turn implies that
the inverse image $\Sigma(A)$ of $\Phi(\Sigma(A))$ is compact. This completes the proof.
\end{proof}
Our work thus yields, by Proposition \ref{abtinssatsomuniform}, a characterization of the spectrum, 
$\Sigma(\mathbf{PA}_q(K))$, of the generalized uniform algebra $\mathbf{PA}_q(K),$ and that spectrum is proven to be a compact Hausdorff space with the weak-star topology
(the latter is well-known if one removes the word 'generalized', we have simply verified the proofs in our specific setting).
We mention that the map of Eqn.(\ref{gelftrans}) is called the Gel'fand transform\index{Gel'fand transform}.
If $A$ is a uniform algebra then the algebra $\hat{A}$ of Gel'fand transforms is a uniform algebra on $\Sigma(A).$
Also for each $f\in A$ we have $\sigma(A)=\hat{f}(\Sigma).$

\chapter{Some integral formulas and the maximum principle}

\section{Two-radius characterization}
Let $\Omega\subset\C$ be a domain.  
A version of Moreras theorem holds true with respect to circles,
namely assume $f\in C^\infty(\Omega)$ and $p\in \Omega$ 
and $\int_{\{\abs{z-p}=\epsilon\}} f(\zeta)d\zeta$ for sufficiently small $\epsilon.$
Then by Green's theorem
\begin{equation}
0=\int_{\{\abs{z-p}=\epsilon\}} f(\zeta)d\zeta =2i\int_{\{\abs{\zeta}\leq \epsilon\}} \partial_{\bar{z}} fd\mu(\zeta)
\end{equation}
for $\epsilon$ sufficiently small such that $\{\abs{\zeta}\leq \epsilon\}\subset\Omega.$ Letting 
$\epsilon\to 0,$
we obtain $\partial_{\bar{z}} f(p) =0.$ 
Denote by $J_j$ the $j$:th Bessel function of the first kind
\begin{equation}
J_j(z):=\left\{
\begin{array}{ll}
\sum_{k=0}^\infty \frac{(-1)^k}{2^{2k+j}k!(j+k)!}z^{2k+j}\quad & , j\geq \frac{1}{2}\\
(-1)^j J_{-j}(z), \quad & j=<0\\
\end{array}
\right.
\end{equation}
They are known to satisfy 
\begin{equation}
\sum_{j=-\infty}^\infty J_j(z)t^j =\exp\left( \frac{z}{2}\left(t-\frac{1}{t}\right)\right)
\end{equation}
where the series converges uniformly on compacts in $\C_z\times (\C_t\setminus \{0\}).$ 
It is known that the zeros the Bessel functions are all real. See e.g.\ Watson \cite{watson}.
It is possible to write
\begin{equation}
J_j(z)=\frac{1}{2\pi} \int_0^{2\pi} \exp(i(j\theta -z\sin\theta)) d\theta
\end{equation}
Let $Z_j$ denote the set of ratios of the (nonzero) zeros of $J_j.$
For $r>0,$ denote by $\mu_r$ the measure determined by
\begin{equation}
\int f(rt) d\mu(t)=\int f(z)d\mu_r(t),\quad f\in C^0(\R)
\end{equation}
Set for $k\in \N,$ $n\in \Z,$
\begin{equation}
A_{k,n}:=\int t^{2k+\abs{n}} d\mu(t)
\end{equation}
let $s_n$ be the smallest value of $k$ for which $A_{k,n}\neq 0.$
\begin{equation}
F(z)=\sum_{k=0}^\infty \frac{(-1)^k}{k!(k+\abs{n})!}\left(A_{k,n}\left(\frac{-z}{2}\right)\right)^{2k+\abs{n}}
\end{equation}
and denote by $Z(n,\mu)$ the set of ratios of the (nonzero) zeros of $F.$

Let $\mathcal{E}'(\Rn)$ denote 
the Schwartz space of distributions with compact support.
Recall that the Fourier transform of an element $f\in \mathcal{E}'(\Rn)$ is given by 
\begin{equation}
\hat{f}(\xi):=\langle f_z,\exp(iz\cot \xi) \rangle
\end{equation} 
Denote by $\mathbb{E}'$ the space of Fourier transforms of elements of $\mathcal{E}'$ 
The topology of $\mathbb{E}'$ is defined so that the Fourier transform becomes a topological isomorphism.
Let $K$ be the class of continuous positive functions on $\Cn$ of the form
$k(z) = k_1 (\norm{\re z}_{\ell^1})k_2 (\norm{\im z}_{\ell^1}),$ where $k_1$ and $k_2$ are increasing, 
$k_1$ grows faster than any polynomial, and $k_2$ grows faster than any linear exponential.
Then for $G\in \mathbb{E}',$ $\abs{G(z)}/k(z)\to 0$ as $z\to \infty.$ A neighborhood basis at $0$ 
for $\mathbb{E}'$ is given (see Ehrenpreis \cite{ehrenpreis}, p.156) by the sets
$N_k:=\{G\in \mathbb{E}':\sup_{\Cn} \abs{G(z)}/k(z) <1\},$ where $k\in K$.
The Fourier transform satisfies for $q\in \Z_+,$
\begin{equation}
\widehat{(\partial_z^q f)}=\frac{(-1)^q}{2^q}(z_2+iz_1)^q \hat{f}
\end{equation}
\begin{equation}
\widehat{(\partial_{\bar{z}}^q f)}=\frac{1}{2^q}(z_2-iz_1)^q \hat{f}
\end{equation}
and
\begin{equation}
\widehat{\Delta f}=-(z_2^2+iz_1^2) \hat{f}
\end{equation}
Zalcman \cite{zalcman} proved 
the following.
\begin{theorem}\label{zalcmanthm}
Let $f\in L^1_{\mbox{loc}}(\R^2)$ and suppose that there exists distinct positive $r_1,r_2$ such that
for almost all $z\in \C,$
\begin{equation}
\int_{0}^\infty \int_0^{2\pi} f(z+\rho\exp(i\theta))\exp(iq\theta) d\theta d\mu_{r_j}(\rho) =0, \quad j=1,2
\end{equation}
Then, if $r_1/r_2\notin S(n,\mu)$ and $n\geq 0,$ $f$ is a distribution solution to
\begin{equation}
\partial_{\bar{z}}^n \Delta^s f=0
\end{equation}
\end{theorem}
\begin{proof}
Assume $f$ is $C^\infty$-smooth. The Fourier transform $\hat{\sigma}_j,$
of the measure $\sigma_j:=\exp(in\theta)d\theta d\mu_{r_j}(\rho)$ is given by
\begin{equation}
\hat{\sigma}_j=2\pi(z_2-z_1)^n F(r_j\sqrt{z_1^2 +z_2^2})/(\sqrt{z_1^2 +z_2^2})^n
\end{equation}
Hence the functions
\begin{equation}
F_j=2\pi(z_2-z_1)^n F(r_j\sqrt{z_1^2 +z_2^2})/(\sqrt{z_1^2 +z_2^2})^n (z_1^2 +z_2^2)^s 
\end{equation}
belong to $\mathbb{E}'(\R^2)$ and have, by the conditions of the Theorem, no common zeros.
In the case of $n=1$ (i.e.\ distributions on $\R$), 
if two elements $f_1,f_2\in \mathbb{E}'(\R)$ have no common zeros, then (see Schwartz \cite{schwartz1}) 
the ideal generated by $f_1$ and $f_2$ is all of $\mathbb{E}'(\R).$ 
Hence we can find nets of functions $\{f_{1,\iota}\},$ $\{f_{2,\iota}\}$ in $\mathbb{E}'(\R)$ such that
in $\mathbb{E}'(\R),$
\begin{equation}
h_\iota=f_{1,\iota}(z) F_1(z) +f_{2,\iota}(z) F_2(z)\to 1
\end{equation}
As before let $K$ be the class of continuous positive functions on $\Cn$ of the form
$k(z) = k_1 (\norm{\re z}_{\ell^1})k_2 (\norm{\im z}_{\ell^1}),$ where $k_1$ and $k_2$ are increasing, 
$k_1$ grows faster than any polynomial, and $k_2$ grows faster than any linear exponential.
For each $k\in K$ we have for all $\iota$ of sufficiently high order 
\begin{equation}
\sup \abs{1-h_\iota(z)}/k(z)<1
\end{equation}
For $I(z)\in \mathbb{E}'(\R)$ we set $\tilde{I}(z_1,z_2):=I(\sqrt{z_1^2 +z_2^2})\in \mathbb{E}'(\R^2).$
This yields
\begin{equation}
\abs{1-\tilde{h}_\iota(z_1,z_2)}/k_1\left(\abs{\re \sqrt{z_1^2+z_2^2}}\right)k_2\left(\abs{\re \sqrt{z_1^2+z_2^2}}\right)<1
\end{equation}
Now for $(z_1,z_2)\in C^2$ we have 
\begin{equation}\label{dhow}
\abs{\re \sqrt{z_1^2+z_2^2}}\leq \norm{\re (z_1,z_2)}, \abs{\im \sqrt{z_1^2+z_2^2}} \leq \norm{\re (z_1,z_2)}
\end{equation}
where $\norm{\re (z_1,z_2)}=\abs{\re z_1}+\abs{\re z_2}.$
To see this note that if $(a +ib)^2=z_1^2+z_2^2,$ for $a,b\in\R,$ an $(z_1,z_2)=(x_1,x_2)+i(y_1,y_2),$ then
$ab=(x_1,x_2)\cdot (y_1,y_2),$
and
$b^2=\abs{(y_1,y_2)}^2-\abs{(x_1,x_2)}^2 +a^2.$ If 
$b^2>\abs{(y_1,y_2)}^2$ then by the Cauchy-Schwarz inequality $\abs{a}\leq \abs{(x_1,x_2)}\abs{(y_1,y_2)}/\abs{b} 
< \abs{(x_1,x_2)},$
which inserted into $b^2=\abs{(y_1,y_2)}^2-\abs{(x_1,x_2)}^2 +a^2,$ yields $b^2\leq \abs{(y_1,y_2)}^2,$ which is a contradiction.
Hence $\abs{b}^2\leq \abs{(y_1,y_2)}^2$ and hence $\abs{a}^2\leq \abs{(x_1,x_2)}^2.$
Since $k_1,k_2$ are increasing, Eqn.(\ref{dhow}) implies that for all $\iota$ of sufficiently high order
\begin{equation}
\abs{1-\tilde{h}_\iota (z_1,z_2)}/k_1(\norm{\re (z_1,z_2)})k_2(\norm{\im (z_1,z_2)})<1
\end{equation}
This implies that
\begin{equation}
\tilde{h}_\iota (z_1,z_2)=f_{1,\iota}(\sqrt{z_1^2+z_2^2}) F_1(\sqrt{z_1^2+z_2^2})+f_{2,\iota}(\sqrt{z_1^2+z_2^2}) F_2(\sqrt{z_1^2+z_2^2})\to 1
\end{equation}
in $\mathbb{E}'(\R^2)$ so the closure of the ideal generated by $F_1$ and $F_2$ is all of $\mathbb{E}'(\R^2).$
Let $T_j$ be the element in $\mathcal{E}'$ whose Fourier transform is $(-1)^s 2^n F_j.$ Then the system
\begin{equation}
T_1 *g=0,\quad T_2 * g=0,\quad g\in C^\infty(\R^2)
\end{equation}
has only the trivial solution $g\equiv 0.$
Since $\hat{\sigma}_j(z_1,z_2)=(z_2-iz_1)^n(z_1^2+z_2^2)^s F_j(z_1,z_2),$
we must have $\partial_{\bar{z}}^n (\Delta^s T_j)=\sigma_j.$ Since $f*\sigma_j=0,$ $j=1,2,$
we have $T_j*\partial_{\bar{z}}^n (\Delta^s f)=0,$ $j=1,2,$ hence $\partial_{\bar{z}}^n (\Delta^s f)=0.$
This proves the result under the assumption that $f$ is $C^\infty$-smooth.
If $f\in L^1_{\mbox{loc}}(\R^2)$ use a smooth approximation of identity, i.e.\ $\phi_\epsilon\in C^\infty_c(\R^2),$
such that $\phi_\epsilon\to 1$ and set $f_\epsilon:=f*\phi_\epsilon \in C^\infty(\R^2)$. Then 
$f_\epsilon \to f$ locally in $L^1(\R^2)$ thus $\partial_{\bar{z}}^n (\Delta^s f_\epsilon) \to 
\partial_{\bar{z}}^n \Delta^s f$ in distribution sense. Hence $f$ is a distribution solution to 
$\partial_{\bar{z}}^n \Delta^s f=0.$ This completes the proof.
\end{proof}

We immediately obtain the following consequence, originally proved by Reade \cite{reade}.
\begin{theorem}[Two-radius theorem for $q$-analytic functions]
Let $f\in L^1_{\mbox{loc}}(\R^2)$ and suppose that there exists distinct positive $r_1,r_2$ such that
for almost all $z\in \C,$
\begin{equation}
\int_0^{2\pi} f(z+r_j\exp(i\theta))\exp(iq\theta) d\theta =0, \quad j=1,2
\end{equation}
Then, if $r_1/r_2\notin Z_q$ and $q\geq 0,$ $f$ 
agrees almost everywhere with a $q$-analytic function.
\end{theorem}
\begin{proof}
If in Theorem \ref{zalcmanthm} one chooses $\mu=\delta_1,$ by which we mean the point mass at $1$,
then $\mu_{r_j}=\delta_{r_j},$ so
the integral conditions reduce to
\begin{equation}
\int_0^{2\pi} f(z+r_1\exp(i\theta))\exp(iq\theta) d\theta =0, \quad j=1,2
\end{equation}
Furthermore, $A_{k,q}$ reduces (for point mass measure $\delta_1$) for $q\in \Z_+,$ to $A_{k,q}=z^{2k+q}.$
But then
be the smallest value of $k$ for which $A_{k,q}\neq 0$ is $s=0$ and
\begin{equation}
F(z)=J_q(-z)
\end{equation}
So $Z(q,\delta_1)=Z_q$ (i.e.\ the set of ratios of the nonzero zeros of the Bessel function).
Now the differential operator $\partial_{\bar{z}}^q$ is elliptic for $q\in \Z_+,$ thus by the Elliptic regularity theorem
the homogeneous solutions agree almost everywhere with smooth functions.
This completes the proof.
\end{proof}
Let $z=x+iy=\rho\exp(i\varphi)$ $\rho:=\abs{z},$ $\varphi\in [0,2\pi).$ Denote by 
\begin{equation}
W_j(z):=J_j(z)(\rho)\exp(ij\varphi)
\end{equation}
The following is known (see Volchkov \cite{volchkov}, Lemma 3) 
\begin{lemma}\label{korenovlemma}
For $j\in \Z,$ $k\in \N,$ $\lambda >0$ and any $z\in \{\abs{z} <R-r\},$ $0<r<R,$ we have
\begin{equation}
\int_{\abs{\zeta}\leq r} \zeta^k W_m(\lambda(\zeta +z))d\mu(\zeta)=(-1)^k \frac{2\pi r^{k+1}}{\lambda} J_{k+1}(\lambda r)W_{m+k}(\lambda z)
\end{equation}
and
\begin{equation}
\int_{\abs{\zeta}= r} \zeta^{k-1} W_m(\lambda(\zeta +z))\zeta=(-1)^k 2\pi r^{k} J_{k}(\lambda r)W_{m+k}(\lambda z)
\end{equation}
\end{lemma}
\begin{proof}
In the case $m=0$, we have $W_0(z)=J_0(\rho)$ 
and the result follows from the so-called Neumann addition theorem (see Korenov \cite{korenov}, p.39). 
For the case $m>0$ the result follows from the case $m=0$ by applying $\partial_{\bar{z}}^m$ to both sides of 
the equation. 
\end{proof}
Volchkov \cite{volchkov} proved the following.
\begin{theorem}
For $z\in \C$ and $r>0,$ denote $K(z,r):=\{\abs{z-p}\leq r\}$ and let $q\in \Z_+.$ Let $r_1,r_2\in \R_+.$
Let $f\in C^0(K(0,R))$, for a number $R>0,$ such that for all $z\in \C,$
\begin{equation}
\int_{K(z,r_j)} (\zeta-z)^q f(\zeta)d\mu(\zeta) =0,\quad j=1,2
\end{equation}
If $r_1+r_2 <R$ and $\frac{r_1}{r_2}\notin Z_{q+1},$ then $f$ is $q$-analytic.
\end{theorem}
\begin{proof}
Extend $f$ to be zero on $\{\abs{z}\geq R$ and let $\phi\in C^\infty_c(\C)$ such that $\phi=1$ on $\{\abs{z}\leq R\}.$ 
It is then sufficient to prove the result for the $C^\infty$-smooth regularization by $\phi,$ i.e.\ the convolution $f*\phi.$ 
We shall thus keep the notation $f$ but assume $f\in C^\infty(K(0,R))$. 
First we show that if for all $\abs{z}<R-r$ 
\begin{equation}\label{ekvationm1}
\int_{K(z,r)} (\zeta-z)^q f(\zeta)d\mu(\zeta) =0
\end{equation}
then each term of the Fourier series
\begin{equation}\label{fffour}
f(\rho \exp(i\varphi))=\sum_{k=-\infty}^\infty f_k(\rho)\exp(ik\varphi),\quad \rho>0
\end{equation}
defined at $z=0$ by continuity, is also $C^\infty$-smooth
and satisfies for all $\abs{z}<R-r$ 
\begin{equation}
\int_{K(z,r)} (\zeta-z)^q f_k(\rho)\exp(ik\varphi)d\mu(\zeta) =0
\end{equation}
The $C^\infty$-smoothness follows from that of $f$ together with the formula
\begin{equation}
f_k(\rho)\exp(ik\varphi)=\frac{1}{2\pi}\int_{-\pi}^\pi f(x\cos \alpha -y\sin\alpha,x\sin\alpha +y\cos\alpha) \exp(-ik\alpha)d\alpha
\end{equation}
Furthermore, writing Eqn.(\ref{ekvationm1}) for $f(z\exp(i\alpha)),$ multiplying it by $\exp(-ik\alpha)$,
and integrating over $[-\pi,\pi],$ we obtain the wanted result.
Let $\{\nu_k\}_{k\in \Z_+}$ be an increasing sequence of all positive zeros of $J_{n+1}.$
\begin{proposition}[Volchkov \cite{volchkov1991}, see also Lemma 5.6, Volchkov \cite{volchkovbok}]\label{captturre}
Let $R>0$ and $f\in C^\infty (K(z,R)),$ where $K(z,R):=\{\abs{z}<R\}.$ 
Then $f$ satisfies for a fixed $r<R,$ and all
$K(z,r)\subset K(z,R)$
\begin{equation}
\int_{K(z,r)} (\zeta-z)^n f(\zeta)d\mu(\zeta) =0
\end{equation}
if and only if each coefficient $f_k(\rho)$ of the Fourier series
in Eqn.(\ref{fffour}) has the form 
\begin{equation}
f_k(\rho)=\sum_{m=1}^\infty a_{mk} J_k\left(\frac{\nu_m\rho}{r}\right) +\sum_{\substack{0\leq p\leq n-1}{p+k \geq 0}} c_{pk} 
\rho^{2p+k}
\end{equation}
where $a_{mk}=O(m^{-\alpha})$ as $m\to \infty$ for any $\alpha >0.$
\end{proposition}
By Proposition \ref{captturre} (with $r=r_1$) together with Lemma \ref{korenovlemma} (with $r=r_2$) we have
for all $\rho\leq r_1$
\begin{equation}\label{shutyomo}
\sum_{m=1}^\infty \frac{a_{mk}}{\nu_m} J_{n+1}\left(\frac{\nu_m r_2}{r_1}\right) J_{k+n}\left(\frac{\nu_m \rho}{r_1}\right)=0 
\end{equation}
Now if we show that each term in the Fourier series given by Eqn.(\ref{fffour}) is $n$-analytic then 
$f$ is $n$-analytic. From the definition of $n$-analyticity the coefficients of such terms can then be written 
$f_k(\rho)=\sum_{\substack{0\leq p\leq n-1}{p+k \geq 0}} c_{pk} 
\rho^{2p+k}.$
For $k=1$ the formulas for the coefficients of the Fourier-Bessel expansion (see e.g.\ \cite{tolstov}, p.270) 
imply that all $a_{m1}=0,$ for all $m$,
hence each term in the Fourier series given by Eqn.(\ref{fffour}) is $n$-analytic.
If $k>1$, then differentiating both sides of Eqn.(\ref{shutyomo}) decreases the index of the Bessel function 
(see \cite{vladimirov}, p.350) 
and this can be done in a way to transform into a Fourier-Bessel series in order to obtain $a_{mn}=0$ for all $m.$
This shows that each term in the Fourier series given by Eqn.(\ref{fffour}) is $n$-analytic. Hence $f$ is $n$-analytic.
This completes the proof.
\end{proof}

\section{Cauchy formulas and the annular maximum principle}
The results of this subsection are mainly due to 
Dolzhenko \& Danchenko \cite{dolchenkodanchenko2002}, \cite{dolchenkodanchenko2005},
\cite{dolchenkodanchenko2007},\cite{dolchenkodanchenko1998} and Dolzhenko \cite{dolzhenkoensam}.
\begin{definition}
The {\em modulus of continuity}\index{Modulus of conitnuity}, $\omega(f,E,\delta)$, 
of a function $f$ with respect
to a subset $E\subset\C$ is defined as
$\omega(f,E,\delta):=\sup\{\abs{f(z)-f(t)}:z,t\in E,\abs{z-t}\leq \delta\},$ $\delta\geq 0.$
\end{definition}
When it is well-defined (for example when $f$ is bounded on a convex set $E$) the function $\omega(r):=\omega(f,E,r)$ is nonnegative, nondecreasing,
$\omega(0)=0,$ $\omega(a+b)\leq \omega(a)+\omega(b)$ for all nonnegative $a,b.$
\begin{definition}
Let $q\in \Z_+$ and let $\Omega\subset\C$ be a domain.
For $R>0$ and a given function $f$, defined at $z\in \C$ and Lebesgue integrable on the circle
$\{\abs{t-z}=R\}$ define (using polar parametrization of the circle, $t=z+R\exp(i\theta)$)\index{$\Delta_R^q$ ($q$:th Halo difference of radius)}
\begin{equation}\label{dolchsex}
\Delta_R f(z):=\frac{1}{2\pi R} \int_{\{\abs{t-z}=R\}} (f(t)-f(z))\abs{dt}=\frac{1}{2\pi}\int_0^{2\pi}f(z+R\exp(i\theta))d\theta
\end{equation}
(this is called the {\em $q$:th Halo difference of radius})
and set
\begin{equation}
\omega_{\mbox{or}}^{q}(f,\Omega,\delta):=\sup\{\abs{\Delta_R^q f(z)}\},\quad 0<\delta<\infty,\quad \omega_{\mbox{or}}^{q}(f,\Omega,0)=0
\end{equation}
where the supremum is taken over all pairs $(z,R)$ such that $0<R\leq \delta$ and the closed disc of radius $qR$
with center $z$, belongs to $\Omega.$
\end{definition}
Note that if $f(x,y)\in C^2(U)$ for an open neighborhood $U$ of a point $z_0=(x_0,y_0)$ then
$\lim_{R\to 0} R^{-2}\Delta_R f(z_0)=\frac{1}{2}\left(\partial_x^2f(x_0,y_0)+\partial_y^2 f(x_0,y_0)\right)=\Delta f(z_0).$
\begin{proposition}\label{halodiffprop}
Let $q\in \Z_+,$ let $\Omega\subset\C$ be a domain and let $f$ be a $q$-analytic function on $\Omega,$ i.e.\
there exists holomorphic functions $a_k(z)$, $k=0,\ldots,q-1$, on $\Omega$ such that
$f(z)=\sum_{k=0}^{q-1} a_k(z)\bar{z}^k$ on $\Omega.$
Then for $z\in \Omega$ and $0<R<\rho:=\mbox{dist}(z,\partial\Omega)$
\begin{equation}
\Delta_R f(z)= \sum_{p=0}^{q-1} \bar{z}^p \sum_{k=p+1}^{q-1} \binom{k}{p}\frac{R^{2k-2p}}{(k-p)!} \partial_z^{k-p}a_k(z)
\end{equation}
and for $q\geq 2,$ $z\in \Omega,$ $0<R<\frac{\mbox{dist}(z,\partial\Omega)}{q-1}$
\begin{equation}\label{dolchvisatva}
\Delta_R^{q-1} f(z)= (q-1)!R^{2q-2} \partial_z^{q-1}a_{q-1}(z)
\end{equation}
\end{proposition}
\begin{proof}
For $\{t:\abs{t-z}\leq R\}\subset\Omega$ we have by Eqn.(\ref{dolchsex}) together with the representation $f(z)=\sum_{k=0}^{q-1} a_k(z)\bar{z}^k$
\begin{equation}\label{dolchnine}
\Delta_R f(z)= \sum_{k=0}^{q-1} 
\frac{1}{2\pi i} \int_{\{\abs{t-z}=R\}} \frac{\bar{t}^k-\bar{z}^k}{t-z}a_k(t)dt
\end{equation}
Using the relations
\begin{equation}
\bar{t}^k-\bar{z}^k= \sum_{p=0}^{k-1} \binom{k}{p}(\bar{t}-\bar{z})^{k-p}\bar{z}^p
\end{equation}
\begin{multline}
\int_{\{\abs{t-z}=R\}}(\bar{t}-\bar{z})^{q} a_k(t)(t-z)^{-1}dt =R^{2q} 
\int_{\{\abs{t-z}=R\}}a_k(t)(t-z)^{-q-1}dt \\=R^{2q}\frac{2\pi i}{q!} \partial_z^q a_k(z)
\end{multline}
\begin{equation}
\Delta_R f(z)  =\sum_{k=0}^{q-1}\sum_{p=0}^{k-1} \binom{k}{p}\frac{1}{(k-p)!}R^{2(k-p)}\bar{z}^p \partial_z^{k-p} a_k(z)
\end{equation}
we obtain by Eqn.(\ref{dolchnine})
\begin{equation}
\Delta_R f(z)= \sum_{k=0}^{q-1}\sum_{p=0}^{k-1}\binom{k}{p} \frac{1}{(k-p)!} R^{2(k-p)}\bar{z}^p \partial_z^{k-p} a_k(z)
\end{equation}
For $p=q-2,k=q-1$ the sum in the right hand side takes the form
\begin{equation}
\bar{z}^{q-2} \binom{q-1}{q-2} \frac{R^{2(q-1)-2(q-2)}}{((q-1)-(q-2))!} \partial_z^{(q-1)-(q-2)} a_{q-1}(z)=\bar{z}^{q-2}(q-1)R^2\partial_z a_{q-1}(z)
\end{equation}
thus taking the $(q-1)$:th power of $\Delta_R$ yields Eqn.(\ref{dolchvisatva}).
This completes the proof.
\end{proof}
We immediately obtain the following.
\begin{corollary}
For a domain $\Omega\subset\C$, $q \in \Z,$ $q\geq 2$ and $f\in \mbox{PA}_q(\Omega)$ we have
\begin{equation}
\abs{\partial_z^{q-1} a_{q-1}(z)} \leq \frac{q^{2q-2}}{(q-1)!} \rho^{2(1-q)} \omega_{\mbox{or}}^{q-1}(f,\Omega,\rho/q)
\end{equation}
for $z\in \Omega,$ $\rho=\mbox{dist}(z,\partial\Omega).$
\end{corollary}
Let $\Omega\subset\C$ be a domain. For $z\in \Omega,$ $v\in \C,$ and $f\in \mbox{PA}_q(\Omega)$ we
know that $f(z)=\sum_{m=0}^{q-1}\bar{z}^m \phi_m(z)$ for holomorphic $\phi_m,$ and 
we associate to $f$ the function
\begin{equation}
F(v,z):=\sum_{m=0}^{q-1} v^m \phi_m(z)
\end{equation}
For fixed $\zeta\in \Omega,$ $\rho>0$ assume
$f$ is also continuous on $\{ t:\abs{t-\zeta}\leq\rho\},$ i.e.\ $f\in A_q(\{ t:\abs{t-\zeta}\leq\rho\}).$
Define
\begin{equation}\label{danchdolz98a0}
J(k,z):=J(f,\{ t:\abs{t-\zeta}\leq\rho\},z)=\frac{1}{2\pi i} \int_{\{ t:\abs{t-\zeta} =\rho\}}\frac{f(t)}{t-z}
\left(\frac{t-\zeta}{z-\zeta}\right)^{k-1}dt
\end{equation}
Then by Proposition \ref{halodiffprop} (see also Dolzhenko \& Danchenko \cite{dolchenkodanchenko1998}), for $z=\zeta -\rho^2/(\bar{\zeta}-\bar{v})\in \{ t:\abs{t-\zeta}<\rho\}$ and $v$ symmetric with respect to
$\{ t:\abs{t-\zeta} =\rho\}$, we have
\begin{equation}\label{danchdolz98a}
F(\bar{v},z)=J(n,z),\quad F(\bar{v},z)-(\bar{v}-\bar{\zeta})^{q-1}\phi_{q-1}(\zeta)=J(q-1,z)
\end{equation}
By the Sokhotsky formula we have for the nontangential limit as $z\to z_0\in \{ t:\abs{t-\zeta} =\rho\}$
that the singular integral $J(f(\cdot)-f(z_0),\{ t:\abs{t-\zeta} <\rho\},k,z_0),$ is zero for $k\geq q$
and is equal to $-(\bar{z}_0-\bar{\zeta})^{q-1}\phi_{q-1}(\zeta)$ for $k=q-1.$
From the case $k=q-1$ we have that if for some $\alpha\in (0,1)$ and $M>0,$ we
\begin{equation}
\abs{f(z)-f(z_0)}\leq M\abs{z-z_0}^\alpha, \quad z\in \{ t:\abs{t-\zeta} =\rho\}
\end{equation}
then $\abs{\phi_{q-1}(\zeta)}\leq\mbox{const}\cdot M\rho^{\alpha-1}.$ Now let $\epsilon\in (0,1),$
$\rho>0,$ $z\in \C$ and $f\in A_q( \{ t:\abs{t-z}\leq (1+\epsilon)\rho\}).$ Then by
Eqn.(\ref{danchdolz98a0}) and Eqn.(\ref{danchdolz98a}) together with the Cauchy formula for $s=0,\ldots,q-1$
we have
\begin{multline}\label{doldalmult1}
\partial_{\bar{z}}^s f=\partial_v^s F(v,z)|_{v=\bar{z}} =\frac{s!}{2\pi i} \int_{\abs{\bar{v}-\bar{z}}=\frac{\rho}{\epsilon}(1-\epsilon^2)}
F(\bar{z},z)(\bar{v}-\bar{z})^{-s-1}d\bar{v}=\\
\frac{s!}{4\pi^2\epsilon^2(1-\epsilon^2)^s\rho^{2(s+1)}}\int_{\abs{\zeta -z}=\rho\epsilon}(z-\zeta)^{s+2-q}
\int_{\abs{t-\zeta}=\rho} f(t)(t-\zeta)^{q-1}\frac{dt}{t-z}d\bar{\zeta}
\end{multline}
Setting
\begin{equation}
K:=K(z,\rho,\epsilon)=\{t: (1-\epsilon)\rho <\abs{t-z}< (1+\epsilon)\rho  \},\quad r=\abs{t-z}
\end{equation}
we can transform the integral in right hand side of Eqn.(\ref{doldalmult1}) as an area integral over $K$ with respect to the area measure $d\mu(t)$
with respect to $t.$ Using the relations
\begin{equation} 
\abs{dtd\bar{\zeta}}=r\rho^{-1}\abs{\sin^{-1}\alpha}d\mu(t),\quad \alpha=\mbox{arg}((t-\zeta)/(t-z)) 
\end{equation}
\begin{equation}
2r\rho\abs{\sin \alpha}=\sqrt{(\rho^2 (1-\epsilon)^2-r^2)(r^2-\rho^2 (1-\epsilon)^2)}
\end{equation}
it is possible to rewrite (see Dolzhenko \& Danchenko \cite{dolchenkodanchenko1998} for the details) Eqn.(\ref{doldalmult1}) according to
\begin{equation}
\partial_{\bar{z}}^s f=\frac{s!}{2\pi^2}\frac{\epsilon^{s-q+1}}{(1-\epsilon^2)^s}\rho^{1-s}\int_K \frac{f(t)}{r} E(t)\frac{d\mu(t)}{Q(\rho,\epsilon,r)}
\end{equation}
for $s=0,\ldots,q-1,$ where $t=z+r\exp(i\theta),$ 
\begin{equation}
E(t):=\exp(i\alpha_1(t))+\exp(i\alpha_1(t)),\quad Q(\rho,\epsilon,r):=\sqrt{(\rho^2 (1-\epsilon)^2-r^2)(r^2-\rho^2 (1-\epsilon)^2)}
\end{equation}
where $\zeta_j$, $j=1,2$, are points of the points of intersection of the circle $\abs{\zeta-z}=\rho\epsilon$
and $\abs{\zeta-t}=\rho$ respectively,
and $\alpha_j$ is the argument of $(z-\zeta_j)^{s+2-q}(t-\zeta_j)^{q-1}(t-z)^{-1}d\bar{\zeta}_jdt_j.$
Define
\begin{equation}\label{dolcheqargs}
\delta(\rho,\epsilon,r):=\arccos \frac{\rho^2(1+\epsilon^2)-r^2}{2\rho^2\epsilon},\quad \gamma=\gamma(\rho,\epsilon,r):=
\arccos \frac{r^2-\rho^2(1-\epsilon^2)}{2\rho^2\epsilon}
\end{equation}
\begin{equation}
\alpha(s,q,r)=\alpha(\rho,\epsilon,s,n,r)=(s+1)\gamma +q\delta,\quad \delta=\delta(\rho,\epsilon,r)
\end{equation}
This renders using the polar representation $t=z+r\exp(i\theta)$
\begin{equation}
E(t)=(-1)^{s+1}2\exp(is\phi) \cos (\alpha(s,q,r))
\end{equation}
Define for fixed $\rho>0,$ $\epsilon\in (0,1),$ $z\in \C$ 
\begin{equation}\label{dolchj1eq}
J_1(f,q,s,K):=\frac{1}{\pi^2}\int_K \frac{f(t)\exp(is\theta)}{\abs{t-z}}\frac{\cos(\alpha(s,q,\abs{t-z}))}{Q(\rho,\epsilon,\abs{t-z})}d\mu(t)
\end{equation}
Then for $f\in \mbox{PA}_q(\{t:\abs{t-z}<\rho(1+\epsilon)\})\cap C^0(\{t : \abs{t-z} \leq\rho(1+\epsilon)\})$
ans $s=0,\ldots,q-1$
\begin{equation}\label{dolchseven}
\partial_{\bar{z}}^s f=(-1)^{s+1}s! \frac{\epsilon^{s-q+1}}{(1-\epsilon^2)^s}\rho^{1-s} J_1(f,q,s,K)
\end{equation}
Similarly, calculating the convolution in the second part of Eqn.(\ref{danchdolz98a}) and using the relations
\begin{equation}
\bar{v}-\bar{z}=(1-\epsilon^2)(\bar{v}-\bar{\zeta}),\quad \bar{v}-\bar{\zeta}=\frac{\rho^2}{z-\zeta}
\end{equation}
which imply
\begin{equation}
d\bar{v}=(1-\epsilon^2)\rho^2\frac{d\zeta}{z-\zeta}
\end{equation}
yields for $s=0,\ldots,q-1,$ $\rho>0,$ $\epsilon\in (0,1),$ $z\in \C,$ and $f$ as before
\begin{equation}\label{dolchsevenprim}
\partial_{\bar{z}}^s f+(-1)^{q-s} \frac{s!}{(q-s-1)!} \frac{\rho^{2(q-1-s)}}{(1-\epsilon^2)^s} a_{q-1}^{q-1-s} J_1(f,q-1,s,K)
\end{equation}
Eqn.(\ref{dolchseven}) together with Eqn.(\ref{dolchsevenprim}) give
\begin{equation}
\frac{1}{s!}\partial_{\bar{z}}^s a_{q-1}(z)=(-1)^q \rho^{2(q-1-s)}\frac{1}{\epsilon^s} 
\left(J_1(f,q,q-1-s,K)-\epsilon J_1(f,q-1,q-1-s,K)\right)
\end{equation}
Next by Eqn.(\ref{dolcheqargs}) we can (after writing out each side of the equation) verify
\begin{equation}\label{hohoty}
\cos(\gamma +q\delta)=\frac{\epsilon\cos(q\delta)-\cos((q-1)\delta)}{\sqrt{1-\epsilon^2-2\epsilon\cos\delta}}
\end{equation}
Setting $s=0$ in Eqn(\ref{dolchj1eq}) and making the change of variable, for fixed $z$, $r:=\abs{t-z},$
we note that 
\begin{equation}\label{hohoty1}
2rQ^{-1}(\rho,\epsilon,r)=\frac{d\delta(\rho,\epsilon,r)}{dr}\end{equation} 
thus in the polar
representation $t=z+r\exp(i\theta),$ where $r=\rho\sqrt{1-\epsilon^2-2\epsilon\cos\delta},$ we obtain the following.
\begin{theorem}
Let $q\in \Z_+$ and $f\in C^0(\{t:\abs{t-z}<(1+\epsilon)\rho\})\cap\mbox{PA}_q(\{t:\abs{t-z}<(1+\epsilon)\rho\}).$
Then
\begin{multline}
f(z)=\frac{\epsilon^{-q+1}}{2\pi^2} \times\\
\int_0^\pi \frac{\cos((q-1)\delta)-\epsilon \cos(q\delta)}{1-\epsilon^2-2\epsilon\cos\delta}
d\delta \int_0^{2\pi}f(z+\rho\sqrt{1-\epsilon^2-2\epsilon\cos\delta}\exp(i\theta)d\theta
\end{multline}
\end{theorem}
\begin{corollary}\label{maxprincedolch}
For $K:=K(z,\rho,\epsilon)=\{t: (1-\epsilon)\rho <\abs{t-z}< (1+\epsilon)\rho  \}$
we have
\begin{equation}
\abs{f(z)}\epsilon^{q-1}(1-\epsilon)\leq \norm{f}_K,\quad \rho>0,\epsilon \in (0,1)
\end{equation}
\end{corollary}
Analogous to the derivation of Eqn.(\ref{dolchseven}) we have for $v\in \C,\rho>0,\epsilon\in (0,1),$
$s\geq 0$ and $f\in A_q(\{\abs{z-v}<(1+\epsilon)\rho\}),$ taking into account the first equality in 
Eqn.(\ref{danchdolz98a}) we write out the convolution integral
\begin{equation}\label{jjkk11}
\partial_v^s f(v)=\partial_z^s F(\bar{v},z)|_{z=v}=\frac{s!}{2\pi i} 
\int_{\{\abs{z-v}=\rho(1-\epsilon^2)\}} \frac{F(\bar{v},z)}{(z-v)^{s+1}} dz
\end{equation}
Defining, as before with the polar representation $t=z+r\exp(i\theta)$,
\begin{equation}\label{dolchten}
J_2(f,q,s,K):=\frac{1}{\pi^2} \int_K \frac{f(t) \exp(-is\phi)\cos\beta(s,q,\abs{t-z})}{Q(\rho,\epsilon,\abs{t-z})}d\mu(t) 
\end{equation}
with $\beta(s,q,r)=(s+1)\gamma +(s+q)\delta$.
Eqn.(\ref{jjkk11}) can be transformed with the variable $v$ replaced by $z$, to
\begin{equation}\label{dolchnine}
\partial_v^s f(v)=(-1)^{s+1} s!\frac{\epsilon^{-q+1}}{(1-\epsilon^2)^s}\rho^{1-s}J_2(f,q,s,K)
\end{equation}
We fix $z_0\in \C$ and let $q\in \Z_+,$ $\epsilon\in (0,1),\rho>0, 0<d<\rho(1-\epsilon),$
$\lambda:=(\rho(1-\epsilon)-d)/(\rho(1+\epsilon)),$ $\rho_1:=\lambda\rho.$
Let $f\in \mbox{PA}_q(\{t:\abs{t-z_0}<2\rho-d  \})\cap C^0(\{t:\abs{t-z_0}\leq 2\rho-d  \}),$
have the representation $f(z)=\sum_{m=0}^{q-1} \bar{z}^m a_m(z).$ 
Consider the associated function
$F(\bar{v},z)=\sum_{m=0}^{q-1} \bar{v}^m a_m(z).$
Now we have
\begin{equation}
\int_K \frac{d\mu(t)}{Q(\rho,\epsilon,r)}=\pi^2
\end{equation}
In Eqn.(\ref{dolchnine}) we carry out the integration with respect to
$t_1:=z_0+r_1\exp(i\phi_1),$ $(1-\epsilon)\rho<r_1<(1+\epsilon)\rho$, $\phi_1\in [0,2\pi),$
and in Eqn.(\ref{dolchseven}) with respect to $t_2=t_1+r_2\exp(i\phi_2),$ $(1-\epsilon)\rho_1 <r_2<(1+\epsilon)\rho_1,$
$\phi_2\in [0,2\pi)$ for fixed $t_1.$ Noting that the set of values of the argument $t_2$ of the integrable function $f$ belongs to
$K(z_0,q,2\rho-d),$ we obtain for $f\in A_q(\{\abs{t-z_0}<2\rho-d\})$ the inequality 
\begin{equation}\label{rearea}
\abs{\partial_{\bar{v}}^l\partial_z^k F}\leq k!l!\lambda^{-l}\frac{1}{(1-\epsilon)^2}\frac{\epsilon^{l+2(1-q)}}{(\rho(1-\epsilon^2))^{k+l}}\norm{f}_{K(z_0,d,2p-d)}
\end{equation}
Now expand $F(\bar{v},z)$ in the variables $v,z$ in the set
$\{ \abs{\bar{v}-\bar{z}_0}\leq d,\abs{z-z_0}\leq d\}$ according to
\begin{equation}
F(\bar{v},z)=\sum_{m=0}^\infty \sum_{l=0}^{q-1} \frac{1}{l!(m-l)!} (\partial^{(l,m-l)}F)(\bar{v}-\bar{z}_0)^l(z-z_0)^{m-l}
\end{equation}
By Eqn.(\ref{rearea}) the inner sum can be estimated by
\begin{equation}
\lambda^{1-q}\frac{\epsilon^{2(1-q)}}{(1-\epsilon)^2}\frac{\epsilon^q-\lambda^q}{\epsilon -\lambda}\frac{d^m}{\rho^m(1-\epsilon^2)^m}\norm{f}_{K(z_0,d,2\rho -d)}
\end{equation}
which implies, since $d<\rho(1-\epsilon)$, that the series converges uniformly in the closure of $\{ \abs{\bar{v}-\bar{z}_0}\leq d,\abs{z-z_0}\leq d\}$.
This proves the following.
\begin{theorem}
Let $z_0\in \C$, $q\in \Z_+,$ $\epsilon\in (0,1),\rho>0, 0<d<\rho(1-\epsilon),$
$\lambda=\frac{\rho(1-\epsilon)-d}{\rho(1+\epsilon)}$ and $f(z)=F(\bar{z},z)\in A_q(\{\abs{z-z_0}<2\rho-d\}).$
Then for any point $(v,z)\in\{\abs{t-z_0}\leq d\}^2$ we have
\begin{equation}
\abs{F(\bar{v},z}\leq \lambda^{1-q} \frac{\epsilon^{1-2q}}{(1-\epsilon)^2}\frac{\epsilon^q-\lambda^q}{\epsilon -\lambda}
(1+\epsilon)\norm{f}_{K(z_0,d,2p-d)}
\end{equation}
\end{theorem}
Now let $z_0\in \C$, $q\in \Z_+,$ $\epsilon\in (0,1),\rho>0,$ $K:=K(z_0,\rho,\epsilon)$
$=\{t:(1-\epsilon)\rho <\abs{t-z}<(1+\epsilon)\rho\},$ $z=z_0+\rho_0\exp(i\theta_0),$ $t=z_0+r\exp(i\theta),$
$a:=\rho_0(1-\epsilon^2)^{-1}\rho^{-1}$ and 
$f(z)=F(\bar{z},z)\in A_q(\{\abs{t-z_0}<\rho(1+\epsilon)\}).$ 
Define for $\delta=\delta(\rho,\epsilon,r),$ $\delta=\delta(\rho,\epsilon,r)$ given by Eqn.(\ref{dolcheqargs})
\begin{equation}
W(\phi,\phi_0,\delta):=\sum_{s=0}^\infty (-1)^s a^s \exp(-is(\phi-\phi_0)\cos((s+1)\gamma+(s+q)\delta)
\end{equation}
Taylor expansion of $F(\bar{z}_0,z)$ gives for 
$z\in \{\abs{z-z_0}<(1-\epsilon)\rho\}$ and $t=z_0+r\exp(i\phi),$ using Eqn.(\ref{dolchnine}) and Eqn.(\ref{dolchten})
\begin{multline}
F(\bar{z}_0,z)=\sum_{s=0}^\infty \frac{1}{s!}\partial_z^s F|_{z=z_0}(z-z_0)^s=-\epsilon^{1-q}\rho\sum_{s=0}^\infty(-1)^s a^s
\exp(i\theta_0 s)J_2(f,q,s,K)\\
=-\frac{\rho\epsilon^{1-q}}{\pi^2}\int_K\frac{f(t)}{rQ(\rho,\epsilon,r)}W(\phi,\phi_0,\delta)d\mu(t)
\end{multline}
Setting
\begin{equation}
\sigma(\epsilon,\delta)=\sqrt{1+\epsilon^2-2\epsilon\cos\delta}
\end{equation}
and using Eqn.(\ref{hohoty}) and Eqn.(\ref{hohoty1}) with $t=z_0+r\exp(i\phi=z_0+\rho\sigma(\epsilon,\delta)\exp(i\phi)$
we obtain
\begin{equation}\label{dolchfiftheen}
F(\bar{z}_0,z)=-\frac{\epsilon^{1-q}}{2\pi^2}\int_0^{2\pi}d\phi\int_0^\pi \frac{f(z_0+\rho\sigma(\epsilon,\delta)\exp(i\phi)}{\sigma(\epsilon,\delta)}W(\phi,\phi_0,\delta)d\delta
\end{equation}
By the relations
\begin{equation}
\exp(i\alpha= \frac{1-\epsilon\cos\delta+i\epsilon\sin\delta}{\sigma(\epsilon,\delta)},\quad \alpha\in (0,\pi/2)
\end{equation}
\begin{equation}
W(\phi,\phi_0,\delta)=(1-a\exp(-i(\phi-\phi_0)\exp(i\alpha))(1-a\exp(-i(\phi-\phi_0)\exp(-i\alpha)
\end{equation}
we have the following integral formula for $z\in \{z:\abs{z-v}<\rho(1-\epsilon)\},$ $a=\rho(1-\epsilon^2),$ 
$\sigma=\sigma(\epsilon,\delta)=\sqrt{1+\epsilon^2-2\epsilon\cos\delta}$
\begin{multline}\label{dolchfiftheenprim}
F(\bar{v},z)=-\rho(1-\epsilon^2)\frac{\epsilon^{1-q}}{4\pi^2}\int_0^{2\pi} \exp(i\phi)d\phi
\int_0^\pi \frac{f(v+\rho\sigma\exp(i\phi)}{\sigma}\times\\
\left(
\frac{\exp(-i(q-1)\delta}{(z-v)+a\exp(i\phi)(\epsilon\exp(i\delta)-1)/\sigma} +\frac{\exp(i(q-1)\delta}{(z-v)+a\exp(i\phi)(\epsilon\exp(-i\delta)-1)/\sigma}
\right)d\delta
\end{multline}
Since for all real $\phi,\delta$ we have $\abs{\exp(i\phi(\epsilon\exp(i\delta)-1)/\sigma}=1$ define
$\exp(i\psi)= (\epsilon\exp(i\delta)-1)/\sigma$ for a real-valued $\psi=\psi(\epsilon,\delta).$
Since the function $f(v+\rho\sigma\exp(i\phi))/\sigma$ is even with respect to $\delta$ changing in the second integral $\delta$ with $-\delta$
renders
\begin{equation}\label{dolchny2}
F(\bar{v},z)=-\rho(1-\epsilon^2)\frac{\epsilon^{1-q}}{4\pi^2}\int_0^{2\pi}\int_0^{2\pi}\frac{f(v+\rho\sigma\exp(i\phi)}{\sigma}
\frac{\exp(i\phi)\exp(-i(q-1)\delta)d\delta d\phi}{z-v+a\exp(i(\psi+\phi))}
\end{equation}
It is possible to verify (see Danchenko \& Dolzhenko \& \cite{dolchenkodanchenko2007}, p.5189) that this integral vanishes for all $z\in \{\abs{z-v}<(1-\epsilon)\rho\},$ if one 
rewrites
$(z-v+a\exp(i(\psi+\phi)))^{-1}$ as
$J:=(a^2/(\overline{z-v})+a\exp(i(\psi+\phi)))^{-1},$ where we further have
\begin{multline}
\frac{1}{z-v+a\exp(i(\psi+\phi)}=\frac{1}{\frac{a^2}{\overline{z-v}}+a\exp(i(\psi+\phi))}=\\
\frac{\exp(-i(\psi+\phi))}{a}\frac{a^2-\abs{z-v}^2}{\abs{z-v+a\exp(i(\psi+\phi))}^2}
\end{multline}
Using the relation
\begin{equation}
\exp(i\psi)=\frac{\epsilon\exp(i\delta)-1}{\sigma}
\end{equation}
and subtracting the aforementioned zero integral from Eqn.(\ref{dolchny2}) gives
\begin{multline}
F(\bar{v},z)=
\frac{\epsilon^{1-q}}{4\pi^2}\times\\
\int_0^{2\pi}\int_0^{2\pi}\frac{f(v+\rho\sigma\exp(i\phi))}{\sigma^2}
(\exp(-i(q-1)\delta)-\epsilon \exp(-iq\delta))
\frac{(a^2-\abs{z-v}^2)d\delta d\phi}{\abs{z-v+a\exp(i(\psi+\phi))}^2}
\end{multline}
Replacing $q$ by $q+k$ and multiplying by $\epsilon^k$, $k\in \Z_{\geq 0}$, and summing over $k$ renders 
with $\sigma=\sigma(\epsilon,\delta)=\sqrt{1+\epsilon^2-2\epsilon\cos\delta},$
$\exp(i\psi)=(\epsilon\exp(i\delta)-1)/\sigma,$
$\rho>0,\epsilon\in (0,1),$ $a=\rho(1-\epsilon^2),$ $z\in \{\abs{z-v}<(1-\epsilon)\rho\}$,
the following
\begin{multline}\label{dolchny3}
F(\bar{v},z)=\frac{\epsilon^{1-q}}{4\pi^2}\times\\
\int_0^{2\pi}\int_0^{2\pi}\frac{f(v+\rho\sigma\exp(i\phi)}{\sigma^2}
(\exp(-i(q-1)\delta)-\epsilon \exp(-iq\delta))
\frac{(a^2-\abs{z-v}^2)d\delta d\phi}{\abs{z-v+a\exp(i(\psi+\phi))}^2}
\end{multline}
In fact Eqn.(\ref{dolchny3}) holds true in the larger disc $\{\abs{z-v}<(1-\epsilon^2)\rho\}$ due to the fact that the integral is harmonic in this larger disc.
Since the function $f(v+\rho\sigma\exp(i\phi))/\sigma$ is even with respect to $\delta$, setting $z=v$ gives
\begin{equation}
f(v)=\frac{\epsilon^{1-q}}{4\pi^2}(1-\epsilon^2) \int_0^{2\pi}  \frac{\cos(q-1)\delta}{\sigma^2(\epsilon,\delta)}d\delta
\int_0^{2\pi}f(v+\rho\sigma\exp(i\phi))d\phi
\end{equation}
This gives an improved version of Corollary \ref{maxprincedolch}.
\begin{theorem}[Annular maximum principle]
For $q\in \Z_+,$, a given $q$-analytic $f$, and $K:=K(z,\rho,\epsilon)=\{t: (1-\epsilon)\rho <\abs{t-z}< (1+\epsilon)\rho  \}$
we have for $z\in \C$
\begin{equation}
\abs{f(z)}\leq \epsilon^{1-q}\norm{f}_K,\quad \rho>0,\epsilon \in (0,1)
\end{equation}
\end{theorem}
It must be noted that there seems to exists predated similar results on 
analogues of Cauchy estimates and annular maximum principles obtained by Balk \cite{ca1}.
For instance Balk proved the following. 
\begin{theorem}[See Balk \cite{ca1}, Theorem 1.5, p.204]\label{balk204}
Let $B\subset\subset \C$, be a ball centered at $p_0$ of radius $R$ and let $f$ be $q$-analytic function on $B$
with representation $f(z)=\sum_{j=0}^{q-1}f_j(z)\bar{z}^j.$
Let $0<R_0<R.$
Set $M_0:=\sup_B \abs{f}$ and $M_1=\sup_{R_0\leq \abs{z-p_0}<R_1} \abs{f}<\infty$. 
If $M_1<\infty$ then 
there exists a constant $\lambda_{R_0}$ depending only on $q$ and $R_0/R$ such that 
\begin{equation}
\abs{\bar{z}^k f_{k}(z)}\leq \lambda_{R_0} M_0,\quad \abs{z-p_0}<R_0
\end{equation}
and for each $l$ there exists a constant $\lambda_{R_0,k,l}$ depending only on $q,k,l,R_0$ such that 
\begin{equation}
\abs{\partial_{\bar{z}}^k\partial_z^l f_{k}(z)}\leq \lambda_{R_0,k,l} M_0 R_0^{-k-l},\quad \abs{z-p_0}<R_0
\end{equation}
Furthermore, if $M_0<\infty$ then there exists a constant $\lambda_{k,l}$ depending only on $q,k,l$ such that at $p_0$,
\begin{equation}\label{nicolescu}
\abs{\partial_{\bar{z}}^k\partial_z^l f_{k}(p_0)}\leq \lambda_{k,l} M_0 R^{-k-l}
\end{equation}
\end{theorem}
 Balk gives some consequences of this most of which follow from the analogous results (e.g.\ the derivative estimates and the annular maximum principle proved above in this section).
 \begin{theorem}
 Let $U$ be a bounded Jordan domain. Set $U_h:=\{ z\in U:\mbox{dist}(z,\partial U)<h\}$ for $h>0.$
 let $f$ be $q$-analytic function on $U$
 with representation $f(z)=\sum_{j=0}^{q-1}f_j(z)\bar{z}^j.$ Let $V$ be a closed subset of $U.$
 If $\abs{f}\leq M$ on $U_h$ then there exists a constant $L$ depending only on $U,h,q,V,$ such that
 \begin{equation}
 \abs{f(z)}\leq LM,\quad \abs{f_j(z)}\leq LM,\quad j=0,\ldots q-1
 \end{equation}
 on $V.$
 Also for each pair $k,l$ of integers there exists $L_{k,l}$ depending only on $U,h,q,V,k,l$ such that
 \begin{equation}
 \abs{\partial_{\bar{z}}^k\partial_z^l f_{k}(p_0)}\leq L_{k,l} M
 \end{equation}
 on $V.$
 \end{theorem} 
  Obviously this implies     
      \begin{corollary}[See Balk \cite{ca1}, Cor. 5, p.206]\label{balk206b}
         Let $\Omega\subset\C$ be a domain and let $q\in \Z_+.$
         If a sequence $\{f_j\}_{j\in \Z_+}$ of $q$-analytic functions converges uniformly inside $\Omega$
         then the any pair of nonnegative integers $k,l$
          $\partial_{\bar{z}}^k\partial_z^l f_j(z)$ converges uniformly to 
          $\partial_{\bar{z}}^k\partial_z^l f(z)$ on $\Omega.$
         \end{corollary}

We recall the following theorem which shall be used repeatedly and a proof can be found in e.g.\ 
H\"ormander \cite{hormander}, Cor.1.2.5, p.4.
\begin{theorem}[See also Theorem \ref{vitaliporterthm}]\label{hormkonvthm}
	Let $U\subset\Cn$ be a domain. If $\{f_j\}_{j\in \Z_+}$ is a sequence of holomorphic functions on $U$
	such that $f_j\to f$ uniformly on compacts of $U$, then $f$
	is holomorphic on $U$.
\end{theorem}
\begin{proof}
	Let $\Omega\subseteq U$ be a domain whose boundary consists of a finite number of $C^1$-smooth Jordan curves.
		Applying Stokes' formula to $f(\zeta)/(\zeta-z)$ and using the notation $B_\epsilon:=\{\zeta\in \Omega:\abs{\zeta -z}<\epsilon\}$ (where $\epsilon$ is sufficiently small) yields
		\begin{equation}\label{ekvforhorm}
	\int_{B_\epsilon} 
	\frac{\partial_{\bar{\zeta}}h(\zeta)}{\zeta-z}d\mu(\zeta)=
	\int_{\partial \Omega} \frac{h(\zeta)}{\zeta-z} d\zeta+
	\int_0^{2\pi} f(z+\epsilon \exp(i\theta)) id\theta
	\end{equation}
	Since $(\zeta -z)^{-1}$ is integrable over $\Omega$ and $f$ is continuous at $\zeta$, letting $\epsilon\to 0$ yields
	for $z\in \Omega$
	\begin{equation}\label{uu1276}
	f(z)=\frac{1}{2\pi i}\left(
	\int_{\Omega} 
	\frac{\partial_{\bar{\zeta}}f(\zeta)}{\zeta-z}d\mu(\zeta)+
	\int_{\partial \Omega} \frac{f(\zeta)}{\zeta-z} d\zeta
	\right)
	\end{equation}
	Let $K\subset\Omega$ be a compact and let $\phi\in C^\infty_c(\Omega)$ satisfy $\phi=1$ on a neighborhood of $K$.
	Then $\partial_{\bar{z}} (\phi f)=f\partial_{\bar{z}} \phi$ thus by Eqn.(\ref{uu1276}) applied to $\phi f$ we have
	\begin{equation}
	\phi(z)f(z)=\frac{1}{2\pi i}
	\int 
	\frac{\partial_{f(\zeta)\bar{\zeta}}\phi(\zeta)}{\zeta-z}d\mu(\zeta)
	\end{equation}
	Since $\phi=1$ on a neighborhood of $K$ and since $\abs{z-\zeta}$ is bounded from below when $z\in K$ and $\zeta$ in supp$\partial_{\bar{z}}\phi$, differentiation of the last equation gives that there exists constants $C_j$, $j\in \N$ (depending on $K$ and $\Omega$) such that
	\begin{equation}
	\sup_{z\in K}\abs{\partial_z^j f(z)}\leq C_j\norm{f}_{L^1(\Omega)}
	\end{equation}
	Applying this to the sequence of holomorphic functions
	$(f_j-f_m),$ shows that
	$\partial_z f_j$ converges uniformly. Since $\partial_{\bar{z}}f_j=0$ it follows that $\partial_x f_j$ and
	$\partial_y f_j$ converge uniformly on $K$. We obtain that $f\in C^1$ and $\partial_{\bar{z}} f= \partial_{\bar{z}}f_j =0.$ Since $K$ was an arbitrary compact this completes the proof.
\end{proof}

The following is proved in this book using other techniques, see Proposition \ref{ahernbrunakonsekvens}, but seems to have
 been known by Balk \cite{balk69} earlier and can obviously also be proved using the annular maximum principle. 
 \begin{theorem}[See Balk \cite{ca1}, Cor. 3, p.206]
 Let $\Omega\subset\C$ be a domain and let $q\in \Z_+.$
 If a sequence $\{f_j\}_{j\in \Z_+}$ of $q$-analytic functions converges uniformly inside $\Omega$
 then the limit function is also $q$-analytic in $\Omega.$
 \end{theorem}

\section{Classical types of maximum principle}

\begin{definition}
Let $p_0\in \C$. We shall say that a given property holds true {\em for all sufficiently small domains containing $p_0$}, if there is an $r>0$, such that the given property holds true for all domains $D\ni p_0,$ satisfying $D\subset \{\abs{z-p_0}<r\}$ (where $z$ denotes the standard complex coordinate in $\C$).
\end{definition}
\begin{definition}[Local maximum modulus property]\index{Local maximum modulus property}
Let $\Omega\subset\C$ be an open subset, let $f\in C^0(\Omega)$ and 
denote by $\Sigma\subset\Omega$, the set of points at which $\abs{f}$ attains strict local maximum.
We say that $f$ satisfies the {\em local maximum modulus property}, on $\Omega$,
if for each point $p_0\in \Omega$, it holds true that: given any sufficiently small, relatively compact subdomain $D\Subset \Omega,$ with $p_0\in D,$ $f$ satisfies,
\begin{equation}
\max_{z\in \overline{D}}\abs{f(z)}=
\max_{z\in \Sigma\cup \partial D}\abs{f(z)}.
\end{equation}
If every member of a family, $\mathcal{F},$ of continuous functions, satisfies the local
maximum modulus on $\Omega$, then we simply say that $\mathcal{F}$ satisfies the local maximum modulus property with respect to $\Omega.$
\end{definition}

\begin{definition}[Semi-analytic sets]\index{Semi-analytic set}
Let $\Omega\subset \R^N,$ $N\in \Z_+,$ be a real-analytic submanifold. We say that a subset $\mathcal{Z}\subset\Omega$ is {\em semi-analytic} if for every point $p\in \Omega$ there exists a neighborhood $U\ni p$, together with nonnegative integers $n,m,$ and real-analytic functions $f_{j}, g_{jk},$ such that,
\begin{equation}
\mathcal{Z}\cap U =\bigcup_{j=1}^m \{ z\in \Omega\colon f_{j}(z)=0,\, g_{jk}(z)>0, k=1,\ldots, n\}.
\end{equation}
\end{definition}
For an introduction to the theory of semi-analytic sets, see e.g.\
Krantz \& Parks \cite{krantzparks}, Sec.5.4, Lojasiewicz \cite{lojaensemble} (where the above version of the definition appears on p.48), \cite{lojabook}.
\begin{example}
	Let $\Omega\subset \C,$ be a subdomain. If $F:\Omega\to \R,$ is a real-analytic function on $\Omega$, then the critical set,
	$\{\nabla F=0\},$ is a semi-analytic subset of $\Omega,$ because $f_{1}:=\frac{\partial F}{\partial x},$
	$f_{2}:=\frac{\partial F}{\partial y},$ for $x:=\re z,$ $y:=\im z,$ are real-analytic functions and
	$\{\nabla F=0\}=\{ z\in \Omega\colon f_1(z)=0, f_2=0\}.$
\end{example}
In general a $q$-analytic function $f$ will be complex valued and so cannot be subharmonic, but in special cases such as $f(z)=z+\bar{z}$ or
$f(z)=z\bar{z}$, we obtain the local maximum modulus property directly from the maximum principle for subharmonic functions.
\begin{example}
Let $\Omega\subset \C,$ be an open subset and let $f=a(z)+\bar{z}b(z),$
where
$a,b$ are holomorphic on $\Omega.$ 
If $b(z)=\mu$ for some constant $\mu\in \C$, we have 
$\abs{f}^2=\left(\abs{a}^2+\abs{\mu}^2\abs{z}^2\right) +2\re(\bar{\mu}za(z)),$ which is clearly subharmonic since 
$\re(\bar{\mu}za(z))$ is harmonic 
and the expression in the parenthesis is subharmonic, 
thus 
for every relatively compact subdomain $D\Subset \Omega,$ it holds true that,
\begin{equation}
\max_{z\in \overline{D}}\abs{f(z)}=
\max_{z\in \Sigma\cup \partial D}\abs{f(z)},
\end{equation}
where in fact $\Sigma=\emptyset,$
and furthermore one can also assume weak local maximum does not occur
unless $\abs{f}^2$ is constant near that point see e.g.\
Fraenkel \cite{fraenkel}, p.53.
The above is a special case of a function of the form (here $\mathscr{O}(\Omega)$ denotes the set of holomorphic
functions on $\Omega$),
\begin{equation}
f(z)=P(z)+\overline{Q}(z),\quad  P\in \mathscr{O}(\Omega),\quad Q(z)=\sum_{j=0}^{n-1} c_j z^j,
\end{equation}
for a positive integer $n$ and complex constants $c_0,\ldots,c_{n-1}$. Obviously,
$f$ will be $n$-analytic and $\abs{f}^2$ will be subharmonic because
$\re f$ and $\im f$ are both harmonic and using the relation
$\Delta=4\partial\overline{\partial}$, we see that
$4\Delta \abs{f}^2 =2\abs{\partial_x \re f}^2+ 2\abs{\partial_y \re f}^2 +2\abs{\partial_x \im f}^2+2\abs{\partial_y \im f}^2\geq 0,$ where we denote $x=\re z,$ $y=\im z.$
Another easy example of a set of $n$-analytic functions with subharmonic square modulus, is given by functions of the form
$P(z)\cdot \overline{P}(z),$ for a holmorphic polynomial $P$ of order $(n-1).$
\end{example}
\begin{example}\label{example3}
Let $\Omega\subset \C,$ be an open subset and let $f=1-z\bar{z}.$ Obviously 
$f$ is a bianalytic function on $\Omega,$ and attains {\em strict} local maximum
at $p_0=0.$
\end{example}
\begin{example}\label{example4}
Let $\Omega=\{\abs{z}<1\}\subset\C,$ be an open subset and let $f(z):=z-z^2\bar{z}.$ Then
$f$ is bianalytic on $\Omega,$ and,
\begin{equation}
\abs{f}=\abs{z}(1-\abs{z}^2).
\end{equation} 
Setting $r:=\abs{z},$ it becomes clear that, $\abs{f}=\abs{f}(r)=r-r^3, r<1,$
attains weak local maximum for $r=1/\sqrt{3},$ i.e.\ on a circle enclosing $p_0.$
Hence setting $p_0=0,$ and $D=\{\abs{z}<1/2\}$ we obtain a function that disobeys the local boundary
maximum modulus principle, yet does not attain strict local maximum at any point.
\end{example}
Example \ref{example3} shows that, the obviously necessary
exclusion of {\em strict} local maxima, is a non-void condition, 
for any maximum principle for bianalytic functions. Example \ref{example4} shows that {\em locality} is necessary for 
a boundary maximum modulus principle for bianalytic functions.
Since any $n$-analytic function is automatically $(n+1)$-analytic our examples suffice for sustaining the necessity of locality and exclusion of strict local maxima for a maximum modulus principle when $n\geq 2.$
\begin{proposition}
Let $\Omega\subset \C,$ be an open subset, $p_0\in \Omega,$ and let $f$ be a 
complex-valued function on $\Omega,$ whose real and imaginary parts are real analytic. 
Denote by $\Sigma\subset \Omega$, the set of points at which $\abs{f}$ attains strict local maximum.
Then for every sufficiently small domain $D\Subset \Omega,$ containing
$p_0\in D,$ it holds true that,
\begin{equation}\label{condek}
\max_{z\in \overline{D}}\abs{f(z)}=\max_{z\in \Sigma\cup \partial D}\abs{f(z)}.
\end{equation}
\end{proposition}
\begin{proof}
	Let $\re f (z)=:A(z),$ $\im f=: B(z).$
	We introduce some notations regarding our
	real-analytic function $f\bar{f}=\abs{f}^2.$ 
	Denote,
	\begin{equation}\label{denotar}
	F(z):=\abs{f(z)}^2, \quad \mathfrak{C}:=\{z\colon F(z)\mbox{ is a weak local maximum}\}
	\end{equation} 
	\begin{equation}
	\mathcal{Z}:=\{z\in \Omega \colon \nabla F(z)=0\},
	\end{equation}
	Note that $\mathcal{Z}$ is, in particular, a semi-analytic set.
	The following is a special case of a well-known result from real-analytic geometry.
	\begin{lemma}\label{lemma1mp} [Special case of the Curve Selection Lemma, see e.g.\ Lojasiewicz \cite{lojaensemble}, p.103]
	Let $f=A+i B,$ $\Omega,$ $F,$ and $\mathcal{Z},$ be as above, 
	$0\in \Omega$. If $0$ is a point at which $\mathcal{Z}$ accumulates,
	then there exists a real-analytic curve $\kappa\colon (0,1)\to \Omega,$
	$\kappa(0)=p,$ such that $\kappa((0,1))\subset \mathcal{Z}.$
	\end{lemma}
	Note in Lemma \ref{lemma1mp}, that because $\mathcal{Z}$ is closed also $p\in\mathcal{Z},$ and furthermore the conclusion does not exclude that $ \mathcal{Z}$ contains an open neighborhood of $p$ in $\Omega$.
	In particular, we know that in Lemma\ref{lemma1mp}, $F\circ \kappa$ is constant since $\frac{\partial}{\partial t}(F\circ \kappa)(s)=(\nabla F)|_{\kappa(s)}\cdot \frac{\partial\kappa}{\partial t}(s),$ $s\in (0,1).$ We point out that $\mathcal{Z}$ is not itself required to locally be a submanifold, which is clearly seen by the example $\{ x^2-y^2 =0\}$
	which fails to be a submanifold of $\R^2$ at $(0,0).$
	We shall also use the following known result.
	\begin{lemma}[See e.g.\ Massey \& Tr\'ang \cite{massey}, Theorem 6.7]\label{lemma2mp}
	Let $f=A+i B,$ $\Omega,$ $F,$ and $\mathcal{Z},$ be as above, 
	$p_0\in \Omega$.
	Then, there exists an open neighborhood, $U,$ of $p_0,$ in $\Omega,$ such that
	the set $\Omega\cap \mathcal{Z},$ is contained in the preimage set $F^{-1}(F(p_0))$.
	\end{lemma}
	Next we consider the following cases (note that (iii) is complementary to the union of (i) and (ii), thus the three cases cover all possibilities):
	\begin{itemize}
		\item[Case (i)] \textit{The point $p_0$ belongs to $\Sigma$.} Then, by definition, there exist a small domain $D_0\subset D,$
		$p_0\in D_0,$ such that $\max_{z\in \overline{D}_0}\abs{f(z)} =\abs{f(p_0)},$ so $f$ would satisfy equation \ref{condek}
		with $D$ replaced by any subdomain of $D_0$, containing $p_0.$ 
		\item[Case (ii)] \textit{There exists a domain $D_0\subset D,$ $p_0\in D_0$, such that (i) fails and $\mathfrak{C}\cap (D_0\setminus \{p_0\})=\emptyset.$} Then clearly $\abs{f}^2$ does not attain weak local maximum\footnote{If it where to attain weak local maximum at $p_0$ then (i) must hold true because in such case $p_0$ would be an isolated maximum rendering point.} at any point of $D_0$, and because, by continuity, it does attain maximum on $\overline{D}_0,$ the maximum will be attained on the boundary $\partial D_0.$ This can be repeated for any subdomain of $D_0$ containing $p_0.$
		\item[Case (iii)] \textit{Case (i) fails and there exists a sequence $\{z_j\}_{j\in \Z_+},$ $z_j\in \mathfrak{C},$ such that $z_j\to p_0.$}
	\end{itemize}
	It is clear that the result will be proved once we can verify it under the conditions of Case (iii). So let us assume Case (iii) holds true. 
	Obviously, we have in Case (iii),
	\begin{equation}
	p_0\in \overline{\mathcal{Z}}.
	\end{equation}
	Let $U_0$ be as in Lemma \ref{lemma2mp} and let $U\subset U_0$
	be a domain containing $p_0.$ 
	Let $\kappa$ be as in Lemma \ref{lemma1mp}.
	We can choose $U$ sufficiently small such that,
	for every subdomain $D\subset U,$ containing $p_0,$
	we have, 
	\begin{equation}\label{anv2}
	\partial D\cap \kappa\neq \emptyset.
	\end{equation}
	Assume, in order to deduce a contradiction, that there would exist a choice of domain $D\ni p_0$ in $U$, such that,
	\begin{equation}\label{ceq}
	\sup_{z\in D}F(z)>\max_{z\in \partial D}F(z).
	\end{equation}
	Then $F$ must attain the value $\max_{z\in \overline{D}}F(z)$ at a point $p\in D\cap \mathfrak{C}\subset D\cap \mathcal{Z}.$ But by using 
	Lemma \ref{lemma2mp}, we have made sure that 
	$F(\mathcal{Z}\cap U)=\{F(p_0)\},$
	whence $\max_{z\in \overline{D}}F(z)=F(p_0).$
	But by construction $F(z)=F(p_0),$ for all $z\in \kappa,$ thus by
	Eqn.(\ref{anv2}) we arrive at a contradiction to Eqn.(\ref{ceq}).
	We conclude that for all sufficiently small domains $D\ni p_0$,
	\begin{equation}\label{ceq0}
	\max_{z\in \overline{D}}\abs{f(z)}^2\leq \max_{z\in \partial D}\abs{f(z)}^2,
	\end{equation}
	which in turn verifies Eqn.(\ref{condek}). This 
	verifies the conclusion under the conditions of Case (iii), and thus
	completes the proof.
\end{proof}
Since an $n$-analytic function has real-analytic real and imaginary parts, we obtain the following (see Daghighi \& Krantz \cite{daghighikrantz}).
\begin{theorem}\label{mpmainresult}
Let $\Omega\subset \C,$ be an open subset, $p_0\in \Omega,$ and let $f$ be an $n$-analytic function $n\geq 2.$ 
Denote by $\Sigma\subset \Omega$, the set of points at which $\abs{f}$ attains strict local maximum.
Then for every sufficiently small domain $D\Subset \Omega,$ containing
$p_0\in D ,$ it holds true that,
\begin{equation}\label{condek1}
\max_{z\in \overline{D}}\abs{f(z)}=\max_{z\in \Sigma\cup \partial D}\abs{f(z)}.
\end{equation}
\end{theorem}
By Example \ref{example3} and Example \ref{example4} we can also make the following remark.
\begin{remark}
If, in Theorem \ref{mpmainresult}, $\Sigma$ is removed from Eqn.(\ref{condek1}) then, 
for each $n\geq 2,$ there exists an $n$-analytic function $f$ on $\Omega,$ and a subdomain $D$, for which
Eqn.(\ref{condek1}) fails. If instead
the requirement that $D$ be sufficiently small is dropped, then, 
for each $n\geq 2,$ there exists an $n$-analytic function $f$ on $\Omega,$ together with a choice of $D\Subset \Omega$, for which Eqn.(\ref{condek1}) fails.
\end{remark}

Reduced $q$-analytic functions have the following know representation (which is stated without proof in
Balk\cite{ca1}, p.203 Thm 1.4, and also without proof in  Balk \& Zuev \cite{balkzuev}, p.204. It is possible
that a proof could be present in the conference publication Balk \cite{balk71}, but we have not verified this).
\begin{theorem}
	\label{balk71lemma}
Let $g$ be a reduced $q$-analytic function on the disc $\{\abs{z}<R\},$ for a positive $R.$
Let $R_0,R_1,\ldots,R_{q-1}$ satisfy $0<R_0<R_1<\cdots<R_{q-1}<R.$
Then
$g(z)=\sum_{j=0}^{q-1} P_j(\abs{z}^2)g_j(z)$
where each $g_j$ is holomorphic and coincides with $g$ on $\{\abs{z}=R_j\},$
and 
\begin{equation}
P_j(t)=\frac{(R_0^2-t)\cdots(R_{j-1}^2-t)(R_{j+1}^2-t)\cdots(R_{q-1}^2-t)}{(R_0^2-R_j^2)\cdots(R_{j-1}^2-R_j^2)(R_{j+1}^2-R_j^2)\cdots(R_{q-1}^2-R_j^2)}
\end{equation} 
\end{theorem}
\begin{proof}
Set $g_j:=g|_{\abs{z}=R_j}$, let $\tilde{g}_j$ be the unique holomorphic function
on $\{\abs{z}<R\}$ such that $\tilde{g}_j|_{\abs{z}=R_j} =g_j$ and set
$G(z):=\sum_{j=0}^{q-1} P_j(\abs{z}^2) \tilde{g}_j(z).$
It is clear that $f(z)-G(z)$
is a reduced $q$-analytic function on $\Omega$. Furthermore
$(g(z)-G(z))=0$ on the union $\bigcup_{j=0}^{q-1}\{\abs{z}=R_j\}$
because on each circle
$\{\abs{z}=R_j\}$ we have $P_j(R_k^2)=0$ for $k\neq j,$ $P_j(R_j)=1,$
hence $g(z)=P_j(R_j)\tilde{g}_j(z)=G(z),$ for $z\in \{\abs{z}=R_j\}.$
By Proposition \ref{uniqueness203} this implies
that $g-G\equiv 0$ on $\Omega.$ This completes the proof.
\end{proof}
A non-constant holomorphic function is locally open.
This is no longer true in general for $q$-analytic functions when $q>1$.
\begin{example}
If $P:\C\to \C$ is a complex polynomial of order $q$, then $f:=\abs{P}^2=P\bar{P}$
is a $(q+1)-$analytic function which is real-valued thus not open.
\end{example}

\begin{definition}
Let $f=u+iv$ be a differentiable function on $\C,$ where $u,v$ are the real and imaginary parts respectively.
Let $z=x+iy$ denote the complex Euclidean coordinate for $\C.$ A point $p_0$ will be called a {\em critical point}
if $|\mbox{Jac}(f)|(p_0)=0,$ 
where $\mbox{Jac}(f)$ denotes
the Jacobian matrix 
$$\begin{bmatrix}
\frac{\partial u}{\partial x} & \frac{\partial u}{\partial y}\\
\frac{\partial v}{\partial x} & \frac{\partial v}{\partial y}
\end{bmatrix}
$$
\end{definition}
By identifying $\C$ with $\R^2$ we can by the inverse function theorem conclude that, for a point $p_0\in \C$,
a differentiable function $f$ satisfying $|\mbox{Jac}(f)|(p_0)\neq 0$ is locally open at $p_0.$
Here is an example of a $2$-analytic function which is locally open except on a circle. 
\begin{example}
Consider $f(z)=z(1-z\bar{z})$ on the open unit disc. 
It has local weak maximum modulus on the circle $\abs{z}=\frac{1}{\sqrt{3}},$ thus not locally open at each point of that circle 
(and obviously it is zero on the unit circle and the origin).
On the other hand we have $\re f=x(1-x^2-y^2), \im f=y(1-x^2-y^2)$ so that $|\mbox{Jac}(f)|=(1-3x^2-y^2)(1-3y^2-x^2) -4x^2y^2$
$=1-4x^2-4y^2+6x^2y^2+3x^4+3y^4.$
Denoting the standard polar representation of the holomorphic coordinate as $z=r\exp(i\theta),$ we have
$|\mbox{Jac}(f)|=1-4x^2-4y^2+6x^2y^2+3x^4+3y^4=1-4r^2+3r^4.$ The solution to the equation $1-4t+3t^2=0$
are $t=1$ or $t=\frac{1}{3},$ thus the zeros of $|\mbox{Jac}(f)|$ are obtained when
$r^2=1$ or
$r^2=\frac{1}{3}.$ Thus on the open unit disc $|\mbox{Jac}(f)|$ has as zero set the circle $\abs{z}=\frac{1}{\sqrt{3}}$ and by the inverse function theorem it is locally open at all other points within the unit disc. 
\end{example}

Here is an example of a $q$-analytic functions that is locally open but
has non-discrete critical set.
\begin{example}
Consider the $4$-analytic function $f(z,\bar{z})=f(x,y)=\frac{1}{8}(z+\bar{z})^3+\frac{1}{2}(z-\bar{z})=x^3+iy.$
We have 
$|\mbox{Jac}(f)|=3x^2$ which vanishes precisely on the vertical line $\{ x=0\}.$
In particular, $f$ is locally open at all points where $x\neq 0.$ Let us check points on the critical set.
For fixed $y=y_0,$ and $\epsilon>0,$ set $D_{y_0}=\{ x+iy_0\colon \abs{x}<\epsilon\}$.
Then $f(D_{y_0})=\{ (-\epsilon^3,\epsilon^3)+iy_0)\}.$ Hence varying $y$ in an interval $I_{\delta}:=(y_0-\delta,y_0+\delta)$ for $\delta>0,$ 
the sets $\{ f(D_{y})\colon y\in I_\delta\}$ foliate an open neighborhood of $(0,y_0)$. Since $\epsilon$ and $\delta$ where arbitrary 
the sets $\{D_{y}\colon y\in I_\delta\}$ form a neighborhood basis of $(0,y_0)$ (in the sense that any open neighborhood of $(0,y_0)$ contains
a set of the form $\{D_{y}\colon y\in I_\delta\}$) proving that $f$ is locally open at $(0,y_0).$ Since $y_0$ was arbitrary $f$ is locally open.
\end{example}

The following example considers a type of minimum modulus principle.
\begin{example} 
Let $f(z):=C_0+d_0(z) +\sum_{l=1}^{q-1} c_l \bar{z}^l$ where $C_0$ and the $c_l$ are complex constants and $d_0$ 
holomorphic ($d_0(0)\neq 0$ unless $d_0\equiv 0$, i.e.\ all constant terms are included in $C_0$). Note that $f$ is a harmonic (complex valued) function. Assume that $f$ attains weak local minimum modulus at $p_0\in \C.$ Denote the standard polar representation of the holomorphic Euclidean coordinate by $z=r\exp(i\theta).$
By the mean value property for harmonic functions we have for all small $r>0$
$f(p_0)=\frac{1}{2\pi}\int_{\abs{z-p_0}=r} \re f d\theta +i\frac{1}{2\pi}\int_{\abs{z-p_0}=r} \im f d\theta =\frac{1}{2\pi}\int_{\abs{z-p_0}=r} f d\theta$.
If $f(p_0)\neq 0$ then either $\abs{f}|_{\{\abs{z-p_0}=r\}}=\abs{f(p_0)}$ or there exists points $q,q'\in \{\abs{z-p_0}=r\}$ such that
$\abs{f(q)}<\abs{f(p_0)}$, $\abs{f(q')}>\abs{f(p_0)}$. But $\abs{f}$ cannot be constant near any point $z$ such that 
$|\mbox{Jac}(f)(z_0)|\neq 0$, thus if $\abs{f(p_0)}$ is a weak local minimum then it must be zero or $p_0$ will be an accumulation point of 
$\{ |\mbox{Jac}f|=0\}$. That is to say, if $f$ has only isolated critical points then $f(p_0)=0.$
\end{example}

\begin{theorem}[Daghighi \cite{daghighilocopen}]
Let $f$ be a $q$-analytic function on $\C.$
If all critical points of $f$ are isolated then $f$ is locally open.
\end{theorem}
\begin{proof}
Let $p_0$ be a critical point of $f$ in the sense that $|\mbox{Jac}(f)|(p_0)=0.$
Denote by $u$ the real part of $f$ and by $v$ the imaginary part of $f,$ i.e.\ $f=u+iv.$
If there exists a continuous curve $\kappa[0,1)\to \C$, differentiable on $(0,1)$, $\kappa\subset \{ f=f(p_0)\}$ such that $\kappa(0)=p_0$
then for $t\in (0,1),$ $0=\partial_t (u\circ \kappa)(t)=\nabla u \cdot \dot{\kappa}$ and 
$0=\partial_t (v\circ \kappa)(t)=\nabla v \cdot \dot{\kappa}$.
Hence $\nabla u$ is parallel to $\nabla v$ so the Jacobian of $f$ (in the above sense) at each point of $\kappa$ has two linearly dependent columns. Thus
$|\mbox{Jac}(f)|(\kappa(t))\equiv 0$.
This implies that $\kappa\subset\{|\mbox{Jac}(f)|=0\}$, i.e.\ $p_0$ is not an isolated critical point.
The contrapositive statement is that if $p_0$ is an isolated critical point then the restriction of $f$ to any sufficiently small circle centered at
$p_0$ satisfies $f-f(p_0)\neq 0.$
So let $F$ be a $q$-analytic function which has only isolated critical points and let $p_0$ be a critical point. Set $f=F-F(p_0)$, so that $f(p_0)=0$. By what we have already done
there exists a circle $C$, 
centered at $p_0$ such that $f|_C \neq f(p_0)$. Set $2\epsilon=\min_C \abs{f}$. 
 
Let $\eta\in \{\abs{z}<\epsilon\}$.
Then $\abs{f-\eta}>\epsilon$, whereas $\abs{(f-\eta)(p_0)}=\epsilon$. So $\abs{f-\eta}$ attains minimum
in the disc $\hat{C}=\{\abs{z-p_0}<\epsilon\}$.

\begin{lemma}\label{lemma1}
Weak local minimum modulus of a continuous function $G$ at a point $z_0$ implies that $z_0$ is a critical point in the sense that $|\mbox{Jac}(G)|(z_0)=0$.
\end{lemma}
\begin{proof} 
Weak local minimum modulus implies that for an open neighbourhood of $z_0$, $G$ maps the interior point $z_0$ to the boundary of the image, whence
$G$ cannot be locally open. But if $|\mbox{Jac} (G)(z_0)|\neq 0$ then $G$ must necessarily be locally open at $z_0$ due to the inverse function theorem. This completes the proof. 
\end{proof}

Hence 
if $\abs{f-\eta}$ attains local minimum at $p_0$ then 
we can assume  $\abs{f-\eta}$ attains strict local minimum at $p_0$ (since $p_0$ is an isolated critical point).

\begin{lemma}
\label{reducedlemma}
If $g(z)$ is a reduced $q$-analytic function that attains strict local minimum modulus at $0$ then $g(0)=0.$
\end{lemma}

\begin{proof}
Let $g$ have the representation given in
Theorem \ref{balk71lemma} and let each $g_j=c_{j,0}+\sum_{k\geq k_j} c_{j,k}z^k,$ for constants $c_{j,k}$
where either $g_j\equiv c_{j,0}$ or $c_{j,k_j}\neq 0.$
Also let $K_j:=P_j(0)$ (note that $P_j(0)\neq 0$ for all $j$).
Since $0$ is an isolated critical point at least one $g_j$
must be nonconstant.
Set $C_0:=\sum_j K_j c_{j,0}$ and $\nu_0:=\min_{j} k_j.$
Note that $g(0)=C_0.$ If $\nu_0=0$ then $g(r,\theta)=g(r)$ thus $\abs{g(0)}$ can only attain strict local minimum if $g(0)=0.$ So we can assume $\nu_0\geq 1.$ 
Then near $0$ we have
\begin{equation}
\abs{g}^2=\abs{C_0}^2+2\sum_{j=1}^{q-1} \re(\overline{C_0}K_jc_{j,k_j}z^{k_j})+O\left(\abs{z}^{2\nu_0}\right)
\end{equation}
Set $G(z):=C_0+\sum_{j=1}^{q-1}K_j c_{j,k_j}z^{k_j}.$
Then $G$ is holomorphic and satisfies that $\abs{G}^2$
attains strict local minimum at $0$ if $\abs{g}^2$ does,
since near $0$ also
\begin{equation}
\abs{G}^2= \abs{C_0}^2+2\sum_{j=1}^{q-1} \re(\overline{C_0}K_jc_{j,k_j}z^{k_j})+O\left(\abs{z}^{2\nu_0}\right)
\end{equation} 
and the sum must be positive in a punctured neighborhood of $0$ if $\abs{g}^2$ attains strict minimum at $0$.
Since $G$ is locally open this implies that $C_0=0.$ This completes the proof.
\end{proof}

\begin{lemma}\label{lemma2}
If $f$ is $q$-analytic
such that $\abs{f}$ attains strict local minimum at $p_0$ then either $f(p_0)=0$ or  
$p_0$ is a non-isolated critical point (in the sense that it is an isolated zero of $|\mbox{Jac}(f)|$).
\end{lemma}
\begin{proof}
Without loss of generality assume $p_0=0.$
Let $f$ be polyanalytic of exact order $q$ i.e.\
$f(z)=C_0+\sum_{j=0}^{q-1}a_j(z)\bar{z}^j,$ for holomorphic $a_j$ such that
$a_0(0)=0,$ i.e.\ $f(0)=C_0.$ 
Let $z=r\exp(i\theta)$ be the standard polar representation.
It is straightforward (by simply replacing each occurrence of $z\bar{z}$ in the 
standard representation by $r^2$)
to obtain a unique representation according to
\begin{equation}\label{braeq}
f(z)=C_0 +\sum_{j=0}^{q-1} r^{2j}h_j(z)
\end{equation}
\begin{equation}
h_j(z):=d_j(z)+\sum_{l=1}^{q-1-j}c_{j,l}\bar{z}^l 
\end{equation}
where the $d_j$ are holomorphic such that either $d_0\equiv 0$ or $h_0(0)=O(\abs{z})$ (i.e.\ we make sure that $d_0$ has no constant term by including all constant terms into $C_0$). 
Note that for $\sum_{j=0}^{q-1}r^{2j}h_j(z)$ 
the real and imaginary parts can be identified (up to the constant $C_0$ for $h_0$)
as the harmonic components in the representation of $f$ using the Almansi expansion of the polyharmonic
real and imaginary parts respectively.
Since $f$ is assumed to attain strict local minimum modulus at $0$ we automatically know that ($q>1$ and) $f(z)$ is the sum of a constant term, $C_0$, and a term $\sum_{j=0}^{q-1}r^{2j}h_j(z)$, that behaves as $O(\abs{z}^{\nu_0})$, near the origin,
for some positive integer $\nu_0\geq 1.$ 
(To see this note that for complex 
constants $\alpha_k,\beta_k,$ $k=1,\ldots,M,$ we have 
$\sum_{k=1}^M (\alpha_{\nu_0} z^{k-\nu_0} +\beta_k\bar{z}^{k-\nu_0})=r^{\nu_0}\sum_{k=\nu_0}^M (\exp(i\nu_0\theta)\alpha_k z^{k-\nu_0} +\exp(-i\nu_0\theta)\beta_k\bar{z}^{k-\nu_0}),$
where $\nu_0$ is the smallest positive integer satisfying $\abs{\alpha_{\nu_0}}+\abs{\beta_{\nu_0}}\neq 0.$)
Set $\hat{c}_{j,l}:=\exp(2i\mbox{Arg}C_0)\overline{c_{j,l}}$ (note that $\frac{C_0}{\overline{C_0}}=\exp(2i\mbox{Arg}C_0)$). Thus
$\re(\overline{C_0}c_{j,l}\bar{z}^l)=\re(C_0 \overline{c_{j,l}\bar{z}^l})=\re(\overline{C_0}\hat{c}_{j,l}z^l).$
Set 
\begin{equation}
g(z):=C_0 +\sum_{j=0}^{q-1} r^{2j}\left(d_j(z)+\sum_{l=1}^{q-1-j}\hat{c}_{j,l}z^l \right)
\end{equation}
Then near $0$,
\begin{equation}
\abs{g}^2= \abs{C_0}^2 +2\sum_{j=0}^{q-1} r^{2j}\re\left(\overline{C_0}\left(d_j(z)+\sum_{l=1}^{q-1-j}c_{j,l}\bar{z}^l \right)\right)+O\left(\abs{z}^{2\nu_0}\right)
\end{equation}
and 
\begin{equation}
\abs{f}^2= \abs{C_0}^2 +2\sum_{j=0}^{q-1} r^{2j}\re\left(\overline{C_0}\left(d_j(z)+\sum_{l=1}^{q-1-j}c_{j,l}\bar{z}^l \right)\right)+O\left(\abs{z}^{2\nu_0}\right)
\end{equation}
Thus $g$ attains local minimum modulus at $0$ whenever $f$ does. Since $g$ is reduced Lemma \ref{reducedlemma} implies in such case that $g(0)=C_0=0,$ i.e.\ $f(0)=C_0=0.$
This completes the proof.
\end{proof}

By Lemma \ref{lemma2}  $\abs{f(z_0)-\eta}=0$ i.e.\ $f(p_0)=\eta.$ Recall that $p_0\in \hat{C}$ and that $\eta$ is an arbitrary point
in the ball $\{\abs{z}<\epsilon\}$. Hence $f$ is locally open at $p_0.$ This completes the proof.
\end{proof}
This immediately renders a strong maximum principle.
\begin{corollary}[Strong maximum principle]
Let $\Omega\subset\C$ be a domain, let $q\in \Z_+$ and let $f\in \mbox{PA}_q(\Omega).$
If all critical points of $f$ are isolated then $f$ does not attain weak local maximum on $\Omega.$
\end{corollary}

\section{Further notable results}
In this section we collect some noteworthy miscellaneous results, without giving the proofs, that are either directly involve integral representations or heavily dependent upon these in their method of proof. 
\begin{definition}
	Let $p\in \R$ and denote by $q_p$ the maximal integer such that $q_q=p+\epsilon$ for some $\epsilon\in (0,1].$
	Bounded (uniformly continuous) functions in a domain $\Omega\subset\C$ are called bounded
	(uniformly continuous) of order $0$. For $p>0$, a function $f(z) = f(x, y)$ is
	called bounded of order $p$ in $\Omega$ if it possesses all its partial derivatives of order $q_p$ or less
	in $\Omega$, and the derivatives of order less than $q_p$ are Lipschitz continuous of order $1$ 
	and all derivatives of order $q_p$ are Lipschitz continuous of order $(p-q_p)$. Moreover, if all the partial derivatives of
	order $q_p$ exist, and their modulus of continuity in $\Omega$ is of order $o(r^{p-q_p})$ as $r\to 0$, 
	then $f$ is called {\em uniformly continuous of order $p$} in $\Omega$. For $p < 0$, we call a function $f$ bounded
	of order $p$ if it is locally bounded in $\Omega$ and $\abs{f(z)}=O((\mbox{dist}(z,\partial\Omega))^p)$ as 
	$\mbox{dist}(z,\partial\Omega)\to 0.$ If $f$ is continuous in $\Omega$
	and $\abs{f(z)}=o(\mbox{dist}(z,\partial\Omega))^p)$ as $\mbox{dist}(z,\partial\Omega)\to 0$ then $f$
	is called {\em uniformly continuous of order $p$}\index{Uniformly continuous function of order $p$} in $\Omega$.
	The class of all functions bounded of order $p$ is denoted by $M^p(\Omega)$ and the class
	of all functions uniformly continuous of order $p$ is denoted $\mbox{CU}^q(\Omega)$.
\end{definition}
Dolzhenko \cite{dolzhenkoensam} proved using integral representations, the following.
\begin{theorem}
	Let $\Omega\subset \C$ be a domain whose boundary contains a closed Jordan analytic curve $\Gamma,$ such that
	$\mbox{dist}(\Gamma,(\partial\Omega)\setminus \Gamma)>0$. Let $f$ be a polyanalytic function $f$ in $\Omega$ that can be represented
	near $\Gamma$ of the form $f(z)=\sum_{j=0}^{q-1} (\bar{z}-\phi(z))^j a_j(z),$ where the
	$a_j$ are holomorphic and $\Gamma:=\{\phi(z)=\bar{z}\}$. Set $f_j:=(\bar{z}-\phi(z))^j a_j(z).$ Then for all $p\in \R$,
	in some one-sided neigborhood $Q\subset\Omega$ of $\Gamma$, we have the following:
	\begin{equation}
	f\in M^p(Q)\Leftrightarrow f_j\in M^p(Q), 0\leq j\leq q-1 \Leftrightarrow a_j\in M^p(Q), 0\leq j\leq q-1 
	\end{equation}
	\begin{equation}
	f\in \mbox{CU}^p(Q)\Leftrightarrow f_j\in \mbox{CU}^p(Q), 0\leq j\leq q-1 \Leftrightarrow a_j\in \mbox{CU}^p(Q), 0\leq j\leq q-1 
	\end{equation}
\end{theorem}
Trofymenko \cite{trofymenko2008} proves the following theorem for $q$-analytic functions which resembles a generalization of Morera type. 
\begin{theorem}
	Let $\Omega\subset\C$ be a domain and let $q>1.$ Then any $q$-analytic function
	$f$ on $\Omega,$ satisfies for all $z\in \Omega,$
	\begin{equation}
	\int_{l=1}^{q-1} \frac{1}{(l-1)!l!}\left(\partial_z^{l-1}\partial_{\bar{z}}^l\right) f(z) r^{2l} =
	\frac{1}{2\pi i} \int_{\abs{\zeta-z}=r} f(\zeta)) d\zeta
	\end{equation}
\end{theorem}

Hedenmalm \& Haimi \cite{hedenhaimi} proved for $q=2$ the following.
\begin{theorem}
	Let $f$ be a $2$-analytic function on the unit disc $D=\{\abs{z}<1\}$, 
	and let $\psi$ be subharmonic on $D.$
	Then
	\begin{equation}
	\abs{\overline{\partial} f(0)}^2\leq \frac{3}{\pi}\exp(-2\psi(0))\int_D\abs{f}^2\exp(2\psi)d\mu(z)
	\end{equation}
	and
	\begin{equation}
	\abs{f(0)}^2\leq \frac{8}{\pi}(1+6\abs{\psi(0)}^2)\exp(-2\psi(0))\int_D\abs{f}^2\exp(2\psi)d\mu(z)
	\end{equation}
	where $d\mu(z)$ is the standard area measure.
\end{theorem}
Let $\Omega\subset\C$ be a domain. For $f\in L^p(\Omega)$ set 
\begin{equation}
\omega(\Omega,p,f,\delta):=\sup\{\norm{f(\cdot)-f(z)}_{L^p(\{\abs{t-z}<r\})} :r\in (0,\delta], \{\abs{t-z}<r\}\subset\Omega\}
\end{equation}
We define $f\in\mbox{Lip}(\beta,p,\Omega),$ $\beta>0$, if $\omega(\Omega,p,f,\delta)=O(\delta^\beta)$ as $\delta\to 0.$
Ramazanov \cite{ramazanov1} (Theorem 5) proved the following.
\begin{theorem}
	Let $\Omega\subset\C$ be a bounded  domain.
	Let $q,p\in \Z_+,$ let $\beta$ be such that $p^{-1} <\beta\leq 1+2p^{-1},$ and let $f\in \mbox{PA}_q(\Omega)$.
	Set $\rho:=\rho(z,\partial\Omega)=\mbox{dist}(z,\partial\Omega).$
	Then 
	\begin{equation}
	f\in\mbox{Lip}(\beta,p,\Omega) \Leftrightarrow \abs{\partial_z f(z)}+ \abs{\partial_{\bar{z}} f(z)}=O(\rho^{\beta-1-\frac{2}{p}}),\mbox{ as }\rho\to 0
	\end{equation}
\end{theorem}
Ramazanov \cite{ramazanov2004} (Theorem 2) also proved the following.
\begin{theorem}
	Let $\Omega\subset\C$ be a simply connected domain, $0\in \Omega$.
	Let $q\in \Z_+$ and let $u(z)=u(x,y)$ be a $q$-harmonic function on $\Omega.$ Then there is a unique reduced
	$f\in \mbox{PA}_q(\Omega)$, say $f(z)=\sum_{j=0}^{q-1}h_j(z)\abs{z}^{2j},$
	for harmonic $h_j$ such that:\\
	(1) $u(x,y)=\re f(z).$\\
	(2) $u(0,0)=f(0)=h_0(0)$.\\
	(3) $\Delta^j u(x,y)|_{z=0}=4^j (j!)^2 h_j(0),$ $j=1,\ldots,q-1.$ 
	If $u\in L^p(\{\abs{z}<1\},\alpha)$, $\alpha>-1,$ then there is a constant $C(p,q,\alpha)$ depending only on $p,q,\alpha,$ 
	\begin{equation}
	\norm{f}_{p,\{\abs{z}<1\},\alpha}\leq C(p,q,\alpha)\norm{u}_{p,\{\abs{z}<1\},\alpha}
	\end{equation}
	where for a bounded domain $\Omega\subset\C$ we define
	\begin{equation}
	\norm{f}_{p,\Omega,\alpha)}:=\left(\int_{\Omega} (1-\abs{z}^2)^\alpha\abs{f(z)}^p d\mu(z) \right)^{\frac{1}{2}}
	\end{equation}
\end{theorem}
\begin{definition}
	Let $\Omega\subset\C$ be a bounded domain.
	Denote by 
	\begin{equation}
	A_m L_p(\Omega,\alpha):=\{f\in \mbox{PA}_q (\Omega): \norm{f}_{p,\Omega,\alpha)}<\infty\}
	\end{equation}
\end{definition}
\begin{theorem}[Ramazanov \cite{ramazanov1999} and \cite{ramazanov1}]
	Set
	\begin{multline}
	A_k^0(\{\abs{z}<1\},\alpha):=\left\{ f(z)=(1-\abs{z}^2)^{-\alpha}\partial_z^{k-1}(1-\abs{z}^2)^{k-1+\alpha}F(z):\right.\\
	\left.
	F\in \mathscr{O}(\{\abs{z}<1\}),f\in A_k L_p(\{\abs{z}<1\},\alpha)\right\}
	\end{multline}
	The space $A_m L_p(\{\abs{z}<1\},\alpha),$ $\alpha>-1,$ $p\in \Z_+,$ has the following decomposition
	\begin{equation}
	A_m L_p(\{\abs{z}<1\},\alpha) = A_1^0 L_p(\{\abs{z}<1\},\alpha)\oplus \cdots \oplus A_m^0 L_p(\{\abs{z}<1\},\alpha)
	\end{equation}
	Furthermore, for a bounded domain $\Omega\subset\C$, $p>1$ and $s\in\Z_{\geq 0}$ we have
	\begin{multline}
	\abs{\partial_z^s f(z)}\leq \\
	A_1(m,s,p,\alpha)\pi^{-\frac{1}{p}}
	\sum_{k=1}^m \frac{\Gamma(k+s+\alpha)(2k-1+s+\alpha)}{\Gamma(k)\Gamma(k+\alpha)}
	(\mbox{dist}(z,\partial\Omega))^{-\frac{2\alpha+2+sp}{p}} \norm{f}_{p,\Omega,\alpha)}
	\end{multline}
	\begin{multline}
	\abs{\partial_{\bar{z}}^s f(z)}\leq \\
	A_2(m,s,p,\alpha)\pi^{-\frac{1}{p}} 
	\sum_{k=1}^m \frac{\Gamma(k+s+\alpha)(2k-1-s+\alpha)}{\Gamma(k-s)}
	(\mbox{dist}(z,\partial\Omega))^{-\frac{2\alpha+2+sp}{p}} \norm{f}_{p,\Omega,\alpha)}
	\end{multline}
	where for $1/q+1/p=1,$ $p>1$
	\begin{equation}
	A_1(m,s,p,\alpha):=\left(
	\int_0^1 (1-t)^\alpha \abs{t^{-\frac{s}{2}}(1-t)^{-\alpha}(t^{k-1+s}(1-t)^{k-1+\alpha})^{k-1}}^q dt
	\right)^{\frac{1}{q}}
	\end{equation}
	\begin{equation}
	A_2(m,s,p,\alpha):=
	\left(
	\int_0^1 (1-t)^\alpha \abs{t^{-\frac{s}{2}}(1-t)^{-\alpha}(t^{k-1-s}(1-t)^{k-1+\alpha})^{k-1}}^q dt
	\right)^{\frac{1}{q}}
	\end{equation}
\end{theorem}

\chapter{A Cauchy-Riemann complex for certain polyanalytic functions} 
\section{Preliminaries and background}
A {\em chain complex}\index{Chain complex} is a (possibly infinite) sequence
$$
0\overset{d_0}{\longleftarrow} C_0 \overset{d_1}{\longleftarrow} C_1 \overset{d_{2}}{\longleftarrow}\cdots \overset{d_{q+1}}{\longleftarrow} C_{q+1}\longleftarrow \cdots
$$
where the $C_j$ (whose elements are called {\em $j$-chains}\index{Chains}) are abelian groups and the $d_j$, called {\em boundary operators}\index{Boundary operator} are homomorphisms of abelian groups such that
$d_j\circ d_{j+1}=0,$ for $j\geq 0.$ This implies that $B_j(C):=\mbox{Im} d_{j+1}\subseteq \mbox{Ker} d_j=:Z_j(C)$
thus we have a well-defined quotient group $H_j(C):=Z_j(C)/B_j(C)$ called the {\em $j$:th homology group} of the chain complex $C.$
Elements of $B_j$ are called $j$-boundaries and elements of $Z_j$ are called $j$-cycles.
If all homology groups vanish then the sequence is called {\em exact}\index{Exact sequence}.
Homology groups can be regarded as a measure of the failure of a chain complex to be exact.
A {\em cochain complex}\index{Cochain complex} is a sequence
$$
0\overset{d^{-1}}{\longrightarrow} C^0 \overset{d^0}{\longrightarrow} C^1\overset{d^{1}}{\longrightarrow}\cdots \overset{d^{q}}{\longrightarrow} C^{q+1}\longrightarrow \cdots
$$
where the $C^j$ (called {\em cochains}\index{Chains}) are abelian groups and the $d^j$ (called {\em coboundary maps}\index{Coboundary maps} are homomorphisms satisfying $d^{j+1}\circ d^j=0$, $j\geq 0.$
We have $B^j:=\mbox{Im}d^j\subseteq \mbox{Ker}d^{j+1}=:Z^{j+1}$ and the
quotient group (module) $H^j(C):=Z^j/B^j$ is called the $j$:th {\em cohomology group (module)}\index{Cohomology group} of the cochain complex $C$\index{Cochain complex}. 
Elements of $B^j$ are called $j$-coboundaries and elements of $Z^j$ are called $j$-cocycles. Exactness of the cochain complex is defined by vanishing of the
cohomology groups.
Cohomology groups can be regarded as a measure of the failure of a cochain complex to be exact.
Given two chain complexes $C$ and $C'$, a {\em (co)chain map}\index{Chain map (cocahin map)} $f : C \to C'$ is a family $f =(f_j)_{j\geq 0}$ of homomorphisms $f_j:C_j\to C_{'j}$
such that all the squares of the following diagram
commute i.e.\ $f_j\circ d_{j+1} =d'_{j+1} \circ f_{j+1},$ $j\geq 0,$
\begin{equation}
\begin{array}{ccccccccccc}
0 &\overset{d_{0}}{\longleftarrow} & C_0 & \overset{d_1}{\longleftarrow} & C_1 &\overset{d_{2}}{\longleftarrow} & \cdots & \overset{d_{q-1}}{\longleftarrow} & C_{q} & \longleftarrow & \cdots\\
 &                        & \downarrow_{f_0}  &    & \downarrow_{f_1}   &    &        &       &\downarrow_{f_q} &  & \\
0 &\overset{d'_{0}}{\longleftarrow} & C'_0 & \overset{d'_1}{\longleftarrow} & C'_1 &\overset{d'_{2}}{\longleftarrow} & \cdots & \overset{d'_{q-1}}{\longleftarrow} & C'_{q} & \longleftarrow & \cdots
\end{array}
\end{equation}
A chain map induces a homomorphism of homology groups $H_j(f)\colon H_j(C)\to H_j(C').$
Given three chain complexes $C,C',C''$ and two chain maps $f:C\to C'$, $g:C'\to C''$ it holds
$H_j(g\circ f)=H_j(g)\circ H_j(f)$, and $H_j(\mbox{Id}_C)=\mbox{Id}_{H_j(C)}$, $j\geq 0.$
That is to say the map $C\mapsto \{H_j(C)\}_{j\in \N}$ is a functor from the
the category of chain complexes and chain maps to the category of abelian groups and
group homomorphisms, and similarly cochain maps induce such a functor.
Let us look at a first example (see e.g.\ Narasimhan \& Nievergelt \cite{narasimhannieve}, p.28).
Let $\Omega\subset \C$ be an open subset and let $\mathcal{U}=\{U_i\}_{i\in \N}$ be an open cover of $\Omega.$
Let $J\subset I\times I$ be the set of pairs$(i,j)$ such that $U_i \cap U_j \neq \emptyset.$ For an open subset $V\subset\C$ denote by $\C(V)$
the set of locally constant complex valued functions on $V.$
Set $C^1:=\Pi_{(i,j)\in J} \C(U_i \cap U_j)$ (these will correspond to the $1$-cochains
of our example) and $Z^1:=\{ (c_{ij})_{(i,j)\in J}\in C^1 \colon c_{ij}+c_{jk}+c_{ki}=0 \mbox{ on }U_i\cap U_j\cap U_k\mbox{ whenever }
U_i\cap U_j\cap U_k\neq \emptyset\}$ (these will correspond to the $1$-cocycles
of our example).
Set further $C^0:=\Pi_{i\in I} \C(U_i)$ and
\begin{equation}
\delta\colon C^0\to Z^1
\end{equation}
defined as follows: if $(c_{i})_{i\in I}\in C^0, c_i\in \C(U_i)$ then $(\delta c)_{ij}=(c_i|U_i\cap U_j)-(c_j|U_i\cap U_j)=c_i-c_j$ on $U_i\cap U_j$
for $(i,j)\in J.$
Denote by $B^1$ the image of $\delta.$ Then $\delta$ is a $\C$-linear map between complex vector spaces (in particular additive groups) $Z^1, B^1,$
and we obtain the first homology group of our example, with complex coefficients, as the quotient vector space $H^1_{\mathcal{U}}:=Z^1/B^1.$
\begin{theorem}
If $H^1_{\mathcal{U}}=0$ and $f\in \mathscr{O}(\Omega)$ has a primitive, $F_i$, on each $U_i$ (in the sense that there exists a function $F_i$ on $U_i$ satisfying $\partial_{z} F_i=f|_{U_i}$)
then $f$ has a primitive on $\Omega.$
\end{theorem}
\begin{proof}
Set $c_{ij}:=F_i-F_j$ on $U_i\cap U_j$. Obviously, since $U_i\cap U_j$ is connected $c_{ij}$ is constant.
If $U_i\cap U_j\cap U_k\neq \emptyset$ then
\begin{equation}
c_{ij}+c_{jk}+c_{ki}=(F_i-F_j)+(F_j-F_k)+(F_k-F_i)=0
\end{equation}
thus $e:=\{ (c_{ij})_{(i,j)\in J}\}\in Z^1.$ The condition $H^1_{\mathcal{U}}=0$ then implies that there exists $c:=(c_i)_{i\in \N}\in C^0$ such that
$\delta c =e$ which implies
\begin{equation}
F_i-F_j=c_i-c_j\mbox{ on } U_i\cap U_j,\quad (i,j)\in J
\end{equation}
The function defined by $F|_{U_i}=F_i-c_i$ is a primitive of $f$ on $\Omega.$ This completes the proof.
\end{proof}

\begin{theorem}\label{tthhthm}
Suppose given a family of functions $f_{ij}\in \mathscr{O}(U_i\cap U_j)$ satisfying
$f_{ij}+f_{jk}=f_{ik}$ on $U_i\cap U_j\cap U_k$ for $i,j,k\in I.$
Then there is a family $\{f_j\}_{i\in I}$ where $f_i\in \mathscr{O}(U_i)$ such that
$f_i-f_j =f_{ij}$ on $U_i\cap U_j$ for all $i,j\in I.$
\end{theorem}
\begin{proof}
First we find a family $\{\phi_i\}_{i\in I}$ where $\phi_i\in C^\infty_c(U_i)$ (where the subindex $c$ means compact support) 
such that
$\phi_i-\phi_j =f_{ij}$ on $U_i\cap U_j$ for all $i,j\in I.$
Let $\{\eta_j\}_{j\in I}$ be a partition of unity relative to the open cover $\mathcal{U}$ and
\begin{equation}
\psi(z):=\left\{
\begin{array}{ll}
\eta_j(z)f_{ij}(z), & z\in U_i\cap U_j\\
0, & z\in U_i\setminus U_i\cap U_j\\
\end{array}
\right.
\end{equation}
Then $\psi\in C^\infty(U_i),$ and we write $\eta_j f_{ij}$ for this function.
Set $\phi_i:=\sum \eta_j f_{ij}$ on $U_i.$ Since the family $\{\mbox{supp}(\eta_j)\}$ is locally finite 
the sum has only finitely many nonzero terms in a neighborhood of any point of $U_i$ so $\phi_i\in C^\infty(U_i).$
Also $f_{ij}+f_{jk}=f_{ik}$ on $U_i\cap U_j\cap U_k$ so for $i=j=k$ we have $f_{ij}=-f_{ji}$ on $U_i\cap U_j.$
For $k,l\in I$ 
$\phi_k-\phi_l=\sum_j \eta_j(\phi_{kj}-\phi_{lj}).$ Using $f_{kj}-f_{lj}=f_{kj}+f_{jl}$ we get
$\phi_k -\phi_l= \sum_j \eta_j\phi_{kl}=\phi_{kl} \sum_j \eta_j=\phi_{kl}.$
Now taking $\partial_{\bar{z}}$ on both sides yields
$\partial_{\bar{z}}\phi_k = \partial_{\bar{z}}\phi_k.$ on $U_k\cap U_l.$
Hence there is a function $\phi\in C^\infty(\Omega)$ such that $\phi|_{U_j}=\partial_{\bar{z}} \phi_j$ for all $j\in I.$
By the solution to the $\overline{\partial}$-problem in one dimension we can find $u\in C^\infty(\Omega)$ such that $\partial_{\bar{z}} u=\phi$ on $\Omega$ and we 
set $f_i:=\phi_i-u$ on $U_i.$
The each $f_i$ is holomorphic on $U_i$ and $f_i-f_j=(\phi_i-u)-(\phi_j-u)=\phi_i-\phi_j=f_{ij}$ on $U_i\cap U_j.$
This completes the proof.
\end{proof}

\begin{theorem}
Let $\mbox{i}_\Omega:\C\to \mathscr{O}(\Omega)$ be the map sending any complex number $c$ to the constant function $z\mapsto c$ on $\Omega.$
Let $d_\Omega:\mathscr{O}(\Omega)\to \mathscr{O}(\Omega),$ denote the derivative $f\mapsto \partial_z f.$ 
Then the following sequence is exact
\begin{equation}
0{\longrightarrow} \C \overset{\mbox{i}_\Omega}{\longrightarrow} \mathscr{O}(\Omega) \overset{d_{\Omega}}{\longrightarrow}\mathscr{O}(\Omega) \overset{\delta}{\longrightarrow} H^1_\mathcal{U}\longrightarrow 0
\end{equation}
\end{theorem}
\begin{proof}
Obviously $\mbox{i}_\Omega c=0\Rightarrow c=0$ so the map is injective. Also $\partial_z f=0$ implies that $f$ is constant so $\mbox{Im} i_\Omega =\mbox{Ker} d_\Omega.$
If $d_\Omega (F)=f$ then $\delta(f)$ can be defined by taking the family $(F_i)$ of primitives as $F_i=F|_{U_i},$ so $\delta(f)$ is the class of $F_i-F_j=0$ on $U_i\cap U_j$ which implies $\delta(f)=0.$
Thus $\mbox{Im} d_\Omega \subset \mbox{Ker} \delta.$
Conversely if $f\in \mbox{Ker}\delta$ then letting $F_i\in \mathscr{O}(U_i)$ be a primitive of $F$ on $U_i$ and $c_{ij}=F_i-F_j$ on $U_i\cap U_j$.
Since $\delta(f)=0$ there exists $(c_i)_{i\in I}$, for locally constant $c_i$ (thus constants) on $U_i$ such that $c_i-c_j=c_{ij}$ on $U_i\cap U_j$.
Thus $F_i-F_j=c_i-c_j$  on $U_i\cap U_j,$ which implies that there is a function $F$ on $\Omega$ with $F|_{U_i}=F_i-c_i$, which is holomorphic and satisfies
$d_\Omega(F)|_{U_i}=\partial_z(F_i-c_i)=f|_{U_i}$, so $f\in \mbox{Im}d_\Omega.$ Thus $\mbox{Ker} \delta \subset \mbox{Im} d_\Omega $ and thus
$\mbox{Ker} \delta = \mbox{Im} d_\Omega.$
Finally, let $e\in H^1_\mathcal{U}$ and $\{c_{ij}\}\in Z^1$ such that $\delta \{c_{ij}\} =e.$
Since the $c_{ij}$ are locally constant they are holomorphic and by Theorem \label{tthhthm} there exists a family $\{F_i\}_{i\in I}$ with $F_i\in \mathscr{O}(U_i)$ such that $F_i-F_j=c_{ij}$
on $U_i\cap U_j.$
Since $\partial_z F_i -\partial_z F_j=\partial_z c_{ij}=0$ on $U_i\cap U_j$ there is a holomorphic $f$ on $\Omega$ such that $f|_{U-i}=\partial_z F_i.$
We can define $\delta(f)$ by choosing a primitive on $U_i$ to be $F_i$. Thus $\delta(f)$ becomes the class in $H^1_\mathcal{U}$
of $\{ (F_i-F_j)|_{U_i\cap U_j}\},$ i.e.\ of $\{c_{ij}\}.$
Thus $\delta(f)=e.$
This completes the proof.
\end{proof}

We shall be interested in the case of cochain complexes where the cohomology groups $H^j(\Omega,\mathcal{F})$
take values in a sheaf $\mathcal{F}$ on the space $\Omega$. This means that the
modules $\mathcal{F}(U)$, where $U$ is a domain in $\Omega$, may vary with the
$U$. This has been extensively studied for the case of sheaves of germs of holomorphic functions. Let us first recall the definitions.
\begin{definition}
A connected complex manifold $X$ of dimension $1$ having a countable base for its open sets is called a {\em Riemann surface}.\index{Riemann surface}
$X$ is said to be equipped with a {\em projective structure}\index{Projective structure}, if its atlas is given 
by holomorphic charts $\{(\phi_i,U_i)\}_{i\in I}$ such that for any pair $i,j$ with $U_i\cap U_j\neq 0$ the coordinate change
is given by a linear fractional, i.e.\ 
$\phi_j\circ\phi_i^{-1}=\frac{az +b}{cz +d}$ with $ad-bc=1,$ see e.g.\ Gustafsson \& Peetre \cite{peetre}.
\end{definition}
A theorem of Rad\'o (for a proof, see e.g.\ Foster \cite{foster}, 23.1, p.187) states that any connected complex manifold $X$ of dimension $1$ has a countable base.
A typical example is given by the complex projective line\index{Complex projective line $\mathbb{P}^1$} $\mathbb{P}^1$ defined as follows.
Let $X=\C\cup \{\infty\}.$ Define the open sets $U_1=\C$ and $U_2 =\{ z\in \C \colon z\neq 0\}\cup \{\infty\}.$ Then
$X=U_1\cup U_2$. Set $\phi_1:U_1\to \C$ be the identity, and $\phi_2:U_2 \to \C$ the map $\phi_2(\infty)=0,$
$\phi_2(z)=1/z$ if $z\in U_2\setminus \{\infty\}.$ Then $\phi_1,\phi_2$ are complex analytic 
charts in the sense that they are homeomorphisms satisfying $\phi_1\circ \phi_2^{-1}:\phi_2(U_1 \cap U_2)\to \phi_1(U_1\cap U_2)$
is holomorphic. Indeed, it is given by $z\mapsto 1/z$ from $\C\setminus \{0\}$ onto itself. 
We have $\hat{C}\simeq \mathbb{P}^1$ and the space of lines in $\C^2$
i.e.\ $\mathbb{P}^1\simeq (C^2\setminus \{0\})/\C^*$ via the isomorphism $[Z_1,Z_2]\mapsto z=Z_1/Z_2.$
Inverse charts $\C \to \mathbb{P}^1$ are given by $\psi(t)=
[ct + d, at + b],$ for complex constants $a,b,c,d,$ where $a/c$ does not belong to the image. The image omits a=c. The automorphism group of $\hat{C}$ is in fact
the quotient of $SL_2(\C)$ by $\pm I$. This means every automorphism lifts to a linear map on $\C^2$ with determinant
$1$. 
\begin{definition}
A 2 dimensional manifold (surface) is called orientable\index{Orientable surface} iff 
on the overlap of any pair of charts the determinant of the differential of the transition map has 
constant sign.
The {\em genus}\index{Genus of an ortientable surface} of a connected, orientable surface, 
$S,$ is a topological invariant
defined as the maximum number of cuttings by non-intersecting closed simple curves 
that do not disconnected the manifold. For an orientable surface the {\em Euler characteristic}\index{Euler characteristic}
can be defined as $\chi=2-2g.$
Alternatively, the Euler characteristic can be defined via
$\int_S K =2\pi\chi$ where $K$ is the Gaussian curvature of the surface.
\end{definition}
Every Riemann surface is orientable since it is a complex manifold, so the transition maps can be regarded as
maps $U\to V$ for open sets $U,V\subset \R^2$
where the determinant of the Jacobian is given by the absolute value of the
complex derivative of the transition map.
Any compact Riemann surface with genus $g>1$ has finitely many
conformal automorphism groups.
For deeper knowledge on complex analysis on compact Riemann surfaces see e.g.\ the book of Gunning \cite{gunningriemann}.
\begin{definition}
Let $X$ and $S$ be connected manifolds. Recall that two paths $c_i:[0,1]\to S,$ $i=1,2,$ are homotopic if
there exists a continuous function $H: [0,1]\times[0,1]\to S,$
such that
$H(t,0)=c_1(t),$ $t\in [0,1],$
$H(t,1)=c_2(t),$ $H(0,s)=c_1(0)=c_2(0)$, $s\in [0,1],$
$H(1,s)=c_1(1)=c_2(1).$ This clearly gives an equivalence relation
and the equivalence class of a path $c$ is denoted $[c].$
The equivalence classes, with respect to homotopy, of loops with fixed base point $p\in S$ form a group, $\pi_1(S,p),$ called the {\em fundamental group} of $S$ with base point $p.$ 
The group action if given by $[c_1]\cdot[c_2]:=[c_1\cdot c_2]$
for representative loops $c_1,c_2,$ where $c_1\cdot c_2$ is the path that first traverses $c_1$ then $c_2.$
A continuous map $\pi:X\to S$ is called a {\em covering map}\index{Covering map} if for each $p\in M,$
there is a neighborhood $U$ of $p$ in $S$ such that each connected component of $\pi^{-1}(U)$ is mapped homeomorphically onto $U.$
If $p\in M$ and $H$ is a subgroup of $\pi_1(S,p)$, then it is well-known that there is a covering $\pi:X\to S$ (called the {\em universal covering map})
with the property that for every $x\in X$ with $\pi(x)=p,$ we have $\pi_*(\pi_1(X,x))=H.$ Choosing $H=\{1\}$ we obtain a simply connected manifold $\tilde{S}$ called the {\em universal covering} of $S$ (since it has a trivial fundamental group, it is also its own universal cover as well), such that $\pi:\tilde{S}\to S.$
\end{definition}
Let $X,Y$ be topological spaces. 
 A {\em lift} of a map $f:Y\to X$
is a map $\tilde{f}_Y\to X$ such that $\pi\circ \tilde{f}=f.$ Let $p_0\in X$ and $\tilde{p_0}\in \pi^{-1}(p_0).$
It is well-known that the homomorphism between fundamental groups $\pi_*:\pi_1(\tilde{X},\tilde{p}_0)\to \pi_1(X,p_0)$
is injective. Different choices of base point will simply render subgroups of $\pi_1(X,p_0)$ that are conjugate to each other.
(recall that two subgroups $G_1$ and $G_2$ of a group $G$ are conjugate if there is a $g\in G$ such that
$gG_1 g^{-1}=G_2$).
$X$ has a universal cover if it is connected, locally path-connected and if every point 
in $X$ has a neighborhood with the property that every loop in it is nullhomotopic.

\begin{definition}
If $S$ is an orientable compact real surface of genus $g$. A projective structure on $S$ induces a complex structure on $S$ just by pulling
back the complex structure of $\mathbb{P}^1$ by the projective charts. Denote by $\Sigma$ the corresponding compact Riemann surface  so obtained. Denote by $\pi:C\to \Sigma$, the universal covering map. By the Riemann mapping theorem we can assume that $C$ is either $\mathbb{P}^1,$ $\C$ or the unit disc $\{\abs{z}<1\}$ respectively for $g=0,1$ and $> 2$ respectively. 
There is thus an inherited representation of the fundamental group $\pi_1(\Sigma)\to \mbox{Aut}(C)$ whose image $W$ is a subgroup of
the set of linear fractionals (i.e.\ M\"obius transformations or members of $PGL_2(\C)$). 
The atlas defined on $\Sigma$ by all local determinations of $\pi^{-1}:\Sigma\to \mathbb{P}^1$ defines a projective structure on $\Sigma$, compatible with the complex structure. Note that any two determinations of $\pi^{-1}$ differ by left composition with an element of $W.$
In particular, in the case of genus greater than $1$, it makes sense to speak of a projective structure on a Riemann surface, {\em induced
by a covering map of the unit disc.}
\end{definition}
\begin{definition}\label{sheafdefii}
Given a topological space $X$ and a class, $C$, of groups (or modules), a {\em presheaf}\index{Presheaf} on $X$ with values in $C$ consists
of an assignment of a group $\mathcal{F}(U)$ in $C$ to every open subset $U\subset X$ and of a map
$\rho_V^U : \mathcal{F}(U) \to \mathcal{F}(V)$ of the class of structures in $C$ to every inclusion 
of open
subsets $V \subseteq U \subseteq X$, such that
$\rho_W^U = \rho_W^V \circ \rho_V^U$, 
$\mathcal{F}(\mbox{Id}_U) = \mbox{Id}_{\mathcal{F}(U)}$
for any two consecutive inclusions $W\subseteq V \subseteq U$.
A {\em sheaf}\index{Sheaf} on $X$ with values in $C$ is a
presheaf $\mathcal{F}$ on $X$ such that for any open subset $U$ of $X$ and every open cover $\{U_i\}_{i\in I}$ of $U$
we have:
(i) For every family $\{f_i\}_{i\in I}$ with $f_i \in \mathcal{F}(U_i)$, satisfying
\begin{equation}
\rho^{U_i}_{U_i\cap U_j}(f_i)= \rho^{U_j}_{U_i\cap U_j}(f_i),\quad i,j\in I
\end{equation}
then there is some $f\in\mathcal{F}(U)$ such that $\rho^{U}_{U_i}(f)=f_i$ for all $i\in I.$
(ii) For any two $f,g\in \mathcal{F}(U)$, 
\begin{equation}
\rho^{U}_{U_i}(f)= \rho^{U}_{U_i}(g),\mbox{ for all }i\in I
\end{equation}
then $f = g$.
\end{definition}
In general, a presheaf is not a sheaf (consider for example the constant presheaf).
\begin{definition}
If $\mathcal{A}$ is a sheaf of rings (or vector spaces) we say that a sheaf $\mathcal{B}$, is a 
{\em sheaf of $\mathcal{A}$-modules} if in the definition of $\mathcal{B}$ each of the assigned groups $\mathcal{F}(U)$
to open sets $U$, is a $\mathcal{A}(U)$-module.
\end{definition}  

\begin{definition}
Let $X$ be a Riemann surface, and let $p_0\in X$. Consider the pairs $(U, f)$, where $U$ is an open neighbourhood of $p_0$ and $f\in \mathscr{O}(U)$. 
Two pairs $(U, f)$ and $(V, g)$ are called equivalent, and define the same germ of holomorphic function at $p_0$, 
if there exists an open neighbourhood $W$ of $p_0$, $W \subset U\cap V$, such that $f|_W = g|_W$. An equivalence class $[f]_{p_0}$, is called a {\em germ of the holomorphic function $f$ at $p_0$}\index{Germ of holomorphic functions}. 
The value of $[f]_{p_0}$ at $p_0$ is defined as $[f]_{p_0}(p_0)=f(p_0)$ (note that $g(p_0)=f(p_0)$ for any $g\in [f]_{p_0}$).
The set of all germs at $p_0$ is denoted $\mathscr{O}_{p_0}.$ It is clearly a ring (in fact can be turned into a $\C$-algebra).
Define the disjoint union $\mathscr{O}_X:=\bigcup_{p\in X} \mathscr{O}_p$ and the map $\pi:\mathscr{O}_X\to X,$ $\pi(f)=p_0.$
Let $(U,f)$ be a pair defining $[f]_{p_0}$. Set $N(U,f)=\{ [f]_p \colon p\in U\}.$
Define a topology on $\mathscr{O}_X$ by the condition that the sets $\{N(U,f)\}$ form a fundamental system of neighborhoods of $[f]_{p_0}$ when $(U,f)$ runs over
all pairs defining $[f]_{p_0}$. With this topology $\mathscr{O}_X$ is Haussdorff and $\pi$ is a local homeomorphism (see e.g.\ Narasimhan \cite{narasimhancompact}, p.12). $\pi^{-1}(\{p_0\})$
is called the {\em stalk}\index{Stalk} (or fiber) over $p_0$.
\end{definition}
We can thus consider the sheaf of germs of holomorphic functions\index{Sheaf of germs of holomorphic functions} and similarly one can define the sheaf of germs of meromorphic functions.
\begin{definition}
Let $\mathcal{F},\mathcal{G},\mathcal{H}$ be sheaves on a topological space $X$. 
A morphism $\alpha :\mathcal{F}\to \mathcal{G}$ of sheaves on $X$ consists of a morphism $\alpha_U:\mathcal{F}(U)\to \mathcal{G}(U)$ for all open sets $U\subset X$ such that
$$
\begin{array}{lll}
\mathcal{F}(U) &\overset{\alpha_U}{\longrightarrow} &\mathcal{G}(U)\\
\downarrow_{\rho_{V}^U} & & \downarrow_{\rho_{V}^U}\\
\mathcal{F}(V) &\overset{\alpha_V}{\longrightarrow} &\mathcal{G}(V)\\
\end{array}
$$
We say that 
\begin{equation}
0\longrightarrow \mathcal{F} \overset{\alpha}{\longrightarrow} \mathcal{G} \overset{\beta}{\longrightarrow} \mathcal{H} {\longrightarrow} 0
\end{equation}
is exact if for all $U$, the sequence
\begin{equation}
0\longrightarrow \mathcal{F}(U) \overset{\alpha}{\longrightarrow} \mathcal{G}(U) \overset{\beta}{\longrightarrow} \mathcal{H}(U) {\longrightarrow} 0
\end{equation}
is exact and if $h\in \mathcal{H}(U)$ and $z\in U$, there exists an open neighborhood $V\subset U$ of $z$, and $g\in \mathcal{G}$ such that $\beta_V(g)=h|_V.$
The {\em kernel sheaf}\index{Kernel sheaf} of $\alpha$ is defined as $\mbox{Ker}\alpha(U)=\mbox{Ker}(\alpha_U :\mathcal{F}\to \mathcal{G}(U)).$
If 
$$
0\longrightarrow \mathcal{F} \overset{\alpha}{\longrightarrow} \mathcal{G} {\longrightarrow} \mathcal{H} {\longrightarrow} 0
$$
is a short exact sequence of sheaves then we define the cokernel $\mbox{coker}\alpha :=\mathcal{H}.$
This does not imply that $\mathcal{G}(U) {\longrightarrow} \mathcal{H}(U)$ is surjective. In general, there
is a sheaf which is the cokernel of any sheaf morphism but to obtain it explicitly it must be defined via a procedure called sheafification which we shall not do here.
\end{definition}
\begin{definition}
A complex of sheaves
$$
\mathcal{F}_1 \overset{d_1}\longrightarrow \mathcal{F}_2 \overset{d_2}{\longrightarrow} \cdots \overset{d_{q-1}}{\longrightarrow}\mathcal{F}_q \overset{d_q}{\longrightarrow} \cdots
$$
is {\em exact} if
$$
0 \longrightarrow \mbox{Ker}d_j {\longrightarrow} \mathcal{F}_j {\longrightarrow} 
\mbox{Ker} d_{j+1} {\longrightarrow}0
$$
is a short exact sequence for all $j.$
\end{definition}

\begin{definition}\index{Resolution of a sheaf}
A {\em resolution} of a sheaf $\mathcal{F}$ is a complex of sheaves $\mathcal{F}^j$, $j\in \Z_+,$ such that 
$$
0\longrightarrow \mathcal{F} \longrightarrow \mathcal{F}^1 \longrightarrow\mathcal{F}^2 \longrightarrow \cdots
$$
is exact. 
\end{definition}
Recall that the 
set, $\mathcal{A}^{p,q},$ of complex $(p,q)$-forms and the maps 
$\partial :\mathcal{A}^{p,q}\to \mathcal{A}^{p+1,q}$ and $\overline{\partial} :\mathcal{A}^{p,q}\to \mathcal{A}^{p+1,q}$,
where for an open set $U$, one can start from the exterior derivative and define its $\C$-linear extension
$d$ to complex $p$-forms and then denote by $\overline{\partial}$ ($\partial$) the composition of $d$ with projection
to $\mathcal{A}^{p,q+1}$ ($\mathcal{A}^{p+1,q}$). On $\mathcal{A}$ these are just the operators $\partial=\partial_{z}$ and 
$\overline{\partial}=\partial_{\bar{z}}$.
If $\alpha=fdz_I\wedge d\bar{z}_J \in \mathcal{A}^{p,q}$ then
$d\alpha=\partial\alpha +\overline{\partial}\alpha,$
$\partial\alpha:=\sum_j \partial_{z_j} f dz_j\wedge dz_I\wedge d\bar{z}_J$
and $\overline{\partial}\alpha:=\sum_j \partial_{\bar{z}_j} f d\bar{z}_j\wedge dz_I\wedge d\bar{z}_J$.
Due to the alternating property of wedge products we obviously have $\partial^2=0$, $\overline{\partial}^2=0$.
Also the operators $\overline{\partial}$ and $\partial$ commute with restriction since $d$ does, this can be used to realized that they induce
maps of sheaves. 
For a complex manifold $X$ the set, $\mathcal{A}^{p,q},$ of complex $(p,q)$-forms together with the $\overline{\partial}$-operator
yields the following complex of sheaves
\begin{equation}
\cdots \overset{\overline{\partial}_{q-2}}{\longrightarrow} \mathcal{A}^{p,q-1} \overset{\overline{\partial}_{q-1}}{\longrightarrow} \mathcal{A}^{p,q} 
\overset{\overline{\partial}_{q}}{\longrightarrow} \mathcal{A}^{p,q+1} \overset{\overline{\partial}_{q+1}}{\longrightarrow} \cdots
\end{equation}

\begin{definition}
A differential form $\phi\in \mathcal{A}^{p,q}$ is called {\em $\overline{\partial}$-closed}\index{$\overline{\partial}$-closed form} if 
$\overline{\partial}\phi =0$ and it is called {\em $\overline{\partial}$-exact}\index{$\overline{\partial}$-exact form} if 
there exists $\psi\in\mathcal{A}^{p,q-1}$ such that $\overline{\partial}\psi =\phi$. 
\end{definition}

\begin{definition}
The $(p,q)$ Dolbeault cohomology group\index{Dolbeault cohomology group} is defined as
\begin{equation}
H_{\overline{\partial}}^{p,q}(X):=\frac{\mbox{Ker}\overline{\partial}\colon \mathcal{A}^{p,q}\to \mathcal{A}^{p,q+1}}{\mbox{Im}\overline{\partial}\colon \mathcal{A}^{p,q-1}\to \mathcal{A}^{p,q} }
\end{equation}
where $\mbox{Im}\overline{\partial}$ denotes the image of the operator.
By convention negative groups are zero and $H_{\overline{\partial}}^{p,q}(X)$ is the set of holomorphic functions.
\end{definition}

\begin{theorem}\label{dolbeaultlemma}
	Let $K\subset\C$ be compact, $U$ an open neighborhood of $K$ and $g\in C^\infty(U).$
	Then there exists $f\in C^\infty(V)$ for an open neighborhood $V\subset U$ of $K$
	such that
	$\partial_{\bar{z}} f=g$ on $V$. Further if $g$ is $C^\infty$ or holomorphic in some additional  then $f$ 
	can be chosen to have the same properties. 
\end{theorem}
\begin{proof}
	We have
	\begin{equation}
	d\left(f(\zeta)\frac{d\zeta}{\zeta -z}\right)=\partial_{\bar{\zeta}}
	\left(\frac{f(\zeta)}{\zeta -z}\right)d\bar{\zeta}\wedge d\zeta=
	\partial_{\bar{\zeta}} f(\zeta)\frac{d\bar{\zeta}\wedge d\zeta}{\zeta -z}
	\end{equation}
	Let $K\subset W\subset U$ where $W$ is chosen such that it is bounded by a rectifiable closed curve.
	Recall that Stokes formula for a smooth complex function $h$ on $\overline{W}$ gives
	\begin{equation}
	\int_{\partial W} h(z)dz=\int_W dh\wedge dz =
	\int_W \partial_{\bar{z}} h dz\wedge d\bar{z} 
	\end{equation}
	Applying this to $W\setminus B_\epsilon,$ $B_\epsilon:= \{\zeta\colon \abs{\zeta-z}<\epsilon\}$
	for $\epsilon$ sufficiently small such that $\overline{B}_\epsilon\subset W,$
	we get
	\begin{equation}\label{ekvforhorm}
	\int_{\partial W\cup B_\epsilon} \frac{\partial_{\bar{\zeta}}h(\zeta)}{\zeta-z}d\bar{\zeta}\wedge d\zeta=
	\int_{\partial W} \frac{h(\zeta)}{\zeta-z} d\zeta -
	\int_{\partial B_\epsilon} \frac{h(\zeta)}{\zeta-z} d\zeta 
	\end{equation}
	Since $(\zeta -z)^{-1}d\bar{\zeta}\wedge d\zeta$ is a bounded measure in the plane, the integral
	on the right hand side tends to the integral over $W$ as $\epsilon\to 0.$
	Passing to polar coordinates, $z=r\exp(i\theta)$ 
	we have
	\begin{equation}
	\lim_{\epsilon\to 0}\int_{\partial B_\epsilon} \frac{h(\zeta)d\zeta}{\zeta -z}=
	\lim_{\epsilon\to 0}\int_{0}^{2\pi} h(z+\epsilon\exp(i\theta))id\theta = 2\pi i h(z)
	\end{equation}
	thus
	\begin{equation}
	h(z)=\frac{1}{2\pi i}\left(
	\int_{\partial W} \frac{h(\zeta)}{\zeta-z}d\zeta+\int_{W} \frac{\partial_{\bar{\zeta}}h(\zeta)}{\zeta-z} d\bar{\zeta}\wedge d\zeta
	\right)
	\end{equation}
	which when $\partial_{\bar{\zeta}}h\equiv 0$ is the Cauchy formula.
	
	Let $\phi\in C^\infty_c(\C)$ such that $\phi\equiv 1$ on an open neighborhood $V$ of $K$ (which we choose 
	such that $V\subset W$).
	Then $g\phi\in C^\infty$ and
	\begin{equation}
	\int_{U} \frac{g(\zeta)}{\zeta-z}d\zeta=\int_{U} \frac{g(\zeta)\phi(\zeta)}{\zeta-z}d\bar{\zeta}\wedge d\zeta
	+\int_{U} \frac{(1-\phi(\zeta))g(\zeta)}{\zeta-z} d\bar{\zeta}\wedge d\zeta
	\end{equation}
	where the last integrand is not supported in $U$. 
	Hence we only need to prove the result for $g$ with compact support in $\C.$
	For such $g$
	the integral 
	\begin{equation}
	f(z)=\frac{1}{2\pi i}\int_{\C} \frac{g(\zeta)}{\zeta-z}d\bar{\zeta}\wedge d\zeta
	\end{equation}
	is well-defined and $C^\infty$-smooth on $\C$ and $f$ is 
	$C^\infty$ or holomorphic in some additional parameters whenever $g$ 
	is. 
	By change of variables $\zeta\mapsto \zeta +z$
	\begin{multline}
	\partial_{\bar{z}} f(z)=
	\frac{1}{2\pi i}\partial_{\bar{z}}\int_{\C} \frac{\partial_{\bar{\zeta}}g(\zeta+z)}{\zeta} d\bar{\zeta}\wedge d\zeta=\\
	\frac{1}{2\pi i}\int_{\C} \frac{g(\zeta+z)}{\zeta}d\bar{\zeta}\wedge d\zeta=
	\frac{1}{2\pi i}\int_{\C} \frac{g(\zeta)}{\zeta-z}d\bar{\zeta}\wedge d\zeta
	\end{multline}
	By the Cauchy integral formula the right hand side is $g(z)$. This completes the proof.
\end{proof}

\begin{theorem}\label{dolbeaultlemma2}
Let $U\subset\Cn$ be an open subset and let $B=\{\abs{z-p}<r\}$ be a bounded polydisc $\overline{B}\subset U.$ 
For any given $\overline{\partial}$-closed $\alpha\in \mathcal{A}^{p,q}(\tilde{B})$ (i.e.\ $\overline{\partial}\alpha =0$) 
there is a $\beta\in \mathcal{A}^{p,q-1}(\tilde{B})$
satisfying $\overline{\partial}\beta =\alpha$ on $\tilde{B}$, where $\tilde{B}$ is an open neighborhood of $\overline{B}$.
\end{theorem}
\begin{proof}
We use induction in the highest index $m$ such that $d\bar{z}_m$ appears in the representation of $\alpha.$
If $m=0$ then $\alpha=0$
Then $\alpha=d\bar{z}_m\wedge \beta' \wedge  +\alpha'$, where $\alpha'$ and $\beta$ involve 
only the conjugate differentials $d\bar{z}_1,\ldots, d\bar{z}_{m-1}.$ 
Since
$0=\overline{\partial}\alpha=-d\bar{z}_m\wedge \overline{\partial}\beta' +\overline{\partial}\alpha'$
both
$\beta$ and $\alpha'$ are holomorphic in $z_{m+1},\ldots,z_{n}$ since each coefficient is annihilated by $\partial_{z_{m+1}},\ldots, \partial_{z_{n}}.$
Now each coefficient, say $g$, of $\alpha$, is $C^\infty$-smooth with respect to $z_m$ in an open neighborhood of 
$\{\abs{z_m-p_m}\leq r_m\}$ and also 
$C^\infty$-smooth with respect to $z_1,\ldots,z_{m-1}$ 
and holomorphic with respect to $z_{m+1},\ldots,z_{n}$ in $\{\abs{z_{j}-p_{j}}<r_j, j=m+1,\ldots,n\}$.
By Theorem \ref{dolbeaultlemma} there exists a $C^\infty$-smooth $f$ in an open neighborhood of $\overline{B}$
such that $\overline{\partial} f=g$ and such that $f$ is holomorphic with respect to $z_{m+1},\ldots,z_n.$
Let $\psi$ be the differential form obtained by replacing each coefficient $g$ of $\alpha$ by the corresponding coefficient $f$ above.
Then $\overline{\partial}\psi =\delta +d\bar{z}_m\wedge \alpha,$ where $\delta$ is a differential form involving only
the conjugate differentials $d\bar{z}_1,\ldots, d\bar{z}_{m-1}.$
Set $\kappa=\alpha-\overline{\partial}\psi=\alpha'-\delta.$ Then
$\kappa$ involves only the conjugate differentials $d\bar{z}_1,\ldots, d\bar{z}_{m-1}$ and
$\overline{\partial}\kappa=\overline{\partial}f-\overline{\partial}\overline{\partial}\psi=0.$
By the induction hypothesis there exists a $C^\infty$-smooth differential form $\eta$ in an open neighborhood
of $\overline{B}$ such that $\overline{\partial} \eta =\kappa.$ Then $\beta:=\eta + \psi\in \mathcal{A}^{p,q}$ and satisfies
$\overline{\partial} \beta =\alpha.$ This completes the proof.
\end{proof}

\begin{theorem}[Poincar\'e lemma for polydiscs]\label{dolbeaultlemma3}
Let $B$ be a (possibly undbounded) polydisc. 
$H_{\overline{\partial}}^{p,q}(B)=0,$ for all $q\in Z_+.$
\end{theorem}
\begin{proof}
It suffices to show that for any given $\overline{\partial}$-closed $\alpha\in \mathcal{A}^{p,q}(B)$ (i.e.\ $\overline{\partial}\alpha =0$) 
there is a $\psi\in \mathcal{A}^{p,q-1}(B)$
satisfying $\overline{\partial}\psi =\alpha$ on $B$.
Let $\{D_j\}_{j\in \Z_+}$ be a sequence of bounded polydiscs in $\Cn$ such that
$D_{j}\subseteq D_{j+1},$ $\bigcup_j D_j=D,$ and let $\alpha\in \mathcal{A}^{p,q},$ $q\in \Z_+.$
\\
\\
\textit{Case: $q>1$.}
We construct by induction on $j$ a sequence of $(p,q-1)$-forms $\psi_j$ such that
$\psi_j$ is $C^\infty$-smooth on an open neighborhood of $\overline{D}_j$, $\overline{\partial}\psi_j=\alpha$ in that neighborhood and $\psi_{j+1}|_{D_j}=\psi_j.$
Then $\psi$ defined by $\psi|_{D_j}=\psi_j$ yields the wanted result.
For $j=1$, the existence of $\psi_1$ is given by Theorem \ref{dolbeaultlemma2}. Assume the case some fixed $j>1$ holds true. We show that
it holds true also for the case $j+1.$ Suppose $\psi_1,\ldots,\psi_{j}$ are given satisfying the wanted conditions. By
Theorem \ref{dolbeaultlemma2} there exists a $(p,q-1)$-form $\psi_{j+1}'$ which is $C^\infty$-smooth on an open neighborhood of $\overline{D}_{j+1}$ and satisfies $\overline{\partial}\psi_{j+1}=\alpha$ in that neighborhood. Then $\psi_{j+1}-\psi_j$ is $\overline{\partial}$-closed
$(p,q-1)$-form in an open neighborhood of $\overline{D}_j,$ and because $q-1>0$, Theorem \ref{dolbeaultlemma2}
yields a $(p,q-2)$-form $\eta_j$ in an open neighborhood $U_j$ of $\overline{D}_j$ such that
$\overline{\partial} \eta_j =\psi_{j+1}' -\psi_j$ in $U_j.$
Let $\rho_j\in C^\infty(\Cn)$ such that $\rho\equiv 1$ on an open neighborhood $\overline{D}_j$
and $\rho_J\equiv 0$ on an open neighborhood of $D\setminus U_j.$ 
The product $\rho_j\eta_j$ extends to a $C^\infty$-smooth differential form on $D$ that vanishes on $D\setminus U_j,$
we keep the same notation for the extension.
Then $\psi_{j+1}=\psi_{j+1}'-\overline{\partial}(\rho_j\eta_j)$ in an open neighborhood of $\overline{D}_{j+1}$
has the wanted properties. This completes the induction and proves the cases $q>1.$
\\
\\
\textit{Case: $q=1$.}
We construct by induction on $j$ a sequence of $(p,0)$-forms $\psi_j$ such that
$\psi_j$ is $C^\infty$-smooth on an open neighborhood of $\overline{D}_j$, 
$\overline{\partial}\psi_j=\alpha$ in that neighborhood and $\psi_{j+1}-\psi_j=\theta_j$ is 
a $\overline{\partial}$-closed $(p,0)$-form in an open neighborhood of $\overline{D}_j$ with all coefficients of $\theta_j$
being bounded by $2^{-j}$ on $\overline{D}_j.$ Then $\psi(z)=\lim_j\psi_j(z)$ is a well-defined $(p,0)$-form
in $D$ and since for any $j$ the limit can be written as $\psi|_{D_j}=\psi_j|_{D_j} +\sum_{k=j}^\infty \theta_j|_{D_j}$,
where $\sum_{k=j}^\infty \theta_j|_{D_j}$ is a uniformly convergent series of $\overline{\partial}$-closed $(p,0)$-forms, thus $\psi|_{D_j}-\psi_j|_{D_j}$
is a $\overline{\partial}$-closed  $(p,0)$-form and therefore $\psi$ has the wanted properties. We use induction in $j.$
The existence of $\psi_1$ is given by Theorem \ref{dolbeaultlemma2}. Assume that $\psi_1,\ldots,\psi_j$ have been
constructed as wanted. By Theorem \ref{dolbeaultlemma2} there exists a $(p,0)$-form $\psi_{j+1}'$ which is
$C^\infty$-smooth in an open neighborhood of $\overline{D}_{j+1}$ such that $\overline{\partial}\psi_{j+1}'=\alpha$
in that neighborhood.
Then $\theta'_j:=\psi_{j+1}'-\psi_j$
is a $\overline{\partial}$-closed $(p,0)$-form on an open neighborhood of $\overline{D}_j.$
The power series expansion of the coefficients of $\theta'_j$ near the common center of all polydiscs are absolutely and uniformly 
convergent in $\overline{D}_j$, thus $\theta'_j=\theta_j''+\theta_j$, where the coefficients of $\theta''_j$ are polynomials
and the coefficients of $\theta_j$ are bounded by $2^{-j}$ on $\overline{D}_j.$ The form $\psi_{j+1}:=\psi_{j+1}'-\theta''_j$ defined in an
open neighborhood of $\overline{D}_{j+1},$ then satisfies the wanted conditions. This completes the induction.
This completes the proof.
\end{proof}

Note that in the case of $C^1$ and a form $g\in \mathcal{A}^{0,1}(B)$ for a polydisc $B$,
Theorem \ref{dolbeaultlemma2} with $q=1,$ simply yields the existence of a $C^\infty$-smooth {\em function}
$f\in \mathcal{A}^{0,0}(B)$ such that
$\partial_{\bar{z}} f=g$ on $B.$

\begin{remark}
Let $g_1,\ldots,g_n$ be given $C^\infty$-smooth functions on a polydisc $B$ in $\Cn$
such that 
\begin{equation}\label{yy00}
\partial_{\bar{z}_j}g_k=\partial_{\bar{z}_k}g_j,\quad 1\leq j\leq n, 1\leq k\leq n
\end{equation}
Then there exists a $C^\infty$-smooth function $f$ on $B$
such that $\partial_{\bar{z}_j}f=g_j,$ $j=1,\ldots,n.$
Indeed, let $\gamma:=\sum_{j=1}^n g_j d\bar{z}_{j}\in \mathcal{A}^{0,1}(B).$
By Theorem \ref{dolbeaultlemma2}
there exists a $C^\infty$-smooth complex-valued {\em function}, $f\in \mathcal{A}^{0,0}$, such that
$\overline{\partial}f =\gamma.$
But $\overline{\partial}f=\sum_j \partial_{\bar{z}_j} f(z) d\bar{z}_j$ so the conclusion is equivalent 
to the system
\begin{equation}
\partial_{\bar{z}_j} f(z) =g_j(z),\quad j=1,\ldots,n
\end{equation}
On the other hand 
\begin{multline}
\overline{\partial} \gamma=\sum_j \overline{\partial} g_j\wedge d\bar{z}_j =
\sum_{j=1}^n \sum_{k=1}^n (\overline{\partial}_{\bar{z}_k} g_j)\wedge d\bar{z}_k \wedge d\bar{z}_j=\\
\sum_{k<j} \left(\overline{\partial}_{\bar{z}_j} g_k-\overline{\partial}_{\bar{z}_k} g_j-\right)\wedge d\bar{z}_k \wedge d\bar{z}_j
\end{multline}
Thus the conditions of Eqn.(\ref{yy00}) is equivalent to
\begin{equation}
\overline{\partial} \gamma =0
\end{equation}
\end{remark}
We should point out that  
determining the resolution to the chain complex where the $\overline{\partial}$-operator acts on complex differential forms of \index{Bidegree of complex differential form}bidegree $(p,q)$ can be done using functional analysis with certain advantages. 
Here is a short description of how $L^2$-methods\index{$L^2$-method} are used for such a resolution.
Our description glosses over many important details, for a rigorous in depth exposition, see e.g.\ Krantz \cite{krantzflerdim}, Ch. 4 and H\"ormander \cite{hormander}.
Let $H_1,H_2,H_3$ be Hilbert spaces of square integrable differential forms (the $H_j$ are to be thought of as forms of consecutive degree when acted upon by the operator $\overline{\partial}$), with inner products $\langle \cdot,\cdot\rangle_{H_j}, j=1,2,3$.
\begin{definition}[Adjoint]\index{Adjoint operator}
Let $\psi\in H_2$. 
Let $\mathbf{t}$ be a linear operator $H_1\to H_2.$ 
The {\em adjoint}, $\mathbf{t}^*$ of $\mathbf{t}$ is the operator defined via
$\langle \mathbf{t}x,y\rangle_{H_2}=\langle x,\mathbf{t}^* y\rangle_{H_1}$ 
(in particular in the case $H_1=H_2$ and a given basis we have a matrix representation of $\mathbf{t}$ and the adjoint is then represented by the transpose matrix).
Let $X_j\subset H_j$ be dense subspaces and let $T$ be a linear operator $X_1\to H_2.$ We denote $\psi\in X_1^{T^*}$ if there is a constant $C(\psi)>0$ such that $\abs{\langle T\phi,\psi\rangle_{H_2}}\leq C\norm{\phi}_{H_1}$ for all $\phi\in X_1.$
\end{definition}
We state without proof that for each $\psi\in X_1^{T^*}$ there exists a unique element $T^*\psi\in H_1$ such that
$\langle T\phi,\psi\rangle_{H_2}=\langle \phi,Y^*\psi\rangle_{H_1}.$ 
Now let $A\colon H_1\to H_2,$
$B\colon H_2\to H_3,$ be linear closed operators (both to the thought of as the operator $\overline{\partial}$) defined only on dense subspaces (thus their adjoint operators will also be defined on dense domains) such that $AB=0.$ Solving the overdetermined $Bf=g$, we require $Ag=0.$ Define the self-adjoint operator
$L:=A^*A+BB^*\colon H_2\to H_2.$
\begin{proposition}
If $\norm{g}^2\leq C(\norm{B^*g}^2+\norm{Ag}^2)$ for a positive constant $C$, then $L$ is invertible.
\end{proposition}
\begin{proof}
Given the estimate $L$ is injective since $Lg=0$ implies 
$0=C\langle (A^*A+BB^* )g,g\rangle=C(\norm{B^*g}^2+\norm{Ag}^2)\geq \norm{g}^2.$ Since it self-adjoint and taking inverse commutes with taking the adjoint, the operator $L$ is invertible on its domain.
\end{proof}
We can thus decompose
\begin{equation}\label{hopo}
g=(A^*AL^{-1})g+(BB^*L^{-1})g
\end{equation}
Recall that we are choosing $g$ such that $Ag=0$, then by Eqn.(\ref{hopo}) $0=Ag=AA^*AL^{-1} g$ (since $AB=0$). Thus $0=\langle AA^*AL^{-1} g,AL^{-1}g \rangle=\norm{A^*AL^{-1} g}$ which implies $g=B(B^*L^{-1}g).$
Specifying $B$ to be the operator $\overline{\partial}$ (viewed as an operator acting on complex $(p,q)$-forms for a fixed pair of nonnegative integers $(p,q)$) we obtain, given the above estimate, a solution
$f:=\overline{\partial}^* L^{-1} g$ to the equation $\overline{\partial} f=g.$ 
It is well-known that on bounded pseudoconvex domains
with smooth boundary the estimate holds true. Let us look a little closer at how the particular example
can be described.
A densely define linear operator between Hilbert spaces is called closed if its graph
is closed.
For example, for a domain $\Omega\subset\Cn,$ 
the unbounded linear 
operator $\overline{\partial}:L^2_{(0,q)}(\Omega)\to L^2_{(0,q+1)}(\Omega),$ can be verified to be a closed operator.
Its adjoint $\overline{\partial}^*:L^2_{(0,q+1)}(\Omega)\to L^2_{(0,q)}(\Omega),$ is also densely defined, linear and closed.
The orthogonal complement of the closure of the range of $\overline{\partial}$
in $L^2_{(0,q+1)}(\Omega)$ is the kernel of $\overline{\partial}^*$ and thus we have
an orthogonal direct sum. 
$L^2_{(0,q+1)}(\Omega)=\mbox{Ker}\overline{\partial}^* \oplus \overline{\mbox{Ran}\overline{\partial}}$
The $L^2$-version of the inhomogeneous Cauchy-Riemann equation is
to given $g\in L^2_{(0,q+1)}(\Omega)$ such that $\overline{\partial}g=0$, find
$f\in L^2_{(0,q)}(\Omega)$ such that $\overline{\partial}f\in L^2_{(0,q+1)}(\Omega)$, satisfying 
$\overline{\partial}f=g$. Often one is really interested in the orthogonal projection of the solution $f$,
onto the orthogonal complement of the closed subspace $\mbox{Ker}\overline{\partial}$ in $L^2_{(0,q)}(\Omega).$
Note that we can decompose $g=g_1+g_2,$ $g_1\in \mbox{Ker}\overline{\partial}^*,$ $\mbox{Ran}\overline{\partial}$,
in particular $g_2\in \mbox{Ker}\overline{\partial}.$
In the case of a domain $\Omega$ with smooth boundary and functions $f$ smooth up to and including the boundary,
the problem of finding a situation where the sought norm estimate can be obtained could be
handled by introducing the {\em weighted} inner products on $L^2_{(0,q)}(\Omega),$
$\langle f,g\rangle_\phi:=\int_\Omega \langle f,g\rangle\exp(-\phi)d\mu(z),$ for appropriate $\phi.$
with induced Hilbert space norm $\norm{\cdot}_\phi.$ It can be verified that $\overline{\partial}$ is still closed under these
new norms for appropriate $\phi$ and that the adjoint is well-defined.
The point is to arrive at 
$\norm{g}_\phi \leq C(\norm{\overline{\partial}g}_\phi +\norm{\overline{\partial}^*g}_\phi),$ for $g\in L^2_{(0,q+1)}(\Omega)
\cap\mbox{Dom}(\overline{\partial})\cap\mbox{Dom}(\overline{\partial}^*_\phi),$ $q\in N.$
Without giving the details we state that for a defining function $\rho$ for a bounded domain $\Omega$ with smooth boundary and
for appropriate choice of $\phi$ we have
\begin{multline}
\norm{\overline{\partial} f}_\phi^2 +\norm{\overline{\partial}^* f}_\phi^2 =\\
\int_\Omega \sum_J\sum_{j=1}^{n}\abs{\partial_{\bar{z}_j}f_J}^2\exp(-\phi) +
\int_\Omega \sum_{J'}\sum_{j,k=1}^{n}\partial_{\bar{z}_j}f_{j,J'}\overline{f_{k,J'}} \frac{\partial^2 \phi}{\partial z_j\partial \bar{z}_k}\exp(-\phi)+\\
\int_{\partial\Omega} \sum_{J'}\sum_{j,k=1}^{n}\partial_{\bar{z}_j}f_{j,J'}\overline{f_{k,J'}} \frac{\partial^2 \rho}{\partial z_j\partial \bar{z}_k}\exp(-\phi)
\end{multline}
Comparing this to Definition \ref{pseudoconvexdef} and noting that $\sum_j f_{j,J'}\partial_{z_j}\in T_z^{(1,0)}\partial\Omega$, we see that if $\Omega$ is pseudoconvex then for $\phi:=\abs{z}^2$ the second therm on the right hand side becomes $\norm{f}_\phi^2$ and in the last term we know that the matrix defined by $\frac{\partial^2 \rho}{\partial z_j\partial \bar{z}_k}$ is non-negative definite. This can be used to obtain the sought estimate.
One advantage of the above described $L^2$ method of proof is that it puts in place the opportunity to pursue the further investigation of regularity of the solution using the regularity properties of the right hand side (usually additional conditions such as plurisubharmonic exhaustion functions and subelliptic estimates are introduced). We shall not pursue this theory in this text but refer to e.g.\ Krantz \cite{krantzflerdim} and H\"ormander \cite{hormander}.

\section{Bol's Lemma}
We shall need this lemma in order to investigate how linear fractional maps interact with the polyanalytic differential operator.
\begin{theorem}[Bol's lemma]
If $f$ is holomorphic and 
\begin{equation}
F(w):=f(\phi(w))(E(w))^{\mu-1}
\end{equation}
where
\begin{equation}
\phi(w):=\frac{aw+b}{E(w)}
\end{equation}
\begin{equation}
E(w):=cw+d
\end{equation}
for $ad-bc=1.$
Then
\begin{equation}
\partial_w^\mu F=(\partial_w^\mu (\phi(w))(E(w))^{-(\mu+1)}
\end{equation}
\end{theorem}
\begin{proof}
Consider the Taylor expansion in $t$, near $t=0,$ of both sides of 
\begin{equation}
F(w+t):=f(\phi(w+t))(E(w+t))^{\mu-1}
\end{equation}
For the left hand side we have the coefficient for $t^k$ given by
\begin{equation}\label{klux1}
\frac{1}{k!}\partial^\mu_t F(w)
\end{equation}
For the right hand side using the relations
\begin{equation}
\phi(w+t)=\phi(w)+\frac{t}{E(w)E(w+t)}
\end{equation}
\begin{equation}
E(w+t)=E(t)+ct
\end{equation}
we obtain the  coefficient for $t^k$ as
\begin{equation}\label{klux2}
\sum_{j=0}^k \frac{1}{(k-j)!}\partial^{k-j}_t(\phi(w))\binom{\mu-1-k-j}{j}(E(w))^{\mu-1-2k+j}c^j
\end{equation}
Equating the expressions of Eqn.(\ref{klux1}) and Eqn.(\ref{klux2}) gives the desired result. This completes the proof. 
\end{proof}
We mention that an alternative faster proof can be given if one is familiar with group actions of $SL_2(\C)$. Namely,
it is known that the set of polynomials of degree less than $\mu$ is invariant under the action of the special linear group $\mbox{SL}_2(\C)$ action
$f(w)\mapsto f(\phi(w))(E(w))^{\mu-1}.$
Writing the Taylor polynomial, $T$, of $f$ near a point $p_0,$ with remainder $R$, i.e.\
$f=T+R$ near $p_0$ with $R(w)=O((w-p_0)^\mu).$
Then $R(\phi(w))(E(w))^{\mu-1}$ must be the remainder $f(\phi(w))(E(w))^{\mu-1}$ at the point $\phi^{-1}(p_0).$
 Bol's Lemma \cite{ball} can be used to show that the operator
 $\overline{\partial}^q$ for positive integers $q$, is invariant under linear fractional mappings.
 \begin{corollary}[Vasin \cite{vasin}] 
 If $A:\C\to \C,$ $z\mapsto \frac{az+b}{cz+d}$ then for $f\in C^2(\C)$
 \begin{equation}
 \overline{\partial}^q\left(f(Az)\overline{A'(z)}^{\frac{1-n}{2}}\right)=\overline{\partial}^q f(Az)\cdot \overline{A'(z)}^{\frac{1-n}{2}}
 \end{equation}
 \end{corollary}
In the main theorem of this section we shall consider the following special case of $\Omega.$
Let $\Omega$ be a compact Riemann surface of genus $m>1$ with projective structure. 
Assume $\Omega$ has atlas given 
by holomorphic charts $\{(\phi_i,U_i)\}_{i\in I}$, 
such that for any pair $i,j$ with $U_i\cap U_j\neq 0$ the coordinate change 
is not merely given by a linear fractional, but are in fact of the form 
$\phi_j\circ\phi_i^{-1}=\frac{a_{ij}z +b_{ij}}{\bar{b_{ij}}z +\bar{a_{ij}}}$ with $a_{ij}\bar{a_{ij}}-b_{ij}\bar{b_{ij}}=\abs{a_{ij}}^2-\abs{b_{ij}}^2=1.$ 
We set
\begin{equation}
A_{ij}:=\begin{bmatrix}
a_{ij} & b_{ij}\\
\bar{b}_{ij} & \bar{a}_{ij}
\end{bmatrix}
\end{equation} 
For a domain $U\subset\C,$ denote by $\mathcal{C}^{p,q}(U)$ the sheaf of germs of $(p,q)$-forms.
Bol's lemma implies that $\overline{\partial}_z^q$ is a sheaf homomorphism 
\begin{equation}
\mathcal{C}^{p,1-q}\overset{\overline{\partial}_z^q}{\longrightarrow} \mathcal{C}^{p,1+q}
\end{equation}
whose kernel, denoted $\mathcal{P}^{p,1-q}$, is the sheaf of $p,1-q$-forms annihilated by $\overline{\partial}_z^q.$ 

\section{A cochain complex}
Let $M$ be a Riemann surface with atlas $\{(U_\alpha,\phi_\alpha)\}_{\alpha}$.
Denote by $\mathcal{C}^{p,q}$ the sheaf of germs of complex $(p,q)$-forms.
On the intersection of any pair $U_\alpha\cap U_\beta$ the sections of the sheaf $\mathcal{C}^{p,q}$ are related by
\begin{equation}
f_\beta(z_\beta)dz_\beta^p d\bar{z}_\beta^q =f_\alpha(z_\alpha)dz_\alpha^p d\bar{z}_\alpha^q,\quad 2p,2q\in \Z
\end{equation}
The differentials for half-integer order $(p,q)$ are well-defined since the square root of the linear cotangent bundle $\mu^2=k$ is defined.
When $M$ is projective (i.e.\ has a projective covering $\{(U_\alpha,\phi_\alpha)\}_{\alpha}$) then on any
 $U_\alpha\cap U_\beta$ the role of a cocycle can be played by
\begin{equation}
\mu_\alpha\beta =(c_{\beta\alpha} z_\alpha +d\beta\alpha)^{-1}
\end{equation}
where 
\begin{equation}
z_\beta =\frac{a_{\beta\alpha} z_\alpha +b_{\beta\alpha}}{c_{\beta\alpha} z_\alpha +d_{\beta\alpha}}
\end{equation}
Since any $q$-analytic function $f$ has a representation of the form $\sum_{j=0}^{q-1}a_j(z)\bar{z}^j$ for holomorphic $a_j$, it can be identified as belonging to a germ of a $q$-analytic differential in the sheaf 
$\mathcal{P}^{p,1-q}$.
Given a projective coordinate, $(V,z),$ let $\mathscr{O}_z$ denote a fiber of a sheaf $\mathscr{O}$ of holomorphic functions. 
A fiber of the sheaf $\mathcal{P}^{p,1-q}$ has a direct decomposition according to
\begin{equation}
\mathcal{P}_z^{p,1-q}=\mathscr{O}_z\oplus\bar{z}\mathscr{O}_z\cdots\oplus \bar{z}^{q-1}\mathscr{O}_z
\end{equation} 
On a Riemann surface with projective structure this yields an exact sequence
\begin{equation}\label{exactseq}
0\longrightarrow \mathcal{P}^{p,1-q}\overset{\mbox{Id}}{\longrightarrow} \mathcal{C}^{p,1-q}\overset{\overline{\partial}_z^q}{\longrightarrow}
\mathcal{C}^{p,1+q}\longrightarrow 0
\end{equation}
where exactness at the last step follows from solvability of the $q$-analytic equation. 

Now the germ in the sheaf 
$\mathcal{P}^{p,1-q}$ to which a given $q$-analytic function $f$ belongs, in turn belongs to the ring $\mathcal{O}_z[\bar{z}]$ of polynomials. 
A polynomial $P$ over the field of fractions for $\mathcal{O}_z$ (i.e.\ the ring $\mathcal{M}_z$ of germs of meromorphic functions)
can be factorized as
\begin{equation}
P=\Pi_{i=1}^{q-1}(\varphi_i\bar{z}+\psi_i)
\end{equation}
where $\varphi_i,\psi_i,i=1,\ldots,q-1,$ are algebraic over $\mathcal{M}_z$ and are determined up to an invertible element of $\mathcal{O}_z.$
\begin{lemma}\label{vasinlemma1}
Let $P$ be a section of the sheaf $\mathcal{P}^{p,q-1}$ over an open neighborhood $U$. Then in projective coordinates, 
\begin{equation}
*P(z):=\sum_{i}^{q-1} a_i(z)\bar{z}^{q-1-i}
\end{equation}
defines a section of the sheaf $\mathcal{O}^{p+1-q}$ over $U.$ 
\end{lemma}
\begin{proof}
It suffices to verify the transformation law for the
passage from one coordinate neighborhood to another. 
To do this denote by $\mu_{\beta\alpha}=(\bar{b}_{\beta\alpha}z_\alpha +\bar{a}_{\beta\alpha})^{-1}$ 
the cocycle on the intersection $U_\beta\cap U_\alpha$ of two coordinate neighborhoods 
that determines the square root of the linear
cotangent bundle $\mu^2=k.$ Then
\begin{equation}
P_\beta(z_\beta)\mu_{\beta\alpha}^p\bar{\mu}_{\beta\alpha}^{1-q}=P_\alpha(z_\alpha)
\end{equation}
which can be factorized according to
\begin{equation}
\mu_{\beta\alpha}^p\bar{\mu}_{\beta\alpha}^{1-q}=\Pi_{i=1}^{q-1}(\varphi_{i\beta}\bar{z}_\beta +\psi_{i\beta})=
\Pi_{i=1}^{q-1}(\varphi_{i\alpha}\bar{z}_\alpha +\psi_{i\alpha})
\end{equation}
and
\begin{equation}
\mu_{\beta\alpha}^p\Pi_{i=1}^{q-1}(\varphi_{i\beta}(
\bar{a}_{\alpha\beta}\bar{z}_\alpha +b_{\beta\alpha}) +
\psi_{i\beta}(b_{\beta\alpha}\bar{z}_\alpha +a_{\beta\alpha}))=
\Pi_{i=1}^{q-1}(\varphi_{i\alpha}\bar{z}_\alpha +\psi_{i\alpha})
\end{equation}
By the uniqueness of the factorization (see e.g.\ Fritzsche \& Grauert \cite{frigrau}, III, 4) the corresponding factors are equal up to an invertible 
element of $\mathcal{O}(U_\alpha \cap U_\beta),$ i.e.\ there are holomorphic functions
$\xi_{i\beta\alpha},$ $i=1,\ldots,q-1,$ on $U_\alpha\cap U_\beta$ that are cocycles of linear 
holomorphic fiber bundles, $\xi_i$, over $U$ such that $\pi_{i=1}^{q-1} \xi_{i\beta\alpha}=\mu^p_{\beta\alpha}$ and 
\begin{equation}
\xi_{i\beta\alpha}(\varphi_{i\beta}(
\bar{a}_{\alpha\beta}\bar{z}_\alpha +b_{\beta\alpha}) +
\psi_{i\beta}(b_{\beta\alpha}\bar{z}_\alpha +a_{\beta\alpha}))=\varphi_{i\alpha}\bar{z}_\alpha +\psi_{i\alpha}
\end{equation}
Meromorphic linear fiber bundles are trivial bundles (see e.g.\ Gunning \& Rossi \cite{gunningrossi}, p.247) thus there exists for each
$i=1,\ldots,q-1,$
a cochain $\{\xi_{i\alpha}\}_\alpha$ of meromorphic functions such that
\begin{equation}
\xi_{i\beta\alpha} =\frac{\xi_{i\beta}}{\xi_{i\alpha}},\quad \xi_{i\alpha} \in \mathcal{M}(U_\alpha)
\end{equation}
Comparing the coefficients of $\bar{z}_\alpha$ gives
\begin{equation}
\begin{bmatrix}
\xi_{i\beta} \varphi_{i\beta}\\
-\xi_{i\beta} \psi_{i\beta}
\end{bmatrix}=
\begin{bmatrix}
a_{\beta\alpha} & b_{\beta\alpha}\\
\bar{b}_{\beta\alpha} & \bar{a}_{\beta\alpha}
\end{bmatrix}
\begin{bmatrix}
\xi_{i\alpha} \varphi_{i\alpha}\\
-\xi_{i\alpha} \psi_{i\alpha}
\end{bmatrix}
\end{equation}
Hence the family $(\xi_{i\beta} \varphi_{i\beta},-\xi_{i\beta} \psi_{i\beta})^T$ is a meromorphic section of a two-dimensional
vector bundle with cocycle transformation matrix $A_{\beta\alpha}.$
The same transformation matrix holds for the family $(d_\alpha z_\alpha,d_\alpha)^T,$ where
$d$ is a holomorphic $1$-differential over $U$, i.e.\
\begin{equation}
\begin{bmatrix}
\xi_{i\beta} \varphi_{i\beta} & d_\beta z_\beta\\
-\xi_{i\beta} \psi_{i\beta} & d_\beta
\end{bmatrix}=
\begin{bmatrix}
a_{\beta\alpha} & b_{\beta\alpha}\\
\bar{b}_{\beta\alpha} & \bar{a}_{\beta\alpha}
\end{bmatrix}
\begin{bmatrix}
\xi_{i\alpha} \varphi_{i\alpha} & d_\alpha z_\alpha\\
-\xi_{i\alpha} \psi_{i\alpha} & d_\alpha
\end{bmatrix}
\end{equation}
and taking determinants
\begin{equation}
d_\beta \xi_{i\beta} (\varphi_{i\beta} +z_\beta\psi_{i\beta})=d_\alpha \xi_{i\alpha}(\varphi_{i\alpha} +z_\alpha\psi_{i\alpha})
\end{equation}
This yields after multiplication over $i$
\begin{equation}
d_\beta^{q-1}\mu_{\beta\alpha}^p\Pi_{i=1}^{q-1} 
(\varphi_{i\beta} +z_\beta\psi_{i\beta})= d_\alpha^{q-1} \Pi_{i=1}^{q-1}
(\varphi_{i\alpha} +z_\alpha\psi_{i\alpha})
\end{equation}
i.e.\
\begin{equation}
\mu_{\beta\alpha}^{p+q-1}\Pi_{i=1}^{q-1} 
(\varphi_{i\beta} +z_\beta\psi_{i\beta})= \Pi_{i=1}^{q-1}
(\varphi_{i\alpha} +z_\alpha\psi_{i\alpha})
\end{equation}
which recalling the factorization used and returning to the $f_i$, implies
\begin{equation}
\mu_{\beta\alpha}^{p+q-1}\sum_{i=0}^{q-1}f_{i\beta}z_\beta^{q-1-i}=
\sum_{i=0}^{q-1}f_{i\alpha}z_\alpha^{q-1-i}\sum_{i=0}^{q-1}f_{i}z_\beta^{q-1-i}\end{equation} 
is
a holomorphic $(p+1-q)$-differential over $U.$ This completes the proof.
\end{proof}

Let $\lambda^{-1}$ be the metric on $\Omega$ defined in projective coordinates by 
$\lambda^{-1}=(1-\abs{z}^2)^{-1}.$ This metric is a global section of the sheaf $\mathcal{C}^{1,1}.$ If $h(z)$
is a representative of the germ of a differential in the sheaf $\mathcal{P}^{p+1,2-q}$, then $\lambda h(z)$
is a representative of the germ of a differential in the sheaf $\mathcal{P}^{p,1-q}.$
Denote by $\Lambda$ the operator corresponding to multiplication by $\lambda$ and denote by $*$ the operator corresponding to $*$
in Lemma \ref{vasinlemma1}.
\begin{lemma}\label{vasinlemma2}
The sequence
\begin{equation}
0\longrightarrow \mathcal{P}^{p+1,2-q}\overset{\Lambda}{\longrightarrow} \mathcal{P}^{p,q-1}
\overset{*}{\longrightarrow} 
\mathcal{O}^{p+1-n}\longrightarrow 0
\end{equation}
is exact.
\end{lemma}
\begin{proof}
Since $\lambda$ does not vanish on the unit disc, $\Lambda$ is injective. 
 For any $P\in \mathcal{O}_z^{p+1-q},$ the germ $\bar{z}^{q-1}P\in \mathcal{P}_z^{p,1-q}$ satisfies
 $*(\bar{z}^{q-1}P)=P$ which shows that $\Lambda$ is also surjective.
 Now let $h(z):=\sum_{i=0}^{q-2}f_i \bar{z}^i\in \mathcal{P}_z^{p+1,2-q} .$
 Then
 \begin{multline}
 *\Lambda h(z)=*\left( (1-\abs{z}^2)\sum_{i=0}^{q-2}f_i \bar{z}^i\right)=*\left( \sum_{i=0}^{q-2}f_i \bar{z}^i -\sum_{i=0}^{q-2}f_i z\bar{z}^{i+1}\right)=\\
 \sum_{i=0}^{q-2}f_i \bar{z}^{q-1-i} -\sum_{i=0}^{q-2}f_i z\bar{z}^{q-1-(i+1)}=0
 \end{multline}
 Thus $\mbox{Im}\Lambda\subset \mbox{Ker}*.$ On the other hand if $P(z):=\sum_{i=0}^{q-2}f_i \bar{z}^i\in \mathcal{P}_z^{p,1-q}$ satisfies $*P=0$ then
 factorizing $P$ we have
 \begin{equation}
 *P=*\Pi_{i=0}^{q-1} (\varphi_i\bar{z} +\psi_i)=\Pi_{i=0}^{q-1} (\varphi_i +z\psi_i)=0
 \end{equation}
 thus there is an $i_0$ such that $\varphi_{i_0} +z\psi_{i_0}=0,$ which implies
 \begin{equation}
 P(z)=\psi_{i_0}(1-\abs{z}^2)\Pi_{i\neq i_{0}}(\varphi_i\bar{z} +\psi_i)=
 (1-\abs{z}^2)\sum_{i=0}^{q-2} f'_i\bar{z}^i=(1-\abs{z}^2)h(z)
 \end{equation}
 where the coefficients $f_i'$, $i=1,\ldots,q-2,$ 
 are integral combinations over $\mathcal{O}_z$, of the $f_i,$ $i=0,\ldots,q-1.$
 Hence $h(z)\in \mathcal{P}^{p+1,2-q}$, proving that 
 $\mbox{Im}\Lambda\supset \mbox{Ker}*,$ thus by what we have previously done
 $\mbox{Im}\Lambda = \mbox{Ker}*.$ This completes the proof of the lemma.
\end{proof}

If $\{U_j\}_{j\in I}$ is an open cover of a topological Hausdorrf space $X$ and $\mathcal{F}$ a sheaf of abelian groups on $X.$ If $p$ is non-negative integer, denote by $s=(s_1,\ldots,s_p)\in I^{p+1}$ and $U_s=U_{s_0}\cap \cdots\cap U_{s_p}.$
A $p$-cochain $c$ of the given covering with values in $\mathcal{F}$ is the map which assigns to each $s$ a section $c_s\in \Gamma(U_s,\mathcal{F})$ such that $c_s$ is an alternating function of $s$ (i.e.\ $c_s$ changes sign if two indices are permuted). Let $\delta^p$ denotes the coboundary operator with respect to the given open cover (the atlas)
from the set of $p$-cochains to the set of $p+1$-cochains, namely $(\delta^p c)_s
=\sum_{j=0}^{p+1} (-1)^j c_{s_0\cdots \hat{s_j}\cdots s_{p+1}}$, where
$\hat{s_j}$ means that the index $s_j$ is removed.
These maps satisfy $\delta^{p+1}\circ \delta^{p}=0$ for all $p\geq 0,$ thus define a cochain complex
of sheaves. Usually, given a topological space $X$, an open cover $\mathcal{U}=\{U_i\}_{i\in I}$ and values in a sheaf $G$
one call the cophomolgy groups $\mbox{Ker}\delta^p/\mbox{Im}\delta^{p-1},$ $p\geq 0$, the $p$:th {\em $\Check{C}$hech cohomology group} (with respect to the cover $\mathcal{U}$). 
We mention that it is possible to define cohomology groups that do not depend 
upon the cover, by a process called the direct limit. Roughly speaking, it is based upon considering refinements: a cover
$\mathcal{V}=\{V_j\}_{j\in J}$ is called a refinement of $\mathcal{U}$ if there is a function $\iota: J\to I$
such that $V_j\subseteq U_{\iota(j)},$ $j\in J.$ This yields a direct preorder on the set of 
open covers and turns the set of cohomology groups (with values in a sheaf $\mathcal{F}$) 
with respect to all possible open covers into a so called direct mapping family so the
direct limit makes sense and is denoted $\mbox{lim}_{\overrightarrow{\mathcal{U}}} \Check{H}^p(\mathcal{U},\mathcal{F}).$
Let $A=(A^j)_{j\in Z}, B=(B^j)_{j\in \Z}, C=(C^j)_{j\in \Z}$ be three cochain complexes (where the $A_j,B_j,C_j$ are $R$-modules) 
and we have a short exact sequence of chain complexes
\begin{equation}
0\longrightarrow A\overset{f}{\longrightarrow} B\overset{g}{\longrightarrow} C\to 0
\end{equation}
i.e.\ for each $j$ we have exact sequences
\begin{equation}
0\longrightarrow A_j\overset{f_j}{\longrightarrow} B_j\overset{g_j}{\longrightarrow} C_j\to 0
\end{equation}
We shall construct a so-called {\em connecting homomorphism}\index{Connecting homomorphism} $\delta^*:H^j(C)\to H^{j+1}(A).$
Consider the diagram, where we denote by $d_A^j$ the coboundary operator, $A^j\to A^{j+1},$ of the complex $A$ and similarly for $B,C.$
\begin{equation}
\begin{array}{lllllllll}
0 &\longrightarrow & A^{j+1}          &\overset{f^{j+1}}{\longrightarrow} & B^{j+1} &\overset{g_{j+1}}{\longrightarrow} & C^{j+1} &\longrightarrow & 0\\
 &                 & \uparrow_{d_A^j} &                                   & \uparrow_{d^j_B} &                        &  \uparrow_{d^j_C} &  &\\
0 &\longrightarrow & A^j & \overset{f_j}{\longrightarrow} & B^j &\overset{g_j}{\longrightarrow} & C^j &\longrightarrow & 0\\
\end{array}
\end{equation}
Let $[c]\in H^j(C).$ Since $g_j$ is surjective there exists $b\in B^j$ such that $g_j(b)=c.$
Since the diagram is commutative $g_j d^j g_j b=d^j g_j b=d^j c=0$ because $c$ is a cocycle. By exactness
$\mbox{Ker} g_j =\mbox{Im} f_j,$ which implies $d^j b=f_j a$ for some $a\in A^{j+1}.$
By injectivity of $f_j$ $a$ is uniquely determined by $b$ and 
$f_j(d^j a)=d^j(f_j a)=d^jd^j b=0$ so
$d^j a=0$ thus $a$ is a cocycle and defines a cohomolgy class $[a].$ Set $\delta^*[c]=[a]\in H^{j+1}(A).$
This is a well-defined linear map $H^j(C)\to H^{j+1}(A),$ that in fact does not depend on the choices the cocycles $c$ and $b$ made.
Now for a cochain map $f:A\to B$ we have for cocycles $a\in Z^j(A),$ $d_B(f_j(a))=f_j(d_A a)=0$ and for coboundaries $a'\in A^{j-1},$ $f(d_A(a'))=d_B(f(a')).$
Thus $f$ induces a linear map of cohomolgy
\begin{equation}
f_j^*:H^j(A)\to H^j(B),\quad f_j^*[a]=[f(a)]
\end{equation}
Similarly, we have induced maps $g_j^*$ for $g_j.$
Then we claim that the following sequence of cohomology is exact (see e.g.\ H\"ormander \cite{hormander}, Thm 7.3.7, p.197, or Foster \cite{foster}, p.123)
\begin{multline}
\cdots \overset{g^*}{\longrightarrow} H^{j-1}(C)\overset{\delta^*}{\longrightarrow} H^j(B)
\overset{f^*}{\longrightarrow} H^j(B)\overset{g^*}{\longrightarrow}
 H^j(C)\overset{\delta^*}{\longrightarrow} H^{j+1}(A)\\
\overset{f^*}{\longrightarrow} H^{j+1}(B) \overset{g^*}{\longrightarrow} \cdots
\end{multline}
We prove exactness at $H^j(C).$ The exactness at $H^j(A)$ and $H^j(B)$ is obtained analogously.
\\
\textit{$\mbox{\em Im}g_j^*\subset \mbox{\em Ker}\delta^*$:}
Let $[b]\in H^j(B)$. Then $\delta^* g_j^*[b]=\delta^*[g_j(b)],$ and we can choose
the element in $B^j$ that maps to $g_j(b)$ as $b$ and take $d_B\in B^{j+1}$ to obtain
\begin{equation}
\begin{array}{lllll}
0 & \longrightarrow & d_B b & &\\
 & &\uparrow & &\\
 & & b&\longrightarrow & g_j(b)
 \end{array}
\end{equation}
Since $b$ is a cocycle $d_B b=0$, so by commutativity of the above diagram $d_B[g_j(b)]=0,$ hence $g_j^*[b]\in \mbox{Ker}\delta^*.$
\\
\textit{$\mbox{\em Im}g_j^*\supset \mbox{\em Ker}\delta^*$:}
Suppose $\delta^*[c]=[a]=0,$ where $[c]\in H^j(C).$ This means that $a=d_A(a')$ for some $a'\in A^j.$ The calculation of $\delta^*[c]$
can be written, for any element $b\in B^j$ with $g_j(b)=c$ as
\begin{equation}
\begin{array}{lllll}
a & \longrightarrow & d_B b & &\\
\uparrow_{d_A} & &\uparrow_{d_B} & &\\
a' & & b &\overset{g_j}{\longrightarrow} & c
\end{array}
\end{equation}
Thus 
$d_B(b-f_j(a'))=d_B b -d_B f_j(a')=d_B b-f_j d_A(a')=d_B b-f_j(a)=0,$
and $g_j(b-f_j(a'))=g_j(b)-g_j f_j(a')=g_j(b)=c.$
Hence $g_j^*[b-f_j(a')]=[c].$ This proves $\mbox{Im}g_j^*\supset \mbox{Ker}\delta^*$
and thus exactness at $H^j(C).$
\\
\\
Now the exact sequence in Eqn.(\ref{exactseq})
induces (see e.g.\ H\"ormander \cite{hormander}, Thm 7.3.7, p.197) the exact sequence of cohomologies (for a connection homomorphism $\delta^*$)
\begin{equation}
0\longrightarrow H^0(\mathcal{P}^{p,1-q})\overset{\mbox{Id}}{\longrightarrow} H^0(\mathcal{C}^{p,1-q})\overset{\overline{\partial}_z^q}{\longrightarrow}
H^0(\mathcal{C}^{p,1+q})\overset{\delta^*}{\longrightarrow}H^1(\mathcal{P}^{p,1-q})\longrightarrow 0
\end{equation}
\begin{theorem}
Let $\Omega$ be a compact Riemann surface of genus $>1$, with a projective structure satisfying the following conditions:
\begin{itemize}
\item[(i)] Assume $\Omega$ has atlas is given 
by holomorphic charts $\{(\phi_i,U_i)\}_{i\in I}$, 
such that for any pair $i,j$ with $U_i\cap U_j\neq 0$ the coordinate change 
is not merely given by a linear fractional, but are in fact of the form 
$\phi_j\circ\phi_i^{-1}=\frac{a_{ij}z +b_{ij}}{\bar{b_{ij}}z +\bar{a_{ij}}}$ with $a_{ij}\bar{a_{ij}}-b_{ij}\bar{b_{ij}}=\abs{a_{ij}}^2-\abs{b_{ij}}^2=1.$ 
\item[(ii)] $\Omega$ does not allow a family $\{g_\beta(z_\beta)\}_\alpha$
of holomorphic functions such that for any for two charts $(U_\alpha ,z_\alpha)$, $(U_\beta, z_\beta)$
\begin{equation}
\frac{g_\beta(z_\beta)}{2}\frac{\partial z_\beta}{\partial z_\alpha} -\frac{g_\beta(z_\alpha)}{2}=
\frac{\partial^2 z_\beta}{\partial z_\alpha^2}/\frac{\partial z_\beta}{\partial z_\alpha}
\end{equation}
\end{itemize}
Then the following sequences of finite-dimensional spaces for $q\geq 2$, are exact
\begin{equation}
0\longrightarrow H^0(\mathcal{P}^{p+1,2-q})\longrightarrow H^0(\mathcal{P}^{p,1-q})\longrightarrow 0,\quad p\leq q-1
\end{equation}
\begin{equation}
0\longrightarrow H^0(\mathcal{P}^{p+1,2-q})\longrightarrow H^0(\mathcal{P}^{p,1-q})\longrightarrow H^0(\mathcal{O}^{p+1-q})\longrightarrow 0,\quad p\geq q
\end{equation}
\end{theorem}
\begin{proof}
Denote by $\mathcal{D}^{2-p,2-q}$ the sheaf of continuous linear functionals over the sheaf $\mathcal{C}^{p,q}.$
The adjoint of the sheaf homomorphism
\begin{equation}
\mathcal{C}^{p,q-1}\overset{\overline{\partial}_z^q}{\longrightarrow} \mathcal{C}^{p,1+q}
\end{equation}
is given by
\begin{equation}
\mathcal{D}^{p,1-q}\overset{(-\overline{\partial}_z)^q}{\longrightarrow} \mathcal{D}^{2-p,1+q}
\end{equation}
Since $\overline{\partial}^q_z$ is elliptic, $\mbox{Ker}(-\overline{\partial}_z)^q =\mathcal{P}^{2-p,1-p}.$
Now the exact sequence of sheaves in Lemma \ref{vasinlemma2} induce (where as before one uses the snake lemma to obtain the
connecting homomorphisms $\delta^*$), since the $H^2$ cohomologies are trivial, an exact cohomology sequence
\begin{multline}
0\longrightarrow H^0(\mathcal{P}^{p+1,2-q})\overset{\Lambda}{\longrightarrow} H^0(\mathcal{P}^{p,1-q})
\overset{*}{\longrightarrow} H^0(\mathcal{O}^{p+1-q})\overset{\delta^*}{\longrightarrow} \\
\longrightarrow H^1(\mathcal{P}^{p+1,2-q})\overset{\Lambda}{\longrightarrow} H^1(\mathcal{P}^{p,1-q})
\overset{*}{\longrightarrow} H^1(\mathcal{O}^{p+1-q})
\longrightarrow 0
\end{multline}
For $q=2$ this reduces to
\begin{multline}\label{ekvationsex}
0\longrightarrow H^0(\mathcal{C}^{p+1})\overset{\Lambda}{\longrightarrow} H^0(\mathcal{P}^{p,-1})
\overset{*}{\longrightarrow} H^0(\mathcal{O}^{p-1})\overset{\delta^*}{\longrightarrow} \\
\longrightarrow H^1(\mathcal{O}^{p+1})\overset{\Lambda}{\longrightarrow} H^1(\mathcal{P}^{p,-1})
\overset{*}{\longrightarrow} H^1(\mathcal{O}^{p-1})
\longrightarrow 0
\end{multline}
The cohomologies with coefficients in holomorphic sheaves are finite-dimensional (see e.g.\ 
Gunning \& Rossi \cite{gunningrossi}, VIII.A, 10, p. 245) and by induction in the order of the operator
$\overline{\partial}_z^q$ also the cohomologies with coefficients in polyanalytic sheaves are finite-dimensional. This
implies that $H^0(\mathcal{P}^{2-p,1-q}),$ and $H^1(\mathcal{P}^{p,1-q})$ are finite-dimensional.
By applying the proof of the Serre duality theorem as it is presented in Gunning \cite{gunningriemann},
to the current settings Vasin \cite{vasin} deduces that  
\begin{equation}
H^1(\mathcal{P}^{p,1-q})=H^0(\mathcal{C}^{p,1+q})/\overline{\partial}_z^q H^0(\mathcal{C}^{p,1-q})
\end{equation}
has the same dimension as the space $H^0(\mathcal{P}^{2-p,1-q}).$ 
We shall proceed by induction in $q.$
\\
\\
\textit{Base case: q=2.} Consider the sequence in Eqn.(\ref{ekvationsex}). 
For $p<1$ and a given compact Riemann surface we have $H^0(\mathcal{O}^{p+1})=0,$ so we have an exact sequence
\begin{equation}
0\longrightarrow H^0(\mathcal{O}^{p+1})\overset{\Lambda}{\longrightarrow} H^0(\mathcal{P}^{p,-1})
{\longrightarrow} 0
\end{equation}
For $p>1$, the Serre duality theorem (see e.g.\ Gunning \cite{gunningriemann}, Thm 1.18) implies that the $H^1(\mathcal{O}^{p+1})$ is isomorphic to
$H^0(\mathcal{O}^{1-p})=0,$ which yields an exact sequence
\begin{equation}
0\longrightarrow H^0(\mathcal{O}^{p+1})\overset{\Lambda}{\longrightarrow} H^0(\mathcal{P}^{p,-1})
\overset{*}{\longrightarrow} H^0(\mathcal{O}^{p-1})
{\longrightarrow} 0
\end{equation} 
Since $H^0(\mathcal{O}^0)=\C,$ we have for $p=1$ 
\begin{equation}
0\longrightarrow H^0(\mathcal{O}^{2})\overset{\Lambda}{\longrightarrow} H^0(\mathcal{P}^{1,-1})
\overset{*}{\longrightarrow}\C
{\longrightarrow} 0
\end{equation}
If we prove that $\mbox{Im}*=\{0\}$ then we obtain also for $p=1$ 
\begin{equation}
0\longrightarrow H^0(\mathcal{O}^{2})\overset{\Lambda}{\longrightarrow} H^0(\mathcal{P}^{1,-1})
{\longrightarrow} 0
\end{equation}
If $\mbox{Im}*\neq \{0\}$ then there exists, in
a local coordinate $z_\alpha,$ a differential $P_\alpha =f_\alpha \bar{z}_\alpha +g_\alpha,$
such that $f_\alpha +z_\alpha g_\alpha=1.$ For two charts $(U_\alpha ,z_\alpha)$, $(U_\beta, z_\beta)$
we have on $U_\alpha \cap U_\beta$
$P_\beta \mu_{\beta\alpha}\bar{\mu}_{\beta\alpha}^{-1}=P_\alpha,$
thus
\begin{equation}
\frac{g_\beta(z_\beta)}{(\bar{b}_{\beta\alpha}z_\alpha +a_{\beta\alpha})^2}-g_\alpha(z_\alpha)=
-\frac{\bar{b}_\beta}{\bar{b}_{\beta\alpha}z_\alpha +\bar{a}_{\beta\alpha}}
\end{equation}
which can be written as
\begin{equation}
\frac{g_\beta(z_\beta)}{2}\frac{\partial z_\beta}{\partial z_\alpha} -\frac{g_\beta(z_\alpha)}{2}=
\frac{\partial^2 z_\beta}{\partial z_\alpha^2}/\frac{\partial z_\beta}{\partial z_\alpha}
\end{equation}
By condition (ii)
$\Omega$ allows no such family. 
The last statement implies $f_{\alpha}+z_\alpha g_{\alpha}\neq 1,$ which takes care of the case $p=1.$
This proves the theorem for the case of $q=2.$
\\
\\
\textit{Induction: assume the theorem is proved for $1,\ldots, q-1$.}
Since $H^0(\mathcal{O}^{p+1-q})=0$ for $p<q-1,$ we have an exact sequence
\begin{equation}
0\longrightarrow H^0(\mathcal{P}^{p+1,2-q})\overset{\Lambda}{\longrightarrow} H^0(\mathcal{P}^{p,1-q})
{\longrightarrow} 0
\end{equation}
For $p\geq q,$ we have by the duality theorem
$H^0(\mathcal{P}^{p+1,2-q})\simeq H^0(\mathcal{P}^{1-p,2-q})$, and by induction
$H^0(\mathcal{P}^{p+1,2-q})\simeq \cdots \simeq H^0(\mathcal{O}^{q-p-1}),$ which yields the exact sequence
\begin{equation}
0\longrightarrow H^0(\mathcal{P}^{p+1,2-q})\overset{\Lambda}{\longrightarrow} H^0(\mathcal{P}^{p,1-q})
\overset{*}{\longrightarrow} H^0(\mathcal{O}^{p+1-q})\longrightarrow 0
\end{equation}
Letting $p=q-1$ this reduces to (as before using the natural connection homomorphism $\delta^*$)
\begin{multline}
0\longrightarrow H^0(\mathcal{P}^{q,2-q})\overset{\Lambda}{\longrightarrow} H^0(\mathcal{P}^{q-1,1-q})
\overset{*}{\longrightarrow} H^0(\mathcal{O}^{0}) \overset{\delta^*}{\longrightarrow}\\
H^1(\mathcal{P}^{q,2-q}) \overset{\Lambda}{\longrightarrow} H^1(\mathcal{P}^{q-1,1-q})
\overset{*}{\longrightarrow} H^0(\mathcal{O}^{0})
\longrightarrow 0
\end{multline}
By the duality theorem we have for the last three spaces
$H^1(\mathcal{P}^{q,2-q})\simeq H^0(\mathcal{P}^{q-2,2-q})$,
$H^1(\mathcal{P}^{q-1,1-q})\simeq H^0(\mathcal{P}^{3-q,1-q})$ and
$H^0(\mathcal{O}^{0})\simeq H^0(\mathcal{O}^{2}),$ respectively.
By the induction hypothesis
\begin{equation}
H^0(\mathcal{P}^{2-p,2-q})\simeq \cdots \simeq H^0(\mathcal{P}^{-1,-1})\simeq H^0(\mathcal{O}^0)\simeq \C,
\end{equation}
\begin{equation}
H^0(\mathcal{P}^{3-p,1-q})\simeq \cdots \simeq H^0(\mathcal{O}^2)
\end{equation}
So the exact sequence becomes
\begin{multline}
0\longrightarrow H^0(\mathcal{P}^{q,2-q})\overset{\Lambda}{\longrightarrow} H^0(\mathcal{P}^{q-1,1-q})
\overset{*}{\longrightarrow} \C \overset{\delta^*}{\longrightarrow} \C\\
 \overset{\alpha}{\longrightarrow} H^0(\mathcal{O}^2)
\overset{\beta}{\longrightarrow} H^0(\mathcal{O}^{2})
\longrightarrow 0
\end{multline}
The homomorphism $\beta$ is surjective between isomorphic spaces thus $\mbox{Ker}\beta$ is trivial. By 
exactness $\alpha=0$ and $\delta^*$ is surjective. Since $\delta^*$ acts between isomorphic spaces, 
$\mbox{Ker}\delta^*$ is trivial. This implies by exactness that $*$ is the zero homomorphism.
This completes the induction. This completes the proof.
\end{proof}
Vasin \cite{vasin} states that the conditions (i) and (ii) of the theorem are automatic, but as the author of this book cannot verify this we have included them as conditions in the theorem.

\chapter{Cluster values and Picard type theorems}
\section{Preliminaries}
We recall first some preliminaries.
	Given $t\in \{\abs{z}=1\}\subset\C$ and $\theta\in (0,2\pi)$ define the {\em Stoltz angle}
	at $t$ to be the set
	\begin{equation}
	S(t,\theta):=\{ z\colon \abs{z}<1, \abs{\mbox{arg}(1-\bar{t}z)}<\theta\}
	\end{equation}
	A function $f$ on $\{\abs{z}<1\}$ is said to have {\em nontangential}\index{Angular limit} 
	(or {\em angular}) limit $a$ at $t\in \{\abs{z}=1\}$ if for every Stoltz angle $S(t,\theta)$ we have
	\begin{equation}
	a=\lim_{S(t,\theta)\ni z\to t} f(z)
	\end{equation}
	$a$ is then called the boundary value of $f$ at $t,$ and we define $f(t)=a.$ If the
	angular limit exists a.e.\ on $\{\abs{z}=1\}$ then the function
	$f^*:\{\abs{z}=1\} \to \C,$ $t\mapsto f(t)$ is called the {\em boundary function}\index{Boundary function} of $f.$
Recall that the Luzin-Privalov theorem states that if $f\in \mathscr{O}(\{\abs{z}<1\}$ has angular
limit zero on a set of positive Lebesgue measure then $f\equiv 0.$ On the other hand there are nontrivial holomorphic
functions whose {\em radial limit},
$\lim_{r\to 1-} f(rt),$ 
vanish almost everywhere.
Fatou's theorem states that any function $f\in H^p(\{\abs{z}<1\})$ has a boundary function $f^*\in L^p(\{\abs{z}=1\}),$
$p\in \Z_+$ and furthermore the norm of $f$ coincides with the norm of $f^*.$
Given a Jordan domain $\Omega$
it is known that there is a conformal mapping $\Phi:\Omega\to U$
of $\Omega$ to a domain $U$ bounded by finitely many circles. If 
$E\subseteq \partial\Omega$ consists of finitely many circular arcs, then 
there is bounded harmonic
function $\omega$ on $\Omega$ with boundary limit $1$ at points of $E$ and boundary
limit $0$ at points of $\partial\Omega\setminus E$ (except possibly at the endpoints
of the arcs that define up $E$), in fact it is given by the
Poisson integral
of the characteristic function of $E$, which can be
shown to exist by standard arguments of conformal mapping to the disc. 
By the Carath\'eodory theorem $\Phi$ and its inverse extend continuously to the respective
boundaries, thus a homeomorphism between the
closures, and we thus have well-defined pull-back, by $\Phi$, value $\omega(z,\Omega,E)$, at $z\in \Omega$.
Furthermore, note that when $\Omega$ is bounded, then by the maximum principle
the value $\omega(z,\Omega,E)$ is unique (and in the unbounded case uniqueness follows from the Phragmen-Lindel\"of theorem). 

\begin{definition}
	Let $\Omega\subset\C$ be a domain with boundary consisting of finitely
	many Jordan curves. Let
	$E$ be a finite union of arcs in $\partial\Omega$
	The {\em harmonic measure}\index{Harmonic measure} of $E$ at the point
	$z\in \Omega$ with respect to $\Omega$ is the value at $z$ of the (unique solution to the Dirichlet problem) bounded harmonic
	function $\omega$ on $\Omega$ with boundary limit $1$ at points of $E$ and boundary
	limit $0$ at points of $\partial\Omega\setminus E$ (except possibly at the endpoints
	of the arcs that define up $E$). We denote the harmonic measure by
	$\omega(z,\Omega,E)$.
\end{definition}

We shall need an auxiliary result which we state without proof.
\begin{theorem}[See e.g.\ Privalov \cite{privalov}, p.287]\label{fanpriv208}
	Let $f(z)$ be a bounded analytic function on $\{\abs{z}<1\}$ which possesses radial boundary values on a set $E\subset\{\abs{z}=1\}$
	of positive measure. Denote by $E_f$ the set of values of $f$ on $E.$ Let $G$ be an domain containing $E_f$ and $f(z_0)$ for some $z_0,$ $f(z_0)\notin E_f.$ Then
	\begin{equation}
	\omega_i(E_f,f(z_0),G)\geq \omega(E,z_0,\{\abs{z}<1\})
	\end{equation}
\end{theorem}

We can now deduce a version of the Luzin-Privalov uniqueness theorem for meromorphic functions (see, for example, Privalov \cite{privalov}, p.292 and Luzin \& Privalov \cite{luzinprivalov}).
\begin{theorem}\label{luzinprivalovmero}
	Let $f(z)$ be a nonconstant meromorphic function on $\{\abs{z}<1\}$ which possesses (finite or infinite) boundary values
	on a subset $E\subset\{\abs{z}=1\},$ with positive measure. Then
	$E_f:=f(E)$ contains a closed set of positive harmonic measure.
\end{theorem}
\begin{proof}
We can assume $E_f$ is a proper closed and bounded subset and that $E$ a perfect subset. Let $a\in \C\setminus E_f.$
Then the function $f_a:=1/(f(z)-a)$ takes the value $\infty$ on $E$. By a result of Luzin (see e.g.\ Privalov \cite{privalov}, Ch. I, 1.3)
there exists a perfect subset $E'\subset E$ of positive measure on which $f(z)$ is continuous. If we replace $E$ by $E'$ this yields
a set of values of $f$ that is closed and bounded. For each $z_0\in E$ let $S_{z_0,n}$ be a sector with vertex $z_0$, with sides formed by
rays with angle $\pi/4$ from the tangent, and with
radius $1/n$ for $n\in \Z_+.$ Let $M=\max_{E} \abs{f},$ and $n(z_0)\in \Z_+$ sufficiently large such that $\abs{f(z)}<2M$ on $S_{z_0,n}$.
the function $n(z_0)$ takes countably many values for each $E$ of positive measure thus there exists an $N\in \Z_+$ such that
the set of points in $E$ for which $n(z_0)=N$ has positive measure so we can choose a perfect set $P\subset E$ of positive mass
for which $n(z_0)=N$ for all points of $P.$ Let $r$ be such that $r>1-1/N$ such that $f(z)$ has no poles on $\{\abs{z}=r\}.$
Denote by $D$ the union of $\{\abs{z}<r\}$ and all sectors corresponding to points of $P.$ The boundary, $\Gamma$ of $D$ is the union of $\{\abs{z}=r\}$ and parts of rays emanating from end points of subintervals of $P$ at angle $\pi/4$ with respect to the tangent and the corresponding arcs in $\{\abs{z}=r\},$ thus $\Gamma$ will be a rectifiable curve.
	$f(z)$ is bounded on $\Gamma$ and has finitely many poles in $D.$
	We can deform $\Gamma$ so that the bounded set enclosed by $\Gamma$ does not contain any poles, and simultaneously not change $\Gamma$ in a neighborhood
	of $P.$
	Then $f(z)$ is bounded in the bounded set enclosed by $\Gamma.$ Now under the extension of a conformal map of the bounded set enclosed by $\Gamma$ onto the unit disc, $P$ is mapped on a set $P_1$ of positive measure in $\{\abs{z}=1\}$, and $f$ is transformed to a bounded function $\varphi$.
	The function $\phi$ is bounded and possesses a.e.\ on $P_1$ boundary values along all nontangential paths. 
	We choose a perfect subset $P_1'\subset P_1$ of positive measure on which all points $\phi(z)$ have angular boundary values and such that
	the angle is preserved.
	Note that the tangent of the curve $\Gamma$ agrees a.e.\ with on $P$ with the tangent of $\{\abs{z}=1\}.$ 
		Thus
	at a.e.\ point of $P$ the nontangential paths to $\Gamma$ also nontangential paths to $\{\abs{z}=1\}$ and conversely nontangential paths to
	$\{\abs{z}=1\}$ are nontangential also to $\Gamma.$
	Let $E_f'$ be the set of values of $\varphi(z)$ on the set $P'_1.$
	If the values of $\phi$ on $\{\abs{z}<1\}$ belong to $E_f'$ then $E_f'$ has positive harmonic measure. since $\phi$ is nonconstant. Thus $E_f'$ must have an interior point. If for some $p_0$ we have $\varphi(p_0)\notin E_f'$ then by Theorem \label{fanpriv208}, for a sufficiently large $R>0$,
	\begin{equation}
	\omega(E_f',\phi(p_0),\{\abs{z}=R\})\geq \omega(P_1',p_0,\{\abs{z}<1\})>0
	\end{equation}
	Hence $E_f'$ has positive harmonic measure. Since $E_f'\subset E_f$ this completes the proof.
\end{proof}
\begin{corollary}
Let $f(z)$ be a meromorphic function on $\{\abs{z}<1\}$ which takes nontangential boundary value $0$
on a subset $E\subset\{\abs{z}=1\},$ with positive measure. Then $f\equiv 0.$
\end{corollary}
\begin{remark}
	Note that if $\Omega,\Omega'$ are domains bounded by rectifiable Jordan curves and $\phi:\Omega\to \omega'$ a conformal map,
	$E\subset\partial \Omega,$ $E'\subset\partial\omega'$ $z\in \Omega,$ and
	$h$ is the harmonic function $\omega(\phi(z),\Omega',\phi(E))$ then $h\circ \phi=
	\omega(z,\Omega,E),$ i.e.\ harmonic measure is conformally invariant.
	In the case of the plane (but not in higher dimension), $\omega(z,E,\Omega)$ for a simply connected domain $\Omega$
	is the image of normalized Lebesgue
	measure on the unit circle under a conformal map $\Phi:\{\abs{z}<1\}\to \Omega$ with $\Phi(0)=z.$
	It is known 
	(see Krantz \cite{krantzp203}, p. 203, and Riesz \& Riesz \cite{riesz1916})
	that if $\Omega$ is a simply connected domain bounded by a rectifiable Jordan curve then Harmonic measure and linear measure are mutually absolutely continuous.
\end{remark}
We have thus proved a meromorphic version of the Luzin-Privalov theorem.  
\begin{theorem}[Luzin-Privalov theorem for meromorphic functions]\label{luzinprivalovmero1}
	Let $f(z)$ be a nonconstant meromorphic function on $\{\abs{z}<1\}$ which takes angular boundary values
	zero on a subset $E\subset\{\abs{z}=1\},$ with positive Lebesgue measure then
	$f\equiv 0.$
\end{theorem}
Here is another result we shall use in Chapter \ref{bvpsec}. 
\begin{theorem}[Golubev-Privalov, see e.g. Privalov \cite{privalov}, III, p.202 and Golubev \cite{golubev}, \cite{golubev1}]\label{golubev}
	Let $\gamma$ be a rectifiable Jordan curve in the plane that enclose a bounded domain $U$. Let $g$ be
	an $L^1$ function on $\gamma.$
	Then $g$ is the boundary value of a holomorphic function iff
	\begin{equation}\label{golubevek}
	\int_\gamma \zeta^k g(\zeta)d\zeta=0,\mbox{ for all }k\in \N
	\end{equation}
\end{theorem}
\begin{proof}
	Consider the Cauchy-type integral
	\begin{equation}
	G(z)=\frac{1}{2\pi i} \int\frac{g(\zeta)d\zeta}{\zeta -z}
	\end{equation}
	By the (Sokhotsky) jump formula we know that $G$ has nontangential (inner) boundary values, $G^+$, that coincide
	a.e.\ with $g$ if and only if the outer boundary values $G^-$ vanish, which is equivalent to
	$G\equiv 0$ on the complement of $(U\cup \gamma).$
	The result follows from the expansion at infinity
	\begin{equation}
	\frac{1}{2\pi i}\int_\gamma \frac{f(\zeta)d\zeta}{\zeta -z}=-\sum_{k=0}^\infty \frac{1}{2\pi i}\frac{1}{z^{k+1}} \int_\gamma \zeta^k f(\zeta)d\zeta
	\end{equation}
	Conversely every boundary value of a holomorphic function on $U$ satisfies Eqn.(\ref{golubevek}). 
	This completes the proof.
\end{proof}
This also yields and interesting corollary.
\begin{corollary}\label{schwarzapplic1lem1}
	Let $\Omega$ be a simply connected, bounded domain, bounded by a closed analytic curve, $\gamma$.
	A function $t(z)\in C^0(\gamma)$ is the boundary values of a function 
	$\mathscr{O}(\Omega)\cap C^0(\overline{\Omega})$
	if
	\begin{equation}
	\int_\gamma t(z)f(z)dz =0
	\end{equation}
	for all $f\in \mathscr{O}(\Omega)\cap C^0(\overline{\Omega}).$
\end{corollary}
The classical Lindel\"of theorem \cite{lindelof} is the following.
\begin{theorem}
	Let $f(z)$ be a bounded analytic function in an angular domain (sector) $S=\{\abs{\mbox{arg}z}<\theta_0,0<\abs{z}<R_0\}$ 
such that $f(z)$ tends to a unique value $\lambda$ as $z\to 0$ along a continuous path $\ell$ in $S$ terminating at $z=0.$
Then $f(z)$ tends to $\lambda$ uniformly as $z\to 0$ inside any angular subdomain (inner sector)
of opening $-\theta_0+\delta\leq \mbox{arg} z\leq \theta_0-\delta.$
\end{theorem}
\begin{proof}
	Assume first that $\lambda\neq\infty$ and set for $n\in \N$
	\begin{equation}
	S_n:=\left\{ z:\frac{R}{2^{n+2}}<\abs{z}<\frac{3R}{2^{n+2}}, \abs{\mbox{arg} z}<\theta_0\right\}
	\end{equation}
	\begin{equation}
	f_n(z):=f\left(\frac{z}{2^n}\right),\quad z\in S_0
		\end{equation}
		Then $f_n(S_0)=f(S_n),$ $n\in \Z_+.$
		Since $f(z)$ is bounded it omits at least two values thus by Theorem \ref{vitaliporterthm} (Vitali) applies since
		$\lim_{n\to \infty} f_n(z)=\lambda$ on the segment $\mbox{arg}z=\gamma$,
		$R/4<\abs{z}<3R/4.$
		Hence $f_n(z)\to \lambda$ uniformly on compacts of $S_0.$ W.l.o.g.\ we may assume $f_n\to \lambda$ uniformly on compacts of
		\begin{equation}
		K:=\left\{ z:\frac{R}{4}<\abs{z}<\frac{3R}{4}, \abs{\mbox{arg} z}<\theta_0\right\}
		\end{equation}
		For $\epsilon>0$ there exists $n_0=n_0(\epsilon,K)$ such that
		\begin{equation}\label{h897867}
		\abs{f_n(z)-\lambda}<\epsilon,\quad n\geq n_0,\quad z\in K
		\end{equation}
		Hence if $\abs{\zeta}\leq 3R/2^{n_0+2}$, $-\theta_0+\delta<\mbox{arg}z<\theta_0-\delta$
		then Eqn.(\ref{h897867}) implies
		\begin{equation}
		\abs{f_n(\zeta)-\lambda}<\epsilon
		\end{equation}
		If $\lambda=\infty$ we can apply the analogous arguments to $1/f(z)$. This completes the proof.
		\end{proof}
By applying a conformal mapping this can be reformulated as follows.
\begin{theorem}\label{lindelofthm}
Let $f(z)$ be a bounded analytic function in $\{\abs{z}<1\}.$
If $f(z)\to \lambda$ as $z\to \exp(i\theta_0)$ along some arc $\ell$ in the unit 
disc and terminating at $\exp(i\theta_0)$, then $f(z)$ tends to $\lambda$ uniformly, as $z\to \exp(i\theta_0)$ inside any
angular domain (inner sector) in $\{\abs{z}<1\}$ having $\exp(i\theta_0)$ as vertex.
	\end{theorem}
	We now prove a generalization of the classical Lindel\"of theorem.
\begin{theorem}
		Let $f(z)=u(z)+i(z)$ be analytic and bounded in $\{\abs{z}<1\}$ and let $\ell_1,\ell_2$ be two arcs in $\{\abs{z}<1\}$ terminating at
		a point $p\in \{\abs{z}=1\}.$
		If $u(z)\to a$ as $z\to p$ along $\ell_1$
		and $v(z)\to b$ as $z\to p$ along 
		$\ell_2$ then $f(z)\to a+ib$ as $z\to p$ inside any angular domain of opening
		$\pi-\epsilon$ having vertex $p$ and bisected by the radius drawn to $p.$
\end{theorem}
\begin{proof}
	Setting $g=f(z)-a-ib$ we may assume $\re g\to 0$ and $\im g g\to 0$ respectively as $z\to p$
	along $\ell_1$ and $\ell_2$ respectively and w.l.o.g.\ we assume $p=1.$\\ 
	\textit{Case 1:} Suppose $\ell_1\cap \ell_2=\emptyset$ in a neighborhood of $z=1.$
	If necessary, after a conformal mapping, it suffices to show that if $g(z)$ is a bounded analytic map
	on $\{\abs{z}<1\}$ such that the boundary function $\re g(\exp(i\theta))$ ($\im g(\exp(i\theta)$ respectively) 
	has limit $0$ at $z=1$, when approaching the point $z=1$ from above (below respectively), then $g$ has radial limit $0$ at $z=1$.
	\\
	Set $g_1(z):=\bar{g}(\bar{z})$. Then $g_1$ is a bounded analytic function on $\{\abs{z}<1\}$ and so is $gg_1.$ Also 
	$\re g(\exp(i\theta))$ ($\im g(\exp(i\theta))$ respectively) has limit $0$ at $z=1$, which is a point of continuity.	
	\begin{definition}
		A sequence $\{ c_j\}_{j\in\N}$ is called {\em summable by arithmetic means}\index{Summability by arithmetic means}
		if $s_j=\frac{1}{n}(\sum_{j=0}^n c_j)$ tends to a limit as $n\to \infty.$ 
	\end{definition}
	Note that if $s_j\to s$ then also $c_j\to s$
	since if $s_j=s+\delta_j$ then $c_j=s+\frac{1}{n}(\sum_{j=0}^n \delta_j)$ where the sum in the right hand side goes to $0$
	if $\delta_j\to 0.$ However $c_j$ may have a limit even if $s_j$ does not, take e.g.\ the series $1-1+1-1+\cdots.$ Which has
	partial sums alternating between $1$ and $0$ but $s_j\to \frac{1}{2}.$
	By the so-called Fej\'er's theorem (see e.g.\ Titchmarsh \cite{titchmarsh0}, p.414) the Fourier series of $gg_1$ is summable by arithmetic means
	at $z=1$ and its real part has radial limit $0$ at $z=1.$ Hence $\abs{g}^2$ has radial limit $0$ at $z=1$.
	\\
	\textit{Case 2:} Suppose $\ell_1\cap \ell_2$ intersect arbitrary close to $z=1.$
		We can assume that $\ell_1,\ell_2$ are non-self-intersecting polygonals and their paths determine simple
		closed polygons. Let $\ell_j'$ denote that part of $\ell_j$, $j=1,2,$ on one such polygon. Let $z_1,z_2$ be the end points of $\ell_j.$ Let $\abs{u-a}$
		($\abs{v-b}$) be bounded by $c>0$ on $\ell_1$ ($\ell_2$).
		Let $M:=\sup_{\abs{z}<1} \abs{g}.$  
		Now we map the interior of the polygon onto the unit disc such that $z_1,z_2$ respectively are mapped to $1,-1$ respectively. Let $\ell'$ denote the inverse image
		of the diameter through $\pm 1.$  
		Let $f(\xi)$, $\abs{\xi}\leq 1$ denote the transform of $g(z)$ and $f_1:=\bar{f}(\bar{\xi})$. 
			Now
			$\re f <c$ on $\{\abs{z}=1,\im \xi\leq 0\}$, 
			and $\im f<c$ on $\{\abs{xi}=1,\im \xi\geq 0\}$, $\abs{f}<M$ on $\abs{\xi}\leq M.$
			The function $ff_1$ is analytic which implies that $\re ff_1$ is harmonic on the closed disc
			thus satisfies on the diameter through $\pm 1$,
			$\re ff_1<2cM.$ On the other hand we have on that diameter
			$ ff_1=\abs{f}^2=\re ff_1<2cM.$ 
						\begin{equation}
			\abs{g}^2\leq 2cM
			\end{equation}
			on the diameter though $\pm 1.$
			By Theorem \ref{lindelofthm} this completes the proof.
		\end{proof}
\begin{definition}[Concordant functions]
	Let $\Omega\subset\C$ be a domain bounded by an OSCAR (see Definition \ref{oscardef}), $\gamma$, with Schwarz function
	$S(z),$ such that on a neighborhood $U$ of $\gamma$ we have $\gamma=\{S(z)=\bar{z}\},$ and let $f(z)=\sum_{j=0}^{q-1} a_j(z)\bar{z}^j$, for holomorphic $a_j(z),$ be a $q$-analytic function on $\Omega.$ Let $p\in \gamma.$
	The function
	\begin{equation}
	f_p(z):=\sum_{j=0}^{q-1} a_j(z)\bar{p}^j
	\end{equation}
	is called the (holomorphic) function concordant with $f$ at $p,$ and
		the function (defines on $U\cap \Omega$)
		\begin{equation}
		f_\gamma(z):=\sum_{j=0}^{q-1} a_j(z)(S(z))^j
		\end{equation}
		is called the (holomorphic) function concordant with $f$ on $\gamma$.
\end{definition}
\begin{proposition}\label{concordantprop}
	Let $\Omega\subset\C$ be a domain bounded by an OSCAR (see Definition \ref{oscardef}), $\gamma$, with Schwarz function
	$S(z),$ such that on a neighborhood $U$ of $\gamma$ we have $\gamma=\{S(z)=\bar{z}\},$ and let $f(z)=\sum_{j=0}^{q-1} a_j(z)\bar{z}^j$, for holomorphic $a_j(z),$ be a $q$-analytic function on $\Omega.$ Let $p\in \gamma.$ Then
		\begin{equation}\label{utopot1}
	f_p(z):=f(z)+\sum_{j=0}^{q-1}\frac{1}{j!}(\bar{p}-\bar{z})^j \partial_{\bar{z}}^j f(z)
	\end{equation}
	\begin{equation}\label{utopot2}
	f_\gamma(z):=f(z)+\sum_{j=0}^{q-1}\frac{1}{j!}(S(z)-\bar{z})^j \partial_{\bar{z}}^j f(z)
	\end{equation}
\end{proposition}
\begin{proof}
	Set for $(z,w)\in\Omega\times \C$
	\begin{equation}
	F(z,w):=\sum_{j=0}^{q-1} a_j(z)w^j
	\end{equation}
	Taylor Expansion with fixed $z$ near $w=w_0$ gives
	\begin{equation}\label{utopot}
	F(z,w):=F(z,w)+\sum_{j=0}^{q-1}\frac{1}{j!}(w-w_0)^j \partial_{w}^j F(z,w)|_{w=w_0}
	\end{equation}
	Since $\partial_{w}^j F(z,w)|_{w=\bar{z}}=\partial_{\bar{z}}^j f(z)$ we have with $w_0:=\bar{z}$ and $w=\bar{p}$ in Eqn.(\ref{utopot})
	the formula given by Eqn.(\ref{utopot1}). Similarly, setting $w_0=\bar{z}$ and $w=S(z)$ we have Eqn.(\ref{utopot2}).
	This completes the proof.
	\end{proof}

\section{Cluster and boundary values}
Let $H$ be a holomorphic function on $D=\{\abs{z}<1\}$ mapping the disc univalently onto a 
Jordan domain, $\Omega$, with rectifiable boundary, $\gamma$, such that $\abs{H'(z)}$ is bounded below by a positive number.
Let $C^0(\gamma,\C)$ be the Banach space equipped with the uniform norm. Denote by $S_q$ the subspace of
$C^0(\gamma,\C)$ consisting of boundary values of $q$-analytic functions in $\Omega$ that are continuous on $\overline{\Omega}.$
Denote $S=\bigcup_q S_q.$
\begin{theorem}[Mazalov \cite{mazalov}]
If \begin{equation}
\int_D \abs{H''(z)}drd\theta <\infty
\end{equation}
(where $z=re^{i\theta}$) then $S$ has the first category in $C^0(\gamma,\C)$.
\end{theorem}
In particular, since $C^0(\gamma,\C)$ is a Banach space, the Baire Category Theorem implies that 
\begin{equation}
S\neq C^0(\gamma,\C)
\end{equation}
i.e.\ there are continuous functions on $\gamma$ that cannot be identified as the boundary values of any $q$-analytic function
on the domain enclosed by $\gamma.$\\
\\
When $q=1$ the well-known Luzin-Privalov boundary uniqueness theorem implies that a solution (a $1$-analytic function on $D$), must be
uniquely determined by the values of $f$ on a subset of $\gamma$ of positive length.
Already in the case $q=2$ this is no longer true as can be exemplified by the function
$f(z)=(1-z\bar{z})z^N$ for a nonnegative integer $N$, and $\Omega=\{\abs{z}<1\}.$
On the other hand, Mazalov \cite{mazalov} recently proved that
there exists a Jordan domain $\Omega$ bounded by a $C^{\infty}$-smooth curve $\gamma$ which has the following property:
Let $\Omega'\subset\Omega$ be any Jordan domain such that it has a boundary arc $\kappa\subset\gamma$ and let $f$ be $q$-analytic in $\Omega'.$
If $f$ has angular boundary values equal to zero on a subset of $\kappa$ of positive length then $f\equiv 0.$
\begin{definition}
Let $\Omega$ be a domain, let $p\in \partial\Omega$ and let $\gamma^p$ be a continuous path
from $\Omega$ into by which we mean an arc with equation $z=\lambda(t)$, where
$\lambda(t)$ is continuous on $[0,1],$
$\lambda([0,1))\subset\Omega$ $\lambda(1)=p.$ Let $f$ be a function on $\Omega.$ 
The {\em complete cluster set}\index{Complete cluster set} of $f$ at $p$, denoted $C(f,p,\Omega).$
Recall that $C(f,p,S)$ is the set of points $A$ such that
for each $a\in A$ there exists a sequence $\{z_j\}_{j\in \N}$ with $z_j\to p$ as $j\to \infty$ and
$\lim_{j\to \infty} f(z_j)=a\in A.$
We shall by $C(f,p,\gamma^p)$ shall denote that the cluster set
of $f$ at $p$ {\em along} the path $\gamma^p$.
If $C(f,p,\gamma^p)=\{s_0\}$ for some $s_0\in \C\cup\{\infty\}$ the $s_0$ is called an {\em asymptotic value}\index{Asymptotic value}
of $f$ at $p.$ If $\gamma^p$ is a closed regular analytic arc then $s_0$ is called a {\em regular asymptotic value}
of $f$ at $p.$ If there are $k$ mutually non-contiguous and mutually non-tangential, at $p$, regular paths
from $\Omega$ into $p$ such that$f$ has asymptotic value $s_0$ at $p$ along each the $k$ paths, then $f$ is said to have at the point $p$
the {\em regular asymptotic value of order not less than $k$}\index{Order of a regular asymptotic value}.
Let $E\subset\partial \Omega$ and let $\Gamma=\{\gamma^p:p\in E\}$ be a family of continuous paths in $\Omega$ 
leading into points of $E.$ The  family is called {\em uniformly continuous}\index{Uniformly continuous family och paths}
if there exists a function $z=\phi(t,p)$ which defines a path $\gamma^p$ for any fixed $p\in E$ with $t$ varying in $[0,1]$, 
which is uniformly continuous (with respect to $p,t$) for $p\in E$ and $t\in [0,1]$.
\end{definition}
\begin{definition}
Let $S$ be a circular sector with vertex $p$ (we assume the sector is open and with angle not greater than $\pi$ radians). A sector $S'$ (also open) is said to be
{\em inner}\index{Inner sector} to $S$ if it has the same vertex $p$ and its closure $\overline{S}$ belongs to
$S\cup \{p\}.$ The least among the angles between boundary radii of the sector $S'$ and the boundary radii of $S$
is called the {\em inclination}\index{Inclination of a sector} of $S'$ to $S$. Note that if a path $\gamma$ leads from $S$ into $p$ then there exists an OSCAR
$\tilde{\gamma}$ which contains $\gamma$ ($\gamma$ can be analytically continued beyond $p$). 
\end{definition}
A cluster set is nonempty and closed. A cluster set that contains only a single point is called {\em degenerate}\index{Degenerate cluster set}.
Given a function $w=f(z)$. A cluster set of $f$ at a point $f$
that covers the entire $w$-sphere, it is called {\em total}\index{Total cluster set}. 
Perhaps surprisingly, existence of boundary values of polyanalytic functions of arbitrary order can be 
treated only via conditions of first order.

\begin{theorem}[Balk \& Tutschke \cite{balktutschke}]
	Suppose $f$ is bounded and continuous in the open circular sector $S$ with vertex $s_0$ such that 
	$\partial_{\bar{z}} f$ exists a.e.\ (in the Sobolev sense) and
	$\partial_{\bar{z}} f\in L^p(\overline{\Omega}),$ for some $2<p<\infty.$
	Suppose that there exists an open circular sector $S$ with vertex $s_0$ such that 
	there exists two paths $\gamma_1,\gamma_2$ leading to $s_0$ in $S$, satisfying
	$\re f(z)\to a'$ ($\im f(z)\to a''$ respectively) as $\gamma_1 \ni z\to s_0$ ($\gamma_2 \ni z\to s_0$ respectively).
	Then for each choice of sector $S'$ (closed in $S$) inner to $S$ with the same vertex, $s_0,$ we have $f(z)\to a'+ia'',$
	as $S'\ni z\to s_0.$
\end{theorem}
\begin{proof}
		We shall begin by proving a lemma. 
	\begin{lemma}\label{balktuslem}
		Suppose $\Omega$ is a domain bounded by a Jordan curve $\gamma.$
		Let $S\subset\Omega$ such that $\gamma\cap \partial S=\{s_0\}.$ and let 
		$f$ be a continuous function on $\Omega$ such that
		$\partial_{\bar{z}} f$ exists a.e.\ (in the Sobolev sense) and
		$\partial_{\bar{z}} f\in L^p(\overline{\Omega}),$ for some $2<p< \infty.$ Then there exists a
		holomorphic function $h$ in $\Omega$
		and a constant $d$ such that $C(f,s_0,S)=C(h,s_0,S)+d.$
	\end{lemma}
	\begin{proof}
				Consider the Pompieu integral
		\begin{equation}
		\phi(z)=-\frac{1}{\pi}\int_\Omega \frac{\partial_{\bar{\zeta}} f(\zeta)}{\zeta -z}d\zeta
		\end{equation}
		which in well-known to satisfies (see Lemma \ref{lemmanvibalktus}) $\partial_{\bar{z}} \phi =
		\partial_{\bar{z}} f(z)$ a.e.\ on $\Omega.$
		By The H\"older inequality we have a constant $c_1$ such that 
		(see Theorem \ref{vekuaineqthm})
		for any $z_1,z_2\in \C$
		\begin{equation}
		\abs{\phi(z_1)-\phi(z_2)}\leq \abs{z_1-z_2}^{\frac{p-2}{2}}
		\end{equation}
		Hence $\phi$ is continuous on $\C,$
		thus
		\begin{equation}
		\partial_{\bar{z}}(f-\phi)=0\mbox{ a.e.} \Omega
		\end{equation}
		Since $f-\phi$ is continuous, there exists a holomorphic function
		$h$ on $\Omega$ such that $f=h+\phi$ a.e.\ on $\omega.$ Set $d:=\phi(s_0)$ and let 
		$a\in C(h,s_0,S)$ be an arbitrary cluster value. By definition there exists 
		a sequence $\{z_j\}_{j\in \N}$ such that $z_j\to p$ as $j\to \infty$ and
		$\lim_{j\to \infty} h(z_j)=a,$ i.e.\ $\lim_{j\to \infty} f(z_j)= a+d.$
		Hence $C(f,s_0,S)\subset C(h,s_0,S)+d$. Similarly, we can realize that $C(h,p,S)+d\subset C(f,p,S)$.
		This completes the proof.
	\end{proof}
	(The lemma holds also for the case $p=\infty$ and the proof is analogous).
	By Lemma \ref{balktuslem}, $f=h(z)+\phi(z)$ for $z\in S,$ for holomorphic $h$ on $S$ and continuous $\phi$ on $\C.$
	Denote $\phi(s_0)=:d=d_1+id_2.$ 
	Then for each closed inner sector $S'$ to $S$ we have for each sequence $\{z_j\}_{j\in \N}$
	in $S$ that if $z_j\to s_0$ then $h(z_j)\to (a'-d_1)+i(a''-d_2)=a'+ia''-d.$
	Hence $f(z)=h(z_j)+\phi(z_j)\to a'+ia''.$
	This completes the proof.
\end{proof}

\begin{definition}
	Let $\Omega\subset\C$ be a domain bounded by a Jordan curve and let 
	$\gamma$ be a rectifiable arc in $\partial\Omega$. 
	Then there exists a tangent at a.e.\ point on $\gamma$. Let $p$ be such a point. Consider
	the set $\Sigma$ of 
	circular sectors $S$ with vertex $p$ such that
	$S\subset \Omega$ and the boundary radii of $S$ do not belong to the tangent to $\gamma$ at $p.$
	Let $f$ be a function on $\Omega.$ We denote by
	$P(f)$ the set of such points $p$ such that for any choice of $S$ as above
	$\C\cup\{\infty\}\subset C(f,p,S)$. We denote by 
	$F(f)$ the set of such points $p$ such that for any choice of $S$ as above
	$C(f,p,S)$ contains a single point (the same point for any choice of $S$),
	in particular $f$ has angular limit, that is a limit
	along non-tangential paths leading to $p$.
\end{definition}

The definition of $F(f)$ and $P(f)$ can be formulated for a function $f$ on the unit disc
as follows.
If $U$ is a domain in the unit disc, two sides of which consists of distinct
rectilinear segments terminating at a point $\exp(i\theta_0)$, then the
point $\exp(i\theta_0)$ belong to $P(f)$ if
the cluster set $C(f,p,U)$ is total for each angular 
sector $S$ of arbitrarily small angle $v$ between pairs of chords through 
$\exp(i\theta_0)$. A point $\exp(i\theta_0)$ belongs to $F(f)$ if
$\bigcup_{S}C(f,\exp(i\theta_0),S)$ is degenerate and if $\lim f(z)$ exists uniformly as $z\to \exp(i\theta_0)$ in any
angle (sector) $S$ between chords through $\exp(i\theta_0).$
We shall need the following known result.
We shall need the following result of Riesz which we state without proof.
\begin{theorem}\label{rieszlohwater33}
	Let $w=f(z)$ be an analytic function mapping $\{\abs{z}<1\}$ conformally onto a domain 
	$\omega$ bounded by a rectifiable Jordan curve $\gamma$, then under the homeomorphism 
	$w=f(\exp(i\theta))$ of the boundaries induced by $w=f(z),$ any set of measure zero on $\{\abs{z}=1\}$ is mapped onto a set of measure
	zero in $\gamma.$
\end{theorem}

\begin{theorem}[Collingwood \& Lohwater \cite{coolingwoodlohwater},Ch. 6, Thm 8.2]\label{colowat}
	If $f(z)$ is meromorphic on $\{ \abs{z}<1\}$ then a.e.\ point of $\{\abs{z}=1\}\setminus(P(f)\cup F(f))$
	has zero measure.
\end{theorem}
\begin{proof}
	Let $E_1$ denote the set of points $\exp(i\theta)$
	for which $C(f,\exp(i\theta),S)$ is degenerate for all $S$ with vertex $\exp(i\theta)$, i.e.\
	$C(f,\exp(i\theta),S)$ consists of a single point. Note that $\lim_{z\to \exp(i\theta)} f(z)$ exists in any angle at $\exp(i\theta).$
	Let $E_2$ denote the set of points on $\{\abs{z}=1\}$ for which 
	$C(f,\exp(i\theta),S)$ is total for any $S$ with vertex $\exp(i\theta)$ and let $E_3$ denote
	$\{\abs{z}=1\}\setminus(E_1\cup E_2).$
	Then $E_1$ and $E_2$ ae measurable thus $E_3$ is measurable.
	Assume (in order to reach a contradiction) that $E_3$ has positive measure.
	The set of curvelinear triangles on the sphere whose vertices have rational coordinates is countable, denote the set of such triangles by $\tau_1,\tau_2,\ldots.$
	If $\exp(i\theta_0)\notin E_2$ then there is a sector $S$ with vertex $\exp(i\theta_0)$ such that
	$C(f,\exp(i\theta_0,S))$ is not total, in particular there exists a $\tau_k$ such that $\tau_k\cap C(f,\exp(i\theta_0,S)) =\emptyset$.
	Inside $S$ we choose another sector $S_1$ in the form of a circular sector with vertex $\exp(i\theta_0)$
	and radius $1/m$, $m\in \Z_+$ such that $f(z)$ takes no value of $\tau_k$ in $S_1.$
	Choose $S_1$ such that the angular opening of $S_1$ is $\pi/2^n$, $n\in \Z_+$ where the angle bisector of
	$S_1$ at $\exp(i\theta_0)$ forms the angle $p\pi/q$, for integers $p,q$, with the radius of 
	$\{\abs{z}<1\}$ drawn to $\exp(i\theta_0)$. 
	Denote by $E_{n,m,p,q,k}$ the set of all $\exp(i\theta)\in E_3$ for which such a construction
	is possible. Since the five indices are natural numbers this collection is countable. Since each 
	$E_{n,m,p,q,k}$ is defined via means of inequalities, equalities and limiting processes applied to a measurable function, they are all
	measurable.
	If $E_3$ has positive measure there exists at least one member, say
	$E_{n_0,m_0,p_0,q_0,k_0}$ of positive measure that contains a perfect set $E$ of positive measure.
	If $w_0$ is an interior point of $\tau_{k_0}$ there exists $\delta>0$ such that $1/\abs{f(z)-w_0}<\delta$ in any region $S_1$ with vertex $\exp(i\theta)\in E'$.
	Extend the sides of $S_1$ from the point $\exp(i\theta\in E'$ until it either meets a side of another $S_1$ or else meets $\abs{z}=\rho_0>1-1/m_0.$
	The union of $\{\abs{z}<\rho_0\}$ and the interiors of the domains
	$S_1$ forms a simply connected domain, which we denote $D_z$, whose boundary intersects $\{\abs{z}=1\}$ in the set $E'.$
	By the definition of $D_z$ $\partial D_z$ is a rectifiable Jordan curve.
	The function $(f(z)-w_0)^{-1}$ has at most finitely many poles $a_1,\ldots,a_\nu$ in $D_w.$ Since the poles must lie in $\{\abs{z}\leq \rho_0\}$ the function
	$F(z)=(f(z)-w_0)^{-1}\Pi_{j=1}^{\nu}(z-a_j)$ is analytic and bounded in $D_z.$
	By Theorem \ref{rieszlohwater33} $F(z)$ has a unique limit value as $z$ goes to a point of $\partial D_z$
	for a.e.\ point of $\partial D_z.$
	Since $\partial D_z$ has a tangent a.e.\ on $E'$ and since any such tangent must coincide with the tangent to $\{\abs{z}=1\}$, when it exists,
	there exists at least one point $\exp(i\theta_0)$ on $E$ for which $F(z)$ (and therefore $f(z)$) has a unique limit value in any angle of opening $<\pi$
	and bisected by $\theta=\theta_0$.
	This implies that $\exp(i\theta_0)$ belongs to $E_1\setminus E_3$ which is a contradiction.
	This completes the proof.
\end{proof}

\begin{theorem}\label{balktusthm1}
	Suppose $f$ is continuous on a domain $\omega$ bounded by a Jordan curve whose boundary
	contains a rectifiable arc $\gamma.$
	Suppose
	$\partial_{\bar{z}} f$ and $\partial_z f$ exist a.e.\ and
	belong to $L^p(\overline{\Omega}),$ for some $2<p<\infty.$
	Then $\gamma$ can be decomposed into
	\begin{equation}\label{balktusekv}
	\gamma=P(f)\cup F(f)\cup N(f)
	\end{equation}
	where $N(f)$ is a set of measure zero on $\gamma.$
\end{theorem}
\begin{proof}
	In the case of a holomorphic function $h$ on $\Omega$ the result is known,
	see Theorem \ref{colowat}, i.e.\ 
	\begin{equation}
	\gamma=P(h)\cup F(h) \cup N(h)
	\end{equation}
	for a set $N(h)$ of measure zero on $\gamma.$
	Also by Lemma \ref{balktuslem} $f$ can be decomposed as $f=h+\phi$
	for holomorphic $h$ and continuous $\phi.$ Let
	$p\in P(h)$ and set $\phi(p)=:d.$
	Then for each choice of sector $S$ in $\Omega$ with vertex $p$ we have
	\begin{equation}
	C(h,p,S)=\C\cup\{\infty\},\quad C(f,p,S)=C(h,p,S)+d
	\end{equation}
	hence $C(f,p,S)=\C\cup\{\infty\}=\C\cup\{\infty\}.$ Hence $P(h)\subseteq P(f).$ Similarly
	$F(h)\subseteq F(f).$ Hence also $N(f)=N(h),$ and we have Eqn.(\ref{balktusekv}).
	This completes the proof.
	\end{proof}	
	\begin{corollary}\label{balktuscor}
		Suppose $f$ is continuous and bounded on a domain $\Omega$ bounded by a Jordan curve whose boundary
		contains a rectifiable arc $\gamma.$
		Suppose
		$\partial_{\bar{z}} f$ and $\partial_z f$ exist a.e.\ and
		belong to $L^p(\overline{\Omega}),$ for some $2<p<\infty.$ Then $f$ has angular limit a.e.\ on $\gamma.$
		If $f\in H^1(\{\abs{z}<1\})$ such that $\partial_{\bar{z}} f\in L^p(\{\abs{z}\leq 1\}),$ with $2<p<\infty$
		then $f$ has finite angular limit a.e.\ on $\{\abs{z}=1\}.$
	\end{corollary}
	
	\begin{remark}
		Note that Theorem \ref{balktusthm1} and Corollary \ref{balktuscor} applies in particular to
		$q$-analytic functions $f$ for any $q\in \Z_+.$
	\end{remark}
Let us give an alternative proof for the case of polyanalytic functions in domains bounded by a not necessarily rectifiable 
Jordan curve, as given by Balk \& Vasilenkov \cite{balkvasilenkov1992}.
\begin{theorem}[Cf. Theorem \ref{lindelofthm}]\label{balkvasilthm1}
Let $f(z)$ be a bounded $n$-analytic function in a sector $S$ with vertex $\zeta,$ $S$ in a domain $\Omega\subset\C$ bounded by
a not necessarily rectifiable 
Jordan curve, $\Gamma$, containing one or more OSCAR:s. Suppose $\partial_{\bar{z}} f$ is bounded on $S$ and suppose
$f$ tends to a finite limit 
$f(z)\to \lambda$ as $z\to \zeta$ along some continuous path $\gamma$ in $\Omega$. 
Then $f(z)$ tends to $\lambda$ as $z\to \zeta$ interior to any fixed sector $S_1$
inner to $S$.
\end{theorem}
\begin{proof}
Let $S$ be a sector with vertex $\zeta$ whose angle is $\leq \pi$ and let $S_1$ be a sector inner to $S$
of half radius and let $\theta$ be the inclination of $S_1$ to $S$. Suppose $\abs{f(z)}\leq M$ on $S$. Then for any non-negative integers $n,k,m$ $n\geq 1,$ $k+m\geq 1$ we can find a constant $C=C(n,k,m,\theta)$ such that
	\begin{equation}\label{balkvasilek00}
	\abs{z-\zeta}^{k+m}\abs{\partial^m_z\partial_{\bar{z}}^k f}\leq CM,\mbox{ on } S_1
	\end{equation}
	This follows by dropping a perpendicular $[z,t]$ from the point $z\in S_1$ to the nearest to $z$ boundary radius of the sector $S$ and applying 
Theorem \ref{balk204}, Eqn.(\ref{nicolescu}) together with the inequality
\begin{equation}
\abs{z-t}>\abs{z-\zeta}\sin\theta
\end{equation}
Now choose a sector $S_2$with vertex $\zeta$ such that $S_1$ is inner to $S_2$ and $S_2$ is inner to $S.$
W.l.o.g.\ we can assume their radii satisfy $R(S_2)\leq R(S)/2,$ $R(S_1)\leq R(S_2)/2.$ Set $f_1:=\partial_{\bar{z}} f$ which by assumption is bounded in $S$ so there exists $M_1>0$ such that, By Eqn.(\ref{balkvasilek00})
\begin{equation}\label{balkvasilek00}
\abs{z-\zeta}^{k-1}\abs{\partial_{\bar{z}}^{k-1} f}\leq C_{k-1}'M_1,\quad k=1,\ldots,n-1
\end{equation}
where the constants $C_{k-1}'$ are independent of $f,M_1.$ 
Thus for \begin{equation}\lim_{S_2\ni z\to \zeta} \abs{z-\zeta}^{k-1}\abs{\partial_{\bar{z}}^{k-1} f} =0\end{equation} 
We then have for the concordant function
$f_\zeta(z)\to \alpha$ for $\gamma \ni z\to\zeta.$
Also by the same argument using that $f$ is bounded we have that $f_\zeta$ is bounded on $S_2$. By Theorem \ref{lindelofthm} we obtain that
$f_\zeta(z)\to \alpha$ as $z\to \zeta$ interior to $S_2.$ This implies that also $f(z)\to \alpha$ as $S_1\ni z\to \zeta.$
This completes the proof.
\end{proof}
Note that if $f(z)$ is a bounded $n$-analytic function in a domain $\Omega\subset\C$ whose boundary contains
an OSCAR, $\gamma$ and if $\partial_{\bar{z}} f$ is bounded on $\Omega$ then we can locally express
$\gamma$ in terms of a Schwarz function $A(z)$ as $\gamma\cap V=\{A(z)=\bar{z}\}$ for a domain $V$ containing a point $\zeta\in \gamma.$
By Theorem \ref{balkvasilthm1} $f$ has a finite angular limit a.e.\ $\zeta$ on $\gamma.$
By Eqn.(\ref{balkvasilek00}) we obtain for a.e.\ $\zeta\in \gamma,$ that for $k+m\geq 1$
\begin{equation}
(A(z)-\bar{z})^{k+m}\partial^m_z \partial_{\bar{z}}^k f \to 0\mbox{ as }z\to \zeta
\end{equation}
This implies that
for a.e.\ $\zeta\in \gamma$ and any nonnegative integers $k,m$ such that $k+m\geq 1$ we have
\begin{equation}
\partial^m_z \partial_{\bar{z}}^k f=o\left(\frac{1}{(A(z)-\bar{z})^{k+m}}\right)
\end{equation}
as $z\to \zeta$ along any path in $\Omega$ which is not tangential to $\gamma.$
This immediately yields the following corollary using the obvious Schwarz function $A(z)=1/z$, for the unit circle.
\begin{corollary}
	If $f(z)$ is a bounded $n$-analytic function in the unit disc then for
	any nonnegative integers $k,m$ such that $k+m\geq 1$ we have
	\begin{equation}
	\partial^m_z \partial_{\bar{z}}^k f=o\left(\frac{1}{(1-\abs{z}^2)^{k+m}}\right)
	\end{equation}
	\end{corollary}
\begin{theorem}
	\label{balkvasilthm3}
	Let $f(z)$ be a bounded $n$-analytic function in a sector $S$ with vertex $\zeta,$ in the boundary of a domain $\Omega\subset\C$ bounded by
	a not necessarily rectifiable 
	Jordan curve, $\Gamma$. Suppose $\partial_{\bar{z}} f$ is bounded on $S$ and suppose
	$\alpha$ is a regular asymptotic value of order $\geq n$ of $f$ at $\zeta.$ Then 
	$f$ tends to a finite limit 
	$f(z)\to \lambda$ as $z\to \zeta$ along some continuous path $\gamma$ in $\Omega$. 
	Then $f(z)$ tends to $\lambda$ as $z\to \zeta$ interior to any fixed sector $S_1$
	inner to $S$.
\end{theorem}
\begin{proof}
	W.l.o.g.\ assume $\alpha=0,$ $\zeta=0.$ Let $\gamma_0,\ldots,\gamma_{n-1}$ be regular paths along which $f(z)\to \alpha=0$
	as $\gamma_k\ni z\to \zeta,$ $k=0,\ldots,n-1.$
	Let $\gamma_1=\{z:z=\lambda(t),t\in [a,b]\}.$ Let $S_2$ be a sector inner to $S$ such that $S_1$ is inner to $S_2.$
	W.l.o.g.\ we can assume $R(S_1)<R(S_2)/2<R(S)/4.$ Suppose
	$f(z)=\sum_{j=0}^{n-1}h_j(z)\bar{z}^k,$ for holomorphic $h_j$ and suppose $\abs{f}<M,$ for a constant $M>0.$
	Then $\gamma_1$ belongs to an OSCAR $\gamma=\{z:z=\lambda(t),t\in(\alpha',\beta')\}$ for $\alpha'<a<b<\beta'.$
	Let $A(z)$ be the Schwarz function of $\gamma,$ i.e.\ $A(z)=\bar{z}$ defines $\gamma$ for some holomorphic $A(z)$ on a neighborhood of $\gamma.$Let $\delta(0,\rho)$ be a small disc centered at $0$ on which $A(z)$ is holomorphic and bounded and since $z=0$ belongs to $\gamma$ we have $A(0)=0.$ We assume $t$ is the arc length parameter of $\gamma$ from some initial point (with initial parameter value $t_0$), on $\gamma,$ and let $\zeta=$ correspond to $\lambda(t_0)=0.$
	We have for a real parameter $\theta$
	\begin{equation}
	A'(0)\lambda'(t_0)=\overline{\lambda'(t_0)},\quad \abs{A'(0)}=1,\quad A'(0)=\exp(i\theta)
	\end{equation}
	Also $\varphi:=\mbox{arg}\lambda'(t_0)$ is the angle between the arc $\gamma$ and the real axis, measured at $\zeta=0.$ Thus $A'(0)=\overline{\lambda'(t_0)}/\lambda'(t_0)=\exp(-2i\varphi)$ thus modulo $2\pi$ we have $\theta=-2\varphi.$ By Proposition \ref{concordantprop} the concordant function $f_\gamma$ satisfies
	 	\begin{equation}
	 \abs{f_\gamma}\leq \abs{f(z)}+\sum_{j=1}^{n-1}\frac{1}{j!}\abs{\bar{z}-A(z)}^j\abs{\partial_{\bar{z}}^j f}
	 \end{equation}
In $\delta(0,\rho),$ $A(z)$ can be written for $z\to 0$ as
		\begin{equation}
	A(z)=A(0)+A'(0)z+z\cdot o(z)=z(\exp(i\theta)+o(z)))
	\end{equation}
	Thus on $S_2\cap \delta(0,\rho)$ we have for sufficiently small $\rho>0$ that there is a constant $C$ independent of $f,M$ such that
	\begin{equation}
\abs{\bar{z}-A(z)}\leq \abs{z}+\abs{A(z)}\leq \abs{z}+\abs{z}\abs{\exp(i\theta)+o(z)}<3\abs{z}
	\end{equation}
\begin{equation}
\abs{f_\gamma}= M+\sum_{j=1}^{n-1}\abs{\partial_{\bar{z}}^j f}3^j \abs{z}^j\leq MC
\end{equation}
	Since $f_\gamma(z)=f(z)$ on $\gamma$ we have $f_\gamma(z)\to 0$ as $\gamma\ni z\to 0.$ By Theorem \ref{lindelofthm} this implies $f_\gamma(z)\to 0$ as $S_1\ni z\to 0$ for $S_1$ inner to $S_2$ thus	
		\begin{equation}
	\abs{f_\gamma}\leq \sum_{j=0}^{n-1}h_j(z)z^j\exp(i\theta)+o(\abs{z}),\quad S_1\ni z\to 0
	\end{equation}
	Repeating this for $\gamma_k,$ $k=0,\ldots,n-1$ and denoting in each case the OSCAR $\gamma=\gamma^{(k)}$ and the concordant function $f_k:=f_{\gamma_k}$ we have, for mutually distinct numbers $\theta_0,\ldots,\theta_{n-1}$, between $0$ and $2\pi,$ a system
	\begin{equation}
	f_k(z)=\sum_{j=0}^{n-1}h_j(z)z^j\exp(i\theta_k)+o(\abs{z}),\quad S_1\ni z\to 0,\quad k=0,\ldots,n-1
	\end{equation}
	and furthermore $f_k(z)\to 0$ as $S_1\ni z\to 0$. 
	The determinant, $V(z)$, of this system is a Vandermonde determinant and we have as $z\to 0$ that
	$V(z)\to$\\
	$V(\exp(i\theta_0),\ldots,\exp(i\theta_{n-1}))\neq 0.$ Hence for sufficiently small $z$ we obtain from the system $n$ functions $h_j(z)z^j,$ $j=0,\ldots,n-1,$. Since $f_k(z)\to 0$ as $S_1\ni z\to 0$, $k=0,\ldots,n-1$, this holds true for all the functions $h_j(z)z^j$, hence also for $f(z).$ This proves theorem \ref{balkvasilthm3}	
\end{proof}

\begin{theorem}
	\label{balkvasilthm6}
	Let $f(z)$ be a bounded $n$-analytic function in a domain $\Omega\subset\C$ bounded by
	a not necessarily rectifiable 
	Jordan curve, whose boundary contains $n$ pairwise non-continuous OSCAR:s $\Gamma_1,\dots,\Gamma_n.$\\
	(1) Let $\Gamma_k$ contain a set $M_k$ of second category and let $G$ contain a uniformly continuous family $L_{1,k}$ of paths $L_{1,k}^\zeta$
	conducting from $\Omega$ to the points $\zeta$ of $M_k$, where $M_k$ and $L_{1,k}$ satisfy that for any $\zeta\in M_k$ 
	\begin{equation}
	\infty\notin C(f,\zeta,L_{1,k}^\zeta),\quad \infty\notin C(\partial_{\bar{z}}f,\zeta,L_{1,k}^\zeta),\quad k=1,\ldots,n
	\end{equation}
	(2) Let $\Gamma_k$ contain a (relatively) metrically dense set $E_k$ and let $\Omega$ contain a family $L_{2,k}$ of paths $L_{2,k}^\zeta$,
	nontangential to $\partial\Omega$, conducting from $\Omega$ to the points $\zeta$ of $E_k$ such that for each $\zeta\in E_k$
	\begin{equation}
	0\in C(f,\zeta,L_{2,k}^\zeta),\quad\quad k=1,\ldots,n
	\end{equation}
	Then $f(z)\equiv 0.$
\end{theorem}
\begin{proof}
	Set $\phi(z)0\abs{f(z)}+\abs{\partial_{\bar{z}}f(z)}.$ Then $\phi$ maps continuously the domain $\Omega$
	into a Riemannian sphere and it is a known result of Dolzhenko \cite{dolchenko267} that
		the set of all points $\zeta\in M_k$ such that $C(\varphi,\zeta,L_{1,k}^\zeta)\neq C(\varphi,\zeta,\Omega)$ is a set of second category. But $M_k$ itself is a set of second category thus there exists a point $\zeta_0$
	such that $C(\varphi,\zeta_0,L_{1,k}^{\zeta_0})\neq C(\varphi,\zeta_0,\Omega)$, hence by (1) the set 
	$C(\varphi,\zeta_0,\Omega)$ is bounded by a constant, say $M.$
	This implies (which can be realized by contradiction) that there exists a disc $K(\zeta_0,\rho)$, centered at $\zeta_0$ with radius $\rho>0$
	such that $\varphi(z)<M$ for all $z\in \Omega\cap K(\zeta_0,\rho).$
	We then have that $\abs{f(z)}<M$ and $\abs{\partial_{\bar{z}} f}<M$ on $\Omega\cap K(\zeta_0,\rho).$
	Let the OSCAR $\Gamma_k$ have Schwarz function $A_k(z)$ on a neighborhood $\Delta_k$ of
	$\Gamma_k.$ Let $F_k(z)$ be the concordant function with $f$ on $\Gamma_k$,
	\begin{equation}
	F_k(z)=\sum_{j=0}^{n-1} h_j(z)A_k^j(z)
	\end{equation}
	Since $f$ has angular limit zero for each $\zeta\in E_k$ the same holds true for $F_k(z).$ Since $E_k$ has positive relative measure we have (by the Luzin-Privalov theorem) $F_k(z)\equiv 0$ on $\Omega\cap \Delta_k.$
	
	This implies that there exists a discrete set $d_k$ in $\Omega$ such that$A_k(z)$ can be analytically continued from $\Omega\cap \Delta_k$ to any simply connected domain $D$ belonging to $\Omega\setminus d_k.$ Since this is true for any $k=1,\ldots,n$, we have, setting $d=\bigcup_{k=1}^n d_k$, that each $A_k(z)$, $k=1,\ldots,n$ can be analytically continued to any simply connected domain $D\subset\Omega\setminus d.$
	Since the functions $A_j(z),$ $j=1,\ldots,n$ are mutually non-continguous and the
	system
	\begin{equation}
	\sum_{j=0}^{n-1} h_j(z)A_k^j(z)\equiv 0,\quad k=1,\ldots , n
	\end{equation}
	has corresponding matrix that is a Vandermonde matrix $V(A_1,\ldots,A_n)$, thus 
	all $h_j(z)\equiv 0$ on $D$ and by connectedness on all of $\Omega$, $j=0,\ldots,n-1.$
	Hence $f\equiv 0$. This completes the proof of Theorem \ref{balkvasilthm6}
\end{proof}

\section{Picard type theorems}

A direct consequence of Theorem \ref{fundamentalthmbalk} is the following.
\begin{proposition}
	Suppose $P(z,\bar{z})$ is a polynomial
	such that the exact degree of $P(v,w)$ as a joint polynomial in complex variables $v,w$ is $n\in \N$
	such that $n> 2\max\{n_1,n_2\}$ where $n_1$ is the exact degree of $P$ with respect to $v$ and 
	$n_2$ is the exact degree of $P$ with respect to $w$. 
	Then $P$ takes all complex values on $\C.$
\end{proposition}

On the other hand let $n_2\leq n\leq 2n_2.$ Then for example the polynomial $P(z):=(z+\bar{z})^{2n_2-n}(z\bar{z})^{n-n_2}$
omits all values with nonzero imaginary part.
\begin{theorem}\label{transcendental001}
	An entire transcendental polyanalytic function, $f(z)=\sum_{j=0}^{n} a_j(z)\bar{z}^j$ for entire $a_j(z)$, takes all complex values except at most one exception.
\end{theorem}
\begin{proof}
	Suppose $f$ assumes two values (w.l.o.g.\ we can assume these to be $0$ and $1$) only on a bounded set $M$.
	The theorem follows if we show that $f$ must be a polynomial.
	Now there exists $r_0$ such that $M\subset\{\abs{z}<r_0\}.$ this implies that there exists constants $p_1,q_1$ such that for $t\geq r_0$, $\lambda>2$ 
	\begin{equation}
	\Psi f(z):=\frac{1}{2\pi}\int_{\{\abs{z}=\lambda t\}} d\mbox{arg} f(z)=p_1, \Psi(f(z)-1)=q_1
	\end{equation}
	Now setting
	\begin{equation}
	\Phi(z,\lambda t):=\sum_{\nu}^n \lambda^{2\nu}t^{2\nu} \frac{\phi_\nu(z,t)}{(z-t)^\nu}
\end{equation}
where $\phi_\nu(z)=\sum_{k=\nu}^n C_{k,\nu} a_k(z)t^{k-\nu},$ $\nu=0,\ldots,n,$
we have constants $C_{k,\nu}$ such that
\begin{equation}
\Phi(z,\lambda t)=f(z),\quad z\in \{\abs{z}=\lambda t\}
\end{equation}
This implies that $\Psi \Phi(z,\lambda t)=p_1$, $\Psi (\Phi(z,\lambda t)-1)=q_1.$ This shows that $\Phi$
assumes the value $0$ ($1$) on $\{\abs{z-t}<\lambda t\}$, $p$ ($q$) times, where
$p=p_1+n,$ $q=q_1+n.$ Furthermore, the function
$\phi:=\Phi(tz,\lambda t)$ takes the value $0$ ($1$) not more than $p$ ($q$) times
on $\{\abs{z}<1\}.$ Set $\lambda=t^s,$ $s:=2n+p$ and $\phi(z):=\sum_{k=0}^\infty b_k z^k.$
Then $b_k=o(\lambda^{3n})$ as $t\to \infty,$ $k=0,\ldots,p.$
We shall need the following lemma which we state without proof
\begin{lemma}[Saxer \cite{saxer}, p.215]\label{saxerlemma}
	If the function $\phi(z)=\sum_{j=0}^\infty b_j z^j$ is regular in $\{\abs{z}<1\}.$ and assumes the values $0$ ($1$)
	not more than $p$ ($q$) times, and if $\abs{b_k}<K_k,$
	$k=0,\ldots,p,$ $K_0>e,$ then for $\abs{z}<1-\theta$ ($\theta\in (0,1)$)
	there is a constant $C>0$ independent of $\theta,\phi(z)$, $K_0,\ldots,K_p,$ such that
	\begin{equation}\label{saxerek4}
	\abs{\phi(z)}<(\sum_{k=0}^p K_k)^{\frac{C}{\theta}}
	\end{equation}
\end{lemma}
Setting in Lemma \ref{saxerlemma}, $\theta=1/2$ we obtain by Eqn.(\ref{saxerek4}) that there exists an $m$ such that
$\phi(z)=o(t^m)$ as $t\to \infty$ for $\abs{z}\leq 1/2.$ This implies that $\psi(z):=(z-1)^n\phi(z)$ also satisfies $\psi(z)=o(t^m)$
for $\abs{z}\leq 1/2$ as $t\to \infty.$
Let $a_\nu(z):=\sum_{\mu=0}^\infty c_{\mu,\nu} z^\mu,$ $\nu=0,\ldots,n.$ By the Cauchy inequalities for $\psi$
we have for each $\mu$, as $t\to \infty$ 
\begin{equation}
\lambda^{2n}t^{n+\mu}\abs{o(1)+c_{\mu,n}}=o(t^m)
\end{equation}
This implies that $c_{\mu,n}=0$ for all sufficiently large $\mu$, hence $a_n(z)$ is a polynomial. Similarly, we obtain by iteration that
$a_\nu(z)$ is a polynomial for $\nu=n-1,\ldots,0.$ This completes the proof.
\end{proof}
\begin{corollary}
Let $P(x,y)$ and $Q(x,y)$ be arbitrary polynomials in the real variables $x,y$ such that 
$P(x,y)\not\equiv 0$. Let $f(z)$ be an entire analytic function with respect to the complex variable $z=x+iy.$ If 
$P(x,y)\cdot f(z)+Q(x,y)$ assumes both of the two complex numbers $a$ and $b$ ($a\neq b$) on bounded sets,
then $f(z)$ is a polynomial.
\end{corollary}
\begin{proof}
	This an immediate consequence of Theorem \ref{transcendental001} together with
	the fact that $P(x,y)\cdot f(z)+Q(x,y)$ is an entire polyanalytic function.
\end{proof}
\begin{definition}
	Let $\Omega\subset\C$ be a domain and let $f(z)$ be a polyanalytic function
	on $\Omega$.
	A point $p\in \Omega$ is called an {\em isolated singularity}\index{Isolated singularity} of a polyanalytic function 
	$g$ on some subset $U\subset\Omega$ if there is a punctured neighborhood $V$ of $p$
	such that $V\subset U.$
	A point in $\Omega$ is called an {\em essential singularity}\index{Essential singularity}
	of $f(z)$ if it is an essential singularity of at least one of the analytic components of $f(z).$
	For a point $p\in \Omega$ and a punctured
	neighborhood $U$ of $p$, we denote $M(a,U):=(U\setminus \{p\})\cap f^{-1}(a)$.
	A function $a(z)$ is meromorphic at a point $p$ if it is holomorphic in some punctured neighborhood of
	$p$ and has no essential singularity at $p.$
	Note that on a domain $\Omega$ a polyanalytic function $f$ has a unique representation
	$f(z)=\sum_{j=0}^{q-1} a_j(z)\bar{z}^j$ for holomorphic $a_j(z),$ $j=0,\ldots,q-1,$ which implies that if $p$ is an isolated singularity
	of $f$ then it is an isolated singularity for each $a_j$.
	A point $p$ is called a {\em nonessential singularity}\index{Nonessential singularity point}
of $f$ if it is a nonessential isolated singularity of each $a_j(z).$
\end{definition}
Let $f$ be a polyanalytic function with an isolated singularity at a point $p\in \C.$ Recall that a number $a\in\C$
is called an exceptional value for $f$ and $p$ if there exists a neighborhood $U$ of $p$
such that $f(z)-a$ never vanishes on $U\setminus \{p\}.$

\begin{definition}
A polyanalytic function $f(z)$ on a domain is called meromorphic at a point $p$
if all of its analytic components are meromorphic at $p.$ $f(z)$ is called
reducible on a domain if it can be expressed in that domain as a product of two polyanalytic 
nonholomorphic functions. Otherwise $f$ is called {\em irreducible}\index{Irreducible polyanalytic function} on the given domain.
\end{definition}
\begin{definition}
For a continuous function $f(z)$ and a real positive number $\rho>0$ such that
$f(z)$ is not zero on $\{\abs{z}=\rho\}$ we denote by $\Delta_\rho f$, $1/2\pi$ times the change in argument of $f$ around the positively oriented circumference $\{\abs{z}=\rho\}$.	
\end{definition}

\begin{theorem}\label{boschkraj70thm}
Let $f(z),g(z),h(z)$ be polyanalytic functions
on an annular neighborhood $A$ of a point $z_0\in \C.$, finite or infinite,
such that $g$ and $h$ do not have an essential singularity at $z_0$ and $g-h$ is not identically zero in $A.$ 
If $f-g$ and $f-h$ never vanish in $A$ then $z_0$ is not an essential singularity of $f.$
\end{theorem}
\begin{proof}
W.l.o.g.\ we can assume $h\equiv 0.$ First assume $z_0=\infty.$ Let $A=\{R<\abs{z}<\infty\},$ some $R>0.$ Since $f$ and $f-g$ 
have no zeros on in $A$ there exists integers $r,s$ such that $\Delta_\rho f=r$
and $\Delta_\rho (f-g)=s$ for all $\rho>R.$ 
For a function $h$ on $\{\abs{z}\geq R\}$ that is an $n+1$-analytic function $h(z)=\sum_{j=0}^{n}b_j(z)\bar{z}^j$ on $A$, denote by $h(z,\rho):=\rho^{2j} b_j(z)/z^j$ for $\abs{z}>R.$
Then $h(z,\rho)=h(z)$ on $\{\abs{z}= R\}$ and holomorphic on $\{\abs{z}> R\}$.
Let $\delta,\sigma>0$ and $u,v$ be integers such that
$\delta>\sigma >R$ and $\Delta_\sigma f(z,\rho)=u,$
$\Delta_\sigma (f(z,\rho)-g(z,\rho))=v$ for all $\rho>\delta.$ This implies that
the functions $f(z,\rho)$ and $f(z,\rho)-g(z,\rho)$ respectively have
$r-u$ zeros and $s-v$ zeros respectively on $\{\sigma<\abs{z}<\rho\}$ for all $\rho>\delta.$
Hence $f(\rho z,\rho)$ and $f(\rho z,\rho)-g(\rho z,\rho)$ respectively have
at most $r-u$ zeros and $s-v$ zeros respectively on $\{\sigma/\delta <\abs{z}< 1\}$ for all $\rho>\delta.$
\begin{lemma}\label{boschkraj70lem3}
Let $f(z)\not\equiv 0$ be polyanalytic in $\{R<\abs{z}<\infty\}.$ Suppose $z_0$ is not an essential singularity of $f.$ Then there exists a finite subset $F$
of $\{0<\abs{z}<\infty\}$ such that if $A$ is any closed bounded subset of $\{0<\abs{z}<\infty\}$ which has empty intersection with $F$
then there exists positive constants $K,L$ and a nonnegative integer $t$ such that $L\rho^t\leq\abs{f(\rho z,\rho)}\leq K\rho^t$
for all $z\in A$ and all sufficiently large $\rho.$ 
\end{lemma}
\begin{proof}
There exists nonnegative integers $s,t$ and a polynomial $P(z)\not\equiv 0$ such that
$f(\rho z,\rho)/\rho^t \to P(z)/z^s$ uniformly for $z$ in a closed and bounded subset of $\{0<\abs{z}<\infty\}$
as $\rho\to +\infty.$ The conclusion of the 
lemma follows by letting $F$ denote the set of nonzero roots of $P(z).$
This proves Lemma \ref{boschkraj70lem3}.
\end{proof}
By Lemma \ref{boschkraj70lem3} there exists some $\sigma/\delta <a<b<1$ and some $\mu>\delta$ such that $L\rho^t\leq \abs{g(\rho z,\rho)}\leq K\rho^t$
for $a<\abs{z}<b$ and $\rho>\mu$ where $L,K$ are positive constants and $t\in \N.$
Let $\{\rho_j\}_{j\in \N},$ $\rho_j>\mu$, be a sequence diverging to $+\infty.$ For $j\geq 1$
define on $\{ z:a<\abs{z}<b\}$, 
$H_j(z):=f(\rho z,\rho_j)/g(\rho_j z,\rho_j)$, in particular each $H_j$ is analytic in $\{ z:a<\abs{z}<b\}$.
Also each $H_j$ takes the value $0$ and $1$ at most $r-u$ and $s-v$ times respectively in $\{ z:a<\abs{z}<b\}$.
This implies that the family of the functions $H_j$ is quasinormal of order $q$
on $\{ z:a<\abs{z}<b\}$ where $0\leq q\leq \max\{r-u,s-v\}$ (see Proposition \ref{montelpropenquasinormal} or e.g.\ Montel \cite{montel}, p.57).
First consider the case $q\geq 1.$ Then there exists a subsequence $\{H_{j_v}\}_{v\in \Z}$ of $\{H_j\}_{j\in \Z_+}$
and $q$ points $z_1,\ldots,z_q$ in $\{a<\abs{z}<b\}$ such that
$H_{j_v}(z)\to\infty$ almost uniformly on $\{a<\abs{z}<b\}\setminus\{z_1,\ldots,z_q\}$ as $v\to \infty.$
We can then choose $\lambda\in (a,b)$ such that $H_{j_v}(z)\to \infty$ uniformly on $\abs{z}=\lambda$ as $v\to \infty.$
Hence there exist a an integer $v_0$ such that $\abs{H_{j_v}(z)}\geq 1$ for $\abs{z}=\lambda$ and $v\geq v_0.$
Thus
$\abs{f(\rho z,\rho)}\geq \abs{g(\rho z,\rho)}\geq L\rho^t \geq L\mu^t$ for $\abs{z}=\lambda$, $\rho=\rho_j$ and $v\geq v_0.$
Also $\Delta_{\lambda_\rho} f(z,\rho)\leq \Delta_\rho f(z,\rho)=r$ for $\rho=\rho_j$ and
$v\geq v_0.$

\begin{lemma}\label{boschkraj70lem2}
	Suppose $f$ is polyanalytic on $\{R<\abs{z}<\infty\}$, for some $R>0$,
	and that there exists some $\lambda\in (0,1]$, a positive constant $K,$
	a nonnegative integer $s$, and some sequence $\{\rho_j\}_{j\in \N}$, $\rho_j>R/\lambda$ 
	diverging
	to $+\infty$ such that $\abs{f(\rho z,\rho)}> K$ for $\abs{z}=\lambda$ and $\rho=\rho_j.$
	If the sequence of integers $\{\Delta_{\lambda_{\rho_j}} f(z,\rho)\}_{j\in \N},$
	is bounded then $z_0=\infty$ is not an essential singularity of $f.$
\end{lemma}
\begin{proof}
	Suppose $f$ is polyentire and let $f(z)=\sum_{j=0}^n f_j(z)\bar{z}^j$. Let $g(z)=z^n f(z)$ for $z\in \C$,
	in particular
	$g(z,\rho)=z^n f(z,\rho)$ for $\rho\geq 0,$ $\abs{z}>0.$ Then $g$
	satisfies the conditions of the lemma. W.l.o.g.\ we can assume
	$\abs{g(\rho z,\rho)}\geq 1$ for $\abs{z}=\lambda$ and $\rho=\rho_j$.
	The function $g(z,\rho)$ is entire for all $\rho\geq 0.$ 
	There exists an integer $p$ such that
	$\Delta_{\lambda_\rho} g(z,\rho)\leq p$ for $\rho=\rho_j$.
	By the Jensen formula we have
	\begin{multline}\label{boschkraj70ekv2}
	\frac{1}{2\pi} \int_0^{2\pi} \log\abs{ g(\lambda\rho \exp(i\theta),\rho)}d\theta =\\
	\sum_{j=1}^q \log \frac{\lambda\rho}{\abs{a_j}}+\log \frac{\abs{g^{(s)}(0,\rho)}}{s!} + s\log \lambda\rho
		\end{multline}
		for $\rho=\rho_j,$ where $a_1,\ldots,a_q$ are the roots of $g(z,\rho)$ in 
		$0<\abs{z}<\lambda \rho$ with multiplicity
		and $s\geq 0$ is the multiplicity of the root $z=0$ of $g(z,\rho).$
		Now there exists some $R_0>R/\lambda$ such that $s$ is independent of $\rho$ for $\rho>R_0.$
		Thus $g^{(s)}(0,\rho)$ is a fixed polynomial in $\rho$ for $\rho>R_0.$ 
		Furthermore, $q=\Delta_{\lambda_\rho} g(z,\rho)-s\leq p-s$ for $\rho=\rho_j>R_0.$
		Using elementary estimates on $g(z,\rho)$ together with the fact that $z=0$ is a root of multiplicity $s$ of $g(z,\rho)$ for $\rho>R_0$,
		we obtain that there exists some $\sigma>0$, a nonnegative integer, $t$, and
		some $R_1>R_1+1$ such that $\abs{a_1},\ldots,\abs{a_q} >\sigma/\rho^t$ for $\rho>R_1.$
		By Eqn.(\ref{boschkraj70ekv2}) there thus exists a positive integer $K$ together with a 
		positive integer $\mu$ such that for all $\rho=\rho_j>R_1$
		\begin{equation}\label{boschkraj70ekv3}
		\frac{1}{2\pi} \int_0^{2\pi} \log \abs{g(\lambda\rho\exp(i\theta),\rho)}d\theta \leq \log K\rho^\mu
		\end{equation}
		For $z=r\exp(i\theta))$ where $0\leq r <\lambda\rho$ we have by the Poisson-Jensen formula
		\begin{multline}\label{boschkraj70ekv4}
		\log\abs{g(z,\rho)}=\frac{1}{2\pi}\int_0^{2\pi} \log \abs{g(\lambda\rho \exp(i\theta),\rho)}
		\frac{\lambda^2\rho^2 -r^2}{\lambda^2\rho^2 -2\lambda \rho r\cos (\theta -\theta)+r^2}d\theta\\
		+\sum_{l=1}^q \log \abs{\frac{\lambda\rho(z-a_l)}{\lambda^2 \rho^2-\bar{a}_l z}} +s\log \frac{r}{\lambda\rho}
		\end{multline}
		for $\rho=\rho_j>R_1$ whenever $g(z,\rho)\neq 0$. For $\abs{z}=\lambda \rho$ and $\rho=\rho_j.$
		$\abs{g(z,\rho)}\geq 1$. By Eqn.(\ref{boschkraj70ekv4}) we thus have for $\abs{z}=r<\lambda\rho$ and $\rho=\rho_j>R_1$
		and whenever $g(z,\rho)\neq 0$
		\begin{equation}\label{boschkraj70ekv5}
		\log\abs{g(z,\rho)}\leq \frac{\lambda \rho +r}{\lambda \rho -r}\frac{1}{2\pi} \int_0^{2\pi} \log \abs{g(\lambda \rho\exp(i\theta),\rho)}d\theta
		\end{equation}
		By Eqn.(\ref{boschkraj70ekv3}) and Eqn.(\ref{boschkraj70ekv5}) with $r=\lambda\rho/2$, we have that there exists a constant $M>0$ and $\nu\in \Z_+$ such that
		$\abs{g(\rho z,\rho)}\leq M\rho^\nu$ for $\abs{z}=\lambda/2$ and $\rho=\rho_j >R_1.$
		By Lemma \ref{boschkraj70lem1} $z_0$ is not an essential singularity of $f.$
		This proves Lemma \ref{boschkraj70lem2} in the case that $f$ is polyentire.
		Now suppose $f$ is polyanalytic in $\{0<\abs{z}<\infty\}$ and $z_0=0$ is not an essential singularity of $f.$ Then
		there exists $s\in \N$ and a polyentire function $g$ such that $g(z)=z^s f(z)$ for $0<\abs{z}<\infty.$
		By what we have done $z_0=\infty$ is not an essential singularity of $g$ and thus not an essential singularity of $f.$ 
		Next suppose that
		$f$ is polyanalytic in $\{R<\abs{z}<\infty\}$.
		Let $f_j(z)=\sum a_{k,j} z^k$ for $R<\abs{z}<\infty$ be the Laurent expansion at $z_0=\infty.$
		$g_j(z)=\sum_{k=-n}^\infty a_{k,j} z^k$ for $0<\abs{z}<\infty$ and let $g(z)=g_j(z)\bar{z}^j$ for 
		$0<\abs{z}<\infty$. Then $g$ is polyanalytic on $0<\abs{z}<\infty$ and $z_0=0$ is not an essential singularity
		of $g.$ Now $f(\rho z,\rho)-g(\rho z,\rho)\to 0$ uniformly for $\abs{z}=\lambda$ as $\rho\to \infty.$ Hence
		$z_0=\infty$ is not an essential singularity of $g$ and thus $z_0=\infty$ is not an essential singularity of $f.$
		This completes the proof of Lemma \ref{boschkraj70lem2}.	
	\end{proof} 

By Lemma \ref{boschkraj70lem2} $z_0$ is not an essential singularity of $f.$ Next consider the case $q=0.$
Then the family of functions $H_j$, defined above 
on $\{ z:a<\abs{z}<b\}$, 
$H_j(z):=f(\rho z,\rho_j)/g(\rho_j z,\rho_j)$, is normal in $\{ z:a<\abs{z}<b\}.$
If there is a subsequence $H_j$ which diverges to $\infty$ almost uniformly on $\{ z:a<\abs{z}<b\}$
then we can as above verify that $z_0=\infty$ is not an essential singularity of $f$.
Suppose there exists a subsequence $H_{j_v}$ of $H_j$ which converges
almost uniformly on $\{ z:a<\abs{z}<b\}$.
Let $\lambda\in (a,b).$
So there exists a constant $M>0$ and an integer $v_0$ such that $\abs{H_{j_v}(z)}\leq M$
Hence $\abs{f(\rho z,\rho)}\leq M\abs{g(\rho z,\rho)}\leq MK\rho^t$ 
for $\abs{z}=\lambda$ and $v\geq v_0.$
for $\abs{z}=\lambda$ and $\rho=\rho_{j_v}$ and $v\geq v_0.$

\begin{lemma}\label{boschkraj70lem1}
	Suppose $f$ is polyanalytic on $\{R<\abs{z}<\infty\}$, for some $R>0$,
	and that there exists some $\lambda>0$, a positve constant $K,$
	a nonnegative integer $s$, and some sequence $\{\rho_j\}_{j\in \N}$, $\rho_j>R/\lambda$ 
	diverging
	to $+\infty$ such that $\abs{f(\rho z,\rho)}\leq K\rho^s$ for $\abs{z}=\lambda$ and $\rho=\rho_j.$
	Then $z_0=\infty$ is not an essential singularity of $f.$
\end{lemma}
\begin{proof}
	Let $f(z)=\sum_{j=0}^n f_j(z)\bar{z}^j$ where 
	$f_j(z)=\sum a_{j,k} z^k$, is the Laurent expansion of $f_j$ at $z_0=\infty.$ Then the coefficient $b_k$ of $z^k$ in the Laurent  expansion of $f(z,\rho)$ at $z_0=\infty$ is
	$b_k=\sum \rho^{2k} a_{j,k+j}$ By the Cauchy inequalities we have
	\begin{equation}
	\abs{\sum \rho^{2k} a_{j,k+j}}\leq K/\lambda^k
	\end{equation}
	which holds true for $\rho=\rho_j.$ This implies that the coefficients of positive powers of $\rho$
	are zero so $a_{j,l}=0$ for $l>s-j$. Thus $z_0=\infty$ is not an essential singularity of $f_j$. This completes the proof of Lemma
	\ref{boschkraj70lem1}.
\end{proof} 
By Lemma \ref{boschkraj70lem1} $z_0=\infty$ is not an essential singularity of $f.$
This proves Theorem \ref{boschkraj70thm} for the case when $z_0=\infty$.
By inversion with respect to the center $z_0$ the case when $z_0$ is finite follows from the case when $z_0$ is not finite.
This proves Theorem \ref{boschkraj70thm}.
\end{proof}
\begin{corollary}
Let $\Omega\subset\C$ be a domain and let $f(z)$ be a polyanalytic function
on $\Omega$. Let $p\in \Omega$ be an isolated singularity of $f(z).$ If $f$ omits two distinct finite values in some annular neighborhood of $p$, then $p$ is not an essential isolated singularity of $f.$
\end{corollary}
Krajkiewicz \cite{boschkraj1974} also proves the following.
\begin{theorem}
	Let $f$ be polyanalytic on $.$
	with the representation $f(z)=\sum_{j=0}^n f_j(z)\bar{z}^j,$ where $f_n\equiv 1$ on the set, $A(\infty,R)$, of
	all finite complex numbers $z$ such that $1/R\leq \abs{z}<+\infty.$
	If $f$ admits the exceptional value $0$ at the point $\infty$ then the
	functions
	$f_0,\ldots,f_{n-1}$ have nonessential isolated singularity at $\infty.$ 
\end{theorem}

\chapter{Boundary value problems for polyanalytic functions}\label{bvpsec}
\section{Background}\label{bvpintro}
The theory of boundary value problems is possibly the most natural branch of partial differential equations and is accordingly deep and multi-faceted.
Riemann \cite{riemann} and Hilbert \cite{hilbert} (see e.g.\ the exposition, in terms of historical notes, given in Gakhov \cite{gakhov}, Ch. II, parag.19) initiated the analysis of some particularly natural settings in complex analysis of one variable.
Let us begin by considering the so-called Riemann boundary value problem (see e.g.\
Gakhov \cite{gakhov}, sec. 14.3).
\begin{definition}\label{riemanproblemdef}
	Given a simple closed contour $\gamma$ enclosing the domain $D^+$ and setting $D^-=\hat{\C}\setminus (D^+\cup \gamma)$, 
	and two H\"older continuous functions $h,g$ on $\gamma$, $h\neq 0$,
	the so-called {\em Riemann boundary value problem associated to $\gamma$}\index{Riemann boundary value problem} is that of finding two functions
	$f^\pm\in \mathscr{O}(D^\pm),$ including $z=\infty,$
	satisfying on the boundary
	$f^+=h f^-$ (This is the homogeneous case of the problem). Replacing the boundary condition by 
	$f^+=h f^- +g$ is the inhomogeneous version of the problem.
	The {\em index}\index{Index of the Riemann problem}
	of the problem 
	denoted by $\kappa$, is defined as the integer $\frac{1}{2\pi}[\mbox{arg} g(t)]_\gamma:=\frac{1}{2\pi}\int_{\gamma} d\arg g(x)$ (i.e.\ the change in argument
	of $g$ when traversing $\gamma$ in the positive direction, passing every point once) which is known to be the sum of the number of zeros of $f^\pm$
	in $D^\pm$. 
	More generally, we may consider the case of a domain $D^+$ bounded by smooth contours $\gamma_0,\ldots,\gamma_p$, 
	mutually non-intersecting such that $\gamma_0$ encloses all others (in the case of an unbounded $D^+$, $\gamma_0$ may be absent), and denote by $\gamma$
	the union $\gamma0\cup_{j=0}^p \gamma_j,$ (the positive orientation on $\gamma$ will be such that $D^+$ lies to the left when $\gamma$ is described
	in that direction). $D^-$ will again denote $D^-=\hat{\C}\setminus (D^+\cup \gamma)$. We may then consider
	the same boundary value problem for such more general domains, and we use the same name.
\end{definition}
Taking logarithms in 
\begin{equation}\label{muskesh351inhomo}
f^+(t)=h(t) f^-(t) +g(t) 
\end{equation}
with $g\equiv 0$ yields
\begin{equation}
\ln f^+(t)-\ln f^-(t)= \ln h(t) 
\end{equation}
As $t$ moves in the positive direction along $\gamma_j$ (i.e.\ counter-clockwise for $\gamma_0$ and clockwise for the remaining
$\gamma_j$, $j=1,\ldots,p$) $\log h(t)$ increases by multiples of $2\pi i$ so that
\begin{equation}
\frac{1}{2\pi i}[\ln h(t)]_{\gamma_j}=\frac{1}{2\pi}[\mbox{arg} h(t)]_{\gamma_j}=:\lambda_j,\quad j=0,\ldots,p
\end{equation}
where $\lambda_j$ are integers.
We define the index $\kappa$ of the boundary value problem (or of $h(t)$ with respect to $\gamma$) in this more general case to be the sum
\begin{equation}
\kappa=\sum_{j=0}^n\lambda_j =\frac{1}{2\pi i}[\ln h(t)]_\gamma=\frac{1}{2\pi}[\mbox{arg}h(t)]_\gamma
\end{equation}
Let $a_1,\ldots,a_p$ be fixed points in the corresponding regions $D^+_1,\ldots,D_p^+$ and set 
$\Pi(z)=\Pi_{j=1}^p (z-a_j)^\lambda_j$. Then $h_0(t):=t^{-\kappa}\Pi(t)h(t)$ will return to its initial value after traversing any circuit of the contours $\gamma_j,$
thus
the function $\ln h_0(t)$ is a one-valued, continuous and H\"older continuous function on $\gamma,$
on each of the $\gamma_j$ we fix a branch of this function.
If $f(z)$ is a solution to the homogeneous Riemann boundary value problem then $\Pi(z)f(z)$ is holomorphic on $D^+$ and $z^\kappa f(z)$ is holomorphic on $D^-$ (except possibly at $\infty$).
Set
\begin{equation}
\psi(z):=\left\{
\begin{array}{ll}
\Pi(z)f(z) & , z\in D^+\\
z^{\kappa}f(z) &, z\in D^-
\end{array}
\right.
\end{equation}
and multiplying Eqn.(\ref{muskesh351inhomo}) with $g\equiv 0$, by $\Pi(z)$
we transform the boundary condition to
\begin{equation}
\psi^+(t)=h_0(t)\psi^-(t)
\end{equation}
i.e.\
\begin{equation}\label{muskeshtransformed}
\ln \psi^+(t)-\ln \psi^-(t)=\ln h_0(t)
\end{equation}
Assuming $\ln \psi(z)$ to be one-valued, sectionally holomorphic and vanishing at infinity, we have
\begin{equation}
\ln \psi(z)=\frac{1}{2\pi i}\int_\gamma \frac{\ln h_0(t)dt}{t-z},\quad \psi(t)=\exp(\tilde{\Gamma}(z))
\end{equation}
where
\begin{equation}
\tilde{\Gamma}(z)=\frac{1}{2\pi i} \int_\gamma \frac{\ln h_0(t)dt}{t-z}
\end{equation}
Then $\psi(z)$ 
is sectionally holomorphic, equal to infinity at infinity and everywhere nonzero and
for an arbitrary $t_0\in \gamma,$ $\tilde{\Gamma}^+(t_0)-\tilde{\Gamma}^-(t_0)=\ln h_0(t_0)$,
so $\psi$ is a particular solution to the transformed homogeneous problem of Eqn.(\ref{muskeshtransformed})
thus we obtain a particular solution to the problem Eqn.(\ref{muskesh351inhomo}) with $g\equiv 0$
according to
\begin{equation}
f_0(z):=\left\{
\begin{array}{ll}
\frac{\exp(\tilde{\Gamma}(z)}{\Pi(z)} & , z\in D^+\\
z^{\kappa}\exp(\tilde{\Gamma}(z) &, z\in D^-
\end{array}
\right.
\end{equation}
By the Sokhotsky-Plemelj formula
\begin{equation}
\tilde{\Gamma}^\pm(t_0)=\pm \frac{1}{2}\ln h_0(t_0) +\tilde{\Gamma}(t_0)
\end{equation}
so that
\begin{equation}
f_0^+= \frac{\exp(\frac{1}{2}\ln h_0(t_0))\exp(\tilde{\Gamma}(t_0))}{\Pi(t_0)} =\frac{\exp(\tilde{\Gamma}(t_0))
	\sqrt{h(t_0)}}{t_0^{\frac{\kappa}{2}}\sqrt{\Pi(t_0)}}
\end{equation}
\begin{equation}
f_0^-= t_0^{-\kappa}\exp(-\frac{1}{2}\ln h_0(t_0))\exp(\tilde{\Gamma}(t_0))=
\frac{\exp(\tilde{\Gamma}(t_0))}{t_0^{\frac{\kappa}{2}}\sqrt{\Pi(t_0)}\sqrt{h(t_0)}}
\end{equation}
where $f_0$ is nonzero on the finite part of the plane (if $\kappa >0,$ $f_0(\infty)=0$)
including the boundary values $f_0(t)$ and $f_0(z)$ has degree at infinity $-\kappa.$
Now let $f(z)$ be an arbitrary solution. Since $f_0$ is nonzero
\begin{equation}
\frac{f^+(t)}{f^{-1}(t)}=\frac{f_0^+(t)}{f_0^{-1}(t)}
\end{equation}
thus the function $f(z)/f_0(z)$ is holomorphic but has finite degree at infinity thus is a polynomial.
Hence we conclude that $f(z)=f_0(z)P(z)$ for a complex polynomial and H\"older continuity of $h(t)$ implies the same property for $f^\pm (t).$
If $k$ is the degree of $P(z)$ then the degree of $f(z)$ at infinity is $-\kappa +k$ thus the degree
of a solution is not less than $-\kappa$, the degrees of $f(z)$ and $f_0(z)$ at infinity are equal only 
when $P(z)$ is constant nonzero. 
Now any solution vanishing at infinity is given by $f(z)=f_0(z)P_{\kappa-1}(z)$
where
\begin{equation}
P_{\kappa-1}(z)=\sum_{j=0}^{\kappa-1} c_jz^{\kappa-1-j}
\end{equation}
for constants $c_j,$ $j=0,\ldots,\kappa -1.$
In particular, if $\kappa\leq 0$ then the homogeneous Riemann boundary value problem has no nontrivial 
solution, vanishing at infinity; if $\kappa>0$ it has precisely $\kappa$ linearly independent solutions, vanishing at infinity
given by $f_0(z),zf_0(z),\ldots,z^{\kappa -1} f_0(z).$
\\
Now consider the inhomogeneous Eqn.(\ref{muskesh351inhomo}), i.e.\ with $g\not\equiv 0$ (the index of the problem is defined to be the index
$\kappa$ of the homogeneous version).
Setting $h(t)=f_0^+(t)/f^-(t)$ gives
\begin{equation}
\frac{f^+(t)}{f^+(t)} -\frac{f^-(t)}{f^-(t)}=\frac{g^(t)}{f^+(t)}
\end{equation}
Since $f(z)/f_0(z)$ has finite degree at infinity we have
\begin{equation}
\frac{f(t)}{f_0(t)}=\frac{1}{2\pi i} \int_\gamma \frac{g(t)dt}{f_0^+(t)(t-z)} +P(z) 
\end{equation}
for a polynomial $P(z),$ thus we obtain the general solution
\begin{equation}\label{muskeshinhomoekv}
f(z)=\frac{1}{2\pi i} \int_\gamma \frac{g(t)dt}{f_0^+(t)(t-z)} +f_0(z)P(z) 
\end{equation}
Since $-\kappa$ is the exact degree of $f_0$ at infinity the solution $f(z)$
vanishes at infinity for $\kappa\geq 0$ if and only if the degree of $P(z)$ is $\leq \kappa -1$;
for $\kappa =0$, we must use $P(z)\equiv 0.$ For $\kappa <0$ $P\equiv 0$
and
the coefficients of $z^{-1},\ldots, z^{-\kappa}$ in the expression
\begin{equation}
\frac{1}{2\pi i}\int_\gamma \frac{g(t)dt}{f^+(t)(t-z)}=-\frac{z^{-1}}{2\pi i}\int_\gamma \frac{g(t)dt}{f_0^+(t)}-
\frac{z^{-2}}{2\pi i}\int_\gamma \frac{tg(t)dt}{f_0^+(t)}-\cdots
\end{equation}
vanish.
We conclude that the general solution vanishing at infinity is given by Eqn.(\ref{muskeshinhomoekv}) with $P(z)=P_{\kappa-1}(z)$
of degree $\leq \kappa -1$ (where $P_{\kappa -1}\equiv 0$ for $\kappa =0$); if $\kappa<0$, then
given the necessary and sufficient conditions
$\int_\gamma \frac{t^j g(t)dt}{f_+^+(t)}=0,$ $j=0,\ldots,-\kappa-1$ for a solution, the solution is given by 
\begin{equation}
f(z)=\frac{f_0(z)}{2\pi i}\int_\gamma \frac{g(t)dt}{f^+(t)(t-z)}
\end{equation}
For $\kappa <0$ the solution is unique, for $\kappa=0$ there is only one solution vanishing at infinity and for $\kappa >0$
there is a parameter family of solutions (i.e.\ infinitely many).
We shall later use the solutions to the Riemann boundary value problem when considering some techniques regarding boundary value problems for polyanalytic functions.
\\
\\ 
  We mention (without giving the details) also the following problem as it is central and classic in the theory of boundary value 
  problems, although we shall not directly use it in this book. 
The so-called {\em Riemann-Hilbert problem}\index{Riemann-Hilbert problem} in its classical formulation, for a domain $\Omega\subset\C$ whose boundary, $\gamma,$ is a simple contour, consists in finding, for any given triplet $a(t),b(t),c(t)$ of real-valued continuous functions on $\gamma$, a function $f(z)$ that is holomorphic on $\Omega$ and satisfying
\begin{equation}\label{riemahilbert}
\re (a(t)+ib(t)) f^+(t)=a(t)\re f^+(t)-b(t)\im f^+(t)=c(t)
\end{equation}
on $\gamma,$ where $f^+$ denotes the boundary values of $f$, with limit taken from inside the domain $\Omega.$
Note that any linear combination, with real coefficients, of solutions $f_1,\ldots,f_n$, to the homogeneous problem ($c(t)\equiv 0$) will again be a solution to the homogeneous problem.
\\
\\
In the case of the half-plane we use the following replacement for the usual H\"older condition for the point at $\infty$. Let $\gamma$ be a straight line in $\C$ parametrized by $t$ and let $\varphi(t)$ be a H\"older continuous function  that takes the value $\varphi(\infty)$ as $t\to \infty$ or $-\infty$. The H\"older condition (with coefficient $\alpha>0$) at $\infty$ is that for some constant $c>0$ we have for all sufficiently large $\abs{t}$,
\begin{equation}\label{holderkondissss}
\abs{\varphi(t)-\varphi(\infty)}\leq \frac{c}{\abs{t}^\alpha}
\end{equation}
Then the integral
$\frac{1}{2\pi i}\int_\gamma\frac{\varphi(t)dt}{t-z}$ is defined as the limit over the segment $(r,s)$ where $r,s$ are points tending to $\infty$ on either side of a fixed point. 
The Sokhotsky-Plemelj formula remains in some sense valid for $\gamma$.
If (the line) $\gamma$ splits the plane into $\Omega^\pm$ then we have for $z\notin \gamma$, $z\in \Omega^\pm$ 
\begin{multline}
\frac{1}{2\pi i} \int_\gamma\frac{(\varphi(t)-\varphi(\infty))dt}{t-z}+\frac{\varphi(\infty)}{2\pi i} \int_\gamma\frac{dt}{t-z}=\\
\frac{1}{2\pi i} \int_\gamma\frac{(\varphi(t)-\varphi(\infty))dt}{t-z}\pm \frac{1}{2}\varphi(\infty)
\end{multline}
where the integral in the right hand side converges absolutely due to the condition of Eqn.(\ref{holderkondissss}). See also Section \ref{plemeljsec} of the appendix. Such a H\"older condition at $\infty$ (in the sense of Eqn.(\ref{holderkondissss})) is usually added when considering the Riemann-Hilbert problem for e.g.\ the upper half-plane and furthermore that
$a^2(t)+b^2(t)\neq 0$ on $\gamma$ including $t=\infty.$
\\
\\
Using the solution to the inhomogeneous Riemann boundary value problem of Eqn.(\ref{muskeshinhomoekv}) 
described above one can characterize the solvability for the case of the {\em unit disc}, $D=\{\abs{z}<1\},$ 
$\gamma=\partial D=\{\abs{z}=1\},$ as follows.
The boundary conditions may obviously be written 
$2c=(a+ib)f^+(t)+(a-ib)\overline{f^+}(t)$ on $\gamma.$
Let $f_0$ be the solution to the homogeneous Riemann boundary value problem 
defined via 
\begin{equation}
(a+ib)f^+(t) +(a-ib)f^-(t)=2c
\end{equation}
i.e.\
\begin{equation}
h(t)=-\frac{a-ib}{a+ib},\quad g(t)=\frac{2c}{a+ib}
\end{equation}
It can be written
$f_0(z)=c\exp(\tilde{\Gamma}(z))$ for $\abs{z}<1$ and $f_+(z)=c z^{-\kappa}\exp(\tilde{\Gamma}(z))$ for $\abs{z}>1$ where 
$c$ is an arbitrary nonzero constant.
We may extend $\tilde{\Gamma}(z)$ on $D^+$ to $D_-$ according to $\tilde{\Gamma}_*(z)=\overline{\tilde{\Gamma}}(1/z)$ (in the case of the unit disc).
Here
\begin{equation}
\Gamma(z)=\frac{1}{2\pi i}\int_\gamma \frac{\ln (t^{-\kappa}h(t))dt}{t-z}=\frac{1}{2\pi}\int_\gamma \frac{\mbox{arg}\left(-t^{-\kappa}\frac{a-ib}{a+ib}\right)dt}{t}
\end{equation}
where the numerator in the integral is a continuous function on $\gamma.$
Then
$\tilde{\Gamma}_*(z)=\tilde{\Gamma}(z)-i\alpha$
where $\alpha$ is the real constant
\begin{equation}
\alpha=\frac{1}{2\pi i}\int_\gamma \frac{\mbox{arg}\left(-t^{-\kappa}\frac{a-ib}{a+ib}\right)dt}{t}=\frac{1}{2\pi} \int_0^{2\pi}\mbox{arg}\left(-\exp(is)^{-\kappa}\frac{a-ib}{a+ib}\right)ds
\end{equation} 
and we have used $t=\exp(is).$
Set $f_{0,*}(z)=\bar{c}\exp(\tilde{\Gamma}_*(z))$ for $\abs{z}>1$ and
$f_{0,*}(z)=\bar{c}z^\kappa\exp(\tilde{\Gamma}_*(z))$ for $\abs{z}<1.$
Choosing $c=\exp\left(i\frac{i\alpha}{2}\right)$ we obtain $f_{*,0}(z)=z^\kappa f_{0}(z).$
By what we have already done, if $\kappa\geq 0$ then the solutions $f(z)$ bounded at infinity are
of the form $f(z)=P(z)f_0(z)$ for 
$P(z)=\sum_{j=0}^\kappa c_jz^{\kappa-j}.$ Thus $f(z)$ is a solution to the Riemann-Hilbert problem for the disc 
if and only if $f_*(z)=f(z),$ i.e.\
\begin{equation}
f_{0,*}(z)=z^\kappa f_0(z),\quad P_*(z)=\bar{P}(1/z)
\end{equation}
so that $z^\kappa\bar{P}(1/z)=P(z)$ which yields the condition
\begin{equation}
c_\kappa =\bar{c}_{\kappa-j},\quad j=0,\ldots,\kappa
\end{equation}
Setting $c_j=a_j+ib_j, j=0,\ldots,\alpha/2$ for real $a_j,b_j,$ and $b_{\kappa/2}=0$ we have
$\kappa$ different constants
$c_k=a_{\kappa-k}+ib_{\kappa-k}, k=\kappa/2,\ldots,\kappa$.
Thus for $\kappa\geq 0$ the general solution takes the form
$f(z)=\sum_{j=0}^\kappa c_j \phi_j(z)$
for arbitrary constants $c_j$, $j=0,\ldots,\kappa$, and linearly independent (over the real numbers) particular solutions $\phi_j,$
and each $\phi_j(z)$ will have the form $\tilde{f}_0(z)z^k$ for some $k=0,\ldots,\kappa$,
i.e.\ we have exactly $\kappa +1$ linearly independent solutions and
\begin{equation}
f(z)=\tilde{f}_0(z)\sum_{j=0}^\kappa c_{\kappa-j}z^j
\end{equation} 
For $\kappa\leq -2$ the homogeneous Riemann boundary value problem has no nontrivial 
solution (as we have proved) and thus also the Riemann-Hilbert problem will have no nontrivial solution.

\subsection{Reduction to the unit disc}
Let $\Omega\subset\C$ be a domain bounded by a H\"older continuous simple contour $\gamma.$ Let $z=\phi(\zeta),$ $\zeta=\phi^{-1}(z)=:\psi(z)$, be a conformal mapping of $\Omega$ onto the unit disc $\{\abs{\zeta}<1\}.$ It is known that both derivatives $\phi'(\zeta)$ and $\psi'(z)$ have continuous extension to the boundary.
It follows that if $F(t)$ is a H\"older continuous (with H\"older coefficient say $\alpha$) function on $\gamma$ then this will hold true
under the conformal transformation onto the unit disc i.e.\ for the function $(F\circ \phi)(\psi(z))=(F\circ \phi)(\zeta))$ with the same index, and the same is true if the roles of $\gamma$ and the unit circle are interchanged. Hence the Riemann-Hilbert problem, for $\Omega$ given by Eqn.(\ref{riemahilbert}), can be transformed to a
Riemann-Hilbert problem for the unit disc with boundary condition
\begin{equation}\label{riemahilbert}
\re (a(\phi(t))+ib(\phi(t))) f^+(\phi(t))=a(\phi(t))\re f^+(\phi(t))-b(\phi(t))\im f^+(\phi(t))=c(\phi(t))
\end{equation}
for $t$ varying in the unit circle.
Once we have a solution $f(\zeta)$ on the unit circle then we may obtain a solution to the problem for $\Omega$
according to $(f\circ \psi)(z).$
The index of the problem is obtained by
\begin{equation}
\kappa=\frac{1}{2\pi i}\left[\log\frac{a-ib}{a+ib}\right]_\gamma=\frac{1}{2\pi }\left[\mbox{arg}\frac{a-ib}{a+ib}\right]_\gamma
\end{equation}
where
$\frac{1}{2\pi }[\mbox{arg}g]_\gamma:=\frac{1}{2\pi}\int_{\gamma} d\arg g(t)$ (i.e.\ the change in argument
of $g$ when traversing $\gamma$ in the positive direction, passing every point once).
Note that 
once $a(t),b(t),c(t)$ and the boundary $\gamma,$ are H\"older continuous (with some H\"older coefficient) then
the same will be true for the boundary value $f^+$ of any solution since it is well-known that this holds true for the case of the unit disc, and as we have seen the solution for the more general domain
as above is obtained by conformally mapping to the disc, and such will preserve the H\"older continuity.
\\
\\
Some special cases
of the higher order generalizations of such boundary value problems have in more recent literature been more popular than others.
We shall begin by describing some of the most common problems that appear in the recent literature concerning $q$-analytic functions,
namely the Schwarz problem, the Dirichlet problem and the Neumann problem, and we shall mainly consider the problems for the unit disc, but we will shall, in sections that follow, give some information about more arbitrary domains and we also conclude with some remarks on further generalizations. 

\section{The Dirichlet problem for the disc}
In the Riemann-Hilbert problem (Eqn.(\ref{riemahilbert})), for the unit disc $D=\{\abs{z}<1\},$ setting $a(t)\equiv 1,b(t)\equiv 0$ and $c(t)=g(t)$ we obtain the problem of finding a holomorphic function $f(z)$ on the unit disc such that it has prescribed boundary values $f(z)=g(z)$ for $z\in \partial D.$
\\
The Pompieu operator (see Eqn.(\ref{pompdef0})) is a useful tool in solving boundary value problems for first order equations. 
In the case of $q=1$, for a bounded domain $\Omega$ with continuous rectifiable boundary $\gamma$, given continuous functions $f$ on $\gamma$ and $g$ on $\Omega$, the function
\begin{equation}
F(z):=\frac{1}{2\pi i}\int_\gamma f(\zeta)\frac{d\zeta}{\zeta -z}-\frac{1}{\pi}\int_\Omega g(\zeta)\frac{d\mu(\zeta)}{\zeta -z}
\end{equation}
(where $\zeta=\eta+i\nu$, $d\mu(\zeta)$ denotes $d\eta d\nu:=d\eta\wedge d\nu$) is a solution in the weak sense to the equation $\partial_{\bar{z}} F=g$ on $\Omega$, by the properties of the Pompieu operator. 
For convenience we shall in the remainder of the section simply use $\mu(z)$
for the area measure in $\C$ which for $z=x+iy,$ upon the identification $\R^2\simeq \C$ means $d\mu(z)=dxdy=dx\wedge dy$.
However, the boundary values of $F$ will in general not coincide with the given $f.$ This can be remedied as shown below.
For $q>2$ (indeed, for the treatment of the more general model equation $\partial_z^m\partial_{\bar{z}}^q F=0$) Begehr \& Hile \cite{begehrhile} introduce the 
higher order Pompieu operators $T_{p,q}$ given by Eqn.(\ref{pomieueq1}) and Eqn.(\ref{pomieueq2}). These satisfy,
\begin{equation}
\partial_z^\mu\partial_{\bar{z}}^\nu T_{m,n}f=T_{m-\mu,n-\nu}f,\quad \mu+\nu\leq m+n
\end{equation}
in the weak sense. If $0<m+n$ then $T_{m,n}$ are weakly singular
and for $m+n=0<m^2+n^2$ the operators are strongly singular of Calderon-Zygmund type to be interpreted as Cauchy principal integrals.
Here we shall initially focus on presenting some known results for the {\em unit disc} in the complex plane, but we shall follow up with comments and examples on more general cases. 
\begin{theorem}[Dirichlet problem for the homogeneous $\overline{\partial}$-equation]\label{dirichlethomo}
Given a continuous boundary function $\gamma\in C^0(\partial D)$,
$D=\{\abs{z}<1\},$ the {\em Cauchy integral}\index{Cauchy integral} 
\begin{equation}\label{hoppa}
f(z)=\frac{1}{2\pi i}\int_{\abs{\zeta}=1}\frac{\gamma(\zeta) d\zeta}{\zeta-z}
\end{equation} 
solves the Dirichlet problem $\partial_{\bar{z}} f=0$ on $D$, $f|_{\partial D}=\gamma,$
(i.e.\ $f$ is assumed to be continuous on $\overline{D}$) as soon as
\begin{equation}\label{lobo00}
\frac{1}{2\pi i}\int_{\abs{\zeta}=1}\frac{\bar{z} d\zeta}{1-\bar{z}\zeta}=0, \quad \abs{z}<1
\end{equation} 
\end{theorem}
\begin{proof}
The last condition stems from the 
fact that the Cauchy integral of Eqn.(\ref{hoppa}) gives one holomorphic function in $D$ and one in $\hat{\C}\setminus D$
where by the Sokhotzki-Plemelj formula, given 
any $z_0\in \{\abs{z}=1\}$
\begin{equation}
\lim_{z\to z_0,\abs{z}<1} f(z)-\lim_{z\to z_0,\abs{z}>1} f(z)=\gamma(z_0)
\end{equation} 
Namely, we have for $z$ such that $\abs{z}>1$ (i.e.\ $\frac{1}{\bar{z}}\in D$),
\begin{equation}\label{loboekv}
f\left(\frac{1}{\bar{z}}\right)= -\frac{1}{2\pi i}\int_{\abs{\zeta}=1}\frac{\bar{z} d\zeta}{1-\bar{z}\zeta}
\end{equation} 
We obtain
\begin{multline}
f(z)-f\left(\frac{1}{\bar{z}}\right)= \frac{1}{2\pi i}\int_{\abs{\zeta}=1}\gamma(\zeta)
\left(\frac{\zeta}{\zeta-z}+\frac{\overline{\zeta}}{\overline{\zeta-z}}-1\right)\frac{d\zeta}{\zeta}=\\
\frac{1}{2\pi i}\int_{\abs{\zeta}=1}\gamma(\zeta)
\re\left(\frac{2\zeta}{\zeta-z}-1\right)\frac{d\zeta}{\zeta}=\\
\frac{1}{2\pi i}\int_{\abs{\zeta}=1}\gamma(\zeta)
\re\left(\frac{\zeta+z}{\zeta-z}\right)\frac{d\zeta}{\zeta}
\end{multline} 
Now the real part of $(\zeta+z)/(\zeta-z)$
is (recall from Section \ref{realandimagsec})
the Poisson kernel
\begin{equation}
\frac{\abs{\zeta}^2-\abs{z}^2}{\abs{\zeta-z}^2}
\end{equation}
So the limit $z\to z_0, \abs{z}>1,$ for $z_0\in \{\abs{z}=1\},$ exists in Eqn.(\ref{loboekv}) 
by the properties of the Poisson kernel and clearly $\lim_{z\to z_0} z=\lim_{z\to z_0} \frac{1}{\bar{z}}.$
In particular, given Eqn.(\ref{lobo00}) we have
\begin{equation}
\lim_{z\to z_0,\abs{z}>1} f(z)=-\lim_{z\to z_0,\abs{z}>1}\frac{1}{2\pi i}\int_{\abs{\zeta}=1}\frac{\bar{z} d\zeta}{1-\bar{z}\zeta}=0
\end{equation}
proving sufficiency of Eqn.(\ref{lobo00}).
To prove necessity, note that since we are assuming the $f(z)$ given on $\overline{D}$ solves the Dirichlet problem
we already know $\lim_{z\to z_0,\abs{z}<1} f(z)=\gamma(z_0)$ thus 
\begin{equation}
\lim_{z\to z_0,\abs{z}>1} f(z)=-\lim_{z\to z_0,\abs{z}>1}\frac{1}{2\pi i}\int_{\abs{\zeta}=1}\frac{\bar{z} d\zeta}{1-\bar{z}\zeta}=0
\end{equation}
Since on $\hat{\C}\setminus D$ we have $f(\infty)=0$ the maximum principle implies that $f\equiv 0$ for $\abs{z}>0.$
This proves necessity of Eqn.(\ref{lobo00}).
This completes the proof.
\end{proof}
This immediately yields the solution to the inhomogeneous problem.
\begin{theorem}[Dirichlet problem for the inhomogeneous $\overline{\partial}$-equation]\label{dirichletholomorphic}
Given a continuous boundary function $\gamma\in C^0(\partial D)$ and a function $g\in C^0(D),$
$D=\{\abs{z}=1\},$ then the function 
\begin{equation}\label{hoppa1}
f(z)=\frac{1}{2\pi i}\int_{\abs{\zeta}=1}\frac{\gamma(\zeta) d\zeta}{\zeta-z} -
\frac{1}{\pi}\int_{\abs{\zeta}=1}\frac{g(\zeta) d\mu(\zeta)}{\zeta-z}
\end{equation} 
solves the Dirichlet problem $\partial_{\bar{z}} f=g$ on $D$, $f|_{\partial D}=\gamma,$
(i.e.\ $f$ and $\partial_{\bar{z}} f=$ are assumed to be continuous on $\overline{D}$) as soon as
\begin{equation}\label{lobo001}
\frac{1}{2\pi i}\int_{\abs{\zeta}=1}\frac{\bar{z} d\zeta}{1-\bar{z}\zeta}-
\frac{1}{\pi}\int_{\abs{z}<1} \frac{g(\zeta) d\mu(\zeta)}{\zeta-z} =0
\end{equation} 
\end{theorem}
\begin{proof}
The function $\varphi:=f(z)-Tg$ is a weak solution to the homogeneous equation $\overline{\partial}\varphi=0,$ thus a strong solution by elliptic regularity. 
Applying the condition given by Eqn.(\ref{lobo00}) of Theorem \ref{dirichlethomo} to the boundary function
$\gamma-Tg$, associated to $\varphi$, on $\partial D$, yields using the Cauchy formula via
\begin{multline}
\frac{1}{2\pi i}\int_{\abs{\zeta}=1}\frac{1}{\pi}\int_{\tilde{\zeta}} f(\tilde{\zeta})\frac{d\mu(\tilde{\zeta})}{\tilde{\zeta}-\zeta}\frac{\bar{z}d\zeta}{1-\bar{z}\zeta} 
=\\
-\frac{1}{\pi}\int_{\abs{\tilde{\zeta}}<1}f(\tilde{\zeta}) 
\frac{1}{2\pi i}\int_{\abs{\zeta}=1} \frac{\bar{z}}{1-\bar{z}\zeta} \frac{d\zeta d\mu(\tilde{\zeta})}{\zeta-\tilde{\zeta}} =
-\frac{1}{\pi}\int_{\abs{\tilde{\zeta}}=1}f(\tilde{\zeta})\frac{\bar{z}}{1-\bar{z}\zeta}d\mu(\tilde{\zeta})
\end{multline}
which yields Eqn.(\ref{lobo001}). 
\end{proof}

\begin{remark}\label{greenrem}
Let us point out also how the Dirichlet problem for holomorphic functions on the unit disc can be solved by use of the Green function in accordance with Section \ref{realandimagsec}.
Let $\Omega\subset \C$ be a simply connected bounded domain, let 
 $g(z,z_0)$ be the Green function given by Eqn.(\ref{greeneq}) and let 
 \begin{equation}\label{greenremekv}
 \eta(z,z_0)=\int_{p_0}^z \frac{\partial g(\zeta,z_0)}{\partial n}d\mu(\zeta) +c
 \end{equation}
 where $c$ is a real constant.
 Then the function $f(z):=u(z)+iv(z),$ defined by
 \begin{equation}\label{holomorphdirichlet}
 f(z)=\int_{\partial \Omega} u(\zeta) \frac{\partial }{\partial n}(g(\zeta,z)+i\eta(\zeta,z))d\mu(\zeta)
 \end{equation}
 is holomorphic in $\Omega$ and is clearly uniquely determined by the boundary values of its real part, i.e.\ solves the classical Dirichlet problem for holomorphic functions.
 Note that $\eta(z,z_0)$ is a 
 harmonic conjugate of $g$. The function $g(z,z_0)+i\eta(z,z_0)$
 is called the {\em complex Green function} (with pole at $z_0$) and we have for the Riemann mapping, $\omega$, of $\Omega$ onto the disc, 
 mapping $z_0$ to $0$, such that $\omega(z_0,z_0)=0$ and 
 $\partial_\zeta\omega(z_0,z_0)\neq 0$
 then
 $\omega(z,z_0)=\exp(-G(z,z_0))$.
 The Green function $g$ and the complex Green function are invariant under conformal mappings.
\end{remark}

\subsection{The Dirichlet problem for the inhomogeneous $q$-analytic equation, $q>1$, for the disc.}\label{dirichletpolysecc}
Let $\Omega$ be a Jordan domain with rectifiable boundary, $\gamma$. 
Given $g\in C^1(\Omega,\C)$ and $f_j\in C^0(\gamma,\C),$ $j=0,\ldots ,q-1,$ find $F\in C^0(\overline{\Omega}),$ such that $\overline{\partial}^q F=g$ on $\Omega$ and $\overline{\partial}^j F=f_j$ on $\gamma,$ $j=0,\ldots ,q-1.$
\\
\\
 The Dirichlet problem for the homogeneous polyanalytic equation was solved by
 Theodorescu, see Theorem \ref{theodorescuthm}, 
   which immediately shows that $q$-analytic functions are determined by the boundary values of the
  functions $\frac{\partial^j f}{\partial \bar{z}^j},$ $j=0,\ldots,q-1.$
    Recall that the theorem states that for a $q$-analytic function $f$ on $\overline{\Omega}$ (for $q\in \Z_+$ and a domain $\Omega$ bounded by a rectifiable closed contour $\Gamma$)
  \begin{equation} 
  f(z)= 
  \frac{1}{2\pi i}\sum_{j=0}^{q-1}\int_{\Gamma} \frac{1}{j!(t-z)}(\bar{t}-\bar{z})^j \frac{\partial^j f}{\partial \bar{t}^j} dt, \quad z\in \Omega 
  \end{equation}
  This result will be included in the result for the inhomogeneous case which we shall prove below, see e.g.\ Begehr \cite{begehrbok}. 
  %
   \begin{theorem}[Dirichlet problem for the inhomogeneous polyanalytic equation]\label{dirichletpolyettan}
   Given $g\in L^p(D)$, $p>2,$ $\gamma_j\in C^0(\partial D,\R),$ $j=0,\ldots ,q-1,$ the problem 
   \begin{equation}
   (\partial_{\bar{z}})^q f=g,\mbox{ on } D, \quad (\partial_{\bar{z}})^j f=\gamma_j,\mbox{ on } \partial D,\quad j=0,\ldots ,q-1
   \end{equation}
   is uniquely solvable for $f\in L^p(D,\C),$ iff
   \begin{multline}\label{solvcond}
   0=
   \frac{(-1)^{q-l}\bar{z}}{\pi} \int_{\abs{\zeta}<1} \frac{g(\zeta)}{1-\bar{z}\zeta}
   \frac{(\overline{\zeta}-\overline{z})^{q-1-l}}{(q-1-l)!}d\mu(\zeta)+\\
   \sum_{l=j}^{q-1} \frac{\bar{z}}{2\pi i} \int_{\abs{\zeta}<1}  \frac{(-1)^{l-j}\gamma_l(\zeta)}{1-\bar{z}\zeta}
   \frac{(\overline{\zeta}-\overline{z})^{l-j}}{(l-j)!}d\zeta  
   \end{multline}
   for $j=0,\ldots,q-1,$
   and then the solution is given by
   \begin{multline}\label{diricheq}
   f(z)=
   \frac{(-1)^q}{\pi} \int_{\abs{\zeta}<1}  \frac{g(\zeta)}{(q-1)!}
   \frac{(\overline{\zeta}-\overline{z})^{q-1}}{\zeta-z}d\zeta  
   +\\
   \sum_{j=0}^{q-} \frac{(-1)^j}{2\pi i} \int_{\abs{\zeta}=1} \frac{\gamma_j(\zeta)}{j!}
   \frac{(\overline{\zeta}-\overline{z})^j}{\zeta-z}d\zeta 
   \end{multline}
   \end{theorem}
   \begin{remark}
   The last sum in Eqn.(\ref{diricheq}) is precisely the solution to the homogeneous Dirichlet problem defining
   $q$-analytic functions, as was proved by Theodorescu \cite{theodorescu}, see Theorem \ref{theodorescuthm}.
   \end{remark}
   \begin{proof}
   We shall use induction in $q.$ The case $q=1$ is precisely Theorem \ref{dirichletholomorphic}. So assume $q\geq 2$ and that 
   the result holds true for the $(q-1)$-analytic case.
   The problem can be reformulated as the system
   \begin{equation}
   (\partial_{\bar{z}})^{q-1} f=h,\mbox{ on } D, \quad (\partial_{\bar{z}})^j f=\gamma_j,\mbox{ on } \partial D, j=0,\ldots ,q-2
   \end{equation}
   \begin{equation}
   \partial_{\bar{z}} h=g,\mbox{ on } D, \quad \partial_{\bar{z}} h=\gamma_{q-1},\mbox{ on } \partial D
   \end{equation}
   By the induction hypothesis this gives for $j=0,\ldots,q-2,$ the condition of Eqn.(\ref{solvcond})
   \begin{multline}\label{solvcond001}
   0=
   \frac{(-1)^{q-l}\bar{z}}{\pi} \int_{\abs{\zeta}<1} \frac{h(\zeta)}{1-\bar{z}\zeta}
   \frac{(\overline{\zeta}-\overline{z})^{q-1-l}}{(q-1-l)!}d\mu(\zeta)+\\
   \sum_{l=j}^{q-1} \frac{\bar{z}}{2\pi i} \int_{\abs{\zeta}=1}  \frac{(-1)^{l-j}\gamma_l(\zeta)}{1-\bar{z}\zeta}
   \frac{(\overline{\zeta}-\overline{z})^{l-j}}{(l-j)!}d\zeta  
   \end{multline}
   together with the condition
   \begin{equation}
   0=   -\frac{\bar{z}}{\pi} \int_{\abs{\zeta}<1} \frac{\gamma_{q-1}(\zeta)d\mu(\zeta) }{1-\bar{z}\zeta} +
   \frac{\bar{z}}{2\pi i} \int_{\abs{\zeta}=1}  \frac{g(\zeta)}{\zeta -z}d\zeta
   \end{equation}
   This implies that for $j=0,\ldots,q-2,$
   \begin{multline}
   0=\frac{(-1)^{q-1-l}\bar{z}}{\pi} \int_{\abs{\zeta}<1} 
   \frac{h(\zeta)}{1-\bar{z}\zeta}
   \frac{(\overline{\zeta}-\overline{z})^{q-2-l}}{(q-2-l)!}d\mu(\zeta)+\\
   \sum_{l=j}^{q-2} \frac{\bar{z}}{2\pi i} \int_{\abs{\zeta}=1}  \frac{(-1)^{l-j}\gamma_l(\zeta)}{1-\bar{z}\zeta}
   \frac{(\overline{\zeta}-\overline{z})^{l-j}}{(l-j)!}d\zeta  =\\
   \frac{(-1)^{q-l}\bar{z}}{\pi} \int_{\abs{\zeta}<1} 
   \frac{g(\zeta)}{1-\bar{z}\zeta}
   \frac{(\overline{\zeta}-\overline{z})^{q-1-l}}{(q-1-l)!}d\mu(\zeta)+\\
   \sum_{l=j}^{q-2} \frac{\bar{z}}{2\pi i} \int_{\abs{\zeta}=1}  \frac{(-1)^{l-j}\gamma_l(\zeta)}{1-\bar{z}\zeta}
   \frac{(\overline{\zeta}-\overline{z})^{l-j}}{(l-j)!}d\zeta 
   \end{multline}
   which proves Eqn.(\ref{solvcond}) for the case $q.$
   Next note that
   \begin{multline}
   \frac{1}{\pi} \int_{\abs{\zeta}<1} h(\zeta) \frac{(\overline{\zeta}-\overline{z})^{q-2}}{(q-2)!(\zeta-z)}d\mu(\zeta)=\\
   \frac{1}{2\pi i} \int_{\abs{\tilde{\zeta}}=1} \gamma_{q-1}(\tilde{\zeta})\Psi_{q-1}(\tilde{\zeta},z)d\tilde{\zeta}
   -\frac{1}{\pi} \int_{\abs{\tilde{\zeta}}<1} g(\tilde{\zeta})\Psi_{q-1}(\tilde{\zeta},z)d\mu(\tilde{\zeta})
   \end{multline}
   where
   \begin{multline}
   \Psi_{q-1}(\tilde{\zeta},z):=-\frac{1}{\pi} \int_{\abs{\tilde{\zeta}}<1}
   \frac{(\overline{\zeta}-\overline{z})^{q-2}}{(q-2)!(\zeta -z)}  \frac{1}{\tilde{\zeta}-\zeta}d\mu(\zeta)
   =\\
   -\frac{1}{\pi} \int_{\abs{\tilde{\zeta}}<1}
   \frac{(\overline{\zeta}-\overline{z})^{q-2}}{(q-2)!(\tilde{\zeta} -z)}
   \left(
   \frac{1}{\tilde{\zeta}-\zeta}-\frac{1}{\zeta -z}
   \right)
   d\mu(\zeta)
   =\\
   \frac{(\overline{\tilde{\zeta}}-\overline{z})^{q-1}}{(q-1)!(\tilde{\zeta} -z)}
   -\frac{1}{2\pi i}\frac{1}{(p-1)!(\tilde{\zeta}-z)} \frac{1}{\pi} \int_{\abs{\tilde{\zeta}}<1} 
   \left(
   \frac{(\overline{\tilde{\zeta}}-\overline{z})^{q-1}}{\zeta -\tilde{\zeta}}-
   \frac{(\overline{\tilde{\zeta}}-\overline{z})^{q-1}}{\zeta -z}
   \right) d\mu(\zeta)
   =\\
   \frac{(\overline{\tilde{\zeta}}-\overline{z})^{q-1}}{(q-1)!(\tilde{\zeta} -z)}
   -
   \frac{1}{2\pi i}\frac{1}{(p-1)!(\tilde{\zeta}-z)} \frac{1}{\pi} \int_{\abs{\tilde{\zeta}}=1} 
   (\overline{\zeta}-\overline{z})^{q-1}
   \left(
   \frac{(\overline{\tilde{\zeta}}-\overline{z})^{q-1}}{1-\tilde{\zeta}\bar{\zeta}}-
   \frac{(\overline{\tilde{\zeta}}-\overline{z})^{q-1}}{1-z\tilde{\zeta}}
   \right)\frac{d\bar{\zeta}}{\bar{\zeta}}\\
   =\frac{(\overline{\tilde{\zeta}}-\overline{z})^{q-1}}{(q-1)!(\tilde{\zeta} -z)}
   \end{multline}
   This implies that
   \begin{multline}
   f(z)=\sum_{l=j}^{q-2} \frac{(-1)^j}{2\pi i} \int_{\abs{\zeta}=1}  \frac{\gamma_j(\zeta)}{j!}
   \frac{(\overline{\zeta}-\overline{z})^{j}}{\zeta -z}d\zeta +
   \\ \frac{(-1)^{q-1}}{\pi}\int_{\abs{\zeta}<1} \frac{h(\zeta)}{(q-2)!}
   \frac{(\overline{\zeta}-\overline{z})^{q-2}}{\zeta -z}d\mu(\zeta)=\\ 
   \sum_{l=j}^{q-2} \frac{(-1)^j}{2\pi i} \int_{\abs{\zeta}=1}  \frac{\gamma_j(\zeta)}{j!}
   \frac{(\overline{\zeta}-\overline{z})^{j}}{\zeta -z}d\zeta 
   + \frac{(-1)^{q-1}}{\pi}\int_{\abs{\zeta}<1} \frac{g(\zeta)}{(q-1)!}
   \frac{(\overline{\zeta}-\overline{z})^{q-1}}{\zeta -z}d\mu(\zeta)
   \end{multline}
   This completes the proof.
   \end{proof}
  
\section{The Schwarz problem for the disc}
\begin{definition}
 An operator $S\colon C^0(\partial \Omega,\R)\to \{f\in \mathscr{O}(\Omega)\cap C^0(\overline{\Omega},\C)\colon 
\re Sf=\re f\}$
 is called a {\em Schwarz operator}.\index{Schwarz operator}. In particular, $S$ is uniquely determined up to an imaginary term,
 i.e.\ there is a unique $S$ satisfying the condition $\im (Sf)(p_0)=0$ for a point $p_0\in \Omega.$
 \end{definition}
 When $\Omega$ has a Green function $g$ and associated harmonic conjugate $h$
 then
 the Schwarz operator of $\Omega$ is given by
 $(Sf)(z):=-\frac{1}{2\pi}\int_{\partial \Omega} u(\zeta)\frac{\partial}{\partial n}(g(\zeta,z)+ih(\zeta,z))d\zeta$,
 where $u=\re f$ and $\mu(\zeta)$ denote the Lebesgue surface measure.
 For $\Omega=D$ the Schwarz operator is given by
 \begin{equation}
 (Sf)(z):=-\frac{1}{2i\pi}\int_{\partial D} u(\zeta)\frac{\zeta +z}{z-\zeta}\frac{d\zeta}{\zeta}+i\im(Sf)(0),\quad z\in D
 \end{equation}

 Finding explicit Schwarz operators and reproducing kernels have long been popular approaches to solving many  boundary value problems in complex analysis.
 The analysis of boundary value problems involves of course some theory
regarding existence and uniqueness, but much of the analysis concerns
 searching for closed form solutions that are domain-dependent, for example, 
 solving the Schwarz problem for a half-disc
 and for the unit disc respectively, in closed form, are two different results.
Oftentimes, once a boundary value problem has been solved for the unit disc there are methods to derive solutions
 for more general domains. Gakhov \cite{gakhov}, p.259, provides (citing an unpublished work of M.P. Ganin) 
 for example a method to treat certain 
 cases where the domains are obtained as the image, under rational functions, of the unit disc.
 If we let $\Omega\subset\C$ be domain bounded by a 
 simple $C^q$-smooth closed contour, with complex coordinate $z,$ then there is a conformal mapping $\phi:\Omega\to D:=\{\abs{z}<1\},$
 such that $\phi$ and its  derivatives up $q$:th order, to can be extended continuously to $\partial\Omega$ mapping $\partial\Omega$ to $\partial D.$ The Green function $g(z,\zeta)$ of $\Omega$
  (cf. Section \ref{realandimagsec}) can be written
 \begin{equation}\label{ekv761}
 g(z,\zeta)=\log\abs{\frac{1-\phi(z)\overline{\phi(\zeta)}}{\phi(z)-\phi(\zeta)}}=-\log\abs{z-\zeta}+\omega(z,\zeta)
 \end{equation}
 where 
 \begin{equation}\label{ekv762}
 \omega(z,\zeta):=\log\abs{\frac{z-\zeta}{\phi(z)-\phi(\zeta)}(1-\phi(z)\overline{\phi(\zeta)})}
 \end{equation}
 Recall from Section \ref{realandimagsec} that we have explicit expressions
 for certain Green functions $g(z,\zeta)$ of $\Omega$, in particular for the unit disc, see Eqn.(\ref{greeneqdisken}). 
 In this text we shall begin by focusing on presenting some known results for the {\em unit disc} in the complex plane which we believe will render the disposition easier to digest.
 We shall however follow up with more general cases.
Let us start with the Schwarz problem
for the inhomogeneous $\overline{\partial}$-equation.
 By this we shall mean the following. Let $\Omega\subset\C$ be domain bounded by a 
 simple smooth closed contour, with complex coordinate $z.$ We seek a holomorphic function on $\Omega$
 continuous on $\overline{\Omega}$ such that $\re \partial_{\bar{z}} f(z)=0,$ 
 on $\partial\Omega,$ $\im f(0)=0.$ 
 We stick to the case $\Omega=\{\abs{z}<1\}=:D.$  
 Consider the following operator 
 \begin{equation}
 (S_1f)(z):=-\frac{1}{2i\pi}\int_{D} \frac{\zeta +z}{\zeta -z}\frac{f(\zeta)}{\zeta}+
 \frac{1+z\bar{\zeta}}{1-z\bar{\zeta}}\frac{\overline{f(\zeta)}}{\bar{\zeta}}d\mu(\zeta),\quad z\in D
\end{equation}
 For $f\in L^p(\overline{D})$, $S_1$ satisfies (see e.g.\ Begehr \cite{begehrbok}, p.233)
 \begin{equation}
 \partial_{\bar{z}}(S_1f)(z)=f(z), \quad z\in D,\quad
 \quad \re (S_1f)=0 \mbox{ on } \partial D, \im (S_1f)(0)=0
 \end{equation}

\begin{definition}
 For $f\in L^p(\overline{D})$, $p,j\in \Z_+,$ define
 \begin{equation}
 (S_jf)(z):=\frac{(-1)^j}{2\pi(j-1)!}\int_{D}(2\re (\zeta-z))^{j-1}\left(
  \frac{\zeta +z}{\zeta-z}\frac{f(\zeta)}{\zeta}+
 \frac{1+z\bar{\zeta}}{1-z\bar{\zeta}}\frac{\overline{f(\zeta)}}{\bar{\zeta}}
 \right)d\mu(\zeta)
 \end{equation}
 \end{definition}

We shall need the following formula.
 \begin{theorem}\label{formulalem}
 Any complex-valued $f\in C^1(D)\cap C^0(\overline{D})$
 can be represented as
 \begin{multline}
 f(z)=\frac{1}{2\pi i} \int_{\partial D} \re f(z)\frac{\zeta +z}{\zeta-z}\frac{d\zeta}{\zeta}
 +\frac{1}{2\pi}\int_{\partial D}\im f(\zeta)\frac{d\zeta}{\zeta}-\\
 \frac{1}{\pi}\int_D \left(
 \frac{\partial_{\bar{\zeta}} f(\zeta)}{\zeta-z}+\frac{z\overline{\partial_{\bar{\zeta}} f(\zeta)}}{1-z\bar{\zeta}}
 \right)d\mu(\zeta),\quad z\in D
 \end{multline}
 \end{theorem}
 \begin{proof}
 For $\Omega\subset\C$ be domain bounded by a simple smooth closed contour, with complex coordinate $z,$ 
 recall that the so-called {\em Gauss theorem}\index{Gauss theorem} implies that for a complex valued $f\in C^0(\overline{\Omega})\cap C^1(\Omega)$
 \begin{equation}\label{gausseq}
 \int_\Omega \partial_{\bar{z}} f(z)d\mu(z)=\frac{1}{2i}\int_{\partial\Omega} f(z)dz
 \end{equation}
 Let $\Omega=D$, $p_0\in D$ and $\epsilon>0$ sufficiently small such that $D_{\epsilon,p_0}:=\{\abs{z-p_0}<\epsilon\}\subset D.$
 By Eqn.(\ref{gausseq}) 
 \begin{equation}
 f(z)=\frac{1}{2\pi i}
 \int_{\partial(D\setminus D_{\epsilon,p_0})}  \frac{f(\zeta)d\zeta}{\zeta -p_0}-
 \int_{ D\setminus D_{\epsilon,p_0}}  \frac{\partial_{\bar{\zeta}}f(\zeta)}{\zeta -p_0}d\mu(\zeta)=0
 \end{equation}
 In polar coordinates $z=r\exp(i\theta)$ we have
 \begin{equation}
 \int_{ \overline{ D_{\epsilon,p_0}}}  \frac{\partial_{\bar{\zeta}}f(\zeta)}{\zeta -p_0}d\mu(\zeta)
 =\int_{0}^\epsilon \int_0^{2\pi} \partial_{\bar{z}}f(z_0+r\exp(i\theta)) dr d\theta
\end{equation}
 and 
 \begin{equation}
 \int_{D} \frac{ f(\zeta)}{\zeta -p_0}d\mu(\zeta)=
 \int_{D\setminus D_{\epsilon,p_0}} \frac{ f(\zeta)}{\zeta -p_0}d\mu(\zeta)-
 \int_{ \overline{ D_{\epsilon,p_0}}}  \frac{\partial_{\bar{z}}f(\zeta)}{\zeta -p_0}d\mu(\zeta)
 \end{equation}
 implies that
 \begin{equation}
 \lim_{\epsilon\to 0} \int_{D\setminus D_{\epsilon,p_0}}  \frac{f(\zeta)}{\zeta -p_0}d\mu(\zeta)=
 \int_{D}  \frac{f(\zeta)}{\zeta -p_0}d\mu(\zeta)
 \end{equation}
 Since
 \begin{equation}
 \int_{\partial D_{\epsilon,p_0}}  \frac{f(\zeta)d\zeta}{\zeta -p_0}=i\int_0^{2\pi} f(p_0+\epsilon\exp(i\theta))d\theta
 \end{equation}
 and
 \begin{equation}
 \int_{\partial(D\setminus D_{\epsilon,p_0})}  \frac{f(\zeta)}{\zeta -p_0}d\mu(\zeta)=
 \int_{\partial D}  \frac{f(\zeta)d\zeta}{\zeta -p_0}-i\int_0^{2\pi} f(p_0+\epsilon\exp(i\theta))d\theta
 \end{equation}
 we obtain
 \begin{equation}
 \lim_{\epsilon\to 0} \int_{\partial(D\setminus D_{\epsilon,p_0})}  \frac{f(\zeta)}{\zeta -p_0}d\mu(\zeta)=
 \int_{\partial D}  \frac{f(\zeta) d\zeta}{\zeta -p_0} -2\pi if(p_0)
 \end{equation}
 This proves that
 \begin{equation}\label{eq7}
 f(z)=\frac{1}{2\pi i}
 \int_{\partial D} f(\zeta) \frac{d\zeta}{\zeta -z}-\frac{1}{\pi}
 \int_{D}  \frac{\partial_{\bar{z}}f(\zeta)}{\zeta -z}d\mu(\zeta)
 \end{equation}
 Now by Eqn.(\ref{gausseq}) we have for fixed $z\in D$
 \begin{equation}
 \frac{1}{2\pi i}\int_{\partial D} f(\zeta) \frac{\bar{z}d\zeta}{1 -\bar{z}\zeta}
 -\frac{1}{\pi}
 \int_{D}  \frac{\bar{z}\partial_{\bar{\zeta}}f(\zeta)}{1 -\bar{z}\zeta}d\mu(\zeta)=0
 \end{equation}
 Taking the complex conjugate and adding to Eqn.(\ref{eq7}) and replacing $\bar{\zeta}=-\zeta \overline{d\zeta},$ for 
 $\abs{\zeta}=1$, gives
 \begin{equation}\label{slutklameq}
 f(z)=\frac{1}{2\pi i}\int_{\partial D} \left(\frac{\zeta f(\zeta)}{z-\zeta}+\frac{z\overline{f(\zeta)}}{z -\zeta}\right)
 \frac{d\zeta}{\zeta}
 -\frac{1}{\pi }\int_{D}\left(\frac{\partial_{\bar{\zeta}}f(\zeta)}{\zeta -z}+\frac{z\overline{\partial_{\bar{\zeta}}f(\zeta)}}{1 -z\bar{\zeta}}\right)d\mu(\zeta)
 \end{equation}
 for $z\in D.$
This completes the proof. 
 \end{proof}

 By subtracting $i\im f(0)$ from Eqn.(\ref{slutklameq}) we obtain the following corollary.
 \begin{corollary}[Modified Cauchy-Pompieu formula]\label{modcauchypompformula}
 Any complex-valued $f\in C^1(D)\cap C^0(\overline{D})$
 can be represented as
 \begin{multline}
 f(z)=\frac{1}{2\pi i} \int_{\partial D} 
 \re f(z)\frac{\zeta +z}{\zeta-z}\frac{d\zeta}{\zeta}-\\
 \frac{1}{2\pi}\int_{D}\left(
 \frac{\partial_{\bar{\zeta}} f(\zeta)}{\zeta}\frac{\zeta +z}{\zeta-z}+
 \frac{\overline{\partial_{\bar{\zeta}} f(\zeta)}}{\bar{\zeta}}\frac{1+z\bar{\zeta}}{1-z\bar{\zeta}}
 \right)d\mu(\zeta)+i\im f(0),\quad z\in D
 \end{multline}
 \end{corollary}

If $f$ in Corollary \ref{modcauchypompformula} is holomorphic then
$f(z)=\frac{1}{2\pi i} \int_{\partial D}\re f(\zeta)\left(\frac{2\zeta}{\zeta-z}-1\right)\frac{d\zeta}{\zeta} +i\im f(0)$,
where 
\begin{equation}
\frac{2\zeta}{\zeta-z}-1=\frac{\zeta+z}{\zeta-z}
\end{equation}
is a Schwarz kernel
\begin{equation}
Sg(z)=\frac{1}{2\pi i}\int_{\partial D}g(\zeta)\frac{\zeta+z}{\zeta-z}\frac{d\zeta}{\zeta}
\end{equation}
for $g\in C^0(\partial D,\R),$ and satisfies $\lim_{D\ni z\to\zeta} \re(Sg)(z)=g,$ $\zeta\in\partial D.$

 \begin{theorem}\label{sjthm0}
\begin{equation} 
S_j f(z)=S^j_1 f(z),\quad j\in \Z_+
\end{equation}
 \end{theorem}
 \begin{proof}
 We use induction in $j\in \Z_+.$ The case $j=1$ is trivial, so assume $j>1$ and that
 the statement holds true for $1,\ldots,j-1.$
 We have
 \begin{equation}\label{presdefeq}
 \partial_{\bar{z}} S_j f(z)=
 \frac{(-1)^{j-1}}{2\pi(j-2)!}\int_{D}(2\re (\zeta-z))^{j-2}\left(
  \frac{\zeta +z}{\zeta-z}\frac{f(\zeta)}{\zeta}+
 \frac{1+z\bar{\zeta}}{1-z\bar{\zeta}}\frac{\overline{f(\zeta)}}{\bar{\zeta}}
 \right)d\mu(\zeta)=S_{j-1}f(z)
 \end{equation}
 Furthermore
 \begin{equation}
 S_j f(0)=
 \frac{(-1)^{j}}{2\pi(j-1)!}\int_{D}(2\re \zeta)^{j-1}\left(
  \frac{f(\zeta)}{\zeta}+
 \frac{\overline{f(\zeta)}}{\bar{\zeta}}
 \right)d\mu(\zeta)
 \end{equation}
 implies \begin{equation}\label{preeqsdef3}
 \im S_j f(0)=0.\end{equation}
Also on $\partial D$
 \begin{multline}
 S_j f(z)=
 \frac{(-1)^j}{2\pi(j-1)!}\int_{D}(2\re (\zeta-z))^{j-1}\left(
  \frac{\zeta +z}{\zeta-z}\frac{f(\zeta)}{\zeta}+
 \frac{\overline{z+\zeta}}{\overline{z-\zeta}}\frac{\overline{f(\zeta)}}{\bar{\zeta}}
 \right)d\mu(\zeta)=\\
 \frac{(-1)^ji}{2\pi(j-1)!}\int_{D}(2\re (\zeta-z))^{j-1}\im\left(
  \frac{\zeta +z}{\zeta-z}\frac{f(\zeta)}{\zeta}
 \right)d\mu(\zeta)
 \end{multline}
 which implies \begin{equation}\label{preeqsdef}
 \re S_j f(z)=0 \mbox{ on }\partial D
 \end{equation}
 Now $S_1(S_{j-1}f)$ can be written
 \begin{multline}
 -\frac{1}{2\pi}
 \int_{D}\left(
  \frac{\zeta +z}{\zeta-z}\frac{(S_{j-1}f)(\zeta)}{\zeta}+
 \frac{\overline{z+\zeta}}{\overline{z-\zeta}}\frac{\overline{(S_{j-1}f)(\zeta)}}{\bar{\zeta}}
 \right)d\mu(\zeta)=\\
 \frac{(-1)^ji}{2\pi(j-1)!}\int_{D}(2\re (\zeta-z))^{j-1}\im\left(
  \frac{\zeta +z}{\zeta-z}\frac{f(\zeta)}{\zeta}
 \right)d\mu(\zeta),\quad z\in D
 \end{multline}
 On the other hand we obtain the same expression by applying Corollary \ref{modcauchypompformula}
 to $S_{j}.$ 
 This completes the proof.
 \end{proof}
 
 \begin{theorem}\label{sjthm1}
 For $f\in L^p(\overline{D},\C),$ $1\leq p, j\in \N,$ $S_kf$ satisfies
 \begin{equation}\label{sdefeq1}
 (\partial_{\bar{z}})^l S_j f(z)=S_{j-l} f(z),\quad 1\leq l\leq j, \mbox{ if }S_0f:=f
 \end{equation}
  \begin{equation} \label{sdefeq2}
  \re((\partial_{\bar{z}})^l S_j f)=0,\mbox{ on }\partial D, \quad 1\leq l\leq j-1
  \end{equation}
  \begin{equation} \label{sdefeq3}
    \im((\partial_{\bar{z}})^l S_j f)(0)=0, \quad 1\leq l\leq j-1
    \end{equation}
  \end{theorem}
  \begin{proof}
  Eqn.(\ref{sdefeq1}) is proved by Eqn.(\ref{presdefeq}) above.
  The result of Theorem \ref{sjthm0} together with Eqn.(\ref{preeqsdef}) (Eqn.(\ref{preeqsdef3}) respectively) yields Eqn.(\ref{sdefeq2}) (Eqn.(\ref{sdefeq3}) respectively)
  This completes the proof.
  \end{proof}
 The following can be found in Begehr \cite{begehrbok}, p.235, Theorem 55.
 \begin{theorem}[Schwarz problem for the inhomogeneous polyanalytic equation]
 The equation 
 \begin{equation}\label{dirichletpolyek1}
 \frac{\partial^q f}{\partial\bar{z}}=g(z)
 \end{equation}
 on the unit disc $D,$
 has a unique solution satisfying 
 \begin{equation}
 \re \frac{\partial^j f}{\partial\bar{z}}=0\mbox{ on }\partial D, \im \frac{\partial^j f(0)}{\partial\bar{z}}=0,\quad
  0\leq j\leq q-1 
 \end{equation}
 The solution is given by $f=S_q g.$
 \end{theorem}
 \begin{proof}
 The general solution to Eqn.(\ref{dirichletpolyek1}) is given by
 \begin{equation}\label{dirichletpolyek1}
f(z)=\sum_{j=0}^{q-1}\phi_j\bar{z}^j +S_q g(z)
 \end{equation}
 for holomorphic $\phi_j.$
 Indeed, the homogeneous version of Eqn.(\ref{dirichletpolyek1}) is simply the defining equation for 
 $q$-analytic functions, thus the general solution is given in the form
 $\sum_{j=0}^{q-1}\phi_j\bar{z}^j,$ for holomorphic $\phi_j(z)$ and 
  by
   Theorem \ref{sjthm1}
   $S_q g$ is a particular solution. We state without proof that all particular solutions are of this
   form.
 Also by
  Theorem \ref{sjthm1}
  $S_k g$ satisfies 
  \begin{equation} 
    \re((\partial_{\bar{z}})^{j-1} f=\re\left((j-1)!\phi_{j-1}+S_1f\right)
    =(k-1)!\re\phi_{j-1}=0\mbox{ on }\partial D
    \end{equation} 
   Hence $ \phi_{j-1}(z)=ic_{j-1}$ for real constants $c_{j-1}.$
   Since we also have $\im \frac{\partial^j f(0)}{\partial\bar{z}}=0,\quad
     0\leq j\leq q-1$, this implies that $c_{j-1}=0.$ Thus $\im \phi_{j-1}(0)=0.$ By iteration we can prove that
     $\phi_l(z)=0$ for all $l=0,\ldots,j-1.$
     This completes the proof.
 \end{proof}

 \section{The Neumann problem for the disc}
 Let $\nu$ denote the  outward normal derivative at the boundary of a simply connected
 bounded domain $\Omega$ with smooth boundary. In the case where the boundary is a circle $\{\abs{z-p_0}=R\},$ $p_0\in C$
 the outward normal is $\nu=\frac{z-p_0}{R},$ and in polar coordinates $z=r\exp(i\theta),$
 $\partial_\theta$ vanishes, giving the normal derivative
 \begin{equation}
 \partial_\nu=\partial_r=\frac{z}{R}\partial_z +\frac{\bar{z}}{R}\partial_{\bar{z}}
 \end{equation}
 For the unit disc $D=\{\abs{z}=1\},$ this yields
 \begin{equation}
 \partial_\nu=z\partial_z +\bar{z}\partial_{\bar{z}}
 \end{equation}
 The {\em Neumann problem for the homogeneous $\overline{\partial}$-equation} on the unit disc is
 to find for a given set of data $g_0\in C^0(\partial D,\C)$, a function $f\in C^0(\overline{D})$ such that $f$ is
 analytic on $D$ and
 satisfying $f=g_0$ on $\partial D,$ and $\partial_\nu f(0)=c_0,$ for a complex constant $c_0.$
 \\ \index{Neumann problem}
 The {\em Neumann problem for the inhomogeneous $\overline{\partial}$-equation} on the unit disc is
 to find, for a given set of data $g_0\in C^0(\partial D,\C)$, $g_1\in C^0(D),$ a function $f\in C^0(\overline{D})$ such that 
 $\partial_{\bar{z}}) f=g_1$ on $D$ and
 satisfying $f=g_0$ on $\partial D,$ and $(\partial_\nu) f(0)=c_0,$ for a complex constant $c_0.$
 \begin{theorem}[Neumann problem for the homogeneous $\overline{\partial}$-equation]
 The Neumann problem for the homogeneous $\overline{\partial}$-problem is solvable iff for $\abs{z}<1$
 \begin{equation}\label{kkk0}
 \frac{1}{2\pi i}\int_{\partial D} g_0(\zeta)\frac{d\zeta}{(1-\bar{z}\zeta)\zeta} =0
 \end{equation}
 The solution is given by
 \begin{equation}\label{neumaneq}
 f(z)=c_0-\frac{1}{2\pi i}\int_{\partial D} g_0(\zeta)\log(1-z\bar{\zeta})\frac{d\zeta}{\zeta}
 \end{equation}
 \end{theorem}
 \begin{proof}
 Using the assumed holomorphy of $f$ on $D$, the boundary condition can be rewritten as
 \begin{equation}
 zf'(z)=g_0(\zeta),\quad \abs{z}=1
 \end{equation}
 This becomes a Dirichlet condition for the $\overline{\partial}$-equation for the function $zf'(z)$ which 
 can be solved by
 \begin{equation}
 zf'(z)=\frac{1}{2\pi i}\int_{\partial D} g_0(\zeta)\frac{d\zeta}{\zeta-z}
 \end{equation}
 iff for $\abs{z}<1$
 \begin{equation}\label{kkk1}
 \frac{1}{2\pi i}\int_{\partial D} g_0(\zeta)\frac{\bar{z}d\zeta}{1-\zeta\bar{z}} =0
 \end{equation}
 Since $zf'$ is holomorphic and vanishes at $z=0$ we also have
 \begin{equation}\label{kkk2}
 \frac{1}{2\pi i}\int_{\partial D} \frac{g_0(\zeta)}{\zeta-0}d\zeta =0
 \end{equation}
 Since
 \begin{equation}
 \frac{1}{(1-\bar{z}\zeta)\zeta} =\frac{\bar{z}}{(1-\bar{z})}\frac{\zeta}{\zeta}+\frac{1}{\zeta}\frac{1-\bar{z}\zeta}{1-\bar{z}\zeta}
 \end{equation} 
 Eqn.(\ref{kkk0}). 
 Also we obtain
 \begin{equation}
 f'(z)=\frac{1}{2\pi i}\int_{\partial D} \frac{g_0(\zeta)}{\zeta(\zeta-z)}d\zeta 
 \end{equation}
 and since $\zeta\bar{\zeta}=1$ on $\partial D,$ integration yields
 \begin{equation}
 f(z)=c_0-\frac{1}{2\pi i}\int_{\partial D} g_0(\zeta)\log
 \underbrace{\left(\frac{\zeta-z}{\zeta}\right)}_{=(1-z\bar{\zeta})}\frac{d\zeta}{\zeta} 
 \end{equation}
 for a complex constant $c_0,$
 which proves Eqn.(\ref{neumaneq}).
 This completes the proof.
 \end{proof}
 
 \begin{theorem}[Neumann problem for the inhomogeneous $\overline{\partial}$-equation, see Begehr \cite{begehrboletinI}]
 The Neumann problem for the inhomogeneous $\overline{\partial}$-equation is solvable with solution $f$ in $L^p(D)$, 
 iff for $\abs{z}<1$
 \begin{equation}
 \frac{1}{2\pi i}\int_{\abs{\zeta}=1} g_0(\zeta)\frac{d\zeta}{(1-\bar{z}\zeta)\zeta}+
 \frac{1}{2\pi i}\int_{\abs{\zeta}=1} g_1(\zeta)\frac{\overline{d\zeta}}{1-\bar{z}\zeta}+
 \frac{1}{\pi }\int_{\abs{\zeta}<1} \frac{\bar{z}g_1(\zeta)}{(1-\bar{z}\zeta)^2}d\mu(\zeta)=0
 \end{equation}
 The solution is then given by
 \begin{equation}
 f(z)=c_0-\frac{1}{2\pi i}\int_{\abs{\zeta}=1} (g_0(\zeta)-\bar{\zeta}g_1(\zeta))\log(1-z\bar{\zeta})\frac{d\zeta}{\zeta}
 -\frac{1}{\pi}\int_{\abs{\zeta}<1} \frac{zf(\zeta)}{\zeta(\zeta-z)}d\mu(\zeta)
 \end{equation}
 \end{theorem}
 \begin{proof}
 Set $\varphi:=f(z)-Tg_1.$ 
 A known property (see Calderon-Zygmund \cite{calderonzygmund}) of the Pompieu operator $T$ is that
 $\partial_{z}Tg_1(z)=-\frac{1}{\pi}\int_D g_1(\zeta)\frac{d\mu(\zeta)}{(\zeta-z)^2}$
 holds in distribution sense a.e.\ for $z\in D$ when $g_1\in L^p(D,\C),$ $p>1,$
 and the integral is understood as a Cauchy principal value
 \begin{equation}
 \int_D g_1(\zeta)\frac{d\mu(\zeta)}{(\zeta-z)^2}=\lim_{\epsilon\to 0}
 \int_{D\setminus D_{\epsilon,z}}g_1(\zeta)\frac{d\mu(\zeta)}{(\zeta-z)^2}
 \end{equation}
 This implies that
 $\partial_{\bar{z}} \varphi=0$ on $D,$ $\partial_\nu \varphi=g_0-z(\partial_{z}Tg_1(z))-\bar{z}g_0$ on $\partial D$
 and $\varphi(0)=c_0-Tg_1(0).$
 As pointed out in Section \ref{areolarsec}, Vekua proved (see e.g.\ Vekua \cite{vekua}, Ch.1, parag.8) that
 $\Pi g_1\in C^\alpha(\overline{D})$ whenever $g_1\in C^\alpha(\overline{D}),$
 where $\Pi g_1=\partial_z Tg_1.$
 We thus have a homogeneous Neumann problem for $\varphi$, hence there is a solution 
 \begin{multline}\label{neumanmiddle}
 \varphi(z)=c_0-Tg_1(0)-\\
 \frac{1}{2\pi i}\int_{\partial D} (g_0(\zeta)-\zeta(\partial_{z}Tg_1(z))(\zeta)-\bar{\zeta}g_1(\zeta))\log(1-z\bar{\zeta})\frac{d\zeta}{\zeta}
 \end{multline}
 iff
 \begin{equation}\label{neumanmiddle2}
 \frac{1}{2\pi i}\int_{\partial D} (g_0(\zeta)-\zeta(\partial_{z}Tg_1)(\zeta)-\bar{\zeta}g_1(\zeta))\frac{d\zeta}{(1-\bar{z}\zeta)\zeta}=0
 \end{equation}
  Now up to the factor $-1$, the middle term of Eqn.(\ref{neumanmiddle}) can be written as
 \begin{multline}
 \frac{1}{2\pi i}\int_{\abs{\zeta}=1} (\partial_{z}Tg_1)(\zeta)\log (1-z\bar{\zeta})\frac{d\zeta}{\zeta}=\\
 -\frac{1}{\pi }\int_{\abs{\tilde{z}}<1} g_1(\tilde{\zeta})\frac{1}{2\pi i}
 \int_{\abs{\zeta}=1}\frac{\log(1-z\bar{\zeta})}{(\zeta-\tilde{\zeta})^2}d\zeta d\mu(\tilde{\zeta})
 =\\
 \frac{1}{\pi }\int_{\abs{\tilde{z}}<1} g_1(\tilde{\zeta})\frac{1}{2\pi i}
 \int_{\abs{\zeta}=1}\frac{\log(1-z\bar{\zeta})}{(1-\tilde{\zeta}\bar{\zeta})^2}\overline{d\zeta} d\mu(\tilde{\zeta})
 =0
 \end{multline}
 and up to the factor $-1$, the middle term of Eqn.(\ref{neumanmiddle2}) can be written as
 \begin{multline}
 \frac{1}{2\pi i}\int_{\abs{\zeta}=1} (\partial_{z}Tg_1)(\zeta)\frac{d\zeta}{1-\bar{z}\zeta}=
 -\frac{1}{\pi }\int_{\abs{\tilde{\zeta}}<1} g_1(\tilde{\zeta})\frac{1}{2\pi }
 \int_{\abs{\zeta}=1}\frac{1}{(\zeta-\tilde{\zeta})^2}\frac{d\zeta}{1-\bar{z}\zeta}d\mu(\tilde{\zeta})
 =\\
 -\frac{1}{\pi }\int_{\abs{\tilde{\zeta}}<1} g_1(\tilde{\zeta})\partial_\zeta 
 \frac{1}{1-\bar{z}\zeta}|_{\zeta=\tilde{\zeta}} d\mu(\tilde{\zeta})
 =
 -\frac{1}{\pi }\int_{\abs{\zeta}<1} g_1(\zeta) 
 \frac{\bar{z}}{(1-\bar{z}\zeta)^2} d\mu(\zeta)
 \end{multline}
 This completes the proof.
 \end{proof}
Note that there is no single natural way to generalize the Neumann conditions to higher order. One possible way is to introduce powers of $\partial_\nu$ whereas another way is to consider $\partial_\nu(\partial_{\bar{z}})^j,$ $j=0,\ldots,q_1$ in the side conditions. We shall therefore not consider the Neumann problem in the case of higher order here.

 One way the theory of complex boundary value problems\index{Mixed boundary value problems} can be enriched is by introducing mixed-boundary conditions. We give
 here one such example, see Begehr \cite{begehrboletinI}. 
 \begin{theorem}[The Dirichlet-Neumann problem for the inhomogeneous bianalytic equation]
 Let $D=\{\abs{z}<1\}\subset\C$ and $g_0,g_1\in C^0(\partial D,\C),c_0\in \C.$ Consider the equation
 \begin{equation}
 (\partial_{\bar{z}})^2f=g \mbox{ on }D, \quad f=g_0,\quad \partial_\nu (\partial_{\bar{z}}f)=g_1\mbox{ on }\partial D,\quad \partial_{\bar{z}}f(0)=c_0 
 \end{equation}
 for $g$ belonging to the subset of $L^1(D,\C)$ functions which are continuous on the boundary $\partial D.$
 Then there exists a solution $f$ iff for $z\in D$
 \begin{equation}\label{forstacond}
 c_0-\frac{1}{2\pi i}\int_{\partial D} g_0(\zeta)\frac{d\zeta}{1-\bar{z}\zeta}+\frac{1}{\pi}\int_{D}g(\zeta)\frac{1-\abs{\zeta}^2}{\zeta(1-\bar{z}\zeta)}d\mu(\zeta)
 =0
 \end{equation}
 and
 \begin{equation}\label{andracond}
 \frac{1}{2\pi i}\int_{\partial D} (g_1(\zeta)-\bar{\zeta}g(\zeta))\frac{d\zeta}{\zeta(1-\bar{z}\zeta)}+
 \frac{1}{\pi}\int_{D} \frac{\bar{z}g(z)}{(1-\bar{z}\zeta)^2}d\mu(\zeta)
 =0
 \end{equation}
 The solution is then
 \begin{multline}
 f(z)=c_0\bar{z}+\frac{1}{2\pi i}\int_{\partial D} g_0(\zeta)\frac{d\zeta}{\zeta-z}+
 \frac{1}{2\pi i}\int_{\partial D} (g_1(\zeta)-\bar{\zeta} g(\zeta))\frac{1-\abs{z}^2}{z}\log(1-z\bar{\zeta})\frac{d\zeta}{\zeta}
 +\\
 \frac{1}{\pi }\int_{D} g(\zeta)\frac{\abs{\zeta}^2-\abs{z}^2}{\zeta(\zeta-z)} d\mu(\zeta)
 \end{multline}
 \end{theorem}
 \begin{proof}
 The problem can be reformulated as the system 
 \begin{equation}
 \partial_{\bar{z}}f=h \mbox{ on }D, \quad f=g_0\mbox{ on }\partial D 
 \end{equation}
 \begin{equation}
 \partial_{\bar{z}}h=g \mbox{ on }D, \quad \partial_\nu h=g_1\mbox{ on }\partial D,\quad h(0)=c_0 
 \end{equation}
 of an inhomogeneous Dirichlet problem and an inhomogeneous Neumann problem.
 Hence the conditions for existence of unique solutions according to what we have previously done is
 \begin{equation}
 \frac{1}{2\pi i}\int_{\abs{\zeta}=1} g_0(\zeta)\frac{d\zeta}{1-\bar{z}\zeta}=\frac{1}{\pi}\int_{D} \frac{h(\zeta)}{1-\bar{z}\zeta}d\mu(\zeta)
 \end{equation}
 and 
 \begin{equation}
 \frac{1}{2\pi i}\int_{\abs{\zeta}=1} (g_1(\zeta)-\bar{\zeta}g(\zeta))\frac{d\zeta}{(1-\bar{z}\zeta)\zeta}
 +\frac{1}{\pi}\int_{\abs{\zeta}<1} \frac{\bar{z}g(\zeta)}{(1-\bar{z}\zeta)^2}d\mu(\zeta)
 =0
 \end{equation}
 respectively, 
 and the unique solutions are
 \begin{equation}
 f(z)=\frac{1}{2\pi i}\int_{\abs{\zeta}=1} g_0(\zeta)\frac{d\zeta}{\zeta-z}
 +\frac{1}{\pi}\int_{\abs{\zeta}<1}\frac{ h(\zeta)}{\zeta-z}d\mu(\zeta)
 \end{equation}
 and
 \begin{equation}
 h(z)=c_0-\frac{1}{2\pi i}\int_{\abs{\zeta}=1} (g_1(\zeta)-\bar{\zeta}g(\zeta))\log(1-z\bar{\zeta})\frac{d\zeta}{\zeta}
 -\frac{1}{\pi}\int_{\abs{\zeta}<1} \frac{zg(\zeta)}{\zeta(\zeta-z)}d\mu(\zeta)
 \end{equation}
 respectively.
 Using the identities
 \begin{equation}
 \frac{1}{\pi}\int_{\abs{\zeta}<1}\frac{1}{1-\bar{z}\zeta}d\mu(\zeta)=1
 \end{equation}
 \begin{equation}
 \frac{1}{\pi}\int_{\abs{\zeta}<1} \frac{\log(1-\zeta\bar{\zeta})}{1-\bar{z}\zeta}d\mu(\zeta)=
 \frac{1}{\pi}\int_{\abs{\zeta}=1} \log(1-\zeta\bar{\zeta})\frac{d\zeta}{(1-\bar{z}\zeta)\zeta}=0
 \end{equation}
 and
 \begin{equation}
 \frac{1}{\pi}\int_{\abs{\zeta}<1} \frac{1}{1-\bar{z}\zeta}\frac{\zeta}{\tilde{\zeta}-\zeta}d\mu(\zeta)
 =\frac{\abs{\tilde{\zeta}}^2-1}{\tilde{\zeta}(1-\bar{z}\tilde{\zeta})}
 \end{equation}
 the condition of Eqn.(\ref{forstacond}) follows.
 Similarly, using the identities
 \begin{equation}
 -\frac{1}{\pi}\int_{\abs{\zeta}<1}\frac{d\mu(\zeta)}{\zeta-z}=\bar{z}
 \end{equation}
 \begin{equation}
 \frac{1}{\pi}\int_{\abs{\zeta}<1} \frac{\log(1-\zeta\bar{\tilde{\zeta}})}{\zeta-z}d\mu(\zeta)=
 \frac{1-\abs{z}^2}{z}\log(1-z\bar{\tilde{\zeta}})=0
 \end{equation}
 and
 \begin{equation}
 \frac{1}{\pi}\int_{\abs{\zeta}<1} \frac{1}{\zeta-z}\frac{\zeta}{\bar{\tilde{\zeta}}-\zeta}d\mu(\zeta)=
 \frac{\abs{\bar{\tilde{\zeta}}}^2 -\abs{z}^2}{\bar{\tilde{\zeta}}-z}
 \end{equation}
 the condition of Eqn.(\ref{andracond}) follows.
 This completes the proof.
 \end{proof}

\section{The Dirichlet problem for Jordan domains}{The Dirichlet problem 
	for Jordan domains with rectifiable boundary}\label{newsec5}
It is a completely natural step (clearly motivated by the Caratheodory theorem for Jordan domains, or other generalizations of the Riemann mapping theorem) 
to consider domains which are obtained as a conformal image of the unit disc, after having handled the unit disc (we already touched upon this in Remark \ref{greenrem}).
\\
A successful strategy to obtain generalizations based upon our previously described techniques
is to study, 
in the representation formulas for solutions or solvability conditions, where and how the dependence of the domain of the given boundary value problem appears. 
It turns out that in several cases the representations can be written such that the dependence on the domain 
appears in terms of a version of the Green function. Examples of recent results which implement the strategy of using representations of solutions which
appropriately feature the domain dependence can be found in
Aksoy, Begehr \& Celebi \cite{begaksoy1}, \cite{begaksoy3} and Begehr \& Shupeyeva \cite{begaksoy2}.
Here we only give one short simplified example from Begehr \& Shupeyeva \cite{begaksoy2}, for illustration, omitting many details of the proof. 
\\
First let us explain why we believe that considering Jordan domains
suffices to press the methodological point we wish to make.
\begin{remark}\label{dirichharmrem}
	Consider the classical Dirichlet problem for harmonic functions 
	with respect to a Jordan domain, $\Omega$, with boundary $\partial\Omega=:\gamma$ given by a rectifiable curve. 
	Applying the Caratheodory theorem,
	some $\phi$ is chosen as a conformal mapping of the unit disc, $D$, onto $\Omega,$ with continuous extension to the boundary, mapping $\partial D$ onto $\gamma$,
	such that $\phi(\exp(i\theta))$ parametrizes $\gamma$.
	Let $f$ be a Borel function on $\gamma$ such that $f\circ \phi$ is integrable on $\partial D.$ If $w:=\phi^{-1}(z)$ then
	the function $u_f(z):=\frac{1}{2\pi}\int_{0}^{2\pi} f\circ\phi(\exp(i\theta))\frac{1-\abs{w}^2}{\abs{\exp(i\theta)-w}^2} d\theta$
	is harmonic on $\Omega$, satisfying $\lim_{\Omega\ni z\to \zeta}=f(\zeta)$, whenever $\zeta\in \phi^{-1}(\zeta)\in 
	\partial D$ is a point of continuity of $f$. 
	It follows that if $f$ is continuous on $\gamma$ then $u(z)$ is a solution to the Dirichlet problem for Laplace equation on $\Omega.$ 
As usual we may define the function 
	$g(z,\zeta):=\ln\abs{z-\zeta}+h(z,\zeta)$ 
	such that $h(z,\zeta)$ as a function of $z$ coincides with
	$-\ln\abs{z-\zeta}$ for $z\in \gamma$, and such $h$ exists by the above existence of solution for the Laplace equation, where $g$
	is the Green function
	(see Remark \ref{greenrem} and Eqn.(\ref{ekv761}), Eqn.(\ref{ekv762})).
	\\
	Define for each Green function, $g(z,\zeta)$, the function
		\begin{equation}\label{begehrekv}
	h_1(z,\zeta):=\ln\abs{z-\zeta}^2+g(z,\zeta) \, (=\ln\left((z-\zeta)(\bar{z}-\bar{\zeta})\right)+g)
	\end{equation}
	\\
	The strategy 
	is now to formulate an appropriate general representation for solvability conditions for the solution of the polyanalytic Dirichlet problem
	whose dependence of the domain appears primarily in terms of $h_1(z,\zeta)$ and will continue to be a representation for solutions independent of the particular choice of simple, continuous curve $\gamma$, thus independent of $h_1$ (as long as $h_1$ is obtained with respect to a Jordan domain) and ultimately independent of the particular choice of
	Green function.
\end{remark}

\begin{theorem}[Dirichlet problem for the inhomogeneous polyanalytic equation on Jordan domains with rectifiable boundary]
	Let $\Omega$ be a Jordan domain with rectifiable boundary.
	Given $u\in L^p(\Omega)$, $p>2,$ $f_j\in C^0(\partial \Omega,\R),$ $j=0,\ldots ,q-1,$ the problem 
	\begin{equation}
	(\partial_{\bar{z}})^q f=u,\mbox{ on } \Omega, \quad (\partial_{\bar{z}})^j f=f_j,\mbox{ on } \partial \Omega,\quad j=0,\ldots ,q-1
	\end{equation}
	is uniquely solvable for $f\in L^p(\Omega,\C),$ iff for each $0\leq j\leq q-1,$
	\begin{multline}\label{solvcondny}
	\frac{1}{\pi}\int_\Omega u(\zeta_{q-j})\left(\frac{1}{\pi}\int_\Omega\right)^{q-j-1}\partial_{\zeta_1}h_1(z,\zeta_1)\displaystyle{\prod_{l=1}^{q-j-1}}\frac{d\xi_l d\eta_l}{\zeta_l-\zeta_{l+1}}d\xi_{q-j}d\eta_{q-j}=
	\\
	\frac{1}{2\pi i} \int_{\partial \Omega}  
	f_j(\zeta)\partial_{\zeta_1} h_1(z,\zeta)d\zeta+\\
	\sum_{m=j+1}^{q-1} \frac{1}{2\pi i} 
	\int_{\partial \Omega}
	f_m(\zeta_{m-j+1})\left(\frac{1}{\pi}\int_\Omega\right)^{m-j}\partial_{\zeta_1} h_1(z,\zeta_1)
	\displaystyle{\prod_{l=1}^{m-j}}\frac{d\zeta_l d\eta_l}{\zeta_l-\zeta_{l+1}}d\zeta_{m-j+1}
	\end{multline}
	for $j=0,\ldots,q-1,$
	and then the solution is given by
	\begin{equation}\label{diricheqny}
	f(z)=
	\frac{(-1)^q}{\pi} \int_{\Omega}  \frac{u(\zeta)}{(q-1)!}
	\frac{(\overline{\zeta}-\overline{z})^{q-1}}{\zeta-z}d\zeta  
	+\\
	\sum_{j=0}^{q-1} \frac{(-1)^j}{2\pi i} \int_{\partial\Omega} \frac{\gamma_j(\zeta)}{j!}
	\frac{(\overline{\zeta}-\overline{z})^j}{\zeta-z}d\zeta 
	\end{equation}
\end{theorem}
\begin{proof}
	The proof is by induction in $q$. First consider the base case $q=1$ ($f=f_0$). 
Suppose we have the condition given by Eqn.(\ref{solvcondny}), namely
$0=\frac{1}{2\pi i}\int_{\partial\Omega} f_0(\zeta)
\partial_{\zeta} h_1(z,\zeta)d\zeta$
$-\frac{1}{\pi}\int_\Omega$
$u(\zeta)\partial_{\zeta} h_1(z,\zeta)d\xi d\eta=:\tilde{f}$.
Then we may add this to the expression (which is the right hand side off Eqn.(\ref{diricheqny}) for $q=1$) 
	$\frac{1}{2\pi i}\int_{\partial\Omega}f_0(\zeta)\frac{d\zeta}{\zeta -z}-\frac{1}{\pi}\int_\Omega u(\zeta)\frac{d\zeta d\nu}{\zeta-z}$ 
	in order to obtain
	\begin{equation}
	\frac{1}{2\pi i}\int_{\partial\Omega}f_0(\zeta)\left(\frac{1}{\zeta -z}+\partial_\zeta g(z,\zeta)\right)d\zeta
	-\frac{1}{\pi}\int_\Omega u(\zeta)\left(\frac{1}{\zeta -z}+\partial_\zeta g(z,\zeta)\right)d\xi d\eta
	\end{equation}
	and writing
	$\frac{1}{\zeta -z}+\partial_\zeta h(z,\zeta)$
	as $\partial_\zeta(\ln\abs{\zeta-z}+h(z,\zeta))=\partial_\zeta g(z,\zeta)$
	renders a solution
	to the boundary value problem. 
	Conversely, if a solution 
	$f(z)$ is given, 
	it necessarily has a representation according to the Cauchy-Pompieu formula
  i.e.\ be written in the form, see e.g.\ Eqn.(\ref{cauchypompformulaendim}),
$f(z)=\frac{1}{2\pi i}\int_{\partial\Omega}f_0(\zeta)\frac{d\zeta}{\zeta -z}-\frac{1}{\pi}\int_\Omega u(\zeta)\frac{d\zeta d\nu}{\zeta-z}$, and again 
adding the side condition will
also give a solution, 
and the difference has zero boundary values
thus condition of Eqn.(\ref{solvcondny}) follows. 
Uniqueness of the solution with respect to boundary values follows from the uniqueness with respect to boundary values for the homogenous Dirichlet problem which in turn follows from Theorem \ref{theodorescuthm},
since the difference of two solutions to the Dirichlet problem for the inhomogeneous polyanalytic equation on Jordan domains with rectifiable boundary, is a $q$-analytic function with zero boundary values.
	For the induction step, assume $q\geq 2$ and that 
	the result holds true for the case up to and including $q$, 
	and consider the case $q+1$.
	As in the case of the proof of Theorem \ref{dirichletpolyettan} 
	the problem can be reformulated as the system
	\begin{equation}
	(\partial_{\bar{z}})^{q} f=w,\mbox{ on } \Omega, \quad (\partial_{\bar{z}})^j f=f_j,\mbox{ on } \partial \Omega, j=0,\ldots ,q-1
	\end{equation}
	\begin{equation}
	\partial_{\bar{z}} w=f,\mbox{ on } \Omega, \quad \partial_{\bar{z}} w=\gamma_{q-1},\mbox{ on } \partial \Omega
	\end{equation}
	Writing $w(z)$ according to
	\begin{equation}
	w(z)=\frac{1}{2\pi i}\int_{\partial\Omega} f_q(\zeta_{q-j+1})\frac{d\zeta_{q-j+1}}{\zeta_{q-j+1}-z}-\frac{1}{\pi}\int_\Omega u(\zeta_{q-j+1})
	\frac{d\xi_{q-j+1}d\eta_{q-j+1}}{\zeta_{q-j+1}-z}
	\end{equation}
	and applying the solvability condition
	 \begin{equation}\label{tolvankal}
	 \frac{1}{2\pi i}\int_{\partial\Omega} f_q(\zeta)\partial_\zeta h_1(z,\zeta) d\zeta=\frac{1}{\pi}\int_\Omega u(\zeta)h_1(z,\zeta)d\xi d\eta
	 \end{equation}
	 we may right the left hand side of Eqn.(\ref{solvcondny}), with $u$ replaced by $w$, according to
	 \begin{multline}%
	 \frac{1}{\pi}\int_\Omega\left(
	 \frac{1}{2\pi i} \int_{\partial \Omega}  
	 f_q(\zeta_{q-j+1})\frac{d\zeta_{q-j+1}}{\zeta_{q-j+1}-\zeta_{q-j}}\right)\times\\
	 \left(\frac{1}{\pi}\int_\Omega \partial_{\zeta_1} h_1(z,\zeta_1)\right)^{q-j-1}
	 \displaystyle{\prod_{l=1}^{q-j-1}}\frac{d\xi_l d\eta_l}{\zeta_l-\zeta_{l+1}}d\xi_{q-j}d\xi_{q-j}d\eta_{q-j}-\\
	 \frac{1}{\pi}\int_\Omega\left(\frac{1}{\pi}u(\zeta_{q-j+1})\frac{d\xi_{q-j+1} d\eta_{q-j+1}}{\zeta_{q-j+1}-\zeta_{q-j}}\right)\times
	 \\
	 \left(\frac{1}{\pi}\int_\Omega\right)^{q-j-1} \partial_\zeta h_1(z,\zeta_1)
	 \Pi_{l=1}^{q-j-1}\frac{d\zeta_l d\eta_l}{\zeta_l-\zeta_{l+1}}d\xi_{q-j}d\eta_{q-j}
	\end{multline}
	 The latter in turn can be rewritten as
	 \begin{multline}\label{sistaskiten}
	 -\frac{1}{2\pi i} \left(\int_{\partial \Omega}  
	 f_q(\zeta_{q-j+1})\left(\frac{1}{\pi}\int_\Omega\right)^{q-j}\partial_{\zeta_1} h_1(z,\zeta_1)\times\right.
	 \\
	  \left.\displaystyle{\prod_{l=1}^{q-j-1}}\frac{d\xi_l d\eta_l}{\zeta_l-\zeta_{l+1}}\frac{d\xi_{q-j} d\eta_{q-j}}{\zeta_{q-j}-\zeta_{q-j}}\right)  d\zeta_{q-j+1}
	 +\\
	 \frac{1}{\pi}\int_\Omega\left(f(\zeta_{q-j+1})\left(\frac{1}{\pi}\int_\Omega\right)^{q-j}\partial_{\zeta_1} h_1(z,\zeta_1)\times\right.\\
	 \left.
	 \displaystyle{\prod_{l=1}^{q-j-1}}\frac{d\xi_l d\eta_l}{\zeta_l-\zeta_{l+1}}
	 \frac{d\xi_{q-j} d\eta_{q-j}}{\zeta_{q-j}-\zeta_{q-j+1}}\right)d\xi_{q-j+1}d\eta_{q-j+1}
	 \end{multline} 
	 Replacing in the left hand side of Eqn.(\ref{solvcondny}) the area integral, with $u$ replaced by $w$,
	 with the use of Eqn.(\ref{sistaskiten}) as a substitution for Eqn.(\ref{tolvankal}), renders the
	 condition of Eqn.(\ref{solvcondny}) for the case $q+1$. This completes the proof sketch.	 
\end{proof}

\section{Domains with rational boundaries}
Let $\Omega$ be a bounded simply connected domain bounded by a closed simple rectifiable Jordan curve, $\gamma$. Assume 
$0\in \Omega.$ Let $q\in \Z_+$ and denote by $\mathcal{A}(\overline{\Omega})$ the set of complex analytic functions on $\Omega$ that 
have bounded angular boundary values on $\gamma.$ For $p>0$ denote by $E_p(\overline{\Omega})$ the {\em Smirnov class}\index{Smirnov class} of
functions holomorphic on $\Omega$ such that there exists a sequence
of closed rectifiable Jordan curves
$\gamma_k(f)\subset\Omega$, $k\in \Z_+$, where each $\gamma_k(f)$ bounds a bounded domain $V_k$ such that:\\
(1) $V_1(f)\subset \cdots V_k(f)\subset\Omega,$ $\bigcup_{k=1}^\infty V_k(f)=\Omega.$\\ 
(2) $\sup_k\left(\int_{\gamma_k(f)}\abs{f(z)}^p\abs{dz}\right)<\infty.$
This definition of $E_p$ is equivalent to the condition that there is a conformal onto 
map $\phi:\{\abs{z}<1\}\to \Omega$ such that $\phi'(w)$ is an {\em outer function}\index{Outer function}
which by definition means that it can be written in the form
$\phi'(w)=C\exp\left(\int_{\{\abs{\zeta}=1\}}\frac{\zeta+w}{\zeta-w}\log v(\zeta)d\zeta\right),$
for some $C,v$ with $\abs{C}=1,$ $v>0$ a.e.\ on $\{\abs{\zeta}=1\}$ and $\log v\in L^1(\{\abs{\zeta}=1\}),$
see e.g.\ Garnett \cite{garnettbok}.
We say that a function $f$ that is $q$-analytic on $\Omega$ belongs to the class $E_p^{q-1}(\overline{\Omega})$
if each analytic component of $f$ belongs to $E_p(\overline{\Omega}).$ If each analytic component of $f$
belongs to $\mathcal{A}(\overline{\Omega})$ then we say that $f$ belongs to the class $\mathcal{A}^{q-1}(\overline{\Omega})$.
Note that 
\begin{equation}
E_p^{q-1}(\overline{\Omega})\subset \mathcal{A}^{q-1}(\overline{\Omega})
\end{equation}
A curve $\lambda_1$ given by $z=\lambda_1(t)=\lambda_{11}(t)+i\lambda_{12}(t),$ parametrized by  $\alpha\leq t\leq \beta$, 
will be called an {\em analytic curve} if the function $\lambda_1 (t)$ is an analytic
function of the variable $t$ such that $\lambda_1'(t)\neq 0,$ $\alpha\leq t\leq \beta$.
Recall that for every simple analytic arc, $\gamma$, there exists a (unique) Schwarz function $G(z)$ 
that is analytic in some open neighborhood $U_\lambda$ of $\lambda$ such that
$\bar{z}=G(z)$ on $\lambda$ (see Proposition \ref{schwarzfuncprop}).
If a closed analytic curve $\gamma$ is the boundary of a finite simply-connected domain $\Omega$, then the Schwarz function $\gamma$ 
inside the domain $\Omega$ must have at least one singular point, see Section \ref{schwarzsec}.
\begin{definition}
Suppose that the boundary $\gamma$ of the domain $\Omega$ contains a 
simple analytic arc $\gamma_1$ with Schwarz function $G_1(z)$ (in particular $\gamma_1$ is defined by $\bar{z}-G_1(z)$, where 
$G_1$ is analytic on a neighborhood $U$ of $\gamma_1$), and let $f(z)$ be a $q$-analytic function on $\Omega.$
By definition $f$ can be represented in the form $\sum_{j=0}^{q-1} a_j(z)\bar{z}^j$ on $\Omega.$ 
Then we have on $\Omega\cap U$ a holomorphic function
\begin{equation}
f_{\gamma_1} (z) = \sum_{j=0}^{q-1} a_j(z)(G_1(z))^j
\end{equation}
that coincides with $f$ on $\gamma_1.$
$f_{\gamma_1}$ will be called {\em coherent with $f(z)$ on $\gamma_1$} (see Balk \cite{ca1}, Section 4.1). 
\end{definition}
Let $z=\lambda(w),$ $w=\xi+i\eta$, be a conformal univalent map of the unit disc in $\C$ to $\Omega$ such that
$\lambda(0)=0,$ $\lambda'(0)>0$. 
It is known that if $\Omega$ is bounded by a simple closed analytic curve with Schwarz function $G(z)$, 
then in the notation $z=\lambda(w),$ $w=\lambda^{-1}(z),$ we have
\begin{equation}\label{rasulovsix}
G(z)=\overline{\lambda(1/\lambda^{-1}(z))}
\end{equation}
\begin{lemma}\label{rasulovlem1}
Let $\Omega$ be a bounded simply connected domain bounded by a simple closed analytic curve $\gamma$ 
with Schwarz function $G(z)$. If $G (z)$ is meromorphic in $\Omega$ then 
$z+\lambda(w)$, which is a conformal mapping of $\{w:\abs{w}<1\}$ onto $\Omega$, is rational. 
\end{lemma}
\begin{proof}
Suppose $G(z)$ is meromorphic on $\Omega$. Then $G(z)=\phi(z)/\varphi(z)$, for holomorphic $\phi,\varphi$ on $\Omega.$ 
By Eqn.(\ref{rasulovsix}), we have for $\abs{w}=1$,
\begin{equation}\label{rasulovseven}
\overline{\lambda(w)}=\frac{\phi(\lambda(w))}{\varphi(\lambda(w))}=\overline{\lambda(1/w)}
\end{equation}
Set
\begin{equation}
B(w):=\left\{
\begin{array}{ll}
\frac{\phi(\lambda(w))}{\varphi(\lambda(w))} & ,\abs{w}\leq 1\\
\overline{\lambda(1/w)} & ,\abs{w}>1
\end{array}
\right.
\end{equation}
where $\frac{\phi(\lambda(w))}{\varphi(\lambda(w))}$ is
meromorphic in $\{\abs{w}<1\}$ and $\overline{\lambda(1/w)}$ is holomorphic in $\{ 1 < \abs{w} \leq \infty \}$ and vanishes at infinity 
since $\lambda(0)= 0$. Hence $B$ has only a finite number of poles, thus is a rational function, 
and so there exist polynomials $P_l(w),Q_m(w)$, $\mbox{deg}(P_l)=l,$ $\mbox{deg}(Q_m)=m,$
such that $B$ can be written as $P_l/Q_m$. Since $B$ vanishes at infinity we have $l<m.$ 
By the definition of $B$ together with Eqn.(\ref{rasulovseven}), there exists polynomials $c_l,d_m$ such that on the circle $\{\abs{w}=1\}$
we have
\begin{equation}
\lambda(w)=\frac{\overline{P}_l(1/w)}{\overline{Q}_m(1/w)}=w^{m-l}\frac{c_l(w)}{d_m(w)}
\end{equation}
This implies that for $\abs{w}<1$ we have $\lambda(w)=w^{m-l}c_l(w)/d_m(w).$
This completes the proof.
\end{proof}
\begin{remark}
Note that the proof shows that if the Schwarz function $G(z)$ of the curve $\gamma$ is itself rational, i.e.\ $G(z)=P(z)/Q(z)$ for polynomials $P,Q,$
then necessarily $\mbox{deg}P<\mbox{deg}Q.$ On the other hand, it is known (see Section \ref{schwarzsec}) that if $G(z)$ is meromorphic on $\Omega$, 
then the curve $\gamma$ is algebraic, that is, it can be implicitly given by the equation $R(z,\bar{z})= 0$, 
for a polynomial $R$ in $z,\bar{z}.$ 
\end{remark}
We now consider the following Dirichlet problem (see Rasulov \cite{rasulov1995}):\\
$(\star)$ {\em 
Let $q\geq 2$ and assume $\gamma$ is  a simple analytic curve with Schwarz function $R(z)$, such that $R$ is rational in the sense that
$R(z)=P(z)/Q(z),$ $P,Q$ polynomials that have no common divisors and $Q(z)$ has no zeros outside $\Omega.$ 
Given a function $g\in L^1(\gamma)$, find all functions $f\in \mathcal{A}^{q-1}(\overline{\Omega})$ such that
$f=g$ a.e.\ on $\gamma.$}
\\
since the problem is linear, the general solution will have the form $f(z) = f_0(z) + f_1(z)$, where $f_0(z)$ 
is the general solution of the corresponding homogeneous problem ($g\equiv 0$), and $f_1(z)$ 
is any particular solution of the inhomogeneous problem. By Lemma \ref{rasulovlem1}
we know that $\Omega$ is the rational image of the unit disc.
\begin{lemma}\label{rasulovlem21}
The homogeneous problem $(\star)$ (i.e.\ with $g\equiv 0$) has an infinite set of linearly independent solutions 
of the form
\begin{equation}\label{jojorasulov}
f_0(z)=(Q(z)\bar{z}-P(z))f_{q-1}(z), \quad f_{q-1}\in  \mathcal{A}^{q-1}(\overline{\Omega})
\end{equation}
\end{lemma}
\begin{proof}
Let $q\in \Z_+.$
It is clear that a.e.\ on $\gamma$ we have $Q(z)/P(z) =\bar{z},$ hence any function
in Eqn.(\ref{jojorasulov}) is a solution to the homogeneous version of the problem $(\star).$ Suppose that
$\tilde{f}$ is a $q$-analytic function on $\Omega$ with representation 
\begin{equation}\label{rasulov13}
\tilde{f}(z)=\sum_{j=0}^{q-1} \tilde{a}_j(z)\bar{z}^j
\end{equation}
such that $\tilde{f}\in \mathcal{A}^{q-1}(\overline{\Omega})$ is a solution to the 
homogeneous version of the problem $(\star),$ i.e.\ $\tilde{f}$ has a.e.\ zero angular boundary values.
Since $\gamma$ has Schwarz function $R(z)$ the function {\em coherent} with $\tilde{f}(z)$ on $\gamma$ has the form
\begin{equation}
\tilde{f}_\gamma(z)=\sum_{j=0}^{q-1} \tilde{a}_j(z)(R(z))^j =\sum_{j=0}^{q-1} \tilde{a}_j(z)\left(\frac{P(z)}{Q(z)}\right)^j
\end{equation}
Hence $\tilde{f}_\gamma$ is holomorphic on $\Omega$ except possibly at the finitely many zeros of $Q(z)$, i.e.\
$\tilde{f}_\gamma$ is meromorphic on $\Omega.$ 
Since the functions $\tilde{f}$ and $\tilde{f}_\gamma$ have the same angular boundary values a.e.\ on $\gamma$ 
(see Balk \cite{ca1}, Sec 4.3), we have $\tilde{f}=\tilde{f}_\gamma=0$ a.e.\ on $\gamma.$ 
By the Luzin-Privalov uniqueness theorem for 
meromorphic functions (Theorem \ref{luzinprivalovmero1})
this implies 
$\tilde{f}_\gamma(z)\equiv 0$ on $\Omega$ which in turn implies
\begin{equation}\label{rasulov16}
\tilde{a}_0(z)\equiv -\sum_{j=1}^{q-1} \tilde{a}_j(z)\left(\frac{P(z)}{Q(z)}\right)^j
\end{equation}
Since $P,Q$ have no common divisors, this implies by iteration for $\nu=1,\ldots,q-1$
\begin{equation}\label{rasulov17}
\sum_{j=1}^{q-1} \tilde{a}_j(z)\left(\frac{P(z)}{Q(z)}\right)^j =Q(z)\psi_\nu(z) ,\quad z\in \Omega
\end{equation}
where the $\psi_\nu$ are holomorphic functions on $\Omega$ which belong to $\mathcal{A}^0(\overline{\Omega}).$
Now Eqn.(\ref{rasulov16}) inserted into Eqn.(\ref{rasulov13}) gives 
\begin{multline}\label{rasulov18}
\tilde{f}(z)=\\
(Q(z)\bar{z}-P(z))\left(\frac{1}{Q(z)}\sum_{k=0}^{q-1}\sum_{j=0}^{k-1} \bar{z}^{k-j-1}\left(\frac{P(z)}{Q(z)}\right)^j \tilde{a}_k(z)\right)
\end{multline}
The second factor on the right hand side is a $(q-1)$-analytic function on 
$\Omega$ belonging to $\mathcal{A}^{q-1}(\overline{\Omega}).$ Thus the $q$-analytic
function $\tilde{f}$ in Eqn(\ref{rasulov18}) has the form given by Eqn.(\ref{jojorasulov}).
This completes the proof.
\end{proof}
Hence in order to solve the inhomogeneous version of $(\star)$ (supposing the existence of a solution), 
it suffices, by Lemma \ref{rasulovlem21}, 
to find a particular solution.  
Suppose that $g\in L^1(\gamma)$ and that a $q$-analytic function, $F(z)=\sum_{j=0}^{q-1}\bar{z}^j a_j(z)$, on $\Omega$
is a solution of the inhomogeneous problem in the sense that $F=g$ a.e.\ $\gamma$ and $F\in \mathcal{A}^{q-1}(\overline{\Omega}).$   
As before introduce the function $F_\gamma$ ({\em coherent} as previously defined with $f$ with respect to $\gamma$) according to
\begin{equation}
F_\gamma(z)=\sum_{j=0}^{q-1}\left(\frac{P(z)}{Q(z)}\right)^j a_j(z)
\end{equation}
Note that $\left(Q(z)\right)^{q-j} F_\gamma(z)$ is holomorphic on $\Omega$ and belongs to 
$\mathcal{A}^0(\overline{\Omega}).$ It is further known (see Balk \cite{ca1}, Sec 4.3) that
$\lim_{z\to t\in \gamma} F(z)=\lim_{z\to t\in \gamma}F_\gamma (z)$, hence $F=g$ a.e.\ $\gamma$ gives (taking the limit
$\lim_{z\to t\in \gamma}$ to mean that $z$ approaches $t$ nontangentially)
\begin{equation}
\lim_{z\to t\in \gamma}((Q(z))^{q-1}F_\gamma(z))= (Q(z))^{q-1}g(t)
\end{equation}
This implies, by Theorem \ref{golubev}, that the function $(Q(z))^{q-1}g(t)$ is the boundary function of
a holomorphic function thus 
\begin{equation}
\int_\gamma (Q(t))^{q-1}g(t)t^pdt=0,\quad p\in \Z_{\geq 0}
\end{equation}

\begin{theorem}
Let $q\geq 2$ and assume $\gamma$ is  a simple analytic curve with Schwarz function $R(z)$, such that $R$ is rational in the sense that
$R(z)=P(z)/Q(z),$ $P,Q$ polynomials that have no common divisors and $Q(z)$ has no zeros outside $\Omega.$ 
For the solvability of the inhomogeneous problem $(\star)$, i.e.\ that given a function $g\in L^1(\gamma)$, finding 
all functions $f\in \mathcal{A}^{q-1}(\overline{\Omega})$ such that
$f=g$ a.e.\ on $\gamma$, it is necessary and sufficient that
\begin{equation}\label{rasulov23}
\int_\gamma (Q(t))^{q-1}g(t)t^pdt=0,\quad p\in \Z_{\geq 0}
\end{equation}
and the latter is equivalent to
\begin{equation}\label{rasulov25}
g(t)=G^+(t)+ (Q(t))^{1-q}\sum_{k=0}^\nu c_kt^k
\end{equation}
where $G^+$ is the boundary value of an arbitrary analytic function on $\Omega$ of class $E_1$, $c_0,\ldots,c_\nu$ are 
arbitrary complex constants, $\nu=m(q-1)-1$ and $m$ is the degree of the polynomial $Q (t).$
\end{theorem}
(Note that the number of solutions of the homogeneous problem is not necessarily
finite and neither are the number of solvability conditions)
\begin{proof}
By what we have previously done Eqn.(\ref{rasulov23}) is necessary for the solution of $(\star)$.
Next note that Eqn.(\ref{rasulov23}) is equivalent to the condition that for $t\in \gamma$
\begin{equation}\label{rasulov24}
-\frac{(Q(z))^{q-1}g(t)}{2}+\frac{1}{2\pi i} \int_\gamma \frac{(Q(t))^{q-1}g(t)}{\tau -t}d\tau =0
\end{equation}
This can be written
\begin{multline}
(Q(z))^{q-1}\left( -\frac{g(t)}{2}+\frac{1}{2\pi i} \int_\gamma \frac{g(t)}{\tau -t}d\tau+\right.\\
\left.\frac{(Q(z))^{1-q}}{2\pi i} \int_\gamma \frac{(Q(\tau))^{q-1}-(Q(t))^{q-1}}{\tau -t}d\tau\right)
\end{multline}
Since $(Q(z))^{q-1}$ is a polynomial that is zero-free on $\Omega$, the last equation is equivalent to
the so-called {\em degenerate} singular integral equation 
(see Gakhov \cite{gakovartikel}; the corresponding characteristic equation is obtained by removing the term with $K(t,\zeta)$, and this equation can be solved by the
Plemelj-Sokhotsky-formula for limits of Cauchy integrals at points 
on the path of integration, see Gakhov \cite{gakhovartikel2})  
\begin{equation}\label{rasulov26}
-\frac{g(t)}{2}+\frac{1}{2\pi i} \int_\gamma \frac{g(t)}{\tau -t}d\tau+\\
\frac{1}{2\pi i} \int_\gamma K(t,\tau)g(\tau)d\tau =0
\end{equation}
for
\begin{multline}\label{rasulov27}
K(t,\tau)=(Q(z))^{1-q}\left( \frac{(Q(\tau))^{q-1}-(Q(t))^{q-1}}{\tau -t}\right) =\\
(Q(z))^{1-q}\sum_{k=1}^{m(q-1)}q_k \sum_{p=1}^{k-1} t^{k-p-1}\tau^p
\end{multline}
where the $q_k$ are the coefficients of the polynomial $(Q(t))^{q-1}$. This implies that for $\nu=m(q-1)-1$, we have
\begin{equation}\label{rasulov28}
-\frac{g(t)}{2}+\frac{1}{2\pi i} \int_\gamma \frac{g(t)}{\tau -t}d\tau+\\
(Q(z))^{1-q}\sum_{k=0}^{\nu} c_k t^k =0
\end{equation}
where
\begin{equation}
c_k:=\sum_{p=k}^\nu q_{p+1}a_p, a_p:=\frac{1}{2\pi i}\int_\gamma \tau^p g(\tau)d\tau
\end{equation}
Since the degree of the polynomial $(Q(z))^{q-1}$ is $m(q-1)$ where $m$ is the degree of $Q(t)$, $\nu<m(q-1)$, the expression
\begin{equation}
-\frac{g(t)}{2}+\frac{1}{2\pi i} \int_\gamma \frac{g(t)}{\tau -t}d\tau
\end{equation} 
is the boundary value of a holomorphic function, $F^-(z)$, on $\Omega^-:=\hat{C}\setminus\overline{\Omega}$
thus a.e.\ on $\gamma$ we have
\begin{equation}\label{rasulov30}
(Q(z))^{1-q} \sum_{k=0}^{\nu} c_k t^k =F^-(t),\quad t\in \gamma
\end{equation} 
so Eqn.(\ref{rasulov28}) takes the form
\begin{equation}\label{rasulov31}
-\frac{g(t)}{2}+\frac{1}{2\pi i} \int_\gamma \frac{g(t)}{\tau -t}d\tau=-F^-(t)
\end{equation}
This singular integral equation is solvable (see 
e.g.\ \cite{gakovartikel},
\cite{gakhov}; the solution can be obtained by a certain regularization method, which consists of an application of a suitable operator to the given equation in order to yield a Fredholm equation).
and the general 
solution $F^-(t)$ takes the form
\begin{equation}
\tilde{g}(t)=G^+(t)+F^-(t),\quad t\in \gamma
\end{equation}
where $G^+$ is given by the angular boundary values of an arbitrary analytic function on $\Omega$, belonging
to $E_1.$ Replacing $g$ by $\tilde{g}$ in Eqn.(\ref{rasulov26}) gives
\begin{equation}\label{rasulov33}
-F^-(f)+\frac{1}{2\pi i}\int_\gamma K(t,\tau)(G^+(\tau)+F^-(\tau))d\tau=0
\end{equation}
By the expression for $K(t,\tau)$ given in Eqn.(\ref{rasulov26})
this yields
\begin{equation}
\frac{1}{2\pi i}\int_\gamma K(t,\tau)G^+(\tau)d\tau=0
\end{equation}
\begin{equation}
H^-(t):=\frac{1}{2\pi i}\int_\gamma K(t,\tau)F^-(\tau)d\tau=0
\end{equation}
where $H^-$ is a function of the form given by Eqn.(\ref{rasulov30}) for some $c_0,\ldots,c_\nu.$
This implies that Eqn.(\ref{rasulov33}) takes the form 
$-F^-(t)+H^-(t)=0,$ $z\in \gamma$
hence by boundary uniqueness 
\begin{equation}
-F^- (t)+H^- (t)\equiv 0,\mbox{ on } \Omega
\end{equation}
This proves that Eqn.(\ref{rasulov24}) is equivalent to Eqn.(\ref{rasulov23}) in the theorem.
It remains to show that if $g$ satisfies Eqn.(\ref{rasulov25})
then the problem $(\star)$ has at least one nonzero solution. Given Eqn.(\ref{rasulov23})
suppose we solve the following equation for a polynomial $\psi_{q-1}(z)$
\begin{equation}\label{rasulov35}
(P(z))^{q-1}\psi_{q-1}(z)=\sum_{k=0}^\nu c_k z^k
\end{equation}
where $P$ is the numerator of a rational function $R(z)=P(z)/Q(z).$ This renders a $q$-analytic
function on $\Omega$ of the form
\begin{equation}
F_1(z)=\frac{1}{2\pi i}\int_\gamma \frac{G^+(\tau)}{\tau -z}d\tau +\bar{z}^{q-1}\psi_{q-1}(z)
\end{equation}
where we have, since $R$ is a Schwarz function, that for each $t\in \gamma$, $\bar{t}=P(t)/Q(t)$ thus Eqn.(\ref{rasulov35})
gives a.e.\ on $\gamma$
\begin{multline}
\lim_{z\to t\in \gamma}\left(\frac{1}{2\pi i}\int_\gamma \frac{G^+(\tau)}{\tau -z}d\tau +\bar{z}^{q-1}\psi_{q-1}(z)\right)=\\
=G^+(t)+\bar{t}^{q-1}\psi_{q-1}(t)=G^+(t)+\bar{t}^{q-1}(Q(z))^{1-q}\sum_{k=0}^\nu c_k t^k
\end{multline}
Hence $F_1$ is a nontrivial solution to $(\star)$. Since the $\mbox{deg}(P (z))<\mbox{deg}(Q (z))$
we have $\mbox{deg}(P)^{q-1})<\mbox{deg}(Q )^{q-1}$ and consequently
for any polynomials $P (z), Q (z)$ with this property,
there is always a polynomial $\psi_{q-1}$ satisfying Eqn.(\ref{rasulov35}) 
This completes the proof.
\end{proof}

On the other hand, already in the case of $q=2$ it is known that the boundary being rational is, in some sense, necessary.

\begin{definition} 
Let $\Omega\subset\C$ be a simply connected domain
Given a continuous function $f\in \partial\Omega$, consider the problem of finding solutions $u$ to the 
Dirichlet problem $\partial_{\bar{z}}^2 u=0,$ $u|_{\partial\Omega}=f$, such that $u\in E_2.$
The problem is called {\em Hausdorff} if
the condition that for every solution $v\in E_2$ to the homogeneous problem (i.e.\ with $f$ replaced by $0$)
\begin{equation}
\int_\gamma u\partial_{\bar{z}} \bar{v} dz =\int_\gamma f \partial_{\bar{z}} \bar{v} dz=0
\end{equation}
is necessary and sufficient for the solvability of the given inhomogeneous Dirichlet problem.
\end{definition}
\begin{definition}
A simple continuous curve, $\gamma$, in the plane is called a {\em Lyapunov curve}\index{Lyapunox curve} if it has a well-defined tangent at each point and
there exist $c>0,$ $\alpha\in (0,1]$ such that for any pair $x,y\in \gamma,$ we have
$\theta(x,y)<c\abs{x-y}^\alpha$, where $\theta(x,y)$ denotes the angle between $x$ and $y.$
\end{definition}
Hop \cite{hop} proved the following.
\begin{theorem}[Hop \cite{hop}]
Let $\Omega\subset\C$ be a simply connected domain and let $\phi$ denote a conformal mapping of 
the unit disc onto $\Omega$ such that $\partial\Omega$ is given by an analytic curve. Consider
the Dirichlet problem of finding $u\in E_2$ such that $\partial_{\bar{z}}^2 u=0,$ $u|_{\partial\Omega}=f$, for a given continuous $f$. 
If $f\equiv 0$ then the problem has a nonzero solutions if and only if $\phi$ is rational.
In the inhomogeneous case when $\gamma$ is given by a Lyapunov curve the problem is Hausdorff if and only if $\phi$ is rational.
\end{theorem}
  
\section{Some comments on related problems and generalizations}\label{concludrem} 

\subsection{An example of a case of a multiply connected domains}
Yet another way the theory of complex boundary value problems can be enriched is by considering
multiply connected domains. 
See the survey of Aksoy \& Celebi \cite{aksoy} for a generalization of results for the unit disc to multiply connected domains. 
Here we only 
give an example and related technique.
Let $S^+$ be an $m$-connected domain bounded by $m$ closed disjoint Lyapunov curves
$L_1,\ldots,L_m,$ where $L_1$ contains all the $L_j$, $j=2,\ldots,m$ and denote $L=\cup_{j=0}^m.$
Let $\Omega_j^+$ be bounded simply connected
domains bounded by non-intersecting contours $\gamma_j,$ $j=1,\ldots,m,$ respectively.
Let $w(z,\bar{z})$ be polyanalytic on $\Omega^+$
and $w_j$ polyanalytic on $\Omega_j^+,$ $j=1,\ldots,m,$ for some $m\in \Z_+.$
We suppose that $L$ and $\gamma_j$ are traversed in the 
positive direction with respect to their interiors $S^+,$ $\Omega_j^+.$
The complement of $\Omega^+\cup\Gamma_j$ is denoted $\Omega^-$,
for $j=2,\ldots,m,$ the simply connected domain enclosed by $L_j$ is denoted $S_j^-$
and for $j=1$ $S_1^-$ denotes the complement of the union of $L_1$ with the
bounded domain it encloses.
\\
Let $\alpha_{-1}^j(t)$ denote the sense preserving functions on $\gamma_j$ 
such that it transforms $\gamma_j$ homeomorphically into some closed contour belonging to $L$;
$\alpha_{-1}^j(\gamma_j)\cap \alpha_{-1}^k(\gamma_k)\neq \emptyset$, $j,k=1,\ldots,m$ and
such that each $a_{-1}^j$ has continuous nonzero derivatives at points of $\gamma_j.$
\\
Define
\begin{equation}
\omega(t,\gamma_k)=\delta_{k,j},\quad t\in \gamma_j,k,j=1,\ldots,m
\end{equation}
where $\delta_{j,k}=1$ for $j=k$ and $\delta_{j,k}=0$ otherwise.
Define on $L$ 
\begin{equation}
\alpha(t)=\omega(t,\gamma_j)\alpha^j(t),\quad t\in L_j,j=1,\ldots,m
\end{equation}
Consider the following boundary value problem:
\\
\\
$(\star)$ Find $w$ polyanalytic in $S^+$ and $w_j$ polyanalytic in $\Omega_j^+$ respectively such that on the corresponding contours
\begin{equation}
\partial_z^k w_j(\alpha(t),\overline{\alpha(t)})=A_k(t)\overline{\partial_z^k w(t,\bar{t})}+B_k(t)
\end{equation}
$k=0,\ldots,q-1,$ $j=1,\ldots,m,$
where $A_k(t),B_k(t)$ are H\"older continuous on $L$ and $A_k(t)\neq 0$ on $L.$
\\
\\
We
make the Ansatz
\begin{equation}
w(z,\bar{z})=\sum_{\nu}^{q-1} \bar{z}^\nu U_\nu(z),\quad z\in S^+
\end{equation}
\begin{equation}
w_j(z,\bar{z})=\sum_{\nu}^{q-1} \bar{z}^\nu U_{\nu,j}(z),\quad z\in D_j^+
\end{equation}
for holomorphic functions $U_{\nu}(z),U_{\nu,j}(z),$
The boundary conditions of problem $(\star)$ can be rewritten as the system
\begin{equation}
U_{p,j}(\alpha(t))=A_p(t)\overline{U_p(t)} +g_{p,j}(t)
\end{equation}
for $p=0,\ldots,q-1,$ $j=1,\ldots,m$, where
\begin{multline}
g_{p,j}(t)=\\
\sum_{i=1}^{n-p-1} \binom{p+i}{j} 
\left(A_p(t)tî\overline{P_{p+i}(t)}-\overline{\alpha(t)}^i U_{p+i,j}
(\alpha(t)) \right)+\frac{B_p(t)}{p!}
\end{multline}
Assume $0\in S^+$ and define $A_{0,p}(t):=\bar{t}^{\kappa_p} A_p(t)$, $t\in L,$
and denote by $X_{0,p,j}(z), X_{0,p}(z)$ 
the solution to the boundary value problem
\begin{equation}
X_{0,p,j}(\alpha(t))=A_{0,p}(t) \overline{X_{0,p}(t)}
\end{equation}
so that we can represent
\begin{equation}
A_{p}(t)=\frac{X_{0,p,j}(\alpha(t))}{\overline{X_{0,p}(t)}\bar{t}^{\kappa_p}},\quad t\in L
\end{equation}
By the ansatz we have
\begin{equation}
\frac{U_{p,j}(\alpha(t))}{X_{0,p,j}(\alpha(t))}= 
\frac{\overline{U_{p}(t)}}{ \bar{t}^{\kappa_p}\overline{X_{0,p,j(t)}}}+G_{p,j}(t)
\end{equation}
\begin{equation}
G_{p,j}(t)=\frac{g_{p,j}(t)}{X_{0,p,j}(t)}
\end{equation}
In the case of positive index $\kappa_p>0$ for all $p$ we can set
\begin{equation}
\frac{U_{p}(z)}{z^{\kappa_p} X_{0,p}(z)}= 
\sum_{i=1}^{\kappa_p} \frac{e_{i,p}}{z^i}+\psi_p(z)
\end{equation}
where the
$\psi_p(z)$ are indefinite analytic functions on $S^+$ and $C_{i,p}$ are complex constants, $p=0,\ldots,q-1.$
Damjanovi\'c \cite{damjanovic} showed that there exists $\xi$ such that when we 
set
$B_{2k-1,p}=\re C_{k,p'},$ $B_{2k,p}=\re C_{k,p'},$
\begin{equation}
g_{2k-1}(t)=\frac{1}{\bar{t}^k},\quad g_{2k}(t)=\frac{-i}{\bar{t}^k},\quad t\in L 
\end{equation}
\begin{equation}
M_{2k-1}(z)=\frac{1}{z^k}-\frac{z-a_0}{2\pi i}\int_L \frac{\overline{\phi_{2k-1}(t)}dt}{(t-a_0)(t-z)}
\end{equation}
\begin{equation}
M_{2k}(z)=\frac{i}{z^k}-\frac{z-a_0}{2\pi i}\int_L \frac{\overline{\phi_{2k}(t)}dt}{(t-a_0)(t-z)}
\end{equation}
$k=1,\ldots,\kappa_{p'},$ $a_0\in S^-=\cup_{j=1}^m S_j^-,$ $a_j\in D_j^-,$ 
$j=1,\ldots,m,$
$M_{2\kappa_p +1}(z)=1,$ $M_{2\kappa_p +2}(z)=i,$ $z\in S^+,$
\begin{equation}
N_{k,j}(z):=\frac{z-a_j}{2\pi i}\int_L \frac{\phi_{k}(\alpha_{-1}^j(t))dt}{(t-a_j)(t-z)}
\end{equation}
$k=1,\ldots,2\kappa_p,$ $j=1,\ldots,m,$ and
$N_{2\kappa_p +1,j}(z)=1,$ $N_{2\kappa_p +1,j}(z)=-i,$ $z\in D_j^+$, where
$\phi_i(t)$ are
the solutions to the Fredholm integral equations $F_{\phi_i}(t)=g_i,$
$i=1,\ldots,2\kappa_{p'},$
\begin{equation}
F_{\phi_i}(t)=\phi(t)+\frac{1}{2\pi i}\int_L \left(\frac{K(\tau,t)\alpha'(\tau)}{\alpha(\tau)-\alpha(t)}
-\frac{(t-a_0)\overline{\xi}}{(\tau -a_0)(\tau-\bar{t})}\phi(\tau)
\right)d\tau
\end{equation}
\begin{equation}
K(\tau,t)=\left\{
\begin{array}{ll}
\frac{\alpha(t)-a_j}{\alpha(\tau) -a_j} & ,\tau,t\in L_j, j=1,\ldots,m\\
0 & ,\tau\in L_j,t\in L_k,j\neq k,j,k=1,\ldots m
\end{array}
\right.
\end{equation}
then we obtain that $(\star)$ is solved by
\begin{equation}\label{damj4a}
w(z,\bar{z})=\sum_{\nu}^{q-1} \bar{z}^\nu U_\nu(z),\quad z\in S^+
\end{equation}
\begin{equation}\label{damj4b}
w_j(z,\bar{z})=\sum_{\nu}^{q-1} \bar{z}^\nu U_{\nu,j}(z),\quad z\in D_j^+
\end{equation}
\begin{equation}
U_\nu:=z^{\kappa_\nu} X_{0,\nu}(z) \left(
\sum_{i=1}^{2\kappa_\nu +2}B_{i,\nu} M_i (z)^j -\frac{z-a_0}{2\pi i}
\int_L \frac{\overline{\lambda_\nu(t)}dt}{(t-a_0)(t-z)}
\right)
\end{equation}
\begin{equation}
U_{\nu,j}:=X_{0,\nu,j}(z)\left(
\sum_{i=1}^{2\kappa_\nu +2}B_{i,\nu} N_{i,j} (z) -\frac{z-a_j}{2\pi i}
\int_L \frac{\lambda_\nu(t) dt}{(t-a_j)(t-z)}
\right)
\end{equation}
(where the $\lambda_\nu(t)$ are the solutions of the Fredholm integral equations
$(F\lambda_\nu)(t)=G_{\nu,j}(t),$ $t\in L_j,$ $\nu=0,\ldots,q-1,$ $j=1,\ldots,m$ and
$B_{i,\nu}$ complex constants).
If however all $\kappa_p<0,$
$p=0,\ldots,q-1,$ then the problem is uniquely solvable with solution given by
\begin{equation}\label{damj5a}
w(z,\bar{z})=\sum_{\nu=0}^{q-1} \bar{z}^\nu U_\nu(z), z\in S^+
\end{equation}
\begin{equation}\label{damj5b}
w_j(z,\bar{z})=\sum_{\nu=0}^{q-1} \bar{z}^\nu U_{\nu,j}(z), z\in D_j^+
\end{equation}
\begin{equation}
U_\nu:=z^{\kappa_\nu} X_{0,\nu}(z) \left(
C_\nu -\frac{z-a_0}{2\pi i}\int_L \frac{\overline{\lambda_\nu(t)}dt}{(t-a_0)(t-z)}\right),\quad z\in S^+
\end{equation}
\begin{equation}
U_{\nu,j}:=z^{\kappa_\nu} X_{0,\nu,j}(z) \left(
\overline{C}_\nu +\frac{z-a_j}{2\pi i}\int_L \frac{\lambda_\nu(t)dt}{(t-a_j)(t-z)}\right),\quad z\in D_j^+
\end{equation}
(where the $\lambda_\nu(z)$ are the solutions of the integral equation
$(F_{\lambda_\nu})(t)=G_{\nu,j}(t)$, $t\in L_j,$ and $C_\nu$ are arbitrary complex constants,
$\nu=0,\ldots,q-1,$ $j=1,\ldots,m$).
For the solvability conditions we expand the function
\begin{equation}
C_\nu =\frac{z-a_0}{2\pi i}\int_L \frac{\overline{\lambda_\nu(t)}dt}{(t-a_0)(t-z)}
\end{equation}
expanding it in a Taylor series about $z=0$, $U_p(z)$ is analytic in $S^+$ if all the coefficients
at $z^0,\ldots,z^{-\kappa_p -1}$ in the expansion vanish.
Choosing
$C_p=\frac{1}{2\pi i}\int_L \frac{\overline{\lambda_\nu(t)}dt}{(t-a_0)(t-z)}$ and separating the real and imaginary
parts in the conditions
\begin{equation}
\int_L \frac{\overline{\lambda_\nu(t)}dt}{(t-a_0)(t-z)}=0,\quad j=1,\ldots,-\kappa_p -1, p=1,\ldots,q-1
\end{equation}
we obtain $-2\sum_{p=0}^{q-1}(\kappa_p +1)$ real conditions for solvability (in the case
$\kappa_p<0,$
$p=0,\ldots,q-1$)
\begin{equation}\label{damj6}
\re \int_L \frac{\overline{\lambda_p(t)}}{t^{j+1}}dt=0,
\im \int_L \frac{\overline{\lambda_p(t)}}{t^{j+1}}dt=0
\end{equation}
$j=1,\ldots,q-1.$
This yields the following theorem;
however as this result is not central to the scope of the book, we have assumed known the 
general solution of certain Fredholm integral equations and the solution certain standard
boundary value problems depending on the indices $\kappa_\nu$ (see
Gakhov \cite{gakhov} and Damjanovi\'c \cite{damjanovic}), so the proof is not complete however we believe
that the method of proof is clear.
\begin{theorem}
Assume for all $p\in\{0,1,\ldots,q-1\}$ the indices
$\kappa_p=\frac{1}{2}\mbox{arg} A_p(t)|_L \geq 0.$
Then $(\star)$ is solvable and its solution is given by
Eqn.(\ref{damj4a}) and Eqn.(\ref{damj4b})
Thus the solution contains $2\sum_{p=0}^{q-1} (k_p +1)$ arbitrary real constants. 
\\
\\
If however all $\kappa_p<0,$
$p=0,\ldots,q-1,$ then the problem is uniquely solvable with solution given by
Eqn.(\ref{damj5a}) and Eqn.(\ref{damj5b})
if and only 
the $-2\sum_{p=0}^{q-1} (\kappa_p +1)$ real conditions of solvability of Eqn.(\ref{damj6}) are satisfied.
\end{theorem}
We mention that Gabrinovi\u{c} \cite{gabrinovicbvp} solved a boundary value problem of Carleman type for certain polyanalytic functions
on multiply-connected domains.

\subsection{Several complex variables}
We mention briefly that some of the techniques describes can be adapted to several complex variables, see e.g.\ Begehr \& Gilbert \cite{begehrgilbert}. 
By iteration in each variable of the derivation of formula for functions $f\in \mathscr{O}(D)\cap C^0(\overline{D})$
\begin{equation}
f(z)=\frac{1}{2\pi i}\int_{\abs{\zeta}=1}\re f(\zeta)\frac{\zeta+z}{\zeta-z}\frac{d\zeta}{\zeta}+\frac{1}{2\pi}\int_{\abs{\zeta}=1}\im f(\zeta)\frac{d\zeta}{\zeta},\quad \abs{z}<1
\end{equation}
used to solve the Schwarz problem for $q=1$, 
yields for the unit polydisc $D^n\subset \Cn,$ the representation 
\begin{equation}
f(z)=\frac{1}{(2\pi i)^n}\int_{\partial D^n}\re f(\zeta)\left(2\frac{\zeta}{\zeta-z}-1\right)\frac{d\zeta}{\zeta}+
\im f(0),\quad z\in D^n
\end{equation}
where 
\begin{equation}
\frac{\zeta}{\zeta-z}:=\Pi_{j=1}^n \frac{\zeta_j}{\zeta_j-z_j},\quad 
\frac{d\zeta}{\zeta}:=\Pi_{j=1}^n \frac{d\zeta_j}{\zeta_j}
\end{equation}
Recalling how the Pompieu operator $Tf:=-\frac{1}{\pi}\int_D f(\zeta)\frac{d\mu(\zeta)}{\zeta-z},z\in \C$, was used to find particular solutions to the inhomogeneous Cauchy-Riemann equation, one can find particular solutions for the overdetermined inhomogeneous system
\begin{equation}
\partial_{\bar{z}_j} f=f_j,\quad \partial_{\bar{z}_k} f_j=\partial_{\bar{z}_j} f_k,\quad 1\leq j,k\leq n
\end{equation}
in $\Cn$ by iterating the Pompieu operators $T_j$ with respect to the variable $z_j$. Namely when the $f_j$ have mixed partial derivatives with respect to the variables $z_k,$ $k=1,\ldots, n$, a particular solution 
is given by 
\begin{equation}
f(z)=\sum_{\nu=0}^{n}(-1)^{\nu+1}\sum_{1\leq j_1<\cdots<j_\nu \leq n} T_{j_\nu}T_{j_{\nu-1}}\cdots T_{j_1} \partial_{\bar{z}_{j_2}}\partial_{\bar{z}_{j_3}}\cdots \partial_{\bar{z}_{j_\nu}} f_{j_1}
\end{equation}
where the second sum is taken over multi-indices $(j_1,\ldots,j_\nu)$
satisfying $1\leq j_1<\cdots<j_\nu \leq n$.
Using this technique one can prove the following.
\begin{theorem}[Schwarz problem on $D^n$]
Let $f_k,$ $k=1,\ldots,n$, have mixed partial derivatives of first order with respect to the variables $z_l,$
$l=1,\ldots,n$ and their complex conjugates $\bar{z}_l,$ $l\neq k$, in $L^1(\overline{D}^n)$ and satisy the compatibility conditions
\begin{equation}
\partial_{\bar{z}_l} f_k=\partial_{\bar{z}_k} f_l,\quad 1\leq k,l\leq n,\quad l\neq k,\mbox{ in }D^n
\end{equation}
Let $\gamma$ be real-valued and continuous on $\partial D^n$ and
\begin{multline}
\re\{\sum_{k=1}^{n-1}\frac{1}{(2\pi i)^n}\int_{\partial D^n} \gamma(\zeta)\left(\Pi_{\nu=1}^k\frac{\zeta_\nu}{\zeta_\nu-z\nu}-1\right)\times\\
\frac{\bar{z}_{k+1}}{\bar{\zeta}_{k+1}-\bar{z}_{k+1}}
\Pi_{\nu=k+2}^n \left(  
\frac{\zeta_{\nu}}{\zeta_{\nu}-z_{\nu}}-
\frac{\bar{\zeta}_{\nu}}{\bar{\zeta}_{\nu}-\bar{z}_{\nu}}
\right)-\\
\sum_{\nu=2}^{n}\sum_{\lambda=1}^{\nu-1}
\sum_{\stackrel{1\leq k_1<\cdots <k_\lambda\leq n}{1\leq 1\leq k_{\lambda+1}<\cdots <k_\nu\leq n}}\frac{(-1)^\nu}{\nu}\int_{D_{k_1}}\cdots \int_{D_{k_\nu}} \partial_{\bar{\zeta}_{k_2}}\cdots \partial_{\bar{\zeta}_{k_\lambda}}\partial_{\zeta_{k_{\lambda+1}}}\partial_{\zeta_{k_\lambda}} f_{k_1}\times\\
 \frac{\bar{z}_{k_1}}{1-\bar{z}_{k_1}\zeta_{k_1}}\cdots
\frac{\bar{z}_{k_\lambda}}{1-\bar{z}_{k_\lambda}\zeta_{k_\lambda}}
\frac{z_{k_{\lambda+1}}}{1-z_{k_{\lambda+1}}\bar{\zeta}_{k_{\lambda+1}}}\cdots\\
\frac{z_{k_\lambda}}{1-z_{k_\lambda}\bar{\zeta}_{k_\lambda}}d\mu(\zeta_{k_1},\ldots,\zeta_{k_{\lambda+1}})\}=0
\end{multline}
Then the Schwarz problem
\begin{equation}
\re f=\gamma \mbox{ on }\partial D^n
\end{equation}
for the inhomogeneous Cauchy-Riemann system
\begin{equation}
\partial_{\bar{z}_k} f=f_k,\quad 1\leq k\leq n, \mbox{ on }\partial D^n
\end{equation}
is solvable. The solution is
\begin{multline}
f(z)=\frac{1}{(2\pi i)^n}\int_{\partial D^n} \gamma(\zeta)\left(2\Pi_{\nu=1}^n\frac{\zeta_\nu}{\zeta_\nu-z\nu}-1\right)\Pi_{\nu=1}^n\frac{d\zeta_\nu}{\zeta_\nu} +\\
\sum_{\nu=1}^{n}\sum_{1\leq k_1<\cdots <k_\lambda\leq n}
\frac{(-1)^\nu}{\pi^\nu}\times\\
\int_{D_{k_1}}\cdots \int_{D_{k_\nu}}\{ \overline{
\partial_{\bar{\zeta}_{k_2}}\cdots \partial_{\bar{\zeta}_{k_\nu}}f_{k_1}}
\frac{z_{k_1}}{1-z_k\bar{\zeta}_{k_1}}\cdots-
\frac{z_{k_\nu}}{1-z_\nu\bar{\zeta}_{k_\nu}}-\\
(-1)^\nu\partial_{\bar{\zeta}_{k_2}}\cdots \partial_{\bar{\zeta}_{k_\nu}}f_{k_1}
\frac{z_{k_1}}{1-z_k\bar{\zeta}_{k_1}}\cdots
\frac{z_{k_\nu}}{1-z_k\bar{\zeta}_{k_\nu}}
\}
d\mu(\zeta_{k_1},\ldots,\zeta_{k_{\nu}})+
\\
+\sum_{\nu=2}^{n}\sum_{\lambda=1}^{\nu-1}\sum_{\stackrel{1\leq k_1<\cdots <k_\lambda\leq n}{1\leq 1\leq k_{\lambda+1}<\cdots <k_\nu\leq n}}\frac{(-1)^\nu}{\pi^\nu}\times\\
\int_{D_{k_1}}\cdots \int_{D_{k_\nu}}
\overline{
\partial_{\bar{\zeta}_{k_2}}\cdots \partial_{\bar{\zeta}_{k_\lambda}}\partial_{\zeta_{k_{\lambda+1}}}\cdots\partial_{\zeta_{k_\nu}} f_{k_1}}\frac{z_{k_1}}{1-z_k\bar{\zeta}_{k_1}}\cdots\\
\frac{z_{k_\lambda}}{1-z_\lambda\bar{\zeta}_{k_\lambda}}
\frac{1}{\zeta_{k_{\lambda+1}}-z_{k_{\lambda+1}}}\cdots
\frac{1}{\zeta_{k_{\nu}}-z_{k_{\nu}}}
d\mu(\zeta_{k_1},\ldots,\zeta_{k_{\nu}}) +ic
\end{multline}
for an arbitrary real constant $c.$
\end{theorem}
For a proof see Begehr \& Gilbert \cite{begehrgilbert}.

\subsection{An example on using growth conditions}
Let $H^+:=\{ z\in \C\colon \im z>0\},$ $H^-:=\{ z\in \C\colon \im z<0\},$ $H_0:=\{\im z=0\}.$

Obviously, when considering higher orders of the Pompieu operator for the half plane
\begin{equation}
T_{0,q} f(z):=\int_{\im \zeta >0} f(\zeta)\frac{(\bar{z}-\bar{\zeta})^{q-1}}{z-\zeta} d\mu(\zeta)
\end{equation}
the behaviour at infinity must be controlled. 

\begin{definition}
For a domain $\Omega\subset\Rn$ denote by $W^{k,p}(\Omega)$ the {\em Sobolev space}\index{Sobolev space} which consists of 
all $L^p(\Omega)$ functions such that for each multi-index $\abs{\alpha}\leq k,$ the weak derivative $D^\alpha$
exists and belongs to $L^p(\Omega),$ and the space is equipped with the Sobolev norm
\begin{equation}
\norm{f}_{W^{k,p}}:=\left(\sum_{\abs{\alpha}\leq k}\int_\Omega \abs{D^\alpha f}^p dx\right)^{\frac{1}{p}},\quad 1\leq p<\infty
\end{equation}
and
\begin{equation}
\norm{f}_{W^{k,\infty}}:=\sum_{\abs{\alpha}\leq k}\mbox{ess-sup}_{x\in \Omega}\abs{D^\alpha f}
\end{equation}
The space of locally $W^{k,p}(\Omega)$ functions is denoted $W^{k,p}_{\mbox{loc}}(\Omega).$
\end{definition}

If $\abs{f(z)}=O(\abs{z}^{-q-\delta}),\mbox{ as }z\to \infty$ then it is known (see Chaudhary \& Kumar \cite{chaudury} and the references
therein) that $T_{0,q} f(z)$ exists as a Lebesgue integral, is continuous in the upper half plane and satisfies
in the weak sense
\begin{equation}\label{cchow}
\partial_{\bar{z}} T_{0,q} f=T_{0,q-1} f
\end{equation}
For example this is the case for all functions in $W^{q,1}(\{\im z>0\},\C)$ such that
\begin{equation}
\abs{z}^{q-2}\partial^q_{\bar{z}} f(z)\in L^1(\{\im z>0\},\C) 
\end{equation}
and
\begin{equation}
\lim_{R\to \infty} R^j\left(\max_{\abs{z}=R} \abs{\partial^j_{\bar{z}} f(z)}\right)=0
\end{equation}

They consider among other things the following Dirichlet problem. 
 For $g\in L^p(H^+,\C),$ $\gamma_j(t)\in L^p(\R,\C),$ $p>2,$
 find $f\in W^{q,1}(H^+,\C)$ such that
 \begin{equation}
 \partial_{\bar{z}} f=g \mbox{ on }H^+,\quad \partial_{\bar{z}}^j f=\gamma_j\mbox{ on }H_0,\quad j=0,\ldots,q-1
 \end{equation}
 with $\bar{z}^{q-2}\partial_{\bar{z}}^q f\in L^1(H^+,\C)$
 and $\lim_{R\to \infty} R^j\left(\max_{\abs{z}=R,0\leq \im z} \abs{\partial^j_{\bar{z}} f(z)}\right)=0$,
 and where $t^j\gamma_j(t)\in L^p(\R,\C)\cap C(\R,\C),$ $j=0,\ldots,q-1$. 
 They give a necessary and sufficient condition, namely  
 \begin{multline}\label{chowdeq}
T_{0,q-\nu}g(\bar{z})=\frac{(-1)^{q-\nu}}{(q-\nu-1)!\pi}\int_{H^+}\partial_{\bar{\zeta}}^q f(\zeta)(\bar{\zeta}-z)^{q-\nu-1}\frac{d\mu(\zeta)}{\zeta-\bar{z}}=\\
-\sum_{j=\nu}^{q-1}\frac{(-1)^{j-\nu}}{(q-\nu)!2\pi i}\int_{-\infty}^\infty \partial_{\bar{\zeta}}^\nu f(t)(t-z)^{j-\nu}\frac{dt}{t-\bar{z}}
\end{multline}
$0\leq \nu\leq q-1,$
for solvability, that makes use of the growth conditions.
 Their theorem also involves the fact that when the problem is solvable the solution is, similar to the bounded case, given by
\begin{multline}
 f(z)=\\
 \frac{(-1)^q}{(q-1)!\pi}\int_{H^+} g(\zeta)(\bar{\zeta}-\bar{z})^{q-1}\frac{d\mu(\zeta)}{\zeta -z}+
\sum_{j=0}^{q-1}\frac{1}{2\pi i}\frac{(-1)^j}{j!}\int_{-\infty}^\infty \gamma_j(t)(t-\bar{z})^j\frac{dt}{(t-z)}=\\
T_{0,q} g(z)+\sum_{j=0}^{q-1}\frac{1}{2\pi i}\frac{(-1)^j}{j!}\int_{-\infty}^\infty \gamma_j(t)(t-\bar{z})^j\frac{dt}{(t-z)} 
 \end{multline}
 But we are here focused on the growth conditions.
 First note that for $q=1$, Eqn.(\ref{chowdeq}) reduces to
 \begin{equation}\label{gg00}
 \frac{1}{\pi}\int_{H^+} f(\zeta)\frac{d\mu(\zeta)}{\zeta-\bar{z}}=
 \frac{1}{2\pi i}\int_{-\infty}^\infty  f(t)\frac{dt}{t-\bar{z}}
  \end{equation}
which is an unbounded version of the complex form of the Gauss theorem, see e.g.\ Begehr \cite{begehrboletinI}, (6).  
Let the contour $C$ consist of the boundary of the intersection with $H^+$, of the circle, centered at $0$, with radius $R>0$. Denote by $C^+$ the bounded domain enclosed by $C.$
When $q=1$, $f$ is holomorphic thus 
\begin{equation}
 \frac{1}{\pi}\int_{C^+} f(\zeta)\frac{d\mu(\zeta)}{\zeta-\bar{z}}=
 \frac{1}{2\pi i}\int_{\theta=0}^\pi  f(t)\frac{dt}{t-\bar{z}},\quad t=R\exp(i\theta)
\end{equation}
The conditions $\bar{z}^{-1}\partial_{\bar{z}} f\in L^1(H^+,\C)$
 and $\lim_{R\to \infty} R\left(\max_{\abs{z}=R} \abs{f(z)}\right)=0$ 
 mean in particular that the integrand on the left hand side satisfies that $\abs{ \frac{f(\zeta)}{\zeta-\bar{z}}}$ on $C\cap H^+$ goes to zero as $R\to \infty$
 thus what remains will be the integral on the real axis, hence Eqn.(\ref{gg00}).
By induction in $q$ and using the inequality
\begin{equation}
\abs{\frac{1}{2\pi i}\int_{\abs{\zeta}=R,0<\im \zeta}\partial_{\bar{\zeta}}^j f(\zeta)
\frac{(\bar{\zeta}-z)^{j-\nu}d\zeta}{\zeta -\bar{z}} } \leq (R+\abs{z})^{j-\nu-1} 
\left(\max_{\abs{z}=R,0\leq \im z} \abs{\partial^j_{\bar{z}} f(z)}\right)R
\end{equation}
where the right hand side tends to $0$ as $R\to \infty,$ they obtain Eqn.(\ref{chowdeq}). Further, given a solution $f$ to the problem and using the fact that under the circumstances
$T_{0,q} f(z)$ exists as a Lebesgue integral, is continuous in the upper half plane and satisfies
in the weak sense Eqn.(\ref{cchow}) 
they obtain
\begin{multline}\label{chowchow}
\partial_{\bar{z}}^j f(z)=\frac{1}{\pi}\int_{-\infty}^\infty \gamma_j(t)\frac{ydt}{\abs{t-z}^2}+\\
\sum_{\nu=j+1}^{q-1}\frac{(-1)^{\nu-j}}{(\nu-j)2\pi i}\int_{-\infty}^\infty \gamma_\nu(t)\left(\frac{(t-\bar{z})^{\nu-j+1}-(t-z)^{\nu-j+1}}{\abs{t-z}^2}\right)dt+\\
T_{0,q-j}(g(z)-g(\bar{z})
\end{multline}
Given the conditions on $\gamma_j$ it is possible to verify 
\begin{equation}
\lim_{\overline{H^+}\ni z\to t_0\in H_0}\frac{1}{\pi}\int_{-\infty}^\infty \gamma_j(t)\frac{ydt}{\abs{t-z}^2} =\gamma_j(t_0)
\end{equation}
and that this is the only terms that remains in Eqn.(\ref{chowchow}) as $z\to t_0\in H_0,$
i.e.\
\begin{equation}
\lim_{z\to t_0}\partial_{\bar{z}}^j f(z)=\gamma_j(t_0), t_0\in H_0
\end{equation}

\subsection{An example using the technique of reduction to systems of Riemann boundary value problems}
Let $H^+:=\{ z\in \C\colon \im z>0\},$ $H^-:=\{ z\in \C\colon \im z<0\},$ $H_0:=\{\im z=0\}.$ 
\begin{definition}
Given a simple closed contour $\gamma$ enclosing the domain $D^+$ and setting $D^-=\hat{\C}\setminus (D^+\cup \gamma)$, 
and two H\"older continuous functions $h,g$ on $\gamma$, $h\not\equiv 0$,
the so-called {\em Riemann boundary value problem associated to $\gamma$}\index{Riemann boundary value problem} is that of finding two functions
$f^\pm\in \mathscr{O}(D^\pm),$ 
satisfying on the boundary
$f^+=h f^-$ (this is the homogeneous case of the problem). Replacing the boundary condition by 
$f^+=h f^- +g$ is the inhomogeneous version of the problem.
The {\em index}\index{Index of the Riemann problem}
of the problem 
denoted by $\kappa$, is defined as the integer $\frac{1}{2\pi}\int_{\gamma} d\arg g(x)$ (i.e.\ the change in argument
of $g$ when traversing $\gamma$ in the positive direction, passing every point once) which is known to be the sum of the number of zeros of $f^\pm$
in $D^\pm$. 
\end{definition}
It is known that for the homogeneous Riemann problem to be solvable it is necessary that $\kappa\geq 0$. If the index is positive the problem has
$\kappa +1$ linearly independent solutions $f^\pm(z)=\exp(\tilde{\Gamma}^\pm(z))z^{-\kappa} P_\kappa(x)$ where 
$\tilde{\Gamma}=\frac{1}{2\pi i}\int_\gamma \frac{\ln t^{-\kappa}h(t)}{t-z}dt$
$P_\kappa$ is a complex polynomial of degree no higher than $\kappa$ (see e.g.\
Gakhov \cite{gakhov}, sec. 14.3, and Section \ref{bvpintro}).
\\
\\
Suppose we are given an integrable function $g$ satisfying the H\"older condition on the real axis $H_0$
such that  the following (Cauchy principal value integral) is well-defined
\begin{equation}
Sg(z):=\left\{
\begin{array}{ll}
\frac{1}{2\pi i}\int_{-\infty}^\infty \frac{g(t)}{t-z}dt,& \im z\neq 0\\
\frac{1}{\pi i}\int_{-\infty}^\infty \frac{g(t)}{t-z}dt,& z\in H_0
\end{array}
\right.
\end{equation}
with the additional condition
\begin{equation}
\lim_{z\in H_0,\re z=x\to \infty} Sg(x)=0
\end{equation}
By the properties of the Cauchy principal value integral $Sg\in \mathscr{O}(H^+\cup H^-)$.
Define also 
\begin{equation}
S^{\pm}g(z):=\left\{
\begin{array}{ll}
Sg(z), & \im z\neq 0\\
\pm\frac{1}{2}g(z)+\frac{1}{2}(Sg)(z),& z\in H_0
\end{array}
\right.
\end{equation}
The function $Sg$ is in particular a so-called {\em sectionally analytic}\index{Sectionally analytic function} function with jump on the real axis in the sense that
it satisfies the jump condition
\begin{equation}
Sg(x)=S^+ g(x)-S^- g(x) ,\quad x\in H_0
\end{equation}
Here
$S^{\pm}g\in  \mathscr{O}(H^\pm)$ and $S^{\pm}g(\infty)=\frac{1}{2}g(\infty)$,
see e.g.\ Du \& Wang \cite{wang} and Section \ref{plemeljsec} of the appendix. 
\\
\\
\textbf{Problem.}
For a given integrable function $g$ as above on the real axis $H_0$, 
find a $q$-analytic function 
$f\in \mbox{PA}_q(H^+\cup H^-)$ such that it is {\em sectionally polyanalytic}\index{Sectionally polyanalytic} 
with jump on the real axis, meaning that we have two functions $f^\pm$ defined on $H^\pm \cup H_0$, $f^\pm=f$ on $H^\pm,$
and satisfying on the boundary
\begin{equation}
\left\{
\begin{array}{ll}
f^+(x)=C_0(x)f^-(x)+\varphi_0(x),\quad x\in H_0\\
\left(\partial_{\bar{z}}f\right)^+ (x)=C_1(x)\left(\partial_{\bar{z}}f\right)^- (x)+\varphi_1(x)\\
\cdots\\
\left(\left(\partial_{\bar{z}}\right)^{q-1}f\right)^+ (x)=C_1(x)\left(\left(\partial_{\bar{z}}\right)^{q-1}f\right)^- (x)+\varphi_{q-1}(x)\\
\end{array}
\right.
\end{equation}
with
\begin{equation}
\limsup_{H^+\cup H^-\ni z\to \infty} \abs{z^{1-q} \tilde{f}(z)}<\infty
\end{equation}
where 
\begin{equation}
\tilde{f}(z):=
\left\{
\begin{array}{ll}
f(z), &z\in H^+\cup H^-\\
f^{\pm}(z),& z\in H_0
\end{array}
\right.
\end{equation}
Note that the requirements imply
\begin{equation}\label{pluto0}
\limsup_{\overline{H^\pm}\ni z\to \infty} \abs{z^{1-q} f^\pm(z)}<\infty
\end{equation}
To solve this problem Du \& Wang \cite{wang} make the Ansatz $f(z)=\sum_{j=0}^{q-} \bar{z}^j f_j(z)$
where the $f_j(z)$ are sectionally analytic functions.
Supposing $f(z)$ solves the problem, then each of the distributions
$(\partial_{\bar{z}})^j f(z),$ $j=0,\ldots,q-1$ defines a sectionally $(q-j)$-analytic function with respect to the real axis
as soon as $f$ belongs to $W^{q,1}(\Omega)$ solving the additional growth condition
\begin{equation}
\limsup_{\overline{H^\pm}\ni z\to \infty} \abs{z^{1-q} \left((\partial_{\bar{z}})^j f(z)\right)^\pm(z)}<\infty
\end{equation} 
The Ansatz reduces the boundary conditions to a system
\begin{equation}\label{klf}
M(x)F^+(x)=C(x)M(x) F^-(x) +\Phi(x)
\end{equation}
where
\begin{equation}
F^\pm (x):=
\begin{bmatrix}
f_0^\pm\\
\vdots\\
f_{q-1}^\pm
\end{bmatrix}
,\quad 
\Phi (x)=
\begin{bmatrix}
\varphi_0^\pm\\
\vdots\\
\varphi_{q-1}^\pm
\end{bmatrix}
\end{equation}

\begin{equation}
M (x):=
\begin{bmatrix}
1 & x & x^2 \cdots x^{q-1}\\
0 & 1! & 2x \cdots (q-1)x^{q-2}\\
0 & 0 & 2! \cdots (q-1)(q-2)x^{q-3}\\
\vdots & \vdots &\vdots & \ddots & \vdots\\
0 & 0 & 0& \cdots (q-1)!
\end{bmatrix}
\end{equation}
\begin{equation}
C(x):=I_{q\times q}
\begin{bmatrix}
C_0(x)\\
C_1(x)\\
\vdots \\
C_{q-1}(x)
\end{bmatrix}
\end{equation}
(here $I_{q\times q}$ denotes the identity matrix).

If $F$ is a solution to Eqn.(\ref{klf}) then for $MF=:G=(g_0,\ldots,q_{q-1})$ we have
\begin{equation}\label{klf}
G^+(x)=C(x)G^-(x) +\Phi(x)
\end{equation}
\\
\\
As a generalization to the index to the case of the real axis define
\begin{equation}
\kappa_j:=\frac{1}{2\pi}\int_{H_0}d\arg C_j(x),\quad j=0,\ldots, q-1
\end{equation}
Set for $j=0,\ldots, q-1$
\begin{equation}
Y_j:=\left\{
\begin{array}{ll}
\exp(\Gamma_j(z)), & z\in H^+\\
\left(\frac{z+i}{z-i}\right)^{\kappa_j}\exp(\Gamma_j(z)), & z\in H^-
\end{array}
\right.
\end{equation}

\begin{equation}
X_j(z):=(z+i)^{-\kappa_j} Y_j(z)
\end{equation}
where
\begin{equation}
\Gamma_j(z):=\frac{1}{2\pi i}\int_{-\infty}^\infty \frac{\ln G'_j(x)}{x-z}dx,\quad G'_j(x):=\left(\frac{z+i}{z-i}\right)^{\kappa_j} G_j(x)
\end{equation}
By the properties of the Cauchy type integrals defining $S^\pm$ we have
for complex constants $a_j$ and $X_j^+$ based upon $X_j,$ defined in the same way as $f^+$ based upon $f,$ $j=0,\ldots,q-1$
\begin{equation}
\lim_{H^+\ni z\to\infty} Y_j^+(z)=\lim_{H^+\ni z\to\infty} z^{\kappa_j} X_j^+(z)=a_j\neq 0 
\end{equation}
The solution to the individual Riemann boundary value problems are (see Wang \& Du \cite{wang} and the references
therein) for $j=0,\ldots,q-1$
\begin{equation}
g_j(z)=\frac{Y_j(z)}{2\pi i}\int_{-\infty}^\infty \frac{\varphi_j(x)}{Y_j^+(x)(x-z)} dx +X_j(z)h_j(z)
\end{equation}
where the $h_j$ are entire functions.

Now $M ^{-1}(x)=[b_{lk}(z)]_{q\times q}$ where
\begin{equation}
b_{lk}(z):=
\left\{
\begin{array}{ll}
\frac{(-1)^{k+l}}{(l-k)!(k-1)!}z^{l-k}, & l\geq k\\
0, & l<l
\end{array}
\right.
\end{equation}
Setting $P(z):=(1,\bar{z},\cdots,\bar{z}^{q-1})$ the solution to the system defined by Eqn.(\ref{klf})
is for the case when all $\kappa_j$ are non-negative
\begin{equation}\label{pluto}
f(z)=P(z)M^{-1}(z) G(z)=\sum_{j=0}^{q-1} \frac{(\bar{z}-z)^j}{j!}g_j(z)
\end{equation}
Du \& Wang \cite{wang} further prove that
the functions $h_j$ can be replaced by polynomials
of order no more than $\kappa_j+q-j-1,$
and they give further remarks also for other signatures $(\kappa_0,\ldots,k_{q-1})$.
Chaudhary \& Kumar \cite{chaudury} prove via other techniques similar results for also the Schwarz problem and a mixed boundary value problem for polyanalytic
functions on the half plane.
   \subsection{The reduced case}
   We shall briefly, for the sake of exposition, sketch (skipping over many details and proofs) the treatment of one 
   type of boundary value problem 
  for {\em reduced} polyanalytic functions,
   see Gakhov \cite{gakhov}, parag.\ 32, and the references therein. 
   Let $\gamma$ be a Hölder continuous simple closed contour in the complex plane parametrized
   by arc length  
   according to
   $\gamma(s)=x(s)+iy(s)$, $s\in \R,$ which has derivatives up to order $(2q+2)$. 
   Denote by $D^+$ the bounded domain enclosed by $\gamma.$
   Assume we are given real
   functions $a_j(s),b_j(s),$ of the 
   arc length parameter, $s,$ for $\gamma$, which satisfy the H\"older condition up to the derivatives of order $2q-j-2$ and
   $a_j^2(s)+b_j^2(s)\neq 0.$
   We consider the problem of finding a reduced $q$-analytic function
   $f(x,y)=u(x,y)+iv(x,y)$ on $D^+$ which is continuous on $\gamma$ satisfying 
   \begin{equation}
   a_j(s)\frac{\partial^{q-1} u}{\partial x^{q-j}\partial y^{j-1}} +
   b_j(s)\frac{\partial^{q-1} v}{\partial x^{q-j}\partial y^{j-1}}=c_j(s)
   \end{equation}
   We consider only the case when $\gamma$ is the unit circle.
   The boundary conditions can be proven to take the form
   \begin{equation}\label{recond}
   \re \left((a_j(s)-ib(s))\frac{\partial^{q-1} f}{\partial x^{q-j}\partial y^{j-1}}\right)=c_j(s),\quad j=1,\ldots,q
   \end{equation}
   Using $\partial_x=\partial_z+\partial_{\bar{z}}$ and $\partial_y=-i\partial_z+i\partial_{\bar{z}}$
   we have
   \begin{multline}
   \frac{\partial^{q-1} }{\partial x^{q-j}\partial y^{j-1}}=i^{j-1}\left( \partial_z+\partial_{\bar{z}}\right)^{q-j}\left( \partial_{\bar{z}}-\partial_z\right)^{q-1}=\\
   i^{j-1} \sum_{l=0}^{q-j}\sum_{m=0}^{j-1} (-1)^m C_{q-j}^l C_{j-1}^m \frac{\partial^{q-1} }{\partial z^{q-l-m-1}\partial \bar{z}^{l+m}}
   \end{multline}
   Because $f$ is assumed to be reduced, there exists holomorphic functions $\phi_j$, $j=0,\ldots,q-1$
   such that
   $f(z)=\sum_{j=0}^{q-1} (z\bar{z})^j\phi_j(z),$ whence 
   \begin{multline}
   \frac{\partial^{q-1} }{\partial x^{q-j}\partial y^{j-1}}f=   i^{j-1} \sum_{l=0}^{q-j}\sum_{m=0}^{j-1} (-1)^m C_{q-j}^l C_{j-1}^m \sum_{\nu=l+m}^{q-1}\frac{\nu!}{(\nu-l-m)!}\bar{z}^{\nu-l-m}\times \\
   (z^\nu \phi_\nu(z))^{q-l-m-1}
   \end{multline}
   Introduce for $j=1,\ldots,q,$ the holomorphic functions
   \begin{multline}
   \Phi_j(z):=\sum_{l=0}^{q-j}\sum_{m=0}^{j-1} (-1)^m C_{q-j}^l C_{j-1}^m \sum_{\nu=l+m}^{q-1}\frac{\nu!}{(\nu-l-m)!}
   z^{q-\nu+l+m-1}\times \\
   (z^\nu \phi_\nu(z))^{q-l-m-1}
   \end{multline}
   Then on $\gamma$ we have
   \begin{equation}
   \frac{\partial^{q-1} f}{\partial x^{q-j}\partial y^{j-1}}=\frac{i^{j-1}}{z^{q-1}}\Phi_j(z),\quad z\in \gamma
   \end{equation}
   and Eqn.(\ref{recond}) takes the form
   \begin{equation}\label{recond1}
   \re \left(i^{j-1} \frac{(a_j(s)-ib(s))}{z^{q-1}}\Phi_j(z)\right)=c_j(s),\quad z\in \gamma,\quad j=1,\ldots,q
   \end{equation}
   This is a system of $q$, so-called {\em Hilbert boundary value problems}\index{Hilbert boundary value problem} for determining the holomorphic functions $\Phi_j,$
   which then in turn can be used to determine the functions $\phi_k,$ $k=0,\ldots,q-1$, using the relations:
   \begin{equation}
   \bar{z}=\frac{1}{z},\quad z\in \gamma=\{\abs{z}=1\}
   \end{equation}
   \begin{equation}
   \frac{\partial^{q-1} }{\partial x^{q-l-m-1}\partial y^{l+m}}f|_{\bar{z}=\frac{1}{z}}=\frac{i^{m+l}}{z^{q-1}}\Phi_{m+l-1}(z)
   \end{equation}
  \begin{multline}
   \frac{\partial^{q-1} }{\partial z^{q-j}\partial \bar{z}^{j-1}}f|_{\bar{z}=
   \frac{1}{z}}=\\
   \left(\frac{1}{2}(\partial_x-i\partial_y)\right)^{q-j}\left(\frac{1}{2}(\partial_x+i\partial_y)\right)^{j-1}(f(z)) |_{\bar{z}=\frac{1}{z}}=\\
   \frac{1}{2^{q-1}}\sum_{p=0}^{q-j}\sum_{m=0}^{j-1} (-1)^{l} i^{l+m}C_{q-j}^l C_{j-1}^m \frac{\partial^{q-1} }{\partial x^{q-l-m-1}\partial y^{l+m}}f|_{\bar{z}=\frac{1}{z}}=\\
   \frac{1}{2^{q-1}z^{q-1}} \sum_{p=0}^{q-j} C_{q-j}^l \sum_{m=0}^{j-1} (-1)^m C_{j-1}^m \Phi_{l+m+1}(z)
   \end{multline}
   see Gakhov \cite{gakhov}, parag.\ 32, and the references therein where it is also described how the same
   boundary value problem can be solved for domains that are mapped onto the unit disc by rational functions. 
 
 \subsection{Concluding comments} 
   We mention that sometimes the term Robin boundary value problem (after
      Gustave Robin (1855-1897) \index{Robin condition} originally concerning the heat equation $\partial f-\Delta f=0$ on $\R_+\times \Omega,$
      $\Omega\subset \Rn,$ a bounded Lipschitz domain and 
      $\partial_\nu f+cf=0$ on $\R_+\times \partial\Omega$, i.e.\ a mixed type condition)
      appears in the literature regarding polyanalytic functions, and for example in
      case of functions $f\in \mathscr{O}(D)\cap C^0(\overline{D})$ it involves the condition $f+zf=\gamma$ on $\partial D,$ see e.g\ Begehr \& Harutyunyan \cite{begehrharu}. 
      \\
      \\
 We believe that the examples given in this text includes most of the basic concepts that are used to solve more general kinds of boundary value problems. We believe that in general the reduction of higher order problems to lower order systems and the use of modified kernel operators (or rather appropriate integral representations) and modification of kernels to explicitly given boundaries (where particular attention must be made to e.g.\ corners or rather the end points of analytic arcs in the boundary
   and in the unbounded case 
 growth conditions included) 
 are the main implemented methods.  
 For this reason we feel that the handful of results given on the topic of boundary value problems can be considered sufficient as an introduction.

\chapter{$q$-analyticity in hypercomplex analysis}\label{hypercomplexsec}
 \section{The associative division algebra $\mathbb{H}$}
 \begin{definition}
 Denote by $\mathbb{H}$ the algebra of real quaternions with standard multiplication $\times$ and componentwise addition \lq$+$\rq,
 recall that $q\in \mathbb{H}$ if $q=q_0+q_1 i+q_2 j+q_3 k,$ where $1,i,j,k$ are basis elements and
 the multiplication is the one induced by $i^2=j^2=k^2=ijk=-1$ 
 (and satisfies $ij=k=-ji,$ $jk=i=-kj,$ $ki=j=-ik$).
 The identity is given by the real quaternion $1$ and the real quaternions form the
 center of the algebra, i.e.\ they are precisely the elements that commute with all members of the algebra.
 For $q\neq 0,$ we have an inverse $q^{-1}$ given by $q^{-1}=(q_0^2+q_1^2+q_2^2+q_3^2)(q_0+q_1 i+q_2 j+q_3 k),$ i.e.\ $\mathbb{H}$ is a non-commutative but associative 
 division algebra over $\R.$ Note that $q=(q_0 +q_1i)+(q_2+q_3i)j=:z+wj,$ with 'conjugate' given by $\bar{q}:=q_0-q_1 i-q_2 j-q_3 k=(q_0-q_1i)-j(q_2-q_3i)=\bar{z}-j\bar{w}$
 and $q\bar{q}=\bar{q}q=\abs{q}^2=q_0^2+q_1^2+q_2^2+q_3^2=\abs{z}^2+\abs{w}^2.$
 \end{definition}
 Note that the $2$-dimensional subspace in $\mathbb{H}$ generated by $\{1,i\}$ is isomorphic to $\C,$ 
 (in fact 
 $\C=\mbox{span}\{1,i\}=\{x:ix=xi\},$ whereas $\mbox{span}\{j,k\}=\{x:ix=-xi\}$) 
 thus $\mathbb{H}$ can be identified as a $\C$-vector space using left multiplication
 for scalars. Let $A$ be a finite dimensional division algebra over $\R$. The fundamental theorem of algebra 
 states that if
 $D$ is {\em commutative} (i.e.\ a field) then it is isomorphic to either $\R$ or $\C.$
 It turns out that $\mathbb{H}$ renders the natural setting for noncommutative complex analysis in the same way
 that $\C$ does for the commutative case (this was proved by Frobenius in 1877, see Palais \cite{palais} for an elementary proof).
  \begin{theorem}
 Let $A$ be a finite dimensional division algebra over $\R$. 
 Then $A$ is isomorphic to precisely on of $\R, \C$ or $\mathbb{H}$
 (in particular, if $A$ is non-commmutative then it is isomorphic to $\mathbb{H}$).
 \end{theorem} 
 \begin{proof}
 Let $a\in A,$ $a\notin \R$. 
 \begin{lemma}\label{palaislem}
 $\R+a\R$ is a maximal commutative subset of $A$, consisting of
 all elements of $A$ commuting with $a$, and it is a field isomorphic to $\C.$
 \end{lemma}
 \begin{proof}
 Let $B\subseteq A$ be a subset of maximal dimension such that $\R+a\R\subseteq F$
 and such that $F$ is commutative. If $x\in A$ is an element that commutes with all elements of $A$ then $F+\R x$
 is commutative, thus equal to $B$, so $x\in B,$ hence $F$ is a maximal commutative
 subset of $A.$ If $0\neq x\in F$ then $x^{-1}$ commutes with all elements $y\in B$
 then $xy=yx\Rightarrow yx^{-1}=x^{-1}y$ so $x^{-1}\in B,$ which implies
 that $B$ is a field, so by the fundamental theorem of algebra is isomorphic $\C$ (since it has dimension $2$).
 If $x\in A$ commutes with $a$ then $x$ commutes with $\R +a\R=B$ hence $x\in B.$ 
 This proves Lemma \ref{palaislem}.
 \end{proof}
 By Lemma \ref{palaislem} we can pick an element $i\in A$ such that $i^2=-1$ and
 $\R +i\R\simeq \C,$ so that $A$ can be viewed as a $\C$-vector space where the scalar operations are given by left multiplication.
 Since right multiplication by $i$ is a $\C$-linear transformation on the $\C$-vector space $A$
 given by $Tx=xi,$ $T^2=-\mbox{Id}$ the only possible values of $T$ are $\pm i.$ Denote 
 $A^{\pm}:=\{x\in A:xi=\pm ix\}$.  
 Since each $x$ can be written
 $x=1/2(x-ixi)+1/2(x+ixi)$ we have $A=A^+\oplus A^-.$ By Lemma \ref{palaislem} $A^+=\C$.
 Since $A^+\cap A^- =\{0\}$ we have $x,y\in A^-\Rightarrow xy\in A^+$.
 If $A^-=\{0\}$ then we immediately have $A\simeq \C$. Suppose $A^- \neq \{0\}.$
 Let $c\in A^-.$ Right multiplication by $c$ is a $\C$-linear transformation on $A$ with inverse given by right
 multiplication ny $c^{-1}$, i.e.\ the transformation is nonsingular and thus
 $x,y\in A^-\Rightarrow xy\in A^+$ implies that 
 right multiplication by $c$ interchanges $A^+$ and $A^-$ so $\mbox{dim}_\C A^-=\mbox{dim}_\C A^+=1,$ 
 thus $A^+\simeq \C$ and $A=A^+\oplus A^-$ implies that
 the real dimension of $A$ must be $4$ thus its complex dimension $2.$
 Since $\R+c\R$ is a fields we have $c^2\in \R+c\R.$ 
 $A^+\simeq \C$ and $x,y\in A^-\Rightarrow xy\in A^+$ implies that
 $c^2\in \C$ thus $c^2\in \C\cap (\R+c\R)=\R.$ If $c^2>0$ it would have two roots in $\R$ thus 3 roots in $\R+c\R,$  which is impossible since $\R+c\R\simeq \C$.
 Hence we have, $c^2\in \R$ and $c^2<1.$
 Let $j$ be a multiple of $c$ such that $j^2=-1$ and set $k:=ij.$ Since 
 $\mbox{dim}_\C A^-=\mbox{dim}_\C A^-=1,$ $\{j,k\}$ forms a basis for $A^-$ over $\R$ thus  
 $A^+\cap A^- =\{0\}$ implies that $\{1,i,j,k\}$ forms a basis for $A$ over $\R.$ Since
 $j,k\in A^-$, they satisfy $ji=-ij,$ $ki=-ik,$ which together with $i^2=j^2=-1,k=ij,$ proves that $\{1,i,j,k\}$
 satisfy the multiplication table for $\mathbb{H}.$
 This completes the proof.
 \end{proof}
 
 \section{$k$-monogenic functions and basic properties}
 The defining equations for polyanalytic functions have natural counterparts in the more general theory of Clifford algebras, but
 we shall here begin with the specific case of $\R^4$ so that we can  give an account of the original pioneering work of
 Brackx \cite{brackx1}. Once this is done we shall present the more general theory.
 \begin{definition}
 A {\em $k$-cell}\index{$k$-cell}, $c$, in $\Rn$, is the image of a differentiable map $\phi:[0,1]^k\to \R^n$, with a specified orientation. The 
 same cell with opposite orientation is denoted $-c.$ 
 By convention a $0$-cell is a isolated point (oriented), i.e.\ a map of $\{0\}$ to a point $c(0)$.
 A {\em $k$-chain}, $C$, is a formal linear combination of $k$-cells. 
 Clearly, we have for the set of $k$-chains, a free additive Abelian group structure (generated by the $k$-cells) 
 using that for $k$-chains $c,d$ and integers $a,b,$ we have $c-c=\emptyset,$
 $ac+bc=(a+b)c$, $c+d=d+c.$ 
 For the standard cube, $I^k$, by which we mean that $\phi_{I^k}$ is given by the inclusion map, $\phi_{I^k}: [0,1]^k\to \Rn$, define
 $\iota_{i,w}(x_1,\ldots,x_{k-1})=$
 $(x_1,\ldots,x_{i-1},w,x_{i},\ldots,x_{k-1})$ and set $\phi_{i,w}:[0,1]^{k-1}\to \Rn,$ $\phi_{i,w}=\phi\circ \iota_{i,w},$
 denote the image of $\phi_{i,w}$ by $I_{i,w}$ (this is a $(k-1)$-cell called the {\em $(i,w)$-face} of $I^k$).
 The boundary $\partial I^k$ is defined as $\sum_{i=1}^k \sum_{w=0}^1 (-1)^{i+w} I_{i,w}.$
 The $(i,w)$-face of a $k$-cell $c$ is defined as the composition $c_{i,w}:=c\circ I_{i,w}$ and the boundary 
 $\partial c :=\sum_{i=1}^k \sum_{w=0}^1 (-1)^{i+w} c_{i,w}.$ 
 Let $C=\sum_{j=1}^N a_j c_j$ be a $k$-chain, for some $N\in \Z_+,$ $a_j\in \Z,$
 and $c_j$ $k$-chains. The boundary of $C$ is defined as $\partial C:=\sum_{j=1}^N a_j \partial c_j.$
 Since  taking the $i$-face of a $j$-face is the same as taking the $j$-face of an $i$-face we have for any $k$-chain $C$, $\partial(\partial C)=0.$
 Let $\omega$ be a $k$-form. Then the integral of $\omega$ over $C$
 is defined by
 \begin{equation}\label{eqnchainint}
 \int_C \omega=\sum_{j=1}^N a_j \int_{c_j} \omega
 \end{equation}
 Here use that for a basis $k$-form, $\omega$, on $[0,1]^k$ we can write $\omega=fdx_1\wedge \cdots \wedge dx_k,$ for a unique function
 $f$ so that
 $\int_{[0,1]^k} \omega =\int_{[0,1]^k} f=$
 $\int_{[0,1]^k} fdx_1\wedge \cdots \wedge dx_k=\int_{[0,1]^k} fdx_1\cdots dx_k$,
 so for a $k$-cell $c$ we can set
 $\int_c\omega=\int_{[0,1]^k} c^*\omega,$ and hence Eqn.(\ref{eqnchainint}) is well-defined (for a $0$-form we define the integral to be $\omega(c(0))$).
 \end{definition}
 \begin{theorem}[Stokes' theorem for real-valued functions]
 Let $\omega\subset \Rn$ be an open subset, $\omega$ a differential $k-1$-form on $\Omega$ and $c$ a $k$-chain on $\Omega.$ Then
 \begin{equation}
 \int_c d\omega =\int_{\partial c} \omega
 \end{equation}
 \end{theorem}
 \begin{proof}
 It suffices to prove the theorem for the case that $\omega=fdx_1\wedge \cdots \wedge \widehat{dx_i} \wedge \cdots\wedge dx_k$
 where $\widehat{dx_i}$ denotes the removal of the $i$:th component.
 First let $\omega$ be a $(k-1)$-form on $[0,1]^k$. 
 Then $\int_{[0,1]^{k-1}} I_{j,w}^* (fdx_1\wedge \cdots \wedge \widehat{dx_i} \wedge \cdots \wedge dx_k)$
 is equal to $0$ for $i\neq j$ and equal to
 $\int_{[0,1]^{k}} I_{j,w}^* f(x_1,\ldots,w,\ldots,x_k)dx_1 \cdots dx_k$ for $i=j.$
 Hence 
 \begin{multline}
 \int_{\partial [0,1]^k} fdx_1\wedge \cdots \widehat{dx_i} \wedge dx_k=
 \sum_{j=1}^k\sum_{w=0}^1 (-1)^{j+w} \int_{[0,1]^{k}} I_{j,w}^* fdx_1\wedge \cdots \wedge \widehat{dx_i} \wedge \cdots \wedge dx_k\\
 =(-1)^{i+1}\int_{[0,1]^k}\int f(x_1,\ldots,1,\ldots,x_k)dx_1 \cdots dx_k +\\
 (-1)^{i}\int_{[0,1]^k}\int f(x_1,\ldots,0,\ldots,x_k)dx_1 \cdots dx_k
 \end{multline}
 For the left hand side, $\int_c d\omega$ we have, by Fubini's Theorem and the one dimensional Fundamental Theorem of Calculus
 \begin{multline}
 \int_{\partial [0,1]^k} d(fdx_1\wedge \cdots \widehat{dx_i} \wedge dx_k)=
  \int_{[0,1]^{k}} \partial_{x_i} fdx_1\wedge \cdots \widehat{dx_i} \wedge \cdots\wedge dx_k\\
 =(-1)^{i-1} \int_0^1\cdots\left(\int_0^1 \partial_{x_i} f(x_1,\ldots,x_k)dx_i\right) dx_1\cdots \widehat{dx_i} \cdots dx_k=\\
 (-1)^{i-1} \int_0^1\cdots \int_0^1 (f(x_1,\ldots,1,\ldots,x_k)-f(x_1,\ldots,0,\ldots,x_k)) dx_1\wedge \cdots \widehat{dx_i} \wedge \cdots\wedge dx_k=\\
 (-1)^{i-1}\int_{[0,1]^k} \partial_{x_i} f(x_1,\ldots,1,\ldots,x_k)dx_1 \cdots dx_k +\\
 (-1)^{i}\int_{[0,1]^k} \partial_{x_i} f(x_1,\ldots,0,\ldots,x_k)dx_1 \cdots dx_k\\
 \end{multline}
 which proves the result for $[0,1]^k.$ By definition we have for a $k$-chain $c$, $\int_{\partial c} \omega =\int_{\partial [0,1]^k} c^*\omega$
 thus
 \begin{multline}
 \int_c d\omega = \int_{\partial [0,1]^k} c^*(d\omega)=\int_{\partial [0,1]^k} d(c^*\omega)=\int_{\partial c}\omega
 \end{multline}
 Using additivity for the general case of $c=\sum_j a_j c_j$, i.e.\ $\int_c d\omega=\sum_j a_j\int_{c_i} d\omega$ yields the wanted result.
 This completes the proof.
 \end{proof}
 
 \begin{definition}
 Let $M$ be a $4$-dimensional, differentiable, oriented manifold, 
 $M\subset \Omega$ for an open non-empty subset $\Omega\subset\R^4.$ Let $p\in \{1,2,3,4\},$ and let
 $C$ be a $p$-chain on $M.$ Let $\omega_\alpha$, $\alpha=0,1,2,3,$ be real $p$-forms on $M$ i.e.\
 $\omega_\alpha =\sum_h \eta_{\alpha,h} dx^h,$ where $h=(h_1,\ldots,h_p)\in \{0,1,2,3\}^p,$ $0\leq h_1 <\cdots<h_p\leq 3,$
 $dx^h\in \Lambda^p w,$ where $w$ is the $4$-dimensional vector space with basis
 $\{dx_0,dx_1,dx_2,dx_3\}$ and
 $\eta_{\alpha,h}\in C^r(\Omega),$ $\eta_{\alpha,h}:\R^4 \to \R,$ for all $\alpha,h.$
 Let $\{e_0,e_1,e_2,e_3\}$ denote the basis of the algebra $\mathbb{H}$. Then we can decompose each
 {\em quaternion-$p$-form}\index{Quaternion-$p$-form} $\omega$ as 
 $\omega=\sum_{\alpha=0}^3 e_\alpha \omega_\alpha.$ Define
 \begin{equation}
 \int_C \omega =\sum_{\alpha=0}^3 e_\alpha \int_C \omega_\alpha
 \end{equation}
 A function $f:\R^4\to \mathbb{H}$ can be represented by
 \begin{equation}
 f=\sum_{\alpha=0}^3 e_\alpha f_\alpha,\quad x=(x_1,x_2,x_3,x_4)\mapsto \sum_{\alpha=0}^3 e_\alpha f_\alpha(x)
 \end{equation}
 where $f_\alpha$, $\alpha=0,1,2,3,$ are real-valued.
 Define
 \begin{equation}
 D:=\sum_{\beta=0}^3 e_\beta \frac{\partial}{\partial x_\beta}
 \end{equation}
 so that 
 \begin{equation}
 Df=\sum_{\alpha,\beta =0}^3 e_\beta e_\alpha \frac{\partial f_\alpha}{\partial x_\beta},\quad (fD)=\sum_{\alpha,\beta=0}^3  e_\alpha e_\beta \frac{\partial f_\alpha}{\partial x_\beta}
 \end{equation}
 A function $f:\R^4\to \mathbb{H}$ is called {\em left (right) $k$-monogenic}\index{$k$-monogenic function} on $\Omega$ if
 (i) $f_\alpha \in C^k(\Omega), \alpha=0,1,2,3$ (this shall be denoted $f\in C^k(\Omega)$) and (ii)
 \begin{equation}
 D^kf=D(D^{k-1} f)=0 \mbox{ in } \Omega \quad ((fD^k=0)\mbox{ in }\Omega)
 \end{equation}
 Define further 
 \begin{equation}\overline{D}:=\sum_{\beta=0}^3 \epsilon_\beta e_\beta \frac{\partial}{\partial x_\beta},\quad \epsilon_0=1,\epsilon_\beta =-1,\beta=1,2,3
 \end{equation}
 \end{definition}
 Not only are the defining equations analogous to those of $q$-analytic functions, but 
 the relation 
 \begin{equation}D\overline{D}=\overline{D}D=\Delta e_0, \qquad
 \Delta:=\sum_{\beta=0}^3 e_\beta \frac{\partial^2}{\partial x^2_\beta}
 \end{equation}
 immediately implies 
 \begin{equation}
 D^k\overline{D}^k f= \overline{D}^k (D^K f)=\Delta^k f
 \end{equation} 
 thus we obtain the
 following counterpart to the fact that $q$-analytic functions have $q$-harmonic real and imaginary parts, thus are complex-valued
 $q$-harmonic functions.
 \begin{proposition}
 If $f\in C^{2k}(\Omega)$ is left (right) $k$-monogenic on $\Omega$ then $f$ is $k$-harmonic in the sense that
 $\Delta^k f=(\overline{D}D)^kf=\sum_{\alpha=0}^3 \frac{\partial^2}{\partial x_\alpha^2}f=0$ on $\Omega.$
 \end{proposition}
 We denote $d\hat{x}_0:=dx_{1}\wedge dx_{2} \wedge dx_{3},$
 $d\hat{x}_1:=dx_{0}\wedge dx_{2} \wedge dx_{3},$ $d\hat{x}_2:=dx_{0}\wedge dx_{1} \wedge dx_{3},$
 $d\hat{x}_3:=dx_{1}\wedge dx_{1} \wedge dx_{2},$
 $d\sigma_x:=\sum_{\alpha =0}^3 (-1)^\alpha e_\alpha d\hat{x}_{\alpha}.$
 \begin{theorem}\label{stokesthmbrackx}
 Let $\Omega \subset\R^4$ be a nonempty open subset, let $M\subset \omega$ be a $4$-dimensional differentiable
 orientable manifold. Let $C$ be a $4$-chain on $M$ and $f,g\in C^1(\Omega).$ Then
 \begin{equation}
 \int_{\partial C} fd\sigma_x g=\int_C \left((fD) g+f(Dg)\right) dx_0\wedge dx_{1}\wedge dx_{2} \wedge dx_{3}
 \end{equation} 
 \end{theorem}
 \begin{proof}
 The left hand side can be written as
 \begin{equation}
 \int_{\partial C} fd\sigma_x g=\int_{\partial C} \sum_{\alpha,\beta,\gamma=0}^3 
 (-1)^\alpha e_\beta e_\alpha e_\gamma f_\beta g_\gamma d\hat{x}_\alpha
 \end{equation}
 By Stokes' theorem for real-valued functions 
 \begin{equation}
 \int_{\partial C} f_\beta g_\gamma d\hat{x}_\alpha =(-1)^\alpha \int_{C} 
 \partial_{x_\alpha} (f_\beta g_\gamma) dx_0\wedge dx_{1}\wedge dx_{2} \wedge dx_{3}
 \end{equation}
 which implies
 \begin{multline}
 \int_{\partial C} fd\sigma_x g = \\
 \sum_{\alpha,\beta,\gamma=0}^3 e_\beta e_\alpha e_\gamma \int_C \left(
 (\partial_{x_\alpha} f_\beta)g_\gamma +f_\beta \partial_{x_\alpha} g_\gamma 
 \right) dx_0\wedge dx_{1}\wedge dx_{2} \wedge dx_{3}\\
 =\int_C ((fD)g +f(Dg))dx_0\wedge dx_{1}\wedge dx_{2} \wedge dx_{3}
 \end{multline}
 This completes the proof.
 \end{proof}
 
 \begin{theorem}[A quaternion Green's theorem]\label{brackxthm22}
 Let $\Omega \subset\R^4$ be a nonmepty open subset, let $M\subset \omega$ be a $4$-dimensional differentiable
 orientable manifold. Let $C$ be a $4$-chain on $M$ and $f,g\in C^k(\Omega).$ Then
 \begin{equation}
 \int_{\partial C} \sum_{j=0}^{k-1} (-1)^j fD^{k-1-j}d\sigma_x D^j g=
  \int_C  (fD) g+(-1)^{k-1} f(D^k g) dx_0\wedge dx_{1}\wedge dx_{2} \wedge dx_{3}
 \end{equation} 
 \end{theorem}
 \begin{proof}
 Substituting in Theorem \ref{stokesthmbrackx} $f,g$ with $fD^{k-1-j},D^jg$ and summing over $j=0,\ldots,k-1$ we obtain
 \begin{multline}
 \int_{\partial C} \sum_{j=0}^{k-1} (-1)^j fD^{k-1-j}d\sigma_x D^j g =\\
 \sum_{j=0}^{k-1} (-1)^j \int_C (fD^{k-1-j})( D^j g) dx_0\wedge dx_{1}\wedge dx_{2} \wedge dx_{3}+\\
 \sum_{j=1}^k (-1)^{j-1} \int_C (fD^{k-1-j})(D^j g) dx_0\wedge dx_{1}\wedge dx_{2} \wedge dx_{3}=\\
 \int_C (fD^{k}) g dx_0\wedge dx_{1}\wedge dx_{2} \wedge dx_{3}+\\
 \sum_{j=1}^k ((-1)^j+(-1)^{j-1})\int_C (fD^{k-j})(D^j g) dx_0\wedge dx_{1}\wedge dx_{2} \wedge dx_{3}+\\
 (-1)^{k-1}\int_C f(D^kg) dx_0\wedge dx_{1}\wedge dx_{2} \wedge dx_{3}=\\
 \int_C f(D^kg)+(-1)^{k-1}f(D^kg) dx_0\wedge dx_{1}\wedge dx_{2} \wedge dx_{3}
 \end{multline}
 This completes the proof.
 \end{proof}
 \begin{definition}
 For a function $g\in L^1_{\mbox{loc}}(\R^4,\mathbb{H})$ we define the associated
 {\em quaternion distribution} as follows:
 For each real-valued test function $\phi\in C^\infty_c(\R^4)$, set
 \begin{equation}
 g[\phi]:=\int g(x)\phi(x) dx_0\wedge dx_{1}\wedge dx_{2} \wedge dx_{3}
 \end{equation}
 \index{Quaternion distribution}
 Denote by $\rho_x^2:=\abs{x}^2=\sum_{j=0}^3 x_j^2,$
 $x=\sum_{j=0}^3 e_jx_j,$ $\overline{x}=\sum_{j=0}^3 \epsilon_j e_j x_j.$ 
  For a given surface $\Sigma$, denote by 
 $d\hat{u}_j:=(-1)^j n_j dS,$ $j=0,1,2,3,$ where $n_j$ is the $j$:th component of the external surface normal to $\Sigma$
 and $dS$ the 3-form that is the elementary surface element. For the unit sphere the 3-form can be written $d\sigma_u=\sum_{j=0}^3 e_j u_j dS=udS.$
 \end{definition}
 \begin{theorem}
 The functions 
 \begin{equation}
 g_k(x)=\frac{1}{2\pi^2}\frac{\overline{x}}{\rho_x^4}\frac{x_0^{k-1}}{(k-1)!},\quad x\in \R^4\setminus\{0\}, k\in \Z_+
 \end{equation}
 (where $2\pi^2$ should be interpreted as the surface area of the $4$-dimensional unit sphere)
 satisfy
 \begin{equation}
  D^k g_k =g_k D^k =\delta
 \end{equation}
 where $\delta$ is the Dirac distribution,
 i.e.\ $g_k$ is a two-sided fundamental solution for the operator $D^k.$ 
 \end{theorem}
 \begin{proof}
 Note that each component of $g_k$ is analytic on $\R^4\setminus \{0\}$ and for
 $k>1, k\in \Z_+,$ we have
 $D g_k=g_k D=g_{k-1}$ on $\R^4\setminus \{0\}$ which implies that
 \begin{equation}
 D^k g_k =g_kD^k =0\mbox{ on }\R^4\setminus\{0\}
 \end{equation}
 Taking the distribution derivatives we have for each test function $\phi\in C^\infty_c(\R^4)$
 \begin{multline}
 D^k g_k[\phi] = (-1)^k [D^k\phi]g_k=\\
 (-1)^k\int (D^k\phi(x))g_k(x) dx_0\wedge dx_{1}\wedge dx_{2} \wedge dx_{3}
 \\
 =\lim_{\epsilon \to +0}(-1)^k \int_{\{ \abs{x}> \epsilon\}} (\phi D^k)g_k dx_0\wedge dx_{1}\wedge dx_{2} \wedge dx_{3}
 \end{multline}
 By Theorem \ref{brackxthm22} we have
 \begin{equation}
 D^k g_k[\phi]=\\
 \lim_{\epsilon \to +0} \int_{\{\abs{x}=\epsilon\}} \sum_{j=0}^{k-1} (-1)^{k+1-j} (\phi(x)D^{k-1-j})
 d\sigma_x D^jg_k(x)
 \end{equation}
 Setting $x=\epsilon u$ this gives
 \begin{multline}
 D^k g_k[\phi]=\\
 \lim_{\epsilon \to +0} \int_{\{\abs{x}=1\}} \sum_{j=0}^{k-2} (-1)^{k+1-j} (\phi(\epsilon u)D^{k-1-j}) \epsilon^3
 d\sigma_u g_{k-j}(\epsilon u) +\\
 \lim_{\epsilon \to +0} \int_{\{\abs{x}=1\}} \phi(\epsilon u) \epsilon^3 d\sigma_u g_{1}(\epsilon u)
 \end{multline}
 Since $\phi\in C^\infty_c(\R^4)$, each $\abs{\phi D^{k-1-j}}$ is bounded on $\{\abs{x}<1\}$ by say $M_j>0$ for $j=0,\ldots,k-2.$
 This gives
 \begin{equation}
 0\leq \abs{\int_{\{\abs{x}=1\}}\phi(\epsilon u) D^{k-1-j}\epsilon^3 d\sigma_u g_{k-j}(\epsilon u)}\leq 
 \frac{ \epsilon^{k-1-j} M_j}{(k-j-1)!}
 \end{equation}
 Since for sufficiently small $\epsilon$ we have
 $\phi(\epsilon u)=\phi(0)+Q(\epsilon)$ for some $Q$ satisfying
 $\lim_{\epsilon\to +0} Q(\epsilon)=0$, we deduce that
 \begin{equation}
 \int_{\{\abs{x}=1\}}\phi(\epsilon u) \epsilon^3 d\sigma_u g_1(\epsilon u)=\frac{1}{2\pi^2}
 \int_{\{\abs{x}=1\}} (\phi(0)+Q(\epsilon))dS
 \end{equation}
 and
 \begin{equation}
 \lim_{\epsilon \to +0}\int_{\{\abs{x}=1\}}\phi(x)d\sigma_x g_1(x)=\phi(0)
 \end{equation}
 This renders
 \begin{equation}
 D^k g_k[\phi]=(-1)^k \int (D^k\phi(x))g_k(x) dx_0\wedge dx_{1}\wedge dx_{2} \wedge dx_{3}=\phi(0)
 \end{equation}
 This implies $g_kD^k[\phi]=\phi(0).$ This completes the proof.
 \end{proof}
 \begin{theorem}[Quaternion version of the Cauchy integral formula]\label{brackscauchyformula}
 Let $\Omega\subset\R^4$ be an open subset, let $f$ be a left-k-monogenic
 function on $\Omega$ and let $S\subset\Omega$ be a $4$-dimensional, compact, differentiable, oriented
 manifold with boundary. Then for each $x\in \mathring{S}$
 \begin{multline}
 f(x)=\int_{\partial S} \sum_{j=0}^{k-1} (-1)^j g_{j+1}(u-x)d\sigma_u D^jf(u)=\\
 \frac{1}{2\pi^2} \int_{\partial S}\sum_{j=0}^{k-1} (-1)^j \frac{\bar{u}-\bar{x}}{\abs{u-x}^4}\frac{(u_0-x_0)^j}{j!}d\sigma_u D^j f(u)
 \end{multline}
 where $\abs{u-x}$ denotes the Euclidean distance in $\R^4,$ $\mbox{dist}(x,u).$
 \end{theorem}
 \begin{proof}
 For fixed $x\in \mathring{S}$ pick $t\in \R^+$ such that $\{\abs{x}\leq t\}\subset\mathring{S}.$
 A consequence of Theorem \ref{brackxthm22} (which is known to hold true also for the case of manifolds with boundary, see
 e.g.\ Ward \cite{ward})
   is that for a right-$k$-monogenic function $f$ and a left-$k$-monogenic function $g$ on $\Omega,$
 we have for any $4$-chain $C$ on $M\subset\Omega$
 \begin{equation}
 \int_{\partial C} \sum_{j=0}^{k-1} (-1)^j fD^{k-1-j} d\sigma D^j g=0
 \end{equation}
 Applying this for $g_k$ and $f$ with the domain replaced by $S\setminus \{\abs{x}\leq t\}$
 yields
 \begin{multline}\label{brackx41ekv}
 \int_{\partial S} \sum_{j=0}^{k-1} (-1)^j g_{j+1}(u-x)d\sigma_u D^j f(u)=\\
 \frac{1}{2\pi^2} \int_{\{\abs{x} = t\}} \sum_{j=0}^{k-1} (-1)^j 
 \frac{\bar{u}-\bar{x}}{\abs{u-x}^4}\frac{(u_0-x_0)^j}{j!}d\sigma_u D^j f(u)=\\
 \frac{1}{2\pi^2} \sum_{j=0}^{k-1} \frac{(-1)^j}{j!t^4} \int_{\{\abs{x} = t\}} \sum \epsilon_\gamma e_\gamma e_i e_{\beta_1}\cdots e_{\beta_j} e_\alpha (-1)^i \times\\
 \frac{\partial^j f_\alpha}{\partial u_{\beta_1}\cdots \partial u_{\beta_j}}(u_\gamma -x_\gamma)
 (u_0-x_0)^j d\hat{u}_i
 \end{multline}
 By Theorem \ref{stokesthmbrackx} we have
 \begin{multline}
 (-1)^i \int_{\{\abs{x} = t\}} \frac{\partial^j f_\alpha}{\partial u_{\beta_1}\cdots \partial u_{\beta_j}}(u_\gamma -x_\gamma)
 (u_0-x_0)^j d\hat{u}_i =\\
 \int_{\{\abs{x} < t\}} \frac{\partial^{j+1} f_\alpha}{\partial u_i \partial u_{\beta_1}\cdots \partial u_{\beta_j}}(u_\gamma -x_\gamma)
 (u_0-x_0)^j du_0\wedge du_{1}\wedge du_{2} \wedge du_{3}+\\
 \int_{\{\abs{x} < t\}} \frac{\partial^{j} f_\alpha}{\partial u_{\beta_1}\cdots \partial u_{\beta_j}}(u_\gamma -x_\gamma)j
 (u_0-x_0)^{j-1} \delta_{0i} du_0\wedge du_{1}\wedge du_{2} \wedge du_{3} +\\
 \int_{\{\abs{x} < t\}} \frac{\partial^{j} f_\alpha}{\partial u_{\beta_1}\cdots \partial u_{\beta_j}}
  \delta_{i\gamma}(u_0-x_0)_j du_0\wedge du_{1}\wedge du_{2} \wedge du_{3}
 \end{multline}
 For $f\in C^k(\Omega)$ we have the estimates on $S$
 \begin{equation}
 \abs{\frac{\partial^{j+1} f_\alpha(x)}{\partial u_i\partial u_{\beta_1}\cdots \partial u_{\beta_j}}}\leq
 M_{\alpha,i,\beta_1,\ldots,\beta_j},\mbox{ for all }\alpha,i,\beta_1,\ldots,\beta_j
 \end{equation}
 \begin{equation}
 \abs{\frac{\partial^{j} f_\alpha(x)}{\partial u_{\beta_1}\cdots \partial u_{\beta_j}}}\leq
 M_{\alpha,\beta_1,\ldots,\beta_j},\mbox{ for all }\alpha,\beta_1,\ldots,\beta_j
 \end{equation}
 Setting $M_1:=\max M_{\alpha,i,\beta_1,\ldots,\beta_j},$ $M_2:=\max M_{\alpha,\beta_1,\ldots,\beta_j},$
 yields
 \begin{multline}\label{haarek0n2}
 0\leq \\
 \abs{\int_{\{\abs{x} < t\}}\frac{\partial^{j+1} f_\alpha(x)}{\partial u_i\partial u_{\beta_1}\cdots \partial u_{\beta_j}}
 (u_\gamma-x_\gamma)(u_0-x_0)^j
 du_0\wedge du_{1}\wedge du_{2} \wedge du_{3} }\\
 \leq M_1 \int_{\{\abs{x} < t\}} \abs{u_\gamma-x_\gamma}\abs{u_0-x_0}^j
 du_0\wedge du_{1}\wedge du_{2} \wedge du_{3}\\
 \leq M_1t^{j+5}\frac{1}{4}\frac{1}{2\pi^2}
 \end{multline}
 \begin{multline}\label{haarek0n1}
 0\leq \\
 \abs{\int_{\{\abs{x} < t\}}\frac{\partial^{j} f_\alpha(x)}{\partial u_{\beta_1}\cdots \partial u_{\beta_j}}
 (u_\gamma-x_\gamma)j(u_0-x_0)^{j-1}
 du_0\wedge du_{1}\wedge du_{2} \wedge du_{3} }\\
 \leq jM_2 \int_{\{\abs{x} < t\}} \abs{u_\gamma-x_\gamma}\abs{u_0-x_0}^{j-1}
 du_0\wedge du_{1}\wedge du_{2} \wedge du_{3}\\
 \leq M_2 \int_{\{\abs{x} < t\}} \abs{u_\gamma-x_\gamma}\abs{u_0-x_0}^{j-1}
 du_0\wedge du_{1}\wedge du_{2} \wedge du_{3} \leq jM_2 t^{j+4}\frac{2\pi^2}{4}
 \end{multline}
 
 \begin{multline}\label{haarek0n}
 0\leq \\
 \abs{\int_{\{\abs{x} < t\}}\frac{\partial^{j} f_\alpha(x)}{\partial u_{\beta_1}\cdots \partial u_{\beta_j}}
 (u_\gamma-x_\gamma)du_0\wedge du_{1}\wedge du_{2} \wedge du_{3} }\\
 \leq M_2 \int_{\{\abs{x} < t\}} \abs{u_0-x_0}^{j}
 du_0\wedge du_{1}\wedge du_{2} \wedge du_{3}\\
 \leq M_2 t^{j+4}\frac{2\pi^2}{4}
 \end{multline}
 where $\frac{1}{4}2\pi^2$ is to be interpreted as the volume of the $4$-dimensional unit sphere.
 For $t>0$ sufficiently small we have
 $f_\alpha(u)=f_\alpha(x)+Q_\alpha(t),$ for some $Q_\alpha(t)$
 such that $\lim_{t\to +0}Q_\alpha(t)=0,$ $\alpha=0,1,2,3.$
 Thus for $j=0$ in Eqn.(\ref{haarek0n}) we have
 \begin{equation}
 \int_{\{\abs{x} < t\}} f_\alpha du_0\wedge du_{1}\wedge du_{2} \wedge du_{3} =t^4 \frac{2\pi^2}{4} f_\alpha +t^4
 \frac{2\pi^2}{4} Q_\alpha (t)
 \end{equation}
 Letting $t\to +0$ in Eqn.(\ref{brackx41ekv}) and applying Eqn.(\ref{haarek0n}), Eqn.(\ref{haarek0n1}),
 Eqn.(\ref{haarek0n2}) we obtain
 \begin{multline}
 \int_{\partial S} \sum_{j=0}^{k-1} (-1)^j g_{j+1}(u-x)d\sigma_u D^j f(u)=\\
 \frac{4}{2\pi^2} \sum_{\alpha,i=0}^{3} \epsilon_i e_i e_i e_{\alpha} \frac{2\pi^2}{4} f_\alpha(x)=\\
  \sum_{\alpha,i=0}^{3} e_{\alpha} f_\alpha(x)=f(x)
 \end{multline}
 This completes the proof.
 \end{proof}
 
 \begin{theorem}
 Let $\Omega\subset\R^4$ be an open subset, let $f$ be a left-$k$-monogenic
 function on $\Omega$. Note that any closed ball centered at $a\in \Omega$, with radius $t>0,$ such that 
 $\{\abs{x-a}\leq t\}\subset\Omega$ defines  
 a $4$-dimensional, compact, differentiable, oriented
 manifold with boundary. Then
 \begin{equation}
 f(a)=\frac{1}{t^4 2\pi^2}\int_{\{\abs{x-a}<t\}}\sum_{j=0}^{k-1} (-1)^j \frac{(u_0 -a_0)^j}{j!} D^j f(u) 
 du_0\wedge du_{1}\wedge du_{2} \wedge du_{3}
 \end{equation}
 \end{theorem}
 \begin{proof}
 By Theorem \ref{brackscauchyformula}
 \begin{equation}
 f(a)=\frac{1}{2\pi^2}\int_{\{\abs{x-a}=t\}}\sum_{j=0}^{k-1} (-1)^j \frac{\bar{u} -\bar{a}}{\abs{u-a}^4}
 \frac{(u_0 -a_0)^j}{j!} d\sigma_u D^j f(u) 
 \end{equation}
 By Theorem \ref{stokesthmbrackx} we obtain
 \begin{multline}
 t^4 2\pi^2 f(a)= 
 \sum_{j=0}^{k-1} (-1)^j \times\\
  \left( 
 (\bar{u} -\bar{a})\frac{(u_0 -a_0)^j}{j!} D^{j+1} f(u) +4 \frac{(u_0 -a_0)^j}{j!}D^j f(u)
 \right) du_0\wedge du_{1}\wedge du_{2} \wedge du_{3} \\
 +\sum_{j=0}^{k-1} (-1)^j  \int_{\{\abs{x-a}<t\}} \frac{(\bar{u} -\bar{a})(u_0-a_0)^{j-1}}{(j-1)!} (D^j f(u))
 du_0\wedge du_{1}\wedge du_{2} \wedge du_{3}
 \end{multline}
 Hence
 \begin{equation}
 f(a)=\frac{4}{t^42\pi^2}\int_{\{\abs{x-a}<t\}}\sum_{j=0}^{k-1} (-1)^j \frac{(u_0 -a_0)^j}{j!} D^j f(u) 
 du_0\wedge du_{1}\wedge du_{2} \wedge du_{3}
 \end{equation}
 This completes the proof.
 \end{proof}
 
 \begin{definition}
 A {\em homogeneous left-$k$-monogenic polynomial}\index{Left-$k$-monogenic polynomial} of degree $m$, in a neighborhood 
 $\Omega_0$ of the origin, is a function
 \begin{equation}
 Q_m^{(0)}:\R^4\to \mathbb{H},\quad x\mapsto \sum_{j=0}^3 e_j Q_{m,j}^{(0)}(x), j=0,1,2,3
 \end{equation}
 where each $Q_{m,j}^{(0)} :\R^4\to \R$ is a homogeneous polynomial of degree $m$ in $x_0,x_1,x_2,x_3,$
 satisfying
 $D^k Q_m^{(0)} =0$ on $\Omega_0 .$
 If $a=(a_0,a_1,a_2,a_3)\neq 0$ we use the notation
 $Q_m^{(a)}$ for the quaternion polynomial 
 $Q_m^{(0)}(x)=\sum_{j=0}^3 e_j Q_{m,j}^{(0)}(x),$ where for each $j,$
 $Q_{m,j}^{(a)} :\R^4\to \R$ is a homogeneous polynomial of degree $m$ in
 $(x_0-a_0),(x_1-a_1),(x_2-a_2),(x_3-a_3).$ For practical reason we shall
 work always near the origin $a=0$, thus we shall drop the superscript and write simply 
 $Q_m(x):=Q_m^{(0)}(x),$ and we assume $m<k.$
 A homogeneous left-$1$-monogenic polynomial is called {\em regular}. Denote by 
 $P_{\mbox{Fueter}}^{m-l}$ the set of regular homogeneous polynomials of the form
 \begin{equation}
 p_{r_1,\ldots,r_{m-l}}(x) =\frac{1}{(m-l)!} \sum_{\pi(r_1,\ldots,r_{m-l})} z_{r_1}\cdots z_{r_{m-l}}
 \end{equation}
 where $\pi(r_1,\ldots,r_{m-l})$ denotes the set of all permutations of $\{r_1,\ldots,r_{m-l}\}$ 
 and $z_{r},$ $r=1,2,3$ are the so-called hypercomplex variables defined by
 \begin{equation}
 z_r =x_re_0 -x_0 e_r,\quad r=1,2,3
 \end{equation}
 \end{definition}
 Brackx \cite{brackx1} proves the following results on polynomials and quaternion Taylor series expansion.
\begin{theorem}
If $P_n(x)$ is a homogeneous left-$k$-monogenic polynomial of degree $q$ in a neighborhood $U$
of the origin then on $U$
\begin{equation}
P_q(x)=\sum_{j=0}^{k-1} \sum_{(r_1,\ldots,r_{q-j})} \frac{x_0^j}{j!} p_{r_1,\ldots,r_{m-l}} \cdot \frac{\partial^{q-j} D^j P_n}{\partial x_{r_1}\cdots \partial x_{r_{q-j}}}
\end{equation}
where the second sum runs over all combinations of $(1,2,3)$ in sets of $(q-j)$ elements.
Furthermore, any polynomial of the form $\frac{1}{j!} p_{r_1,\ldots,r_{m-l}},$
$(r_1,\ldots,r_{m-l})\in \{1,2,3\}^{q-j}$ of degree $q$ is left- and right-$(j-1)$-monogenic on $\R^4.$
\end{theorem}
As a consequence we have the following.
\begin{theorem}
If $P_q$ is a homogeneous polynomial of degree $q$ in
$x_0,x_1,x_2,x_3$ is left-$k$-monogenic on an open neighborhood $U$ of the origin then it is left-$k$-monogenic on $\R^4.$
\end{theorem}
\begin{proof}
Since $D^k P_q$ is homogeneous of degree $(q-k)$ it can be written
\begin{equation}
D^kP_q =  \sum_{(r_1,\ldots,r_{q-k})}  x_{r_1}\cdots x_{r_{q-k}} \cdot c_{r_1,\ldots,r_{q-k}}
\end{equation}
where the sum is over all possible combination over $(0,1,2,3)$ with repetition in the sets of
$(q-k)$ elements. Since $P_q$ is left-$k$-monogenic 
$\sum_{(r_1,\ldots,r_{q-k})}  x_{r_1}\cdots x_{r_{q-k}} \cdot c_{r_1,\ldots,r_{q-k}}=0$
on $U$ thus $ c_{r_1,\ldots,r_{q-k}}=0$ thus $D^k P_q =0$ on $\R^4.$
This completes the proof.
\end{proof}
\begin{proposition}
Denote by ${}^k[P_q]$ the set of left-$k$-monogenic functions on $\R^4.$
Note that $({}^k[P_q],+,\mathbb{H},+,x)$ is a right vector space over $\mathbb{H}.$
Then the space
\begin{multline}
{}^k[B_q]:=\\
\{\frac{1}{j!}x_0^j p_{r_1,\ldots,r_{m-l}}:j=0,\ldots,k-1;(r_1,\ldots,r_{m-l}), r_i\in \{1,2,3\}\}
\end{multline}
is a set of right-$\mathbb{H}$-free generators for the vector space ${}^k[P_q]$,
in particular it is a basis for the right vector space 
${}^k[B_q]$ for each $q\in \N.$
\end{proposition}
\begin{definition}
A function $\R^4 \to \mathbb{H}$ is called analytic in some open nonempty subset
$\Omega\subset\R^4$ iff the four components $f_j,$ $j=0,1,2,3,$ are real analytic on $\Omega$
which in turn is equivalent to $f\in C^\infty(\Omega)$ such that
for each compact $K\subset\Omega$ there exists $C_1(K),C_2(K)$ satisfying that for all $q\in \N,$
\begin{equation}
\sup_{x\in K} \abs{\frac{\partial^q f_j(x)}{\partial x_0^{q_0}\cdots \partial x_3^{q_3}} }\leq q! C_1(K)C_2(K)^q,\sum_{i=0}^3 q_i=q,j=0,1,2,3
\end{equation}
\end{definition}
Hence for an analytic $f$ near $0$, setting
\begin{equation}
\frac{\partial^q f(x)}{\partial x_0^{r_0}\cdots \partial x_3^{r_3}}|_{x=0} :=\sum_{i=0}^3 e_j 
\frac{\partial^q f_j(x)}{\partial x_0^{r_0}\cdots \partial x_3^{r_3}}|_{x=0}
\end{equation}
gives an absolutely and uniformly local {\em Taylor expansion} near $0$ according to
\begin{equation}
f(x)= \sum_{m=0}^\infty \frac{1}{m!} \sum_{r_1=0}^3 \cdots \sum_{r_m=0}^3  x_{r_1}\cdots x_{r_{m}}
\frac{\partial^m f(x)}{\partial x_{r_1}\cdots \partial x_{r_m}}|_{x=0}
\end{equation}
\begin{theorem}
If $f$ is left-$k$-analytic in an open nonempty subset $U\subset\R^4,$ then $f$ is analytic on $U$ and for each $p\in U$
there exists an open $U_p\ni p$ on which we have a uniformly convergent expansion
\begin{multline}
f(x)= \sum_{n=0}^\infty \sum_{m=0}^{k-1}  \frac{1}{m!(n-m)!} \sum_{r_1=1}^3 \cdots \sum_{r_{n-m}=1}^3
(x_0-p_0)^m\times\\
(z_{r_1}-p_{r_1}')\cdots (z_{r_{n-m}}-p_{r_{n-m}}')
\frac{\partial^{n-m} D^m f}{\partial x_{r_1}\cdots \partial x_{r_{n-m}}}|_{x=p}
\end{multline}
where $z_j -p_j':=(x_j-p_j)e_0 -(x_0-p_0)e_1,$ $j=1,2,3.$
This expansion of a left-$k$-monogenic function is unique.
\end{theorem}
\begin{theorem}
If $f$ is a function in an open nonempty subset $U\subset\R^4,$ such that for each $p\in U$
there exists an open $U_p\ni p$ on which we have a uniformly convergent expansion
\begin{multline}
f(x)= \sum_{n=0}^\infty \sum_{m=0}^{k-1}  \frac{1}{m!(n-m)!} \sum_{r_1=1}^3 \cdots \sum_{r_{n-m}=1}^3
(x_0-p_0)^m\times\\
(z_{r_1}-p_{r_1}')\cdots (z_{r_{n-m}}-p_{r_{n-m}}')c^{(m)}_{r_1,\ldots,r_{q-k}}
\end{multline}
where each $c^{(m)}_{r_1,\ldots,r_{q-k}}\in\mathbb{H},$ then $f$
is left-$k$-monogenic on $U.$
\end{theorem}

\section{Clifford algebras}
Let us now describe the more general situation.
Recall that the {\em tensor product}\index{Tensor product}, $V\otimes W$ of two vector spaces $V,W$ over a field $K$, is a
pair $(V\otimes W,\otimes)$ where $V\otimes W$ is a $K$-vector space and $\otimes:V\times W\to V\otimes W$
a bilinear map such that for every $K$-vector space, $Z$ and every bilinear map $\varphi:V\times W\to Z$
there exists a unique map $\varphi_\otimes:V\otimes W\to Z$ such that
$\varphi(x,y)=\varphi_\otimes (x\otimes y),$ for all $x\in V,y\in W.$ Note that $\otimes$ depends upon $K$.
The vector space $V\otimes W$ is defined up to isomorphism and the vectors $x\otimes v$ generate $V\otimes W.$
If $V,W$ are $K$-algebras with multiplication $m_V,m_W$ respectively we have associated linear maps
$m_V':V\otimes V\to V,$ $m_W':W\otimes W\to W$ which in turn render
an associated linear map 
\begin{equation}
m_V'\otimes m_W':(V\otimes V)\otimes (W\otimes W)\to (V\otimes W)
\end{equation}
Using associativity and the canonical isomorphism $(V\otimes W)\simeq (W\otimes V)$, we have an isomorphism  
$(V\otimes V)\otimes (W\otimes W)\to (V\otimes W)\otimes (V\otimes W)$ which renders a linear map
\begin{equation}
m_{V\otimes W}: (V\otimes W)\otimes (V\otimes W)\to V\otimes W 
\end{equation}
which defines a multiplication on $V\otimes W$ under which $V\otimes W$ is a $K$-algebra using
\begin{equation}
(x\otimes x') m_{V\otimes W} (y\otimes y')=(xx')\otimes (yy')
\end{equation}
for all $x,x'\in V,$ $y,y'\in W.$ We state without proof the existence of the universal {\em tensor algebra} $T(V)$, for a $K$-vector space $V$,
where $T(V)$ is a $K$-algebra together with a map $\iota:V\to T(V)$ such that given any $K$-algebra $A$ and any linear map
$\varphi:V\to A$, there exists a unique $K$-algebra homomorphism, $\tilde{\varphi}: T(V)\to A$ such that $\varphi=\tilde{\varphi}\circ\iota.$
$T(V)$ can be constructed by setting 
\begin{equation}
T(V):=\oplus_{i\geq 0} \overbrace{V\otimes\cdots \otimes V}^{i\mbox{-times}}
\end{equation}
where $V^0:=K,$ see e.g.\ Atiyah and Macdonald \cite{atiyah}.
\begin{definition}[Clifford algebra]\index{Clifford algebra}
Let $K=\R$, let $V$ be a finite-dimensional $K$-vector space and let $\xi:V\times V\to K$ be a symmetric
bilinear form with associated quadratic form $Q:V\to K$, $Q(x)=\xi(x,x).$ 
A Clifford algebra, $\mbox{Cliff}(V,Q)$, is a unital associative algebra 
together with a linear map $\nu:V\to \mbox{Cliff}(V,Q)$ satisfying:
\\
(1) $(\nu(x))^2=Q(x)\cdot 1,$ (where $1$ denotes the multiplicative identity in the algebra) for all $x\in \mbox{Cliff}(V,Q).$
\\
(2) Given any unital associative algebra $A$ over $K$ and any linear map $\eta : V\to A$ such that 
$(\eta(x))^2 = Q(x)\cdot 1_A$ for all $x\in V$ (where $1_A$ denotes the multiplicative identity in $A$), 
there is a unique algebra homomorphism $\phi :\mbox{Cliff}(V,Q)\to A$
(i.e.\ $\phi$ is a linear map that is also a ring homomorphism such that
$\phi(1_{\mbox{Cliff}}=1_A$) such that $\phi\circ \nu = \eta$.
\end{definition}
It is classic that any two Clifford algebras associated to $V$ and $Q$ must be isomorphic. The existence of 
$\mbox{Cliff}(V,Q)$ is provided by noting that if we set
$\mathcal{J}$ to be the (multiplicative) ideal of $T(V)$ generated by all elements of the form
$x\otimes x-Q(x)\cdot 1,$ $x\in V$, then the map $\nu:V\to \mbox{Cliff}(V,Q)$
is the composition 
\begin{equation}
V\stackrel{i}{\longrightarrow} T(V)\stackrel{\pi}{\longrightarrow} T(V)/\mathcal{J}
\end{equation}
where $i$ is the inclusion map and $\pi$ the natural quotient map.
Any Clifford algebra is associative but not necessarily commutative. 

\begin{example}
Let $V$ be an $n$-dimensional $\R$-vector space and let $\{e_1,\ldots,e_n\}$ be a standard basis for $V.$
Suppose $\nu(x_1 e_1)=x_1^2.$ Then $\mbox{Cliff}(V,Q)$ is spanned by $(1,e_1)$ and we have
$e_1^2=-1$. Hence $\mbox{Cliff}(V,Q)$ is isomorphic to $\C$ with isomorphism given by $e_1\mapsto i.$
If instead $n=2$ and $\{e_1,e_2\}$ the standard basis for $V\simeq \R^2$. Suppose
$Q(x_1e_1+x_2e_2)=-(x_1^2+x_2^2).$ Then $\mbox{Cliff}(V,Q)$ is spanned by
$(1,e_1,e_2,e_1e_2).$ Furthermore, we have 
$e_2e_1=-e_1e_2,$ $e_1^2=-1,$ $e_2^2=-1$ and 
$(e_1e_2)^2=-1.$ Hence under the bijection
\begin{equation}
e_1\mapsto i,\quad e_2\mapsto j,\quad e_1e_2\mapsto k
\end{equation}
we obtain $\mbox{Cliff}(V,Q)\simeq \mathbb{H}.$
\end{example}

\begin{remark}
A quadric form, $Q$, on an $n$-dimensional $\R$-vector space, $V,$ can be represented as
$Q(x)=\sum_{j=1}^n\sum_{k=1}^n q_{jk}x_jx_k,$ which can be identified in matrix notation 
as $x^T [Q] x,$ where $x^T$ denotes the transpose vector and $[Q]:=[q_{jk}]_{jk}$ a symmetric matrix.
The signature (the triplet of the number of positive,zero,negative eigenvalues respectively), 
$g,$ of $[Q]$ is a parameter upon which the Clifford algebra $\mbox{Cliff}(V,Q)$ depends.
A quadric form is called non-degenerate if $[Q]$ has no zero eigen-values. Note that necessarily the sum of 
the elements of the signature is $n.$
In this text we shall consider only the case of so-called {\em standard Clifford algebras} over a finite 
dimensional 
$\R$-vector space, $V,$ with orthonormal basis $\{e_1,\ldots,e_n\},$ 
by which we mean that $\mbox{Cliff}(V,Q)$ is generated
by the basis elements and the non-degenerate quadric form $Q$ of signature $(0,0,n)$ 
chosen to induce 
satisfies the conditions $e_k^2=-1,$ $e_ke_j=-e_je_k,$ $k\neq j.$
Any element $x$ in the standard Clifford algebra can be written
$x=\sum_J c_J e_J$ where $J=\{i_1,\ldots,i_k\}$ is any subset (possibly empty) of 
$\{1,\ldots,n\}$ and $i_1<\cdots<i_k,$ and $e_J:=e_{i_1}\cdots e_{i_k},$ $e_\emptyset :=1,$
the identity element.
So the dimension over $\R$ is $\sum_{k=0}^n\binom{n}{k}=2^n.$
In the case of non-degenerate $Q$, sometimes authors 
denote by $\R_{0,n}$ the standard Clifford algebra associated to
$\R^n$ and $Q$ a nondegenerate quadric form with signature $(0,0,n)$.
\end{remark}

\section{$H^*$-algebras and approximation}
 \begin{definition}
 	A complex algebra (not necessarily commutative) $A$ is called a $H^*$-algebra if 
 	it is equipped with an inner product $(,)$ and an involution $x\to \overline{x}$ such that 
 	(i) $(xy,z)=(y,\overline{x}z).$ (ii) $(yx,z)=(y,z\overline{x})$ for all $x,y,z\in A$ and (iii) $A$
 	is a Banach algebra for the norm $\norm{\cdot}_0$ induced by the inner product.
 	\end{definition}
 Recall that an $H^*$-algebra is called {\em proper} if it contains no elements that annihilate the whole algebra.
 	Ambrose \cite{ambrose} proved that if $A$ is a proper $H^*$-algebra then $A$ is an involutive Banach algebra.
  \begin{proposition}
 	If $A$ is an $H^*$-algebra and $a\in A$ then $aA=(0)$ iff $Aa=(0).$
 \end{proposition}
 \begin{proof}
 	Let $b,c\in A$. Since $ab=0$ we have
 	\begin{equation}
 	0=(ab,c)=(a,c\bar{b})=(\bar{c}a,b)
 	\end{equation}
 	Since $b$ was arbitrary we have $\bar{c}a=0$ and since $c$ was arbitrary this implies $Aa=(0).$
 	\end{proof}
 	Hence an $H^*$-algebra is proper iff $aA=(0)\Rightarrow a=0$ (and this is true iff $Aa=(0)\Rightarrow a=0$).  
 	
 	\begin{proposition}
 		If $A$ is an $H^*$-algebra then $A$ is proper iff each element $a\in A$ has a unique adjoint.
 	\end{proposition}
 	\begin{proof}
 		Let $b,c\in A$ and let $a', a''$ be two adjoints of $a.$ 
 		If $A$ is proper the 
 		\begin{equation}
 		0=(ab,c)=(b,a'c)=(b,a''c)
 		\end{equation}
 		thus
 		\begin{equation}
 		(b,(a-a'')c)=0,\quad \forall b,c\in A
 		\end{equation}
 		Hence $(a'-a'')c=0$ for all $c\in A$ which yields $a'-a''=0.$
 		If instead $A$ is not proper then there exists $d\in A$, $d\neq 0$ such that $dA=Ad=(0).$ If $a$ is any element 
 		and $a'$ any adjoint of $a$also $a'+d$ is an adjoint of $a$ This completes the proof.
 		\end{proof}
 
 	Let $A$ be a Clifford algebra over a quadric $n$-dimensional real vector 
 	space with orthogonal basis $\{ e_1,\ldots,e_n\},$ 
 	$e_i^2=-e_0$, $i=1,\ldots,n$ where $e_0$ is the identity in $A.$
 	Let furthermore $e_\alpha=e_{i_1}e_{i_2}\cdots e_{i_h}$ be an arbitrary basis
 	element of $A$ where $i_1<i_2\cdots <i_h,$ $i_j\in \{1,\ldots,n\}, j=1,\ldots,h,$
 	$\alpha=\{i_1,\ldots,i_n\}.$ We shall simply denote $\alpha\in p_h(\{1,\ldots,n\})$.
 	We define an inner product $(,)_0$ and norm $\abs{\cdot}_0$ and an involution on $A$ that turns $A$ into a $H^*$-algebra\index{$H^*$-algebra}.
 	For each $\lambda=\sum_{\alpha} \lambda_\alpha e_\alpha\in A$ we use the norm
 	$\abs{\lambda}_0^2:=2^n \sum_\alpha \lambda_\alpha^2$ on $A.$
 	Now recall that $\R_{0,m}$ has a multivector structure where each element $a\in \R_{0,m}$ takes the form
 	\begin{equation}
 	a=\sum_\alpha a_\alpha e_\alpha=\sum_{k=0}^m \sum_{\abs{\alpha}=k} a_\alpha e_\alpha,\quad a_\alpha\in \R
 	\end{equation}
 	which yields the decomposition
 	\begin{equation}
 	\R_{0,m}=\R^0_{0,m}\oplus \cdots\oplus \R^m_{0,m}
 	\end{equation}
 	\begin{equation}
 	\R^k_{0,m}=\mbox{Span}_\R \{a_\alpha :\alpha=(\alpha_1,\ldots,\alpha_k),1\leq \alpha_1<\cdots <\alpha_k\leq m\}
 	\end{equation}
 	Hence $a$ can be written as
 	\begin{equation}
 	a=[a]_0+[a]_1+\cdots+[a]_m
 	\end{equation}
 	where $[a]_k$ denotes the projection on the space $\R^k_{0,m}$ of so-called {\em $k$-vectors}\index{$k$-vectors}.
 	We can introduce involutions on the (standard) Clifford algebra $\R_{0,m}$ in different ways, and it suffices to
 	introduce these on the basic elements $e_\alpha$ and then extend by linearity to the Clifford algebra.
 	We use the conjugation $a\mapsto \bar{a}$ defined by
 	\begin{equation}
 	\overline{e}_j=-e_j,\quad j=1,\ldots,m,\quad \overline{ab}=\bar{b}\bar{a}
 	\end{equation}
 	We mention that other examples of involutions are given by (reversion) $a\mapsto a^*,$
 	defined by $e_j^*=e_j,$ for the basis elements $j=1,\ldots,m$, and $(ab)^*=b^*a^*$ for $a,b$ in the Clifford algebra $R_{0,m}$.
 	It follows that $e^*_\alpha=e_{\alpha_k}\cdots e_{\alpha_1}$ when $e_\alpha=e_{\alpha_1}\cdots e_{\alpha_k}$
 	and $e_\alpha^*=(-1)^{\frac{k(k-1)}{2}} e_\alpha$ whenever $\abs{\alpha}=k.$
 		Another example is $a\mapsto \tilde{a}$ given by $\tilde{e}_j=-e_j$ on the basis elements and $\tilde{ab}=\tilde{a}\tilde{b}.$
 		It follows that $\tilde{e_\alpha}=(-1)^ke_\alpha$ whenever $\abs{\alpha}=k.$
 		This clearly extends to the complexification $\C_m:=\C\otimes \R_{0,m}$. Recall that as an associative $\C$-algebra
 		$\C_m$ has dimension $2^m$ and each element $\lambda\in \C_m$ can be written
 		$\lambda=\sum_\alpha \lambda_\alpha e_\alpha,$ $\lambda_\alpha\in \C,$
 		or $\lambda=a+ib,$ $a,b\in \R_{0,m}.$ Once can define the so-called Hermitian conjugation $^\dagger$ on $\C_m$
 		according to
 		$\lambda^\dagger=\bar{a}-i\bar{b}$ for $\lambda\in \C_m$ or 
 		$\lambda^\dagger=\sum_\alpha \overline{\lambda_\alpha} \bar{e}_\alpha$ where $\overline{\lambda_\alpha}$ denotes the usual conjugation in $\C.$
 		\\
 		Clearly, $a\mapsto \overline{a}$ is the composition of $a\mapsto \tilde{a}$ with $a\mapsto a^*$ hence
 		\begin{equation}
 		\overline{e}_\alpha=(-1)^{\frac{k(k+1)}{2}} e_\alpha,\mbox{ whenever }\abs{\alpha}=k
 			\end{equation} 
 			For any $a\in \R_{0,m}$ we have
 			\begin{equation}
 			\bar{a}=[a]_0-[a]_1-[a]_2+[a]_3+[a]_4-[a]_5-\cdots
 			\end{equation}
 			The conjugation $a\mapsto \bar{a}$ induces an inner product and associated norm on $\R_{0,m}$ according to 
 			\begin{equation}
 			(a,b)=[\bar{a}b]_0=[b\bar{a}]_0,\quad \abs{a}=\sqrt{[a\bar{a}]_0}=\sqrt{\sum_\alpha a_\alpha^2}
 				\end{equation}
 				for $a,b\in \R_{0,m}$.
 				The Hermitian conjugation on $\C_m$ yields correspondingly
 				\begin{equation}
 				(\lambda,\mu)=[\lambda^\dagger]_0=[\lambda \mu^\dagger]_0,\quad \abs{\lambda}=\sqrt{[\lambda\lambda^\dagger]_0}=\sqrt{\sum_\alpha \abs{\lambda_\alpha}^2}
 					\end{equation}
 					for $\lambda,\mu\in \C_{m}$.
 	We shall (following Brackx, Delanghe \& Sommen \cite{brackxbook}) modify the above inner product by a scalar $2^n$ as follows.
  	\begin{definition}
 		Let $A$ be a finite dimensional (standard universal) Clifford algebra with basis $\{e_1,\ldots,e_n\}.$
 		The subspace spanned by the $\binom{n}{p}$ products $e_\alpha,$ cardinality of $\alpha$ equal to $p,$ the subspace will be
 		denoted $A_p.$
 		Let $\lambda=\sum_\alpha e_\alpha \lambda_\alpha$ be a Clifford number. 
 		The coefficient $\lambda_\alpha$ of $e_\alpha$ will be denoted $[\lambda]_0$
 		The number $[\lambda]_0$ is called the {\em scalar part}\index{Scalar part of a Clifford number}
 		of $\lambda.$
 		An inner product in $A$ is defined by putting for all $a,b\in A$
 		\begin{equation}
 		(a,b)_0=2^n [a,\bar{b}]_0=2^n\sum_\alpha \lambda_\alpha b_\alpha
 		\end{equation}
 		\end{definition}
 		\begin{remark}\label{clifforrem0101}
 		Note that 
 		\begin{equation}
 		(a,b)_0=(b,a)_0=(\bar{a},\bar{b})_0=(\bar{b},\bar{a})_0
 		\end{equation}
 		\begin{equation}
 		\abs{\lambda}_0=\sqrt{(\lambda,\lambda)_0}=2^{\frac{n}{2}}\sqrt{[\lambda,\bar{\lambda}]_0}
 		2^{\frac{n}{2}}\sqrt{\sum_\alpha \lambda_\alpha^2}
 		\end{equation}
 		In this way $A$ is a real Hilbert space and a Banach algebra with $\abs{ab}\leq \abs{a}_0\abs{b}_0.$
 		Note that $e_\alpha\bar{e}_\alpha =\bar{e}_\alpha e_\alpha =e_0$, 
 		$\abs{e_\alpha}_0=2^{\frac{n}{2}} \neq 1$ and
 		if $\lambda=\lambda_0e_0+\sum_{i=1}^n \lambda_i e_i \in A_0\oplus A_1$,
 		$\mu=\mu_0e_0+\sum_{i=1}^n \mu_i e_i \in A_0\oplus A_1$ 
 		then $\lambda\bar{\lambda}=\bar{\lambda}\lambda =\sum_{i=0}^n (\lambda)^2$
 		and $\abs{\lambda\mu}_0=2^{\frac{n}{2}}\abs{\lambda}_0\abs{\mu}_0,$
 		$\abs{\lambda}_0^2=2^n \sum_{i=0}^n \lambda_i^2,$
 		$\lambda^{-1}=2^n \frac{\bar{\lambda}}{\abs{\lambda}_0^2}$ for $\lambda\neq 0.$
 			\end{remark}
 	
 	\begin{proposition}
 		The finite dimensional Clifford algebra $A$ equipped with inner product $(\cdot,\cdot)$ and involution $a\mapsto \overline{a}$
 		is a finite dimensional $H^*$-algebra.
 	\end{proposition}
 	\begin{proof}		
 		It suffices to prove that for all $\lambda,\mu,\nu\in A$
 		\begin{equation}
 		(\lambda\mu,\nu)_0=(\lambda,\bar{\mu}\nu)_0=(\lambda,\nu\bar{\mu})_0
 		\end{equation}
 		Now
 		\begin{equation}
 		(\lambda\mu,\nu)_0=2^n[\lambda\mu\bar{\nu}]_0=2^n[\lambda\overline{\nu\bar{\mu}}]_0
 		=(\lambda,\nu\bar{\mu})_0
 		\end{equation}
 		Also we have
 		\begin{equation}
 		(\mu,\bar{\lambda}\nu)_0=2^n[\mu\overline{\bar{\lambda}\nu}]_0=2^n[\mu\bar{\nu}\lambda]_0
 		=\\(\overline{\nu\bar{\mu}},\bar{\lambda})_0=(\nu\bar{\mu},\lambda)_0
 		=(\lambda,\nu \bar{mu})_0
 		\end{equation}
 		This completes the proof.
 	\end{proof}
 	When $A$ is a proper finite dimensional algebra then it has an identity $e$ with $\norm{e}_0\geq 1.$ 
 	In this case (see Bonsail \& Goldie \cite{bonsai}) the so-called trace $\tau(a,b)=\re(b,\bar{a}), a,b\in A$
  	is nonsingular which means that if for fixed $a\in A$
 	$\tau_a(b)=\tau(ab)=0$ for all $b\in A,$ then $a=0.$ To see that the trace is nonsingular note that if the linear functional 
 	$\langle \tau_a,b\rangle=\tau(ab),$ satisfies that for all $b\in A$
 	$\langle \tau_a,b\rangle=0$ then
 	\begin{equation}
 	\langle \tau_a,\bar{a}\rangle=\tau(a\bar{a})
 	=\re(\bar{a},\bar{a})=\re(a,a)=\re\abs{a}^2=0
 	\end{equation}
 	Hence $a=0.$
 	This can be used to show that for any real linear functional $\phi$ on $A$ there exists a unique $a\in A$ such that $\phi(x)=\tau_a(x)=\tau(ax),$ $x\in A.$
 	In particular, for $a=e_\alpha$ we have
 	$\langle \tau_{e_\alpha},x\rangle=2^n(-1)^{n(\alpha)(n(\alpha)+1)/2} a_\alpha,$
 	where $n(\alpha)$ denotes the cardinality of $\alpha.$
 	\begin{theorem}[Delanghe \cite{delanghe1980}, Thm. 1.1]
 		If $X$ is a right $A$-module (in particular we assume $xe=x$ for all $x\in X$ so that $X$ is a real or complex vector space) 
 		and $T$ a right $A$-linear functional on $X$, then $\tau_e T$ is a real linear functional on $X$ and $T=0$ iff $\tau_e T=0.$
 	\end{theorem}
 	\begin{proof}
 		$\R$-linearity of $\tau_e T$ is immediate and also that $T=0$ implies $\tau_e T=0$ from the definition. Assume
 		$\tau_e T=0$ and let $x\in X.$ The for each $\lambda\in A$ we have
 		\begin{equation}
 		\tau_{T(x)}(\lambda)=\tau(T(x)\lambda)=\tau_e(T(x)\lambda)=\tau_e(T(x\lambda))=0
 		\end{equation}
 		Since $\tau$ is a nonsingular trace operator on $A$ this implies $T(x)=0.$
 		This completes the proof.
 	\end{proof}
 	Let $m\leq n,$ $m\neq 0$ and let $\Omega\subset\R^{m+1}$ be a nonempty open subset.
 	Denote by $M_k(\Omega,A)$\index{$M_k(\Omega,A)$, $k$-monogenic functions with values in an $H^*$-algebra $A$} the set of functions $f\in C^k(\Omega,A)$ such that $D^k f=0$ on $\Omega,$
 	where $k\in \N$ and $D:=\sum_{i=0}^m e_i\partial_{x_i}.$
 	$D$ is the hypercomplex generalization of the Cauchy-Riemann operator and when $n=m=1$ the solutions 
 	to $D^k f=0$ on $\Omega$ are precisely the $k$-analytic functions.
 	For general $n,m,k$ the space $M_k(\Omega,A)$ is the subspace of the set of $\R^{2^n}$-valued $k$-harmonic functions.
 	It is known that $M_k(\Omega,A)$ equipped with the topology of uniform convergence on compacts, is a right $A$-Fr\'echet module.
 	Since dim$A=2^n$ the equation $D^k f=0$ is equivalent to a system of $2^n$ linear partial differential equations, each of order $k$, in the unknown real valued functions $f_\alpha.$
 	If the basis elements $e_\alpha$ of $A$ are ordered in a certain way, then the left regular representation
 	of $A$ allows us to associate to each $\lambda\in A$ a $2^n\times 2^n$ real matrix $\Theta(\lambda).$ Since $A$ has an identity
 	this representation is an isomorphism.
 	Setting $\overline{D}=\sum_{i=0}^m \overline{e}_i\partial_{x_i}=e_0\partial_{x_0}-\sum_{j=1}^m e_j\partial_{x_j}$ and $\Delta=\sum_{i=0}^m \partial^2_{x_i}$
 	we have
 	$D\overline{D}-\overline{D}D=\Delta e_0$.
 	\begin{proposition}
 		The system of differential equations associated to the hypercomplex differential operator $D^k$ is strongly elliptic.
 	\end{proposition}
 	\begin{proof}
 		We have $\Theta(D^k)=(\Theta(D))^k.$ Since power of a strongly elliptic operator is strongly ellitic we only need
 		to prove the result for the first order $\Theta(D).$ The multiplication rules
 		for $A$ imply that for each $a,b\in PN$ the element $x_{a,b}$ at the $b$:th row and $a$:th column of $\Theta(e_j$ equals $-x_{a,b}$.
 		Hence $(\Theta(e_j))^T+\Theta(e_j)=0$ or
 		$(\Theta(e_j))^T=\Theta(-e_j)=\Theta(\overline{e_j}),$ $j=1,\ldots,n$ where
 		$(\Theta(e_j))^T$ denotes the transpose.
 		Since $\Theta(e_0)=\delta$ is the identity matrix we have $\Theta(D)^T=\Theta(\overline{D}),$ which implies
 		\begin{equation}
 		\Theta(\overline{D}D)=\Theta(\overline{D})\Theta(D)=\Delta \delta
 		\end{equation}
 		Hence
 		\begin{equation}
 		(\mbox{det}\Theta(D))^2=\Delta^{2n}
 		\end{equation}
 		This completes the proof.
 	\end{proof}
 \begin{definition}
 	Let $A$ be a finite dimensional (standard universal) Clifford algebra and let $X$ be a right $A$-module. A family $\mathscr{P}$ of functions $p:X\to \R$
 	is called a {\em proper system of semi-norms} on $X$\index{Proper system of semi-norm in hypercomplex analysis}
 	if:
 	\begin{itemize}
 		\item[(i)] For any finite number of elements $p_1,\ldots,p_k\in \mathscr{P}$ there exists 
 		$p\in \mathscr{P}$ and $C>0$ such that for all $f\in X$
 		\begin{equation}
 		\sup_{i=1,\ldots,k} p_i(f)\leq Cp(f)
 		\end{equation}
 		\item[(ii)] If $p(f)=0$ for all $p\in\mathscr{P}$ then $f=0$
 		\item[(iii)] For any $p\in \mathscr{P},$ $\lambda\in A$ and $f,g$ we have
 		\begin{equation}
 		p(f\lambda)\leq \abs{\lambda}_0 p(f) \mbox{ and }p(f\lambda)=\abs{\lambda}p(f)\mbox{ if }\lambda\in \C
 		\end{equation}
 		\begin{equation}
 		p(f+g)\leq p(f)+p(g)
 		\end{equation}
 		A function $\norm{\cdot}$ is called a {\em proper norm}\index{Proper norm} on $X$ if it satisfies (iii) and
 		$\norm{f}=0\Rightarrow f=0.$
 		A function $T:X\to A$ is called a right-linear functional on $X$ if for all $\lambda,\mu\in A$
 		and $f,g\in X$ we have
 		$T(f\lambda +g\mu)=T(f)\lambda +T(g)\mu.$
 		If $X$ is a right $A$-module with a proper system of semi-norms $\mathscr{P}$ the a right $A$-linear functional $T$ on $X$ is called {\em bounded}
 		if there exists $C>0$ and $p\in \mathscr{P}$ such that $\abs{T(f)}_0\leq Cp(f)$ for all $f\in X.$ 
 	\end{itemize}
 \end{definition}
 
 \begin{definition}
 	Let $\Omega\subset\R^{m+1}$ be a nonempty open subset and $M_k(\Omega,A)$ as above.
 	Let $\{E_j\}_{j\in \N}$ be a sequence of compacts of $\Omega$
 	such that $E_j\subset \mbox{int}E_{j+1}$ and the $\mbox{int}E_j\uparrow\Omega$. Then the functions 
 	$p_j:M_k(\Omega,A)\to \R$, $p_j(f)=\sup_{x\in E_j}\abs{f(x)}_0$ define a proper countable system of semi-norms on $M_k(\Omega,A)$
 	and we denote by $M_k(\Omega,A,\mathscr{P})$ the space $M_k(\Omega,A)$ equipped with this system.
 \end{definition}
 
 \begin{proposition}
 	$M_k(\Omega,A,\mathscr{P})$ is a right Frech\'et $A$-module.
 \end{proposition}
 \begin{proof}
 	Let $\{f_j\}_{j\in \N}$ be a Cauchy sequence in $M_k(\Omega,A,\mathscr{P})$.
 	Then for all $\epsilon>0$ and any compact $K\subset\Omega$ there exists $N(\epsilon,K)\in \N$
 	such that
 	\begin{equation}
 	\sup_{x\in K}\abs{f_r(x)-f_s(x)}_0\leq \epsilon ,\quad r,s\geq N(\epsilon,K)
 	\end{equation}
 	Passing to the components $f_{j,\alpha}$ of the $f_j$ we have for all $\alpha$
 	\begin{equation}
 	\sup_{x\in K}\abs{f_{r,\alpha}(x)-f_{s,\alpha}(x)}\leq \epsilon ,\quad r,s\geq N(\epsilon,K)
 	\end{equation}
 	Since each $f_j\in M_k(\Omega,A)$ we have $\Delta^k f_j =0$ on $\Omega$
 	or for each $\alpha$, $\Delta^k f_{i,\alpha}=0$ on $\Omega.$
 	This implies that there exists functions
 	$f_\alpha\in C^\infty(\Omega,\R)$ such that
 	$\{f_{j,\alpha}\}_{j\in\N}$ together with all its derived sequences
 	$\{f^{(\alpha)}_{j,\alpha}\}_{j\in\N}$, $\alpha\in \N^{m+1}$
 	we have
 	$f^{(\alpha)}_{j,\alpha}\to f^{(\alpha)}$ uniformly on compacts of $\Omega.$
 	Hence setting $f=\sum_\alpha f_\alpha e_\alpha$ we have $f\in C^\infty(\Omega,\R)$ and for any $\alpha\in \N^{m+1},$
 	$f^{(\alpha)}_j\to  f^{(\alpha)}$ uniformly on compacts of $\Omega$. 
 	In particular, $D^j f_j\to D^j f$ uniformly on compacts of $\Omega$ 
 	so $D^k f=0$i.e.\ $f\in M_k(\Omega,A).$
 	Hence $f_j\stackrel{\mathscr{P}}{\longrightarrow} f$. This completes the proof.
 		This completes the proof.
 	\end{proof}
 
\begin{definition}
	Let $\Omega\subset\R^{m+1}$ be an open subset and $A$ a (universal) Clifford algebra as before.
	Let $K\subset\R^{m+1}$ be a compact. We denote by $M_k(K,A)$ the set of functions $f$ for which there exists an open
	neighborhood $U$ of $K$ such that $f\in M_k(U,A).$
\end{definition}
Brackx \& Delanghe \cite{brackxdelanghe1980} prove several approximation theorems, here mention just as an example the following.
\begin{proposition}
	Let $K\subset\R^{m+1}$ be a compact with connected complement. 
	Then for each $f\in M_k(K,A)$ and each $\epsilon>0$ there exists $g\in M_k(\R^{m+1},A)$ such that $\sup_{x\in K}\abs{f(x)-g(x)}_0\leq \epsilon$
	(in other words $M_k(\R^{m+1},A)$ is uniformly dense in $M_k(K,A)$).
\end{proposition}
 We believe that we have provided the reader with sufficient introductory knowledge in order to be able to 
 further explore the cited literature with respect to polyanalyticity from the perspective of hypercomplex analysis.

\chapter{Reproducing kernels spaces of polyanalytic functions}\label{reproducingsec}
Recall that a Hilbert space is a vector space with an inner product space such that
with respect to the induced metric, it becomes a complete metric space.
A Banach space is a complete normed vector space. 
Let us begin by a standard example.
Let $I\subset \R$ be a closed interval and let $k\in C^0(I\times I, \R),$
and $g\in C^0(I,\R).$ The {\em Fredholm integral equation} with data $k,g$ is
the fixed point equation
\begin{equation}\label{aaa}
Tf=f, Tf(x):=g(x)+\int_I k(x,y)f(y)dy
\end{equation}
Since $C^0(I,\R)$ can be realized as a complete (the uniform limit of continuous functions is continuous) Banach space under the sup-norm, the Banach fixed point theorem (contraction mapping theorem) implies that
there exists a unique fixed point (i.e.\ a unique solution $f\in C^0(I,\R)$
to $Tf=f$) whenever $T$ is a contraction.
Clearly, it suffices that $\sup_{x\in I}\int_I k(x,y)f(y)dy <1$ 
since for $f_1,f_2\in C^0(I),$ $\norm{Tf_1-Tf_2}\leq \norm{f_1-f_2}\sup_{x\in I}\int_I k(x,y)f(y)dy.$
It can be proved that the map $\tilde{k}:g\mapsto \int_I k(x,y)g(y)dy,$
satisfies $\lim_{j\to\infty} \tilde{k}^j f_0 =0.$
Furthermore, the Banach fixed point theorem
states that the fixed point $f$ can be obtained as $\lim_{j\to \infty} T^j f_0$ where $f_0$ is an arbitrary starting member
of the Banach space. Hence $f=\lim_{j\to\infty}  (T^j f_0+\sum_{l=0}^{j}\tilde{k}^j g)=\sum_{j=0}\tilde{k}^j g.$
The Banach fixed point theorem applies to any complete metric space so one can
generalize the above situation to complex analysis using so called reproducing kernels and appropriately normed spaces.
In the case of analytic and harmonic functions, Bergman \cite{bergman} was a pioneer (see also Bergman \& Schiffman \cite{bergmanshiff}).
with regards to introducing reproducing kernels on normed subspaces of $L^2(\C)$ and use these in the theory of boundary value problems.
If a class of functions $\mathcal{F}$ on $\Cn$, forms a complex Hilbert space the a function $k(x,y)$ is called a reproducing kernel of $\mathcal{F}$
if for each $y\in \Cn$, $k(x,y)\in \mathcal{F}$ and for each $y\in \Cn, f\in \mathcal{F},$ $f(y)=\langle f(x),k(x,y)\rangle_x$
where the subscript $x$ indicated that the scalar product applies to functions of $x.$
Note that if $k'(x,y)$ is another reproducing kernel then by the reproducing property we have
$\norm{k(x,y)-k'(x,y)}^2=\langle k-k',k\rangle -\langle k-k',k'\rangle=0$.
The quadric form $\sum_{i,j=1}^n k(y_i,y_j)\bar{\xi}_i\xi_j$ is non-negative for all $y_1,\ldots,y_n$ in $\Cn,$ (see e.g.\ Aronzajn \cite{aronzajn})
in particular, $k(x,x)\geq 0,$ $k(x,y)=\overline{k(y,x)},$ $\abs{k(x,y)}^2\leq k(x,x)k(y,y).$
Conversely, to each positive matrix $k(x,y)$
there corresponds a unique class of functions
with a uniquely determined quadratic form in it, forming a Hilbert space
and admitting $k$ as a reproducing kernel. 
If $\mathcal{F}$ is a sub-Hilbert space of a larger Hilbert space $\tilde{\mathcal{F}}$ then 
$f(y)=\langle h(x),k(x,y)\rangle_x$ gives the projection $f$ of the element $h$ in $\tilde{\mathcal{F}}$.
If $\{f_j\}_{j\in \N}$ is an orthonormal system in $\mathcal{F}$ then for each sequence $\{ a_j\}_{j\in \N}$ 
satisfying $\sum_{j=1}^{\infty} \abs{a_j}^2<\infty,$
we have $\sum_{j=1}^{\infty} \abs{a_j}\abs{f_j(x)}\leq (k(x,x))^{\frac{1}{2}}\left(\sum_{j=1}^{\infty} \abs{a_j}^2\right)^{\frac{1}{2}}.$ 
For a fixed $y$ the Fourier coefficients of $k(x,y)$ for the system $\{f_j\}_{j\in \N}$
are $\langle k(x,y), f_j(x)\rangle =\overline{\langle f_j(x),k(x,y)\rangle} =\overline{f_j(y)}.$ Hence
$\sum_{j=1}^{\infty} f_j(y)^2\leq  \langle k(x,y),k(x,y)\rangle_x=k(y,y).$
and thus $\sum_{j=1}^{\infty} \abs{a_j}\abs{f_j(x)}\leq k(x,x)^{\frac{1}{2}}\left(\sum_{j=1}^{\infty} \abs{a_j}^2\right)^{\frac{1}{2}}.$
Recall that for a given linear differential operator $L:C^\infty(\Omega)\to C^\infty(\Omega)$ on a domain $\Omega$ in Euclidean space,
a {\em fundamental solution}\index{Fundamental solution}, $E$, for $L$ is defined as a weak solution to the equation $LE=\delta_x$, for
each fixed $x\in \Omega,$ where $\delta_x$ denotes the Dirac distribution, i.e.\ $\delta_x *f = f(x)$ where $*$ denotes convolution.
Given a fundamental solution $E$ we can solve the equation $Lf(x)=g(x)$, using $E*g$ since $L(E*g)=(LE)*g=g(x).$
If $k(x,y)$ is a reproducing kernel 
and the inner product is defined in terms of integration then in $f(y)=\langle f(x),k(x,y)\rangle_x$ the right hand side can be viewed as 
a convolution $k*f$. This closely relates fundamental solutions to reproducing kernels.
In the introduction are given several reproducing kernels generalizing those for the Laplacian and the Cauchy-Riemann operator
respectively.

\section{Some Fredholm theory}\label{fredholmsec}
Fredholm theory has many applications in partial differential equations. Also as we shall see
it is necessary to have some basic knowledge about Fredholm theory if one is to study the literature on 
the Toeplitz operators.
In terms of an eigenvalue problem we can identify Eqn.(\ref{aaa}) as the case $\lambda=1$ of the operator equation
\begin{equation}
(\tilde{k}-\lambda I)(f)=g
\end{equation} 
where $\tilde{k}$ is a linear operator from a Banach space to itself and $I$ the identity operator.
If $\tilde{k}$ is bounded the homogeneous case ($g=0$) 
has solutions for at most countably many eigenvalues $\lambda=\lambda_j,$ $j\in\N$, and the
corresponding solutions are called eigenfunctions. The eigenvalues are real if
the kernel is symmetric. When $g\not\equiv 0$ then the equation has solutions whenever $\lambda$ is not an eigenvalue.
Recall that for $\abs{\lambda}< \norm{\tilde{k}}^{-1}$ the resolvent $(\tilde{k}-\lambda I)^{-1}$ can be expressed via the Neumann series
i.e.\ a powers series in $\tilde{k},$ as we have seen in the above example. 

More generally a bounded linear operator between two Banach spaces $T:B_1\to B_2$ is called {\em Fredholm}
if (i) dim$(\mbox{Ker}(T))<\infty$. (ii) The range, Ran$(T)\subset B_2$, is closed in $B_2$. (iii) dim$(B_2/$Ran$(T))<\infty$. 
Recall that the set of compact operators is closed in the set of bounded linear operators, $L(B_1,B_2).$ 
For any $T\in L(B_1,B_2)$, Ker$(T)$ is always closed and we always have a bijection $\hat{T}:B_1/$Ker$(T)\to$Ran$(T)$, and if Ran$(T)$ is closed then $\hat{T}^{-1}$ is bounded.
If dim$(\mbox{Ker}(T))<\infty$ then there is a closed subspace $B_1'\subset B_1$ such that $B_1=\mbox{Ker}(T)\oplus B_1'.$
Similarly, if (iii) holds then there is a closed subspace $B_2'\subset B_2$ such that $B_2=\mbox{Ran}(T)\oplus B_2'.$
Hence one could regard a Fredholm operator as a kind of isomorphism modulo finite dimensional subspaces. In particular, the conditions (i),(ii),(iii)
render a kind of well-posedness for the equation $Tf=g.$ 

If $T$ is a Fredholm operator, then the solvability of the equation $Tf = g$ for a given
$g$ is equivalent to determining whether $g$ is orthogonal to the finite-dimensional
subspace Ker$(T^*)$. The space of solutions of the equation $Tf = g$ for a given
$g$ is finite dimensional.
We mention in passing that it turns out that the
difference of  dim Ker$(T)-$dim Ker$(T^*)$ (called the classical {\em index})\index{Index of Fredholm operator}
plays an important role in the study of Fredholm operators as it is invariant under compact perturbations.
It turns out a bounded linear operator $T:H_1\to H_2$ between Hilbert spaces is Fredholm iff it is 
of the form $T=T_1-T_2$ for an isomorphism $T_1$ and a finite rank (in particular compact) operator $T_2$, see e.g.\ Ramm \cite{ramm}.
Many theorems exists on the solvability of equations of the form $Tf=g$ when $T$ is Fredholm.
For example, for Hilbert spaces, we have the following. 
\begin{theorem}[See e.g.\ Hsiao \& Wendland \cite{hsiao}, Thm 5.3.7]
Let $T:H_1\to H_2$ be a bounded linear operator which is Fredholm, where
$H_1,H_2$ are Hilbert spaces (see our reference to e.g.\ Ramm \cite{ramm}).
Then precisely one of the following holds true:\\
(i) Ker$(T)=\{0\}$ and Ran$(T)=H_2$\\
\\
(ii)
$0<$dim Ker$(T)=$dim ker$(T^*)<\infty$ and the equations
$Tf=g, T^*u=v$ have solutions $f\in H_1$ and $u\in H_2$ respectively iff the corresponding right hand side satisfies the finitely many orthogonality conditions
\begin{equation}
\langle g,u_0\rangle=0,\quad \forall u_0\in\mbox{Ker}(T^*) 
\end{equation}
or
\begin{equation}
\langle v, f_0\rangle=0,\quad \forall f_0\in\mbox{Ker}(T) 
\end{equation}
respectively.
If these conditions are satisfied then the general solution is of the form
\begin{equation}
f=f^*+\tilde{f_0}\quad (u=u^*+\tilde{u_0}\mbox{ respectively})
\end{equation}
where $f^*$ ($u^*$ respectively) are particular solutions depending continuously on $f$ ($u$ respectively)
and $\tilde{f_0}\in$Ker$(T)$ ($\tilde{u_0}\in$Ker$(T^*)$ respectively).
\end{theorem}
See also Begehr \cite{begehrbok} p.227, regarding the application of Fredholm theory to singular integral equations that are highly relevant in the theory of 
higher order complex equations.
For the case of Banach space one could formulate an analogue of the theorem by saying that
for a Fredholm operator between two Banach spaces
the equation $Tf=g$ is solvable for every $g\in$Ran$(T)$ and
the solution is unique modulo Ker$(T).$  
Here are some ways to derive the Fredholm property for bounded linear operators.
\begin{theorem}[See e.g.\ Rudin \cite{rudin} Thm. 4.25]\label{rudinthm0}
Let $K$ belong to the set of bounded linear operators, $L(B_1,B_2)$, between the Banach spaces $B_1,B_2.$
If $K$ is compact then $I+K$ is Fredholm.
\end{theorem}
\begin{proof}
Let $B$ be the unit ball in Ker$(I + K )$. Then $B = K (B)$ thus the image of a bounded set under a compact operator hence is 
precompact. Since it is closed it is compact. By Riesz’s lemma Ker$(I + K )$ is finite dimensional. 
Let $\{x_i\}_{i\in \N}$ be a bounded sequence such that $x_i + K_i x_i$ converges to $y\in B_2.$ 
Since $\{x_i\}_{i\in \N}$ is bounded there is a subsequence $\{x_{i_j}\}$ such that $\{K x_{i_j} \}$ converges, which in turn implies that 
$\{x_{i_j} \}$ converges. Hence Ran$(I + K )$ is closed.
Applying the same argument to the adjoint $I + K^*$ renders the result. 
\end{proof}
\begin{theorem}[See e.g.\ Abramovic \& Aliprantis \cite{abramovic}, p.158]\label{abramovicthm}
If $F:B_1\to B_2$, $G:B_2\to B_3$ are two Fredholm operators on Banach spaces then $GF$ is Fredholm.
\end{theorem}
\begin{lemma}
\label{rudinthm1}
If $S\in L(B_1,B_2)$ and $K:B_2\to B_3$ is compact, for Banach spaces $B_1,B_2,B_3.$
$KS$ is compact. If $B_3=B_1$ then $SK$ is compact.
\end{lemma}
\begin{proof}
Let $\{x_i\}_{i\in \N}$ be a bounded sequence in $B_1$. Since $S$ is bounded, the 
sequence $\{Sx_i\}_{i\in \N}$ is bounded in $B_2$. Since
$K$ is compact the sequence $\{Sx_i\}_{i\in \N}$ has a convergent subsequence, thus $KS$ is compact.
If $B_3=B_1$ and $\{y_i\}_{i\in \N}$ is a bounded sequence in $B_2$ then we extract a convergent subsequence
$\{Ky_{i_k}\}_{k\in \N}.$ Since $S$ is a bounded linear operator, it is continuous
thus $SKy_{i_k}$ converges to an element in $B_1$ as $k\to \infty.$ Hence $SK$ is compact. 
\end{proof}
The following generalization which is a standard characterization of Fredholm operators.
\begin{theorem}\label{fredhm00}
A bounded linear operator between two Banach spaces $T:B_1\to B_2$ is
Fredholm iff there exists $S\in L(B_2,B_1)$ and compact operators (in fact they can be chosen as finite rank projections) $K_1\in L(B_1,B_1),K_2\in L(B_2,B_2)$ such that
$ST=I-K_1$ and $TS=I-K_2$.
\end{theorem}
\begin{proof}
To see $(\Leftarrow)$ note that dim$(\mbox{Ker}(T))\leq$dim$(\mbox{Ker}(I-K_1))<\infty$
and that dim($B_2/\mbox{Ran}(T)$)$\leq$dim($B_2/\mbox{Ran}(I-K_2)$)$<\infty.$
It suffices to show that $\mbox{Ran}T$ is a closed and complemented subspace of $B_2.$
Let $\{ y_j +R(T)\}_{j=1}^N$ be a basis for $B_2/\mbox{Ran}(T)$ and set $Z:=\mbox{span}\{y_j\}_j\subset B_2.$ 
Denote by $\hat{\oplus}$ the {\em extrinsic direct sum} defined for two Banach spaces $A,B$ as
$A\hat{\oplus} B:=\{(a,b)\in A\times B\colon a\in A, b\in B\}.$ Recall that the well-known 
bounded inverse theorem states that a bounded bijective linear operator has bounded inverse.
Indeed, if $P\in L(A,B)$ is a bounded linear operator between Banach spaces $A,B,$ with inverse $P^{-1}$ then
$P^{-1}$ is known to be a linear operator from $B$ to $A$ since for $a_1,a_2\in A,$ and scalars $c_1,c_2$,
$P^{-1}(c_1a_1+c_2a_2)=P^{-1}(P(c_1P^{-1}(a_1) +c_2P^{-1}(a_2)))=c_1P^{-1}(a_1) +c_2P^{-1}(a_2).$ 
Since $P$ is continuous it is closed. 
This implies that $P^{-1}$ is also closed. To see this let $b_j\to b$ in $B.$ and $P^{-1}(b_j)\to a$ in $A.$Letting
$a_j=P^{-1}(b_j)$ we have $a_j\to a$ in $A$ and $P(a_j)\to b$ in $B.$ Since $P$ is closed, $b=P(a)$ thus $a=P^{-1}(b).$
By the closed graph theorem $P^{-1}$ must be continuous, 
thus bounded.
In particular, if $A$ and $B$ are complementary closed subspaces of $A\oplus B$ then, by the
bounded inverse theorem, the map $A\hat{\oplus} B\to A\oplus B\colon (a,b)\mapsto a+b$ is a
topological isomorphism and there is a constant $C<\infty$ such that for all $a\in A,b\in B,$
$\norm{a+b}_A\leq \norm{a}_A +\norm{b}_B\leq C\norm{a+b}_A.$
Now the map $B_1/\mbox{Ker}(T) \hat{\oplus} Z \to B_2:$
$(x+\mbox{Ker}(T),z)\mapsto Tx+z$ is a bounded linear bijective map. 
So there exists a $C<\infty$
such that $\norm{x+\mbox{Ker}(T)}_{B_1} +\norm{z}_Z \leq C\norm{Tx +z}_{B_2}.$
Setting $z=0$ shows that $B_1/\mbox{Ker}(T)\to\mbox{Ran}(T)$ has a bounded inverse. Thus $\mbox{Ran}(T)$ is closed.
\\
For $(\Rightarrow)$ recall that $B_1=\mbox{Ker}(T)\oplus B_1'$ and $B_2=\mbox{Ran}(T)\oplus B_2'$
for closed subspaces $B_j'\subset B_j,$ $j=1,2.$ 
Note that $T_1:=T|_{B_1'}$ is an isomorphism from$B_1'$ into $\mbox{Ran}(T)$.
Recall that a linear map $P$ from a Banach space $A$ to itself is called a projection if $P^2=P.$
Then $I-P$ is also a projection and we have a splitting $A=\mbox{Ran}(P)+\mbox{Ker}(P)$, 
$\mbox{Ran}(P)\cap\mbox{Ker}(P)=\{0\}.$ Conversely if $A=A'+A''$ and $A'\cap A''=\{0\},$ for subspaces $A',A''$
then for each $a\in A$ there are unique $a'\in A',a''\in A''$ such that $a=a'+a'',$ so the linear map given by $a\mapsto a'$
is a projection, and we call it the projection onto $A'$ {\em along} $A''.$ 
Let $K_1$ be the projection onto Ker$(T)$ along $B_1'.$
and let $K_2$ be the projection onto $B_2'$ along Ran$(T)$. Let $\hat{T}:B_1'\to $Ran$(T)$ be the restricted isomorphism and set 
$S:=\hat{T}^{-1}(I-K_2).$ Then $K_1,K_2$ are finite rank projections and $ST=\hat{T}^{-1}(I-K_2)T=\hat{T}^{-1}T=I-K_1,$ and
$TS=T\hat{T}^{-1}(I-K_2)=I-K_2.$ This completes the proof.
\end{proof}
For a Banach space $B$ let $\mathcal{K}(B,B)$ denote the set of compact operators.
Then $\mathcal{K}(B,B)$ is a two-sided ideal in $L(B,B)$ thus
the quotient space $L(B,B)/\mathcal{K}(B,B)$ is a Banach algebra. Let $[T]\in L(B,B)/\mathcal{K}(B,B).$
Then the previous theorem implies that $T\in L(B,B)$ is Fredholm iff $[T]$ is invertible in $L(B,B)/\mathcal{K}(B,B)$.
We can now derive a result that can be applied to the case of a general complex equation as presented in the introduction (see e.g.\ Douglas \cite{douglas} for the case of Hilbert spaces).
\begin{corollary}
The sum $T+K$ of a Fredholm operator $T:B_1\to B_2$ and a compact operator $K:B_1\to B_2$ on Banach spaces $B_1,B_2,$ is Fredholm. 
\end{corollary}
\begin{proof}
According to Theorem \ref{fredhm00} there exists compact operators $K_1\in$\\ $L(B_1,B_1),K_2\in L(B_2,B_2)$ and $S\in L(B_2,B_1)$
such that $ST=I-K_1$ and $TS=I-K_2$. This implies that
$S(T+K)=I+(K_1+SK)$ and $(T+K)S=I+(K_2+KS).$ Since $K$ is compact, Lemma \ref{rudinthm1} implies that $FK$ and $KF$ are compact.
Clearly, the sum of two compact operators is compact (the sum of two totally bounded sets of a normed space is totally bounded).
Hence $(K_2+KS)$ and $(K_2+KS)$ are compact operators.
Another application of Theorem \ref{fredhm00} implies that $T+K$ is Fredholm.
\end{proof}
Now we mention that in PDE theory on sometimes uses the notion of a parametrix.
\begin{definition}[Parametrix]\index{Parametrix}
Let $L:X\to Y$ be a bounded linear operator between Banach spaces. An bounded operator $P:Y\to X$
is called a {\em left parametrix} for $L$ if 
\begin{equation}
PL=I_1+K_1
\end{equation}
where $I_1$ is the identity operator and $K_1$ is a compact operator.
$P$
is called a {\em right parametrix} for $L$ if 
\begin{equation}
LP=I_2+K_2
\end{equation}
where $I_2$ is the identity operator and $K_2$ is a compact operator.
A {\em two-sided paramterix} is one that is both a left and a right parametrix.
\end{definition}
A parametrix can be viewed as a sort of substitute for the inverse of a given operator that does not necessarily have
an inverse.
This is also true for the Fredholm property, since being Fredholm implies that the kernel and range both have direct complements
(any finite dimensional subspace of a Banach space is complemented and the quiotient with any finite dimensional subspace
of a Banach space is complemented) and
we have seen that if $L:X\to Y$ is a Fredholm operator 
then the map $X_1\to \mbox{Ran}$, where $X_1$ is the direct complement of $\mbox{Ker}X$, is one-to-one and 
continuous and as we have previously proved the inverse mapping is also continuous, i.e.\
$L$ yields a continuous isomorphism between $X_1$ and its range in $Y.$
We have previously proved that any Fredholm operator has a two sided parametrix.
Let us give an alternative proof of the Fredholm property characterization.
\begin{theorem}
If an operator $L:X\to Y$ between Banach spaces has a left parametrix, then $\mbox{Ker}(L)$ is finite dimensional.
If $L$ has a right parametrix, then $\mbox{Ran}(L)$ is closed and the cokernel is finite dimensional.
In particular, if $P$ is a two-sided parametrix for $L$ then $L$ is Fredholm (and $P$ is also Fredholm since $L$ is a paramtrix for $P$).
\end{theorem}
\begin{proof}
We have compact operators $K_1,K_2$ such that 
(because $P$ is a left parametrix)
\begin{equation}\mbox{Ker}(L)\subset \mbox{Ker}(LP)=\mbox{Ker}(I_1 +K_1) \end{equation}
and (because $P$ is a right parametrix)
\begin{equation}
\mbox{Ran}(L)\supset \mbox{Ran}(LP)=\mbox{Ran}(I_2 +K_2) 
\end{equation}
and
\begin{equation}
Y=\mbox{Ran}=\mbox{Ran}(I_2+K_2)+Q
\end{equation}
where $Q$ is a finite dimensional subspace of $Y$ thus
$Y=\mbox{Ran}(L)+Y_1$ where the cokernel of $L$, $Y_1,$ is a finite dimensional subspace of $Q.$
\end{proof}

We can in particular cases obtain the following very useful estimates 
\begin{theorem}
If an operator $L:X\to Y$ between Banach spaces is Fredholm
such that $X$ is continuously embedded in a Banach space $\tilde{X}.$
Then there exists a constant $C>0$ such that for all $x\in X,$
\begin{equation}
\norm{x}_X\leq C(\norm{Lx}_Y +\norm{x}_{\tilde{X}})
\end{equation}
If $\mbox{Ker}L=\{0\}$ then the last term in the parenthesis can be dropped. 
\end{theorem}
\begin{proof}
On the direct complement $X_1$ of $\mbox{Ker}L$, the estimate holds true without the last term in the parenthesis
since $L:X_1\to \mbox{Ran}L$ is invertible. If $\mbox{Ker}L\neq \{0\}$ then $\norm{x}_X\leq C_1\norm{x}_{\tilde{X}}$ on $\mbox{Ker}L.$
Decompose $x=x_1+x_2,$ $x_1=Px$ where $P$ is the projection on $\mbox{Ker}L$, and $x_2=(I-P)x$ (where $(I-P)$ is the projection on $X_1.$
Then $\norm{x}_X\leq\norm{x_1}_X +\norm{x_2}_{X}$ and for a constant $C_1>0,$
\begin{equation}
\norm{x_1}_X\leq C_1\norm{x_1}_{\tilde{X}}\leq C_1(\norm{x_1}_{\tilde{X}}+\norm{x_2}_{\tilde{X}})
\end{equation}
Since the embedding of $X$ in $\tilde{X}$ is continuous we have $\norm{x_1}_{\tilde{X}}\leq C_2\norm{x_2}_X,$ for a constant $C_2>0.$
Hence for a constant $C_3>0,$
\begin{equation}
\norm{x_2}_X=\norm{x_2}_{X_1}\leq C_3\norm{Lx_2}_{\mbox{Ran}(L)}= C_3\norm{Ax_2}_{Y}=C_3\norm{Ax}_{Y}
\end{equation}
This completes the proof.
\end{proof}

Here is one way Fredholm theory can be applied to the theory of elliptic equations.
\begin{theorem}[See e.g.\ Agranovich \cite{agranovich}, Thm 6.2.1]\label{agranovthm}
Let $L$ be an $m$:th order partial differential operator with $C^\infty$-smooth coefficients
on a closed (compact) smooth manifold $M.$ Then the following are equivalent:
\begin{itemize}
\item[(1)] $L$ is elliptic.
\item[(ii)] $L$ defines for each $s,$ a Fredholm operator $H^s(M)\to H^{s-m}(M).$
\item[(iii)] $L$ has a two-sided parametrix $P$ acting boundedly from $H^{s-m}(M)$ to $H^s(M)$ such that
$K_1=PL-I$ and $K_2=LP-I$ are bounded operators $H^s(M)$ to $H^{s+1}(M)$ for any $s.$
\item[(iv)] There is a constant $C_s$ such that
\begin{equation}
\norm{f}_{H^s(M)}\leq C_s(\norm{Lf}_{H^{s-m}(M)} +\norm{f}_{H^{s-1}(M)})
\end{equation}
\end{itemize}
\end{theorem}
One can give an alternative proof of the elliptic regularity theorem (see the Appendix, Section \ref{ellipticapp}) using Fredholm theory. 

\section{Bergman spaces} 
A Hilbert space $H$ of functions on a set $\Omega$, is called a {\em reproducing kernel Hilbert space}
\index{Reproducing kernel} if for each $z_0\in \Omega$, there exists a
positive number $C_{z_0},$ such that $\abs{f(z_0)}\leq C_{z_0} \norm{f}_H,$ for all $f\in H.$
By identifying $f(z_0)$ with the evaluation functional $f\mapsto f(z_0)$ 
it follows from the Riesz representation theorem
that for all $z_0\in \Omega$, there exists a unique $K_{z_0}\in H$
such that $f(z_0)=\langle f,K_{z_0} \rangle.$
The {\em reproducing kernel} of $H$ is defined by 
\begin{equation}
K\colon \Omega\times\Omega\ni (x,y)\mapsto 
\langle K_{x},K_{y} \rangle
\end{equation}
Given a domain $\Omega\subset\C$
and a positive integer $p>0,$
the {\em Bergman space}, $A_1^p(\Omega)$, is defined by the set of holomorphic $L^p(\Omega)$ functions,
i.e.\ $f\in A_1^p(\Omega)$ iff $f\in \mathscr{O}(\Omega)$ and having finite norm according to
\begin{equation}
\norm{f}_{A_1^p(\Omega)}=\langle f,f \rangle:=\left(\frac{1}{\pi}\int_{\Omega} \abs{f(x+iy)}^p  dx\wedge dy\right)^{\frac{1}{p}}<\infty
\end{equation}
where we interpret $dxdy =dx\wedge dy$ as the standard (Lebesgue) area measure on $\C.$ 
For $p=2$ these spaces are Banach spaces with respect to the norm $\norm{\cdot }_{A_1^2(\Omega)}$
as a consequence of the fact that uniform limits of holomorphic functions
on a compact are holomorphic in the interior. 
If $f(z)$ is holomorphic with expansion $f(z)=\sum_{j=0}^\infty c_j z^j$
on the disc then we have
\begin{equation}
\int_{\{\abs{z}<1\}} \abs{f(z)}^2 \frac{d\mu(z)}{\pi} =\sum_{j=0}^\infty \frac{\abs{c_j}^2}{j+1}
\end{equation}
Thus $f$ belongs to the Bergman space iff the right hand side is finite.
For $g=\sum_j d_j z^z,$ the scalar product of $f$ and $g$ can be defined by the formula 
$\langle f,g\rangle=\sum_{j=0}^\infty \frac{c_j \bar{d}_j}{j+1}$ and the set
$\{z^j\sqrt{j+1}\}_{j\in \N}$ is an orthonormal basis for the Bergman space. 
The polynomials will obviously be dense in 
$A_1^2(\{\abs{z}<1\}).$ 
\begin{proposition}\label{evalenglis}
Let $\Omega\subset\Cn$ be a bounded domain. For each $q\in \Omega,$, the linear functional $A^2(\Omega)\ni f\mapsto
f(q)$ is bounded. As a consequence $f(q)=\langle f,g_q\rangle$ for some $g_q\in A_1^2(\Omega)$ and furthermore
\begin{equation}
\norm{g_q}\leq \left(\frac{1}{\omega_n (\mbox{\em dist}(q,\Cn\setminus \Omega)^{2n})}\right)^{\frac{1}{2}}
\end{equation}
where $\omega_n$ denotes the Lebesgue volume of the unit ball in $\Cn.$
\end{proposition}
\begin{proof}
Define
\begin{equation}
F_q(z):=
\left\{
\begin{array}{ll}
0 & , \abs{z-q}\geq \mbox{dist}(q,\Cn\setminus \Omega)\\
(\omega_n (\mbox{dist}(q,\Cn\setminus \Omega)^{2n}))^{-1} & ,\abs{z-q} <\mbox{dist}(q,\Cn\setminus \Omega)
\end{array}
\right.
\end{equation}
Then
\begin{multline}
\int_\Omega \abs{F_q(z)}^2 d\mu(z)=\omega_n (\mbox{dist}(q,\Cn\setminus \Omega)^{2n}\left(
\frac{1}{\omega_n (\mbox{dist}(q,\Cn\setminus \Omega))^{2n}}
\right)^2\\
=\frac{1}{\omega_n (\mbox{dist}(q,\Cn\setminus \Omega))^{2n}}<\infty
\end{multline}
thus $F_q\in L^2(\Omega).$ 
Denote by $P$ the orthogonal projection of $L^2(\Omega)$ onto $A_1^2(\Omega).$
Then for any $f\in A_1^2(\Omega)$
\begin{multline}
\langle f,P( F_q)\rangle_{A_1^2(\Omega)}= \langle f,P(F_q)\rangle_{L^2(\Omega)}=\\
\frac{1}{\omega_n (\mbox{dist}(q,\Cn\setminus \Omega))^{2n}}
\int_{\{ \abs{z-q}<\mbox{dist}(q,\Cn\setminus \Omega)\}} f(z)
\end{multline}
which by the mean value theorem is equal to $f(z).$ Hence choosing $g_q=P( F_\lambda)$
we have
\begin{equation}
\norm{g_q}_2^2 \leq \norm{F_q}_2^2 =\frac{1}{\omega_n (\mbox{dist}(q,\Cn\setminus \Omega))^{2n}}
\end{equation}
This completes the proof.
\end{proof}
For the unit disc the kernel function can be written 
\begin{equation}
g_\zeta(z)=\frac{1}{(1-\bar{\zeta}z)^2} =\sum_{j=0}^\infty (j+1)\bar{\zeta}^j z^j
\end{equation}
Proposition \ref{evalenglis} can be used to realize that point evaluations are bounded thus plays a fundamental role in the analysis
of reproducing kernel spaces (in particular Bergman spaces). For the $q$-analytic case we shall later prove an analogous result.

\section{Polyanalytic Bergman spaces}
\subsection{Introduction}
Let $\Omega\subset\C$ be a domain. For a continuous function $\phi\colon \Omega\to \R_+$ we define the {\em weighted $q$-polyanalytic Bergman space} $A^2_{q,\phi}(\Omega)$
\index{Polyanalytic Bergman space} 
to be the set of $q$-analytic functions belonging to the weighted $L^2$ space with norm
\begin{equation}
\norm{f}_{A_{\phi}^2(\Omega)}:=\left(\pi^{-1}\int_{\Omega} \abs{f(z)}^2\phi(z) d\mu(z)\right)^{\frac{1}{2}}
\end{equation}
where $z=x+iy$ denotes the complex coordinate and $\mu(z)$ the standard area measure on $\C$.
We mention that when $\Omega=\C$ and $\phi(z)$ is chosen as $\exp(-\abs{z}^2)$ we have an inner product space
with inner product 
\begin{equation}
\langle f,g \rangle:=\frac{1}{\pi}\int_{\C} f\bar{g}(z)(z)\exp(-\abs{z}^2)d\mu(z)
\end{equation}
called the {\em Segal-Bargmann space}. 
\begin{proposition}\label{closednessthm}
$A^2_{q,\phi}(\Omega)$ is a closed subspace 
of the weighted space $L^2_\phi(\Omega),$
with norm $\norm{f}=\left(\frac{1}{\pi}\int_{\Omega} \abs{f(z)}^2\phi(z) d\mu(z)
\right)^{\frac{1}{2}}.$
\end{proposition}
\begin{proof}
The fundamental solution to 
$\partial_{\bar{z}}$ is $\frac{1}{\pi z}$
and by Leibniz rule for distributions $\frac{\bar{z}^{q-1}}{z\pi (q-1)!}$
is a fundamental solution for $\partial_{\bar{z}}^q.$
For $f\in A^2_q(\Omega)$ and a test function $\xi\in C_c^\infty(\Omega),$ such that $\xi=1$
on a neighborhood of $z$,
we have
\begin{multline}
\partial_{\bar{z}}^j f(z)=\frac{(-1)^{q-1}}{(q-1)!}\int_\Omega \partial_{\bar{\zeta}}^q
(\partial_{\bar{\zeta}}^j f(\zeta)\xi(\zeta))\frac{(\bar{\zeta}-\bar{z})^{q-1}}{\zeta -z}d\mu(\zeta)
=\\
\frac{(-1)^{q-1}}{(q-1)!}\int_\Omega \sum_{k=0}^{q-j-1}\partial_{\bar{\zeta}}^{k+j}f(\zeta)
\partial_{\bar{\zeta}}^{q-k} \xi(\zeta) \frac{(\bar{\zeta}-\bar{z})^{q-1}}{\zeta -z}d\mu(\zeta)
\\=\int_{\mbox{supp}\partial_{\bar{\zeta}}\xi} f(\zeta)U_{j,\xi}(\zeta,z)d\mu(\zeta)
\end{multline}
where 
\begin{equation}
U_{j,\xi}(\zeta,z):=\frac{1}{(q-1)!} \sum_{k=0}^{q-j-1}\partial_{\bar{\zeta}}^{k+j}
\left(
\partial_{\bar{\zeta}}^{q-k}\xi(\zeta)\frac{(\bar{\zeta}-\bar{z})^{q-1}}{\zeta -z}
\right)
\end{equation}
and $U_{j,\xi}\in C^\infty_c(\Omega).$
Hence
\begin{equation}
\abs{\partial_{\bar{z}}^j f(z)} \leq \norm{ f}_{A^2_{\phi}(\Omega)}\norm{ U_{j,\xi}}_{A^2_{\phi}(\Omega)}
\end{equation}
since $\phi$ was arbitrary these bounds are locally uniform in $z$.
This completes the proof.
\end{proof}
An alternative proof of closedness of $A_q^2$ in $L^2$ can, in the case of the unit disc be given as follows.
Let $\varphi_z(\zeta):=(z-\zeta)/(1-\zeta\bar{z}),$ and $E(z,r):=\{\zeta\colon \abs{\varphi_z(\zeta)}<r\},$
$0<r<1.$
\begin{proposition}
\begin{equation}
E(z,r)=\{ \abs{z-c}<R\},\quad c=\frac{z(1-r^2)}{1-r^2\abs{z}^2},\quad  R=\frac{r(1-\abs{z}^2)}{1-r^2\abs{z}^2}
\end{equation}
\end{proposition}
\begin{proof}
For $\zeta,z\in D,$ 
\begin{multline}
\abs{\frac{z-\zeta}{1-z\bar{\zeta}}}<r \Leftrightarrow \abs{z-\zeta}^2<r^2\abs{1-z\bar{\zeta}}^2\Leftrightarrow\\
\abs{z}^2 +\abs{\zeta}^2-2\re(z\zeta)<r^2(1+\abs{z}\abs{\zeta}^2 -2\re(z\bar{\zeta}))\Leftrightarrow\\
(1-r^2\abs{z}^2)\abs{\zeta}^2 -2\frac{1-r^2}{1-\abs{z}^2r^2}\re(z\bar{\zeta})<\frac{r^2-\abs{z}^2}{1-\abs{z}^2r^2}
\Leftrightarrow\\
\abs{\zeta -\frac{1-r^2}{1-\abs{z}^2r^2}z}^2 <\frac{r^2-\abs{z}^2}{1-\abs{z}^2r^2} +\frac{(1-r^2)^2\abs{z}^2}{(1-\abs{z}^2r^2)^2}
\Leftrightarrow \\
\abs{\zeta -\frac{1-r^2}{1-\abs{z}^2r^2}z}^2 <\frac{(1-\abs{z}^2)^2r^2}{(1-\abs{z}^2r^2)^2} \Leftrightarrow \\
\abs{\zeta -\frac{1-r^2}{1-\abs{z}^2r^2}z}^2 <\frac{1-\abs{z}^2}{1-\abs{z}^2r^2}r
\end{multline}
This completes the proof.
\end{proof}
\begin{proposition}[Pavlovi\'c \cite{pavlovic}]\label{pavlovic}
There exists $C>0$ such that for any $z_0\in D$, $s\in (0,1-\abs{z_0}),$ and any $f$ satisfying $\Delta^q f=0$ (in particular
this includes all $f\in A^2_q(D)$) 
\begin{equation}
\abs{f(z_0)}^2\leq \frac{C}{t^2}\int_{\abs{\zeta -z_0}<s}\abs{f(\zeta)}^2 d\mu(\zeta)
\end{equation}
\end{proposition}
Set $R'=R-\abs{c-z}.$ For $s=R'$ in Proposition \ref{pavlovic} and using $\{ \zeta:\abs{\zeta -z}<R'\}\subset E(z,t),$
we have
\begin{multline}
\abs{h(z)}^2 \leq \frac{C(1+t\abs{z})^2}{t^2(1-\abs{z}^2)^2}\int_{\abs{\zeta-z}<R'} \abs{f(\zeta)}^2d\mu(\zeta)\leq \\
 \frac{C(1+t^2)^2}{t^2(1-\abs{z}^2)^2}\int_{E(z,t)} \abs{f(\zeta)}^2d\mu(\zeta)=
\frac{C'}{(1-\abs{z}^2)^2}\int_{E(z,t)} \abs{f(\zeta)}^2d\mu(\zeta)
\end{multline}
where $C'$ is independent of $f$ and $z.$ This proves that $A^2_q$ is closed in $L^2(D).$
\\
\\
As a consequence $A^2_{q,\phi}(\Omega)$ is a Hilbert space. Since the point evaluations $z\mapsto 
\overline{\partial}_z^j f(p_0)$ are bounded, the Riesz representation theorem implies that for each 
$p_0\in \Omega$ there exists $K_{q,\phi}(\zeta,p_0)\in A^2_{q,\phi}(\Omega)$
such that $f(z)=\langle f, K_{q,\phi}(\zeta,p_0)\rangle,$ i.e.\
$A^2_{q,\phi}(\Omega)$ is a reproducing kernel space.
Analysis of polyanalytic versions of Bergman spaces appeared already in the work of
Ko\u{s}elev \cite{koselev} (1977) and has recently seen increasing interest and activity,
see e.g.\ Vasilevski \cite{vasilevski}, \cite{vasilevski2}, \cite{vasilevskiquaternion}, \cite{vasilevskibok} (and the references therein), Wolf \cite{wolf},  Ramazanov \cite{ramazanov1}, \u{C}u\u{c}kovi\'c \& Le \cite{cuckovic} and Karlovich \cite{karlovich}.
Let $\phi:\Omega\to \R_+$ be a continuous weight function and consider the Bergman space $A^2_{q,\phi}(\Omega).$
Let $\{\psi_j\}_{j\in \Z_+}$ be an orthonormal basis for $A^2_{q,\phi}(\Omega).$
Then $K_{q,\Omega,\phi}(z,\zeta):=\sum_{j=1}^{\infty} \psi_j(z)\overline{\psi(\zeta)}, (z,\zeta)\in \Omega\times \Omega,$
satisfies $f(\zeta)=\langle f,K_{q,\Omega,\phi}(\cdot,\zeta)\rangle_{L^2_\phi(\Omega)},\zeta\in \Omega$ for all $f\in 
A^2_{q,\phi}(\Omega).$ Thus it makes sense to call the $K_{q,\Omega,\phi}$ Bergman kernels\index{Bergman kernel}.
Let us look at the basic case $q=1$ closer. For a domain, if $K_\Omega(z,\zeta)$ is the Bergman kernel function 
satisfying $f(z)=\int_\Omega K_\Omega(z,\zeta)f(\zeta)d\mu(\zeta)$ for $f\in A^2_1(\Omega)$
then
the function $k_\zeta(z):=\overline{K_\Omega(\zeta,z)}$ belongs to $A^2_1(\Omega)$
and
\begin{multline}
\overline{K_\Omega(\zeta,z)}=\int_\Omega K_\Omega(z,\tau)\overline{K_\Omega(\zeta,\tau)}d\mu(\tau)=\\
\overline{\int_\Omega K_\Omega(z,\tau)\overline{K_\Omega(\zeta,\tau)}d\mu(\tau)}=
\overline{\overline{K_\Omega(z,\zeta)}}=K_\Omega(z,\zeta)
\end{multline}
If $\{\psi_j\}_{j\in\in \Z_+}$ is an orthonormal basis in $A^2_1(\Omega)$ then
$k_\zeta(z)=\overline{K_\Omega(\zeta,z)}=K_\Omega(z,\zeta)$ is analytic and has an expansion
\begin{equation}
K_\Omega(z,\zeta)=\sum_{j=1}^\infty c_j(\zeta)\psi_j(z)
\end{equation}
where the coefficients are 
\begin{multline}
c_j(\zeta)=\langle k_\zeta(z),\psi_j(z)\rangle =\int K_\Omega(z,\zeta)\overline{\psi_j(z)}d\mu(z)=\\
\overline{\int_\Omega K_\Omega(\zeta,z)\psi_j(z)d\mu(z)}=\overline{\psi_j(\zeta)}
\end{multline}
which yields 
\begin{equation}
\int_\Omega K_\Omega(z,\zeta)=\sum_{j=1}^\infty \psi_j(z)\overline{\psi_j(\zeta)}
\end{equation}
For the unit disc $\{\psi_j\}_{j\in \Z_+}$, $\psi_j:=\sqrt{\frac{j}{\pi}} z^{j-1}$, forms
an orthonormal basis of $A^2_1(\{\abs{z}<1\})$ thus we obtain the well-known formula
\begin{equation}
\int_\Omega K_\Omega(z,\zeta)=\sum_{j=1}^\infty \frac{j}{\pi} (z\bar{\zeta})^{j-1}=\frac{1}{\pi(1-z\bar{\zeta})^2}
	\end{equation}
	
	\begin{theorem}
		If $\omega=\phi(z)$ is a biholomorphic map of the domain $D,$ onto a domain $\Omega$ then
		\begin{equation}
		K_{D}(z,\zeta)=
		K_{\Omega}(w,\omega)\phi'(z)\overline{\phi'(\zeta)}
		\end{equation}
	\end{theorem}
	\begin{proof}
		Let $\{g_j(w)\}_{j\in \Z_+}$ be an orthonormal basis for $A^2_1(\Omega).$ Then
		\begin{equation}
		\delta_{j,k}-\int_\Omega g_j(w)\overline{g_k(w)}d\mu(w)=
		\int_D g_j(\phi(z))\overline{g_k(\phi(z))} \phi'(z)\overline{\phi'(z)} d\mu(z)
		\end{equation}
		which implies that $\psi_j:=g_j(\phi(z))\phi'(z)$, $j\in \Z_+$ defines an orthonormal basis for $A^1_1(D).$ Thus
		\begin{multline}
		\int_\Omega K_{D}(z,\zeta)=\sum_{j=1}^\infty \psi_j(z)\overline{\psi_j(\zeta)}=
		\sum_{j=1}^\infty g_j(\phi(z))\phi'(z)\overline{g_j(\phi(\zeta))\phi'(\zeta)}=\\
		\sum_{j=1}^\infty g_j(\phi(z))\overline{g_j(\phi(\zeta))}\phi'(z)\overline{\phi'(\zeta)}=
		\int_\Omega K_{G}(w,\omega)\phi'(z)\overline{\phi'(\zeta)}
		\end{multline}
		This completes the proof.
	\end{proof}
	Obviously for the upper half-plane, $H^+=\{\im z >1\}$, the a biholomorphism onto the unit disc 
	$D:=\{\abs{z}<1\}$,
	is given by M\"obius transformation
	\begin{equation}
	z=\frac{w-i}{1-iw}
	\end{equation}
	Thus the Bergman kernel is given by
	\begin{equation}\label{orthokernelekvenbra}
	K_{H^+}(w,\omega)= K_{D}(z,\zeta)=\\
	\frac{1}{\pi\left(1-\frac{w-i}{1-iw}\frac{\bar{\omega}+i}{1+i\bar{\omega}} \right)} \frac{2}{(1-iw)^2}\frac{2}{(1-i\bar{\omega})^2}
	=\\
	\frac{1}{\pi(w-\bar{\omega})^2}
	\end{equation}
	As, for a domain $\Omega,$ $A^2_1(\Omega)$ is a closed subspace of the Hilbert space $L^(\Omega)$ there exists an orthogonal projection of $L^2(\Omega)$ onto
	$A^2_1(\Omega).$ 
	If we denote this projection by
	$B_\Omega$, and call it the {\em Bergman orthogonal projection}\index{Bergman orthogonal projection} then we can verify that
	for $f\in L^2(\Omega)$
	\begin{multline}
	(B_\Omega f)(z)=\langle B_\Omega ,f \rangle - \langle f, B_\Omega k_z \rangle -\langle f,k_z\rangle=\\
	\int_\Omega f(\zeta)\overline{k_z(\zeta)}d\mu(\zeta)=\int_\Omega K_\Omega (z,\zeta)f(\zeta)d\mu(\zeta)
	\end{multline}
	It is further known (although we shall not need this) that if $g(z,\zeta)$ it the Green function of a domain $\Omega$ then
	\begin{equation}
	K_\Omega(z,\zeta)=-\frac{2}{\pi}\Delta g(z,\zeta),\quad \Delta=4\partial_z\partial_{\bar{z}}
\end{equation}
 In the case $\Omega=\{\abs{z}<1\}$ 
Ko\u{s}elev \cite{koselev} 
showed that
\begin{equation}\label{tututytu}
K_q(z,\zeta)
=\sum_{j=1}^{q-1}\sum_{m=0}^{\infty} e_{mj}(z)\overline{e_{mj}(\zeta)}
\end{equation}
\begin{equation}\label{koselevbasekv}
e_{mj}(z):=\sqrt{\frac{m+j+1}{\pi}}\frac{1}{(m+j)!}\frac{\partial^{m+j}}{\partial z^j\partial \bar{z}^m}
(\abs{z}^2-1)^{m+j}
\end{equation}
where the $\{e_{mj}(z)\}_{m=0}^\infty,$ $j=0,\ldots,q-1,$ forms a complete orthonormal system in $A^2_q(D),$
for the unit disc $D.$
In the sums that follow we use for convenience the condition
\begin{equation}
\binom{n}{-m}=0,\quad n\in \Z ,m=\in \Z_+
\end{equation}
\begin{equation}
\binom{n}{n+m}=0,\quad n\in \N ,m=\in \Z_+
\end{equation}
Let us denote for the unit disc $D$, $K_{q,D}(z,\zeta)=:K_q(z,\zeta).$
%
	Since $\partial_{\bar{z}}(z\bar{z}-1)=z$ we have
	for $k=1$ 
	\begin{equation}
	\partial_{\bar{z}}^{1+j}(\abs{z}^2-1)^{1+j}=
	z(1+j)(\abs{z}-1)^{j}
	\end{equation}
	and by induction in $k$ we deduce
	\begin{equation}
	\partial_{\bar{z}}^{k}(\abs{z}^2-1)^{k+j}=\frac{(k+j)!}{j!}
	z^k(\abs{z}^2-1)^{j}
	\end{equation}
	This yields
	\begin{multline}
	\frac{\partial^{k+j}}{\partial z_j\partial \bar{z}_k}
	(\abs{z}^2-1)^{k+j}=
	\partial_{z}^{j}(z^k(\abs{z}^2-1)^{j})=\partial^j_z \left(\sum_{i} (-1)^i\binom{j}{i}z^k(z\bar{z})^{j-i}\right)=\\
	\sum_{i} (-1)^j\binom{j}{i}\partial^j_z \left(z^{k+j+i}\bar{z}^{j-i}\right)=\sum_{i} (-1)^j\binom{j}{i}\frac{(k+i)!}{j!} z^{k+i}\bar{z}^{j-i}
	\end{multline}
%
%
This can be used togheter with Eqn.(\ref{koselevbasekv}) and Eqn.(\ref{tututytu}) to deduce in the case of the unit disc (see Ko\u{s}elev \cite{koselev})
\begin{equation}\label{koseleveq}
\begin{split}
K_{q}(z,\zeta) & = \sum_{j=0}^{q-1}\sum_{k=0}^\infty e_{kj}(z)\overline{e_{kj}(\zeta)}\\
 & =q\sum_{j=0}^{q-1}(-1)^j\binom{q}{j+1}\binom{q+j}{q}\frac{(1-z\bar{\zeta})^{q-j-1}\abs{z-\zeta}^{2j}}{(1-z\bar{\zeta})^{q+j+1}}\\
& =\frac{q}{(1-\bar{\zeta}z)^{2q}}\sum_{j=0}^{q-1}(-1)^j\binom{q}{j+1}\binom{q+j}{q}\abs{1-z\bar{\zeta}}^{2(q-j-1)}
\abs{z-\zeta}^{2j}
\end{split}
\end{equation}
and 
\begin{equation}
K_{q}(z,z)=\frac{q^2}{\pi(1-\abs{z}^2)^2}
\end{equation}
Note that $K_q(z,z)=\sum_{j=1}^\infty\abs{\psi(z)}^2,$ hence $\Delta\log K_q(z,z)\geq 0$ i.e.\ $K_q(z,z)$ is subharmonic.
Setting $C_{l,j}:=l(l-1)\cdots (l-1+1)/j!$, the property
\begin{equation}
\frac{1}{m}C_{q,1}C_{q,q}+\frac{1}{m+1}C_{q,2}C_{q+1,q}\cdots+ +\frac{1}{m+q-1}C_{q,q}C_{2q-1,q}=\left\{
\begin{array}{ll}
\frac{1}{m} &,m=1\\
0 &, m=2,3,\ldots
\end{array}
\right.
\end{equation}
can be used to verify that 
\begin{equation}
(\abs{z}^{2(m-1)},K_q(z,0))=\left\{
\begin{array}{ll}
1 & m=1\\
0 &, m=2,3,\ldots,q
\end{array}
\right.
\end{equation}
which in turn together with the standard representation pf $q$-analytic functions, and orthogonality of $z^l,z^j$ for $i\neq j$, can be used to very $(f(z),K_q(z,0))=f(0).$
\\
For the half-space $H^+=\{\im z>0\},$ the
reproducing kernel of the Hilbert spaces $A^2_{q,(\lambda +1)(2y)^\lambda 
}(H^+),$ $z=x+iy,$ $\lambda\in (-1,\infty),$ have been
calculated using Fourier transforms (see Ramírez-Ortega \cite{ortega}). In particular, the so-called {\em poly-Bergman orthogonal projection} 
\begin{equation}
\mathcal{B}_{q,\lambda} :L^2(H^+,(\lambda +1)(2y)^\lambda)\to 
A^2_{q,(\lambda +1)(2y)^\lambda}(H^+)\end{equation}
takes the form
\begin{equation}
\mathcal{B}_{q,\lambda}f(z):=\int_{H^+} f(\zeta) K_{q,\lambda}(z,\zeta)(\lambda +1)(2y)^\lambda d\mu(z)
\end{equation}
for
\begin{multline}
K_{q,\lambda}(z,\zeta):=((z-\bar{z})-(\zeta-\bar{\zeta}))^{-1}(i\bar{\zeta}-iz)^{-\lambda-2}\times\\
\sum_{j,k=0}^{q-1} \alpha_{jkq}\left(
\frac{z-\bar{z}}{\lambda +j+1} -\frac{\zeta-\bar{\zeta}}{\lambda +k+1}
\right)\frac{(z-\bar{z})^j(\zeta-\bar{\zeta})^{k}}{(z-\bar{\zeta})^{j+k}}
\end{multline}
and 
\begin{equation}
\alpha_{jkq}:=\frac{(-1)^{j+k}\Gamma(q)\Gamma(\lambda +q)\Gamma(\lambda +j+k+2)(q+\lambda)}{\pi(\lambda+1)\Gamma(\lambda +j+1)j!(q-1-j)!\Gamma(\lambda +k+1)k!(q-1-k)!}
\end{equation}
 Recall that the {\em orthogonal difference}\index{Orthogonal difference} of two closed linear subspaces spaces $X,Y$ of a Hilbert space 
is given by $X\ominus Y:=X\cap Y^{\perp}.$
Reproducing kernels for the so called {\em true weighted poly-Bergman spaces}
\begin{equation}
A^2_{(q),(\lambda +1)(2y)^\lambda }(H^+):=A^2_{q,(\lambda +1)(2y)^\lambda }(H^+) \ominus A^2_{q-1,(\lambda +1)(2y)^\lambda }(H^+)
\end{equation}
where $\ominus$ denotes the orthogonal difference, have also been calculated, see Ramírez-Ortega \cite{ortega2}.
The weights $\phi$ of a weighted Bergman space can be viewed as the association of 
Riemannian metric to $\Omega,$ with respect to a special the two dimensional\index{Riemannian metric}
volume form,
see
Bergman \cite{bergman0}. 
For example
in the case $\Omega=\{\abs{z}<1\}$ 
Bergman considered the (hyperbolic) metric $ds^2:=\frac{1}{2}\Delta\ln K_\Omega(z,\bar{z})=\frac{1}{\pi}(1-\abs{z}^2)^{-2}(dx^2+dy^2)$ (where $\Delta=4\partial_z\partial_{\bar{z}}$ is the Laplace operator).
Recall that a metric, $g$ on $\Omega$, gives the
length of tangent vectors $v\in T_p\C$, $p\in \Omega$ according to
$\abs{v}^2_{g}=g\abs{v}^2,$ thus the length of a continuous curve $\gamma:[0,1]\to \C$, piecewise $C^1$ on $(0,1)$, is given by 
\begin{equation}
\tilde{\ell}(\gamma)=\int_{0}^1\abs{\frac{\partial \gamma(t)}{\partial t}}_{g} dt
\end{equation}
and the distance between two points $p,q$ is defined by \begin{equation}
\inf\{\tilde{\ell}(\gamma)\colon = \gamma\in C^0([0,1]) \mbox{ and piecewise }C^1, \gamma(0)=p,\gamma(1)=q\}
\end{equation}
The length element of the Bergman metric, where one uses $ds^2=\frac{1}{2}\Delta K_\Omega(z,\bar{z})(dx^2+dy^2)$
is given by $d\ell =\left(\frac{1}{2}\Delta \ln K_\Omega (z,\bar{z})\right)^{\frac{1}{2}}dl$ where $dl$ is the Euclidean length element. The
Bergman distance $\delta_\Omega(z_1,z_2),$ between two points $z_1,z_2$, for the case of the unit disc $\Omega=\{\abs{z}<1\}$ is
\begin{equation}
\delta_\Omega(z_1,z_2)=\frac{1}{2}\ln\left(\frac{\abs{1-z_1\bar{z}_2}+\abs{z_1-z_2}}{{1-z_1\bar{z}_2}-\abs{z_1-z_2}}\right)
\end{equation}
It is well-known that the metric induced by the volume form $ds^2$ above is invariant under biholomorphic mappings, see e.g.\ Krantz \cite{krantzflerdim}.
This means that if $\phi(z)$ is a biholomorphism of the unit disc onto a domain $U$ then
$\delta_\Omega(z_1,z_2)=\delta_U(z_1,z_2).$

\subsection{Vasilevski's method of projections}
As before denote the upper half-space by $H^+:=\{z\in \C:\im z>0\}.$
 Denote by $\chi_+(u)$ the characteristic function of the upper half line.
Consider the Bergman projection
\begin{equation}\label{9103}
(B_{H^+} f)(z)=\frac{1}{\pi} \int_{H^+} \frac{f(\zeta)}{(z-\bar{\zeta})^2}d\mu(\zeta),\quad L^2(H^+)\to A^2_1(H^+)
\end{equation}
which (see the Eqn.(\ref{orthokernelekvenbra}) and for the more general case we shall later prove Theorem \ref{ortokerneltheorem}) is the orthogonal projection of $L^2(H^+)$ onto the Bergman space $A^2_1(H^+)$.
Denote the kernel of this operator by 
\begin{equation}
K_{H^+}(z,\zeta)=-\frac{1}{\pi}\frac{1}{(z-\bar{\zeta})^2}=\frac{1}{\pi}\frac{1}{(i\bar{\zeta}+\eta -i\xi)^2}
\end{equation}
where we denote $\xi=\re z,$ which has Fourier transform (see Bateman \cite{bateman}) with respect to $\xi$ given by
\begin{equation}
\sqrt{\frac{2}{\pi}}\chi_+(\xi)\xi\cdot \exp((\eta+i\bar{\zeta})\xi)
\end{equation}
The space $A^2_{(q)}(H^+)$ of {\em true-$q$-analytic functions on $H^+$} \index{True-$q$-analytic functions on the upper half-plane}
is defined as 
\begin{equation}
A^2_{(q)}(H^+)=A^2_q(H^+)\ominus A^2_{q-1}(H^+), \quad A^2_{(1)}(H^+)=A^2_1(H^+)
\end{equation}
where $A^2_q(H^+)$ is the unweighted poly-Bergman space.
Vasilevski defines the unitary operator
$U_1=F\otimes I: L^2(H^+)=L^2(\R)\otimes L^2(\R_+)\to L^2(\R)\otimes L^2(\R_+),$
where $F$ is the Fourier transform $\hat{f}(x):=Ff(x)=\frac{1}{2\pi}\int_\R \exp(-ix\cdot\xi)f(\xi)d\xi$.
Note that the image $A_1^2=U_1(A^2_1(H^+))$ is the closed subspace of $L^2(H^+)$ consisting of all functions $f$ satisfying
\begin{equation}
(F\otimes I)2\partial_{\bar{z}} (F^{-1}\otimes I)f=i(x+\partial_y)f=0
\end{equation}
which has general solution $f(x,y)=\phi(x)\exp(-xy).$ 
The orthogonal projection $B_1:L^2(H^+)\to A_1^2$ takes the form
$B_1=(F\otimes I) B_{H^+} (F^{-1}\otimes I).$
We have
\begin{multline}
(F^{-1}\otimes I)B_1 f=\langle f,(F\otimes I)K_{H^+}(z,\bar{\zeta})\rangle=\\
\sqrt{\frac{2}{\pi}}\int_{H^+} f(\xi,\eta)\chi_+(\xi)\xi \exp(-(\eta-i\zeta)\xi)d\xi d\eta=\\
	\sqrt{\frac{2}{\pi}}\int_{H^+} f(\xi,\eta)\chi_+(\xi)\xi \exp(-(\eta+y)\xi)\exp(ix\xi)d\xi d\eta\\
\frac{1}{\sqrt{2\pi}}\int_{\R}\left(\int_{\R_+} 2\xi\chi_+(\xi)f(\xi,\eta)\exp(-(\eta+y)\xi)d\eta\right)\exp(ix\xi)d\xi=\\
(F^{-1}\otimes I)\left(2x\chi_+(x)\exp(-xy)\int_{R_+}f(x,\eta)\exp(-x\eta)d\eta\right)
\end{multline}
which yields
\begin{equation}
(B_1 f)(x,y)=2x\chi_+(x)\exp(-xy)\int_{R_+}f(x,\eta)\exp(-x\eta)d\eta
\end{equation}
Thus each $f\in L^2(H^+)$ takes the form
$f(x,y)=\chi_+(x)\sqrt{2y} \phi(x)\exp(-xy),$ $\phi(x)\in L^2(\R).$
Conversely any $f$ of this form satisfies
\begin{multline}
(B_1 f)(x,y)=2\chi_+(x)\exp(-xy)\int_{\R_+} f(x,\eta)\exp(-\eta x)d\eta=\\
2x\chi_+(x)\exp(-xy)\sqrt{2x}\phi(x)\int_{\R_+} \exp(-2\eta x)d\eta=\\
\chi_+(x)\sqrt{2x}\phi(x)\exp(-xy)\int_{\R_+} 2x\exp(-2\eta x)d\eta=\\
\chi_+(x)\sqrt{2x}\phi(x)\exp(-xy)=f(x,y)
\end{multline}
Thus $A_1^2$ consists of all functions $f$ having the form
\begin{equation}
f(x,y)=\chi_+(x)\sqrt{2y} \phi(x)\exp(-xy),\quad \phi(x)\in L^2(\R)
\end{equation}
and $\norm{f(x,y)}_{A_1^2}=\norm{\phi(x)}_{L^2(\R_+)}.$	
Now let us introduce the variable (which will allows us to revert to the usual $x,y$ after applying certain operators) by $w=u+iv.$ Define
\begin{equation}
(U_2 f)(u,v):=\frac{1}{\sqrt{2\abs{x}}}f\left(x,\frac{y}{2\abs{x}}\right)
\end{equation} and 
\begin{equation}
U:=U_2U_1:
L^2(\R)\otimes L^2(\R_+)\to L^2(\R)\otimes L^2(\R_+)
\end{equation}
We have the inverse $U_2^{-1}=U_2^*:L^2(\R)\otimes L^2(\R_+)\to L^2(\R)\otimes L^2(\R_+)$
\begin{equation}
U_2^{-1}:f(x,y)\mapsto \sqrt{2\abs{u}} f(u,2\abs{u}v)
\end{equation}
Denote by $A_2^2:=U_2 (A_1^2).$ Then 
the operator $B_2:=U_2 B_1 U_2^{-1}$ is the orthogonal projection of $L^2(H^+)$ onto $A_2^2.$
As before we use $w=u+iv$ for the variable before applying, in this case the transformation $U_2$, 
in order to  revert to $x,y$ in the image (recall that 
for example a transformation like the Fourier transform require new notation of the variable).
We have
\begin{multline}
(B_2f)(x,y)=U_2\left(2u\chi_+(u)\exp(-uv)\int_{R^+}\exp(-u\eta)\sqrt{2u} f(u,2u\eta)d\eta\right)\\
=U_2\left(\sqrt{2u}\chi_+(u)\exp(-uv)\int_{R^+} f(u,\nu)\exp\left(-\frac{\nu}{2}\right)d\nu\right)=\\
\chi_+(x)\exp\left(-\frac{y}{2}\right)\int_{R^+} f(u,\nu)\exp\left(-\frac{\nu}{2}\right)d\nu
\end{multline}
Next recall that the Laguerre polynomials\index{Laguerre polynomial} of degree $j\in \N,$ are defined according to
\begin{equation}
L_j(y):=L_j^0(y)=\frac{\exp(y)}{j!}\partial_y^j(\exp(-y)y^n)=\sum_{k=0}^j\frac{j!}{k!(j-k)!}\frac{(-y)^k}{k!}, y\in \R_+
\end{equation}
and the system of functions $\ell_j(y):=\exp(-y^2/2)L_j(y),$ $j\in \N,$ form an orthonormal basis for $L^2(\R_+).$
Let $L_j$ be the one-dimensional subspace of $L^2(\R_+)$ generated ny $\ell_j(y).$
Set 
\begin{equation}
L_n^\oplus:=\oplus_{k=0}^n L_k
\end{equation}
In particular, the one-dimensional projection $P_0$ of $L^2(\R_+)$ onto $L_0$ has the form
\begin{equation}
P_0(\psi)(y)=\langle \psi,\ell_0\rangle\cdot \ell_0=\exp\left(-\frac{y}{2}\right)\int_{\R_+}\psi(t)\exp\left(-\frac{t}{2}\right)dt
\end{equation}
which means $B_2=\chi_+(x)I\otimes P_0.$
The results described render the following theorem, see Vasilevski \cite{vasilevskibok}, Theorem 3.1.1, p.69.
\begin{theorem}\label{vasilevskithm311}
	The unitary operator $U=U_2 U_1$ defines isometric isomorphisms
	\begin{equation}
	U:A^2_{1}(H^+)\to L^2(\R_+)\otimes L_{0}
	\end{equation} 
	and under $U$ the Bergman projection is unitarily equivalent to $\chi_+ I\otimes P_0$ i.e.\
	\begin{equation}
U B_{H^+} U^{-1}=\chi_+ I\otimes P_0
	\end{equation} 	
\end{theorem}

Denote by $H^2_+(\R)$ the Hardy space on the upper half-plane, i.e.\ the set of $L^2(\R)$ functions that have analytic continuation $f(z)$ to $H^+$
such that
\begin{equation}
\sup_{y>0}\int_{\R+iy}\abs{f(x+iy)}^2dx<\infty
\end{equation}
Denote by $P_\R^+$ the so-called {\em Szeg\"o transform}, which is the orthogonal projection of 
$L^2$ onto $H^2_+(\R)$.
Then the Fourier transform $F$ gives an isometric isomorphism of $L^2(\R)$ under which
the Hardy space $H^2_+(\R)$ is mapped onto $L^2(\R_+)$ and  $P_\R^+:L^2(\R)\longrightarrow H^2_+(\R)$
is unitarily equivalent to
$FP^+_\R F^{-1}=\chi_+ I.$ Furthermore, if $H_-^2(\R)$ denotes the Hardy space on the lower half-plane and $P^-_\R$ the orthogonal projection
of $L^2(\R)$ onto $H^2_-(\R)$ then
$L^2(\R)=H^2_+(\R)\otimes H^2_-(\R)$ and we have $F(H^2_-(\R))=L^2(\R_-)$ and $FP_{\R}^- F^{-1} =\chi_- I.$ 
By Theorem \ref{vasilevskithm311} this renders the following relation to Hardy spaces (Vasilevski \cite{vasilevskibok}, Thm. 3.2.1).
\begin{theorem}
The unitary operator $W=(F^{-1}\otimes I)U_2(F\otimes I)$ gives an isometric isomorphism of $L^2(H^+)=L^2(\R)\otimes L^2(\R_+)$ 
such that
\begin{equation}
W(A^2_1(H^+))=H^2_+(\R)\otimes L_0
\end{equation}
and
\begin{equation}
W B_{H^+} W^{-1}=P^+_{\R}\otimes P_0
\end{equation}
\end{theorem}
\begin{definition}
For a domain $\Omega\subset\C$, define $\tilde{A}^2_q(\Omega)$ to be the closed subspace of $L^2(\Omega)$ 
consisting of all functions anti-$q$-analytic on $\Omega$
(i.e.\ $\partial_z^q f=0$ on $\Omega$ and $f\in L^2(\Omega)$), $\tilde{A}^2_{(1)}(\Omega)=\tilde{A}^2_1(\Omega)$, and (analogous to the definition of the spaces
$A^2_{(1)}(\Omega)$) for $q>1$,
$\tilde{A}^2_{(q)}(\Omega)=\tilde{A}^2_q(\Omega)\ominus \tilde{A}^2_{q-1}(\Omega).$
It is well-known (see e.g.\ Dzhuraev \cite{dzhuraevbok}) that the orthogonal projection $\tilde{B}_\Omega$ of $L^2(\Omega)$ onto $\tilde{A}^2_q(\Omega)$
is given by 
\begin{equation}
(\tilde{B}_\Omega f)(z)=\int_\Omega K_\Omega(\zeta,z)f(\zeta)d\mu(\zeta)
\end{equation}
\end{definition}
\begin{theorem}
The unitary operator $U=U_2U_1$ gives an isometric isomorphism of $L^2(H^+)=L^2(\R)\otimes L^2(\R_+)$ 
under which \\
(1) $\tilde{A}_1^2(H^+)$ is mapped onto $L^2(\R_)\otimes L_0$.
\\
(2) $\tilde{B}_{H^+}$ is unitarily equivalent to $\chi_-I\otimes P_0$ (where $\chi_-$ denotes the characteristic function
of the negative half-line) i.e.\
\begin{equation}
U\tilde{B}_{H^+}U^{-1}=\chi_- I\otimes P_0
\end{equation}
\end{theorem}
\begin{proof}
	This is obtained by following all the steps of Theorem \ref{vasilevskithm311}.
	The $L^2(H^+)$-solutions of $(u+\partial_v)f=0$ have the form
	\begin{equation}
	f(u,v)=\chi_-(u)\sqrt{2\abs{u}}f(u)\exp(uv),\quad f(u)\in L^2(\R)
	\end{equation}
	Furthermore, with $w=u+iv,$ the Fourier transform of $\overline{K_{H^+}(w,\zeta)}=-\frac{1}{\pi}\frac{1}{(\bar{w}-\zeta)^2}$
	with respect to $\xi=\re w$ is given by (see Bateman \cite{bateman})
	\begin{equation}
	-\sqrt{\frac{2}{\pi}}\xi \chi_-(\xi)\exp((\eta -iw)\cdot\xi)
	\end{equation}
	The projection $\tilde{B}_1=U_1\tilde{B}_{H^+} U^{-1}$ onto $\tilde{A}_1^2 =U_1(\tilde{A}_1^2(H^+))$ takes the form
	\begin{equation}
	(\tilde{B}_1 f)(u,v)=2\abs{u}\chi_-(u)\exp(uv)\int_{\R_+} f(u,v)\exp(u\eta)d\eta
	\end{equation}
	The proof follows from this in the same manner as for Theorem \ref{vasilevskithm311}.
\end{proof}
\begin{corollary}
	The unitary operator $W=(F^{-1}\otimes I)U_2(F\otimes I)$ gives an isometric isomorphism of $L^2(H^+)=L^2(\R)\otimes L^2(\R_+)$ 
	such that
	\begin{equation}
	W(\tilde{A}^2_1(H^+))=H^2_-(\R)\otimes L_0
	\end{equation}
	and
	\begin{equation}
	W \tilde{B}_{H^+} W^{-1}=P^-_{\R}\otimes P_0
	\end{equation}
\end{corollary}
Now clearly
\begin{equation}
L^2(\R)=H^2_+ \oplus H^2_-(\R)
\end{equation}
and the orthogonal complement $A_1^2(H^+)\oplus \tilde{A}^2_1(H^+)$ 
\begin{equation}
(A_1^2(H^+)\oplus \tilde{A}^2_1(H^+))^\perp =L^2(H^+)\ominus (A_1^2(H^+)\oplus \tilde{A}^2_1(H^+))
\end{equation}
satisfies
\begin{equation}
W\left((A_1^2(H^+)\oplus \tilde{A}^2_1(H^+))^\perp\right) =L^2(\R)\otimes L_0^\perp
\end{equation}

Introduce the isometric embedding $R_0:L^2(\R_+)\longrightarrow  L^2(\R)\otimes L^2(\R_+)$ according to
\begin{equation}
(R_0f)(x,y):=\chi_+(x)f(x)\ell_0(y)
\end{equation}
here the function $f(x)$ is extended to an element of $L^2(\R)$ by setting $f(x)\equiv 0$ for $x<0.$ Note
that the image of $R_0$ is $A_2^2.$ We have that the adjoint $R_0^*:L^2(H^+)\longrightarrow  L^2(\R_+)$ is given by
\begin{equation}
(R^*_0 f)(x)=\chi_+(x)\int_{\R_+}f(x,y)\ell_0(\eta)d\eta
\end{equation}
and
$R^*_0 R_0=I :L^2(\R+)\longrightarrow  L^2(\R_+),$
$R^*_0 R_0=B_2 :L^2(H^+)\longrightarrow  A_2^2=L^2(\R_+)\otimes L_0$.
Then 
\begin{equation}
R:=R_0^* U
\end{equation}
maps $L^2(H^+)$ onto $L^2(\R_+)$ and its restriction to $A_1^2(H^+)$ is an isometric isomorphism
\begin{equation}
R|_{A_1^2(H^+)} \longrightarrow  L^2(H^+)
\end{equation}
Furthermore, its adjoint is an isometric isomorphism of $L^2(\R_+)$ onto $A_1^2(H^+)$
\begin{equation}
R^*=U^*R_0: L^2(H^+)\longrightarrow A_1^2(H^+)  
\end{equation}
and
$R^* R=I :L^2(\R+)\longrightarrow  L^2(\R_+),$
$R^* R=B_2 :L^2(H^+)\longrightarrow  A_1^2(H^+)$.
\begin{proposition}
$R^*=U^*R_0:L^2(\R_+)\longrightarrow A_1^2(H^+)$ is given by
\begin{equation}
(R^*f)(z)=\frac{1}{\sqrt{\pi}}\int_{\R_+}\sqrt{\xi}f(\xi)\exp(iz\cdot\xi)d\xi
\end{equation}
\end{proposition}
\begin{proof}
We have
\begin{multline}
(R^*f)(z)=(U_1^*U_2^*R_0 f)(z)=\\
(F^{-1}\otimes I)\left(\chi_+(\xi)f(\xi)\sqrt{2\xi}\exp(-xy)\right)=\\
\frac{1}{\sqrt{2\pi}}\int_{\R}\chi_+(\xi)f(\xi)\sqrt{2\xi}\exp(-xy)\exp(ix\xi)d\xi
=\frac{1}{\sqrt{\pi}}\int_{\R_+}\sqrt{\xi}f(\xi)\exp(iz\cdot\xi)d\xi
\end{multline}
This completes the proof.
\end{proof}
This yields the inverse isomorphism
\begin{equation}
R:A_2^2\longrightarrow L^2(H^+)
\end{equation}
according to
\begin{equation}
(Rf)(x)=\sqrt{x}\frac{1}{\sqrt{\pi}}\int_{H^+}f(w)\exp(-i\bar{w}\cdot x)d\mu(w)
\end{equation}
\begin{theorem}
	The unitary operator $U=U_2U_1 : L^2(\R)\otimes L^2(\R_+)\longrightarrow L^2(\R)\otimes L^2(\R_+)$
	maps $A_n^2(H^+)$ onto $L^2(\R_+)\otimes L^\oplus_{n-1}.$
\end{theorem}
\begin{proof}
	The space $U(A^2_n(H^+))$ coincides with the set of functions in $L^2(H^+)=L^2(\R)\otimes L^2(\R_+)$ satisfying
	\begin{multline}\label{vasseekvven}
	U(\partial_u +i\partial_v U^{-1} f=i^n U_2(u+\partial_v)^n U_2^{-1} f=\\
	i^n\abs{x}^n\mbox{sign}(x)+2\partial_y)^nf=0
	\end{multline}
	The intersection of the general solution of Eqn.(\ref{vasseekvven}) with $L^2(\R)\otimes L^2(\R_+)$ consists of 
	the set of functions of the form
	\begin{equation}
	\sum_{k=0}^{n-1} \chi_+(x)\psi_k(x)y^k\exp\left(-\frac{y}{2}\right),\quad \psi(x)\in L^2(\R),\quad k=0,\ldots,n-1
	\end{equation}
	which upon rearranging polynomials with respect to $y$ coincides with the functions of the form
	\begin{equation}
	\sum_{k=0}^{n-1} \chi_+(x)g_k(x)L_k(y)\exp\left(-\frac{y}{2}\right)=\sum_{k=0}^{n-1} \chi_+(x)g_k(x)\ell_k(y)
	\end{equation}
	where $g_k(x)\in L^2(\R),$ $k=0,\ldots,n-1.$
	This completes the proof.
\end{proof}
Note that by the definition of true-$q$-polyanalytic functions we have the decomposition
\begin{equation}
A_n^2=\oplus_{k=0}^n A_{(k)}^2,\quad \tilde{A}_n^2=\oplus_{k=0}^n \tilde{A}_{(k)}^2
\end{equation}
which renders the following corollary.
\begin{corollary}
	The unitary operator $U=U_2U_1 : L^2(\R)\otimes L^2(\R_+)\longrightarrow L^2(\R)\otimes L^2(\R_+)$
	maps $A_{(n)}^2(H^+)$ onto $L^2(\R_+)\otimes L_{n-1}.$
\end{corollary}
\begin{theorem}
	The unitary operator $U=U_2U_1 : L^2(\R)\otimes L^2(\R_+)\longrightarrow L^2(\R)\otimes L^2(\R_+)$
	maps $\tilde{A}_{n}^2(H^+)$ onto $L^2(\R_+)\otimes L^{\oplus}_{n-1}.$
\end{theorem}
\begin{proof}
	The space $U(\tilde{A}^2_n(H^+))$ coincides with the set of functions in $L^2(H^+)=L^2(\R)\otimes L^2(\R_+)$ satisfying
	\begin{multline}\label{vasseekvven1}
	U(\partial_u -i\partial_v U^{-1} f=i^n U_2(u-\partial_v)^n U_2^{-1} f=\\
	i^n\abs{x}^n\mbox{sign}(x)-2\partial_y)^n f=0
	\end{multline}
	The intersection of the general solution of Eqn.(\ref{vasseekvven1}) with $L^2(\R)\otimes L^2(\R_+)$ consists of 
	the set of functions of the form
	\begin{equation}
	\sum_{k=0}^{n-1} \chi_-(x)\psi_k(x)y^k\exp\left(-\frac{y}{2}\right),\quad \psi(x)\in L^2(\R),\quad k=0,\ldots,n-1
	\end{equation}
	which upon rearranging polynomials with respect to $y$ coincides with the functions of the form
	\begin{equation}
	\sum_{k=0}^{n-1} \chi_-(x)g_k(x)L_k(y)\exp\left(-\frac{y}{2}\right)=\sum_{k=0}^{n-1} \chi_-(x)g_k(x)\ell_k(y)
	\end{equation}
	where $g_k(x)\in L^2(\R),$ $k=0,\ldots,n-1.$
	This completes the proof.
\end{proof}
\begin{corollary}
	The unitary operator $U=U_2U_1 : L^2(\R)\otimes L^2(\R_+)\longrightarrow L^2(\R)\otimes L^2(\R_+)$
	maps $\tilde{A}_{(n)}^2(H^+)$ onto $L^2(\R_-)\otimes L_{n-1}.$
\end{corollary}
The above results immediately render the following decomposition due to Vasilevski \cite{vasilevski1999}, Thm .4.5, p.482. 
\begin{theorem}
	The operator $U=U_2 U_1$ defines isometric isomorphisms
	\begin{equation}
	W:A^2_{(j)}(H^+)\to H^2_+(\R)\otimes L_{j-1},\quad W:\tilde{A}^2_{(j)}(H^+)\to H_-^2(\R)\otimes L_{j-1}
	\end{equation} 
	\begin{equation}
	W:A^2_{j}(H^+)\to H_+^2(\R)\otimes \left(\oplus_{k=0}^{j-1} L_{k}\right),\quad W:\tilde{A}^2_{(j)}(H^+)\to H_-^2(\R)\otimes \left(\oplus_{k=0}^{j-1} L_{k}\right)
	\end{equation} 
	\begin{equation}
	W:\oplus_{k=1}^{\infty} A^2_{(j)}(H^+)\to L^2(\R_+)\otimes H_+^2(\R)
	\end{equation}
	\begin{equation}
	 W:\oplus_{k=1}^{\infty} \tilde{A}^2_{(j)}(H^+)\to H_-^2(\R)\otimes L^2(\R_+)  
	\end{equation}
	and we have a decomposition
	\begin{equation}
	L^2(H^+)=\oplus_{k=1}^{\infty} (A^2_{(j)}(H^+)\oplus \tilde{A}^2_{(j)}(H^+))=
	\left(\oplus_{k=1}^{\infty} A^2_{(j)}(H^+)\right)\oplus \left(\oplus_{k=1}^{\infty} \tilde{A}^2_{(j)}(H^+)\right)
	\end{equation} 
\end{theorem}
\begin{definition}
	Let $H^+$ denote the upper half-plane and let $B_{H^+,(n)}$ ($B_{H^+,n}$) denote the {\em poly-Bergman projection}\index{Poly-Bergman projection}
		where the Bergman projection for $A^2_{1}(H^+)$ is given by Eqn.(\ref{orthokernelekvenbra}) 
	i.e.\ it is given by a reproducing kernel function 
	\begin{equation}
	K_{H^+}(z,\zeta)=-\frac{1}{\pi}\frac{1}{(z-\bar{\zeta})^2}
	\end{equation}
	and for $n\geq 1$ we define the kernels
	\begin{equation}
	K^{(n)}_{H^+}(z,\zeta):=K_{H^+}(z,\zeta)\sum_{j=0}^{n-1}\sum_{k=0}^{n-1} \kappa_{j,k}^{n-1}
	\left(\frac{z-\bar{z}}{z-\bar{\zeta}}\right)^j\left(\frac{\zeta-\bar{\zeta}}{z-\bar{\zeta}}\right)^k
	\end{equation}
	where
	\begin{equation}
	\kappa_{j,k}^{n-1}=(-1)^{j+k}\left(\frac{n!}{j!k!}\right)^2 \frac{(j+k+1)!}{(n-j)!(n-k)!}
		\end{equation}
		and define
		\begin{equation}
		(B_{H^+,(n)}f)(z)=\int_{H^+} K^{(n)}_{H^+}(z,\zeta)f(\zeta)d\mu(\zeta) 
		\end{equation}
	\end{definition}
	\begin{theorem}\label{ortokerneltheorem}
		$B_{H^+,(n)}$ is the orthogonal projection of $L^2(H^+)$ onto $A^2_{(n)}(H^+)$.
		\end{theorem}
		\begin{proof}
		The orthogonal projection $B_2^{(n)}=UB_{H^+}U^{-1}:L^2(H^+)\to L^2(\R_+)\otimes_{n-1}$ is given by
		\begin{equation}
		B_2^{(n)}=\chi_+(x)\otimes P_{n-1}
		\end{equation}
		where	
		\begin{equation}
		(P_{n-1}\psi)(y)=\langle \psi,\ell_{n-1}\rangle\ell_{n-1}=\ell_{n-1}(y)\int_{\R_+}\psi(t)\ell_{n-1}(t)dt
		\end{equation}
		is the orthogonal projection of $L^2(\R_+)$ onto $L_{n-1}$, i.e.\
		\begin{equation}
		(B_2^{(n)}f)(x,y)=\chi_+(x)\ell_{n-1}(y)\int_{\R_+}f(x,t)\ell_{n-1}(t)dt
		\end{equation}
		Let $w=u+iv,$ $\zeta=\xi+i\eta$. For $B_1^{(n)}=U_1^{1}B_2^{(n)}U_2$ we have
		\begin{multline}
		(B_1^{(n)}f)(u,v)=\sqrt{2\abs{u}}\chi_+(u)\ell_{n-1}(\abs{u}v)\int_{\R_+}\frac{1}{\sqrt{2\abs{u}}} f\left(u,\frac{t}{2\abs{u}}\right)\ell_{n-1}(t)dt\\
		=2u\exp(-uv)\chi_+(u)L_{n-1}(2uv)\int_{\R_+} f(u,\eta)\ell_{n-1}(2u\eta)d\eta
		\end{multline}
		This implies
		\begin{multline}
		(F^{-1}\otimes I)B_1^{(n)}f=\\
		\frac{1}{\sqrt{2\pi}}\int_\R
		\left(\int_{\R_+} \chi_+(\xi)2\xi\exp(-v\xi)L_{n-1}(2\xi v)L_{n-1}(2\xi\eta)f(\xi,\eta)\exp(-\eta\xi)d\eta\right)\exp(iu\xi)d\xi\\
		=\frac{1}{\sqrt{2\pi}}\int_{H^+} f(\xi,\eta)\chi_+(\xi)\xi L_{n-1}(2\xi v)L_{n-1}(2\xi\eta)\exp(-(\eta-iw)\xi)d\xi d\eta
			\end{multline}
			Since $B_2^{(n)}=U_2(F\otimes I)B_{H^+,(n)}(F^{-1}\otimes I)U_2^{-1}$ we have
			\begin{equation}
			B_1^{(n)}=U_2^{-1}B_2^{(n)}U_2 =(F\otimes I)B_{H^+,(n)}(F^{-1}\otimes I)
			\end{equation}
			and
			\begin{multline}
			(F^{-1}\otimes I)B_1^{(n)}f= B_{H^+,(n)}(F^{-1}\otimes I)f=\\
			\langle (F^{-1}\otimes I)f,\overline{K^{(n)}_{H^+}(w,\zeta)}\rangle =
			\langle f,(F\otimes I)\overline{K^{(n)}_{H^+}(w,\zeta)}\rangle
			\end{multline}
			which yields
			\begin{equation}
			(F\otimes I)\overline{K^{(n)}_{H^+}(w,\bar{\zeta})}=\frac{1}{\sqrt{2\pi}}\chi_+(\xi)\xi
			L_{n-1}(2\xi v)L_{n-1}(2\xi\eta)\exp(-(\eta-i\bar{w})\xi)
			\end{equation}
			Now by the definition of $L_n(y)$ we have
			\begin{equation}
			L_{n-1}(y)=\sum_{k=0}^{n-1}\frac{n!}{k!(n-k)!}\frac{(-y)^k}{k!}=\sum_{k=0}^{n-1}\lambda_k^{n-1}y^k,\quad 
			\lambda_k^{n-1}:=\frac{(-1)^k n!}{k!(n-k)!k!}
			\end{equation}
			thus
			\begin{equation}
			(F^\otimes I)\overline{K^{(n)}_{H^+}(w,\zeta)}=
			\sum_{j=0}^{n-1}\sum_{k=0}^{n-1}\lambda_j^{n-1}\lambda_k^{n-1} 2^{j+k}v^j\eta^k\left(\frac{1}{\sqrt{2\pi}}\chi_+(\xi)\xi^{j+k+1}\exp(\eta +i\bar{w})\xi\right)
			\end{equation}
			Now the Fourier transform of 
			\begin{equation}
			\frac{1}{\pi}\frac{m!}{(i\bar{w}+\eta -i\xi)^{m+1}}
				\end{equation}
				is given by (see Bateman \cite{bateman})
				\begin{equation}
				\frac{2}{\sqrt{\pi}}\chi_+(\xi)\xi^m\exp(-(\eta+i\bar{w})\xi)
				\end{equation}
				This yields
				\begin{equation}
				\overline{K^{(n)}_{H^+}(w,\zeta)}=\frac{1}{\pi} \sum_{j=0}^{n-1}\sum_{k=0}^{n-1}\lambda_j^{n-1}\lambda_k^{n-1} 2^{j+k}v^j\eta^k
				\frac{(j+k+1)!}{(i\bar{w} +\eta -i\xi)^{j+k+2}}
				\end{equation}
				thus
				\begin{equation}
				K^{(n)}_{H^+}(w,\zeta)=-\frac{1}{\pi}\frac{1}{(z-\bar{\zeta})^2} \sum_{j=0}^{n-1}\sum_{k=0}^{n-1}\lambda_j^{n-1}\lambda_k^{n-1}2^{j+k}v^j\eta^k
				\frac{(iv)^j(i\eta)^k}{(w-\bar{\zeta})^{j+k}}
				\end{equation}
				Introducing
				\begin{equation}
				\kappa_{j,k}^{n-1}=\lambda_j^{n-1}\lambda_k^{n-1}(j+k+1)!
				=
				(-1)^{j+k}\left(\frac{n!}{j!k!}\right)^2 \frac{(j+k+1)!}{(n-j)!(n-k)!}
				\end{equation}
				gives
				\begin{multline}
				K^{(n)}_{H^+}(w,\zeta)=-\frac{1}{\pi}\frac{1}{(z-\bar{\zeta})^2} 
				\sum_{j=0}^{n-1}\sum_{k=0}^{n-1}\kappa_{j,k}^{n-1}\left(\frac{2iv}{w-\bar{\zeta}}\right)^j
				\left(\frac{2i\eta}{w-\bar{\zeta}}\right)^k =\\
				-\frac{1}{\pi}\frac{1}{(z-\bar{\zeta})^2} 
				\sum_{j=0}^{n-1}\sum_{k=0}^{n-1}\kappa_{j,k}^{n-1}\left(\frac{w-\bar{w}}{w-\bar{\zeta}}\right)^j
				\left(\frac{\zeta -\bar{\zeta}}{w-\bar{\zeta}}\right)^k
				\end{multline}
				This completes the proof.
		\end{proof}
Let us also mention the following popular operators.
\begin{equation}\label{dzhuraevop1}
(Sf)(z):=-\frac{1}{\pi}\int_{H^+} \frac{f(\zeta)}{(\zeta -z)^2} d\mu(\zeta)
\end{equation}
with its adjoint given by
\begin{equation}\label{dzhuraevop2}
(S^*f)(z)=-\frac{1}{\pi}\int_{H^+} \frac{f(\zeta)}{(\bar{\zeta} -\bar{z})^2} d\mu(\zeta)
\end{equation}
Vasilevski (see \cite{vasilevskibok} Cor.3.5.6, p.87) proves that any true $j$-polyanalytic function $f$ can be written as
\begin{equation}
f(z)=S^{j-1}\phi, \mbox{ some }\phi\in A^2_{1}(H_+)
\end{equation}
and any $f\in \tilde{A}_2^{(j)}$ can be written as
\begin{equation}
f(z)=(S^*)^{j-1}\varphi, \mbox{ some }\varphi\in A^2_{1}(H_+)
\end{equation}
Vasilevski applies the same method as described above for the half-plane, for obtaining a decomposition for the case of the unit disc, 
see e.g.\ Vasilevski \cite{vasilevskibok}, Ch.4,
and as we shall see below the same methods can be repeated for studying so-called poly-Fock spaces, which are more or less 
a special case of weighted poly-Bergman spaces that, due to their applications often get separate treatment.

\subsection{Vasilevski's decomposition of poly-Bergman spaces for the upper half-space over complex quaternions}
Vasilevski \cite{vasilevskiquaternion} also performed analogous studies on poly-Bergman spaces in hypercomplex analysis.
As before let $H^+$ denote the upper half-plane $\{\im z>0\}$ and $H_-$ the lower half-plane.
Denote by $\mathbb{H}(\C):=\mathbb{H}\otimes \C$, called the complex quaternions.
Let $\Omega\subset \R^4$ be a domain.
If $\psi=\{\psi^0,\psi^1,\psi^2,\psi^3\}$ is a set of orthonormal real quaternions and $\psi^0=1$ then $\psi^1,\psi^2,\psi^3$
are pure imaginary quaternions.
On the set $C^1(\Omega,\mathbb{H})$ (or $C^1(\Omega,\mathbb{H}(\C))$) the quaternionic Cauchy-Riemann operator is defined
by
\begin{equation}
^\psi D=\sum_{k=0}^3 \psi^k \partial_{x_k}
\end{equation}
And functions annihilated by the $q$:th power of this operator are called $q$-left-$\psi$-hyper-holomorphic (where we shall
oftentimes simply call them $q$-$\psi$-hyper-holomorphic).
Denote 
\begin{equation}
\sigma_{\psi,x}:=\sum_{k=0}^3 (-1)^k\psi^k d\hat{x}_k
\end{equation}
(where $d\hat{x}_k$ denotes $dx_0\wedge \cdots\wedge dx_3$ with $dx_k$ removed). Denote the fundamental solution
\begin{equation}
\mathcal{K}_\psi(x)=\frac{1}{2\pi^2\abs{x}^4}\sum_{k=0}^3 \bar{\psi^k}\cdot x_k
\end{equation}

\begin{definition}
	Let $\Pi$ denote the upper half-space in $\R^4$.
	Let $L^2(\Pi,\mathbb{H}(\C))$ denote the complex quaternion valued $L^2$ functions.
	Let $\psi=\{\psi^+,\psi^1,\psi^2,\psi^3\}$ be the set of orthonormal real quaternions, $\psi^0=1$.
	For each $n\in \N$ denote by $^\psi\mathcal{A}^2_n(\Pi,\mathbb{H}(\C))$ ($^{\bar{\psi}}\mathcal{A}^2_n(\Pi,\mathbb{H}(\C))$)
	the closed subspace of $L^2(\Pi,\mathbb{H}(\C))$ consisting of 
	$n$-$\psi$-hyper-holomorphic ($n$-$\bar{\psi}$-hyper-holomorphic)\index{$n$-$\psi$-hyper-holomorphic function} functions, 
	in the sense that they
	satisfy
	\begin{equation}
	\left(\sum_{k=0}^3 \psi^k \partial_{x_k}\right)^n f=0,\quad \left(\left(\sum_{k=0}^3 \bar{\psi}^k \partial_{x_k}\right)^n f=0\right)
	\end{equation}
	
	Define the set of {\em true-$n$-$\psi$-hyper-holomorphic functions} ({\em true $n$-$\bar{\psi}$-hyper-holomorphic function})
	by \index{True-$n$-$\psi$-hyper-holomorphic function}
	\begin{equation}
	^\psi\mathcal{A}^2_{(n)}(\Pi,\mathbb{H}(\C))\ominus ^\psi\mathcal{A}^2_{(n-1)}(\Pi,\mathbb{H}(\C))
	\end{equation}
	\begin{equation}
	^{\bar{\psi}}\mathcal{A}^2_{(n)}(\Pi,\mathbb{H}(\C))\ominus ^{\bar{\psi}}\mathcal{A}^2_{(n-1)}(\Pi,\mathbb{H}(\C))
	\end{equation}
	Introduce the operators
	\begin{equation}
	^\psi Q_+=^\psi q_+(\xi)I,\quad ^\psi Q_-=^\psi q_-(\xi)I
	\end{equation}
	\begin{equation}
	^\psi q_\pm (\xi):= \frac{1}{2}\left(1\pm i\sum_{k=1}^3 \psi^k \frac{\xi_k}{\abs{\xi}}\right),\quad \xi\in \R^3
	\end{equation}
	Then $^\psi Q_+$ and $^\psi Q_-$ are mutually orthogonal projections on $L^2(\R^3,\mathbb{H}(\C))$
	and $^\psi Q_+^\psi Q_-=I$.
	Let $\{\ell_j(y)\}_{j=0}^\infty$ be the orthonormal base in $L^2(\R_+)$ given by 
	$\ell_j(y):=\exp(y/2)L_j(y)$ where $L_j(y)$ is the Laguerre polynomial of degree $j.$ Denote by $L_j$
	the one-dimensional subspace of $L^2(\R_+)$ generated by $\ell_j(y)$ and denote by $P_j$ the orthogonal projection of $L^2(\R_+)$
	onto $L_j.$
	Set
	\begin{equation}
	L_j^\oplus \oplus_{k=0}^j L_k
	\end{equation}
\end{definition}

See Chapter \ref{hypercomplexsec} for an introduction to the analogues of poly-analytic functions
in hypercomplex analysis. In particular, The following Cauchy integral formula is known, cf. Theorem \ref{brackscauchyformula}.
\begin{theorem}
	Let $\Omega\subset\R^4$ be a bounded domain with smooth boundary $\Gamma$ and let 
	$f$ be a $\psi$-hyper-holomorphic function on $\Omega$
	and continuous on $\overline{\Omega}$. Then 
	\begin{equation}
	\int_\Gamma \mathcal{K}_\psi(\tau-x)\cdot\sigma_{\psi,\tau}\cdot f(\tau)=
	\left\{
	\begin{array}{ll}
	f(x) & , x\in \Omega\\
	0 & ,x\in \R^4\setminus \overline{\Omega}
	\end{array}
	\right.
	\end{equation}
\end{theorem}
For a bounded domain $\Omega\subset\R^4$ with smooth boundary $\Gamma$, the following 
singular integral is bounded on $L^2(\Gamma,\mathbb{H}(\C))$
\begin{equation}
(^\psi S_\Gamma f)(t):=2\int_\Gamma \mathcal{K}_\psi (\tau-t)\cdot \sigma_{\psi,\tau}\cdot f(\tau)
\end{equation}
and furthermore it is self-adjoint and satisfies $^\psi S^2_\Gamma =I.$
For $\Gamma=\{(0,x_1,x_2,x_3)\}\subset\R^4$ we have
\begin{multline}\label{vasilq21}
(^\psi S_\Gamma f)(t):=\frac{1}{\pi^2}\int_\R^3 \frac{1}{\abs{\tau-t}^4} \sum_{k=1}^3 \bar{\psi}^k(\tau_k-t_k)f(\tau)d\tau=\\
\sum_{k=1}^3 \bar{\psi}^k \frac{1}{\pi^2}\int_\R^3 \frac{\tau_k-t_k}{\abs{\tau-t}^4} f(\tau)d\tau=
\sum_{k=1}^3 \bar{\psi}^k R_k =-\sum_{k=1}^3 \psi^k R_k
\end{multline}
where the $R_k$ are called the Riesz operators and are given by
\begin{equation}\label{vasilq22}
(R_k f)(t):=\frac{1}{\pi^2} \int_\R^3 \frac{\tau_k-t_k}{\abs{\tau-t}^4} f(\tau)d\tau,\quad k=1,2,3
\end{equation}
Let us (in this section solely, in order to conform with the cited literature) 
Denote the Fourier transform and inverse Fourier transform respectively of a function $f\in L^2(\R^3,\mathbb{H}(\C))$ by
\begin{equation}\label{vasilq23}
(Ff)(t)=(2\pi)^{-\frac{3}{2}}\int_{\R^3} \exp(-it\cdot\xi)\cdot f(\xi)d\xi,\quad 
(F^{-1}f)(t)=(2\pi)^{-\frac{3}{2}}\int_{\R^3} \exp(it\cdot\xi)\cdot f(\xi)d\xi
\end{equation} 
we have
\begin{equation}\label{vasilq24}
R_k=F^{-1}i\frac{\xi_k}{\abs{\xi}} F
\end{equation}
Denote by $^\psi H^2_\pm \subset L^2(\R^3,\mathbb{H}(\C)))$ the Hardy space of functions which admit $\psi$-hyperholomorphic extension to 
$\{(x_0,x_1,x_2,x_3)\in \R^4: \pm x_0 >0\}.$
Then the so-called {\em Szeg\"o} orthogonal projections of $L^2(\R^3,\mathbb{H}(\C))$ onto
$^\psi H^2_\pm \subset L^2(\R^3,\mathbb{H}(\C)))$ respectively are given by
\begin{equation}\label{vasilq25}
^\psi P_{\R^3}^+ :=\frac{1}{2}(I-^\psi S_{\R^3}),\quad ^\psi P_{\R^3}^- :=\frac{1}{2}(I+^\psi S_{\R^3})
\end{equation}
Set for $x\in \R^3$
\begin{equation}\label{vasilq26}
^\psi q_+(\xi) :=\frac{1}{2}\left(1+i\sum_{k=1}^3 \psi^k\frac{\xi_k}{\abs{\xi}}\right),\quad 
^\psi q_-(\xi) :=\frac{1}{2}\left(1-i\sum_{k=1}^3 \psi^k\frac{\xi_k}{\abs{\xi}}\right)
\end{equation}
Then the following operator are mutually orthogonal projections
\begin{equation}\label{vasilq27}
^\psi Q_+ :=^\psi q_+(\xi) I,\quad ^\psi Q_- :=^\psi q_-(\xi) I
\end{equation}
and we have
\begin{equation}\label{vasilq28}
^\psi Q_+ + ^\psi Q_- = I
\end{equation}
Denote by $^\psi L^2_\pm:=^\psi L^2_+(\R^3,\mathbb{H}(\C))$ the images of $^\psi Q_\pm$ respectively.
Applying the Fourier transform to the projections of Eqn.(\ref{vasilq25}), using then representation of 
Eqn.(\ref{vasilq21}) and then applying Eqn.(\ref{vasilq24}),Eqn.(\ref{vasilq27}),Eqn.(\ref{vasilq28}) yields the following
result.
\begin{theorem}\label{vasilqthm21}
	\begin{equation}
	F:L^2(\R^3,\mathbb{H}(\C))\to L^2(\R^3,\mathbb{H}(\C))
	\end{equation}
	is an isometric isomorphism under which:
	(1) $^\psi H^2_\pm (\R^3,\mathbb{H}(\C))$ are mapped onto $^\psi L^2_\pm$.\\
	(2) $^\psi P_{\R^3}^\pm$ are unitary equivalent to $^\psi Q_\pm$
	\begin{equation}
	F^\psi P^\pm F^{-1} =^\psi Q_\pm
	\end{equation}
	and
	\begin{equation}
	L^2(\R^3,\mathbb{H}(\C))=^\psi H^2_+ (\R^3,\mathbb{H}(\C))\oplus ^\psi H^2_- (\R^3,\mathbb{H}(\C))
	\end{equation}
\end{theorem}
Define the unitary operator
\begin{equation}
U_1=I\otimes F
\end{equation}
acting on $L^2(\Pi,\mathbb{H}(\C))=L^2(\R_+)\otimes L^2(\R^3,\mathbb{H}(\C))$. The image of $^\psi\mathcal{A}^2_1=U_1(^\psi\mathcal{A}^2_1(\Pi,\mathbb{H}(\C))$
can be identified as the closed subspace of $L^2(\Pi,\mathbb{H}(\C))$ consisting of functions $f=f(x)=f(x_0,x')$ such that
\begin{equation}\label{vasilqekv31}
(I\otimes F)^\psi D(I\otimes F^{-1})f=\left(\partial_{x_0} +i\sum_{k=1}^3 \psi^k x_k\right)f=0
\end{equation}
The general solution to Eqn.(\ref{vasilqekv31}) is given by
\begin{equation}
\varphi(x_0,x')=\exp(-\abs{x'}\cdot x_0)\frac{1}{2}\left(1+i\sum_{k=1}^3 \psi^k \frac{x_k}{\abs{x'}}\right)h(x')
\end{equation}
where for $h(x'),$ the condition $\varphi\in L^2(\Pi,\mathbb{H}(\C))$ implies that $^\psi \mathcal{A}_1^2$ consists of the functions
\begin{equation}
f(x_0,x')=\exp(-\abs{x'}\cdot x_0)\frac{1}{2}\left(1+i\sum_{k=1}^3 \psi^k \frac{x_k}{\abs{x'}}\right)\sqrt{2\abs{x'}}g(x')
\end{equation}
for $g(x')\in L^2(\R^3,\mathbb{H}(\C))$ and furthermore
\begin{equation}
\norm{f(x_0,x')}_{^\psi \mathcal{A}^2_1}=\norm{(^\psi Q_+ g)(x')}_{L^2(\R^3,\mathbb{H}(\C))}
\end{equation}
Define the unitary operator acting on $L^2(\Pi,\mathbb{H}(\C))$ according to
\begin{equation}
U_2: f(x_0,x')\mapsto \frac{1}{2}f\left(\frac{x_0}{2\abs{x'}},x'\right)
\end{equation}
with inverse given by $U_2^{-1}=U_2^*$
\begin{equation}
U_2^{-1}: f(x_0,x')\mapsto \sqrt{2\abs{x'}} f(2\abs{x'}\cdot x_0,x')
\end{equation}
Then the space $^\psi\mathcal{A}_2^2=U_2(^\psi\mathcal{A}_1^2)$ consists of functions $f$ of the form
\begin{multline}
f(x_0,x')=U_2\left( \exp(-\abs{x'}\cdot x_0)\frac{1}{2}\left(1+i\sum_{k=1}^3 \psi^k \frac{x_k}{\abs{x'}}\right)\sqrt{2\abs{x'}}g(x')\right)\\
=\exp\left(-\frac{x_0}{2}\right)\frac{1}{2}\left(1+i\sum_{k=1}^3 \psi^k \frac{x_k}{\abs{x'}}\right)g(x')
\end{multline}
Let $\ell_0(x_0):=\exp\left(-\frac{x_0}{2}\right)\in L^2(\R_+)$, which satisfies $\norm{\ell_0(x_0)}=1.$
Denote by $L_0$ the one-dimensional subspace of $L^2(\R_+)$ generated by $\ell_0(x_0)$. The one-dimensional projection $P_0$ of $L^2(\R_+)$
onto $L_0$ is given by
\begin{equation}
(P_0\psi)(x_0)=\langle \psi,\ell_0\rangle\cdot\ell_0=\exp\left(-\frac{x_0}{2}\right)\int_{\R_+}\psi(t)\exp\left(-\frac{t}{2}\right)dt
\end{equation}
The orthogonal projection $^\psi B_2$ of $L^2(\Pi,\mathbb{H}(\C))$ onto $^\psi\mathcal{A}_2^2$ is given by
\begin{equation}
^\psi B_2 =P_0\otimes ^\psi Q_+
\end{equation}
Denote by $^\psi B_\Pi$ the Bergman orthogonal projection of 
$L^2(\Pi,\mathbb{H}(\C))$ onto $^\psi \mathcal{A}^2_1(\Pi,\mathbb{H}(\C))$.
We may gather the above remarks into the following result.
\begin{theorem}
	The unitary operator $U:=U_2 U_1$ gives an isometric isomorphism of $L^2(\Pi,\mathbb{H}(\C))=L^2(\R_+)\otimes L^2(\R^3,\mathbb{H}(\C))$
	under which:\\
	(1) $^\psi\mathcal{A}_1^2(\Pi,\mathbb{H}(\C))$ is mapped onto $L_0\otimes \im ^\psi Q_+$\\
	(2) $^\psi B_\Pi$ is unitarily equivalent to $P_0\otimes ^\psi Q_+$according to
	\begin{equation}
	U ^\psi B_\Pi U^{-1} =P_0\otimes ^\psi Q_+
	\end{equation}
\end{theorem}

\begin{theorem}\label{vasilqthm51}
	The unitary operator $U$ acting on $L^2(\R_+)\otimes L^2(\R^3,\mathbb{H}(\C))$
	maps $^\psi \mathcal{A}_n^2(\Pi,\mathbb{H}(\C))$ of
	$n$-$\psi$-hyperholomorphic functions, onto the space $L_{n-1}^\oplus \otimes ^\psi L_2^+$
	(recall that $^\psi L_2^+=^\psi L_2^+(\R^3,\mathbb{H}(\C))$ denotes the image of the projection $^\psi Q_+$).
\end{theorem}
\begin{proof}
	The space $U(^\psi\mathcal{A}_n^2(\Pi,\mathbb{H}(\C))$ coincides with the set of functions $f\in L^2(\Pi,\mathbb{H}(\C))=
	L^2(\R_+)\otimes L^2(\R^3,\mathbb{H}(\C))$ satisfying
	\begin{multline}\label{vassillekv}
	U\left( \sum_{k=1}^3 \psi^k \partial_{x_k}\right)^n U^{-1} f=U_2
	\left(\partial_{x_0} \sum_{k=1}^3 \psi^k x_k\right)^n U_2^{-1} f=\\
	\abs{x'}^n \left(2\partial_{x_0} +i\sum_{k=1}^3 \psi^k \frac{x_k}{\abs{x'}}\right)^n f=0
	\end{multline}
	Now the intersection of the general solution of Eqn.(\ref{vassillekv}) with the space 
	$L^2(\R_+)\otimes L^2(\R^3,\mathbb{H}(\C))$ consists of fall functions of the form
	\begin{equation}\label{vassillekv0}
	\sum_{m=1}^3 x_0^m \exp\left(-\frac{x_0}{2}\right)\frac{1}{2} \left(1+i\sum_{m=1}^3 \psi^k \frac{x_k}{\abs{x'}}\right)h_m(x')
	\end{equation}
	where $h_m(x')\in L^2(\R^3,\mathbb{H}(\C)$, $m=0,\ldots,n-1.$
	By reordering the polynomials with respect to $x_0$ these functions in Eqn.(\ref{vassillekv0}) take the form
	\begin{multline}
	\sum_{m=1}^3 L_m(x_0) \exp\left(-\frac{x_0}{2}\right)\frac{1}{2} \left(1+i\sum_{m=1}^3 \psi^k \frac{x_k}{\abs{x'}}\right)g_m(x')=\\
	\sum_{m=1}^3 \ell_m(x_0)\frac{1}{2} \left(1+i\sum_{m=1}^3 \psi^k \frac{x_k}{\abs{x'}}\right)g_m(x')
	\end{multline}
	for $g_m(x')\in L^2(\R^3,\mathbb{H}(\C))$, $m=0,\ldots,n-1.$ This completes the proof.
\end{proof}
Now the space 
$^\psi\mathcal{A}^2_{(n)}(\Pi,\mathbb{H}(\C)=^\psi\mathcal{A}^2_{n}(\Pi,\mathbb{H}(\C)\ominus ^\psi\mathcal{A}^2_{n-1}(\Pi,\mathbb{H}(\C)$ 
of true-$n$-$\psi$-hyperholomorphic functions
satisfies
\begin{equation}
^\psi\mathcal{A}^2_{n}(\Pi,\mathbb{H}(\C)=\oplus_{k=1}^n {}^\psi\mathcal{A}^2_{(k)}(\Pi,\mathbb{H}(\C)
\end{equation}
and we have the following.
\begin{corollary}
	The unitary operator $U$ maps $^\psi\mathcal{A}^2_{(n)}(\Pi,\mathbb{H}(\C)$ onto
	$L_{n-1} \otimes ^\psi L_2^+.$
\end{corollary}

\begin{theorem}
	The unitary operator $U$ acting on $L^2(\R_+)\otimes L^2(\R^3,\mathbb{H}(\C))$
	maps ${}^{\bar{\psi}} \mathcal{A}_n^2(\Pi,\mathbb{H}(\C))$ onto the space $L_{n-1}^\oplus \otimes^\psi L_2^-$
	(where ${}^\psi L_2^-=^\psi L_2^+(\R^3,\mathbb{H}(\C))$ denotes the image of the projection $^\psi Q_-$).
\end{theorem}
\begin{proof}
	We may follow the proof of Theorem \ref{vasilqthm51}.
	The space $U(^{\bar{\psi}}\mathcal{A}_n^2(\Pi,\mathbb{H}(\C))$ coincides with the set of functions $f\in L^2(\Pi,\mathbb{H}(\C))=
	L^2(\R_+)\otimes L^2(\R^3,\mathbb{H}(\C))$ satisfying
	\begin{multline}\label{vassillekvhobo}
	U\left( \sum_{k=1}^3 \bar{\psi}^k \partial_{x_k}\right)^n U^{-1} f=U_2
	\left(\partial_{x_0} -i\sum_{k=1}^3 \psi^k x_k\right)^n U_2^{-1} f=\\
	\abs{x'}^n \left(2\partial_{x_0} -i\sum_{k=1}^3 \psi^k \frac{x_k}{\abs{x'}}\right)^n f=0
	\end{multline}
	The intersection of the general solution of Eqn.(\ref{vassillekvhobo}) with the space 
	$L^2(\R_+)\otimes L^2(\R^3,\mathbb{H}(\C))$ consists of fall functions of the form
	\begin{equation}\label{vassillekv0hobo}
	\sum_{m=1}^3 x_0^m \exp\left(-\frac{x_0}{2}\right)\frac{1}{2} \left(1-i\sum_{m=1}^3 \psi^k \frac{x_k}{\abs{x'}}\right)h_m(x')
	\end{equation}
	where $h_m(x')\in L^2(\R^3,\mathbb{H}(\C)$, $m=0,\ldots,n-1.$
	By reordering the polynomials with respect to $x_0$ these functions in Eqn.(\ref{vassillekv0}) take the form
	\begin{multline}
	\sum_{m=1}^3 L_m(x_0) \exp\left(-\frac{x_0}{2}\right)\frac{1}{2} \left(1-i\sum_{m=1}^3 \psi^k \frac{x_k}{\abs{x'}}\right)g_m(x')=\\
	\sum_{m=1}^3 \ell_m(x_0)\frac{1}{2} \left(1-i\sum_{m=1}^3 \psi^k \frac{x_k}{\abs{x'}}\right)g_m(x')
	\end{multline}
	for $g_m(x')\in L^2(\R^3,\mathbb{H}(\C))$, $m=0,\ldots,n-1.$ This completes the proof.
\end{proof}
\begin{corollary}
	The unitary operator $U$ maps $^{\bar{\psi}}\mathcal{A}^2_{(n)}(\Pi,\mathbb{H}(\C))$ onto
	$L_{n-1} \otimes ^\psi L_2^-.$
\end{corollary}
As a consequence to what we have established in this section, we have in particular the following special case of the result given by Vasilevski \cite{vasilevskiquaternion}.
\begin{theorem}[Vasilevski \cite{vasilevskiquaternion}]
	It holds true that
	\begin{equation}
	\begin{split}
	L^2(\Pi,\mathbb{H}(\C)) & = \oplus_{j=1}^\infty \left({}^\psi\mathcal{A}^2_{(j)}(\Pi,\mathbb{H}(\C))\oplus
	{}^{\bar{\psi}}\mathcal{A}^2_{(j)}(\Pi,\mathbb{H}(\C))\right)
	\\
	& = \left(\oplus_{j=1}^\infty {}^\psi\mathcal{A}^2_{(j)}(\Pi,\mathbb{H}(\C))\right)\oplus\left( \oplus_{j=1}^\infty
	{}^{\bar{\psi}}\mathcal{A}^2_{(j)}(\Pi,\mathbb{H}(\C))\right)
	\end{split}
	\end{equation}
\end{theorem}

\section{Toeplitz operators}\label{cuckovic}
Let $\Omega\subset\C$ be a domain. As we have seen 
$A_q^2(\Omega)$ is  closed in $L^2$ (Theorem \ref{closednessthm}), 
and using the Riesz representation theorem
there is a reproducing kernel $K_{q,\Omega}(z,\zeta)$, for $(z,\zeta)\in \Omega\times \Omega$ such that $f(z)=\langle f,K_{q,\Omega}(z,\cdot)\rangle,$ $z\in \Omega.$
For example in the case of the unit disc $\Omega=\{\abs{z}<1\}$
Eqn.(\ref{koseleveq}) shows that the poly-Bergman kernels are given by
\begin{multline}
K_{q,\Omega}(z,\zeta)=\frac{q}{(1-\zeta\bar{z})^{2n}}\sum_{j=0}^{q-1}(-1)^j\binom{q}{j+1}\binom{q+j}{q}\abs{1-\zeta \bar{z}}^{2(q-1-j)}
\abs{\zeta -z}^{2j}=\\
\frac{q}{(1-\zeta\bar{z})^{2q}}\sum_{j=0}^{q-1}(-1)^j\binom{q}{j+1}\binom{q+j}{q}\abs{\varphi_z(\zeta)}^{2j}
\end{multline}
\begin{definition}[Toeplitz operator with bounded symbol]\index{Toeplitz operator}\label{toeplitzdef}
	Let $\Omega\subset\C$ be a domain with continuous boundary. Let $q\in \Z_+.$
	The Toeplitz operator $T_f=T_{q,f}: A^2_q(\Omega)\to A^2_q(\Omega)$ {\em generated} by the function $f\in L^\infty(\Omega)$  
	is defined as the operator
	\begin{equation}
	T_f g=P(fg),\quad g\in A^2_q(\Omega)
	\end{equation}
	where $P$ denotes the orthogonal projection of $L^2(\Omega)$ onto $A^2_q(\Omega).$
	Sometimes the function $f$ is called the {\em symbol} ot the Toeplitz operator.\index{Symbol of a Toeplitz operator} When it is clear from the context what $q$ is, we simply use the notation $T_f$.
\end{definition}
We have explicit expressions, for instance for the upper half-plane (see Theorem \ref{ortokerneltheorem}) for the orthogonal projection of
$L^2(\Omega)$ onto $A^2_q(\Omega)$ via the poly-Bergman kernels and the analogous result holds true e.g.\ for the unit disc, i.e.\ the orthogonal projection is given precisely by the Bergman projections, defined via the poly-Bergman kernels.
For example we have that for a given $f\in C^0(\overline{\Omega})$ 
the Toeplitz operator with respect to $A^2_1(\Omega)$ is given by
\begin{equation}
T_f g(z)=\int_\Omega K_{q,\Omega}(z,\zeta)f(\zeta)g(\zeta),\quad g\in A^2_q(\Omega)
\end{equation}
For the unit disc, $\Omega=\{\abs{z}<1\}$, this gives for example 
\begin{equation}
(T_f g)(z):=\frac{1}{\pi}\int_\Omega \frac{f(\zeta)g(\zeta)}{(1-z\bar{\zeta})^2}\mu(\zeta)
\end{equation}
Clearly, $T_f$ is bounded and $\norm{T_f}\leq \norm{f}_\infty$
and
$T_{c_1f_1+c_2f_2}=c_1T_{f_1}+c_2T_{f_2}$
for constants $c_1,c_2,$ and functions $f_1,f_2\in L^\infty.$ 
Furthermore, $T^*_f =T_{\bar{f}}.$
Also if $\Omega$ is the unit disc and $T_f=0$ then for each $g\in A^2_q(\Omega)$ we have $fg\in (A^2_q(\Omega))^\perp$, in particular for each pair
$n,m\in \N$ we have
\begin{multline}
0=\langle fz^n,z^m\rangle -\int_\Omega f(z)z^n\bar{z}^m d\mu(z)=\int_\Omega f(z)\overline{z^m\bar{z}^n} d\mu(z)\\
=\langle f,z^nz^m\rangle_{L^2(\Omega)}
\end{multline}
and since the linear span of the monomials $z^n\bar{z}^m$ is dense in $L^2(\Omega)$ 
we must have $f=0$ a.e.\ (and obviously the converse also holds true, i.e.\ we may conclude $f=0$ a.e.\ iff $T_f=0$).

\subsection{Characterization of Toeplitz operators on the upper half plane in the case of symbols independent of $\re z$}
\begin{definition}
		Let $f$ be a functions on the upper half plane satisfying $f(z)=f(y)$.
		The Toeplitz operator $T_{n,f}$ ($T_{(n),f}$) generated by (or with symbol) $f$ acting on $\mathcal{A}_n^2(H^+)$ ($\mathcal{A}_{(n)}^2(H^+)$) is defined as the operator
		\begin{equation}
		T_{n,f}:\mathcal{A}_n^2(H^+)\ni g\mapsto B_{H^+,n}(fg)\in \mathcal{A}_n^2(H^+)
		\end{equation}
		\begin{equation}
		T_{(n),f}:\mathcal{A}_{(n)}^2(H^+)\ni g\mapsto B_{H^+,(n)}(fg)\in \mathcal{A}_{(n)}^2(H^+)
		\end{equation}
	\end{definition}

Vasilevski (see \cite{vasilevskibok}, Thm. 10.4.8, p.254) proved the following.
\begin{theorem}\label{vasilgammathm}
	For any $f(y)\in L^\infty(H^+)$ the Toeplitz operator $T_{f}$ acting on 
$\mathcal{A}_{1}^2(H^+)$ is unitarily equivalent to the multiplication operator
$\gamma_f(x)I=R_0 T_f R_0^*$ acting on $L^2(\R_+)$ where
\begin{equation}
\gamma_f(x):=\int_{\R^+} f\left(\frac{y}{2\abs{x}}\right)\exp(-y)dy,\quad x\in \R_+
\end{equation}
\end{theorem}
\begin{proof}
We calculate 
	\begin{multline}
R_{0}T_{f}R_{0}^*=
R_{0} B_{H^+}f B_{H^+}R_{0}^* =\\
R_{0} R_{0}^* R_{0} f R_{0}^* R_{0} R_{0}^*= 
(R_{0} R_{0}^*) R_{0} f R_{0}^* (R_{0} R_{0}^*)=\\
R_{0} f R_{0}^*= R^*_{0} U_2 U_1 f(y) U_1^{-1} U_2^{-1}  R_{0}=\\
R^*_{0} U_2  f(y) U_2^{-1}  R_{0} 
\end{multline}
Now
\begin{equation}
U_2  f(y) U_2^{-1}= f\left(\frac{y}{2\abs{x}}\right)
\end{equation}
So it remains to verify that $\gamma_f(x)I$ equals
\begin{equation}
R^*_{0}   f\left(\frac{y}{2\abs{x}}\right)  R_{0}
\end{equation}
given that $f=f(y).$
Recall that $\ell_j(y)=(-1)^j\exp\left(-\frac{y}{2}\right)L_j(y),$ so in particular
$\ell_0=\exp\left(-\frac{y}{2}\right)$. Furthermore, recall that
\begin{equation}
(R_0 \phi)(x,y)=\chi_+(x)\phi(x)\ell_0(y)
\end{equation}
and 
\begin{equation}
(R_0^* f)(x)=\chi_+(x)\int_{\R_+}f(x,y)\ell_0(\eta)d\eta
\end{equation}
Hence
\begin{multline}
\left(R^*_{0}   f\left(\frac{y}{2\abs{x}}\right)  R_{0}\phi\right)(x)=\int_{\R_+}  f\left(\frac{y}{2\abs{x}}\right)\phi(x)(\ell_0(y))^2 dy=\gamma_f(x)\cdot \phi(x)
\end{multline}
This completes the proof.
\end{proof}
A similar result holds true for the unit disc, see Vasilevski \cite{vasilevskibok}. 
Ramírez-Ortega \& S\'anchez-Nungaray \cite{ortegasanchez} define
\begin{equation}
R_{0,(n)}:L^2(\R_+)\longrightarrow L^2(\R_+)
\end{equation}
\begin{equation}
(R_{0,(n)}f)(x,y)=\chi_+(x)f(x)\ell_{n-1}(y)
\end{equation}
which has adjoint
\begin{equation}
(R_{0,(n)}^*\phi)(x)=\chi_+(x) \int_{\R_+}\phi(x,t)\ell_{n-1}(t)dt
\end{equation}
Since the image $U(A_{(n)}^2(H^+)=L^2(\R_+)\otimes L_{n-1}$ we have
\begin{equation}
R_{0,(n)}^*R_{0,(n)}=I:L^2(\R_+)\longrightarrow L^2(\R_+)
\end{equation}
and 
\begin{equation}
R_{0,(n)}R^*_{0,(n)}=\chi_+(x)I\otimes P_{(n-1)}:L^2(H^+)\longrightarrow L^2(\R_+)\otimes L_{n-1}
\end{equation}
Set also
\begin{equation}
R_{(n)}=R_{0,(n)}^* U
\end{equation}
Similarly, define
\begin{equation}
R_{0,n}:(L^2(\R_+))^n\longrightarrow L^2(\R)\otimes L^2(\R_+)
\end{equation}
\begin{equation}
(R_{0,n}f)(x,y)=\sum_{k=1}^n \chi_+(x)f_k(x)\ell_{k-1}(y)\chi_+(x)[N_n(y)]^Tf(x)
\end{equation}
where  and 
\begin{equation}
f=(f_1,\ldots,f_n)^T,\quad N_n(y)=(\ell_0(y),\ldots,\ell_{n-1}(y))^T
\end{equation}
The adjoint 
\begin{equation}
R_{0,n}^*:(L^2(H^+))^n\longrightarrow (L^2(\R_+))^2
\end{equation}
is given by
\begin{multline}
(R_{0,n}^* \phi)(x)=\\
 \left(\chi_+(x)\int_{\R_+} \phi(x,y)\ell_0(y)dy,\ldots,\chi_+(x)\int_{\R_+} \phi(x,y)\ell_{n-1}(y)dy\right)^T=\\
\chi_+(x)\int_{\R_+} \phi(x,y)N_n(y)dy
\end{multline}
Since the image $U(A_n^2(H^+))=L^2(\R_+)\otimes L_{n-1}^\oplus$ we have
\begin{equation}
R_{0,n}^*R_{0,n}=I:(L^2(\R_+))^n\longrightarrow (L^2(\R_+))^n
\end{equation}
\begin{equation}
R_{0,n}R_{0,n}^*=\chi_+ I\otimes P_{n-1}:L^2(H^+)\longrightarrow L^2(\R_+)\otimes L_{n-1}^\oplus
\end{equation}
The operator $R_n:=R^*_{0,n}U$ maps $L^2(H^+)$ onto $(L^2(\R_+))^n$ and its restriction to $A^2_n(H^+)$ is an isometric isomorphism. Furthermore, $R^*_n=U^*R_{0,n}$ is an isometric isomorphism from $(L^2(\R_+))^n$ onto $A_n^2(H^+)$.
This implies
$R_{n}^* R_n =B_{H^+ ,n}:L^2(H^+)\longrightarrow A^2_n(H^+)$ and
$R_{n} R_n^* =I:(L^2(\R_+))^n\longrightarrow (L^2(\R_+))^n$ and
	Ramírez-Ortega \& S\'anchez-Nungaray \cite{ortegasanchez} generalized Theorem \label{vasilgammathm} as follows.
	\begin{theorem}
		Let $H^+$ denote the upper half plane $H^+$. For any $f(y)\in L^\infty(H^+)$ the Toeplitz operator $T_{(n),f}$ acting on 
	$\mathcal{A}_{(n)}^2(H^+)$ is unitarily equivalent to the multiplication operator
	$\gamma_{(n),f}(x)I=R_{(n)} T_f R_{(n)}^*$ acting on $L^2(\R_+)$ where
	\begin{equation}
	\gamma_{(n),f}(x):=\int_{\R^+} f\left(\frac{y}{2\abs{x}}\right)(\ell_{n-1}(y))^2 dy
	\end{equation}
\end{theorem}
\begin{proof}
	We may calculate
	\begin{multline}
	R_{(n)}T_{(n),f}R_{(n)}^*=
	R_{(n)} B_{H^+,(n)}f B_{H^+,(n)}R_{(n)}^* =\\
	R_{(n)} R_{(n)}^* R_{(n)} f R_{(n)}^* R_{(n)} R_{(n)}^*= 
	(R_{(n)} R_{(n)}^*) R_{(n)} f R_{(n)}^* (R_{(n)} R_{(n)}^*)=\\
	R_{(n)} f R_{(n)}^*= R^*_{0,(n)} U_2 U_1 f(y) U_1^{-1} U_2^{-1}  R_{0,(n)}=\\
	R^*_{0,(n)} U_2  f(y) U_2^{-1}  R_{0,(n)}=R^*_{0,(n)}   f\left(\frac{y}{2\abs{x}}\right)  R_{0,(n)}
		\end{multline}
		which combined with
		\begin{equation}
		\left(R^*_{0,(n)}   f\left(\frac{y}{2\abs{x}}\right)  R_{0,(n)} \phi\right)(x)=\int_{\R_+} f\left(\frac{y}{2\abs{x}}\right)
				\phi(x)(\ell_{n-1}(y))^2dy=\gamma_{(n),f}(x)\cdot f(x)
				\end{equation}
				completes the proof.
			\end{proof}

			\begin{theorem}
				For any $f(y)\in L^\infty(H^+)$ the Toeplitz operator $T_{n,f}$ acting on 
			$\mathcal{A}_{n}^2(H^+)$ is unitarily equivalent to the multiplication operator
			$\gamma^{n,f}(x)I=R_{n} T_{n,f} R_{n}^*$ acting on $(L^2(\R_+))^n$ where
			the matrix-valued function $\gamma^{n,f}(x)I=(\gamma_{ij}^{n,f})$ is given by
			\begin{equation}
			\gamma^{n,f}(x):=\int_{\R^+} f\left(\frac{y}{2\abs{x}}\right) N_n(y)(N_n(y))^T dy
			\end{equation}
			that is for $i,j=1,\ldots,n,$
			\begin{equation}
			\gamma_{ij}^{n,f}(x):=\int_{\R^+} f\left(\frac{y}{2\abs{x}}\right)\ell_{i-1}(y)\ell_{j-1}(y)dy
			\end{equation}
		\end{theorem}
		\begin{proof}
			We may calculate
			\begin{multline}
			R_{n}T_{n,f}R_{n}^*=R_{n}B_{H^+,n} f B_{H^+,n} R_{n}^*=
			R_{n}R_n^* R_n f  R_{n}^* R_n R_{n}^*=\\
			R_n f  R_{n}^*=R^*_{0,n} U_2 U_1 f U_1^{-1} U_2^{-1}  R_{0,n}=R_n f  R_{n}^*=R^*_{0,n} U_2 f U_2^{-1}  R_{0,n}=
			R^*_{0,n}  f\left(\frac{y}{2\abs{x}}\right)   R_{0,n}
			\end{multline}
			For $g=(g_1,\ldots,g_n)^T\in (L^2(\R_+))^n$
			\begin{multline}
			\left(
			R^*_{0,n}  f\left(\frac{y}{2\abs{x}}\right)   R_{0,n} g\right)(x)=R_{0,n}^*\left( f\left(\frac{y}{2\abs{x}}\right)
			\chi_x(x)(N_n(y))^T g(x)\right)=\\
			\chi_+(x)\int_{R_+}\left(
			f\left(\frac{y}{2\abs{x}}\right) (N_n(y))^T g(x)\right)N_n(y)dy=\\
			\chi_+(x)\int_{R_+}\left(
			f\left(\frac{y}{2\abs{x}}\right) N_n(y)(N_n(y))^T g(x)\right)g(y)dy
			=\chi_+(x)\gamma^{n,f}(x)g(x)
			\end{multline}
			and the components of $\gamma^{n,f}$ can now be verified to be given by 
			\begin{equation}
			\int_{\R^+} f\left(\frac{y}{2\abs{x}}\right)\ell_{i-1}(y)\ell_{j-1}(y)dy
			\end{equation}
			This completes the proof.
	\end{proof}

\subsection{A conditions for Fredholmness of Toeplitz operators}
Let $q\in \Z_+$. We denote for a continuous $f\in C^0(\{\abs{z}=1\})$
the associated Toeplitz operator $T_f=T_{q,f}$ with respect to $A^2_q(\{\abs{z}=1\})$, i.e.\ we drop the subindex for dependence of $q$ in $T_f$.
We may define the so-called {\em Toeplitz algebra}, $\mathcal{T},$ as $\{T_f +K:f\in C^0(\{\abs{z}=1\}),K-\mbox{compact}\}.$
\begin{proposition}\label{windingwolfprop}
	If $T_f$ belongs to the Toeplitz algebra $\mathcal{T}$ then $\mbox{ind}T_f=-\mbox{wind}(f)$
	where $\mbox{wind}(f)$ denotes the winding number of $f\in C^0(\abs{z}=1)$ with respect to the origin.
\end{proposition}
\begin{proof} 
	If $R,S$ are bounded operators from a Hilbert space to itself then 
	$\mbox{ind}(S+K)=\mbox{ind}S$ for compact operators $K$, thus
	it suffices to prove the result for $f$ nowhere vanishing on the unit circle. 
	Such $f$ is homotopically equivalent to
	$z^n$ for some $n$.
	If $\{S_t\}$ is a path of Fredholm operators then 
	$\mbox{ind}(S_0)=\mbox{ind}(S_1)$ 
	thus by the homotopy invariance of the winding number
	it suffices to prove the result for $T_{z^n}.$ Now for a Fredholm operator $S$ we have $\mbox{ind}S^*=\mbox{ind}S.$ 
	Using $\mbox{ind}(R+S)=\mbox{ind}(R)+\mbox{ind}(S)$ 
	together with the algebraic properties of the index it suffices to prove the result for
	$T_z$ which is injective, i.e.\ $\mbox{dim ker}T_z=\{0\}.$ Also $\mbox{ker}T_z^*=\mbox{ker}T_{z^{-1}}$ is spanned by the constant $1.$
	Thus $\mbox{dim ker}T_z^*=1.$ Hence
	$\mbox{ind}T_z=0-1=-1=\mbox{wind}(z).$
	If $S$ is bounded operator from a Hilbert space to itself then $S$ is Fredholm with index zero. This completes the proof.
\end{proof}

Let $a=a(\exp(i\theta))$ be a continuous function on the unit circle with Fourier series
$a(\exp(i\theta))\sum_{j\in \Z} a_j \exp(ij\theta).$
If $a\neq 0$ everywhere then the so-called Wiener-Hopf operator\index{Wiener-Hopf operator}
\begin{equation}
W(a):\ell^2(\N)\to \ell^2(\N)
\end{equation}
whose matrix representation in the basis $e_k:=(\delta_{jk})_{j=0}^\infty$
is $(a_{j-k})_{j,k=0}^\infty$ is Fredholm
 and its index is
\begin{equation}
\mbox{ind}W(a)=-\mbox{wind}(a)
\end{equation}
where $\mbox{wind}(a)$ denotes the winding number of $a$ about the origin.
If in addition this winding number is zero 
then (see Wolf \cite{wolf}) the equation
\begin{equation}
W(a)\phi =g,\quad \phi,g\in \ell^2(\N)
\end{equation}
has a unique solution $\phi$ for each given $f$. Furthermore, denoting by $P_n$ the orthogonal projection of $\ell^2(\N)$ onto the linear span of the first $n+1$ vectors $e_j$, $j=0,\ldots,n$,
the equations
\begin{equation}
P_n W(a)P_n\phi^{(n)} =P_n f =f,\quad \phi,f\in \ell^2(\N)
\end{equation}
are uniquely solvable for all sufficiently large $n$ and $\phi^{(n)}$ goes to 
$\phi$ in the $\ell_2(\N)$-norm as $n$ approaches infinity.
We shall denote the unit disc by $D=\{\abs{z}<1\}$ and as before $A^2_q(D)$ 
denotes the space of functions $q$-analytic on the disc
with $\int_D \abs{f(z)}^2 \frac{1}{\pi} d\mu(z),$ where $\mu(z)$ is the standard area measure. An orthogonal basis is given by 
Eqn.(\ref{koselevbasekv}) as we have pointed out.
\begin{definition}
	\index{Asymptotic invertability of a Toeplitz operator}
	A Toeplitz operator $T_f :A^2_q(D)\to A^2_q(D)$ generated by a continuous function $f\in C^0(\overline{D})$
	is called {\em asymptotically invertible} if\\
	(i) $T_f$ is invertible on $A^2_q(D).$\\
	(ii) Each $P_j T_f P_j$ restricted to the image $\mbox{im}P_j$ are invertible for all sufficiently large $j.$\\
	(iii) $(P_j T_f P_j)^{-1} P_j\to T^{-1}_f$ strongly on $A^2_q(D).$
	\end{definition}
	In particular, if $T_f$ is asymptotically invertible then the equation $T_f \phi=g$ for $\phi, g\in A^2_q(D)$
	has a unique solution and the unique solution of the truncated equation
	\begin{equation}
	P_j T_f P_j\phi^{(j)} =P_j g
	\end{equation}
	is in some sense arbitrarily close to the exact solution $\phi$ for sufficiently large $j.$
Denote by $H_j$ the closed linear hull of $\{e_{k,j}:k\in \N\},$ $j=0,\ldots,q-1.$
Then
\begin{equation}
A_q^2(D)=\oplus_{j=0}^{q-1} H_j
\end{equation}
and accordingly each $g\in A^2_q(D)$ can be represented by a vector $g=(g_0,\ldots,g_{q-1})$ where each $g_j\in H_j.$
Let $P_j$ be the orthogonal projection of $A^2_q(D)$ onto $H_j$, $j=0,\ldots,H_{q-1}.$
This decomposes
\begin{equation}
T_f=\sum_{j,k=0}^{q-1} P_k T_f P_j
\end{equation}
where $P_k T_f P_j :H_j\to H_k.$ For $f\in C^0(\overline{D})$ we have
the Fourier coefficients of the restriction to the unit circle given by
\begin{equation}
g_j=\frac{1}{2\pi}\int_0^{2\pi} g(\exp(i\theta))\exp(-iq\theta))d\theta
\end{equation}
Let $W(f)$ be the operator whose matrix representation on $\ell^2(\N)$ is given by 
$(f_{j-k})_{j,k=0}^\infty$. Then $W(f)$ is a discrete Wiener-Hopf operator with symbol $f|_{\{\abs{z}=1\}}$ on $\ell^2(\N)$ thus a bounded operator.
	Denote by $U_j:H_j\to \ell^2(N)$ the unitary operator that sends $e_j$ to $e_k$ where
	$(e_j)_j$ is the standard basis for $\ell^2(\N).$ 	
	\begin{theorem}\label{wolfthm31}
		Let $f\in C^0(\overline{D})$. Then the following operator is compact
		\begin{equation}
		K:=T_f-\sum_{j=0}^{q-1} P_j U_j^* W(f) U_j P_j
		\end{equation}
	\end{theorem}
	\begin{proof}
		Denoting by $\delta_{jk}$ denote the Kroenecker function we may write
		\begin{equation}
		K:=\sum_{j,k=0}^{q-1} P_k (T_f-\delta_{jk} U_k^* W(f) U_j) P_j
		\end{equation}
		Thus it suffices to show that each of the following operators are compact.
		\begin{equation}
		P_k (T_f-\delta_{jk} U_k^* W(f) U_j) P_j
		\end{equation}
		Assume this holds true for the case of $f(z)=z^l\bar{z}^n.$ 
		Choose a sequence of trigonometric polynomials $\{f^{(n)}\}_{n=0}^\infty$ with $f^{(n)}\to f$
		uniformly on $\overline{D}.$
Now
\begin{equation}
\norm{P_k T_f P_j}\leq \norm{T_f}\leq \norm{a}_{L^\infty}(\overline{D})
\end{equation}
and
\begin{equation}
\norm{P_k U_k^* W(f)U_j P_j}\leq \norm{W(f)}\leq \norm{a}_{L^\infty}(\partial{D})
\end{equation}
Thus $P_kT_{f^{(n)}} P_j\to P_k T_f P_j$ and as the uniform convergence of $f^{(n)}$ to $f$
implies the uniform convergence of the restrictions to the unit circle,
$P_k U_k^* W(f^{(n)}) U_j P_j\to P_k U_k^* W(f) U_j P_j$.
However, we know that each
\begin{equation}
P_k (T_{f^{(n)}}-\delta_{jk} U_k^* W(f^{(n)}) U_j) P_j
\end{equation}
is compact hence
\begin{equation}
P_k (T_f-\delta_{jk} U_k^* W(f) U_j) P_j
\end{equation}
is also compact since it is the uniform limit of compact operators.
Let $(a_{j,k}^{\nu,\mu})_{j,k=0}^\infty$ denote the matrix representation of the
Toeplitz operator $P_\mu T_f P_\nu$ in the basis $e_{j,k}$ of Eqn.(\ref{koselevbasekv}) which is an operator
$H_\nu\to H_\mu$ 
\begin{multline}
a_{j,k}^{\nu,\mu}=(P_\mu T_{z^l\bar{z}^n} P_\nu e_{k,\nu},e_{j,\mu})=
(T_{z^l\bar{z}^n} e_{k,\nu},e_{j,\mu})=\\
(P(z^l\bar{z}^n e_{k,\nu}),e_{j,\mu})=(z^l\bar{z}^n e_{k,\nu},e_{j,\mu})=\\
\sqrt{(k+\nu+1)(j+\mu+1)}\sum_{i,s}(-1)^{i+s}\binom{\nu}{i}\binom{\mu}{s}\times\\
\binom{k+\nu +i}{\nu}\binom{j+\mu-s}{\mu} (z^{k+i+\mu-i-s},z^{j+n+\nu-i-s})
\end{multline}
where
\begin{equation}
(z^{k+i+\mu-i-s},z^{j+n+\nu-i-s})=\delta_{k+i+\mu,j+m+\nu}\frac{1}{j+m+\nu-i-s+1}
\end{equation}
which yields with the index transformation $i=\nu-i$ and $s=\mu-s$ in the sum
\begin{equation}
a_{j,k}^{\nu,\mu}=\left\{
\begin{array}{l}
0 \qquad , k+l+\mu \neq j+m+\nu\\
\sqrt{(k+\nu+1)(j+\mu+1)}(-1)^{\mu+\nu}\sum_{i,s}(-1)^{i+s}\binom{\nu}{i}\binom{\mu}{s}\times\\
\binom{k+i}{\nu}\binom{j+s}{\mu} \frac{1}{j+n-\mu+i+s+1} ,\qquad otherwise
\end{array}
\right.
\end{equation}
Next we shall state the following identities the proof of which may be found in Wolf \cite{wolf}, Sec.4,
\begin{equation}\label{wolflem32}
\lim_{j\to \infty} a_{j,k}^{\nu,\mu}=1,\quad k+l+\mu =j+n+\nu,\quad a_{j,k}^{\nu,\mu}=O(j^{-1}),\quad \mu\neq \nu
\end{equation}
By Eqn.(\ref{wolflem32}) the operator
$P_\mu (T_f -\delta_{\nu\mu} U_\mu^* W(f)U_\nu)P_\nu$ is a compact weight shift operator $H_\nu\to H_\mu.$
Thus it is a compact operator on $A_q^2(D).$ This completes the proof.	
	\end{proof}

Note that once the Fredholmness of $T_f$ is concluded 
we have the index $\mbox{ind}T_f =-q\mbox{wind}(f)$ (where the winding number is with respect to the origin).
\begin{theorem}\label{wolfthm33}
	Let $f\in C^0(\overline{D}).$ Then $T_f$ is Fredholm if and only if $f$ does not vanish on the unit circle. 
\end{theorem}
\begin{proof}
	By Theorem \ref{wolfthm31}, $T_f$ is Fredholm if and only if the following operator is Fredholm
	\begin{equation}
	T^0_f:=\sum_{\nu=0}^{q-1} P_\nu U_\nu^* W(f) U_\nu
	\end{equation}
	Now using $P_\nu P_\mu =\delta_{\nu\mu} P_\mu$ we have
	\begin{equation}
	\mbox{ker}T^0_f=\oplus_{\nu=0}^{q-1}\mbox{ker}P_\nu U_\nu^*W(f)U_\nu P_\nu\simeq
	\oplus_{\nu=0}^{q-1}\mbox{ker}W(f)
	\end{equation}
	\begin{equation}
	\mbox{coker}T^0_f=\oplus_{\nu=0}^{q-1}\mbox{coker}P_\nu U_\nu^*W(f)U_\nu P_\nu\simeq
	\oplus_{\nu=0}^{q-1}\mbox{coker}W(f)
	\end{equation}
	Hence $T_f$ is Fredholm if and only if $W(f)$ is Fredholm which in turn is true if and only if $f|_{\abs{z}=1}$ is invertible in the algebra $C^0(\{\abs{z}=1\}).$
	This completes the proof.
\end{proof}
\begin{theorem}
	Let $f\in C^0(\overline{D})$. $T_f$ is asymptotically invertible only if $T_f$ is invertible on $A^2_j(D).$
\end{theorem}
\begin{proof}
	If $T_f$ is invertible then by Theorem \ref{wolfthm33}, $W(f)$ is Fredholm of index zero.
	A consequence of a well-known result called Coburn's theorem (see e.g.\ B\"ottcher \& Silbermann \cite{bottcher}, Cor.2.40, p.71) a Toeplitz operator is invertible if and only if it is Fredholm operator with zero index.
	 This implies that $W(f)$ is invertible. Thus
	$W(a)$ is asymptotically invertible which implies that $T^0_f$ is asymptotically invertible. Since $T_f$ is invertible and by Theorem \ref{wolfthm31} is also a compact perturbation of 
	$T^0_f$ we obtain that $T_f$ is asymptotically invertible. This completes the proof.
\end{proof}
We mention finally, regarding Toeplitz operators, as they appear in the theory of polyanalytic functions, that there exists a field of studies of  $C^*$-algebras of Toeplitz operators on weighted 
Bergman spaces and poly-Bergman spaces, see e.g.\ Vasilevski \cite{vasilevskibok}.
Further recent results on this topic can be found  e.g.\ regarding sectors (see Karlovich \cite{karlovich})
and characterizations for the case of the upper half-plane 
(see Vasilevski \cite{vasilevskibok} for the Bergman space and Ramírez-Ortega \& S\'anchez-Nungaray \cite{ortegasanchez} for the polyanalytic generalization).

\section{Poly-Fock spaces, Bargmann projections and the Berezin transform}\label{vasilevskisec}
In this section we collect short explanations of some miscellaneous notions that appear in the literature
on operator theoretical analysis of polyanalytic functions so that the reader will be able to recognize these notions. We shall in this book not emphasize these notions and we shall not need them in other chapters or sections. These notions are Fock-spaces, Bargmann projections and the Berezin transform.
\subsection{Poly-Fock spaces}
The so-called {\em true-polyanalytic Fock spaces} are basically
just a special case of weighted Bergman spaces which due to their different readership and current popularity often get separate treatment. 
However, as they where recently characterized in a sense by Vasilevski \cite{vasilevski} using the important method of projections, 
we devote a short subsection to this topic.
For holomorphic functions in one variable these Hilbert spaces where studied by Bargmann \cite{bargmann}.
\begin{definition}
	The space $\mathcal{F}_2^1(\Cn)$ is defined as the closed subspace of all entire functions in 
	the weighted space $L^2(\Cn,\pi^{-1}\exp(-\abs{z}^2)d\mu(z)).$ Similarly, the {\em $j$-poly-Fock space}\index{Poly-Fock space} 
	$\mathcal{F}_2^j(\Cn)$ is defined as the space of $q$-analytic functions in $L^2(\Cn,\pi^{-n}\exp(-\abs{z}^2)d\mu(z)).$
\end{definition}
Vasilevski \cite{vasilevski} points out that the true-j-Fock spaces can be written as orthogonal differences
\index{Orthogonal difference}
\begin{equation}
\mathcal{F}^{(j)}_2(\C) =\mathcal{F}^{j}_2(\C) \ominus \mathcal{F}^{j-1}_2(\C), \quad j>1
\end{equation} 
\begin{equation}
\mathcal{F}^{(1)}_2(\C) =\mathcal{F}^{1}_2(\C) 
\end{equation} 
where the $\mathcal{F}^{j}_2(\C)$ are the weighted poly-Fock spaces (i.e.\ without the prefix {\em true}).
Recall that the {\em orthogonal difference} of two closed linear subspaces spaces $X,Y$ of a Hilbert space 
is given by $X\ominus Y:=X\cap Y^{\perp}.$
Vasilevski defines the {\em Bargmann projector}\index{Bargmann projector}
\begin{equation}
P_n:L^2(\Cn)\to \mathcal{F}_2^1(\Cn)
\end{equation} 
\begin{equation}
(P_n f):=\pi^{-n}\int_{\Cn} f(\zeta)\exp(\bar{\zeta}\cdot z)d\mu(\zeta)
\end{equation} 
The decompositions $L^2(\Cn)=L^2(\C^{j})\otimes L^2(\C^{n-j}),$
$\mathcal{F}^1_2(\Cn)=\mathcal{F}^1_2(\C^{j})\otimes \mathcal{F}^1_2(\C^{n-j}),$
induce $P^2(\Cn)=P^2(\C^{j})\otimes P^2(\C^{n-j}).$
Vasilevski defines the unitary operator 
\begin{equation}
U_1\colon L^2(\C)\to L^2(\R^2),\quad U_1\phi(z)=\pi^{-\frac{1}{2}}\exp i(-1/2z\bar{z})\phi(z)
\end{equation}
or $U_1\phi(x,y)=\pi^{-\frac{1}{2}}\exp i(-(x^2+y^2)/2)\phi(x+iy).$
The space $\mathcal{F}_{(1),_j}:=U_1(\mathcal{F}_{2}^j(\C))$ is the set of all smooth $L^2(\R^2)$ functions
satisfying
\begin{equation}
D_j^{(1)} f=U_1 \partial_{\bar{z}}^j U_1^{-1} f=\left( \partial_{\bar{z}} +\frac{z}{2}\right)f=0
\end{equation}
Set 
\begin{equation}
U_2:=I\otimes F,\quad (F\phi)(y):=(2\pi)^{-\frac{1}{2}}\int_\R \exp(-it\cdot y) \phi(t)dt
\end{equation}
Then $U_2$ maps $L^2(\R^2)=L^2(\R,dx)\otimes L^2(\R,dy)$ isometrically onto itself.
The space $\mathcal{F}_{(2),_j}:=(I\otimes F)(\mathcal{F}_{(1),_j})$ is the closure of the set of smooth $L^2$ functions
such that
\begin{equation}
D_j^{(2)} f=(I\otimes F)D_j^{(1)} (I\otimes F^{-1})f=\frac{1}{2^j}\left(\partial_x-y+x-\partial_y\right)^jf=0
\end{equation}
Define also
\begin{equation}
U_3:L^2(\R^2)\to L^2(\R^2),\quad U_3f:=f\left(\frac{x+y}{\sqrt{2}},\frac{x-y}{\sqrt{2}}\right)
\end{equation}
Then $U_3=U_3^*=U_3^{-1},$ and $U_3$ maps $\mathcal{F}_{(2),_j}$ onto $\mathcal{F}_{(3),_j}$ which is the closure
of the set of smooth functions satisfying
\begin{equation}
D_j^{(3)} f=U_3D_j^{(2)} U_3^{-1} f=2^{-k/2}(\partial_y+y)^j f=0
\end{equation}
with general solution given by $f=\exp(-y^2/2)\sum_{k=0}^{j-1} y^k d_k(x),$ 
which can be rewritten as
\begin{equation}
f=\sum_{k=0}^{j-1} (2^k k!\sqrt{\pi})^{-\frac{1}{2}}\tilde{H}_k(y)\exp(-y^2/2) g_k(x)
\end{equation}
where the $\tilde{H}_k(y)$ are the so-called {\em Hermite polynomials} of degree $k$\index{Hermite polynomial}
\begin{equation}
\tilde{H}_k(y):=(-1)^k\exp(y^2)\partial_y^k\exp(-y^2)=k! \sum_{m=0}^{[k/2]} \frac{(-1)^m(2y)^{k-2m}}{m!(k-2m)!}
\end{equation}
$g_k(x),d_k(x)\in L^2(\R,dx).$
The functions
\begin{equation}
h_k(y):=(2^k k!\sqrt{\pi})^{-\frac{1}{2}}\tilde{H}_k(y)\exp(-y^2/2)
\end{equation}
form an orthonormal basis for $L^2(\R).$ Let $H_k$ denote the $1$-dimensional subspace of $L^2(\R)$ generated
by $h_k(y)$ and denote by $Q_k$ the orthogonal projection of $L^2(\R)$ onto $H_k,$ which is given by
\begin{equation}
(Q_k\phi)(y):=h_k(y)\int_\R \phi(t)h_k(t)dt
\end{equation}
Then 
\begin{equation}
\mathcal{F}_{(3),_j}=L^2(\R)\otimes H_j^{\oplus}, \quad H_j^{\oplus}:=\oplus_{k=0}^{j-1} H_k 
\end{equation}
This proves the following.
\begin{theorem}[Vasilevski \cite{vasilevski}]
	The unitary operator $U:=U_3U_2U_1$ is an isometric isomorphism of $L^2(\C)$ onto $L^2(\R,dx)\otimes L^2(\R,dy)$ that maps $\mathcal{F}^j_2(\C)$
	onto $L^2(\R)\otimes H_k^\oplus$.
\end{theorem}
\begin{corollary}
	$U$ maps the true-j-Fock space $\mathcal{F}_2^{(j)}(\C)$ onto $L^2(\R)\times H_{j-1}$ and the true-$j$-Bargmann
	projection $P_{(k)}$ is unitary equivalent to $UP_{(j)}U^{-1}=I\otimes Q_{j-1}.$ Furthermore
	$U(\mathcal{F}^1_2(\C))=L^2(\R)\otimes H_0$ and $UPU^{-1}=I\otimes Q_0.$
\end{corollary}
\begin{corollary}
	$L^2(\C,\pi^{-1}\exp(-z\bar{z})d\mu(z))= \oplus_{j=1}^\infty \mathcal{F}^{(j)}_2(\C).$
\end{corollary}
The operator $(R_j f)(x):=f(x)h_{j-1}(y),$ $L^2(\R,dx)\otimes L^2(\R,dy)$ yields an isometric embedding with adjoint
$R^*_j:L^2(\R^2)\to L^2(\R^2),$ $(R^*_j f)(x)=\int_\R f(x,y)h_{j-1}(y)dy,$ satisfying
$R_j^* R_j=I,$ $L^2(\R,dx)\to L^2(\R,dx),$ and  
$R_j R_j^*=I\otimes Q_j,$ $L^2(\R^2)\to L^2(\R)\otimes H_{j-1}.$
The map $R^*_j U$ maps $L^2(\C,\pi^{-1}\exp(-z\bar{z})d\mu(z))$ onto $L^2(\R,dx)$ and its restriction to $\mathcal{F}_2^{(j)}(\C)$ is an isometric isomorphism
$\mathcal{F}_2^{(j)}(\C)\to L^2(\R,dx).$ Its adjoint is an onto isometric isomorphism $L^2(\R,dx)\to \mathcal{F}_2^{(j)}(\C).$
\begin{theorem}
	$U^*R_j : L^2(\R,dx)\to \mathcal{F}_2^{(j)}(\C)$ is given by
	\begin{equation}
	(U^*R_j f)(z)=(2^{j-1}(j-1)!\sqrt{\pi})^{-\frac{1}{2}}\int_\R \exp(-(z^2+t^2)/2+\sqrt{2}zt)\tilde{H}_{j-1}\left(\frac{z+\bar{z}}{\sqrt{2}}-t\right)f(t)dt
	\end{equation}
\end{theorem}
\begin{proof}
	\begin{multline}
	(U^*R_j f)(z)=(U_1^{-1}U_2^{-1}U_2^{-1} R_j f)(z)=U_1^{-1}(I\otimes F^{-1})f\left(\frac{x+y}{\sqrt{2}}\right)
	h_{j-1}\left(\frac{x+y}{\sqrt{2}}\right)=\\
	\pi^{\frac{1}{2}}\exp((x^2+y^2)/2)(2\pi)^{-\frac{1}{2}}\int_\R \exp(iyt)f\left(\frac{x+t}{\sqrt{2}}\right)
	h_{j-1}\left(\frac{x-t}{\sqrt{2}}\right)dt=\\
	(2^{j-1}(j-1)!\sqrt{\pi})^{-\frac{1}{2}}
	\int_\R \exp\left(\frac{x^2}{2}+\frac{y^2}{2}+iy(\sqrt{2}x-t)-\frac{(\sqrt{2}x-t)^2}{2}\right)\times\\ f(t) \tilde{H}_{j-1}(\sqrt{2}x-t)dt
	=(2^{j-1}(j-1)!\sqrt{\pi})^{-\frac{1}{2}}\times\\
	\int_\R \exp\left(\frac{z^2+t^2}{2}+\sqrt{2}zt\right)
	\tilde{H}_{j-1}\left(\frac{z+\bar{z}}{\sqrt{2}}-t\right)f(t)dt
	\end{multline}
	This completes the proof.
\end{proof}

For example 
\begin{equation}
(U^*R_1 f)(z)=\pi^{-\frac{1}{4}}\int_\R \exp\left(\frac{z+\bar{z}}{\sqrt{2}}+\sqrt{2}zt\right) f(t)dt
\end{equation}
is an isomorphism $L^2(\R,dx)\to F_2^1(\C).$ Note also that 
\begin{equation}
(R_j^*U)(R_j^*U)^*=I,\quad L^2(\R)\to L^2(\R)
\end{equation}

\begin{definition}[True-$j$-Bargmann projection]\index{True-Bargmann projection}
	The {\em true-$j$-Bargmann projection} $P^{(j)}:L^2(\C,\pi^{-1}\exp(-\abs{z}^2)d\mu(z))\to \mathcal{F}_2^{(j)}(\C)$
	is defined as the operator $P^{(j)}:=(U^*R_j)^*(U^*R_j).$
\end{definition}

By definition of the Hermite polynomials $\tilde{H}_k(y)$, Taylor expansion near $y=0$
gives
\begin{equation}\label{generating000}
\exp(2yt-t^2)=\sum_{k=0}^\infty \frac{H_k(y)}{k!}t^k
\end{equation} 
which implies
\begin{equation}\label{analoghermit}
\partial_y H_k(y)=2k H_{k-1}(y),\quad H_k(y)=2H_{k-1}(y)-\partial_y H_{k-1}(y)
\end{equation}

\begin{proposition}
	The family $\{h_k\}_k$ form an orthonormal basis for $L^2(\R,dy).$
\end{proposition}
\begin{proof}
	We do this by first showing that if we set 
	$\tilde{h}_k(y):=\tilde{H}_k(y)\exp(-y^2/2)=h_k(y)(2^k k!\sqrt{\pi})^{\frac{1}{2}},$ then the family
	$\{\tilde{h}_k\}_k$ form an
	orthogonal basis, and further proving the relation
	\begin{equation}\label{jkaspo1}
	\int_\R \tilde{h}_k(y)^2 dy =2^k k!\sqrt{\pi}
	\end{equation}
	First note that 
	\begin{equation}
	\int_\R \tilde{h}_k(y)=\frac{1}{\sqrt{\pi}}(-2i)^k\exp(\frac{y^2}{2})
	\int_\R t^k\exp(-t^2+2iyt)dt 
	\end{equation}
	where one uses $\exp(-y^2)=\frac{1}{\sqrt{\pi}}\int_\R t^k\exp(-t^2+2iyt)dt$.
	This implies that for $\abs{t}<1,$
	\begin{multline}\label{jkaspo}
	\sum_{k=0}^\infty \frac{\tilde{h}_k(x)\tilde{h}_k(y)}{2^k k!}t^k=\\
	\frac{1}{\pi}\exp\left(\frac{x^2+y^2}{2}\right)\sum_{k=0}^\infty \int_\R \int_\R \frac{(-2uvt)^k}{k!}
	\exp\left(-u^2-v^2+2ixu+2ivy\right)dudv=\\
	\frac{1}{\pi}\exp\left(\frac{x^2+y^2}{2}\right)\int_\R \int_\R \exp\left(-u^2-v^2+2ixu+2ivy-2uvt\right)dudv=\\
	\frac{1}{\sqrt{\pi}}\exp\left(\frac{x^2+y^2}{2}\right)\int_\R \exp\left(-(1-t^2)u^2+2i(x-yt)u\right)du=\\
	(1-t^2)^{-\frac{1}{2}}\exp\left(\frac{x^2+y^2}{2}\right)\exp\left(\frac{(x-yt)^2}{1-t^2}\right)
	\end{multline}
	Setting $x=y$ Eqn.(\ref{jkaspo}) gives
	\begin{equation}
	\sum_{k=0}^\infty \frac{\tilde{h}^2_k(y)}{2^k k!}t^k =(1-t^2)^{-\frac{1}{2}}\exp\left(-\frac{1-t}{1+t}y^2\right)
	\end{equation}
	Integrating both sides gives
	\begin{multline}
	\sum_{k=0}^\infty \left(\int_\R \tilde{h}_k(y)^2 dy\right)\frac{t^k}{2^k k!}=\\
	(1-t^2)^{-\frac{1}{2}}\exp\left(-\frac{1-t}{1+t}y^2\right)dy=\sqrt{\pi}(1-t^2)^{-\frac{1}{2}}
	\left(\frac{1+t}{1-t}\right)^{\frac{1}{2}} =\sqrt{\pi}\sum_{k=0}^\infty t^k.
	\end{multline}
	This proves Eqn.(\ref{jkaspo1}).
	Now the analogous relations to Eqn.(\ref{analoghermit}) for the $\tilde{h}_k(y)$ become 
	\begin{equation}
	(\partial_y +y)\tilde{h}_k(y) =2k\tilde{h}_{k-1}(y),\quad (-\partial_y +y)\tilde{h}_k(y) =\tilde{h}_{k+1}(y)
	\end{equation}
	which in turn gives the recursion formula
	\begin{equation}
	2y\tilde{h}_k(y) =\tilde{h}_{k+1}(y)+2k\tilde{h}_{k-1}(y)
	\end{equation}
	These relations imply that
	\begin{equation}
	(\partial_y^2+y^2)(\tilde{h}_k)=(2k+1)\tilde{h}_k
	\end{equation}
	which in turn gives
	\begin{equation}\label{integreraekv}
	2(k-j)\tilde{h}_k(y)\tilde{h}_j(y)=\tilde{h}_k(y)\partial_y^2\tilde{h}_j(y) -\tilde{h}_j(y)\partial_y^2\tilde{h}_k
	\end{equation}
	In particular, the $\tilde{h}_k$ are eigenfunctions of the so-called {\em Hermite operator}\index{Hermite operator} $H:=\partial_y^2+y^2.$
	Integrating the relation in Eqn.(\ref{integreraekv}) shows that the $\tilde{h}_k$ form an orthogonal family in $L^2(\R,dy).$
	Hence by Eqn.(\ref{jkaspo1}) the $h_k$ form an orthogonal family in $L^2(\R,dy).$
\end{proof}
Note that Eqn.(\ref{generating000}) implies
$H_k(-y)=(-1)^kH_k(y)$ so $h_{2j}$ is even and $h_{2j+1}$ is odd. In particular, $h_{2j+1}(0)=0$ and
$h_{2j}(0)$ can be obtained from
\begin{multline}
\sum_{k=0}^\infty \tilde{h}_k(y)^2 t^{2k}=\pi^{-\frac{1}{2}}(1-t^2)^{-\frac{1}{2}}=\\
\sum_{k=0}^\infty \tilde{h}_k(y)^2 t^{2k}=\pi^{-\frac{1}{2}} \sum_{k=0}^\infty
\frac{\Gamma\left(k+\frac{1}{2}\right)}{\Gamma\left(k+1\right)\sqrt{\pi}}t^{2k}
\end{multline}
which gives
\begin{equation}
h_{2j}^2(0)=\frac{1}{\pi}\frac{\Gamma\left(k+\frac{1}{2}\right)}{\Gamma\left(k+1\right)\sqrt{\pi}}
\end{equation}
Set 
\begin{equation}
\mathbf{a}:=\frac{1}{\sqrt{2}}(\partial_y +y),\quad \mathbf{b}:=\frac{1}{\sqrt{2}}(-\partial_y +y)
\end{equation}
\begin{equation}
\mathfrak{a}:=I\otimes \mathbf{a},\quad \mathfrak{b}:=I\otimes\mathbf{b}
\end{equation}
Their action on Hermite functions $h_k$ satisfy
\begin{equation}
\mathbf{a}h_j(y)=\sqrt{j} h_{j-1}(y),\quad \mathbf{b}=\sqrt{j+1} h_{j+1}(y)
\end{equation}
which gives an isometric isomorphisms 
$\frac{1}{\sqrt{j}}\mathfrak{a}|_{\mathcal{F}^{(j+1)}}:\mathcal{F}^{(j+1)}\to \mathcal{F}^{(j)},$
and their inverses
$\frac{1}{\sqrt{j}}\mathfrak{b}|_{\mathcal{F}^{(j)}}:\mathcal{F}^{(j)}\to \mathcal{F}^{(j+1)}.$
This proves the following lemma.
\begin{lemma}\label{lemmapsoas}
	For natural numbers $k,l$ such that $k<l$ we have isometric isomorphisms
	$\sqrt{\frac{(k-1)!}{(l-1)!}}\mathfrak{a}^{l-k}|_{\mathcal{F}^{(l)}}:\mathcal{F}^{(l)}\to \mathcal{F}^{(k)},$
	and their inverses
	$\sqrt{\frac{(k-1)!}{(l-1)!}}\mathfrak{b}^{l-k}|_{\mathcal{F}^{(k)}}:\mathcal{F}^{(k)}\to \mathcal{F}^{(l)}.$
\end{lemma}
Define on the weighted space $L^2(\C,\pi^{-1}\exp(-\abs{z}^2))$ the operators
\begin{equation}
\mathbf{A}:=\partial_z +\bar{z},\quad \mathbf{B}:=\partial_{\bar{z}}
\end{equation}

\begin{theorem}
	For natural numbers $k,l$ such that $k<l$ we have isometric isomorphisms
	$\sqrt{\frac{(k-1)!}{(l-1)!}}\mathbf{A}^{l-k}|_{\mathcal{F}^{(l)}_2(\C)}:\mathcal{F}^{(l)}_2(\C)\to \mathcal{F}^{(k)}_2(\C),$
	and their inverses
	$\sqrt{\frac{(k-1)!}{(l-1)!}}\mathbf{B}^{l-k}|_{\mathcal{F}^{(k)}_2(\C)}:\mathcal{F}^{(k)}_2(\C)\to \mathcal{F}^{(l)}_2(\C).$
\end{theorem}
\begin{proof}
	For $k\in \N$ we have $\mathcal{F}_2^{(k)}(\C)=U^{-1}(\mathcal{F}_2^{(k)}(\C))$ and using
	the relations $\mathbf{A}=U^{-1} \mathfrak{a} U,$ $\mathbf{B}=U^{-1} \mathfrak{b} U,$ in Lemma \ref{lemmapsoas}
	completes the proof.
\end{proof}

We thus obtain, similar to the case of $q$-analytic functions on bounded domains the following
representation theorem, decomposing the function with holomorphic components.
\begin{corollary}\label{psoascorren}
	For any element $f\in \mathcal{F}_2^{(q)}(\C)$ there is a holomorphic $\phi \in \mathcal{F}_2^{(1)}(\C)=\mathcal{F}_2^{1}(\C),$ such that
	$f$ has a representation
	\begin{equation}
	f(z)=\sum_{j=0}^q(-1)^j\frac{\sqrt{(q-1)!}}{j!(q-1-j)!}\bar{z}^{q-1-j} f^{(j)}(z)
	\end{equation}
	where $f^{(j)}$ is the $j$:th derivative of $\phi.$
	Furthermore, $\norm{f}_{\mathcal{F}_2^{(q)}(\C)} =\norm{\phi}_{\mathcal{F}_2^{1}(\C)}.$
\end{corollary}
Define 
\begin{multline}
e_z^j(\zeta):=\frac{1}{(j-1)!}\mathbf{B}_\zeta^{j-1} \overline{\mathbf{B}}_z^{q-1}\exp(\zeta\bar{z})=\\
\frac{1}{(j-1)!}\left( -\partial_\zeta +\bar{\zeta} \right)^{j-1}\left( -\partial_{\bar{z}} +z \right)^{j-1}\exp(\zeta\bar{z}), \quad \zeta\in \C
\end{multline}
Note that
\begin{multline}
g(\zeta):=\frac{1}{(j-1)!}\left(-\partial_{\bar{z}} +z \right)^{j-1}\exp(\zeta\bar{z})=\\
\sum_{k=0}^{j-1} (-1)^{k} \frac{\sqrt{(j-1)!}}{k!(j-1-k)!} z^{j-1-j}\zeta^j \exp(\zeta\bar{z})=
(z-\zeta)^{j-1}\exp(\zeta\bar{z})
\end{multline}
is holomorphic in $\zeta$ and belongs to the appropriately weighted $L^2(\C)$ space thus to $\mathcal{F}_2(\C),$ so
by Lemma \ref{lemmapsoas} 
\begin{equation}
e_z^j(\zeta)=\frac{1}{\sqrt{(j-1)!}}\left(-\partial_{\bar{z}} +z \right)^{j-1}g(\zeta)
\end{equation}
belongs to $\mathcal{F}_2^{(j)}(\C).$
Furthermore, 
\begin{equation}
\left( -\partial_\zeta +\bar{\zeta} \right)\left( -\partial_{\bar{z}}+z\right)\exp(\zeta\bar{z})=\exp(\zeta\bar{z})((z-\zeta)(\bar{\zeta}-\bar{z})+1)
\end{equation}
\begin{multline}
\left( -\partial_\zeta +\bar{\zeta} \right)\left( -\partial_{\bar{z}} +z \right)\exp(\zeta\bar{z})
(z-\zeta)^k(\bar{\zeta}-\bar{z})^k=\\
\exp(\zeta\bar{z})\left((z-\zeta)^{k+1}(\bar{\zeta}-\bar{z})^{k+1}+2(k+1)(z-\bar{\zeta})^k(\bar{\zeta}-\bar{z})^k +k^2(z-\zeta)^{k-1}(\bar{\zeta}-\bar{z})^{k-1}\right)
\end{multline}
Hence there are polynomials $p_{j}$, with real coefficients, such that
\begin{equation}
e_z^j(\zeta)=\exp(\zeta\bar{z}) p_{j-1}((z-\bar{\zeta})(\bar{\zeta}-\bar{z}))
\end{equation}
and 
\begin{equation}
e_z^j(\zeta)=\overline{e_z^j(\zeta)}
\end{equation}

\begin{theorem}[Representation of the orthogonal Bargmann projection]
	The operator
	\begin{equation}
	(P_{(j)} f)(z)=\langle f(\zeta),e_z^j(\zeta)\rangle =\int_\C f(\zeta)e_z^j(\zeta) (\pi^{-1}\exp(-\abs{\zeta})d\mu(\zeta))
	\end{equation}
	is the orthogonal Bargmann projection of $L^2(\C,\pi^{-1}\exp(-\abs{\zeta})d\mu(\zeta))$ onto the
	true-$j$-Fock space $\mathcal{F}_2^{(j)}(\C).$ 
\end{theorem}
\begin{proof}
	By the remark previous to the theorem, $e_z^j(\zeta) \in \mathcal{F}_2(\C),$ hence for $f\in (\mathcal{F}_2(\C))^\perp,$
	we have 
	\begin{equation}
	(P_{(j)} f)(z)=\langle f(\zeta),e_z^j(\zeta)\rangle =0
	\end{equation}
	i.e.\ $(\mathcal{F}_2(\C))^\perp \subset \mbox{Ker}P_{(j)}.$
	Next note that
	\begin{multline}
	(-\partial_z+\bar{z})^{j-1}\exp(-\abs{\zeta})=\sum_{k=0}^{j-1} (-1)^k \frac{(j-1)!}{k!(j-1-k)!}\bar{z}^{j-1-k}\bar{\zeta}^k\exp(-\abs{\zeta})=\\
	(\bar{z}-\bar{\zeta})^{j-1}\exp(-\abs{\zeta})
	\end{multline}
	which can be used to show that
	in Corollary \ref{psoascorren}, we can rewrite
	\begin{multline}
	f(z)=\sum_{k=0}^{j-1} (-1)^k \frac{(j-1)!}{k!(j-1-k)!}\bar{z}^{j-1-k}\phi^{(k)}(z)=\\
	\frac{1}{\sqrt{(j+1)!}}\langle \phi(\zeta),\overline{\mathbf{B}}^{j-1}\exp(-\abs{\zeta})\rangle
	\end{multline}
	Where $\phi\in \mathcal{F}^1_2(\C)$.
	By Corollary \ref{psoascorren} and Lemma \ref{lemmapsoas}, this implies
	\begin{equation}
	(P_{(j)} f)(z)=\langle f(\zeta),e_z^j(\zeta)\rangle = \langle \phi(\zeta),\frac{1}{\sqrt{(j-1)!}}\overline{\mathbf{B}}_z^{j-1}
	e_z^j(\zeta)\rangle =f(z)
	\end{equation}
	Hence $P_{(j)}|_{(\mathcal{F}_2(\C))^\perp}=0$ and $P_{(j)}|_{\mathcal{F}_2(\C)}=I.$
\end{proof}

Now all the arguments for the one-dimensional case have direct natural generalization to the case of $\Cn.$ We shall not give the details.
Let $dm_1:=$\\
$\pi^{-1}\exp(-\abs{\zeta}^2)d\mu(\zeta)$, $\zeta\in \C,$ and more generally $dm:=\pi^{-n}\exp(-\abs{z}^2)d\mu(z)$ where $d\mu$ is the standard measure on $\Cn$.
The poly-Fock spaces $\mathcal{F}_2^j(\Cn)$ are defined via tensor products 
$L^2(\C,dm)=L^2(\C,dm_1)\otimes \cdots \otimes L^2(\C,dm_1),$
$\mathcal{F}_2^1(\Cn):=\mathcal{F}_2^1(\C)\otimes \cdots \otimes \mathcal{F}_2^1(\Cn),$
For multi-indices $\alpha\in \N^n,$ we define the {\em true-$\alpha$-poly Fock spaces over $\Cn$}
\begin{equation}
\mathcal{F}_2^{(\alpha)}(\Cn):=\otimes_{j=1}^n \mathcal{F}_2^{(\alpha_j)}(\Cn)
\end{equation}
Denote by $\mathbf{P}_{(\alpha)}$ the orthogonal projection of $L^2(\Cn,dm)$ onto
$\mathcal{F}_2^{(\alpha)}(\Cn).$
The operators $U_1,U_2,U_3$ can be generalized naturally according to
$U_1:L^2(\C,dm)\to L^2(\R^{2n},\mu(x,y)),$ $(Uf)(z)=\pi^{-\frac{n}{2}}\exp(-\abs{z})f(z),$
$U_2:=I\otimes F,$ $L^2(\R^{2n},\mu(x,y))\to L^2(\R^{n},\mu(x))\otimes L^2(\R^{n},\mu(y)),$
$(Ff)(y):=(2\pi)^{n/2}\int_{\Rn}\exp(-i t\cdot y)f(t)dt,$ and $U_3=U_3^*=U_3^{-1}:L^2(\R^{2n})\to L^2(\R^{2n}),$
$(U_3f)(x,y):=f\left(\frac{x+y}{\sqrt{2}},\frac{x+y}{\sqrt{2}}\right).$
Then the functions $\mathbf{h}_\alpha(y)=\Pi_{k=1}^n h_{\alpha_k}(y_k)$ form an orthonormal
basis in $L^2(\Rn).$
Denote by $\mathbf{H}_\alpha$ the one-dimensional subspace of $L^2$ generated by $\mathbf{h}_\alpha(y).$
Then $\mathbf{Q}_\alpha=\otimes_{k=1}^n Q_{\alpha_k},$ where $Q_{\alpha_k}$ is the one-dimensional orthogonal projection
acting on the $k$:th component
of $L^2(\Rn)$ onto $\mathbf{H}_\alpha$, according to
$(\mathbf{Q}_\alpha f)(y):=\mathbf{h}_\alpha(y)\int_{\Rn}f(t)\mathbf{h}_\alpha(t)dt.$
Then (Vasilevski \cite{vasilevski}, \cite{vasilevski1999}) 
$U:=U_3U_2U_1$ is an isometric isomorphism of $L^2(\Cn,dm)$ onto $L^2(\Rn)$ under which $\mathcal{F}_2^{(\alpha)}(\C)$ is mapped to
onto $L^2(\Rn)\otimes \mathbf{H}_{\alpha-(1,\cdots,1)}$; 
$\mathbf{P}_{(\alpha)}$ is unitary equivalent to $U\mathbf{P}_{(\alpha)}U^{-1} =I\otimes \mathbf{Q}_{\alpha-(1,\cdots,1)},$
and we have the following decomposition
\begin{equation}
L^2(\Cn,dm)=\otimes_{\abs{\alpha}=n}^{\infty} \mathcal{F}_2^{(\alpha)}(\Cn)
\end{equation}
The multi-dimensional analogues of $\mathbf{A}$ and $\mathbf{B}$ are
\begin{equation}
\tilde{\mathbf{A}}:=(\partial_{\bar{z}_1},\ldots,\partial_{\bar{z}_n})
,\quad \tilde{\mathbf{B}}:=(-\partial_{z_1}+z_1,\ldots,-\partial_{z_n}+z_n)
\end{equation}
and we the system defined by
\begin{equation}
\mathbf{e}_z^\alpha(\zeta):=\frac{1}{(\alpha-(1,\ldots,1))!}(\tilde{\mathbf{A}}\tilde{\mathbf{B}})^{\alpha-(1,\ldots,1)}\exp(\zeta\cdot \bar{z}), \zeta\in \C
\end{equation}
parametrized by $z\in \Cn$, forms an orthonormal system in $\mathcal{F}_2^{(\alpha)}(\Cn).$ Furthermore, any function $f$ in 
$\mathcal{F}_2^{(\alpha)}(\Cn)$ can be written in the form
\begin{equation}
f(z)=\sum_{\abs{\beta}<\alpha-(1,\ldots,1)} (-1)^{\abs{\beta}}\frac{\sqrt{(\alpha-(1,\ldots,1))!}}{\beta!(\alpha-\beta-(1,\ldots,1))!}\abs{z}^{\alpha-\beta-(1,\ldots,1)}
\partial^\beta \phi(z)
\end{equation}
where $\phi(z)\in \mathcal{F}_2(\C)$ such that $\norm{f}_{\mathcal{F}_2^{(\alpha)}(\Cn)}=\norm{\phi}_{\mathcal{F}_2(\Cn)}$
The operator
\begin{equation}
(\mathbf{P}_{(\alpha)} f)(z)=\langle f(\zeta),\mathbf{e}_z^\alpha(\zeta)\rangle=\int_{\Cn} f(\zeta)\mathbf{e}_z^\alpha(\zeta)dm(\zeta)
\end{equation}
is the orthogonal Bargmann projection of onto $\mathcal{F}_2^{(\alpha)}(\Cn)$.
The {\em polyanalytic Fock space} $\mathbf{F}_2^{(q)}(\C)$ is defined as
\begin{equation}
\mathbf{F}_2^{(q)}(\C)=\oplus_{j=1}^{q} \mathcal{F}_2^{(j)}(\C)
\end{equation}
The reproducing kernels
of the polyanalytic Fock spaces have been computed using several different methods.
We mention that
using the so called Laguerre polynomials can also be written
\begin{equation}
L_k^\alpha=\sum_{i=0}^k(-1)^i\binom{k+\alpha}{k-i}\frac{x^i}{i!}
\end{equation}
and are known to satisfy $\sum_{k=0}^{q-1} L_k^\alpha =L_{k-1}^{\alpha+1},$
the reproducing kernel for $\mathcal{F}_2^{(q)}(\C)$ can be written
\begin{equation}
\mathcal{K}^q(z,w)=\pi L_{q-1}^0(\pi\abs{z-w}^2)\exp(\pi z\bar{w})
\end{equation}
and the reproducing kernel for $\mathbf{F}_2^{(q)}(\C)$ can be written
\begin{equation}
\mathbf{K}^q(z,w)=\pi L_{q-1}^1(\pi\abs{z-w}^2)\exp(\pi z\bar{w})
\end{equation}

If for true poly-Fock spaces $\mathcal{F}_2^{(q)}(\C)$ one replaces the weight $\phi=\exp(-\abs{z}^2)$ 
with a weight depending on an integer $j$, according to 
$\phi=\phi_j(z):=\exp(-j\abs{z}^2)$, then it turns out that the resulting modified true poly-Fock spaces 
are
the $L^2$-eigenspaces
\begin{equation}
\{f\in L^2(\C)\colon \mathfrak{m}_q f=qf\}
\end{equation}
where
\begin{equation}
\mathfrak{m}_j :=-\frac{\partial^2}{\partial z\partial\bar{z}} +j\bar{z}\frac{\partial}{\partial\bar{z}}
\end{equation}
the latter is a so-called magnetic Schr\"odinger operator
(See e.g.\ Askour et. al\ \cite{askour}).
Sometimes the eigenspaces (or the eigenvalues) are called {\em Landau levels.}\index{Landau level} See e.g.\ Shigekawa \cite{shigekawa} for
more on information on such eigenspaces. 
We have chosen in this book, not to put any emphasis on results that border on applied mathematics. 
\subsubsection{The Berzein transform}
In this section we begin by discussing the poly-Bergman spaces 
on the unit disc, $D,$ 
and this will be denoted as before $A^2_q(D).$ 
Let \begin{equation}
\varphi_z(\zeta):=\frac{z-\zeta}{1-\zeta\bar{z}}
\end{equation}
As we have noted Eqn.(\ref{koseleveq}) shows that the Bergman kernels for the disc are given by
\begin{multline}
K_q(z,\zeta)=\frac{q}{(1-\zeta\bar{z})^{2n}}\sum_{j=0}^{q-1}(-1)^j\binom{q}{j+1}\binom{q+j}{q}\abs{1-\zeta \bar{z}}^{2(q-1-j)}
\abs{\zeta -z}^{2j}=\\
\frac{q}{(1-\zeta\bar{z})^{2q}}\sum_{j=0}^{q-1}(-1)^j\binom{q}{j+1}\binom{q+j}{q}\abs{\varphi_z(\zeta)}^{2j}
\end{multline}
Set $K_z(\zeta):=K_q(z,\zeta).$
We have $\norm{K_z}^2= K(z,z)=\frac{q^2}{(1-\abs{z}^2)^2}.$ Set
\begin{multline}
k(z,\zeta):=\frac{1-\abs{z}^2}{q}K(z,\zeta)=\\
\frac{(1-\abs{z}^2)\abs{1-\zeta\bar{z}}^{2(q-1)}}{(1-\zeta\bar{z})^{2q}}
\sum_{j=0}^{q-1}(-1)^j\binom{q}{j+1}\binom{q+j}{q}\abs{\varphi_z(\zeta)}^{2j}
\end{multline}
and $k_z(\zeta):=K_z(\zeta)/\norm{K_z(\zeta)}.$
Using 
\begin{equation}
1-\abs{\varphi_z(\zeta)}^{2}=\frac{(1-\abs{z}^2)(1-\abs{\zeta}^2)}{\abs{1-\zeta\bar{z}}^2}
\end{equation}
we have
\begin{multline}\label{ekvationmulticuck}
\abs{k(z,\zeta)}=\frac{1-\abs{z}^2}{\abs{1-\zeta\bar{z}}^{2}} \abs{\sum_{j=0}^{q-1}(-1)^j\binom{q}{j+1}\binom{q+j}{q}\abs{\varphi_z(\zeta)}^{2j}}=\\
\frac{1-\abs{\varphi_z(\zeta)}^2}{1-\abs{\zeta}^2} \abs{\sum_{j=0}^{q-1}(-1)^j\binom{q}{j+1}\binom{q+j}{q}\abs{\varphi_z(\zeta)}^{2j}}=\\
\frac{1-\abs{\varphi_z(\zeta)}^2}{1-\abs{\zeta}^2} \abs{q+\sum_{j=0}^{q-1}(-1)^j\binom{q}{j+1}\binom{q+j}{q}\abs{\varphi_z(\zeta)}^{2j}}
\end{multline}

\begin{lemma}\label{lemmashosho}
	There exists $c>0$ such that for any $f\in A^2_q$ and any $z\in D,$
	\begin{equation}
	\abs{f(z)}^2\leq c\int_D\abs{f(\zeta)}^2\abs{k(z,\zeta)} d\mu(\zeta)
	\end{equation}
\end{lemma}
\begin{proof}
	By Eqn.(\ref{ekvationmulticuck}), $(1-\abs{w}^2)\abs{k_z(w)}\to q$ as $\abs{\varphi_z(w)}\to 0$ i.e.\ for any $\epsilon>0$ there exists $\delta>0$ such that
	$\abs{\varphi_z(w)}<\delta$ for all $z,w.$ Hence there exists $r\in (0,1)$ such that
	$(1-\abs{w}^2)\abs{k_z(w)}\geq 1/2,$ whenever $\abs{\varphi_z(w)}<r.$ Since
	$\abs{k_z(w)}>\frac{1}{2(1-\abs{w}^2)}$ for all $z\in E(w,r)$ we obtain 
	\begin{equation}
	\int_D\abs{f(z)}^2\abs{k_z(w)}^2 d\mu(z)\geq \frac{1}{4(1-\abs{w}^2)^2}\int_{E(w,r)}\abs{f(z)}^2 d\mu(z)
	\geq \frac{\abs{f(w)}^2}{4C_r}
	\end{equation}
	For $c=4C_r$ this completes the proof.
\end{proof}
In this section we shall deal with the slightly more general situation of {\em symbols} that are merely $L^1$ (instead of previously being required to be bounded).
Hence
for any $f\in L^1(D),$ where $D=\{\abs{z}<1\},$ the Toeplitz operator $T_f :A^2_q(D)\to A^2_q(D)$ 
is defined as
\begin{equation}
(T_f g)(z):=\int_D f(\zeta)g(\zeta)P(fg),\quad g\in A^2_q(\Omega)
\end{equation}
for $g\in A^2_q(D)$ such that the integral is well-defined for all $z\in D,$
where $P$ denotes the orthogonal projection of $L^2(\Omega)$ onto $A^2_q(\Omega).$ 
When $f$
is bounded $T_fg$ is defined for all $g\in A^2_q(D),$
and is given by
$T_f g = P(fg),$ where $P$ is the orthogonal projection from $L^2(D)$ onto $A^2_q(D).$ Since $\overline{K}_z$
is a $q$-analytic function of $z$, $T_f(z)$ is $q$-analytic in $D.$ Furthermore, the domain of $T_f$ is the space, $W$, of all linear combinations of kernel functions $K_z$ which is dense in
$A^2_q(D).$ 
We have as pointed out before, in the case of the unit disc $D$
(see Eqn.(\ref{orthokernelekvenbra}))
\begin{equation}
(T_f g)(z):=\frac{1}{\pi}\int_D \frac{f(\zeta)g(\zeta)}{(1-z\bar{\zeta})^2}\mu(\zeta)
\end{equation}
The operator $T_f$ is called bounded on $A^2_q(D)$ if there is a constant $c>0$ such that $\norm{T_f g}\leq c\norm{g},$ for all $g\in
A^2_q(D).$ In such case $T_f$ extends uniquely to a bounded operator on $A^2_q(D)$, and for $h\in A^2_q(D)$ and $g:=\sum_{j=1}^m a_jK_{z_j}\in W,$
we have
\begin{multline}\label{fyratvaanp}
\langle T_f h,g\rangle =\sum_{j=1}^m \abs{a}_j \langle T_f h, K_{z_j} \rangle= \sum_{j=1}^m \bar{a}_j (T_f h)(z_j)=\\
\sum_{j=1}^m \bar{a}_j \int_D f(\zeta)h(\zeta)K_{z_j}(\zeta)d\mu(\zeta)
=
\int_D f(\zeta)h(\zeta)\bar{g}(\zeta)
\end{multline} 
If $f\in L^2(D)$ and $h\in A_q^2(D)$ is bounded then $(fh)\in L^2(D)$ thus
\begin{equation}
P(fh)=\langle P(fh),K_z \rangle =\langle fh, K_{z} \rangle=
\int_D f(\zeta)h(\zeta)K_{z}(\zeta) 
=(T_f h)(z)\end{equation} 
If $f$ is bounded then $\norm{T_f}\leq \norm{f}_\infty$ and $T_f^*=T_{\bar{f}}.$
If $f$ is holomorphic then 
\begin{equation}
\langle T_{\bar{f}} K_z, K_\zeta\rangle =\langle K_z, fK_\zeta\rangle =\bar{f}(z)\bar{k}_\zeta(z)=
\langle \bar{f}(z) K_z, K_\zeta\rangle
\end{equation}
and since the linear span of $\{K_\zeta :\zeta\in D\}$ is dense in $A^2_q(D)$ we have
$T_{\bar{f}} K_z=\bar{f}(z)K_z,$ $z\in D$ (and $T_f g=fg$). 
Now if $u\in L^2(D)$ is harmonic then there are holomorphic $f,g$ such that $u=f+\bar{g}$ so for all $z,\zeta\in D,$
\begin{equation}
\langle T_{u} K_z, K_\zeta\rangle =\langle \bar{g}K_z, K_\zeta\rangle +\langle K_z, \bar{f}K_\zeta\rangle =(f(\zeta)+\bar{g}(z))\langle 
K_z, K_\zeta\rangle
\end{equation}
In particular, $\langle T_u K_z, K_z\rangle =u(z)\norm{K_z}^2,$ $z\in D.$
If $T_u$ is bounded then $\norm{T_u}\leq \norm{u}_\infty,$ ans since the converse inequality holds for bounded $u$
we have $\norm{T_u}\leq \norm{u}_\infty.$
\begin{definition}
	Denote by $\mathcal{B}(A_q^2(D))$ the space of bounded linear operators $A_q^2(D)\to C^\omega(D)$ (where $C^\omega(D)$ denotes the set of real-analytic functions).
	The {\em Berezin transform}\index{Berezin transform}, $B(T),$ of an operator $T$ on $A_q^2(D)$ is defined as
	\begin{equation}
	B(T)(z)=\norm{K_z}^{-2}\langle TK_z,K_z\rangle =\langle Tk_z,k_z\rangle
	\end{equation}
\end{definition}
Note that $B(\mathcal{B}(A^2_q(D))\subseteq C^\omega(D)$ and $\norm{B(T)}\leq \norm{T}$, $T\in \mathcal{B}(A^2_q(D).$
The formula for $k_z$ gives for $z\in D,$
\begin{multline}\label{berzeinekv}
B(T_f)(z)=\langle T_f k_z,k_z\rangle =\int_D f(\zeta)\abs{k_z(\zeta)}^2d\mu(\zeta)=\\
\int_D f(\zeta)\frac{(1-\abs{z}^2)^2}{\abs{1-\zeta\bar{z}}^4} \left(
\sum_{j=0}^{q-1}\binom{q}{j+1}\binom{q+j}{q}\abs{\varphi_z(\zeta)}^{2j}
\right) 
\end{multline}
which is defined for all $f\in L^1(D)$ and $z\in D,$ independent of boundedness of $T_f.$
\begin{definition}
	$Bf(z)$ given by Eqn.(\ref{berzeinekv}) is called the
	{\em Berzein transform} of $f\in L^1(D).$
\end{definition}

\begin{lemma}
	\label{lemmafyra1}
	There is a constant $c > 0$ such that for any non-negative function
	$f \in L^1(D)$ satisfying $g= Bf\in L^1(D)$ and that $T_g$ is bounded on $A^2_q(D),$
	we
	have $\norm{T_f}\leq c\norm{T_g}$.
\end{lemma}
\begin{proof}
	Let $h$ be a linear combination of kernel functions. Then
	by Eqn.(\ref{fyratvaanp}) together with Lemma \ref{lemmashosho}
	\begin{multline}
	\langle T_g h,h\rangle =\int_D g(z)\abs{h(z)}^2d\mu(z)=\\
	\int_D \int_D f(\zeta)\abs{k_z(\zeta)}^2 \abs{h(z)}^2d\mu(z)=
	\int_D \left(\int_D f(\zeta)\abs{k_z(\zeta)}^2 \abs{h(z)}^2d\mu(z)\right)f(\zeta)d\mu(\zeta)\geq\\
	\frac{1}{c}
	\int_D \int_D\abs{h(\zeta)}^2 f(\zeta)d\mu(\zeta) =
	\frac{1}{c}
	\int_D \langle T_f h,h\rangle
	\end{multline}
	Since the set of linear combination of kernel functions is dense in $A^2_q(D)$ and since $T_g$ is bounded we have $T_f\leq cT_g.$ This completes the proof.
\end{proof}

\begin{definition}
	A bounded linear operator $T:X\to Y$ between Hilbert spaces is called a {\em Hilbert–Schmidt} operator if 
	there exists an orthonormal basis $\{e_j\}_{j\in \Z_+},$ in $X$ such that $\sum_{j}\norm{Te_j}^2$ is convergent.
\end{definition}
\begin{proposition}
	A Hilbert-Schmidt operator $T$ is compact.
\end{proposition}
\begin{proof}
	Let $\{e_j\}_{j\in \N},$ be an orthonormal basis such that $\sum_{j}\norm{Te_j}^2 <\infty.$
	Define the finite rank (thus compact) operators $T_j(f)=\sum_{k=1}^j \langle f,e_k\rangle Te_k$ which belongs to the span of $\{Te_1,\ldots,Te_j\}.$
	We have 
	\begin{multline}
	\norm{(T-T_j)f}\leq \sum_{k=j+1}^\infty \norm{Te_j}\abs{\langle f,e_k\rangle}\leq\\
	\left(\sum_{k=j+1}^\infty \abs{\langle f,e_k\rangle}^2 \right)^{\frac{1}{2}}
	\left(\sum_{k=j+1}^\infty \norm{Te_k}^2 \right)^{\frac{1}{2}}
	\end{multline}
	For $f$ such that $\norm{f}\leq 1,$ we thus have
	\begin{equation}
	\norm{(T-T_j)f}\leq \left(\sum_{k=j+1}^\infty \norm{Te_k}^2 \right)^{\frac{1}{2}}\to 0,\mbox{ as }j\to 0
	\end{equation}
	The limit is compact as the operator limit of finite rank operators.
	This completes the proof.
\end{proof}
An integral operator $T$ on on $L^2(X)$ for a separable Hilbert space $X$, with kernel $K\in L^2(X\times X)$ is a Hilbert-Schmidt operator.
To see this, set $K_z(\zeta)=K(z,\zeta),$ let $\{e_j\}_{j\in \N},$ be an orthonormal basis for $L^2(X)$, 
and note that $(Te_j)(z)=\int_X K_z(\zeta)e_j(\zeta)d\mu(\zeta)=$
$\langle K_z,\bar{e}_j\rangle.$ By dominated convergence
$\sum_{j=1}^\infty \norm{Te_j}^2 =$\\
$\int_X\sum_{j=1}^\infty \abs{\langle K_z,\bar{e}_j\rangle}^2 d\mu(\zeta).$
Furthermore, $\{\bar{e}_j\}_{j\in \N},$ is also an orthonormal basis hence
$\sum_{j=1}^\infty \norm{Te_j}^2 $
$=\int_X \sum_{j=1}^\infty \norm{K_z}^2 d\mu(\zeta)=$
$\int_X\int_X \sum_{j=1}^\infty \abs{K(z,\zeta)}^2 d\mu(\zeta)d\mu(z)<\infty ,$ since we assume $K\in L^2.$

\begin{lemma}\label{lemmafyra2}
	If $f$ is a bounded function on $D$ such that $\lim_{D\ni z, \abs{z}\to 1} f(z)=0$ then $T_f$ is compact on $A^2_q(D).$ 
\end{lemma}
\begin{proof}
	If $f$ has compact support then
	$T_f$ is an integral operator with bounded kernel $f(\zeta)\overline{K}_z(\zeta),$ $(z,\zeta)\in D\times D,$ 
	thus is a Hilbert-Schmidt operator. 
	hence compact. Also by a standard diagonal argument we have that if a
	a sequence of compact operators $\{T_k\}_{k\in \Z}$ converges to an operator $K$ in the norm topology, i.e.\ 
	$\norm{T_j -K}\to 0$, then $K$ is compact. Since 
	we can approximate $f$ by multiples of $f$ with appropriate test functions supported in $D$ we obtain that $T$ is compact.
	This completes the proof.
\end{proof}

\begin{theorem}[\u{C}u\u{c}kovi\'c \& Le \cite{cuckovic} Thm 4.3]\label{cuckovicthm}
	If $\in L^1(D)$ is non-negative then $T_f$ is compact on $A_q^2(D)$ iff $Bf(z)\to 0$ as $D\ni z,\abs{z}\to 1.$
\end{theorem}
\begin{proof}
	Using the fact that $\norm{K_q(z,z)}\to \infty$ as $D\ni z,\abs{z}\to 1,$ and that $q$-analytic polynomials are dense in $A^q_2(D)$ we have
	for $f\in A^2_q(D),$
	\begin{equation}
	\lim_{D\ni z, \abs{z}\to 1} (1-\abs{z}^2)f(z)=q \lim_{D\ni z,\abs{z}\to 1} \langle f(z),k(z,\cdot)\rangle =0
	\end{equation}
	In particular, $k_z\to 0$ in the weak sense as $D\ni z, \abs{z}\to 1.$ If $T_f$ is compact then
	\begin{equation}
	Bf(z)=B(T_f)(z)=\langle T_f k_z,k_z\rangle\to 0,\mbox{ as }D\ni z, \abs{z}\to 0
	\end{equation}
	For the converse assume $g=Bf\to 0$ as $D\ni z,\abs{z}\to 1.$
	Then $g$ is bounded. 
	By Lemma \ref{lemmafyra2} $T_{g}$ is compact.
	By Lemma \ref{lemmafyra1} $T_f$ is bounded. This completes the proof. 
\end{proof}

Axler \& Zheng \cite{axlerzheng}  prove the following.
\begin{proposition}[Axler \& Zheng \cite{axlerzheng}, Cor. 2.5]
	If $f$ is a bounded function on $\Omega$ then the Toeplitz operator $T_f$ is compact on $A^2_1(D)$ iff $B_0 f(z)\to 0,$
	$z\in D,\abs{z}\to 1.$
	Here $B_0 f(z)=\int_D f(\varphi_z(\zeta))d\mu(\zeta).$
\end{proposition}
\u{C}u\u{c}kovi\'c \& Le \cite{cuckovic} pose the question on whether the "if" direction can be generalized to the cases $A^2_q(D),$ $q\geq 2.$
In such a setting the integral kernel used in the definition of the Toeplitz operator $T_f$ and the Berezin transform $B$ depends upon $q\geq 2.$

\section{$A^2_q(\Omega)$-removability}
\begin{definition}[See Axler, Conway \& Mcdonald \cite{axlerconway}]
Let $\Omega\subset\C$ be a domain. 
A point $p_0\in \partial \Omega$ is called {\em $A^2_q(\Omega)$-removable} with respect to 
if there is an open neighborhood $V$ of $p_0$ such that any $A^2_q(\Omega)$ function
can be extended to a member of $A^2_q(\Omega\cup V).$ The set of $A^2_q(\Omega)$-removable points will be denoted
$\partial_{2-r}^q\Omega.$ The {\em $q$-Bergman essential boundary of $\Omega$}, $\partial_{2-e}^q\Omega$
is the set of all points of $\partial \Omega$ which are not $A^2_q(\Omega)$-removable.
\end{definition}
Note that $\partial_{2-r}^q\Omega$ is relatively open in $\partial\Omega,$ thus $\partial_{2-e}^q\Omega$ is compact and a measurable subset of $\C.$
\begin{proposition}\label{inrepunktremovethm}
Let $p_0\in \partial\Omega$ is an isolated point of $\partial \Omega$ then $p_0\in \partial_{2-r}^1\Omega.$
In particular, if there is an open neighborhood $U\subset\C$ of a point $p_0$ such that $(U\setminus \{p_0\})\subset \Omega$
then $p_0$ is $A_1^2(\Omega)$-removable.
\end{proposition}
\begin{proof}
There exists a $\delta>0$ such that $W:=\{ 0<\abs{z-p_0}<\delta\}\subset\Omega.$ If $f\in A^2_1(\Omega)$ then
the Laurent expansion near $p_0,$ 
\begin{equation}
f(z)=\sum_{-\infty}^\infty c_j(z-p_0)^j
\end{equation} 
yields with the polar representation $z=r\exp(i\theta)$,
\begin{equation}
\infty >\int_W \abs{f}^2 =\int_{0}^\delta \int_0^{2\pi} \abs{f}^2 drd\theta =\\ 2\pi \sum_{-\infty}^\infty \abs{c_j}^2 
\int_0^{\delta} r^{2j+1} dr
\end{equation} 
Since $\int_0^{\delta} r^{2j+1} dr$ diverges for $j<0$, $a_j=0,$ $j<0,$ so $f$ is holomorphic and thus has 
a holomorphic extension to $W\cup\{p_0\}$. Hence $p_0\in \partial_{2-r}^1\Omega.$ This completes the proof.
\end{proof}

\begin{proposition}\label{bddholoprop}
Let $K\subset \C$ be a set that can be covered by a family of discs $\{D_i\}_i$ where each $D_i$ has radius $r_i>0.$ 
If $\sum_i r_i <C$ for a constant $C>0,$ then the set of bounded holomorphic functions on $(\C\cup\{\infty\})\setminus K$ 
(denoted $H^\infty(\C\cup\{\infty\})\setminus K$))
consists only of constant functions. 
\end{proposition}
\begin{proof}
Denote the boundary $\Gamma_C:=\partial\left(\bigcup_{j} D_j\right).$ Let $f\in H^\infty(\C\cup\{\infty\})\setminus K),$
$f(\infty)=0.$
By the Cauchy theorem
\begin{equation}
f(z)=\frac{1}{2\pi i}\int_{\Gamma_C} \frac{f(\zeta)d\zeta}{z-\zeta},\quad z\notin \overline{\left(\bigcup_{j} D_j\right)}
\end{equation}
Hence
\begin{equation}
\abs{f(z)}\leq \frac{C\norm{f}}{2\pi\mbox{dist}(z,\Gamma_C)}
\end{equation}
Repeating the process for a sequence of $C$ going to zero we obtain $\abs{f(z)}=0$ for all $z\in \C\setminus K.$
This completes the proof.
\end{proof}

\begin{proposition}
$\partial_{2-r}^1\Omega$ has zero area.
\end{proposition}
\begin{proof}
Assume there exists a compact subset $K\subset \partial_{2-r}^1\Omega$ such that $\mu(K)>0$ where $\mu$ is the standard area measure on $\C.$  
Set 
\begin{equation}
F(z):=\int_E \frac{d\mu(\zeta)}{\zeta-z}
\end{equation}
Then $F$ is holomorphic on $\C\setminus K$ and it is continuous on $\C$ since it is 
the convolution of the locally integrable function $-\frac{1}{\zeta}$ and the characteristic function of a bounded set.
Since $\lim_{z\to \infty} F(z) =0$, $F$ has continuous extension to $\C$ and since 
$\lim_{z\to \infty} zF(z) =-\mu(E)$, $F$ is nonconstant. 
In particular, $F|_\Omega \in A^2_1(\Omega)$ extends to a non-constant bounded analytic
function defined on all of $\C$. This contradicts the Liouville theorem.
Hence each compact $K\subset \partial_{2-r}^1\Omega$ must have zero area and thus 
$\mu(\partial_{2-r}^1\Omega)=0.$ This completes the proof. 
\end{proof}

\begin{proposition}
Let $p_0\in \partial\Omega.$ If the connected component of $\partial \Omega$ containing $p_0$ contains more than one
point then $p_0\in \partial_{2-e}^1\Omega.$
\end{proposition}
\begin{proof}
Assume that the connected component of $\partial \Omega$ 
contains $p_0$ and also contains more than one
point but that $p_0\in \partial_{2-r}^1\Omega.$ Since
$\partial_{2-r}^1\Omega$ is relatively open in $\partial\Omega$ there is a $\delta>0$ such that
$M:=\{ \abs{z-p_0}\leq \delta\}\cap \partial\Omega \subset \partial_{2-r}^1\Omega.$ Let 
$K$ be the connected component of $M$ containing $p_0.$ Then $K$ contains more than one point
so there exists a conformal map $f$ of $(\C\cup\{\infty\})\setminus K$ onto the open unit disc such that 
$f|_\Omega \in A^2_1(\Omega)$ and since $K\subset \partial_{2-r}^1\Omega$ $f$ extends to a non-constant bounded
holomorphic function on $\C.$ This is a contradiction to Liouvilles theorem. 
This completes the proof.
\end{proof}

\begin{proposition}\label{frageteckprop}
$\partial\overline{\Omega}\subset \partial_{2-e}\Omega.$ Furthermore
$\Omega\cup \partial_{2-e}\omega$ is an open subset of $\C.$
\end{proposition}
\begin{proof}
Let $p_0\in \partial \overline{\Omega}$ and let $U$ be an open neighborhood of $p_0.$
Then $U$ contains a point $q_0\notin \overline{\Omega}.$
Then the function
$\frac{1}{z-q_0}|_\Omega$ belongs to $A^2_1(\Omega).$ but has
no holomorphic extension to $\Omega\cup U,$ thus
$q_0\in \partial_{2-e}\Omega.$ This proves
$\partial\overline{\Omega}\subset \partial_{e-q}\Omega.$
Now let $p_0\in \partial_{2-e}\Omega.$ Since 
$\partial_{2-e}\Omega$ is relatively compact in $\partial\Omega$
there exists $\epsilon>0,$ such that
$\{\abs{z-p_0}<\epsilon\}\cap \partial\Omega \subset \partial_{2-r}\Omega.$
Since $\partial\overline{\Omega}\subset \partial_{e-q}\Omega$,
$p_0$ belongs to the interior of $\overline{\Omega}.$
Hence there is a $\delta>0$ such that
$\{\abs{z-p_0}<\delta\}\subset \overline{\Omega}=\Omega \cup\partial\Omega$ so
$\{\abs{z-p_0}<\delta\}\subset \Omega\cup \partial_{2-r}\Omega$ and thus
$\Omega\cup \partial_{2-r}\Omega$ is open in $\C.$ This completes the proof.
\end{proof}

For what follows the reader may want to see Section \ref{ellipticapp} for some preliminary definitions from distribution theory.
Here we shall, in line with e.g.\ Friedlander \& Joshi \cite{friedlander}, use the notation
$\langle u,\phi \rangle$ for defining the distribution $u\in \mathcal{D'}$ acting on $\phi\in \mathcal{D}.$ In the literature an alternative notation
for this is $(u,\phi)$.
\begin{definition}[Special case of Theorem 1.4, Harvey \& Polking \cite{harveypolking}]
Let $\mathbf{F}\subset\mathcal{D}'(\Omega)$ be a class of distributions on an open subset $\Omega\subset\Rn$.
Let $P(x,D)$ be a linear partial differential operator with $C^\infty$-smooth coefficients defined on $\Omega.$ 
The {\em $L^p(\Omega)$-capacity\index{$L_P$-capacity} of a subset $A\subset \Omega$ with respect to $P(x,D)$}
(denoted $L^p-\mbox{Cap}_P(A,\Omega)$)
is defined as
\begin{equation}
\sup\{\abs{\langle P(D)f,1 \rangle} \colon f\in L^p(\Omega),\norm{f}\leq 1,\mbox{supp}(P(D)f)\subset A\mbox{ is compact}\}
\end{equation}
If $P(D)$ has constant coefficients on $\Rn$ and has a fundamental solution (in the sense that there exists $E\in \mathcal{D}'(\Rn)$ 
such that $P(D)E=\delta$).
Then the {\em $L_P^E$-capacity of a subset $A\subset \Rn$ with respect to $E$}
is defined as
\begin{equation}
\sup\{\abs{\langle u,1 \rangle} \colon u\in \mathcal{E}'(\Rn),\norm{E*u}\leq 1,\mbox{supp}(u)\subset A\}
\end{equation}
A relatively closed subset $A\subset\Omega$ is called {\em removable for $\mathbf{F}$ with respect to $P(D)$}
\index{Removable set with respect to linear partial differential operators} if
each $f\in \mathbf{F}$ which satisfies $P(x,D)f=0$ in $\Omega\setminus A$ also satisfies $P(x,D)f=0$ in $\Omega.$
\end{definition}
If a set $A\subset \Omega$ is removable for $\mathbf{F}$ with respect to $P(x,D)$ then the $L_P^p(\Omega)$-capacity of $A$
is zero. Denote by $P^{(\alpha)}(\xi)$ the map $\partial_\xi^\alpha P(\xi),$ for a polynomial $P(\xi).$
\begin{theorem}[Harvey \& Polking \cite{harveypolking}]\label{harveypolking1}
Let $\Omega \subset\subset \Rn$ be a relatively compact open subset and let $A\subset\subset \Omega$ be a relatively closed subset.
Let $P(D)$ be a linear partial differential operator with $C^\infty$-smooth coefficients defined on $\Omega.$
Suppose $P(D)$ has a fundamental solution $E$ satisfying for each $\alpha,$ 
\begin{equation}\label{konditionseekk}
P^{(\alpha)}(D)E\in \mathcal{M}_{\mbox{loc}}(\Rn)  
\end{equation}
where $\mathcal{M}_{\mbox{loc}}(\Rn)$ denotes the space of regular Borel measures on $\Omega.$
Then $A$ is removable for $L^p_{\mbox{loc}}(\Omega),$ $1<p\leq \infty$, iff $L^p-\mbox{Cap}_P(A,\Omega)=0.$
\end{theorem}
For the readers convenience we reproduce here the proof with the notations for the particular case of interest.
\begin{proof}
As we have already pointed out, if a set $A\subset \Omega$ is removable for $L^p_{\mbox{loc}}(\Omega),$ with respect to $P(D)$, 
then the $L_{\mbox{loc}}^p(\Omega)$-capacity of $A$ is zero. We prove the converse.
Assume $L^p-\mbox{Cap}_P(A,\Omega)=0$ and that $f\in L_{\mbox{loc}}^2(\Omega)$ such that $P(D)f=0$ on $\Omega\setminus A.$
Let $\phi\in C^\infty_c(\Omega)$ and set $g:=E*(\phi P(D)f).$ Then
\begin{equation}
P(D)g=\phi P(D)f
\end{equation}
Let $\psi\in C^\infty_c(\Rn).$ From Leibniz generalized rule for distributions
\begin{multline}
\langle \phi P(D)f,\psi\rangle =\langle f, P^*(D) (\phi\psi)\rangle=\\
\sum ((-1)^{\abs{\alpha}}/\alpha!) \langle f, D^\alpha \phi^* P^{(\alpha)}(D)\psi \rangle=\\
\langle \sum ((-1)^{\abs{\alpha}}/\alpha!) P^{(\alpha)}(D)(fD^\alpha \phi),\psi \rangle
\end{multline}
Hence we have adjoint formula 
\begin{equation}
\phi P(D)f= \sum_{\abs{\alpha}\geq 0} ((-1)^{\abs{\alpha}}/\alpha!)P^{(\alpha)}(D)(fD^\alpha \phi)
\end{equation}
which implies
\begin{equation}\label{ocksa}
g=E*(\phi P(D)f)= \sum_{\abs{\alpha}\geq 0} ((-1)^{\abs{\alpha}}/\alpha!)(P^{(\alpha)}(D)E)*(fD^\alpha \phi)
\end{equation}
Since $P^{(\alpha)}(D)E\in \mathcal{M}_{\mbox{loc}}(\Rn)$ and $\Omega\subset\subset\Rn$, 
convolution by $\partial^\alpha P(D)E$ followed by restriction to $\Omega$ is a continuous linear mapping $L^2_{\mbox{loc}}\to L^2_{\mbox{loc}},$
Since each $fD^\alpha \phi\in L^2_{\mbox{loc}}$ Eqn.(\ref{ocksa}) implies $g\in L_{\mbox{loc}}^2(\Omega).$ Since by assumption $L^p-\mbox{Cap}_P(A,\Omega)=0,$
\begin{equation}
\langle P(D)f,\phi \rangle= \langle \phi P(D) f,1 \rangle=\langle P(D)g,1\rangle =0 
\end{equation}
i.e.\ $P(D)f=0$ on $\Omega.$ This completes the proof.
\end{proof}
Harvey \& Polking \cite{harveypolking} also prove the following.
\begin{theorem}\label{harveypolking2}
Suppose $P(D)$ has a fundamental solution $E$ such that $P^{(\alpha)}(D)E\in \mathcal{M}_{\mbox{loc}}(\Rn)$ for all $\alpha.$
Let $A$ be a relatively closed subset of an arbitrary open set $\Omega\subset\Rn.$
Let $\psi\in C^\infty(\Rn)$ such that $\psi$ vanishes on a neighborhood of the origin and is identically one in a neighborhood
of infinity. If $\psi P^{(\alpha)}(D)E\in L^2(\Rn)$ for all $\alpha,$ then
$A$ is removable for $L^2_{\mbox{loc}}(\Omega)$ iff $L^2-\mbox{Cap}_E(A)=0.$
\end{theorem}
\begin{proof}
The proof follows from the same arguments as the proof of Theorem \ref{harveypolking1}
except that we must verify that each $P^{(\alpha)}(D)E*(fD^\alpha \phi)\in L^2(\Rn).$ Now
$(1-\psi)P^{(\alpha)}(D)E$ is a compactly supported Borel measure and $fD^\alpha \phi \in L^2(\Rn).$
Hence $(1-\psi)P^{(\alpha)}(D)E*fD\phi \in L^2(\Rn).$ By hypothesis
$\psi P^{(\alpha)}E\in \in L^2(\Rn),$ and $fD^\alpha \phi\in L^1(\Rn).$ Hence $(\psi P^{(\alpha)}E)*fD^\alpha \phi\in L^2(\Rn).$
This completes the proof.
\end{proof}
\begin{remark}
Harvey \& Polking \cite{harveypolking} point out that if $P$ and $Q$ are two elliptic operators of order $m$ then
$L^2-\mbox{Cap}_{P}(A) = 0$ iff $L^2-\mbox{Cap}_{Q}(A) = 0.$ Furthermore
they define the $\mbox{Lip}^\delta -\mbox{Cap}_P(A,\Omega)$ analogous to $L^2-\mbox{Cap}_P(A,\Omega)$ but 
replacing $L^2_{\mbox{loc}}(\Omega)$
with the Lipschitz space with norm 
\begin{equation}
\norm{f}_{\mbox{Lip}^\delta(\Omega)}:=\norm{f}_\infty +\sup_{x,y\in\Omega}\frac{\abs{f(x)-f(y)}}{\abs{x-y}^{\delta}}
\end{equation}
They prove that the relatively closed subset $A\subset \Omega$ is removable with respect to $P$ in
in $\mbox{Lip}^\delta(\Omega)$ iff $\mbox{Lip}^\delta -\mbox{Cap}_P(A)=0,$ and they point out that 
the latter happens iff $\Lambda_{n-m+\delta}(A)=0$, where $\Lambda_{k}(A)$ denotes the $k$-dimensional 
Hausdorff measure.  
\end{remark}
Theorem \ref{harveypolking1} cannot be used directly for $A^2_p(\Omega)$ when $q>1$, but as we shall see below it could be used indirectly, but first it is interesting to analyze the reason why it cannot be used directly.
Indeed, for $q=1$, poles of meromorphic functions render functions that are not $L_{\mbox{loc}}^2(\Omega)$, and it is natural that isolated points are 
$L_{\mbox{loc}}^2(\Omega)$-removable. For $q>1$ there exists $q$-analytic functions whose 
holomorphic components could be realized as restrictions of meromorphic functions, and thus (due to the uniqueness of the representation
of holomorphic components) those poles will in general not be $A^2_q(\Omega)$-removable. For example if $\Omega=\{0<\abs{z}<1\},$
then the function $\bar{z}^j z^{-j},$ $j\in \Z_+$, is clearly $A^2_p(\Omega)$ but has no extension to $A^2_p(\{\abs{z}<1\}).$
Let us look at the functions which render positive $L^2-\mbox{Cap}_{\partial_{\bar{z}}^q}$-capacity for the origin as an isolated point.
So we are looking for functions $f$, $\norm{f}\leq 1,$ such that $\abs{\langle \partial_{\bar{z}}^q f,1\rangle}>0.$
We know that the normalization, $\tilde{E},$ of the fundamental solution $E_q:=\frac{\bar{z}^{q-1}}{z\pi (q-1)!}$ satisfies this.
We shall see below that the span of functions of the simple form $\frac{\bar{z}^{k}}{z^l},$
$k=1,\ldots,j-1,l=1,\ldots ,k$ and functions in $A^2_j(\Omega)$ can be used to specify all
members of $A^2(\Omega\setminus \{0\}).$
It is known that if the conditions of Theorem \ref{harveypolking2} are satisfied then the
$L^2-\mbox{Cap}_{\partial_{\bar{z}}^q}(A)$ is zero iff the
$L^2-\mbox{Cap}_{E_q}(A)$ is zero.
The above observations makes it interesting to consider the following. 
\begin{definition}
For a domain $\Omega\subset\Cn,$ and a multi-index $\alpha\in\Z_+^n$, denote by $WA^2_{\alpha,\mbox{loc}}(\Omega)$, the intersection of the space of $\alpha$-analytic functions
on $\Omega$ ($\mbox{PA}_\alpha(\Omega)$) with the
Sobolev spaces $W^{\max_j\alpha_j,2}_{\mbox{loc}}(\Omega).$
\end{definition}
\begin{proposition}\label{lemprop1}
Let $q\in\Z_+,$ let $\Omega\subset\C$ be a domain and let $f\in WA^2_{q,\mbox{loc}}(\Omega)$. Then each analytic component of
$f$ is $L^2_{\mbox{loc}}(\Omega).$
\end{proposition}
\begin{proof}
If $\Omega\subset \C$ is a domain, and $f$ a $q$-analytic function on $\Omega$ then $f$ has, on $\Omega,$ the representation
$\sum_{j=0}^{q-1} a_j(z)\bar{z}^j.$ If $f$ is also $W^{q-1,2}_{\mbox{loc}}(\Omega),$
then each $\partial_{\bar{z}}^j f\in L^2_{\mbox{loc}}(\Omega),$ $j=0,\ldots,q-1.$ In particular
$(q-1)!a_{q-1}(z)=\partial_{\bar{z}}^{q-1} f\in L^2_{\mbox{loc}}(\Omega)$.
Suppose that
the components $a_{q-1},\ldots,a_{q-1-k},$ are $L^2_{\mbox{loc}}(\Omega)$ for some $k\in \{1,\ldots,q-1\}$
(the case $k=0$ is already handled). 
It follows that each $\bar{z}^l a_j(z)\in L^2_{\mbox{loc}}(\Omega),$ for $j=q-1-k,\ldots q-1,$
and $l\in \Z_+.$
Hence
$\partial_{\bar{z}}^{q-2-k} f(z) =(q-2-k)! a_{q-2-k}(z) +\varphi(z) $ for a function $\varphi\in L^2_{\mbox{loc}}(\Omega),$
and so $a_{q-2-k}=(f-\varphi) \in L^2_{\mbox{loc}}(\Omega).$ 
This completes the proof.
\end{proof}
\begin{proposition}\label{harvpolkabtinprop}
Let $q\in \Z_+,$ let $\Omega\subset\C$ be a domain and let $A\subset\Omega$ be a subset such that
$L^2-\mbox{Cap}_{\partial_{\bar{z}}}(A)=0$ and such that $\Omega\setminus A$ is connected. Then 
\begin{equation}
WA^2_{q,\mbox{loc}}(\Omega\setminus A)=WA^2_{q,\mbox{loc}}(\Omega)
\end{equation} 
\end{proposition}
\begin{proof}
Let $f\in WA^2_{q,\mbox{loc}}(\Omega\setminus A).$ Let $p_0\in \Omega.$
Then on $\Omega\setminus A$ 
$f$ has the representation
$f(z)=\sum_{j=0}^{q-1} a_j(z)\bar{z}^j,$ for functions $a_j$ holomorphic.
By Proposition \ref{lemprop1}, each $a_j\in L^2_{\mbox{loc}}(\Omega\setminus A).$ Since 
$L^2-\mbox{Cap}_{\partial_{\bar{z}}}(A)=0$, Theorem \ref{harveypolking1} 
implies that each $a_j$ has an extension to some $\tilde{a}_j\in A^2_{1,\mbox{loc}}(\Omega).$ 
This implies (since due to locality the factor $\bar{z}^j$ cannot turn the product into a non-member of
$L^2_{\mbox{loc}}(\Omega)$) that each $\bar{z}^j\tilde{a}_j\in A^2_{j+1,\mbox{loc}}(\Omega).$
Hence for any $p_0\in \Omega$ there is an open neighborhood $U$ of $p_0$ in $\Omega$ such that
\begin{equation}
\norm{\sum_{j=0}^{q-1} \tilde{a}_j(z)\bar{z}^j}_{L^2(U)}\leq \sum_{j=0}^{q-1} \norm{\tilde{a}_j(z)\bar{z}^j}_{L^2(U)}<\infty
\end{equation}
Hence the $C^\infty$-smooth function $\tilde{f}(z):=\sum_{j=0}^{q-1} \tilde{a}_j(z)\bar{z}^j$ defines an extension of $f$ to
$WA^2_{q,\mbox{loc}}(\Omega).$
\end{proof}
\begin{definition}\label{flerdimharvpolkdef}
Let $(z_1,\ldots,z_n)$ denote Euclidean coordinates for $\Cn$ and let $\Omega\subset\Cn$ be a domain.
For each $c\in \C^{n-1}$, denote $M_{c,j}:=\{z\in \Cn\colon z_k=c_k, j\neq k \}$.
Let $\Omega\subseteq\Cn$ be a domain and let $A\subset\Omega$ be a subset.
$A$ is said to have {\em separately zero 
 $L^p(\Omega)$-capacity with respect to the Cauchy-Riemann operator}
(denoted $L^p-\mbox{Cap}_{\overline{\partial}}(A,\Omega)=0$)
if for each $c\in \C^{n-1}$, and each $j=1,\ldots ,n,$ we have
$L^2-\mbox{Cap}_{\partial_{\bar{z}_j}}(A\cap M_{c,j})=0$ (i.e.\ the $L^2$-capacity of $A\cap M_{c,j}$, calculated on the complex one-dimensional Euclidean copy $M_{c,j}$, of $\C$, with respect to the operator
$\partial_{\bar{z}_j}$, is zero).
\end{definition}
\begin{proposition}
Let $\alpha\in \Z_+^n,$ let $\Omega\subset\Cn$ be a domain and let $A\subset\Omega$ be a subset such that
$\Omega\setminus A$ is connected and $L^p-\mbox{Cap}_{\overline{\partial}}(A,\Omega)=0$. Then 
\begin{equation}
WA^2_{\alpha,\mbox{loc}}(\Omega\setminus A)=WA^2_{\alpha,\mbox{loc}}(\Omega)
\end{equation} 
\end{proposition}
\begin{proof}
Let $f(z)\in WA^2_{\alpha,\mbox{loc}}(\Omega\setminus A).$
For each $c\in \C^{n-1}$ denote\\ $\hat{z}_j:=(c_1,\ldots,c_{j-1},z_j,c_{j},\ldots,c_{n-1}),$ for $j=2,\ldots,n-1,$
$\hat{z}_1:=(z_1,c_1,\ldots,c_{n-1}),$ and $\hat{z}_n:=(c_1,\ldots,c_{n-1},z_n).$
For each $j=1,\ldots ,n,$ set $f_{j,c}(z_j):=f(\hat{z}_j).$
Then each $f_{j,c}(z_j)\in WA^2_{\alpha_j,\mbox{loc}}(M_{c,j}\cap(\Omega\setminus A)).$
Since $L^p-\mbox{Cap}_{\overline{\partial}}(A,\Omega)=0$, Proposition \ref{harvpolkabtinprop} together with Definition \ref{flerdimharvpolkdef} implies that $f_{j,c}(z_j)$ has an extension $\tilde{f}_{j,c}(z_j)\in WA^2_{\alpha_j,\mbox{loc}}(M_{c,j}\cap(\Omega)).$
Since this holds true for all $c\in\C^{n-1}$ and all $j=1,\ldots,n$,
$f$ has a separately $\alpha$-analytic extension, $\tilde{f}$ to $\Omega$. 
By Theorem \ref{hartog1} this implies that
$\tilde{f}$ is locally a.e.\ equal to a $C^\infty$-smooth function with respect to $(\re z_1,\im z_1,\ldots,\re z_n,\im z_n)$ that is
jointly $\alpha$-analytic on $\Omega,$ in particular each derivative of order $\leq \max_j \alpha_j$ of $\tilde{f}$ is a member of $L_{\mbox{loc}}^2(\Omega)$ hence $\tilde{f}\in WA^2_{\alpha,\mbox{loc}}(\Omega).$ 
This completes the proof.
\end{proof}

\begin{definition}
Let $\Omega\subseteq\C$ be a domain and let $X\subset \Omega$ be a subset. $X$ 
is called {\em $WA^2_q$-removable} with respect to $\Omega,$
if $WA^2_{q,\mbox{loc}}(\Omega\setminus X)=WA^2_{q,\mbox{loc}}(\Omega).$
The {\em $q$-essential boundary with respect to $WA^2_{q,\mbox{loc}}(\Omega)$}, denoted $\partial_{2-r}^{q,W}$, is
defined as the set of $WA^2_q$-removable points with respect to $\Omega.$
\end{definition}
Since the Harvey-Polking theorem is an 'if and only if' we obtain a generalization regarding
the essential boundary.
\begin{proposition}
Let $\Omega\subset\C,$ be a bounded, open, connected,
non-empty subset. Then for each $q\in \Z_+$ we have
\begin{equation}
\partial_{2-r}^{1}\Omega =\partial_{2-r}^{q,W}\Omega
\end{equation}
\end{proposition}
Clearly, the fact that isolated points are not automatically removable for $A^2_q$, $q>1$,
implies that for a domain $\Omega\subset \C,$ the essential $q$-Bergman boundary $\partial_{2-e}^q \Omega$ 
does not (in contrast to the case $q=1$) necessarily belong to
$\partial (\Omega\cup\partial\Omega).$ For this reason, some authors (see Pessoa \cite{pessoaess})
introduce alternative definitions, replacing both the notion of removability and the notion
essential boundary with other things. We shall not do this here. 
Pessoa \cite{pessoa2013} shows an estimation similar to that of Theorem \ref{balk204} 
based upon previous results of Vekua \cite{vekua} and Dzhuraev \cite{dzhuraevbok}.
\begin{definition}
For a non-empty finitely connected domain $\Omega \subset \C,$ 
define for $j\in \Z_+,$
\begin{equation}
S_{\Omega,j} f(z) := -\frac{(-1)^j j}{\pi}\int_\Omega \frac{(\zeta -z)^{j-1}}{(\bar{\zeta} -\bar{z})^{j+1}}f(\zeta)d\mu(\zeta),
\end{equation}
\begin{equation}
\tilde{S}_{\Omega,j} f(z) := -\frac{(-1)^{-j} (-j)}{\pi}\int_\Omega \frac{(\zeta -z)^{-j-1}}{(\bar{\zeta} -\bar{z})^{-j+1}}f(\zeta)d\mu(\zeta),
\end{equation}
Let $K_{\Omega,j}$ denote the reproducing kernel of $A^2_j(\Omega).$
\end{definition}
Note that $S_1$ can be considered as the weak derivative of the Pompieu operator 
\begin{equation}
Tf=-\frac{1}{\pi}\frac{f(\zeta)}{\zeta -z}d\mu(\zeta)
\end{equation}
Note also that 
\begin{equation}
S_{\Omega,j}^*f(z) := \tilde{S}_{\Omega,j}
\end{equation}
Vekua \cite{vekua} (see p.61) proved that for $f\in C^{\infty}(\Omega)\cap L^2$
\begin{equation}\label{vekuaeq}
\partial_{\bar{z}} S_{\Omega,1} f =\partial_z f,\quad \partial_{z} S^*_{\Omega,1} f =\partial_{\bar{z}} f
\end{equation}
\begin{proposition}\label{propositionpessoa2}
Let $U\subset \C$ be a bounded domain and $p_0\in U.$ Then for each $f\in A^2_j(U)$ and each $k\in \Z_+$
there exists a positive constant $M$ depending on $j,k$ such that
\begin{equation}
\abs{\partial_{\bar{z}}^k f_{k}(p_0)}\leq M\frac{\norm{f}}{(\mbox{dist}(z,\partial U))^{j+k+1}}
\end{equation}
\end{proposition}
\begin{proof}
Assume w.l.o.g.\ $p_0=0$ and let $B$ be a ball centered at $0$ and relatively compact in $U.$
Let $D:=\{\abs{z}<1\}.$
Let $f(z)\in A^2_1(tD),$ $t>1$, have the representation $f(z)=\sum_{l=0}^{j-1} f_l(z)\bar{z}^l.$ Taylor expansion
of each $f_l(z)$ near $0$ gives the uniformly convergent series
\begin{equation}
f(z)=\sum_{l=0}^{j-1} \sum_{m=0}^\infty \frac{1}{m!}\partial_z^m f_l(0)z^m\bar{z}^l
\end{equation}
in $D.$
the orthogonality condition
\begin{equation}
\langle z^m,z^n\rangle =\frac{\pi}{j+1} \delta_{m,n},\quad n,m\in \N,
\end{equation}
where $\delta_{n,m}=1$ for $n=m$, and $\delta_{n,m}=1$ for $n\neq m$.
This implies 
\begin{equation}
\langle f(z),z^n\rangle =\frac{\pi \partial_z^n f(0)}{n+1},\quad f\in A_1^2(tD), n\in \N, t>1
\end{equation}
Hence for $j=1$ the polynomial $p_{j,1}(z,\bar{z}):=\frac{(n+1)!}{\pi}z^n$
satisfies $\partial_z^n f(0)=\langle f,p_{n,j}\rangle,$ $n\in \N.$
If $f\in A^2_1(D)$ and $0<t<1,$ set $f_t(z):=f(tz)\in A^2_1(t^{-1} D).$ We have in $L^2(D),$ 
$\lim_{t\to 1-} f_t =f.$ Which implies for $f\in A^2_1(D),$
\begin{equation}
 \partial_z^n f(0) =\lim_{t\to 1-} \partial_z^n f_r(0) =\lim_{t\to 1-} \langle f_r,p_{n,1}\rangle=\langle f,p_{n,1}\rangle
\end{equation}
We show by induction in $j\geq 1,$ that there exists a polynomial $p_{n,j}(z,\bar{z})$
such that for $f\in A^2_j(D)$,
\begin{equation}
\langle f(z),p_{n,j}\rangle = \partial_z^n f(0)
\end{equation}
The case $j=1$ is proved so let $j>1$ and assume the claim holds true for $A^2_k(D),$ $k=1,\ldots ,j-1.$
Let $f\in A^2_j(tD),$ $t>1$ and write near $0$ gives the uniformly convergent series in $D,$
\begin{equation}
f(z)=\sum_{l=0}^{j-1} \sum_{m=0}^\infty \frac{1}{m!l!}\partial_z^m \partial_{\bar{z}}^m f(0)z^m\bar{z}^l
\end{equation}
For a $j$-analytic $f$ on the unit disc, $D$, Eqn.(\ref{vekuaeq}) implies $(S_D^*)^j f=0$  
For $g_k:=(S_D^*)^{k} f=0,$ $k=1,\ldots, j-1,$ we have $(S_D^*)^{j-k} g=0,$ i.e.\ $g_k\in A^2_{j-k}(D).$
The conditions
\begin{equation}
\langle z^{l_1}\bar{z}^{s_1},z^{l_2}\bar{z}^{s_2}\rangle =\langle z^{l_1+s_2},z^{l_2+s_1}\rangle,
\quad l_k,s_k\in \N,k=1,2
\end{equation}
together with the Vekua derivation formulas gives
\begin{multline}
\langle f,z^n\rangle =\frac{\partial_z^n f(0)}{n!} \langle z^{n},z^{n}\rangle =\sum_{k=1}^{j-1}
\frac{\partial_z^{n+k}\partial_{\bar{z}}^{k} f(0)}{(n+k)!k!} \langle z^{n+k}\bar{z}^k,z^n\rangle =\\
\frac{\pi \partial_z^n f(0)}{(n+1)!} +\sum_{k=1}^{j-1} 
\frac{\partial_z^{n+2k} g(0)}{(n+k)!k!} \langle z^{n+k},z^{n+k}\rangle=\\
\frac{\pi \partial_z^n f(0)}{(n+1)!} +\pi \sum_{k=1}^{j-1} 
\frac{\langle (S_D^*)^k f,p_{n+2k,j-k}\rangle}{(n+k+1)!k!} =\\
\frac{\pi \partial_z^n f(0)}{(n+1)!} +\pi \sum_{k=1}^{j-1} 
\frac{\langle f, S_D^^k p_{n+2k,j-k}\rangle}{(n+k+1)!k!} 
\end{multline}
By Karlovich \& Pessoa \cite{karlovichpess}, Lemma 2.2, it holds that for each $n,m\in \Z_+,$ and each $z\in D,$
\begin{equation}
(S_D (\bar{z}^n z^m))(z) =\frac{m}{n+1}\bar{z}^{n+1} z^{m-1} +\frac{\min\{0,n+1-m\}}{m+1}z^{m-n-2}
\end{equation}
which implies that for each $k\in \N$ and each polynomial $p(z,\bar{z}),$
$(S_D)^kp$ is a polynomial in $\bar{z}$ and $z.$ Hence for each $f\in A^2_j(tD),$ $t>1,$ $n\in \N,$
\begin{equation}
\partial_z^n f(0) =\langle f,p_{n,j}\rangle
\end{equation}
where $p_{n,j}(z,\bar{z})$ is defined recursively by $p_{n,1}=\frac{(n+1)!}{\pi}z^n,$
and 
\begin{equation}
p_{n,j}:=p_{n,1}-(n+1)!\sum_{k=1}^{j-1}\frac{(S_D^*)^k p_{n+2k,j-k}}{(n+k+1)!k!}
\end{equation}
As before, setting $f_t(z)=f(tz),$ for an $f\in A^2_j(D),$ $0<t<1,$ we obtain by repeating the arguments of the case $j=1$,
$\partial_z^n f(0) =\langle f,p_{n,j}\rangle$ for $f\in A^2_j(D).$
Thus the functional $f\mapsto \partial_z^n f(0)$ has norm $\norm{p_{n,j}},$ which implies for $f\in A^2_j(D)$ and $n\in \N,$ 
\begin{equation}
\abs{\partial_z^n f(0)}\leq \norm{p_{n,j}}\norm{f}
\end{equation}
Let $w\in \C$ and $r>0$, $f\in A^2_j(\{\abs{z-w}<r\}),$ and $g(\xi):=f(w+\xi r),$ for $\xi\in D.$
Then $g\in A^2_j(D)$, $(S_D^*)^k g\in A^2_{j-k}(D)$, and on $D$,
\begin{equation}
\partial_{\bar{z}}^k g = \partial_z^k(S_D^*)^k g
\end{equation}
Hence
\begin{multline}
\abs{\partial_{\bar{z}}^k f(w)} =\frac{\abs{\partial_{\bar{z}}^k g(0)}}{r^{k}}=\frac{\abs{z^k (S_D^*)^k g(0)}}{r^{k}}=\\
\frac{\abs{\langle g, S_D^k p_{k,j-k}\rangle}}{r^{k}}\leq \frac{\norm{ S_D^k p_{k,j-k}}}{r^{k}}\norm{g}_D =
\frac{\norm{ S_D^k p_{k,j-k}}}{r^{k+1}}\norm{f}_{\{\abs{z-w}<r\}}
\end{multline}
If $f\in A^2_j(U)$ and $w\in U,$ then $f\in A^2_j(\{\abs{z-w}<r\})$ which implies
\begin{equation}
\abs{\partial_{\bar{z}}^k f(w)} \leq \frac{\norm{ S_D^k p_{k,j-k}}}{(\mbox{dist}(w,\partial U))^{k+1}}\norm{f}_{\{\abs{z-w}<r\}}
\leq \frac{\norm{ S_D^k p_{k,j-k}}}{(\mbox{dist}(w,\partial U))^{k+1}}\norm{f}
\end{equation}
This completes the proof.
\end{proof}
\begin{proposition}[Pessoa \cite{pessoa2016}]
Let $\Omega\subset \C$ be a bounded domain and $p_0\in \Omega.$ 
Then for $j>1,$
\begin{equation}
A_j^2(\Omega\setminus\{p_0\})=\mbox{span}\left\{
\psi,\frac{(\bar{z}-\bar{\xi})^k}{(z-\xi)^l}\colon \psi\in A^2_j(\Omega),k=1,\ldots,j-1;l=1,\ldots ,k\right\}
\end{equation}
\end{proposition}
\begin{proof}
Without loss of generality assume $\xi=0.$ The method of proof follows the idea of the proof of Proposition \ref{inrepunktremovethm}.
There exists $f_j\in A_j^2(U\setminus\{0\})$ such that
$f(z)=\sum_{k=0}^{j-1} f_k(z)\bar{z}^k, z\in U\setminus\{0\}.$
By Proposition \ref{propositionpessoa2}
there exists $\delta >0, M>0$ such that for $\abs{z}<\delta,$
\begin{equation}
\abs{f_{j-1}(z)} \leq =\frac{\abs{\partial_{\bar{z}}^{j-1} f(z)}}{(j-1)!}\leq 
M\abs{z}^{-j}
\end{equation}
Furthermore, if near the origin $\abs{f_{k+l}(z)}\leq M_l\abs{z}^{-k-l-1}$ for $l=1,\ldots, j-k-1,$ for constants $M_l$ there 
are  
positive constants $N_1,N_2$ 
such that near the origin
\begin{multline}
\abs{f_{k}(z)}\leq 
N_1\abs{z}^{-k-1}+N_2\abs{z}^{-k-1}
\leq (N_1+N_2)\abs{z}^{-k-1}
\end{multline}
Hence near the origin
\begin{equation}
\abs{f_{k}(z)}\leq M_k\abs{z}^{-k-1} \mbox{ for }l=1,\ldots, j-k-1,
\end{equation}
The Laurent series of $f_k$ thus takes the form
\begin{equation}
f_{k}(z)=\sum_{l=-(k+1)}^\infty a_{l,k}z^{l}
\end{equation}
for $l=1,\ldots, j-k-1.$
Then there exists a $r>0$ and a bounded measurable function $g_k$ such that
\begin{equation}
\bar{z}^k f_k(z)=a_{-(k+1),k}\frac{\bar{z}^k}{z^{k+1}} +g_k(z), \quad \abs{z}<r
\end{equation} 
This implies that the function
\begin{equation}
h(z):=f(z)-\sum_{k=0}^{j-1} g_k=\frac{1}{z}\sum_{k=0}^{j-1}a_{-(k+1),k}\frac{\bar{z}^k}{z^{k}}, \quad \abs{z}<r
\end{equation} 
belongs to $A_q^2(\{\abs{z}<r\})$.
If for some $k=0,\ldots,j-1,$ $a_{-(k+1),k}\neq 0$ then there exists $M>0$ and a sector $\sigma\subset \{\abs{z}<r\}$
such that
$\abs{h(z)}\geq M/\abs{z},$ $z\in \sigma,$ which can be realized by rewriting
$\sum_{k=0}^{j-1}a_{-(k+1),k}\frac{\bar{z}^k}{z^{k}}=\sum_{k=0}^{j-1}a_{-(k+1),k}\exp(-2i\theta k),$ $z=\abs{z}\exp(i\theta).$
This would contradict the $L^2$-property of $h.$
This implies that $a_{-(k+1),k}=0$ for $k=0,\ldots,j-1,$ i.e.\
\begin{equation}
f_{k}(z)=\sum_{l=-k}^\infty a_{l,k}z^{l}
\end{equation}
which together with the representation
$f(z)=\sum_{k=0}^{j-1} f_j(z)\bar{z}^k, z\in U\setminus\{0\}$
yields
\begin{equation}
f(z)=\sum_{k=0}^{j-1} \bar{z}^k f_k(z)=\sum_{k=1}^{j-1}\sum_{l=0}^{k-1} a_{l-k,k}\frac{\bar{z}^k}{z^{k-l}} +\sum_{k=0}^{j-1}\bar{z}^k \phi_k(z)
\end{equation}
for holomorphic $\phi_k(z)=\sum_{l=0}^\infty a_{l,k}z^l,$ $k=0,\ldots ,j-1.$
Since $U$ is bounded $\psi:=\sum_{k=0}^{j-1}\bar{z}^k \phi_k(z)$ belongs to $A^2_q(U).$
This completes the proof.
\end{proof}
The usual definition of logarithmic capacity usually involves a limit of normalized so-called energy integrals.
Sometimes the following definition is used in the literature (see e.g.\ 
Axler, Conway \& Mcdonald \cite{axlerconway})
\begin{definition}
Given a metric space $X$, a $\sigma$-algebra, $C(X),$ on $X$ is a collection of subsets $B$ such that  
$X\in B,$ and the collection is closed with respect to complements and countable 
unions. 
Denote by $B(X)$
the smallest $\sigma$-algebra in $X$ that contains all open subsets of $X$. 
When $X$ is separable $B(X)$ is $\sigma$-algebra
generated by the open (or closed) balls of $X$.
A {\em measure} $\mu$ on a $X$ is a non-negative map on $C(X)$
such that $\mu(\emptyset)=0$ and $\mu$ is countably additive, i.e.\ 
for any countable collection of pairwise disjoint sets $\{ A_i\}_{i\in I}$ in $B(X)$, 
we have $\mu\left(\bigcup_{i=1}^\infty\right)=\sum \mu(A_i).$
$\mu$ is called a normalized measure if $\mu(X)=1.$
A measure $\mu$ on $X$ is called a {\em Borel regular measure} if for every 
set $A'\in B(X)$ and every $\mu$-measurable set $A\subset X$ 
we have $\mu(A\cap A')+\mu(A\setminus A'),$ and 
there exists a set $A\in B(X),$ such that $\mu(A')=\mu(A).$
A compact set $K\subset\C$ is said to have zero logarithmic capacity if 
\begin{equation}
\sup_{z\in \C}\left\{ \int_K\log \frac{1}{\abs{z-\zeta}}dm(\zeta)\right\}=\infty
\end{equation}
for every normalized measure $m$ supported in $K.$
\end{definition}
F\'ekete \cite{fekete} introduced the notion of so-called transfinite diameter and proved that for a {\em compact set}
that notion can equivalently be used to obtain the same number as the logarithmic capacity. Taking the supremum over all compacts in a nonempty set can then be used as an equivalent definition as follows.
\begin{definition} 
For a nonempty compact subset $K\subset\C$ define 
\begin{equation}
c(K):=\lim_{m} \max_{p_1,\ldots, p_n\in K} \Pi_{1\leq j<k\leq l}^m \abs{p_k-p_j}^{\frac{2}{n(n-1)}}
\end{equation} 
The {\em logarithmic capacity}\index{Logarithmic capacity}, $c(E),$ of a nonempty set $E\subset\C$ is defined as
\begin{equation}
c(E):=\sup \{ c(K)\colon K\subset E, K\mbox{ compact}\}, \quad c(\emptyset):=0
\end{equation} 
\end{definition}
We mention also the Ahlfors (or analytic) capacity.
\begin{definition}\index{Ahlfors capacity}
For a compact set $K\subset \C$ the (analytic) Ahlfors capacity of $K$ is defined as
\begin{multline}
\gamma(K):=\\\sup\{\abs{\lim_{z\to \infty} z(f(z)-f(\infty))}\colon f\in H^\infty(\C\setminus K),
\norm{f}_\infty =1, f(\infty)=\lim_{z\to \infty} f(z)=0\}
\end{multline}
For an arbitrary set $E$ we define
\begin{equation}
\gamma(E):=\sup_{K\subset E, K\mbox{ compact}} \gamma( K)
\end{equation} 
\end{definition}
Ahlfors proved that
a set $E\subset\C$ is removable for $H^\infty$ if and only if $\gamma(E)=0$, in particular
there exists an open neighborhood $\Omega$ of $E$ such that $H^{\infty}(\Omega)=H^\infty(\Omega\setminus E),$
and that implies that $H^\infty(\C\setminus E)=\{0\}.$

\begin{theorem}[See e.g.\ Carleson p.73]\label{carlesonthm}
Let $\Omega$ be a connected domain in $\C$ containing the point at infinity.
Then $A^2_q(\Omega)=\{0\}$ iff $c(\Omega^c)=0,$ where $\Omega^c$ denoted the compact complement of $\Omega.$
\end{theorem}
\begin{proof}
Assume $c(\Omega^c) >0.$ Then $\Omega^c$ contain non-intersecting closed subsets
$E_j,$ such that $c(E_j)>0,$ $j\in \Z_+.$ Let $\mu_j$ be distributions of unit mass on $E_j,$ with bounded logarithmic potentials 
and set $\mu:=\mu_1-\mu_2$ and
\begin{equation}
f(z)=\int_{\Omega^c}\frac{d\mu(\zeta)}{\zeta -z} =\frac{2}{z^2} +\cdots,\quad \abs{z}>R
\end{equation} 
Then for sufficiently large $R$
\begin{equation}
f(z)=\int_{\abs{z}>R}\abs{f(z)}^2 dx\wedge dy <\infty
\end{equation} 
and for a constant $k>0$
\begin{multline}
\int_{\abs{z}\leq R}\abs{f(z)}^2 dx\wedge dy \leq \int_{\Omega^c} \int_{\Omega^c}  \abs{d\mu(\zeta)} d\mu(\zeta')
\int_{\abs{z}\leq R} \frac{dx\wedge dy}{\abs{z-\zeta'}\abs{z-\zeta}}\\
\leq k \int_{\Omega^c} \int_{\Omega^c} \log\abs{\frac{3R}{\zeta -\zeta'}}\abs{d\mu(\zeta)}\abs{d\mu(\zeta')}<\infty
\end{multline}
Let $\gamma$ denote a system of analytic curves enclosing $E_1.$ Then
\begin{equation}
\frac{1}{2\pi i}\int_\gamma f(z)dz =-\mu_1(E_1)=-1
\end{equation}
thus $f(z)\not\equiv 0.$ This proves that $A^2_q(\Omega)\neq \{0\} \rightarrow c(\Omega^c)=0.$
For the converse suppose that
$c(\Omega^c)=0.$ Choose $\Omega_j,$ $\Omega_j\subset \Omega_{j+1},$ $j\in \Z_+,$ such that $\Omega_j\to \Omega,$
where each $\Omega_j$ is bounded by a finite number of analytic curves. Let
$g_j$ denote the Green function of $\Omega_j$ with pole at $\infty$ and let $h_j$ denote the conjugate function
of $g_j.$ 
For a compact set $E$ denote by $\mathcal{M}(E)$ the collection of all
positive unit Borel measures $\mu$ supported on $E.$
The Green function $g_E$ is defined as $V_E -U^{\mu_E}(z)$ where for a positive normalized Borel measure $\mu,$
$U^\mu:=\int \log\frac{1}{z-\zeta} d\mu(\zeta),$ 
\begin{equation}
V_E:=\inf \{ I(\mu) \colon \mu\in \mathcal{M}(E)\},\quad I(\mu):=\int U^\mu d\mu=\int\int \log\frac{1}{\abs{z-\zeta}} d\mu(\zeta)d\mu(z)
\end{equation}
and where $\mu_E$ denotes the measure for which infimum is attained in the definition of $V_E$
The existence of $\mu_E$ can be realized as follows.
By the Banach-Alaoglu Theorem $\mathcal{M}(E)$ is compact in the weak star topology (this is also a version of the Helly's selection theorem).
Let $\{\mu_k\}_{k\in \Z_+}$ be a sequence in $\mathcal{M}(E)$ such that $\lim_{k\to \infty} I(\mu_k)=V_E$ and let 
$\mu_E$ denote the weak-star cluster point of the sequence
$\{\mu_k\}_{k\in \Z_+}.$ 
Since $1/\abs{z-\zeta}$ is lower semi-continuous with respect to $\zeta$, we have
by the monotone convergence theorem
\begin{equation}
\int_E f d\mu_E \leq \liminf_{k\to\infty} \int_E f d\mu_k
\end{equation} 
So $U^{\mu_E}(z)\leq \liminf_{k\to\infty} U^{\mu_k}(z)$
and since also $\mu_k\times \mu_k$ converges weak-star to $\mu_E\times \mu_E$ we have
$I(\mu_E)\leq \liminf_{k\to\infty} I(\mu_k).$ Hence there exists $\mu_E$ such that $I(\mu_E)=V_E.$
We point out that sometimes the logarithmic capacity of a compact set $E$ is in fact defined as 
$c(E):=\exp(-V_E),$ this can be shown to be equivalent to the
the definition via the transfinite diameter.
If $V_E =\infty$, we have $\exp(-V_E)=0,$ such sets $E$ are called {\em polar sets.}\index{Polar set}
The Green function satisfies:\\
(a) $g_j$ is harmonic on $\Omega \setminus \{\infty\}.$\\ 
(b) $g_j(z)-\log\abs{z}=O(1),$ as $z\to \infty.$\\
(c) $g_j(z)\to \infty$ as $\Omega_j \ni z\to \zeta \in \partial \Omega_j$ for all $\zeta\in \partial\Omega_j\setminus S$, for 
a set $S$ of zero capacity. \\
The condition (b) implies that
\begin{equation}
g_j(z)-\log\abs{z} \to V_{\Omega_j} \to \log \frac{1}{c(\Omega_j)},\mbox{ as }z\to \infty
\end{equation}
One can deduce that
$c(\Omega^c)=0$ if and only if $g_j(z)\to \infty,$ as $j\to \infty,$ for each $z\in \Omega,$ uniformly on inside domains.
Introduce the coordinates $\zeta:=\xi +i\eta=g_j+ih_j$ in $\Omega_j$
and set $U_j:=\{ z: 0<g_j <1\}.$ For $f\in A^2_1(\Omega)$ we have
\begin{equation}
\int_{U_j} \abs{f(z)}^2 dx\wedge dy =\int_{0}^{2\pi} \frac{\abs{f}^2 d\eta}{(\partial_\eta g_j)^2}=\epsilon_j \to 0, j\to \infty
\end{equation} 
Also
\begin{multline}
\int_{0}^{1} d\xi \left( \int_{g_j=\xi} \abs{f(z)}\abs{dz}\right)^2\leq
\int_{0}^{1} d\xi \left( \int_{g_j=\xi} \frac{\partial g_j}{\partial n}\abs{dz} \cdot 
\int_{g_j=\xi} \frac{\abs{f}^2}{\partial g_j/\partial n}\abs{dz} \right)\\
=2\pi \int_{0}^{1}\int_{0}^{2\pi} \frac{\abs{f}^2}{(\partial_\eta g_j)^2}d\eta d\xi =2\pi \epsilon_j
\end{multline} 
Hence there exists $\xi_j\in (0,1)$ such that
\begin{equation}
\int_{g_j=\xi_j} \abs{f(z)}\abs{dz} \leq \sqrt{2\pi \epsilon_j}
\end{equation} 
The curves $g_j =\xi_j$ approach $\Omega^c$ as $j\to \infty.$ Let $p_0\in \Omega.$ Since
$f(\infty)=0,$ we have $j_0,$ and a constant $c_1$, such that for $j>j_0$
\begin{equation}
f(p_0)=\frac{1}{2\pi i}  \int_{\{g_j=\xi_j\}} \frac{f(z)dz}{z-p_0}\leq c_1 \sqrt{2\pi \epsilon_j}\to 0
\end{equation} 
Hence $f\equiv 0.$
This completes the proof.
\end{proof}
We can rephrase the last theorem by saying that if $K\subset \C$ is a compact subset 
then $A^2_1(\C\setminus K)$ is trivial iff $c(K)=0.$

\begin{proposition}
Let $K\subset\C$ be a compact subset and let $U\subset\C$ be an open neighborhood of $K.$
Then $K$ has zero logarithmic capacity iff
every function in $A^2_1(U\setminus K)$
has an analytic extension to $U.$
\end{proposition}
\begin{proof}
($\Leftarrow$) 
By Theorem \ref{carlesonthm} $c(K)=0$ iff $A^2_1(\C\setminus K) =\{0\}.$
Assume each $f\in A_1^2(\C\setminus K)$ has an extension to $A^2_1(\C).$
Let $g\in A^2_1(\C)$ and write $g(z)=\sum_j c_j z^j.$
we have $A_1^2(\C)=\{0\}$, so by Theorem \ref{carlesonthm}
$c(K)=0.$\\ 
($\Rightarrow$)
Suppose $A^2_1(\C\setminus K) =\{0\}.$ Any $h\in A^2_1(\C\setminus K)$ can be written as $h=h_1+h_2$ where $h_1\in \mathscr{O}(U)$,
$h_2\in \mathscr{O}((\C\cup\{\infty\})\setminus K),$ $h_2(\infty)=0.$
We show that $h\in \mathscr{O}(U).$
Let $R>0$ such that $K\subset \{\abs{z}<R\}$ and let $g\in \mathscr{O}(\{\abs{z}>R\}\cup\{\infty\}.$ 
For $\zeta\in \C\setminus K$ we have
\begin{equation}
\lim_{z\to \infty} z\left(\frac{h_2(z)-h_2(\zeta)}{z-\zeta}\right) =-h_2(\zeta)
\end{equation}
If $h_2\not\equiv 0$, choose $p_0,q_0\in (\C\setminus K)\cap \{ h_2\neq 0\}$ and set
\begin{equation}
g(z)=\frac{1}{h_2(p_0)}\left(\frac{h_2(z)-h_2(p_0)}{z-p_0}\right) -\frac{1}{h_2(q_0)}\left(\frac{h_2(z)-h_2(q_0)}{z-q_0}\right)
\end{equation}
Then $g(\infty)=\lim_{z\to \infty} g(z)=0,$ and $g'(\infty)=\lim_{z\to \infty} zg(z)=0.$
The power series of $g$ at $\infty$ takes the form $g(\infty)=c_0 +\sum_{j=1}^\infty c_j z^{-j},$ for complex constants $c_j,$
so $g(\infty)=0$ and  $g'(\infty)=0$ implies that $g\in A^2_1(\{\abs{z}>R\}),$
since $c_0=g(\infty)$ and $c_1=g'(\infty).$ 
Thus $g\in A^2_1(\{\abs{z}>R\})$ if $K\subset \{ \abs{z}<R\}.$ If $V\Subset U$ is an open neighborhood of $K$
(where $\Subset$ denotes relatively compact) then 
$h_2=h-h_1\in A_1^2(V\setminus K),$ thus $g\in A_1^2(V\setminus K).$ This implies $g\in A_1^2(\C\setminus K)$
which by the hypothesis forces $g\equiv 0.$
This implies
\begin{equation}
\frac{1}{h_2(p_0)}\left(\frac{h_2(z)-h_2(p_0)}{z-p_0}\right) =\frac{1}{h_2(q_0)}\left(\frac{h_2(z)-h_2(q_0)}{z-q_0}\right)
\end{equation}
Solving for $h_2(z)$ shows that $h_2(z)$ is rational with precisely one pole, denote the pole by $p_1.$
If $p_1\notin U$ then $h$ is holomorphic on $U$. If $p_1\in U$ then $h$ is holomorphic on $U\setminus \{p_1\}$ with pole at $p_1.$
Since $h\in A_1^2(U\setminus K)$, $p_1\in K$, and since $K$ has zero area $h\in A_1^2(U\setminus \{p_1\}).$
By Proposition \ref{frageteckprop},
any isolated point $p_1$ of $\partial U$ belongs to $\partial_{2-r} U,$ hence $h$ extends to a holomorphic function on $U.$ 
This completes the proof.
\end{proof}

\begin{theorem}[Kouchekian \cite{kouchekian}]\label{kouchekthm}
For a domain $\Omega\subseteq\C$,
$A_1^2(\Omega)$ is trivial iff $c(\C\setminus \Omega) = 0$.
\end{theorem}
\begin{proof}
Let $K\subset \C\setminus \Omega$ be a compact subset. Since $L^2(\C\setminus K)\subset L^2(\Omega)$ we have that
$L^2(\C\setminus K)$ is trivial. By Theorem \ref{carlesonthm} this happens iff $c(K)=0$.
Since $K$ was arbitrary $c(\C\setminus \Omega)=0.$
Conversely assume $c(\C\setminus \Omega)=0.$ In particular, $\Omega$ is unbounded.
Let $U$ be a neighborhood of infinity. Then $c(\C\setminus \Omega)=0$
implies that $c(\mbox{cl}U\setminus \Omega)=0.$
Since $c(K)=0$ then $L^2(\C\setminus K)$ is trivial. This completes the proof.
\end{proof}
This extends partially 
to the polyanalytic case.
\begin{theorem}\label{kouchech1}
Let $\Omega\subseteq\C$ be a domain.
Then 
$A_q^2(\Omega)$ is trivial only if $c(\C\setminus \Omega) = 0$.
\end{theorem}
\begin{proof}
Assume $c(\C\setminus \Omega) > 0$.
In particular, $\C\setminus\Omega\neq \emptyset$ since isolated points have zero logarithmic capacity.
Let $p_0\in \C\setminus\Omega.$
By Theorem \ref{kouchech1}
there exists a nonconstant function $a_0(z)\in A^2_1(\Omega).$
Setting $f(z)=(\bar{z}-\bar{p_0})^{q-1} \frac{1}{(z-p_0)^{q-1}}a_0(z),$
we have $\norm{f}_{L^2(\Omega)}=
\norm{a_0}_{L^2(\Omega)}<\infty.$
Hence the function $f$
is nonconstant and belongs to $A_q^2(\Omega)$.
The contrapositive of what we have proved is the following.
\begin{equation}
(A_{q}^2(\Omega)\mbox{ is trivial}) \Rightarrow (c(\C\setminus \Omega) = 0)
\end{equation}
 This completes the proof. 
\end{proof}
For more specific necessary conditions, see Theorem 2.8 of Pessoa \cite{pessoa2016}.
The following shows that a set, $E$, of zero logarithmic capacity must have zero $L^2$-capacity with respect to
$\partial_{\bar{z}},$ i.e.\ $\mbox{Cap}-L^2_{\partial_{\bar{z}}}(E)=0.$
\begin{theorem}
For a domain $\Omega\subset\C,$ a point $p_0\in \partial_{2-r}^1 \Omega$ if and only if
there exists a $\delta>0$ such that $\{ \abs{z-p_0}\leq \delta\}\setminus \Omega$ has zero logarithmic capacity.
\end{theorem}
\begin{proof}
For $p_0\in \partial_{2-r}^1 \Omega$ we have by Proposition \ref{frageteckprop} 
 that there exists a $\delta>0$ such that $\{ \abs{z-p_0}\leq \delta\}
\subset \Omega.$ Since $\partial_{2-r}^1 \Omega$ is relatively open in $\partial\Omega$ we have for sufficiently small 
$\delta$ that 
$\{ \abs{z-p_0}\leq \delta\}\cap\Omega \subset \partial_{2-r}^1 \Omega,$ thus
$\{ \abs{z-p_0}\leq \delta\} \setminus  \Omega \cup \partial^1_{2-r} \Omega.$  
So $\{ \abs{z-p_0}\leq \delta\} \setminus  \Omega \subset \partial^1_{2-r} \Omega.$
In particular, $\{ \abs{z-p_0}\leq \delta\} \setminus  \Omega$ is compact in 
$\Omega \cup \partial^1_{2-r} \Omega.$
Then each $f\in A_1^2\left((\Omega \cup \partial^1_{2-r} \Omega)\setminus(\{ \abs{z-p_0}\leq \delta\} \setminus  \Omega)\right)$
has analytic extension to $\Omega \cup \partial^1_{2-r} \Omega.$ Hence
$c(\{ \abs{z-p_0}\leq \delta\} \setminus  \Omega)=0.$
For the converse assume $\delta>0$ is such that $c(\{ \abs{z-p_0}\leq \delta\} \setminus  \Omega)=0.$
In particular, 
each bounded analytic function on
the bounded $\Omega \cup \partial^1_{2-r} (\{ \abs{z-p_0}\leq \delta\} \setminus  \Omega)$ 
has holomorphic extension to $\Omega \cup \partial^1_{2-r} \Omega$,
implying that $H^\infty(\C\setminus (\{ \abs{z-p_0}\leq \delta\} \setminus  \Omega))$ is trivial.
Next note that if $J\subseteq L$ is a connected component of compact set $L$, 
then the conformal map of $\C\setminus J$ onto the interior of the unit disc,
becomes a nonzero bounded analytic function on $\C\setminus K$
(see e.g.\ Gamelin \cite{gamelinbok}, p.198).
This implies that $\{ \abs{z-p_0}\leq \delta\} \setminus  \Omega$ is totally disconnected 
This implies that there exists an open neighborhood $V\subset \{ \abs{z-p_0}\leq \delta\}$ of $p_0$ such that
$V\cap \{\abs{z-p_0}\leq \delta\} \setminus  \Omega$ is compact.
Then
$V\cap ( \{\abs{z-p_0}\leq \delta\} \setminus  \Omega)$ is compact in 
$\Omega \cup V$ and the difference $(\Omega\cup V)\setminus (V\cap (\{ \abs{z-p_0}\leq \delta\} \setminus  \Omega))=\Omega.$
Since $V\cap ( \{\abs{z-p_0}\leq \delta\} \setminus  \Omega)$ has zero logarithmic capacity each function
in $A^2_1(\Omega)$ has analytic extension to $\Omega\cup V,$ hence $p_0\in \partial_{2-r}^1 \Omega.$
This completes the proof.
\end{proof}

\chapter{Polyentire functions that allow special factorization}\label{polyentiresec}

\section{Conditions on the zero set}
With regards to preliminaries from Nevanlinna theory in this chapter, see also Section \ref{nevanlinnasec}.
We shall use the following lemma which we state without proof 
(to obtain the stated representation one can, 
staring from the collection $\varphi_0(z),\ldots,\varphi_{n-1}(z)$ 
repeatedly go through the collection, and 
each time that
one $\varphi_k(z)$ turns out to be a linear combination of
the functions $\varphi_{k+1}(z),\ldots,\varphi_{n-1}(z)$,
replace this function 
appropriately and repeat the process until, if necessary after relabeling,
one obtains the wanted representation).
\begin{lemma}\label{balkentirepalem1}
	Any entire reduced polyanalytic function $f(z)=\sum_{j=0}^{n-1}\abs{z}^{2j}\varphi_j(z)$ 
	has a representation of the form
	\begin{equation}\label{balkentireaekv1}
	f(z)=\sum_{k=1}^s p_k(\abs{z}^2)\psi_k(z)\end{equation}
	for entire $\psi_j,$
	such that $\psi_j(z)$ $j=0,\ldots,n,$ linearly independent and
	where each $\psi_k$ coincides with some $\varphi_j$ and $\varphi_{n-1}=\psi_s$, $\psi_k\not\equiv 0$, $k=1,\ldots,s$
	and where $p_k(t)$ is a polynomial in $t$ of degree
	$\alpha_k$ such that
	\begin{equation}
	0\leq \alpha_1< \cdots <\alpha_s =n-1
	\end{equation}
\end{lemma}
We further assume that \begin{equation}
\tau(\psi_k)=1,\quad k=1,\ldots,s
\end{equation}
where $\tau(\psi_k)$ is the first nonzero coefficient in the expansion of $\varphi_k$
with respect to $z.$
\\
For a reduced entire polyanalytic function with representation given by Eqn.(\ref{balkentireaekv1})
we set
\begin{equation}B=B(z)=\{\psi_1,\ldots,\psi_s\}
\end{equation}
\begin{equation}
T(r,B):=\frac{1}{2\pi}\int_0^{2\pi} u(r\exp(i\varphi))d\varphi
\end{equation}
where
\begin{equation}\label{manuilovekv6}
u(z)=u(z,B):=\max_k \{\ln \abs{\psi_k(z)}\}
\end{equation}
Then it is well-known from Nevanlinna theory (the Jensen formula)
that
\begin{equation}\label{jensenformek7}
\frac{1}{2\pi}\int_0^{2\pi}\ln \abs{f(\exp(i\varphi))}d\varphi=N(r,0,f)-N(r,\infty,f)+\ln \abs{\tau(f)}
\end{equation}
where $f(z)$ is a function meromorphic on $\{\abs{z}<r\}$ and
$\tau(f)$ the first nonzero coefficient in the Laurent expansion at $z=0,$
$N(r,a,f)$ the Nevanlinna function
\begin{equation}
N(r,a,f)=\int_0^r\frac{n(t,a)-n(0,a)}{t}dt +n(0,a)\ln r
\end{equation}
For the holomorphic function $\psi(z)$ on $\{\abs{z}<r\}$ we have
$N(r,\infty,\psi)=0$ and by (the Jensen formula) Eqn.(\ref{jensenformek7})
we have 
\begin{equation}\label{manuilovekv9}
\frac{1}{2\pi}\int_0^{2\pi} \ln \abs{\psi(r\exp(i\varphi))}d\varphi =N(r,0,\psi)+\ln \abs{\tau(\psi)}
\end{equation}
Let $\psi_k(z)$ and $\psi_s(z)$ be entire functions and define the meromorphic function
\begin{equation}\label{manuilovekv10}
f_{ks}(z)=\frac{\psi_k(z)}{\psi_s(z)}
\end{equation}
By Eqn.(\ref{manuilovekv6}) we obtain for the meromorphic curve $B_{ks}=\{\psi_k,\psi_s\}$
\begin{equation}\label{manuilovekv12}
u_{ks}(z)\equiv u(z,B_{ks})=\ln^+\abs{f_{ks}(z)}+\ln\abs{\psi_s(z)}
\end{equation}
Now the Nevanlinna characteristic satisfies
\begin{equation}\label{manuilovekv13}
T(r,f)=m(r,\infty,f)+N(r,\infty,f)
\end{equation}
where
\begin{equation}\label{manuilovekv14}
m(r,\infty,f)\equiv m(r,f)=\frac{1}{2\pi}\int_0^{2\pi} \ln^+ \abs{f(r\exp(i\varphi))}d\varphi
\end{equation}
By Eqn.(\ref{manuilovekv12}) and Eqn.(\ref{manuilovekv9})
we have
\begin{multline}\label{manuilovekv15}
T(r,f_{ks})=\frac{1}{2\pi}\int_0^{2\pi} u_{ks}(r\exp(i\varphi))d\varphi -\frac{1}{2\pi}\int_+^{2\pi}\ln \abs{\psi(r\exp(i\varphi))}d\varphi
+\\
N(r,\infty,f_{ks})=T(r,B_{ks})+N(r,\infty,f_{ks})-N(r,0,\psi_s)-\ln\abs{\tau(\psi_s)}
\end{multline}
Denote by $\psi_{ks}$ the analogue of a greatest common divisor of $\psi_k$ and
$\psi_s$, in the sense that
it is an entire function such that there exists
\begin{equation}
\psi_k(z)=\sigma_k(z)\psi_{ks}(z),\quad \psi_s(z)=\sigma_s(z)\psi_{ks}(z)  
\end{equation}
for $\sigma_k(z),$ $\sigma_s(z)$ without common zeros, and where we further assume the normalization
$\tau(\psi_{ks})=1.$
By Eqn.(\ref{manuilovekv10}) we have
\begin{equation}\label{manuilovekv16}
N(r,0,\psi_s)=N(r,\infty,f_{ks})+N(r,0,\psi_{ks})
\end{equation}
By Eqn.(\ref{manuilovekv15}) and Eqn.(\ref{manuilovekv16}) we obtain
\begin{equation}\label{manuilovekv17}
T(r,f_{ks})\leq T(r,B_{ks})-N(r,0,\psi_{ks})-\ln\abs{\tau(\psi_{ks})}
\end{equation}
By Eqn.(\ref{manuilovekv6}) we have
\begin{equation}
u(r,B_{ks})\leq u(r,B)
\end{equation}
which yields
\begin{equation}\label{manuilovekv19}
T(r,B_{ks})\leq T(r,B)
\end{equation}
For $r\geq 1$ Eqn.(\ref{manuilovekv17}) and Eqn.(\ref{manuilovekv19}) imply
\begin{equation}\label{manuilovekv20}
T\left(r,\frac{\psi_k}{\psi_s}\right)\leq T(r,B).
\end{equation}
\begin{lemma}\label{haymanIIlem}
	If $T(r)$ is continuous, increasing and $T(r)\geq 1$ for $r_0\leq r <\infty,$ then
	\begin{equation}
	T\left(r+\frac{1}{T(r)}\right)<2T(r)
	\end{equation}
	outside a set $E_0$ of $r$ which has linear measure at most $2.$
\end{lemma}
\begin{proof}
	Let $r_1$ be the lower bound of all $r>r_0$ such that 
	\begin{equation}\label{haylem41andraekv}
	T\left(r+\frac{1}{T(r)}\right)\geq 2T(r)
	\end{equation}
	If there does not exist such $r$ we are done. We define by
	iteration a sequence of numbers $r_n$ as follows: given $n$, let $r'_n=r_n +1/T(r_n)$
	and define $r_{n+1}$ as the lower bound of all $r\geq r_n'$ for which
	Eqn.(\ref{haylem41andraekv}) holds true. This yields a sequence $\{r_n\}_{n}$
	such that, by continuity of $T(r)$, Eqn.(\ref{haylem41andraekv}) holds true for
	$r=r_n,$ $n\in \Z_+.$ Thus $r_n$ belongs to $E_0.$ By definition there are no points of $E_0$ in $(r'_n,r_{n+1})$
	so the set of closed intervals $[r_n,r']$ contains $E_0$. 
	If there exists infinitely many $r_n$
	then $r_n$ cannot tend to a finite limit $r$ since then $r_n<r_n'\leq r_{n+1},$ 
	so also $r_n'\to r$. However, since $T(r)$ is increasing we have
	for all $n$, $r'_n-r_n=1/T(r_n)\geq 1/T(r)>0$, which gives a contradiction.
	Since $r_0$ belongs to $E_0$ we have
	\begin{equation}
	T(r'_n)=T(r_n+1/T(r_n))\geq 2T(r_n)
	\end{equation}
	Since $T(r)\geq 1$ we have
	\begin{equation}
	T(r_{n+1})\geq 2T(r_n)\geq\cdots\geq 2^nT(r_1)\geq 2^n
	\end{equation}
	Thus
	\begin{equation}
	\sum_{n=1}^\infty(r'_n -r_n)=\sum_{n=1}^\infty \frac{1}{T(r_n)}\leq \sum_{n=1}^\infty 2^{1-n}=2
	\end{equation}
	This completes the proof.
\end{proof}
We are now ready to prove the following theorem.
\begin{theorem}\label{entirepathm00}
	If $f(z)=\sum_{j=0}^{n-1}f_j(z)\abs{z}^{2j}$ is a reduced polyentire function with bounded  
	zero set then there exists an entire function $g(z)$ and a polynomial 
	$P(z,\abs{z}^2)$ such that $f(z)=\exp(g(z))P(z,\bar{z}).$
\end{theorem}
\begin{proof}
	By Proposition \ref{balkunboundedzerolem} we know that if
	If $f_{n-1}(z)$ has unbounded zero set then so does $f(z).$
	Since by assumption $f(z)$ has bounded zero set we have by the contrapositive version of Lemma \ref{balkentirepalem1} that
	$f_{n-1}(z)$ is an entire function with bounded zero set thus 	
	there exists a complex polynomial $\phi_{n-1}(z)$ and an entire $g(z)$ such that 
	\begin{equation}\label{manuilovekv23}
	f_{n-1}(z)=\phi_{n-1}(z)\exp(g(z))
	\end{equation}
	Set
	\begin{equation}\label{manuilovekv25}
	f_j(z)\exp(-g(z))=\varphi(z),\quad j=0,\ldots,n-2
	\end{equation}
	which implies that the following function has unbounded zero set
	\begin{equation}\label{manuilovekv26}
	\Psi(z):=f(z)\exp(-g(z))=\sum_{j=0}^{n-1} \abs{z}^{2j}\varphi_j(z)
	\end{equation}
	By Lemma \ref{balkentirepalem1}
	we have
	\begin{equation}\label{balkentirepaekv3}
	f(z)=\sum_{k=1}^s p_k(\abs{z}^2)\psi_k(z)
	\end{equation}
	for entire $\psi_j,$
	such that $\psi_j(z)$ $j=0,\ldots,n,$ linearly independent and
	where each $\psi_k$ coincides with some $\psi_j$ and $\varphi_{n-1}=\psi_s$, $\psi_k\not\equiv 0$, $k=1,\ldots,s$
	and where $p_k(t)$ is a polynomial in $t$ of exact degree
	$\alpha_k$ such that
	$0\leq \alpha_1< \cdots <\alpha_s =n-1.$
	Suppose (in order to reach a contradiction) that at least one of the $\psi_k(z)$ is not a polynomial.
	Let $B=B(z)=\{\psi_1,\ldots,\psi_s\}$ as above.
	We introduce the notation
	\begin{equation}\label{manuilovekv28}
	G(z)=G(z,c):=\{g_1(z),\ldots,g_s(z)\}
	\end{equation}
	\begin{equation}\label{manuilovekv29}
	g_k(z)=g_k(z,c):=p_k(c^2)\psi_k(z),\quad k=1,\ldots,s
	\end{equation}
	Setting
	\begin{equation}\label{manuilovekv31}
	U(z):=\max_k\{\ln \abs{g_k(z)}\}
\end{equation}
we have
\begin{equation}\label{manuilovekv32}
T(r,G)=\frac{1}{2\pi}\int_0^{2\pi} U(r\exp(i\theta))d\theta
\end{equation}
and using Eqn.(\ref{manuilovekv6}) we have for sufficiently large $c$ that there is a constant 
$A$ independent of $c$ such that
\begin{multline}
u(z)=\max_k\{\ln \abs{\psi_k(z)}\}=\max_k \{\ln \abs{g_k(z)}-\ln \abs{p_k(c^2)}\}\\
<\max_k \ln\abs{g_k(z)} +A
\end{multline}
so that $u(z)<U(z)+A.$
Thus for sufficiently large $c$ there is a constant $A$ independent of $c$ such that
\begin{equation}\label{manuilovekv30}
T(c,B)<T(c,G)+A
\end{equation}
Now for sufficiently large $c$
we obtain
\begin{equation}\label{manuilovekv33}
T(c,G)<A\ln T(c,B)+A\ln c
\end{equation}
According to the definitions and remarks given preceding the theorem we obtain
for sufficiently large positive $c$ there is a constant $A$ independent of $c$ such that
\begin{equation}
T(c,B)<A\ln T(c,B)+A\ln c
\end{equation}
By the conditions of the theorem we may choose a sufficiently large $\rho$ such that for $\gamma=\{\abs{z}=\rho\}$ the function 
$\Psi(z)$ is nonzero for $\abs{z}\geq \rho.$
Set
\begin{equation}\label{manuilovekv34}
B_{p_\gamma}\Psi(z)=\frac{1}{2\pi}\Delta_\gamma \mbox{arg}\Psi(z)=\frac{1}{2\pi}\int_\gamma d\mbox{arg}\Psi(z)
\end{equation}
and denote
\begin{equation}
B_{p_\gamma}\Psi(z)=h
\end{equation}
Then for $c\geq \rho,$ we have with $\Gamma:=\{\abs{z}=c\}$
\begin{equation}\label{manuilovekv35}
B_{p_\Gamma}\Psi(z)=h
\end{equation}
Now on $\Gamma$ the function $\Psi(z)$ coincides with the analytic function
\begin{equation}\label{manuilovekv36}
\Phi(z)=\Phi(z,c^2):=\sum_{k=1}^s p_k(c^2)\psi_k(z)
\end{equation}
which yields
\begin{equation}
B_{p_\Gamma}\Phi(z)=h
\end{equation}
Thus $\Phi(z)$ has inside $\Gamma$, $h$ zeros including multiplicity. 
Define
\begin{equation}
M:=\begin{bmatrix}
a_1^{\alpha_1} & a_1^{\alpha_2} &\cdots & a_1^{\alpha_s} \\ 
\vdots & \vdots & \ddots& \vdots\\
a_{s+1}^{\alpha_1} & a_{s+1}^{\alpha_2} &\cdots & a_{s+1}^{\alpha_s} \\ 
\end{bmatrix}
\end{equation}
where $a_1,\ldots,a_{s+1}$ are positive numbers such that $1\leq a_1<\ldots <a_{s+1}.$
Removing from $M$, one row, and denoting the remaining rows $\beta_1,\ldots,\beta_s$, we set
\begin{equation}\label{manuilovekv39}
\sigma[\beta_1,\ldots,\beta_s]:=\begin{bmatrix}
a_{\beta_1}^{\alpha_1} & a_{\beta_1}^{\alpha_2} &\cdots & a_{\beta_1}^{\alpha_s} \\ 
\vdots & \vdots &\ddots & \vdots\\
a_{\beta_s}^{\alpha_1} & a_{\beta_s}^{\alpha_2} &\cdots & a_{\beta_s}^{\alpha_s} \\ 
\end{bmatrix}
\end{equation}
Consider $\sigma_\beta'(\beta_1):=\sigma[\beta_1,\beta']$ as a function of $\alpha_{\beta_1}$ with
$\beta'=(\beta_2,\ldots,\beta_s)$ fixed. This has at most $s$ nonzero coefficients and by Decartes rule of signs it has at most $s-1$
positive roots. But there also exists $s-1$ positive roots since each $a_{\beta_2},\ldots,a_{\beta_s}$ is a root of this polynomial. 
If there is a positive number $a_{\beta_1}$ distinct from the $a_{\beta_2},\ldots,a_{\beta_s}$ then
the matrix $M$ has nonzero determinant. Since the $a_j$ are positive, real, distinct, nonzero the determinant of this square matrix is nonzero 
as a consequence of the Decartes rule of signs.
Setting
\begin{equation}
R:=\max\left\{ \frac{1}{\abs{\sigma[\beta_1,\ldots,\beta_s]}}\right\}
\end{equation}
we have $0<R<\infty.$
Define
\begin{equation}\label{manuilovekv41}
\Phi_j(z):=\Phi(z,c^2a_j),\quad j=1,\ldots,s+1
\end{equation}
each of which has oat most $h$ zeros on $d:=\{\abs{z}\leq c\}.$
By a known theorem of Cartan \cite{cartanentirepa} (see also Goldberg \cite{goldberg1960}, p.295, with $p=s,$ $q=s+1$)
we have
\begin{equation}\label{manuilovekv43}
T(c,G)<\frac{1}{2\pi}\int_0^{2\pi} v(c\exp(i\theta))d\theta +\ln K
\end{equation}
where $K$ is a constant depending only on $s$ and independent of $c, \psi_k(z),$ and where
\begin{equation}
v(z):=\max {1<\beta <s+1} \{ \ln \abs{\Phi_\beta (z)}\}
\end{equation}
Now we have the following relation for the Wronskinan $w$
\begin{equation}
w[\Phi_{\beta_1},\ldots,\Phi_{\beta_s}]\equiv C[\beta_1,\ldots,\beta_s,c]\cdot w[g_1,\ldots,g_s]
\end{equation}
where
\begin{equation}\label{manuilovekv45}
C[\beta_1,\ldots,\beta_s,c]:=
\begin{bmatrix}
\frac{p_1(c^2 a_{\beta_1})}{p_1(c^2)} & \cdots & \frac{p_s(c^2 a_{\beta_1})}{p_s(c^2)}\\ 
\vdots & \ddots  & \vdots\\
\frac{p_1(c^2 a_{\beta_s})}{p_1(c^2)} & \cdots & \frac{p_s(c^2 a_{\beta_s})}{p_s(c^2)}
\end{bmatrix}
\end{equation}
By the properties of the $p_k(z)$ this determinant has a finite limit as $c\to \infty$ given by Eqn.(\ref{manuilovekv39}) 
\begin{equation}\label{manuilovekv46}
\lim_{c\to\infty} C[\beta_1,\ldots,\beta_s,c]=\sigma[\beta_1,\ldots,\beta_s]
\end{equation}
Set
\begin{equation}\label{manuilovekv47}
H(z)=H(z,c):=\frac{\Pi_{j=1}^{s+1}\Phi_j(z)}{w[g_1,\ldots,g_s]}
\end{equation}
\begin{equation}\label{manuilovekv48}
\Theta(z;\beta_1,\ldots,\beta_s;c):=
\frac{1}{C[\beta_1,\ldots,\beta_s;c]}\times
\begin{bmatrix}
1 & \cdots & 1\\
\frac{\Phi'_{\beta_1}}{\Phi_{\beta_1}} & \cdots & \frac{\Phi'_{\beta_s}}{\Phi_{\beta_s}}\\
\frac{\Phi''_{\beta_1}}{\Phi_{\beta_1}} & \cdots & \frac{\Phi''_{\beta_s}}{\Phi_{\beta_s}}\\
\vdots & \ddots  & \vdots\\
\frac{\Phi^{(s-1)}_{\beta_1}}{\Phi_{\beta_1}} & \cdots & \frac{\Phi^{(s-1)}_{\beta_s}}{\Phi_{\beta_s}}\\
\end{bmatrix}
\end{equation}
If $\beta$ is not among the $\beta_1,\ldots,\beta_s$ then
\begin{equation}\label{manuilovekv49}
\Phi_\beta(z)=\Theta(z;\beta_1,\ldots,\beta_s;c)\cdot H(z)
\end{equation}
Set
\begin{equation}\label{manuilovekv50}
W(z):=\max\{\ln \abs{\Theta(z;\beta_1,\ldots,\beta_s;c)}\}
\end{equation}
where the maximum is taken over fixed $z,c$ for all possible combinations of $\beta_1,\ldots,\beta_s.$
By a known result of Cartan \cite{cartanentirepa} (see also Goldberg \cite{goldberg1960}, p.296)
\begin{equation}\label{manuilovekv51}
v(z)=\ln \abs{H(z)}+W(z)
\end{equation}
Using Eqn.(\ref{manuilovekv43}) and Eqn.(\ref{manuilovekv51}) we will obtain a bound for $T(c,B)$
if we can bound from above, for sufficiently large $c$, the integrals
\begin{equation}\label{manuilovekv52}
\frac{1}{2\pi}\int_0^{2\pi}\ln \abs{H(c\exp(i\theta))}d\theta
\end{equation}
\begin{equation}\label{manuilovekv53}
 \frac{1}{2\pi}\int_0^{2\pi}\ln W(c\exp(i\theta))d\theta
\end{equation}
\begin{equation}\label{manuilovekv54}
\psi_s(z)=z^\lambda\chi(z),\quad  \chi(0)=1
\end{equation}
Choose $r_0$ sufficiently small such that $\chi(z)$ is nonzero on $\delta:=\{\abs{z}\leq r_0\}$  
By Rouch\'e's theorem we have
$p_s(c^2) \psi_s(z)$ and $\sum_{k=0}^{s-1} p_k(c^2)\psi_k(z)$
and we can make sure that 
the function $\Phi(z;c^2)$ has precisely $\lambda$ zeros (with multiplicity) on $\delta,$
and may vanish at $z=0$.
Furthermore, the multiplicity of a zero of $\Phi$ at $z=0$ will be $\mu.$
We may assume that for sufficiently large $c$ the number $\mu$ does not depend upon $c.$
So for $\{c_m\}_{m\in \N}$ with $c_m\to \infty$ as $m\to \infty$ we have that the function
$\Phi(z,c^2)$ has at $z=0$, for sufficiently large $m$, a zero of order $\mu_1$ for a constant
a constant $\mu_1\leq \lambda$ and we may after relabeling assume the sequence begins with sufficiently large index
such that 
\begin{equation}\label{manuilovekv55}
\Phi(z;c^2_m)=z^{\mu_1} E(z;c^2)
\end{equation}
where $E(z;c^2)$ is regular and nonzero at $z=0$.
By Eqn.(\ref{manuilovekv36} we obtain
\begin{equation}\label{manuilovekv56}
\frac{z^{\mu_1} E(z;c_m^2) -z^\lambda\chi(z)p_s(c_m^s)}{p_{s-1}(c_m^2)}=\psi(z)+
\sum_{k=1}^{s-2}\frac{p_k(c_m^2)}{p_{s-1}(c_m^2)}\psi_k(z)
\end{equation}
where $\psi_{s-1}(z)\not\equiv 0.$ For $c_m\to \infty$ the right hand side converges on $\delta$ to $\psi_{s-1}(z).$
Thus for large $c_m$ the right hand side has a zero at $z=0$ of order 
not greater than the order of $\psi_{s-1}(z).$
On the other hand the right hand side of Eqn.(\ref{manuilovekv56}) has for large $m$ a zero at 
$z=0$ of order not lower than $\mu_1$ and consequently
$\psi_{s-1}(z)$ has a zero at 
$z=0$ of order not lower than $\mu_1$.
By iterating these arguments we obtain that each $\psi_k(z),$ $k=1,\ldots,s,$
are divisible by $z^{\mu_1}$.
By Eqn.(\ref{manuilovekv36})
this implies that $\Phi(z;c^2)$ is divisible by $z^{\mu_1}$ for all sufficiently large $c$ (i.e.\ not only those
of the sequence $\{c_m\}_m$). This shows that the order of the zero at $z=0$ for large $c$ is idependent of $c.$
We need to consider the following two cases:\\
(I) $\mu=\lambda$:
Then for large $c$ the function $\Phi(z;c^2)$ does not have zeros in $\{0<\abs{z}\leq r_0\}$ and the same is true for $\Phi_j(z),$
$j=1,\ldots,s+1.$
\\
(II) $\mu<\lambda$:
\\
In this case we claim that there exists $A>0$ and $\alpha\geq 0$ (independent of $c$)
such that for sufficiently large $c$ the functions $\Phi(z;c^2)$ do not have zeros in the annulus
\begin{equation}\label{manuilovekv57}
0<\abs{z}\leq Ac^{-\alpha}
\end{equation}
and the same is true for $\Phi_j(z),$
$j=1,\ldots,s+1.$\\
\\
To prove the claim note that by Eqn.(\ref{manuilovekv36}) all the zeros of 
$\Phi(z;c^2)$
$0<\abs{z}\leq r_0$
must satisfy
\begin{equation}\label{manuilovekv58}
z^{\lambda-\mu}+\sum_{k=1}^{s-1}\frac{p_k(c^2)}{p_{s}(c^2)}\chi_k(z)=0
\end{equation}
$\chi_k(z):=\psi_k(z)/(z^\mu \chi(z)),$ where $\chi_k(z),$ $k=1,\ldots,s-1,$ are regular and at least one of them is nonzero at $z=0.$
Introducing $\zeta:=c^{-2}$ and
$\pi_k(\zeta)=\frac{p_k(c^2)}{p_s(c^2)}$ we have by continuity $\pi_k(0)=0$, $k=1,\ldots,s-1.$
By Eqn.(\ref{manuilovekv59})
By Eqn.(\ref{manuilovekv58}) we have
\begin{equation}\label{manuilovekv59}
z^{\lambda-\mu}+\sum_{k=1}^{s-1}\pi_k(\zeta)\chi_k(z)=0
\end{equation}
By the Weierstrass preparation theorem this is equivalent (as an equation in the variable $z$) to
an equation of the form
\begin{equation}\label{manuilovekv60}
z^{\lambda-\mu}+\sum_{j=0}^{\lambda-\mu -1} A_j(\zeta)z^j=0
\end{equation}
for functions $A_j(\zeta)$ analytic near $\zeta=0$ and satisfying $A_j(0)=0,$
$j=0,\ldots,\lambda-\mu-1.$
Decomposing the pseudopolynomial, near $(0,0)$, 
in left hand side into a product of irreducible factors it is possible to deduce
(see Theorem 4.12 of
Fuks \cite{fuks1962}, p.88-89, on representation of analytic surfaces) 
that Eqn.(\ref{manuilovekv59}) determines in a neighborhood of $z=0$ finitely many functions of the form
\begin{equation}
z=d_{1,j}\zeta^{\frac{1}{m_j}}+d_{2,j}\zeta^{\frac{2}{m_j}}+\cdots
\end{equation} 
where $m_j$ are natural numbers such that $\sum_j m_j=\lambda-\mu$
Thus for sufficiently large $c$ the zeros of the functions $\Phi(z;c^2)$
lying in $0<\abs{z}\leq r_0$ satisfy that there exists constants $A,\alpha$ independent of $c$ such that
\begin{equation}\label{manuilovekv61}
\abs{z}\geq \frac{A}{c^\alpha}
\end{equation} 
and the same is true for $\Phi_j(z),$
$j=1,\ldots,s+1,$ and we may choose $A$ such that Eqn.(\ref{manuilovekv60}) to holds simultaneously
for all these finitely many functions with the same $A.$
This proves the claim and 
in fact we have that in both
cases I and II will satisfy Eqn.(\ref{manuilovekv57}), namely in the case I we can use
$\alpha=0,$ $A=r_0.$
\\
\\
Thus if $\psi_s(z)$ has a zero at $z=0$ of order $\lambda$ then there exists an integer $0\leq \mu\leq \lambda$
together with constants $A>0$ and $\alpha\geq 0$ such that for sufficiently large $c$
the functions $\Phi_j$, $j=1,\ldots,s+1$ each have a zero at $z=0$ of order $\mu$
and the other zeros of these functions satisfy Eqn.(\ref{manuilovekv61}).
By Eqn.(\ref{jensenformek7}) and Eqn.(\ref{manuilovekv52}) we have
\begin{equation}
\frac{1}{2\pi} \int_0^{2\pi} \ln \abs{H(c\exp(i\theta))}d\theta =
N(c,0,H)-N(c,\infty,H)+ \ln\abs{\tau(H)}
\end{equation} 
thus
\begin{equation}\label{manuilovekv62}
\frac{1}{2\pi} \int_0^{2\pi} \ln \abs{H(c\exp(i\theta))}d\theta \leq \ln
\abs{\tau(H)}+N(c,0H)
\end{equation} 
But
\begin{equation}
\abs{\tau(H)} =\frac{\abs{\tau(\Phi_1)}\cdots\abs{\tau(\Phi_{s+1})}}{\abs{\tau(w[\psi_1,\ldots,\psi_s])}\cdots\abs{p_1(c^2)\cdots p_s(c^2)}} <Ac^{(82s+2)\alpha_s}
\end{equation} 
thus for sufficiently large $c$ we have
\begin{equation}\label{manuilovekv63}
\ln\abs{\tau(H)}\leq A\ln c
\end{equation} 
Denote by $\bar{n}_j(r,0)$ the number of zeros of $\Phi_j(z)$ in $\{\abs{z}\leq r\}$
and each zero is counted with multiplicity determined by the mimum of the order of the zero and the number $s-1.$
We denote by $\bar{N}_j(r,0)$ the expression
\begin{equation}\label{manuilovekv64}
\int_0^r \frac{\bar{n}_j(t,0)-\bar{n}_j(0,0)}{t}dt +\bar{n}_j(0,0)\ln r
\end{equation} 
By a result of Cartan \cite{cartanentirepa} (see also Goldberg \cite{goldberg1960}, p.296)
we have
\begin{equation}\label{manuilovekv65}
N(r,0,H)\leq \sum_{j=1}^s \bar{N}_j(r,0)
\end{equation}
Indeed, for sufficiently large $c$ we have
\begin{equation}
\bar{N}_j(c,0)\leq \int_{Ac^{-\alpha}}^c \frac{\bar{n}_j(t,0)}{t}dt
+\mu\ln c \leq\\
h\ln \frac{c}{Ac^{-\alpha}} +\lambda \ln c\leq A\ln c
\end{equation}
By Eqn.(\ref{manuilovekv65})
\begin{equation}\label{manuilovekv66}
N(c,0,H)\leq A\ln c
\end{equation}
By Eqn.(\ref{manuilovekv62}), Eqn.(\ref{manuilovekv63}) and
By Eqn.(\ref{manuilovekv66}) we obtain
\begin{equation}\label{manuilovekv67}
\frac{1}{2\pi} \int_0^{2\pi} \ln \abs{H(c\exp(i\theta))}d\theta \leq A\ln c
\end{equation}
By Eqn.(\ref{manuilovekv46}), Eqn.(\ref{manuilovekv48}), Eqn.(\ref{manuilovekv50}),
and Eqn.(\ref{manuilovekv14}) we have
\begin{equation}\label{manuilovekv68}
\frac{1}{2\pi} \int_0^{2\pi} W(c\exp(i\theta))d\theta <
A+\sum_{j=1}^{s+1}\left( m\left( c,\frac{\Phi'_j}{\Phi_j}\right)+\cdots+
m\left( c,\frac{\Phi^{(s-1)}_j}{\Phi_j}\right) \right)
\end{equation}
But for any $r,j,k$ we have
\begin{equation}\label{manuilovekv69}
m\left( r,\frac{\Phi^{(k)}_j}{\Phi_j}\right)=
m\left( r,\Pi_{l=1}^k \frac{\Phi^{(l)}_j}{\Phi_j^{(l-1)}}\right)\leq 
\sum_{l=1}^{k}m\left(r,\frac{\Phi^{(l)}_j}{\Phi_j^{(l-1)}} \right)
\end{equation}
We thus obtain an estimate form above for the integral in Eqn.(\ref{manuilovekv53}) via
\begin{equation}\label{manuilovekv70}
m\left( r,\frac{\Phi^{(k)}_j}{\Phi_j^{(k-1)}}\right),\quad j=1,\ldots,s+1; k=1,\ldots,s-1
\end{equation}
We shall now give two inequalities that will be used in order to estimate the expressions in Eqn.(\ref{manuilovekv70})
and we denote these by (a) and (b) respectively.
\\
(a): 
By a theorem of Nevanlinna \cite{nevanlinna29}, p.61, we have for a meromorphic function $f(z)$ satisfying $\tau(f)=c_\lambda$ and with $0<r<\rho$ that
\begin{multline}\label{manuilovekv71}
m\left( r,\frac{f'}{f}\right)<34+5\ln^+\abs{\lambda}+3\ln^+ \frac{1}{\abs{c_\lambda}}
+\\
7\ln^+ \frac{1}{r} +4\ln^+\rho +
3\ln^+\frac{1}{\rho-r}+4\ln^+ T(\rho,f)
\end{multline}
Now we shall apply Eqn.(\ref{manuilovekv71}) with the assumption that the meromorphic function $f(z)$ depends on
a parameter $c>0$ (it will depend on $\Phi_j^{(k)}(z)$), where $c$ is sufficiently large that 
$\abs{\tau(f)}$ is greater than a positive constant, a finite number $\abs{\lambda}$ and $r$ we will let be a large positive number, and under these assumptions Eqn.(\ref{manuilovekv71})
can be written
\begin{equation}\label{manuilovekv72}
m\left( r,\frac{f'}{f}\right)<A +4\ln^+\rho +
3\ln^+\frac{1}{\rho-r}+4\ln^+ T(\rho,f)
\end{equation}
(b):\\
In this case we claim
\begin{equation}\label{manuilovekv73}
T(r,f')\leq 2 T(r,f)+m\left( r,\frac{f'}{f}\right)
\end{equation}
To see this note that
\begin{multline}
T(r,f')=N(r,f')+m(r,f')=2N(r,f)-\\
N_1(r,f)+m\left( r,f\cdot\frac{f'}{f}\right)\leq
N(r,f)+(N(r,f)+m(r,f))+m\left( r,\frac{f'}{f}\right)\leq\\
2T(r,f)++m\left( r,\frac{f'}{f}\right)
\end{multline}
where $N_1(r,f)$ denotes
\begin{equation}
\int_0^r \frac{n_1(t,a)-n_1(0,a)}{t}dt +n_1(0,a)\ln r
\end{equation}
and $n_1(t,a)$ denotes the set of $a$-points of $f(z)$ in $\{\abs{z}\leq t\}$, each $a$-point of order $\lambda$ is counted $\lambda-1$ times.
By Eqn.(\ref{manuilovekv72}) and Eqn.(\ref{manuilovekv73}) we have
\begin{equation}\label{manuilovekv74}
T(r,f')< A+4\ln\rho +\ln \rho +3\ln^+\frac{1}{\rho -r} ++2T(r,f)+ 4\ln^+ T(\rho,f)
\end{equation}
Thus
\begin{equation}\label{manuilovekv75}
T(r,f')< A+\ln \rho +3\ln^+\frac{1}{\rho -r} +6T(\rho,f)
\end{equation}
Let
$f(z)=\Phi_j^{(k-1)}(z).$
$\rho_0,\ldots,\rho_s,$
satisfying
$1<c=\rho_0<\rho_1<\cdots <\rho_s.$
These numbers will be specified later. By
Eqn.(\ref{manuilovekv71}) and Eqn.(\ref{manuilovekv75})
together with the properties of the positive logarithm and the inequalities
$\ln x< x$ for $x>0$ and $T(\rho_k)\leq T(\rho_s)$ for $\rho_k\leq \rho_s$, we obtain
\begin{equation}\label{manuilovekv76}
m\left( c,\frac{\Phi^{(k)}_j}{\Phi_j^{(k-1)}}\right)<
A+A\ln \rho_1 +A\ln^+\frac{1}{\rho_1 -\rho_0} +A\ln^+ T(\rho_1,\Phi_j^{(k-1)})
\end{equation}
for $l=1,\ldots , k-1,$
\begin{equation}\label{manuilovekv76a}
T(\rho_l,\Phi_j^{(k-l)})< 
A+A\ln \rho_{l+1} +A\ln^+\frac{1}{\rho_{l+1} -\rho_l} +A\ln^+ T(\rho_{l+1},\Phi_j^{(k-l)})
\end{equation}
This implies
\begin{equation}
T(\rho_1,\Phi_j^{(k-1)})< 
A+A\ln \rho_{s} +A \sum_{l=0}^{s-1} \ln^+ \frac{1}{\rho_{l+1}-\rho_l}
+A\ln^+ T(\rho_s,\Phi_j)
\end{equation}
\begin{equation}\label{manuilovekv77}
m\left( c,\frac{\Phi^{(k)}_j}{\Phi_j^{(k-1)}}\right)<
A+A\ln \rho_{s} +A \sum_{l=0}^{s-1} \ln^+ \frac{1}{\rho_{l+1}-\rho_l}
+A\ln^+ T(\rho_s,\Phi_j)
\end{equation}
By Eqn.(\ref{manuilovekv36}), Eqn.(\ref{manuilovekv41}), Eqn.(\ref{manuilovekv77})
\begin{multline}\label{manuilovekv78}
m\left( c,\frac{\Phi^{(k)}_j}{\Phi_j^{(k-1)}}\right)<
A+A\ln \rho_{s} +A \sum_{l=0}^{s-1} \ln^+ \frac{1}{\rho_{l+1}-\rho_l}+\\
A\ln c +A \sum_{l=0}^{s} \ln^+ T(\rho_s,\psi_l)
\end{multline}
By Eqn.(\ref{manuilovekv68}), Eqn.(\ref{manuilovekv69}), Eqn.(\ref{manuilovekv78})
we conclude that for sufficiently large $c$ we have
\begin{multline}\label{manuilovekv79}
\frac{1}{2\pi}\int_0^{2\pi} W(c\exp(i\theta))d\theta<
A+A\ln c +A\ln \rho_s + \\
A \sum_{l=0}^{s} \ln^+ T(\rho_s,\psi_l)
+A \sum_{l=0}^{s} \ln^+ \frac{1}{\rho_l -\rho_{l-1}}
\end{multline}
Here we may assume that
$T(c,\phi_l)\geq 1$, $l=1,\ldots,s.$
We now specify the numbers $\rho_k$ as follows
for $k=0,\ldots,s$
\begin{equation}\label{manuilovekv80}
\rho_k=c+\frac{k}{s\sum_{j=1}^s T(c,\psi_j)}
\end{equation}
For $0\leq k\leq s$ we then have
\begin{equation}\label{manuilovekv81}
\rho_s< c+\frac{1}{T(c,\psi_k)}
\end{equation}
By Eqn.(\ref{manuilovekv79}) and Eqn.(\ref{manuilovekv80})
we obtain
\begin{multline}\label{manuilovekv82}
\frac{1}{2\pi}\int_0^{2\pi} W(c\exp(i\theta))d\theta<
A+A\ln c +\\
A \sum_{k=0}^{s} \ln T(\rho_s,\psi_k)
+A \sum_{k=0}^{s} \ln T\left(c+\frac{1}{T(c,\psi_k)},\psi_k\right)
\end{multline}
By Lemma \ref{haymanIIlem} 
applied to $h(c)$ which is monotone increasing in $c$ we have
\begin{equation}\label{manuilovekv83}
h\left(c+\frac{1}{h(c)}\right)<2h(c)
\end{equation}
for all $c$ except possible a set of finite measure.
For $h(c)=T(c,\psi_k)$ we have by Eqn.(\ref{manuilovekv82})
that for all $c$ except possible a set of finite measure
\begin{equation}\label{manuilovekv84}
\frac{1}{2\pi}\int_0^{2\pi} W(c\exp(i\theta))d\theta<
A+A\ln c +
A \sum_{k=0}^{s} \ln T(c,\psi_k)
\end{equation}
By Eqn.(\ref{manuilovekv43}), Eqn.(\ref{manuilovekv50})-(\ref{manuilovekv53}), Eqn.(\ref{manuilovekv67})
and Eqn.(\ref{manuilovekv84}) we have for $c>0$ (except possible a set of finite measure)
\begin{equation}\label{manuilovekv85}
T(c,G)<A+A\ln c +A \sum_{k=0}^{s} \ln T(c,\psi_k)
\end{equation}
where
\begin{equation}
T(c,\psi_k)= T\left( c,\frac{\psi_k}{\psi_s}\cdot\psi_s\right)\leq
T\left( c,\frac{\psi_k}{\psi_s}\right)+T(c,\psi_s)
\end{equation}
By Eqn.(\ref{manuilovekv20}) and for polynomial $\psi_s$ we have
\begin{equation}
T(c,G)<A\ln T(c,B)+A\ln c
\end{equation}
which is the formula of Eqn.(\ref{manuilovekv33}).
As we have previously pointed out, by Eqn.(\ref{manuilovekv30}) and Eqn.(\ref{manuilovekv33}) 
this gives
\begin{equation}
T(c,B)<A\ln T(c,B)+A\ln c
\end{equation}
thus
\begin{equation}\label{manuilovekv86}
1<A\frac{\ln T(c,B)}{T(c,B)}+A\frac{\ln c}{T(c,B)}
\end{equation}
for all $c$ except possible a set of finite measure.
Suppose that at least one of the functions
$f_{ks}=\psi_k/\psi_s$, $k=1,\ldots,s-1,$
is transcendental.
It then follows (see Wittich \cite{vuttuch},  p.35) 
\begin{equation}
\lim_{c\to \infty} \frac{\ln c}{T(c,f_{ks})}=0
\end{equation}
so by Eqn.(\ref{manuilovekv20})
\begin{equation}\label{manuilovekv87}
\lim_{c\to \infty} \frac{\ln c}{T(c,B)}=0
\end{equation}
and since $T(c,B)$ is increasing to $\infty$ as $c\to \infty$, by Eqn.(\ref{manuilovekv20}),
we have as $c\to \infty$
\begin{equation}
\lim_{c\to \infty} \frac{\ln T(c,B)}{T(c,B)}=0
\end{equation}
Letting in Eqn.(\ref{manuilovekv86}) $c$ increase to infinity (while staying outside the exceptional intervals, which are of finite length),
we arrive at a contradiction. Hence all the functions $f_{ks}$ ar rational and since $\psi_s$ is a polynomial
all the $\psi_k$, $k=1,\ldots,s-1$ are polynomials.
By Eqn.(\ref{balkentirepaekv3}) and Eqn.(\ref{manuilovekv26}) this completes the proof.
	\end{proof}
The following is a polyentire version of Theorem \ref{balkdiscthm}.
\begin{theorem}[Balk \cite{balkentirepa}]\label{entirepathm}
	If $f(z)$ is a polyentire function with bounded zero set then there exists an entire function $g(z)$ and a polynomial 
	$P(z,\bar{z})$ such that $f(z)=\exp(g(z))P(z,\bar{z}).$
\end{theorem}
\begin{proof}
The reduced $n$-analytic function
\begin{equation}
F(z):=z^{n-1}f(z)
\end{equation}
satisfies the conditions of Theorem \ref{entirepathm00}
thus there exists polynomials $\psi_j$, $j=0,\ldots,n-1$ together with an entire $g(z)$
such that
\begin{equation}
F(z)=\left(\sum_{j=0}^{n-1} \psi_j(z)\bar{z}^j \right)\exp(g(z))
\end{equation}
This implies
\begin{equation}
\sum_{j=0}^{n-1} \psi_j(z)\bar{z}^j =\sum_{j=0}^{n-1} z^{n-1}\bar{z}^j f_j(z)\exp(-g(z))
\end{equation}
which by the uniqueness of analytic components gives for $j=0,\ldots,n-1$
\begin{equation}
\psi_j(z) = z^{n-1}f_j(z)\exp(-g(z))
\end{equation}
Hence there are polynomials $Q_j(Z)$ defined via
\begin{equation}
\psi_j(z)=z^{n-1}Q_j(z),\quad 
\end{equation}
such that
\begin{equation}
f(z)=\left(\sum_{j=0}^{n-1} Q_j(z)\bar{z}^j\right) \exp(g(z))
\end{equation}
This completes the proof.
\end{proof}

With regards to Theorem \ref{boschkraj70thm} we mention that Balk \cite{balk75factor} proved the following.
\begin{theorem}\label{balk75factorthm}
	Let $f(z)$ be an $n$-analytic function on a punctured domain $\Omega\subset \C$ of a point $z_0\in \C.$
	If $f$ has an isolated essential singularity at $z_0$ and $z_0$ is not a limit point of zeros of $f$, then $f$  
	may be written as $f(z)=p(z,\bar{z})\exp(G(z))$, where $G(z)$ is holomorphic near $z_0$ with an isolated essential singularity at $z_0$, while $p$ is polyanalytic of order $n$  
	near $z_0$ and does not have an essential singularity there.
\end{theorem}

Balk \& Goldberg \cite{balkgoldberg76}, prove, by a long argument based upon contradiction, using several 
results outside the scope of this book e.g.\ on algebroid functions, the Iversen property and complete covers,
the following result.
\begin{theorem}\label{balkgoldberg76factor1}
	Let $K=K(R):=\{R\abs{z}<\infty\}$ and let $f(z)=\sum_{j=0}^{q} a_j(z)\bar{z}^j,$ 
	for holomorphic $a_j(z),$ on $K$, $a_{q}\not\equiv 0.$
	Assume $f$ is irreducible in $K'=K'(R'):=\{R'<\abs{z}<\infty\}$ for each $R'$,
	$R\leq R'<\infty.$ If the set $E$ of all nonisolated zeros of $f$ in $K$ is unbounded then
	\begin{equation}
	f(z)=H(z)M(z)
	\end{equation}
	for holomorphic $H(z)$ on $K$ and a polyanalytic $M(z)$ such that $M(z)$ is meromorphic at $\infty.$
\end{theorem}
\begin{theorem}\label{balkgoldberg76factor}
	Let $R>0$ and let $f(z)$ be an $n+1$-analytic function on the punctured neighborhood 
	$U:=\{R<\abs{z}<\infty\}$ of $\infty.$ 
	Let $E$ be the set of nonisolated zeros of $f(z).$ If $E$ is unbounded then $f(z)$ can be expressed in some punctured neighborhood
	$U'$ of $\infty$ as $f(z)=I(z)M(z)$ where $I(z)$ is polyanalytic with only isolated zeros in $U'$, and $M(z)$ is a polyanalytic
	function that is meromorphic at $\infty.$
\end{theorem}
\begin{proof}
	Let $f(z)=\sum_{j=0}^n a_j(z)\bar{z}^j$ for $a_n(z)\not\equiv 0$ and $n\in \Z_+$ 
	or holomorphic $a_j(z),$ $j=0,\ldots,n$ on $\delta.$
	If $f(z)$ is reducible in some neighborhood $\delta_1K'(R_1)$ then it can be written as a product $F_1(z)F_2(z)$ for polyanalytic $F_1,F_2$ that are not analytic. 
	If there
	does not exist $R'>R_1$ such that $F_1,F_2$ are reducible in
	$\{R'<\abs{z}<\infty\}$ then 
	$f(z)$ can be written as the product
	\begin{equation}
	f(z)=\Pi_{k=1}^m f_k(z)
	\end{equation}
	for polyanalytic functions $f_k(z)$ which are irreducible at $\infty.$
	If at least one of $F_1,F_2$ is reducible in $K'(R_2)$ for some $R_2>R_1$
	we repeat the process 
	and after finitely many steps $f(z)$ will be expressed (see Mori \cite{mori}) as a product of polyanalytic functions that are not analytic in a domain of type $K'(R_2).$
		The process cannot be iterated indefinitely 
	since at each step the order of analyticity of the factors decreases and it cannot vanish.
	Hence
	we obtain that
	in a neighborhood $\delta_0=\{z:R_0<\abs{z}<\infty\}$, for some $R_0\geq R$ $f(z)$ can be written as the product
	\begin{equation}\label{balkgoldprodek}
	f(z)=\Pi_{k=1}^m f_k(z)
	\end{equation}
	for polyanalytic functions $f_k(z)$ which are irreducible at $\infty.$
	By assumption at least one of the $f_k(z)$ has unbounded set of nonisolated zeros.
	By Theorem \ref{balkgoldberg76factor1} there exists a holomorphic $H_k(z)$ and a polyanalytic $M_k(z)$
	that is meromorphic at $\infty$ such that
	\begin{equation}
	f_k(z)=H_k(z)M_k(z)
	\end{equation}
	Since the assertion of Theorem \ref{balkgoldberg76factor} holds true for
	functions with bounded set of nonisolated zeros the theorem holds true for each factor of the
	finite product given by Eqn.(\ref{balkgoldprodek}). This completes the proof of Theorem
	\ref{balkgoldberg76factor}.
\end{proof}
\begin{theorem}\label{balkgoldberg76factor2} 
	Let $f(z)$ be a polyanalytic function on the punctured neighborhood 
	$U$ of a point $c$ and also have essential singularity at $c$. Suppose there exists a polyanalytic function $a(z)$, meromorphic at $c$, such that $c$ is not a limit point of $\{f(z)-a(z)=0\}.$ Then for any $b(z)$, $b(z)\not\equiv 0$, that is polyanalytic on $U$ and meromorphic at $c$, there exists a neighborhood $U'$ of $c$ such that the set 
	$M(b,b'):=\{f(z)-b(z)=0\}\cap U'$ has $c$ as its unique limit point.
	\end{theorem}
\begin{proof}
	W.l.o.g.\ we assume $a(z)\equiv 0,$ $b(z)\not\equiv 0,$ $c=\infty$ assuming $f(z)\neq 0$ on $\delta.$ By 
	the Picard theorem for polyanalytic functions (see Bosch \& Krajkiewicz \cite{boschkraj70} and Theorem \ref{boschkraj70thm}) the set of roots to $f(z)=b(z)$ has $c$ as unique limit point.
	We show that there exists a neighborhood $\delta'$ of $c=\infty$ on which $M(b,\delta')$ has no limit points other than $\infty.$
	Assume (in order to reach a contradiction) that this doen not hold true.
	Then $\infty$ must be a limit point of the set of nonisolated zeros of $f(z)-b(z)$.
	By Theorem
	\ref{balkgoldberg76factor} there exists a neighborhood $\delta_0=K'(R_0)$ a polyanalytic $g(z)$ on $\delta_0$
	with only isolated zeros on $\delta_0$ and a polyanalytic $\pi(z)$ that is meromorphic at $\infty$ 
	such that on $\delta_0$
	\begin{equation}\label{balkgoldprodekl11}
	f(z)-b(z) =g(z)\pi(z)
	\end{equation}
	Since $\infty$ is a limit point of the set of nonisolated zeros of $f(z)-b(z)$ Eqn.(\ref{balkgoldprodekl11})
	$\infty$ is a limit point of the set of nonisolated zeros of $\pi(z).$ In particular, $\pi(z)$ is not holomorphic. However, since $c=\infty$ is not a limit point for the set of zeros of $f(z)$ we have 
	by the theorem on factorization in a neighborhood of isolated singularites (see Theorem \ref{balk75factorthm}) that
		there exists an entire $G(z)\neq 0$ and a polyanalytic $p(z)\neq 0$ on $\delta_0$ where $p(z)$ is meromorphic at $\infty$
	such that
	\begin{equation}\label{balkgoldprodekl12}
	f(z)=G(z)p(z)
	\end{equation}
	By Eqn.(\ref{balkgoldprodekl11}) and Eqn.(\ref{balkgoldprodekl12}) we have
	\begin{equation}\label{balkgoldprodekl13}
	G(z)p(z)-g(z)\pi(z)-b(z)\equiv 0
	\end{equation}
	We write
	$\pi(z)=\sum_{j=0}^n \pi_j(z)\bar{z}^j,$
	$p(z)=\sum_{j=0}^s p_j(z)\bar{z}^j,$
	$g(z)=\sum_{j=0}^k g_j(z)\bar{z}^j,$
	$b(z)=\sum_{j=0}^s b_j(z)\bar{z}^j$
	for holomorphic $\pi_j,p_j,g_j,b_j$ on $\delta_0$
	such that $\pi_j,p_j,b_j$ are meromorphic at $c=\infty$,
	$\pi_n(z)\not\equiv 0$, $n\geq 1,s=n+k.$ So Eqn.(\ref{balkgoldprodekl13}) becomes
	\begin{equation}\label{balkgoldprodekl14}
	G(z)\sum_{j=0}^s p_j(z)\bar{z}^j -\sum_{j=0}^k g_j(z)\bar{z}^j\sum_{j=0}^n \pi_j(z)\bar{z}^j-\sum_{j=0}^s b_j(z)\bar{z}^j\equiv 0
	\end{equation}
	Identifying the coefficients of powers of $\bar{z}$
	we obtain, setting $\pi_j=0$ for $j<0$
	\begin{equation}\label{balkgoldprodekl15}
	G(z) p_j(z) - \sum_{k=0}^{s-j} g_{j-n+k}(z) \pi_{n-k}(z)- b_j(z) = 0
	\end{equation}
	Setting $j=s$ we get
	\begin{equation}
	\pi_n g_k=GP_k-B_k
	\end{equation}
	where $P_k$ and $B_k$ are holomorphic in $\delta_0$ and meromorphic at $\infty.$
	Setting $\nu=s-1$ we get
	\begin{equation}
	\pi_n g_{k-1}=Gp_{s-1}-\pi_{n-1}g_k -b_{s-1}
	\end{equation}
	which implies
	\begin{equation}
	\pi_n^2 g_{k-1}=G\pi_n p_{s-1}-\pi_{n-1}(GP_k-B_k)-\pi_nb_{s-1}=GP_{k-1}-B_{k-1}
	\end{equation}
	where $P_{k-1}$ and $B_{k-1}$ are holomorphic in $\delta_0$ and meromorphic at $\infty.$ Iterating this
	we obtain
	\begin{equation}
	\pi_n^3 g_{k-2}=GP_{k-2}-B_{k-2},\ldots, \pi_n^{k+1} g_{0}=GP_{0}-B_{0}
	\end{equation}
	where $P_{k-2},B_{k-2},\ldots,P_0B_0$ are holomorphic in $\delta_0$ and a meromorphic at $\infty.$ This renders
	\begin{equation}\label{balkgoldprodekl16}
	\pi_n^{k+1}g(z)=GP(z)-B(z)
	\end{equation}
	where
	\begin{equation}
	P(z)=\sum_{j=0}^k P_j\pi_n^j\bar{z}^j,\quad B(z)=\sum_{j=0}^k B_j\pi_n^j\bar{z}^j
	\end{equation}
	By Eqn.(\ref{balkgoldprodekl13}) and Eqn.(\ref{balkgoldprodekl16}) we obtain
	\begin{equation}\label{balkgoldprodekl17}
	G(\pi_n^{k+1} p(z)-P(z)\pi(z)) (\pi_n^{k+1}b(z)-B(z)\pi(z))\equiv 0
	\end{equation}
	Since $G(z)$ has an essential singularity at $\infty$, Eqn.(\ref{balkgoldprodekl17}) holds only if
	$\pi^{k+1}_n(z)\equiv P(z)\pi(z).$ Since $\pi(z)$ has unbounded set of nonisolated zeros and $p(z)\neq 0$ in 
	$\delta_0$ this is a contradiction. This
	verifies that there exists a neighborhood $\delta'$ of $c=\infty$ on which $M(b,\delta')$ has no limit points other than $\infty.$
	This completes the proof of Theorem \ref{balkgoldberg76factor2}.
\end{proof}
For a compact $K\subset\R^2$, we shall in this section denote by $\mbox{Cap}_2(K)$
the logarithmic capacity of $K$ (See Ronkin \cite{ronkinbok}, p.51, see also Chapter \ref{reproducingsec}). 
If $A\subset \R^2$ is an arbitrary set then we define the {\em inner capacity}\index{Inner capacity}
$\underline{\mbox{Cap}}(A)$ as the least upper bound of the capacities of all compacts subsets $K$ of $A$ i.e.\
\begin{equation}
\underline{\mbox{Cap}}(A):=\sup_{K\subset A} \mbox{Cap}(K)
\end{equation}
we define the {\em outer capacity}\index{Outer capacity}
$\underline{\mbox{Cap}}(A)$ as the greatest lower bound of the inner capacities of all open supersets of $A$ i.e.\
\begin{equation}
\overline{\mbox{Cap}}(A):=\inf_{A\subset G} \underline{\mbox{Cap}}(G)
\end{equation}
\begin{definition}[See Ronkin \cite{ronkinbok}, p.86]
	For each $z'\in C^{n-1}$ we define $M_{z'}^{n-k}:=\{z\in \Cn: z_j=z_j',j=1,\ldots,n-k\},$ $k=1,\ldots,n-1.$
	For $E\subset\Cn$, $n>1$, we define the {\em $\Gamma$-projection}\index{$\Gamma$-projection}, $\Gamma_n^{n-1}(E)$, of $E$
	on $\C^{n-k}$ to be the set of all points $(z_1,\ldots,z_{n-1})$ such that 
	\begin{equation}
	\underline{\mbox{Cap}}_2 (E\cap M_{z}^{n-1})>0
	\end{equation}
	(where $\underline{\mbox{Cap}}_2$ denotes the inner capacity)
	and 
	\begin{equation}\Gamma_n^{l}(E)=\Gamma_{l+1}^{l}(\Gamma_n^{l+1}(E))
	\end{equation}
	in particular 
	$\Gamma_n^{1}(E)=\Gamma_{2}^{1}(\Gamma_n^{2}(E))$. For $E\subset \C$ we define $\Gamma_1^1(E)=E.$
	When $E$ is a Borel set we may replace $\underline{\mbox{Cap}}_2$ by $\mbox{Cap}_2$.
	It is known that if $E$ is a Borel set such that $\underline{\mbox{Cap}}_2 \Gamma_n^1(E)=0$
	then the $2n$-dimensional Lebesgue measure $\lambda_{2n}(E)=0$ (see Ronkin \cite{ronkinbok}, p.90).
	The {\em $\Gamma$-capacity}\index{$\Gamma$-capacity, $\Gamma$-Cap$(E)$} of $E\subset\Cn$, $\Gamma-\mbox{Cap}(E)$
	is defined as
	\begin{equation}
	\Gamma-\mbox{Cap}(E):=\sup_{\alpha\in T} \overline{\mbox{Cap}}_2 \Gamma_n^1 (\alpha E)
		\end{equation}
		where $\underline{\mbox{Cap}}_2$ denotes the outer capacity, and $T$ is the class of unitary transformations on $\Cn$,
		that is $z\mapsto \zeta,$ $\zeta_j=\langle z,a^{(j)}\rangle$, 
		$j=1,\ldots,n$ where the vectors $a^{(j)}\in \Cn$ satisfy $\abs{z}=\abs{\zeta}$ for all $\zeta\in \Cn.$ 
	The $\Gamma$-capacity is invariant under 
	unitary transformations of $\Cn.$
\end{definition}
The following is known (see Ronkin \cite{ronkinbok}, p.92)
\begin{proposition}
	Let $E\subset \Cn,$ $E\neq \Cn$ such that the intersection of $E$ with an analytic plane $\{z_i=a_iw+b_i,i=1,\ldots,n,\abs{w}<\infty\}$ which is
	not a subset of $E$ has inner capacity zero. Then the $\Gamma$-capacity of $E$ is zero. 
	The $\Gamma$-capacity of any analytic set (the common zeros of analytic functions in $\Cn$)
	is zero.  
\end{proposition}
Ronkin \cite{ronkinarticle}, in a study of sets $\{w\in \C^m:f(z,w)=0\}$ for entire $f(z,w),$ where $(z,w)\in\C^{n+m}$, and conditions for when
the latter set coincides with a set of the form $\{w\in \C^m:p(w)=0\}$ for a polynomial $p(w)$,
proves the following.
\begin{theorem}\label{ronkinthm}
	Let $(z,w)\in\C^{n+m}$ denote the complex coordinates, let $f(z,w)$ be an entire function
	and let $E\subset \C^n$ be a set of positive $\Gamma$-capacity. 
	If for each $z\in E$ the set $\{w:f(z,w)=0\}$ coincides with the zero set of a polynomial in $w$
	then there exists an entire function $h(z,w)$ and function $P(z,w)$ that is a polynomial with respect to $w$ with coefficients that
	are entire with respect to $z$, such that
	$f(z,w)$ has the representation 
	\begin{equation}
	f(z,w)=P(z,w)\exp(h(z,w))
	\end{equation}
\end{theorem}
Balk \cite{ca1} refers to a local publication (which we have not been able to obtain) by Manuilov \cite{manuilov} with regards to the following result.
\begin{theorem}\label{manuilovthm}
	Let $(z,w)\in\C^{n+m}$ denote the complex coordinates. Let $P_1(z,w)$ be a pseudopolynomial in $z$ with entire coefficients with respect to $w$
	and Let $P_2(z,w)$ be a pseudopolynomial in $w$ with entire coefficients with respect to $z$.
	If there exists an entire function $G(z,w)$ such that $P_1(z,w)\equiv Z_2(z,w)\exp(G(z,w))$
	then there exists a polynomial $P(z,w)$ together with entire functions $g_1(z),$ $g_2(z)$ such that
	\begin{equation}
	P_1\equiv P\exp (g_1),\quad P_2\equiv \exp (g_2) 
	\end{equation}
\end{theorem}
Balk \& Manuilov \cite{balkmanuilov} prove a generalization of Theorem \ref{entirepathm}.
\begin{theorem}
	If $f(z_1,z_2))$ is a polyentire function on $\C^2$, with bounded 
	zero set then there exists an entire function $g(z_1,z_2)$ and a polynomial 
	$P(z_1,z_2,\bar{z}_1,\bar{z}_2)$ such that $f(z)=\exp(g(z))P(z,\bar{z}).$
	\end{theorem}
	\begin{proof}
	Let $(n_1,n_2)\in Z_+^2$ and set
	\begin{equation}
	F(z_1,z_2,z_3,z_4):=\sum_{k_1=0}^{n_1-1}\sum_{k_2=0}^{n_2-1} a_{k_1,k_2}(z_1,z_2)z_3^{k_1}z_4^{k_2}
	\end{equation}
	For each choice of $z_2=z_2^0$, $z_4=\bar{z}_2^0$ the function $F$ is a pseudopolynomial with respect to $z_3$
	with coefficients that are entire functions of $z_3.$ This restricted function has bounded zero set in its domain
	thus by Theorem \ref{entirepathm}
	there exists an entire function $g(z_1,z_2^0,\bar{z}_2^0)$, with respect to $z_1$, and a polynomial $\pi(z_1,z_3,z_2^0,\bar{z}_2^0)$ with respect to 
	$z_1,z_3$ such that
	\begin{equation}\label{balkmanuilov9}
	F(z_1,z_2^0,z_3,\bar{z}_2^0)=\pi(z_1,z_3,z_2^0,\bar{z}_2^0)\exp(g(z_1,z_2^0,\bar{z}_2^0))
	\end{equation}
	Denote by $E$ the subset in the two-dimensional complex space with respect to the variables $z_2,z_4$ defined by 
	$z_4=\bar{z}_2$. By Eqn.(\ref{balkmanuilov9}) we have for each choice of $z_2,z_4$ that
	$E$ has positive $\Gamma$-capacity thus by Theorem \ref{ronkinthm} there exist an entire function
	$\varphi(z_1,z_2,z_3,z_4)$ and a function $P_{13}(z_1,z_2,z_3,z_4)$ that is a polynomial 
	with respect to $(z_1,z_3)$ with coefficients that are entire with respect to $z_2,z_4$
	such that
	\begin{equation}\label{balkmanuilov10}
	F(z_1,z_2^0,z_3,\bar{z}_2^0)\equiv P_{13}(z_1,z_2,z_3,z_4)\exp(\varphi(z_1,z_2,z_3,z_4))
	\end{equation}
	Analogously we obtain an entire function
	$\psi(z_1,z_2,z_3,z_4)$ and a function\\ 
	$P_{24}(z_1,z_2,z_3,z_4)$ that is a polynomial 
	with respect to $(z_2,z_4)$ with coefficients that are entire with respect to $z_1,z_3$
	such that
	\begin{equation}\label{balkmanuilov11}
	F(z_1,z_2,z_3,z4)\equiv P_{24}(z_1,z_2,z_3,z_4)\exp(\psi(z_1,z_2,z_3,z_4))
	\end{equation}
	By Eqn.(\ref{balkmanuilov10}) and Eqn.(\ref{balkmanuilov11}) we obtain
	\begin{equation}
	P_{24}=P_{13}\exp \eta,\quad \eta=\varphi-\psi
	\end{equation}
	By Theorem \ref{manuilovthm} there exists a polynomial $P(z_1,z_2,z_3,z_4)$ together with entire functions
	$\varphi_1(z_2,z_4)$, $\psi_1(z_1,z_3)$ such that
	\begin{equation}
	P_{13}\equiv P\exp(\varphi_1),\quad P_{24}\equiv P\exp \psi_1
	\end{equation}
	This implies that $F\equiv P\exp \Phi$ where $\Phi=\varphi +\psi.$
	Comparing the growth of the functions $F/P$ and $\exp(\Phi)$ shows that
	$\Phi$ is independent of $z_3,z_4.$ Hence the restricted function
	$f(z_1,z_2):=F(z_1,z_2,\bar{z}_1,\bar{z}_2)$ has a representation of the wanted form. This completes the proof.
	\end{proof}

	\section{Lelong's indicator of growth}
	As is pointed out in Balk \cite{ca1}, Theorem 2.6, p.213, it is possible, in the case of
	a reduced polyentire function, $f$, to obtain factorization conditions based upon the growth of
	$N(r,0,E_r)$, where $E_r$ is a holomorphic function coinciding with $f$ on $\{\abs{z}=r\},$ 
	with namely $N(r,0,E_r)=O(\log r),$ $r\to \infty$ will be a sufficient condition.
	We shall here look at a slightly different approach which predates these observations, and is mainly due
	to the works of Lelong.
	\begin{definition}
		By a {\em current}\index{Positive current} of bidegree $(p,q)$ we simply mean the dual of the space of
		$C^\infty$-smooth differential forms of bidegree $(p,q)$ forms on $\Cn.$  
		A current $\theta$ of bidegree $(n-p,n-p)$ is called {\em positive} if for all  
		collections of $p$, $C^\infty$-smooth, differential 
		forms $\alpha_1,\ldots,\alpha_p$ of bidegree $(1,0)$, with compact support, the distribution
		\begin{equation}
		\theta\wedge (i\alpha_1\wedge)\wedge \cdots\wedge (i\alpha_p\wedge)
		\end{equation} 
		is positive.
		A differential form is called positive if its natural associated current is positive.
		If $X$ is an analytic subvariety of $\Cn$ of dimension $p$ and $\phi$ a differential $(p,p)$-form on $\Cn$ then we associate to $\theta$
		the current 
		\begin{equation}\label{lelskodbasicekv}
		\langle \theta,\phi\rangle =\int_X \phi
		\end{equation}
		For a differential form $\beta$ and a positive integer $n$ we denote
		\begin{equation}
		\beta^n:= \overbrace{\beta\wedge\cdots\wedge \beta}^{n-\mbox{times}}
		\end{equation}	
	\end{definition}
	Lelong proved (see e.g.\ Skoda \cite{skoda}, Proposition 1.1) that the current defined by Eqn.(\ref{lelskodbasicekv})
	is closed and positive.
	\begin{example}
		Consider the current $\theta$ of bidegree $(1,1)$
		\begin{equation}
		\theta=i\sum_{j,k} \theta_{jk} dz_j\wedge d\bar{z}_k
		\end{equation}
		where $\theta_{jk}$ are $0$-currents or distributions on $\Cn$, $1\leq j,k\leq n.$ Then the condition for positivity
		is
		\begin{equation}
		\sum_{j,k} \theta_{jk}\lambda_j\bar{\lambda}_k \geq 0,\mbox{ for all }\lambda_j\in \C
		\end{equation}
		If$\theta$ is a positive current, $\omega$ is a differential form of degree $(1,1)$ and $\rho$ is a $C^\infty$-smooth function with compact support 
		then the regularization $\rho*\theta$ and the current $\theta\wedge \omega$ are both positive. 
	\end{example}
	\begin{definition}
		If $\theta$ is a positive current on $\Cn$ and 
		\begin{equation}
		\beta:=i\sum_{j=1}^n dz_j\wedge d\bar{z}_j
		\end{equation}
		we associate to $\theta$ the so called {\em trace measure}
		\begin{equation}
		\sigma:=\frac{2^{-n}}{n!} \theta \wedge \beta^n
		\end{equation}
		and the {\em projecive measure}\index{Projective measure}
		\begin{equation}
		\nu:=\frac{1}{(2\pi)^n} \theta \wedge \alpha^n,\quad \alpha=i\Delta \log\abs{z}^2
		\end{equation}
		where $\Delta=\partial\overline{\partial}$.
		We define the measure supported in the ball with radius $r$ according to
		\begin{equation}
		\sigma(r):=\int_{\abs{x}<r} d\sigma(x)
		\end{equation}
		\end{definition}
		Lelong \cite{lelong68} proved the following. 
			\begin{proposition}
		If $\theta$ is a closed positive current then the function
		$r^{-2n}\sigma(r)$, $r>0$ is increasing in $r$.
		Let 
		\begin{equation}
		\nu(0)=\lim_{r\to 0} \pi^{-n} n! r^{-2n} \sigma(r)
		\end{equation}
		Then the $\int_{0\leq \abs{x}<r} d\nu(x)$ is a convergent integral and
		\begin{equation}
		\nu(0)+\int_{0<\abs{x}<r} d\nu(x)=\pi^{-n} n! r^{-2n}\sigma(r)
		\end{equation}
		\end{proposition}
		The number $\nu(z)$ defined by
		\begin{equation}
		\nu(z)=\lim_{r\to 0} \pi^{-n} n! r^{-2n} \int_{\abs{x-z}<r}d\sigma(x)
		\end{equation}
		is called the {\em Lelong number}\index{Lelong number} of the current $\theta$ at $z\in \Cn$.
		We have
		\begin{equation}
		\nu(z)+(2\pi)^{-n}r^{-2n} \int_{\abs{x-z}<r} \beta^n \wedge \theta
		\end{equation}
		and it can be showed that $\nu(z)$ is an integer $\geq 1$ (see e.g.\ Skoda \cite{skoda}, p.362).
		\begin{definition}
		If $\theta$ is a positive closed current on $\Cn$ we say that 
		$\theta$ has {\em order} $\rho$ if\index{Order of a current}
		\begin{equation}
		\limsup_{r\to \infty} \frac{\log \nu(r)}{\log r} =\rho
		\end{equation}
		\end{definition}
		\begin{definition}
		A real valued function $\phi,$ $-\infty \leq \phi(x)<\infty$, is called {\em subharmonic}\index{Subharmonic function}
		on a domain $U\subset\R^m$ if $\phi$ is upper semi-continuous and $\not\equiv -\infty$ in $U$
		such that \begin{equation}
		\phi(x)\leq \lambda(x,r,\phi)\equiv \frac{1}{\omega_m^{-1}}\int \phi(x+r\alpha)d\omega_m(\alpha)
		\end{equation}
		for $r<\mbox{dist}(x,U)$,
		where $\omega_m$ is the Lebesgue measure on the unit sphere $S^{m-1}$ in $\R^m$ and $\lambda$ is the 
		mean value of $\phi$ on $S^{m-1}$ relative to the Haar measure $d\omega_m(\alpha).$
		For a domain $\Omega\subset\Cn$ A real valued function $\phi,$ $-\infty \leq \phi(x)<\infty$, is called {\em plurisubharmonic}
		on $\Omega$ if $\phi$ is upper semi-continuous and $\not\equiv -\infty$ in $\Omega$
		such that 
		\begin{equation}
		\phi(x)\leq l(z,w,r,\phi)= \frac{1}{2\pi}\int_0^{2\pi} \phi(z+wr\exp(i\theta))d\theta
		\end{equation}
		for all $w,r$ such that $z+uw\in \Omega,$ $u\in \C,$ $\abs{u}\leq r.$
		If $p(z)$ is a subadditive, real-valued function on $\Cn$ and $a(z):\Cn\to \R_+$ we define
		\begin{equation}
		M_{a,p}(r)=\sup_{p(z)\leq r} a(z)
		\end{equation}
	\end{definition}
	Let $\phi(t)$ be an increasing function of $t$ for $t\geq 0$ such that $\limsup_{t\to\infty} \frac{\phi(tu)}{\phi(t)}<\infty$,
	for $u\geq 0$
	and $\lim_{t\to \infty}\phi(t)=\infty$.
	If
	$f(z)$ is an entire (holomorphic) function on $\Cn$ such that for a constant $C=C(f)>0$ depending on $f$
	\begin{equation}\label{lelonggrumanekv1}
	\log \abs{f(z)}\leq \phi(\abs{z}) +C(f)
	\end{equation}
	then
	\begin{equation}
	\limsup_{t\to\infty}\frac{\log \abs{f(tz)}}{\phi(t)}
		\end{equation}
		is a measure of the asymptotic growth with respect to the {\em weight} factor\index{Weight factor for asymptotic growth}
		on the real lines through the origin.
		Similarly, one can consider the expressions of the form $\limsup_{i\in I} C_i \log \abs{f_i}$ for 
		entire $f_i$ satisfying Eqn.(\ref{lelonggrumanekv1}) with $C,f$ replaced by $C_i,f_i$. 
		It is well-known that if $\phi$ is plurisubharmonic then $\{z\in \Omega :\phi(z)=-\infty\}$ has Lebesgue measure zero and
		the functions $r\mapsto l(z,w,r,\phi)$ and $r\mapsto \sup_{z'\in\{z+uw,u\in \C,\abs{u}\leq r\}} \phi(z)$ are increasing convex functions of $\log r$ and hence
		$\lambda(z,r,\phi)$ is also (see e.g.\ Lelong \& Gruman \cite{lelonggruman}).
		Consider the functions  and $M(z,r,\phi)=\sup_{z'-z}\leq r\phi(z')$
		and $m_{\phi}(z,z',r):=\sup_{\abs{u}\leq r} \phi(z+uz').$
		Recall that if $p$ and $q$ are two norms on $\Cn$ there exists positive finite constants $C_1$ and $C_2$ such that
		\begin{equation}
		0<C_1\leq \frac{p(z)}{q(z)}\leq C_2
		\end{equation}
		So if $a(z):\Cn\to \R_+$ and $M_{a,p}(r)=\sup_{p(z)\leq r}a(z)$ then there exists $0<C_1<C_2<\infty$ depending only on $p(z),q(z)$ 
		such that for each $a(z)$
		\begin{equation}\label{ekplsh0022}
		M_{a,p}(C_1r)\leq M_{a,p}(r)\leq M_{a,p}(C_2r)
		\end{equation}
		Any subharmonic function on a domain $U\subset \R^2$ is locally integrable and the set of plurisubharmonic functions on a domain $\Omega\subset
		\Cn$
		belong to the set of subharmonic functions on the corresponding domain $\Omega\subset \R^{2n}.$
		\begin{proposition}\label{grumanbokpro14}
			If $\phi(z)$ is plurisubharmonic on $\Cn$ then
			$m_{\phi}(z,z',r):=\sup_{\abs{u}\leq r} \phi(z+uz')$ is $\equiv -\infty$ or an increasing convex function of $\log r$
			and if $p(z)$ is a complex norm then $M_{\phi,p}(r)$ is an increasing convex function of $\log r.$
		\end{proposition}
		\begin{proof}
			\begin{equation}
			\end{equation}
			The first assertion follows from the fact that the function $\phi(z+uz')$ is $-\infty$ for all $u\in \C$ or
			a subharmonic function of the variable $u:=\alpha +i\beta$ in $C\simeq \R^2.$ 
			For the second assertion note that $M_{\alpha,p}(r)=\sup_{z\in p^{-1}(1)}\left(\sup_{\abs{u}\leq r} \phi(uz)\right)$
			and $\sup_{\abs{u}\leq r} \phi(uz)$ is an increasing convex function of $\log r$ or $\equiv -\infty,$
			but is not $\equiv  -\infty$ for all $z.$ This completes the proof.
		\end{proof}
		\begin{theorem}\label{nulimkvotthm}
			(i) Let $p(z)$ be a norm and $\phi(z)$ a plurisubharmonic function on $\Cn.$ Then $C=\lim_{r\to \infty} \frac{M_{\phi,p}(r)}{\log r}$
				and $C(z,z')=\lim_{r\to \infty} \frac{m_{\phi,p}(z,z',r)}{\log r}$
					exist, either finite or infinite such that:
					\\
					(a) $C\geq 0$ and $C(z,z')\geq 0$ except possibly $C(z,z')=-\infty$ in which case $\phi(z+uz')\equiv -\infty$ for $u\in \C.$\\
					(b) $C(z,uz')\equiv C(z,z')$ for each $u\in \C,$ $u\neq 0.$
					\\
					\\
					(ii) If $p(z)$ is a norm on $\Cn$ and $\phi(z+uz')\not\equiv -\infty$, $u\in \C$ then 
					$\frac{\partial}{\partial \log r} M_{\phi,p}(r)$ and $\frac{\partial}{\partial \log r} m_{\phi,p}(z,z',r)$ exist except possibly
					for a countable set of $r$ values and
					\begin{equation}
					\lim_{r\to \infty}\frac{\partial}{\partial \log r} M_{\phi,p}(r)=\lim_{r\to \infty} \frac{M_{\phi,p}(r)}{\log r} 
					\end{equation}
					\begin{equation}
					\lim_{r\to \infty}\frac{\partial}{\partial \log r} m_{\phi,p}(r)=\lim_{r\to \infty} \frac{m_{\phi,p}(r)}{\log r} 
					\end{equation}
					\end{theorem}
					\begin{proof}
					There exists $r_0>0$ such that $M_{\phi,p}(r_0)>-\infty$ and since
					it is an increasing convex function of $\log r$, $M_{\phi,p}(r)>-\infty$ for $r\geq r_0.$ This proves that $C\geq 0$. If $\phi(z+uz')\equiv -\infty$
					then for $r>1$ we have
					$(\log r)^{-1} m_\phi(z,z',r)=-\infty$ hence $C(z,z')=-\infty.$ Otherwise $\phi(z+uz')$ is subharmonic with respect to $u$ 
					on $\R^2$ thus $C(z,z')\geq 0.$ By direct calculation we have $C(z,z')=C(z,uz')$ for $u\neq 0.$
					The assertion (ii) is direct consequence of Proposition \ref{grumanbokpro14}.
				\end{proof}
		\begin{definition}
			If $a(z):\Cn\to \R_+$ and $p(z)$ is a norm the {\em order}\index{Order of a positive real-valued function with respect to a given norm}, $\rho$, 
			of $a(z)$ (with respect to $p$) is defined as
			\begin{equation}
			\rho:=\limsup_{r\to \infty} \frac{\log M_{a,p}(r)}{\log r}
			\end{equation}
			If $\rho<\infty$, then $a(z)$ is said to be of maximal, normal, or minimal type respectively if the positive number
			\begin{equation}
			\sigma=\limsup_{r\to \infty}\frac{M_{a,p}(r)}{r^\rho}
			\end{equation}
			is $+\infty$, finite, zero respectively. $\sigma$ is said to be the {\em type} of $a(z)$ with respect to $p(z)$.
			For a plurisubharmonic $\phi$ we define the order $\rho$ of $\phi$ by using
			$a(z):=\phi^+(z)=\sup\{\phi(z),0\}.$ If $f$ is an entire function of order $\rho$ we say that $f$ has order $\rho$
			\index{Order of an entire function} if the plurisubharmonic function $\log \abs{f}$ has order $\rho.$ 
				\end{definition}
		The order and type are in some sense invariant with respect to the choice of norm $p$ and with respect to translation, thus depend only on the topology of $\Cn.$
		\begin{remark}
			Let $\phi$, be a plurisubharmonic function on $\Cn$, and $\lambda(\phi,z_0,r)$, the
			spherical mean of $\phi$ with respect to the sphere with center $z_0$ and radius $r.$ As we have seen
			it is known that this is increasing with respect to $\log r$ (see Skoda  \cite{skoda}, Proposition I.17) and after a 
			verifying existence of partial derivative except possibly a countable set off $r$ values, with respect to $\log r$, it is possible to prove (see Skoda \cite{skoda}, p.45)
			that if $t=i\partial \overline{\partial} \phi$ then $\sigma_t =\frac{1}{2\pi}\Delta V$ and
			\begin{equation}
			\nu_t(z_0)=\lim_{r\to 0} \frac{\partial}{\partial \log r} \lambda(z_0,r,\phi)=\lim_{r\to 0}\frac{\lambda(z_0,r,\phi)}{\log r}
			\end{equation}
				\end{remark}
		Let $\alpha\in \Z_+^n,$ and 
	let $f(z)=\sum_{j=0}^n\sum_{\gamma_j=0}^{\alpha_j-1} a_\gamma(z)\bar{z}^\gamma,$
	for entire $a_\beta(z),$ be an $\alpha$-analytic function on $\Cn.$
	Set (see Skoda \cite{skoda}, p.46)
	\begin{equation}
	\phi_f:=\frac{1}{2}\log \left(\sum_{j=0}^n\sum_{\gamma_j=0}^{\alpha_j-1} \abs{a_\gamma}^2\right)
	\end{equation}
	Then $\phi_f$ has the associated positive current (see Lelong \cite{lelong64})
	\begin{equation}
	\theta:=\frac{i}{\pi}\sum_{i,j} \partial_{z_i}\partial_{\bar{z}_j} \phi_f dz_i\wedge d\bar{z}_j
	\end{equation}
		We set
	\begin{equation}
	\beta:=\frac{i}{2}\sum_{J}  dz_i\wedge d\bar{z}_j,\quad \beta_{n-1}:=\frac{\beta^{n-1}}{(n-1)!}
	\end{equation}
	The trace of this current\index{Trace of a current}
	is the positive measure
	\begin{equation}
	\sigma=\theta\wedge \beta_{n-1} =\frac{1}{2\pi} \Delta \phi_f \beta_n
	\end{equation}
	Denote by $\sigma(r)$ the measure supported in the ball centered at the origin with radius $r.$
	We have the increasing function with respect to $r$
	\begin{equation}
	\nu(r)=\frac{(n-1)!}{\pi^{n-1}}\frac{\sigma(r)}{r^{2n-2}}
	\end{equation}
\begin{definition}
	Let $\alpha\in \Z_+^n$ and let $f(z)=\sum_{j=1}^n\sum_{\gamma_j<\alpha_j}^q a_{\gamma}(z)\bar{z}^\gamma$ be a polyanalytic function, i.e.\ the $a_\gamma(z)$ are entire functions.
	Let us denote for $\gamma\in \N^n,$ by $\gamma< \alpha$ that $\gamma_j< \alpha_j$, $j=1,\ldots,n.$
	Let $E:=\{z:a_\gamma(z)=0,\gamma< \alpha\}$. Let $z_0\in E$ and denote by $n_0$ the lowest sum-norm of the multi-integers corresponding to the multi-degrees
	of the $a_\gamma(z)$ about $z_0$, $\gamma< \alpha.$
	If $E$ is a proper nonempty subset then $n_0$ is positive.
	When $E$ is algebraic, i.e.\ when each $a_\gamma(z)$ is a polynomial then
	we denote by $N$ the highest sum-norm of the multi-integers corresponding to the multi-degrees of the $a_\gamma(z)$,$\gamma< \alpha.$
	Let $\lambda(\phi_f,z_0,r)$ denote the spherical mean of the plurisubharmonic function $\phi_f$ and set, for $z_0$ not necessarily the origin
	$\nu(z_0,r)$ analogous to above but where we replace the origin by $z_0$
	(it is possible  
	to identify this as $\nu(z_0,r)=\frac{\partial\lambda}{\partial \log r}$ at $z_0$, we know this is the case off a countable set of points, and for our purposes it suffices to define the value as the left limits of the difference quotients)
	This is a positive function that has a limit as $r\to \infty$ as noted above. 
\end{definition}
\begin{proposition}\label{traoreprop31}
	Let $\alpha\in \Z_+^n$ and let $f(z)=\sum_{\gamma\in N^n,\gamma<\alpha} a_{\gamma}(z)\bar{z}^\gamma$ be a polyanalytic polynomial,
	i.e.\ the $a_j(z)$ are polynomials. Let $E:=\{z:a_j(z)=0,\gamma<\alpha\}$ and let $z_0\in E$.
	Then 
	\begin{equation}
	n_0\leq \nu(z_0,r)\leq N,\quad \lim_{r\to 0}\nu(z_0,r)=n_0,\quad \nu(z_0,r)=N
	\end{equation}
\end{proposition}
\begin{proof}
	W.l.o.g.\ suppose $z_0=0.$ Let $\xi\in \Cn$, 
	$\abs{\xi}=1,$ $z:= s\cdot \xi,$ $s\in \C$, $\abs{s}=r.$
	Then for complex coefficients $a_{\gamma,\beta}\in \C$ we have
	\begin{equation}
	a_{\gamma}(z)=\sum_{\beta_\gamma\leq \beta\leq \mu_\gamma} a_{\gamma,\beta}s^{\abs{\beta}} \xi^\beta,
	\end{equation}
	$\beta\leq \mu_\gamma$
	where $\mu_\gamma$ is the degree of $a_\gamma,$ 
	and where in the sum $\beta_\gamma \leq \beta \leq \mu_\gamma$ for some $\beta_\gamma\in \N^n,$ $\abs{\beta_\gamma}>0.$
	Then since $\gamma-\beta_\gamma\geq 0,$
	\begin{equation}
	a_\gamma(z)=s^{\abs{\beta_\gamma}} \sum_{\beta_\gamma\leq \beta\leq \mu_\gamma} a_{\gamma_j,\beta_j} s^{\abs{\beta-\sigma_\gamma}}\xi^\beta
	=:s^{\beta_\gamma} A_\gamma(\xi,s)
	\end{equation} 
	where $A_\gamma(\xi,0)\neq 0.$
	Now $n_0=\inf_\gamma \abs{\beta_\gamma}$ which gives
	\begin{equation}
	\phi_f(z)=n_0\log r +\frac{1}{2}\log \sum \abs{A_\gamma^2(\xi,s)}=n_0\log r + G(\xi,s)
	\end{equation}
	As we have previously pointed out (see e.g.\ Lelong \cite{lelong64}) it is known that $\lambda(\phi_f,r)$ is convex with respect to $\log r$ and 
	$\lim_{r\to 0} \nu(r)$ is the limit of $\frac{\lambda(\phi_f,r)}{\log r}$ hence
	\begin{equation}
	\lambda(\phi_f,r)=n_0\log r +\lambda(G,1)
	\end{equation}
	However, $G(\xi,s)$ is continuous and bounded for $s=0$. Thus
	$\frac{\lambda(\phi_f,r)}{\log r}$ goes to $n_0$, hence $\nu(0)=n_0.$
	Now given $N=\sup_\gamma \abs{\mu_\gamma}$ we have
	\begin{equation}
	a_\gamma(s,\xi)=s^N\sum_\beta a_{\gamma,\beta}\frac{\xi^\beta}{s^{N-\abs{\beta}}}=s^N A_\gamma(\xi,s)
		\end{equation}
		which yields
		\begin{equation}
		\phi_f(z)=\frac{1}{2}\log \sum \abs{a_j}^2 =N\log r +\frac{1}{2}\log
		\sum \abs{A_j(\alpha,s)}^2 =N\log r+G(\alpha,s)
		\end{equation}
		For each parameter value $s$ such that $\abs{s}=r>0$,
		$G(\alpha,s)$ is a plurisubharmonic function of $\alpha.$
		This defines a family $F$, for $r\geq r_0>0$, of plurisubharmonic functions
		such that $F$ is bounded on compacts and there does not exists a sequence of functions
		in $F$ that converge uniformly on compacts to $-\infty.$	
			It follows from a known result of Lelong \cite{lelong50} that there exists a constant independent od $s,$ $\abs{s}<r_0$, such that $\lambda(G,1)\geq a.$
		This implies that
		\begin{equation}
		\lambda(\phi_f,r)=N\log r +\lambda(G,1),\quad\abs{s}=r
		\end{equation}
		\begin{equation}
		\frac{\lambda(\phi_f,r)}{\log r}=N +\frac{\lambda(G,1)}{\log r}
		\end{equation}
		where the right hand side goes to $N$ as $r\to \infty,$ since $\lambda(G,1)$ is bounded. This completes the proof.
	\end{proof}
\begin{proposition}\label{traoreprop32}
	Let $f(z)=\sum_{\beta<\alpha} a_\beta(z)\bar{z}^\beta$ be an $\alpha$-analytic function near a point $z_0\in \Cn$.
	If $z_0$ is not a common zero the analytic components $a_\beta(z)$ then
	$\frac{\sigma(z_0r)}{r^{2n}}$ is bounded as $r\to \infty.$ 
\end{proposition}
\begin{proof}
	We have that$ \sigma$ is, in terms of distributions, up to a constant,
	$\frac{1}{2\pi} \Delta\phi.$ Let $a=\{a_\beta\}_\beta$ be the array of analytic components and $\abs{a}=\sum_\beta \abs{a_\beta}^2.$
	At a point $z_0$ that does not belong to the singular support of $\sigma$ we have
	\begin{equation}
	2\partial_{z_i}\phi =\frac{1}{\abs{a}^2}\sum_\beta \bar{a}_\beta \partial_{z_i} a_\beta
	\end{equation}
	\begin{multline}
	2\partial_{z_i}\partial_{z_i}\phi =\frac{1}{\abs{a}^2} \sum_\beta \partial_{\bar{z}_i}\bar{a}_\beta \partial_{z_i} a_\beta
	-\frac{1}{\abs{a}^4}\left( \sum_\beta \bar{a}_\beta  \partial_{z_i} a_\beta\right)
	\left( \sum_\beta a_\beta  \partial_{\bar{z}_i} \bar{a}_\beta\right)=\\
	\frac{1}{\abs{a}^2} \sum_\beta \abs{\partial_{z_i}a_\beta}^2
	-\frac{1}{\abs{a}^4} \abs{ \sum_\beta a_\beta  \overline{\partial_{z_i} a_\beta}}^2
		\end{multline}
Next recall Lagrange's identity for $m$-dimensional $v,w$
\begin{multline}
\abs{v\wedge w}^2=\sum_{i=1}^{m}\sum_{j=i+1}^m \abs{v_iw_j-v_jw_i}^2=\frac{1}{2}
\sum_{i=1}^{m}\sum_{j=1}^{m} \abs{v_iw_j-v_jw_i}^2\\
=
\frac{1}{2}
\sum_{i=1}^{m}\sum_{j=1}^{m} (v_iw_j-v_jw_i)(\bar{v}_i\bar{w}_j-\bar{v}_j\bar{w}_i)=\\
\frac{1}{2}
\sum_{i=1}^{m}\sum_{j=1}^{m} \left(\abs{v_i}^2\abs{w_j}^2 -2\re (v_jw_i\bar{v}_i\bar{w}_j)+\abs{v_j}^2\abs{w_i}^2\right)=\\
=\left( \sum_{i=1}^{m} \abs{v_i}^2 \right)\left( \sum_{j=1}^{m} \abs{w_i}^2 \right)-\abs{\sum_{i=1}^{m} v_i\bar{w}_i}^2
=\abs{v}^2\abs{w}^2-\abs{v\cdot w}^2
\end{multline}
By the Lagrange identity we have
\begin{equation}
\abs{a\partial_{z_i} a}^2 +\abs{a\wedge \partial_{z_i} a}^2 =\abs{a}^2 \abs{\partial_{z_i} a}^2
\end{equation}
This implies
\begin{equation}
2\partial_{\bar{z}_i} \partial_{z_i} \phi_f =\frac{1}{\abs{a}^4} \abs{a\wedge \partial_{z_i} a}^2 =\frac{1}{\abs{a}^4}\sum_{\gamma,\beta}
	\abs{a_\gamma \partial_{z_i} a_\beta -a_\beta \partial_{z_i} a_\gamma}^2
	\end{equation}
	Obviously, summing over $i=1,\ldots,n$ yields a multiple of $\Delta \phi_f.$
	If $\abs{a(z_0)}\neq 0$ then there exists a neighborhood of $z_0$ on which $\abs{a}\geq \delta >0.$
	Thus, for constants $K,C,M$ we have
	\begin{equation}
	\sigma(z_0,t)=\int_{\abs{z-z_0}\leq t} K\frac{1}{\abs{a}^4} \sum_i 
	\sum_{\gamma,\beta}
	\abs{a_\gamma \partial_{z_i} a_\beta -a_\beta \partial_{z_i} a_\gamma}^2d\mu(z)\leq C\delta^{-4} Mt^{2n}
	\end{equation}
	Hence $\sigma$ is a positive measure. 
	This completes the proof.
		\end{proof}
Traour'e \cite{traore94} proved the following.
\begin{theorem}
Let $f(z)$ be a polyentire function on $\Cn$. Then $\nu(r)$ is bounded by $\nu_0$ if and only if
there exists a polyanalytic polynomial whose multi-degree, $\tau$, with respect to $z$ satisfies $\abs{\tau}\leq \nu_0$, and an entire function $g$
such that $f(z)=P(z,\bar{z})\exp(g(z)).$
\end{theorem}
\begin{proof}	
	Suppose $\phi_f(0)\neq -\infty.$ By Proposition \ref{traoreprop31}
	we have sufficiency of the given factorization, i.e.\ if $f(z)$ has the factorization 
	$P(z,\bar{z})\exp(g(z))$ as in the theorem then necessarily $\nu(r)\leq \nu_0.$
	For the converse implication suppose $\nu(r)\leq \nu_0.$ Let $\theta$ be the current associated to $\phi_f$, in particular
	its support is the support of the measure $\sigma$. 
	We have that if the support of $\theta$ does not contain the origin then
	\begin{equation}
	2i \partial\overline{\partial} V=\theta
	\end{equation}
	for (see Skoda \cite{skoda}, p.76, Theorem 3.26 and Lelong \cite{lelong64})
	\begin{equation}
	V(z):=C_n \int_{\Cn} \left(\frac{1}{\norm{u}^{2n-2}} -\frac{1}{\norm{u-z}^{2n-2}}\right) d\sigma(u)
	\end{equation}
	for a coefficient $C_n$ depending on $n.$ 
	However, it suffices instead of not containing the origin, that $\frac{\nu(t)}{t}$ is integrable on a
	neighborhood of the origin, see Skoda \cite{skoda}, p.164, and the latter is guaranteed by Proposition \ref{traoreprop32}.
	There thus exists a pluriharmonic function $h$ such that 
	\begin{equation}\label{ttrraa02}
	\phi_f=2\pi V +h
	\end{equation}
	In order to find a bound for 
	$V(z)$ we decompose
	\begin{multline}
	C^{-1}_n V(z)= \int_{\abs{u}\leq \abs{z}} \left(\frac{1}{\norm{u}^{2n-2}} -\frac{1}{\norm{u-z}^{2n-2}}\right) d\sigma(u)+
	\\
	\int_{\abs{u}> \abs{z}} \left(\frac{1}{\norm{u}^{2n-2}} -\frac{1}{\norm{u-z}^{2n-2}}\right) d\sigma(u)
	=:I_1 +I_2
	\end{multline}
	By the mean value theorem applied to $\abs{u}^{2-2n}$ (Lagrange's formula for finite increments)
	we have that $I_2$ is bounded by
	\begin{equation}
	I_2\leq (2n-2)\int_{\abs{u}>\abs{z}} \frac{\abs{z}d\sigma(u)}{\abs{u}^{2n-1}}
		\end{equation}
		Setting $\abs{z}=r$ we have by integration by parts
		\begin{equation}
		(2n-2)r\int_r^\infty \frac{d\sigma(u)}{t^{2n-2}}=(2n-2)r \frac{\sigma(t)}{t^{2n-2}}|_{r}^\infty
		+k(2n-2)r \int_r^\infty (2n-1)\frac{\sigma(t)}{t^{2n}}dt
		\end{equation}
		By the hypothesis on $\nu$ the integral in the right hand side is bounded by a constant independent of $r.$
		This bounds $I_2$.
		For sufficiently large $r$ we have, since $\phi_f(0)\neq 0$ yields $\nu(0)=0$ (where we denote by
		$\nu(z_0)$ the limit of $\nu(z_0,r)$ as $r\to 0$)
		for $I_1$ 
		\begin{equation}
		I_1\leq \int_{\abs{u}\leq r} \frac{d\sigma(u)}{\abs{u}^{2n-2}} =
		\int_{\abs{u}\leq 1} \frac{d\sigma(u)}{\abs{u}^{2n-2}}+
		\int_{1<\abs{u}\leq r} \frac{d\sigma(u)}{\abs{u}^{2n-2}}
		\end{equation}
		\begin{multline}
		\int_{\abs{u}\leq 1} \frac{d\sigma(u)}{\abs{u}^{2n-2}} = \frac{\sigma(t)}{t^{2n-2}}|_0^1 +(2n-2)\int_0^1\frac{\sigma(t)}{t^{2n-1}}dt
			\\
			=\frac{\pi^{n-1}}{(n-1)!}\left(\nu(1)+(2n-2)\int_0^1 \frac{\nu(t)}{t}dt\right)
			\end{multline}
			By Proposition \ref{traoreprop32} $\frac{\nu(t)}{t}$ is integrable on $[0,1]$, hence
			\begin{multline}
			I_1\leq O(1)+\int_{1\leq \abs{u}\leq r} \frac{d\sigma(u)}{\abs{u}^{2n-2}} =O(1)+\frac{\sigma(t)}{t^{2n-2}}|_1^r+\\
			(2n-2)\int_1^r \frac{\sigma(t)}{t^{2n-1}}dt\leq O(1)+\frac{\pi^{n-1}}{(n-1)!}(2n-2)\nu_0 \log r
			\end{multline}
			where $\nu_0$ is the bound for $\nu(r).$ This renders
			\begin{equation}
			V\leq C_n\left(O(1)+\frac{(2n-2)\pi^{n-1}}{(n-1)!}\nu_0\log r\right)
			\end{equation}
			where the coefficient $C_n=(n-2)!/(2\pi^{n-1})$ can be found in Skoda \cite{skoda}, p.76,
			to be a constant depending only on $n$, hence 
			\begin{multline}
			V\leq O(1)+\frac{(n-2)!}{2\pi^{n-1}} \frac{(2n-2)\pi^{n-1}}{(n-1)!}\nu_0\log r =\\
			O(1)+ \frac{1}{n-1}\frac{2(n-1)}{2}\nu_0\log r\leq
			 O(1) +\nu_0\log r  
			\end{multline}
			By Eqn.(\ref{ttrraa02}) this implies $\phi_f-h\leq O(1)+\nu_0\log r.$
			Let $\epsilon>0$. Since $h$ is pluriharmonic there exists an entire function $g$ such that $h=\re g$ which yields for
			small $\epsilon>0$
			\begin{equation}
			\frac{1}{2}\log \sum\abs{a_j}^2 -\log \exp(\re g)\leq (\nu_0 +\epsilon)\log r
			\end{equation}
			\begin{equation}
			\log\left( \frac{\left(\sum\abs{a_j}^2\right)^{\frac{1}{2}} }{\abs{\exp(g)}}\right)
				\leq (\nu_0+\epsilon)\log r
				\end{equation}
				which gives for $r\geq \epsilon$
				\begin{equation}
				\left(\sum\abs{a_j}^2\right)^{\frac{1}{2}}\abs{\exp(-g)}\leq r^{\nu_0 +\epsilon}
					\leq (\nu_0+\epsilon)\log r
					\end{equation}
					Hence for each $j$, $\abs{a_j \exp(-g)}\leq r^{\nu_0 +\epsilon}$ which gives $a_j=P_j\exp(g)$ for a complex polynomial in $z$ of multi-order $\kappa$
					with $\abs{\kappa}\leq \nu_0.$ This yields the wanted conclusion under the assumption $\phi_f(0)\neq -\infty.$
					If $\phi_f(0)=-\infty$ let $z_0$ be such that $\phi_f(z_0)\neq -\infty.$
					Setting $z=z_0+w$, the function $\psi(z)=\phi_f(z_0+w)$ is plurisubharmonic and bounded at $w=0.$ This reduces the proof to the first case.
					This completes the proof.	
	\end{proof}
\begin{corollary}
	Let $f(z)$ be a polyentire function on $\Cn$. 
	Suppose that the analytic components of $f(z)$ have no common zeros on $\{\abs{z}\geq r_0\}$ for some $r_0>0.$ If there exists a constant $C$ such that 
	$\abs{\mbox{Jac} \phi_f} \leq Cr^{-1}$ then $f(z)$ is a polyanalytic polynomial.
\end{corollary}

\begin{corollary}
	Let $f$ be a polyentire function on $\Cn.$ Then $f$ has the factorization $f(z)=P(z,\bar{z}) \exp(g(z))$ for an entire $g(z)$ and a polyanalytic polynomial $P(z,\bar{z})$,
	if and only if $\lim_{r\to \infty} \frac{\lambda(\phi_f,r)}{\log r}<\infty.$
	\end{corollary}

\section{Further notable results}
See also Chapter \ref{liouvillesec} for Liouville type theorems.
Balk \cite{balk1969} proved the following. 
\begin{theorem}
	If $f(z)$ is an entire reduced polyanalytic function of order $q$, $q\in \Z_+,$ without zeros
	then it has a representation of the form $f(z)=P(z,\bar{z})\exp(g(z))$, 
	where $g(z)$ is an entire function and $P(z,\bar{z})$ is a polyanalytic polynomial (which is necessarily polyanalytic of order $q$ i.e.\ of
	degree at most $q-1$ with respect to the
	$\bar{z}$ component). 
\end{theorem}
The following appears in a review of a local publication Balk \cite{balksmolensk69} (that we have not been able to obtain).
\begin{theorem}\label{balk1977}
		Let $f(z)=P(z,\bar{z})\exp(g(z))$ be a polyentire function of order $n$ with bounded zero set on $\C$, for a 
	polynomial $P$ and an entire holomorphic $g.$
	If $f(z)$ has no zeros, then the degree of $P$ is less than $2n-1$. 
		If $g(z)=\sum_{j=0}^m p_j(z,\bar{z})\exp(\alpha_j\bar{z}),$
	$m\geq 2,$ $\alpha_j\neq \alpha_k,$ for polynomials $p_j(z,\bar{z})$, then for each complex value $a$  
	the set $\{z:g(z)=a\}$ is unbounded. 
\end{theorem} 

\chapter{Approximation results}\label{approxsec}
\begin{definition}
For a compact $X\subset \Cn,$ 
by $C^0(X)$ the space of continuous functions on $X$ equipped with the sup-norm
$\norm{f}_\infty :=\sup_{x\in X} \abs{f}.$
Set $A_q(X)=C^0(X)\cap \mbox{PA}_q(\mbox{int}(X)),$ where 
$\mbox{int}(X)$ denotes the interior of $X.$ Denote by $P_q(X)$ the set of uniform limits, on $X$, of  
$q$-analytic polynomials and for a compact subset $Y\subseteq X,$ denote by $R_q(X,Y)$
the set of uniform limits, on $X$, of functions of the form $g(z)=\sum_{j=0}^{q-1}a_j(z)\bar{z}^j,$
such that $a_j(z)$ are rational functions with poles outside of $Y.$
We shall for the cases $q=1$ denote $A(X)=A_1(X)$, $P(X)=P_1(X)$
and $R(X)=R_1(X,\Cn)$. 
\end{definition}

Let us begin with an observation by Runge.
\begin{proposition}\label{rungesobs}
Let $X\subset\C$ be a compact subset and denote by $\mathscr{O}(X)$ the algebra 
of functions that are holomorphic in an open neighborhood of $X.$
Denote by $\overline{\mathscr{O}(X)}$ the closure of
$\mathscr{O}(X)$ with respect to the uniform norm on $X.$ Then
$R(X)=\overline{\mathscr{O}(X)}$
\end{proposition}
\begin{proof}
It is clear that $R(X)\subseteq\overline{\mathscr{O}(X)}$. For the converse,
let $U$ be an open neighborhood of $X$, let $f\in \mathscr{O}(X)$ and let
$\kappa$ be a rectifiable contour which has winding number $1$ at each point of $X.$
Then
\begin{equation}
f(z)=\frac{1}{2\pi i}\int_\kappa \frac{f(\zeta)}{\zeta -z} d\zeta,\quad z\in X
\end{equation}
The Riemann integral in the right hand side is for a fixed $z$ a limit of sums of the form 
\begin{equation}
\frac{1}{2\pi i}\sum_{j=1}^N c_j \frac{1}{\zeta_j -z} ,\quad \zeta_j\in \kappa
\end{equation}
for positive integers $N$ and constants $c_j.$
Since $X$ is contained in the set enclosed by $\kappa$ these sums converge uniformly on $X$.
Since $f$ was arbitrarily chosen in $\mathscr{O}(X)$
this shows $\mathscr{O}(X)\subset R(X).$ This completes the proof.
\end{proof}
Three early theorems have become so popular that they can be considered common knowledge so we shall
state them without proof, namely for the case $q=1$ we have Runge's theorem stating that $P_1(X)=R_1(X,X)$ if and only if $\C\setminus X$ is connected. We also have by the classical theorem of Mergelyan \cite{mergelyan}
stating that
$A_1(X)=P_1(X)$ if and only if $\C\setminus X$ is connected. The Oka-Weil\index{Oka-Weil theorem} approximation theorem states that
if $f_1,\ldots,f_r$ are entire functions on $\Cn$ then
$\Omega:=\{ z\in \Cn \colon \abs{f_j(z)}<1, 1\leq j\leq r\}$ is a Runge domain (any holomorphic function on $\Omega$  
can be uniformly approximated by complex polynomials).
The theory of rational approximation is wide and deep and many generalizations and variations of the above results exists in the literature, we refer
to Boivin \& Gauthier \cite{boivingauthier} and the
references therein, and shall not attempt to survey them here. Many results have original proof
utilizing the theory of uniform algebras. Such techniques will in general fail for the polyanalytic situation
as the spaces $A_q(X),$ $P_q(X)$ and $R_q(X)$ are no longer uniform algebras.
\section{The Malgrange-Lax theorem}
One of the earliest results on polyanalytic approximation is the following (which more generally concerns elliptic equations, thus contains 
powers of the Cauchy-Riemann operator in one variable, as a special case).
The result was found by Malgrange \cite{malgrange56} (and for order two also Lax \cite{lax56}) in 1956.
We present a proof by Narasimhan (see \cite{narasimhan}).
See Section \ref{ellipticapp2} in the appendix for definitions 
and preliminaries required for the statement and the proof.
We shall need the following result.
\begin{theorem}[See Narasimhan \cite{narasimhan}, p.197]\label{sobolevlemma}
Let $\Omega\subset \Rn$ be an open subset and $m>n/p$. Then for any compact $K\subset\Omega$ there exists a constant 
$C_K$ such that
for any $C^\infty$-smooth function $f:\Omega\to \Cn$ 
we have $\sup_{x\in K} \abs{f(x)}\leq C_{K}\abs{f}_{m,p}.$
\end{theorem}
\begin{proof}
Let $K'$ be a compact such that for any $y\in K,$ $g(y):x\mapsto f(x+y)\in C^\infty_c(K').$
Hence it suffices to prove that for $\Omega=\Rn$
and for each $f\in C^\infty_c(K')$ there is a $C$ such that $\abs{f(0)}\leq C\abs{f}_{m,p}.$
Let $\eta:\R_+\times S^{n-1}\to \Rn\setminus\{0\}$, $\eta(t,x):=tx.$
If $x_1,\ldots,x_n$ are the restrictions to $S^{n-1}$ of the Euclidean coordinates in $\Rn$ then the $(n-1)$ form
$\omega:=\sum_{k=1}^n x_kdx_1\wedge\cdots \wedge d\hat{x}_k\wedge\cdots \wedge dx_n$
(where $\hat{}$ denotes that the component is omitted) satisfies
$\eta^*(dy\wedge\cdots \wedge dy_n)=t^{n-1} dt\wedge \omega,$ $y\in \Rn\setminus \{0\}.$
Since $\int_U dy\wedge\cdots \wedge dy_n$ is positive over  $U\subset\Rn\setminus\{0\}$, $U=\eta(\{\frac{1}{2}<t<1\}\times S^{n-1}),$ we have $\int_{S^{n-1}} \omega \neq 0.$
If $f\in C^\infty(\Rn)$ set $g=f\circ\eta.$ Then there are homogeneous polynomials $q_\alpha$ of degree $m$ on $\Rn$ such that
\begin{equation}\label{gmref}
\partial_t^m g(t,x)=\sum_{\abs{\alpha}=m} q_\alpha\left(\frac{y}{\norm{y}}\right)\partial^\alpha f(y),\quad y=\eta(t,x)=tx
\end{equation}
in particular, the $q_\alpha\left(\frac{y}{\norm{y}}\right)$ are bounded.
For sufficiently large $M>0$ and $x\in S^{n-1}$ we have
\begin{equation}
f(0)= -\int_0^M \partial_t f(tx)dt=\frac{(-1)^m}{(m-1)!}\int_0^M t^{m-1}  \partial_t^m g(tx) dt 
\end{equation}
Multiplying by $\omega$ and integrating over $S^{n-1}$ we have
for a constant $C_m$
\begin{equation}
f(0)\int_{S^{n-1}}\omega=C_m \int_{S^{n-1}} \int_0^M t^{m-1}  \partial_t^m g(tx) dt \wedge \omega
\end{equation}
Since $\int_{S^{n-1}}\omega\neq 0$ this yields for a constant $C_m'$
\begin{multline}
f(0) = C_m' \int_{S^{n-1}}\int_0^M t^{m-n} \partial_t^m g(tx) dt\wedge \omega=\\
=C_m' \int_{\norm{y}\leq M}\int_0^M t^{m-n} \partial_t^m g(y) dy,\quad t=\norm{y}
\end{multline}
where $\partial_t^m g$ is given by Eqn.(\ref{gmref}). 
For $p=1$, $m\geq n,$ this gives $\abs{f(0)}\leq C_k\abs{f}_{m,1}$ for a constant $C_k.$ So assume $p>1.$
By H\"older's inequality
we have 
\begin{multline}
\abs{f(0)}\leq \\
C'\left(\int_{\{\abs{y}<M\}} t^{(m-n)p'}dy \right)^{\frac{1}{p'}}
\left(\int_{\{\abs{y}<M\}} \abs{\frac{\partial^m g(y)}{\partial t^m}}^p dy\right)^{\frac{1}{p}}
\end{multline}
Since $m>n/p$, $(m-n)p'+n-1>-1$ 
\begin{equation}
\int_{\norm{y}\leq M}t^{(m-n)p'}dy=\int_{S^{n-1}}\int_0^M t^{(m-n)p'}t^{n-1} dt\wedge \omega<\infty
\end{equation}
and since
$\partial^m g(tx)/\partial t^m =\sum_{\abs{\alpha}\leq m} q_\alpha(y)\partial^\alpha f(y)$ for bounded $q_\alpha$,
there exists a constant $C''>0$ such that
\begin{equation}
\int_{\norm{y}<M} \abs{\partial_t^m g(y)}^p dy \leq C'' \abs{f}_{m,p}^p
\end{equation}
Hence exists a constant $C_K>0$ depending on $K$ such that$\abs{f(0)}\leq C_K\abs{f}_{m,p}$. This proves the result for $p>1,$ for the case that $f$ is supported in $K.$
Finally, this can be applied to $\phi f$ for a $\phi\in C^\infty_c(\Omega)$ such that $\phi(x)=1$ on $K.$
This completes the proof.
\end{proof}

\begin{definition}
For an open subset $Y\subset\C$, denote by $h(Y)$ the so-called {\em hull of $Y$}\index{Hull of an open set in $\C$}, by which we mean the sum
of $Y$ and all connected components of $\C\setminus Y.$ More generally for a $C^\infty$-smooth manifold, $V$, with countable basis
we define the hull $h_V(Y)$ of a subset $Y\subset V$ to be the union of $Y$ and all connected components
of $V\setminus Y.$
\end{definition}
\begin{remark}\label{narasimrem}
It is clear that $Y\to h(Y)$ is increasing (i.e.\ If $Y_1\subset Y_2$ then $h(Y_1)\subset h(Y_2)$).
To see this note that if $C$ is a relatively compact connected component of $V\setminus Y_1$
and if $C'$ is a connected component of $V\setminus Y_2$ not contained in $C$ and if $C\cap C'\neq \emptyset,$ 
then $C\cup C'$ is a connected subset of $V\setminus Y_1$ properly containing $C$ which is impossible since $C$ is a component.
Thusn $C\setminus Y_2$ is a union of connected components of $V\setminus Y_2.$ Since $C$ is relatively compact $C\subset h_V(Y_2).$
This proves also that $h_V(h_V(Y))=h_V(Y).$
If $U$ is a relatively open subset containing a compact $Y$, let $\Omega_1,\ldots,\Omega_n$ be open connected sets with 
$\partial U\subset \cup_i \Omega_i,$ 
$\Omega_i \cap Y=\emptyset.$ Then there exists at most $n$ relatively compact connected components of $V\setminus Y$
which are not contained in $U$, thus $h(Y)$ is relatively compact and since $V\setminus h(Y)$ is the union 
of unbounded components of $V\setminus Y$, $V\setminus h(Y)$ is open. Hence $h(Y)$ is compact whenever $Y$ is compact..
If $Y$ is a compact subset of an open set $\Omega$ and if
$V\setminus \Omega$ has no compact components then $h(Y)\subset \Omega.$ To see this note that
if $W$ is a bounded component of $V\setminus Y$, not contained in $\Omega,$ we can let $a\in W$, $a\notin \Omega..$
If $V_a$ is the connected component of $V\setminus \Omega$ containing $a$ then we have $V_a\subset W,$ so
$V_a$ is relatively compact which yields a contradiction. Hence $h(Y)$ must belong to $\Omega.$
If $Y$ is compact and $Y=h(Y)$, then $Y$ has a fundamental system of open (compact) neighborhoods $U(W)$ such that $U=h(U)$ $(W=h(W)).$
It is known also that taking the hull preserves closedness. Indeed, if $Y$ is closed in $V$ any component of $V\setminus Y$ is open so $V\setminus h(Y)$ being a union of the components of $V\setminus Y$ that are relatively compact, is open.
\end{remark}

By Runge's theorem all holomorphic functions on $Y$ can be uniformly approximated by complex polynomials
if and only if $Y=h(Y).$

\begin{lemma}\label{narasimlem1}
Let $V$ be an oriented $C^\infty$-manifold of real dimension $n$, and let $E,F$ be 
$C^\infty$-vector bundles on $V$.
Let $L:E\to F$ be an elliptic operator of order $m$ with $C^\infty$ coefficients. 
Let $\Omega\subset V$ be an open subset and let $Lf=0,$ $Lf_\nu=0$ on $V$. 
Then the following are equivalent:\\
(i) $f_\nu\to f$ in $L^2$ locally on $\Omega$
\\
(ii) $f_\nu\to f$ uniformly on compact subsets of $V$
\\
(iii) $f_\nu$ and $\partial^\alpha f_\nu$ for every $\alpha,$ converge to $f$ and $\partial^\alpha f$ respectively, uniformly
on compact subsets of $V$ (note that because of elliptic regularity $f$ and $f_\nu$ are $C^\infty$).
\end{lemma}
\begin{proof}
(iii)$\Rightarrow$(ii) is trivial and so is (ii)$\Rightarrow$(i). Thus it suffices to prove (i)$\Rightarrow$(iii).
Let $K\subset V$ be a compact subset and let $r>0.$
We can assume $E,F$ are trivial and $V\subset \Rn$ is an open subset. Let $U,U'$ be open subsets such that
$K\subset U\Subset U'\Subset \Omega$.
By Proposition \ref{friedrichs} 
there exists a constant $C'$ such that for any $g\in C^\infty(\Omega)$
\begin{equation}
\abs{g}^U_{m+r} \leq C'( \abs{Lg}^U_r +\abs{g}_0^{U'})
\end{equation}
If $f_\nu \to f$ in $L^2(U')$ and $Lf_\nu=0$ on $V$ and $Lf=0$ on $\Omega$ this implies
\begin{equation}\label{doeroe}
\abs{f_\nu -f}_{m+r}^U\leq C'\abs{f_\nu -f}_0^{U'}
\end{equation}
Now $r>0$ was arbitrary thus
application of Theorem \ref{sobolevlemma} 
shows that $\norm{f_\nu -f}_{r}^K$ is bounded by a multiple of an expression of the form apprearing in the left hand side of Eqn.(\ref{doeroe}). 
Hence (i)$\Rightarrow$(iii). This completes the proof.
\end{proof}

\begin{theorem}[Malgrange-Lax]\label{malgrangelaxthm}
Let $V$ be an oriented real-analytic manifold and let $E,F$ be 
analytic vector bundles on $V$ with rank$(E)=$rank$(F)$.
Let $L:E\to F$ be an elliptic operator of order $m$ with analytic coefficients. 
Let $\Omega\subset V$
be an open subset such that $V\setminus \Omega$ has no compact connected components. 
Then
any solution $f\in C^\infty(\Omega,E)$ to $Lf=0$ on $\Omega$ 
can be approximated uniformly on compacts of $\Omega$ by
solutions $s\in C^\infty(V,E)$ to $Ls=0.$
\end{theorem}
\begin{proof}
Let $K\subset\Omega$ be a compact.
As we have pointed out in Remark \ref{narasimrem}, $K$ is compact then the hull $h(K)$ is compact and $h(K)\subset \Omega.$ 
We can suppose (if necessary by replacing $K$ with $h(K)$ that $K=h(K),$ and by Remark \ref{narasimrem}.
Let $K'\subset V$ be a compact such that $h(K)\subset \mbox{int}(K')$, and let
\begin{equation}
A(K'):=\{ f\in H^0(K',E)\colon Lf=0\mbox{ on }\mbox{int}(K')\}
\end{equation}
Let $S(K)$ be the restrictions of sections $f\in C^{\infty}(N,E)$ to $H_0(K,E)$ such that
\begin{equation}
Lf=0,\mbox{ some open neighborhood }N\mbox{ of} h(K)\end{equation}
Define the map $\eta:H^0(V,E)\to H^0(h(K),E)$ given by $\eta(s)=s$ on $h(K)$ and $\eta(s)=0$ outside $h(K).$ 
Then $\eta(A(K'))\subset \eta(S(K)).$ 
We show that is dense in $\eta(S(K)).$
Let $l$ be a continuous linear functional on $H^0(h(K),E)$ such
that $l(s)=0$ for $s\in \eta(A(K'))$. 
Recall that $H^m(K,E)$ denotes the set of sections $s:V\to E$ which are locally in $H^m$ with supp$s\subset K.$
We shall need the following.
\begin{lemma}\label{narasim9prop1}
For any continuous linear functional $l$ on $H^0(K, E)$, there 
exists a unique $s'\in H^0(K, E')$ such that $l(s) = \langle s,s'\rangle$, for any
$s'\in H^0(K, E)$.
\end{lemma}
\begin{proof}
If there exists $s'\in H^0(K,E')$ such that $l(s)=\langle s,s'\rangle$ for all $s\in H^0(K,E')$,
then $s'$ is unique. Let $U$ be a coordinate neighborhood such that 
$E|_U$ 
is trivial. 
Let $N$ be a compact subset of $U$. Then there is an $s'\in H^0(N,E')$ such that
$l(s)=\langle s's\rangle_E,$ for all $s\in H^0(N,E).$
Let $\tau:E_{U}\to U\times \Cn$ be an isomorphism and $\tau^*:E^*_{U}\to U\times \C^n$ the corresponding isomorphism of $E^*_U.$
Let $\psi:U\to \Omega$ be an isomorphism of $U$ onto an open subset $\Omega\subset \Rn$
and let 
\begin{equation}
(s_1,\ldots,s_r)=\tau(s)\circ\psi^{-1},\quad s\in H^0( ,E)
\end{equation}
Then $s_j\in L^2(\Omega)$ and vanish off $\psi(N).$
Now $L^2(U)$ is a Hilbert space so by the Riesz theorem there exists $t=(t_1,\ldots,t_n)\in L^2(U)$ that vanish off $\psi(N)$ such that
\begin{equation}
l(s)=(s,t)=\int_U\sum_{i=1}^n s_i t_i dx_1\wedge \cdots \wedge dx_n
\end{equation}  
Setting 
\begin{equation}
s':=(\tau^*)^{-1}(t_1\circ\psi,\ldots,t_r\circ \psi)\otimes \psi^*(dx_1\wedge \cdots \wedge dx_n)\end{equation}
we have $s'\in H^0(N,E')$ and $l(s)=\langle s',s\rangle_E,$ for all $s\in H^0(N,E).$
This proves Lemma \ref{narasim9prop1}.
\end{proof}
By Lemma \ref{narasim9prop1} 
there exists $u\in H^0(h(K),E')$ 
such that $l(s) =\langle s, u \rangle$ for $s \in H^0(K, E)$. Then $\langle s, u \rangle = 0$ for
every $s$ with $Ls = 0$ on
$K'$. 
We shall use the following result (the proof is based upon 
the fact that for $V$ a compact oriented $C^\infty$-smooth manifold, $E, F,$ $C^\infty$-smooth
vector bundles on $V$, and a differential
operator $L: C^\infty(V, E) \to C^\infty(V, F)$, we have that by definition $L^*s' = 0$ (where $L^*$ denotes the transpose operator) for $s'\in H^0(V,F')$ means precisely 
$\langle Lu,s'\rangle =0$ for all $u\in C^\infty(V,E)$).
\begin{proposition}\label{narasim9prop2prim}
Let $V$ be an oriented $C^\infty$-smooth manifold, $E, F,$ $C^\infty$-smooth
vector bundles on $V$. Let $L$ be an elliptic differential
operator $L: C^\infty_c(V, E) \to C^\infty(V, F)$. Let $K$ be a compact subset of $V$ and $s\in H^0(K,F)$ be such that
$\langle s, s' \rangle = 0$ for every $s' \in H^0(K, F')$ such that
$L^*s' = 0$ on $\mbox{int}(K).$ Then there exists $\sigma\in H^m(K,E)$ such that $L\sigma =s.$
\end{proposition}
\begin{proof}
	Let
	\begin{multline}
	N:=\{t\in H^0(K,E):\langle t',t\rangle_F=0,\forall t'\in H^0(K,F')\\
	\mbox{ s.t.} L^*t'=0\mbox{ on int}K\}
	\end{multline}
	Now $L^*t'=0$ on int$K$ means that $\langle Ls,t'\rangle_{F'}=0$ for all $s\in C^\infty_0(\mbox{int}K,E).$ Let $l$
	be a continuous linear functional on $H^0(K,E)$ that vanishes on $L_K(H^m(K,E)).$
	By Lemma \ref{narasim9prop1} there exists $t_0'\in H^0(K,F')$ such that $l(t)=\langle t_0',t\rangle_E$ for all $t\in H^0(K,E).$
	Since $l$ vanishes on $L_K(H^m(K,E))$ we have
	$l(Ls)=\langle t_0',Ls\rangle_E=0$ for all $s\in C^\infty_c(\mbox{int}K,E),$ such that$L^*t_0'=0$ on
	$\mbox{int}K.$ By the definition of
	$N$
	\begin{equation}
	\langle t_0',t\rangle_E=l(t)=0,\quad \mbox{ for all } t\in N
	\end{equation}
	Since $L_K(H^m(K,E))$ is closed and $l(t_0)=0$ this completes the proof.
\end{proof}

By Proposition \ref{narasim9prop2prim}, 
there exists
$v\in H^m(K', F')$ such that $L^*v = u$ (where $L^*$ denotes the transpose operator).
Since supp$u\subset h(K)$ i.e.\ $L^*v = 0$ on $V\setminus h(K),$ 
by Theorem \ref{analyticity}, $v$ is analytic on $V\setminus h(K).$ Since
supp$v\subset K'$ and $V\setminus h(K)$
has no relatively compact
connected components we have $v=0$ on $V\setminus h(K)$ i.e.\ $v\in H^m(h(K),F')$.
For $s\in S(K)$, let $U$ be a neighborhood of $h(K)$ so that $s$ is defined and
$Ls = 0$ on $U$. Then  $\langle \eta(s), u \rangle =\langle s, u \rangle =$
$\langle s, L^*u \rangle_U=$
$\langle Ls, u \rangle_U= 0,$
because supp$v\subset h(K)$, i.e.\ $l(s)=0$ for any $s \in \eta(S(K))$.
By the Hahn-Banach Theorem this proves that $\eta(A(K'))$ is dense in $\eta(S(K)).$
Thus if $Lf=0$ in a neighborhood of $h(K)$ then there exists a sequence of functions
$\{ f_\nu\}$ in $A(K')$ such that $f_\nu\to f$ in $H^0(h(K),E)$ such that by Lemma \ref{narasimlem1}, $f_\nu \to f$
uniformly on compact sets in $\mbox{int}(h(K))$.
Let $\{K_r\}$ be a sequence of compact subsets such that $\cup K_r = V$
and $h(K)\subset \mbox{int}(K_1),$ $h(K_1)\subset \Omega,$ $h(K_r)\subset \mbox{int}(K_{r+1}),$ for $r\geq 1.$
So if $L f = 0$ in a neighborhood
of $h(K_1)$ there exists $f_1\in A(K_2)$ such that
$\norm{ f - f_1}_{h(K)} <\epsilon/2.$
By induction there exists a sequence $\{f_\nu\}$ in $A(K_{\nu+1})$ such that 
$\norm{ f_\nu - f_{\nu+1}}_{h(K_\nu)} <\epsilon/2^\nu,$ with $f_\nu\in C^\infty\mbox{int}(K_{\nu+1}).$
Set $g_r:=f_r +\sum_{s=r+1}^\infty(f_s-f_{s-1})$ on $K_r$
and $g$ on $V$ defined by $g:=g_r$ on $K_r.$ The series is $\lim_{s\to \infty} f_s$ and the converges uniformly on compacts of $V.$
Furthermore, $g_r=g_{r+1}$ on $K_r$
and for any section $u\in C^\infty_c(V,F')$ we have $(Lg)(u) =\langle g, L^*u \rangle =$
$\lim_{s\to \infty} \langle  f_s, L^*u \rangle=$
$\lim_{s\to \infty} \langle  Lf_s, u \rangle=0,$
where $\langle  f_s, L^*u \rangle=\langle  Lf_s, u \rangle$ if supp$u \subset K_s$. 
Then $\norm{ f- g}_{h(K)} <\epsilon$ and by Lemma \ref{narasimlem1} $Lg=0.$
This completes the proof.
\end{proof}
Combining this 
Now $\partial_{\bar{z}}^q$ being elliptic on $\R^2\simeq \C$ yields
after combining Theorem \ref{malgrangelaxthm} with Lemma \ref{narasimlem1}
the following.
\begin{corollary}
Let $q\in \Z_+$ and let $\Omega\subset\C$ be an open subset such that $\C\setminus \Omega$ has no compact connected components. 
Then for any $q$-analytic function $f$ on $\Omega$ there exists a sequence $\{f_j\}_{j\in \Z_+}$ of
entire $q$-analytic functions such that for all $k,l\in \Z_{\geq 0},$ 
$\partial_z^k\partial_{\bar{z}}^l f_j\to \partial_z^k\partial_{\bar{z}}^l f$
uniformly on compacts of $\Omega$.
\end{corollary}
We mention that Gauthier \& Zwick \cite{gauthierzwick}, p.364. prove the following (we state only a special case of their results for sheaves).
\begin{theorem}
	Let $L$ be an elliptic operator on an open subset $\Omega\subset \Rn.$
	The following are equivalent:\\
	(a) If $U\subset\Omega$ is open and $K\subset U$ compact such that each precompact component $M$ of $\Omega\setminus K$ meets $\Omega\setminus U$, then each
	.then any solution to $Lg=0$ on $K$, is the uniform limit on $K$ of solutions to $Lu=0$ on $U.$ \\
	(b) For any two open sets $U\subset V$ in $\Omega$, if $V\setminus U$ has no compact components, then any solution to 
	$Lu=0$ on $U$ may be uniformly approximated on compacts by restrictions $p_j=\tilde{p}_j|_U$ of solutions to $L\tilde{p}_j=0$ on $V.$\\
	\\
	If $K\subset V$ is a compact where $V\subset\Omega$ is open, and if $K=\widehat{K}_V$ (by which we denote the hull of $K$ with respect to $V$),
	and if $Lg=0$ on $K$, then $g$ is the uniform limit on $K$ of solutions to $Lu=0$ on $V.$ 
	\end{theorem}

\section{Vitushkin's theorem}

The following remark points out a technique which is often used for reducing an approximation proof to the case of smooth functions.

\begin{remark}\label{introapproxbeviser}
	Regularization (see Appendix \ref{elliptch}) is often used in proofs on approximation in the theory of polyanalytic functions. A recurring technique is to prove the result
	for $C^\infty$-smooth functions, and then use regularization of a $C^1$-smooth solution, and use the approximation property
	on each member, and round off by a standard diagonal argument in order to obtain the wanted result also for $C^1$-smooth solutions.
	The argument goes as follows. Usually one considers an elliptic operator $L$ on an open subset $\Omega$ of $\Rn,$
	and one works on Sobolev spaces, in the vector-valued situation we have
	$(f_1,\ldots,f_n):\Omega\to \C^q$, and the norms given by $\abs{f}_{m,p}^\Omega$.
	Let
	$\phi\in C^\infty_c(\Rn),$
	$\phi\geq 0,$ $\mbox{supp}(\phi)\subset \{x\in \Rn:\abs{x}<1\}$, $\int_{\Rn} \phi(x)dx=1.$ 
	Set $\phi_\epsilon(x):=\epsilon^{-n} \phi(x/\epsilon).$ Then 
	for any $f\in L^p(\Omega,q)$ that has weak derivatives up to order $m$ in $L^2$
	we consider the convolution (regularization)
	\begin{equation}
	f_\epsilon:=\phi_\epsilon *f(x)=\int_\Omega \phi_\epsilon(x-y)f(y)dy
	\end{equation}
	It is clear that $f_\epsilon$ has compact support and $\partial^\alpha \phi_\epsilon$ is bounded and goes to $0$ at every point if $\abs{\alpha}\geq 1,$
	and $\phi_\epsilon\to 1.$
	Thus
	\begin{equation}
	\partial^\alpha f_\epsilon =\sum_{\beta\leq \alpha}\binom{\alpha}{\beta} \partial^\beta \phi_\epsilon \partial^{\alpha-\beta} f
	\end{equation}
	implies that
	\begin{equation}
	\partial^\alpha f_\epsilon \to \partial^\alpha f\mbox{ in }L^2(\Rn,q)
	\end{equation}
	It also converges in $H^m(\Rn).$
	Let us look at the details for the case of $f\in L^p(\Omega),$ i.e.\ $f$ is complex-valued.
	Extend $f$ to $\Rn$ by setting $f=0$ on the complement of $\Omega.$
	we have
	\begin{equation}
	(f_\epsilon -f)(x)=\int_{\abs{y}\leq \epsilon} \phi_\epsilon(-y)(f(y+x)-f(x))dy
	\end{equation}
	By H\"older's inequality we have for $p>1$
	\begin{equation}
	\abs{(f_\epsilon -f)(x)}^p\leq \left(\int_{\abs{y}\leq \epsilon}\abs{\phi_\epsilon(y)}^{p'}\right)^{\frac{p}{p'}}
	\int_{\abs{y}\leq \epsilon} \abs{f(y+x)-f(x)}^p dy
	\end{equation}
	where $1/p+1/p'=1$, for $p=1$ we use $\sup \abs{\phi_\epsilon(y)}=\epsilon^{-n}\sup\abs{\phi(y)}.$
	Integrating over $\Rn$ and taking the $p$:th roots gives
	\begin{equation}
	\norm{(f_\epsilon -f)(x)}_{L^p}\leq \epsilon^{\frac{n}{p}}
	\left(\int_{\abs{y}\leq \epsilon}dy\right)^{\frac{1}{p}}\norm{\phi}_{L^{p'}}\sup_{\abs{y}\leq \epsilon}
	\left(
	\int_{\abs{y}\leq \epsilon} \abs{f(y+x)-f(x)}^p dy
	\right)^{\frac{1}{p}}
	\end{equation}
	and the right hand side goes to $0$ as $\epsilon\to 0$ for $f$ with compact support, and for $f\in L^p(\Omega)$
	the same holds true since continuous functions with compact support are dense.
	Also note that if $\{g_j\}_{j\in \Z_+}$ is a sequence in $C^\infty_c(\Omega)$ converging to $f$ in $L^2(\Omega)$
	where $f$ has weak derivatives up to and including order $m$
	then
	\begin{multline}
	\partial^\alpha f_\epsilon =\int_\Omega \partial^\alpha \phi_\epsilon(x-y)f(y)dy 
	=
	\lim_{j\to \infty} \int_\Omega \partial^\alpha \phi_\epsilon(x-y)g_j(y)dy
	=\\
	\lim_{j\to \infty} \int_\Omega  \phi_\epsilon(x-y) \partial^\alpha f(y)dy
	=(\phi_\epsilon *(\partial^\alpha f))(x)
	\end{multline}
\end{remark}
One can define the analytic (Ahlfors) capacity $\gamma,$ and continuous capacity $\alpha$\index{Continuous capacity} as follows: 
For a compact $K\subset \C$, call a function $f$, $K$-admissable if $f$ is holomorphic on $(\C\cup\{\infty\})\setminus K$ such that $f(\infty)=0.$
Denote by $B_0$ the class of $K$-admissable\index{$K$-admissable function} functions in the unit ball of the space of bounded Borel functions on $\C\setminus K.$
Denote by $B_1$ the class of $K$-admissable functions in the unit ball of the space of continuous functions on $\C\cup \{\infty\}.$
We define
\begin{equation}
\gamma(K):=\sup_{f\in B_0} \left\{ \lim_{z\to \infty} \abs{zf(z)}\right\}=\sup_{f\in B_0}\abs{f'(\infty)}
\end{equation}
\begin{equation}
\alpha(K):=\sup_{f\in B_1} \left\{ \lim_{z\to \infty} \abs{zf(z)}\right\}=\sup_{f\in B_1}\abs{f'(\infty)}
\end{equation}
Here $f'(\infty):=\lim_{z\to \infty} z(f(z)-f(\infty)),$
$f(\infty):=\lim_{z\to \infty} f(z).$
For an arbitrary bounded set $V$ define the analytic (Ahlfors) capacity $\gamma(V)$ of $V$ and the continuous capacity $\alpha(V)$ of $V$ by
\index{Continuous capacity}
\begin{equation}
\gamma(V):=\sup\{\gamma(K):K\subset V,K\mbox{ compact}\}
\end{equation}
\begin{equation}
\alpha(V):=\sup\{\alpha(K):K\subset V,K\mbox{ compact}\}
\end{equation}

The following important results of Vitushkin must be mentioned because some of the methods of proof
of these results has influenced generalization results on approximation to the polyanalytic case.

\begin{theorem}[Vitushkin \cite{vitushkin59}]\label{vitushkins}
Let $X\subset \C$ be a compact subset. The following are equivalent:\\
(i) $R(X)=C(X)$.\\
(ii) For every open disc $D$, $\gamma(D)=\gamma(D\setminus X).$
\end{theorem}

\begin{theorem}[Vitushkin \cite{vitushkin66}]\label{vitushkins0}
Let $X\subset \C$ be a compact subset. The following are equivalent:\\
(i) $R(X)=A(X)$.\\
(ii) For every open disc $D$, $R(K\cap \overline{D})=A(K\cap \overline{D}),$
where $\overline{D}$ denotes the closure of $D.$
\end{theorem}
Theorem \ref{vitushkins} can equivalently be stated as follows.
\begin{theorem}[Vitushkin \cite{vitushkin67}]\label{vituskin1}
Let $X\subset \C$ be a compact subset. The following are equivalent:\\
(i) $R(X)=A(X)$.\\
(ii) For every open disc $D$, $\alpha(D\setminus\mbox{int}(X))=\alpha(D\setminus X).$
\end{theorem}
We shall give a proof of Theorem \ref{vituskin1}, which has been chosen in part because it describes a scheme
for obtaining rational approximation that is also described in Gamelin \& Garnett \cite{gamelingarnett}. The proof follows Zalcman \cite{zalcmanbok}, see also Section \ref{vitushkinapp} in the Appendix, for
auxiliary results.
Let us begin by presenting a less complicated (compared to Vitushkin's capacity theorem which we shall later prove) preliminary result of Vitushskin which gives a taste of some of the methodology
employed in the theory of approximation in terms of capacity conditions.
\begin{definition}\index{$C(S,m)$}
	For an arbitrary set $S\subset \C$ and $m>0$ denote by 
	$C(S,m)$ the set of all continuous functions $f$ on $\hat{\C}$ 
	such that $f(\infty)=0,$ $f$ is analytic off some compact subset of $S$ and
	$\norm{f}_\infty \leq 1.$
\end{definition} 
Let $g(z)$ be a non-negative $C^\infty_c(\{\abs{z}<1\})$ function such that $\int_\C g(z)d\mu(z)=1$.  
Divide the plane into equal squares, $\{\sigma_j\}_{j\in \N},$ each of side $1$, the family being 
without pairs of common interior points, and let $g_k(\zeta):=\int_{\sigma_k}g(z-\zeta)d\mu(z)$.
Then $g_k$ is bounded, $C^\infty$-smooth and non-negative. We have $\sum_{k=1}^\infty g_k(z)\equiv 1$, since $\int g(z-\zeta)d\mu(z)=1.$
Setting $g_k(z,\delta):=g_k(\frac{6z}{\delta})$ we obtain (see Vitsuhkin \cite{vitushkin67}) with
	$\{g_k(z,\delta)\}_{k\in \N}$ a family of non-negative real $C^\infty$-smooth functions such that
	for each $\delta>0$,
	$\sum_{k=1}^\infty g_k(z,\delta)\equiv 1$; $\max_z \abs{ D g_k(z,\delta)}\leq \frac{\lambda}{\delta}$ for a constant $\lambda>0$ (independent of $k,\delta$);
	for each $\delta,k$ the support $v_k^\delta$ of $g_k(z,\delta)$ is convex with diameter $<\delta$;
	for each $\delta$ and each circle $\xi$ of radius $\delta$ 
	the number of values of different $k$ such that $\xi\cap v^\delta_k\neq\emptyset$, is $\leq \lambda.$
	Let for $\delta>0$ and a continuous function $\omega(f,\delta):=\sup_{\abs{z-\zeta}\leq\delta}\abs{f(z)-f(\zeta)}$, which is
	called the modulus of continuity.
	\begin{proposition}[Vitushkin \cite{vitushkin67}, Ch.2, parag. 4 Lemma 1]\label{vitusanvandicarmona}
		Let $f(z)$ be a continuous function with modulus of continuity $\omega(f,\delta)$ such that
		$f(z)$ has a representation of the form 
		$f(z)=\sum_{k=1}^\infty f_k^\delta(z)$ where
		\begin{equation}
		f_k^\delta(z)=\frac{1}{2\pi i}\int \frac{g_k(\zeta,\delta)}{\zeta -z}\partial_{\bar{\zeta}} f d\zeta\wedge d\zeta
		\end{equation}
		Let $\varphi(z)$ be a function of the form $\varphi(z)=\sum_{k=1}^\infty \varphi_k(z)$
	such that for all $k$, $\varphi_k\in C(v^{r,\delta}_k,m\omega(f,\delta))$ and
	\begin{equation}\label{vitsuekvnnn1}
	\lim_{z\to \infty} z^2(f^\delta_k(z)-\varphi_k(z))=0
	\end{equation}
	where $v^{r,\delta}_k$ is the set of points of $v^\delta_k$ within radius $r\delta,$ for $r\geq 0.$
	Then $\abs{\varphi(z)-f(z)}<\epsilon\leq c(m,r)\omega(f,\delta)$, where $c(m,r)$ depends only upon $m$ and $r.$
		\end{proposition}
		\begin{proof}
		Set $f_k(z)=f_k^\delta (z)-\varphi_k(z).$ For fixed $z_k\in v_k^\delta$ we write the Laurent expansion
		\begin{equation}
		f_k(z)=\sum_{n=1}^\infty a_{kn}(z_k-z)^{-n}
		\end{equation}
		By Eqn.(\ref{vitsuekvnnn1}) we have $a_{k1}=a_{k2}=0.$ 
		By the condition of the lemma
		$\abs{\varphi_k(z)}\leq m\omega(f,\delta)$.
		Also for fixed $z_1$ in the closure of the support, $U$, of $\partial_{\bar{z}} f$ we have
		\begin{equation}
		f^\delta_k(z_1)=\frac{1}{2\pi i}\int_U \frac{g_k(\zeta,\delta)}{\zeta -z_1}\partial_{\bar{\zeta}}(f(\zeta)-f(z_1))
		\end{equation}
		so that
		\begin{equation}
		\abs{f^\delta_k(z)}\leq \frac{1}{2\pi}\omega(f,\delta)\max_z \abs{\partial_{\bar{z}}g_k}\int_U\abs{\frac{d\zeta\wedge d\bar{\zeta}}{\zeta -z_1}}
		\end{equation}
		and since
		\begin{equation}
		\int_U\abs{\frac{d\zeta\wedge d\bar{\zeta}}{\zeta -z_1}}=\int_U\frac{2dx\wedge dy}{\abs{\zeta -z_1}}\leq 2\int_{\abs{\zeta}\leq \delta}
		\abs{\frac{dx\wedge dy}{\zeta}}=4\pi\delta
		\end{equation}
		we have 
		that
		$\abs{f^\delta_k(z)}\leq 2\lambda\omega(f,\delta).$
		Thus
		\begin{equation}
		\abs{f_k(z)}\leq \abs{f_k^\delta (z)}+\abs{\varphi_k(z)}\leq (m+2\lambda)\omega(f,\delta)=m_1\omega(f,\delta)
		\end{equation}	
		\begin{lemma}\label{prelemmetvit}
			If $F\in C(U,m)$ for a subset $U\subset\C$ and positive $m$, such that $F(z)$ has Laurent expansion about $z_0$ given by 
			$F(z)=\sum_{k=1}^\infty$ then
			$\abs{a_n}\leq 4m\alpha(U)(2r)^{n-1},$ where $r=\sup_{z\in U} \abs{z-z_0}.$
		\end{lemma}
		\begin{proof}
			Set
			\begin{equation}
			F_1(z):=\frac{1}{2\pi i}\int_{\partial U} \frac{f(\zeta)(\zeta -z_0)^{n-1}}{\zeta -z}d\zeta,\quad F_2(z):=\frac{1}{2\pi i}\int_{\abs{\zeta -z_0}=2r} \frac{f(\zeta)(\zeta -z_0)^{n-1}}{\zeta -z}d\zeta
			\end{equation}
			and let $\psi(z)$ be defined as the function which coincides with $F_1(z)$ off $U$ and coincides with
			$F(z)(z-z_0)^{n-1}-F_2(z)$ on $\{\abs{z-z_0}<2r.$ Then $\psi\in C(U,\mu)$ with $\mu=\max\abs{\psi(z)}.$
			For $\abs{z-z_0}\leq r$ we have
			\begin{equation}
			\abs{F_2(z)}=\abs{\frac{1}{2\pi i}\int_{\abs{\zeta-z_0}=2r} \frac{F(\zeta)(\zeta -z_0)^{n-1}}{\zeta -z}d\zeta}\leq
			\frac{1}{2\pi}\int \frac{m(2r)^{n-1}}{r}\abs{d\zeta}\leq m2^nr^{n-1}
			\end{equation}
			For $\abs{z-z_0}\leq 2$ we have $\psi(z)=F(z)(z-z_0)^{n-1}-F_2(z)$ and
			\begin{equation}
			\mu\leq \max_{\abs{z-z_0}\leq r}\abs{F(z)(z-z_0)^{n-1}}+\max_{\abs{z-z_0}\leq r}\abs{F(z)}\leq mr^{n-1}+m2^{2n}r^{n-1}\leq 4m(2r)^{n-1}
			\end{equation}
			Next note that  
			an examination of the definition shows that
			for $\theta\in C(U,1)$,
			$\alpha(U)=$ whereas for $\theta\in C(U,\max\abs{\theta})$ we have
			$\frac{\gamma(U,\theta)}{\max\abs{\theta}}\leq \alpha(U).$ 
				Since $\psi(z)=F_1(z)$ on $U$ we have $\gamma(U,F_1)=\gamma(U,\psi)$ which implies
				\begin{multline}
				\abs{a_n}= \abs{\frac{1}{2\pi i}\int_{\partial U} F(z)(z -z_0)^{n-1}dz}=
				\abs{\frac{1}{2\pi i}\int_{\partial U} (F_1(z)+F_2(z))dz}=\\
				\abs{\frac{1}{2\pi i}\int_{\partial U} F_1(z)dz}=\abs{\gamma(U,F_1)}=\gamma(U,\psi)\leq \alpha(U)\mu\leq 4m\alpha(U)(2r)^{n-1}
					\end{multline}
					This proves Lemma \ref{prelemmetvit}.
					\end{proof}
					Recall that $a_{k1}=a_{k2}=0$. 
					By Lemma \ref{prelemmetvit} we have $\abs{a_{kn}}\leq 4m_1\omega(f,\delta)(m_2\delta)^n$ for $\abs{z_k-z}\geq 2m_2\delta$ we have
					\begin{multline}
					\abs{f_k(z)}\leq \sum_{n=3}^\infty \abs{a_{kn}(z_k -z)^{-n}}\leq 4m_1\omega(f,\delta)\sum_{n=3}^\infty (m_2\delta)^n\abs{z_k -z}^{-n}\leq \\
					4m_1\omega(f,\delta)\delta^3\abs{z_k-z}^{-3}\sum_{n=3}^\infty (m_2\delta)^{n-3}(2m_2\delta)^{-n+3}=m_3\omega(f,\delta)\delta^3\abs{z_k -z}^{-3}
					\end{multline}
					By the property of the number of different $k$ for which the $v^\delta_k$ intersect the circle, we have
					\begin{multline}
					\abs{f(z)-\varphi(z)}\leq \sum_{k} \abs{f_k(z)}=\sum_{k:\abs{z_k-z}\leq 4m\delta}\abs{f_k(z)}+
					\sum_{k:\abs{z_k-z} > 4m\delta}\abs{f_k(z)}\leq \\
					m_4\omega(f,\delta)+\sum_{k:\abs{z_k-z} > 4m\delta} m_3\omega(f,\delta)\cdot \delta\cdot \delta^2\abs{z_k -z}^{-3}\leq\\
					m_4\omega(f,\delta)+m_5\omega(f,\delta)\delta\int_{\abs{\zeta -z}\geq 4m\delta} \abs{\zeta -z}^{-3}dx\wedge dy\leq \\
					m_4\omega(f,\delta)+m_6\omega(f,\delta)=c(m,r)\omega(f,\delta)
					\end{multline}
					This completes the proof of Proposition \ref{vitusanvandicarmona}.
				\end{proof}
\begin{definition}
Let $X$ be a compact set in the plane. Let $f\in C^0(X)$ and $\mu$ a complex regular Borel measure on $X.$
We say that $\mu$ is orthogonal to $R(X)$ if for each $f\in R(X)$ we have
$\int_X fd\mu=0$, we denote the latter condition by $\mu \perp f.$
\end{definition}
Let $f\in C^0(X).$ Since $R(X)$ is a closed subspace of the Banach space $C^0(X)$ under the sup-norm, $f\in R(X)$ iff $\mu \perp f$
for any $\mu$ orthogonal to $R(X).$ Hence $R(X)=A(X)$
if and only if they have the same set of orthogonal measures. Also $R(X)=C^0(X)$ iff only the zero measure is orthogonal
to $R(X).$
For a complex regular Borel measure, $\mu$, on $X$
denote $\hat{\mu}:=\int (\zeta -z)^{-1} d\mu(\zeta).$ 
\begin{lemma}
$\hat{\mu}=0$ on $\C\setminus X$ iff $\mu \perp R(X).$
\end{lemma}
\begin{proof}
Necessity is clear. Suppose $\hat{\mu}=0$ on $\C \setminus X.$ Then $\mu\perp (\zeta -z)^{-1}$ for 
every $z\in \C\setminus X.$ By fractional decomposition $\mu$ annihilates any rational function with poles off $X.$ Since any element of $R(X)$ can be 
approximated by such functions this completes the proof.
\end{proof}
\begin{theorem}
Let $f\in A(X)$ such that $f$ can be extended to a continuously differentiable function in a neighborhood of $X.$
Then $f\in R(X).$
\end{theorem}
\begin{proof}
Suppose $f$ is continuously differentiable function with compact support. Then
\begin{equation}
f(z)=-\frac{1}{\pi}\int_\C \frac{\partial_{\bar{\zeta}}f(\zeta)}{\zeta -z} d\zeta
\end{equation}
Suppose $\mu \perp R(X).$ Then $\hat{\mu}=0$ on $\C\setminus X.$ Since $f$ is analytic 
on $\mbox{int}(X),$ we have by continuity that $\partial_{\bar{\zeta}}f=0$ on $X,$
hence $\hat{\mu}\cdot \partial_{\bar{\zeta}}f =0$ everywhere which implies
\begin{equation}
f(z)=\int_X f(z)d\mu(z)=\frac{1}{\pi} \hat{\mu}(\zeta)\partial_{\bar{\zeta}}f(\zeta) =0
\end{equation}
This completes the proof.
\end{proof}

Let as usual $\mu$ denote the standard Lebesgue measure on $\C,$
in the sense that $d\mu(z)=dx\wedge dy,$ $z=x+iy.$
We recall, for the readers convenience the following well-known result,
see e.g.\ Eqn.(\ref{tgreffen}) in Section \ref{areolarsec}.
\begin{lemma}\label{browderlemma} 
Let $\Omega\subset \C$ be a bounded domain with smooth positively oriented boundary $\kappa.$
Let $f\in C^1(U)$ for a domain $U$ containing $\overline{\Omega}$.
Then for all $z\in \Omega,$
\begin{equation}
f(z)=\frac{1}{2\pi i}\int_\kappa \frac{f(\zeta)}{\zeta -z} -\frac{1}{\pi} \int_\Omega \frac{\partial_{\bar{z}} f}{\zeta -z}d\mu(\zeta)
\end{equation}
\end{lemma}
\begin{proof}
Let $\epsilon>0$ be sufficiently small such that
$\{\abs{z}\leq \epsilon\}\subset \Omega.$ Set $\Omega_\epsilon 
:=\Omega \setminus \{\abs{z}\leq \epsilon\}.$
Then $\Omega$ is a bounded domain with smooth boundary
where we assume $\{\abs{z}\leq \epsilon\}$ has positive orientation.
By the Green formula
\begin{multline}
\int_\kappa \frac{f(\zeta)}{\zeta -z}d\zeta -\int_{\{\abs{z}= \epsilon\}} \frac{f(\zeta)}{\zeta -z}d\zeta
=\int_{\Omega_\epsilon} \partial_{\bar{\zeta}}\left(\frac{f(\zeta)}{\zeta -z}\right)d\mu(\zeta)=
2i\int_{\Omega_\epsilon}  \frac{\partial_{\bar{\zeta}}f(\zeta)}{\zeta -z}d\mu(\zeta)
\end{multline}
Also
\begin{equation}
\int_{\{\abs{z}= \epsilon\}} \frac{f(\zeta)}{\zeta -z}d\zeta =i\int_0^{2\pi} f(\zeta +\epsilon\exp(i\theta))d\theta
\stackrel{\epsilon\to 0}{\longrightarrow} 2\pi if(\zeta)
\end{equation}
and
\begin{equation}
\int_{\Omega_\epsilon} \frac{ \partial_{\bar{\zeta}}f(\zeta)}{\zeta -z}d\zeta
\stackrel{\epsilon\to 0}{\longrightarrow}
 \int_{\Omega}\frac{ \partial_{\bar{\zeta}}f(\zeta)}{\zeta -z}d\zeta
\end{equation}
This completes the proof.
\end{proof}
A direct consequence is that if $f\in C^1_c(\C)$ where the subindex $c$ denotes compact support)
the by replacing $\kappa$ with large circles, we have for all $z\in \C,$
\begin{equation}\label{vitusanvekm}
f(z)=\frac{1}{\pi}\int_\C \frac{ \partial_{\bar{\zeta}}f(\zeta)}{\zeta -z}d\mu(\zeta)
\end{equation}

\begin{lemma}\label{browderlemma0} 
Let $K\subset \C$ be a measurable set.
Then for all $z\in \C$
\begin{equation}
\int_K \frac{d\mu}{\abs{\zeta -z}} \leq 2\sqrt{\pi \mu(K)}.
\end{equation}
\end{lemma}
\begin{proof}
If $\mu(K)=\infty$ or $0$ the conclusion is trivial so assume $\mu(K)<\infty.$
We have
\begin{equation}
\int_{\{\abs{z}< \sqrt{\mu(K)/\pi}\}} \frac{d\mu(\zeta)}{z-\zeta}
=\int_0^{2\pi}\int_0^{\sqrt{\mu(K)/\pi}} \frac{1}{r} rdrd\theta =2\sqrt{\pi \mu(K)}
\end{equation}
We have $\mu(K)=\mu(\{\abs{z}< \sqrt{\mu(K)/\pi}\})$ so 
Write 
\begin{equation}
K=(K\cap \{\abs{z}< \sqrt{\mu(K)/\pi}\})\cup (K\setminus \{\abs{z}< \sqrt{\mu(K)/\pi}\})
\end{equation}
and note that
\begin{equation}
\mu(K\setminus \{\abs{z}< \sqrt{\mu(K)/\pi}\})=\mu(\{\abs{z}< \sqrt{\mu(K)/\pi}\}\setminus K)
\end{equation}
Since
$1/\abs{z-\zeta}\geq 1/\sqrt{\mu(K)/\pi}$ on $\{\abs{z}< \sqrt{\mu(K)/\pi}\}\setminus K$
and $1/\abs{z-\zeta}\leq 1/\sqrt{\mu(K)/\pi}$
on $K\setminus \{\abs{z}< \sqrt{\mu(K)/\pi}\}$
we obtain
\begin{equation}
\int_{K\setminus \{\abs{z}< \sqrt{\mu(K)/\pi}\}} \frac{d\mu}{\abs{\zeta -z}} \leq \int_{\{\abs{z}< \sqrt{\mu(K)/\pi}\}\setminus K} \frac{d\mu}{\abs{\zeta -z}}
\end{equation}
Hence
\begin{multline}
\int_{K\setminus \{\abs{z}< \sqrt{\mu(K)/\pi}\}} \frac{d\mu}{\abs{\zeta -z}} \leq \\
\left( \int_{K\cap \{\abs{z}< \sqrt{\mu(K)/\pi}\}} +\int_{K\setminus \{\abs{z}< \sqrt{\mu(K)/\pi}\}}\right)\frac{d\mu}{\abs{\zeta -z}} \leq \\
\left( \int_{K\cap \{\abs{z}< \sqrt{\mu(K)/\pi}\}} +\int_{\{\abs{z}< \sqrt{\mu(K)/\pi}\}\setminus K}\right)\frac{d\mu}{\abs{\zeta -z}}
=\int_{\{\abs{z}< \sqrt{\mu(K)/\pi}\}} \frac{d\mu}{\abs{\zeta -z}}
\end{multline}
Hence
\begin{equation}
\int_{K} \frac{d\mu}{\abs{\zeta -z}} \leq \int_{\{\abs{z}< \sqrt{\mu(K)/\pi}\}}\frac{d\mu}{\abs{\zeta -z}}
\end{equation}
This completes the proof.
\end{proof}
In particular, $1/z$ is summable on compacts.
\index{Vitushkin's localization operator}
Define for $f\in C^0(\C)$ and $\psi$ be continuously differentiable with compact support in $\C$
\begin{equation}
T_\psi (f)(z):=\frac{1}{\pi}\int\frac{f(\zeta)-f(z)}{\zeta-z}\partial_{\bar{z}} \psi(\zeta)d\mu(\zeta), z\in \C
\end{equation}
The following is known (see Gamelin \cite{gamelinbok}, p.210).
\begin{proposition}\label{vitushkinopprops} 
$T_\psi (f)$ is continuous on $\C$ and $T_\psi (f)(\infty)=0$; $\partial_{\bar{z}} T_\psi (f)=
\partial_{\bar{z}} f\cdot \psi$ in the weak sense; If diam(supp$\psi$)$\leq \delta$ then
\begin{equation}
\norm{T_\psi (f)}_\infty \leq 2\delta \norm{\partial_{\bar{z}} \psi }_\infty \omega(f,\delta)
\end{equation}
where  $\omega (f,\cdot)$ denotes the {\em modulus of continuity}\index{Modulus of continuity} of $f$,
given by $\omega(f,\epsilon):=\sup_{\abs{z-\zeta}\leq \epsilon} \abs{f(z)-f(\zeta)}.$
\end{proposition}
\begin{proof}
Clearly
\begin{multline}
\abs{\frac{1}{\pi}\int \frac{f(\zeta)-f(z)}{\zeta -z} \partial_{\bar{z}} \psi d\mu(\zeta)}\leq
\frac{1}{\pi}\sup_{z,\zeta\in \mbox{supp}\psi}\abs{f(\zeta)-f(z)} \norm{\partial_{\bar{z}} \psi}_\infty
 \int \frac{1}{\abs{\zeta -z}}  d\mu(\zeta)
\end{multline}
and we also have for $z_0$ such that $\abs{T_\psi (f)(z_0)}=\norm{T_\psi (f)}_\infty,$
\begin{equation}
\int_{\mbox{supp}\psi} \frac{1}{\abs{\zeta -z}}  d\mu(\zeta) \leq 
\int_{\{\abs{\zeta}<\delta\}} \frac{d\mu(\zeta)}{\abs{\zeta}}  =2\pi\delta
\end{equation}
so the stated inequality follows.
Note that $T_\psi (f)(z)$ can, by applying the Green formula to $\phi$
 be written
\begin{equation}\label{hjpipi}
\frac{1}{\pi}\int \frac{f(\zeta)-f(z)}{\zeta -z} \partial_{\bar{z}} \psi d\mu(\zeta)=
\psi (z)f(z) + \frac{1}{\pi}\int \frac{f(\zeta)}{\zeta -z} \partial_{\bar{z}} \psi d\mu(\zeta)
\end{equation}
where the last term is the convolution of $1/z$ with a bounded function with compact support and
continuous with respect to $\zeta.$ Hence $T_\psi (f)$ is measurable and continuous
as soon as $f$ is continuous. Since $\psi$ has compact support $f \psi$ vanishes at infinity and
hence $T_\psi (f)$ vanishes at $\infty.$
Applying $\partial_{\bar{z}}$ to the equality  completes the proof.
\end{proof}
Note that right hand side of Eqn.(\ref{hjpipi}) also implies that 
\\
(i) $f-T_\psi (f)(z)$ is analytic on the interior 
of $\psi^{-1}(1)$.
\\
(ii) $T_\psi (f)(z)$ is analytic on the complement of $X.$
\begin{proposition}
$T_\psi (f)(z)$ is analytic whenever $f$ is analytic.
\end{proposition}
\begin{proof}
If $f$ is analytic on a disc $D$ then $(f(z)-f(\zeta))(z-\zeta)$ is jointly analytic in $z$ and $\zeta$ on $D\times D,$
so the integral
\begin{equation}
\int_D \frac{f(z)-f(\zeta)}{z-\zeta}\partial_{\bar{z}} \psi d\mu(\zeta)
\end{equation}
is analytic on $D,$ thus $T_\psi f$ is analytic on $D.$ This completes the proof.
\end{proof}

\begin{proposition}\label{partitionofunutylemmprop}
Let $\delta_j,$ $j\in \Z_+,$ be a sequence going to $0$. For each $j$ there exists a 
$C^\infty$-smooth partition of unity $\{\phi_{k,j}\}_{k\in \Z_+},$ $\phi_{k,j}\in C_c^\infty(\C)$ satisfying
together with an associated sequence of points $z_{k,j}$ such that
\begin{equation}
0\leq \phi_{k,j}\leq 1, \sum_{k}\phi_{k,j}=1 
\end{equation}
\begin{equation}
\phi_{k,j}=0\mbox{ off }\{\abs{z-z_{k,j}}\leq \delta_j\}
\end{equation}
\begin{equation}
\norm{\partial_{\bar{z}}\phi_{k,j}}_\infty \leq \frac{20}{\delta_j}
\end{equation}
and such that
no point $z$ is contained in more than $25$ of the discs $\{\abs{z-z_{k,j}}\leq \delta_j\}$
\end{proposition}
\begin{proof}
Let $\phi\in C_c^\infty(\C)$ such that $0\leq \phi(z) \leq 1,$
$\phi(z)=0$ for $\abs{z}\geq 1,$ $\int \phi=1$ and 
$\norm{\partial_{\bar{z}}\phi}_\infty \leq 10.$ 
Let $Q_k,$ $k\in \Z_+,$ denote squares with mutually disjoint interiors such 
that $\C=\cup_k Q_k.$
Define
\begin{equation}
\phi_k(z):=\int_{Q_k} \phi(\zeta-z) d\mu(\zeta)
\end{equation}
Then each $\phi_k$ is nonnegative and
\begin{equation}
\sum_{k=1}^\infty \phi_k(z)=\sum_{k=1}^\infty \int_{Q_k} \phi(\zeta-z) d\mu(\zeta)=\int \phi(\zeta -z)d\mu(\zeta)=1
\end{equation} 
If $\mbox{dist}(z,Q_k)\geq 1$ then $\phi_k(z)=0,$ so each $\phi_k$ is supported in a disc $\{ \abs{z-z_k}<2\},$
where $z_k$ denotes the centroid of $Q_k.$ Hence no $z$ belongs to more than $25$ of the
$\{ \abs{z-z_k}<2\}.$
For fixed $j$, the function $\phi_{k,j}(z):=\phi_k(2z/\delta_j)$ is supported in 
a disc $\{ \abs{z-z_{k,j}}\leq \delta_j\}$
and satisfies
$\norm{\partial_{\bar{z}}\phi_{k,j}}_\infty \leq 10\frac{2}{\delta_j}.$
Also $0\leq \phi_{k}\leq 1$, $\sum_{k}\phi_{k}=1$ implies that 
$0\leq \phi_{k,j}\leq 1$, $\sum_{k}\phi_{k,j}=1$. This completes the proof.
\end{proof}
We shall now give a proof of Theorem \ref{vituskin1} (it utilizes a scheme
for obtaining rational approximation that is also described in Gamelin \& Garnett \cite{gamelingarnett}). The proof follows Zalcman \cite{zalcmanbok}, see Section \ref{vitushkinapp} for technical lemmas, definitions and
auxiliary results used in the proof.
\begin{proof}(of Theorem \ref{vituskin1})
\\
$(i)\Rightarrow (ii)$:
Let $\Omega\subset\C$ be a bounded open subset, 
let $\epsilon>0$, $K\subset (\C\setminus \mbox{int}X)\cap \Omega$ a closed subset and $f\in C(K,1)$
such that 
\begin{equation}f'(\infty)>\alpha ((\C\setminus \mbox{int}X)\cap \Omega)-\epsilon \end{equation}
Since $f\in A(X)$, $(i)$ implies $f\in R(X).$ By Corollary \ref{bishopcorr}
$f\in R(X\cup (\C\setminus \Omega)).$ Let $f_n$ be rational such that $f_n\in C((\C\setminus X)\cap \Omega,1)$, $f_n\to f$ uniformly
on $X\cup(\C\setminus \Omega).$ Then $f'_n(\infty)\to f'(\infty).$ Since
$f'_n(\infty)\leq \alpha(\Omega\cap\C\setminus X)$ we have
\begin{equation}
\alpha(\Omega\cap (\C\setminus \mbox{int}X))-\epsilon \leq f'(\infty)\leq \alpha(\Omega\cap (\C\setminus X))
\end{equation}
Letting $\epsilon\to 0$ gives $\alpha(\Omega\cap\C\setminus \mbox{int}X)\leq \alpha(\Omega\cap\C\setminus X)$
and obviously $\alpha(\Omega\cap\C\setminus \mbox{int}X)\geq \alpha(\Omega\cap\C\setminus X)$, so we actually have equality.
In particular, this implies that for every $D_{z,r}:=\{\abs{\zeta-z}<r\}$ (denoting also $\overline{D}_{z,r}:=\{\abs{\zeta-z}\leq r\}$)  
\begin{equation}%
\alpha(\overline{D}_{z,r}\cap (\C\setminus \mbox{int}X))=\alpha(\overline{D}_{z,r}\cap (\C\setminus X)).
\end{equation}
$(ii)\Rightarrow (i)$:
Fix $r>1$. For $\delta>0$
\begin{equation}
\frac{\alpha(\overline{D}_{z,\delta}\cap (\C\setminus \mbox{int}X))}{\alpha(\overline{D}_{z,r\delta}\cap (\C\setminus X))} \leq
\frac{\alpha(D_{z,\frac{\delta(1+r)}{2}}\cap (\C\setminus \mbox{int}X))}
{\alpha(D_{z,\frac{\delta(1+r)}{2}}\cap (\C\setminus X))} =1
\end{equation}
Letting $\delta\to 0$ we obtain that for each $r>1$
\begin{equation}\label{fyrrrannn}
\lim_{\delta\to 0} \frac{\alpha(\overline{D}_{z,\delta}\cap\C\setminus \mbox{int}X)}{\alpha(\overline{D}_{z,r\delta}\cap\C\setminus X)} \leq
1
\end{equation}
Next suppose that
there exists $r\geq 1,$ $m\geq 1$ such that for each $z$ and $\delta$
\begin{equation}\label{iiekl}
\alpha(\overline{D}_{z,\delta}\cap\C\setminus \mbox{int}X) \leq m\alpha (\overline{D}_{z,\delta r}\cap\C\setminus X)
\end{equation}
Let $f\in A(X)$ and extend $f$ as a function with compact support in $\C.$
Let $\phi_{k,n}$ be the partition of unity constructed in Proposition \ref{partitionofunutylemmprop}.
For simplicity omit the subindex $n$, i.e.\
set $\delta:=\delta_n,\phi_k:=\phi_{k,n},$ $z_k:=z_{k,n}.$
Then
\begin{equation}
f_k(z)=f(z)\phi_k(z)+\frac{1}{\pi}\int f(\zeta)\partial_{\bar{\zeta}} \phi_k \frac{1}{\zeta -z} d\mu(\zeta)
\end{equation}
is holomorphic on $(\C\setminus \overline{D}_{z_k,\delta}) \cup \mbox{int}X$
and $\sum_k f_k \equiv f.$
Set $X_k:=(\C\setminus X)\cap \overline{D}_{z_k,r\delta}$, $O_k:=O(X_k),$
\begin{equation}
f_k(z)=\sum_{s=1}^\infty \frac{a_{sk}}{(z-O_k)^s}
\end{equation}
By Proposition \ref{partitionofunutylemmprop} and the properties of operator $T$ in
Proposition \ref{vitushkinopprops} applied to $f_k$,
$\norm{f_k}_\infty \leq$
$2\cdot 2\delta\omega(f,2\delta)20/\delta$
$=80\omega(f,2\delta)$ which implies by Eqn.(\ref{iiekl})
\begin{multline}
\abs{a_{lk}}=\abs{f'(\infty)}\leq \\
80\omega(f,2\delta) \alpha(\overline{D}_{z,\delta}\cap\C\setminus \mbox{int}X)
\leq 80m \omega(f,2\delta)\alpha(X_k)
\end{multline}
By the coefficient estimate given in Proposition \ref{coeffestimateprop}
\begin{equation}
\abs{a_{2k}}\leq m_1\omega(f,2\delta)\alpha(X_k)\beta(X_k)
\end{equation}
Set $m_2=20\max(80m,m_1).$ By Lemma \ref{extraskrapppp}
there exists $g_k\in C^0(X_k,m_2\omega(f,2\delta))$ such that $g'_k(\infty)=f'_k(\infty)$
and $\beta(X_k,O_k,g_k)=\beta(X_k,O_k,f_k).$ $\norm{f_k}_\infty \leq m_2\omega(f,2\delta)$ implies
$\norm{f_k -g_k}_\infty \leq 2m_2\omega(f,2\delta)$ which in turn gives
\begin{equation}
\frac{\abs{z-z_k}^3}{r^3\delta^3}\abs{f_k(z)-g_k(z)}\leq 2m_2\omega(f,2\delta),\quad \abs{z-z_k}=r\delta
\end{equation}
Set $m_3:=2m_2r^3.$ By the maximum principle 
\begin{equation}
\abs{f_k(z)-g_k(z)}\leq 2m_2 r^3\delta^3\omega(f,2\delta) \min\{ 1,\delta^3/\abs{z-z_k}^3\}
\end{equation}
For a fixed $w$ every $\overline{D}_{z_k,\delta}$ intersects at least one and at most two of the circles
$C_s:=\{\abs{z-w}=s\delta\},$ $s\in \Z_+.$ Denote by $N_s$ the number of discs that intersect $C_s.$ By Proposition \ref{partitionofunutylemmprop}, $N_s\pi \delta^2 \leq 25$
$\mbox{area}(\{(s-1)\delta \leq \abs{z-w}\leq (s+1)\delta\}),$ which implies $N_s\leq 100s.$
Set $m_4=100m_3(1+\sum_{s=2}^\infty s(s-1)^{-3}).$
For $s\geq 2$ and $\overline{D}_{z_k,\delta} \cap C_s \neq \emptyset,$ 
$\abs{w-z_k}\geq (s-1)\delta,$ gives
$\abs{f_k(w)-g_k(w)}\leq m_3\omega(f,2\delta)/(s-1)^3.$ Hence
\begin{multline}
\sum_{k=1}^\infty \abs{f_k(w)-g_k(w)} \leq m_3\omega(f,2\delta)N_1 +\\
m_3\omega(f,2\delta)\sum_{s=2}^\infty \frac{N_s}{(s-1)^{3}} \leq m_4 \omega(f,2\delta)
\end{multline}
Hence
\begin{multline}
\norm{f-\sum_{k=1}^\infty g_k} =\norm{\sum_{k=1}^\infty f_k-\sum_{k=1}^\infty g_k}=\\
\norm{\sum_{k=1}^\infty (f_k-g_k)}\leq \sup_{z} \sum_{k=1}^\infty \abs{f_k(z)-g_k(z) }\leq m_4\omega(f,2\delta)
\end{multline}
Now there are only finitely many nonzero $f_k,$ say $f_{k_j}$, $j=1,\ldots ,N.$
We have $f=\sum_{j=1}^N f_{k_j}$ and
$\norm{\sum_{j} (f_{k_j} -g_{k_j})}_X\leq m_4\omega(f,2\delta).$ Also each $g_k$ is 
holomorphic on a neighborhood of $(\C\setminus X_k)=X\cup (\C\setminus \overline{D}_{z_k,r\delta})$
hence $g=\sum_j g_{k_j}$ is holomorphic on a neighborhood of $X.$
By Proposition \ref{rungesobs} $g\in R(X).$
Since $m_4$ depends only on $m_1$ and $\delta=\delta_n\to 0$
as $n\to\infty$, this proves (i) given Eqn.(\ref{iiekl}).
\\
\\
It remains to show that there exists $r\geq 1,$ $m\geq 1$ such that for each $z$ and $\delta$
\begin{equation}
\alpha(\overline{D}_{z,\delta}\cap(\C\setminus \mbox{int}X) )\leq m\alpha (\overline{D}_{z,\delta r}\cap (\C\setminus X))
\end{equation}
Assume (in order to reach a contradiction) that Eqn.(\ref{iiekl}) does not hold.
Pick $z_1,$ $\delta_1$ such that
$\alpha(\overline{D}_{z_1,\delta_1}\cap\C\setminus \mbox{int}X) \leq$
$\alpha (\overline{D}_{z_1,3\delta_1}\cap\C\setminus X)$. Since we have proved $(i)\Rightarrow (ii)$ we obtain
$A(X\cap \overline{D}_{z_1,3\delta_1})\neq R(X\cap \overline{D}_{z_1,2\delta_1}).$
Set $X_1=X\cap \overline{D}_{z_1,2\delta_1}.$ Since we have proved that Eqn.(\ref{iiekl}) is sufficient
for (i) we know that there exists
$z_2,\delta_2$ such that 
$\alpha(\overline{D}_{z_2,\delta_2}\cap (\C\setminus \mbox{int}X_1)) >$
$2\alpha (\overline{D}_{z_2,5\delta_2}\cap (\C\setminus X_1))$. If
$(\C\setminus X_1 \cap \overline{D}_{z_2,5\delta_2}$ contained a disc of diameter $\geq \delta_2$ that would imply
$\alpha(\overline{D}_{z_2,\delta_2}\cap (\C\setminus \mbox{int}X_1)) >$
$2\alpha (\overline{D}_{z_2,5\delta_2}\cap (\C\setminus X_1))$, so we can assume
$\overline{D}_{z_2,5\delta_2}\subset \overline{D}_{z_1,2\delta_1+\delta_2}.$
Also if there existed $x\in \in \overline{D}_{z_2,2\delta_2}\setminus \overline{D}_{z_1,2\delta_1}$
and $x+2\exp(i\theta)\delta_2 \in \overline{D}_{z_2,5\delta_2}$ for all real $\theta.$ Since there exists $\theta$ such that
$x+2\exp(i\theta)\delta_2 \notin \overline{D}_{z_1,2\delta_1+\delta_2}$. So we can assume
$\overline{D}_{z_2,2\delta_2} \subset \subset D_{z_1,2\delta_1}$.
Thus $\overline{D}_{z_2,\delta_2}\cap \mbox{int} X_1 = \overline{D}_{z_2,\delta_2}\cap \mbox{int} X$
which implies
$\overline{D}_{z_2,\delta_2}\cap (\C\setminus \mbox{int} X_1) = \overline{D}_{z_2,\delta_2}\cap (\C\setminus \mbox{int}X),$
thus
\begin{equation}
\alpha(\overline{D}_{z_2,\delta_2}\cap (\C\setminus \mbox{int} X_1))>2\alpha(
\overline{D}_{z_2,5\delta_2}\cap (\C\setminus \mbox{int}X))
\end{equation}
Since $5\delta_1 \leq 2\delta_1 +\delta_2$ we have $\delta_2 \leq \delta_1/2.$
Repeating this procedure renders sequences $\{z_n\}_{n\in \Z_+},$ $\{r_n\}_{n\in \Z_+},$
$\{\delta_n\}_{n\in \Z_+},$ such that $r_n \to \infty,$ $\delta_n\to 0$ and
$\overline{D}_{z_{n+1},2\delta_{n+1}}\subset \overline{D}_{z_n,2\delta_n},$
$\delta_{n+1}\leq \delta_n/2$, 
$\alpha(\overline{D}_{z_2,\delta_2}\cap(\C\setminus \mbox{int}X)) >$
$n\alpha (\overline{D}_{z_2,5\delta_2}\cap(\C\setminus X))$.
Set $z_0:= \cap \overline{D}_{z_n,2\delta_n}.$ Then $\abs{z_n -z_0}<2\delta_n$ implies
$\alpha(\overline{D}_{z_0,3\delta_n}\cap (\C\setminus \mbox{int} X))\geq$
$\alpha(
\overline{D}_{z_n,\delta_n}\cap (\C\setminus \mbox{int}X))$ and
$\alpha(\overline{D}_{z_0,(r_n-2)\delta_n}\cap (\C\setminus X))\geq$
$\alpha(\overline{D}_{z_n,r_n\delta_n}\cap (\C\setminus X)) \geq n$. Hence
\begin{equation}
\frac{\alpha(\overline{D}_{z_0,3\delta_n}\cap(\C\setminus \mbox{int}X))}{\alpha(\overline{D}_{z_0,(r_n -2)\delta_n}\cap (\C\setminus X))} 
\geq
\frac{\alpha(\overline{D}_{z_n,\delta_n}\cap (\C\setminus \mbox{int}X))}{\alpha(D_{z_n,r_n\delta_n}\cap (\C\setminus X))}\geq n
\end{equation}
For $n\to \infty$ we have $r_n\to \infty$ thus for $r\geq 1,$ $\epsilon_n:=3\delta_n$ and sufficiently large
$n$ we have
\begin{equation}
\frac{\alpha(\overline{D}_{z_0,\epsilon_n}\cap(\C\setminus \mbox{int}X))}
{\alpha(\overline{D}_{z_0,(r\epsilon_n}\cap (\C\setminus X))} 
\geq
\frac{\alpha(\overline{D}_{z_0,\epsilon_n}\cap (\C\setminus \mbox{int}X))}{\alpha(\overline{D}_{z_0,\frac{(r_n-2)\epsilon_n}{3}}
\cap 
(\C\setminus X))}\geq n
\end{equation}
contradicting Eqn.(\ref{fyrrrannn}).
This proves Eqn.(\ref{iiekl}). 
This completes the proof of Theorem \ref{vituskin1}).
\end{proof}

\section{Some modern results}
Let us now look at some results in the polyanalytic setting that are inspired by Vitushkin's methods.

\begin{theorem}[Carmona \cite{carmona1985}]
Let $q\in \Z_+$ and let $X\subset\C$ be a compact subset such that $\C\setminus X$ has finitely many connected components.
Let $g$ be a three times continuously differentiable complex function in a neighborhood of $X$ which 
satisfies $\partial_{\bar{z}} g(z)\neq 0$ for all $z\in X.$ Then any $f\in C^0(X),$ 
satisfying $\partial_{\bar{z}}\left(\frac{\partial_{\bar{z}}f}{\partial_{\bar{z}}g}\right)=0$
on $\mbox{int}(X)$ in the weak sense 
can be uniformly approximated on $X$ by functions $\psi$ 
that satisfy $\partial_{\bar{z}}\left(\frac{\partial_{\bar{z}}f\psi}{\partial_{\bar{z}}g}\right)=0$ on an open neighborhood of $X.$
\end{theorem}
\begin{proof}
We begin with the following lemma.
\begin{lemma}\label{carmonalemmatva} 
For any $p_0\in \mbox{int}(X),$
\begin{equation}
\abs{\frac{\partial_{\bar{z}}f\psi}{\partial_{\bar{z}}g} (p_0)}\leq c
\frac{\omega(f,\mbox{dist}(p_0,\C\setminus X))}{\mbox{dist}(p_0,\C\setminus X)}
\end{equation}
\end{lemma}
\begin{proof}
A function that satisfies $\partial_{\bar{z}}\left(\frac{\partial_{\bar{z}}f\psi}{\partial_{\bar{z}}g}\right)=0$
can be decompose as $f=h+gk,$ for $h,k\in \mathscr{O}(\mbox{int}(X)),$ by defining 
$k=\partial_{\bar{z}} f/\partial_{\bar{z}} g.$
Let $p_0\in \mbox{int}(X).$
Let $r>0$ such that $2r<\mbox{dist}(p_0,\C\setminus K).$ Let $\rho_1\in \C^\infty_c([-4,4],\R)$ (here
the subindex $c$ denotes compact support in $[-4,4]$) such that
$0\leq \rho_1 \leq 1$ and $\rho_1=1$ on $[-1,1].$ Set $\rho:=\rho_1(\abs{z-p_0}^2/r^2)$ and set
$\psi:=\rho/\partial_{\bar{z}}\in C^2_c(\{ \abs{z-p_0} <2r\}).$ This implies
\begin{equation}
\norm{T_\psi(f)}_\infty \leq c\omega(f,2r)r\norm{\partial_{\bar{z}} \psi}_\infty
\end{equation}
Since
\begin{equation}
\norm{\partial_{\bar{z}} \psi}_\infty =\norm{\partial_{\bar{z}}\left(\frac{1}{\partial_{\bar{z}} g}\right)\rho
+\partial_{\bar{z}} \rho \frac{1}{\partial_{\bar{z}} g}}_\infty \leq \frac{c}{r}
\end{equation}
we have
\begin{equation}\label{ggg0}
\abs{\int_{\abs{z-p_0}=r} T_\psi(f)(z)dz}\leq cr\omega(f,2r)
\end{equation}
Using the property that if diam(supp$\psi$)$\leq \delta$ then
\begin{equation}
\norm{T_\psi (f)}_\infty \leq 2\delta \norm{\partial_{\bar{z}} \psi }_\infty \omega(f,\delta)
\end{equation}
Stokes theorem yields
\begin{multline}
\int_{\abs{z-p_0}=r} T_\psi(f)(z)dz =\int_{\abs{z-p_0}<r} \partial_{\bar{z}}f\psi d\bar{z}\wedge dz =\\
\int_{\abs{z-p_0}<r} k d\bar{z}\wedge dz =cr^2 k(p_0)
\end{multline}
Eqn.(\ref{ggg0}) then implies
\begin{equation}
\abs{k(p_0)}\leq c\frac{\omega(f,2r)}{r}
\end{equation}
Letting $r\to \frac{1}{2} \mbox{dist}(p_0,\C\setminus X)$ gives the wanted result.
This proves Lemma \ref{carmonalemmatva}.
\end{proof}
Now we extend $f$ and $g$ (keeping the same notation) to $\C$ such that $\mbox{supp}f\subset U,$ $U:=\{\partial_{\bar{z}} g\neq 0\}.$
Define in the standard way an approximation of unity by letting $\rho\in C^\infty_c(\{ \abs{z}<1\}),$
$0\leq \rho\leq 1,$ $\int\rho =1$ and set $\rho_\epsilon (x)=(1/\epsilon^2)\rho(x/\epsilon).$
Then $\int\rho_\epsilon =1,$ $\int \partial_{\bar{z}}^j\rho_\epsilon =0,$
$\int \partial_{\bar{z}}^j \abs{\rho_\epsilon} \leq c/\epsilon^2$, $j=1,2$ and $\mbox{supp} \rho_\epsilon\subset \{\abs{z}\leq \epsilon\}.$
Set $f_\epsilon(x):=f*\rho_\epsilon(x)=\int f(x-t)\rho_\epsilon(t)d\mu(t),$ $f_\epsilon\in C^\infty(\C)$ and
$\mbox{supp}f_\epsilon \subset U,$ for sufficiently small $\epsilon$.
Here $c$ is the bound from Lemma \ref{carmonalemmatva}.
Let $m$ be a measure on $X$ which is orthogonal to $R(X,g).$ Define for $z\in \C$
\begin{equation}
\Check{m}(\zeta):=\int K(z,\zeta) dm(z),\quad  K(z,\zeta):=\int\frac{g(z)-g(\zeta)}{z-\zeta}
\end{equation} 
\begin{lemma}\label{carmonlem001}
$\Check{m}(z)$ is well-defined on $\C$ and continuous except possibly on the countable set $\{\abs{m}(\{\zeta\})>0\}$
and is $\mu$-locally integrable.
\end{lemma}
\begin{proof}
By the mean value theorem $K(\cdot,\zeta)$ is bounded on $\mbox{supp}\mu$. Let $\zeta$ be such that 
$\abs{\mu}(\{\zeta\})=0$ and pick $\zeta_j$ such that $\zeta_j\to \zeta.$ Then $K(z,\zeta_j)\to K(z,\zeta)$
except possibly for $z=\zeta,$ so we have $\mu$-a.e.\ convergence. The mean value theorem implies that
$K(\cdot,\zeta_j)$ is uniformly dominated on $\mbox{supp}\mu$ thus
$\Check{\mu}(\zeta_j)\to \mu(\zeta)$, hence $\Check{\mu}$ is continuous at $\zeta.$ 
Also
\begin{equation}\label{assjack}
\abs{\Check{\mu}(\zeta)}\leq (\norm{g}_{\mbox{supp}\mu} +\abs{g(\zeta)})\int\frac{d\abs{\mu}(z)}{\abs{z-\zeta}}
\end{equation}
Now define (see Browder \cite{browder}, p.153)
\begin{equation}
m'(\zeta):=\int\frac{d\abs{\mu}(z)}{\abs{z-\zeta}}
\end{equation}
By Lemma \ref{browderlemma0} we have for any compact $X\subset \C,$ and $a>0,$
\begin{equation}
\int_X d\abs{m}(z)\int_{\{\abs{z}<a\}} \frac{d\mu(\zeta)}{\abs{z-\zeta}}\leq 2\pi a\abs{m}(X)
\end{equation}
where $1/\abs{z_1-z_2}$ is measurable with respect to the product measure $\abs{m}\times \mu$ on
$X\times \{\abs{z}<a\}.$ So by Fubinis theorem 
$m'$ is $\mu$-summable and $m' <\infty$ a.e.\ $\mu.$ Hence by Eqn.(\ref{assjack})
$\Check{m}$ is locally integrable.
This completes the proof of Lemma \ref{carmonlem001}.
\end{proof}
With the bound $c$ from Lemma \ref{carmonalemmatva} we have $\abs{\Check{m}}\leq c$.
Hence
\begin{equation}
\int f_\epsilon d\mu =\frac{1}{\pi} \int_\C \partial_{\bar{z}}\left( \frac{f_\epsilon}{\partial_{\bar{z}}g}\right)
\Check{m} d\mu =
\frac{1}{\pi} \int_X \partial_{\bar{z}}\left( \frac{f_\epsilon}{\partial_{\bar{z}}g}\right)
\Check{m} d\mu
\end{equation}
Set 
\begin{equation}
X_1:=\{z\in X :\mbox{dist}(z,\C\setminus X)\geq 2\epsilon\}, X_2 :=\{z\in X :\mbox{dist}(z,\C\setminus X)< 2\epsilon\},
\end{equation}
and
\begin{equation}
A_1:= \int_{X_1} \partial_{\bar{z}}\left( \frac{f_\epsilon}{\partial_{\bar{z}}g}\right)\Check{m}d\mu,\quad 
A_2:= \int_{X_2} \partial_{\bar{z}}\left( \frac{f_\epsilon}{\partial_{\bar{z}}g}\right)\Check{m}d\mu
\end{equation}
Then 
\begin{equation}
\pi \int f_\epsilon d\mu =A_1+A_2
\end{equation}
For $z\in X_1$ we have for $\abs{\zeta}<\epsilon,$ $z-\zeta\in \mbox{int}(X)$ and $\mbox{dist}(z-\zeta,\C\setminus X)\geq 0$ and
\begin{equation}
\partial_{\bar{z}}\left( \frac{f_\epsilon}{\partial_{\bar{z}}g}\right)=
\int \partial_{\bar{z}}\left( \frac{f_\epsilon (z-\zeta)}{\partial_{\bar{z}}g(z)}\right) \rho_\epsilon (\zeta)d\mu(\zeta)
\end{equation}
Inserting $f=h+gk$ yields
\begin{equation}
\partial_{\bar{z}}\left( \frac{f_\epsilon}{\partial_{\bar{z}}g}\right)=
\int \partial_{\bar{z}}\left( \frac{\partial_{\bar{z}} g(z-\zeta)}{\partial_{\bar{z}}g(z)}\right)k(z-\zeta)\rho_\epsilon (\zeta)d\mu(\zeta)
\end{equation}
Since $g$ is differentiable, and as we have seen, $\int\abs{\partial_{\bar{z}} \rho_\epsilon d\mu} \leq c/\epsilon^2$, 
application of Lemma \ref{carmonalemmatva}
gives with the bound $c$ from Lemma \ref{carmonalemmatva}
\begin{multline}
\abs{\partial_{\bar{z}}\left( \frac{f_\epsilon}{\partial_{\bar{z}}g}\right)}\leq \\
\frac{c}{\epsilon}\sup_{\substack{x\in X}{\abs{\zeta}<\epsilon}} \left\{ \abs{ \partial_{\bar{z}}\left( 
\frac{\partial_{\bar{z}} g(z-\zeta)}{\partial_{\bar{z}}g(z)}\right)}\omega(f,\mbox{dist}(z-\zeta,\C\setminus X))\right\}\leq
\epsilon c(z)
\end{multline}
Then by the dominated convergence theorem  we obtain $\lim_{\epsilon\to 0} A_1=0.$
Next let $p_0\in X_2.$ $\int \partial_{\bar{z}}^j \abs{\rho_\epsilon} \leq c/\epsilon^2,$ $j=1,2$,
implies that
\begin{equation}
\partial_{\bar{z}}^i f_\epsilon (z)=\int (f(z-\zeta)-f(z))\partial_{\bar{z}}^i \rho_\epsilon (\zeta) d\mu(\zeta),\quad i=1,2
\end{equation}
thus for $j=1,2$
\begin{equation}
\abs{\partial_{\bar{z}}\left( \frac{f_\epsilon}{\partial_{\bar{z}}g}\right)(z)}\leq \abs{\partial_{\bar{z}}^i f_\epsilon (z)}\leq 
c\frac{\omega(f,\epsilon)}{\epsilon^2}
\end{equation}
Hence
\begin{equation}\label{skunda5}
\abs{A_2}\leq c\frac{\omega(f,\epsilon)}{\epsilon^2} \int_{X_1} \abs{\Check{m}} d\mu
\end{equation}
Since $\C\setminus X$ has finitely many component
we can fo sufficiently small $\epsilon,$ cover $X_2$ by discs $D_1,\ldots,D_n$ of radius $4\epsilon$ with centers off $X.$
For each $D_j$ let $E_j\subset D_j\setminus X$ with diam$(E_j)\geq \epsilon.$
\begin{lemma}\label{camonasistalem}
There exists functions $Q_j(u,\cdot)\in R(X,g)$ such that for $z\in \C\setminus E_j,$ and $u\in D_j$
\begin{equation}\label{skunda1}
\abs{Q_j(u,z)}\leq c
\end{equation}
\begin{equation}\label{skunda2}
\abs{Q_j(u,z) -\frac{g(z)-g(u)}{z-u}} \leq c\frac{\epsilon^3}{\abs{z-u}^3}
\end{equation}
\end{lemma}  
\begin{proof}
Assume $D_j$ has center $z_0$ and radius $r>0.$ Let $f_1$ be a conformal representation of the connected set
$(\C\setminus E_j)\cup\{\infty\}.$ in $\{ \abs{z}<1\}$ such that $f_1(\infty)=0.$
Denote by $\gamma(E)$ the Ahlfors capacity. 
The compact connected set $E_j$ satisfies 
\begin{equation}
\gamma(E)\leq \mbox{diam}(E_j)\leq 4\gamma(E_j).
\end{equation}
To see this note that for a disc $\{\abs{z-p_0}<R\},$ the function $G(z)=R/(z-p_0)$ maps 
$(\C\setminus \{\abs{z-p_0}<R\})\cup\{\infty\}$ conformally to the unit disc and $f'(\infty)=R.$
So $\gamma(\{\abs{z-p_0}<R\})=R.$ Hence $\gamma(E_j)\leq \mbox{diam}(E_j)$
(since $E_j$ is contained in a disc with radius $\mbox{diam}(E_j)$.  
On the other hand the function $G_1(z)=(z+1/z)(L/4)$ maps 
the unit disc conformally to $(\C\cup\{\infty\})\setminus [-L/2,L/2]$ and satisfies $G_1(0)=\infty.$
Hence $\gamma([-L/2,L/2])=\lim_{z\to \infty} \abs{zG_1^{-1}(z)}=\lim_{z\to 0} \abs{zG_1(z)}=L/4.$
Letting $U$ be the unbounded component of $(\C\setminus E_j)\cup\{\infty\}$ consider the conformal map 
$G_2:U\to B(0,1)$ with $G_2(\infty)=0.$ Let $z_1,z_2\in E_j$ such that
$\abs{z_1-z_2}=\mbox{diam}({E_j}).$ Then the function
\begin{equation}
G_3(z)=\frac{\gamma(E_j)}{G_2^{-1}(z)-z_1}
\end{equation}
is a univalent analytic function $\{\abs{z}<1\} \to \C$ with $G_2(0)=0,$ $G'(0)=1.$ 
Thus by the $\frac{1}{4}$-Koebe theorem we have $\{\abs{z}<1/4\}\subset G_2(\{\abs{z}<1\})$ which
shows the wanted bound.
\\
\\
Now set $a:=f'(\infty).$ Then $\abs{a}=\gamma(E)\geq 4/r.$
The function $F:=f_1/a,$ $\Omega\to \{\abs{z}<1/\abs{a}\}$ is holomorphic such that

\begin{equation}\label{carmonatooo}
F(\infty)=0,\quad F'(\infty)=1, \quad \norm{F}_\infty \leq 4/r
\end{equation}
Set
\begin{equation}
b=\frac{1}{2\pi i} \int_{\{\abs{z-p_0}=r\}} (z-z_0)f(z)dz,\quad b':=\frac{1}{2\pi i} \int_{\{\abs{z-p_0}=r\}} (z-p_0)^2 F(z)dz
\end{equation}
Expand in a neighborhood of $\infty$, for fixed $u$ and for $\abs{z-u}>2r$
\begin{equation}\label{caromaexpand1}
F(z)=\frac{1}{z-u} +\frac{\lambda_2 (u)}{(z-u)^2} +\frac{\lambda_3 (u)}{(z-u)^3}+\cdots
\end{equation}
\begin{equation}\label{caromaexpand2}
F^2(z)=\frac{1}{(z-u)^2} +\frac{2\lambda_2(u)}{(z-u)^3} +\cdots
\end{equation}
\begin{equation}\label{caromaexpand3}
F^3(z)=\frac{1}{(z-u)^3} +\frac{3\lambda_2(u)}{(z-u)^4} +\cdots
\end{equation}
Let $\Gamma$ be a positively oriented circle centered at $p_0$ with sufficiently large radius such that
\begin{equation}
\lambda_2(u)=\frac{1}{2\pi i} \int_\Gamma (z-u)F(z) dz,\quad \lambda_3(u)=\frac{1}{2\pi i} \int_{\Gamma}(z-u)^2 F(z) dz,
\end{equation}
Then $1=F'(\infty)=(1/2\pi i)\int_{\Gamma} F(z)dz$ implies
\begin{equation}
\lambda_2(u)=b-u-p_0,\quad \lambda_3(u)=b'-2(u-p_0)b+(u-p_0)^2
\end{equation}
Define 
\begin{equation}
Q_j(u,z):=(g(z)-g(u))H(u,z)
\end{equation}
\begin{equation}
H(u,z):=F(z)+(u-p_0-b)f^2(z)-(b'-2(u-p_0)b-2b^2-(u-p_0)^2)f^3(z)
\end{equation}
Since $H(u,\cdot)$ is holomorphic we have $Q(u,\cdot)$ belong to
the closure in $C^0(X)$ of the
the space of holomorphic functions on $\Omega$. By Eqn.(\ref{carmonatooo})
we have $\abs{b}\leq 4r$ and $\abs{b'}\leq 4r^2$ which implies for $z\in \Omega,$ $u\in D_j$
\begin{multline}\label{ffuuffuu}
\abs{H(u,v)}\leq \frac{4}{r} +(\abs{u-p_0} +\abs{b})\frac{16}{r^2} +\\
(\abs{b'} +2\abs{u-p_0}\abs{b} +2\abs{b}^2 +\abs{u-p_0}^2)\frac{64}{r^3}\leq \frac{c}{r}.
\end{multline}
So for fixed $u\in D_j,$ $H(u,z)(z-u)\in \mathscr{O}(\Omega)$ is bounded by $c$ (independently of $u$) as soon as $\abs{z-u}\leq 2r.$
By the maximum principle $\abs{H(u,z)}\abs{z-u}\leq c$ for $u\in D_j,$ $z\in \omega.$
By the mean value theorem we obtain $\abs{Q_j(u,z)}\leq c.$
Eqn.(\ref{caromaexpand1})-Eqn.(\ref{caromaexpand3}) imply that there exists a holomorphic function $h(u,\cdot)$ on $\{\abs{z-u}>2r\}\cup\{\infty\},$ such that
\begin{equation}
\abs{H(u,z)-\frac{1}{z-u}}=\frac{1}{\abs{z-u}^4} \abs{h(u,z)},\quad z\in \omega,u\in D_j
\end{equation}
Thus for $\abs{z-u}\leq 2r,$ $z\in \Omega,$ Eqn.(\ref{ffuuffuu})  implies that
\begin{equation}
\abs{h(u,z)}=\abs{H(u,z)-\frac{1}{z-u}}\abs{z-u}^4\leq cr^3
\end{equation}
By the maximum principle this holds for all $z\in \Omega.$ Hence for $z\neq u$, $z,u\in \Omega$ 
\begin{multline}
\abs{Q(u,v) -\frac{g(z)-g(u)}{z-u}}\leq c\abs{z-u}\abs{H(u,z)-\frac{1}{z-u}}  \leq c\frac{r^3}{\abs{z-u}^3}
\end{multline}
This completes the proof of Lemma \ref{camonasistalem}.
\end{proof}
Now set $L_1:=D_1\cap X_2$ and $L_j:=D_j\cap X_2\setminus (L_1\cup\cdots \cup L_{j-1}),$
$1<j\leq n.$ Then we have a disjoint union $X_2=\bigcup_j L_j.$ For each
$u\in X_2$ there exists a unique $j$ such that $L_j$ thus
\begin{equation}
\abs{\Check{m}(u)}\leq \int_X\abs{Q_j(u,z)-\frac{g(z)-g(u)}{z-u}}d\abs{m}(z)
\end{equation}
We have
\begin{equation}\label{skunda3}
\sum_j \int_{L_j} \abs{\Check{m}(u)} d\mu(u)\leq \int_X d\abs{\Check{m}}(z)  \sum_j \int_{L_j} 
\abs{Q_j(u,z)-\frac{g(z)-g(u)}{z-u}}d\mu(u)
\end{equation}
For fixed $z\in X,$
\begin{multline}
\sum_j \int_{L_j} \abs{Q_j(u,z)-\frac{g(z)-g(u)}{z-u}}d\mu(u) \leq \\
\sum_j \int_{L_j\cap\{\abs{u-z}<4\epsilon\}} \abs{Q_j(u,z)-\frac{g(z)-g(u)}{z-u}}d\mu(u) +\\
\sum_j \int_{L_j\cap\{\abs{u-z}\geq 4\epsilon\})} \abs{Q_j(u,z)-\frac{g(z)-g(u)}{z-u}}d\mu(u)
\end{multline}
Estimating in the right hand side, the first sum using the Eqn.(\ref{skunda2}), and
the second sum using the Eqn.(\ref{skunda2}), we obtain
\begin{multline}
\sum_j \int_{L_j\cap\{\abs{u-z}<4\epsilon\}} \abs{Q_j(u,z)-\frac{g(z)-g(u)}{z-u}}d\mu(u) +\\
\sum_j \int_{L_j\cap\{\abs{u-z}\geq 4\epsilon\})} \abs{Q_j(u,z)-\frac{g(z)-g(u)}{z-u}}d\mu(u) \leq\\
\sum_j \int_{L_j\cap\{\abs{u-z}<4\epsilon\}} cd\mu(u) +\sum_j \int_{L_j\cap\{\abs{u-z}\geq 4\epsilon\})} 
c\frac{\epsilon^3}{\abs{z-u}^3}d\mu(u) \leq\\
c\epsilon^2 +c\int_{4\epsilon}^\infty \int_0^{2\pi} \rho \frac{\epsilon^3}{\rho^3} d\rho d\theta =c\epsilon^2
\end{multline}
This together with Eqn.(\ref{skunda3}) and Eqn.(\ref{skunda5})
yields
\begin{equation}
\abs{A_2}\leq c\frac{\omega(f,\epsilon)}{\epsilon^2} \abs{m}(X)\epsilon^2 =O(\omega(f,\epsilon))
\end{equation}
Hence $\lim_{\epsilon\to 0} (A_1+A_2) =0.$
This completes the proof.
\end{proof}
In particular, since $\C\setminus X \neq\emptyset$ we can let $p_0\in \C\setminus X$
and $g(z)=\bar{z}-p_0.$ Then the condition $\partial_{\bar{z}}\left(\frac{\partial_{\bar{z}}f\psi}{\partial_{\bar{z}}g}\right)=0$
becomes $\partial_{\bar{z}}\left(\partial_{\bar{z}}f\right)=0,$ i.e.\ $f\in A_2(X),$
and in the conclusion the members of the approximating sequence satisfy 
$\partial_{\bar{z}}\left(\partial_{\bar{z}}\psi\right)=0$ on an open neighborhood of $X.$
If we further require that each derivative of $f$ of order $\leq q$
has continuous extension to $X$ (we shall denote this by $f\in C^q(X)$) then, since $g\in C^\infty(\C),$ the theorem can then be repeated for
$f$ replaced by $\partial_{\bar{z}}^j f$ for $j=1,\ldots,q-1$ in order to yield the following corollary.

\begin{corollary}\label{carmonacorextrabra}
Let $q\in \Z_+$ and let $X\subset\C$ be a compact subset such that $\C\setminus X$ is connected.
Then $A_q(X) \cap C^{q-1}(X)=P_q(X).$
\end{corollary}
\begin{proof}
Any $f(z)\in C^{q-1}(X)$ that takes the form $f(z)=\sum_{j=0}^{q-1} a_j(z)\bar{z}^j,$
for $a_j\in A_1(X)$ belongs to $A_q(X)\cap C^{q-1}(X)$. Conversely let
$f\in A_q(X)\cap C^{q-1}(X)$. On int$X$, $f$ takes the form $f(z)=\sum_{j=0}^{q-1} a_j(z)\bar{z}^j,$
for $a_j\in \mathscr{O}(\mbox{int}X).$ Since $f \in C^{q-1}(X)$ each $\partial_{\bar{z}}^j f,$ $j=0,\ldots,q-1,$ has continuous extension to $X.$
This implies that each $a_j$ has an extension $\tilde{a}\in A_1(X).$ 
Indeed, $\partial_{\bar{z}}^{q-1} f=a_{q-1}$ so $a_{q-1}$ has a continuous extension to $X$. Therefore
the process can be repeated for $f-\bar{z}^{q-1}a_{q-1}$ calculating $\partial_{\bar{z}}^{q-2} (f-\bar{z}^{q-1}a_{q-1})$ and so on.
Thus $f$ has an extension
$\sum_{j=0}^{q-1} \tilde{a}_j(z)\bar{z}^j\in A_q(X)\cap C^{q-1}.$ The representation is unique since the boundary values are determined by continuity.
So let $f\in A_q(X)\cap C^{q-1}(X)$. We have seen that this implies that $f$ can be written
as $\sum_{j=0}^{q-1} a_j(z)\bar{z}^j,$
for $a_j\in A_1(X)$. Since
$\C\setminus X$ is connected each $a_j$ can be uniformly approximated on $X$ by
polynomials $P_{j,k},$ in $z.$ Hence $f_k:=\sum_{j=0}^{q-1} P_{j,k}(z)\bar{z}^j$
is a sequence in $P_q(X)$ converging to $f.$ The inclusion $P_n(X)\subset A_q(X)\cap C^{q-1}(X)$ is trivial.
\end{proof}
See Section \ref{noproof} for further consequences.
Let $X\subset\C$ be a compact subset. For $w\in C$ and $r>0$ denote by $d(w,r,X)$ the least upper bound of the diameters of all connected components of 
$\{ \abs{z-w}<r\}\setminus X.$ Set
\begin{equation}
\Theta(X):=\inf \left\{ \frac{d(w,r,X)}{r}\colon z\in \partial X,r>0\right\}
\end{equation}
For a subset $Y\subset\C$ and a function $g:Y\to \C,$ 
set $\omega_Y(g,\delta):=\sup \{\abs{g(z)-g(w)}\colon z,w\in Y,\abs{z-w}\leq \delta\}.$
When $Y$ is $\C$ we drop the subindex.
For $j\in \Z_+$ and $\phi\in C^m_{\mbox{loc}}(\C)$ set
\begin{equation}
\nabla^j\phi:=(\partial_{z}^{0}\partial_{\bar{z}}^{j-0}\phi,\ldots, \partial_{z}^{j-1}\partial_{\bar{z}}^{1}\phi,\partial_{z}^{j}\partial_{\bar{z}}^{0}\phi)
\end{equation}
\begin{equation}
\norm{\nabla^j\phi}_Y:=(\partial_{z}^{0}\partial_{\bar{z}}^{j-0}\phi,\ldots, \partial_{z}^{j-1}\partial_{\bar{z}}^{1}\phi,\partial_{z}^{j}\partial_{\bar{z}}^{0}\phi)
\end{equation}
\begin{equation}
\omega_Y(\nabla^j\phi,\delta):=\max_{0\leq s\leq m}\omega_Y(\partial_{z}^{s}\partial_{\bar{z}}^{j-s}\phi,\delta)
\end{equation}
where $\norm{\cdot}_Y$ denotes the sup-norm.
Fedorovskiy \cite{fedorovsky2011} proved the following.
\begin{proposition}\label{fedorovskyproppen}
Let $X\subset\C$ be a compact subset such that $\Theta(X)>0.$ 
Then for each $f\in C^{n-1}(\C)\cap  \mbox{PA}_n(\mbox{int}X)$ there exists a sequence
$\{f_k\}_{k\in \Z_+},$
$f\in \mbox{PA}_n(X)$ (where the latter means $n$-analytic in a neighborhood of $X$)
such that for all $s,t\in \Z_+$ with $s+t\leq n-1,$
\begin{equation}
\partial_{z}^{s}\partial_{\bar{z}}^{t} f_k\to \partial_{z}^{s}\partial_{\bar{z}}^{t} f\mbox{ uniformly as } k\to \infty
\end{equation}
\end{proposition}
\begin{proof}
Without loss of generality we can assume $f$ has compact support, the general case follows by a standard regularization argument (see e.g.\ Remark \ref{introapproxbeviser} and the paragraph following Lemma \ref{carmonalemmatva}).
Recall that the fundamental solution to $\partial_{\bar{z}}^n$ is given by $E:=\bar{z}^{n-1}((n-1)!\pi z)^{-1}.$
For $g\in \mbox{PA}_q(\{\abs{z}<r\})$ we have for $\abs{w-a}<r,$
\begin{equation}
g(w)=\sum_{s=0}^{n-1}\sum_{t=0}^\infty \frac{\partial_{z}^{s}\partial_{\bar{z}}^{t} g(a)}{s!t!}(\bar{w}-\bar{a})^s(w-a)^t
\end{equation}
with converges in $C^\infty(\{\abs{w}<r\}).$
If $T$ is a distribution with compact support in $\{\abs{w}<r\}$ then $f:=E_n*T$ can be written, for $\abs{z-a}>r$
\begin{equation}
f(z)=\sum_{s=0}^{n-1}\sum_{t=0}^\infty c_t^s(f,a) (\bar{w}-\bar{a})^s(w-a)^t\partial_{z}^{s}\partial_{\bar{z}}^{t} E_n(z-a)
\end{equation}
with uniform convergence in $C^\infty(\{\abs{z}> r\})$ and
where
\begin{equation}
c_t^s(f,a):=\frac{(-1)^{n-s-t}}{s!t!}\langle T(w), (\bar{w}-\bar{a})^s(w-a)^t \rangle
\end{equation}
For $\varphi\in C^\infty_c(\{\abs{z}<r\})$ define
$V_\phi \colon (C^\infty_c(\C))' \to (C^\infty_c(\C))',$ $V_\varphi f:=E_n *(\varphi \partial_{\bar{z}}^n f).$
So let $f\in C_c^{n-1}(\C)\cap \mbox{PA}_n(\mbox{int}X).$ Let $R>2$ such that supp$(f)\cup X\subset\{\abs{z}<R/2\}.$
For $\delta\in (0,1)^2$ let $\{ \phi_j\}_{j\in \Z^2}$ be a partition of unity with 
$\phi_j \in C_c^\infty(\{\abs{z-j\delta}<\delta\}),$
such that $0\leq \phi_j\leq 1,$ $\norm{\nabla^k\phi_j}\leq A\delta^{-k}$ for $k=0,\ldots,n$ 
and $\sum_{j\in \Z^2} \phi_j =1.$
(The use of vector notation for $\delta$ and $j$ comes from the fact that the proof
follows that of Mazalov, Paramonov \& Fedorovskiy \cite{mazalovfedorovskiy}, 
which concerns elliptic operators $L$ on $\R^2$, where
of course $L=\partial_{\bar{z}}$ becomes a special case. We could consider $\delta$ as $\delta\cdot (1,1)$ when necessary
for the notation).
For $j\in \Z^2$ set $f_j:=E_n*(\phi_j\partial_{\bar{z}}^n f).$
\begin{lemma}\label{fedorolem23}
Let $f\in C_c^{n-1}(\C),$ $\phi\in C_c^\infty(\{\abs{z-a}<\delta\})$ for all $a\in \C,$ $\delta>0.$
The function $F:=V_\phi f$ satisfies for an absolute constant $A$:\\
(i) $F\in C^{n-1}_{\mbox{loc}}(\C)$, $\partial_{\bar{z}}^n F=0$ on $\C\setminus (\mbox{supp}(\partial_{\bar{z}}^nf)
\cap \mbox{supp}(\phi)).$\\
(ii) $\norm{\partial_z^s\partial_{\bar{z}}^k F}_{\{\abs{z-a}<q\delta\}}\leq A \delta^n\omega(\nabla^{n-1}f,\delta)
\norm{\nabla^{s+t+1}\phi}$ for $s,t\in \Z_+,$ $s+t\leq n-1.$\\
(iii) $\abs{c_t^s(F,a)}\leq \frac{1}{s!t!}A\delta^{s+t+2}\omega(\nabla^{n-1}f,\delta)\norm{\nabla\phi}$ for all $s=0,\ldots,n-1,$
$t\in \Z_+.$
\end{lemma}
\begin{proof}
We give the proof for $f\in C_c^\infty(\C),$ a standard regularization argument (see e.g.\ Remark \ref{introapproxbeviser} and the paragraph following Lemma \ref{carmonalemmatva})
provides the proof for the general case. The proof of (i) is essentially the same as the proof
in the case of the properties of the Vitsushkin localization operator, see Proposition \ref{vitushkinopprops}. 
Setting
\begin{equation}
f_a(w):=f(w)-\sum_{k=0}^{n-1}\frac{\partial_z^s\partial_{\bar{z}}^{k-s} f(a)}{s!(k-s)!}(\bar{w}-\bar{a})^s
(w-a)^{k-s}
\end{equation}
we can write
\begin{equation}
F(z)=\langle E_n(w-z),\phi(w)\partial_{\bar{z}}^n f(w)\rangle =\langle E_n(w-z),\phi(w)\partial_{\bar{z}}^n f_a(w)\rangle
\end{equation}
Let $s,t\in \Z_+$ such that $s+t\leq n-1.$ Now
\begin{multline}\label{usdsa}
\abs{\partial_z^s\partial_{\bar{z}}^t F(z)}=
\abs{\langle\partial_z^s\partial_{\bar{z}}^t E_n(w-z),\phi(w)\partial_{\bar{z}}^n f_a(w)\rangle}\leq \\
\abs{\langle \partial_{\bar{z}}(\phi(w)\partial_z^s\partial_{\bar{z}}^t E_n(w-z)),\partial_{\bar{z}}^{n-1}f_a(w)\rangle}=\\
\abs{\langle \partial_{\bar{z}}(\phi(w)\partial_z^s\partial_{\bar{z}}^t E_n(w-z)+
\phi(w)\partial_z^{s+1}\partial_{\bar{z}}^t E_n(w-z),\partial_{\bar{z}}^{n-1}f_a(w)\rangle}\leq\\
\abs{\partial_{\bar{z}} \phi(w)\partial_z^{s}\partial_{\bar{z}}^t E_n(w-z),\partial_{\bar{z}}^{n-1}f_a(w)\rangle} +
\abs{\phi(w)\partial_z^{s+1}\partial_{\bar{z}}^t E_n(w-z),\partial_{\bar{z}}^{n-1}f_a(w)\rangle}
\end{multline}
If $s+t\leq n-2$ then integrating the terms in the last expression in polar coordinates with origin $z$ yields
for absolute constants $A_1,A_2,$
\begin{equation}
\norm{\partial_z^{s}\partial_{\bar{z}}^t F}_{\{\abs{\zeta-a}<q\delta\}}\leq\\
A_1\omega(\nabla^{n-1}f,\delta)\norm{\nabla\phi}\delta^{n-s-t}+A_2\omega(\nabla^{n-1}f,\delta)\norm{\phi}\delta^{n-s-t-1}
\end{equation}
Since also $\norm{\nabla^l\phi}\leq A\delta^{m+1-l}\norm{\nabla^{m+1}\phi},$ $l=0,\ldots,m,$ $m\in\Z_+,$
this proves (ii) whenever $s+t\leq n-2.$
Writing
\begin{multline}
\phi(z)\partial_z^s\partial_{\bar{z}}^t E_n(w-z)=\\
\partial_z^t\left(\phi(w)\partial_{\bar{z}}^s E_n(w-z)\right)-
\sum_{j=0}^{t-1}\binom{t}{j} \partial_z^s\partial_{\bar{z}}^j E_n(w-z)\partial_z^{t-j}\phi(w)
\end{multline}
gives
\begin{multline}
\abs{\partial_z^s\partial_{\bar{z}}^t F(z)}\leq \abs{\langle \partial_z^t(\phi(w)\partial_{\bar{z}}^s E_n(w-z),\partial_{\bar{z}}^n f_a(w)\rangle}+\\
\sum_{j=0}^{t-1}\binom{t}{j} \abs{\langle \partial_z^s\partial_{\bar{z}}^j  E_n(w-z)\partial_z^{t-j}\phi(w),\partial_{\bar{z}}^n f_a(w)\rangle}
\end{multline}
For the first term in the right hand side we have
\begin{multline}
\abs{\langle \partial_z^t(\phi(w)\partial_{\bar{z}}^s E_n(w-z),\partial_{\bar{z}}^n f_a(w)\rangle}=\\
\abs{\langle \partial_{\bar{z}}^{n-s}(\phi(w)\partial_{\bar{z}}^s E_n(w-z),\partial_z^s\partial_{\bar{z}}^t f_a(w)\rangle} \leq \phi(z)(\partial_z^s\partial_{\bar{z}}^t f(z)-
\partial_z^s\partial_{\bar{z}}^tf(a))+\\
\sum_{j=1}^{n-s} \binom{n-s}{j}\abs{\langle \partial_{\bar{z}}\phi(w)\partial_{\bar{z}}^{n-j}E_n(w-z),\partial_z^s\partial_{\bar{z}}^t f_a(w)\rangle }
\end{multline}
thus for some absolute constants $A_j,$
\begin{multline}
\sup_{z\in \{\abs{z-a}<q\delta\}} \abs{\langle \partial_z^t(\phi(w))\partial_{\bar{z}}^s E_n(w-z),\partial_{\bar{z}}^n f_a(w)\rangle} \leq A_0\omega(\nabla^{n-1}f,\delta)\norm{\phi}
+\\
\sum_{j=1}^{n-s} \binom{n-s}{j} A_j\omega(\nabla^{n-1}f,\delta)\norm{\nabla^j \phi}\delta^j\leq A\delta^{s+t+1}
\norm{\nabla^{s+t+1}\phi}\omega(\nabla^{n-1}f,\delta)
\end{multline}
Also we have
\begin{multline}
\abs{\langle \partial_z^s\partial_{\bar{z}}^j E_n(w-z)\partial_z^{t-j} \phi(w),
\partial_{\bar{z}}^n f_a(w)\rangle} \leq\\
\abs{\partial_z^{s+1}\partial_{\bar{z}}^j E_n(w-z)\partial_z^{0}\partial_{\bar{z}}^{t-j}\phi(w),\partial_{\bar{z}}^{n-1} f_a(w)\rangle}+\\
\abs{\langle \partial_z^{s}\partial_{\bar{z}}^j E_n(w-z),\partial_{\bar{z}}^{n-1} f_a(w)\rangle}
\end{multline}
Since $s+j\leq n-2$ for $j=1,\ldots ,t-1,$ both terms in the right hand side can be estimated analogous to
the terms $\abs{\langle \partial_{\bar{z}}\phi(w)\partial_z^{s}\partial_{\bar{z}}^t E_n(w-z),\partial_{\bar{z}}^{n-1}f_a(w)\rangle}$
and $\abs{\langle \phi(w)\partial_z^{s+1}\partial_{\bar{z}}^t E_n(w-z),\partial_{\bar{z}}^{n-1}f_a(w)\rangle}$ in 
Eqn.(\ref{usdsa}). This proves (ii). Finally, for $s,t\in \Z_+,$ $s\leq n-1$ we have
\begin{multline}
\abs{c_t^s(F,a)}s!t! =\abs{\langle \phi(w)\partial_{\bar{z}}^nf(w),(\bar{w}-\bar{a})^s(w-a)^t\rangle}=\\
\abs{\langle \partial_{\bar{z}}^{n-1} f_a(w),\partial_{\bar{z}}(\phi(\bar{w}-\bar{a})^s(w-a)^t\rangle}\leq A\delta^{s+t+2}\omega
(\nabla^{n-1}f,\delta)\norm{\nabla\phi}
\end{multline}
which proves (iii). This completes the proof of Lemma \ref{fedorolem23}.
\end{proof}
By Lemma \ref{fedorolem23} each $f_j\in C_{\mbox{loc}}^{n-1}(\C)$, $f_j$ is $n$-analytic on $\mbox{int}X$ and
on the complement of $\{\abs{z-a}<q\delta\},$ and satisfies
for $s,t\in \Z_+,$ with $s+t\leq n-1,$
\begin{equation}\label{tva5an}
\norm{\partial_z\partial_{\bar{z}}^t f_j}_{\{\abs{z-a_j}<q\delta\}}\leq A\delta^{n-1-s-t}\omega(\nabla^{n-1}f,\delta)
\end{equation}
and for $s=0,\ldots,n-1,$ $t\in \Z_+,$
\begin{equation}\label{tva7an}
\abs{c_t^s(f_j,a_j)}\leq \frac{A}{s!t!}\delta^{s+t+1}\omega(\nabla^{n-1}f,\delta)
\end{equation}
Also the sum $f(z)=\sum_{j\in \Z^2} f_j(z)$ is finite since $f_j\equiv 0$ if 
supp$(\phi_j)\cap$supp$(\partial_{\bar{z}}^n f)=\emptyset .$
Since $f_j\equiv 0$ if $\{\abs{z-a_j}<\delta\}$ and $f_j$ is $n$-analytic in a neighborhood of $X$ whenever $\{\abs{z-a_j}<\delta\}\cap X=\emptyset,$
it suffices to approximate only the functions $f_j$ such that $j\in J:=\{j\in \Z^2:\{\abs{z-a_j}<\delta\}\cap
\partial X\neq \emptyset\}.$
Let $\delta$ be such that for all $j\in J, \{\abs{z-a_j}<2\delta/\theta\}\subset \{\abs{z}<R\}.$
Then there is $a_j^*\in \{\abs{z-a_j}<2\delta/\theta\}\setminus X$ and a Jordan curve $\gamma_j\subset \{\abs{z-a_j^*}\leq \delta\}$
with initial point
$a_j^*$ and end point in $\{\abs{z-a^*_j}=\delta\}.$ 
Pick a constant $c_0$ and polynomial $p_0$ independent of $\gamma/\delta$ and a holomorphic
branch of the square root $\sqrt{z(z-a/\delta)}$ outside $\gamma/\delta$, such that $\lim_{z\to \infty} zh(z)=1$ with
\begin{equation}
h(z):=c_0(z^n(z-a/\delta)^n \sqrt{z(z-a/\delta)}-p_0(z))
\end{equation}
\begin{equation}
g_\gamma(z):=\frac{\bar{z}^{n-1}h(z/\delta)}{\pi\delta(n-1)!}
\end{equation}
For $j\in J$ set $\gamma_j^*:=\gamma_j-a^*_j$ and $g_j^*(z):=g_{\gamma_j}^*(z-a_j^*).$ Then
for $s=0,\ldots,n-1,$ $t=0,\ldots,n,$
\begin{equation}\label{tva8an}
\norm{\partial_z^s\partial_{\bar{z}}^t g_j^*}_{\abs{z-a_j}<q\delta\}\setminus \gamma_j^*}\leq A\delta^{n-2-s-t}
\end{equation}
and for $\abs{z-a_j^*}>\delta,$ 
\begin{equation}
g_j^* (z)=\sum_{t=0}^\infty d_{j,0}^* \partial_z^t E_n(z-a_j^*) 
\end{equation}
for $d_{j,t}^*$ such that $d_{j,0}^*=1,$ $\abs{d_{j,t}^*}\leq A\delta^t/t!,$ $t\in \Z_+ .$
Set for $s=0,\ldots,n-1,$ $l=0,\ldots,n-s,$
\begin{equation}
\beta_{j,s,0}:=c_0^s(f_j,a_j^*),\quad \beta_{j,s,l}:=c_l^s(f_j,a_j^*)-\sum_{k=0}^{l-1}\beta_{j,s,k}d^*_{j,l-k}
\end{equation}
Setting
\begin{equation}
g_j(z)=\sum_{s=0}^{n-1}\sum_{l=0}^{n-s}\beta_{j,s,l}\partial_z^s\partial_{\bar{z}}^l g_j^*(z)
\end{equation}
Thus 
\begin{equation}
g_j(z)=\sum_{s=0}^{n-1}\sum_{l=0}^{\infty} c_t^s(g_j,a_j^*)\partial_z^s\partial_{\bar{z}}^t E_n(z-a_j^*)
\end{equation}
converges in $C^\infty(\abs{z-a_j^*}>q\delta)$ with
\begin{equation}
c_t^s(g_j,a_j^*)=\sum_{l=0}^{\min\{t,n-s\}}\beta_{j,s,l}d^*_{j,t-l}
\end{equation}
for $s=0,\ldots,n-1,$ $t\in \Z_+.$ 
By Eqn.(\ref{tva8an}) and Eqn.(\ref{tva5an}) together with $\abs{\beta_{j,s,l}}\leq A\delta^{s+l+1}\omega(\nabla^{n-1}f,\delta)$ gives
for $s,t\in \Z_+,$ such that $s+t\leq n-1,$
\begin{equation}\label{tva10an}
\norm{\partial_z^s\partial_{\bar{z}}^t g_j}_{\{\abs{z-a_j}<q\delta\}\setminus \gamma_j^*} \leq A\delta^{n-1-s-t}\omega(\nabla^{n-1}f,\delta)
\end{equation}
and Eqn.(\ref{tva7an}) together with the inequalities $\abs{d_{j,t}^*}\leq A\delta^t/t!,$ $t\in \Z_+ $ given above 
we obtain for $s=0,\ldots,n-1,$ $t\in \Z_+,$
\begin{equation}\label{tva11an}
\abs{c_t^s(g_j,a_j^*)}\leq \frac{A}{s!t!}\delta^{s+t+1}\omega(\nabla^{n-1}f)t^n
\end{equation}
Since each $g_j$ is $n$-analytic in a neighborhood of $X$ as it is $n$-analytic outside of $\gamma_j$, it suffices to
show for $s,t\in \Z_+,$ such that $s+t\leq n-1,$
\begin{equation}
\norm{\sum_{j\in J}( \partial_z^s\partial_{\bar{z}}^t (f_j-g_j)}_X \leq A\delta^{n-1-s-t}\omega(\nabla^{n-1}f,\delta)
\end{equation}
By definition of $\beta_{j,s,l}$ we have for $s=0,\ldots,n-1,$
$t=0,\ldots,n-s,$ $c_t^s(f_j,a_j^*)=c_t^s(g_j,a^*_j)$ which for $s,t\in \Z_+,$ such that $s+t\leq n-1,$ implies together with 
Eqn.(\ref{tva7an}) and Eqn.(\ref{tva11an})
that for $\abs{z-a_j^*}>q\delta,$ 
\begin{equation}
\abs{\partial_z^s\partial_{\bar{z}}^t (f_j-g_j)}\leq \frac{A\delta^{n+2}\omega(\nabla^{n-1}f,\delta)}{\abs{z-a_j^*}^{3+s+t}}
\end{equation}
This implies together with Eqn.(\ref{tva5an}) and Eqn.(\ref{tva10an}), that for $\abs{z-a_j^*}>q\delta,$ $z\notin \gamma_j^*$
\begin{equation}
\abs{\partial_z^s\partial_{\bar{z}}^t (f_j-g_j)}\leq A\delta^{n-1-s-t}\omega(\nabla^{n-1}f,\delta)
\end{equation}
Since the number of indices $j$ such that
$k\delta\leq \abs{z-a_j^*}<(k+1)\delta,$ $k\geq 1,$ is $\leq Ak$
we have
\begin{multline}
\abs{\sum_{j\in J} \partial_z^s\partial_{\bar{z}}^t (f_j-g_j)}\leq\\
\sum_{j\in J} \abs{\partial_z^s\partial_{\bar{z}}^t (f_j-g_j)}\leq\\
\sum_{j:\abs{z-a_j^*}<q\delta} \abs{\partial_z^s\partial_{\bar{z}}^t (f_j-g_j)}\leq\\
\sum_{k=[q]}^\infty \sum_{k\delta\leq \abs{z-a_j^*}<(k+1)\delta}  \abs{\partial_z^s\partial_{\bar{z}}^t (f_j-g_j)}\leq\\
\delta^{n-1-s-t}\omega(\nabla^{n-1}f,\delta)\left( A_1+A_2\sum_{k=[q]}^\infty k^{-2-s-t}\right)\leq A\delta^{n-1-s-t}\omega(\nabla^{n-1}f,\delta)
\end{multline}
where $[q]$ denotes the least upper integer.
So a function in $C^{n-1}(\C)$, coinciding with $\sum_{j\in J} g_j(z)+\sum_{j\notin J} f_j(z)$ in a neighborhood of $X$
yields the sought approximation of $f.$ This completes the proof of Proposition \ref{fedorovskyproppen}.
\end{proof}
As a consequence of Proposition \ref{fedorovskyproppen} together with 
the Malgrange-Lax theorem incorporation Lemma \ref{narasimlem1}
(the method can be found in e.g.\ Narasimhan \cite{narasimhan}, parag 3.10,
and applies for solutions to elliptic equations), Fedorovskiy \cite{fedorovsky2011} proved
the following.
\begin{theorem}
Let $X\subset\C$ be a compact subset. Then the following are equivalent.\\
(i) If $f\in C^{n-1}(\C)\cap \mbox{PA}_q(\mbox{int}X)$ then there exists a sequence $\{p_k\}_{k\in \Z_+}$ of $n$-analytic polynomials such that
$\partial_z^s\partial_{\bar{z}}^t p_k \to \partial_z^s\partial_{\bar{z}}^t f$ uniformly on $X$, for all $s,t\in \Z_+,$ such that $s+t\leq n-1,$ 
as $k\to \infty$.\\
(ii) $\C\setminus X$ is connected.
\end{theorem}
\begin{proof}
(i)$\Rightarrow$(ii). Assume $\C\setminus X$ is not connected and let $U$ be an unbounded connected component and
$p_0\in U$ and let $R=$diam$(X).$ Let $f(z)=E_n(z)$ for $\frac{1}{2}\mbox{dist}(p_0,U)<\abs{z}<R$ and extended to a function in $C^{n-1}(\C)$ and assume
there exists a sequence $\{p_j\}_{j\in \Z_+}$ such that
$\partial_z^s\partial_{\bar{z}}^t p_j\to \partial_z^s\partial_{\bar{z}}^t f$ uniformly on $X$ for all $s,t\in \Z_{\geq 0},$ such that $s+t\leq n-1.$
Then the sequence $\{ \partial_{\bar{z}}^{n-1} p_j \}_{j\in \Z_+}$ converges uniformly on $X$ to 
$\partial_{\bar{z}}^{n-1} E_n=1/(\pi z),$ hence $1/z$ can be approximated uniformly on $\partial U$ by complex polynomials which
yields a contradiction. This proves (i)$\Rightarrow$(ii).\\
(ii)$\Rightarrow$(i). Suppose $\C\setminus X$ is connected and let $f\in C^{n-1}(\C)\cap \mbox{PA}_q(\mbox{int}(X)).$ By
Proposition \ref{fedorovskyproppen} there exists a sequence $\{f_j\}_{j\in \Z_+}$ of functions in $A_n(X)$ such that
$\partial_z^s\partial_{\bar{z}}^t f_k \to \partial_z^s\partial_{\bar{z}}^t f$ uniformly on $X$, for all $s,t\in \Z_+,$ such that $s+t\leq n-1,$ 
as $k\to \infty$. Suppose $f_k$ is $n$-analytic on a neighborhood $U_k$ of $X$ and let $V_k\subset U_k$ be a neigbhorhood
of $X$ such that each point of $\C\setminus U_k$
can be joined to $\infty$ by a curve lying outside the closure of $V_k.$ 
By e.g.\ the proof of Narasimhan \cite{narasimhan}, Theorem 3.10.7 (recalling that the operator $\partial_{\bar{z}}^n$ is an elliptic operator) each $f_n$ and its partial derivatives up to order $n-1$ can be uniformly approximated
by $n$-analytic polynomials on the closure of $V_k$. 
This proves (ii)$\Rightarrow$(i).
This completes the proof.
\end{proof}
\begin{definition}[Nevanlinna domain]\index{Nevanlinna domain}
Let $\Omega\subset\C$ be a bounded simply connected domain.
Let $H^\infty(\Omega)$ denote the set of bounded holomorphic functions on $\Omega.$
$\Omega$ is called a {\em Nevanlinna domain} if there exists two functions $u,v\in H^\infty(\Omega)$, $v\not\equiv 0,$ such that
\begin{equation}
\bar{\zeta}=\frac{u(\zeta)}{v(\zeta)}, \mbox{ for almost every }\zeta\in \partial\Omega 
\end{equation}
in the sense that
the following equality of angular values holds
\begin{equation}
\overline{\varphi(\xi)}=\frac{(u\circ\varphi)(\xi)}{(v\circ\varphi)(\xi)}, \mbox{ for almost every }\xi\in \{\abs{z}=1\}
\end{equation}
where $\varphi$ is some conformal mapping from $\{\abs{z}<1\}$ onto $\Omega.$
By the Luzin-Privalov theorem the ratio $u/v$ is uniquely defined in $\Omega$ for a Nevanlinna domain $\Omega$. 
Denote by $\mathcal{N}$ the class of Nevanlinna domains.
\end{definition}
\begin{definition}
A bounded domain $\Omega\subset\C$ is called a {\em Carath\'eodory domain}\index{Carath\'eodory domain}
if $\partial\Omega=\partial\Omega_\infty$, where $\partial\Omega_\infty$ denotes the unbounded connected component of the set
$\hat{\C}\setminus \overline{\Omega},$ $\hat{\C}=\C\cup\{\infty\}.$
For a compact $X\subset \C,$ denote by $\hat{X}$ the union of $X$ and all bounded connected components of the set 
$\C\setminus X.$ $X$ is called a {\em Carath\'eodory compact} if $\partial X=\partial\hat{X}.$ $\hat{X}$ is precisely the 
{\em polynomial hull of $X$}.
\end{definition}
Each Carath\'eodory domain $\Omega$ is simply connected and satisfies $\Omega=\mbox{int}(\overline{\Omega}).$
Denote the Cauchy transform by
\begin{equation}
\hat{\mu}(z):=\frac{1}{2\pi i}\int \frac{d\mu(\zeta)}{\zeta -z}
\end{equation}
and recall that 
this is known to satisfy that $\hat{\mu}$ is holomorphic off supp$(\mu)$
and $\overline{\partial} \hat{\mu} =\frac{i}{2}\mu$ in distribution sense.
The following is yet an example where the techniques (scheme of approximation) of 
Vitushkin has inspired a result in the polyanalytic case (see Carmona, Paramonov \& Fedorovskiy \cite{cpf}, Proposition 2.5). 
\begin{theorem}\label{vitushkincopythm}
If $X\subset \C$ is a Carath\'eodory compact then $R_n(X,X)=A_n(X).$ 
\end{theorem}
\begin{proof}
	We note that 
	this is a generalization of a known result for holomorphic functions, see Vitushkin \cite{vitushkin67}, Ch.2, parag. 5, Thm.1.
Recall that the fundamental solution to $\partial_{\bar{z}}^n$ is given by $E:=\bar{z}^{n-1}((n-1)!\pi z)^{-1}.$
We shall follows the method (scheme) in the proof of
Proposition \ref{fedorovskyproppen} (see also Mazalov, Paramonov \& Fedorovskiy \cite{mazalovfedorovskiy} which in turn is the blueprint used in the original proof
of Proposition \ref{fedorovskyproppen}), hence we give only the main steps, see the aforementioned proof for details.
Let $f\in A_n(X)$ and we keep the same notation for an extension of $f$ to a continuous function with compact support in $\C.$
For a fixed $\delta>0$ let $\{\phi_j\}_{j\in \Z^2}$ be a partition of unity
$\phi_j\in C^\infty_c(\{\abs{z-a_j}\}<\delta\}),$ $a_j=j\delta,$
$0\leq \phi_j\leq 1 ,$ $\sum_j \phi_j=1,$ $j=(j_1,j_2)$ and set 
\begin{equation}
J:=\{j\in \Z^2\colon \{\abs{z-a_j}<\delta\} \cap \partial X\neq \emptyset\} 
\end{equation}
For each $j\in J$ there exists $a^*_j \in \{\abs{z-a_j}<\delta\} \cap \partial X$ and a Jordan curve $\gamma_j$ in 
$\{\abs{z-a_j}\leq \delta\} \cap \partial X$ with initial point $a_j^*$ and end point in $\{\abs{z-a_j}= \delta\}.$
For $j\notin J$ set $a_j^*=a_j.$ Set $f_j:=E*(\phi_j \partial_{\bar{z}}^n f),$ so 
in distribution sense
$\partial_{\bar{z}}^n f_j=\phi_j\partial_{\bar{z}}^n f$ and (as in the proof of
Proposition \ref{fedorovskyproppen}, see also Mazalov, Paramonov \& Fedorovskiy \cite{mazalovfedorovskiy}, p.295) 
we have for constants $c^s_m$ depending on $j$
\begin{equation}
f_j(z)=\sum_{s=0}^{n-1}\sum_{m=0}^\infty c_m^s \partial_{\bar{z}}^s \partial_z^m E(z-a_j^*)
\end{equation}
with uniform convergence in $C^\infty(\{ \abs{z-a_j^*}>2\delta\})$ 
such that, for a constant $A$ depending only on $n$, and all $s=0,\ldots,n-1,$ $m\in \Z_{\geq 0},$
\begin{equation}
\norm{f_j}_{\{\abs{z-a_j^*}<3\delta\}}\leq A\omega_{\{\abs{z-a}<\delta\}}(f,\delta),\quad \abs{c_m^s}\leq\frac{A\omega_{\{\abs{z-a}<\delta\}}(f,\delta)(2\delta)^{2-n+m+s}}{m!s!}
\end{equation}
Let $\gamma$ be a Jordan curve in $\{\abs{z}\leq \delta\}$
with initial point
$0$ and end point in $a\in\{\abs{z}=\delta\}.$
Pick a constant $c_0$ and polynomial $p_0$ independent of $\gamma/\delta$ and a holomorphic
branch of the square root $\sqrt{z(z-a/\delta)}$ outside $\gamma/\delta$, such that $\lim_{z\to \infty} zh(z)=1$ with
\begin{equation}
h(z):=c_0(z^n(z-a/\delta)^n \sqrt{z(z-a/\delta)}-p_0(z))
\end{equation}
\begin{equation}
g_\gamma(z):=\frac{\bar{z}^{n-1}h(z/\delta)}{\pi\delta(n-1)!}
\end{equation} 
As in the proof of
Proposition \ref{fedorovskyproppen} we have
for $s=0,\ldots,n-1,$ $m=0,\ldots,n,$ and $\abs{z}>\delta,$
\begin{equation}
\norm{\partial_z^m\partial_{\bar{z}}^s g_j^*}_{\{\abs{z}<3\delta\}\setminus \gamma}\leq A\delta^{n-2-s-m}
\end{equation}
\begin{equation}
g_\gamma \sum_{t=0}^\infty d_m \partial_z^m E_n(z-a)
\end{equation}
for $d_m$ such that $d_0=1,$ $\abs{d_m}\leq A\delta^m/m!,$ $m\in \Z_{\geq 0} .$
For $j\in J$ set $\gamma_j^*=\gamma_j-a_j^*$ and $g_j^*(z):=g_{\gamma_j^*}(z-a_j^*)$ and
\begin{equation}
g_j(z):=\sum_{s=0}^{n-1}\sum_{l=0}^{n-s}\beta_{l}^s \partial_{\bar{z}}^s\partial_{z}^l g_j^*(z)=:
\sum_{s=0}^{n-1}\sum_{m=0}^{\infty}b_m^s \partial_{\bar{z}}^s\partial_{z}^m E_n(z-a_j^*)
\end{equation}
and the series converges in $C^\infty(\abs{z-a_j^*}>2\delta)$ with
$c_m^s=\beta_m^s$ for
for $s=0,\ldots,n-1,$ $m=0,\ldots,n-s.$
Let
\begin{equation}
g_j^*(z)=\sum_{m=0}^\infty d_m^* \partial_z^m E_n(z-a_j^*)
\end{equation}
for $d_0^*=1,$ $\abs{d_m^*}\leq A\delta^m/m!,$ $m\in \Z_{\geq 0},$ which implies
\begin{equation}
b_m^s=\sum_{l=0}^{\min\{m,n-s\}} \beta_{l}^s d^*_{m-l}
\end{equation}
for $s=0,\ldots,n-1,$ $m=0,\ldots,n-s,$ and using
$\abs{\beta_l^s}\leq A\omega_{\{\abs{z-a}<\delta\}}(f,\delta)\delta^{2-n+l+s}$ also
\begin{equation}
\norm{g_j}_{\{\abs{z-a_j^*}<3\delta\}} \leq A\omega_{\{\abs{z-a}<\delta\}}(f,\delta)
\end{equation}
\begin{equation}
\abs{b_m^s}\leq \frac{A\omega_{\{\abs{z-a}<\delta\}}(f,\delta)(2d)^{2-n+m+s} m^n}{s!m!}
\end{equation}
for $s=0,\ldots,n-1,$ $m\in \Z_{\geq 0}.$
Hence for $z\in \C,$ and using on $\C\setminus\left\{\abs{z-a_j^*}\geq 3\delta\right\}$, the above estimates for $c_m^s$, $b_m^s,$
$\partial_{\bar{z}}^s \partial_z^m E_n(z-a_j^*)$ (and the corresponding Laurent-expansions at $a_j^*$) we obtain
\begin{equation}
\abs{f_j(z)-g_j(z)}\leq A\omega_{\{\abs{z-a}<\delta\}}(f,\delta) \min\left(1,\delta^3\abs{z-a_j^*}^{-3}\right)
\end{equation} 
The inequality holds true also if $a_j^*$ is replaced by $a_j$ in the right hand side, thus by Proposition \ref{vitusanvandicarmona}
$f$ can be uniformly approximated on $X$ by
$\sum_{j\in J} g_j+\sum_{j\notin J} f_j$ as $\delta\to 0.$ Since the $g_j$ are $n$-analytic in a neighborhood
of $X$ this implies $f_j\equiv 0$ whenever $\{\abs{z-a_j}<\delta\}\subset\mbox{int}(X)$ and
$\partial_{\bar{z}}f_j=0$ on $\{\abs{z-a_j}\geq \delta\}.$ Finally, applying Runge's theorem completes the proof
of Theorem \ref{vitushkincopythm}.
\end{proof}
The following was proved by Carmona, Paramonov \& Fedorovskyi \cite{cpf}.
\begin{theorem}
Let $q\geq 2$ be an integer.
\\
(1) If $\Omega\subset \C$ is a Carath\'eodory domain then
$C^0(\partial\Omega)=R_q(\partial\Omega,\overline{\Omega})$ if and only if $\Omega$ is not a Nevanlinna domain.\\
(2) If $X\subset \C$ is a  Carath\'eodory compact, then $A_q(X)=P_q(X)$ if and only if
$C^0(\partial\Omega)=R_q(\partial\Omega,\overline{\Omega})$ for each bounded connected component $\Omega$ of the set 
$\C\setminus X$ and, therefore, if and only if each bounded connected component of the set $\C\setminus X$ is not a Nevanlinna domain.
\end{theorem}
\begin{proof}
(i) $(\Rightarrow)$:\\
We prove the contrapositive.
Let $\Omega\subset\mathcal{N}.$ We show that this implies $C^0(\partial\Omega)\neq R_n(\partial\Omega,\overline{\Omega}).$
Let $k$ be a conformal map of $\{\abs{z}<1\}$ onto $\Omega$, and let $k(t)$ denote the boundary values, in particular by definition
there exists bounded holomorphic functions $u,v$,$v\not\equiv 0$ on $\Omega$
such that 
\begin{equation}
\overline{k(t)}=\frac{u(k(t))}{v(k(t))}
\end{equation}
and the ratio $u/v$ is uniquely determined in $\Omega.$
Let $z_0\in \Omega$ such that $\abs{u(z_0)+\bar{z}_0v(z_0)}\neq 0.$
Assume (in order reach a contradiction) that 
\begin{equation}\label{carm2002contra}
\frac{(\bar{z}-\bar{z}_0)^{n-1}}{z-z_0}|_\partial\Omega \in R_n(\partial\Omega,\overline{\Omega})
\end{equation}
Then for any $\delta>0$ there exists rational functions $f_0,\ldots,f_{n-1}$ with poles off $\overline{\Omega}$
such that for $\zeta\in \partial\Omega$
\begin{equation}
\abs{ \sum_{s=0}^{n-1} f_s(\zeta)-\frac{(\bar{\zeta}-\bar{z}_0)^{n-1}}{\zeta-z_0} }<\delta
\end{equation}
Setting $w_0:=k^{-1}(z_0)\in \{\abs{z}<1\}$ gives for almost all $t\in \{\abs{z}=1\}$ (here $k(t) \in \partial \Omega$)
\begin{equation}
\abs{ \sum_{s=0}^{n-1} \bar{\zeta}^s f_s(k(t))\frac{u(k(t))^s}{v(k(t))^s} -
\frac{(u(k(t))-\overline{k(w_0)}v(k(t)))^{n-1}}{v(k(t))^{n-1} (k(t)-k(w_0))} }<\delta
\end{equation}
This implies that for almost all $t\in \{\abs{z}=1\}$
\begin{multline}
\abs{ \sum_{s=0}^{n-1} f_s(k(t))\frac{u(k(t))^s}{v(k(t))^{s+1-n}}(k(t)-k(w_0))-(u(k(t))-\overline{k(w_0)}v(k(t)) )^{n-1} }\\
\leq \delta \sup_{\tau\in \{\abs{z}=1\}} \abs{v(k(\tau))^{n-1}(k(\tau)-k(w_0)) }
\end{multline}
Since all involved functions are boundary values of bounded holomorphic functions we have by the maximum principle
\begin{equation}
\abs{ u(k(w_0))- \overline{k(w_0)}v(k(w_0)) )^{n-1} } \leq \delta \sup_{\tau\in \{\abs{z}=1\}} \abs{v(k(\tau))^{n-1}(k(\tau)-k(w_0)) }
\end{equation}
which for sufficiently small $\delta$ implies that
Eqn.(\ref{carm2002contra}) does not hold. This proves
\begin{equation}
\frac{(\bar{z}-\bar{z}_0)^{n-1}}{z-z_0}|_\partial\Omega \notin R_n(\partial\Omega,\overline{\Omega})
\end{equation}
The contrapositive to what we have proved is precisely (i) $(\Rightarrow)$.\\
\\
(i) $(\Leftarrow)$:\\
Again we prove the contrapositive. So assume $C(\partial\Omega)\neq R_n(\partial\Omega,\overline{\Omega}),$
we can assume $C(\partial\Omega)\neq R_2(\partial\Omega,\overline{\Omega}).$
Then there exists a measure $m$ on $\partial\omega,$ $m\not\equiv 0$ such that $m\perp R_1(\partial\Omega,\overline{\Omega})$
such that 
$\bar{\zeta}m\perp R_1(\partial\Omega,\overline{\Omega}).$ 
If $\Omega$ is a Caratheodory domain then $m$ and $\bar{\zeta}m$ vanish on isolated points
(see Bishop \cite{bishoplemmanine}, Lemma 9).
Set $F(z):=\hat{m}(z),$ $G(z):=\widehat{\bar{\zeta}m},$ $f(w):=F(k(w))k'(w),$
$g(w):=G(k(w))k'(w).$
\begin{lemma}\label{atomlemma}
$f$ and $g$ belong to the Hardy class $H^1(\{\abs{z}<1\})$ of $L^1$ analytic functions.
\end{lemma}
\begin{proof}
Let $z_0=k(0)$. We can assume $k'(0)>0.$ Let $\{T_l\}_{l\in \Z_+}$ be a sequence of rectifiable contours,
where each $T_l$ encloses the domain $D_m$ such that
$D_l\subset D_{l-1}$ and such that the $T_l$ converge to $\partial\Omega$. 
Let $\phi_l:D_l\to \{\abs{z}<1\}$ be conformal maps satisfying
$\phi_l(z_0)=0,$ $\phi'_l(z_0)>0.$ 
Then $\pi_m$ converge to $k^{-1}$ uniformly on compacts of $\Omega$
(see Carath\'eodory \cite{carateodoryboken}, p.74-80).
For a fixed $w\in \{\abs{z}<1\}$ set $z_l:=\phi_l^{-1}(w)$ and
\begin{equation}
h(\zeta)=
\left\{
\begin{array}{ll}
\frac{1}{\phi_l(\zeta) -\phi_l(z_l)} - \frac{1}{\phi_l'(\zeta)(\zeta-z_l)}, &  \zeta \neq z_l \\
-\frac{\phi_l''(\zeta)}{\phi'_l(z_l)}, & \zeta=z_l 
\end{array}
\right.
\end{equation}
Then $h\in \mathscr{O}(D_l)$ thus $h\in R_1(\partial\Omega,\overline{\Omega}),$
so $m\perp R_1(\partial\Omega,\overline{\Omega})$ implies
\begin{equation}
\frac{1}{2\pi i} \int \frac{dm(\zeta)}{\phi_l(\zeta) -\phi_l(z_l)} =\frac{\hat{m}(z_l)}{\phi_l'(z_l)}
\end{equation}
Set $m_l(S):=m(\phi_l^{-1}(S\cap \{\abs{z}<1\}))$ for Borel subsets $S$ of $\C.$
We have supp$m_l\subset\{\abs{z}<1\}$, $m_l\perp P_1$, $\norm{m_l}\leq \norm{m}$ and
\begin{equation}\label{oooooooooo}
\frac{1}{2\pi i} \int \frac{dm_l(\zeta)}{\tau -w} = \hat{m}(\phi_l^{-1}(w))(\phi_l^{-1})'(w)
\end{equation}
Let $m_0$ be a limit point of the sequence $m_l$ in the weak-star topology in the space of measures.
Then supp$m_0\subset\{\abs{z}=1\}$ and $m_0\perp P_1(\{\abs{z}=1\}).$
By the Riesz representation theorem there exists $h\in H^1(\{\abs{z}<1\})$ such that $m_0=h(t)dt|_{\{\abs{z}=1\}}$, where $h(t)$ denote
the boundary values of $h.$ 
Letting $l\to \infty$ in Eqn.(\ref{oooooooooo}) gives for $w\in \{\abs{z}<1\}$
\begin{multline}
\hat{m}_0(w)=\hat{m}(k(w))k'(w)=\frac{1}{2\pi i} \int_{\{\abs{z}=1\}}\frac{dm_0(t)}{t -w} =\\ 
\frac{1}{2\pi i} \int_{\{\abs{z}=1\}}\frac{h(t)dt}{t -w}
=h(w)
\end{multline}
Thus $\hat{m}(k(w))k'(w)\in H^1(\{\abs{z}<1\}).$
The same arguments can be repeated for $m$ replaced by $\bar{\zeta}m$.
This completes the proof of Lemma \ref{atomlemma}.
\end{proof}
Set $\Phi(z):=\bar{z}F(z)-G(z)=\frac{1}{2\pi i}\int \frac{\bar{z}-\bar{\zeta}}{\zeta -z}dm(\zeta).$
Since $\frac{\bar{z}-\bar{\zeta}}{\zeta -z}$ is bounded and $m$ vanishes on isolated points we have $\Phi\in C^0(\C).$
For fixed $z\in \C\setminus \overline{\Omega}$ we have $\frac{\bar{z}-\bar{\zeta}}{\zeta -z}\in R_2(\partial\Omega,\overline{\Omega}),$ so $\Phi=0$
on $\C\setminus \overline{\Omega}$ which implies
$\bar{z}F(z)-G(z)\to 0$ as dist$(z,\partial\Omega)\to 0,$ $z\in\Omega.$ Hence as $\{\abs{\zeta}<1\}\ni w,$ $\abs{w}\to 1$
where$z=k(w)$ we have
\begin{equation}
\frac{f(w)\overline{k(w)}-g(w)}{k'(w)}=F(k(w))\overline{k(w)} -G(k(w))=\bar{z}F(z)-G(z)\to 0
\end{equation}
Now $f(w)\overline{k(w)} -g(w)$ has angular boundary values a.e.\ on $\{\abs{\zeta}=1\}$ because $f,g\in H^1(\{\abs{z}<1\}).$
Now the function $1/k'(w)$ cannot have zero angular values on a set of positive measure 
in $\{\abs{z}=1\}$, by the Luzin-Privalov theorem. 
Hence for almost all $t\in \{\abs{z}=1\}$ and each Stoltz angle $S(t,\alpha)$ (with vertex $t$ and bisector passing through the origin)
there exists a sequence $\{w_l\}\in S(t,\alpha)$ such that $w_l\to t$ as $l\to \infty.$
Thus for almost all $t\in \{\abs{z}=1\}$ we have $f(t)k(t)-g(t)=0.$
Now $f(w)=\hat{m}(k(w))k'(w)$ so if $f\equiv 0$ then 
$\frac{1}{2\pi i}\int \frac{dm(\zeta)}{\zeta -z} \equiv 0$ so $m\perp R_1(\partial \omega,\C\setminus\partial\Omega)$. By Mergelyan's theorem
any continuous function on $\partial\Omega$ is uniformly approximable by such rational functions 
so $m\equiv 0.$ The contrapositive of this is that
$m\not\equiv 0$ implies $f\not\equiv 0.$ Since $H^1(\abs{z}<1\})$ is a subset of the Nevanlinna class and
(see Privalov \cite{privalov}, Ch.2, parag. 2.1)
each function in the Nevanlinna class is the ratio of two functions in $H^\infty(\{\abs{z}<1\})$, there exists $f_1,g_1\in H^\infty(\{\abs{z}<1\})$
and $u,v\in H^\infty(\Omega)$ such that for almost all $t\in \{\abs{z}=1\}$
\begin{equation}
\overline{k(t)}=\frac{g(t)}{f(t)}=\frac{g_1(t)}{f_1(t)}=\frac{u(k(t))}{v(k(t))}
\end{equation}
Hence $\Omega\in \mathcal{N}.$ The contrapositive to what we have proved is precisely (i) $(\Leftarrow)$.\\
\\
(ii) $(\Rightarrow)$ We prove the contrapositive. Assume there exists a bounded component $\Omega\in \mathcal{N}$
of $\C\setminus X.$ By (i) there exists a measure $m\not\equiv 0$ with supp$(m)\subset\partial\Omega$
such that $m\perp R_n(\partial\Omega,\overline{\Omega})$ thus $m\perp P_n(X).$ 
Since $\hat{m}(z)=0$ $z\in \C\setminus \overline{\Omega}$
$m\not\equiv 0$ implies (by the argument given in the proof of (i)) 
that there exists a point $z_0\in \Omega$ such that $\hat{m}(z_0)\neq 0.$ Hence
$m \not \perp (\frac{1}{z-z_0}\in A_n(X).$ This implies $A_n(X)\neq P_n(X)$. This proves $(\Rightarrow)$.\\
\\
(ii) $(\Leftarrow)$.
\begin{lemma}\label{carmonrunglem1}
Let $C\subset\C$ be a compact balanced subset. Let $U:=\mbox{int}C$ and let $V$ be a component of $U.$
Then there exists a sequence $\{q_j\}$ of functions analytic in an open neighborhood of $C$ such that:\\
(1) $q_n\to 1$ uniformly on each compact of $U\setminus V.$\\
(2) $q_n\to 0$ uniformly on each compact of $V.$\\
(3) $\norm{q_j}_{C}\leq 1$.
\end{lemma}
\begin{proof}
Let $z_0\in V$ and let $\{J_n\}$ be any sequence of simple closed curves surrounding $C$ and converging to 
some compact set containing the boundary, $B$, of $C.$ Let $\varphi_n$ be the Riemann
map of the enclosed set of $J_n$ onto $\{\abs{z}<1\}$ such that $\varphi_n(z_0)=0$ and $\varphi_n'(z_0)>0.$
Let $\varphi$ be the Riemann map of $V$ onto $\{\abs{z}<1\}$ with $\varphi(z_0)=0$ and $\varphi'(z_0)>0,$
such that $\varphi_n\to \varphi$ uniformly on compacts of $V.$ 
Then the $\varphi_n$ are uniformly bounded on each component $W$ of $U$ 
so if necessary we can replace $\varphi_n$ (keeping the same notation)
by a subsequence so that $\varphi_n\to \varphi$ uniformly on compacts of $W$ for each component $W$ of $U.$
If $W\neq V$ then $\varphi_n(W)$ is disjoint from $\varphi_n(V)$ and since $\varphi_n(V)$ eventually fills $\{\abs{z}<1\}$
the limit $\lim_n \varphi_n$ on $W$ has modulus $1$ everywhere, the limit is a constant $\alpha_W$ of modulus $1$.
Given a finite set $\Gamma$ of components of $U$ not containing $V$, $K>0$ and $\epsilon>0,$
there exists an integer $N>K$ such that 
$\abs{\alpha_w^N-1}<\epsilon$ for all $W\in \Gamma$ (see Bishop \cite{bishoplemmanine}, p.336).
Let $S$ be a compact subset of $U.$
Then there exists $L>0$ and an integer $N$ such that $\abs{(\varphi_n(z))^N-1}<\epsilon$
all $n>L$ and all $z\in S\cap (\cup_{W\in \Gamma} W).$
Since $\sup_{z\in S\cap V}\abs{\varphi(z)} <1$ we can choose $K,n$
sufficiently large such that $\abs{(\varphi_n(z))^N}<\epsilon$ for all $z\in S\cap U$. This proves the existence 
of a sequence $\{q_n\}$ satisfying (1), (2) and that $q_n\to 1$ uniformly on compacts of $\cup_{W\in \Gamma}W.$
Since this holds true for every choice of $\Gamma$ and since $U$ has a countable set of components
the proof follows from a diagonal argument.
\end{proof}
By Runge's theorem, Lemma \label{carmonrunglem1} yields.
\begin{lemma}\label{ggeerr}
Let $X$ be a Carath\'eodory compact, let $\Omega$ be a bounded component of its complement and let 
$U:=(\mbox{int}\hat{X})\setminus\Omega.$
Let $\hat{X}$ denote the polynomially convex hull of $X.$
Then there exists a sequence $\{q_j\}\subset P_1$ such that:\\
(1) $q_j\to 1$ uniformly on each compact of $\Omega.$\\
(2) $q_j\to 0$ uniformly on each compact of $U.$\\
(3) $\norm{q_j}_{\hat{X}}\leq C,$ for an absolute constant $C>0.$
\end{lemma}
Let $m\not\equiv 0$ be a measure on $X$ such that $m\perp P_n(X)$ and let $\Omega$ be a fixed
bounded component of $\C\setminus X$ so $\Omega\notin \mathcal{N}.$
Let $\{q_j\}, q_j\in P_1,$ be a sequence as in Lemma \ref{ggeerr} and set $m_j:=q_jm.$ Let $m_\Omega$
be a limit point of
$\{m_j\}$ in the weak-star topology in the space of measures, so there exists $\{j_l\}$
such that $j_l\to \infty,$ $\mu_{jl}$ converges to $m_\Omega$ in the weak-star topology as $l\to \infty.$
Denoting by 
$\Omega_\infty'$ the unbounded component of $\C\setminus X.$
Since supp$m_\Omega\subset \partial X$ for all $s=0,\ldots,n-1,$ and $z_0\notin\mbox{supp}m$ we have
\begin{multline}
\widehat{\bar{\zeta}^s m_\Omega}(z_0)=\frac{1}{2\pi i}\int \frac{\bar{z}^s dm_\Omega(z)}{z-z_0}=\lim_l \frac{1}{2\pi i}
\int \frac{\bar{z}^s q_{jl}(z)dm(z)}{z-z_0}=\\
\lim_l\left(
\frac{1}{2\pi i}\int \frac{(q_{jl}(z)-q_{jl}(z_0)\bar{z}^s dm(z)}{z-z_0}+\frac{q_{jl}(z_0)}{2\pi i}
\int \frac{\bar{z}^s dm(z)}{z-z_0}
\right)=\\
\lim_l \frac{q_{jl}(z_0)}{2\pi i}\int \frac{\bar{z}^s dm(z)}{z-z_0}=\left\{
\begin{array}{ll}
\widehat{\bar{\zeta}^sm}(z_0) & , z_0\in \Omega\\
0 & ,z_0\in (U\setminus\mbox{supp}(m))\cup \Omega_\infty'
\end{array}
\right.
\end{multline}
It follows that $m_\Omega\perp P_n(X).$ 
This implies that 
$\widehat{\bar{\zeta}^s m_\Omega}(z_0)=0$ on $U\cup \Omega_\infty'.$
We can
assume supp$m\subset \partial X$ (if necessary we can repeat the arguments replacing $m$ with $m_\Omega$) thus
$m\perp P_1(X).$ Since $m$ vanishes on isolated points the function 
\begin{equation}
f_\Omega(w)=\frac{1}{\pi}\int \frac{(\bar{z}-\bar{w})dm_\Omega(z)}{z-w}
\end{equation}
is continuous on $\C$ and zero on $U\cup\Omega_\infty'$ (see Bishop \cite{bishoplemmanine}, Lemma 9)
in particular $f_\Omega=0$ off $\overline{\Omega}.$
Furthermore, $\partial_{\bar{z}}^2 f_\Omega =m_\Omega$ in distribution sense. Hence supp$m_\Omega\subset \partial \Omega$
which implies $m_\Omega \perp R_n(\partial\Omega,\overline{\Omega}).$ Since $\Omega\notin \mathcal{N}$
(1) implies $m_\Omega =0.$ This shows that
for each $\Omega$, $\widehat{\bar{\zeta}^s m_\Omega}(z)=0,$ for $z\in \Omega,$ $s=0,\ldots,n-1.$ Hence $m\perp R_n(X,X).$
By Theorem \ref{vitushkincopythm} $R_n(X,X)=A_n(X)$ so that
$m\perp A_n(X).$
This proves (ii) $(\Leftarrow)$. This completes the proof.
\end{proof}
We believe that the handful of results given here present the central types of questions
and the major techniques of proof that appear. However, 
some notable results on approximation are also presented in Section \ref{noproof}.
See also Mazalov, Paramonov \& Fedorovskiy \cite{mazalovfedorovskiy}, Boivin \& Gauthier \cite{boivingauthier},
Fedorovskiy \cite{fedorovskynevanlinibok}, \cite{fedorovskysurvey0} and \cite{fedorovskysurvey} for further results of this kind.

\section{Further notable results}\label{noproof}
In this section we collect some noteworthy miscellaneous results 
without giving the proofs.\\
\\
As pointed out in Bagby \cite{bagby}, the Harvey-Polking theorem applied to the Cauchy-Riemann operator
has applications to approximation theory for example the following.
\begin{theorem}[Special case of Bagby \cite{bagby}, Theorem 2.1]
	Let $X\subset\Rn$ be a compact set of positive Lebesgue measure. Then the following are equivalent.
	\\
	(a) If $f\in L^2(X)$ and $\partial_{\bar{z}} f=0$ on $\mbox{int}(X)$ then for each $\epsilon>0$ there exists a
	solution $u$ to $\partial_{\bar{z}} f=0$ on a neighborhood of $X$ such that $\norm{u-f}_{L^2(X)} <\epsilon.$\\
	(b) If $G$ and $\Omega$ are open relatively compact subsets of $\Rn$ such that $G\subset\subset \Omega$ then
	\begin{equation}
	L^2-\mbox{Cap}_{\partial_{\bar{z}}}(G\setminus \mbox{int}(X),\Omega)=L^2-\mbox{Cap}_{\partial_{\bar{z}}}(G\setminus X,\Omega)
	\end{equation}
\end{theorem}
Carmona \& Fedorovsky \cite{carmonafedorovsky} prove the following regarding approximation on compacts in $\C$.
\begin{theorem}
	For each integer $q\geq 1$, there exists a compact set $X\subset\C$ such that
	$P_{2n}(X) = C(X) \neq P_n(X).$
\end{theorem}
In the proof the set $X$ is chosen as a union of confocal ellipses.
Carmona \cite{carmona1982} proved the following.
\begin{theorem}
	Let $X\subset \C$ be a compact subset with empty interior, and
	let $g\in C^2(U)$ where $U$ is an open neighborhood of
	$X$. Denote by $R_0(X)$ the algebra of rational functions with poles off $X$ and $Z := \{ x \in X \colon \partial_{\bar{z}} g(z) = 0\}$. Then, $R_0(X) + gR_0(X)$ 
	is uniformly dense in
	$C(X)$ if and only if $R(Z)=C(Z).$
\end{theorem}

Carmona \cite{carmona1985} gives a stronger version of Corollary \ref{carmonacorextrabra} by removing the $C^{q-1}$ condition.
\begin{theorem}[Carmona \cite{carmona1985}]
	Let $q\in \Z_+$ and let $X\subset\C$ be a compact subset such that $\C\setminus X$ has finitely many connected components.
	Any $f\in A_q(X),$ can be uniformly approximated on $X$ by functions that are $q$-analytic on an open neighborhood of $X.$
\end{theorem}
\begin{corollary}
	Let $q\in \Z_+$ and let $X\subset\C$ be a compact subset such that $\C\setminus X$ is connected.
	Then $A_q(X)=P_q(X).$
\end{corollary}

Let $X\subset \C$ be a comactt, let $\mathcal{P}$ denote the space of complex polynomials and let
$\mathcal{R}(X)$ be the space of all rational functions with no poles on $X.$
Set for positive integers $k_1,\ldots,k_m,$ such that $k_1<\cdots <k_m$,
\begin{equation}
\mathcal{P}(\bar{z}^{k_1},\ldots,\bar{z}^{k_m}):=\{ p_0+\bar{z}^{k_1},\ldots,\bar{z}^{k_m} :p_0,\ldots,p_m\in \mathcal{P}\}
\end{equation}
Denote by $A(X,\bar{z}^d)$ the space of functions $f\in C^0(X)$ such that their restriction to $\mbox{int}(X)$ takes the form $f_0+\bar{z}^d f_1$ for functions $f_0,f_1$ holomorphic in 
$\mbox{int}(X).$
Recently Baranov et al. 
\cite{baranov2016} 
introduced so-called {\em $d$-Nevanlinna domains}.
\begin{definition}
	Let $d\geq 1$ be a fixed integer.
	A bounded simply connected domain $\Omega$ is called a
	{\em $d$-Nevanlinna domain}
	if there exists two functions $u,v\in H^\infty(\Omega)$ such that a.e.\ on $\Omega$ we have in the sense of conformal mappings
	\begin{equation}
	\bar{z}^d =\frac{u(z)}{v(z)}
	\end{equation}
	\index{$d$-Nevanlinna domain}
	If $varphi$ is a fixed conformal mapping of the unit disc onto $\omega$ then this means that we have for the angular boundary values the equality
	\begin{equation}
	\overline{\varphi(\xi)}^d =\frac{(u\circ\varphi)(\xi)}{v\circ\varphi(\xi)},\quad \mbox{a.e.\ on }  \{\abs{\xi}=1\}
	\end{equation}
\end{definition}
Baranov et al. 
\cite{baranov2016} proved the following.
\begin{theorem}
	If  $d$ is a positive integer and $X\subset \C$ is a compact subset 
	then
	$A(X,\bar{z}^d)=R(X,\bar{z}^d)$.
	\end{theorem}
	The result is rather quickly obtained for the special case when $\partial X =\widehat{X}$ (where $\widehat{X}$ denotes the union of $X$ 
	and all bounded connected components of $\C\setminus X$).
	\begin{theorem}
	If $\Omega$ is a $d$-Nevanlinna domain, $d\geq 2$, then 
	$C^0(\partial\Omega)\neq R(\partial\Omega,\overline{\Omega},\bar{z}^d).$
		If $U\subset \C$ is a Carath\'eodory domain (i.e.\
		$\partial U =\partial U_\infty$ where $\partial U_\infty$ is the unbounded connected
		component of the set $\hat{\C}\setminus\overline{U}$)
		$d\geq 2$ an integer and if $C^0(\partial\Omega)\neq R(\partial\Omega,\overline{\Omega},\bar{z}^d)$, then
			$U$ is a $d$-Nevanlinna domain.
			\end{theorem}	
				\begin{theorem}
					If $X\subset \C$ is a compact subset such that $\partial X =\widehat{X}$
					and if $d\geq 2$ is an integer then 
					$A(X,\bar{z}^d)=R(X,\bar{z}^d)$ if and only if each bounded connected component of $\C\setminus X$ is not a $d$-Nevanlinna domain.
				\end{theorem}

\chapter[Generalization of Laurent series]{Generalization of Laurent series and Liouville's theorem for polyanalytic functions}\label{liouvillesec}

\section{Krajkiewicz's generalization of countably analytic functions}
We point out that Krajkiewicz \cite{kraj1973} introduced a generalization of countably analytic functions, in terms of analogues of Laurent series, as follows.
\begin{definition}
	Let $z_0$ be a complex number (finite or infinite) and let $0<R<+\infty.$ If $z_0$ is finite, then
	denote by $A(z_0,R)$ the set of all finite $z$ such that $0<\abs{z-z_0}<R.$ 
	If $z_0=\infty$ then $A(z_0,R)=A(\infty,R)$ denotes the set of all finite complex numbers $z$ such that $R<\abs{z}< +\infty.$
	A function $f:A(z_0,R)\to \C$ is called {\em multianalytic} on $A(z_0,R)$ 
	if there exists a nonnegative integer $n$ and a sequence $\{f_j\}_{j=-n}^\infty$
	of holomorphic functions on $A(z_0,R)$ such that for all $z\in A(z_0,R)$ we have
	\begin{equation}
	f(z)=\sum_{j=-n}^\infty (\bar{z} -\bar{z_0})f_j(z)
	\end{equation}
	\index{Multianalytic function}
	
	When $f(z)$ is a function that is multianalytic at $z_0=\infty$
	there exists an $R>0$ and a nonnegative integer $n\geq 0$ together with a sequence $\{f_j\}_{j\geq -n}$ 
	of functions analytic on
	$A(\infty,R)$ for all $j\geq -n$ such that
	\begin{equation}\label{multikrajekv001}
	f(z)=\sum_{j=-n}^\infty \frac{f_j(z)}{\bar{z}^j}
	\end{equation}
	Let the Laurent series of each $f_j(z)$ with center $z_0=\infty$ be
	\begin{equation}\label{multikrajekv002}
	f_j(z)=\sum_{j=-\infty}^\infty a_{j,k} z^k,\quad \abs{z}>R
	\end{equation}
	Define for each $\rho>R$
	\begin{equation}\label{multikrajekv003}
	f(z,\rho):=\sum_{j=-n}^\infty z^j\frac{f_j(z)}{\rho^{2j}},+quad R\leq \abs{z}\leq \rho
	\end{equation}
	This is an analytic function on $R<\abs{z}<\rho$ and continuous on $R<\abs{z}\leq \rho$ and $f(z,\rho)=f(z)$ 
	on $\abs{z}=\rho>R.$
	\begin{proposition}
		Let $f(z)$ be multianalytic on $A(z_0,R)$ and assume $f$ is represented on $A(\infty,R)$
		according to
		\begin{equation}
		f(z)=\sum_{j=-n}^\infty (\bar{z}-\bar{z}_0)^j f_j(z)
		\end{equation}
		for analytic functions $f_j(z)$ on $A(z_0,R).$ Then the functions
		$f_j(z)$ are uniquely determined on $A(z_0,R)$ by $f(z).$
	\end{proposition}
	\begin{proof}
		W.l.o.g.\ we assume $z_0=\infty$. Let $f$ have on $A(\infty,R)$ the representation
		given by Eqn.(\ref{multikrajekv001}) and Eqn.(\ref{multikrajekv002}).
		It suffices to show that if $f\equiv 0$ on $A(\infty,R)$ then
		$f_k\equiv 0$ on $A(\infty,R)$ for all $k\geq -n$. Now $f\equiv 0$ on $A(\infty,R)$ implies
		$f(z,\rho)=0$ for all $\abs{z}=\rho>R.$ Then according to Eqn.(\ref{multikrajekv003}) we have for all $\rho>R$
		and each integer $k$ that
		\begin{equation}
		\sum_{j=-n}^\infty \frac{a_{j,k-j}}{\rho^{2j}}=\frac{1}{2\pi i} \int_{\abs{z}=\rho} \frac{f(z,\rho)dz}{z^{1+k}}=0
		\end{equation}
		This implies $a_{j,k-j}=0$
		for $j\geq -n$ and every integer $k.$
		Hence $a_{j,l}=0$ for all $j\geq -n$ and every integer $l.$ This implies that $f_j\equiv 0$ on $A(\infty,R)$ for each $j\geq -n.$ 
		This completes the proof.
	\end{proof}
	
	Krajkiweicz \cite{kraj1973} 
	defines
	the order $d_j=d(f_j)$ at $z_0$ as follows:
	If $f_j\equiv 0$ then $d_j=-\infty.$ If $z_0$ is an isolated singularity of $f_j$ of 
	then $d=+\infty.$ Suppose $z_0$ is not an isolated singularity of $f_j$ of 
	$f_j\not\equiv 0$. If $z_0=\infty$ then there is a unique integer $\nu$ such that $f_j(z)/z^\nu$ is holomorphic
	and nonzero at $z_0.$ One then defines $d_j=\nu.$ If $z_0\neq\infty$ 
	there is a unique integer $\nu$ such that $(z-z_0)^\nu f_j(z)$ is holomorphic
	and nonzero at $z_0.$ One then defines $d_j=\nu.$
	We define $d=d(f)$ to be the least upper bound of $d_j-j,$ for $j\geq -n.$ So $d$ is either an integer or $\pm \infty.$
	The point $z_0$ is called an essential isolated singularity
	of $f$ (or a pole of $f$) of order $d$ if and only if $d= +\infty.$  
\end{definition}
Note that if $z_0$ is an essential isolated singularity of $f$ 
then for each nonzero complex number $\alpha$ the value $f(z_0)-\alpha$ is an exceptional value for $f$ at $z_0$. We mention (although it is not directly related to Liouville type theorems, which are focus of this chapter, it is worth noting as it is currently one of the few results appearing in the literature on multianalytic functions) the following result proved by Krajkiweicz \cite{kraj1973} for multianalytic functions.
\begin{theorem}[Krajkiweicz \cite{kraj1973}]
	Let $f$ be a function which is multianalytic at a
	point $z_0$, finite or infinite. If $z_0$ is an essential isolated singularity
	of $f$, then $f$ can have at most one exceptional value at $z_0$.
\end{theorem}

\subsection{Application for obtaining a classical type of Liouville theorem for polyanalytic functions}
As an application of multianalytic functions we give the following.
\begin{lemma}
	Let $f(z)$ be multianalytic on $A(\infty,R)$ (in particular it is multianalytic at $\infty$
	so Eqn.(\ref{multikrajekv001}) applies) and let $\lambda\in (0,1]$ be fixed and
	let 
	$\{\rho_j\}_{j\in \N}$ be a sequence of numbers $\rho_j>R$ diverging to $+\infty.$ Suppose there exists a constant $K>0$ together with an integer $t$ such that
	\begin{equation}
	\abs{f(\rho z,\rho)}\leq K\rho^t,\quad \abs{z}=\lambda,\rho=\rho_j
	\end{equation} 
	Then $z_0=\infty$ is a pole of $f$ of order $d\leq t.$
\end{lemma}
\begin{proof}
	By Eqn.(\ref{multikrajekv001}) there exists
	a nonnegative integer $n\geq 0$ together with a sequence $\{f_j\}_{j\geq -n}$ 
	of functions analytic on
	$A(\infty,R)$ for all $j\geq -n$ such that
	\begin{equation}
	f(z)=\sum_{j=-n}^\infty \frac{f_j(z)}{\bar{z}^j}
	\end{equation}
	and we suppose 
	the Laurent series of each $f_j(z)$ with center $z_0=\infty$ is given by
	\begin{equation}
	f_j(z)=\sum_{j=-\infty}^\infty a_{j,k} z^k
	\end{equation}
	We have for $\rho>R$ and each integer $k$ that
	\begin{equation}
	\sum_{j=-n}^\infty a_{j,k-j}\rho^{k-2j} =\frac{1}{2\pi i} \int_{\abs{z}=\lambda} \frac{f(\rho z,\rho)dz}{z^{1+k}}
	\end{equation}
	This yields for each $\rho=\rho_m$ and each integer $k$ that
	\begin{equation}\label{multikrajekv00aa}
	\abs{ \sum_{j=-n}^\infty \frac{a_{j,k-j}}{\rho^{2j+t-k}}}\leq \frac{K}{\lambda^t}
	\end{equation}
	In the infinite sum in Eqn.(\ref{multikrajekv00aa}) the coefficients of positive powers of $\rho$ must vanish, hence $a_{j,k-j}=0$ for $k\geq -n$ and $2k+t-j\leq 1.$
	This implies that $a_{j,l}=0$ for
	$k\geq -n$ and $l\geq j+t+1.$ Hence $z=\infty$ is a removable singularity or a pole of $f_j$ of order $d\leq t.$ 
	This completes the proof.
\end{proof}
\begin{corollary}
	Let $f(z)$ be multianalytic on $A(\infty,R)$ and
	let 
	$\{\rho_j\}_{j\in \N}$ be a sequence of numbers $\rho_j>R$ diverging to $+\infty.$ 
	Suppose there exists a constant $K>0$ together with an integer $t$ such that
	\begin{equation}
	\abs{f(z)}\leq K\rho^t,\quad \abs{z}=\lambda,\rho=\rho_j
	\end{equation} 
	Then $z_0=\infty$ is a pole of $f$ of order $d\leq t.$
\end{corollary}
As a special case of the corollary we obtain the following is a generalization of the more standard Liouville theorem.
\begin{theorem}
	If $f(z)$ is a polyentire function such that $\abs{f(z)}$ is bounded then $f(z)$ is a complex constant.
\end{theorem}

\section{Liouville's theorem for polyharmonic functions}
Recall that a real valued function on a domain in the plane is polyharmonic of order $q\in \Z_+$ if and only if it is the real part of
a $q$-analytic function.
Recently generalizations of Liouville's theorem for polyharmonic functions have been
established, see e.g.\ Nakai \& Tada \cite{nakai}, \cite{nakai0},
Armitage \cite{armitage} and Futamara \cite{futamara}.
In particular, we have the following due to Nakai \& Tada \cite{nakai}.
\begin{theorem}
	Let $u$ be a polyharmonic function of order $m$ on $\Rn$. Let $s>2(m-1).$ Then $u$ is a polyhamonic polynomial
	of degree less than $s$ if and only if there exists an increasing divergent sequence $\{r_j\}_{j\in \Z_+}$,
	$r_j>0, j\in \Z_+,$ such that
	\begin{equation}\label{cond3takai}
	\liminf_{i\uparrow\infty} \left(\min_{\abs{x}=r_i} \frac{u(x)}{\abs{x}^s}\right)\geq 0
	\end{equation} 
\end{theorem}
\begin{proof}
	Introduce the polar coordinate representation $x=r\xi$ for points $x\in\Rn$ where $\xi=x/\abs{z}\in S^{n-1}$
	for $x\neq 0$ and $\xi=(1,0,\ldots,0)\in S^{n-1}$ for $x=0$.
	Let $\{S_{jk},j=1,\ldots,N(k)\}$ be a fixed orthonormal basis for the subspace of all spherical harmonics of degree $k$ of $L^2(S^{n-1},d\sigma)$ where $d\sigma$ is the area element on 
	$S^{n-1}$ such that $\{S_{jk}\}_{k\in\N}$ is a complete orthonormal system in $L^2(S^{n-1},d\sigma).$
	By the so-called addition theorem (see Axler, Boudron \& Ramsey \cite{axlerbok}, 5.11)
	we have
	\begin{equation}
	\sum_{j=1}^{N(k)} S_{jk}(\xi)^2=\frac{N(k)}{\sigma_n}
	\end{equation}
	where $\sigma_n$ is the surface area $\sigma(S^{n-1})$ of $S^{n-1}.$
	Here we have (see Axler, Boudron \& Ramsey \cite{axlerbok}, 5.17)
	$N(0)=1$ and
	\begin{equation}
	N(k)=(2k+n-2)\frac{\Gamma(k+n-2)}{\Gamma(k+1)\Gamma(n-1)},\quad k\in \Z_+
	\end{equation} 
	This implies that for each $\xi\in S^{n-1}$
	\begin{equation}
	\abs{S_{jk(\xi)}}\leq \sqrt{\frac{N(k)}{\sigma_n}}=:A_k,\quad j=1,\ldots,N(k),k\in \N
	\end{equation} 
	For each harmonic $h(r\xi)$ we have the expansion
	\begin{equation}
	h(r\xi)=\sum_{k=0}^\infty\left(\sum_{j=1}^{N(k)} a_{jk}S_{kj}(\xi)\right)r^k
	\end{equation}
	for constants $a_{jk0},$ $j=1,\ldots,N,k\in \N$ and the series converges uniformly in $\xi\in S^{n-1}$ for any fixed $0<r<\infty.$
	Since $2m-2<s$ there exists a unique integer $n$ such that 
	\begin{equation}\label{cond3takai22}
	2m-2\leq n<s\leq n+1
	\end{equation}
	It suffices to show that
	given Eqn.(\ref{cond3takai}) any $m$-harmonic $u$ is a polynomial of degree at most $n$.
	Now by Eqn.(\ref{cond3takai})there exists, for each $\epsilon>0$, a positive integer $i(\epsilon)$
	such that for each $i\geq i(\epsilon)$ we have $u(x)/\abs{x}^s\geq -\epsilon$ for $\abs{x}=r_i$ and $r_i>1$
	and since by Eqn.(\ref{cond3takai22}) $\abs{x}^s\leq \abs{x}^{n+1}$ for $\abs{x}\geq 1$, 
	we have for each $i\geq i(\epsilon)$
	\begin{equation}\label{cond3takai23}
	u(x)\geq -\epsilon \abs{x}^{n+1},\quad \abs{x}=r_i 
	\end{equation}
	Let the Almansi expansion of $u$ be
	\begin{equation}
	u(x)=\sum_{p=1}^m \abs{x}^{2p-2} h_p(x)
	\end{equation}
	for harmonic $h_p$, $p=1,\ldots,m$ which in the coordinates $x=r\xi$ take the form
	\begin{equation}
	h_p(r\xi)=\sum_{k=0}^\infty \left(\sum_{j=1}^{N(k)} a_{pkj}S_{kj}(\xi)\right)r^k,\quad p=1,\ldots,m
	\end{equation}
	for constants $a_{pjk},$ and the series converges uniformly in $\xi\in S^{n-1}$ for any fixed $0<r<\infty.$ Hence
	\begin{equation}\label{takai2punkt4}
	u(r\xi)=\sum_{p=1}^m r^{2p-2} \left(\sum_{k=0}^\infty \left(\sum_{j=1}^{N(k)} a_{pkj}S_{kj}(\xi)\right)r^k\right)
	\end{equation} 
	converges uniformly in $\xi\in S^{n-1}$ for $0<r<\infty.$
	By Eqn.(\ref{cond3takai23}) this gives for any $\xi\in S^{n-1}$
	\begin{equation}\label{takai25}
	\sum_{p=1}^m r^{2p-2} 
	\left(\sum_{k=0}^\infty \left(\sum_{j=1}^{N(k)} a_{pkj}S_{kj}(\xi)\right)r^k\right)
	+\epsilon r_i^{n+1}\geq 0,\quad i\geq i(\epsilon)
	\end{equation} 
	\begin{lemma}\label{takailem26}
		Given Eqn.(\ref{takai25}) it holds true that
		\begin{equation}
		a_{pkj}=0,\quad (2p-2)+k\geq n+1
		\end{equation}
	\end{lemma}
	\begin{proof}
		First we consider the case $(2p-2)+k>n+1.$
		Assume (in order to reach a contradiction) that the lemma is false for this case.
		Then there exists $1\leq p\leq m,$
		$0\leq k<\infty$, $1\leq j\leq N(k)$ such that
		$(2p-2)+k>n-1$ but $a_{pkj}\neq 0.$ 
		Let $p_1$ be the greatest $p$ satisfying $a_{pkj}\neq 0.$
		Choose $\eta\in \{\pm 1\}$ such that
		\begin{equation}\label{takai7}
		\eta a_{p_1 kj}>0
		\end{equation}
		Multiplying both sides of Eqn.(\ref{takai25}) by $A_k-\eta S_{jk}(\xi)\geq 0$ and integrating with respect to$d\sigma$ on $S^{n-1}$ yields
		\begin{equation}\label{takai8}
		\underbrace{\sigma_n A_k\left(\sigma_d^{-\frac{1}{2}}\sum_{p=1}^m a_{p01} r_i^{2p-2} +\epsilon r_i^{n+1}\right)}_{=:P(r)}
		-\underbrace{\sum_{p=1}^{p_1}\eta a_{pkj} r_i^{(2p-2)+k} }_{=:Q(r)}\geq 0
		\end{equation}
		for $i\geq i(\epsilon).$
		Now $P(r)$ is a polynomial of degree $n+1$ whereas $Q(r)$ is a polynomial of degree $(2p_1-2)+k>n+1$
		by Eqn.(\ref{takai7}). This implies that $\mbox{deg}P<\mbox{deg}Q$, the leading coefficient of $Q$ 
		is strictly positive by Eqn.(\ref{takai7})
		whereas Eqn.(\ref{takai8}) yields $P(r_i)\geq Q(r_i)$ for $i\geq i(\epsilon)$, which is a contradiction.
		This proves the lemma for the last case $(2p-2)+k>n+1.$
		\\
		\\
		Next we prove $a_{pkj}=0$ for the final case $(2p-2)+k=n+1.$ Assume (in order to reach a contradiction) that this 
		does not hold true.
		Then there exists $1\leq p_1\leq m,$
		$0\leq k<\infty$, $1\leq j\leq N(k)$ such that
		$(2p-2)+k=n-1$ but $a_{pkj}\neq 0.$ Note that $p_1=m$ or $p_1<m.$ 
		If $p_1<m$ then $a_{pj}=0$ for $p_1<p\leq m$ by what we have already done since $(2p-2)+k>(2p_1 -2)+k=n+1.$
		Choose $\eta\in \{\pm 1\}$ such that
		\begin{equation}\label{takai9}
		\eta a_{p_1 kj}-\sigma_n A_k\epsilon >0
		\end{equation}
		Multiplying both sides of Eqn.(\ref{takai25}) by $A_k-\eta S_{jk}(\xi)\geq 0$ and 
		integrating with respect to $d\sigma$ on $S^{n-1}$ yields
		\begin{equation}\label{takai10}
		\underbrace{\sigma_n^{\frac{1}{2}} A_k \sum{p=1}^m a_{p01} r_i^{2p-2}}_{=:P(r)}
		-\underbrace{(\eta a_{p_1 kj}-\sigma_n A_k\epsilon)r^{n+1}_i+ \sum_{p=1}^{p_1-1}\eta a_{pkj} r_i^{(2p-2)+k} }_{=:Q(r)}\geq 0
		\end{equation}
		Now $P(r)$ is a polynomial of degree at most $2m-2<n+1$ whereas $Q(r)$ is a polynomial of degree $n+1$.
		This implies that $\mbox{deg}P<\mbox{deg}Q$, the leading coefficient of $Q$ is strictly positive by 
		Eqn.(\ref{takai9})
		whereas Eqn.(\ref{takai10}) yields $P(r_i)\geq Q(r_i)$ for $i\geq i(\epsilon)$, which is a contradiction.
		This proves the lemma for the last case $(2p-2)+k=n+1.$
		This proves Lemma \ref{takailem26}.
	\end{proof}
	By Lemma \ref{takailem26} we may rewrite Eqn.(\ref{takai2punkt4}) for each $\xi\in S^{n-1}$ as
	\begin{equation}\label{takai11}
	u(r\xi)=\sum_{p=1}^m \left(\sum_{(2p-2)+k\leq n} \left(\sum_{j=1}^{N(k)} a_{pkj}r^{2p-2}r^k S_{kj}(\xi)\right)\right)
	\end{equation}
	Since $S_{kj}(x):=r^kS_{kj}(\xi)$ ($x=r\xi$) is a homogeneous 
	harmonic polynomial in $x$ of degree $k$, Eqn.(\ref{takai11}) gives that
	\begin{equation}
	u(x)=\sum_{p=1}^m \left(\sum_{(2p-2)+k\leq n} \left(\sum_{j=1}^{N(k)} a_{pkj}\abs{x}^{2p-2} S_{kj}(x)\right)\right)
	\end{equation}
	is a polynomial in $x$ of degree at most $n$. This completes the proof.
\end{proof}
This can be reformulated as follows.
\begin{corollary}[See e.g.\ Armitage \cite{armitage}, Thm. A, (iii)]\label{armitageliouvillethm}
	Let $u$ be a polyharmonic function of order $m$ on $\Rn$. Let $s>2(m-1).$ Then $u$ is a polyhamonic polynomial
	of degree less than $s$ if and only if
	\begin{equation}
	\liminf_{r\to \infty} r^{-s} \max\{ u^+(z):\abs{z}=r\}=0
	\end{equation}
\end{corollary}

\begin{corollary}\label{abbeliouvillenollanthm}
	If $f(z,\bar{z})$, $z=x+iy,$ is a function on $\C$ such that $u(x,y)=\re f$ and $v(x,y):=\im f$ are $q$-harmonic, some $q\in \Z_+,$ 
	such that for a constant $C>0$ and some $s\geq 2(q-1)$ we have for sufficiently large $\abs{z}$ (say $\abs{z}>R_0$ for some $R_0>0$)
	\begin{equation}
	\frac{\abs{f(z)}}{\abs{z}^s}\leq C
	\end{equation}
	then
	$f(z,\bar{z})$ is a polynomial of joint order $2(q-1)$ i.e.\ 
	$f(z,\bar{z})=\sum_{\abs{\beta}\leq 2(q-1)} z^{\beta_1}\bar{z}^{\beta_2}$.
\end{corollary}
\begin{proof}
	Note that the conditions imply that
	\begin{equation}
	\lim_{\abs{z}\to \infty} \frac{\abs{f(z)}}{\abs{z}^{2q-1}}=0
	\end{equation}
	Let $s\geq 2q-2.$
	Let $z=x+iy$ and $u(x,y):=\re f(x,y),$ $v(x,y):=\im f(x,y)$. Then
	$\abs{f(z)}\leq C\abs{z}^s$ for all sufficiently large $\abs{z}$ implies
	\begin{equation}
	\abs{f}^2\leq C^2 \abs{z}^{2s} \Leftrightarrow \frac{\abs{u}^2+\abs{v}^2}{\abs{z}^{2s}}\leq C^2
	\Leftrightarrow \abs{\frac{\abs{u}}{\abs{z}^{s}}}^2+ \abs{\frac{\abs{v}}{\abs{z}^{s}}}^2\leq C^2
	\end{equation}
	Hence
	\begin{equation}
	\frac{\abs{u}}{\abs{z}^{s}}\leq C,\quad  \frac{\abs{v}}{\abs{z}^{s}}\leq C
	\end{equation}
	In particular, 
	\begin{equation}
	\liminf_{r\to \infty} r^{-s-1} \max\{ u^+(z):\abs{z}=r\}=0,\quad \liminf_{r\to \infty} r^{-s-1} \max\{ v^+(z):\abs{z}=r\}=0
	\end{equation}
	By Corollary \ref{armitageliouvillethm} this implies that
	$u(x,y)$ and $v(x,y)$ are both polynomials of degree $2(q-1)$ jointly with respect to $x,y$ (i.e.\
	the monomial $x^{\alpha_1}y^{\alpha_2}$ appearing in the representation satisfies
	$\abs{\alpha}=\alpha_1+\alpha_2 =2(q-1).$
	Hence
	\begin{equation}
	f(z,\bar{z})=(u+iv)\left(\frac{z+\bar{z}}{2},\frac{z-\bar{z}}{2i}\right)
	\end{equation}
	is a polynomial that has degree of at most $2(q-1)$ with respect to $\bar{z}$ and $z.$
\end{proof}	
As a consequence we obtain the following. 
\begin{theorem}
	Let $f(z)$ be a polyentire function of order $q$, $q\in\Z_+.$
	If there exists a constant $C>0$ such that for sufficiently large $\abs{z}$ (say $\abs{z}>R_0$ for some $R_0>0$) we have $\abs{f(z)}\leq C\abs{z}^{2(q-1)}$,
	then $f(z)$ is a polyanalytic polynomial whose analytic components $a_0(z),\ldots,a_{2(q-1)}(z)$ 
	satisfy that each $a_j(z)$ is a polynomial of degree $\leq 2q-j-1.$
	\begin{equation}
	f(z)=\sum_{j=0}^{2q-1} a_j(z)\bar{z}^j
	\end{equation}
	then $a_j(z)$ is a polynomial of degree $\leq 2q-j-1.$ (In particular if $q=1$, i.e.\ $f(z)$ is holomorphic, then $f(z)$ must be constant).
\end{theorem}
\begin{proof}	
	A polyanalytic function of order $q$ satisfies that $u:=\re f$ and $v=\im f$ are both
	$q$-harmonic functions.
	By Theorem \ref{abbeliouvillenollanthm} we thus have 	
	\begin{equation}
	f(z,\bar{z})=(u+iv)\left(\frac{z+\bar{z}}{2},\frac{z-\bar{z}}{2i}\right)
	\end{equation}
	is a polynomial that has degree of at most $2(q-1)$ with respect to $\bar{z}$ and $z.$
	Hence there are constants $c_\beta\in \C$ such that
	\begin{equation}
	f(z)=\sum_{\abs{\beta}\leq 2(q-1)} c_\beta (z+\bar{z})^{\beta_1}(z-\bar{z})^{\beta_2} 
	\end{equation}
	Hence if $a_0(z),\ldots,a_{2(q-1)}(z)$ are the analytic components of $f(z)$ i.e.\
	\begin{equation}
	f(z)=\sum_{j=0}^{2q-1} a_j(z)\bar{z}^j
	\end{equation}
	then $a_j(z)$ is a polynomial of degree $\leq 2q-j-1.$
	This completes the proof.
\end{proof}

We have the following Liouville type theorem for polyanalytic functions (it appears without proof in Balk \& Zuev \cite{balkzuev} (2), p.212, see also Balk \cite{balkliouville}).
\begin{theorem}
	If $f(z)$ is a polyanalytic function on $\C$ such that for some $s\in \N$ and $C>0$ we have
	for all $z\in \C$ with sufficiently large $\abs{z}$ (say $\abs{z}>R_0$ for some $R_0>0$)
	\begin{equation}\label{hjasqwtrekv}
	\frac{\abs{f(z)}}{\abs{z}^s}\leq C
	\end{equation}
	Then $f(z)$ is a polynomial of order $s$ (irrespective of the order of polyanalyticity of $f$).
\end{theorem}
\begin{proof}
	Let $q\in \Z_+$ and suppose $f(z)$ is $q$-analytic and let $z=x+iy$ be the complex coordinate in $\C.$
	If $s\geq 2(q-1)$ then $f(z)$ is a polynomial in $x,y$ of degree $2(q-1).$
	by Corollary \ref{abbeliouvillenollanthm}.
	If $s<2(q-1)$ 
	then
	for sufficiently large $\abs{z}$
	\begin{equation}
	\frac{\abs{f}}{\abs{z}^s}\leq \frac{\abs{f}}{\abs{z}^{2(q-1)}}
	\end{equation}
	so if $f(z)$ satisfies 
	\begin{equation}
	\frac{\abs{f(z)}}{\abs{z}^s}\leq C
	\end{equation}
	thus again $f$ satisfies the conditions of Theorem \ref{abbeliouvillenollanthm}.
	Hence for all $s\in \N$ we have that the conditions of the corollary implies that
	\begin{equation}
	f(z)=P(x,y)=\sum_{\abs{\beta}\leq 2(q-1)} c_\beta x^\beta_1y\beta_2
	\end{equation}
	for complex coefficients $c_\beta$. In polar coordinates $z=r\exp(i\theta)$ this is
	\begin{equation}
	f(r,\theta)=\sum_{\abs{\beta}\leq 2(q-1)} c_\beta r^{\abs{\beta}} \cos^{\beta_1}\theta \sin^{\beta_2}\theta
	\end{equation}
	and the condition of Eqn.(\ref{hjasqwtrekv}) is 
	\begin{equation}
	\abs{\sum_{\abs{\beta}\leq s}\frac{r^{\abs{\beta}}}{r^s}  c_\beta  \cos^{\beta_1} \theta \sin^{\beta_2} \theta} \leq C
	\end{equation}
	for sufficiently large $r>0$ and all $\theta\in \R.$
	This implies that
	\begin{equation}
	f(r,\theta)=\sum_{\abs{\beta}\leq s} c_\beta r^{\abs{\beta}} \cos^{\beta_1}\theta \sin^{\beta_2}\theta
	\end{equation}
	This completes the proof.
\end{proof}
In particular, if a polyanalytic function $f$ is essentially bounded (i.e. $f\in L^\infty(\C)$) then 
$f$ is constant, so if also $f\in L^2(\C)$ then $f\equiv 0.$ 

\section[Application of generalized Laurent series]{Application of generalized Laurent series with respect to a real variable}
Balk \cite{balk58liouville} proved the following which generalized a theorem of Hadamard.
\begin{proposition}\label{balkhadamard}
	Let $f(z)=\sum_{n=-\infty}^\infty p_n(r)\exp(in\theta)$ where $p_n(r)$ polynomials in $r$ and $r^{-1}$.
	Suppose that the series
	\begin{equation}
	\sum_{n=0}^\infty p_n(r)\exp(in\theta),\quad \sum_{n=0}^\infty p_{-n}(r)\exp(-in\theta)
	\end{equation}
	converge uniformly on $\{\abs{z}=r\}$. If $J(r)=\sup\im f(z)$ on $\{\abs{z}=r\}$
	then for $n\in \Z_+\setminus \{0\}$
	\begin{equation}
	2\im p_0(r)+\abs{p_n(r)-\overline{p_{-n}(r)}}\leq \max \{4J(r),0\}
	\end{equation}
\end{proposition}
\begin{proof}
	We have
	\begin{multline}
	\im f(r\exp(i\theta))=\im p_0(r)+\sum_{n=1}^\infty (\im p_n(r) \cos n\theta+\re p_n(r)\sin n\theta)+\\
	\sum_{n=1}^\infty (\im p_{-n}(r) \cos n\theta+\re p_n(r)\sin n\theta)
	\end{multline}
	Furthermore
	\begin{equation}\label{balkliouvillekev01}
	\frac{1}{2\pi} \int_0^{2\pi} \im f(z)d\theta =\im p_0(r)
	\end{equation}
	\begin{equation}\label{balkliouvillekev02}
	\frac{1}{2\pi} \int_0^{2\pi} \im f(z) \sin n\theta d\theta =\re p_n(r)-\re p_{-n}(r)
	\end{equation}
	\begin{equation}\label{balkliouvillekev03}
	\frac{1}{2\pi} \int_0^{2\pi} \im f(z) \cos n\theta d\theta =\im p_n(r)+\im p_{-n}(r)
	\end{equation}
	Also
	\begin{equation}\label{balkliouvillekev3}
	\frac{1}{2\pi} \int_0^{2\pi} \im f(z)(\sin n\theta+i \cos n\theta) d\theta = p_n(r)-\overline{ p_{-n}}(r)
	\end{equation}
	By Eqn.(\ref{balkliouvillekev01})-Eqn.(\ref{balkliouvillekev3}) we have
	\begin{equation}
	\abs{p_n(r)-\overline{p_{-n}}(r)}+\im p_0(r)\leq
	\frac{1}{2\pi} \int_0^{2\pi} (\abs{\im f(z)}+\im f(z)) d\theta\leq \max\{4J(r)\}
	\end{equation}
	This completes the proof of Proposition \ref{balkhadamard}.
\end{proof}
Balk \cite{ca1}, p.211, states without proof that a consequence of this is 
the following (cf. Theorem \ref{abbeliouvillenollanthm}).
\begin{theorem}
	If $f(z)$ is a polyanalytic function of order $n$ on $\C$ such that for some $s\in \N$ and $C>0$ we have
	for all $z\in \C$ with sufficiently large $\abs{z}$ (say $\abs{z}>R_0$ for some $R_0>0$)
	\begin{equation}\label{hjasqwtrekv}
	\abs{\im f(z)}\leq C \abs{z}^s
	\end{equation}
	then $f(z)=R(z,\bar{z})+S(z,\bar{z})$ where $R,S$ are polynomials, $\im R\equiv 0$
	$\mbox{deg}_z R\leq n-1,$ $\mbox{deg}_{\bar{z}} R\leq n-1,$ $\mbox{deg}_{z,\bar{z}} S\leq s.$
\end{theorem}
Utilizing Proposition \ref{balkhadamard}, Balk \cite{balk58liouville} proves the following.
\begin{theorem}\label{balkliouvillethm1}
	Let $P(x,y)$ be a real polynomial of degree $s$ and let $D^+$ ($D^-$ respectively), denote the set of points in the $(x,y)$ 
	plane where $P(x,y)>0$ ($P(x,y)<0$ respectively). If $\Phi(z)$ 
	is a univalent analytic function in the $z=x+iy$ plane such that $\im \Phi(z)$  is bounded from above (below respectively)
	in $D^+$ ($D^-$ respectively), then every isolated singularity of $\Phi(z)$ is a pole of order $s$ 
	and lies on the curve $P(x,y)=0$. In addition any pole of multiplicity $>1$  
	can lie only in a singular point of $P(x,y)=0$. 
	\end{theorem}
\begin{proof}
	We shall need the following lemma.
	\begin{lemma}\label{balkliouvillelem3}
		Let $P(x,y)$ be a real valued polynomial of degree $s$ with respect to $x$ and $y$ jointly,
		and let $z=x+iy=r\exp(i\theta).$ Then
		\begin{equation}\label{balkliouvilleekv4}
		P(x,y)=\sum_{n=-s}^s q_n(r)\exp(in\theta)
		\end{equation}
		where $q(r)$ is a polynomial of degree not higher than $s$.
		Furthermore, for $n\in \Z_+\setminus \{0\}$ the polynomials $q_n(r)$ and $q_{-n}(r)$ are divisble 
		by $r^n$ and satisfy
		\begin{equation}\label{balkliouvilleekv5}
		q_{-n}(r)=\overline{q_n}(r),\quad n=\Z_+\setminus \{0\}
		\end{equation}
	\end{lemma}
	\begin{proof}
		Let
		\begin{equation}
		P(x,y)=\sum_{m=0}^s \sum_{k=0}^m A_k^{(m)} x^k y^{-k}
		\end{equation}
		where $\{A^{(m)}_k\}$ are constants.
		For $x=\frac{1}{2}r(\exp(i\theta)+\exp(-i\theta)),$
		$y=\frac{1}{2i}r(\exp(i\theta)-\exp(-i\theta))$ we have
		\begin{multline}
		P(x,y)=\sum_{m=0}^s r^m \sum_{k=0}^m B_k^{(m)} (\exp(i\theta)+\exp(-i\theta))^k
		(\exp(i\theta)-\exp(-i\theta))^{m-k}\\
		=\sum_{m=0}^s r^m \sum_{k=0}^m C_k^{(m)} \exp(ik\theta)
		\end{multline}
		for constants $\{B^{(m)}_k\}$ and $\{C^{(m)}_k\}$.
		Collecting the terms containing $\exp(in\theta)$ we obtain a representation for $P(x,y)$ in the form of 
		Eqn.(\ref{balkliouvilleekv4}) where the degree of each $q_n(r)$ cannot be higher than $s$
		(however at least one of the polynomials $q_n(r)$ has degree at least $s$ if $P(x,y)$ has degree at least $s$).
		This proves Eqn.(\ref{balkliouvilleekv4}).
		\\
		\\
		For $n\in \N,$ we collect the terms containing $\exp(in\theta)$ 
		from the expression \\
		$\sum_{m=0}^s r^m \sum_{k=0}^m C_k^{(m)} \exp(ik\theta)$
		yields
		\begin{equation}\label{balkliouvilleekv6}
		r^m \sum_{k=-m}^m C_k^{(m)} \exp(ik\theta)
		\end{equation}
		Obviously, for $m<n$, $\exp(in\theta)$ does not appear in the sum.
		Thus in the ratio of $q_n(r)$ with $\exp(in\theta)$, the members containing $r^m$, for $m<n$, 
		do not appears. Bu then $r^n$ divides $q_n(r)$. Similarly, we deduce that
		$r^n$ divides $q_{-n}(r).$ 
		If $P(x,y)$ is a real valued $\im P(x,y)\equiv 0$ means
		\begin{equation}
		\im q_0+\sum_{n=0}^s ((q_n+\im q_{-n})\cos n\theta +(\re q_n-re q_{-n})\sin n\theta)
		\equiv 0
		\end{equation}
		Thus $\im q_n+\im q_{-n}=0,$ $\re q_n-\re q_{-n}=0$ i.e.\
		$q_{-n}(r)=\overline{q_n}(r).$
		This completes the proof of Lemma \ref{balkliouvillelem3}.
	\end{proof}
	Let $\im \Phi(z)\leq A_1$ on $D^+$ and
	$\Phi \geq A_2$ on $D^-$ for constants $A_1,A_2$ and we assume $A_1>A_2.$
	Let $a=a_1+ia_2$ be an arbitrary isolated singularity of $\Phi(z).$
		Suppose $P(a_1,a_2)> 0$. 
	Then $\im \Phi(z)\leq A_1$ on $\delta.$
	But an analytic function cannot be of bounded modulus in a 
	neighborhood of one of its poles (and this is true also at the point at infinity). 
	Furthermore, the so-called Sokhotsky-Weierstrass theorem states that if $\Phi$ is a holomorphic function
	on a punctured neighborhood of a point $a$ then either $\lim_{z\to a}\Phi(z)$ exists in the extended plane or the cluster set of $\Phi$ at $a$ is 
	all of the extended complex plane (i.e.\ $a$ is an essential singularity).
	This  
	implies that $\im \Phi(z)\leq A_1$ on $\delta$ is not possible.
	Similarly, we obtain a contradiction by assuming that $P(a_1,a_2)< 0$.
	We therefore conclude that $(a_1,a_2)$ belongs to $\{P(x,y)=0\}.$\\
	\\
	Let 
	$f(z)=\Phi(z)-iA_2,$ $c=A_1-A_2$, where $c>0$ since $A_1>A_2$. Then $\im f\leq c$ on $D^+$ and $\im f
	\geq 0$ on $D_-.$ Then everywhere, except singular points of $\Phi$, we have
	$P(x,y)\im f(z)\leq c\abs{P(x,y)}$
	thus
	\begin{equation}\label{balkliouvilleekv7}
	\im (P(x,y)f(z))\leq c\abs{P(x,y)}
	\end{equation}
	Let $a$ be an isolated singularity of $f(f)$. In a neighborhood of $a$ the function $f(z)$ can be represented as
	$f(z)=\sum_{n=-\infty}^\infty c_n(z-a)^n$ for constants $\{c_n\}_{n\in \Z}.$
	W.l.o.g.\ we assume $a=0$ so that
	$f(z)=\sum_{n=-\infty}^\infty c_nz^n =\sum_{n=-\infty}^\infty c_n r^n\exp(in\theta)$
	By Eqn.(\ref{balkliouvilleekv4}) we have
	\begin{equation}
	P(x,y)f(z)=\sum_{n=-\infty}^\infty p_n(r) \exp(in\theta)
	\end{equation}
	where
	\begin{equation}\label{balkliouvilleekv8}
	p_{n}(r)=\sum_{m=-s}^s c_{-m+n} r^{-m+n}
	\end{equation}
	By Lemma \ref{balkliouvillelem3} we have for $m=0,\ldots,s$ that $q_m(r)/r^m$
	is a polynomial in $r,$ of degree $s-m,$ we denote it by $Q_m.$ Since $P(x,y)\not\equiv 0$
	we know that $Q_m\not\equiv 0.$ So there exists an integer $\nu$ such that $r^{\nu+1}$
	divides each of the polynomials $Q_0,\ldots,Q_s.$ By Lemma \ref{balkliouvillelem3}
	it follows that $q_n$ and $q_{-n}$ are both divisible by $r^{m+\nu}.$
	By Eqn.(\ref{balkliouvilleekv8}) we have for $n>0$
	\begin{equation}\label{balkliouvilleekv9}
	p_{-n}(r)=\frac{1}{r^{n-\nu -2}}\sum_{m=1}^s \frac{q_{-m}(r)}{r^{m+\nu}} c_{m-n} r^{2m+2}
	+\frac{1}{r^{n-\nu}}\sum_{m=0}^s \frac{Q_{-m}}{r^{\nu}} c_{-m-n}
	\end{equation}
	Denote the free term of $Q_m/r^\nu$ by $A_m$, $m=0,\ldots,s.$
	By the choice of $\nu$ there exists at least one $m$ such that $A_m\neq 0.$
	Considering the behaviour as $r\to 0$ in Eqn.(\ref{balkliouvilleekv9}) we have
	\begin{equation}\label{balkliouvillelem10}
	p_{-n}(r)=\frac{1}{r^{n-\nu}}\sum_{m=0}^s A_m c_{m-n} 
	+O\left(\frac{1}{r^{n-\nu-1}}\right)
	\end{equation}
	By Eqn.(\ref{balkliouvilleekv8}) we have for $n\geq n$
	\begin{equation}
	p_{n}(r)=\sum_{m=1}^s \frac{q_{-m}(r)}{r^{m+\nu}} c_{m+n} r^{2m+n+\nu}
	+\sum_{m=0}^s \frac{q_{m}}{r^{m+\nu}} c_{-m-n}r^{n+\nu}
	\end{equation}
	thus as $r\to 0$ we have
	\begin{equation}\label{balkliouvilleekv11}
	p_{n}(r)=O\left(r^{n+\nu}\right)
	\end{equation}
	in particular
	\begin{equation}\label{balkliouvilleekv12}
	p_{0}(r)=O\left(r^{\nu}\right)
	\end{equation}
	Denote by $A_\lambda$ the smallest nonzero among $A_0,\ldots,A_s$ (we have seen that there exists at least one nonzero)
	i.e.\ $A_\lambda\neq 0,$ $A_m=0$, $m<\lambda.$ The polynomial $Q_i$ has degree at most $s-i$ and
	since it is divisible by $r^\nu$ we have $s-i<\nu,$ $Q_i\equiv 0$ so that
	$A_i=0$ for 
	$s-i<\nu.$ Thus $A_\lambda\neq 0$ for $\lambda +\nu\leq s.$ 
	Let $T(r):=\max_{\abs{x+iy}=r}\abs{P(x,y)}.$ Then $T(r)=O(1)$ as $r\to 0.$
	By Eqn.(\ref{balkliouvilleekv7}) together with Proposition \ref{balkhadamard}
	\begin{equation}\label{balkliouvilleekv13}
	\abs{ \frac{1}{r^{n-\nu}} \sum_{m=0}^s A_m c_{-m-n} +O\left(\frac{1}{r^{n-\nu-1}}\right)
		-O\left(r^{n+\nu}\right)}+O(r^\nu)\leq O(1)
	\end{equation}
	For $n>\nu$ this is only possible if
	$\sum_{m=0}^s A_m c_{-m-n}=0$. 
	Setting
	$n=\nu+k$, $k=1,\ldots,s$, 
	the fact that $m<\lambda A_m=0,$ $A_\lambda \neq 0,$
	implies that there exists constants
	$\{c_{-n}:n=\lambda+\nu+1,\lambda+\nu+2,\ldots\}$
	satisfying
	\begin{equation}\label{balkliouvilleekv14}
	A_\lambda c_{-(\lambda+\nu +k)} +A_{\lambda+1}c_{-(\lambda +\nu+k+1)}+\cdots +A_s c_{-(+\nu+k)}=0
	\end{equation}
	It follows that the function
	$\psi(z)=\sum_{n=0}^\infty c_{-n} z^{-n}$ is a rational function of $z.$
	This implies that $f(z)$ cannot have an essential singularity at $z=0$
	This completes the proof of Theorem \ref{balkliouvillethm1}. So beginning from some $n$, we have $c_{-n}=0.$ Since $A_\lambda \neq 0$ we have by 
	Eqn.(\ref{balkliouvilleekv14}) that $c_{-n}=0$ for $n=\lambda +\nu +k,$ $k\in \Z_+,$ i.e.\ for all $n>\lambda +\nu.$
	But since $s\geq \lambda+\nu$ we have $c_{-n}=0$ for $n>s.$ This implies that
	$f(z)$ (and $\Phi(z)$) has a pole of order not higher than $s$ at $z=0.$
	This completes the proof of Theorem \ref{balkliouvillethm1}.
\end{proof}
We mention that Goldberg \cite{goldbergliouville}
gave a generalization of the above theorem in the sense that the domains $D^+,D^-$ can be replaced by two families of domains
given by more involved conditions.
\begin{corollary}
	Let $a$ be an essential singularity of an analytic function $\Phi(z)$ and let $G$ be a neighborhood of $a.$
	Let $P(x,y)$ be a polynomial and denote by $G_+$ ($G_-$) the set $\{P(x,y)>0\}$ ($\{P(x,y)<0\}$).
	For any two lines $a_1 x+b_1y +c_1$ and $a_2 x+b_2y +c_2$, it cannot hold that one of them does not contain
	any values from $\Phi|_{G_+}$ and simultaneously the other does not contain
	any values from $\Phi|_{G_-}$.
\end{corollary}
\begin{proof}
	Two lines $a_1 x+b_1y +c_1$ and $a_2 x+b_2y +c_2$ will either be parallel or intersect. 
	By rotating by an angle $\alpha$ we can find constants
	$A_1,A_2$ such that for the half planes defined by $\im z< A_1$, $\im z>A_1$, and 
	$\im z <A_2,$ $\im z>A_2$, the function $\exp(i\alpha)\Phi(z)$ is bounded from above or below on $G_+,$ $G_-.$
	so that a point, $a$, cannot be an essential singularity for $\exp(i\alpha)\Phi(z)$ and as a consequence 
	it cannot be an essential singularity for $\Phi(z).$ This completes the proof.
\end{proof}
\begin{corollary}
	Let $a$ be a singular point of a function $\Phi(z)$ satisfying the conditions of Theorem \ref{balkliouvillethm1}.
	Suppose that at $(x,y)=(a_1,a_2)$ at least one of $P'_x$ or $P_y'$ is nonzero. Then $a$ is a simple pole
	for $\Phi(z)$ 
	and the residue of $\Phi(z)$ at $a$ (considered as a vector) is collinear with the tangent to the curve $\{P(x,y)=0\}$ at $(a_1,a_2).$
\end{corollary}
\begin{proof}
	Suppose w.l.o.g.\ $P_y'\neq 0.$ By the implicit function theorem we can for some function $\Psi$, near $(a_1,a_2)$ express the coordinate $y$ on the curve
	$\{P(x,y)=0\}$
	as $y=\Psi(x)$. Let $\tau$ be the tangent to the curve $\{y=\Psi(x)\}$ at $(a_1,a_2)$ and let $\alpha$ be the angle of this tangent to the $x$-axis,
	$\alpha\in [0,\pi).$
	With $z=x+iy$ sufficiently close to the curve $\{P(x,y)=0\}$ its distance to 
	$\abs{z-a}$ will be made arbitrary small. It follows that for sufficiently small $\epsilon>0$ there
	exists $\delta>0$ such that the intersection of $\gamma:=\{y=\Psi(x)\}$ with $\{\abs{z-a}<\delta\}$ belongs to
	\begin{equation}
	\{\alpha -\epsilon <\mbox{arg}(z-a)<\alpha +\epsilon,\quad \alpha-\epsilon +\pi<\mbox{arg}(z-a)<\alpha +\epsilon+\pi\}
	\end{equation}
	and the end points of $\gamma$ will lie in two arcs of a circle $\abs{z-a}=\delta$ enclosed within these angles. 
	Consider two such sectors, $\sigma_1$, $\sigma_2$
	of the circle $\abs{z-a}<\delta,$ where $\sigma_1$ is defined by $\alpha+s<\mbox{Arg}(z-a)<\alpha-\epsilon +\pi,$
	$\sigma_2$ is defined by $\alpha+\epsilon-\pi<\mbox{Arg}(z-a)<\alpha-s.$ Since $a$ is a pole of $\Phi(z)$ by Theorem 
	\ref{balkliouvillethm1}, the proof of the theorem shows that 
	we cannot have that the whole disc $\abs{z-a}<\delta$ belongs to one of the sets $D_+,$ $D_-.$
	This implies that $\gamma$ decomposes the disc into two parts, one belonging to $D_+$ and the other to $D_.$We may suppose $\sigma_1\subset D_-$
	and $\sigma_2\subset D_+.$ Using the notation of the proof of Theorem \ref{balkliouvillethm1}
	we have 
	\begin{equation}\label{balkliouvilleekv15}
	\im f(z)\geq 0\mbox{ on }\sigma_1
	\end{equation}
	Let $a$ be a pole of order $k$ for $\Phi(z)$ (and thus also for $f(z)$).\\
	We claim that $k=1.$
	\\
	Assume (in order to reach a contradiction) that $k>1.$ We may assume that $\delta$ is sufficiently small so that on $\abs{z-a}<\delta$
	$f(z)$ has the representation $f(z)=\sum_{n=-k}^\infty c_n(z-a)^n$ for constants $c_{-n}$, $c_{-k}\neq 0.$
	For $z\to a$ and $k\mbox{arg}(z-a)\neq\mbox{arg}c_{-k}+n\pi$ we have
	$\im f(z)=(1+o(1)\im \frac{c_{-k}}{(z-a)^k}.$
	Put $c_{-k}=\abs{c_{-k}}\exp(i\phi),$ $z-a=r\exp(i\varphi),$ $0\leq \beta <2\pi,$
	$\alpha +\epsilon<\varphi <\alpha -\epsilon +\pi.$
	If $z\in \sigma_1$ then $\varphi$ may take any value between $\alpha +\epsilon$ and $\varphi <\alpha -\epsilon +\pi$.
	For sufficiently small $\epsilon$ we have for $k\geq 2$ that $k\pi-2k\epsilon>\pi.$ 
	Thus there exists a value $\varphi_0$ for the angle $\varphi$ such that
	$\alpha +\epsilon<\varphi_0 <\alpha -\epsilon +\pi$ and $\sin(\beta-k\varphi_0)<0.$ This implies that
	we can choose $r$ sufficiently small so that for $z_0=r\exp(i\varphi_0)$ we have $\im f(z)<0$ for 
	$z_0\in \sigma_1$ which contradicts Eqn.(\ref{balkliouvilleekv15}). We conclude that $k=1.$
	\\
	So when $k=1$, $\im f\geq 0$ on $\sigma_1$ it is necessary that for any $\varphi$ in the interval $(\alpha+\epsilon,\alpha-\epsilon+\pi)$
	we have $\sin(\beta-\varphi)\geq 0$ i.e.\ $2m\pi +\varphi \leq \beta \leq (2m+1)\pi +\varphi$ for integers $m.$
	For each fixed $m$ the inequality $\beta\geq \varphi +2m\pi$ is possible with $\varphi\in 
	(\alpha+\epsilon,\alpha-\epsilon+\pi)$ only if $\beta\geq \alpha-\epsilon +\pi +2m\pi.$
	Analogously, the inequality $\beta\leq \varphi +2(m+1)\pi$ is possible with $\varphi\in 
	(\alpha+\epsilon,\alpha-\epsilon+\pi)$ only if $\beta\leq \alpha+\epsilon +2(m+1)\pi.$
	Since $\epsilon$ can be chosen arbitrarily small we have $\beta=\alpha+\pi+2m\pi.$ Choosing $\beta\in [0,2\pi]$
	we get $\beta=\alpha +\pi.$ This implies that the vector $c_{-1}$ is collinear to the tangent of $\{P(x,y)=0\}$ at $(a_1,a_2)$.  
	Similarly, we obtain for
	$\sigma_2\subset D_-$ that $\beta=\alpha.$ This completes the proof.
\end{proof}

\section{Application of closeness functionals}
It is known that once we have established a class of functions satisfying Liouville's theorem in the classical form
then we can obtain the Liouville property for functions lying {\em close} to the given class, in a certain sense.
Hence our Liouville property for polyanalytic functions can be applied to more general classes.
\begin{definition}\label{tharkclosedef}
	Let $F$ be a class of mappings $\Rn$ into $\R^k$. We assume that for any open $U\subset\Rn$, $U$ is the domain of definition of at least one $f\in F.$
	We denote by $F(U)$,
	the set of mappings $U\to \R^k$ that belong to $F.$ Consider the following properties.
	\begin{itemize}
		\item[(i)] For each open $U\subset\Rn$ the members of $F(U)$ are bounded on compacts of $U$.
		\item[(ii)] If $T_1:\Cn\to \C,$ $T_1 (z)= a_1 z+z_0,$
		and $T_2:\C\to \C,$ $T_2 (\zeta)= a_1 \zeta+\zeta_0,$ for constants $z_0\in \Cn,\zeta_0\in \C$
		and real constants $a_1,a_2>0$, then for all $f\in F(U)$, for an open $U\subset\Cn$
		we have that $T_2\circ f\circ T_1\in \mbox{PA}_\alpha(T_1^{-1}(U)).$
		\item[(iii)] For an open $U\subset\Cn$, any uniformly bounded family in $ F(U)$
		is equicontinuous on compact subsets of $U.$
		\item[(iv)] $F$ is closed under uniform convergence on compact subsets of the domains of definition.
		\item[(v)] For an open $U\subset\Cn$, if $f\in  F(U)$ then
		$ F(V)$ for any open $V\subset U.$
		\item[(vi)] If $f:U\to \R^k$ such that for each $x\in U$ there exists a neighborhood $U_x\subset U$ of $x$ such that
		$f|_{U_x}\in F$ then $f\in F.$
	\end{itemize}
\end{definition}
\begin{remark}
	Tarkhanov \cite{tarkhanov}, p.413, points out the following:
	Let $L$ be a homogeneous differential operator ($\Rn\times \R^k \to \Rn\times \R^l$) 
	with constant coefficients and injective symbol in $\Rn$ (in particular if $L$ is elliptic).
	Denote by $S$ the sheaf of solutions to $Lf=0$ on open subsets of $\Rn.$
	Then $S$ has the properties (i)-(vi). 
	For example the class of holomorphic mappings from open sets in $\C$ into $\C$ has the properties (i)-(vi), if one identifies
	the mappings with those of domains in $\R^{2}$ into $\R^{2}.$ Since also $\partial_{\bar{z}}^q$ is elliptic for all $q\in \Z_+$
	we have that $\mbox{PA}_q(\C)$ satisfies (i)-(vi). 
\end{remark}
\begin{definition}
	Let $F$ be a family of locally bounded maps
	$\Rn\to \R^k$ and for an open subset $V\subset\Rn$ define $F(v)$ as above.
	$U\subset \Rn$ be a domain and let $f:U\to\R^k$
	be a locally bounded map. For $\theta\in (0,1)$ and arbitrary $R>0$ such that
	$B(x,R):=\{y\in \abs{y-x}<R\}\subset U$ define
	\index{$\mathfrak{d}_{\theta}(f,F)$, Closeness functional}
	\begin{equation}
	\mathfrak{d}_{\theta}(f,F):=\left\{
	\begin{array}{ll}
	\inf_{u\in F(B(x,R))}\left(\sup_{y\in B(x,\theta R)}\frac{\abs{f(y)-u(y)}}{\mbox{diam} f(B(x,R))}\right) & ,\mbox{ if diam}f(B(x,R))\neq 0,\infty\\
	0 & , \mbox{otherwise}
	\end{array}
	\right.
	\end{equation}
		For a ball $B\subset\Rn$, a map $f:\Rn\to \R k$ and a real number $\delta>0$
	we denote in this section 
	\begin{equation}
	m_{B}(f)(r):=\sup_{\stackrel{x,y\in B}{\abs{y-x}<\delta}} \abs{f(y)-f(x)}
	\end{equation}
\end{definition}
We shall need the following lemma due to Kopylov \cite{kopylov}, which we state without proof.
\begin{lemma}[Kopylov \cite{kopylov}]\label{thark9117}
	Suppose $F$ satisfies (ii),(iii),(v) of Definition \ref{tharkclosedef} and let $f:B(x_0,R)\to \R^k$ be a bounded map
	on $B(x_0,R)\subset\Rn$ for some $R>0,$ such that $\mathfrak{d}_\delta(f,F)<\frac{1}{2}$ for some $\theta\in (0,1).$
	Then for each $t\in (0,1)$ and each $\delta\in (0,(1-t)\epsilon\theta R)$ we have
	\begin{multline}\label{8888333311}
	m_{B(x_0,tR)}(f)(\delta)\leq \\
	( (1+2\mathfrak{d}_\theta (f,F))
	\sup_{\stackrel{u\in F(B(0,1))}{\abs{u}<1}} m_{B\left(0,\frac{1}{2}\right)}(u)(\epsilon)+2\mathfrak{d}_\theta(f,F) )^\nu
	\mbox{diam}f(B(x_0,R))
	\end{multline}
	where 
	\begin{equation}
	\nu:=\frac{\log(1+2\frac{1}{2}\delta(1-t)(1-\frac{1}{2}\epsilon\theta)R)}{-\log\left(\frac{1}{2}\epsilon\theta\right)}-1
	\end{equation}
\end{lemma}
Note (iii) of Definition \ref{tharkclosedef} implies that $m_{B(0,\frac{1}{2})}(u)(\epsilon)$ becomes arbitrary small uniformly in $u\in F(B(0,1))$, as $\epsilon\to 0$, with $\abs{u}<1$ and 
$2\mathfrak{d}_\theta(f,F)<1$ as in the hypothesis of Lemma \ref{thark9117}, i.e.\
$\epsilon$ can be chosen such that in Eqn.(\ref{8888333311})
\begin{equation}
( (1+2\mathfrak{d}_\theta (f,F))\sup_{\stackrel{u\in F(B(0,1))}{\abs{u}<1}} m_{B\left(0,\frac{1}{2}\right)}(u)(\epsilon)+2\mathfrak{d}_\theta(f,F) )
<1
\end{equation}
\begin{theorem}[See Tarkhanov \cite{tarkhanov}, Thm 9.1.18]
	Let $F$ be a class of functions satisfying (ii), (iii), (v). 
	Let $f:\R^n\to \R^k$ be a bounded mapping such that $\mathfrak{d}_\theta (f,F)<\frac{1}{2}$ for some $\theta\in (0,1].$ Then
	$f$ is a constant mapping.
\end{theorem}
\begin{proof}
	Since $\mathfrak{d}_\theta (f,F)\leq \mathfrak{d}_1 (f,F)$ for all $\theta\in (0,1)$ and each $f:U\to \R^k$, it suffices to
	prove the result for the case $\theta\in (0,1).$ Assume $\abs{f(z)}\leq C$ for some $C>0$ and all $x\in \Rn.$
	Let $t=1/2$ and $\epsilon\in (0,1/2)$ such that
	\begin{equation}
	v:=(1+2\mathfrak{d}_\theta (f,F))\sup_{\stackrel{u\in F(B(0,1))}{\abs{u}<1}}m_{B\left(0,\frac{1}{2}\right)}(u)(\epsilon)+2\mathfrak{d}_\theta(f,F)<1
	\end{equation}
	For each pair $x,y\in\Rn$ let $N_0\in \N$ be fixed such that
	$N_0>2\frac{1}{\epsilon\theta}\abs{y-x}$. Each ball $B(x,N/2)$, $N\geq N_0$, $N\in \N$, contains $y.$					
	By Lemma \ref{thark9117} applied to $f|_{B(x,N)}$ for $N\geq N$ we have
	\begin{equation}
	\abs{f(y)-f(x)}\leq m_{B(x,\frac{N}{2})}(f)(\abs{y-x})\leq v^\nu \mbox{diam}f(B(x,N))
	\end{equation}
	\begin{equation}
	\nu=\frac{\log\left(1+\frac{1}{\abs{x-y}}(1-\frac{1}{2}\epsilon\theta)N\right)}{-\log\left(\frac{1}{2}\epsilon\theta\right)}-1
	\end{equation}
	But $f(B(x,N))\leq 2C$ for large $N$ so that $\lim_{N\to \infty} v^\nu \mbox{diam}f(B(x,N))=0$
	which implies $f(x)=f(y).$
	This completes the proof.
\end{proof}
\begin{corollary}
	Let $q\in \Z_+$ and 
	let $f:\C\to \C$ be a bounded mapping such that 
	$\mathfrak{d}_\theta (f,\mbox{PA}_q(\C))<\frac{1}{2}$ for some $\theta\in (0,1].$ Then
	$f$ is a constant mapping.
\end{corollary}

\chapter{$q$-polyanalytic functions on subsets of $\Z[i]$}
\section{Introduction}
Some of the pioneers of the investigation of analogues of complex analytic functions on $\Z[i]$ were  
Isaacs \cite{isaacs1}, \cite{isaacs2}, Ferrand \cite{ferrand1} and Duffin \cite{duffin1},\cite{duffin2}.
Notable modern contributions to the field include the works of Kiselman \cite{kiselman1}, \cite{kiselman2}.
In Isaacs \cite{isaacs1} the {\em monodiffric functions of the first kind} on  
the discrete complex plane,
where defined square-wise as those that where annihilated by a certain first order linear difference operator, in particular 
a complex-valued function $f$ on $\Z[i]$, is monodiffric of the first kind 
on a square with vertices\index{Monodiffric functions of the second kind}\index{Monodiffric functions of the first kind}
$\{z,z+1,z+i,z+i+1\},$ whose lower left point is $z\in \Z[i],$ 
if and only if $f$
satisfies, 
\begin{equation}\label{firstkind}
f(z+1)-f(z)=\frac{f(z+i)-f(z)}{i}
\end{equation} 
We shall say that $f$ is monodiffric of the first kind {\em at} $z$
if and only if $f$
satisfies equation \ref{firstkind}.
In this paper we shall say that a function $f$ in the discrete complex plane is {\em monodiffric  of the second kind} at $z\in \Z[i],$ if and only if
\begin{equation}\label{secondkind}
\frac{f(z+1+i)-f(z)}{i+1}=\frac{f(z+i)-f(z+1)}{i-1}
\end{equation}
Ferrand \cite{ferrand1} (who uses a discrete
version of Moreras theorem) 
used the term {\em preholomorphic} for the monodiffric functions of the second kind. 
In this paper we shall say that a function $f$ in the discrete complex plane is {\em monodiffric functions of the third kind} at $z\in \Z[i],$
\begin{equation}\label{thirdkind}
f(z+1)-f(z-1)=\frac{f(z+i)-f(z-i)}{i}
\end{equation}
The monodiffric functions of the third kind (these where also introduced by Isaacs \cite{isaacs1}, p.179) appear less 
frequently in the literature, and then they are not referred to as 
monodiffric functions of the third kind. 
we shall be interested in powers of the operators in equation \ref{firstkind}, \ref{secondkind} and \ref{thirdkind} respectively.
To avoid confusion we point out that in 
Kurowski \cite{kurowski2}, the functions that we here call monodiffric of the third kind, are called 
monodiffric of the second kind. We have chosen not to adapt that terminology and instead use the
terminology used by e.g.\ Kiselman \cite{kiselman1}, \cite{kiselman2}, regarding monodiffric functions of the first and second kind, see also Daghighi \cite{daghighidiscrete}.
\begin{definition}[$q$-polyanalytic functions (polyanalytic functions of order $q$)]\label{qanaldef0}
	Define for complex-valued functions $f$ on $\Z[i]$,
	\begin{equation}
	L_1 f(z):=f(z+1)-f(z)+i(f(z+i)-f(z))
	\end{equation} 
	\begin{equation}
	L_2 f(z):=f(z)+if(z+1)-f(z+1+i)-if(z+i)
	\end{equation}
	\begin{equation}
	L_3 f(z):=f(z+1)-f(z-1)+i(f(z+i)-f(z-i))
	\end{equation} 
	We define, for a given positive integer $q$, and a fixed $j\in \{1,2,3\},$ a complex-valued function $f\colon \Z[i]\to \C$ to be:\\
	{\em $q$-polyanalytic} (or {\em polyanalytic of order $q$}) {\em of the first kind} at $z\in \Z[i]$
	if
	\begin{equation}
	L_1^q f(z)=0
	\end{equation} 
	{\em $q$-polyanalytic} (or {\em polyanalytic of order $q$}) {\em of the second kind} at $z\in \Z[i]$
	if
	\begin{equation}
	L_2^q f(z)=0
	\end{equation}
	and 
	{\em $q$-polyanalytic} (or {\em polyanalytic of order $q$}) {\em of the third kind} at $z\in \Z[i]$
	if
	\begin{equation}
	L_3^q f(z)=0
	\end{equation}
	If the condition holds true at each point of a subset $S\subseteq \Z[i]$ where the defining operator is defined, then we say that $f$ is
	$q$-polyanalytic (or polyanalytic of order $q$) of the first, second or third kind, respectively
	on $S$ and when it is clear from the context what $S$ is we simply say that 
	$f$ is $q$-polyanalytic (or polyanalytic of order $q$) of the first, second or third kind respectively.
\end{definition}
Kiselman \cite{kiselman1}, Sec 3, pointed out that at the level of ideas, the operators defined by equation \ref{secondkind}
and \ref{thirdkind} are quite similar. Kurowski \cite{kurowski1}, p.1, expresses a point that we interpret to be in the same vein as that of Kiselman.
It is clear however that the solution spaces defined by the operators $L_2,L_3$
are not, in any formal rigorous way, equivalent. This can easily be displayed by example. 
See Proposition \ref{inequiv} for a more comprehensive result regarding any pair
among $L_1,L_2,L_3$.
\begin{example}\label{anv}
	Let $z\in \Z[i],$ and let
	$V:=\{z,z+1,z+2,z+2+i,z+2+2i,z+1+2i,z+2i,z+1,z+1+i\}$. Now relative to $V$,
	the points $z,z+1,z+1+i,z+i,$ are the only points where $L_2$ is defined whereas $z+1+i$ is the only point where
	$L_3$ is defined.  
	To this end, set $f(z+2+i)=1, f(z+i)=f(z+1+2i)=f(z+1)=0$ (making sure that $L_3 f(z+1+i)\neq 0$) and set
	$f(z+1+i)=0.$ Then there are four undefined values $f(z),f(z+2),f(z+2+2i),f(z+2i),$ and we invoke four conditions
	$L_2 f(z)=L_2 f(z+1)=L_2 f(z+1+i)=L_2 f(z+i)=0$, giving,
	$$
	\begin{bmatrix}
	1 & 0 & 0 & 0\\
	0 & i & 0 & 0\\
	0 & 0 & -1 & 0\\
	0 & 0 & 0 & -i
	\end{bmatrix}
	\begin{bmatrix}
	f(z)\\
	f(z+2)\\
	f(z+2+2i)\\
	f(z+2i)
	\end{bmatrix}
	=$$
	$$
	\begin{bmatrix}
	-if(z+1)+f(z+1+i)+if(z+i)\\
	-f(z+1)-if(z+2+i)+if(z+1+i)\\
	-f(z+1+i)-if(z+2+i)+if(z+1+2i)\\
	-f(z+i)-if(z+1+i)+f(z+1+2i)
	\end{bmatrix}
	=
	\begin{bmatrix}
	0\\
	1\\
	1\\
	0\\
	\end{bmatrix}
	$$ 
	Which gives,
	$$
	\begin{bmatrix}
	f(z)\\
	f(z+2)\\
	f(z+2+2i)\\
	f(z+2i)
	\end{bmatrix}=
	\begin{bmatrix}
	1 & 0 & 0 & 0\\
	0 & i & 0 & 0\\
	0 & 0 & -1 & 0\\
	0 & 0 & 0 & i
	\end{bmatrix}
	\begin{bmatrix}
	0\\
	1\\
	1\\
	0\\
	\end{bmatrix}
	=
	\begin{bmatrix}
	0\\
	i\\
	-1\\
	0\\
	\end{bmatrix}
	$$
	Thus, we have a complex-valued function, $f$, on $V$
	such that $L_2 f(z)=L_2 f(z+1)=L_2 f(z+1+i)=L_2 f(z+i)=0$ but $L_3 f(z+1+i)\neq 0.$ \\
	We choose $z=0.$
	Set $f(w)=0,$ for all $w\in V'$ where
	$V'=\{ z=x+iy\in \Z[i]\setminus V\colon x\cdot y=0\}.$
	This uniquely defines an extension of $f$ that satisfies $L_2 f=0$  
	at all points of $\Z[i]$. Indeed, consider the point 
	$(2+k)+i$ where $k=1$. The value $f(2+k+i)$ is uniquely determined by the value at
	$2+k, 2+k-1,(2+k-1)+i$ which are all known for $k=1.$ Once the value
	at $f(2+k+i)$ is determined the process can be repeated by replacing 
	$k$ with $k+1$, thus uniquely determining $f$ on the set
	$V\cup \{ y=1, x\geq 3\}.$ Obviously we can then iteratively repeat the process row-wise until
	we have $f$ uniquely determined on the upper right quadrant. Analogously, we have the
	unique extension of $f$ that satisfies $L_2f=0,$ to each the three remaining quadrants by the same iterative process
	(See Section \ref{solutionsec} for a formalization of this process).
	This yields an extension of $f$ that is an entire $1$-polyanalytic of the second kind, and clearly
	no extension of $f$ can be $1$-polyanalytic of the third kind at $z+1+i.$
	Proposition \ref{inequiv} below gives a more comprehensive result on these matters.
\end{example}

\section{Integrals}\label{cauchysec}
A {\em polygon}, $\Gamma$, in the complex plane $\C$ consists of a set of $N$ edges $[a_0, a1],
[a_1, a_2],$\\
$\ldots,[a_{N-1}, a_0]$, where $a_j, j=0,\ldots ,N$ are given points in $\C$. 
In particular, it is determined by the ordered set of vertices
$(a_0,\ldots, a_{N})\in\C^N$.
A function $f$ defined on $\Gamma$ is called {\em piecewise affine} if $f$ is affine on each segment $[a_j , a_{j+1}]$
with the possible exception of the points that belong to two or more segments. $\Gamma$ is called {\em closed} if $a_0=a_{N}.$
Let $f$ be piecewise affine on a polygon $\Gamma$ determined by
$(a_0,\ldots, a_{N})\in\C^N$. 
Then (see e.g.\ Kiselman \cite{kiselman2}, p.2)
\begin{equation}\label{jo}
\int_{\Gamma} f(z)dz = \frac{1}{2}\sum_{j=1}^{N}(f(a_{j})+f(a_{j-1}))(a_{j}-a_{j-1}) 
\end{equation}
\begin{definition}[Integral of complex functions along polygons]
	Let $f$ be a complex-valued function on a polygon $\Gamma$ determined by
	$(a_0,\ldots, a_{N})\in\C^N$. We define the integral of $f$ along $\Gamma$ as 
	\begin{equation}
	\int_{\Gamma} f(z)dz = \frac{1}{2}\sum_{j=1}^{N}(f(a_{j})+f(a_{j-1}))(a_{j}-a_{j-1}) 
	\end{equation}
\end{definition}
Let $f$ be a complex-valued function on $\Z[i]$. Let $p_0\in \Z[i],$ 
and denote by $\Gamma_{p_0}$ the closed polygon defined
by the ordered set of vertices $(a_1,a_2,a_3,a_4):=(p_0, p_0 +1,p_0+1+i,p_0+i)$, i.e.\ moving counter-clockwise.
It is easy to verify that $f$ is $1$-polyanalytic of the second kind at $p_0$ iff 
\begin{equation}\label{samre}
\int_{\Gamma_{p_0}} \tilde{f}(z)dz=0
\end{equation}
where $\tilde{f}$ is the unique piecewise affine function on the closed polygon $\Gamma_{p_0}$ 
such that $\tilde{f}(z)=f(z)$ for $z\in \{ p_0, p_0 +1,p_0+i,p_0+(1+i)\}.$
Indeed, 
$\int_{\Gamma_{p_0}} f= (f(p_0+1+i) +f(p_0+1))i +
(f(p_0+i) +f(p_0+1+i))(-1) + 
(f(p_0) +f(p_0+i))(-i)+$
$(f(p_0+1) +f(p_0))=
(1-i)f(p_0) 
+(1+i)f(p_0+1)
+(-1-i)f(p_0+i)
+(i-1)f(p_0+1+i)
=$
$(1-i)(f(p_0) + if(p_0+1) -f(p_0+1+i) -if(p_0+i))=$
$
(1-i)L_2 f(p_0)=0.$
(In fact a stronger result holds true, see Remark \ref{dubbrem}).
\\
As we have pointed out the condition of equation \ref{secondkind} at $p_0$ does not involve all four adjacent points to $p_0$ but it actually
involves the point $p_0+i+1$ which is not adjacent to $p_0$ in the usual sense.
Considering the operator $L_3$ which at $z$ involves precisely the four points $z\pm i,z\pm 1,$ we have the following.

\begin{observation}\label{obs1}
	Let $f$ be a complex-valued function on $\Z[i]$, let $p_0\in \Z[i],$
	and denote by $\Gamma_{p_0}$ the closed polygon defined
	by the ordered set of vertices $(p_0+1, p_0 +i,p_0-1,p_0-i)$.
	Then,
	\begin{equation}\label{battre}
	L_3 f(p_0)=0 \Leftrightarrow \int_{\Gamma_{p_0}} \tilde{f}(z)dz =0
	\end{equation}
	where $\tilde{f}$ is the unique piecewise affine function on the closed polygon $\Gamma_{p_0}$ 
	such that $\tilde{f}(z)=f(z)$ for $z\in \{p_0+1, p_0 +i,p_0-1,p_0-i\}.$ 
\end{observation}
\begin{proof}
	Using the notation from equation \ref{samre} together with equation \ref{jo}, we have
	\begin{multline}
	2\int_{\Gamma_{p_0}} f(z)dz =0 \Leftrightarrow 
	(f(p_0+i) +f(p_0+1))(i-1) +\\
	(f(p_0-1) +f(p_0+i))(-1-i) + 
	(f(p_0-i) +f(p_0-1))(1-i)+\\
	(f(p_0+1) +f(p_0-i))(1+i)
	=\\
	-2f(p_0+i)
	+2if(p_0+1)
	-2if(p_0-1)
	2f(p_0-i) 
	=\\
	-\frac{2}{i}(i(  f(p_0+i)  -f(p_0-i)) +
	f(p_0+1)
	-f(p_0-1))
	=\\-\frac{2}{i} L_3 f (p_0)=
	0\Leftrightarrow L_3 f(p_0) =0 
	\end{multline}
\end{proof}

\begin{definition}[Zig-zag polygons]\label{zigdef}
	A polygon determined by the ordered set of (possibly infinite) points
	$(a_0,\ldots,a_N),$(or $(a_0,\ldots)$) $a_j\in \Z[i],$ $j=0,\ldots,N$ (or $j=0,1,\ldots$) is called a {\em zig-zag polygon} if $a_j-a_{j-1}\in \{1\pm i,-1\pm i\},$ 
	$j=1,\ldots,N$ (or $j=0,1,\ldots$) some positive integer $N.$ It is {\em non-self-intersecting} if $a_k\neq a_l$ for $k\neq l$ except possibly for $(k,l)\in \{(0,N),(N,0)\}.$
	It is further called {\em closed} if it has $N-1$ points where $a_0=a_N$.
	A point of any subset $\omega\subset \Z[i]$ is called an {\em interior point}
	if and only if all four of its adjacent points of first order belong to $\omega.$ The set of interior points is denoted
	$\mathring{\omega}.$ 
	A {\em subset with zig-zag boundary}, $\omega\subset\Z[i]$ is the union of a possibly infinite set of points 
	together with their sets of adjacent points of first order, such that the set of non-interior points can be ordered to yield a 
	non-self-intersecting
	zig-zag polygon.
	Such $\omega$ is called a {\em domain with zig-zag boundary} if each pair of interior points can be connected by a 
	zig-zag polygon contained in $\mathring{\omega}$ and it is called a {\em simple domain with zig-zag boundary} if
	each pair of non-interior points can be connected by a 
	zig-zag polygon in the set of non-interior points.
	A point is called {\em zig-zag even} ({\em odd}) if it can be connected by a zig-zag polygon to $0$ ($1$). 
	We also call $0$ zig-zag even and we call $1$ zig-zag odd.
\end{definition}
Obviously no zig-zag even point can be connected to a zig-zag 
odd point by a zig-zag polygon, in particular a zig-zag polygon
consists either entirely of zig-zag even points or entirely of zig-zag odd points.
\begin{proposition}\label{pathind}
	Let $f$ be a complex-valued function on $\Z[i]$. Let $\Omega\subset\Z[i]$ be a simple domain with zig-zag boundary.
	The function $f$ satisfies $L_3 f(z)=0$ on $\mathring{\Omega}$ if and only if
	for any closed non self-intersecting zig-zag polygon $\gamma\subset \Omega$ defined
	by an ordered set of vertices
	\begin{equation}\label{battre1}
	\int_{\gamma} f(z)dz=0
	\end{equation}
\end{proposition}
\begin{proof}
	First of all, a function on a domain with zig-zag boundary satisfies
	$L_3f(z)=0$ on $\mathring{\Omega}$ if and only if 
	$L_3f(z)=0$ on $\mathring{\omega}$ for any domain with zig-zag boundary $\omega\subseteq\Omega$, with finitely
	many elements. Hence we only need to prove the statement for the case of finitely many interior points.
	We use induction in the number, $n$, of interior points (see Definition \ref{zigdef}) of the discrete 
	domain, $\Omega$, with $\Gamma$ as its zig-zag boundary, where $\Gamma$ is defined by the ordered set of points $(a_0,\ldots,a_{N-1})$.
	The case $n=1$ is precisely equation \ref{battre}.
	Assume $n>1$ and that the result holds true for the case of $n-1$ interior points. Since $n$ is finite
	we can find 
	an interior point, $z_0+i=x_0+iy_0$, such that
	$y_0$ is minimal and finite.  
	In particular, the points $z_0-1+i,z_0,z_0+1+i$ belong to $\Omega$ but are not interior points and we can assume
	$a_{n_0}=z_0$. 
	We can assume that the boundary is traversed counter-clockwise as vanishing of the integral will be independent
	with reversed direction.
	First consider when also the point $z_0+2i$ is not an interior point of $\Omega,$
	i.e.\ $\Gamma$ contains the ordered subsequence $z_0+2i,z_0-1+i,z_0,z_0+1+i.$
	The set $\Omega':=\Omega\setminus \{z_0,z_0+i,z_0-1+i\},$ is again a domain with zig-zag boundary, $\Gamma'$, defined by
	$(a_0,\ldots,a_{n_0-3},z_0+2i,z_0+1+i,a_{n_0+2},\ldots,a_{N-1}).$
	The only domains in $\Omega$ with boundary $\gamma$ that is a
	closed non self-intersecting zig-zag polygon, such that $\gamma$ does not also lie in $\Omega'$, is either
	simply the polygon, 
	$\hat{\gamma}_{z_0},$ defined by
	$(z_0+2i,z_0-1+i,z_0,z_0+1+i)$ 
	or 
	$\gamma$
	is 
	defined by 
	$(b_0,\ldots,z+2i,z_0-1+i,z_0,z_0+1+i,\ldots,b_{M-1})$ for some positive integer $M$ and points $b_j$ in $\Omega'$, in particular does not contain
	$z_0+i$.
	(Note that the case where $\gamma$ contains $z_0+i$ cannot occur because the point $z_0+2i$ is assumed not to be an interior point).
	The first case is handled by a translated version of equation \ref{battre}. 
	So assume the second case.
	Then letting $\gamma'$ be the closed non self-intersecting zig-zag polygon $\gamma'\subseteq \Omega'$
	defined by $(b_0,\ldots,z_0+2i,z_0+1+i,\ldots,b_{M-1})$ (i.e.\ the two points $z_0-1+i,z_0$ are removed) we have
	\begin{multline}
	2\left(\int_{\gamma'}-\int_{\gamma}\right) f(z)dz=
	(f(z_0+2i)+f(z_0+1+i))(i-1)
	-\\
	[(f(z_0+2i)+f(z_0-1+i))(i+1)
	+
	(f(z_0-1+i)+f(z_0))(-1+i)
	+\\
	(f(z_0)+f(z_0+1+i))(-1-i)
	]
	=\\
	-2 f(z_0+2i)
	+2 f(z_0)
	+2i f(z_0+1+i)
	-2if(z_0-1+i)
	=\\
	-\frac{2}{i}(if(z_0+2i))-if(z_0) +f(z_0+1+i)-f(z_0-1+i))\\
	=-\frac{2}{i}L_3 f(z_0+i)=0
	\end{multline}
	But $\Omega'$ has $n-1$ interior points ($z+i$ is not part of the set) thus by the induction hypothesis
	$\int_{\gamma'} f(z)dz=0$ for any closed non self-intersecting zig-zag polygon $\gamma'\subseteq \Omega'$.
	This takes care of the case when $z_{0}+2i$ is not an interior point.
	
	Now assume that $z_{0}+2i$ is an interior point of $\Omega$ (in particular this implies that $z+3i, z+2i\pm 1$ belong to $\Omega$).
	
	Then the set $\Omega':=\Omega\setminus \{z_0,z_0+i\}$ is again a domain with zig-zag boundary, $\Gamma'$, defined by
	$(a_0,\ldots,z_0-1+i,z_0+2i,z_0+1+i,\ldots,a_{N-1})$,
	but $\Omega'$ has $n-1$ interior points thus by the induction hypothesis
	$\int_{\gamma'} f(z)dz=0$ for any closed non self-intersecting zig-zag polygon $\gamma'\subseteq \Omega'$.
	
	It is easy to see that the only closed non self-intersecting zig-zag polygon $\gamma\subseteq \Omega$
	which does not also lie in $\Omega'$, is one defined by either
	$(b_0,\ldots,z_0-1+i,z_0,z_0+1+i,\ldots,b_{M-1})$ 
	or $(b_0,\ldots,z_0-1+2i,z_0+i,z_0+1+2i,\ldots,b_{M-1})$
	for some positive integer $M$ and points $b_j$ in $\Omega'.$
	
	For each such $\gamma$ define the associated closed non self-intersecting zig-zag polygon $\gamma'\subseteq \Omega'$
	by $(b_0,\ldots ,z_0-1+i,z_0+2i,z_0+1+i,\ldots ,b_{M-1})$ in the first case and by 
	$(b_0, \ldots,z_0-1+2i,z_0+3i,z_0+1+2i,\ldots ,b_{M-1})$ in the second case.
	In the first case we have
	\begin{multline}\label{ko}
	2\left(\int_{\gamma'}-\int_{\gamma}\right)
	f(z)dz=
	(f(z_0+2i) +f(z_0-1+i))(1+i)+\\
	(f(z_0+1+i) +f(z_0+2i))(1-i)
	-[
	(f(z_0) +f(z_0-1+i))(1-i)+\\
	(f(z_0+1+i) +f(z_0))(1+i)
	]
	=
	(f(z_0+2i) +f(z_0-1+i))(1+i)+\\
	(f(z_0+1+i) +f(z_0+2i))(1-i)
	+
	(f(z_0) +f(z_0-1+i))(i-1)+\\
	(f(z_0+1+i) +f(z_0))(-1-i)
	=
	-2f(z_0)
	+2if(z_0-1+i)
	-2if(z_0+1+i)
	+
	\\
	2f(z_0+2i)=
	\frac{2}{i}(i(f(z_0+2i)-f(z_0)) +f(z_0+1+i)-f(z_0-1+i))=\\
	\frac{2}{i} L_3 f(z_0+i)=0
	\end{multline}
	
	In the second case, if $z_0+3i$ is not an interior point then we can repeat the procedure applied to the case
	when $z_0+2i$ was not an interior point. If $z_0+3i$ is an interior point the same calculations as in equation \ref{ko}, but translated one step in the $\im z$-directions  yields
	\begin{equation}
	2\left(\int_{\gamma'}-\int_{\gamma}\right) f(z)dz=\frac{2}{i} L_3 f(z_0+2i)=0
	\end{equation}
	
	This proves the induction step. 
	This completes the proof.
\end{proof}
\begin{remark}\label{dubbrem}
	Kiselman \cite{kiselman2} defines a polygon determined by the ordered set
	$(a_0,\ldots,a_N),$ $a_j\in \Z[i],$ $j=0,\ldots,N$ to be a {\em 4-curve} if $a_j-a_{j-1}\in \{\pm 1,\pm i\},$ 
	$j=1,\ldots,N$ 
	and it is a well-known result see e.g.\ Isaacs \cite{isaacs1}, p.183,
	that if $f$ is a monodiffric function of the first kind then 
	\begin{equation}\label{dubb}
	\int_{\gamma} f(z)dz=0
	\end{equation}
	for each closed (non-self-intersecting) 4-curve $\gamma$. The corresponding result
	for monodiffric functions of the second kind also holds true (see e.g.\ Duffin \cite{duffin2}, Corollary 2.1.1).
\end{remark}

\section{Characterization of $q$-polyanalytic functions}
\label{solutionsec}
In this section we shall obtain the kernels of the powers of the operators $L_1,L_2,L_3$ which in turn 
give the defining difference equations for $q$-polyanalytic functions on $\Z[i]$. 
For background on solving finite difference equations see e.g.\ Mickens \cite{mickens} and Jordan \cite{jordan}. For this particular section we shall 
in the interest of conformity with previous literature use some special notations.

\subsection{First kind, $q=1$}
For a complex-valued function $f$ on $\Z\times \Z$ (or $\Z[i]$ in which case we shall still write $f(k,l)$ instead of $f(k+il)$ where
$(k,l)\in \Z^2$), we use the notation 
\begin{equation}
E_1f(k,l):=f(k+1,l)-f(k,l),\quad E_2 f(k,l):=f(k,l+1)-f(k,l)
\end{equation}
This implies that a function is monodiffric of the first kind on $\Z[i]$ if at each $k+il\in \Z$, we have
\begin{equation}
E_1 f(k,l)=(iE_2-I+i)f(k,l)
\end{equation}
where $I$ denotes the identity operator.
For $k\geq 0,$ this can be is solved by the symbolic method of Boole (see e.g.\ Jordan p.616)
to yield,
\begin{equation}
f(k,l)=(iE_2-I+i)^k \phi(l)=\sum^k_{j=0} i^j\binom{k}{j} E^j_2 \phi(l)
\end{equation}
where $\phi$ is an arbitrary function. Hence $f(k,l)$ is monodiffric on $\Z[i]$ if and only if at each point $k+il\in \Z[i], k\geq 0,$ we have,
\begin{equation}\label{gensolfirst}
f(k,l)=\sum^l_{j=0} i^j\binom{k}{j} \phi(l+j)
\end{equation}
for an arbitrary $\phi.$
Obviously, given a function $f$ defined for $\{ z\in \Z[i]\colon \re z\geq 0 \},$ and $1$-polyanalytic of the first kind, to 
obtain an extension to $\Z[i]$ that is $1$-polyanalytic of the first kind it is sufficient to 
know the values of $f$ on $\{ z\in \Z[i]\colon \re z<0 \}.$ 
\begin{observation}\label{firstkind1}
	Let $f$ be $1$-polyanalytic of the first kind on $\Z[i]$, let $z=x+iy$ denote the standard coordinate in $\Z[i],$ and let $p_0\in \Z[i].$ Then
	$f$ is uniquely determined by its values on $D_0:=\{ z\in \Z[i]\colon \re z=\re p_0 \}\cup \{ z\in \Z[i]\colon (\re z<\re p_0) \wedge (\im z=\im p_0)\}$ 
	Also any proper subset of $D'\subset D$ is not a set of uniqueness (in the sense that there are
	two different $1$-polyanalytic functions on $\Z[i]$ that agree on $D'$). 
\end{observation}
\begin{proof}
	We can use the same procedure
	as that in Example \ref{anv} but applied to $L_1$ instead of $L_2$ namely start from the set, $D_0$. Define iteratively the sets $D_j, j\in \N,$
	as follows: Let $D_{j+1}$ be the set of all points $s\in \Z[i]$ 
	satisfying
	$s\in \{ w,w+1,w+i\}$ for some $w\in \Z[i]$ such that
	precisely one point of the set
	$\{ w,w+1,w+i\}$ does not belong to $D_j.$ Then $\bigcup_{j} D_j=\Z[i]$. $f$ can be iteratively extended to each $D_j$ by assigning the value of $f$
	at $s\in D_j\setminus D_{j-1}$ to be determined by the equation $f(w+1)+if(w+i)-(1+i)f(w)=0.$ 
	Obviously replacing the value of $f$ at a point of $D_0$ yields a different extension to $\Z[i].$
	This completes the proof.
\end{proof}

\begin{remark}[Uniqueness of extension]\label{unique}
	In the proof of Observation \ref{firstkind1}, \ref{secondkind1},\ref{thirdkind1} respectively, the procedure for obtaining the extension of $f$ from its values on the set $D_0$ can only be done in one way.
	In particular, the only extension of a function that vanishes on $D_0$ is the 
	function that vanishes identically on $\Z[i].$
	In other words, if two functions $F_1,F_2$ on $\Z[i],$ satisfy for some $j=1,2,3,$
	$L_j F_1\equiv 0,$ and $L_j F_2\equiv 0$ then
	\begin{equation}
	F_1|_{D_0}=F_2|_{D_0} \Leftrightarrow F_1\equiv F_2
	\end{equation}
	\begin{definition}
		For a fixed $j\in \{ 1,2,3\},$ we say that a set $D\subset \Z[i]$
		is a {\em set of uniqueness} (with respect to $L_j$) if
		any function $f$ satisfying $L_j f\equiv 0$ on $\Z[i]$
		and $f|_D\equiv 0$ must vanish identically. It is a {\em minimal}
		set of uniqueness if it does not properly contain any other set of uniqueness.
	\end{definition}
\end{remark}

\subsection{Second kind, $q=1$}\label{secete}
From a purely theoretical perspective we can find two parametrized independent solutions in the kernel of $L_2$ as follows. 
We use Lagrange's method (see e.g.\ Mickens Section 5.3, p.186) in order to find a particular solution that depends on a parameter. 
Set $\tilde{f}(k,l):=f(k-1,l-1),$ $(k,l)\in \Z^2.$ Obviously, $\tilde{f}$ is $1$-polyanalytic if and only if $f$ is $1$-polyanalytic.
we can write the defining equations of $f$ being $1$-polyanalytic of the second kind according to,
\begin{equation}
(E_1^2E_2-I+iE_2^2E_1-iI)\tilde{f}(k,l)=0,\quad \forall (k,l)\in \Z^2
\end{equation}
Set $\phi(E_1,E_2):=E_1^2E_2-I+iE_2^2E_1-iI.$
To find a particular solution we look for those of the form $\lambda^k\mu^l,$ and we consider the equation,
$\phi(\lambda,\mu)=0,$ i.e.,
\begin{equation}
\lambda^2\mu +i\mu^2\lambda-1-i=0
\end{equation}
Denote by $\lambda_j(\mu)$ the two roots of this equation, which yields the two particular solutions
$(\lambda_j(\mu))^k\mu^l, j=1,2.$
By linearity the sum of all such expressions for all possible values of $\mu$ will also be solutions.
Let
$D_j(\mu)$ be arbitrary functions of $\mu,$ $j=1,2.$
This yields two independent solutions ($1$-polyanalytic functions of the second kind)
\begin{equation}\label{trak}
\tilde{f}_j(k,l) :=\int_{(\re\mu,\im \mu)\in \R^2} D_j(\mu) (\lambda_j(\mu))^k\mu^l d\re \mu d\im \mu, \quad j=1,2
\end{equation}
However, from a practical point of view, equation \ref{trak} are rather intractable and
intangible and more work is required to determine whether every $1$-polyanalytic function of the second kind
corresponds to such a sum.
We can instead, as in the case of $1$-polyanalytic functions of the first kind,
consider minimal determining sets (minimal sets of uniqueness) in order to describe the kernel of $L_2.$

\begin{observation}\label{secondkind1}
	Let $f$ be $1$-polyanalytic of the second kind on $\Z[i]$, let $z=x+iy$ denote the standard coordinate in $\Z[i],$ and let $p_0\in \Z[i].$ Then
	$f$ is uniquely determined by its values on $D_0:=\{ z\in \Z[i]\colon (\re z-\re p_0)(\im z-\im p_0)=0\}$ 
	Also any proper subset of $D'\subset D$ is not a set of uniqueness (in the sense that there are
	two different $1$-polyanalytic functions on $\Z[i]$ that agree on $D'$). 
\end{observation}
\begin{proof}
	We can easily use the procedure
	in Example \ref{anv} as follows. Start from the set, $D_0$. Define iteratively the sets $D_j, j\in \N,$
	by letting $D_{j+1}$ be the set of all points $s\in \Z[i]$ 
	satisfying
	$s\in \{ w+1,w-1,w+i,w-i\}$ for some $w\in \Z[i]$ such that
	precisely 1 point of the set
	$\{ w,w+1,w+1+i,w+i\}$ does not belong to $D_j.$ Then $\bigcup_{j} D_j=\Z[i]$. $f$ can be iteratively extended to each $D_j$ by assigning the value of $f$
	at $s\in D_j\setminus D_{j-1}$ to be determined by the equation $f(w+1)+if(w+1)-f(w+1+i)-if(w+i)=0.$ 
	Obviously replacing the value of $f$ at a point of $D_0$ yields a different extension to $\Z[i].$
	This completes the proof.
\end{proof}

\subsection{Third kind, $q=1$}
Let $f(k,l)=u(k,l)+iv(k,l)$ be a $1$-polyanalytic function of the third kind (where $u$ and $v$ are the real and imaginary parts of $f$ respectively). 
We will start by explaining from a theoretical point of view how one could go about solving for the $1$-polyanalytic function of the third kind, however as will be clear this approach can be non-tractable in practice which is why we then give a description of the kernel of $L_3$ also in terms of minimal sets of uniqueness.
\\
For the theoretical perspective we start by noting that
$L_3f(k,l)=0$ is equivalent to the pair of equations
\begin{equation}\label{oone}
u(k+1,l)-u(k-1,l)=v(k,l+1)-v(k,l-1)
\end{equation}
\begin{equation}\label{twoo}
u(k,l+1)-u(k,l-1)=v(k-1,l)-v(k+1,l)
\end{equation}

Now equation \ref{twoo} yields
\begin{equation}\label{twoo10}
v(k+1,l)=v(k-1,l) -u(k,l+1)+u(k,l-1)
\end{equation}
If we replace $(k,l)$ by $(k-1,l+1)$ and $(k-1,l-1)$ respectively in equation \ref{twoo10} we get 
the pair of equations
\begin{equation}\label{twoo2}
v(k,l+1)=v(k-2,l+1) -u(k-1,l+2)+u(k-1,l)
\end{equation}
\begin{equation}\label{twoo3}
v(k,l-1)=v(k-2,l-1) -u(k-1,l)+u(k-1,l-2)
\end{equation}

Now equations \ref{twoo2} and \ref{twoo3} combined with equation \ref{oone} yield
\begin{multline}
u(k+1,l)-u(k-1,l)=v(k-2,l+1)-
v(k-2,l-1)+\\
u(k-1,l-2)-
u(k-1,l+2)+2u(k-1,l)
\end{multline}
i.e.\
\begin{multline}\label{kjh00}
u(k+1,l)-u(k-1,l)=u(k-3,l)-
u(k-1,l)+\\
u(k-1,l-2)-u(k-1,l+2)+2u(k-1,l)
\end{multline}
Hence a partial difference equation for a real-valued function in two variables.
To solve it we shall use Laplace's method of generating
functions (see Jordan \cite{jordan}, p.607).
First let 
\begin{equation}\label{backeq}
\tilde{u}(k,l):=u(k-3,l-2)
\end{equation}
so that equation \ref{kjh00} becomes
\begin{multline}\label{kjh}
\tilde{u}(k+4,l+2)-
\tilde{u}(k+2,l+2)-
\tilde{u}(k,l+2)+
\tilde{u}(k+2,l+2)-\\
\tilde{u}(k+2,l)+\tilde{u}(k+2,l+4)-2\tilde{u}(k+2,l+2)=0
\end{multline}
Let $a_{r,s}$ denote the coefficient of $\tilde{u}(k+r,l+s)$ in equation \ref{kjh}, e.g.\
$a_{4,2}=1,a_{2,2}=-1$, etc.
Denote the generating function of $\tilde{u}$ by, 
\begin{equation}
U(t,t_1):=\sum_{l=0}^{\infty} \sum_{k=0}^{\infty} \tilde{u}(k,l) t^k t_1^l
\end{equation}

Then we can deduced (see Jordan \cite{jordan}, p.608) that, using the notation 
\begin{equation}
w(k,t_1):=\sum_{l=0}^{\infty} \tilde{u}(k,l) t_1^l
\end{equation}
we have 
\begin{multline}\label{shoto}
\sum_{s=0}^{4} \sum_{r=0}^{4} a_{r,s} t_1^{4-s} \left(w(k+r,t_1)-t_1^0
\tilde{u}(k+r,0)-
t_1^1\tilde{u}(k+r,1)-\cdots
\right.
-\\
\left.
t_1^{s-1}\tilde{u}(k+r,s-1)  \right) =0
\end{multline}
Now equation \ref{shoto} is a linear difference equation in $k$ (for the function $w(k,t_1)$ with $t_1$ fixed) with constant coefficients (in the sense that they are
independent 
of the variable $k$) of order $4$ and it contains already $4$ arbitrary functions 
of $k$ 
\begin{equation}\label{shoto1}
\phi_{j}(k):=\tilde{u}(k,j-1),\quad j= 1,2,3,4
\end{equation}
To be clear set
\begin{multline}
K(t_1,k):=\sum_{s=0}^{4} \sum_{r=0}^{4} a_{r,s} t_1^{4-s} \left(t_1^0
\tilde{u}(k+r,0)+t_1^1\tilde{u}(k+r,1)+\cdots
\right.
+\\
\left.
t_1^{s-1}\tilde{u}(k+r,s-1)  \right)
\end{multline}
\begin{equation}
A(t_1):=\sum_{s=0}^{4} a_{4,s} t_1^{4-s}, \quad 
B(t_1,r):=-\sum_{s=0}^{4} a_{r,s} t_1^{4-s} 
\end{equation}
and write equation \ref{shoto} as
\begin{equation}
A(t_1) w(k+4,t_1)+ \sum_{r=0}^{3} 
B(t_1,r) w(k+r,t_1) =K(t_1,k)
\end{equation}
The general solution to such an equation is given by
the sum $S_{h,t_1}+S_{p,t_1}$ where $S_{h,t_1}$ is the general solution to the homogeneous problem
($K$ replaced by $0$) and $S_{p,t_1}$ is a particular solution.
Since the coefficients of the homogeneous equation are independent of $k$, the 
$S_{h,t_1}$ can be obtained 
via the roots of the characteristic equation
\begin{equation}
\theta(\lambda):=A(t_1) \lambda^4 + \sum_{r=0}^{3} 
B(t_1,r) \lambda^r =0
\end{equation}
If there are $\kappa$ different roots, say $\lambda_1,\ldots,\lambda_{\kappa}$ of multiplicities
$\tau_1,\ldots,\tau_{\kappa}$ then
\begin{equation}
S_{h,t_1}(k)=\sum_{\sigma=1}^{\kappa} 
\left(H_{\sigma,1}(t_1)+H_{\sigma,2}(t_1)k+\cdots +H_{\sigma,\tau_{\kappa}}(t_1)k^{\tau_{\kappa}-1} \right)
\end{equation}
where the $H_{\sigma,\tau}(t_1)$ are arbitrary functions of $t_1$.
A particular solution is usually found by an ansatz e.g.\ using the
method of Section \ref{secete}.
The expansion of $w(k,t_1)$ into a power-series in $t_1$ will yield $\tilde{u}(k,l)$, the arbitrary functions of $t_1$ will
after expansion yield $4$ arbitrary functions of $l$.
In this way one could determine all possible solutions, $\tilde{u}(k,l),$ to equation \ref{kjh}.
This in turn yields the starting real part $u$ via equation \ref{backeq}.
Finally, given $u$ we have the following.
\begin{observation}
	Let $f$ be $1$-polyanalytic of the third kind, $u:=\re f, v:=\im f.$ Then
	$f$ is uniquely determined by the values of $u$ on $\Z[i]$ 
	together with the 
	values of $v$ on a set of the form $\{p_0,p_0+1,p_0+i,p_0+1+i\},$ for some point $p_0\in \Z[i].$ 
\end{observation}
\begin{proof}
	Let $z$ denote the standard coordinate in $\Z[i].$ 
	By equation \ref{oone} and \ref{twoo} we have,
	\begin{equation}\label{oone1}
	v(q_0+i)=-v(q_0-i)-(u(q_0+1)-u(q_0-1))
	\end{equation}
	\begin{equation}\label{twoo1}
	v(q_0+1)=v(q_0-1)-(u(q_0+1)-u(q_0-1))
	\end{equation}
	Hence having the two values of $v$ at $p_0=q_0-1$ and $p_0+1=q_0$, we can obtain $v$
	on the set $S_1:=\Z[i]\cap \{ \im z=\im p_0\},$ via equation \ref{twoo1}. Similarly, we obtain $v$
	on the set $S_1:=\Z[i]\cap \{ \im z=\im p_0 +1\},$ via the values of $v$ at $p_0+i$ and $p_0+1+i.$
	Analogously, given $v$ at the two adjacent points $p_0+i$ and $p_0,$
	equation \ref{oone1} yields $v$ on the set $S_3:=\Z[i]\cap \{ \re z=\re p_0\},$
	whereas the two adjacent points $p_0+1+i$ and $p_0+1,$
	yields $v$ on the set $S_4:=\Z[i]\cap \{ \re z=\re p_0+1\}.$
	This process can now be iterated for each subset of $\bigcup_{j=1}^4 S_j,$ of the form 
	$\{w_0,w0+1,w_0+i,w_0+1+i\},$ for some point $w_0\in \Z[i].$
	This completes the proof.
\end{proof}

Now, for the sake of practicality, we can also for the $1$-polyanalytic functions of the third kind,
consider sets of uniqueness in order to describe the kernel of $L_3.$

\begin{observation}\label{thirdkind1}
	Let $f$ be $1$-polyanalytic of the third kind on $\Z[i]$, let $z=x+iy$ denote the standard coordinate in $\Z[i],$ and let $p_0, q_0\in \Z[i]$ 
	such that $p_0$ ($q_0$) is zig-zag even (odd).
	Then
	$f$ is uniquely determined by its values on $D_0:=D_0^{\mbox{even}}\cup D_0^{\mbox{odd}},$
	where $D_0^{\mbox{even}}$ can be either 
	$\{p_0\}\cup \{p_0+\sum_{j=0}^k ((-1)^j +i), k\in \Z\}$ or $\{p_0\}\cup \{p_0+\sum_{j=0}^k (1+(-1)^j i), k\in \Z\}$
	and $D_0^{\mbox{odd}}$ can be either 
	$\{q_0\}\cup \{q_0+\sum_{j=0}^k ((-1)^j +i), k\in \Z\}$ or $\{q_0\}\cup \{q_0+\sum_{j=0}^k (1+(-1)^j i), k\in \Z\}$
\end{observation}
\begin{proof}
	Start from the set, $D_0$. Define iteratively the sets $D_j, j\in \N,$
	by letting $D_{j+1}$ be the set of all points $s\in \Z[i]$ satisfying
	$s\in \{ w+1,w-1,w+i,w-i\}$ for some $w\in \Z[i]$ such that
	precisely 3 points of the set
	$\{ w+1,w-1,w+i,w-i\}$ belongs to $D_j.$ $\bigcup_j D_j=\Z[i]$ and $f$ can be iteratively extended to each $D_j$ by assigning the value of $f$
	at $s\in D_j\setminus D_{j-1}$ to be determined by the equation $f(w+1)-f(w-1)+if(w+i)-if(w-i)=0.$ 
	Obviously replacing the value of $f$ at a point of $D_0$ yields a different extension to $\Z[i].$
	This completes the proof.
\end{proof}

\section{Some natural multiplications and pseudo-polynomials}\label{polynom}

We believe that from an algebraic perspective it makes sense that 
some notion of multiplication is used such that the corresponding notion of (pseudo-)polynomial will be $1$-polyanalytic, 
and that the multiplication would need to be given
by a binary relation, distributive over addition. 
Obviously being both left and right distributive as well as associative and abelian would be satisfying properties as well, however
such requirements are overly restrictive given the circumstances.
Multiplicative structures that are non-associative do occur in modern research but non-distributive ones seem to be rare
which is why we have
required, at minimum, distributivity.
Natural analogues of polynomials in the function spaces that satisfy our above requirements, more or less exist in the literature, for the case of 
monodiffric functions of the first and second kind already. We shall similarly introduce such analogues for 
$1$-polyanalytic functions of the third kind.

\subsection{First kind}
The analogues of polynomials that we endorse, when it comes to multiplication in the function spaces, are finite linear combinations of what is called pseudo-powers, $z^{(j)},$ $j\in \N.$ 
In  the case of $1$-polyanalytic functions of the first kind, these where introduced by Isaacs \cite{isaacs2} using a multiplication, 
between the coordinate function $x+iy\mapsto x+iy,$ and another complex function $f\colon \Z[i]\to \C,$ according to,
\begin{equation}
(x+iy)\odot_1 f(x+iy) :=xf((x-1)+iy)+iyf(x+i(y-1)) 
\end{equation}
and $c\odot_1 f:= cf$ for constants $c\in \C.$

For two complex-valued functions $f,g$ on $\Z[i]$ we define
\begin{multline}
(g\odot_1 f)(x+iy) :=\re g(x+iy) f((x-1)+iy)+\\i\im g(x+iy)f(x+i(y-1)) 
\end{multline}
and $c\odot_1 f:= cf$ for constants $c\in \C.$

\begin{example}
	We have $z^2=(x^2-y^2) +2ixy,$ so that $z\odot z^2=x((x-1)^3-(x-1)y^2 +2iy(x-1))+iy(x^2-(y-1)^2)+2ix(y-1))$
	whereas
	$z^2\odot_1 z=(x^2-y^2)(x-1+iy)+i2xy(x+i(y-1))$ which implies
	$z\odot_1 z^2(2)=2-4i
	\neq 4 = z^2\odot_1 z.$
	Hence the multiplication $\odot_1$ is not commutative.
	Furthermore, note that
	$z\odot_1 z =x(x-1+iy) +iy(x+i(y-1))=
	(x^2-x+y-y^2) +i2xy,$
	thus $(z\odot_1 z)\odot_1 z =(x^2-x+y-y^2)(x-1 +iy) +i2xy(x+i(y-1))$
	whereas 
	$z\odot_1 (z\odot_1 z)=x((x-1)^2-x+1+y-y^2) +i2(x-1)y) +iy((x^2-x+y-1-(y-1)^2) +i2x(y-1)),$
	giving
	$(z\odot_1 z)\odot_1 z (-1)=-4\neq
	-1(4+1+1) +i2=z\odot_1 (z\odot_1 z)(-1).$
	Hence $\odot_1$ is non-associative.
\end{example}

The multiplication $\odot_1$ is obviously distributive over addition but as we have seen not abelian, and therefore yields two different kinds of pseudo-powers which we shall call 
{\em left} and {\em right} polynomials respectively.
The left pseudo-monomials are defined recursively according to
\begin{equation}
z^{0,\odot_1,l}:=1,\quad z^{j+1,\odot_1,l}=z \odot_{1} z^{j,\odot_1,l}, \quad j\in \Z_+
\end{equation}
whereas the right pseudo-powers are defined according to
\begin{equation}
z^{0,\odot_1,r}:=1,\quad z^{j+1,\odot_1,r}= z^{j,\odot_1,r} \odot_{1} z, \quad j\in \Z_+
\end{equation}

By distributivity we obtain natural extension to left (right) {\em pseudo-polynomials} $P_{1,\mbox{left}}$ ($P_{1,\mbox{right}}$) of degree $N$, according to
\begin{equation}
P_{1,\mbox{left}}(z)=\sum_{j=0}^{N} c_j z^{j,\odot_1,l},\quad  (P_{1,\mbox{right}}(z)=\sum_{j=0}^{N} c_jz^{j,\odot_1,r})
\end{equation}
where the $c_j$ are complex constants.

We call complex multiples of pseudo-powers, {\em pseudo-monomials}.
\begin{proposition}
	Let $f$ be a $1$-polyanalytic function of the first kind on $\Z[i]$. Then
	$z\odot_1 f(z)$ is a $1$-polyanalytic function of the first kind on $\Z[i].$
\end{proposition}
\begin{proof}
	Denote $x=\re z, y=\im z.$
	We have $z\odot_1 f(z)=xf(x-1+iy)+iyf(x+i(y-1)).$ Thus
	\begin{multline}
	L(z\odot_1 f(z))=
	(x+1)f(z)-f(x-1+iy)+iyf(x+1+
	i(y-1))-\\
	iyf(x+i(y-1)
	+i[xf(x-1+i(y+1)-xf(x-1+iy)+
	i(y+1)f(z)-
	\\
	iyf(x+i(y-1)]=
	xf(z)+f(z)-xf(z-1)+iyf(f+1-i)-
	\\
	iyf(z-i)+
	ixf(z-1+i)-ixf(z-1)-yf(z)-
	f(z)+\\
	yf(z-i)=
	y(-f(z)+f(z-i)+i(f(z+1-i)-f(z-i))
	+\\
	x(f(z)-f(z-1)+
	i(f(z-1+i)-f(z-1))=\\
	-\frac{y}{i} L_1 f(z-i) +x L_1 f(z-1) =0
	\end{multline}
	This completes the proof.
\end{proof}
Since the pseudo-powers $z^{j,\odot_1,l}$ are defined iteratively and obviously $f(z)=z$ is $1$-polyanalytic of the first kind, we immediately have the following.
\begin{corollary}
	The pseudo-power $z^{j,\odot_1,l}$ is, for each $j\in \N,$ a $1$-polyanalytic function of the first kind.
\end{corollary}

\subsection{Second kind}
Also in the case of $1$-polyanalytic functions of the second kind, there exists in previous literature 
a natural multiplication (see e.g.\ Duffin \& Petersson \cite{duffin2}, p.626) that we endorse.  
We define the following binary relation, $\odot_2,$ on the space of complex valued functions on $\Z[i]$, 
\begin{equation}
(g \odot_2 f) (z) := \int_{\gamma_{g(z)}}  f(w) dw
\end{equation}
where $\gamma_{g(z)}$ is a 4-curve with initial point $a_0=0$ and end point $a_N=g(z)$ for some positive integer $N.$
As a direct consequence of equation \ref{dubb} (path-independence), the number
$(g \odot_2 f) (z)$ is independent of the choice of different, non-self-intersecting 4-curves, $\gamma,$ that share initial and end point.
Clearly, $\odot_2$ is distributive over addition. It is however not abelian or associative.
\begin{example}
	That $\odot_2$ is not abelian take $f(z)=z^2,g(z)=z$ and calculate
	$2(g \odot_2 f) (-1)=(1+0)(-1-0)=-1,$
	whereas $2(f \odot_2 g) (-1)=(1-0)(1-0)=1.$
	It is not associative which can be seen by 
	setting $A(z):=2(g \odot_2 g)(z),$ and
	noting that $A(-1)=(-1+0)(-1-0)=1,$ $A(0)=0,$ and
	$((g \odot_2 g) \odot_2 g) (-1)=(1-0)(1-0)=1,$
	whereas
	$(g \odot_2 (g \odot_2 g)) (-1)=(A(-1)+A(0))(-1-0)=-1.$
\end{example}

The multiplication $\odot_2$ therefore yields two different kinds of candidates for so called {\em pseudo-monomials},
namely we introduce the left pseudo-powers
\begin{equation}
z^{(0),l}:=1,\quad z^{(j+1),l}=z \odot_{3} z^{(j),l}, j\in \Z_+
\end{equation}
whereas the right pseudo-powers are defined according to
\begin{equation}
z^{(0),r}:=1,\quad z^{(j+1),r}= z^{(j),r} \odot_{3} z, j\in \Z_+
\end{equation}
Up to multiplication by $j$, the left pseudo-monomials can be found in e.g.\ Duffin \& Petersson \cite{duffin2}, 
p.626, there denoted $z^{(j)}$ (which means that $jz^{(j+1),l}=:z^{(j)},$ for $j>0$) and they are known to be $1$-polyanalytic of the second kind,
see e.g.\ Theorem 2.6, Duffin \& Petersson \cite{duffin2} 
(note that the
the reason the multiple $j$ does not appear in our definition is that our
definition arises as a consequence of a more general multiplication whereas the $z^{(j)}$ are stand-alone definitions).
It is precisely the {\em left} pseudomonimials that we endorse (as natural analogues of powers of $z$ in the case of holomorphic functions)
in the case of $1$-polyanalytic functions of the second kind.

\subsection{Third kind}
Let $v_1$ be a zig-zag even point and let $v_1$ ($v_2$) be a zig-zag even (odd) point in $\Z[i].$
Let  $\Gamma_{0,v_1}(\Gamma_{1,v_2})$ be a zig-zag polygon staring at $0 (1)$ with end point $v_1(v_2).$
Obviously for any zig-zag polygon $\Gamma^+$($\Gamma^{-}$) from a point $a$(b) to a point $b$(a)
we have for any complex-valued function $f$ on $\Z[i],$
$\int_{\Gamma^-} f(z) dz = -\int_{\Gamma^+} f(z) dz.$
As a consequence of Proposition \ref{pathind} (zig-zag path-independence) the numbers
\begin{equation}
\mbox{(even)}\int_{0}^{v_1} f(w)dw:=\int_{\Gamma_{0,v_1}} f(w) dw ,\, \mbox{(odd)}\int_{1}^{v_2} f(w)dw:= \int_{\Gamma_{1,v_2}} f(w)dw 
\end{equation}
are independent of the choice of zig-zag polygon  $\Gamma_{0,v_1}$ ($\Gamma_{1,v_2}$) as long as they are zig-zag polygons 
starting at $0 (1)$ with end point $v_1(v_2).$
Define, for each function $f\colon \Z[i] \to \C,$ the number 
\begin{equation}
\int_{\mbox{zig-zag},z} f(w) dw :=\left\{
\begin{array}{cr}
\mbox{(even)}\int_{0}^z f(w) dw & ,\mbox{ if } z \mbox{ is zig-zag even}\\
\mbox{(odd)}\int_{1}^z f(w) dw  & ,\mbox{ if } z \mbox{ is zig-zag odd}
\end{array}
\right.
\end{equation}
Obviously, by zig-zag path-independence we can now define integration along an arbitrary zig-zag polygon $\Gamma$, with starting point $a_0$ 
and end point $a_N$ (for a positive integer $N$)
namely 
$\int_{\Gamma}  f(w) dw = \int_{\mbox{zig-zag},a_N}  f(w) dw-\int_{\mbox{zig-zag},a_0}  f(w) dw.$

\begin{proposition}\label{hofr}
	The function $F(z):=\int_{\mbox{zig-zag},z} f(w)dw$ is $1$-polyanalytic of the third kind whenever $f$ is.
\end{proposition}
\begin{proof}
	Without loss of generality, assume $z$ is zig-zag even.
	we chose paths from $z-1$ to $z+1$ namely $(z-1,z- i,z+1),$ and from $z-i$ to $z+i$ 
	$(z-i,z- 1,z+i),$ and as a consequence of Proposition \ref{pathind} (zig-zag path-independence) we can write
	\begin{multline}
	2 \int_{\mbox{zig-zag},z}  f(w) dw = 
	2\left(\int_{\mbox{zig-zag},z+1}  f(w) dw
	-\int_{\mbox{zig-zag},z-1}  f(w) dw\right)
	+\\
	i2\left(\int_{\mbox{zig-zag},z+i}  f(w) dw
	-\int_{\mbox{zig-zag},z-i}  f(w) dw\right)
	\\
	=2\left(\int_{\mbox{zig-zag},z+1}  f(w) dw
	-\int_{\mbox{zig-zag},z-1}  f(w) dw\right)
	+\\
	i2\left(\int_{\mbox{zig-zag},z+i}  f(w) dw
	-\int_{\mbox{zig-zag},z-i}  f(w) dw\right) 
	=\\
	(f(z+1)+f(z-i))((z+1)-(z-i))+
	(f(z-i)+f(z-1))((z-i)-(z-1))
	+\\
	i(f(z+i)+f(z-1))((z+i)-(z-1))+
	i(f(z-1)+f(z-i))((z-1)-(z-i))
	=\\
	(f(z+1)+f(z-i))(1+i)+
	(f(z-i)+f(z-1))(1-i)
	+\\
	i(f(z+i)+f(z-1))(1+i)+
	+i(f(z-1)+f(z-i))(-1+i)
	=\\
	(1+i)f(z+1)
	+(-i-1)f(z-1)
	+(i-1)f(z+i)
	+(1-i)f(z-i))
	=\\
	(1+i)(f(z+1)-f(z-1))
	+i(1+i)f(z+i)
	-i(1+i)f(z-i)=
	\\
	(1+i)\cdot L_3(f(z)) =0
	\end{multline}
	This completes the proof.
\end{proof}

We define the following binary relation, $\odot_3,$ on the space of complex valued functions on $\Z[i]$, 
\begin{equation}
(g \odot_3 f) (z) := \int_{\mbox{zig-zag},g(z)}  f(w) dw
\end{equation}
Again, $\odot_3$ is distributive over addition but not abelian and not associative (see example \ref{hors}, and therefore yields two different 
(see Example \ref{hors}) kinds of monomials,
\begin{equation}
z^{[0],l}:=1,\quad z^{[j+1],l}=z \odot_{3} z^{[j],l}, j\in \Z_+
\end{equation}
whereas the right monomials are defined according to
\begin{equation}
z^{[0],r}:=1,\quad z^{[j+1],r}= z^{[j],r,1} \odot_{3} z, j\in \Z_+
\end{equation}
By Proposition \ref{hofr} the pseudo-polynomials obtained via the left-monomials are all $q$-polyanalytic of the third kind.
\begin{example}\label{hors}
	To see that $\odot_3$ is non-abelian take $g(z)\odot_3 f(z)$ at the point $z=1+i,$ where $f(z)=z^2,g(z)=z.$
	We have $g(z)\odot_3 f(z) = \int_{\mbox{zig-zag},z}  w^2 dw,$ whereas
	$f(z)\odot_3 g(z) = \int_{\mbox{zig-zag},z^2}  w dw.$ Hence at $z=1+i$ we have $z^2=2i$ so that
	$f(z)\odot_3 g(z) = ((1+i)^2+0)(1+i-0)=(1+i)^2=2i,$ and 
	$g(z)\odot_3 f(z) = \int_{\mbox{zig-zag},2i}  w dw=
	(2i+(1+i))(2i-(1+i))+((1+i)+0)((1+i)-0)=
	(2i+1)(i-1)+(1+i)^2=-3+i
	\neq f(z)\odot_3 g(z).$
	
	To see that $\odot_3$ is non-associative, 
	set $A(z):=g(z)\odot_3 g(z)=\int_{\mbox{zig-zag},z}  w dw .$ 
	Then $(g(z)\odot_3 g(z)) \odot_3 g(z)
	=\int_{\mbox{zig-zag},A(z)}  w dw
	=(A(a_N)+A(a_{N-1}))(a_N-a_{N-1})+\cdots +(A(a_1)+A(a_{0}))(a_1-a_{0}),$
	for a zig-zag polygon with ordered set or vertices $(a_N,\ldots , a_0),$
	where $a_{j+1}\in \{a_j \pm ( 1+i), a_j\pm (1-i)\}.$
	Choose $z=1+i.$
	Then we have a path of integration with only two vertices namely $a_1=1+i, a_0=0.$
	Now
	$g(z)\odot_3 (g(z)  \odot_3 g(z))
	=\int_{\mbox{zig-zag},z}  A(w) dw.$
	Let $\Gamma_{w}:=(b_{M,w},\ldots,b_{0,w}),$ be a zig-zag polygon of minimal length such that
	$b_{M,w}=w.$
	Then 
	$A(w)=(A(b_{M,w})+A(b_{M-1,w}))(b_{M,w}-b_{M-1,w})+\cdots +(A(b_{1,w})+A(b_{0,w}))(b_{1,w}-b_{0,w}),$
	Now we have $A(0)=0, A(1+i)=\int_{\mbox{zig-zag},1+i}  w dw =(1+i)^2=2i.$
	Hence
	$(g\odot_3 (g  \odot_3 g))(1+i)=\int_{\mbox{zig-zag},1+i}  A(w) dw=(A(1+i)+A(0))(1+i-0)=2i(1+i)=2i-2,$
	whereas
	$((g\odot_3 g)  \odot_3 g)(1+i)=\int_{\mbox{zig-zag},A(1+i)}  w dw=
	\int_{\mbox{zig-zag},2i}  w dw=((1+i)+0)((1+i)-0)+(2i+(1+i))(2i-(1+i))=
	(1+i)^2+(3i+1)(i-1)=-4 \neq (g\odot_3 (g  \odot_3 g))(1+i).$
\end{example}

\begin{remark}
	Note the slight difference in notation for the pseudo-monomials $z^{(j),l}$ and $z^{[j],l}$ used to separate between 
	the case of the second and third kind respectively. 
\end{remark}

\subsection{The set of $q$-polyanalytic functions when $q>1$}

Now that we have seen for each of the kernels of $L_1,L_2,L_3$ how to obtain any member
by choosing appropriate values on a minimal set of uniqueness we can use these members in order to solve for
the members of the kernels of $L_1^q,L_2^q,L_3^q$ for $q>1.$ 
Define the convolution of two complex-valued functions $f,g$ on $\Z[i]$ by
\begin{equation}
(f*g)(z) := \sum_{w\in\Z[i]}
f(w)g(z - w), z \in Z[i],
\end{equation}

\begin{observation}\label{wanted}
	Assume that we have an operator 
	$\mathcal{E}_j$ acting on complex-valued functions on $\Z[i]$ and
	satisfying $L_j (\mathcal{E}_j g)(z)= g(z).$
	If 
	$f$ is a $q$-polyanalytic function of the $j$:th kind on $\Z[i]$,
	we have $L_j^q f=0\Rightarrow L_j (L_j^{q-1} f)=0$ so
	$g_1:=L_j^{q-1} f$, is a $1$-polyanalytic function of the $j$:th kind.
	Since $L_j (\mathcal{E}_j g_1)=g_1$
	this means that the function $G_2:=\mathcal{E}_j g_1 +g_2$ is 
	$2$-polyanalytic of the $j$:th kind
	satisfying $L_j^{2} G_2 =g_1,$
	for any function $g_2$ that is $1$-polyanalytic of the $j$:th kind.
	Continuing in this fashion
	we obtain a function
	$G_q = \mathcal{E}_j^{q-1} g_1 +\mathcal{E}_j^{q-2} g_2 +\cdots + \mathcal{E}_j g_{q-1}+ g_q,$
	where the $g_k$ are arbitrary $1$-polyanalytic functions of the $j$:th kind, such that
	$L^q_{j} G_q =0.$
	Hence we can obtain any $q$-polyanalytic function of the $j$:th kind starting from
	$1$-polyanalytic functions of the $j$:th kind and we have in the previous section
	described how all $1$-polyanalytic functions of the $j$:th kind can be determined. 
\end{observation}

Hence in order to determine the $q$-polyanalytic functions of the $j$:th kind for $q>1$ we only need to find
the appropriate operators $\mathcal{E}_j$ associated to the operators $L_j,$ $j=1,2,3$. 

We start with
$\mathcal{E}_1$. It was resolved by Isaacs \cite{isaacs1}, p.194. Define
$B_{+}:=\{z\in \Z[i] \colon 1-\re z\leq \im z\leq 0\}$,
$B_{-}:=\{z\in \Z[i] \colon 1-\im z\leq \re z\leq 0\}$.
Define the operator $Q_1$ as follows:
$Q_+(z) = 0$ when $\re z\leq  0,$ $Q_{+}(1) = 1,$
$Q_{+}(1+iy) = 0$ for $y\neq 0$ and then define recursively $Q_{+}(z)$ for $\re z = p+1, p > 1$,
by $Q_{+}(p + 1 + iy)=(1 + i)Q_{+}(p + iy)-iQ_{+}(p + i(y + 1)).$
Analogously define $Q_{-}$ such that $Q_{-}(z)=0$ for $\im z\leq 0$, $Q_{-}(i)=1,$  
$Q_{-}(x+i) = 0$ for $x\neq 0$ and then
recursively $Q_{-}(z)$ for $\im z = p+1, p>1$,
by 
$Q_{-}(x + i(p+1))=(1-i)Q_{-}(x + ip)+iQ_{-}(x+1+ip).$
Kiselman \cite{kiselman1}, showed in the proof of Theorem 4.2,
that the following operator satisfies the wanted conditions
in Observation  \ref{wanted},
\begin{equation}\label{fund1}
\mathcal{E}_1\colon f\mapsto Q_+ + *(\chi  f) + Q_{-}* (1 - \chi)f
\end{equation}
In the case of the second kind
Let $\chi$ denotes the characteristic function
of the set $A_+:=\{\re z + \im z \geq 0\}$, and set
\begin{equation}
S_{+} (x +iy) = i^{y-x}d(x,y),\quad x +iy \in Z[i]
\end{equation}
where $d(x,y), (x,y) \in \Z^2,$ is defined as $0$
when $x \leq -1$ or when $y \leq -1,$ as $1$ when $(x,y) = (0,0)$, and
for $(x,y) \in \N^2 \setminus \{(0;0) \}$ by the recursion formula
$d(x,y) = d(x -1,y)+d(x -1,y -1)+d(x,y -1).$
Kiselman \cite{kiselman2}, showed that $S_+$ is a fundamental solution supported in
$A_+:=\{(\re z \geq 1) \wedge ( \im z \geq 1) \}$, and pointed out that
there is a natural analogue of $S_+$ but whose support is 
$A_-:=\{(\re z \leq 0) \wedge ( \im z \leq 0) \}$ (instead of $A_+$)
and the existence of which is proved similarly. 
Kiselman \cite{kiselman2}, showed in the proof of Theorem 4.2,
that the following operator satisfies the wanted conditions
in Observation  \ref{wanted},
\begin{equation}\label{fund2}
\mathcal{E}_2\colon f\mapsto S_+ *(\chi  f) + S_{-}* (1 - \chi )f
\end{equation}
\\
\\
Finally, in the case of the third kind, 
define for a complex-valued $f$ on $\Z[i],$
\begin{equation}
L_3' f(z):=f(z + 1) - f(z - 1) - if(z + i) + if(z - i) 
\end{equation}
and denote by
$\Delta_{(2)}$ the operator which acts 
according to 
\begin{multline}
\Delta_{(2)} f(z):=(L_3 \circ L_3') f(z) = f(z + 2) + f(z - 2) +\\
f(z + 2i) + f(z - 2i) - 4f(z)
\end{multline}
We call this operator the {\em two-step} discrete Laplacian.
Define
\begin{equation}\label{fund3}
\mathcal{E}_3\colon f\mapsto f*(L_3'\mathscr{G})
\end{equation}
where we choose $\mathscr{G}$ 
to be the following 
\begin{equation}
\mathscr{G}(m,n):=\frac{1}{2\pi}\int_0^{2\pi} d\psi\left(1-\exp(i(\pm m\phi \pm n\psi))\right)\left(4-2\cos\phi-2\cos \psi \right)^{-1}
\end{equation}
satisfying (see van der Pol \cite{vanderpol})
\begin{equation}
\Delta_{(2)}\mathscr{G}(m+in)=\delta_{(m,n)}
\end{equation}
where $\delta_{(m,n)}=1$ if $m=n=0$ and $\delta_{(m,n)}=0,$ otherwise.
Then in light of $L_3(L_3'\mathscr{G})=\Delta_{(2)}\mathscr{G}$, 
we have $L_3 (\mathcal{E}_3f)(z)= L_3(f*(L_3'\mathscr{G}))(z)= 
\sum_{w\in\Z[i]}
L_3 f(w)(L_3'\mathscr{G})(z - w)
=\sum_{w\in\Z[i]} f(w) (L_3(L_3'\mathscr{G})(z-w))$,
$z\in Z[i].$
Now given a fixed $w\in \Z[i],$ we have 
$(L_3(L_3'\mathscr{G})(z-w))=\Delta_{(2)} \mathscr{G}(z-w)$
and the two-step Laplacian commutes with translation $\tau_w\colon z\mapsto z-w$,
i.e.\ $\Delta_{(2)} \tau_w\circ \mathscr{G}(z) =\tau_w\circ (\Delta_{(2)} \mathscr{G})(z)$.
Thus $\Delta_{(2)} \mathscr{G}(z-w)$
equals $1$ if $z=w$ and $0$ otherwise.
Hence $\mathcal{E}_3$ satisfies the wanted conditions
in Observation \ref{wanted}.
\begin{remark}
	We are considering homogeneous equations
	of the form $L_j^q f=0,$ where $j\in \{1,2,3\}$ and $q$ a fixed positive integer.
	We mention, as a remark, that Tu \cite{shih2}, p.46, presented the following statement:
	Let $n$ be a positive integer, let $c_0,\ldots ,c_{n-1}$ be arbitrary constants and let
	$a_{1}, \cdots, a_{n}$ be distinct roots of
	$a^{n}+c_{n-1}a^{n-1}+\cdots+c_{1}a+c_{0}=0$.
	Then the general solution to $\sum_{j=1}^n c_{j-1} \frac{1}{2^j}L_1^j F=0$ is
	$F(z)= \sum_{i=1}^{n}B_{i}e^{a_{i},z},$
	$B_{i}, i=1,\cdots, n,$ are arbitrary constants and 
	$e^{a,x}=(1+a)^{x}(1+ia)^{y},$ for $z=x+iy$, and $a\in \C.$
	\\
	Clearly
	any equation of the form
	$\frac{1}{2^n}L_1^n F=0$ would have associated to it, 
	$c_{n-1}=1, c_j=0, j<n-1,$ i.e.\
	the equation $a^{n} +1\cdot a^{n-1}=a^{n-1}(a+1)=0$,
	which has $n$ distinct roots only in the case of $n=1$ or $n=2.$
	If $n=1,c_0=0$ then $a_1=0$ is the only root of 
	$a +c_0=0$, and we are considering (up to multiplication by the constant $1/2$) 
	precisely the equation for monodiffric functions of the first kind,
	$L_1 f=0.$ However, the function $\sum_{i=1}^{1}B_{1}e^{0,z}=B_1  1^{x+y}=B_1.$ 
	If $n=2$ then the roots are $a_1=0$ and $a_2=-1,$ thus
	giving the function $F(z)=B_{1}e^{0,z} +B_{2}e^{-1,z}$
	which, for $x\neq 0,$ can be evaluated as
	$B_{1}1^{x+y} +B_{2}\cdot 0\cdot (1-i)^{y}=B_1.$  
	Obviously  
	a complex constant is not the general solution to $L_1f =0$ or $L_1^2 f=0.$
	We conclude that the statement in Tu \cite{shih2} is not meant to apply to equations of the form $L_1^q f=0$ for any positive integer $q$. 
\end{remark}

\section{Pairwise inequivalence of the three kinds}
Here is a proposition that illustrates one way that the minimal sets of uniqueness from
Section \ref{solutionsec} can be useful.

\begin{proposition}\label{inequiv}
	Denote for $j=1,2,3,$ by Ker$(L_j)$ the set of complex-valued functions on $\Z[i]$ that are annihilated by $L_j.$
	Then Ker$(L_k)\setminus$Ker$L_l\neq \emptyset$ for $k\neq l.$
\end{proposition}
\begin{proof}
	By Observation \ref{firstkind1} the set $D=\{ z\colon 0\leq \re z \leq 1\}$ is a minimal determining set
	for functions annihilated by $L_3$.
	Define $f$ on $D$
	according to 
	$f(0)=1, f(1)=-i, f(i+1)=f(i)=0,$
	and $f(z)=0$ otherwise.
	Then we know that $f$ has a unique extension (see Remark \ref{unique}) to $\Z[i]$ which is $1$-polyanalytic of the third kind.
	However, $L_2 f(0)=2\neq 0$ and $L_1 f(0)=-2i-1\neq 0.$
	This takes care of the cases
	Ker$(L_3)\setminus$Ker$(L_1)$
	and Ker$(L_3)\setminus$Ker$(L_2).$
	By Observation \ref{secondkind1}
	we know that
	$D=\{ z\in \Z[i]\colon \re z \im z =0\}$ is a minimal set of uniqueness
	for any function $f$ satisfying $L_2 f\equiv 0.$ 
	On the other hand, by Observation \ref{firstkind1},
	we know that
	$D':=D\setminus \{ z\in \Z[i]\colon (\re z=0)\wedge (\im z > 0) \}$
	is a minimal determining set
	for any function $f$ satisfying $L_1 f\equiv 0.$ 
	Defining a function $g$ on $D$ according to $g(z)=f(z)$ for $z\in D'$
	and $g(z)=f(z)+1$ for $z\in D\setminus D'$, we know that
	$g$ determines uniquely an extension to $\Z[i]$ that is annihilated by 
	$L_2$. On the other hand $g$ cannot be annihilated by $L_1$ because it does not coincide with the unique extension
	from the set $D'$ for functions that are annihilated by $L_1.$
	Furthermore, define a function $h$
	on $D'$ according to $h(-i)=1$ and $h(z)=0$ otherwise.
	Then by Observation \ref{firstkind1} $h$ has unique (see Remark \ref{unique}) extension to $\Z[i]$ that is annihilated by $L_1$.
	However, $L_3 h(0)=-i\neq 0.$
	Furthermore, we know that the extension of $h$ satisfies $g(-i+1)=1+i,$
	hence $L_2 h(-i)=1 +i(1+i)=i\neq 0.$
	This takes care of the cases
	Ker$(L_2)\neq$Ker$(L_1)$, Ker$(L_1)\neq$Ker$(L_3)$ and 
	Ker$(L_1)\neq$Ker$(L_2)$ respectively.
	Finally
	define the function 
	$G(z)$ on $D$ according to $G(1)=1,G(-1)=-1, G(i)=0, G(-i)=0$ and $G(z)=0$ otherwise.
	By Observation \ref{secondkind1}, $G$ determines uniquely an extension to $\Z[i]$ that is annihilated by 
	$L_2$ but by construction $L_3 G(0)=2\neq 0.$ This takes care of the case
	Ker$(L_2)\neq$Ker$(L_3)$,
	This completes the proof.
\end{proof}

\section[Motivation for the third kind]{A motivation for the definition of $q$-polyanalytic functions of the third kind}\label{nysec}

We believe that adjacancy is useful in motivating the defined operators and furthermore some structure such as multiplication is a priori 
not required. 
For this reason we formalize our work using so called
Gaussian structures. Recall that the notation $\Z[i]$ usually implies the {\em ring} of Gaussian integers 
(in particular with a priori given multiplication) and with no graph structure (i.e.\ no adjacancy).

\begin{definition}[Gaussian structure]\label{predef}
	Let $G$ be an additive abelian group. 
	Also equip $G\times G$ with the additive group structure 
	\begin{equation}
	(p_1,p_2)+(q_1,q_2):=(p_1+q_1,p_2+q_2)
	\end{equation}
	$(p_1,p_2),(q_1,q_2)\in G\times G.$
	Define for each $(v_1,v_2)\in G\times G,$ 
	\begin{equation}\mathcal{J}:=(v_1,v_2)\mapsto (-v_2,v_1).
	\end{equation}
	
	Let $G$ also be a directed graph with 
	adjacancy relation $\sim_{G}$. 
	Define an extension, $\sim,$ of the adjacancy relation $\sim_{G},$ by defining for any pair of points
	$p,q\in G\times G$ such that, $p=(p_1,p_2),$ $p\neq q$: $q\sim p$ $\Longleftrightarrow$ $p=q+\mathcal{J}^j((s_1-p_1,0))$ 
	for some $j\in \Z_{\geq 0},$ and some $s_1\sim_{G} p_1.$
	The structure, $\mathcal{G},$ so obtained is called the {\em Gaussian structure induced by $G$}.
	When $G=\Z$, we shall 
	denote the Gaussian structure by $\mathcal{G}_{\Z}.$
\end{definition}

It is clear that letting $G=\Z$ with adjacancy being determined by 
($q\sim p$, $q\neq p$) $\Longleftrightarrow$ ($p\in \{q\pm 1\}$),
and 
$\Z^2$ assumed to have the natural addition induced by
$\Z$, we obtain a Gaussian structure
which aside from its graph properties, can, when equipped with the usual multiplication, be identified with $\Z[i].$
Indeed, we have $G\times G=\Z^2$, and the map $\mathcal{J}$ can be identified with 90 degree clockwise rotation in the plane.
However, we are introducing graph properties (which are not a priori part of the definition of the Gaussian structure induced by $\Z$) 
which in the particular example of $\Z[i],$ implies
$z\sim w,$ and $z\neq w,$ then $z_1=w_1\pm i$ or $z_2=w_2\pm 1$ and
each point has precisely four adjacent points except itself. 
We may obviously introduce multiplication $(z_1,z_2)\cdot (w_1,w_2):=(z_1w_1-z_2w_2,z_1w_2+z_2w_1),$
and thus be able to identify the
Gaussian structure, $\mathcal{G}_{\Z},$ induced by $G=\Z$ with $\Z[i]$ but with additional graph structure as above. 
Note however that we have not introduced a multiplication in our definitions. 

\begin{definition}[$q$-polyanalytic functions of the third kind on Gaussian structures]\label{qanaldef}
	
	Let $q\in \Z_+$ and let $\mathcal{G}$ be a Gaussian structure induced by a group $G$ (in particular we have
	an adjacancy relation $\sim$ on $G\times G$). Since $G$ is directed we can assign to each ordered pair of adjacent points, $s,t,$ $\lambda_{s,t}=1$ 
	($\lambda_{s,t}=-1$)
	if the ordered pair is of positive (negative) direction.
	We define a complex-valued function $f\colon \mathcal{G}\to \C$ to be 
	{\em $q$-polyanalytic of the third kind} at $p\in \mathcal{G}$
	if and only if,
	$L_3 f(p)=0,$ where $L_3f(p)=i\sum_{q\sim p,q_2\neq p_2}f(q)\cdot\lambda_{p_2,q_2}+\sum_{q\sim p,q_1\neq p_1}f(q)\cdot\lambda_{p_1,q_1}.$
\end{definition}

In practice, we shall be working in the case where the inducing group $G$ is $\Z$ and in such cases the other two kinds of
$1$-polyanalytic functions have equally natural formulations.

\begin{definition}
	Let $q\in \Z_+$ and let $\mathcal{G}_{\Z}$ be the Gaussian structure induced by $\Z$. 
	We define a complex-valued function $f\colon \mathcal{G}\to \C$ to be {\ $q$-polyanalytic of the $j$:th kind} at $z\in \mathcal{G}$
	if and only if,
	$L_j^q f(z)=0$, $j=1,2,3,$ 
	where $L_1 f(z):=f(z+1)-f(z)+i(f(z+i)-f(z))$, $L_2 f(z):=f(z+1)-f(z-1)+i(f(z+i)-f(z-i))$, $L_3 f(z):=f(z+1+i)-f(z)+if(z+i)-if(z+1).$
	If the condition holds true at each point of a subset $S\subseteq \mathcal{G}_{\Z}$ where the defining operator is defined, then we say that $f$ is
	$q$-polyanalytic of the $j$:th kind on $S$ and when it is clear from the context what $S$ is we simply say that 
	$f$ is $q$-polyanalytic of the $j$:th kind.
\end{definition}

\begin{definition}[Order of adjacancy]
	Let $G$ be  
	a graph, with
	the adjacancy relation $\sim$. Let $p\in G.$ 
	Denote $\mbox{adj}(p,0):=\{p\},$ and define $\mbox{adj}(p,1)$ as the set of points 
	$\{r\in G\colon r\sim p\}\setminus \{p\}.$
	Iteratively define for each $k\in \Z_+,$ 
	$\mbox{adj}(p,k+1)=\{ z\in G\colon z\sim q \mbox{ for some }q\in \mbox{adj}(p,k)\}\setminus
	\bigcup_{j=1}^{k-1} \mbox{adj}(p,j).$
	The set $\mbox{adj}(p,k)$ will be called the set of points that are {\em 
		adjacent of order $k$ to $p$.}
\end{definition}

From the perspective of graph theory, it may be notable that when applied to Gaussian structures,
the defining operator for $1$-polyanalytic functions of the second kind invokes second order adjacancy when defining a discrete 
analogue of a first order operator and we note that the definition of $1$-polyanalytic functions of the 
first kind does not does not use all first order adjacent point.
In both cases, we find ourselves with rather skewed powers of the given operators in the sense that
the $q$:th power of the operator at a point $z$, will involve points which lie unsymmetrically about $z$. 
This is not the case for the operator appearing in equation \ref{thirdkind}. 
We shall now give yet another motivation
for $1$-polyanalytic functions of the third kind, from the perspective of differential geometry.
\\
\\
Let $M$ be a complex one-dimensional manifold, and let $f\colon M\to \C$ be a differentiable function.
It is well-known that $f$ is holomorphic on $M$ if and only if $df$ is $\C$-linear.
Let $z=x+iy,$ denote the standard complex coordinate for $\C,$ and let $p\in M.$
If $J$ is the complex structure map on $M$ then a basis for $T_{p} M$ is given by $v=\frac{\partial}{\partial x},$ 
$Jv=\frac{\partial}{\partial y}.$ 
Obviously, if $df$ is $\C$-linear then $d_p f(Jv)=id_p f(v),$ for $v=\frac{\partial}{\partial x}.$
Conversely, if 
$d_p f(iv)=id_p f(v)$  
then $d_p f(\lambda v)=\lambda d_p f(v),$ and identifying the complex structure map $J$ with multiplication by $i$,
we see that for all $w\in T_p M,$ 
$d_p f(\lambda w)=\lambda d_p f(w),$
i.e.\ $d_p f$ is $\C$-linear. 

Obviously, the real-linearity of $d_p f$ together with the above implies that
$f$ satisfies the Cauchy-Riemann equations at $p$ if and only if 
\begin{equation} %
d_p f(v)+id_p f(iv)=0, \quad \forall v\in T_p M 
\end{equation}
and by definition $d_pf$ is $\R$-linear so that,
\begin{equation}\label{star1}
2d_p f(v)=d_p f(v)-d_p f(-v), \quad \forall v\in T_p M 
\end{equation}
Hence
\begin{multline}\label{starst}
d_p f\mbox{ is $\C$-linear}\Leftrightarrow \\ 2d_p f(v)+i2d_p f(iv)=0 \Leftrightarrow 
(d_pf(v)-d_p f(-v))+i(d_p f(iv)-d_p f(-iv))=0
\end{multline}
It is an analogue of these equations that we shall use to define 
a symmetric discrete operator whose $q$:th powers will be analogous to the powers $\overline{\partial}^q$.
\\
\\
Recall that if $M$ is an $n$-dimensional smooth real manifold and $p\in M,$ then we can define the set of tangent vectors at
$p$ (or {\em tangent space at $p$}) as the set of vectors $v$ such that
there exists a differentiable curve $\gamma\colon (-\epsilon,\epsilon)\to M,$ some $\epsilon>0,$ $\gamma(0)=p,$
such that $v=\frac{\partial\gamma}{\partial t}(0),$ and acts on the set of differentiable functions, defined on a neighborhood of $p$, according to
$v(g):=\frac{\partial (f\circ\gamma)}{\partial t}(0),$ for differentiable $f\colon U\to \C,$ $p\in U$, $U$ 
an open neighborhood of $p$ in $M.$ 
The tangent space at $p$ is denoted $T_p M.$ Also for differentiable  $f\colon M\to \C,$ we define the differential 
map $d_p f\colon T_p M\to \C,$ as 
$d_{\gamma(0)}f(\frac{\partial\gamma}{\partial t}(0))=\frac{\partial (f\circ\gamma)}{\partial t}(0).$

\begin{definition} 
	\label{ptdef}
	Let $\mathcal{G}$ be a graph and let $p\in \mathcal{G}$. A {\em path $\Gamma$ through $p$ in $G$} 
	is an ordered
	set of points $\Gamma(j)=z_j\in G,$ $j=-m_1,\ldots,m_2,$ for nonnegative integers $m_1,m_2,$ 
	such that $z_j\sim z_{j+1},$ $j=-m_1,\ldots m_2-1$, and
	$p\in \{ \Gamma(j),j=-m_1+1,\ldots,m_2-1\}.$ When the base point is
	not essential to the argument being made we shall 
	simply use the term {\em path in $\mathcal{G}$.} 
	For each $p\in \mathcal{G},$ denote $T_p \mathcal{G}=\{ v\in G\colon v=q-p, \, q\sim p\}.$ This is the set of {\em tangents.}
	Obviously, the cardinality of $T_p \mathcal{G}$ may vary dependent upon the base point $p.$
	Let $f$ be a map $\mathcal{G}\to \mathcal{D},$ 
	for an additive abelian group $\mathcal{D}.$
	For each $p\in \mathcal{G},$
	we have a map $d_p f\colon T_p G\to \mathcal{D},$
	according to $v=(q-p)\mapsto f(q)-f(p).$
	So there exists a path $\Gamma$ containing $p$ and $q$ such that
	$d_p f(v)=f(\Gamma(j_0+1))-f(p)$ where  $\Gamma(j_0)=0.$
\end{definition}
\begin{definition}[$1$-polyanalytic functions of the third on Gaussian structures]\label{solutiondef}
	Let $\mathcal{G}$ be the Gaussian structure induced by $G,$ where $G$ is an additive group. 
	Let $R$ be an additive abelian group 
	and let $f$ be a function $\mathcal{G}\to R^2,$ where $R^2$ is equipped with the componentwise addition.
	$f$ is called a {\em $1$-polyanalytic function of the third kind 
		(with respect to the Gaussian structure $\mathcal{G},$
		at $p$),} if (using the notation of Definition \ref{ptdef}) we have
	\begin{equation}
	d_p f(v)-d_p f(-v)+\mathcal{J}' (d_p f(\mathcal{J} v)-d_p f(-\mathcal{J} v))=0, v\in T_p \mathcal{G}
	\end{equation}
	Where $\mathcal{J}',$ is defined by
	$\mathcal{J}'(A,B)=(-B,A),$ and 
	$\mathcal{J}(v_1,v_2)=(-v_2,v_1)$. 
\end{definition}
From the definitions it is clear that this coincides with the case of $1$-polyanalytic functions of the third kind from Definition \ref{qanaldef},
when e.g.\ $R=\R,$ $G=\Z$. 
	Note that in defining our natural discrete analogue ($L_3$) of the Cauchy-Riemann operator, we have not needed to introduce 
	a multiplicative structure on the domain space (the Gaussian structure), it has been sufficient with 
	a group structure where we on the other hand have required that there exist adjacancy (i.e.\ an additional graph structure).

\chapter{$q$-analyticity in infinite dimensional complex analysis}
A complex polynomial $P$ in $\Cn$ is called homogeneous of degree $m$ if $P(\lambda z)=\lambda^m P(z)$
for all $\lambda\in \C$ $z\in \Cn$. The local (uniformly convergent) power series of any holomorphic function $f$
can be rewritten as a series of $m$-homogeneous polynomials $f(z)=\sum_{m=0}^\infty P_m(z),$
where $P_m$ is homogeneous of degree $m$, and this is called the {\em homogeneous expansion}\index{Homogeneous exapnsion}, 
see e.g.\ Rudin \cite{rudin1980}. If $\phi$ is a linear transformation on $\Cn$
then $P_m\circ \phi$ is again homogeneous of order $m$ and the homogeneous expansion of $f$
is given by $\sum_m (P_m\circ \phi)(z).$ Furthermore, we have for a function $f\in \mathscr{O}(\{\abs{z}<1\}),$
having a homogeneous expansion $f(z)=\sum_{m=0}^\infty P_m(z)$ that for $\zeta\in \{\abs{z}=1\}$
and the function of the variable $\lambda\in \C,$ $\abs{\lambda}<1$ given by
$f_\zeta(\lambda)=\sum_{m=0}^\infty P_m(\zeta)\lambda^m$
is holomorphic with induced coefficients $P_m(\zeta).$ If $\Omega\subset\Cn$ is an open subset, 
$\{P_m(z)\}_{m\in \N}$ a sequence 
where each $P_m$ is $m$-homogeneous and if $\sup_m \abs{P_m(z)}<\infty$ for each $z\in \Omega$ then the series
$\sum_{m=0}^\infty P_m(z)$ converges uniformly on compacts of $\Omega.$
 The notion of homogeneous expansion and the idea of the
restrictions to each complex line being holomorphic have nice generalization to complex Banach spaces.
\\
\\
There is a well-established field of research on infinite dimensional holomorphy, see e.g.\ the books of Dineen \cite{din3} and Mujica \cite{muj}, Soo Bong Chae \cite{chae}
 and  Herv\'e \cite{hervier}.
Three basic examples of properties which have been investigated are the problem of approximation of a holomorphic function by entire functions (e.g.\ it is known that a holomorphic function on the unit ball centered at the origin in a complex Banach manifold can be approximated by entire functions on smaller balls centered at the origin, see Lempert \cite{lempapp08}) and the Levi problem, see overview in Dineen \cite{din3} and the $\overline{\partial}$-problem, see Lempert \cite{lempertI}, \cite{lempertII}, \cite{lempertIII}. 
In this text we shall always assume
any function which is called holomorphic
to be {\em locally bounded}. 
We begin by giving the basic notions of infinite dimensional holomorphy.

\section{Some background on infinite dimensional holomorphy}\label{infoloapp1}
\subsection{Differentiability and polynomials}
Let $X$ be a complex Banach space. By a {\em domain} $\Omega\subset X$ we mean an open connected set.
\begin{definition}[Continuously differentiable mapping]\index{$C^1$ mapping between locally convex spaces}
Let $X,Y$ be two locally convex spaces (e.g.\ Banach spaces), $\Omega\subset X,\Omega$ open, and $f:\Omega \to Y.$ We say that $f\in C^1$ (or $f\in C^1(\Omega ,Y)$) if
\begin{equation}
df(x,v) =\lim_{\R\ni t\to 0} \frac{f(x+tv)-f(x)}{t},
\end{equation}
exists for all $(x,v)\in \Omega\times X$ (in particular $df$ is $\R$-linear).
\end{definition}
We recall the definition of {\em Fréchet derivative}\index{Fréchet derivative}
\begin{definition}[See e.g.\ Mujica \cite{muj}, p.99] 
Let $X,Y$ be separable Banach spaces and $\Omega\subset X,$ an open set and let $h\in X.$. A mapping $u:X\to Y,$ is called {\em Fr\'echet differentiable} at $p\in X$ 
if there is a bounded linear operator
$A :X\to Y$ such that
\begin{equation}
\lim_{\norm{x-p}\to 0} \frac{\norm{u(p )-u(p)-A(x-p)}_Y}{\norm{x-p}_X} =0.
\end{equation}
Sometimes the notation $Df$ is used in the literature for
this so-called Fr\'echet derivative.
\end{definition}
We denote by $cs(X)$ the set of all continuous semi-norms on a topological vector space $X$.
Let $X$ and $Y$ be locally convex vector spaces. Denote by $\mathscr{L}(^m X,Y)$
the space of $m$-linear mappings from $X^m$ (the product space) to $Y,$
and we denote by 
$\mathscr{L}_s (^m X,Y)$ the vector space of all mappings in
$\mathscr{L} (^m X,Y)$ which are symmetric. To every $\phi\in
\mathscr{L} (^m X,Y)$ (where we do not assume continuity, thus when $Y$ is a scalar field, this is a subset of the algebraic dual) we associate
a mapping $\hat{\phi}$ defined by $\hat{\phi}:=\phi\cdot x^m ,$ and call
$\hat{\phi}$ the $m$-homogeneous polynomial associated to $\phi .$ 
Denote by $\mathcal{P}(^m X,Y)$ the sub-vector space of continuous $m$-homogeneous polynomials. Then the linear mapping from the subspace of 
continuous functions $\phi\in \mathscr{L} (^m X,Y)$ to $\mathcal{P}(^m X,Y) ,$
defined by $\phi\mapsto\hat{\phi},$ is surjective. Furthermore
the linear mapping from the subspace of 
continuous functions in $\mathscr{L}^s(^m X,Y)$ to $\mathcal{P}(^m X,Y) ,$
defined by $\phi\mapsto\hat{\phi},$ is bijective.
\begin{definition}\label{jack1}
By a polynomial, we mean a finite sum of elements in 
$\bigcup_m
\mathcal{P}(^m X,Y),$ and we will be considering mainly $Y=\C,$ and
and the set of ($\C$-valued continuous) polynomials on $X$ is denoted $\mathcal{P}(X).$ 
\end{definition}
One usually introduces the norm $\mathcal{P}(^m X,Y)$ to be given by $\norm{P}:=\sup_{\norm{x}\leq 1} \norm{P(x)}.$
The correspondence $P\leftrightarrow \Check{P}$ establishes an isometric isomorphism.
Note that for each $k\in \N,$ if the differential $d^k f(p)$ exists then corresponds to it a natural 
$k$-homogeneous polynomial which we denote $\widehat{d^k f(p)}.$
Specifically for $m$-homogeneous polynomials we denote for $j=0,\ldots,m$
\begin{equation}
\frac{\widehat{d^j P}}{j!}(x): X\ni y\mapsto \binom{m}{j}\Check{P}(x^{m-j},y^j)\in Y
\end{equation}
\begin{equation}
\frac{\widehat{d^j P}}{j!}: X\ni x\mapsto \frac{\widehat{d^j P}}{j!}(x)\in \mathcal{P}(^m X,Y)
\end{equation}
where we have
$\frac{\widehat{d^j P}}{j!}(x)\in \mathcal{P}(^m X,Y),$ $\frac{\widehat{d^j P}}{j!}\in 
\mathcal{P}(^{m-j} X,\mathcal{P}(^j X,Y))$
and
\begin{equation}
P(x+y)=\sum_{j=0}^m \left(\frac{\widehat{d^j P}}{j!}(x)\right) (y)
\end{equation}
is just the Taylor expansion of $P$ about $x.$
Understanding polynomials in infinite dimensional complex analysis is more involved than in the finite dimensional
case. 
\begin{definition}
When the $k$:th Fr\'echet derivative $D^k f(a)$ exists the $k$:th differential of $f$
at $a$ is defined as $d^kf(a)=\widehat{D^k f(a)}$ i.e.\ the $k$-homogeneous polynomial corresponding to the
$k$:th derivative. For evaluation at $v\in X$ we use the notation
$d^k f(a)[v]:=D^k f(a)[^k v].$
For instance we have for $P\in \mathcal{P}(^m X,Y)$ that $P$ is $C^\infty$-smooth and
\begin{equation}
d^k P(x)[v]=m(m-1)\cdots(m-k+1)\Check{P}(^{m-k}x,^kv)
\end{equation}
Thus $D^k P:X\to \mathcal{L}^s(^k X,Y)$ is a polynomial map
from $\mathcal{P}(^{m-k} X,\mathcal{L}^s(^k X,Y)),$ whereas
$d^k P[v]:X\to \mathcal{P}(^k X,Y)$ is a polynomial map
from $\mathcal{P}(^{m-k} X,\mathcal{P}(^k X,Y)).$
\end{definition}
Note that if a function $f:X\to Y$ between Banach spaces over $\mathbb{K}$ ($=\R$ or $\C$), is $k$-times {\em Frech\'et
differentiable} (in the sense that
$D^kf$ exist) then it is $k$-times {\em G\^ateaux differentiable}\index{G\^ateaux differentiable}
in the sense that $\delta^k f\in \mathcal{L}^s(^k X,Y)$ exist in some neighborhood of $x$, where
\begin{equation}
\delta^j f(x)[v_1,\ldots,v_k]=\lim_{t\to 0,t\in \mathbb{K}}\frac{\delta^{j-1} f(x+tv_1)[v_2,\ldots,v_k] -\delta^{j-1} f(x+tv_1)[v_2,\ldots,v_k]}{t}
\end{equation}
exist for all $v_1,\ldots,v_k\in X$ and is a symmetric bounded
$k$-linear mapping, and furthermore
$\delta^k f(x)=D^k f(x),$
see e.g.\ H\`ajek \& Johanis \cite{hajekjohanis}, p.41.
We have the following Taylor formula.
\begin{theorem}
Let $X$ be a normed linear space and $Y$ a Banach space. Let $U\subset X$ be a convex set, $k\in \N$ and $f\in C^k(U,Y).$ Then for any
$x\in U$ and $v\in X$ such that $x+v\in U$ 
\begin{equation}
f(x+v)=\sum_{j=0}^{k-1} \frac{1}{j!} d^j f(x)[v] +\left( \int_0^1 \frac{(1-t)^{k-1}}{(k-1)!}d^k(x+tv)\right)[v]
\end{equation}
and
\begin{multline}
\norm{ f(x+v)-\sum_{j=0}^{k} \frac{1}{j!} d^j f(x)[v]}\leq \\
\frac{1}{k!}\left(
\sup_{t\in [0,1]} \norm{d^k f(x+tv)-d^k f(x)}
\right)\norm{v}^k
\end{multline}
\end{theorem}
\begin{definition}
Let $X$ be a normed linear space and $Y$ a Banach space.
Consider a homogeneous expansion $\sum_{j=0}^\infty P_j(x)$,
$P_j\in \mathcal{P}(^j X,Y).$ The set 
\begin{equation}
\mbox{Int}\{ x\in X: \sum_{j=0}^\infty P_j(x)\mbox{ converges}\}
\end{equation} 
is called the {\em domain of convergence}\index{Domain of convergence} of the series.
The radius of convergence (with respect to the variable $\lambda$) of the series 
$\sum_{j=0}^\infty \norm{P_j}\lambda^j$
is called the {\em radius of norm convergence}\index{Radius of norm convergence} of the series.
\end{definition}
It is known that the domain of convergence is either empty or a balanced neighborhood
and each $P_j\in \mathcal{P}(^j X,Y)$ converges in its domain of convergence absolutely locally
 uniformly and for each $a$ in the domain of convergence there are a neighborhood $V_a$ and $\lambda\in (0,1)$,
$M>0$ such that $\sup_{x\in V_a} \norm{P_j(x)}\leq M\lambda^j$ for each $j\in \N.$
We have the following known result, see H\`ajek \& Johanis \cite{hajekjohanis}, p.67.
\begin{theorem}
Let $X$ be a normed space and $Y$ a Banach space.
If a homogeneous expansion $\sum_{j=0}^\infty P_j(x)$, $P_j\in \mathcal{P}(^m X,Y),$ has positive 
radius of norm convergence $R$ at $x=0,$ 
the map $f:\{x:\norm{x-0}<R\}\to Y$ is $C^\infty$-smooth and
\begin{equation}
d^k f(x)= \sum_{j=0}^\infty d^k P_j(x)
\end{equation}
for $\abs{x}<R$ and $k\in \N.$ Furthermore, the latter series
also has radius of norm convergence $R$ and
\begin{equation}
P_j=\frac{1}{j!} d^j f(0)
\end{equation}
\end{theorem}
Note that if $Q\in \mathcal{P}(^m X,\C)$ for a complex Banach space $X$
with Shauder basis, then the conjugate $\overline{Q(x)}$
satisfies that for each holomorphic function $h$ we have 
$\overline{D} (h\overline{Q(x)})=h \overline{D} (h\overline{Q(x)}).$

\subsection{Bases of spaces of $m$-homogeneous polynomials}
\begin{definition}
Let $X$ and $Y$ be vector spaces. One way to define $m$-homogeneous polynomials is as a map $p:X\to Y$ 
such that there exists an $m$-linear map $M_p:X^m\to Y$
satisfying $p(x)=M_p(x,\ldots,x).$ 
Recall that a sequence is said to converge unconditionally if it convergence independent of reordering.
A {\em Schauder basis} $\{b_j\}_{j\in\N}$ for a Banach space $X$ is a sequence in $X$ such that
for any $x$ there exists a unique set of scalars $\{a_j\}_{j\in \N}$ such that 
 $\lim_{N\to \infty} \norm{x-\sum_{j=0}^N a_j b_j }=0,$
i.e.\ $x=\sum_{j=0}^\infty a_j b_j.$ 
Let $X$ be a Banach space with Schauder basis $\{x_i\}_{i\in \N}.$ The basis
is called {\em shrinking}\index{Shrinking basis} if for all $\phi\in X'$, $\norm{\phi}_n:=
\norm{\phi|_{\overline{\{x_{n+i}:i\in\N\}}} }\to 0$ as $n\to \infty.$ Equivalently a basis is shrinking if
if the sequence of coordinate functionals is a Schauder basis for $X'.$
An element $p\in P(^m X)$ (the space of $m$-homogeneous polynomials) is called a {\em monomial} of degree $m$
if there exists $(i_1,\ldots,i_m),$ $i_j\geq i_{j+1},$ such that $p=p_{i_1,\ldots,i_m}:=x_{i_1}'\cdots x_{i_m}',$
where $\{x_i'\}_{i\in \N}$ denotes the sequence of coordinate functionals.
Denote by $P_f(^m X)$ the set of $m$-homogeneous polynomials of the form $p=\phi^n,$ for $\phi\in X'.$
Denote by $P_w(^m X)$ the closed
subspace of $P(^m X)$ consisting of all $m$-homogeneous polynomials weakly continuous on bounded sets.
An $m$-linear map $p$ on $X$ is called symmetric if $p(x_1,\ldots,x_m)=p(x_{\pi_{(1)}},\ldots,x_{\pi_{(m)}})$
where $\pi$ denotes any permutation of $\{1,\ldots,m\}.$ The map 
\begin{equation}
p\mapsto \Check{p}:=\frac{1}{m!}\sum_{\pi} p(x_{\pi_{(1)}},\ldots,x_{\pi_{(m)}})
\end{equation}
is called the symmetrization map and we call $\Check{p}$ the symmetric $m$-linear map associated to $p.$
The basis $\{x_i\}_{i\in \N}$ is said to be {\em polynomially shrinking}\index{Polynomially shrinking basis}
if for all $p\in P(^m X)$
\begin{equation}
\sup_{\norm{x}\leq 1} \abs{\check{p}\left(\sum_{i=n+1}^\infty x_i'(x)x_i,x^{m-1}\right)}\to 0,\mbox{ as }n\to \infty
\end{equation}
\end{definition}
Note that given vector spaces $X$ and $Y$ and a $m$-homogeneous polynomial
$p:X\to Y$, the $m$-linear map $M_p,$ that gives rise to $p$ is not uniquely determined, in fact
$\Check{p}$ always gives rise to $p.$ However, we have the following so-called polarization formula.
\begin{theorem}[See Bohenblust \& Hille \cite{bohenhille} and Moulis \cite{moulis}]
Let $X$ and $Y$ be vector spaces. For each $m$-homogeneous polynomial $p:X\to Y$ there exists a unique symmetric 
$m$-linear map $\Check{p}$ such that $p(x)=\Check{p}(x,\ldots,x)$ and it satisfies
\begin{equation}
\Check{p}(x,\ldots,x)=\frac{1}{2^{m}m!}\sum_{\epsilon_j=\pm 1} \epsilon_1\cdots \epsilon_m p\left(a+\sum_{j=1}^m \epsilon_j x_j\right)
\end{equation}
where $a\in X$ can be chosen arbitrary. Moreover, if $X,Y$ are normed and $p$ is bounded then
$\norm{p}\leq \norm{\Check{p}}\leq \frac{m^m}{m!}\norm{p}.$
Also for each $n>m$ we have
\begin{equation}
\sum_{\epsilon_j=\pm 1} \epsilon_1\cdots \epsilon_n p\left(a+\sum_{j=1}^n \epsilon_j x_j\right)=0
\end{equation}
\end{theorem}
Pellegrini \cite{pellegrini}, p.274, Thm 2.2, proves the following.
\begin{theorem}
Let $X$ be a Banach space with Schauder basis $\{x_i\}_{i\in \N}.$ The following are equivalent:\\
(i) The basis $\{x_i\}_{i\in \N}$ is polynomially shrinking.\\
(ii) for all $m\in \N$ the monomials of degree $m$ form a Schauder basis for $P(^m X).$\\
(iii) The basis $\{x_i\}_{i\in \N}$ is a shrinking basis and $\overline{P_f(^m X)}=P(^m X),$ for all $m\in \N.$\\
(iv) The basis $\{x_i\}_{i\in \N}$ is a shrinking basis and $P_w(^m X)=P(^m X),$ for all $m\in \N.$
\end{theorem}
See also Dimant \& Zaluendo \cite{dimantzaluendo2}, where they prove an analogous relation between a norm and the
existence of a monomial basis for spaces of multilinear forms over Banach spaces.
We mention also the following notion.
An $m$-homogeneous polynomial is said to be of {\em finite type} if it is a finite linear combination of
$m$-homogeneous polynomials of the form $x\mapsto \phi(x)^m y,$
$\phi\in X',$ $y\in y.$ An $m$-homogeneous polynomial $P$ is called {\em nuclear} if
there exists a bounded sequence $\{\phi_j\}_{j\in \Z_+}\subset X'$ and a sequence $\{c_j\}_{j\in \Z_+}$
such that $P(x)=\sum_j c_j\phi_j(x)^m.$

\subsection{Some notions of holomorphy}
\begin{definition}[See Mujica \cite{muj}, p.33]\label{holomujdef}
Let $\Omega\subset X$ be open and nonempty, $X$ locally convex.
A function $u:X\to Y$ is called {\em holomorphic} 
if
$\forall a\in\Omega, \exists $ a neighborhood $V\subset U,$ and a sequence of polynomials $\{A_m\}_{m\in \N},A_m\in \mathcal{P}(^m U,Y)$
such that,
\begin{equation}
u(x)=\sum_{m=0}^{\infty} A_m(x-a),
\end{equation}
converges uniformly for $x\in V.$ 
\end{definition}
The $A_m$ are of course members of $\mathscr{L}_s (^m X,Y)$ and if
$Y$ is Hausdorff they are uniquely determined by $u.$\\
\\
Now there also exists analogues of the Cauchy-Riemann operator in infinite dimensional holomorphy, see e.g.\
Soraggi \cite{soraggi}, Lempert \cite{lemp1} and Dineen \cite{din3}.
\begin{definition}
Let $X$ be a separable topological complex Banach space, let $\Omega\subset X$ be an open subset
For a map $f:\Omega\to Y$, where $Y$ is a separable Banach space we define for 
$p\in \Omega$ and $v:=(v_1,\ldots, v_k),$ $v_j\in X,$ $j=1,\ldots,k$
\begin{equation}
\overline{D}_{f} (p,v):=\frac{df(p,v)+idf(p,iv)}{2}
\end{equation}
where $df(p,v_j)=\lim_{\R\ni \epsilon\to 0} (f(p+\epsilon v_j)-f(p))/\epsilon.$
For a function $f\in C^1$ (i.e.\ $df(p,v)$ exists for all $p\in \Omega,$ and each vector $v$
and $(p,v)\mapsto df(p,v)$ is continuous), the function is called {\em Fr\'echet holomorphic} if $\overline{D} f\equiv 0$ on $\Omega.$
\end{definition}
It is well-known that these functions locally have Frech\'et differentials of all orders at each point
and they posses locally convergent homogeneous expansion. Furthermore
the homogeneous polynomials in $x-p_0$ in the power series at a point $p_0$,
which form the terms of the power series, are expressible as multiples of the Fr\'echet differentials
at $p_0$ according to the appropriate generalization of the Taylor series, see e.g.\ Chae \cite{chae}, p.392.
Indeed, one has near each point $p_0$, a uniformly convergent (locally near $p_0$)
series $f(x)=\sum_{j=0}^{\infty} P_j(x-\xi)$, $P_j\in \mathcal{P}(^j X,\C),$ and furthermore
(see e.g.\ Chae p.168) $f$ is
$C^\infty$ (i.e.\ Frech\'et differentials of all orders exists) on the ball centered at $p_0$ with radius equal to
the radius of uniform convergence, and
$P_j:=\frac{1}{j!} d^j f(p_0)$, $j\in \N.$

\begin{definition}
A complex valued functional is {\em G\^{a}teaux holomorphic}\index{G\^{a}teaux holomorphic ($G$-holomorphic)} 
(or $G$-holomorphic)
in a domain $\Omega\subset X$ of a complex Banach space $X$ if it is single
valued and its restriction to an arbitrary analytic plane $\{z : z=p_0 +\zeta a\}$
($p_0\in \Omega$, $a\in X$, $\lambda$ a complex parameter) is a holomorphic function of $\zeta$ in
the intersection of the plane with $\Omega$.
\end{definition}

\begin{remark}\label{gateauxremark}
In an infinite dimensional Banach space a $G$-holomorphic function is not necessarily locally bounded.
When the 
difference quotients along each possible direction converge uniformly
then G\^{a}teaux holomorphy at a point implies the existence of the
Fr\'echet derivative at that point. It is known that a function is Fr\'echet differentiable at each point of an open subset
if and only if it is continuous and G\^{a}teaux differentiable at each point, see e.g.\ Chae \cite{chae}, p.392. Also it is known that a function is holomorphic in the sense of Definition \ref{holomujdef} if and only if
it is continuous and $G$-holomorphic.
In fact local boundedness is sufficient, see Bremermann \cite{brem2} p.811.
This can be used to realize that uniform limit of (Fr\'echet) holomorphic functions, in the topology of uniform convergence on compacts
is again holomorphic, see e.g.\ Kim \& Krantz \cite{kimkrantz}, p.196.
\end{remark}
As is pointed out in Dineen \cite{din3}, Hartogs' theorem (in finitely many dimensions) gives that separate holomorphy 
and local boundedness implies holomorphy.
We also point out that there is an extensive theory of smooth analysis on Banach Spaces,
see H\`ajek \& Johanis \cite{hajekjohanis}. For example we have a natural definition of the set of 
real-analytic functions (denoted $C^\omega$)
defined by locally convergent homogeneous expansion (unique Taylor series, where 
homogeneous components can be identified as differentials $d^kf$, see H\`ajek \& Johanis \cite{hajekjohanis}, p.70,
and it is proved that $C^\infty$ functions that are $G$-holomorphic (and equivalently Fr\'echet holomorphic) 
are holomorphic in the
usual sense, i.e.\ have complex homogeneous expansion, see Theorem 160, p.71).

\subsection{Topologies}
Here we give the very basics on topologies mentioned in the text.
First of all the finite topology simply means that open set are the ones satisfying that any finite dimensional slice is open.
\begin{definition}
Let $\{p_{\beta}\}_{\beta}$ be a family of semi-norms (i.e.\ a scalar valued, sub-additive, nonnegative functions) on a vector space $X.$ Then the sets,
\begin{equation}
\{
y\in X |p_{\beta}(x-y)<r
\} ,x\in X,
\end{equation}
{\em generate a topology}\index{Topology generated by semi-norms} (i.e. 
the topology is the smallest topology containing the given sets)
on $X,$ 
which is called the
{\em topology induced by} the family $\{p_{\beta}\}_{\beta} ,$ of semi-norms.
\end{definition}
\begin{definition}
A {\em base for a topology $\tau$}\index{Base for a topology} is a collection $W\subset\tau$ such that,
\begin{equation}
\forall U\in\tau,\forall x\in U\quad \exists V\in W
\mbox{ such that } x\in V\subseteq U.	
\end{equation}
$X$ is called {\em locally convex} if it has a basis consisting of convex sets.
\index{Locally convex}
\end{definition}
We have the following result on semi-norm generated topologies.
\begin{theorem}
If $X$ is a vector space with topology induced from a
family of semi-norms, then $X$ is a locally convex topological vector
space.
\end{theorem}
\begin{definition}
Let $X,Y$ be locally convex spaces and $U\subset Y,$ be open.
{\em The topology on $\mathscr{O}$ of uniform convergence on compact sets,}\index{Topology of uniform convergence on compacts} is the topology generated by
the semi-norms,
\begin{equation} 
p_{\beta, K}(f)=:\norm{f}_{\alpha ,K}=\sup \alpha(f(x)),\quad f\in \mathscr{O}(X),
\end{equation}
where $K$ ranges over the compact subsets of $U$ and $\alpha$ ranges over continuous semi-norms on $Y.$
A basis for this topology is given by the following collection of sets,
\begin{equation} 
B_{g, K,\epsilon } =\{ f\in\mathscr{O} (\Omega ) | \sup_{z\in K} |f(z)-g(z)| <\epsilon \} .
\end{equation}
\end{definition}

\subsection{Some preliminaries on complexification}\label{appapp}
Here we quickly give some facts on the interplay between real and complex vector space structure.
Let $V$ be a real Banach space. Then $V$ is in particular a vector space over $\R,$
and can be complexified, by which we mean $\C\otimes_{\R} V.$ 
If we first consider $\C\simeq \R^2$ with standard basis $\{e_1,e_2\},$ scalar multiplication
is given by $(a+ib)(x\otimes e_1 +y\otimes e_2):=(ax-by)\otimes e_1 +(bx+ay)\otimes e_2.$
$v\mapsto v\otimes e_1$ is an injective real linear map $V\to \R^2 \otimes V,$ thus $V$ is a real subspace of $\C\otimes V.$
One can prove that $\C\otimes V =V\oplus iV,$ which is the reason for the notation, 
$x\otimes e_1 +y\otimes e_2 =:x+iy.$
Complexification induces a conjugation $\overline{av}=\overline{a}v =\overline{a}\otimes v.$ 
The norm of $V$ can be extended in a non-unique way, and two natural requirements on the extension are $\norm{x}_{\C\otimes V} = \norm{x}_{V},\forall x\in V,$ and $\norm{x-iy}_{\C\otimes V} = \norm{x+iy}_{\C\otimes V},\forall x,y\in V.$
Conversely if we start from a complex Banach space, we can decompose it into $V\oplus iV,$ via a projection, $\pi:\C\otimes V\to \C\otimes V,$ satisfying,
$\pi(u+v)=\pi(u)+\pi(v),\pi(au)=a\pi(u),\pi^2(u)=u.$
If $E$ is a linear space $\pi$ a projection on $E,$ such that $\{z\in E:\pi z =z\},\{z\in E\pi(z)=0\}$ are complex linear manifolds then $\{z\in E:\pi (z) =z\}\cap \{z\in E: \pi(z)=0\}=\emptyset,$ and every $z\in E,$ has a unique representation
$z=z_1+z_2,z_1\in \{z\in E:\pi (z) =z\},z_2\in \{z\in E: \pi(z)=0\}$ (see Taylor \cite{taylor}).
Complexifications render under certain conditions well-behaved polynomial extension,
see e.g.\ H\`ajek \& Johanis \cite{hajekjohanis}, p.67.
\begin{proposition}\label{hajekcomplexlem}
Let $X$ be a real normed linear space and let $Y$ be a real Banach space. Let $P_j\in
\mathcal{P}(^j X,Y),$ for $j\in \N$ and let $\tilde{P}(^j \tilde{X},\tilde{Y})$
be their complexifications.
If $R>0$ is the radius of norm convergence of the power series $\sum_{j=0}^\infty P_j$
then the radius of norm convergence of the power series $\sum_{j=0}^\infty \tilde{P}_j$
is at least $R/2.$ If $\Omega$ is the domain of convergence of
$\sum_{j=0}^\infty P_j$ then the series $\sum_{j=0}^\infty \tilde{P}_j$
converges on an open set $\tilde{\Omega}\subset\tilde{X}$ containing $\Omega.$
\end{proposition}

\subsection{Integral representation}
The Cauchy kernel for $\{\abs{z}<1\}\subset\Cn$ is given by $C(z,\zeta)=(1-\langle z,\zeta\rangle)^{-n}$
for all $(z,\zeta)\in \C^n\times\C^n$ satisfying $\langle z,\zeta\rangle\neq 1.$
The Cauchy transform of a function $f\in L^1(\{\abs{z}=1\})$ is defined by 
\begin{equation}
C[f](z):=\int_{\{\abs{\zeta}=1\}} C(z,\zeta)f(\zeta)d\sigma(\zeta)
\end{equation}
where $\sigma$ denotes the standard surface area measure.
For any $f\in C^0(\{\abs{\zeta}\leq 1\})\cap \mathscr{O}(\{\abs{\zeta}< 1\})$
we have
\begin{equation}
f(z)=C[f](z):=\int_{\{\abs{\zeta}=1\}} C(z,\zeta)f(\zeta)d\sigma(\zeta),\quad z\in \{\abs{\zeta}< 1\}
\end{equation}
Furthermore, if $u=\re f$ and $\im f(0)=0$ then
\begin{equation}
f(z)=C[f](z):=\int_{\{\abs{\zeta}=1\}} (2 C(z,\zeta)-1)u(\zeta)d\sigma(\zeta),\quad z\in \{\abs{\zeta}< 1\}
\end{equation} 
This is derived in the one dimensional case from the usual 
formula $f(z)=\frac{1}{2\pi i} \int_{\abs{\zeta}=1} \frac{f(\zeta)}{\zeta -z}d\zeta$
by a change of variables, so in fact for $z$ with $\abs{z}<1,$
\begin{equation}
f(z)=\frac{1}{2\pi i}\int_{\abs{\zeta}=1} \frac{1}{1-\zeta z}\frac{f(1/\zeta)}{\zeta}d\zeta=
-\frac{1}{2\pi i}\int_{\abs{z}=1} \frac{1}{1-\zeta z}f(1/\zeta)\frac{d\zeta}
{2\pi i \zeta}
\end{equation}
so we can use a
normalized Lebesgue measure on the circle and note that on the unit circle (i.e.\ $\abs{\zeta}=1$), we have $\bar{\zeta}=1/\zeta.$
So one could also write
\begin{equation}
f(z)=\int_{\C} \frac{f\left(\zeta\right)}{1-z}\frac{\bar{\zeta}}{\abs{\zeta}}dG(\zeta)
\end{equation}
where $G$ is a Gaussian measure on $\C$.
There exists a natural counterpart for a $G$-holomorphic function $f$ on a Banach space $X$, 
of the one-dimensional integral formula, 
namely that for each $x$ and $\lambda<1$
\begin{equation}
f(\lambda x)=\frac{1}{2\pi i}\int_{\{\abs{\zeta}=1\}} \frac{f(\zeta x)}{\zeta -\lambda}d\zeta
\end{equation} 
The usual procedure of deriving integral formulas however often involves integration of the dual $X'.$
For example, for a reflexive  complex Banach space $X$, and denoting by 
$\mathscr{P}(^m X',\C)$ the set of $m$-homogeneous complex polynomials normed by $\norm{P}=\sup_{\norm{x}\leq 1}\abs{P(x)},$
a member $P\in \mathscr{P}(^m X')$, is called an $m$-homogeneous integral polynomial of degree $m$ 
on $X'$ if there exists a Radon measure $\nu$ on $(B_1,\sigma(X,X'))$ (where $B_1$ denotes the unit ball and the topology
is given by $\sigma(X,X')$) such that
(see Dineen \cite{dineen1971}, p.273)
\begin{equation}
P(\phi)=\int_{\norm{x}<1} \langle\phi,x\rangle^m d\mu(x), \mbox{ for all }\phi\in X'
\end{equation}
The counterpart to the Cauchy type integral  
for a function $f:\{x\in X:\abs{x}<1\}\to \C,$ will with this procedure be
the condition that for every $x$ in the ball
\begin{equation}
f(x)=\int_{\{\gamma\in X':\norm{\gamma}<1\}} \frac{1}{1-\gamma(x)} d\mu(\gamma)
\end{equation}
where $\mu$ is a regular Borel measure on $(B_1,weak-star)$ (where $B_1$ denotes the unit ball in $X'$
and weak-star denotes that the topology is the weak-star topology) and we integrate over the unit ball in $X'$ which is the weak-star closure of the unit sphere.
Such functions are called {\em integral holomorphic functions}\index{Integral holomorphic functions} and
the {\em integral norm}\index{Integral norm} of such a function $P$ is defined
as $\norm{P}_I:=\inf\{ \abs{\mu}:\mu \mbox{ represents }P\},$ and the following is known, 
see Dimant et al.\ \cite{dimant}, Proposition 2,
\begin{proposition}
If $f$ is {\em integral} holomorphic (with associated Borel measure $\mu$) on $\{ x\in X:\abs{x}<1\}$ then $f$ is holomorphic in the sense that
it has a homogeneous expansion, and $m$-homogeneous component in its Taylor series
expansion of $f$ about a point $p_0\in \{ x\in X:\abs{x}<1\}$ is given by
\begin{equation}
\int_{\{\gamma\in X':\norm{\gamma}<1\}} \frac{(\gamma(x))^k}{(1-\gamma(p_0))^{k+1}} d\mu(\gamma)
\end{equation}
Furthermore, the set of {\em integral} holomorphic functions on $\{ x\in X:\abs{x}<1\}$
is a Banach space under the norm $\norm{\cdot}_I.$
\end{proposition}
Most attempts to directly generalize integral formulas
in the infinite dimensional case break down (see e.g.\ Pinasco \& Zalduendo \cite{pinasco}).
The following Cauchy formula for the infinite dimensional complex case is known (see
H\`ajek \& Johanis \cite{hajekjohanis}, p.72).
\begin{theorem}
Let $X,Y$ be Banach spaces, let $U\subset X$ be an open subset and let $f:U\to Y$ be
a holomorphic function. Let $a\in U,$ $v\in X$ and $r>0$ such that
$a+\zeta v\in U$ for all $\zeta$ such that $\abs{\zeta}<r$. Set $\gamma(t)=r\exp(it),$
$t\in [0,2\pi].$ Then
\begin{equation}
d^k f(a+\tau v)[v]=\frac{k!}{2\pi i}\int_{\gamma} \frac{f(a+\zeta v)}{(\zeta -\tau)^{k+1}}d\zeta
\end{equation}
for every $\abs{\tau}<r$ and $k\in \N.$
\end{theorem}
\begin{proof}
Let $\phi\in Y'$ and $g(\zeta):=\phi \circ f(a+\zeta v).$ Then
$g$ is holomorphic on a neighborhood of $\{\zeta:\abs{\zeta}\leq r\}$ thus by the Cauchy formula in one variable
\begin{equation}
\partial_\tau^k g(\tau)= \frac{k!}{2\pi i}\int_{\gamma} \frac{g(\zeta)}{(\zeta -\tau)^{k+1}}d\zeta
\end{equation}
Now 
the functional $\phi$ commutes with differentiation and integration,
which implies
\begin{multline}
\phi\left(d^k f(a+\tau v)[v]\right)= d^k(\phi\circ f)(a+\tau v)[v] = \partial_\tau^k g(\tau)
\\ =\phi\left(\frac{k!}{2\pi i}\int_{\gamma} \frac{f(a+\zeta v)}{(\zeta -\tau)^{k+1}}d\zeta\right)
\end{multline}
Since $Y'$ separates the points of $Y$ this renders the wanted formula.
This completes the proof.
\end{proof}
From this we see also that if $f$ is bounded by $M$ on $\{\abs{x-a}<r\}\subset U$
then for each $k\in \N$ we have $\norm{d^k f(a)}\leq Mk!r^{-k}.$

\subsection{$q$-analytic functions on Banach spaces}
\begin{remark}
In this text we shall always assume holomorphy without additional specification mean holomorphy in the sense of Definition \ref{holomujdef}, the space of holomorphic functions on an open set $U$ is denoted $\mathscr{O}(U)$. In particular, they will always be $G$-holomorphic and conversely whenever we have a continuous $G$-holomorphic function it will be holomorphic in the sense of Definition \ref{holomujdef}.
\end{remark}
We are now ready to introduce a higher order generalization of infinite dimensional holomorphy, this shall be done in terms of a generalization of
$q$-analytic functions to an infinite-dimensional setting. 
For this we 
use monomials in conjugate variables. See e.g.\ Daghighi \cite{daghighinfolo}. 
\begin{definition}
Let $X$ be a complex Banach space and let $\mathcal{P}(^m X,\C)$ denote the space of $m$-homogeneous polynomials in the
sense of Definition \ref{jack1}. 
Denote,
\begin{equation}
\mathcal{P}_{\mbox{conj}}(^m X,\C)=\{\psi:X\to \C\mbox{ such that }\psi=\overline{\phi(z)},\mbox{ some }\phi\in \mathcal{P}(^m X,\C) \}.
\end{equation}
\end{definition}

\begin{definition}[Absolute $q$-analytic functions]\label{absolutqanaldef} 
Let $X$ be a complex Banach space with unconditional (Schauder) basis.
A function $f\in C^{q-1}(X,\C)$
is called {\em polyanalytic of order $q$} or {\em  $q$-analytic} at the origin, if 
in a neighborhood, $U_0$, of the origin in $X,$ $f$ has the representation in terms of a uniformly convergent series
\begin{equation}
f(z)=\sum_{m=0}^{q-1}\sum_{\mathfrak{m}_m\in B_m} a_{\mathfrak{m}_m}(z)\mathfrak{m}_{m},a_{\mathfrak{m}_m}\in \mathscr{O}(U_0) 
\end{equation}
where $B_m$ is a linearly independent basis for $\mathcal{P}_{\mbox{conj}}(^m X,\C).$
$f$ is called {\em countably analytic} at $0$ 
if it has
the representation,
\begin{equation}
f(z)=\sum_{m=0}^{\infty}\sum_{\mathfrak{m}_m\in B_m} a_{\mathfrak{m}_m}(z)\mathfrak{m}_{m}
\end{equation}
and if the required local representation but with translation of the origin holds near every point we simply call $f$ $q$-analytic or countably analytic respectively.
\end{definition}

\begin{remark}
Let $X$ be a complex Banach space with countable unconditional basis $\{e_j\}_{j\in \N}$
(we use here countable basis synonymous with Schauder basis, in particular $X$ is a separable Banach space with a Schauder basis, 
recall that Enflo \cite{enflo} constructed a separable Banach space without a
Schauder basis). So any $x\in X$ has a unique representation $x=\sum_j c_{j,x} e_j$
and this also gives the coordinate functionals $e_j':X\to \C$, $e_j'(x):=c_{j,x}.$
For each $m\in \Z_{\geq 0},$ let $\tilde{B}_m$ be a basis for $\mathcal{P}(^m X,\C)$ (recall that we always assume
our homogeneous polynomials to be continuous). By this we mean that
for each element $p\in \mathcal{P}(^m X,\C)$ there exists a sequence $\{c_j\}_{j\in \N}$
such that we have a (norm) convergent representation 
$p(x)=\sum_{j\in \N} c_j b_j,$ where each $b_j\in \tilde{B}_m,$ and this representation is unique up to 
reordering (and since the basis is unconditional reordering does not affect the convergence).
We note the following.
If $Q\in \tilde{B}_m$ then, in particular, $Q=Q(x)$ is a function satisfying
that $(\mathbf{c}(Q))(x):=\overline{Q(x)}$ is a function satisfying that
$(\mathbf{c}(Q))(\bar{x})$ is $m$-homogeneous with respect to the variable $\bar{x}$ (where $\bar{x}$ denotes the complex conjugate
$(\bar{x}_1,\bar{x_2},\ldots)$), so we could say that $\mathbf{c}(Q)(x)$ it {\em conjugate}-$m$-homogeneous with respect to $x$. 
This is true independent of the choice of basis $B_m.$
Consider a function of the form $p_m(x)=\sum_{Q\in B_m} a_Q(x) (\mathbf{c}(Q))(x)$
where the sum is uniformly convergent and the $a_Q$ are holomorphic.
It is clear that if $p_m$ has such a representation then it will have the analogous representation with any other choice of basis
$B_m$. Furthermore, we can define $P_m(x,y):X\times X \to \C$ according to
$P_m(x,y):=\sum_{Q\in B_m} a_Q(x) (\mathbf{c}(Q))(\bar{y})$. Then $P_m$ is separately holomorphic
thus by the infinite dimensional version of Hartogs theorem it is a holomorphic function, and satisfies
$P_m(x,\bar{x})=p_m(x).$ More specifically, for each fixed $x$, $P(\cdot,y)$ is $m$-homogeneous
with respect to $y$.
We shall call a function $P(x,y)$ a {\em pseudopolynomial of degree $m$} if
it can be written as $P(x,y)=\sum_{j=0}^{m} P_j(x,y)$
for some $P_j(x,y):=\sum_{Q\in B_j} a_{Q,j}(x) (\mathbf{c}(Q))(\bar{y})$,
where the sum is uniformly convergent and the $a_{Q,j}$ are holomorphic. This can be compared to the finite 
dimensional case,
see e.g.\ Fritzche \& Grauert \cite{frigrau}, p.124 (and Definition \ref{pseudopolynomdef}). Notice that $P$ is uniquely determined by the function $p:=\sum_{j=0}^m p_j$ and vice versa.
\end{remark}
\begin{remark}\label{infoloequivrem}
Evidently, any function having the representation given in Definition \ref{absolutqanaldef}, will have such a
representation independent of the choice of bases $B_m,$ and furthermore one can equivalently
define $f$ to be {\em $q$-analytic} if there exists a pseudopolynomial, $F(x,y)$, of degree $m$ with respect to $y$, 
satisfying $F(x,\bar{x})=f(x).$
\end{remark}

\subsection{A characterization based upon restriction} 
Obviously any $q$-analytic function is a $(q+1)$-analytic function.
If $X$ is a complex Banach space and
$V\subset X$ a finite dimensional complex subvector space
then $P\in \mathcal{P}_{\mbox{conj}}(^m X,\C)$ implies that
the restriction $P|_V \in \mathcal{P}_{\mbox{conj}}(^m V,\C).$
If we for instance denote by $z:=(z_1,\ldots ,z_n)$ the complex variables in $V$
then any holomorphic $a_{\beta}$ on $X$ restricts to a function which is holomorphic in $z$
on the finite dimensional complex Euclidean manifold canonically obtained from the complex vector space $V.$
\begin{observation}\label{obs1}
By the homogeneity properties of any member $\mathfrak{m}_m\in \mathcal{P}_{\mbox{conj}}(^m X,\C)$ its restriction
to $V$ has the form of 
a sum of monomials $\bar{z}_1^{\alpha_1}\ldots \bar{z}_n^{\alpha_n},\abs{\alpha}=\sum_j \alpha_j =m,$
i.e.\ 
is a sum of elements in $\bigcup_{j=0}^{q-1}\mathcal{P}_{\mbox{conj}}(^{j} V,\C).$ 
Whence the restriction to $V$ of a $q$-analytic function (near the origin)
has the form $\sum_{\abs{\beta}\leq q-1} a_\beta (z) \bar{z}^{\beta}$ (for some holomorphic $a_{\beta}$ near the origin).
\end{observation}
\begin{observation}\label{chap}
Fixing the variables $z_j,j\neq k,1\leq j\leq n,$ any function of the form $\sum_{\abs{\beta}\leq q-1} a_\beta (z) \bar{z}^{\beta}$
reduces to a $q$-analytic function in the variable $z_k.$ Because the restriction of any $q$-analytic function on 
$X$ restricts to a $q$-analytic function on any finite $n$-dimensional $V$ and by Observation \ref{obs1} we see that
being $q$-analytic on $X$ implies being $q$-analytic (in the one dimensional sense of Definition \ref{endim}) 
along each one complex dimensional slice. 
\end{observation}
Note that we use simply $q$ for the order ans whenever the dimension is finite we assume this is a positive integer, 
i.e.\ absence of conditions specified to the separate components of the variable.
\begin{example}\label{hjlm}
Let $X=\C^2$ which we canonically identify as a Euclidean complex manifold and denote $z=(z_1,z_2)\in \C^2$
the complex coordinates.
The two functions $f_1(z)=\bar{z}_1^4\bar{z}_2^3$ and $f_2(z)=\bar{z}_1^3\bar{z}_2^4$ both belong to
$\mathcal{P}_{\mbox{conj}}(^7 X,\C)$ and one can easily see that they are
$8$-analytic.  
However, $\left(\frac{\partial}{\partial \bar{z}_2}\right)^4f_1\equiv 0\not\equiv \left(\frac{\partial}{\partial \bar{z}_2}\right)^4f_2.$ 
Thus, if one uses a positive integer $q$ to denote the order in higher dimension then one automatically is referring to the slice version, whereas the usual order will be a multi-index e.g.\ $f_1$ can also be said to be of order $\alpha=(5,4)$
whereas $f_2$ then should be of order $\alpha=(4,5).$
\end{example}
\begin{proposition}\label{sats1}
Let $X$ be a complex Banach space with unconditional Schauder basis (which can be viewed as a complex Banach manifold with open unit ball and a single chart)
let $U\subset X$ be open and let $f\in C^{q-1}(U,\C).$ 
Then $f$ is $q$-analytic on $U$ iff 
the restriction of $f$ to any one dimensional complex slice is 
$q$-analytic the sense of Definition \ref{endim}.
\end{proposition}
\begin{proof}
The 'only if'-direction follows from Observation \ref{chap}.
So assume the restriction, $f_\lambda(z,\bar{z})$, of $f$ to any one dimensional complex slice, $\lambda$, is $q$-analytic.
For the converse we use induction in $q\in \Z_+.$ When $q=1$, $f$ is holomorphic and we are done. Thus assume $q\geq 2$ and that the assertion holds true for $(q-1)$-analytic functions.
For a fixed $p\in X$ let $U\subset X$ be an open neighborhood of $p$ on which 
$f$ satisfies that $f|_{\lambda\cap U}$ is $q$-analytic for each one-dimensional
complex slice. Let $a\in U,$ $a\neq p$, $r>0$ such that $U_r:=\{\abs{x-a}<r\}\subset U$ 
and set $g(\zeta)=f(a+\zeta(x-a)),$ $\zeta\in \C,$ and $\abs{\zeta}<r/\norm{x-a}.$
Since $f\in C^{q-1},$ there exists a function 
$\overline{D} f(a)=\overline{D} g(1),$ $C^{q-2}$
with respect to $a$, which (since $\overline{D}$ is defined slice-wise on each one-dimensional complex restriction)
will be $(q-1)$-analytic on each $\lambda\cap U_r.$ By the induction hypothesis
$\overline{D} f$ has a local uniformly convergent representation of the form
$\sum_{j=0}^{q-2}\sum_{\mathfrak{m}_j\in B_j} a_{\mathfrak{m}_j}(x) \mathfrak{m}_j,$
where $B_j$ is basis for the conjugate space $\mathcal{P}_{\mbox{conj}}(^j X,\C).$
Hence (see Remark \ref{infoloequivrem}) there exists a pseudopolynomial, $F(x,w)$, of homogeneous degree $q-1$ (with respect to $w$)
such that for $x\in U_r$, $F(x,\bar{x})=\overline{D}_z f.$ 
Now each $\mathfrak{m}_j(\bar{x})$ is a sum
$\sum_{k} b_{j,k}(x)$ of holomorphic monomials of degree $j$, each of which have a primitive
$\tilde{b}_{j,k}(x)$, i.e.\ $d^1\tilde{b}_{j,k}(x)=b_{j,k}(x)$ on $U_r.$
Hence
which satisfies $(d_w \tilde{F})(x,\bar{x}) =\overline{D} f(x).$ 
Setting $h(x):=\tilde{F}(x,\bar{x})-f(x)$ we have $\overline{D}(\tilde{F}(x,\bar{x})-f(x))=0,$
i.e.\ $h$ is holomorphic, $\tilde{F}(x,\bar{x})$ is  $q$-analytic, thus
$f(x)=\tilde{F}(x,\bar{x})-h(x)$ is $q$-analytic on $U_r$. 
This completes the induction. This completes the proof.
\end{proof}
A consequence of Proposition \ref{sats1} is the following result on zero sets of $q$-analytic functions.
\begin{proposition}\label{radopropinf}
Let $X$ be a complex Banach space with countable basis and let
$f$ be a $C^{q}(U,\C)$-smooth function on an open neighborhood $U$ of $0$ in $X$
which is $q$-analytic on $U\setminus f^{-1}(0)$.
Then $f$ is $q$-analytic on $U.$
\end{proposition}
\begin{proof}
Let $\lambda$ be an arbitrary one dimensional complex slice in $U$ (by which we mean
that $\lambda$ is the intersection with $U$ of the complex span in $X$ of some vector in $X$).
If we assume $f\in C^{q}(U,\C)$ such that
$f$ is $q$-analytic on $U\setminus f^{-1}(0),$ then we know that
$f|_{\lambda}$ is a $q$-analytic function on $\lambda\cap (U\setminus f^{-1}(0))$.
If $f\in C^{q}(U,\C)$ and $f$ is $q$-analytic on $U\setminus f^{-1}(0)$ then Theorem \ref{frabsmooth} together with Proposition \ref{sats1}
implies that 
$f|_{\lambda}$ is a $q$-analytic function on $\lambda\cap U$, so invoking Proposition \ref{sats1} we obtain that 
$f$ is $q$-analytic on $U.$ This completes the proof of Proposition \ref{radopropinf}.
\end{proof}
\section{Tangential solutions}
Let $X$ be a complex Hilbert manifold and $M\subset X$ a submanifold. $T_p X$ itself can be given the structure of a complex Hilbert space (it can be identified with $X$) 
namely via the linear map $J_p:T_p X\to T_p X$ 
i.e.\ $J_p^2 v=-v,\forall v\in T_p X.$ 
Any vector subspace of $T_p X$ which is closed under the application of $J_p$ can then be identified as a complex vector space (with induced complex structure from $X$).
Let $T^{\C}_p M$ denote the largest 
vector subspace of $T_p M$ which is invariant under
the application of $J_p$ i.e.\ the largest 
vector subspace of $T_p M$ which under the induced complex structure is a complex vector subspace of $X.$ 
Recall that for a $C^1$-smooth function $f$ the decomposition into $\C$-linear and $\C$-antilinear parts, $df=\partial f +\overline{\partial} f$
implies that $f$ is holomorphic on an open $U\subset\Cn$ iff $df_p$ is $\C$-linear on $T_p \Cn,\forall p\in U $ 
(in this case $T_p \Cn$ can canonically be equipped with a complex structure).
Let $X$ be a complex Banach space and $M\subset X$ a subspace both with induced topology and differential structure. $T_p X$ itself can be given the structure of a complex Banach space (it can be identified with $X$) 
namely via the linear map $J_p:T_p X\to T_p X$
i.e.\ $J_p^2 v=-v,\forall v\in T_p X.$ 
Any vector subspace of $T_p X$ which is closed under the application of $J_p$ can then be identified as a complex vector space (with induced complex structure from $X$).
Let $H_p M$ (in some literature this is denoted $T^{\C}_p M$ or $T^c_p M$) denote the largest 
vector subspace of $T_p M$ which is invariant under
the application of $J_p$ i.e.\ the largest 
vector subspace of $T_p M$ which under the induced complex structure is a complex vector subspace of $X.$
\begin{example}
Let $X$ be a complex Hilbert space with open unit ball and unit sphere denoted by $M\subset X$. The unit sphere $M$ is then a real-analytic submanifold with (real) tangent space
$T_p M=\{z\in X\colon \re (\langle z,p\rangle ) =0 \},$ where
$\langle \cdot ,\cdot \rangle$ denotes the inner product. 
The maximal complex linear subspace of $X$ contained in $T_p M$ is $H_p M=T_p M\cap iT_p M=$
$\{z\in X\colon \langle z,p\rangle  =0 \}.$ 
\end{example}
Kaup \cite{kaup} (2004) introduced what can be interpreted as solutions to tangential Cauchy-Riemann equations in an infinite dimensional setting, in terms of uniform limits of ambient holomorphic functions. 
\begin{definition}\label{tangentialkaupdef}
Let $X$ be a complex Banach manifold and $M\subset X$ a smooth submanifold. A 
function $f\colon M \to \C$ is said to satisfy the {\em tangential Cauchy-Riemann equations on $M$} if for all $p\in X,$ the differential
$df_p\colon T_p M \to \C$ is complex linear on the subspace $H_p M\subseteq T_p M.$ A continuous function $M\to \C$ is 
to satisfy the tangential Cauchy-Riemann equations on $M$ if it is locally the uniform limit of a sequence of smooth functions that 
satisfy the tangential Cauchy-Riemann equations on $M$.
\end{definition}
It is well-known that the implicit function theorem has the following Banach space version.
\begin{theorem}[See e.g.\ Lang \cite{lang}, Theorem 5, p.18]
Let $U,V$ be open subsets of Banach spaces $X,Y$ respectively. Let $f:U\times V\to Z$
be a $C^k$-smooth mapping. Let $(a,b)\in U\times V$ and assume that the differential of $f$
with respect to the second component (i.e.\ the $Y$-component) is an isomorphism $Y\to Z.$ Let $f(a,b)=0.$ Then there exists a continuous map $g: U_0\to V$ defined on an open neighborhood $U_0$ of $a$ such that $g(a)=b$ and $f(x,g(x))=0$
for all $x\in U_0.$ If $U_0$ is taken to be a sufficiently small ball, then $g$ is uniquely determined, and is also of class $C^k$.
\end{theorem}
\begin{definition}\label{hyperdef}
Let $M\subset X$ be a submanifold of a complex Hilbert manifold $X$ such that $T^{\C}_p M\subset X$ is a complex hypersurface for each $p\in M.$ Then
$X$ is the direct product of $T^{\C}_0 M$ and $N:=\{z\in X:\langle z,w\rangle=0,\forall w\in T^{\C}_0 M\},$ where
$T^{\C}_0 M$ is a complex Hilbert space and $N$ is a complex line, i.e.\ there is an orthonormal basis (of complex vectors) $\mathcal{B}=\{e_i\}_{i\in \N}$ for $X$ 
such that Span$_{\C}e_1=N.$ We denote the coordinates
in such a basis by $(z,w)\in N\otimes T^{\C}_0 M.$ 
$M$ is called an hypersurface of class $C^k$  {\em graphed over its tangent space} near $0$ (in the real analytic case we shall write $C^{\omega}$) if for such a basis we have that near $0$ 
there is an open neighborhood $0\in U_0\subset X$ together with $\rho(z,w)\in C^{\omega}(U_0,\R)$
such that $M\cap U_0 =U_0\cap \{\rho=0\}$
and $d_z \rho(0)\neq 0$ (by which we mean that the differential is applied with respect to the first component $z$, spanned by $e_1$). In such case there exists, for sufficiently small $U_0$ (by the Banach space version of the implicit function theorem)
a $C^k$ (real analytic, $C^\omega$) map $h(\mbox{Re}z,\re w,\im w)$ 
such that
$M\cap U_0 =U_0\cap \{\mbox{Im}z=h(\mbox{Re}\zeta,\re w,\im w)\}.$
\end{definition}
We shall need an auxiliary result.
Recall that a {\em Baire space}\index{Baire space} is a topological space such that
the intersection of countably many dense open subsets is dense. By the Baire category theorem 
any Banach space (any complete metric space) is a Baire space. 
For a proof theoretical point of view we account here for a standard application of the Baire category theorem
in Hartogs' type theorems.
\begin{theorem}
Let $X$ be a Baire space, $Y$ a metric space and $F$ a normed space. If a function $f:X\times Y\to F$ is separately continuous then $f$ is bounded on a non-empty open subset of $X\times Y.$
\end{theorem}
\begin{proof}
For a given $m\in \N$ and $y\in Y$ denote $A_y(m):=\{x\in X:\norm{f(x,y)}\leq m\}.$
For $b\in Y$ and $n\in \N$ define $B_n:=\{y\in Y:\mbox{dist}(b,y)<1/n\}.$
Then the set
$A_{m,n}=\cap \{ A_y(m) :y\in B_n\}$ is closed.For each $x$ the function $y\mapsto f(x,y)$ is continuous at $b$
thus there exists $m,n$ such that $\norm{f(x,y)}\leq m$ for $y\in B_n.$ This implies that $x\in A_{m,n}.$
Thus $X=\cup_{m,n\in \N} A_{m,n}.$ Since $X$ by definition cannot be the countable union of
nowhere dense sets, there exists a nonempty open set $U\subset X$ such that
$U\subset A_{m_0,n_0}$ for some pair of $m_0,n_0\in \N.$ This implies that $f$ is bounded on $U\times B_n$. This completes the proof.
\end{proof}
Now the following generalization of Hartogs' theorem is known (see e.g.\ Bonet et al.\ \cite{bonet}, 
Thm 1.1), see also Chae \cite{chae}, Thm 14.27, p.212, for a specialized version).
\begin{theorem}\label{bonetthm}
Let $E$ be a Baire space, $G$ a metrizable space, $F$ a normed space and $U\subset E\times G$ an open subset.
Let $f:U\to F$ be a separately holomorphic mapping which is bounded on the subsets of $U$ of the form
$K\times L$, for every
finite-dimensional compact $K\subset E$ and every compact $L\subset G.$ Then $f$ is holomorphic.
\end{theorem}
\begin{proof}
Let $U=U_1\times U_2\subset E\times G$ be such that $U_1$ is 
finite-dimensional compact and $U_2$ is compact.
Since $\C\times G$ is metrizable it is known (see Barroso et al.\ \cite{barroso1977}, p.36) that for each $\xi\in U_1$ and $x\in E$, the map $g_{\xi,x}(\lambda,y):=f(\xi +x\lambda,y)$
defined on $\tilde{U}:=\{\lambda\in \C:\xi +\lambda x\in U_1\}\times U_2$ is holomorphic.
For fixed $(a,b)\in U$ let $W_1\times W_2\subset E\times G$ be an open absolutely convex subset (i.e.\
for any any pair $z,w\in W_1\times W_2$ and any pair $\lambda_1,\lambda_2$ such that $\abs{\lambda_1}+\abs{\lambda_2}\leq 1$ 
we have $\lambda_1 p_1 +\lambda_2 p_2\in W_1\times W_2$) such tat $(a,b)\in W_1\times W_2\subset U.$ If $\{V_j\}$ is a decreasing basis of 
absolutely convex neighborhoods of the origin in $G$ such that $V_1\subset W_2.$
Define
\begin{equation}
B_j:=\{ x\in W_1:\norm{\frac{1}{m!} \widehat{d^m f(a,b)}(x,y)} \leq j,\forall y\in V_j,\forall m\in \N\}
\end{equation}
For fix $y\in W_2$ define
\begin{equation}
\tilde{f}:W_1 \times\{\lambda y:\lambda\in \C,\lambda y\in W_2\}\to F
\end{equation}
\begin{equation}
 \tilde{f}(x,\lambda y):=f((a,b)+(x,\lambda y))
\end{equation}
\begin{lemma}\label{valdivia8}
Let $V$ and $W$ be topological spaces. Let $\{G_j\}_{j\in \Z_+}$ be a sequence of dense
open subsets of $V\times W.$ 
Set for each $p\in V,$ $G_j(p):=\{z\in W:(p,z)\in G_j\}.$
If $W$ is metrizable and separable there exists a subset $A$
of $V$ which is of first category such that $G_j(y)$ is dense in $W$ for $j\in \Z_+$ and each $y\in V\setminus A.$
\end{lemma}
\begin{proof}
There is a sequence $\{U_j\}_{j\in \Z_+}$ of nonempty open subsets of $W$
such that every nonempty open subset of $W$ is a union of elements in the sequence. The projection $H_j$ from
$(V\times U_j)\cap G_1\cap \cdots G_j,$ $j\in \Z_+$ onto $V$ is dense in $V$. Thus $V\setminus H_j$ is closed with nonempty interior. Hence
\begin{equation}
A=V\setminus \cap_{j=1}^\infty H_j=\cup_{j=1}^\infty (V\setminus H_j)
\end{equation}
is a subset of $V$ of first category (i.e.\ a set which can be written as a countable union of nowhere dense subsets,
where nowhere dense means that the intersection with any given open set is not dense in the given open set). Let $y\in V\setminus A=\cap_{j=1}^\infty H_j$
and $n\in \Z_+.$ If $B\subset W$ is a nonempty open subset we can find a positive integer $m>n$ such that
$U_m\subset B.$ Then $y\in H_m$ thus $B\cap G_n(y)$ is nonempty and hence $G_n(y)$ is dense.
This completes the proof.
\end{proof}

\begin{lemma}[Valdivia \cite{valdivia}, I parag. 1.2]\label{valdivialem}
Let $V$ and $W$ be Baire spaces. If $W$ is metrizable and separable then $V\times W$ is a Baire space.
\end{lemma}
\begin{proof}
Let $\{G_j\}_{j\in \Z_+}$ be a sequence of dense open subsets of $V\times W.$ 
Set for each $p\in V,$ $G_j(p):=\{z\in W:(p,z)\in G_j\}.$
Let $S$ be a nonempty open subset of $V\times W.$
Let $A$ be the set of first category of $V$ constructed in Lemma \ref{valdivia8}.
Let $P$ denote the projection of $S$ in $V.$
Since $V$ is a Baire space we have $P\cap (V\setminus A)\neq\emptyset.$ Let $y\in P\cap (V\setminus A).$
By Lemma \ref{valdivia8}, $\{G_j(y)\}_{j\in \Z_+}$ is a sequence of dense open subsets of $W.$ Since $W$
is a Baire space the intersection $G:= \cap_{j\in \Z_+} G_j(y)$ is dense in $W.$ Now $S(y)\neq\emptyset,$
$S(y)\in W,$ thus there exists $z\in G\cap S(y)$ which implies that $(y,z)\in S\cap(\cap_{j\in \Z_+} G_j).$
This proves Lemma \ref{valdivialem}.
\end{proof}
By Lemma \ref{valdivialem} $E\times \{\lambda y:\lambda\in \C\}$ is a Baire space.
This implies (see Matos \cite{matos1970}, Theorem 2.3) that $\tilde{f}$ is holomorphic, 
hence $(1/m!) \widehat{d^m \tilde{f}(0,0)}$ is continuous for each $m\in \N,$ and
for each $x\in W_1$ we have
\begin{equation}
\frac{1}{m!} \widehat{d^m \tilde{f}(0,0)}(x,y)=\frac{1}{m!} \widehat{d^m f(a,b)}(x,y)
\end{equation}
This implies that $\frac{1}{m!} \widehat{d^m f(a,b)}(\cdot,y)$ is continuous on $W_1$, thus $B_j$ is closed in $W_1.$
\end{proof}
The following is an infinite dimensional version of a special case of a result of Baouendi \& Treves \cite{bt}.
\begin{proposition}\label{sats2}
Let $X$ be a complex Hilbert space 
(in particular a Banach manifold with a single chart) and $U\subset X$ an open subset.
Let $M\subset U$ be a
real-analytic
hypersurface graphed over its tangent space the sense of Definition \ref{hyperdef}. 
Then a function $f$ belongs to the space
\begin{equation}
\left\{ g\in C^{1}(M,\C): dg_p\mbox{ is }\C\mbox{-linear on }T^{\C}_p M,\forall p\in M\right\}
\end{equation}
if and only if $f$
can be locally realized as the uniform limit of ambient holomorphic functions.
\end{proposition}
\begin{proof}
It is clear that the differential of the local uniform limit of a sequence of holomorphic functions has $\C$-linear differential
on subspaces of $T_p M$ for any $p\in M$ which with the induced complex structure become complex vector subspaces.
For the converse inclusion 
let \begin{equation*}f\in
\left\{ g\in C^{1}(M,\C): dg_p\mbox{ is }\C\mbox{-linear on }T^{\C}_p M,\forall p\in M\right\}\end{equation*}
Let $\{ e_j\}_{j\in \Z_+}$ be an unconditional basis for $X$ 
and denote by $x=(x_j)_{j\in \Z_+}$ the coordinate representation of $x\in M$ with respect to that basis.
We work locally near $x=0,$ say on $U\cap M$ where $U\subset X$ is a small open neighborhood of $0$
in $X.$
First consider the case of hypersurfaces $M\subset X,$ 
satisfying that there exists $j_0\in \Z_+$ such that $M=\{ \im x_{j_0}=0\}.$

\begin{lemma}\label{flatinfololem}
Suppose that $M\subset X,$ 
satisfies that there exists $j_0\in \Z_+$ such that $M=\{ \im x_{j_0}=0\}.$
\end{lemma}
\begin{proof}
Without loss of generality assume $j_0=1,$ and denote $x_1=z,$ $w:=(x_2, x_3,\ldots),$
so that $(z,w)$ denotes the coordinate for $X.$
Then the complex tangent space $T^{\C}_0 M$ can be identified with
$\{ (z,w)\in X: z=0\}$ and $M$ can be assumed to be a graph over $\C\times T^\C_0 M.$
A function $f:M\cap U\to \C$ then satisfies the tangential Cauchy-Riemann equations
if and only if $df$ is $\C$-linear on $\{ (z,w)\in X: z=0\}$ thus if and only if
$f$ is holomorphic with respect to the $w$ variable (independent on  its dependence on $\re z$). 
Hence $f(\re z,w)=f_{\re z}(w)=\sum_{j=0}^\infty P_j^{\re z} (w),$
where $P_j^{\re z}(w)$ is a $j$-homogeneous polynomial with respect to $w$ that is $C^1$-smooth
with respect to $\re z.$ Let $B_j$ denote an unconditional basis for
$\mathcal{P}(^j T^\C_0  X,\C)$ so that each $P_j^{\re z}(w)$
can be written as $P_j^{\re z}(w)=\sum_{\mathfrak{m}\in B_j} a_\mathfrak{m} (\re z) \mathfrak{m}$
where each $a_m(\re z)$ is $C^1$-smooth. 
Now $f$ is $C^1$-smooth, in particular locally bounded and we know
that for sufficiently small $\epsilon>0,$
the sequence $\sum_{j=0}^\infty \sum_{\mathfrak{m}\in B_j} a_\mathfrak{m} (\re z) \mathfrak{m}$
converges uniformly for $\abs{w}\leq \epsilon ,$ $\abs{\re z}<\epsilon.$ 
For each $m$ we have
for sufficiently small $\delta>0,$ and $\abs{\re z}<\delta$,
that there exists an approximating sequence of analytic functions $\tilde{a}_{\mathfrak{m},k}(z)$
converging uniformly on $\{\im z=0, \abs{\re z}<\delta \}$ to $a_\mathfrak{m}$, in particular
we can make sure that for sufficiently large $k\in \N,$ 
$\sup_{\abs{\re z}<\delta} \abs{\tilde{a}_{\mathfrak{m},k}(z) -a_\mathfrak{m}(\re z)}<\epsilon^{j+1}$
(in particular $\abs{\tilde{a}_{\mathfrak{m},k}(z)\mathfrak{m}}\leq$ 
$\abs{a_\mathfrak{m}(\re z)\mathfrak{m}} +\epsilon^{j+1}\abs{\mathfrak{m}}$).
Hence for each sufficiently large $k_0\in \Z_+$, we have for $k>k_0$ that the function 
$\tilde{f}_k(z,w):=\sum_{j=0}^\infty \sum_{\mathfrak{m}\in B_j} \tilde{a}_{\mathfrak{m},k} (z) \mathfrak{m}$ defines a separately 
analytic ($G$-holomorphic) function
function near $0$ (in $X$) such that the sequence $\{\tilde{f}_k\}_{k\geq k_0}$ converges uniformly in $M$, near $0$, to $f.$
Now for sufficiently small sets of the form $L\times K\subset \C \times T^\C_0 M,$
where $K$ is finite dimensional and $L$ compact the restriction of $\tilde{f}_k$ will be bounded.
In Theorem \ref{bonetthm} we can set $F=\C,$ $G=\C$ and $E=T^\C_0 M=\{(z,w)\in X:z=0\}$ and thus obtain that
each $\tilde{f}_k(z,w)$, $k\geq k_0,$ is (jointly) holomorphic on an open neighborhood of $0$,
such that the sequence $\{\tilde{f}_k\}_{k\geq k_0}$ converges uniformly in $M$, near $0$, to $f.$ 
This completes the proof of Lemma \ref{flatinfololem}.
\end{proof}
Now suppose that $M=\{ (\re z+ih(\re z,\re w,\im w),w)\}$ for $z=x_1,$ $w=(x_2,x_3,\ldots),$ and a real-analytic $h$
(which graphs $M$ locally near $0$ over its tangent space).
Introduce the new real variables $\nu,\xi,\eta$ and replace each occurrence of $\re z, \re w,\im w$ respectively by
$\tau_1=\re z+i\nu,$ $\tau_2:=\re w+i\xi,$ $\tau_3:=\im w+i\eta$ respectively
which gives a map $\tilde{h}(\tau_1,\tau_2,\tau_3)=h(\tau_1,\tau_2,\tau_3)$ that is 
by Proposition \ref{hajekcomplexlem} holomorphic with respect to
$\tau_1,\tau_2,\tau_3$ for sufficiently small $\abs{(\re z,\re w,\im w)},$
near $0$ and $\tau$ near the origin in $X\times X\times X.$
Also replace the expression $\re z +h(\re z,\re w,\im w)$
by $(\re z+i\nu) +i\tilde{h}.$
Then $M:=\{ \im \tau =0\}\subset X^3$. 
Define the map $\pi_0(\tau)=(\re z,\re w,\im w)$,
$\pi(\tau)=(\re z,\tau_2,\tau_3).$
Set $\tilde{M}:=\{\nu=0\}.$
Hence we have, for an open neighborhood $W$ of $0$, a map $H_1:W\to \C,$
$\tau\mapsto H_1(\tau)=\tau_1+i\tilde{h}(\tau).$ Define the map $H:W\to X$, $H=(H_1,H_2)$
where $H_2(\tau)=(\tau_2,\tau_3),$ and define the restriction map $H_0:=H|_{\{\nu=0\}}.$ 
Note that $\pi|_{\tilde{M}}= H_0^{-1}.$ 
Define $\tilde{f}(\re z,\re w+i\xi,\im w+i\eta):=f(\re z,\re w,\im w)$.
Obviously $\tilde{f}$ is $C^1$-smooth whenever $f$ is $C^1$-smooth
and since $M$ is embedded in $\tilde{M}$ we can identify $T^\C M$ as a subspace of $T^\C \tilde{M}$.
Also since $\tilde{f}$ is independent of $\xi$ and $\eta$ we obtain that
the sections of $T^\C \tilde{M}$ annihilate $\tilde{f}$ whenever
the sections of $T^\C M$ annihilate $f.$
Hence 
$\tilde{f}$ is a local $C^1$-extension of $f$ near the origin to an open neighborhood
of the origin in $\tilde{M}$ such that $d\tilde{f}$ is $\C$-linear on $T^\C \tilde{M}.$
Set $f_0:= \tilde{f}\circ H_0,$ which is a map $\{ \tau\in X^3: \im \tau_1=0\}\to \C.$
Not that $f_0\circ \pi_0^{-1}=f$ on $M$ near the origin.
By Lemma \ref{flatinfololem} there exists a there exists a sequence of ambiently holomorphic functions 
$\tilde{P}_k(\tau)$ near $0$ in $X^3$ such that $\tilde{P}_k\circ \pi\to f_0$ uniformly near $0$ in $\tilde{M}$.
This implies that the $\tilde{P}_k\circ \pi_0^{-1}$ converge uniformly near $0$ in $M$ to $f_0\circ\pi_0^{-1}=f$. 
Finally, the restriction of a holomorphic function in $X^3$ to
an open neighborhood, $U$, of the origin in $\{ (\re z+i\nu,w)\in X\}$ defines a holomorphic
function on $U$ thus the restriction of each $\tilde{P}_k|_U$ is holomorphic.
This completes the proof.
\end{proof}
The proof uses the ideas of Baouendi, Jacobowits \& Treves \cite{bjt}, p.373.
Hence for the case of $C^1$-functions on certain kinds of hypersurfaces we have that a $C^1$-smooth function satisfies the conditions of Definition
\ref{tangentialkaupdef} if and only if it can locally realized as the uniform limit of ambient holomorphic functions. In the finite dimensional case we know
by the original theorem of Baouendi \& Treves \cite{bt}, that
such an equivalence is true for any generic smooth submanifold $M\subset\Cn$.
In line with this we shall introduce the following definition.
\begin{definition}\label{polycrfunctions}
Let $X$ be a complex Banach manifold and $M\subset X$ a smooth submanifold. A 
function $f\colon M \to \C$ is 
called {\em $q$-pseudo-analytic} on $M$ if 
it can locally realized as the uniform limit of ambient $q$-analytic functions.
\end{definition} 

\section{Meta-analytic functions on Hilbert spaces}
A generalization of polyanalytic functions of order $q$ are the so called {\em meta-analytic functions}\index{Meta-analytic function}.
\begin{definition}[See e.g.\ Balk \cite{ca1}]\label{metanaldef}
Let $\Omega\subset \C$ be a domain and 
let $S(t)=s_0+s_1t+\cdots +s_{n-1}t^{q-1} +t^q$ be a polynomial with complex coefficients.
let $z=x+iy$ denote holomorphic coordinates in $\Cn$. 
A function $f\in C^{q}(\Omega,\C)$ 
is called {\em $S$-meta-analytic}\index{Meta-analytic function} if 
it satisfies on $\Omega$ the equation $S\left(\frac{\partial}{\partial\bar{z}}\right) f=0.$
\end{definition}
The following representation is known.
\begin{theorem}[See Balk \cite{ca1}, p.239, and references therein]\label{239thm}
Let $\Omega\subset \C$ be a domain.
If $S$ is a complex polynomial with roots $a_1,\ldots,a_p$ with the multiplicities $m_1,\ldots,m_p$
then a function $f$ is $S$-meta-analytic in $\Omega$ iff $f(z)=\sum_{k=1}^p P_k(z)\exp(a_k\cdot \bar{z})$
on $\Omega,$ where each $P_k$ is a polyanalytic function (with global representation) of order $m_k.$
\end{theorem}

For $q$-analytic functions there are some known differences regarding the properties of zero sets 
between the cases $q>1$ and $q=1$. 
\begin{example}\label{senare} A set $E\subset\C$ which has a condensation point\footnote{Recall that $p$ is a condensation point
if for each open neighborhood $U$ of $p$ the set $U\cap E$ is uncountable.}
is not necessarily a set of uniqueness when $q>1$ (in contrast to the case $q=1$), see e.g.\ Balk~\cite[p.\,207]{ca1}.
\end{example}
However, the identity principle remains valid when passing over to meta-analytic functions.
\begin{proposition}\label{identity}
Let $\Omega\subset\Cn$ be a domain. Let $f$ and $g$ be two $S$-meta-analytic functions on $\Omega$.
If $f=g$ on an open subset $E\subset \Omega,$ then $f\equiv g$ on $\Omega.$
\end{proposition}
\begin{proof}
The easy case of polyanalytic functions in one variable is contained in Example~\ref{senare} (more specifically Balk~\cite[p.\,207]{ca1}).
By Theorem~\ref{hartog1} this immediately yields the result for $\alpha$-analytic functions (in several variables).
Now let $F(z)$ and $G(z)$ be an $S_1$-meta-analytic functions in one complex variable $z$ (defined on the intersection of the domains of $F$ and $G$, denoted $\omega$, which w.l.o.g.\
is assumed to contain the origin) where
$S_1(t_1)$ is a polynomial in $t_1$ with roots $a_1,\ldots,a_p$ of associated multiplicities $m_1,\ldots,m_p.$
By Theorem~\ref{239thm}, $F$ and $G$ have, on $\omega$, representations
\begin{equation}\label{eq:fog}
F(\zeta)=\sum_{j=1}^{p} P_j(\zeta)\exp (a_j\bar{z}),\quad 
G(\zeta)=\sum_{j=1}^{p} Q_j(\zeta)\exp (a_j\bar{z}),
\end{equation}
where each $P_j$ and $Q_j$ is polyanalytic of order $m_j.$
Now,
\begin{equation}\label{eq:fog0}
S_1(\partial_{\bar{z}})(F-G)=0\Rightarrow H(z)=(F-G)(z)=\sum_{j=1}^{p} R_j(z)\exp (a_j\bar{z}),
\end{equation}
for some $R_j$, polyanalytic of order $m_j.$ Assuming $F=G$ on an open subset $E\subset \omega,$ we have by Eqn.\eqref{eq:fog},
\begin{equation}
\sum_{j=1}^{p} R_j(z)\exp (a_j\bar{z})=0 \Rightarrow R_p(z)=-\sum_{j=1}^{p-1} R_j(z)\exp (\hat{a}_j\bar{z}) \mbox{ on }E,
\end{equation}
where $\hat{a}_j=a_j-a_p,$ $j=1,\ldots, p-1.$
But $R_p(z)$ is polyanalytic of order $m_p$ so by the identity principle of polyanalytic functions
the right hand side must also be polyanalytic of order $m_p.$
However, by definition $a_\iota \neq a_\nu$ for $\iota \neq \nu$ 
thus the openness of $E$ implies that
$R_p\equiv 0$ on $E,$ hence $R_p\equiv 0.$ So Eqn.(\ref{eq:fog0}) reduces to
$(F-G)(z)=\sum_{j=1}^{p-1} R_j(z)\exp (a_j\bar{z}),$ and iteration of the arguments become straight-forward, yielding
$R_j\equiv 0$ for $j=1,\ldots p,$ hence $F-G\equiv 0$. This proves Proposition~\ref{identity}.
\end{proof}
A corollary to Theorem \ref{frabsmooth} is the following.
\begin{corollary} 
\label{metacor}
Let $\Omega\subset \C$ be a domain
and let $S(t)$ be a complex polynomial of the form $(a-t)^m,$ for a complex constant $a$ and positive integer $m$.
Then any function $f\in C^{m}(U,\C)$ 
which $S$-meta-analytic on $\Omega\setminus f^{-1}(0)$ is
$S$-meta-analytic
on $\Omega.$
\end{corollary}
\begin{proof}
By Theorem~\ref{239thm} we have the representation $f(z)=P(z)\exp(a\cdot \bar{z}),$ where each $P$ is polyanalytic of order $m$ on $\Omega\setminus f^{-1}(0).$ In particular, $ \left(P^{-1}(0)\right)\cap (\Omega\setminus f^{-1}(0))=\emptyset .$
Since $f$ is $C^{m}$-smooth, $\exp(-a\cdot \bar{z})f(z)$ is also $C^{m}$-smooth, in particular 
$P(z)$ has $C^{m}$  
extension, $\tilde{P}$, 
to $\Omega$ by defining it to be zero on $f^{-1}(0).$ 
The function $\tilde{P}$ satisfies the conditions with respect to $\Omega$ of
Theorem \ref{frabsmooth} 
and therefore
 defines a polyanalytic function of order $m$ on all of $\Omega.$
This in turn implies that $P(z)\cdot \exp(a\cdot \bar{z})$
extends to the $S$-meta-analytic function
$\tilde{P}(z)\cdot \exp(a\cdot \bar{z})$ on $\Omega.$
This completes the proof.
\end{proof}
We shall now define an analogue of meta-analytic functions on complex Hilbert manifolds.
Let $X$ be a complex Hilbert manifold with inner product denoted $\langle\cdot ,\cdot \rangle$ (in particular 
$\langle v , \zeta z \rangle=\overline{\langle \zeta z,v \rangle}=$
$\overline{\zeta}\langle  v,z \rangle,$ $v,z\in X,$ $\zeta\in \C$).
For each fix $v,w\in X$ we have an anti-linear
functional $\langle v,\cdot\rangle :X\to \C$
whose restriction to the complex line $\lambda =\{\zeta w:\zeta\in \C\}\subset X,$ defines a bianalytic function
$\zeta\mapsto \overline{\zeta}\langle v,w\rangle.$
\begin{definition}
Let $X$ be a complex Hilbert manifold with inner product denoted $\langle\cdot ,\cdot \rangle$.
Let $A=\{a_1,\ldots,a_n\}$ be a set of points in $X$
and let $m=(m_1,\ldots,m_n)\in \Z_+^n.$
A function $f:X\to \C$ is called {$(A,m)$-meta-analytic} at $p_0\in X$
if there exists an open neighborhood $U$ of $p_0$ in $X$
on which $f$ has the representation,
\begin{equation}
f(z)=\sum_{j=1}^n f_j(z)\exp \langle z, a_j\rangle,
\end{equation}
where $f_j(z)$ is a $j$-analytic function on $U,$
$j=1,\ldots ,n.$
\end{definition}
A corollary to the fact that$q$-analytic functions are $q$-analytic along each one-dimensional
complex slice is the following.
\begin{corollary}[To Proposition \ref{sats1}]\label{korl2}
Let $X$ be a complex Hilbert manifold. Let $A=\{a_1,\ldots,a_k\}$ be a set of points in $X$
and let $m=(m_1,\ldots,m_k)\in \Z_+^k$ and let $f$ be a function which is
$(A,m)$-meta-analytic at $p_0\in X.$ 
Then $f$ is an $S_{v}$-meta-analytic function along every one dimensional complex slice 
Span$_{\C} v,$
where $S_{v}$ is a complex polynomial in one variable, of degree $\sum_j m_j,$
such that $\langle v,a_j \rangle\in \C$ is a root of $S_{v}$ with multiplicity $m_j.$
\end{corollary}
The following is an immediate consequence of Corollary \ref{metacor} together with Corollary \ref{korl2}).
\begin{corollary}
Let $X$ be a complex Hilbert manifold and let $U\subset X$ be a bounded domain. Let $a\in X,$
$m \in \Z_+$ and let 
$f\in C^{m}(U)$ be a function which is
$(a,m)$-meta-analytic on $U\setminus f^{-1}(0)$ 
Then $f$ is automatically $(a,m)$-meta-analytic on $U.$
\end{corollary}
Proposition \ref{identity}, immediately yields the infinite dimensional version.
\begin{corollary}\label{identityinf}
$X$ be a complex Hilbert manifold and
let $f$ and $g$ be two $S$-meta-analytic functions on a domain $U\subset X$.
If $f=g$ on an open subset $E\subset U,$ then $f\equiv g$ on $U.$
\end{corollary}

\chapter{Relation to hypoanalytic theory}\label{hypoanalsec}
The study of the tangential $CR$ equations is (even in the case of hypersurfaces) not always part of the core curriculum in undergraduate studies. Therefore we have chosen to collect some relevant preliminaries on the topic in Section \ref{crapp}. Here shall attempt to apply the same thought patterns to the theory of $\alpha$-analytic functions. The archetypical $CR$ submanifold is a $C^\infty$-smooth hypersurface, and that will be our starting point.
$\alpha$-analytic functions on domains in $\Cn$ can always be identified as the restriction to a real-analytic $CR$ submanifold, $M$ of $\C^{n+1},$ of so-called
{\em pseudopolynomials} on $\C^{n+1},$ see Definition \ref{pseudopolynomdef}.
Recalling also the Weierstrass division theorem (see Theorem \ref{weierdivsion}) makes it natural to associate to $q$-analytic functions
certain holomorphic functions in additional complex variables.
\begin{example}
Let $f(z)=\sum_{j=0}^{q-1} a_j(z)\bar{z}^{j}$ be a $q$-analytic function on $\Cn.$
Set $M_1:=\{(z,w)\in \C^2\colon w=\bar{z}\}$.
Then $f$ can be identified as the restriction to $M_1$ of the (holomorphic) pseudopolynomial $F(z,w)=\sum_{j=0}^{q-1} a_j(z)w^{j}$ on $\C^2.$
Now set $M_2:=\{(z,w)\in \C^2\colon w=i\im z\}$. Rewriting $f(z)=\sum_{j=0}^{q-1} a_j(z)(z-2i\im z)^{j}$ we see that
$f$ can be identified as the restriction to $M_2$ of the pseudopolynomial $G(z,w)=\sum_{j=0}^{q-1} a_j(z)(z-2w)^{j}.$ 
\end{example}
This observation opens up the study of the analysis of not only restrictions of polyanalytic functions from ambient complex manifolds but also the study of
$q$-analytic functions in the context of $CR$ manifolds and more generally to be able to relate the study of polyanalytic functions to that of hypoanalytic functions. 
For example the function $\tilde{f}(z,t)=\sum_{j=0}^{q-1} a_j(z)(z-2it)^{j}$ is a $CR$ function on the Levi flat hypersurface
$\tilde{M}_2:=\{ (z,w)\in \C\times R \colon \re w=0\},$ for $t=\im w.$ On the other hand if $g(z)$ is a reduced $q$-analytic function then 
it takes the form $\sum_{j=0}^{q-1} a'_j(z)\abs{z}^{2j}$ for holomorphic $a_j'(z).$ Hence the function $\tilde{g}(z,t)=\sum_{j=0}^{q-1} a'_j(z)t^{j}$
on the strictly pseudoconvex hypersurface $\tilde{M}_3:=\{ (z,w)\in C^2\colon \im w=\abs{z}^2\}$ renders $g$ on the complex submanifold $M_3=\{(z,t)\in \tilde{M}_3\colon t=\im w\}.$

\section{Extension to $CR$ geometry}
We have gathered some preliminaries on $CR$ geometry in Appendix \ref{crapp}. Here we recall that any $CR$ submanifold $M\subset\Cn$ has an appropriate the of 
{\em basic solutions}\index{Basic solution}
see e.g.\ Baouendi, Ebenfeldt \& Rothschild \cite{barot}, p.36, we shall return to this fact later in the section on hypoanalytic structures also.
\begin{definition}\label{hanvisacrhypoanal}
	Let $M\subset\Cn$ be a submanifold of dimension $N$, and denote by $J$ the complex structure map ($J^2=-\mbox{Id}$) defined fiberwise on $T_p\Cn$ and set for $p\in M,$ $T_p^c M=\{X\in T_p M:J(X)\in T_p M\},$ 
	where $T_p^c M=T_p M\cap J(T_p M).$
	Then $J$ defines a complex structure on $T^c_p M.$ $M$ is called a $CR$ submanifold if the (real) dimension of $T_p^c M$ is constant and it is called generic if $T_p M+JT_p M=T_p \Cn$ for all $p\in \Cn.$ $J$ can be extended to a $\C$-linear map defined fiberwise $\C\otimes T_p\Cn\to \C\otimes T_p\Cn.$ We have a natural decomposition
	$\C\otimes T_p^c M=\mathcal{V}_p\oplus \overline{\mathcal{V}}_p,$
	where $\mathcal{V}_p=\{X\in \C\otimes T_p M:J(X)=-iX\}=\{X+iJ(X):X\in T_p^C M\},$ and
	$\re \mathcal{V}_p=T_p^c M,$ where $\re \mathcal{V}_p=\{ X+\overline{X}:X\in \mathcal{V}_p\}.$
	More generally let $M$ be a $C^k$-smooth manifold of dimension $N$. A subbundle $\mathcal{V}\subset\C\otimes TM$ of dimension $m$ is called integrable if for any $p\in M$ there exists $n=N-m$ $C^k$-smooth complex-valued $Z_1,\ldots,Z_n$ defined on an open neighborhood $U$ of $p$ in $M$ with $\C$-linearly independent differentials $dZ_1,\ldots,dZ_n$ such that $LZ_j=0$, $j=1,\ldots,n$, for any section $L$ of $\mathcal{V}$ over $M$. For $p_0\in M$, any such set of functions $Z_1,\ldots,Z_n$, vanishing at $p_0$, is called a family of {\em basic solutions} in $U.$\index{Basic solutions} $M$ is called a $CR$ structure if $\mathcal{V}\cap \overline{\mathcal{V}}=\{0\}$. If $(M,\mathcal{V})$ is integrable then the basic solutions $Z_1,\ldots,Z_n$ render an immersion of a neighborhood of $p_0$ into $\Cn$ and the image is a manifold of real dimension $N$ and its real codimension coincides with the $CR$ codimension of $M$, $N-2m.$ For an integrable $CR$ structure $Z(U)$ can be shown (see Baouendi, Ebenfeldt \& Rothschild \cite{barot}, p.37) to be an embedded real submanifold of $\Cn$, for sufficiently small $U$  and in particular $\mathcal{V}^\perp$ is spanned by $dZ_1,\ldots,dZ_n$
	and $\overline{\mathcal{V}}^\perp$ is spanned by $d\overline{Z}_1,\ldots,d\overline{Z}_n.$
\end{definition}

We note that the techniques described in e.g.\ Baouendi, Ebenfeldt \& Rothschild \cite{barot}, Proposition 1.3.6 and p.37,
render (if necessary after an ambient biholomorphic coordinate change)
near any $p\in M,$ a parametrization of $M$, given by an embedding 
$Z:\Omega\to \Cn,$ where $\Omega\subset\R^{2m}\times\R^{(n-m)}$ is an open subset containing $0$, $Z(0)=p,$ and $Z=(Z_1,\ldots,Z_n)$ is a system of basic solutions, i.e.\ have $\C$-linearly independent differentials $dZ_1,\ldots,dZ_n$,
$LZ_j=0,$ $j=1,\ldots,n,$ for all $L\in \C\otimes TM$ and 
and $d(Z_1)_p,\ldots,d(Z_n)_p,$ spans the annihilator of the $(1,0)$-subbundle of $\C\otimes T^c_p M$ (the latter has complex dimension $m\leq n$) in $\C\otimes T^*_p M,$
satisfying 	
\begin{equation}
\left\{
\begin{array}{ll}
Z_j(x,y)=x_j+iy_j & , j=1,\ldots,m\\
Z_j(x,y)=x_{j}+i\phi_{j-m}(x,y) & ,j=m+1,\ldots,n
\end{array}
\right.
\end{equation}
where $(x,y)\in \R^n\times\R^m$ denote Euclidean coordinates on $\Omega$ and $\phi_j$ denote $C^q$-smooth, real-valued functions, $\phi(0)=\mbox{Jac}\phi(0)=0$ (and where
$Z_1,\ldots,Z_N$ can be identified as the restrictions of ambient complex Euclidean coordinates $x+i(y,s)$ to $\{s=\phi(x,y)\}$, see Theorem \ref{localgraphlemma}).
\begin{remark}\label{chartcrfix}
As we shall see later, in the general theory of hypoanalytic structures, the same locally integrable subbunle $\mathcal{V}\subset \C\otimes T^c M$ can locally underlie different (for instance there are case where some choices can be smooth non-real analytic and simultaneously other choices can be real-analytic) choices of families of basis solutions $Z_1,\ldots,Z_n.$ In $CR$ geometry however, the convention is that when we speak of a generic $CR$ submanifold 
$M\subset\Cn$, then it has a fixed so-called hypoanalytic atlas (see Definition \label{hypoanaldef}) consisting of local charts $(U,Z)$ where $U\subset \R^{\mbox{dim}(m)}$ is an open subset and the $Z$ are systems of local basic solutions with respect to the $CR$ subbundle $\overline{\mathcal{V}}:=H^{0,1} M$,
where $\C\otimes T^c M=H^{1,0} M \oplus H^{0,1} M,$ 
where we require that for any $(Z',U'), (Z'',U'')$ such that $U'\cap U''\neq\emptyset,$
there is for each point $p\in U'\cap U''$ a local ambient biholomorphism $\varphi_{U',U''}$
such that near $Z'(p)=Z''(p)$ we have $Z'=\varphi_{U',U''}\circ Z''.$
We do also allow {\em biholomorphic} equivalence, in the sense
that we speak of the same $M$ even if we make a local ambient biholomorphic coordinate change (with corresponding change of the local basic solutions (local hypoanalytic charts). This is all usually implicitly understood in textbooks on $CR$ geometry but explicitly pointed out in textbooks on hypoanalytic theory, see e.g.\ Treves \cite{treves}.   	
\end{remark}

It is possible to define a $CR$ submanifold $M\subset\Cn$ using so-called defining functions\index{Local defining function for $CR$ submanifold}.
A real submanifold $M\subset\Cn$ of codimension $d$, is a subset such that for each $p_0\in M$ there is an open $U$ neighborhood of $p_0$ and a vector valued function
$\rho=(\rho_1,\ldots,\rho_d): U\to \R^d$ such that $M\cap U=\{Z\in U:\rho(z,\overline{z})=0\}$ and such that the differentials $d\rho_1,\ldots,d\rho_d$ are linearly independent in $U,$ and $z=z(x,y)$ are ambient Euclidean complex coordinates.
The map $\rho$ is called a local defining function for $M$ near $p_0.$ $M$ is generic iff near each $p\in M$ there is a local defining function $\rho$ such that the complex differentials
$\partial\rho_1,\ldots,\partial\rho_d$ are $\C$-linearly independent near $p$ (here $d=\overline{\partial}+\partial$). A generic submanifold of codimension $d$ that is generic, is necessarily a $CR$ submanifold of $CR$ dimension $m=n-d.$
\begin{remark}\label{vectorfieldrem}
	If $M\subset\Cn$ is a $CR$ submanifold of codimension $d$ and with local defining function $\rho$ near $p_0$ on $M$
	then (Baouendi, Ebenfeldt \& Rothschild \cite{barot}, p.15) any 
	any section, $L$, of the complex $CR$ bundle $\overline{\mathcal{V}}$ of $M$ (i.e.\ the $(0,1)$-subbundle of $\C\otimes T^c M$)
	can be written
	\begin{equation}
	L=\sum_{j=1}^n c_j(p) \frac{\partial}{\partial\overline{z}_j}
	\end{equation}
	where $c_j(p)$ are smooth functions for $p$ near $p_0$ such that
	\begin{equation}\label{barothdude}
	\sum_{j=1}^n c_j(p)\frac{\partial \rho_k}{\partial\overline{z}_j}(p)=0,\quad 1\leq k\leq d
	\end{equation}
	A function $f$ if a $CR$ function near $p_0$ on $M$ if and only if it is annihilated by all sections, $L$, of $\mathcal{V}$ near $p_0.$
\end{remark}

\begin{remark}
There are many equivalent ways to define differentiable $CR$ functions on a generic embedded $CR$ submanifold $M$ of $\Cn$.
Denote by $J$ the complex structure map (i.e.\ the linear map on $TM$ such that the holomorphic tangent bundle of $M$ is defined fiberwise according to
$T_p M\cap JT_p M$ $p\in M$ (where $J^2=-I$).
Let $f\in C^1(M)$. Then the following are equivalent: (i) $df$ is $\C$-linear on $H^{0,1}M$. (ii) $df$ is $J$-linear on 
$T^c M=T M\cap JT M$. (iii) For all sections $L$ of the holomorphic tangent bundle we have $Lf\equiv 0$ on $M$
(here we can either require that $L$ be any section of the real subbundle $T^c M$ or that it is any section of the complex subbundle $H^{0,1}M$).
(iv) $f$ can be locally uniformly approximated by entire functions (this is due to the Baouendi-Treves approximation theorem).  
\end{remark}
Property (iv) has a natural analogue in the $\alpha$-analytic case, and when the $CR$-dimension is $1$, so does property (iii).
In higher $CR$ dimension there appears the problem of different choices of basis for the holomorphic tangent bundle so there is no natural decomposition  
of a section $L$ of $T M\cap JT M$.
So from the perspective of property (iii) we shall use a generalized version of the notion of
$q$-analyticity in the sense of Ahern-Bruna
(see Definition \ref{ahernbrunadef}). 
For the reader's convenience we repeat it here. See also Daghighi \cite{daghighicr}
\begin{definition}[$q$-analytic functions in the sense of Ahern-Bruna in ambient space]\label{absqanaldef0}
Let $\Omega\subset\Cn$ be a domain. A function $f:\Omega\to \C$ is called {\em $q$-analytic in the sense of Ahern-Bruna}
if $\partial_{\bar{z}_1}^{\alpha_1}\cdots \partial_{\bar{z}_1}^{\alpha_n} f=0$ on $\Omega$ for all multi-integers $\alpha\in \N^n,$ such that $\sum_{1\leq j\leq n} \alpha_j =q.$
\end{definition}
\begin{proposition}\label{invariantforsta}
	Let $q\in \Z_+.$ Let $z$ and $z'$ be two biholomorphically equivalent ambient Euclidean complex coordinates for $\Cn$, and let
	$\varphi$ be the local biholomorphism near a point $p_0$ in $\Omega$ such that $z'=\varphi\circ z.$
	Let $\Omega\subset\C$ be an open subset.
	Let $f\in C^{q}(\Omega)$. If $f$ satisfies \begin{equation}\left(\frac{\partial}{\partial\overline{z}}\right)^\beta f:=\left(\frac{\partial}{\partial\overline{z}_1}\right)^{\beta_1}\cdots \left(\frac{\partial}{\partial\overline{z}_n}\right)^{\beta_n} f(z)\equiv 0\end{equation} for
	all $\beta\in \N^n,$ such that $\abs{\beta}:=\sum_j \beta_j=q$, then 
	\begin{equation}\left(\frac{\partial}{\partial\overline{z}'}\right)^\beta f:=\left(\frac{\partial}{\partial\overline{z}'_1}\right)^{\beta_1}\cdots \left(\frac{\partial}{\partial\overline{z}'_n}\right)^{\beta_n} f(z')\equiv 0\end{equation} for
	all $\beta\in \N^n,$ such that $\abs{\beta}:=\sum_j \beta_j=q$.
\end{proposition}
\begin{proof}
	We can view $z,z'$ as maps $\R^{2n}\to \Cn$, $z'(p)=\varphi\circ z(p).$ Now $dz'=\partial z' +\overline{\partial}z'$ and the chain rule implies
	$dz'_p=d\varphi_{z(p)}\circ d z=d\varphi_{z(p)}\circ (\partial z+\overline{\partial} z)|_p.$ Here
	composition means multiplication of the associated matrices
	so that by linearity
	$dz'_p=d\varphi_{z(p)}\circ dz_p =
	d\varphi_{z(p)}\circ \partial z_p + d\varphi_{z(p)}\circ \overline{\partial} z_p$,
	where the fact that $\varphi$ is holomorphic implies that $d\varphi$ is $\C$-linear thus
	also $	d\varphi_{z(p)}\circ \partial z_p$ is $\C$-linear whereas
	$d\varphi_{z(p)}\circ \overline{\partial} z_p$ is anti-$\C$-linear.
	Hence the anti-$\C$-linear part of $dz'$ can be identified with
	$d\varphi_{z(p)}\circ\overline{\partial} z|_p$ where the matrix associated to of $d\varphi_{z(p)}$ has holomorphic coefficients with respect to $z$.
	So there exists a matrix $A=[a_{ij}]_{ij}$ (where $i$ denotes row index and the entries are holomorphic functions with respect to $z$) such that
	(because each
	$\frac{\partial}{\partial\overline{z}_j}$ annihilates each of the holomorphic functions $a_{ij}$)
	\begin{multline}
	\left(\frac{\partial}{\partial\overline{z}'_1}\right)^{\beta_1}\cdots \left(\frac{\partial}{\partial\overline{z}'_n}\right)^{\beta_n} f(z')=\\
	\left(\left(\sum_{j=1}^n a_{1j} \left(\frac{\partial}{\partial\overline{z}_j}\right)\right)^{\beta_1}
	\cdots \left(\sum_{j=1}^n a_{nj} \left(\frac{\partial}{\partial\overline{z}_j}\right)\right)^{\beta_n}\right) f(\varphi^{-1}(z))=\\
	\sum_{\abs{\gamma}=q} c_{j,\gamma}
	\left(\frac{\partial}{\partial\overline{z}_1}\right)^{\gamma_1}\cdots \left(\frac{\partial}{\partial\overline{z}_n}\right)^{\gamma_n} f(\varphi^{-1}(z))
	\end{multline}
	for holomorphic functions $c_{j,\gamma}=c_{j,\gamma}(A)$. 
	This completes the proof.
	\end{proof}
One advantage of $q$-analyticity in the sense of Ahern-Bruna (except the fact that it implies $q$-analyticity along each complex line, see Ahern \& Bruna \cite{ahernbruna}, p.132) is the following.

\begin{proposition}\label{mellanstegeh}
Let $q\in \Z_+,$ 
and let $\Omega\subset\C$ be an open subset.
Then a function $f\in C^{q}(\Omega)$ is $q$-analytic in the sense of Ahern-Bruna iff
$L_1\cdots L_q f=0$ for any collection $\{L_1,\ldots,L_q\}$ of sections of $H^{0,1} \Omega,$
such that each $L_j$ is a linear combination of $\frac{\partial}{\partial \bar{z}_{1}},\ldots, \frac{\partial}{\partial \bar{z}_{n}}$ with holomorphic coefficients. 
\end{proposition}
\begin{proof}
Sufficiency is obvious since $L_1\cdots L_q f=0$ for all sections $L$ of $H^{0,1} \Omega$ implies in particular that
we can choose $L_j=\frac{\partial}{\partial \bar{z}_{j_k}},$ $j=1,\ldots,q,$
for any subset $\{j_1,\ldots,j_q\}$ of $\{1,\ldots,n\}.$ To prove necessity, we note that
each of the sections $L_k$ of $H^{0,1} \Omega$ is a linear combination
$\sum_{j=1}^n a_{kj}(z)\frac{\partial}{\partial \bar{z}_{j}}$ for holomorphic $a_{kj}(z).$ 
Hence 
	\begin{multline}
L_1\cdots L_q=
\left(\left(\sum_{j=1}^n a_{1j} \left(\frac{\partial}{\partial\overline{z}_j}\right)\right)
\cdots \left(\sum_{j=1}^n a_{qj} \left(\frac{\partial}{\partial\overline{z}_j}\right)\right)\right) =\\
\sum_{\abs{\gamma}=q} c_{j,\gamma}
\left(\frac{\partial}{\partial\overline{z}_1}\right)^{\gamma_1}\cdots \left(\frac{\partial}{\partial\overline{z}_n}\right)^{\gamma_n}
\end{multline}
for holomorphic functions $c_{j,\gamma}.$
This completes the proof.
\end{proof}

	\begin{corollary}
		$q$-analyticity in the sense of Ahern-Bruna is preserved under local biholomorphism.
	\end{corollary}
\begin{proof}
	Let $\varphi$ be a local holomorphism near a point $p\in \Cn.$
	Then $d\varphi$ is a $\C$-linear map given by a square matrix with holomorphic entries 
	acting on the complex span of $\frac{\partial}{\partial\overline{z}_1},\ldots,\frac{\partial}{\partial\overline{z}_n}.$
	In other words
	each of the vector fields $\frac{\partial}{\partial\overline{z}_j},$
	is transformed to a sum 
	linear in the $\frac{\partial}{\partial\overline{z}_k}$, $k=1,\ldots,k$ 
	with holomorphic coefficients. We can thus apply the proof of Proposition \ref{mellanstegeh}. This completes the proof.
	\end{proof}	
This notion of $q$-analytic functions in the sense of Ahern-Bruna (as opposed to $\alpha$-analyticity in higher dimension than one) generalizes naturally to certain generic $CR$-submanifolds.

\begin{definition}[$q$-analytic functions in the sense of Ahern-Bruna]\label{absqanaldef}
Let $M\subset\Cn$ be a $C^q$-smooth generic $CR$ submanifold of $CR$ dimension $m$. Let $q\in \Z_+.$
 A function $f:M\to \C$ is called {\em $q$-analytic}\index{$q$-analyticity on $CR$ submanifolds} near $p$
if there is a local basis for $L_1,\ldots,L_m$ for the set of sections of $H^{0,1}M$ 
(we will call this a local system of $CR$ vector fields)
near $p\in M$, such that we have
on an open neighborhood of $p$ 
\begin{equation}L_1^{\alpha_1}\cdots L_m^{\alpha_m} f=0,\quad \forall \alpha\in \N^m\mbox{ such that }\sum_{1\leq j\leq m} \alpha_j =q\end{equation}
\end{definition}
It is clear that when $M$ is a complex manifold then Definition \ref{absqanaldef} reduces precisely to
the definition of $q$-analytic functions in the sense of Ahern-Bruna in ambient space.
It remains to verify that this definition is 
locally invariant under local ambient biholomorphic coordinate change (see Proposition \ref{crlemmaforsta}).
Next, from the perspective of property (iv) we introduce the following.
\begin{definition}[$q$-pseudoanalytic function]\label{polycrenlichfunctions}
Let $M$ be a $C^q$-smooth generic $CR$ submanifold in $\C^n$. Let $q\in \Z_+.$
A $C^q$-smooth function $f:M\to \C,$ is called {\em  $q$-pseudoanalytic}\index{$q$-pseudoanalytic function} at $p_0\in M$
if it can be realized, near $p_0$ in $M$, as the local uniform limit of ambient $q$-analytic functions (in the sense of Ahern-Bruna), $F_k(z),$ $j\in \Z_+,$ 
i.e.\ the $F_k$ are defined in an ambient neighborhood of $p_0$ in $\C^n.$ $f$ is called $q$-pseudo-analytic on a relatively open subset $U\subseteq M$ if $f$ is $q$-pseudo-analytic at each point of $U.$
\end{definition}
We shall look at characterization for the function spaces in
Definition \ref{polycrenlichfunctions}
and Definition \ref{absqanaldef}. But first we account for some existing similar notions in the literature.

\subsection{Some related notions occurring in the literature}
Englis \& Zhang \cite{englis} and Boutet de Monvel \& Guillemin \cite{boutetguillemin}
consider the following higher order generalization of 
$CR$ functions on the smooth boundary (hypersurface) $\partial\Omega$ of a bounded domain $\Omega.$
As usual let $\overline{\partial}$ denote the anti-$\C$-linear part of $d$,  
sending $(p,q)$-forms to $(p,q+1)$-forms, in particular for $(0,0)$-forms (complex functions $f$)
$f\mapsto \overline{\partial} f=\sum_{j=1}^n \partial_{\bar{z}_j} f(z)dz_j$.
Let $\rho:\Cn\to \R$ be a defining function for $\Omega$ such that 
$\Omega:=\{ \rho>0\}$ and $\partial\Omega=\{\rho=0\},$
$\abs{\nabla \rho}>0$ on $\partial \Omega.$ The restriction, $\overline{\partial}_b$,
of $\overline{\partial}$ to $\partial\Omega$ is defined by
the restriction of $du$ to the space $T^{0,1}(\partial\Omega)$ of all vectors $X$
on $\partial\Omega$ of the form
$X=\sum_{j=1}^n X_j\partial_{\bar{z}_j},$ $X_j\in \C,$ which are tangent to $\partial\Omega,$
in the sense that $X\rho=0$ so $T^{0,1}(\partial\Omega)$ is spanned by the linearly
independent vector fields
\begin{equation}
L_{jk}:=(\partial_{\bar{z}_j}\rho)\partial_{\bar{z}_k} - (\partial_{\bar{z}_k}\rho)\partial_{\bar{z}_j},\quad 1\leq jzk\leq n 
\end{equation}
In particular, $\overline{\partial}_b u=0$ if and only if $L_{jk}u=0$ for all $j,k.$
Boutet de Monvel \& Guillemin \cite{boutetguillemin} consider the subspaces $\mathcal{B}_q$\index{$\mathcal{B}_q$}
of $L^2(\partial\Omega)$, $q\in \Z_+,$ of functions $u$ that are annihilated by $\overline{\partial}_b^q$
in the sense that 
\begin{equation}
L_{j_1 k_1}\cdots L_{j_q k_q} u=0\mbox{ for all }j_1,\ldots,j_q,k_1,\ldots,k_q
\end{equation}
These spaces are not invariant under biholomorphisms. For example setting 
$\Omega:=\{z\in \C^2:\abs{z}<1\}$, the function $\bar{z}_2$ belongs to $\mathcal{B}_2=\mbox{Ker}L^2_{12}$.
If $\phi_a$ denotes the automorphism of $\Omega$ interchanging $0$ with $a\in \Omega$
then $\bar{z}_2\circ\phi_a\notin \mathcal{B}_2$ if $a\neq 0.$
Englis \& Zhang \cite{englis} introduce the
so-called so-called {\em higher Cauchy-Riemann spaces}\index{Higher Cauchy-Riemann spaces} $\mathcal{C}^q$ for $q\in \Z_+$ which will
be invariant under biholomorphisms.
Any bounded domain has (K\"ahler \cite{kahler}) a natural so-called K\"ahler metric, namely the Bergman metric, which is invariant
under biholomorphic mappings. Let us recall the definition of a K\"ahler metric.
A real $1$-form, $\omega=\sum_{j=1}^{2n} \omega_j dx_j,$ in coordinates $x_1,\ldots,x_{2n}$, 
on a complex $n$-dimensional manifold, $M$, can, with the notation $z_j:=x_{j}+ix_{n+j},$ be expressed as
\begin{equation}
\omega=\sum_{j=1}^n \frac{\omega_{j}-i\omega_{n+j}}{2}dz_j+
\sum_{j=1}^n \frac{\omega_{j}+i\omega_{n+j}}{2}d\bar{z}_j
=:\omega^{(1,0)}+\omega^{(0,1)}
\end{equation}
Recall that a map $H\colon V\times V\to \C$
on a complex vector space $V$ is called a Hermitian form if it is linear in the first coordinate and such that
$H(u,v)=\overline{H(v,u)}.$ Every Hermitian form
has an associated Hermitian matrix, $A,$ such that $H(u,v)=u A \overline{v}^T$.
A Hermitian metric $g$ on $M$
is a tensor field 
of the form
\begin{equation}
\sum_{i,j=1}^n g_{i\bar{j}}(z)dz_i\otimes d\bar{z}_j
\end{equation}
such that $g^{i\bar{j}}=\overline{g^{j\bar{i}}}$ (i.e.\ the matrix defined by $g^{i\bar{j}}(z)$ is Hermitian), the matrix defined by $g^{i\bar{j}}(z)$ is positive definite for each $z$ ($u[g^{i\bar{j}}(z)]_{ij}u^T>0$ for all $u\in \Cn$), and such that it depends smoothly on $z$.
	To any Hermitian metric $g$ we associate 
	a $(1,1)$-form (called the {\em K\"ahler form})
	\begin{equation}
	\omega_g:=\frac{i}{2}\sum_{i,j=1}^n g^{i\bar{j}}dz_i\wedge d\bar{z}_j
	\end{equation}
	This is a real form in the sense that $\overline{\omega}_g=\omega_g$.
	If $d\omega_g =0$ (recall that we are assuming $M$ is a complex manifold) then $\omega_g$ is called a {\em K\"ahler metric}\index{K\"ahler metric}.
	Now suppose that for each $z$ there exists an open $U\ni z$ together with  a function $F:U\to\R$ such that $\omega_g=\partial\overline{\partial}F.$
	Then
	\begin{equation}
	d(\partial\overline{\partial}F)=(\partial+\overline{\partial})\partial\overline{\partial}F-\partial\overline{\partial}\overline{\partial}F=0
	\end{equation}
	which implies $d\omega_g=0.$ Conversely suppose $d\omega_g=0$.
	For each $z$ there exists (see e.g.\ Proposition \ref{primitihorm}) an open $U\ni z$ together with a $1$-form $\eta$
	such that $d\eta=\omega_g.$ Since $\omega_g$ is a $(1,1)$-form we have a decomposition 
	the $1$-form $\eta$,
	$\eta=\eta^{(1,0)}+\eta^{(0,1)}$ such that $\overline{\partial}\eta^{(0,1)}=0=\partial\eta^{(1,0)}.$
	Thus for sufficiently small $U\ni z$ there exists functions $\alpha,\beta$ such that
	$\eta^{(0,1)}=\overline{\partial}\alpha,$ $\eta^{(1,0)}=-\overline{\partial}\beta.$
	Since $\partial\overline{\partial}=-\overline{\partial}\partial$ this yields
	$\omega_g=\partial\overline{\partial}(\alpha+\beta).$ Since $\omega_g$ is real ($\overline{\omega}_g=\omega_g$)
	$F:=\alpha +\beta$ is real.
	 Hence $d\omega_g=0$ is equivalent to the condition that that for each $z$ there exists an open $U\ni z$ together with  a function $F:U\to\R$ such that $\omega_g=\partial\overline{\partial}F.$ 
	Furthermore, the condition $d\omega_g =0$ is in local coordinates precisely 
\begin{equation}
\partial_{z_i} g^{j\bar{k}}(z)=\partial_{z_j} g^{i\bar{k}}(z) 
\end{equation}
or equivalently
\begin{equation}
\partial_{\bar{z}_i} g^{j\bar{k}}(z)=\partial_{\bar{z}_j} g^{j\bar{i}}(z)
\end{equation}
Any complex manifold admits a Hermitian metric but the corresponding form $\omega_g$ is not always a 
K\"ahler metric. 
A complex manifold equipped with a K\"ahler metric is called a K\"ahler manifold.
We mention that, alternatively, one may define K\"ahler metrics via complex structure maps 
(see Section \ref{crapp} in the 
appendix for more on almost complex and complex structure maps).
Let $M$ be a real $C^\infty$-smooth $2n$-dimensional manifold. 
A Riemannian metric $g$ on $M$ is a smooth section of 
$T^*M\otimes T^* M$ defining a positive
symmetric bilinear form on $T_x M$ for each $x\in U.$ 
In local coordinates $x_1,\ldots x_{2n},$
and the natural basis $\partial_{x_1},\ldots , \partial_{x_{2n}}$ for $TM$, $g$ can be represented
by a smooth matrix-valued function $[h^{ij}]_{ij}$, $h^{ij}=g(\partial_{x_i},\partial_{x_j})$ 
where $[h^{ij}]_{ij}$ is positive definite. The doublet $(M,g)$ is called a Riemannian manifold.
An almost complex structure $J$ on $M$ is a bundle automorphism of the tangent bundle $TM$ such that $J^2=-\mbox{Id}.$
$J$ is called a complex structure if $0=[u,v]+J[Ju,v]+J[u,Jv]-[Ju,Jv],$ for all $C^\infty$-smooth
sections
$u,v$ of $TM.$ In such case $J$ 
corresponds to the induced complex multiplication in $\C\otimes TM$ and $T^{0,1}M$ is integrable.
An almost complex structure $J$ on $M$ is called {\em compatible} with the metric $g$ on $M$
if $g(u,v)=g(Ju,Jv).$
When $J$ is a complex structure ($M$ is a complex manifold), the {\em K\"ahler metric}, $\omega_g$, associated to $g$ is the $2$-form
\begin{equation}
\omega_g(u,v):=-g(u,Jv)
\end{equation}
For example if $M=\R^{2n}\simeq \Cn$ with identification via the complex 
structure map
$J$ satisfying $J(\partial_{x_j})=\partial_{x_{n+j}},$ 
$J(\partial_{x_{n+j}})=-\partial_{x_{j}},$ (where clearly $J^2=-\mbox{Id}$), 
and $z_j=x_j+ix_{n+j}$ the we also have, with 
$\partial_{z_j}=\frac{1}{2}(\partial_{x_{j}}-i\partial_{x_{n+j}}),$
$\partial_{\bar{z}_j}=\frac{1}{2}(\partial_{x_{j}}+i\partial_{x_{n+j}}),$ that
$J(\partial_{z_{j}})=i\partial_{z_{j}},$ $J(\partial_{\bar{z}_{j}})=
-i\partial_{\bar{z}_{j}}.$
If $g$ is the (Euclidean) metric $\sum_{j=1}^n dz_j\otimes d\bar{z}_j$ then the associated K\"ahler form becomes
\begin{equation}
\omega_g=\frac{i}{2}\sum_{j=1}^n dz_j\wedge d\bar{z}_j =\sum_{j=1}^n dx_j\wedge d x_{n+j}
\end{equation}
Consider a K\"ahler form $g^{j\bar{k}}d\bar{z}_l\wedge dz_k$, on $\Omega$, we shall denote
the associated K\"ahler metric by $g^{j\bar{k}}.$
Here we use the Einstein summation convention that any variable
appearing in both upper and lower indices will be summed automatically.
Denote by $g^{\bar{l}j}$ the inverse matrix to
$g_{j\bar{k}}$.
Peetre \& Englis \cite{peetreenglis} 
defined the operator (although it seems to have appeared earlier in Peetre, Peng \& Zhang \cite{peetrepengzhang}
and Peetre \& Zhang \cite{peetrezhang1993}).
\begin{equation}
\overline{D} f:=g^{\bar{l}k}\partial_{\bar{z}_l} f
\end{equation}
First let us look at the earliest forms first. 
Denote by $SU(1,n)$ the set of block matrices
\begin{equation}
g=\begin{bmatrix}
A & B\\
C & D
\end{bmatrix}
\end{equation}
where $D\in \C, A,B,C$ are $n\times n$, $n\times 1,$ $1\times n$ matrices respectively with complex entries such that $g$ is unitary with respect to
the indefinite metric $\abs{z_1}^2+\abs{z_2}^2+\cdots +\abs{z_n}^2-\abs{z_{n+1}}^2$ on $\C^{n+1}$ and $\mbox{det}g=1.$
The action of $SU(1,n)$ on the unit ball is given by
$gz=(Az+B)(Cz+D)^{-1}$ and its action on $L^2(\{\abs{z}=1\})$ by
$T_g f(z)=f(gz)(\mbox{det} g'(z))^{\frac{\nu}{n+1}}.$
In Peetre, Peng \& Zhang \cite{peetrepengzhang} is defined for elements $\psi\in SU(1,n)$ 
such that
$\psi 0=z,$ the operators
\begin{equation}
\overline{D}^j f(z)=\left(\sum_{k=0}^n \partial_{w_k} \psi_j(w) \partial_{\bar{w}_k} (f\circ\psi)(w)\right)|_{w=0}
\end{equation}
where
they show, for $R:=\sum_{i=0}^n z_i\partial_{z_i},$ that
\begin{equation}
\overline{D}^j f(z)=(1-\bar{z}^2)\left( \partial_{\bar{z}_j}f -z_j\bar{R} f\right),\quad j=1,\ldots,n
\end{equation}
Now consider the operator from Peetre \& Engli\u{s} \cite{peetreenglis}.
More precisely let $\Omega$ be a complex $n$-dimensional manifold with a K\"ahler metric
with K\"ahler form $g_{j\bar{k}}d\bar{z}_l\wedge dz_k$ on $\Omega$ and let $E$ be a
Hermitian vector bundle over $\Omega.$ One defines a convariant operator $\overline{D}=\overline{D}_E$ (thus dependent upon $E$) locally according to
\begin{equation}
\overline{D}_E (f^\alpha e_\alpha )=g^{\bar{j}k} (\partial_{\bar{z}_j} f^\alpha) e_\alpha \otimes \partial_{z_k} 
\end{equation}
where $\{ e_\alpha\}_\alpha$ is a system of local trivializing sections for $E$ and the
$\partial_{z_k}$, form the standard basis of the holomorphic tangent space at $z$, and $g^{\bar{j}k}$ is the matrix inverse of $g_{k\bar{j}}.$
In particular, $\overline{D}_E$ sends sections of $E$ to sections of the tensor product of $E$
with a symmetric tensor product of the holomorphic tangent bundle $T^{1,0}\Omega$, namely denoting by
$\odot^q T^{1,0}\Omega$ the symmetric tensor subbundle of $\otimes^q T^{1,0}\Omega$
we have that $\overline{D}^q$ maps $C^\infty(E)$ into $C^\infty(E\otimes(\odot^q T^{1,0}\omega)),$ see 
Peetre \& Zhang \cite{peetrezhang1998}.
Invariance is meant by the fact that if $\psi$ is a biholomorphic map $\Omega$ into $\Omega$ then
$\overline{D}(f\circ \psi)=\overline{D}f\circ \psi.$
Define for $q\in \Z_+,$
\begin{equation}
\overline{D}^q f:=g^{\bar{l}_q k_q}\partial_{\bar{z}_{l_q}} g^{\bar{l}_{q-1} k_{q-1}}\partial_{\bar{z}_{l_q-1}}\cdots g^{\bar{l}_{1} k_1}\partial_{\bar{z}_{l_1}}
\end{equation}
The above cited result of Peetre \& Zhang \cite{peetrezhang1998} shows that
$(\overline{D}^q f)^{k_q,\ldots,k_1}$ is symmetric in the indices $k_q,\ldots,k_1.$ 
\begin{definition}
The $q$:th Cauchy-Riemann space $\mathcal{C}_q$ is defined as
\begin{equation}
\mathcal{C}^q:=\{f: \overline{D}^q f=0\mbox{ on }\Omega\}
\end{equation}
\end{definition}
Clearly, $\mathcal{C}^1$ coincides with the holomorphic functions and is independent of the metric $g_{j\bar{k}},$
whereas for $q>1$ $\mathcal{C}^q$ will depend upon the metric. Choosing $\Omega$ to be a bounded strictly pseudoconvex domain with smooth boundary
and the metric $g_{j\bar{k}}$ to be e.g.\ the Bergman metric (or the Poincare metric or the Szeg\"o metric), $\mathcal{C}^q$ will be invariant under biholomorphisms, see
Englis \& Zhang \cite{englis}.
We mention that there also exists a notion, see Peetre \& Zhang \cite{peetrezhang1993} (see also Englis \& Peetre \cite{peetreenglis}) of
the {\em higher Laplace-Beltrami operators}.

\section{Characterization of $q$-analytic functions in $CR$ geometry}
\begin{proposition}\label{crpolyprop1}
Let $M$ be a $C^q$-smooth generic $CR$ submanifold in $\C^n$. Let $q\in \Z_+.$ Then a $C^q$-smooth function $f$ is $q$-pseudoanalytic at $p_0\in M$
if and only if there exists differentiable $CR$ functions $a_\beta,$ $\beta\in \N^n,$ $\sum_j \beta_j< q,$ such that $f$ has, near $p_0$, the representation
\begin{equation}\label{crforstaekv}
f(z)=\sum_{\abs{\beta}<q} a_\beta(z) \bar{z}^\beta
\end{equation}
$z\in \Cn \cap M$ near $p_0$ (here $z=(z_1,\ldots,z_n)$ and the expressions $f(z),a_\beta(z)$ are thus only defined when $z$ lies near $p_0$ on $M$ in $\Cn$).
\end{proposition}
\begin{proof}
Assume $f$ has the representation given by Eqn.(\ref{crforstaekv}). By the Baouendi \& Treves approximation theorem
each $a_\beta$ can be uniformly approximated by entire functions near $p_0$, thus $f$ is $q$-pseudo-analytic at $p_0\in M$.
Conversely, suppose $f$ is $q$-pseudoanalytic at $p_0\in M$ and let $E_j(z)=\sum_{\abs{\beta}<q} E_{j,\beta}(z)\bar{z}^\beta$
be a sequence of ambient $q$-analytic functions in the sense of Ahern-Bruna, in $\Cn$, uniformly converging to $f$ near $p_0$ in $M.$
We claim that there exists a subsequence $\{E_{j_k}\}_{k\in \Z_+}$ such that for each $\beta$ this implies that $\{E_{j_k,\beta}(z)\}_{k\in \Z_+}$ converges uniformly near $p_0$ in $M$, say on an open neighborhood $U$ of $p_0$ in $M$. Since the local uniform limit on $M$
of a sequence of entire functions is a $CR$ function this suffices to prove wanted the result.  W.l.o.g.\ assume $p_0=0.$
We proceed by induction. For $q=1$, each $E_{j}$ is holomorphic so the conclusion follows from the fact that local uniform limits of holomorphic functions on generic $CR$ submanifolds are $CR$ functions.
Assume $q\geq 2$ and that the conclusion holds true for $(q-1)$-pseudoanalytic functions.
Recall that for complex functions $c(z),d(z)$ and natural numbers $N_1,N_2\in \N$, we have
\begin{equation}
\abs{z^{N_1}c(z)+z^{N_2}d(z)}=\abs{z}^{N_1}\abs{c(z)}^2+\abs{z}^{N_2}\abs{d(z)}^2+2\abs{z}^{N_1+N_2}\re(c(z)\overline{d(z)})
\end{equation}
For each assignment 
\begin{equation}
\gamma: \{\alpha\in \N^n:\abs{\alpha}=q-1\}\to \{\alpha\in \N^n:\abs{\alpha}<q-1\}
\end{equation}
pick integers $N_1,N_2\in \N$ such that for $\abs{z}<r,$ where $r>0$ is sufficiently small, we have
\begin{multline}
\abs{z^{(N_1,\ldots,N_1)}\left(\sum_{\abs{\beta}<q-1} E_{j,\beta}(z)\bar{z}^\beta \right) + z^{(N_2,\ldots,N_2)}\left(\sum_{\abs{\beta}=q-1}E_{j,\beta}(z)\bar{z}^{\gamma(\beta)}\right)}\leq \\
C\abs{\sum_{\abs{\beta}<q} E_{j,\beta}(z)\bar{z}^\beta}
\end{multline}
for a constant $C>0.$
Now the sequence $\{F_j\}_{j\in \Z_+},$ with
\begin{equation}F_j(z)=z^{(N_1,\ldots,N_1)}\left(\sum_{\abs{\beta}<q-1} E_{j,\beta}(z)\bar{z}^\beta \right) + z^{(N_2,\ldots,N_2)}\left(\sum_{\abs{\beta}=q-1}E_{j,\beta}(z)\bar{z}^{\gamma(\beta)}\right)
\end{equation}
has only $(q-1)$-analytic (in the sense of Ahern-Bruna) members. 
By the induction hypothesis, there is a subsequence $\{F_{j_k}\}_{k\in \Z_+}$ such that
each sequence of analytic components corresponding to the same order, i.e.\ with members $z^{(N_1,\ldots,N_1)}E_{j_k,\beta}(z)$ or $z^{(N_1,\ldots,N_1)}E_{j_k,\beta}(z)+z^{(N_2,\ldots,N_2)}E_{j_k,\gamma(\beta)}(z)$ respectively, $\abs{\beta}<q$, converges
uniformly near $p_0$ in $M$. By varying over all possible $\gamma$ we obtain,
Since the $E_{j_k,\beta}(z)$ are holomorphic (in particular they have no poles at $0$ in $\Cn$), that each of the sequences 
$\{E_{j_k,\beta}(z)\}_{k\in \Z_+}$, $\abs{\beta}<q$, converges
uniformly near $p_0$ in $M$. Since each $E_{j_k,\beta}(z)$ is an analytic function the local uniform limit of such functions, on $M$, near $p_0$, is a $CR$ function. 
This completes the induction. This completes the proof.
\end{proof}
\begin{proposition}\label{crlemmaforsta}
Let $q\in \Z_+.$ Let $M\subset\Cn$ be a $C^q$-smooth generic $CR$ submanifold of $CR$-dimension $m$. Let
If $f\in C^{q}(M)$ and let
$p_0\in M$.
\\
(i) The property of $f$ being $q$-analytic (in the sense Ahern-Bruna) near $p_0$ in $M$ in invariant with respect to local biholomorphism.\\. 
(ii) Let $L=(L_1,\ldots,L_m)$ be a local system of $CR$ vector fields.
If $f$ satisfies \begin{equation}
L^\beta f:=L_1^{\beta_1}\cdots L_1^{\beta_m}
f\equiv 0\end{equation} for
all $\beta\in \N^n,$ such that $\abs{\beta}:=\sum_j \beta_j=q$, then 
for any $m\times m$ matrix $A$ whose coefficients are differentiable $CR$ functions we have
\begin{equation}
(L')^\beta f:=(L'_1)^{\beta_1}\cdots (L'_m)^{\beta_m}
f\equiv 0\end{equation}
where $L'=AL.$
\end{proposition}
\begin{proof}
	The claim is local and
	w.l.o.g.\ we assume $p_0=0\in M.$
	Let $\varphi$ be a local biholomorphism near $0.$
	First note that a holomorphic map sends $CR$ manifolds to $CR$ manifolds and holomorphic tangent spaces onto holomorphic tangent spaces,
	and it will map basis of $CR$ vector fields onto a basis of $CR$ vector fields. We denote by
	$\varphi_*X$ the push-forward of a complex tangent vector field $X$ on $M$ under $d\varphi,$
	where we use this notation for the restriction of $\varphi$ to $M.$
	Under the coordinate change denote the transformed function by $f(z)=(f\circ\varphi)(\varphi^{-1}\circ z)=:\tilde{f}(\varphi^{-1}(z)).$ 
	Hence the equations 
	\begin{equation}
	(\varphi^{-1}_*L_1)^{\beta_1}
	\cdots (\varphi^{-1}_* L_m)^{\beta_m} \tilde{f}(\varphi^{-1}(z))=0,\quad z\mbox{ near $0$ in }M
	\end{equation}
	for all $\beta\in \N^m$ such that $\abs{\beta}=q$, show, because $\{ \varphi^{-1}_*L_1,\ldots, \varphi^{-1}_*L_m\}$ form a local system of $CR$ vector fields in the new coordinates, that $\tilde{f}$ is $q$-analytic in the sense of Ahern-Bruna.
	This proves $(i).$\\
The claim $(ii)$ follows from the fact that if $L'=AL,$ for a matrix $A=[a_{ij}]_{ij}$, where $i$ denotes row index and the entries are differentiable $CR$ functions,
 then (because each
$L_j$ annihilates $CR$ functions)
\begin{multline}
(L_1')^{\beta_1}\cdots (L_m')^{\beta_m} f=\left(\sum_{j=1}^m a_{1j} L_1^{\beta_1}\right)
\cdots \left(\sum_{j=1}^m a_{mj} L_m^{\beta_m}\right)f=\\
\sum_{\abs{\gamma}=q} c_{\gamma}
L_1^{\gamma_1}\cdots L_m^{\gamma_m}f
\end{multline}
for differentiable $CR$ functions $c_{\gamma}$. This proves $(ii).$ This completes the proof.
\end{proof}
\begin{remark}
	Note that being $q$-analytic in the sense of Ahern-Bruna will thus be independent of the choice of overlying chart $Z_1,\ldots,Z_m,$
	in the sense that if $M\subset\Cn$ is a generic $CR$ submanifold that is is parametrized by some open subset $\Omega\subset\R^{\mbox{dim}(M)}$, $p\in \Omega$
	and if there is an open neighborhood $U$ of $p$ in $\Omega$
	such that two local systems of basic solutions $Z,Z'$ represent $M$ then in part $Z(U)=Z'(U)$
	and there exists, by definition an open neighborhood $V$ of $Z(p)=Z'(p)$ in $\Cn$ and a local biholomorphism
	$\varphi$ on $V$ such that $Z=\varphi Z'.$
\end{remark}
As the choice of nomenclature implies, we have the following.	
\begin{proposition}[$q$-analytic version of the Baouendi \& Treves approximation theorem]\label{crpolyprop2}
Let $M\subset\Cn$ be a $C^q$-smooth generic $CR$ submanifold of $CR$-dimension $m$, let $q\in \Z_+$ and  
let $f\in C^{q}(M)$.
If $f$ is $q$-analytic (in the sense of Ahern-Bruna) on $M$ then $f$ is $q$-pseudoanalytic on $M$.
\end{proposition}
\begin{proof}
The conclusion is local so it suffices to prove the result for a fixed point $p_0\in M.$ 
Also by Proposition \ref{crlemmaforsta} it suffices to prove the result in any one
local system of basic solutions near $p_0.$
W.l.o.g.\ we assume $p_0=0.$
Let $Z:\Omega\to Z(\Omega) \subset \Cn,$ be a local parametrization of $M$ near $0$, where $\Omega\subset\R^{2m}\times\R^{(n-m)}$ is an open subset containing $0$, $Z(0)=0,$ and $Z=(Z_1,\ldots,Z_n)$ a system of basic solutions, 
satisfying
\begin{equation}
\left\{
\begin{array}{ll}
Z_j(x,y)=x_j+iy_j & , j=1,\ldots,m\\
Z_j(x,y)=x_{j}+i\phi_{j-m}(x,y) & ,j=m+1,\ldots,n
\end{array}
\right.
\end{equation}
where $(x,y)\in \R^n\times\R^m$ denote Euclidean coordinates on $\Omega$ and $\phi_j$ denote $C^q$-smooth, real-valued functions, $\phi(0)=\mbox{Jac}\phi(0)=0$, and for some small neighborhood $W$ of $0$ in $\Cn,$
$W\cap M=Z(\Omega).$ We assume this is such that $Z_1,\ldots,Z_n$ can be identified as the restrictions of ambient complex Euclidean coordinates $x+i(y,s)$, by introducing the real Euclidean coordinate $s$ for an open $U\subset \R$, $0\in U,$ to $\{s=\phi(x,y)\}$, see Theorem \ref{localgraphlemma}.
Set $d=n-m.$ If $m=0$ then $M$ is totally real and any differentiable function on $M$ is a $CR$ function.
If $m\geq 1.$
we define for each real $c$ sufficiently near $0$, a local submanifolds $M_c\subset\Cn$ foliating a neighborhood of $0$ in $M$, 
by setting $M_c:=\tilde{Z}_{c}(\Omega),$ for $\tilde{Z}_{c}:\Omega\to \Cn$ given by
\begin{equation}
\left\{
\begin{array}{ll}
\tilde{Z}_{c,1}(x,y)=x_{1}+ic &\\
\tilde{Z}_{c,j}(x,y)=x_j+iy_j & , j=2,\ldots,m\\
\tilde{Z}_{c,j}(x,y)=x_{j}+i\phi_{j-m}(x,y) & ,j=m+1,\ldots,n
\end{array}
\right.
\end{equation}
We see that on a sufficiently small neighborhood $V$ of $0$ in $\Cn$ we have in $\Cn$, $M_c\cap V=\{y_1=c\}\subset\tilde{M}$ and, Letting $(x+i(y,s))$ denote holomorphic coordinates, centered at $0$ in $\Cn$
where $s$ are Euclidean coordinates for $\R^{d},$ $d=n-m$, 
each $M_c\subset\Cn$ is a generic $CR$ submanifold
of codimension $(d-1)$ (this is because it has, near $0$, the local defining function $\rho^{M_c}=(\rho_1^{M_c},\ldots,\rho_{d+1}^{M_c}),$ $\rho_1^{M_c}:=y_1-c,$
$\rho_{j}^{M_c}:=s_j-\phi_{j}$, $j=1,\ldots,d,$ and by construction by $\mbox{Jac}\phi(0)=0$, so that $\rho$ has rank $d+1$ near $0$), and it obviously has $CR$ dimension $m-1$. The local defining function for $M$ is accordingly $\rho:=s-\phi.$
By Remark \ref{vectorfieldrem}
we know that, given a local defining function, $\rho$ for $M$ near $0$, any holomorphic vector field takes the form
\begin{equation}
L=\sum_{j=1}^n c_j(p) \frac{\partial}{\partial\overline{z}_j}
\end{equation}
where $c_j(p)$ are smooth functions for $p$ near $p_0$ such that
\begin{equation}
\sum_{j=1}^n c_j(p)\frac{\partial \rho_k}{\partial\overline{z}_j}(p)=0,\quad 1\leq k\leq d
\end{equation}
Let $L_1,\ldots,L_m$ be a local basis for the sections of $H^{0,1} M$ near $0$ and for $p$ near the origin
\begin{equation}
L_j=\sum_{k=1}^n c_{k,j}(p) \frac{\partial}{\partial\overline{z}_k}
\end{equation}
for smooth functions $c_{k,j}$ satisfying
\begin{equation}
\sum_{k=1}^n c_{k,j}(p)\frac{\partial \rho_l}{\partial\overline{z}_k}(p)=0,\quad 1\leq l\leq d
\end{equation}
in other words 
\begin{equation}
\begin{bmatrix}
\frac{\partial \rho_1(x,y)}{\partial\overline{z}_1} & \cdots & \frac{\partial \rho_1(x,y)}{\partial\overline{z}_n}\\
\vdots & \ddots & \vdots\\
\frac{\partial \rho_d(x,y)}{\partial\overline{z}_1} & \cdots & \frac{\partial \rho_d(x,y)}{\partial\overline{z}_n}\\
\end{bmatrix}
\begin{bmatrix}
c_{1,1} & \cdots &c_{1,m}\\
\vdots & \ddots & \vdots\\
c_{n,1} & \cdots &c_{n,m}\\
\end{bmatrix}=0
\end{equation}
and
\begin{equation}
\begin{bmatrix}
c_{1,1} & \cdots &c_{1,m}\\
\vdots & \ddots & \vdots\\
c_{n,1} & \cdots &c_{n,m}\\
\end{bmatrix}
\end{equation}
has rank $m$. 
The associated local defining function for $M_c$ renders a local basis  $\tilde{L}_1,\ldots,\tilde{L}_{m-1}$ for the sections of $H^{0,1}M_c$ near $(x,y_1,y')=(0,c,0),$ $y'=(y_2,\ldots,y_m),$ according to 
\begin{equation}
\tilde{L}_j=\sum_{k=1}^n d_{k,j}(p) \frac{\partial}{\partial\overline{z}_k}, j=1,\ldots,m-1
\end{equation}
\begin{equation}\label{autoab}
\begin{bmatrix}
\frac{\partial \rho_1^{M_c}(x,y)}{\partial\overline{z}_1} & \cdots & \frac{\partial \rho_1^{M_c}(x,y)}{\partial\overline{z}_n}\\
\frac{\partial \rho_1(x,y)}{\partial\overline{z}_1} & \cdots & \frac{\partial \rho_1(x,y)}{\partial\overline{z}_n}\\
\vdots & \ddots & \vdots\\
\frac{\partial \rho_d(x,y)}{\partial\overline{z}_1} & \cdots & \frac{\partial \rho_d(x,y)}{\partial\overline{z}_n}\\
\end{bmatrix}
\begin{bmatrix}
d_{1,1} & \cdots &d_{1,m-1}\\
\vdots & \ddots & \vdots\\
d_{n,1} & \cdots &d_{n,m-1}\\
\end{bmatrix}=0
\end{equation}
and
\begin{equation}
\begin{bmatrix}
d_{1,1} & \cdots &d_{1,m-1}\\
\vdots & \ddots & \vdots\\
d_{n,1} & \cdots &d_{n,m-1}\\
\end{bmatrix}
\end{equation}
has rank $m-1$.
Clearly, we can use
\begin{equation}
\begin{bmatrix}
d_{1,1} & \cdots &d_{1,m-1}\\
\vdots & \ddots & \vdots\\
d_{n,1} & \cdots &d_{n,m-1}\\
\end{bmatrix}=
\begin{bmatrix}
c_{1,1} & \cdots &  c_{1,m-1}   \\
\vdots & \ddots &  \vdots  \\
c_{n,1} & \cdots &    c_{n,m-1} \\
\end{bmatrix}
\end{equation}
to obtain that the $\tilde{L}_j$ are restrictions, for $p\in M_c$, of
$L_j$ for $j=1,\ldots,m-1.$ 
so any function that is annihilated by
$L^\gamma$ for all $\gamma\in \N^n$, with $\abs{\gamma}=q$ 
will in particular have restriction to $M_c$ that will be, near $(0,c,0)$, annihilated by
$(\tilde{L}_1,\ldots,\tilde{L}_{m-1})^\gamma$. In other words $f|_{M_c}$ is $q$-analytic (in the sense of Ahern-Bruna) near $(0,c,0)$ for each $c$ near $0$.
 \begin{lemma}\label{crabbelemma}
 	Suppose $f$ is a $q$-analytic function near $0$ on $M,$ $q\in \Z_+,$ $q>1$. If each $L_jf$ is $(q-1)$-pseudoanalytic near $0$ then $f$ is $q$-pseudoanalytic near $0.$ (note that the case $q=1$ is trivial since then $f$ is a $CR$ function and local $q$-pseudoanalyticity follows from the Baouendi \& Treves approxiamtion theorem).	
 \end{lemma}
 \begin{proof}
 	We shall use double induction in 
 	the $CR$ dimension $m$ 
 When $m=0$ then $M$ is totally real, i.e.\ $H^{0,1}M=\{0\}$ hence any differentiable function on $M$ is a $CR$ function. So assume also that $m\geq 1$ and that the result holds true for $CR$ dimension $j=1,\ldots,m-1.$	
 Since each $f|_{M_c}$ is $q$-analytic (in the sense of Ahern-Bruna) near $0$ for each $c$ near $0$, and since the restriction of the $q-1$-pseudoanalytic functions $L_jf$ 
 are $q-1$-pseudoanalytic and for $j=1,\ldots,m-1$ coincide with $\tilde{L}_j$ near $0$, we have by the induction hypothesis (with respect to $m$) that each $f|_{M_c}$ is $q$-pseudoanalytic. By proposition \ref{crpolyprop1}
 there exists differentiable $CR$ functions $a^c_\beta,$ $\beta\in \N^n,$ $\sum_j \beta_j< q,$ such that $f|_{M_c}$ has, near $(0,c,0)$, the representation
 \begin{equation}\label{glueabbecr}
 f|_{M_c}(z)=\sum_{\abs{\beta}<q} a^c_\beta(z) \bar{z}^\beta
 \end{equation}
 and also there exists $CR$ functions $a'_{j,\beta},$ $\beta\in \N^n,$ $\sum_j \beta_j< q-1,$ such that $(L_jf)$ has, near $0$, the representation
  \begin{equation}\label{gluetwo}
 (L_jf)(z)=\sum_{\abs{\beta}<q-1}a'_{j,\beta}(z) \bar{z}^\beta,\quad j=1,\ldots,m
 \end{equation}
 But as we have seen the $\tilde{L}_j$ are local restrictions, for $p\in M_c$, near $(0,c,0),$ of
 $L_j$ for $j=1,\ldots,m-1,$ which implies (since the $a^c_\beta$ are $CR$ functions) that  
 \begin{equation}
 \sum_{\abs{\beta}<q} a^c_\beta(z) L_j\bar{z}^\beta=\sum_{\abs{\beta}<q-1}a'_{j,\beta}(z) \bar{z}^\beta,\quad j=1,\ldots,m-1
 \end{equation}
 for $z\in M_c$ near $(0,c,0).$ Since $L_j\bar{z}^\beta$, $j=1,\ldots,m-1,$ is independent of $c$
 the functions 
 $a_\beta^c(z),$ combine to a smooth function $\tilde{a}_\beta$ of $z$ near $0$ in $M,$ hence the functions given by Eqn.(\ref{glueabbecr}) define a smooth representation for $f$ according to
 \begin{equation}
f(z)=\sum_{\abs{\beta}<q} \tilde{a}_\beta(z) \bar{z}^\beta, \mbox{ near $0$ in $M$}
 \end{equation}
 But for each $p$ near $0$ in $M$ we know that $\im p_1=c_p$ for some $c_p,$ in other words
 each $p$ near $0$ in $M$ belongs to some $M_{c_p}$, 
  \begin{equation}
 L_j\tilde{a}_\beta(p) =L_j|_{M_{c_p}} \tilde{a}_\beta(p)= L_j a^{c_p}_\beta(p),\quad j=1,\ldots,m
 \end{equation} 
 Also since each $\tilde{L}_j$ is the restriction of $L_j$ for $j=1,\ldots,m-1,$ we have 
  \begin{equation}\label{fortisju}
 L_j\tilde{a}_\beta(p) =L_j|_{M_{c_p}} \tilde{a}_\beta(p)= L_j a^{c_p}_\beta(p)=0,\quad j=1,\ldots,m-1
 \end{equation}  
 But the $a^{c_p}_\beta$ are $CR$ functions on a generic $CR$ submanifold thus by the Baouendi \& Treves approximation theorem there exists for each $\beta,$ a sequence $\{E^c_{k,\beta}\}_{k\in \Z_+}$ of entire functions such that
 $E_{k,\beta}\to a^{c_p}_\beta$ uniformly near $0$ in $M_{c_p}.$
 This implies
 \begin{multline}\label{lafrth}
 L_m \partial_{\bar{z}}^\kappa f(p)=L_m\partial_{\bar{z}}^\kappa\sum_{\abs{\beta}<q} a^{c_p}_\beta(p) \bar{z}^\beta
 =\\
 L_m\lim_{k\to \infty} \partial_{\bar{z}}^\kappa\sum_{\abs{\beta}<q}E_{k,\beta}(p) ((\bar{z}^\beta)(p))=
 L_m\sum_{\abs{\beta}<q}a^{c_p}_\beta(p) (\partial_{\bar{z}}^\kappa\bar{z}^\beta)(p)
 \end{multline} 
 Also the $a'_{j,\beta}$ in Eqn.(\ref{gluetwo}) are $CR$ functions thus the approximation argument repeated yields
  \begin{equation}
 \partial_{\bar{z}}^\kappa \sum_{\abs{\beta}<q-1}a'_{m,\beta}(z) \bar{z}^\beta=
 \sum_{\abs{\beta}<q-1}a'_{m,\beta}(z) (\partial_{\bar{z}}^\kappa \bar{z}^\beta)
 \end{equation} 
 Choosing $\kappa\in \N^n$ with $\abs{\kappa}=q-1$ Eqn.(\ref{lafrth}) yields together with Eqn.(\ref{gluetwo})
  \begin{equation}
 L_m a^{c_p}_\kappa(p) = 
 \partial_{\bar{z}}^\kappa (L_m f)(p)=
 \sum_{\abs{\beta}<q-1}a'_{m,\beta}(p) (\partial_{\bar{z}}^\kappa\bar{z}^\beta)(p))=
 0
 \end{equation} 
 Hence we have pointwise for $p$ near $0$ in $M$
 \begin{equation}
 L_m \tilde{a}_\kappa(p) =L_m a^{c_p}_\kappa(p) = 0,\mbox{ for all }\kappa\in \N^m,\abs{\kappa}=q-1
 \end{equation} 
 which together with Eqn.(\ref{fortisju}) implies that each 
 $\tilde{a}_\kappa,$ $\kappa\in \N^m,\abs{\kappa}=q-1$, is a $CR$ function near $0$ in $M.$
 We denote $f=f_1+\sum_{\abs{\beta}=q-1}\tilde{a}_{\beta}(z) \bar{z}^\beta$.
 Replacing $f$ by $f_1$, $q$ by $q-1$, and repeating the arguments  
 and letting $\kappa\in \N^m$ run over all multi-indices with $\abs{\kappa}=q-2,$ we obtain that
 \begin{equation}
 L_m \tilde{a}_\kappa(p) = 0,\mbox{ for all }\kappa\in \N^m,\abs{\kappa}=q-2
 \end{equation} 
 and we set $f=f_2+\sum_{\abs{\beta}\geq q-2}\tilde{a}_{\beta}(z) \bar{z}^\beta$.
 Hence each 
 $\tilde{a}_\kappa,$ $\kappa\in \N^m,\abs{\kappa}=q-2$, is a $CR$ function near $0$ in $M.$
 This can be it iterated each time replacing $f_r$ by $f_{r+1}$, $q-r$ by $q-r-1$, 
 in order to obtain that
 each 
 $\tilde{a}_\beta,$ $\beta\in \N^m,\abs{\beta}<q$, is a $CR$ function near $0$ in $M.$
 This shows, by Proposition \ref{crpolyprop1}, that $f$ is $q$-pseudoanalytic near $0$ in $M$. This completes the induction. This completes the proof of Lemma \ref{crabbelemma}. 	
 \end{proof}
We can now proceed by induction in $q$. If $q=1$ then $f$ is a $CR$ function and local $q$-pseudoanalyticity follows from the Baouendi \& Treves approximation theorem. So assume $q>1$ and that the result holds true for $k$-analytic functions, $k=1,\ldots,q-1.$ By the induction hypothesis each $L_jf$ is $q-1$-pseudoanalytic near $0$, $j=1,\ldots,m$. By Lemma \ref{crabbelemma} this implies that $f$ is $q$-pseudoanalytic near $0$. This completes the proof of Proposition  \ref{crpolyprop2}.
\end{proof}
\begin{corollary}
	Let $M$ be a $C^q$-smooth generic $CR$ submanifold in $\C^n$. Let $q\in \Z_+.$ Then $f$ is
	$q$-analytic on $M$ only if it locally has a representation of the form given by Eqn.(\ref{crforstaekv}).
\end{corollary}
In the proof of Proposition \ref{crpolyprop2} we have considered $CR$ hypersurfaces of larger generic $CR$ submanifolds and used $CR$ restrictions of $CR$ functions, however for very special cases. We mention that 
that if $M\subset\tilde{M}$ is a hypersurface such that both $M$ and $\tilde{M}$ are $CR$ submanifolds, then the following is known.

Suppose $M$ is not a complex manifold near $0$. Since $M$ is a hypersurface in $\tilde{M}$, we have $T_p M\neq T^c_p \tilde{M}$ for $p$ near $0$ (this is sometimes called
that $M$ is a {\em noncharacteristic hypersurface}\index{Noncharacteristic hypersurface} in $\tilde{M}$.
If $\tilde{f}$ is a $CR$ function on $\tilde{M}$ which vanishes near a point $p$ in $M$ then $\tilde{f}$ vanishes near $p$ in $\tilde{M}$  (see Hunt, Polking \& Strauss \cite{hunt}, Corollary 3.2, p.441).
\begin{remark}
We stress that the representation in Eqn.(\ref{crforstaekv}) is {\em not} necessarily unique, as it is well-known that
generic $CR$ submanifolds are not sets of uniqueness for $\alpha$-analytic functions in the non-holomorphic case.
\end{remark}
\begin{example}
	In $\C^2$ we consider the flat generic submanifold
	$M:=\{z\in \C^2:\im z_2=0\}$ that the functions $f(z_1,z_2):=(z_1+z_1^2z_2)-z_1^2\bar{z}_2$, $g(z_1,z_2):=(z_1+z_2)-\bar{z}_2$. Then $f,g$ are both $q$-analytic in the sense of Ahern-Bruna with $q=2,$ and $f|_M=g|_M=z_1$.
	\end{example}

\subsection{The method of analytic discs}
In the study of holomorphic extension an old but effective technique consists of the following:
Given a domain $\Omega\subset\Cn$, assume that there exists a bounded one-dimensional complex submanifold
$\Lambda\subset \Cn$ such that there is a relative domain $\lambda\subset\Lambda$ with relative boundary $\partial \lambda\subset \Lambda\cap K$ for some compact $K\subset \Omega.$
Now for each $f\in \mathscr{O}(\Omega)$ there exists a sequence of entire functions, say $E_j,$
such that $E_j\to f$ uniformly on $K.$ 
Denote $H_j:=E_j|_\Lambda.$ These are again holomorphic.
By the maximum principle for holomorphic functions, we know that $\max_{\overline{\lambda}}\abs{H_j}\leq
\max_{\partial\lambda}\abs{H_j},$ thus the uniform convergence of the $H_j$ on $\partial\lambda$
implies that the $H_j$ converge uniformly with respect to the sup-norm, on $\overline{\lambda}.$
The bounded uniform limit of analytic functions on a domain $\lambda$ defines an holomorphic function on $\lambda.$
Hence, if we are able to find a family $\{\lambda_\iota\}_{\iota\in I}$ (where $I$ is some index set) each with the described property such that
the union $U:=\cup_{\iota\in I} \lambda_\iota$ contains an open set $V\subseteq U.$ Then we know that 
the $E_j$ converge uniformly on $V$. Thus the uniform limit defines a holomorphic extension of 
$f$ to $V$ (if $V\cap \Omega\neq \emptyset$). The idea can be used for example in proving the Hartog's extension theorem where one can cover a compact properly lying inside a given domain with a family of one-dimensional complex submanifolds of $\Omega.$
\\
\\
The idea extends to the case of $CR$ submanifolds, and we describe this in Section \ref{crapp}, but some additional effort is then required.
In particular, $\Omega$ is no longer ambiently open (it has complex dimension strictly lower than $n$). One usually needs some generalized notion of convexity and a condition with respect to that definition in order make sure that there exists $V$ such that $V\cap \Omega\neq \emptyset$. Furthermore, the existence of the approximating sequence $E_j$ is a result for ambient holomorphic functions. Perhaps surprisingly, there exists a local version, however, for $CR$ functions (see Section \ref{crapp}) due to a result of Baouendi \& Treves \cite{bt}.
The conclusion of the following proposition (i.e.\ uniqueness) holds true for general holomorphic extension from generic $CR$ submanifolds, but as we shall see, has a particularly simple proof 
in the case when the extension is obtained via the method of analytic discs.
\begin{proposition}
Holomorphic extension obtained via the method of analytic discs, from generic $CR$ submanifolds of $\Cn$, is always unique.
\end{proposition}
\begin{proof}
Assume that $M\subset \Cn$ is a generic $CR$-submanifold, let $p_0\in M$. Let $f$ be a continuous $CR$ 
function on a neighborhood $U$ of $p_0$ in $M.$ By the Baouendi \& Treves approximation theorem 
see Section \ref{crapp} there exists an open neighborhood $V\subset\Cn$ of $p_0$, with $W\cap M\subset U,$ 
together with a sequence of entire functions, $E_j$, $j\in\Z_+,$ such that $E_j|_{W\cap U}\to f|_{W\cap U},$ uniformly. 
Assume now that existence of holomorphic extension via the analytic disc method holds true at $p_0.$
By this we mean the following: 
Assume that there is a family $\{\lambda_\iota\}_{\iota\in I}$ (for an index set $I$) 
of analytic discs whose boundaries are attached to some open neighborhood $C\subseteq W\cap U$ of $p_0$ in $M$.
Assume also that there exists an open $V\subset \cup_{\iota\in I} \lambda_\iota$ in in the ambient $\Cn.$
These two assumptions imply that $f$ has a holomorphic extension $\tilde{f}$ to $V$.
That extension is obtained (by the arguments of the method of analytic discs described above)
as the uniform limit of the $E_j$ on $V$. Now assume $g$ is another continuous $CR$ function
on $U$. Then it has a holomorphic extension $\tilde{g}\mathscr{O}(V)$ given by the uniform limit of entire function say $F_j$, $j\in \Z_+.$
Then $(f-g)|_{U\cap W} =0$ and has in turn via the method of analytic discs
a holomorphic extension, $\tilde{h}$, to $V$, given by the uniform limit of $(E_j-F_j).$ So if $f-g=0$ on $U\cap W$,
then (by the maximum principle for each $(E_j-F_j)$ the uniform limit $\tilde{h}$ is zero, i.e.\ $\tilde{f}=\tilde{g}.$
This completes the proof. 
\end{proof}

\begin{proposition}
Suppose $M\subset \Cn$ is a generic $CR$-submanifold, let $p_0\in M$ and suppose that
any differentiable $CR$ function on a neighborhood $U$ of $p_0$ in $M,$ 
has holomorphic extension, near $p_0$ to some open $V\subset\Cn.$
Let $f$ be a differentiable $CR$ on $U$ and 
denote the holomorphic extension, $\tilde{f}$, to  
$V$. Then $\tilde{f}(V)\subseteq f(U).$
\end{proposition}
\begin{proof}
If $h$ is a nowhere zero differentiable $CR$ function on a generic $CR$ submanifold $N\subset\Cn$, 
then $1/h$ is a differentiable $CR$ function. To see this note that for any $CR$ vector field
$L$ we have (since $h$ is nowhere zero) $L(1/h) =(-1/h^2)Lh=0$. Hence $1/h$ is annihilated by (all sections of) $T^{0,1}N.$
Suppose there exists $q_0\in V$ such that $\tilde{f}(q_0)\notin f(U).$ Then the function
$g(z):=1/(f(z)-\tilde{f}(q_0))$ is a continuous $CR$ function on $U$, but has no holomorphic extension to $V$.
Hence by contradiction $\tilde{f}(V)\subseteq f(U).$ 
\end{proof}
The method of analytic discs is perhaps the most powerful technique of holomorphic extension and in some sense availability 
of certain versions of it can be viewed upon as 
necessary for holomorphic extension and propagation of extension, in higher complex dimension. 
For this reason we have put some effort into describing it in this book.
The method of analytic discs has refinements that can be considered as pioneering techniques in this regard due to Alexander Tumanov (see e.g. Tumanov \cite{tumanov}, \cite{tumanov1}) who used so-called defective analytic discs attached to a 
$CR$ sumbanifold, as propagators of extendability.
Other notable early works are due to Hanges \& Treves \cite{hanges} and Trepreau \cite{trepreau}.

\subsection{Induced $q$-analytic extension}
Obviously, the results of the method of analytic discs (but not really the technique itself) extends in the {\em ambient} space, to the case of $\alpha$-analytic functions, in the following sense:
given a domain $\Omega\subset\Cn$, assume there exists a non-empty open set $V\subset\Cn$
such that $V\cap(\C\setminus \Omega)\neq \emptyset$ satisfying that there exists a subdomain
$\omega\subseteq \Omega$ such that
any holomorphic function on $\Omega$ has holomorphic extension to an open subset $U\subset V.$
Then any function $f$ that is $\alpha$-analytic on $\Omega$
extends (not necessarily uniquely) to an $\alpha$-analytic function $F$ on $\Omega \cup V.$
In other words whenever holomorphic extension holds true to some open set 
then $\alpha$-analytic
extension (but not necessarily a unique one) also holds true and to the same set.
Of course this is because for a domain $\Omega$ we can assume $f(z)$ has the representation $\sum_{j=1}^n\sum_{0\leq \beta_j< \alpha_j} a_\beta(z)\bar{z}^\beta,$
for holomorphic $a_\beta(z)$ and by assumption each $a_\beta(z)$ has a holomorphic extension $\tilde{a}_\beta(z)$ to $\Omega\cup V.$
Furthermore, the functions $\bar{z}^\beta$ are entire $(\beta+(1,\ldots,1))$-analytic, hence the function
$\tilde{f}(z)=\sum_{j=1}^n\sum_{0\leq \beta_j< \alpha_j} \tilde{a}_\beta(z)\bar{z}^\beta,$ yields the sought $\alpha$-analytic extension.
The situation is similar in the case $CR$-submanifolds.
As before, for the method of analytic discs to apply we would need boundary uniqueness which obviously
fails for $\alpha$-analytic functions in general. 
However, we have the following. 
\begin{proposition}[Induced extension for $q$-analytic functions in $CR$ geometry]\label{crpolyprop3}
Let $M$ be a $C^q$-smooth generic $CR$ submanifold in $\C^n$. Let $q\in \Z_+$ and let $0\in M.$
Assume that simultaneous local holomorphic extension of all continuous $CR$ functions near $0$ holds true to some open subset in $\Cn$. 
Then 
any function $f$ that is $q$-pseudoanalytic on $M$ 
extends (not necessarily uniquely) to function $F$, a $q$-analytic in the sense of Ahern-Bruna, on an (ambiently in $\Cn$) open domain $W\subset \Cn.$ 
\end{proposition}
\begin{proof}
By Proposition \ref{crpolyprop2} (Proposition \ref{crpolyprop1} respectively) any $f$ satisfying the conditions of Proposition \ref{crpolyprop3}
has the property that near each $0\in M$
there
exists differentiable $CR$ functions $a_\beta,$ $\beta\in \N^n,$ $\sum_j \beta_j< q,$ such that $f$ has, near $0$, the representation
\begin{equation} 
f(z)=\sum_{\abs{\beta} <q} a_\beta(z) \bar{z}^\beta
\end{equation}
where $\bar{z}^\beta$ denotes the restriction to $M$
of $\bar{z}^\beta.$
By assumption each of the (finitely many) $a_\beta(z)$ has simultaneous local holomorphic extension $\tilde{a}_\beta(z)$ to an open subset of $\Cn$. 
Furthermore, the functions $\bar{z}^\beta$ are entire $(\beta+(1,\ldots,1))$-analytic, hence the function
$F(z)=\sum_{j=1}^n\sum_{\abs{\beta} <q} \tilde{a}_\beta(z)\bar{z}^\beta,$ yields the sought extension.
This completes the proof.
\end{proof}
Again, we stress that the main difference compared to the theory of $CR$ functions is that the $q$-analytic extension is not necessarily unique.
There are many known conditions for the generic $CR$
submanifold to allow holomorphic extension of continuous $CR$
functions, and as we have mentioned, the method of analytic discs usually appears in extension proofs. For instance when $M$ is a smooth hypersurface it is sufficient for holomorphic extension near a point $p_0\in M$ that $M$ be strictly pseudoconvex near $p_0.$ In higher codimension this condition must be replaced with conditions on the so-called Levi-cone, see e.g.\ Boggess \cite{b4}.
We can also introduce the following notion.
\begin{definition}[Cf. Definition \ref {hypocomplexdef}] 
Let $M\subset\Cn$ be a smooth generic $CR$ submanifold.
$M$ is called $q$-hypocomplex\index{$q$-hypocomplexity} at $p_0\in M$ if for each $q$-pseudoanalytic function, $f$, near $p_0$, there exists an open neighborhood $U$ of $p_0$ in $\C^n$, a multi-index $\alpha\in \N^n,$ $\sum_j \alpha_j=q,$ and a local ambient $q$-analytic function
$F$ on $U$
such that $F|_{M\cap U} =f|_{M\cap U}.$
\end{definition}
Obviously we have the following (since local holomorphic extension of $CR$ functions is precisely the case of $q$-hypocomplexity for $q=1$).
\begin{corollary}
[To Proposition \ref{crpolyprop3}]
Let $M\subset\Cn$ be a smooth generic $CR$ submanifold. Then
$M$ is $q$-hypocomplex at a point $p_0\in M$ if and only if
any $CR$ function near $p_0$ has local holomorphic extension to an ambient neighorhood of $p_0$ in $\Cn$
(Cf. Definition \ref {hypocomplexdef}).
\end{corollary}
We have given the results on hypocomplexity in terms of conditions of local holomorphic extension of $CR$ functions. We have not reviewed the vast theory of sufficient geometric conditions for the existence of such extension, but we believe that with the knowledge we have provided on the method of analytic discs, the interested reader will be able to rather quickly commence the reading of the cited literature and explore this field. 
Here we give just one example. 
As the following is somewhat peripheral to the central topic of this book, we shall use here some notions without giving all the details, they are however described more closely in Section \ref{crapp}, Section \ref{lewyapp}.

\begin{theorem}[See Boggess \cite{b4}, p.202]
Let $M\subset \Cn$ be a $C^l,l\geq 4,$ generic embedded $CR$
submanifold and let $p\in M$ be a point such that the Levi cone satisfies $\Gamma_p=N_p(M)$. Then for each
neighborhood $\omega$ of $p$ in $M$ there exists an open $\Omega\subset\Cn$, $p\in \Omega$ such that
each $C^1$-smooth $CR$ function on $\omega$ is the restriction of a unique holomorphic function defined on $\Omega$. 
\end{theorem}
Sometimes in the literature sufficient conditions for local holomorphic extension to a full ambient neighborhood is stated in terms of the directional Levi form at a given point $p\in M$. Namely, it is sufficient that for all nonzero codirections, $\xi\in \chi_p$ with
$\xi\neq 0$, the Hermitian form $\mathcal{L}^{\xi}(\cdot )$ has at least one positive and one negative eigenvalue.
This is particularly illustrative in the case of hypersurfaces in $\Cn$, $n\geq 2,$ (a result due to Lewy \cite{lewy}, see Boggess \cite{b4}, p.199). 
In the case of one-sided holomorphic extension from a strictly pseudoconvex hypersurface $M\subset\Cn$, local convexification (see Section \ref{localconv}) immediately renders the wanted family of analytic discs locally attached to the hypersurface, namely we can, assuming $0$ belongs to $M$, after local convexification write the hypersurface near $0$ in some holomorphic coordinates $z$ for $\Cn$ as
$\{ \im z_n=\abs{(z_1,\ldots,z_{n-1})}^2 +\abs{\re z_n}^2\}.$ Hence each $\{\re z_n=c\}$
is a complex manifold $(n-1)$ dimensional manifold foliated by analytic discs attached to $M$ near $0.$
\section{Relation to hypoanalytic theory}
Now $CR$ geometry can be viewed upon as a special branch of hypoanalytic theory or the
theory of locally integrable structures (see  
the books of Treves \cite{treves}, Berhanu, Cordaro \& Hounie \cite{ber2}, Baouendi, Ebenfelt \& Rothschild \cite{barot} and
Treves \& Cordaro \cite{trevescordaro}). Here one studies distribution solutions are the ones annihilated
by the sections of a locally integrable subbundle of the complex tangent bundle of a real smooth manifold, i.e.\ 
they are distribution solutions to a system of vector fields. 
The theory of polyanalytic functions and the theory of locally integrable structures converge in the case of the holomorphic functions, 
which in one variable are precisely the
$1$-analytic functions and simultaneously the solutions to the locally integrable structure generated over a complex submanifold of $\C$, by
the Cauchy-Riemann vector field $L=\frac{1}{2}\left(\frac{\partial}{\partial x}+i\frac{\partial}{\partial y}\right)$, where $x+iy$ is the standard holomorphic coordinate for $\C$. 
It therefore makes sense to attempt to characterize the solutions in one theory in terms of 
the other theory, possibly allowing a deeper understanding of both subjects. Daghighi \cite{daghighihypoanal} gave such a characterization which we shall describe in this section. To do this we need some basic preliminaries on locally integrable structures.

\subsection{Preliminaries on locally integrable structures}

\begin{definition}[Locally integrable structure, see e.g.\ Berhanu, Cordaro \& Hounie \cite{ber2} p.19] \label{hgh}
Let $\Omega$ be a real $N$ dimensional manifold and let $L\subset \C\otimes T\Omega$ be a subbundle of rank $l.$ 
Set $L_p^{\perp}:=\{ \xi\in \C\otimes T_p^* \Omega: \xi=0\mbox{ on }L_p\}.$ 
$L$ is called a {\em locally integrable structure} on $\Omega$ (for clarity we shall say that
$(\Omega,L)$ is a locally integrable structure)
if for each $p_0\in \Omega$ there is an open neighborhood $p_0\in U\subset \Omega$
together with functions $Z_{1},\ldots ,Z_{n}\in C^{\infty}(U),$
$l=N-n,$ 
such that span$\{dZ_{1p},\ldots ,dZ_{np}\}=L_p^{\perp},$  
$\forall p\in U.$ 
\end{definition}

\begin{definition}[Hypoanalytic structure, see e.g.\ Baouendi, Chang \& Treves \cite{changtreves}]
Let $\Omega$ be a real smooth manifold of dimension $n+l.$ 
A {\em hypo-analytic structure} on 
$\Omega$ is a collection 
of pairs $(U_k,Z_k)$ where $U_k\subset \Omega$ is open and
$Z_k:=(Z_{k1},\ldots,Z_{kn}):U_k\to \C^n$ a smooth map  
s.t.:
\begin{itemize}
\item[(i)] $\{U_k\}$ is an open cover of $\Omega.$
\item[(ii)] $dZ_{k1},\ldots,dZ_{kn}$ are $\C$-linearly independent at each point of $U_k.$ 
\item[(iii)] If $k\neq k'$ and $p\in U_k\cap U_{k'}$ there exists a holomorphic map $F_{k',p}^k$ of an open neighborhood of $Z_k(p)$ in $\C^n$
  into $\Cn$ such that $Z_{k'}=F^k_{k',p}\circ Z_k$ in a neighborhood of $p$ in $U_k\cap U_{k'}.$ 
\end{itemize}
The span of the differentials is denoted $T'$ and the underlying locally integrable structure 
(obtained with respect to the duality between tangent vectors and cotangent 
vectors) 
will be denoted $L=(T')^{\perp}$. 
If $\Omega$ and $Z$ are real-analytic then the hypoanalytic structure is called real-analytic (see Treves \cite{treves}, p.122).
\end{definition}
We shall be working locally and concern ourselves with a fixed hypoanalytic chart and hence will not make any reference to the maps
$F_{k',p}^k$ in our proofs.
Given a locally integrable structure $(\Omega,L),$ it locally near any point $p$, 
say on an open $U,$ $p\in U,$ underlies at least one hypoanalytic structure, which for sufficiently small $U$ can be given with a single chart map $Z$ (this is immediate from Definition \ref{hgh}).
It should be noted however that the same locally integrable structure can underlie significantly different hypoanalytic structures, see e.g.\
Baouendi, Chang \& Treves \cite{changtreves}, p.335.
Once the local hypoanalytic chart map $Z$ (whose differentials span $L^{\perp}$) is known and fixed, one can for clarity use the notation
$(U,L,Z)$ for the given local hypoanalytic structure however the underlying locally integrable structure
$L$ will be locally uniquely determined by the hypoanalytic chart, thus is sufficient to simply denote
the local hypoanalytic structure (with only one chart map) by $(U,Z).$
Even though a hypoanalytic structure determines a unique locally integrable structure, 
namely the orthogonal of the cotangent structure bundle spanned by the differentials of the components of the hypoanalytic chart maps, the choice of hypoanalytic chart is not necessarily unique, i.e., different hypoanalytic structures can yield the same underlying locally integrable structure. 

\begin{definition}[Solution with respect to a locally integrable structure]\label{solss}
Let $\Omega$
be a smooth $(n+l)$-dimensional real manifold
and let $L\subset\C\otimes T\Omega$ be a locally integrable structure
of rank $l$.
Let $L_1,\dots ,L_l,$ define a local basis near a reference point $p_0$
for $L.$ 
A distribution $u$ defined on an open $V\ni p_0$
is called a {\em solution} (with respect to the system induced by $L$) on $V$ if
$L_ju=0,$ on $V,$ $1\leq j\leq l$. If $\Omega$ is real-analytic, and the $L_1,\dots ,L_l,$ are real-analytic (in the sense that
their expressions have real-analytic coefficients with respect to the coordinates of $\Omega$) then we shall call a real-analytic function $u$ a {\em real-analytic 
solution} (with respect to the given locally integrable structure).
\end{definition}

\begin{definition}[Hypoanalytic function, see Treves \cite{treves}, p.123]\label{hypoanfunc}
Let $\Omega$
be a smooth $(n+l)$-dimensional real manifold
equipped with a hypoanalytic structure and let $p\in \Omega.$
A complex-valued function $f$ defined in some neighborhood of $p,$
is said to be {\em hypoanalytic at $p$}, if there exists a hypoanalytic chart $(U,Z),$
$p\in U, Z=(Z_1,\ldots,Z_n)$
(as part of the hypoanalytic structure of $\Omega$)
defined on a neighborhood of $p$ together with a holomorphic function $F$ in
some neighborhood of $Z(p)$ in $\C^n$, 
such that $f=F\circ Z$ in a neighborhood of $p.$
A complex-valued function $f$ is called {\em hypoanalytic} on a subset $S\subset\Omega$, if it is hypoanalytic at each point of $S.$
\end{definition}

\begin{example}\label{exex}
Let us see how a generic $CR$ submanifold of $\Cn$ can be realized from the perspective of locally integrable structures.
Recall that a real smooth submanifold $M\subset\Cn$, of real dimension $n+l$ and real codimension $d=n-l$, is a subset such that for each
$p_0\in M$, there is
an ambient open neighborhood $U\subset\C^n$ of $p_0$ satisfying
$M\cap U=\{\rho=0\},$ where
$\rho: U\to \R^d$ and $d\rho_1,\ldots,d\rho_{n-l}$ are linearly independent in $U.$
Let $p\in M,$ and let $J$ denote the complex structure map on $T_p\C^n$ ($J^2=-\mbox{Id}$). It has $\C$-linear extension to
$\C\otimes T_p\C^n$.
If the complex differentials $\partial \rho_1,\ldots,\partial \rho_d$ are $\C$-linearly
independent at each point, $p,$ then $M$ is called {\em generic}, which in turn is equivalent to (see e.g.\ Baouendi, Ebenfelt \& Rothschild \cite{barot}, p.14)
$T_p M+JT_p M=T_p\C^n$. The complexification $\C\otimes T_p M\cap JT_p M$
can be decomposed with respect to the eigenvalues of $J,$
$\C\otimes T_pM\cap JT_p M=H^{0,1}_p M\oplus H^{1,0}_p M$ (where $H^{0,1}_p M=\{X\in \C\otimes T_p M\colon J(X)=-iX\}$). 
Denote
$L_p:=H_p^{0,1} M,$ and note that $L|_p\cap \overline{L}|_p=\{0\},$ $L|_p+ \overline{L}|_p=T_p\Cn.$ 
Suppose $M\subset\Cn$ is a real-analytic submanifold and that for each $p_0\in M$ there is an open neighborhood $V$ of $p_0$ in $\Cn,$ together with
holomorphic coordinates $(z,w)$ for $V$ centered at $p_0,$
such that $z=(z_1,\ldots,z_l),$ $w=(w_1,\ldots,w_{n-l}),$
and $V\cap M=\{ (z,w)\in V:\im w=\phi(z,\re w)\}$ for an $\R^{n-l}$-valued real-analytic $\phi$ with $d\phi(0)=\phi(0)=0.$
Then the function $\rho=\im w-\phi(z,\re w)$ will satisfy the condition that 
the complex differentials $\partial \rho_1,\ldots,\partial \rho_d$ are $\C$-linearly
independent near $p_0$. The dimension of $T_p M\cap JT_p M$ is constant near $p_0$,
$L:=H^{0,1} (V\cap M)$ is an involutive subbundle (see e.g Boggess \cite{b4}, p.101) of rank $l$, that satisfies 
$L\cap \overline{L}=\{0\}.$ It is known that this implies that this subbundle
is locally integrable see e.g.\ Baouendi, Ebenfelt \& Rothschild \cite{barot}, p.38. $M$ is
what we are used to calling a generic $CR$ submanifold of $\Cn$.
In the proof of the main result of this section we shall be working with local hypoanalytic structures that are associated to 
such manifolds. 
Assume $p_0=0$ and 
denote $z:=x+iy$, $s:=\re w,$ which means that we have for a real-analytic coordinate system 
$(x,y,s)$ on a sufficiently small open $0\in\omega\subset \R^{l+n},$ 
a local representation of $M\cap V,$ via
$Z(x,y,s):\omega\to \Cn$, via
\begin{align} \label{puh}
& Z_j(x,s,y):=x_j+iy_j,\quad  1\leq j\leq l, \\ \label{puh1}
& Z_j(x,s,y):=s_j +i\phi(x,y,s),\quad  l+1\leq j\leq n 
\end{align}
Then $(Z,\omega)$ defines, a local hypoanalytic structure with a single chart, which overlies
the locally integrable structure $L=H^{0,1}(M\cap V)$ and $Z$ will be an embedding \ Baouendi, Ebenfelt \& Rothschild \cite{barot}, p.37. 
Note that we do not necessarily need
$d\phi(0)=0,$ it is sufficient that the components of $(\im w-\phi(z,\re w))$ have $\C$-linearly independent
complex differentials near $0.$
\end{example}

\begin{definition}[Flatness]\label{strict}
Let $M\subset\Cn$ be a local generic submanifold obtained as the image
of a local hypoanalytic chart $Z$. 
$M$ will be called {\em flat} on a domain $V\subset M$, if $V$ is locally foliated by complex submanifolds of
the same dimension as the complex dimension of the locally integrable subbundle
underlying $Z$. 
\end{definition}

\subsection{A characterization linking the solutions from two theories}
\begin{theorem}\label{mainthmhypoanal}
Let $z=x+iy$ be a holomorphic coordinate in $\C$, 
$W\subset\R^2$ an open subset with real-analytic coordinates $(x,y)$ and let
$p\in W.$ Let $q$ be a positive integer. Then there exists a locally integrable structure $(\omega,L)$ where $\omega\subset \R^4$ is open with real-analytic coordinates
$(x,y,s,t)$, $(p,0)\in \omega,$ such that for any
$q$-analytic function $f(z)$ on $z(W),$
the function $\tilde{f}(x,y,s,t):=f(x+iy)$ is a solution near $(p,0)$, with respect to $L$.
Furthermore, there is for sufficiently small $\omega$ a local hypoanalytic chart, $(\omega,Z)$, $Z\colon \omega\to \C^3$, overlying the locally integrable structure, 
such that the image of $Z$ is flat near $Z(p,0)$ (which under the circumstances means that it is foliated by one-dimensional complex manifolds near $Z(p,0)$)
and such that for any
$q$-analytic function $f(z)$ on $z(W),$ the solution $\tilde{f}$ chosen as above is a local hypoanalytic solution near $(p,0)$ (i.e.\ there exists a holomorphic function
$F,$ in the ambient space $\C^3$, near $Z(p,0)$, such that $F\circ Z=\tilde{f}$ near $(p,0)$).
\end{theorem}

\begin{example}\label{extva}
Let $f(z)=z^2+z\bar{z},$ where $z=x+iy$ is the standard complex coordinate in $\C.$ 
Here we have $a(z)=z^2,$ $b(z)=z.$ 
Let $(z,w_1,w_2)$ denote holomorphic coordinates in $\C^3$ and set,
\begin{align*}
Z_1(x,y,\re w_1,\im w_1) & = z \\
Z_2(x,y,\re w_1,\im w_1) & = \re w_1+i(-\re w_2)  \\
Z_3(x,y,\re w_1,\im w_1) & = \re w_2+i(\re w_1-y)
\end{align*}
This implies that $dZ_1,dZ_2,dZ_3$ are $\C$-linearly independent near $0,$ 
thus $Z$ is a hypoanalytic chart with underlying locally integrable structure given by the 
annihilator $\mbox{span}dZ$. Furthermore, there is a neighborhood $\omega$ of the origin such that,
$Z(\omega)=$
$\{ (z,w_1,w_2)\in \C^3 \colon (\im w_1,\im w_2)=(-\re w_2,\re w_1-y)\},$
and by the proof of Theorem \ref{mainthmhypoanal} the image of $Z$ near $Z(0)$ is flat.
Setting 
$\tilde{f}(x,y,\re w_1,\re w_2):=f(x+iy),$
the function 
\begin{equation}
F(z,w_1,w_2):=a+zb-2ib(w_1+iw_2)=2z^2-2iz(w_1+iw_2)
\end{equation}
is holomorphic near $0$ in $\C^3$ and satisfies, $F\circ Z=\tilde{f},$
so that $\tilde{f}$ is a hypoanalytic solution. 
\end{example}
Theorem \ref{mainthmhypoanal} does not exclude that for a given $q$-analytic function $f$, there exist a choice of locally 
integrable structure having the main wanted properties in Theorem \ref{mainthmhypoanal} but 
underlying, near the reference point, a local hypoanalytic chart whose image is not flat.
However, there are cases
where we can conclude that the chosen locally integrable structure cannot be hypocomplex.
\begin{definition}\label{hypocomplexdef}
A locally integrable structure $(\Omega,L),$ is called {\em hypocomplex at $p$} if there is an open neighborhood $U$ of $p$ and
a smooth local chart near the origin $Z=(Z_1,\ldots,Z_n)\to \Cn$ whose components are solutions such that for any
distribution solution $u$ on a neighborhood $U'\subset \Omega$ of $p$ there is a holomorphic function $\tilde{u}$ 
defined on a neighborhood of $Z(p)$ such that 
$u=\tilde{u}\circ Z$
on a neighborhood of $p$ in $U'.$ 
\end{definition}
Recently it was shown that a locally integrable structure is hypocomplex at a given point if and only if it is hypoelliptic at that point, i.e.\ each solution near that point has a $C^\infty$-smooth representation,
see e.g.\ Daghighi \cite{daghighidiss}.
Note that a nonconstant solution defined on a neighborhood of a hypocomplex point $p$, maps open neighborhoods of $p$ to open neighborhoods in $\C$
(see Treves \cite{treves}, p.149). 
For example if $f(z)=1-z\bar{z},$ then $\abs{f}$ attains a strict local maximum at $0$, which
is impossible if there is a hypoanalytic structure together with a solution $\tilde{f}$
such that $|\tilde{f}|$ attains a weak local maximum at $0$ whenever $\abs{f}$ attains a strict local maximum at $0.$
We point out that in the real-analytic category the choice of overlying hypoanalytic structure, once the locally
integrable structure is fixed, is quite restricted.
\begin{theorem}[Treves \cite{treves}, p.127]\label{trevesRA}
Let $\Omega$ be a real-analytic manifold and $L\subset \C\otimes T\Omega$ is an involutive subundle consisting of real-analytic vector fields.
Then there exists a unique (up to local biholomorphism) hypoanalytic structure on $\Omega,$  real-analytic, overlying $L.$
\end{theorem}

\chapter{Special functions}
This chapter treats some interesting subspaces of, and function spaces closely related to, the poly-rational functions. This includes the poly-rational functions, the countably analytic exponential function (which can be used to solve meta-analytic equations and special polynomial type functions.

\section{Poly-rational functions}
In this section we consider for $\alpha\in \Z_+^n,$ functions on domains $\Omega\subset\Cn$
with representation of the form
\begin{equation}
f(z)=\sum_{\beta_j<\alpha_j} r_\beta(z) \bar{z}^\beta
\end{equation}
for rational functions $r_\beta(z)$ with zeros off $\Omega.$ We call these $\alpha$-rational functions.\index{Poly-rational functions ($\alpha$-rational functions)}
In the case of $n=1$ these appear naturally in the analysis of
$q$-analytic functions that have constant modulus. We shall therefore here prove a characterization of such functions (in terms of $q$-rational functions) and also give some applications of this representation. 
First we shall need the following fundamental result. 
The following was originally announced by Weierstrass and first published by Hurwitz \cite{hurwitz}.
\begin{theorem}\label{hurwitzlem} 
A separately rational function $f(u_1,\ldots,u_n)$, in the sense that $f(u_1,\ldots,u_{j-1},a_j,u_{j+1},\ldots,u_n)$ is a rational function with respect to $u_k, k\neq j$ for fixed $a_j$.
Then $f$ is a jointly rational function of $u\in \Cn.$.
\end{theorem}
\begin{proof}
Let $u:=(u_1,\ldots,u_n)$ and without loss of generality assume $f$ is regular at $0$ and $f(0)\neq 0.$ Then for sufficiently small
$\delta>0$ we have $f\neq 0$ and bounded on $G_1:=\{u:\abs{u_j}\leq \delta,j=1,\ldots,n \},$ and has on that set the expansion
$f(u)=\sum_{j=0}^\infty A_j(u) u_1^j$ where the $A_j(u_2,\ldots,u_n)$ are functions, regular near $0$.
For fix $u_2=b_2,$ $\abs{b_2}<\delta,$ let $0<\delta'<\delta$ be such that the function $f(u_1,b_2)$ is a rational function
with respect to $u_1$ and on $G_1':=\{u:\abs{u_1}\leq \delta,\abs{u_2}\leq \delta' \},$ we have for fixed $u_j=a_j$, $j\neq 2,$
\begin{equation}
f(u)=\frac{C_0+C_1(u_2-b_2)+C_2(u_2-b_2)^2\cdots}{C_0'+C_1'(u_2-b_2)+C_2'(u_2-b_2)^2\cdots}
\end{equation}
where the $C_j,C_j'$ are power series in $(u_k-a_k),$ $k\neq 2,$ such that $C_0$ and $C_0'$ do not simultaneously vanish identically. Also if $C_0=0$ and $C'_0\geq 0$ then $f(u_1,b_2,u_3,\ldots,u_n)$ would be independent of $u_1$ whereas if
if $C_0\geq 0$ and $C'_0= 0$ then $f(u_1,b_2,u_3,\ldots,u_n)$ would be unbounded. 
The function $f(u_1,b_2,u_3,\ldots,u_n)$ is rational with respect to $u_1$, satisfies the induction hypothesis and can, for some $r$, be written
\begin{equation}
f(u_1,b_2,u_3\ldots,u_n)=\sum_{j=0}^\infty \bar{A}_j u_1^j=\frac{B_0+B_1u_1+\cdots+B_r u_1^r}{B_0'+B_1'u_1+\cdots+ B_r'u_1^r\cdots}
\end{equation} 
where $\bar{A}_j=A_j(u_1,b_2,u_3,\ldots,u_n)$, the $B_j,B_j'$ are rational functions of $u_3,\ldots,u_n,$ and where at least on of $B_r,B_r'$ is nonzero.
Comparing powers of $u_1$ in 
\begin{equation}
(\sum_j\bar{A}_ju_1^j)(\sum_{k=1}^r B_k' u_1^k)=\sum_{k=1}^r B_k u_1^k
\end{equation}
yields a system of equations
in which the elimination of the $B_k,B_k',$ $k=0,\ldots,r,$
renders that all determinants of degree $(r+1)$ that can be extracted from the following infinite
array vanish
\begin{equation}
\begin{array}{lll}
\bar{A}_1,& \ldots, &\bar{A}_{r+1}\\
\bar{A}_2,& \ldots, &\bar{A}_{r+2}\\
\bar{A}_3,& \ldots, &\bar{A}_{r+3}\\
\vdots &\vdots & \vdots 
\end{array}
\end{equation}
Hence to each $(u_1,b_2,u_3,\ldots,u_n)\in G_1'$ we have a positive integer $r$, the degree of $f(u_1,b_2,u_3,\ldots,u_n)$ with respect to $u_1$. 
If each $r$ only appears for finitely many choices of $b_2$ then $G'$ would consist of countably many points, thus
there exists an $r$ that appears for infinitely many choices of $b_2.$ 
This implies that the determinants
\begin{equation}
\abs{
\begin{array}{lll}
A_1 & \ldots & A_{r+1}\\
A_2 & \ldots & A_{r+2}\\
A_3 & \ldots & A_{r+3}\\
\vdots &\vdots & \vdots 
\end{array}}
\end{equation}
vanish independent of $u_3,\ldots,u_n,$
for infinitely points near any given point with fixed value of $u_2.$ Since such determinants are functions of $u_2,\ldots,u_n$ which are regular at $u_2=\cdots=u_n=0$,
they must vanish identically. Hence 
from
\begin{equation}
\begin{array}{lll}
f_r\bar{A}_1+ & \cdots &+f_0\bar{A}_{r+1}\\
f_r\bar{A}_2+ & \cdots & +f_0\bar{A}_{r+2}\\
f_r \bar{A}_3 + & \cdots, & +f_0\bar{A}_{r+3}\\
\vdots &\vdots & \vdots 
\end{array}
\end{equation}
there exists functions $f_0,\ldots,f_r$, not all simultaneously vanishing. each depending on the variables $u_2,\ldots,u_n,$ regular at $u_2=\cdots=u_n=0$.
Since there exists $\phi_0,\ldots, \phi_r$ which as $f_0,\ldots,f_r$ are functions of $u_2,\ldots,u_n$ regular at $u_2=\cdots=u_n=0,$ such that
\begin{equation}
f(u)=\abs{\sum_{j=0}^r f_j u_1^j}=
(\sum_{j=0}^\infty A_j u_1^j)(\sum_{j=0}^r f_j u_1^j)=
\sum_{j=0}^r \phi_j u_1^j
\end{equation}
we have
\begin{equation}
f(u)=\frac{\sum_{j=0}^r \phi_j u_1^j}{\sum_{j=0}^r f_j u_1^j}
\end{equation}
The function $f(u)$, which we have shown, for fixed $u_2=b_2$ to be rational with respect to $u_j,$ $j\neq 2$,
for $b_2$ sufficiently close to $u_2=0,$ is also, for fixed $u_1=b_1$ sufficiently close to $b_1,$ rational with respect to $u_2,\ldots,u_n.$ 
Let $b_1^{(1)},\ldots,b_1^{(2r+1)},$ be distinct points sufficiently close to $u_1=0$ such that
\begin{equation}
\sum_{j=0}^{r} (f(u)u_1^j)f_j -\sum_{j=0}^{r} u_1^j\phi_j  =0
\end{equation}
and for $i=1,\ldots,2r+1,$
\begin{equation}
\sum_{j=0}^{r} (f(b_1^{(i)},u_2,\ldots,u_n)(b_1^{(i)})^j)f_j -\sum_{j=0}^{r} (b_1^{(i)})^j \phi_j  =0
\end{equation}
From these $2r+2$ equalities we could solve for $f_j\phi_j,$ $j=0,\ldots,r$ (and thus obtain a linear equation for $f(u),$
from which $f$ can be deduced to be rational) if not all determinants of order $2r+1$ extracted from the following array vanish
\begin{scriptsize}
\begin{equation}
\begin{array}{lllllll}
f(b_1^{(1)},u_2,\ldots,u_n), &\ldots & , f(b_1^{(1)},u_2,\ldots,u_n)(b_1^{(1)})^r,& 1 & b_1^{(1)},&\ldots & , (b_1^{(1)})^r\\
\vdots &\vdots & \vdots & \vdots & \vdots & \vdots & \vdots \\
f(b_1^{(2r+1)},u_2,\ldots,u_n), &\ldots & , f(b_1^{(2r+1)},u_2,\ldots,u_n)(b_1^{(2r+1)})^r,& 1 & b_1^{(2r+1)},&\ldots & , (b_1^{(2r+1)})^r\\
\end{array}
\end{equation}
\end{scriptsize}
In the latter case, however, we have that for arbitrary $u$ all determinants of order $2r+1$ extracted from the following array
(which is obtained by replacing $b_1^{(1)}$ in the previous array, by the unknown $u_1$) vanish
\begin{scriptsize}
\begin{equation}
\begin{array}{lllllll}
f(u_1,u_2,\ldots,u_n), &\ldots & , f(u_1,u_2,\ldots,u_n)(b_1^{(1)})^r,& 1 & u_1,&\ldots & , u_1^r\\
f(b_1^{(2)},u_2,\ldots,u_n), &\ldots & , f(b_1^{(2)},u_2,\ldots,u_n)(b_1^{(2)})^r,& 1 & b_1^{(2)},&\ldots & , (b_1^{(2)})^r\\
\vdots &\vdots & \vdots & \vdots & \vdots & \vdots & \vdots \\
f(b_1^{(2r+1)},u_2,\ldots,u_n), &\ldots & , f(b_1^{(2r+1)},u_2,\ldots,u_n)(b_1^{(2r+1)})^r,& 1 & b_1^{(2r+1)},&\ldots & , (b_1^{(2r+1)})^r\\
\end{array}
\end{equation}
\end{scriptsize}
This renders that $f(u)$ is jointly rational whenever
not all determinants of order $2r$ extracted from the array
\begin{scriptsize}
\begin{equation}
\begin{array}{lllllll}
f(b_1^{(2)},u_2,\ldots,u_n), &\ldots & , f(b_1^{(2)},u_2,\ldots,u_n)(b_1^{(2)})^r,& 1 & b_1^{(2)},&\ldots & , (b_1^{(2)})^r\\
\vdots &\vdots & \vdots & \vdots & \vdots & \vdots & \vdots \\
f(b_1^{(2r+1)},u_2,\ldots,u_n), &\ldots & , f(b_1^{(2r+1)},u_2,\ldots,u_n)(b_1^{(2r+1)})^r,& 1 & b_1^{(2r+1)},&\ldots & , (b_1^{(2r+1)})^r\\
\end{array}
\end{equation}
\end{scriptsize}
vanish. In the case that all such determinants vanish we can replace yet 
another $b_1^{(i)}$ by the unknown $u_1$, and so on. 
This completes the proof.
\end{proof}

\begin{proposition}
Let $f(z)$ be separately $\alpha$-rational on an open subset $\Omega\subset\Cn$. Then
$f$ is jointly $\alpha$-rational.
\end{proposition}
\begin{proof}
The function $f=f(z,\bar{z})$ has an associated function given by
$F(z,w)=f(z,w),$ which is separately a polynomial of order $\alpha_j-1$ with respect to $w_j$ whose 
coefficients are rational with respect to 
$z_j$. By Theorem \ref{hurwitzlem} $F(z,w)$ is jointly rational, and since we already know that
it is a polynomial of order $\alpha_j-1$ with respect to $w_j$
it is a polynomial with respect to $w$ of order $\alpha$ whose coefficients are rational in $z.$
Since $F(z,\bar{z})=f(z)$ this completes the proof.
\end{proof}

\begin{lemma}\label{boschlemma}
Let $\Omega\subset\C$ be a domain and let $P,Q\in \C[z,w],$ be relatively prime such that
\begin{equation}
\abs{P(z,\bar{z})}=\abs{Q(z,\bar{z})},\quad z\in \Omega
\end{equation}
There there exists a complex constant $\lambda$ such that $\abs{\lambda}=1$ and
\begin{equation}
P(z,\bar{z})=\lambda \overline{Q(z,\bar{z})},\quad z\in \Omega
\end{equation}
\end{lemma}
\begin{proof}
Define the polynomial
\begin{equation}
F(z,w):=P(z,w)\overline{P(z,w)}-Q(z,w)\overline{Q(z,w)},\quad (z,w)\in \C^2 
\end{equation}
Set $s=\frac{z+w}{2}, t= \frac{z-w}{2i}.$ Then the function $\tilde{F}(s,t)=F(s+it,s-it)$ is holomorphic
with respect to $(s,t)\in C^2$ 
and vanishes in a neighborhood of any point, $p_0,$ belonging to the complex line 
$\{ \im s=\im t=0\}.$ But if $f(\re s,\re t)=\sum_{\abs{\alpha}\geq 0} a_\alpha ((\re s,\re t)-p_0)^\alpha,$
is real-analytic then $\sum_{\abs{\alpha}\geq 0} a_\alpha ((\re s+i\im s,\re t +i\im t)-p_0)^\alpha,$
is the unique holomorphic extension.
Since $F=0$ on $\{ w=\bar{z}\}\cap \Omega\times \Omega$ we have $F\equiv 0$ on $\C^2.$
Since $P,Q$ are relatively prime, $P(z,w)$ divides $Q(z,w)\overline{Q(z,w)},$ i.e.\
there is a polynomial $R(z,w)$ such that $R(z,w)P(z,w)=\overline{Q(z,w)}.$ 
On the other hand $\overline{P},\overline{Q}$ are also relatively prime so 
there is a polynomial $R_1(z,w)$ such that $R_1(z,w)\overline{Q(z,w)}=P(z,w).$
Hence $R R_1=1$ on $\C^2,$ and since both are polynomials (i.e.\ cannot have polynomial inverse unless they are constant) 
this implies that they are both complex constants, i.e.\ 
there is a constant $\lambda$ such that
$R_1=\lambda$ on $\C^2.$ 
Hence $P(z,\bar{z})=\lambda\overline{Q(z,\bar{z})},$ for $z\in \Omega$. At any point of $\{ w=\bar{z}\}$ the equation $P(z,w)\overline{P(z,w)}=Q(z,w)\overline{Q(z,w)}$
yields $\abs{\lambda}=1.$ This completes the proof.
\end{proof}

\begin{theorem}[Bosch \& Krajkiewicz \cite{boschkraj}]\label{boschkrajthm}
Let $\Omega\subset\C$ be a domain. Let $f$ be $q$-analytic on $\Omega$ and let $g$ be $p$-analytic on $\Omega.$
A necessary and sufficient condition that $\abs{f(z)}=\abs{g(z)}$ on $\Omega$
is that there exists a non-identically zero polynomial $P(z,w)\in \C[z,w]$ such that
$P(z,\bar{z})f(z)=\overline{P(z,\bar{z})}g(z)$ on $\Omega.$
Furthermore, if $P$ satisfies the above condition then the degree of $P$ with respect to $z$ can be assumed to be $\leq (q -1)$
and the degree of $P$ with respect to $w$ can be assumed to be $\leq (p -1).$
\end{theorem}
\begin{proof}
Sufficiency is clear thus we only need to prove necessity. We can assume that $f$ and $g$
do not vanish identically on $\Omega$ and that the orders $q,p$ are exact. 
There exists holomorphic $f_k, g_j$, $k=0,\ldots,q-1,$
$j=0,\ldots, p-1,$ such that on $\Omega,$
$f(z)=\sum_{k=0}^{q-1} \bar{z}^k f_k(z),$ $g(z)=\sum_{j=0}^{p-1} \bar{z}^j g_j(z),$ where $f_{q-1}$ and $g_{p-1}$ do not vanish identically.
Hence on $\Omega,$
\begin{equation}
\sum_{k=0}^{q-1} \bar{z}^k f_k(z)\sum_{k=0}^{q-1} z^k \bar{f}_k(\bar{z})=g(z)=\sum_{j=0}^{p-1} \bar{z}^j g_j(z)\sum_{j=0}^{p-1} z^j \bar{g}_j(\bar{z})
\end{equation}
Define for $(z,w)\in X:= \{ (z,w)\colon (z,\bar{w})\in\Omega\times\Omega\}\subseteq\C^2$
\begin{equation}\label{foursums}
F(z,w)=\sum_{k=0}^{q-1} w^k f_k(z)\sum_{k=0}^{q-1} z^k \bar{f}_k(w)=g(z)=\sum_{j=0}^{p-1} w^j g_j(z)\sum_{j=0}^{p-1} z^j \bar{g}_j(w)
\end{equation}
As in the proof of Lemma \ref{boschlemma}, $F|_X =0\Rightarrow F\equiv 0.$ 
Since $f_{q-1}$ and $g_{p-1}$ do not vanish identically there is a $(z_0,w_0)\in X$ and $r>0$ such that
$\{ z: \abs{z-z_0}<r\}\times \{ w: \abs{w-z_0}<r\}\subseteq X$ and such that each of the four sums of Eqn.(\ref{foursums})
are not zero for all $(z,w)\in N_1\times N_2= \{ z: \abs{z-z_0}<r\}\times \{ w: \abs{w-z_0}<r\}.$
Hence for $(z,w)\in N_1\times N_2,$ the following function is holomorphic and nonzero
\begin{equation}
G(z,w):=\frac{\sum_{k=0}^{q-1} w^k f_k(z)}{\sum_{j=0}^{p-1} w^j g_j(z)}=\frac{\sum_{j=0}^{p-1} z^j \bar{g}_j(w)}{\sum_{k=0}^{q-1} z^k \bar{f}_k(w)}
\end{equation}
For each $z\in N_1$, $G$ is a rational function of $w$, and 
for each $w\in N_2$, $G$ is a rational function of $z$.
By Theorem \ref{hurwitzlem} a separately rational function $G(z,w)$, in the sense that $G(z,w_0)$ is a rational function with respect to $z$ for fixed $w_0$
and $G(z_0,w)$ is rational function with respect to $w$ for fixed $z_0,$ then $G(z,w)$
is a rational function of $z$ and $w.$.
Hence 
there exists two non-identically zero relatively prime polynomials $R(z,w)$ and $S(z,w)$ such that $S(z,w)G(z,w)=R(z,w)$ for all $(z,w)\in N_1\times N_2.$ This implies that for all $(z,w)\in N_1\times N_2$
\begin{equation}
S(z,w)\sum_{k=0}^{q-1} w^k f_k(z) =R(z,w)\sum_{j=0}^{p-1} w^j g_j(z)
\end{equation}
\begin{equation}
S(z,w)\sum_{j=0}^{p-1} z^j \bar{g}_j(w)=R(z,w)\sum_{k=0}^{q-1} z^k \bar{f}_k(w)
\end{equation}
Now for all $z\in N_1$ except possibly finitely many points the expressions $R(z, w)$ and $S(z, w)$ considered as
polynomials in $w$ are relatively prime 
(see e.g.\ Bochner \cite{bochneralg}, p.210) 
Hence by what we have done the degree of $R$ with respect to $z$ is $\leq (p-1)$, the degree of $S$ with respect to $z$ is $\leq (q-1)$ and similarly the degree of $R$ with respect to $w$ is $\leq (q-1)$, the degree of $S$ with respect to $w$ is $\leq (p-1)$.
Furthermore, we see that $S(z, \bar{z})f(z)= R(z,\bar{z})g(z)$, for all $z\in \Omega$. Hence $R(z,\bar{z})$ and $S(z,\bar{z})$ are
non-identically zero, relatively prime polynomials in $z$ and $\bar{z},$
with equal modulus on $\Omega$. Hence there is some constant $\lambda$ with
$\abs{\lambda}=1$ such that $R(z,\bar{z})=S(z,\bar{z})$. Writing $\lambda = \exp(-2i\theta)$ for some real $\theta$ we can set $P(z,\bar{z}) =\exp(i\theta)S(z,\bar{z}),$ for
all $z\in \C$, which is a nonidentically zero polynomial in $z$
and $\bar{z}$ such that
$P(z,\bar{z})f(z)= \overline{P(z,\bar{z})}g(z).$
Furthermore, the degree of $P$ with respect to $z$ is $\leq (p-1)$ and the degree of $P$ with respect to $\bar{z}$ is $\leq (q-1)$. This completes
the proof.
\end{proof}
As an immediate corollary we have the following earlier result.
\begin{theorem}[Balk \cite{balk1}]\label{balk1thm} A polyanalytic function $f$, of order $q$ in a domain in $\C,$ has constant modulus
if and only if there is a complex coordinate $z$ such that $f$ is representable in the form $f(z)=\lambda \overline{T(z)}/T(z)$, where $T(z)$ is a polynomial
of degree at most $q-1$, and $\lambda\in \C$ is a constant. 
\end{theorem}
We mention that these results on polyanalytic functions of constant modulus have been generalized to $\alpha$-analytic functions in several complex variables,
see Balk \& Zuev \cite{balkzuev}, p.217, and Sarkisyan \cite{sarkisyan}.
We mention also that Zuev \cite{zuev69b} proved (see e.g.\ Balk \cite{ca1}, Thm 5.3) the following generalization in the case of meta-analytic functions (see Definition \ref{metanaldef} and Section \ref{expmetanalsec}).
\begin{theorem}
Let $M$ be a polynomial with roots $a_j$ of multiplicity $m_j,$ $j=1,\ldots,q.$
Every entire $M$-meta-analytic function $f$ of constant modulus, $c$, on a domain $\Omega$
takes the form
\begin{equation}
f(z)=\lambda \frac{\overline{L(z)}}{L(z)}
\end{equation}
where $\lambda$ is a constant such that $\abs{\lambda}=c,$ and $L(z)=\sum_{j=0}^q P_j(z)\exp(\bar{a}_j z)$
where each $P_j$ is a polynomial of degree $m_j-1.$
\end{theorem}
Bosch \cite{bosch1973} gave necessary and sufficient conditions for meta-analytic functions to have equal modulus.
Gabrinovi\u{c} \cite{gabrinovic} solved a boundary value problem of Carleman type for meta-analytic functions.
One application of Theorem \ref{balk1thm} appears in the study of weak maximum modulus sets.
By a {\em weak local maximum}\index{Weak local maximum} of the absolute value, $\abs{f}$, of a function $f$, at a given point $p,$ we mean that
$\abs{f}(p)$ is a maximum value relative to an open neighborhood of $p$ but $\abs{f}(p)$ is not required to be a strict local maximum.
Recall that the modulus of
a locally open function cannot attain a local maximum. The function $f(z)=z+\bar{z}$ also illustrates the fact that nonconstant $2$-analytic functions, in contrast to nonconstant holomorphic functions, can be locally non-open at each point of an open subset of $\C$. 
Also the boundary maximum modulus principle no longer holds in general for the case $q=2$ exemplified by e.g.\ the function $f(z)=1-z\bar{z}$.
One thing, which plays a role in the presence of a boundary maximum modulus principle, is the properties of the set of local maximum modulus points
(by a local maximum modulus point for a function $f$, we mean a point, $p,$ at which $\abs{f}$ attains a weak local maximum).
If a smooth function $f$ attains weak local maximum modulus, relative to a domain, $U\subset\C$, only at isolated points (i.e.\
the peak set of $f$ relative to $U$ consists of isolated points), then each maximum of 
$\abs{f}$, relative to $U$, is in fact a strict local maximum.
In an example above we pointed out that for general $q$-analytic functions, the zero sets of higher order $q$-analytic functions can (in contrast to the case of 
holomorphic functions) be real-analytic curves. The same is true for weak maximum modulus sets of higher order $q$-analytic functions (in contrast to the case of holomorphic functions). Take e.g.\ $f(z)=1-(z+\bar{z})^2.$ This function is everywhere $3$-analytic and $\abs{f}$ attains weak local maximum value
precisely on the line $\{z\in \C \colon \re z=0 \}.$
Daghighi \cite{daghighibianal} proved the following for the cases where a $2$-analytic function attains weak local maximum at each point of a real-analytic curve in 
$\C$. 
\begin{theorem}\label{prop}
Let $\Omega\subset\C,$ be a domain and let $f(z)=a(z)+\bar{z}b(z),$ where $a,b$ are holomorphic for $z\in \Omega.$
Denote by $\Lambda$ the set of points in $\Omega$ at which $\abs{f}$ attains weak local maximum
and denote by $\Sigma$ the set of points in $\Omega$ at which $\abs{f}$ attains strict local maximum.
We prove the following: \\ 
{\em \bf (i).} At each point, $p\in \Lambda\setminus\Sigma ,$ 
\begin{equation}
\label{braek}
\abs{b(p)}=\abs{\left(\frac{\partial a}{\partial z} +\bar{z}\frac{\partial b}{\partial z}\right)(p)}
\end{equation}
{\em \bf (ii).} If there is a real analytic curve $\kappa:I\to \Lambda\setminus\Sigma$ (where $I$ is a real open interval), if $a,b$ are complex polynomials, and if $f\circ \kappa$ has a complex polynomial extension, 
then either $f$ is constant
or $\kappa$ has constant curvature.
\end{theorem}
The proof uses Theorem \ref{balk1thm}; and also the Curve Selection Lemma (Lojasiewicz \cite{loja}, Proposition 2, p.76 and Cartan \& Bruhat \cite{bruhat}, Theorem 1, see also Wallace \cite{wallace}, Lemma 18.3, p.276, the algebraic version was proved by Milnor \cite{milnor}, p.25)
Theorem \ref{balk1thm}.
Here is an example for $q=2,$ of a one-dimensional weak maximum modulus set which is a circle (and also satisfies the conditions of the second part of Theorem \ref{prop}).
\begin{example}\label{ex1}
Let $\Omega=\{ \abs{z}<1\}\subset \C,$ where $z=x+iy,$ denotes the standard complex coordinate and define,
$f(z):=z-z^2\bar{z}.$ The function $f$ is everywhere a $2$-analytic 
and it attains weak local maximum on each point of the circle $\{\abs{z}=\frac{1}{\sqrt{3}}\}.$ 
\\
It is easy to verify that 
this example satisfies the relation $\abs{b\circ \kappa}$
$=\abs{(\partial a+\bar{z}\partial b)\circ\kappa},$
for $a(z):=z,b(z):=-z^2,$ if $\kappa$ has image $\{z\colon \abs{z}=\frac{1}{\sqrt{3}}\}.$
Indeed, setting $r:=\abs{z},$ we have $\abs{b\circ \kappa}=r^2,$ and
$\abs{(\partial a+\bar{z}\partial b)\circ\kappa}=$
$\abs{\frac{3}{3}-2r^2}.$ Furthermore, (regarding the second part of Theorem \ref{prop}) we have in this example a case where
$a,b$ are both nonconstant complex polynomials and defining $\tilde{f}(z):=2z/3,$ we have
$f\circ \kappa =\tilde{f}\circ \kappa.$
\end{example}
\begin{example}
Here are all the interesting examples where the set of weak maximum modulus points contains an open subset of $\C$.
Theorem \ref{balk1thm} implies that for $q=2,$ there must exist complex constants $\alpha_1,\alpha_2,$ a complex coordinate $z,$ such that, on any open subset where
constant modulus holds true, $\alpha_1 z+\alpha_2\neq 0$, and the function is (up multiplication by a complex constant) representable as,
\begin{equation}
f(z)=\frac{\bar{\alpha_1}\bar{z}+\bar{\alpha_2} }{\alpha_1 z+\alpha_2}, \Rightarrow 
a=\frac{\bar{\alpha_2}}{\alpha_1 z+\alpha_2},\, b=\frac{\bar{\alpha_1}}{\alpha_1 z+\alpha_2}
\end{equation}
We can then verify,
Eqn.\ref{braek} using,
\begin{equation}
\frac{\partial a}{\partial z}=\frac{-\alpha_1\bar{\alpha_2}}{(\alpha_1 z+\alpha_2)^2},
\frac{\partial b}{\partial z}=\frac{-\alpha_1\bar{\alpha_1}}{(\alpha_1 z+\alpha_2)^2}
\end{equation}
and,
\begin{equation}
\frac{\partial a}{\partial z}+ \bar{z} \frac{\partial b}{\partial z}
=\frac{\alpha_1}{\alpha_1 z+\alpha_2}\left(-\frac{\bar{\alpha_1}\bar{z}+\bar{\alpha_2}}{\alpha_1 z+\alpha_2}\right)
\end{equation}
which together with,
\begin{equation}
1=\abs{-\frac{\bar{\alpha_1}\bar{z}+\bar{\alpha_2} }{\alpha_1 z+\alpha_2}}
\end{equation}
implies Eqn.(\ref{braek}).
\end{example}

Here is an example where we show for a $2$-analytic function $f=a+\bar{z}b,$ that $\abs{b}=$
$\abs{\partial a+\bar{z}\partial b},$
along the circular image of a real-analytic curve $\kappa,$ is not a sufficient condition for each point of the image of $\kappa$ to be a weak local maximum.

\begin{example}
Let $\Omega=\{ \abs{z}< 1\}\subset \C,$ where$ z=x+iy,$ denotes the standard complex coordinate and define,
$f(z)=a(z)+\bar{z}b(z), a(z):=1,$ $b(z):=z^2.$
The function $f$ is everywhere a $2$-analytic and  
$\abs{b}=r,$ $\abs{\partial a+\bar{z}\partial b}=2r^2$, where $r:=\abs{z}.$
Thus for each point of the curve $\kappa$ whose image is defined as the circle $\{\abs{z}=\frac{1}{2}\},$ we have
$\abs{b}=\abs{\partial a+\bar{z}\partial b}.$
However, along the curve $\gamma(t):=t+0i,$ we have $\abs{f\circ\gamma}^2(t)=$
$\abs{1+(z^2\bar{z})\circ\gamma}^2(t)=\abs{1+t^3}^2,$ 
which satisfies
$\frac{\partial \abs{f\circ\gamma}^2}{\partial t}|_{t=\frac{1}{2}}>0.$
This implies that the point $z(\gamma(t))=\frac{1}{2}+0i$ (which belongs to the image of $\kappa$) cannot be a weak local maximum for $\abs{f}$.
\end{example}

If $f(z)=u(z)+iv(z)$, for real-valued $u,v,$ is a $2$-analytic function, then it is easy to verify (the function is annihilated by $\overline{\partial}^2$) that, 
\begin{equation}\label{bianal}
\frac{\partial^2 u}{\partial x^2}-\frac{\partial^2 u}{\partial y^2}=2\frac{\partial^2 v}{\partial x\partial y},\quad
\frac{\partial^2 v}{\partial x^2}-\frac{\partial^2 v}{\partial y^2}=-2\frac{\partial^2 u}{\partial x\partial y}
\end{equation}
Now assume that $f$ is independent of one of the variables, say $\partial f/\partial x\equiv 0.$ By Eqn.(\ref{bianal}),
$\frac{\partial^2 u}{\partial y^2}\equiv 0,$ $\frac{\partial^2 v}{\partial y^2}\equiv 0,$ i.e.\
$u$ and $v$ must both be linear functions of $y$, and thus $\abs{f}^2$ is a non-negative polynomial in $y$ of degree at most $2$,
hence has only one critical point (unless it is constant) corresponding to a weak local minimum.
Analogously, if $f$ is independent of $y$ we can confirm that there are no weak local maxima unless $f$ is constant.
Obviously this is no longer true for $q>2,$ e.g.\ for $q=3$ we have the following simple example.

\begin{example}
Let 
$f(z):=1+(z-\bar{z})^2$ where $z,$ denotes the standard complex coordinate. Clearly, $f$ is everywhere a $3$-analytic function  
whose square modulus, $\abs{f}^2,$ 
attains the value $1$ as weak local maximum at each point of the line given by $\{ \im z=0\}.$
\end{example}

We are not aware of an example for $q=2$, where the function does not have constant modulus on any open subset, but where the set of weak maximum modulus points contains 
a curve which is not of constant curvature.

Here is an example for $q=5,$ of a one-dimensional weak maximum modulus set which does not have constant curvature.
\begin{example}
Let $z=x+iy,$ denotes the standard complex coordinate and let 
$f(z):=1-\left(2\frac{z-\bar{z}}{i}-(z+\bar{z})^2\right)^2.$ The function $f$ is everywhere a $5$-analytic 
and it attains the value $1$ on the set $\{x+iy\colon y=x^2\}.$ We can cover any (relatively) closed one-dimensional 
subset $K\subset\{x+iy\colon y=x^2\}$, by a domain $\Omega\subset\C$ such that the peak set of $f$ relative to $\Omega$ contains
$K$.
\end{example}

\subsection{Degenerate polyanalytic functions}
\begin{definition}
	An $n$-analytic function $f(z)=\sum_{j=0}^{n-1} a_j(z)\bar{z}^j$, for holomorphic $a_j(z)$
	on a domain $\Omega\subseteq\C$, is called {\em degenerate}\index{Degenerate polyanalytic transformation} if $f(\Omega)$ has no interior point.
\end{definition}
Note that the Jacobian of a polyanalytic function is countably analytic thus if a polyanalytic function is degenerate
on an open set then it is degenerate on the connected component of that set in its domain of definition.
Now it turns out that the identical vanishing of the Jacobian is precisely the requirement for degeneracy thus
a polyanalytic function $f(z)$ is degenerate if and only if 
\begin{equation}
\abs{\partial_z f}^2\equiv \abs{\partial_{\bar{z}} f}^2
\end{equation}
so known results on polyanalytic functions of equal modulus apply (see e.g.\ Theorem \ref{boschkrajthm}).
Using this Balk \cite{balkbianal64b} proves the following.
\begin{theorem}
	If $f(z)$ is a $2$-analytic (also called bianalytic) function on a simply connected domain $\Omega\subset\C$
	then
	$f(z)$ is degenerate if and only if it has one of the following representations:\\
	(i) $f(z)=A(\exp(i\alpha z)+\bar{z})+B,$ for some constants $A,B\in \C$, $\alpha\in \R$.\\
	(ii) $f(z)=A(z+c)^\gamma(\bar{z}+\bar{c})+B$ for some constants $A,B,c\in \C ,$ $\abs{\gamma}=1.$
\end{theorem}
\begin{proof}
	As we have pointed out $f(z)=u(x,y)+iv(x,y)$ is degenerate on $\Omega$ if and only if the absolute value of the Jacobian vanishes
	identically
	i.e.\
	\begin{equation}
	\abs{\frac{\partial (u,v)}{\partial (x,y)}}^2= \abs{\partial_z f}^2-\abs{\partial_{\bar{z}}f}^2\equiv 0
	\end{equation} 
	Let $f(z)=\phi(z)+\bar{z}\psi(z)$ for analytic $\phi,\psi$ on $\Omega.$ Obviously we only need to prove the result for the case
	$\psi\not\equiv 0.$ So let $a\in \Omega$ and $\rho>0$ be such that $\psi\neq 0$ on $\delta:=\{\abs{z-a}<\rho\}.$  be 
	Set 
	\begin{equation}\label{degenekv21}
	F(\zeta):=f(\zeta+a)\equiv f(z),\quad \zeta:=z-a
	\end{equation}
	We have 
	\begin{equation}\label{degenekv22}
	w=F(\zeta)=\Phi(\zeta)+\bar{\zeta}\Psi(\zeta)
	\end{equation}
	with
	\begin{equation}
	\Psi(\zeta):=\psi(\zeta +a),\quad \Phi(\zeta)=\phi(\zeta +a)+\bar{a}\psi(\zeta+a)
	\end{equation}
	and $F\neq 0$ on $\{\abs{\zeta}<\rho\}.$
	Now $F(\zeta)$ is degenerate if and only if
	\begin{equation}
	\abs{\Phi'(\zeta)+\bar{\zeta}\Psi'(\zeta)}^2-\abs{\Psi(\zeta)}^2 \equiv 0
	\end{equation}
	i.e.\
	\begin{equation}
	\abs{\frac{\Phi'(\zeta)}{\Psi(\zeta)}+\bar{\zeta}\frac{\Psi'(\zeta)}{\Psi(\zeta)}}^2\equiv 1,\mbox{ on }\delta
	\end{equation}
	By Theorem \ref{boschkrajthm} this implies that there exists constants $\lambda,q,$ with $\abs{\lambda}=1$ such that
\begin{equation}
\frac{\Phi'(\zeta)}{\Psi(\zeta)}+\bar{\zeta}\frac{\Psi'(\zeta)}{\Psi(\zeta)}=\lambda \frac{1+\bar{q}\bar{\zeta}}{1+q\zeta}
\end{equation}
Recall further that on a given domain
the analytic components of a polyanalytic function are uniquely determined thus
\begin{equation}\label{degenekv25}
\frac{\Phi'(\zeta)}{\Psi(z)}\equiv \frac{\lambda\bar{q}}{1+q\zeta},\quad \frac{\Psi'(\zeta)}{\Psi(z)}=\lambda \frac{\lambda\bar{q}}{1+q\zeta}
\end{equation}
If $q=0$ then $\Psi'(\zeta)=0$ so $\Psi(\zeta)\equiv A$, a constant, and $\Phi'(\zeta)=\lambda\Psi(\zeta)=\lambda A$ so that $\Phi(\zeta)=\lambda A\zeta +B_1$
for a constant $B_1$. This implies $F(\zeta)=A(\lambda\zeta +\bar{\zeta})+B_1$ 
	so that there is a constant $B$ such that
	\begin{equation}
	f(z)=F(z-a)=A(\exp(i\alpha z)+\bar{z})+B\end{equation}
	Hence assuming $q=0$ we arrive at (i).\\
	If instead $q\neq 0$ we set $b=1/q$ which gives
	\begin{equation}
	\frac{\Psi'(\zeta)}{\Psi(\zeta)}=\frac{\gamma}{\zeta +b},\quad \gamma:=\lambda\frac{\bar{q}}{q},\quad \abs{\gamma}=1
	\end{equation}
	\begin{equation}
	\Psi(\zeta)=A(\zeta+b)^\gamma+B
	\end{equation}
	for a constant $A$. By Eqn.(\ref{degenekv25}) we have
	for a constant $B$
	\begin{equation}
	\Phi(\zeta)=A\bar{b}(\zeta+b)^\gamma +B
	\end{equation}
	\begin{equation}
	w=\Phi(\zeta)+\bar{\zeta}\Psi(\zeta)=A(\zeta+b)^\gamma(\bar{\zeta}+\bar{b})+B
	\end{equation} 
		Reverting back to the variable $z=\zeta +a$ and denoting by $c=b-a$ we obtain (ii). Note that when the imaginary part of 
the function in (ii) is multi-valued we mean here some given branch of the function. 
Let $A\neq 0.$
In order to determine what the given connected domain (not containing $z=-c,\gamma\neq 1$) 
is transformed to under $f$, we need to consider different cases of $\gamma.$
For $\gamma=-1$ we have $f(z)=A\abs{z+c}^2 +B$ so the image coincides with a ray (or part thereof)
$f=At+B,$ $t\geq 0.$
\\
For $\gamma=-1$, we have 
$w=A(\bar{z} +\bar{c})/(z+c) +B$, $w-B=\abs{A}$, so we obtain a circle for the image in the nontrivial case.
\\
Let $\gamma:=\alpha +i\beta,$ $\beta\neq 0$, since $\abs{\gamma}^2=1$ this means $\beta^1=1-\alpha^2$
and $\alpha\in (-1,1).$
Let 
\begin{equation}
z+c=\rho\exp(i\theta),\quad \omega:=\frac{w-B}{A}
\end{equation}
This yields
\begin{equation}
\omega=\exp((\alpha+i\beta)(\ln\rho +i(0+2k\pi)))\rho\exp(-i\theta),\quad k\in \Z
\end{equation}
\begin{equation}
\ln\abs{\omega}=(\alpha +1)\ln\rho -\beta(\theta +2k\pi)
\end{equation}
\begin{equation}
\mbox{Arg} \omega=\beta\ln \rho +(\alpha-1)(\theta +2k\pi)+2k\pi+2m\pi
\end{equation}
\begin{equation}
m\in \Z\omega=\exp((\alpha+i\beta)(\ln\rho +i(0+2k\pi)))\rho\exp(-i\theta),\quad k\in \Z
\end{equation}
Fix $k$ and set $m=-k.$ If $\alpha\neq 1$ then $\beta\neq 0$ so that
\begin{equation}
(1-\alpha)\ln \abs{\omega}=(1-\alpha)\beta\left(\frac{\alpha +1}{\beta}\ln\rho -(\theta +2k\pi)\right)
\end{equation}
\begin{equation}
\beta\mbox{Arg}\omega=(1-\alpha)\beta\left(\frac{\beta}{1-\alpha}\ln\rho -(\theta +2k\pi)\right)
\end{equation}
Since $\beta^2=1-\alpha^2$ we have 
\begin{equation}
\frac{\alpha +1}{\beta} =\frac{\beta}{1-\alpha}
\end{equation}
so that
\begin{equation}
(1-\alpha)\ln\abs{\omega}=\beta\mbox{Arg}\omega
\end{equation}
\begin{equation}
\abs{\omega}=\exp\left(\frac{\beta}{1-\alpha}\mbox{Arg}\omega\right)
\end{equation}
which is the equation of a logarithmic spiral.
Since $w=A\omega +B$, the image can, via solely use of homothety, rotation and shift, be transformed to a logarithmic spiral, and
thus the image will take the form of the image of a function as in (i).
This completes the proof.
This completes the proof.
\end{proof}

\subsection{Quasi-continuity}

\begin{definition}
	Consider, for $\alpha\in \Z_+^n$ and a domain $\omega\subset\Cn$, the set of $L^1_{\mbox{loc}}(\Omega)$ functions, $f$, such that there exists
	meromorphic $a_\beta(z)$, satisfying that, on 
	$\Omega\setminus \sigma$, $f$ has
	a representation of the form
	$\sum_{j=1}^{n}\sum_{\beta_j<\alpha_j} a_\beta(z)\bar{z}^\beta$, where
	$\sigma$ denotes the union of the set of singularities of each $a_\beta$.
	We call such functions {\em $L^1_{\mbox{loc}}$ polymeromorphic}\index{$L^1_{\mbox{loc}}$ polymeromorphic functions}.
	When the $a_\beta$ are rational these are called $L^1_{\mbox{loc}}$ polyrational functions. Note that different (pointwise defined) $L^1_{\mbox{loc}}$ functions
	can have the same representation on $\Omega\setminus \sigma$.
\end{definition}
Obviously, in contrast to the case $q=1$, we have for $q>1$,
existence of bounded $L^1_{\mbox{loc}}$ functions that can be represented in the form $\sum_{j=1}^{q-1} a_j(z)\bar{z}^j$, for meromorphic $a_j$ off the singularities of the $a_j$.
Take e.g.\ the $L^1_{\mbox{loc}}$ function which has a representative, off $0$, given by
$f(z)=g(z)R(\bar{z}/z),$ where $g$ is polyanalytic and $R$ a complex polynomial. 
It is clear that such functions are not necessarily continuous. Take e.g.\ the $L^1_{\mbox{loc}}$ function which has a representative given by
$f(z)=\bar{z}/z,$ for $z\neq 0.$ Note that $\abs{f}$ can be assumed to be bounded (in fact $\equiv 1$).
The representative $f$ can for $z=x+iy$ be written
\begin{equation}
f(z)=\frac{(x^2-y^2)-i2xy}{x^2+y^2}
\end{equation}
which satisfies $\lim_z f(z)\to +1$ as $\{y=0\}\ni (x,y)\to 0$, whereas
$\lim_z f(z)\to -i\neq 1$ as $\{y=x\}\ni (x,y)\to 0$.
Quasi-continuous functions where introduced by Kempisty \cite{kempisty}.
\begin{definition}
	A function $f:X\to Y$ from a topological space $X$ to a metric space $(Y,d)$ is called {\em quasi-continuous}\index{Quasi-continuity} at $p_0\in X$ if $\forall \epsilon >0$ and each neighborhood $\mathscr{U}_{p_0}$ of $p_0$ there exists a nonempty open set $\omega_{p_0}\subset \mathscr{U}_{p_0}$ such that $d(f(x),f(p_0))<\epsilon ,\forall x\in \omega_{p_0}.$
\end{definition}
Martin \cite{martinquasi} proved the following.
\begin{theorem}\label{martinthm}
Suppose $X$ is a Baire space, $Y$ a second countable space, $Z$ a metric space.
Then any function $f:X\times Y\to Z$ which is separately quasi-continuous (in the sense that
$f_{y_0}(x):=f(x,y_0)$ is quasi-continuous for each fixed $y_0$ and $f_{x_0}(x_0,y)$ is quasi-continuous
for each fixed $x_0$) is jointly quasi-continuous.
\end{theorem}
\begin{proposition}
Any bounded $L^1_{\mbox{loc}}$ polymeromorphic function on a domain $\Omega\subset\Cn$ is a.e.\ equal to a quasi-continuous function.
\end{proposition}
\begin{proof}
By Theorem \ref{martinthm} it suffices to show that the function is a.e.\ equal to a function that is quasi-continuous with respect to each variable separately, thus it suffices to show the result for the one-variable case $n=1$.	
So suppose $q\in \Z_+$, $\Omega\subset\C$ and that there 
are $a_j(z)=b_j(z)/c_j(z)$ where $b_j,c_j$ are holomorphic, such that
$f$, takes the form $\sum_{j=0}^{q-1} a_j(z)\bar{z}^j$, off the zeros of $c_j.$
Let $p_0$ be a zero of order $m_j\in \N$ for $c_j,$ $j=0,\ldots,q-1.$ 
Since $c_j$ is holomorphic $p_0$ is an isolated zero of $c_j,$
so suppose that $c_j\neq 0$ on $U\setminus \{p_0\}$ for a domain $U\subset \Omega.$
Without loss of generality we suppose $p_0=0.$
We have holomorphic functions $h_j(z)$ on $U$ such that on $U\setminus \{0\}$
\begin{equation}
f(z)=\sum_{j=1}^{q-1} \frac{1}{\abs{z}^{2m_j}} h_j(z)\bar{z}^{j+m_j}
\end{equation}
Since 
each $h_j(z)\bar{z}^{j+m_j}$, $j=0,\ldots,q-1$ is real analytic on $U$ we know that
the function $\frac{1}{\abs{z}^{2m_j}} h_j(z)\bar{z}^j$ is real-analytic (in particular continuous) on $U\setminus\{0\}.$
Obviously, our sought representative must be pointwise defined, and since we assume boundedness of $f$, we can set (i.e.\ we redefine $f$ in $L^1_{\mbox{loc}}$ such that) for each zero $t$ of any $c_j$,  $\tilde{f}(t):=\lim_{\re z\to \re t} f(\re z+i\im t),$ and $\tilde{f}=f$ for $z\neq t$, $z\in U.$ 
In particular, for each $\epsilon>0$ there exists a $\delta>0$ such that
the preimage $\tilde{f}^{-1}\left(\left\{\abs{\zeta-\tilde{f}(0)}<\epsilon\right\}\right)$ (recall that $\tilde{f}(0)$ is, for this particular choice of representative in $L^1_{\mbox{loc}}$ defined by
$\tilde{f}(t):=\lim_{\re z\to 0} f(\re z+i0)$) contains an open ball in $\{\abs{z}<\delta\}.$ Hence $f$ is quasi-continuous
at the isolated singularity $0$. Since all singularites are isolated, $f$ is locally the sum of quasi-continuous functions, hence quasi-continuous.
\end{proof}

\section{Polyanalytic functions in the sense of \u{C}anak}
Let $\Omega\subset\C$ be a domain.
\u{C}anak \cite{canak1988} considers, for a fixed $p\in \Z_+,$ functions $f(z,\bar{z})$ 
that can be represented, on $\Omega,$ for some $q\in \Z_+,$
in the form 
\begin{equation}
f(z,\bar{z})=\sum_{j=0}^q \left( \frac{z-\bar{z}}{2i}\right)^j f_j(z,\bar{z})
\end{equation}
where each $f_j=u_j+iv_j$ satisfies
\begin{equation}
\partial_x u_j=\frac{1}{p}\partial_y v_j,\quad \partial_y u_j=-\frac{1}{p}\partial_x v_j
\end{equation}
Let us denote the set of such functions by $\mbox{PA}_{\mbox{\u{C}anak},p,q}(\Omega).$
\u{C}anak \cite{canak1992} 
defines $A_c f:=R_c (2\partial_{\bar{z}})R_c^{-1} f$ and
\begin{equation}
A_c^0 f:=R_c R_c^{-1} f=f, A_c^1 f:=A_c f,\cdots, A_c^k f=A_c(A_c^{k-1} f)
\end{equation}
where
\begin{equation}
R_c f:=\frac{(c+1)f -(c-1)\bar{f} }{2c},\quad R_c^{-1} f=\frac{(c+1)f +(c-1)\bar{f} }{2c}
\end{equation}
$R_c^{-1}R_c f=R_c R_c^{-1}f =f$ and furthermore if $g$ is analytic then $h:=R_c g\in 
\mbox{PA}_{\mbox{\u{C}anak},p}(\Omega)$ and conversely
$R_c^{-1}h$ is analytic for a function $h\in \mbox{PA}_{\mbox{\u{C}anak},p_C}(\Omega),$
i.e.\ when the characteristic is $p=c.$
Also
\begin{equation}
R_c^{-1} \left(\sum_{j=0}^q \left( \frac{z-\bar{z}}{2i}\right)^j f_{j_c}\right)
=\sum_{j=0}^q \left( \frac{z-\bar{z}}{2i}\right)^j R_c^{-1} f_{j_c}
\end{equation}
\begin{equation}
R_c \left(\sum_{j=0}^q \left( \frac{z-\bar{z}}{2i}\right)^j \varphi_{j}\right)
=\sum_{j=0}^q \left( \frac{z-\bar{z}}{2i}\right)^j R_c \varphi_{j}
\end{equation}
for holomorphic $\varphi_j.$
Hence
\begin{equation}
A_c \left(\sum_{j=0}^q \left( \frac{z-\bar{z}}{2i}\right)^j f_{j_c}\right)=
\sum_{j=1}^q j\left( \frac{z-\bar{z}}{2i}\right)^{j-1} R_c (i\varphi_{j})
\end{equation}
Denote by $P_c$ the set of $\mbox{PA}_{\mbox{\u{C}anak},p}(\Omega)$ functions with fixed characteristic $p=c$.
Then $\mathscr{O}\subset P_c\subset\C^0.$
For an analytic function $g(z)$ define $\alpha_{g(z)}:C^0\to \mathscr{O}$,
$\alpha_{g(z)} f(z,\bar{z})=f(z,g(z)).$
\u{C}anak \cite{canak1992} proves the following.
\begin{theorem}
Let $f(z,\bar{z})$ be a smooth function with bounded areolar 
derivatives on the strip $\{ \abs{\im z}\leq \delta/2 \}.$ 
Then $f$ can be approximated by members of $P_c$ in the sense that
\begin{equation}
f(z,\bar{z})-\sum_{j=0}^q \left( \frac{z-\bar{z}}{2i}\right)^j R_c \left(\frac{c\alpha_z A^j f}{j!iA^j}\right)=:\mbox{Err}_q 
\end{equation}
satisfies $\abs{\mbox{Err}_q}\to 0$ as $q\to \infty.$
\end{theorem}

\section{The countably analytic exponential}\label{expmetanalsec}
A notable property of the family of all polyanalytic functions (i.e.\ of all possible orders) is that
it is not closed under countable summation (or uniform limits).
\begin{example}
Consider the $(j+1)$-analytic functions $f_j(z):=\frac{1}{j!}\bar{z}^j.$ Then 
$f(z):=\sum_{j=1}^\infty f_j(z)=\exp(i\bar{z})$ which is countably analytic, but not $q$-analytic of any finite order.
\end{example}
Fempl \cite{fempl} defines for a constant $r$ and holomorphic $\phi(z)$, the functions
\begin{equation}
E(z,\bar{z},r):=\phi(z)\exp(r\bar{z})
\end{equation}
If $\varphi(z)$ is an analytic function then the function 
\begin{equation}
f(z):=\varphi(z)+\phi(z)\exp\left(\frac{\bar{z}}{2}\right)=\varphi(z)+E\left(z,\bar{z},\frac{1}{2}\right)
\end{equation}
satisfies that $w(z):=(f(z)-2\partial_{\bar{z}} f(z))$ is analytic.
Fempl \cite{fempl}, p.139, proved the following.
\begin{theorem}
For constants $a_1,\ldots,a_{q}$,
the general solution to the equation
\begin{equation}
2^q \partial_{\bar{z}}^q f(z) +a_1 2^{q-1} \partial_{\bar{z}}^{q-1} f(z)+\cdots +
a_12^{q-1} 2^1 \partial_{\bar{z}}^{1} f(z) +a_q f =0
\end{equation}
takes the form
\begin{equation}
f(z)=\sum_{j=1}^q E_j(z,\bar{z})=\sum_{j=0}^q \phi_j(z)\exp\left( \frac{r_j}{2}\bar{z}\right)
\end{equation}
where $\phi_j$ are holomorphic and the $r_j$, $j=1,\ldots,q,$ are the roots to the equation
\begin{equation}\label{wurzelfempl}
\xi^q+a_1\xi^{q-1}+\cdots +a_{q-1}\xi +a_q =0
\end{equation}
\end{theorem}
\begin{proof}
We have
\begin{equation}
\begin{array}{cc}
2^{1} \partial_{\bar{z}}^{1} f(z) & =\sum_{j=1}^q r_j \phi_j(z) \exp\left( \frac{r_j}{2}\bar{z}\right)\\
2^{2} \partial_{\bar{z}}^{2} f(z) & =\sum_{j=1}^q r_j^2 \phi_j(z) \exp\left( \frac{r_j}{2}\bar{z}\right)\\
\vdots & \vdots\\
2^{q} \partial_{\bar{z}}^{q} f(z) & =\sum_{j=1}^q r_j^q \phi_j(z) \exp\left( \frac{r_j}{2}\bar{z}\right)\\
\end{array}
\end{equation}
Setting $a_0:=1,$ and multiplying each equation by $a_{q-1}, a_{q-2},\ldots,a_1,a_0$ repsectively and
summing we obtain
\begin{multline}
\sum_{j=0}^{q} a_{q-j} 2^j \partial_{\bar{z}}^j f(z) =
\sum_{j=1}^q a_q \phi_j(z) \exp\left( \frac{r_j}{2}\bar{z}\right) +\\
\sum_{j=1}^q a_{q-1} r_j^1 \phi_j(z) \exp\left( \frac{r_j}{2}\bar{z}\right)+\cdots+
\sum_{j=1}^q a_{0} r_j^q \phi_j(z) \exp\left( \frac{r_j}{2}\bar{z}\right)
\end{multline}
We rewrite the right hand side according to
\begin{multline}
\phi_1 \exp\left( \frac{r_j}{2}\bar{z}\right) \sum_{j=0}^{q} a_{j}r_1^{q-j} 
+\phi_2(z) \exp\left( \frac{r_j}{2}\bar{z}\right) \sum_{j=0}^{q} a_{j} r_2^{q-j}+\\
\cdots +
\phi_q(z) \exp\left( \frac{r_j}{2}\bar{z}\right) \sum_{j=0}^{q} a_{j} r_q^{q-j}
\end{multline}
By virtue of Eqn.(\ref{wurzelfempl}) the last equation vanishes identically.
This completes the proof.
\end{proof}
Ke\u{c}ki\'c \cite{keckic}, proved that the formula remains valid  when the $a_j$ are replaced by analytic functions
$a_j(z).$
These results are important in the theory of so-called meta-anlytic functions, whose details lie outside the scope of this book.
We mention here only that they 
We mention that Zuev \cite{zuev69b} proved (see e.g.\ Balk \cite{ca1}, Thm 5.2) the following generalization in the case of meta-analytic functions (see
Definition \ref{metanaldef} and Section \ref{expmetanalsec}).
\begin{theorem}
Every entire meta-analytic function $f$ not growing faster than $\abs{z}^s$ for some constant $s\geq 0,$ can
be represented in the form
\begin{equation}
f(z)=\sum_{j=0}^q P_j(z,\bar{z})\exp(a_j\bar{z}-\bar{a}_jz)
\end{equation}
where $q\in \N,$ the $a_j$ are constants and the $P_j$ a polynomials of degree not higher than $s$ in the variables $z,\bar{z},$
$j=1,\ldots,q.$
\end{theorem}
\section{Minimal polyanalytic polynomials of normal matrices}
Huhtanen \cite{huhtanen} studies the 
\begin{definition}
Let $n\in Z_+.$
An $n\times n$ complex matrix $N$ is called {\em normal}\index{Normal matrix} if
\begin{equation}
NN^T=N^TN
\end{equation}
where $N^T$ denotes the transpose.
The set of complex $n\times n$ normal matrices is denoted $\mathcal{N}.$
\end{definition}
A unitary matrix $M$ is one satisfying $M^T=M^{-1}$, where $M^T$ denotes the transpose.
\begin{proposition}
A matrix $A$ is normal iff if is unitarily diagonalizable (i.e.\ diagonalizable by a nonsingular matrix $U$). 
\end{proposition}
\begin{proof}
Suppose $M,N$ are diagonalizable by the same (nonsingular) matrix $V,$ i.e.\
\begin{equation}
V^{-1}MV=D_M,\quad V^{-1}N V=D_N
\end{equation}
for diagonal matrices $D_M,D_M.$ Then
\begin{equation}
MN=VD_M D_BV^{-1}=VD_N D_M V^{-1}=NM
\end{equation}
Conversely, if $MN=NM$ and $v$ is an eigenvector of $A$ with associated eigenvalue $\lambda$ then
\begin{equation}
MN v=N(Mv)=\lambda Nv
\end{equation} 
thus $Nv$ is also an eigenvector associated to $\lambda,$ hence $Nv=cv$ for a scalar $c,$ which implies that
$v$ is an eigenvector of $N.$ Hence
two matrices $M,N$ commute if and only if there exists a nonsingular $V$
which diagonalizes both $M$ and $N.$
Now suppose $A$ diagonalizable by a nonsingular matrix $U$, i.e.\ $U^{-1}AU=D_A,$
where $U^{-1}=U^T,$ so $AU=UD_A.$ This implies that
\begin{equation}
AA^T=UD_A D_A^T U^T=UD_A^TD_A=A^T A
\end{equation}
so $A$ is normal. Conversely suppose $A$ is normal. Decompose
\begin{equation}
A=B+iC,\quad B=\frac{A+A^T}{2}=B^T,\quad C=-\frac{A-A^T}{2}=C^T 
\end{equation}
then $B,C$ are are Hermitian matrices both diagonalizable by a nonsingular matrix. Since $A$ is normal we have
\begin{equation}
0=AA^T-A^TA=(C+iC)(B-iC)-(B-iC)(B+iC)=2i(CB-BC)
\end{equation}
Hence $B,C$ commute so by our first observation they can be diagonalized by the same matrix $U$.
Then
\begin{equation}
U^TAU=U^T(B+iC)U=D_B+iD_C=D_A
\end{equation}
thus $A$ is unitarily diagonalizable.
This completes the proof.
\end{proof}
Huhtanen \cite{huhtanen} considers the set $\mathcal{P}\mathcal{P}_k$ of functions of the form
\begin{equation}\label{huhteq1}
p(z)=\sum_{j=0}^k h_j(z)\bar{z}^j,\quad h_j\in \mathcal{P}_{k-j},k\in \N
\end{equation}
where $\mathcal{P_k}$ denotes the set of polynomials of degree at most $k.$
If for some $j$ we have deg$(h_j)=k-j$ for $p(z)$ given by Eqn.(\ref{huhteq1})
then the degree of $p$ is defined as $k.$
\index{$\mathcal{P}\mathcal{P}_k$}
Obviously, $\mathcal{P}\mathcal{P}_k$ consists of a special proper subset of what we usually consider to be polyanalytic polynomials,
since the analytic components have to fulfill additional conditions. Functions of the form $z^j\bar{z}^k$ are called polyanalytic monomials.
Set $\mathcal{P}\mathcal{P} =
\cup_{k\in \N} \mathcal{P}\mathcal{P}_k.$
For a normal matrix $N$ the function $p(N)$ is well-defined for any $p\in \mathcal{P}\mathcal{P}$
by identifying $z$ and $\bar{z}$ with $N$ and $N^T$ respectively.
If $j_1 + l_1 > j_2 + l_2$, then $z^{j_1}\bar{z}^{l_1} > z^{j_2}\bar{z}^{l_2}$ and also if
$j_1 + l_1 = j_2 + l_2$ and $j_1 > j_2,$ then again $z^{j_1}\bar{z}^{l_1} > z^{j_2}\bar{z}^{l_2}$. Using this order, 
one can define the {\em minimal polyanalytic polynomial}\index{Minimal polynomial of a Normal matrix} of $N \in \mathcal{N}$ (the set of normal matrices) 
to be the monic $p\in \mathcal{P}\mathcal{P}$, i.e., its leading
term is $z^j\bar{z}^l$, of the least possible degree annihilating $N$. 
Huhtanen \cite{huhtanen} proves the following.
\begin{proposition} Every $N \in\mathcal{N}$ possesses a unique minimal polyanalytic polynomial.
Moreover, unitarily similar $N_1,N_2\in \mathcal{N}$ have the same minimal polyanalytic polynomial.
\end{proposition}
\begin{proof}
Let $p(z) = z^j\bar{z}^l +\mbox{lower terms}$ and $q(z) = z^j\bar{z}^l +\mbox{lower terms}$ be two
minimal polyanalytic polynomials of $N  \in \mathcal{N}$. Then 
$p-q$ yields either a zero polyanalytic polynomial or a smaller nonzero polyanalytic polynomial than $p$
or $q$. In the latter case this $p-q$ would yield an annihilating polyanalytic polynomial
of $N$, which contradicts minimality of $p$ and $q$. This proves the first assertion.
If $p_2$ is the minimal polyanalytic polynomial of $N_2$ and $U$ is a unitary matrix such that $N_1 = U N_2 U^T$ then
$p_2(N_1)=Up_2(N_2)U^T=0$,
thus $p_2$ is also the minimal polyanalytic polynomial of $N_1$. And similarly replacing the roles of
$N_1$ and $N_2$ we obtain the wanted result. This completes the proof.
\end{proof}
The following is clear:
Let $p_{j,l}$ be the minimal polyanalytic polynomial of $N\in \mathcal{N}$.
Then $q_{j,l}:=\overline{p_{j,l}(\bar{z})}$ is the minimal polyanalytic polynomial of $N^T$.
\begin{proposition}
If $N\in \mathcal{N}$ and $p\in \mathcal{P}\mathcal{P},$ then $\sigma(N)\subset
\{ z\in \C:\abs{p(z)}\leq \norm{p(N)}\},$ where $\sigma(N)$ denotes the spectrum of the square matrix $N.$ 
\end{proposition}
\begin{proof}
Since $N$ is normal we have $N^T=q(N)$ for a polynomial $q.$ Inserting this into $p(N)=p(N,N^T)=p(N,q(N))$
gives $\norm{p(N,N^T)}=\max_{\lambda\in \sigma(N)} \abs{p(\lambda,q(\lambda)}.$ This completes the proof.
\end{proof}
\begin{corollary}
Let $p_{j,l}$ be the minimal polyanalytic polynomial of $N\in \mathcal{N}$.
Then the eigenvalues of $N$ are contained in the zero set of $p_{j,l}.$
\end{corollary}

\section{Monogenic functions in the sense of Fedorov}
Monogeneity occurs in different fields of mathematics. 
In algebra a field is called monogenic if its ring of integers is generated by a single element.
In classical complex analysis monogeneity at a point means the existence of the complex differential
\begin{definition}
	A complex valued function $f(z)$ on a domain $\Omega\subset\C$ is called {\em monogenic
		in the classical sense}\index{Monogenic function in the classical sense} if the following limit exists at each point $z\in \Omega$
	\begin{equation}
	\lim_{h\to 0}\frac{f(z+h)-f(z)}{h}
	\end{equation}
\end{definition}
Cauchy proved that if the complex differential of a continuous function 
is continuous at each point of a domain then the function is holomorphic. Goursat \cite{goursat1900} showed that
the continuity condition on the differential can be dropped.
Bohr \cite{bohr1918} proved the following.
\begin{theorem}
	If $w=f(z)$ is a continuous univalent mapping of a domain $\Omega\subset\C$ such that 
	the complex differential is nonzero a.e.\ on $\Omega$ then either $f(z)$ or $\overline{f(z)}$ is holomorphic on $\omega.$
\end{theorem}
This result was generalized by 
Men\u{n}ov \cite{mensov1937}, \label{mensov1936}.
The condition on univalency is necessary, as can be realized by the example
$f(z):=z,$ $\im z\geq 0,$ and $f(z):=\bar{z},$ $\im z<0.$ 
A point $p_0$ is called a $U$-point\index{$U$-point of a continous map} of a continuous map $f(z)$ if there exists two sequences $\{z'_j\}_{j\in \Z_+},$ $\{z''_j\}_{j\in \Z_+}$ 
both converging to $p_0$ such that the semitangents $t',t''$ 
at $p_0$ lie on different straight lines and all the points $w_j'=f(z_j'),$ $w_j''=f(z_j'')$ are different from $w_0=f(p_0).$
The map $w=f(z)$ is called orientation preserving\index{Orientation preserving preserving map at $U$-point}
at the $U$-point $p_0$ if the sequences 
satisfy that if $\{w'_j\}_{j\in \Z_+},$ $\{w''_j\}_{j\in \Z_+}$ have semitangents 
$T',T''$ at $w_0$ such that if
the smallest angle between $t',t''$ lies in $(0,\pi)$ then the smallest angle between $T',T''$ lies in $(0,\pi)$
Trokhimchuk \cite{trokhimchuk} proved the following
\begin{theorem}
	Let $w=f(z)$ be a continuous mapping of a domain $\Omega\subset\C$ 
	with constant stretching a.e.\ on $\Omega$ such that $f(z)$ is orientation preserving at a.e.\ $U$-point.
	Then $f(z)$ is analytic on $\Omega$. If there are no $U$-points of $f$ in $\Omega$ then $f(z)$ is a constant.
	such that 
	the complex differential is nonzero a.e.\ on $\Omega$ then either $f(z)$ or $\overline{f(z)}$ is holomorphic on $\omega.$
	\end{theorem}	
	Petrov \cite{petrov1969} investigated sets of monogeneity of a polyanalytic functions and proved the following.
	\begin{theorem}
		The set of monogeneity of a polyanalytic, non-holomorphic function, $f(z)$ in a domain
		$\Omega\subset\C$ belongs almost everywhere to a circle.
		A necessary and sufficient condition for a set of monogeneity of a polyanalytic
		function $f(z)$ does not depend on the points of $\Omega$ is that $f(z) = c_0 z +c_1\bar{z}$ for constants $c_0,c_1\in \C.$
		In this case, the set of monogeneity is a circle with radius $\abs{c_1}$ and center $c_0.$
	\end{theorem}
We have already seen in Chapter \ref{hypercomplexsec} how $k$-monogenic functions are define in Clifford analysis.
Let $(A,\odot,+)$ be a $\C$-algebra with basis $\{e_0,\ldots,e_m\}.$ In particular, there is
a multiplication table with respect to $\odot$ which expresses each $e_j\odot e_k$ as a linear combination
of the $e_0,\ldots,e_m$ and thus extends to any $a\odot b$ for any $a,b\in A.$ We assume $A$ is commutative, associative
and that multiplication is distributive over addition.
Fedorov \cite{fedorov1958}, \cite{fedorov1960}, introduced the following functions.
\begin{definition}
	Let $\Omega\subset\C$ be a domain.
	A function $f(z)=\sum_{j=0}^m e_j f_j(z),$ $z\in \Omega,$ for complex valued $f_j$ is called {\em monogenic in the sense of Fedorov} with respect to the function
	$t(z)=\sum_{j=0}^m e_j t_j(z),$ for complex valued $t_j$,
	if on $\Omega$ there exists a function $\phi(z)=\sum_{j=0}^m e_j \phi_j(z),$
	such that
	\begin{equation}
	df(z)=\phi(z)dt(z),\quad \mbox{ on }\Omega
	\end{equation}
	which can be written with $z=x+iy$ as 
	\begin{equation}
	\partial_x f(x,y)=\phi(x,y)\partial_x t(x,y),\quad \partial_y f(x,y)=\phi(x,y)\partial_y t(x,y)
	\end{equation}
	The function $\phi$ is called the {\em Fedorov derivative}\index{Fedorov derivative} of $f$ with respect to $t(z)$, and is sometimes denoted $\frac{df}{dt}.$
	The functions $f_j(z),$ are called the $j$:th components of $f,$ $j=0,\ldots,m.$
\end{definition}
We consider the case of $A$ having basis $1,e,\ldots,e^{q-1}$, for an element $e\in A,$ and
the choice of $\odot$ satisfying $e^j\odot e^k=e^{j+k},$ $e^q=0.$
Fedorov points out that any $q$-analytic function $g$
can be regarded as a component of a function $f$ that is monogenic in the sense of Fedorov
with respect to the variable $t=z+e\bar{z}.$
Balk \& Zuev \cite{balkzuev}, p.220, state that Zatulovskaya \cite{zatulovskaya} proved the following. 
\begin{theorem}
	A function $f(z)$ is monogenic in the sense of Fedorov, in the variable $t=z+\bar{z}$ on a domain 
	$\Omega\subset\C$ iff there exists holomorphic functions $a_j(z),$ $j=0,\ldots,q-1$ such that
	all components $f_j(z)$ of $f$ are polyanalytic on $\Omega$
	and can be expressed as
	\begin{equation}
	f_{j}=\sum_{l=0}^{j} \frac{1}{(q-l-1)!} \frac{d^{q-l-1} h_l(z)}{dz^{q-l-1}} \bar{z}^{q-l-1},\quad j=0,\ldots,q-1
		\end{equation}
		For each $q$-analytic $g_{q-1}(z)$ on $\Omega$ there exist functions $g_0,\ldots,g_{q-2}$ such that the function
		$g(z)=\sum_{j=0}^{q-1} e^j g_j$ is monogenic in the sense of Fedorov on $\Omega$ in $t=z+e\bar{z}.$ 
	\end{theorem}
	\begin{theorem}
		If $f_{q-1}$ is a $q$-analytic function on a domain $\Omega\subset\C$ bounded by a rectifiable Jordan curve and if
		$f_0,\ldots,f_{q-2}$ are the associated polyanalytic functions such that 
		$f(z)=\sum_{j=0}^{q-1} e^j f_j$ is monogenic in the sense of Fedorov on $\Omega$ in $t=z+e\bar{z},$
		then 
		\begin{multline}
		\frac{1}{2\pi i} \int_{\partial\Omega} \frac{f_{q-1}(t)dt}{t-z}+\sum_{j=0}^{q-1}
	\int_{\partial\Omega} f_{q-j-1}\frac{(-1)^{j-1}}{j}\frac{d}{dt}\left(
	\frac{(\bar{t}-\bar{z})}{t-z}\right)^jdt=\\
	\left\{
	\begin{array}{ll}
	f_{q-1}(z) & ,z\in \Omega\\
	0 & , z\notin \Omega
	\end{array}
	\right.
	\end{multline}
	\end{theorem}
\begin{appendix} 
\small

\chapter{Some relevant properties of elliptic operators}\label{elliptch}
\section{Elliptic regularity}\label{ellipticapp}
\begin{definition}
Denote by $\mathcal{S}(\Rn)$ the set of $C^\infty$-smooth functions
$f:\Rn\to \C,$ such that
\begin{equation}
x^\alpha D^\beta f(x)\to 0,\mbox{ as }\abs{x}\to \infty
\end{equation}
for every pair of multi-indices $\alpha,\beta\in N^n,$
where $D^\alpha:=\partial_{x_1}^{\alpha_1}\cdots\partial_{x_n}^{\alpha_n}.$
Define the semi-norm estimates
\begin{equation}
\norm{f}_{\alpha,\beta}:=\sup_{\Rn}\abs{x^\alpha D^\beta f(x)}
\end{equation}
A sequence $f_j$ is said to converge in $\mathcal{S}$ to $f$
if $\norm{f_j-f}_{\alpha,\beta}\to 0$ as $j\to \infty,$
for every pair of multi-indices $\alpha,\beta\in N^n.$
The topology in $\mathcal{S}$ is defined via a metric
\begin{equation}
d(f,g):=\sum_{\alpha,\beta\in \N^n}\frac{c_{\alpha,\beta}\norm{f-g}_{\alpha,\beta}}{1+\norm{f-g}_{\alpha,\beta}}
\end{equation}
for any choice of $c_{\alpha,\beta}$ such that $\sum_{\alpha,\beta\in \N^n}c_{\alpha,\beta}<\infty.$
And it can be verified that $\mathcal{S}$ is complete with respect to such a metric.

It is sometimes more convenient to replace the family of semi-norms with $\{p_k\}_{k\in \N},$
where \begin{equation}
p_k(u)=\sum_{\abs{\alpha}\leq k}\sup_{x\in \Rn}(1+\abs{x}^2)^{k/2}\abs{D^\alpha u(x)}\end{equation}
and the distance function then becomes
\begin{equation}
d(f,g):=\sum_{k=0}^\infty 2^{-k}\frac{p_k(f-g)}{1+p_k(f-g)}
\end{equation}

Denote by $\mathcal{D}(\Rn)$ the set of members of $\mathcal{S}(\Rn)$ which have compact support.
Note that if $f\in \mathcal{D}$ and $\phi\in C_c^\infty(\Rn)$ then
$\phi(x/j)f\to f$ as $j\to \infty$, in $\mathcal{S},$ so $\mathcal{D}$ is dense in $\mathcal{S}.$
\end{definition}
$D^\alpha$ is a continuous map on $\mathcal{S}$ since 
$f_j\to f$ in $\mathcal{D}$ implies $D^\alpha f_j\to D^\alpha f$. 
\begin{definition}
A {\em tempered distribution} $u$ on $\Rn$ is a continuous linear functional 
$\mathcal{S}(\Rn)\to \C$, its value is denoted $\langle u,f\rangle$ 
such that there are constants $c,N\geq 0$ 
satisfying $\abs{\langle u,\phi\rangle}\leq c\sum_{\abs{\alpha},\abs{\beta}\leq N} \sup\abs{x^\alpha D^\beta \phi},$ $\phi\in \mathcal{S}(\Rn).$
The set of tempered
distributions is denoted $\mathcal{S}'.$
\index{Tempered distribution}
For $m\in \R$ define $\mathcal{S}_1^m(\Rn)$ as the set of functions $f\in C^\infty(\Rn)$
such that 
\begin{equation}
\abs{D_x^\alpha f(x)}\leq C_\alpha (1+\abs{x}^2)^{(m-\abs{\alpha})/2}, \mbox{ for all }{\alpha}\geq 0
\end{equation}
where $D^\alpha:=D_1^{\alpha_1}\cdots D_n^{\alpha_n},$ $D_j:=-i\partial/\partial x_j.$
\end{definition}

\begin{definition}
The {\em Fourier transform}\index{Fourier transform} $\mathcal{F}:L^1\to L^\infty(\Rn)$ is defined as the map
\begin{equation}f\mapsto \hat{f}(\xi):=(2\pi)^{-\frac{n}{2}}\int f(x)\exp(-ix\cdot \xi)dx\end{equation}
\end{definition}
One verifies directly that the Fourier transform has a restriction satisfying
\begin{equation}\mathcal{F}:\mathcal{S}(\Rn)\to \mathcal{S}(\Rn)\end{equation}
and $\xi^\alpha D_\xi^\beta \mathcal{F} f(\xi)=(-1)^{\abs{\beta}}\mathcal{F}(D^\alpha x^\beta f)(\xi).$
The inversion formula is defined by $\mathcal{F}^* f(\xi)=(2\pi)^{-\frac{n}{2}}\int f(x)\exp(i\cdot \xi)dx$
(the term inversion is used because $\mathcal{F}\mathcal{F}^*=\mathcal{F}^*\mathcal{F}=I$ on $\mathcal{S}(\Rn).$
It is also clear that if $p\in \mathcal{S}'(\Rn)$  is homogeneous of degree $k$ (in fact this extends to the case when $k$ is allowed to be complex), i.e.\
$\langle f,D(t)p\rangle=t^{-k}\langle D(t^{-1}f,p\rangle,$ or simpler $D(t)p=t^k p,$ then 
$\mathcal{F}D(t)=t^{-k} D(t^{-1})\mathcal{F}.$
\begin{proposition}\label{ellipticlemma1}
If $p\in \mathcal{S}_1^k(\Rn)$ then $\hat{p}\in C^\infty (\Rn\setminus \{0\}).$ Furthermore, if $a>0$ and 
$\phi\in C_c^\infty(\Rn)$
such that $\phi(x)=1$ for $\abs{x}<a$ then $(1-\phi)\hat{p}\in \mathcal{S}(\Rn).$
\end{proposition}
\begin{proof}
We have $p\in \mathcal{S}^m_1(\Rn)\Rightarrow D^\beta p\in \mathcal{S}^{m-\abs{\beta}}_1(\Rn)$ and for $\mu <-n,$
\begin{equation}
\mathcal{F}\colon S_1^\mu(\Rn)\to L^\infty(\Rn)\cap C^0(\Rn)
\end{equation}
thus
\begin{equation}
x^\beta\hat{p}=\mathcal{F}(D^\beta p)\in L^\infty\cap C^0,\mbox{ for }\abs{\beta}>m+n
\end{equation}
and more generally $x^\alpha D^\beta p\in \mathcal{S}^{m-\abs{\beta}+\abs{\alpha}}_1(\Rn)$ implies
\begin{equation}
D^\alpha (x^\beta\hat{p})=\mathcal{F}(x^\alpha D^\beta p)\in L^\infty\cap C^0,\mbox{ for }\abs{\beta}>m+n+\abs{\alpha}
\end{equation}
This completes the proof.
\end{proof}

The space $H^s(\Rn)\subset\mathcal{S}'(\Rn)$ is defined as the set of $f\in \mathcal{S}'(\Rn)$
such that the Fourier transform $\hat{f}:=\mathcal{F}(f)$ satisfies
\begin{equation}
\int (1+\abs{\xi}^2)^{1/2} \abs{\hat{f}(\xi)}^2 d\xi <\infty
\end{equation}
It is normed using the inner product
\begin{equation}
\langle f,g \rangle = \int (1+\abs{\xi}^2)^{s} \overline{\hat{f}(\xi)} \hat{g}(\xi) d\xi
\end{equation}

By convention one defines $H^\infty =\cap_{s\in \R} H^s,$ $H^{-\infty} =\cup_{s\in \R} H^s.$
Then $\mathcal{S}\subset H^\infty \subset H^{-\infty}\subset \mathcal{S}'.$

Let $\Omega\subset\Rn$ be an open subset.
The space $C^\infty(\Omega)$ equipped with the family of semi-norms $\norm{f}_{K,\alpha}:=\sup_{x\in K}
\abs{D^{\alpha} f(x)},$ defined for all $K\subset\subset \Omega$ and $\alpha\in \N^n,$ is denoted $\mathcal{E}(\Omega).$
Its dual $\mathcal{E}'(\Omega)\subset\mathcal{D}'(\Omega)$ can be identified as the subspace consisting of
distributions with compact support.

\begin{theorem}\label{regularizationthm}
If $u\in \mathcal{D}'(\Rn)$ and $f\in C_c^\infty(\Rn)$ then $u*f\in C^\infty(\Rn)$
and \begin{equation}\label{suppekvhorm}
\mbox{supp}(u*f)\subset \mbox{supp}(u)+\mbox{supp}(f)
\end{equation}
\end{theorem}
\begin{proof}
A quick way way to realize that $u*f\in C^\infty(\Rn)$ is to note that
$f$ having compact support implies that $u*f$ is a distribution with compact support, and the latter can be identified 
with the dual of $C^\infty(\Rn).$ Indeed, if $\phi\in C^\infty(\Rn)$ and $v$ is a distribution with compact support $K$, then
then let $\psi\in C^\infty_c(\Rn)$ such that $\psi=1$ on $K.$ This means that we can define 
$\langle v,\phi\rangle=\langle v,\psi\phi\rangle +\langle v,(1-\psi)\phi\rangle=$
$\langle v,\psi\phi\rangle$ and we have the semi-norm estimates 
$\abs{\langle v,\phi\rangle}\leq C_{\abs{\alpha}\leq k} \sup_{K'} \abs{\partial^\alpha\phi}$
where $K'$ is the support of $\psi$ and $k,K'$ are constants. Conversely if $v$ is a linear form
on $C^\infty(\Rn)$ that for some compact $K''$ and constants $C',k'$ satisfies the semi-norm estimates
$\abs{\langle v,\phi\rangle}\leq C'_{\abs{\alpha}\leq k'} \sup_{K''} \abs{\partial^\alpha\phi}$
for all $\phi\in C^\infty(\Rn),$ then the restriction, $v',$ of $v$ to $C^\infty_c(\Rn)$ is a distribution with support in the compact $K''.$
Since $\langle v,\phi\rangle =0$ if $K''\cap\mbox{supp}\phi=\emptyset$ we have
$\langle v,\phi\rangle=\langle v',\phi\rangle$ for all $\phi\in C^\infty(\Rn).$
We can however also look explicitly at the definition of derivatives 
\begin{multline}
\frac{(u*f)(x+h)-(u*f)(x)}{h}=\frac{\langle u(y),f(x+h -y)\rangle-\langle u(y),f(x -y)\rangle}{h}=\\
\langle u(y),\frac{f(x+h)-\phi(x)}{h}\rangle
\end{multline}
and since, in the topology of $\mathcal{D},$ for fixed $y$, 
\begin{equation}
\frac{\phi(x+h-y)-\phi(x-y)}{h}\overset{h\to 0}{\longrightarrow} \partial f (x-y)
\end{equation}
we obtain
\begin{equation}
\frac{(u*f)(x+h)-(u*f)(x)}{h}\overset{h\to 0}{\longrightarrow}
\langle u(y),\partial f(x-y) \rangle =\partial (u*f)
\end{equation}

For the second statement of the theorem note that $u*f(x)=0$ unless $x-y\in\mbox{supp} f$ for some $y\in\mbox{supp}u.$ This implies that   
$x\in\mbox{supp}u +\mbox{supp}f.$
This proves Eqn.(\ref{suppekvhorm}).
To prove smoothness note that 
\end{proof}

\begin{definition}
Let $u\in \mathcal{D}'(\Rn).$ 
$u$ is called smooth on an open subset $\Omega\subseteq\Rn$
if there is a $v\in C^\infty(\Omega)$ such that $u=v$ on $\Omega.$
The {\em singular support}\index{Singular support}
of $u$, denoted $\mbox{sing supp} u$ is defined as the complement of the set
$\{x\in \Omega\colon u\mbox{ is smooth on an open neighborhood of }x\}.$
\end{definition}
\begin{theorem}[Elliptic regularity theorem]\label{elliptictheorem}
Let $P(D)=\sum_{\abs{\alpha}\leq m} a_\alpha D^\alpha$
be an operator such that there is a constant $C>0$ satisfying 
$\abs{P(\xi)}\geq C\abs{\xi}^m$ for all sufficiently large $\abs{\xi}$ (this is equivalent to the standard definition of ellipticity of $P(D)$).
If $P(D)u=f$ for some $u,f\in \mathcal{D}'(\Rn)$ then $\mbox{sing supp} u =\mbox{sing supp}f.$
\end{theorem}
\begin{proof}
The inclusion $\mbox{sing supp}f\subset \mbox{sing supp} u$ follows from the property of 
the differential operator $P(D)$, i.e.\ if $u$ is smooth near a point $p_0\in \Rn$,
then of course $P(D)u$ is also smooth near $p_0$ thus $p_0\notin \mbox{sing supp}f.$
We must prove $\mbox{sing supp}f\supset \mbox{sing supp} u$
If $\phi\in C_c^\infty(\Rn)$ satisfies $\phi=1$ for $\abs{\xi}\leq C,$ then 
\begin{equation}
q(\xi):=(1-\phi(\xi))P(\xi)^{-1}\in \mathcal{S}_1^{-m}(\Rn)
\end{equation}
Consider $E:=(2\pi)^{-\frac{n}{2}}\mathcal{F}^* q\in \mathcal{S}'(\Rn).$
By Proposition \ref{ellipticlemma1} $E \in C^\infty (\Rn\setminus \{0\})$ and rapidly decreasing as $\abs{x}\to \infty.$
Set $v:=P(D)E.$ Then $\hat{v}(\xi)=(2\pi)^{-\frac{n}{2}}(1-\phi(\xi)),$ that is 
\begin{equation}
P(D)E=\delta +w,\quad w:=-(2\pi)^{-\frac{n}{2}}\mathcal{F}^*\phi\in \mathcal{S}(\Rn)
\end{equation}
Now $\mbox{sing supp} \delta=\{0\}$ and $\mbox{sing supp} \abs{x}^{2-n}=\{0\}$ for $n\neq 2.$ Furthermore,
if $u$ is a member of $\mathcal{E}'(\Rn)$ satisfying $P(D)u=f$
then $E*f=E*P(D)u=(P(D)E)*u=u+w*u,$ where $w*u\in C^\infty(\Rn).$
On the other hand, for any $h\in \mathcal{E}'(\Rn)$, 
$\mbox{sing supp} E*h\subset \mbox{sing supp} h,$ as soon as $\mbox{sing supp} E\subset \{0\}.$
Since we can multiply the distributions with cut-off functions $\chi\in C^\infty_c(\Rn)$ that are equal to $1$
on an arbitrary large set, we conclude that
that $\mbox{sing supp}(u+w*u)\subset \mbox{sing supp} f.$
By Theorem \ref{regularizationthm} $w*u\in C^\infty(\Rn)$ we have $\mbox{sing supp} (u+w*u)=\mbox{sing supp}u,$ i.e.\
$\mbox{sing supp}u\subset \mbox{sing supp} f.$ This completes the proof.
\end{proof}
An alternative proof can be given using the notion of parametrix and Fredholm theory, see Section \ref{fredholmsec}.

\section{Elliptic analyticity and Friedrich's inequalities}\label{ellipticapp2}
There exists an analytic version of Elliptic regularity which we shall state without proof.
First we gather some definitions, and here we include some definitions that are required in the proof of the Malgrange-Lax approximation theorem.
\begin{definition}
Let $X$ and $E$ be Hausdorff spaces and $\pi:E\to X$ a continuous map. The triplet $(E,\pi,X)$ is called a continuous 
complex (real) vector bundle\index{Vector bundle} of rank $m$
if:\\
(i) For all $x\in X$, $E_x:=\pi^{-1}(x)$ is a vector space of dimension $m$ over $\C$ ($\R$).\\
(ii) For each $x\in X$ there exists a neighborhood $U$ of $x$ and a homeomorphism $h$ of $E_U:=\pi^{-1}(U)$ onto $U\times \C^m$ $(U\times \R^m)$
such that if $\pi$ is the projection of $U\times \C^m$ $(U\times \R^m)$ onto $U$ we have
$\pi(h(y))=x$ whenever $y\in E_x$ and $h|_{E_x}$ is a $\C$-isomorphism ($\R$-isomorphism) of $E_x$ onto $\{x\}\times \C^m$ $(\{x\}\times \R^m$).
If $E$ and $X$ are $C^k$-smooth manifolds and $\pi$ a $C^k$-smooth map and if the
isomorphisms $h|_U$ can be chosen to be a $C^k$-smooth diffeomorphisms, $\pi: E \to X$
is called a $C^k$-smooth bundle.
If $X$ is a real (complex) analytic, a manifold real (complex) analytic
vector bundle can be defined in the same way. Complex analytic bundles
are also called holomorphic vector bundles.
Let $\pi_1:E\to X$, $\pi_2:F\to X$ be two complex vector bundles on $X$. A bundle map or a homomorphism $h:E\to F$ is a continuous hap $h:E\to F$
such that for any $x\in X$,
$h|_{E_x}$ is a $\C$-linear map into $F_x.$ If in addition $h$ is a homeomorphism (so that
$h|_{E_x}$ is an isomorphism onto $F_x$), $h$ is called an isomorphism.
Let $X$ be a $C^k$-smooth manifold and $\pi : E \to X$ a $C^k$-smooth vector bundle 
and $U\subset X$ an open subset. Then a $C^k$-smooth {\em section} $s$ of $E$ on $U$ is a $C^k$-smooth map
$s : U\to  E$ such that $\pi\circ s =\mbox{Id}$ on $U$. 
The set of $C^k$-smooth sections of $E$ over $U$ is denoted $\Gamma^k(U,E).$
The support of a section $s$ defined as the closure of $\{x\colon x\in U,s(x)\neq 0\}$
where $0$ denotes the zero of the vector space $E_x.$ To denote that we have sections of compact support we use the subindex $c$ (for compact),
$\Gamma_c^k(U,E).$
Given a $C^\infty$-smooth manifold $X$ and vector bundles $\pi_1 : E \to X$ and
$\pi_2 : F\to X$, a differential operator $L$ from $E$ to $F$ (written $L:E\to F$)
is an $\R$-linear map $L : \Gamma^\infty_c(X,E)\to\Gamma^0(X, F)$ such that supp$(Ls)\subset$
supp$(s)$ for every $s\in \Gamma_c^\infty(X,E).$
Obviously $L$ gives rise to a linear map as follows: let $x\in X$ and $U$ a relatively compact neighborhood of $x$. Let $\phi:X\to \R$
belong to $C^\infty_c$ and be such that for $\phi=1$ on a neigbhorhood of $x$, and $\phi=0$ outside $U.$ The for any $s\in \Gamma^\infty(X,E)$ set $(Ls)(x):=L(\phi s)(x),$
which is well-defined since $\phi$ has compact support, furthermore $(Ls)(x)$ is independent of $\phi$ because $L$ does not increase supports.
If $E,F$ are $C^\infty$-smooth vector bundles over $X$ of rank $q,r$ respectively, and if,
$U$ is a coordinate neighbourhood of $X$ such that $E_U$ and $F_U$ are trivial
then $\Gamma^\infty_c(U,E)$ can be identified with $C_c^{\infty,q}(U)$ (by which we denote
the set of $q$-tuples of $C^\infty$-smooth functions with compact support in $U$). A linear differential operator
$L$ then defines then an $\R$ linear map $L :C_c^{\infty,q}(U)\to C^{0,r}(U)$. 
The norm on $C^{k,q}(U)$ is denoted by $\norm{f}_m=\sum_{\abs{\alpha}\leq m}\sum_{i=1}^q \sup \abs{\partial^\alpha f_i(x)},$
where $f=(f_1,\ldots,f_q).$ Let $\Omega\subset\Rn$ be an open subset.
Define for $f\in C^\infty(\Omega,\C^q)$, $\abs{f}_{m,l}^\Omega$ via $(\abs{f}_{m,l}^\Omega)^l=\sum_{\abs{\alpha}\leq m}\sum_{i=1}^q \int_\Omega \abs{\partial^\alpha f_i(x)}^l dx ,$
and when the domain $\Omega$ is not relevant we drop the superscript $\Omega.$
If $L:C^{\infty,q}_c(\Omega)\to C^{\infty,r}_c(\Omega)$ is a linear differential operator
then it can be written $Lf=\sum_{\abs{\alpha}\leq m}\sum_{i=1}^q a_\alpha \partial^\alpha f (x)$
for continuous maps $a_\alpha$ of $\Omega$ into the space of $r\times q$ matrices and if
there exists $\alpha$ such that $\abs{\alpha}=m$ and $a_\alpha\not\equiv 0$ on $\Omega$ then
$L$ is said to have order $m$ on $\Omega.$ If for any $\xi\in \Rn\setminus \{0\}$ the map 
$\sum_{\abs{\alpha}=m} \xi^\alpha a_\alpha(x)$ is injective $\C^q\to \C^r$, then $L$ is called elliptic. As before $H^m:=H^{m,2}$ where $H^{m,p}$ 
denote
the Sobolev spaces which are defined as the completion with respect to (Sobolev norm) $\abs{f}_{m,p}^\Omega$
in the space of $C^\infty$-smooth functions on $\Omega$ with finite Sobolev norm.
Let $K$ be a compact set in $V$. We denote by $H^m(K, E)$ the set of
sections $s:V\to E$, which are locally in $H^m$ with supp$s\subset K$. 
Let $U_1',\ldots,U_m'$ be coordinate neighborhoods such that $E|_{U_j'}$ is trivial and let $K\subset \cup U_j,$ Let $U_j\Subset U_j'.$
Let $\tau_j:E|_{U_j'}\to U_j'\times \C^r$ be a $C^\infty$-smooth isomorphism and let $\phi_j\in \C^\infty_c(U_j),$
$\sum \phi_j=1$ on a neighborhood of $K$. Then for $s\in H^m(K,E),$
$\tau_j(\phi_js)$ can be identified as an $r$-tuple of functions on $U_j$. Suppose $U_j'$ is isomorphic to an open subset $\Omega_j'\subset\Rn.$ If $\psi_j:U_j'\to \Omega_j'$ is an isomorphism let $\psi_j(U_j)=\Omega_j.$ Then $\tau_j(\phi_js)$ can be considered as an $r$-tuple $s_j=\tau_j(\phi_js)\circ\psi^{-1}_j$ of functions on $\Omega_j$ and $s_j\in H^m(\Omega_j).$
Let $U=\{U_1,\ldots,U_m\}$ and
\begin{equation}
\abs{s}^2_{m,U}:=\sum_j \abs{s_j}_m^2
\end{equation}
with respect to this norm $H^m(K,E)$ is a Hilbert space. The map $\eta:H^m(K,E)\to \oplus H^m(\Omega_j)$,
$\eta(s)=\oplus s_j$ defines an isometry of $H^m(K,E)$ onto a subspace of $\oplus H^m(\Omega_j).$
Furthermore,
$\tau_j(s|_{U_j})\circ\psi^{-1}_j=s'\in H^m(\Omega_j).$ Setting $\norm{s}_{m,U}^2:=\sum_j\abs{s_j'}_m^2$ the space $H^m(K,E)$ is again a Hilbert space with respect to this norm
and $\abs{s}^2_{m,U}$, $\norm{s}_{m,U}^2$ are equivalent norms on $H^m(K,E)$ and different coverings $U$ and different partitions of unity $\{\phi_j\}$ render equivalent norms.
\end{definition}
If $f=(f_1,\ldots,f_q)\in L^p(\Omega,q),$ by which we mean that each $f_j\in L^p(\Omega)$, and there exists $h^\alpha\in L^k(\Omega,q)$
for $1\leq k\leq \infty,$ $\abs{\alpha}\leq m$ so that for all $g\in C^\infty_c(\Omega,q)$
\begin{equation}
\int_\Omega \langle f(x),\partial^\alpha g(x)\rangle =(-1)^{\abs{\alpha}} \int_\Omega \langle h^\alpha(x),g(x)\rangle dx
\end{equation}
then we say that $f$ has weak derivatives, $h^\alpha$, up to order $m$ in $L^k.$ Let $\{f_\nu\}_\nu$ be a sequence in $C^\infty(\Omega,q)$ which converges in $L^p(\Omega)$ for $\abs{\alpha}\leq m$ in $H_{m,p}(\Omega).$ Then $\{\partial^\alpha f_\nu\}$ converges in $L^p(\Omega)$ for $\abs{\alpha}\leq m$ to a limit $f^{(\alpha)}$.
For each $g\in C^\infty_c(\Omega,q)$ we have 
\begin{equation}
\langle f_\nu,\partial^\alpha g\rangle=(-1)^{\abs{\alpha}}  \langle \partial^\alpha f_\nu,g\rangle\to \langle f^{(\alpha)}_\nu,g\rangle
\end{equation}
Also if two sequences $\{f_\nu\}_\nu$, $\{g_\nu\}_\nu$ define the same element in $H_{m,p}(\Omega)$ then 
\begin{equation}
\lim \partial^\alpha f_\nu=\lim \partial^\alpha g_\nu,\quad \abs{\alpha}\leq m
\end{equation}
thus for $f\in H_{m,p}(\Omega)$ we can define $\partial^\alpha f:=\lim \partial^\alpha f_\nu$ in $H_{m,p}(\Omega),$ $f_\nu\in C^\infty(\Omega,q).$
If $0\leq m'\leq m$ and $f\in H_{m,p}(\Omega)$ then there exists a sequence  $\{f_\nu\}_\nu$, $f_\nu\in C^\infty(\Omega,q),$
that is Cauchy with respect to $\abs{\cdot}_{m,p}$, defining $f$. Since
$\abs{\cdot}_{m',p}\leq \abs{\cdot}_{m,p}$ the sequence is Cauchy with respect to $\abs{\cdot}_{m',p}$ and defines an element $f'\in H_{m',p}(\Omega)$.
independent of the choice of sequence $\{f_\nu\}_\nu$ defining $f.$ 
Set $i(f):=f'.$ Clearly, $i(\mathring{H}^{m,p}(\Omega))\subset \mathring{H}^{m',p}(\Omega).$
We mention the following result which we state without proof.
\begin{theorem}[Rellich Lemma, see Narasimhan \cite{narasimhan}, p.228]\label{relichlemma}
	Let $\Omega\subset\Rn$ be a bounded subset and $0\leq m'<m.$ Then the natural map
	natural map $i:\mathring{H}^{m,p}(\Omega)\to \mathring{H}^{m',p}(\Omega)$ is completely continuous (by which
	is meant that is sends weakly convergent sequences to norm-convergent sequences).
	As a consequence the induced injection
	$i:H^{m}(\Omega,E)\to H^{m-1}(\Omega,E),$ $m\geq 1,$ is completely continuous.
\end{theorem}
\begin{theorem}[Elliptic analyticity theorem]\label{analyticity}\index{Elliptic analyticity theorem}
Let $V$ be an analytic manifold, $E, F$ analytic vector
bundles on $V$ and $L$ an elliptic operator from $E$ to $F$ with analytic
coefficients (i.e.\ for any analytic section $s: U \to E$, $Ls$ is an analytic
section $U\to F$). Then if $s$ is a locally square integrable section such
that $Ls$ is analytic, then $s$ is analytic.
\end{theorem}
The following are a subset of the so-called Friedrichs’ inequalities.\index{Friedrichs’ inequalities}
\begin{theorem}[Friedrichs' inequalities]\label{friedrichs}
	Let $\Omega\subset\Rn$
	be a bounded open subset and let $L$ be 
	an elliptic differential operator on 
	of order $m$, given by
	$L=\sum_{\abs{\alpha}\leq m}a_\alpha \partial^\alpha.$ Let $r\in \Z_{\geq 0}.$
	If $\Omega'',\Omega'$ are relatively compact open subsets of $\Omega$
	with $\Omega''\Subset\Omega'$ then 
	there exists a constant $C_1$ such that for any $f\in C^\infty(\Omega
	,q)$
	\begin{equation}
	\abs{f}^{\Omega'}_{m+r}\leq C_1(\abs{Lf}_r^{\Omega} +\abs{f}_0^{\Omega})
	\end{equation}
	In particular, if $f\in C^\infty_c(\Omega',r)$ then 
	\begin{equation}
	\abs{f}_{m+r}\leq C_1(\abs{Lf}_r +\abs{f}_0)
	\end{equation}
\end{theorem}
\begin{proof}
	Let $\{U_1,\ldots,U_N\}$ be a covering of $\overline{\Omega}'$ such that,
	if supp$f\subset U_i$ then 
	$\abs{Lf}_{k} \geq C''\abs{f}_{m+k}$ for a constant $C''$. 
	Let supp$f\subset \Omega'$ and $\phi_1,\ldots,\phi_N$
	$\phi_i\in C^\infty$ with supp$\phi_i\subset U_i$, $0\leq \phi_i\leq 1,$ $\sum\phi_i^2=1$ on a neighborhood of
	$\overline{\Omega}'.$
	For supp$f\subset\Omega'$ we have for a constant 
	$C''>0,$ $\abs{\abs{f}_{m+k}^2-\sum_i\abs{\phi_i f}_{m+k}^2}
	\leq C''\abs{f}_{m+k-1}$ and
	$\abs{\abs{Lf}_{m+k}^2-\sum_i\abs{L(\phi_i f)}_{k}^2}
	\leq C'\abs{f}_{m+k-1}$. Now $\abs{\phi_i f}_{m+k}\leq C''\abs{L(\phi_i f)}_k$.
	\begin{lemma}\label{showlemo}
	For any $\epsilon>0$ there exists a constant $C(\epsilon)>0$ such that
	\begin{equation}
\abs{f}_{m-1}^2\leq \epsilon\abs{f}_m^2+C(\epsilon)\abs{f}_0^2,\quad f\in \mathring{H}_m(\Omega)
	\end{equation}	
	\end{lemma} 
\begin{proof}
	It suffices to give the proof for $f\in C^\infty_c(\Rn).$
	First note that since $C^\infty_c(\Rn,q)$ is dense in $\mathring{H}(\Rn)$ there exists constants $c_1,c_2,$ such that
	\begin{equation}
	c_1(1+\abs{\xi}^2)^m\leq \sum_{\abs{\alpha}\leq m}\abs{\xi^\alpha}^2\leq c_2(1+\abs{\xi}^2)^m
	\end{equation}
	If $f=(f_1,\ldots,f_q)$ we have
	\begin{multline}
	\abs{f}^2_m = \sum_{\abs{\alpha}\leq m} \sum_{j=1}^q \int\abs{\partial^\alpha f_j(x)}^2 dx=\\
	\sum_{\abs{\alpha}\leq m} \sum_{j=1}^q \int\abs{\widehat{\partial^\alpha f_j}(\xi)}^2 d\xi=
	 \sum_{\abs{\alpha}\leq m} \sum_{j=1}^q \int\abs{\xi^\alpha}^2\abs{\widehat{f_j}(\xi)}^2 d\xi
	\end{multline}
	This implies that there exists a constant $c_2>0$ such that
	\begin{equation}
	\abs{f}^2_m \leq c_2 \int_\Rn (1+\abs{\xi}^2)^{m-1}\abs{\widehat{f}(\xi)}^2 d\xi
	\end{equation}
	Now given $\epsilon>0$ there exists a constant $C(\epsilon)>0$ such that
	\begin{equation}
	(1+\abs{\xi}^2)^{m-1}\leq \epsilon c_2^{-1}(1+\abs{\xi}^2)^{m} +C(\epsilon),\quad \forall \xi\in \Rn
	\end{equation}
	Hence
	\begin{equation}
	\abs{f}^2_{m-1} \leq \epsilon \int_\Rn (1+\abs{\xi}^2)^{m}\abs{\widehat{f}(\xi)}^2 d\xi
	+C(\epsilon)\int_\Rn \abs{\widehat{f}(\xi)}^2 d\xi
	\end{equation}
	This completes the proof of Lemma \ref{showlemo}.
	\end{proof}
	Lemma \ref{showlemo} implies that there is a constant $C'>0$ such that
	\begin{equation}
	\abs{f}_{m+k} \leq C'\left(\abs{Lf}_{k}+\abs{f}_{0}\right),\quad \forall f\in C^\infty_c(\Omega',r)
	\end{equation}
	Now let 
	$\phi\in C_c^\infty(\Omega)$, $0\leq \phi\leq 1,$ $\phi=1$ on a neighborhood of $\overline{\Omega}',$
	and $f\in C^\infty(\Omega,r).$ Let $m'=m+k$. For $\abs{\alpha}\leq m'$ we have
	some functions $C_\beta\in C^\infty_c(\Omega)$ and a constant $A_1,$ such that 
	$\abs{C_\beta}<A_1$ we have
	\begin{equation}
	\partial^\alpha(\phi^{m'} f)=\sum_{\beta\leq \alpha} \binom{\alpha}{\beta}(\partial^\beta \phi^{m'})\partial^{\alpha-\beta}f=
	\phi^{m'}\partial^\alpha f +\sum_{\beta\leq \alpha, \beta\neq 0} C_\beta \phi^{m-\abs{\beta}}	\partial^{\alpha-\beta}f
	\end{equation}{ahernbruna}
	Squaring both sides, applying the Schwarz inequality and summing over $\alpha,$ $\abs{\alpha}\leq m'=m+k,$ gives
	for some constant $A_2>0,$
	\begin{equation}
	\abs{\abs{\phi^{m+k}f}^2_{m+k}-\sum_{\abs{\alpha}\leq m+k}
		\abs{\phi^{m+k})\partial^{\alpha}f}_0^2} \leq A_2\sum_{\abs{\beta}\leq m+k}
	\abs{\phi^{\abs{\beta}}\partial^\beta f}_0^2
	\end{equation}
	
	Since supp$\phi^{m'}f\subset\Omega'$, we have, by what already done there is a constant $A_3$ such that
	\begin{equation}
	\abs{\phi^{m+k}}^2_{m+k}\leq C\left(\abs{L(\phi^{m+k}f)}_r^2+(\abs{f}_{0}^{\Omega})^2\right)
	\end{equation}

	\begin{lemma}\label{narasimlem2}
		Let $\phi\in C^\infty_c(\Omega).$ For any $\epsilon>0$ there is a constant $C(\epsilon)$ such that for $k\geq 1$ and $f\in C^\infty(\Omega)$
		\begin{equation}
		\sum_{\abs{\beta}=k}\abs{\phi^k\partial^\beta f}_0^2 \leq \epsilon \sum_{\abs{\beta}=k+1} \abs{\phi^{k+1}\partial^\beta f}_0^2 +C(\epsilon)
		\sum_{\abs{\beta}=k+1} \abs{\phi^{k-1}\partial^\beta f}_0^2
		\end{equation}
	\end{lemma}
	\begin{proof}
		It suffices to prove that for $k\geq 1,$ $\abs{\beta}=k,$ we have for each $\epsilon>0$ a constant $C(\epsilon)$ such that
		\begin{equation}
		\abs{\phi^k\partial^\beta f}_0^2 \leq \epsilon \sum_{\abs{\alpha}=k+1} \abs{\phi^{k+1}\partial^\alpha f}_0^2 +
		C(\epsilon)\sum_{\abs{\alpha}=k-1}\abs{\phi^{k-1}\partial^\alpha f}_0^2
		\end{equation}
		Let $\beta=\gamma+e$ for some $e$ with $\abs{e}=1.$ Then
		
		\begin{multline}
		\langle \phi^k\partial^\beta f,\phi^k\partial^\beta f \rangle=
		-\langle \partial^\gamma f,\partial^e\phi^{2k}\partial^\beta f \rangle
		=\\
		-\langle \partial^\gamma f,2k\phi^{2k-1}\partial^e \phi \partial^\beta f \rangle
		-\langle \partial^\gamma f,\phi^{2k} \partial^{\beta+e} f \rangle
		=\\
		-2k\langle \partial^{e} \phi\cdot \phi^{k-1}\partial^\gamma f,\phi^{k} 
		\partial^\beta f \rangle -
		\langle \phi^{k-1}\partial^\gamma f,\phi^{k+1} \partial^{\beta+e} f
		\rangle 
		\end{multline}
		Using the inequality
		\begin{equation}
		2\abs{\langle u,v\rangle}\leq \delta \abs{u}_0^2 +\frac{1}{\delta}\abs{v}_0^2,\quad \forall \delta>0
		\end{equation}
		gives for $\delta>0$
		\begin{equation}
		\abs{\phi^{k}\partial^\beta f}_0^2\leq \delta\left( 
		\abs{ \phi^{k}\partial^\beta f}_0^2 + 
		\abs{ \phi^{k+1}\partial^{\beta+e} f}_0^2\right) +C_1(\delta)\abs{ \phi^{k-1}\partial^\gamma f}_0^2  
		\end{equation}
		This yields the wanted inequality. This proves Lemma \ref{narasimlem2}.
	\end{proof}
	
	Application of Lemma \ref{narasimlem2} to $\abs{\phi^{\abs{\beta}}\partial^\beta f}_0^2$ 
	gives
	\begin{equation}
	\sum_{\abs{\beta}\leq m+k}
	\abs{\phi^{\abs{\beta}}\partial^{\beta}f}_0^2 \leq 
	\epsilon \sum_{\abs{\beta}\leq m+k}
	\abs{\phi^{m+k}\partial^{\beta}f}_0^2 +C(\epsilon)\left(\abs{f}_0^{\Omega}\right)^2
	\end{equation}
	which implies that there, for sufficiently small $\epsilon>0,$ is a constant $C_2>0$ such that
	\begin{equation}
	\sum_{\abs{\alpha}\leq m+k}
	\abs{\phi^{m+k}\partial^{\alpha}f}_0^2 \leq C_2 \left(
	\abs{\phi^{m+k}L f}_k^2 +\left(\abs{f}_0^{\Omega}\right)^2\right)
	\end{equation}
	Since $\phi=1$ on $\Omega'$ and supp$\phi\subset\Omega$, this completes the proof.
	This completes the proof.
\end{proof}

\chapter{Some preliminaries from CR geometry}\label{crapp}
A real linear map $J$ on an even dimensional real vector space $V$ is called a complex structure if $J^2=-Id$.
In particular, a complex structure $J$ on $T_p \Cn$ %
is defined as a real linear map $J:T\Cn\to T\Cn$
such that $J^2=-Id,$ specifically $J$ is defined fiberwise on the tangent vector spaces by $\R$-linear maps $J_p: T_p \Cn\to T_p \Cn.$
If $M\subset \Cn$ is a submanifold %
$T^c_p M:=T_p M\cap J_p(T_p M)$ is called
the holomorphic tangent space of $M$ at $p.$ $J$ maps each $T^c_p M$ to itself thus defines a complex structure on $T^c_p M.$
If $T^c_p M$ has constant dimension ($CR$ dimension) for every $p$ then $M$ is called a $CR$ manifold. 
The
$\R$-linear maps $J_p: T_p^c M\to T_p^c M$
have eigenvalues $\pm i.$
$J$ extends to a $\C$-linear map on the bundle $\C\otimes T^cM=\bigcup_{p\in M} \C\otimes T_p^c M $ such that the extension again has eigenvalues $\pm i.$ This decomposes $\C\otimes T^c M=H^{1,0} M\oplus H^{0,1} M$ namely a $\C$-linear and anti-$\C$-linear part, where we denote $H^{0,1} M$ to be the anti-$\C$-linear part. 
\begin{definition}
A differentiable function $f$ on $M$ is called $CR$ if it is annihilated by any $C^1$-section $X$ of $H^{0,1} M$ over $M.$
A distribution $f$ is called $CR$ if 
$Xf=0$ %
in the weak sense
i.e.\ $\langle f,X^{\mbox{adj}}\phi\rangle =0,\forall \phi\in C^{\infty}_c(M),$
where $X^{\mbox{adj}}$ denotes the adjoint.
\end{definition}
Let $M\subset \Cn$ be a generic $CR$ submanifold  
and let $p\in M.$ The Euclidean metric on $\C\otimes T_p \Cn$ induces a metric on 
$\C\otimes T_p M.$
Let $M\subset \Cn$ be a generic $CR$ submanifold  
and let $p\in M.$ The Euclidean metric on $\C\otimes T_p \Cn$ induces a metric on 
$\C\otimes T_p M.$
In particular, the quotient space $\C\otimes T_p M/(H^{1,0}_p M\oplus H^{0,1}_p M)$
can be identified with the orthogonal supplement of 
$H^{1,0}_p M\oplus H^{0,1}_p M$
and furthermore we can identify
$\C\otimes T_p M/(H^{1,0}_p M\oplus H^{0,1}_p M)$
$=\C\otimes (T_p M/T^c_p M).$
Given $X\in H^{1,0}_p M$ we define (the Levi form)\index{Levi form} $\mathcal{L}(X):=\frac{1}{2i}[\overline{\tilde{X}},\tilde{X}]_p$ 
mod $\C\otimes T^c_p M$,
where $\tilde{X}$ is any smooth section of $H^{1,0} M$ extending $X.$
$\mathcal{L}$ is real-valued so the image which lies in $\C\otimes (T_p M/T^c_p M)$, can be identified with 
the real vector space $T_p M/T^c_p M.$ 
Recall that $h\colon V\times V\to \C$
on a complex vector space $V$ is called a Hermitian form if it is linear in the first coordinate and such that
$h(X,Y)=\overline{h(Y,X)}.$ Every Hermitian form
has an associated Hermitian matrix, $A_h$, such that $h(X,Y)=X A_h \overline{Y}^T$.
We consider now
a Hermitian real-valued $\mathcal{L}(X,Y):=\frac{1}{2i}
[\overline{\tilde{X}},\tilde{Y}]_p$ 
mod $\C\otimes T^c_p M.$ Then we can speak of the eigenvalues of the
Hermitian matrix $A_{\mathcal{L}^{\xi}}.$
Let $M\subset \Cn$ be a generic $CR$ submanifold  
and let $p\in M.$
Set $\chi_p:=\{\xi\in T^*_p M:\xi|_{T^c_p M}\equiv 0\}.$ 
The {\em directional Levi form at $p$ in codirection $\xi\in \chi_p$} is defined as $\mathcal{L}^{\xi}
(X):=\langle \xi,\mathcal{L}(X)\rangle.$
A useful fact in higher real codimension is that it is possible to locally describe generic $CR$ manifolds as graphs over their tangent space. In particular, we have the following theorem, see Boggess \cite{b4}, p.103. 
\begin{theorem}[Local graph representation theorem]\label{localgraphlemma}
	Let $M$ is a smooth generic $CR$ submanifold of $\C^n$ of real dimension
	$2n-d,d\leq n$, and let $p\in M.$
	Then there is an open neighborhood $U$ of $p$, 
	a nonsingular complex affine map $A\colon \C^n\to\C^n$
	and a $C^1$-smooth function $h\colon \R^d\times\C^{n-d}\to \R^d$, s.t.,
	\begin{equation}
	A(M\cap U)=\{ (x+iy,w)\in \C^d\times \C^{n-d}:y=h(x,w)  \} ,
	\end{equation}
	where $h(0)=0$, $\mbox{Jac} h(0)=0.$
\end{theorem} 
\begin{proof}
	Assume that $p=0$, by translation. We show that there is a nonsingular
	complex linear map $A\colon \Cn\to \Cn$ that maps $T_{0}(M)$ to $\{ 
	\R^d \ni\im z=0\} \subset \R^{2n} .$ If this holds then the graphing function for $A\{ M\}$ satisfies the criteria for the function $h.$
	We know that $T_{0}M/ (T_{0}M \cap JT_{0}M )$ is a real vector space dimension $d$ and can be identified as a real vector subspace of $T_0 \Cn.$
	We let $b_1 ,\ldots , b_d$ be an orthonormal basis for $T_{0}M/ (T_{0}M \cap JT_{0}M )$ with respect to the induced inner product from
	$T_0 \Cn$. Then
	the $J$-invariant space $(T_{0}M \cap JT_{0}M )$ is orthogonal to $\mathcal{A}_{0}:=T_{0}M/ (T_{0}M \cap JT_{0}M )$ and $J\{ \mathcal{A}_{0} \} .$
	Now $J$ has only the imaginary eigenvalues $\pm i$ and
	$Ja\cdot b=-a\cdot Jb$,
	$a,b\in T_{0}(\R^{2n} ).$
	If $M$ is generic we have an orthonormal basis for $T_{0}(\R^{2n} )
	=T_{0}(M )\oplus J\mathcal{A}_{0} $, 
	consisting of,
	\begin{equation}
	\{ b_1 ,\ldots , b_d , Jb_1 ,\ldots , Jb_d  
	\} .
	\end{equation}
	We set $z=x+iy\in\C^d ,w=u+iv\in\C^{n-d} $, and,
	\begin{align}
	& A(b_j ):=\frac{\partial}{\partial x_j} ,\quad A(Jb_j ):=\frac{\partial}{\partial y_j}, \quad 1\leq j\leq d,\\
	& A(b_j ):=\frac{\partial}{\partial u_{j-d}}, \quad
	A(Jb_j ):=\frac{\partial}{\partial v_{j-d}}, \quad d+1\leq j\leq n.
	\end{align}
	Then $A$ is complex linear, nonsingular and 
	\begin{equation}
	A\{T_{0}M\} =\{
	(x,0,u,v):x\in \R^d, u,v\in \R^{n-d}\} =\{ y=0\} .
	\end{equation}
	This completes the proof.
\end{proof}
More can be said.
Namely (see Boggess \cite{b4}, p.109) given a smooth generic $CR$ submanifold $M\subset \Cn$, $0\in M$, local holomorphic coordinates $(z=x+iy,w)\in \C^d\times \C^m$ such that $M$ has a representation near $0$, given by
	$M=\{ y=h(x,w)\}$,
	$h(0)=\mbox{Jac}( h)(0)=0,$
	such that a basis for $H^{0,1} M$ near the origin is given by vector fields of the form 
	\begin{equation}
	L_j =\frac{\partial}{\partial \bar{w}_j} -2i\sum_{l=1}^d \sum_{k=1}^d \varphi_{lk}
	\frac{\partial h_k}{\partial \overline{w}_j}\frac{\partial}{\partial \overline{z}_k} ,\quad  j=1,\ldots ,m,
	\end{equation}
	where $\varphi_{kl}$ is th $(l,k)$:th entry of the matrix $(I+i(\partial h/\partial x))^{-1}$ where $I$ denotes the identity matrix. 

\section{An approximation theorem and Lewy's theorem}
We shall describe here a proof of Lewy's theorem (see Lewy \cite{lewy},
and a modern generalization due to Boggess \& Polking \cite{bp}) that relies upon a know technique for holomorphic extension, known as the method of analytic discs, which is described in detail in the textbook of Boggess \cite{b4}. We shall also need to be familiar with a known result on local convexification
near strictly pseudoconvex points. 
We cite the theorem here but %
point out that 
it is in fact the method 
of proof (rather than the result itself) of this theorem in Boggess \& Polking \cite{bp} which
is of interest. 
\begin{theorem}[Boggess \& Polking \cite{bp}, textbook form Boggess \cite{b4}, p.200]\label{bpthm0} 
Let $M\subset \Cn$ be a $C^l,l\geq 4,$ generic embedded $CR$ submanifold and let $p\in M$ be a point such that the Levi cone at $p$
has nonempty interior.
Then for every open neighborhood $\omega_1$ of $p$
there is an open neighborhood $\omega_2$ of $p$ and an open $\Omega\subset \Cn$, such that,
\begin{itemize}%
	\item [(a)] $p\in \omega_2\subset \overline{\Omega}\cap M\subset\omega_1$.%
	\item [(b)] For every open cone $\Gamma_1$ smaller than $\Gamma_p$ there is a connected neighborhood $\omega_3$ of $p$ in $M$ and $\epsilon >0$
	such that $\omega_3+\Gamma_1\cap B_{\epsilon}\subset \Omega.$
	\item[(c)] For every continuous $CR$ function $f$ on $\omega_1$ there is a unique holomorphic function $F$ on $\Omega$
	which is continuous on $\overline{\Omega}$ such that $F|_{\omega_2}=f.$
\end{itemize}
\end{theorem}
The proof involves explicit construction of families of analytic discs by solving a Bishops equation, such that the center of these discs pass each point of an open subset of the given normal cone and simultaneously are attached sufficiently close to $p.$ 
For completeness we give a short account for this technique 
and as one of the key elements in the proof is the Baouendi \& Treves \cite{bt} approximation theorem.
(the special case we cite here is formulated in
Boggess \& Polking \cite{bp}, Theorem 2.1, p.761).
\begin{theorem}\label{approthm}
Let $M\subset\Cn$ be a smooth
and $p_0\in M.$
Every open neighborhood $\omega$ of $p_0$ in $M$ contains another open neighborhood $\omega_1$ of
$p_0$ in $M$ such that every $C^1$-smooth $CR$ function on $\omega$ is the uniform limit in $\omega_1$ of holomorphic polynomials.
\end{theorem}
Generalization with respect to regularity can be found in 
Berhanu et al.\ \cite{ber2}, p.53.
\\
\\
Recall that a differential $k$-form $\alpha$ is a smooth section of the vector bundle defined by the $k$:th exterior power of the cotangent bundle, $T^*M,$
i.e.\ $\alpha :\wedge_{j=1}^k T^* M \to M,$ such that $\pi\circ \alpha =Id,$ where $\pi$ is the projection of the vector bundle to $M.$
For complex differential forms use first the splitting $dz_j =dx_j +idy_j, d\bar{z_j} =dx_j -idy_j$ and then 
$dz_{j_1}\wedge \cdots \wedge dz_{j_p} \wedge d\bar{z}_{j_1} \wedge \cdots \wedge d\bar{z}_{j_q}$ is well defined.
\begin{definition}[Defining functions]
Given coordinates $\phi =(\phi_1 ,\ldots ,\phi_n)$ near $p\in M$ where $M$ is a $k$ dimensional submanifold of $X$ which in turn is a real $n$ dimensional manifold, such that $\phi(0)=p$ and
$M\cap U=\{\phi_j =0 ,k+1 \leq j\leq n\} ,$ then the derivative $d(
\{ \phi_{k+1} , \cdots , \phi_n \}),$
has maximal rank which means that the differentials spanning the image, namely, $d\phi_j ,k+1 \leq j\leq n ,$ are linearly independent, which implies that
\begin{equation}
d\phi_{k+1} \wedge \cdots \wedge d \phi_{n} \neq 0 .
\end{equation}
A set of $C^{\infty}$ (or $C^k$) real valued functions 
$\{ \rho_{1} , \cdots , \rho_d \},$ such that
\begin{equation}
M=\{ \rho_{1}=0 , \cdots , \rho_{d} =0 \}, d\rho_{1} \wedge \cdots \wedge d \rho_{d} \neq 0 \mbox{ on }M
\end{equation}
is called a set of {\em defining functions}\index{Defining functions} for $M.$
\end{definition}

\section{The method of analytic discs}\label{methodapp}
Analytic discs are an important tool in the proof of (the generalization of) Lewy's theorem found in Boggess \cite{b4}. 
\begin{definition} An {\em analytic disc}\index{Analytic disc} is a continuous map, 
\begin{equation} 
\mathcal{D} :\overline{D}\to \C^n
\end{equation} 
which is holomorphic on $D.$ Its {\em boundary}\index{Boundary of an analytic disc} is $\Delta|_{\partial D} ,$ and we identify an analytic disc with its image in $\Cn $ (here $D$ denotes the unit disc in $\C$). 
\end{definition}
The so called {\em Bishop equation} is a functional equation that can be used to prescribe
certain characteristic features of an analytic disc,
and in particular it involves the Hilbert transform. 
A major part of the method of analytic discs, whose purpose is to provide $CR$ extension along the interiors of analytic discs, is approximation of $CR$ functions, and this is provided by a theorem due to Baouendi \& Treves \cite{bt}, which we give a more detailed presentation of in Appendix \ref{approxapp}, the original theorem applies to more general smooth first order integrable system of vector fields and provides
local uniform approximation in terms of entire functions.
When applied for the locally integrable system which define the tangential Cauchy Riemann equations (i.e.\ whose homogeneous solutions are $CR$ functions, see Boggess \& Polking \cite{bp} for an interpretation) one obtains the following theorem for continuous $CR$ functions, see e.g.\  Boggess \& Polking \cite{bp}, p.761. 
\begin{theorem}[Baouendi \& Treves \cite{bt}]
Let $p\in \Omega,$ with $\Omega\subset\Cn$ a smooth generic submanifold. Then every neighborhood $U$ of $p$
contains another neighborhood $V$ of $p_{0}$ such that every continuous $CR$ function on $U$
is the uniform limit in $V$ of holomorphic polynomials.
\end{theorem}
We will prove a more general version of this in Theorem \ref{btthm}. 
By applying this version of the approximation theorem of Baouendi \& Treves for continuous $CR$ functions we have that 
there is a sequence of approximating holomorphic polynomials, uniformly convergent on $V\subset \Omega,$ which we will denote $\{ P_j\} .$ 
Next one ensures, by construction, that for a given $p_0\in \partial\Omega ,$ our family of analytic discs, which we associate to $p_0$, will be attached to $\Omega,$ which means that their boundaries (i.e.\ the image, under the analytic disc map, of $\partial D,$ where $D\subset \C$ denotes the unit disc) belong to $\Omega.$ 
In particular, we have uniform convergence of the approximating sequence $\{P_j\}_j ,$ on
the boundary of our analytic discs. By applying the maximum principle for the holomorphic functions, $P_j,$
we can ensure uniform convergence on the interior of the discs. By 
another application of the maximum principle we have uniform convergence on the interior of the disc hull (by disc hull we mean the union of the interior of the analytic discs of our family). 
The limit function, which we denote by $F,$ is holomorphic on a small enough neighborhood of $p_0$ in the disc hull, due to it being the uniform convergence on compacts, of a sequence of holomorphic functions (since the disc hull fills a neighborhood in the ambient space, i.e.\ no real codimension), see Theorem \ref{hormkonvthm}. 
Finally, $F$ serves as a holomorphic extension of $f$ to an open ambient subset consisting of the disc hull.

\section{Local convexification near strictly pseudoconvex points}\label{localconv}
Strict pseudoconvexity is not necessary
for constant modulus on an open set to imply reduction to a constant, but it provides perhaps the simplest affirmative cases because we can use that
a strictly pseudoconvex point of a $C^2$-smooth hypersurface allows local convexification. 
The idea of such convexification can then be generalized to $1$-convex points in order to prove Lewy's theorem for hypersurfaces.
For completeness 
we give here as an appendix a short sketch of the proof for the classical convexification result.
\begin{lemma}[Narasimhan \cite{narasimconvex}]\label{lemmma}
If $M\subset \Cn$ is a $C^2$ real hypersurface which is strictly pseudoconvex at $p_0\in M,$ then $\exists $ an ambient neighborhood $U$ of $p_0$ and a biholomorphism $\phi ,$ such that $\phi(U\cap M)$ is a strictly convex hypersurface in $\phi(U)\subset\Cn .$
\end{lemma}
\begin{proof}
Let $\rho :\Cn\to \R$ be a local defining function for $M$ near $p$ (i.e.\ there exists a domain $p\in U\subset \Cn$ such that 
$M\cap U=\{ \rho=0\}\cap U$ and $d\rho_p\neq 0.$ By strict pseudoconvexity at $p$
the Levi form of $\rho$ is strictly positive definite at a point $p,$
introduce new holomorphic ambient coordinates (recalling that the Levi signature is invariant under biholomorphism) in which the Hessian at
$p$ coincides with the Levi form. 
Let $q_1,\ldots ,q_n$ be a local holomorphic coordinates in $X$ centered at $p$ (i.e.\ $q(p)=0$). After a rotation we can assume that $(1,0,\ldots,0)$ is the outward normal vector at the origin.
Write the second order expansion of $\rho$ near the origin.
\begin{multline}
\rho(Q)=\underbrace{\rho(0)}_{=0} +
\sum_{j=1}^n \frac{\partial\rho}{\partial q_j}(0)q_j +\frac{1}{2}
\sum_{j,k=1}^n \frac{\partial^2\rho}{\partial q_j \partial q_k}(0)Q_j Q_k
+\sum_{j=1}^n \frac{\partial\rho}{\partial \bar{q}_j}(0)\bar{q}_j+\\
\frac{1}{2} \sum_{j,k=1}^n \frac{\partial^2\rho}{\partial \bar{q}_j \partial \bar{q}_k}(0)\bar{Q}_j \bar{Q}_k +
\sum_{j,k=1}^n \frac{\partial^2\rho}{\partial q_j \partial \bar{q}_k}(0)Q_j \bar{Q}_k
+o(\abs{q}^2)=\\
2\mbox{Re}\left(
\sum_{j=1}^n \frac{\partial\rho}{\partial q_j}(0)q_j + \frac{1}{2} \sum_{j,k=1}^n \frac{\partial^2\rho}{\partial q_j \partial q_k}(0)Q_j Q_k
\right) +
\sum_{j,k=1}^n \frac{\partial^2\rho}{\partial q_j \partial \bar{q}_k}(0)Q_j \bar{Q}_k +o(\abs{w}^2).
\end{multline}
By defining new holomorphic coordinates $W=W(Q)$ as follows, 
\begin{equation}
\left\{
\begin{array}{ll}
W_1=Q_1+ \frac{1}{2}\sum_{j,k=1}^n \frac{\partial^2\rho}{\partial q_j \partial q_k}(0)Q_j Q_k \\
W_j=Q_j
\end{array}
\right.
\end{equation}
one obtains,
\begin{equation}
\rho(W)=\mbox{Re}(W_1)+\frac{1}{2}\underbrace{\sum_{j,k=1}^n \frac{\partial^2\rho}{\partial w_j \partial \bar{w}_k}(0)W_j \bar{W}_k}_{=\mbox{Levi form of $\rho$ at $0$}} +o(\abs{W}^2).
\end{equation}
Since the sum in the middle is the Hessian at $0$ and since it coincides with the Levi form at $0$ we have that the Hessian is strictly positive. 
We see that since the Hessian at $0$ is positive (by strict Levi pseudoconvexity) we have by continuity that it is positive on a neighborhood of $0$
thus we have that $\partial \{\rho <0\}$ is in the new coordinates convex near $0.$
\end{proof}

\section[A proof of the approximation theorem]{A proof for the approximation theorem of Baouendi \& Treves in the differentiable case}\label{approxapp}
In this appendix we provide a short proof of the 
the approximation theorem of Baouendi \& Treves \cite{bt}. Several proofs exists in the literature, we shall follow
that of Baouendi, Ebenfeldt \& Rothschild \cite{barot}, an alternative proof can be found in Boggess \cite{b4}, where the theorem is proved for $C^1$-smooth $CR$ functions. The result holds true more generally for $L^p$-distributions.

The result of Baouendi \& Treves is originally one for locally integrable structures. The standard version which appears
in $CR$ gemoetry is the special case where the locally integrable structure is given by the Cauchy-Riemann bundle
(See Section \ref{hypoanalsec} and Definition \ref{hanvisacrhypoanal}). Note that we shall here only prove a local result but it will be apparent that 
the proof will be invariant with respect to the base point (which in our case is chosen to be the origin). Furthermore, we only give the $C^1$-version.

\begin{theorem}[Baouendi \& Treves \cite{bt}]\label{btthm}
	Let $L_1,\ldots,L_m$ be linearly independent $C^\infty$-smooth complex vector fields on an open subset $\Omega\subset\R^{n+m},$ 
	which are locally integrable
in the sense that there exists, in a neighborhood of each point in $\Omega$, $n$ complex-valued solutions (so called basic solutions), denoted $Z_1,\ldots,Z_n$,
satisfying 
\begin{equation}
L_j Z_k=0, \quad j=1,\ldots,m,\quad k=1,\ldots,n 
\end{equation}
such that the differentials $d Z_1,\ldots,dZ_n$ are linearly independent at each point of that neighborhood.
Let $0\in \Omega$ and $Z=(Z_1,\ldots,Z_n)$ a system of basic solutions) as above with $Z(0)=0$
an an open $\Omega'\subset\Omega$ with $0\in\Omega'$. Then there exists an open neighborhood $\Omega''\ni 0$,
such that each $C^1$-solution of the system
\begin{equation}\label{btrevekv1}
L_j h=0,\mbox{ on }\Omega',j=1,\ldots,m  
\end{equation}
is the uniform limit, in $\Omega''$, of a sequence of polynomials with complex coefficients 
in $Z_1,\ldots,Z_n.$  
\end{theorem}
\begin{proof}
	Let $h$ be a $C^1$-solution to Eqn.(\ref{btrevekv1}).
	We shall denote  $Z=(Z_1,\ldots,Z_n).$  
We begin with an appropriate choice of local coordinates.
Let $u=(u_1,\ldots,u_{n+m})$ be local coordinates vanishing at $0$ in $\Omega$ such that
the determinant of the Jacobian matrix
$A(u)=[\frac{\partial Z_j}{\partial u_l}(u)]_{jl}$ is nonzero for $u$ near $0$ in $\R^{n+m}.$
Setting
\begin{equation}
x:= (\re A(0))^{-1}Z(u),\quad y_j:=u_{j+n},\quad j=1,\ldots,m
\end{equation}
we have local coordinates $(x,y)=(x_1,\ldots,x_n,y_1,\ldots,y_m)$ near $0$ in $\Omega$, vanishing at $0$ 
such that
\begin{equation}\label{btekv244}
Z_j(x,y)=x_j+i\phi(x,y),\quad j=1,\ldots,n
\end{equation}
for smooth $\phi_j$, $j=1,\ldots,n$ defined near $0$ in $\R^{n+m}$ such that $\phi(0)=0$ and $\mbox{Jac}_x\phi(0)=0.$
Let $\chi\in C^\infty_c(\R^n)$ such that $\chi(x)\equiv 0$ for $\abs{x}\leq \frac{r}{2}$ and $\chi(x)\equiv 1$ for $\abs{x}\geq r$ for a positive $r$.
Set $K:=\{(x,y):\abs{x}\leq\frac{r}{4},\abs{y}\leq d\}$ for a positive $d.$
For two vectors $(v_1,\ldots,v_n),(w_1,\ldots,w_n)\in \C^n$ we denote $v\cdot w=\sum_{j=1}^n v_jw_j.$
Define the $m$-form $\alpha(x,y)=\alpha_\nu(x,y;z)$ for $z\in \C^n$ and $\nu\in \Z_+$
\begin{equation}
\alpha_\nu(x,y;z):=\left(\frac{\nu}{\pi}\right)^{\frac{n}{2}}\exp(-\nu(z-Z(x,y))\cdot(z-Z(x,y)))\chi(x)h(x,y)dZ(x,y)
\end{equation} 
where $dZ(x,y)=dZ_1(x,y)\wedge\cdots\wedge dZ_n(x,y).$
Fix $y\in \R^m$ such that $0<\abs{y}<d$ and set
\begin{equation}\label{btrevekv6}
D=D_y:=\{(x',y')\in \R^{n+m} :\abs{x'}<r,\abs{y}=ty,t\in (0,1)\}
\end{equation}
By Stokes's theorem we have
\begin{equation}\label{btrevekv7}
\int_D d\alpha_\nu(x',y';z)=\int_{\partial D}\alpha_\nu(x',y';z)
\end{equation} Since $\chi(x')\equiv 0$ for $\abs{x'}\geq 0$ we have
\begin{multline}
\int_{\partial D}\alpha_\nu(x',y')=\\
\int_{\Rn}\left(\frac{\nu}{\pi}\right)^{\frac{n}{2}}\exp(-\nu(z-Z(x',y))\cdot(z-Z(x',y)))\chi(x')h(x',y)d_{x'}Z(x',y)-\\
\left(\frac{\nu}{\pi}\right)^{\frac{n}{2}}\int_{\Rn}\exp(-\nu(z-Z(x',0))\cdot(z-Z(x',0)))\chi(x')h(x',0)d_{x'}Z(x',0)
\end{multline}
where $d_{x'}$ denotes taking the differential with respect to the $x'$ variables.
Let $\mathcal{L}$ denote the locally integrable subbundle spanned by $L_1,\ldots,L_m.$
By definition 
we have for $p$ 
in $\Omega',$ that
$dh(p)$ belongs to the span of $\{dZ_1(p),\ldots,dZ_n(p)\}.$ Since $f(hdZ)=df\wedge dZ$,
where $dZ:=dZ_1(x,y)\wedge\cdots\wedge dZ_n(x,y)$, we have
$d(hdZ)=0$ near the origin. Note that $\exp(-\nu(z-Z(x',y'))\cdot(z-Z(x',y')))$ is holomorphic with respect to the $Z_j(x',y')$ and since 
the product of two solutions is again a solution we have
\begin{multline}
d\alpha_\nu(x',y'):=\\
\left(\frac{\nu}{\pi}\right)^{\frac{n}{2}}\exp(-\nu(z-Z(x',y'))\cdot(z-Z(x',y')))h(x',y')d\chi(x')\wedge dZ(x',y')
\end{multline}
Define the sequence of entire functions in $\Cn$
\begin{multline}
H_\nu(z):=\\
\left(\frac{\nu}{\pi}\right)^{\frac{n}{2}}\int_{\Rn}\exp(-\nu(z-Z(x',0))\cdot(z-Z(x',0)))h(x',0)\chi(x')d_{x'}Z(x',0)
\end{multline}
\begin{lemma}\label{btlemma}
	For any $C^2$-smooth function $G$ on $\Omega$ 
	we have for sufficiently small $r,d$
	\begin{equation}
	\left(\frac{\nu}{\pi}\right)^{\frac{n}{2}}\int_{D}\exp(-\nu(Z(x,y)-Z(x',y'))\cdot(Z(x,y)-Z(x',y')))G(x',y')d\chi(x')\wedge d_{x'}Z(x',y')\to 0
	\end{equation}
	uniformly for $(x,y)\in K$ as $\nu\to \infty.$
\end{lemma}
\begin{proof}
	Since $d_\chi(x)=0$ for $\abs{x}\leq\frac{r}{2}$ and $\abs{x}>r$
	we have 
	\begin{multline}
	\int_{D}\exp(-\nu(z-Z(x',y'))\cdot(z-Z(x',0)))G(x',y')d\chi(x')\wedge d_{x'}Z(x',y')=\\
	\int_{\stackrel{\frac{r}{2}\abs{x'}\leq r}{y'\in[0,y]}}\exp(-\nu(z-Z(x',y'))\cdot(z-Z(x',0)))G(x',y')d\chi(x')\wedge d_{x'}Z(x',y')
	\end{multline}
	By Eqn.(\ref{btekv244}) we have
\begin{equation}
\re ((Z(x,y)-Z(x',y'))\cdot(Z(x,y)-Z(x',y'))=(x-x')^2-(\phi(x,y)-\phi(x',y'))^2
\end{equation}
By the mean value theorem we have
\begin{multline}
\abs{\phi(x,y)-\phi(x',y')}\leq \abs{\phi(x,y)-\phi(x',y)} +\abs{\phi(x',y)-\phi(x',y')}\leq\\
C_1\abs{x-x'}+C_2\abs{y-y'}
\end{multline}
where 
\begin{equation}
C_1:=\sup_{\stackrel{\abs{x'}\leq r}{\abs{y'}\leq d}}\abs{\phi_{x'}(x',y')}
	,\quad
	C_2:=\sup_{\stackrel{\abs{x'}\leq r}{\abs{y'}\leq d}}\abs{\phi_{y'}(x',y')}
		\end{equation}
		Since $D_{x'}\phi(0)=0$ we may choose $r,d$ sufficiently small such that $\abs{a}\leq \frac{1}{4}$
		and choose $d\leq \frac{r}{16A}$ so that
		\begin{equation}
		\abs{\exp(-\nu(Z(x,y)-Z(x',y'))\cdot(Z(x,y)-Z(x',y'))}\leq \exp(-7\nu\left(\frac{r}{32}\right)^2
		\end{equation}
		This proves Lemma \ref{btlemma}
	\end{proof}

\begin{lemma}\label{btmatrixlemma}
	If $B$ is a real $n\times n$ matrix with $\norm{B} <1,$ and $A:=I+iB,$ then,
	\begin{equation}
	\det A\int_{\Rn} e^{-[Ax]^2 } dx =\pi^{n/2} .
	\end{equation}
\end{lemma}
\begin{proof}
	The integral is holomorphic as a function of the entries of $A^T A,$ where 
	\begin{equation}
	[Ax]^2 =A^T Ax \cdot x\Rightarrow e^{-[Ax]^2} =e^{-(A^T A)x\cdot x},
	\end{equation}
	(here the dot product is extended to a $\C$-bilinear form on $\Cn $) where
	\begin{equation}
	\mbox{Re} (A^T A) =I-B^T B,
	\end{equation}
	together with $\norm{B} <1$ implies that
	$A^T A$ has positive definite real part, and the real part is then necessarily nonsingular
	since for $x\in \Cn,$ 
	\begin{equation}
	(A^T A)x=0\Rightarrow 0=\mbox{Re}\langle Ax, \bar{x}\rangle =\langle (\mbox{Re}(A)) x,\bar{x}\rangle\Rightarrow x=0.
	\end{equation}
	Now $\int_{-\infty}^{\infty} e^{-at^2} dt =(x/a)^{1/2} ,a>0.$ If we replace $a$ with a diagonal matrix $C,$ clearly the equality remains.
	For a symmetric positive definite matrix we diagonalize it via an orthonormal transformation. 
	It is known, see e.g.\ H\"ormander \cite{h2}, p.85, that the set of symmetric matrices with positive definite real part, is convex in the vector space of symmetric $n\times n$ matrices, which implies that there is a unique holomorphic branch of the square root $C\mapsto (\det(C))^{1/2},$ satisfying $(\det(C))^{1/2} >0,$ when $C$ is real. Thus for $A^T A$-real matrix, we have, 
	\begin{equation}
	\int_{\Rn} e^{-[Ax]^2} dx =\pi^{n/2}\det(A^T A)^{-1/2} .
	\end{equation}
	Since both sides of the equation are holomorphic when $A^T A$ is symmetric positive definite, equality remains for all such $A^T A,$
	which implies,
	\begin{equation}
	\int_{\Rn} e^{-[Ax]^2} dx =\pi^{n/2}\det(A^T A)^{-1/2}=\pi^{n/2}(\det(A)^2)^{-1/2} .
	\end{equation}
	This proves the lemma.
\end{proof}
\begin{lemma}\label{btsistalemma}
	For any $C^2$-smooth function $G$ on $\Omega$
	we have for sufficiently small $r,d$
	\begin{multline}\label{btsistaekven}
	\lim_{\nu\to \infty} \left(\frac{\nu}{\pi}\right)^{\frac{n}{2}}\int_{x'\in\Rn}
	\exp(-\nu(Z(x,y)-Z(x',y))\cdot(Z(x,y))-Z(x',y))\times\\
	G(x',y)d\chi(x')\wedge d_{x'}Z(x',y)=G(x,y)
	\end{multline}
	uniformly for $(x,y)\in K$ as $\nu\to \infty.$
	\end{lemma}
	\begin{proof}
	We have $d_{x'}Z(x',y)=\mbox{det}[D_{x'}Z(x',y)]dx'$ and by the change of variables $\xi=\sqrt{\nu}(x'-x)$ the integral in Eqn.(\ref{btsistaekven})
	is transformed to
	\begin{multline}
	\pi^{-\frac{n}{2}}\int_{x'\in\Rn}
	\exp\left(
	-\left(\xi+i\sqrt{\nu}\left(\phi\left(x+\frac{\xi}{\sqrt{\nu}},y\right)-\phi(x,y)\right)     \right)\cdot \right.\\
	\left. \left(\xi+i\sqrt{\nu}\left(\phi\left(x+\frac{\xi}{\sqrt{\nu}},y\right)-\phi(x,y)\right)    \right)
	\right)\times\\
		g\left(x+\frac{\xi}{\sqrt{\nu}},y\right)
			\mbox{det}\left[ D_x Z\left(x+\frac{\xi}{\sqrt{\nu}},y\right)\right]d\xi
	\end{multline}
				where $g(x',y'):=\chi(x')G(x',y').$ 
				Choosing $r,d$ sufficiently small such that\\ $C_1=\sup_{\abs{x}\leq r,\abs{y}\leq d} \abs{D_x\phi(x,y)}\leq \frac{1}{2}$
				the integral converges, as $\nu\to \infty$. to
				\begin{equation}
				\pi^{-\frac{n}{2}}\int_{\Rn}
				\exp\left(-(\xi+iD_x\phi(x,y)\xi)\cdot (\xi+iD_x\phi(x,y)\xi)\right)
				G(x,y)\mbox{det}\left[ D_x Z\left(x,y\right)\right]d\xi
					\end{equation}
					By Lemma \ref{btmatrixlemma}
					the right hand side is $G(x,y).$ This proves Lemma \ref{btsistalemma}.
					\end{proof}
					Lemma \ref{btsistalemma} proves that $H_\nu(Z(x,y))\to h(x,y)$ uniformly on $K$.
					Since an entire function is, on any given compact, the uniform limit of complex polynomials and since $K$
					obviously contains an open neighborhood $\Omega''$ of the origin, this proves
					Theorem \ref{btthm}.
\end{proof}

\section{Preliminaries for Lewy's theorem}\label{lewyapp}

Sometimes one identifies the image of the Levi form with a subspace of the normal space $N_p M$ (which in turn can be 
considered as $J(T_p M/T^c_p M)$).
Denote by $\Gamma_p (\subseteq N_pM)$ the cone with constitutes the convex hull of the image of $\mathcal{L}_p$ (see Boggess \& Polking \cite{bp}).
For two subcones $\Gamma_1,\Gamma_2$ of $\Gamma_p$ we say that $\Gamma_1$ is smaller that $\Gamma_2$
if $\Gamma_1\cap S_p$ (where $S_p$ denotes the unit sphere in $N_p M$) is a compact subset of the (relative) interior of $\Gamma_2 \cap S_p.$
Boggess \& Polking \cite{bp} proved a theorem on holomorphic extension of continuous $CR$ functions
near a point such that the Levi cone at $p$
has nonempty interior. The domain of extension has the shape of the product of an open set with a cone.
The proof involves explicit construction of families of analytic discs by solving a Bishops equation, such that the center of these discs pass each point of an open subset of the given normal cone and simultaneously are attached sufficiently close to $p.$ 
For our purposes it is precisely the existence of such families of discs which is useful.
\begin{lemma}[See Boggess \cite{b4}, p.207]\label{bpthm} 
Let $M\subset \Cn$ be a $C^l,l\geq 4,$ generic embedded $CR$ submanifold and let $p\in M$ be a point such that the Levi cone at $p$
has nonempty interior.
Then for every open neighborhood $\omega$ of $p$
and for each cone $\Gamma <\Gamma_p,$ there is a neighborhood $\omega_{\Gamma}\subset\omega$ and a positive number $\epsilon_{\Gamma}$
such that each point in $\omega_{\Gamma} +\{\Gamma\cap B_{\epsilon_{\Gamma}}\}$ is contained in the image of an analytic disc whose boundary image is contained in $\omega.$
\end{lemma}

\section{Brief notes on the proof}
We have already pointed out the role of the Bishop's equation and
Hilbert transform for the construction of analytic discs.
Define $H(v,\mathcal{D}):e^{i\theta}\mapsto h(u(e^{i\theta}),\mathcal{D}(e^{i\theta})$
Then a solution, $u,$ to the so called Bishop's equation, 
\begin{equation}
T\left(H(u(e^{i\theta}),\mathcal{D}(e^{i\theta}))\right)+x=u(e^{i\theta}), \theta\in [0,2\pi),
\end{equation}
(where $T$ denotes the Hilbert transform, recall that the Hilbert transform gives harmonic conjugates)
yields a modified Hilbert transformation of and since $u$ is the boundary value of a harmonic function, denote this by $U(\zeta),$ then in fact the left hand side, given a normalization of the Hilbert transform, this becomes {\em the} harmonic conjugate, thus  $u+iv=u+iH(u,\mathcal{D}):e^{i\theta}\to \C,$ will be the boundary value of an analytic disc 
\begin{equation}
U(\zeta) +iV(\zeta):\overline{D}\to\C^d .
\end{equation}
Then the analytic disc
\begin{equation}
\left( U(\zeta) +iV(\zeta) ,\mathcal{D} (\zeta) \right): \overline{D}\to \C^{d+m} ,
\end{equation}
will satisfy that Im$(U(\zeta) +iV(\zeta))=ih(\mbox{Re}(U(\zeta) +iV(\zeta)), \mathcal{D} ),\abs{\zeta}=1,$ i.e.\ its boundary is contained in $M.$ 
The following is a commutative diagram for an analytic disc, $(U+iV,\mathcal{D}),$ when it is attached to $M.$
\begin{displaymath}
    \xymatrix{
         &  \C^d \times \C^m  \ar[d]^{(V(0),\mbox{Id})} \\
       \R^d \times  \overline{D} \ar[ur]^{(\mathcal{D}(\zeta),U+iV)} \ar[r]^{(x,\mathcal{D}(\zeta))}     & \R^d \times \C^m 
        } , 
\end{displaymath}
$\partial (U+iV,\mathcal{D})(\zeta)\subset M,$
$\zeta\in \overline{D}\subset \C, x+iy\in \C^d.$
\\
For the case of higher real codimension 
and under certain fruitful circumstances the existence of solutions to Bishop's equation 
belonging to certain smoothness classes is known, see e.g.\ Boggess \& Polking \cite{bp}.
Recall that by definition the Hilbert transform of a continuous real valued function $f(e^{i\theta} )$ on $S^1$
is given by
\begin{equation}
Tu:=\frac{1}{2\pi i}\int_{|\zeta | =1} f(\zeta )\mbox{Im}\left(
\frac{\zeta +z}{\zeta -z}
\right)
\frac{d\zeta }{\zeta } .
\end{equation}
We consider $\R^d$-valued maps, so it is important to recall that we require that $T$ acts componentwise. It is part of the details of the proof to ensure that the Hilbert transform is continuous on appropriate Banach function spaces.
Once this is done the strenuous part of the proof is to
be able to prescribe parts of the centers of cleverly chosen
families of discs attached near a reference point.

\subsubsection{Prescribing part of the center}
Let us assume that $0<\alpha <1,$ and that $M$ locally near $0$ has graph representation $M=\{ (x+iy,w)\in\Cn :y=h(x,w)\},$ with $h\in C^{\infty},h(0)=0, \nabla h(0)=0.$
Introduce the parametrized analytic disc $W :\C^m \times \C \to \C^m,$
$(W(w))(\zeta) =w+\beta(\zeta),\beta-$a holomorphic function.
We 
construct an analytic disc, $U(\zeta)+iV(\zeta),\zeta\in \overline{D},$ for 
$H(u,W)(e^{i\theta}):=h(u(e^{i\theta} ),W(e^{i\theta})),$
$U|_{\partial D}=u(e^{i\theta})=(u(x,W))(e^{i\theta}),\theta\in [0,2\pi],$ 
with $u(x,W)\in C^{\alpha} (S^1 ,\R^d ),$ 
by solving the equation,
\begin{equation}
u(e^{i\theta})=-(T(H(u,W))(e^{i\theta}) +x, 
\theta\in [0,2\pi],
\end{equation}
where $T$ denotes the Hilbert transform, which yields the harmonic conjugate (we are assuming that $T$ act componentwise on $\R^d$-valued functions).
We are using the fact that a smooth $u :S^1\to \R^d$ extends to a harmonic function, $U(\zeta),\zeta\in D\subset\C,$ and our equation helps to find its harmonic conjugate. We automatically get (since $u+iv$ is the boundary of $U+iV$) that the disc is attached to $M$, i.e.\
$V(\zeta ) =H(U(\zeta ),W(\zeta )),\abs{\zeta}=1 .$
Note that $V(\zeta)$ is something we get after applying a theorem on harmonic extension from sufficiently smooth boundary values (in fact we are solving a Dirichlet problem), so a priori we do not have explicit description of $V(\zeta)$ on the interior of the unit disc.  
Now the center of the analytic disc becomes:
\begin{equation}
(U+iV,W)(\zeta=0)=(x+i(V(x,w))(\zeta=0),w)=
(x+i(h(u,W))(\zeta=0),w).
\end{equation}
We know that $V(\zeta)$ by construction is harmonic for $\abs{\zeta} <1,$ and continuous to the boundary, thus by the mean value property for harmonic functions we have (since the boundary belongs to $M,$
\begin{equation}
V(\zeta =0)=
\frac{1}{2\pi} \int_{0}^{2\pi} h((U(x,w,t))(e^{i\theta} ),(W(w,t))(e^{i\theta} )d\theta .
\end{equation}
There exists a case, described in detail in e.g.\ Boggess \cite{b4}, where $\beta(\zeta ,t):=\left(\sum_{l=1}^{N} t_l a_l\cdot\zeta^l \right),$
$a_l\in\C^m ,N\in N ,N\geq 1.$
where it is shown that Taylor expansion with respect to $t,$ near $t=0,$ gives: \\
\\
(i) $W(t=0,w)(\zeta)=w\Rightarrow u(t=0,x,w)(\zeta)=x\Rightarrow 
V(x,w,\zeta =0) (t=0)=h(x,w),$ the constant term in the Taylor expansion. 
In the special case that this integral is zero for any $t,$ we get that the center of the constructed analytic disc can be chosen to lie at the origin in $M$ (since $h(0,0)=0$).
\\
\\
(ii)
$\frac{\partial V(x,w,\zeta=0)}{\partial t_j }(t=0),1\leq j\leq N,$ is $0,$ independent of the parameters of $\beta(\zeta),$ so there are no first order terms in the Taylor expansion in $t$ near $t=0:$\\
For the first order terms we have (here $\frac{\partial h }{\partial w } ,\frac{\partial h }{\partial x },$ in general denote matrices
and for codimension one $x$ will be one dimensional so $\frac{\partial h }{\partial x }$ is a scalar function),
\begin{multline}
\frac{\partial V }{\partial t_j } (t=0,x,w)(\zeta =0)=\frac{1}{2\pi} \int_{0}^{2\pi}
\frac{\partial h }{\partial x }(x,w)\cdot \frac{\partial u }{\partial t_j }(e^{i\theta} ) d\theta \vert_{t=0} 
+ \\
\frac{1}{2\pi} \int_{0}^{2\pi}
2\mbox{Re}\left\{
\frac{\partial h }{\partial w }(x,w)\cdot \frac{\partial W }{\partial t_j }(e^{i\theta} ) \right\} d\theta \vert_{t=0}  .
\end{multline}
For the first term $\frac{\partial u }{\partial t_j }$ is harmonic in $D$ because $U(\zeta)$ is, and $u(\zeta=0)=x$ by construction (this is the prescription allowed by the Bishop equation), thus integrating in $\theta$ (i.e.\ independently of $x,w$) a constant times $\frac{\partial u }{\partial t_j }$ gives zero, so the first term is zero. 
Furthermore, 
\begin{equation}
\left(\frac{\partial h}{\partial w_1} ,\ldots , \frac{\partial h}{\partial w_m}\right)\cdot
\left(\frac{\partial W_1}{\partial t_j} ,\ldots , \frac{\partial W_m}{\partial t_j}\right)^T 
=\sum_{k=1} \frac{\partial h}{\partial w_k} \frac{\partial W_k}{\partial t_j} ,
\end{equation}
where $^T$ denotes transpose.
Now for $\beta$ as above, the second term also vanishes because
$\frac{\partial W }{\partial t_j }(\zeta=0)=0,$ and (recalling that $\nabla h(0,0)=0$ so the product is zero) that leaves no first order term.
\\
\\
(iii) \begin{multline}
\frac{t_j t_k}{2}\left(\frac{\partial^2 V(x,w,\zeta=0)}{\partial t_j \partial t_k}(t=0)\right) =\\
\frac{t_j t_k}{2\pi}\int_{0}^{2\pi} \sum_{r,s=1}^m \left(
\frac{\partial^2 h(x,w,\zeta=0)}{\partial w_r \partial \bar{w}_s}\right)
\left(\frac{\partial W_r (e^{i\theta})}{\partial t_j }\right)
\left( \frac{\partial \overline{W}_s (e^{i\theta})}{\partial t_k } \right) d\theta +
\eta_{jk}(t,x,w),
\end{multline}
$\abs{\eta_{jk}(t,x,w)}\leq x\abs{t}^2 .$
This construction can be used (under the assumption that the Levi cone is nonempty) to put the centers in (a subcone of) the convex hull of the Levi cone.\\
The control of the size of the boundary of the parametrized discs are very important.
We have by,
Theorem 1, p.220, Boggess \cite{b4}, that given 
$x\in \R^d,$ $\exists \delta$ such that for $\abs{x}<\delta ,W\in C^{\alpha} (S^1 ,\C^m ),\norm{W}_{\alpha} <\delta ,$  
which states that we can
construct uniquely an analytic disc, $U(\zeta)+iV(\zeta),\zeta\in \overline{D},$ for 
$U(\zeta =0)=x_0,$ and $U|_{\partial D}=u(e^{i\theta})=(u(x,W))(e^{i\theta}),\theta\in [0,2\pi],$ 
with $u(x,W)\in C^{\alpha} (S^1 ,\R^d ),u(W=0,x=0)=0,$ by solving the equation,
\begin{equation}
u(t)=-(T(h(u,W))(e^{i\theta}) +x, \theta\in [0,2\pi],
\end{equation}
and finally the boundary of the analytic disc $A(\zeta):=(U(\zeta)+iV(\zeta),W(\zeta))$ belongs to $M,$ such that $u:C^{\alpha} (S^1 ,\C^m )\times \R^d\to S^1 ,$ depends in a $C^1$ fashion on $x,W.$ The construction forces $V|_{\partial D} =:v(\zeta) =h(u(\zeta),W(\zeta)),\abs{\zeta}=1.$
Note that only the real part of the analytic disc $U+iV$ is prescribed.%

\chapter[Some technical lemmas]{Some technical lemmas for the proof of a theorem of Vitushkin}\label{vitushkinapp}
In this section we collect some definitions and proofs of some technical lemmas used in the proof of
Theorem \ref{vituskin1}. The proof follows Zalcman \cite{zalcmanbok}. 
\begin{definition}
Suppose $S\subset\C$ is a bounded set such that $\alpha(S)>0$. For $f\in C(S,m)$, 
 and $z_0\in \C$
we have
\begin{equation}
f(z)=\frac{a_1}{z-z_0}+\frac{a_2}{(z-z_0)^2}+\cdots
\end{equation}
Set
\begin{equation}
\beta(S,z_0,f):=\frac{a_2}{\alpha(S)}
\end{equation}
Then for a sufficiently large rectifiable contour with winding number $1$ about each point of $S$ we have
\begin{equation}
\beta(S,z_0,f):=\frac{1}{\alpha(S)}\frac{1}{2\pi i}\int_\sigma f(\zeta)(\zeta -z_0)d\zeta
\end{equation}
Define
\begin{equation}
\beta(S,z):=\sup_{f\in C(S,1)} \abs{\beta(S,z,f)}
\end{equation}
and
\begin{equation}
\beta(S):=\inf_{z} \beta(S,z)
\end{equation}
Then $\beta(S,z)$ is a continuous function of $z$ on $\C$ and
$\beta(S,z)\to\infty$ as $z\to \infty.$ Hence there exists a point
$O(S)$ such that $\beta(S,O(s))=\beta(S).$ 
\end{definition}
\begin{lemma}\label{prepreplem}
$f\in A(X)$ and let $\phi$ be a continuously differentiable function
with support in $\{\abs{z-z_0}\leq\delta\}.$ 
Then
\begin{equation}
\abs{\frac{1}{\pi}\int f(\zeta)\partial_{\bar{\zeta}}\phi(\zeta-z_0) d\mu(\zeta)}\leq
4\delta \omega(f,2\delta)\norm{\partial_{\bar{\zeta}}\phi}_\infty
\alpha(\{\abs{z-z_0}\leq \delta\}\cap \C\setminus\mbox{int}(X))
\end{equation}
\end{lemma}
\begin{proof}
Keeping the same notation we consider a continuous extension of $f$ to $\C.$
We have from Proposition \ref{vitushkinopprops} 
that $f_\phi$ is analytic on $(\C\setminus\{\abs{z-z_0}\leq \delta\})\cup\mbox{int}(X)$ and
\begin{equation}
\norm{f_\phi}_\infty 
\leq  4\delta \omega(f,2\delta)
\norm{\partial_{\bar{z}}\phi}_\infty 
\end{equation} 
thus
\begin{equation}
\abs{f_\phi'(\infty)} 
\leq  4\delta \omega(f,2\delta)
\norm{\partial_{\bar{z}}\phi}_\infty \alpha((\C\setminus\{\abs{z-z_0}\leq \delta\})\cup\mbox{int}(X))
\end{equation} 
On the other hand
\begin{equation}
f_\phi'(\infty)=\lim_{z\to \infty} z f_\phi (z)=-\frac{1}{\pi}\int f(\zeta)\partial_{\bar{\zeta}}\phi d\mu(\zeta) 
\end{equation} 
Hence
\begin{equation}
\abs{\frac{1}{\pi}\int f(\zeta)\partial_{\bar{\zeta}}\phi d\mu(\zeta)}\leq
4\delta \omega(f,2\delta)\norm{\partial_{\bar{\zeta}}\phi}_\infty
\alpha(\{\abs{z-z_0}\leq \delta\}\cap (\C\setminus\mbox{int}(X)))
\end{equation}
Repeating the arguments but replacing $\phi$ with
$g(\zeta):=(\zeta -z_0)\phi(\zeta)$ and noting that
$\partial_{\bar{z}} g=(\zeta -z_0)\partial_{\bar{z}} \phi,$
which implies $\norm{\partial_{\bar{z}} g}_\infty\leq \delta \norm{\partial_{\bar{z}} \phi}_\infty,$
completes the proof.
\end{proof}

\begin{lemma}\label{zalcman108}
Let $S\subset \C$ and suppose $\{ S_j\}$ is a 
collection of subsets of $S$ such that each disc of radius $\alpha(S)$
has nonempty intersection with at most $p$ of the sets 
$S_j$. Then
\begin{equation}\label{fjockmokk}
\sum_{j=1}^\infty \alpha(S_j) <400p\alpha(S)
\end{equation}
Furthermore, if $f_j\in C(S_j,1)$ then 
\begin{equation}
\max \sum_{j=1}^\infty \abs{f_j(z)}<200p
\end{equation}
\end{lemma}
\begin{proof}
We can assume $\alpha(S)>0.$ Fix $z_0$ and relabel the
$S_j$ to obtain a sequence with double indices $\{S_{n,k}\}$ such that
\begin{equation}
n\alpha(S)\leq \mbox{dist}(z_0,S_{n,k}) <(n+1)\alpha(S)
\end{equation}
and there are at most $p$ sets $S_{0,k}.$ Since 
$\{n\alpha(S)\leq \abs{z-z_0}\leq (n+1)\alpha(S)\}$
can be covered by $20$ discs of radius $\alpha(S),$
there are, for fixed $n,$ at most $20np$ sets $S_{nk}$.
Let $f_{n,k}\in C(S_{n,k},1).$ 
We claim that for $n\geq 1,$ if $f_{n,k}\in C^0(\{ \abs{z-z_0}\leq R\},1)$ for some $R>0$ then
\begin{equation}\label{seconclaimef0}
\abs{f_{n,k}(z_0)}\leq \gamma(\{ \abs{z-z_0}\leq R\})/(t-R), \quad \abs{z-z_0} \geq t >R
\end{equation}
and if
\begin{equation}
f(z)=\frac{a_1}{z-z_0} +\frac{a_2}{(z-z_0)^2} +\cdots
\end{equation}
we have
\begin{equation}\label{seconclaimef}
\abs{a_n}\leq \exp(1)\alpha(K) R^{n-1}n
\end{equation}
To see this fix $z_1$ such that $\abs{z_1-z_0}\geq t >R$ and set
\begin{equation}
\tilde{f}(z)=\frac{t-R}{z-z_1} \frac{f(z)-f(z_1)}{1-\overline{f(z_1)} f(z)} =-\frac{t-R}{z-z_0}f(z_1)+\cdots
\end{equation}
Then $\abs{\tilde{f}(z)}<1$ for $z\in\in \Omega(K)\cap \{\abs{z-z_0}<R\}$ where $\Omega(K)$
denotes the unbounded component of $\C\setminus K.$ By the maximum principle $\tilde{f}$ is an admissable function for
$K$ so that
$\abs{\tilde{f}'(\infty)}=\lim_{z\to \infty} \abs{z\tilde{f}(z)}=(t-R)\abs{f(z_1)}\leq \gamma(K).$
For Eqn.(\ref{seconclaimef}) use for $n=2,3\ldots,$
\begin{equation}
a_n =\frac{1}{2\pi i} \int_{\abs{z-z_0}=t>R} f(z)(z-z_0)^{n-1} dz
\end{equation}
which gives $\abs{a_n}\leq \gamma(K) t^n(t-R)^{-1}$ where the maximum of the right hand side is obtained for $t=nR/(n-1)<e.$
This proves the claim. 
By this claim we obtain
\begin{equation}
\abs{f_{n,k}(z_0)}\leq \alpha(S_{n,k})/\mbox{dist}(z_0,S_{n,k})\leq \alpha(S_{n,k})/n\alpha(S)
\end{equation}
Next using the inequality (which we state without proof) for $0\leq c_{n,k}\leq 1,$ $k=1,2,\ldots,np,$ $n\in \Z_+$ 
\begin{equation}
\left(\sum_{n,k} \frac{c_{n,k}}{n}\right)^2 \leq 4p \sum_{n,k} c_{n,k}
\end{equation}
we obtain
\begin{equation}
\sum_{j=1}^\infty \abs{f_j(z_0)}\leq p+ \sum_{n,k} \frac{\alpha(S_{n,k})}{n\alpha(S)} <p+9\sqrt{pm}
\end{equation}
Pick $\phi_j\in C(S_j,1)$ such that $\phi_j'(\infty)=\frac{1}{2}\alpha(S_j)$ and set
$\phi :=\sum \phi_j.$ Then $\norm{\phi}_\infty \leq p+9\sqrt{pm}$ thus
\begin{equation}
\phi_j'(\infty)=\frac{1}{2}\sum_{j=1}^\infty \alpha(S_j)\leq \alpha(S)(p+9\sqrt{pm})
\end{equation}
which implies $m\leq 2p+18\sqrt{pm}<400p,$ proving Eqn.(\ref{fjockmokk}).
This in turn implies that 
\begin{equation}
\max \sum_{j=1}^\infty \abs{f_j(z)}\leq p+9\sqrt{400p^2}
\end{equation}
This completes the proof of Lemma \ref{zalcman108}.
\end{proof}

Let $X\subset\C$ be a compact set. 
Suppose there exists constants $m\geq 1,$ $r\geq 1,$ such that
for all $z$ and all $\delta>0$
\begin{equation}\label{starfds}
\alpha((\C\setminus\mbox{int}(X))\cap \{\abs{z-z_0}\leq \delta\})\leq
m \alpha((\C\setminus X)\cap \{\abs{z-z_0}\leq r\delta\})
\end{equation}
Let $\{\phi_{k,n}\}$ be the sequence of partitions of unity in Proposition \ref{partitionofunutylemmprop}
and set $X_{k,n}:=(\C\setminus X)\cap \{\abs{z-z_{k,n}}\leq r\delta_n\}.$

\begin{lemma}\label{elvalemfem}
Let $f\in A(X).$ If Eqn.(\ref{starfds}) holds for all $z$ and all $\delta >0,$ then
\begin{equation}
\abs{\frac{1}{\pi}\int f(\zeta)\partial_{\bar{\zeta}}\phi_{k,n}(\zeta -O(X_{k,n}))d\mu(\zeta)}\leq
m_1 \omega(f,2\delta_n)\alpha(X_{k,n})\beta(X_{k,n})
\end{equation}
where $m_1$ is a constant depending only on $m$ and $r.$
\end{lemma}
\begin{proof}
Use the simplified notation $\phi=\phi_{k,n},\delta=\delta_n,$
$z_0=z_{k,n},$ $\mathcal{O}=O(X_{k,n}).$ Choose $\beta$ such that
\begin{equation}
\beta\leq \beta(X_{n,k})\leq 2\beta, \mbox{ if } \quad\beta(X_{k,n})<\delta
\end{equation}
\begin{equation}\label{betabetabeta0}
\beta\leq \delta, \quad \mbox{ if }\beta(X_{k,n})\geq \delta
\end{equation}
and cover $\{\abs{z-z_0}\leq \delta\}$ by 
a finite number of discs $\{\abs{z-t_i}\leq \beta\}$ such that (using our previous observations)
each disc of radius $\beta$ intersects at most $25r^2$ of the discs 
$\{\abs{z-t_i}\leq r\beta\}$ and such that 
$\{\abs{z-t_i}\leq r\beta\}\subset \{\abs{z-z_0}\leq \beta\}$ for all $i.$
By Lemma \ref{zalcman108}, it follows that
\begin{equation}
\alpha(X_{k,n})\leq (2+r)\beta
\end{equation}
and that each disc of radius $\alpha(X_{n,k})$
has nonempty intersection with at most $p$ of the discs $\{\abs{z-t_i}\leq r\beta\}$  
where $p$ depends only on $r.$
Choose $g_j\in C^\infty_c(\{\abs{z-t_i}\leq r\beta\})$
such that $0\leq g_i\leq 1,$
$\sum_i g_i(z)=1,$ on a neighborhood of 
$\{\abs{z-z_0}\leq \delta\},$ $\norm{\partial_{\bar{z}} g_i}_\infty \leq 20/\beta .$
Then 
\begin{multline}
\frac{1}{\pi}\int f(\zeta)\partial_{\bar{\zeta}}\phi(\zeta -\mathcal{O})
d\mu(\zeta) =\\
\sum_i \frac{1}{\pi}\int f(\zeta)\partial_{\bar{\zeta}}(\phi g_i)(\zeta -t_i)
d\mu(\zeta) + 
\sum_i \frac{1}{\pi}\int f(\zeta)\partial_{\bar{\zeta}}(\phi g_i)(t_i -\mathcal{O})
d\mu(\zeta)
\end{multline}
Since $\beta \leq \delta,$ we have by the proof of Lemma \ref{prepreplem}
\begin{multline}
\abs{ \frac{1}{\pi}\int f(\zeta)\partial_{\bar{\zeta}}(\phi g_i)\mu(\zeta)} 
\leq\\ 
4\beta \omega(f,2\beta)(20/\delta +20/\beta)\alpha(\{\abs{z-t_i}\leq \beta\}\cap(\C\setminus
\mbox{int}(X)))\leq \\
160 \omega(f,2\beta)\alpha(\{\abs{z-t_i}\leq \beta\}\cap(\C\setminus \mbox{int}(X)))
\end{multline}
By Lemma \ref{prepreplem}
\begin{equation}
\abs{ \frac{1}{\pi}\int f(\zeta)\partial_{\bar{\zeta}}(\phi g_i)(\zeta- t_i)\mu(\zeta)}
\leq 
160\beta \omega(f,2\beta)\alpha(\{\abs{z-t_i}\leq \beta\}\cap(\C\setminus\mbox{int}(X)))
\end{equation}
The last three equations imply
\begin{multline}\label{kokolo}
\abs{ \frac{1}{\pi}\int f(\zeta)\partial_{\bar{\zeta}}(\phi)(\zeta- \mathcal{O})\mu(\zeta)}
\leq \\
160\beta \omega(f,2\beta)\sum_i \alpha(\{\abs{z-t_i}\leq \beta\}\cap(\C\setminus
\mbox{int}(X))) +\\
160\beta \omega(f,2\beta)\sum_i \abs{t_i-\mathcal{O}}
\alpha(\{\abs{z-t_i}\leq \beta\}\cap(\C\setminus\mbox{int}(X)))
\end{multline}
Since each disc of radius $\alpha(X_{n,k})$ has nonempty intersection wit at most $p$
of the $\{\abs{z-t_i}\leq \beta\}\subset \{\abs{z-z_0}\leq r\delta\}.$
Eqn.(\ref{starfds}) together with Lemma \ref{zalcman108}, implies
\begin{multline}
\sum_i \alpha(\{\abs{z-t_i}\leq \beta\}\cap(\C\setminus\mbox{int}(X))) \leq \\
m \sum_i \alpha(\{\abs{z-t_i}\leq r\beta\}\cap(\C\setminus X)) \leq m 400p\alpha(X_{n,k})
\end{multline}
By the conditions for the choice of $\beta$ we obtain, for a constant $C_1>0,$ the following bound for 
the first term in the right hand side of Eqn.(\ref{kokolo}) 
\begin{multline}\label{sexkokolo}
160\beta \omega(f,2\beta)\sum_i \alpha(\{\abs{z-t_i}\leq \beta\}\cap(\C\setminus\mbox{int}(X))) \leq \\
C_1 mp\beta\omega(f,2\beta)\alpha(X_{k,n})\leq
C_1 mp\omega(f,2\beta)\alpha(X_{k,n})\beta(X_{k,n})
\end{multline}
For an appropriate bound of the sum in 
Eqn.(\ref{kokolo}) we proceed as follows.
For $t_1\neq \mathcal{O},$ pick $\psi_i\in C^0((\C\setminus X)\cap \{\abs{z-t_i}\leq r\beta\},1)$ such that
\begin{equation}
\psi_i'(\infty)=\frac{1}{2} \frac{\abs{t_i-\mathcal{O}}}{t_i-\mathcal{O}} \alpha((\C\setminus X)\cap \{\abs{z-t_i}\leq r\beta\})
\end{equation}
and set $\psi:=\sum_i \psi_i,$ $\sigma_i:=\{\abs{z-t_i}= r\beta\},$ $\sigma:=\{\abs{z-z_0}= r\beta\}$
By Lemma \ref{zalcman108} $\norm{\psi}_\infty \leq 200p.$ 
$\beta(X_{n,k})=\beta(X_{n,k},\mathcal{O})\geq \abs{\beta(X_{n,k},\mathcal{O},\psi/200p)}$ yields
\begin{multline}
\beta(X_{n,k})\geq (2\pi 200p \alpha(X_{n,k}))^{-1} 
\abs{\int_\sigma \psi(z)(z-\mathcal{O})dz}\geq\\
(400\pi p\alpha(X_{n,k}))^{-1}\abs{\sum_i (t_i -\mathcal{O})\int_{\sigma_i}\psi_i(z)dz+\sum_i \int_{\sigma_i}\psi_i(z)(z-t_i)dz}
\end{multline}
By the expression for $\psi'(\infty)$
we have
\begin{multline}
\sum_i (t_i -\mathcal{O})\int_{\sigma_i}\psi_i(z)dz =
\frac{1}{2}\sum_i \abs{t_i-\mathcal{O}}\alpha((\C\setminus X)\cap \{\abs{z-t_i}\leq r\beta\})
\end{multline}
By Eqn.(\ref{seconclaimef}) and Eqn.(\ref{seconclaimef0}) we have
\begin{equation}
\abs{\int_{\sigma_i} \psi_i(z)(z-t_i)dz} \leq 
2\pi 2\exp(1) r\beta \alpha((\C\setminus X)\cap \{\abs{z-t_i}\leq r\beta\}))
\end{equation}
The last three equations render
\begin{multline}
\sum_i \abs{t_i -\mathcal{O}}\alpha((\C\setminus X)\cap \{\abs{z-t_i}\leq r\beta\})\leq\\
800\pi p \alpha(X_{n,k})\beta(X_{n,k})+8\pi\exp(1) r\beta \sum_i \alpha((\C\setminus X)\cap \{\abs{z-t_i}\leq r\beta\})
\end{multline}
Again applying Lemma \ref{zalcman108} yields for a constant $C_2>0$
\begin{equation}\label{tolvvvvnn}
\sum_i \abs{t_i -\mathcal{O}}\alpha(
(\C\setminus X)\cap \{\abs{z-t_i}\leq r\beta\}))\leq
C_2 r p \alpha(X_{n,k})\beta(X_{n,k})
\end{equation}
Eqn.(\ref{starfds}) together with Eqn.(\ref{tolvvvvnn}) gives for a constant $C_3>0$
\begin{multline}\label{trettvvvnn}
160\omega(f,2\beta) \sum_i \abs{t_i -\mathcal{O}}\alpha((\C\setminus X)\cap \{\abs{z-t_i}\leq r\beta\})\leq\\
160  \omega(f,2\beta) m \sum_i \abs{t_i -\mathcal{O}}\alpha((\C\setminus X)\cap \{\abs{z-t_i}\leq r\beta\})\leq \\
C_3 rpm \omega(f,2\beta) \alpha(X_{n,k})\beta(X_{n,k})
\end{multline}
Eqn.(\ref{trettvvvnn}) together with Eqn.(\ref{sexkokolo}) inserted into Eqn.(\ref{sexkokolo}) renders
\begin{multline}
\abs{\frac{1}{\pi} \int f(\zeta)\partial_{\bar{\zeta}} \phi(\zeta -\mathcal{O})d\mu(\zeta)}\leq \\
C_1 mp \omega(f,2\beta)\alpha(X_{n,k})\beta(X_{n,k})
+C_3 rpm \omega(f,2\beta) \alpha(X_{n,k})\beta(X_{n,k})
\end{multline}
By Eqn.(\ref{betabetabeta0}) we have $\omega(f,2\beta)\leq \omega(f,2\delta)$ thus setting $m_1:=mp(C_1+C_3)$  
completes the proof of Lemma \ref{elvalemfem}.
\end{proof}

We shall also need the following lemma.
\begin{lemma}\label{extraskrapppp}
Let $S\subset\C $ be a bounded set such that $\alpha(S)>0.$
Let $\alpha,\beta\in C$ such that $\abs{\alpha}\leq \alpha(S)$ and $\abs{\beta}\leq \beta(S).$ Then there exists $g\in C(S,20)$ 
such that $g'(\infty)=\alpha$ and $\beta(S,O(S),g)=\beta.$
\end{lemma}
\begin{proof}
By the definition of $\beta(S,z,f)$
\begin{equation}
\beta(S,z,f)=\beta(S,t,f)+ \frac{f'(\infty)}{\alpha(S)}(t-z)
\end{equation}
Let $\phi\in C(S,2)$ such that $\phi'(\infty)=\alpha(S)$, and $z_0:=O(S)+\beta(S,O(S),\phi).$ Then
\begin{equation}
\beta(S,z_0,f)=\beta(S,O(S),\phi)+ \frac{\phi'(\infty)}{\alpha(S)}(O(S)-z_0)=0
\end{equation} 
Let $f_1\in C(S,2)$ such that $\beta(S,z_0,f_1)=\beta(S)$ and $\epsilon$ such that
$f'(\infty)+\epsilon\phi'(\infty)=0$ and $f_2:=f_1 +\epsilon \phi.$ Now $f'_2(\infty)=0$
implies
\begin{equation}
\beta(S,O(S),f_2)=\beta(S,z_0,f_2)=\beta(S,z_0,f_1)+
\epsilon \beta(S,z_0,\phi)=\beta(S,z_0,f_1)=\beta(S)
\end{equation} 
Since $\beta(S,O(S),\phi)=z_0-O(S)$ we have $\abs{O(S)-z_0}\leq 2\beta(S)$ and since $f_1\in C(S,2)$ for $\abs{\epsilon}\leq 2$
we have $f_2\in C(S,6)$. Hence the following function belongs to $C(S,20)$ and satisfies the conclusion of the Lemma
\begin{equation}
f=\frac{\alpha}{\alpha(S)}\phi +\frac{\alpha(O(S)-z_0)}{\alpha(S)\beta(S)}f_2+\frac{\beta}{\beta(S)}f_2
\end{equation} 
This proves Lemma \ref{extraskrapppp}.
\end{proof}

\begin{proposition}[Coefficient estimate]\label{coeffestimateprop}
Let $f\in A(X).$ If Eqn.(\ref{starfds}) holds for all $z$ and all $\delta >0.$ 
Let $\{\phi_{k,n}\}$ be the sequence of partitions of unity in Proposition \ref{partitionofunutylemmprop}
and set $X_{k,n}:=(\C\setminus X)\cap \{\abs{z-z_{k,n}}\leq r\delta_n\}.$
Set
\begin{equation}
g(z)=f(z)\phi_{k,n}+\frac{1}{\pi} \int f(\zeta)\partial_{\bar{\zeta}}\phi_{k,n} \frac{1}{\zeta -z}d\mu(\zeta)
=\sum_{s=1}^\infty \frac{a_s}{(z-O(X_{k,n}))^s}
\end{equation}
Then $\abs{a_2}\leq m_1 \omega(f,2\delta_n) \alpha(X_{n,k})\beta(X_{n,k}).$
\end{proposition}
\begin{proof}
For $\sigma:=\{\abs{z-z_{k,n}}=2\delta\}$ we have
\begin{multline}
a_2=\frac{1}{2\pi i}\int_\sigma g(z)(z-\mathcal{O} )dz =\\
\frac{1}{2\pi i}\int_\sigma \left( \frac{1}{\pi} \int f(\zeta)\partial_{\bar{z}} \phi_{k,n} \frac{1}{\zeta -z}(z-\mathcal{O})d\mu(\zeta)\right) dz=\\
\frac{1}{\pi} \int f(\zeta)\partial_{\bar{z}} \phi_{k,n}\left(\frac{1}{2\pi i} \int_\sigma \frac{z-\mathcal{O}}{\zeta -z} dz\right)d\mu(\zeta)=
\frac{1}{\pi} \int f(\zeta)\partial_{\bar{z}} \phi_{k,n}(\mathcal{O}-\zeta)d\mu(\zeta)
\end{multline}
By Lemma \ref{elvalemfem}, this completes the proof.
\end{proof}

\begin{theorem}\label{bishopthm}
Let $X\subset \C$ be compact and $f\in C^0(\hat{\C}).$
Suppose that for every $z\in X$ there exists a closed neighborhood $K_z:=\{\zeta:\abs{\zeta-z}\leq \delta_z\}$
such that $f|_{X\cap K_z}\in R(X\cap K_z).$ Then $f\in R(X).$
\end{theorem}
\begin{proof}
By compactness there exists finitely many $z_1,\ldots,z_n$ such that
$X\subset \cup_{j=1}^n K_{z_j}.$ Choose $\phi_j\in C^\infty(\hat{C})$ such that
$0\leq \phi_j(z)\leq 1,$ $\phi_j=0$ off $K_{z_j}$ and
$\phi(z)=\sum_{j=1}^n \phi_j(z)=1,$ for $z$ in a closed neighborhood, $V$, of $X,$
such that $V\subset \cup_{j=1}^n K_{z_j}.$ For  
\begin{equation}
f_j(z):=f(z)\phi_j(z) +\frac{1}{\pi} \int f(\zeta)\partial_{\bar{\zeta}}\phi_j \frac{1}{\zeta -z}d\mu(\zeta)
\end{equation}
we have $\sum_{j=1}^n f_j =f+\Phi$ for 
\begin{equation}
\Phi(z)=\frac{1}{\pi}\int f(\zeta)\partial_{\bar{\zeta}}\phi \frac{1}{\zeta-z}d\mu(\zeta)
\end{equation}
Define 
\begin{equation}
C=\max_{1\leq j\leq n} \sup_z \frac{1}{\pi}\int \abs{\partial_{\bar{z}} \phi_j} \frac{1}{\zeta -z}d\mu(\zeta)
\end{equation}
and for $\epsilon >0,$ choose rational functions $h_j$, $1\leq j\leq n,$ with poles off $X$
such that $\norm{f-h_j}_{X\cap K_{z_j}} <\epsilon /8nC.$ Choose a closed neighborhood $U_j$
of $X\cap K_{z_j}$ such that $\norm{f-h_j}_{U_j} <\epsilon /4nC$ and $h_j$ is analytic on $U_j.$
We modify $h_j$ (keeping the same notation) off $U_j$ such that 
$\norm{f-h_j}_\C <\epsilon/4nC$ and set
\begin{equation}
g_j(z)=h_j(z)\phi_j(z) +\frac{1}{\pi}\int h_j(\zeta)\partial_{\bar{\zeta}}\phi_j\frac{1}{\zeta -z}d\mu(\zeta)
\end{equation}
By Proposition \ref{vitushkinopprops}
$g_j$ is analytic on $U_j \cup(\C\setminus K_{z_j})$
which is a neighborhood of $X$. Also
\begin{multline}
\norm{g_j-f_j}_\C \leq \norm{h_j-f}_\C +\norm{\frac{1}{\pi}\int_{K_{z_j}}
(h_j(\zeta)-f(\zeta)\partial_{\bar{\zeta}} \phi_j \frac{1}{\zeta -z}d\mu(\zeta)}_\C\\
\leq \norm{h_j-f}_\C +\norm{h_j -f}_{K_{z_j}} C\leq \frac{\epsilon}{2n}
\end{multline}
Set $g:=\sum g_j.$ Then 
\begin{equation}
\norm{g-\sum f_j}_V \leq \sum_{j=1}^n \norm{g_j-f_j}_V <\frac{\epsilon}{2}
\end{equation}
Since $\partial_{\bar{\eta}} \phi =0$ on $V$, $\Phi(z)$ is analytic on $V.$
By Proposition \ref{rungesobs} there exists $r\in R(X)$ such that $\norm{\Phi -r}_X <\epsilon/4.$ Since 
each $g_j$ is analytic on a neighborhood of $X$ so is $g$ so there exists $q\in R(X)$ such that
$\norm{g-q}_X <\epsilon/4.$ Hence $\norm{f-(r+q)}<\epsilon.$ This completes the proof.
\end{proof}
We immediately have the following corollary.
\begin{corollary}\label{bishopcorr}
Let $X$ be a compact and $U$ a bounded domain. Suppose that $f\in R(X)$ and $f$ analytic
on a neighborhood of $\C\setminus U.$ Then $f\in R(X\cup (\C\setminus U)).$
\end{corollary}

\chapter[Sokhotsky-Plemelj formula and Weierstrass theorems]{The Sokhotsky-Plemelj formula and the Weierstrass preparation theorem and division theorem}\label{plemeljsec}
In this section we shall give the proof of two well-known results which are used repeatedly in the proofs of the book,
namely the Sokhotsky-Plemelj formula and the Weierstrass preparation theorem.
We begin with the Sokhotzki-Plemelj\index{Sokhotzki-Plemelj formula} jump formula, see e.g.\ Muskhelishvili \cite{muskesh}, 17, Eqn.(17.2),  Markushevic \cite{markushevich} or Gakhov \cite{gakhov}. 
Recall that a fucntion $f$ on a continuous arc $\gamma$ is called {\em H\"older continuous}\index{H\"older condition (H\"older continuity)} 
(or to satisfy the {\em H\"older condition} on $\gamma$) if there exists a constant $C>0$ and a constant $\alpha\in (0,1),$
(we also allow $\alpha=1$ which is Lipschitz continuity)
such that for all $x_1,x_2\in \gamma$ we have
\begin{equation}
\abs{f(x_2)-f(x_1)}\leq C\abs{x_2-x_1}^\alpha
\end{equation}
Sometimes $\alpha$ is called the H\"older index. 
The set of such functions 
is denoted $C^{0,\alpha}(\gamma).$
Consider for a H\"older continuous function $\varphi\in C^{0,\alpha}(\gamma)$ on a rectifiable curve $\gamma$
the Cauchy integral
\begin{equation}
f(z):=\frac{1}{2\pi i}\int_\gamma \frac{\varphi(t)dt}{t-z}
\end{equation}
This is single-valued for all $z\notin \gamma$ and for the case that $\gamma$ is bounded we clearly have 
$f(z)=O\left(\frac{1}{\abs{z}}\right)$ as $\abs{z}\to\infty.$ 
The function $\varphi$ is sometimes called the {\em density function}.\index{Density function of a Cauchy type integral}
If $\gamma$ is closed, enclosing the bounded domain $D^+$, and denoting by $D^-$ the exterior domain $\C\setminus (\gamma\cup D^+)$, 
and $\varphi$ has holomorphic extension to $U^+$ then for each $z_0\in \gamma$ we have 
$f_+(z_0):=\lim_{U^+\ni z\to z_0} f(z)=\varphi(z_0),$ whereas $f_-(z_0):=\lim_{D^-\ni z\to z_0} f(z)=0.$ 
We assume the orientation of $\gamma$ is positive in the sense that
an observer moving along $\gamma$ has $D^+$ to the left. In particular, $\varphi(z_0)=f_+(z_0)-f_-(z_0).$
In the more general case, with continuous $\varphi$, we may differentiate under the integral sign so that
$\partial_z^j f(z_0)=\int_\gamma \frac{j!\varphi(t)dt}{(t-z_0)^{j+1}}$ exists for all $j\in \N$ and $z_0\notin \gamma.$
For $z_0\in \gamma$ we must consider the Cauchy principal value
\begin{equation}
\mbox{p.v.}\int_\gamma \frac{\varphi(t)dt}{t-z_0}=\lim_{\epsilon\to 0} \int_{\gamma_\epsilon} \frac{\varphi(t)dt}{t-z}
\end{equation}
where $\gamma_\epsilon:=\gamma\setminus \{\abs{z-z_0}<\epsilon\}.$
\begin{proposition}\label{pvexistsprop}
	Let $\gamma$ be a (non-self-intersecting) rectifiable arc with end points $z_1,z_2$. Then the principal value exists at $z_0\in \gamma$ for any density $\varphi$ which is H\"older continuous near $z_0$ on $\gamma.$
\end{proposition}
\begin{proof}
	We have
	\begin{equation}\label{hardyfyratva}
	\int_{\gamma} \frac{\varphi(t)dt}{t-z_0}
	=\int_{\gamma_\epsilon} \frac{\varphi(t)dt}{t-z_0} 
	+\int_{\{z\in \gamma:\abs{z-z_0}<\epsilon\}} \frac{\varphi(t)dt}{t-z_0}
	=:I_1+
	I_2
	\end{equation}
	For a given $\epsilon>0$ denote by $t_1,t_2$ the end points of the arc segment $\{z\in \gamma:\abs{z-z_0}<\epsilon\}$
	and denote by $z_1,z_2$ the end points of $\gamma$,
	where
	\begin{equation}
	I_2=\ln(z_1-z_0)-\ln(z_2-z_0)-(\ln(t_1-z_0)-\ln(t_2-z_0))
	\end{equation}
	We may assume the condition $\abs{t_1-z_0}=\abs{t_2-z_0}$ is satisfied as taking the limit $\epsilon\to 0$ in order to get
	\begin{equation}
	\lim_{\epsilon\to 0}\left( \ln(t_1-z_0)-\ln(t_2-z_0)\right)=-\pi i
	\end{equation}
	hence
	\begin{equation}
	\lim_{\epsilon\to 0} I_2= \pi i +\ln\frac{z_1-z_0}{z_2-z_0}
	\end{equation}
	For $z\in \gamma$ near $z_0$ we have by the H\"older condition
	a constant $C>0$ and a constant $\alpha\in (0,1),$ 
	such that 
	\begin{equation}
	\abs{f(z)-f(z_0)}\leq C\abs{z-z_0}^\alpha
	\end{equation}
	thus the limit $\lim_{\epsilon\to 0} I_1$ exists near $z_0$ so 
	by Eqn.(\ref{hardyfyratva}) we have
	\begin{equation}
	f(z)=\mbox{p.v.}\int_\gamma \frac{\psi(t)dt}{t-z_0}=
	\frac{1}{2}\varphi(z_0) +\frac{1}{2\pi i}f(z_0)\ln\frac{z_2-z_0}{z_1 -z_0} +\frac{1}{2\pi i}\int_\gamma \frac{\varphi(t)-\varphi(z_0)}{t-z_0}dt
	\end{equation}
	This completes the proof.
\end{proof}
\begin{remark}\label{gakovremmet}
	The proof of Proposition \ref{pvexistsprop} can be repeated for the case of a closed contour $\gamma$ with $z_1=z_2.$
	Note that the result is however local, thus may be applied to the case when $\gamma$ is an open arc for each point
	$z_0$ that is not an end point of $\gamma$ and that we only require that the H\"older condition is {\em locally}
	present for $\varphi$ near $z_0.$
	Given a bounded open rectifiable arc $\gamma$, we may of course complement it to a closed (Jordan) arc by an arbitrary bounded 
	non-self-intersection
	arc whose
	end points coincide with that of $\gamma.$ This will introduce the possibility of associating an orientation to $\gamma$ (and this 
	choice can be done in two different ways).
\end{remark}
Given a domain $\Omega\subset \C$, a rectifiable arc $\gamma\subset\overline{\Omega}$, sometimes one calls a function $f(z)$ defined on $\Omega,$
{\em sectionally holomorphic}\index{Sectionally holomorphic function} with jump along $\gamma$, 
if $f(z)$ is holomorphic on $\Omega\setminus\gamma$ and continuous
on $\gamma$ (with possible exception of the end points of $\gamma$) when approaching $\gamma$ from either side (for a fixed orientation chosen for $\gamma$).
Note that near the end-points (if they lie in $\Omega$) we have by the
H\"older condition
a constant $C>0$ and a constant $\alpha\in (0,1],$
such that for an end point $z_1\in \Omega$
\begin{equation}
\frac{\abs{f(z)}}{\abs{z-z_1}^\alpha}\leq C
\end{equation}
Set for $z_0\in \gamma$
\begin{equation}\label{psiekv}
\psi(z):=\frac{1}{2\pi i}\int_\gamma \frac{\varphi(t)-\varphi(z_0)}{t-z}
\end{equation}
\begin{proposition}\label{gakovpreprop}
	For each path $\kappa\subset \C\setminus \gamma$,
	the limit $\lim_{\kappa\ni z\to z_0} \psi(z)$ exists and is independent of the choice of path $\kappa$. We shall denote the limit
	by $\psi(z_0).$
\end{proposition}
\begin{proof}
	We have setting $\gamma_\epsilon:=\gamma\setminus \{\abs{z-z_0}<\epsilon\}.$
	\begin{multline}
	\psi(z)-\psi(z_0)=\frac{1}{2\pi i} \int_\gamma (z-z_0)\frac{\varphi(t)-\varphi(z_0)}{(t-z)(t-z_0)}dt
	=\\
	\frac{1}{2\pi i} \int_{\gamma\setminus \gamma_\epsilon} (z-z_0)\frac{\varphi(t)-\varphi(z_0)}{(t-z)(t-z_0)}dt
	+\frac{1}{2\pi i} \int_{\gamma_\epsilon} (z-z_0)\frac{\varphi(t)-\varphi(z_0)}{(t-z)(t-z_0)}dt
	=:I_1+I_2
	\end{multline}
	\\	
	Suppose first that $z$ approaches $z_0$ along a path nontangential 
	to $\gamma.$ For sufficiently small $\delta$ the non-obtuse angle, $\omega$
	at $z_0$ between the line segments $z_0z,$ $z_0t$ has a lower bound, say $\omega_0>0.$ By the Sine theorem applies to the triangle $zz_0t$
	we have
	\begin{equation}\label{gakovplem41}
	\frac{z-z_0}{t-z}=\frac{\sin\beta}{\sin\omega}\leq \frac{1}{\sin\omega_0}=K
	\end{equation}
	for a constant $K>0.$
	By the H\"older condition
	there exists a constant $C>0$ and a constant $\alpha\in (0,1)$ such that
	\begin{equation}\label{gakovplem42}
	\abs{\frac{\varphi(t)-\varphi(z_0)}{t-z_0}}<C\abs{t-z_0}^{\alpha -1}=Cr^{\alpha -1},\quad r:=\abs{t-z_0}
		\end{equation} 	
		Also if $s$ denotes the arc-length for $\gamma$ and $r$ the chord length we have that 
		there exists a constant $m>0$ such that
		\begin{equation}
		\abs{\frac{ds}{dr}}\leq m
		\end{equation}
		which implies
		\begin{equation}\label{gakovplem43}
		\abs{t}=\abs{ds}\leq m\abs{dr}
		\end{equation}
		By Eqn.(\ref{gakovplem41})-Eqn.(\ref{gakovplem43})
		we have for $\delta>0$
		\begin{multline}
		\abs{I_1}\leq \int_{\gamma\cap \{\abs{t-z_0}<\delta\}} \abs{\frac{z-z_0}{t-z}}
		\abs{\frac{\varphi(t)-\varphi(z_0)}{t-z_0}}\abs{dt}<KC m\int_{\gamma\cap \{\abs{t-z_0}<\delta\}} r^{\alpha-1}\abs{dr}\\
		=2KCm\int_0^\delta r^{\alpha-1}dr=\frac{2KCm\delta^\alpha}{\alpha}
		\end{multline}
		For a given $\epsilon>0$ we may find $\delta>0$ such that $\abs{I_1}<\frac{\epsilon}{2}.$ 
		Also $I_2$ is continuous (with respect to $z$) on $\gamma_\delta$ (which does not contain $z_0$) so
		for sufficiently small $\abs{z-z_0}$ we may obtain 
		$\abs{I_2}<\frac{\epsilon}{2}$, and we may assume $\delta$ is chosen such that
		simultaneously, $\abs{I_1}<\frac{\epsilon}{2}$ so that
		\begin{equation}
		\abs{\psi(z)-\psi(z_0)}\leq \abs{I_1}+\abs{I_2}<\epsilon
		\end{equation}
		This takes care of the case when $z$ approaches $z_0$ along a path nontangential 
		to $\gamma.$
		\\
		Next suppose that $z$ approaches $z_0$ along a path, $\kappa,$ tangential 
		to $\gamma.$ 
		For $z,z'\in \gamma$ we have
		\begin{equation}\label{gakovkkk0}
		\abs{\psi(z_0)-\psi(z')}\leq \abs{\psi(z_0)-\psi(z)}+\abs{\psi(z)-\psi(z')}
		\end{equation}
		Since the estimate of $\abs{\psi(z)-\psi(z_)}$ above is independent of $z_0$ the limit of $\psi(z)$ to its limit is uniform
		thus the limit function $\psi(z)|_\gamma$ is continuous as both terms of the right hand side of
		Eqn.(\ref{gakovkkk0}) become arbitrarily small for $\abs{z-z'}$ sufficiently small.
		Let $z$ be a point sufficiently close to $z_0$ and draw a through it a curve $\kappa_1$ that intersects $\gamma$ non-tangentially at a point $z_1$
		sufficiently close to $z_0.$ Since $\gamma$ is rectifiable,
		the curve $\kappa_1$ can be chosen such that the lengths of the chord $zz_1$ is bounded up to a constant by the length of the chord $zz_0$.
		Applying the result for approach along non-tangential paths 
		we conclude that $\abs{\psi(z)-\psi(z_1)}$ and $\abs{\psi(z_0)-\psi(z_1)}$ become arbitrarily small
			hence also
			\begin{equation}
			\abs{\psi(z)-\psi(z_0)}\leq \abs{\psi(z)-\psi(z_1)}+\abs{\psi(z_0)-\psi(z_1)}
				\end{equation}
				becomes arbitrarily small. This proves the result also for approach along tangential paths.
				This proves Proposition \ref{gakovpreprop}.	
				\end{proof}
				As in Remark \ref{gakovremmet} we note that the result is local, $\varphi$ is not required to be H\"older continuous except
				near $z_0.$
				Let $\Gamma$ be a smooth arc in the complex plane. 
				and let $f$ belong to the H\"older class
				$C^\alpha(\Gamma),$ $\alpha \in (0,1).$
				If $\Gamma$ is closed let $D^+$ be the bounded domain with $\partial D^+=\Gamma$, let $D^-:=\hat{\C}\setminus (D^+\cup\Gamma),$
				If $\Gamma$ is not closed
				complement $\Gamma$ to a closed arc $\partial D^+(\supset\gamma)$, enclosing a bounded domain $D^+$, let 
				$D^-:=\hat{\C}\setminus (D^+\cup\Gamma),$ and extend $\varphi(z)$ to the extension of $\Gamma$ by defining it to be $0$ on that set.
				This induces an orientation on $\gamma$ (for the case when $\Gamma$ is not closed). Let $\psi$ be given by Eqn.(\ref{psiekv}).
				Denote 
				$f_+(z_0):=\lim_{D^+\ni z\to z_0} f(z),$ $f_-(z_0):=\lim_{D^-\ni z\to z_0} f(z).$ 
				Now we have
				\begin{equation}
				\int_\gamma \frac{dt}{t-z}=
				\left\{
				\begin{array}{ll}
				2\pi i & ,z\in D^+\\
				0 & ,z\in D^-\\
				\pi i & ,z\in \gamma
				\end{array}
				\right.
				\end{equation}
				This yields
				\begin{equation}
				\psi_+(z_0):=\lim_{D^+\ni z\to z_0}\left(\frac{1}{2\pi i}\int_\gamma\frac{\varphi(t)dt}{t-z}-\frac{\varphi(z_0)}{2\pi i}\int_\gamma\frac{dt}{t-z}\right)=f_+(z_0)-\varphi(z_0)
				\end{equation}
				\begin{equation}
				\psi_-(z_0):=\lim_{D^-\ni z\to z_0}\left(\frac{1}{2\pi i}\int_\gamma\frac{\varphi(t)dt}{t-z}-\frac{\varphi(z_0)}{2\pi i}\int_\gamma\frac{dt}{t-z}\right)=f_-(z_0)
				\end{equation}
				\begin{equation}
				\psi(z_0):=\left(\frac{1}{2\pi i}\int_\gamma\frac{\varphi(t)dt}{t-z}-\frac{\varphi(z_0)}{2\pi i}\int_\gamma\frac{dt}{t-z}\right)=f(z_0)-\frac{1}{2}\varphi(z_0)
				\end{equation}
				
				By Proposition \ref{gakovpreprop}
				$\psi(t)$ is continuous which gives $\psi^+(z_0)=\psi^-(z_0)=\psi(z_0)$ i.e.\ 
				\begin{equation}
				f^+(z_0)-\varphi(z_0)=f^{-1}(z_0)=f(z_0)-\frac{1}{2}\varphi(z_0)
				\end{equation}
				All together we have proved the following theorem.
				\begin{theorem}[Sokhotzki-Plemelj]\label{plemsohotthm}
				Let $\Gamma$ be a bounded smooth arc in the complex plane. 
				and let $f$ belong to the H\"older class
				$C^\alpha(\Gamma),$ $\alpha \in (0,1).$
				If $\Gamma$ is closed let $D^+$ be the bounded domain with $\partial D^+=\Gamma$, let $D^-:=\hat{\C}\setminus (D^+\cup\Gamma),$
				If $\Gamma$ is not closed
				complement $\Gamma$ to a closed arc $\partial D^+(\supset\gamma)$, enclosing a bounded domain $D^+$, let 
				$D^-:=\hat{\C}\setminus (D^+\cup\Gamma),$ and extend $\varphi(z)$ to the extension of $\Gamma$ by defining it to be $0$ on that set.
				This induces an orientation on $\gamma$ (for the case when $\Gamma$ is not closed). 
				Then 
				\begin{equation}\label{plemsohot}
				f(z):=\frac{1}{2\pi i} \int_{\Gamma} \varphi(\zeta)\frac{d\zeta}{\zeta -z}
				\end{equation}
				exists as a Cauchy principal integral on $\Gamma$ (except possible at the ends of $\gamma$, in the case of a non-closed $\gamma$)
				and
				defines a holomorphic function in $\hat{\C}\setminus \Gamma$ (where $\hat{\C}$ denotes the Riemann sphere $\C\cup \{\infty\}$),
				vanishing at $\infty$. Furthermore, (for any choice of $\partial D^+,$ in the case of a non-closed $\gamma$) there exists boundary values
				\begin{equation}
				f^+(t)=\lim_{z\to t,z\in D^+} f(z),\quad f^-(t)=\lim_{z\to t,z\in D^-} f(z)
				\end{equation}
				such that for $t\in \Gamma$
				\begin{equation}
				f^+(t)=\frac{1}{2}\varphi(t)+f(t),\quad f^-(t)=-\frac{1}{2}\varphi(t)+f(t)
				\end{equation}
				where $\phi(t)$ is understood as the principal value.
				\end{theorem}
				As in Remark \ref{gakovremmet} we note that the result is local.
				\begin{proposition}\label{gakovindexproppen}
				If $\gamma$ is a smooth closed contour and $\varphi$ a H\"older continuous density with index $\alpha\in (0,1)$ 
				and let $f(z)=\int_\gamma \frac{\varphi(t)dt}{t-z}.$ Then the limit functions $f_+(z)$ and $f_+(z)$ are
				H\"older continuous with the same index. 
				\end{proposition}
				\begin{proof}
				Let $\psi$ be given by Eqn.(\ref{psiekv}).				
				For two arbitrary sufficiently close $z_1,z_2$ we estimate $\abs{\psi(z_2)-\psi(z_1)}.$ 
				Set $s=s(z,t)$ be the length of
				the smaller of the two arcs of $\gamma$ with tend points $z_0,z$ and $z_0,t$ respectively. 
				Since $\gamma$ is smooth we have for a positive constant $m>0$
				\begin{equation}
				s(z_1,z_2)\leq m\abs{z_1-z_2}
				\end{equation}
				Remove an arc $l$ from $\gamma$ where $l$ consists of an arc containing $z_1$ and subarcs of $\gamma$ of lenght $2s(z_1,z_2)$ to both sides of $z_1.$
				The end points of $l$ will be denoted $a,b.$		
				We have
				\begin{multline}
				\psi(z_2)-\psi(z_1)=\frac{1}{2\pi i}\int_l \frac{\varphi(t)-\varphi(z_2)}{t-z_2}dt
		-\abs{\frac{1}{2\pi i}\int_l \frac{\varphi(t)-\varphi(z_1)}{t-z_1}}dt +\\
		\frac{1}{2\pi i}\int_{\gamma\setminus l} \left(\frac{\varphi(t)-\varphi(z_2)}{t-z_2}
	-\frac{1}{2\pi i}\int_l \frac{\varphi(t)-\varphi(z_1)}{t-z_1}\right)dt =\\
	\frac{1}{2\pi i}\int_l \frac{\varphi(t)-\varphi(z_2)}{t-z_2}dt
-\frac{1}{2\pi i}\int_l \frac{\varphi(t)-\varphi(z_1)}{t-z_1}dt
+\\
\frac{1}{2\pi i}\int_{\gamma\setminus l} \frac{\varphi(z_1)-\varphi(z_2)}{t-z_1}dt+
\frac{1}{2\pi i}\int_{\gamma\setminus l} \frac{(\varphi(t)-\varphi(z_2))(z_2-z_1)}{(t-z_1)(t-z_2)}dt\\
=:I_1+I_2+I_3+I_4
\end{multline}
According to the proof of Proposition \ref{gakovpreprop} there is a constant $C'>0$ such that
\begin{multline}
\abs{I_2}\leq \frac{1}{2\pi i} \int_l \abs{\frac{\varphi(t)-\varphi(z_1)}{t-z_1}}\abs{dt}
\leq C'\frac{1}{2\pi i} \int_l \abs{\frac{1}{t-z_1}}^{1-\alpha}\abs{dt}\leq\\
	C''\int_0^{C(\abs{z_2-z_2}}\frac{ds}{s^{1-\alpha}}\leq C'''s^\alpha (z_1,z_2)\leq A_1\abs{z_2-z_1}^\alpha
	\end{multline}
	for positive constants $C'',C''',A_1.$
	Analogously we obtain that there exists a constant $A_2>0$ such that
	\begin{equation}
	\abs{I_1}\leq A_2\abs{z_2-z_1}^\alpha
	\end{equation}
	Now we have a positive constant $C_0$ such that
	\begin{equation}
	\abs{I_3}\leq \frac{\abs{\varphi(z_1)-\varphi(z_2)}}{2\pi}\abs{\int_{\gamma\setminus l} \frac{dt}{t-z_1}}\leq \frac{C_0\abs{z_2-z_1}^\alpha}{2\pi}
		\abs{\int_{\gamma\setminus l} \frac{dt}{t-z_1}}
		\end{equation}
		where 
		\begin{equation}
		\int_{\gamma\setminus l} \frac{dt}{t-z_1}=\ln\frac{a-z_1}{b-z_1}
		\end{equation}
		which is bounded for $z_1\in \gamma$ so there is a constant $A_3>0$ such that 
		\begin{equation}
		\abs{I_3}\leq A_3\abs{t-z_1}^\alpha
		\end{equation}
		By the H\"older condition together with Eqn.(\ref{gakovplem43}) we have for some constants $A',A''>0$
		\begin{multline}
		\abs{I_4}\leq A'\frac{\abs{z_2-z_1}}{2\pi}\int_{\gamma\setminus l}\frac{ds}{\abs{t-z_1}\abs{t-z_2}^{1-\alpha}}\abs{dt}\leq\\
		A''\abs{z_2-z_1}\int_{\gamma\setminus l}\abs{t-z_1}^{\alpha-2}\abs{\frac{t-z_1}{t-z_2}}^{1-\alpha} \abs{dt}
		\end{multline}
		However
		$\abs{t-z_1}-\abs{z_1-z_2}\leq \abs{t-z_2}$ and by the choice of $l$ we have $\abs{t-z_2}\geq k\abs{z_2-z_1}=:\delta$ for constant $k>1$
		so that
		\begin{equation}
		\abs{t-t_1}\leq \frac{k-1}{k}\abs{t-z_2}
		\end{equation}
		This implies that for $R:=\max_{t\in\gamma\setminus l}\abs{t-z_1}$ there is a constant $A'''>0$ such that
		\begin{equation}
		\abs{I_4}\leq A'''\left(\frac{k-1}{k}\right)^{1-\alpha}\abs{z_1-z_2}\int_R^\delta r^{\alpha -2}dr
		\end{equation}
		If $\alpha<1$ this implies
		$\abs{I_4}\leq A_4\abs{z_2-z_1}^\alpha$.
		This completes the proof of Proposition \ref{gakovindexproppen}.
		\end{proof}
		
		Next suppose $\varphi(t)$ is a H\"older continuous function on the real axis $\gamma:=\{\im z=0\}$ for finite $t$ and
		such that $\varphi$ takes a definite limit $\varphi(\infty)$ as $t\to \pm \infty.$ The H\"older condition must be 
		supplemented when the point at infinity is included by the requirement that there exists a constant $C>0$ and $\alpha>0$ such that
		\begin{equation}\label{thusassumeekv}
		\abs{\varphi(t)-\varphi(\infty)}<\frac{C}{\abs{t}}^\alpha
			\end{equation} 
			which we thus assume.
			Set $f(z)=\frac{1}{2\pi i}\int_{-\infty}^\infty \frac{\varphi(t)dt}{t-z},$ for $z\notin X_0.$ 
			Now for $N_1>0,$ $N_2<0$
			we have
			\begin{equation}
			\int_{N_1}^{N_2} \frac{\varphi(t)dt}{t-z}= \int_{N_1}^{N_2} \frac{(\varphi(t)-\varphi(\infty))dt}{t-z}+\varphi(\infty)
			\int_{N_1}^{N_2} \frac{dt}{t-z}
			\end{equation}
			where $\int_{N_1}^{N_2} \frac{(\varphi(t)-\varphi(\infty))dt}{t-z}$ is of order $\abs{t}^{-1-\alpha}$ for large $\abs{t}$ 
			thus converges as $N_1\to \infty$ and simultaneously $N_2\to -\infty$ independently.
			However
			\begin{equation}
			\varphi(\infty)
			\int_{N_1}^{N_2} \frac{dt}{t-z}=\ln(N_2-z)-\ln(N_1 -z)=\ln\frac{\abs{N_2-z}}{\abs{N_1 -z}} \pm iv
			\end{equation}
			where $v$ denotes the angle between the straight lines respectively connecting $z$ with $N_1+i0$ and $N_2+i0$ respectively,
			and where $+iv$ ($-iv$) is used in the case that $z$ belongs to the upper (lower) half plane.
			Then $v\to \pi$ as $N_1\to \infty$ and $N_2\to -\infty$ independently, and in the case that we take the limit symmetrically in the sense that $N_1=-N_2=N\to \infty$
			we obtain
			\begin{equation}
			\lim_{N\to \infty}\ln\frac{\abs{N-z}}{\abs{-N -z}}=\ln 1=0
			\end{equation}
			If the improper integral converges then it equals the principal value integral
			and in this case this yields the Sokhotsky-Plemelj formula for the case of the real axis $\gamma=\{\im z=0\}$
			\begin{equation}
			\frac{1}{2\pi i}\lim_{N\to \infty} \int_{\stackrel{t\in \{\im t=0\}}{\re t\in [-N,N]}} 
				\frac{\varphi(t)dt}{t-z}=  \frac{1}{2\pi i}\int_{-\infty}^\infty \frac{\varphi(t)-\varphi(\infty)}{t-z}dt\pm \frac{1}{2\pi i}\pi i \varphi(\infty) 
				\end{equation}
				where $+$ ($-$) is used in the case that $z$ belongs to the upper (lower) half plane, where we interpret the 
				integral in the right hand in terms of a principal value integral.
				If $\varphi(\infty)=0$ then the improper integral 
				\begin{equation}\label{pvvarphi0}
				\int_{-\infty}^\infty \frac{\varphi(t)dt}{t-z}
				\end{equation}
				converges.
				Since the integral
				\begin{equation}
				\int_{-\infty}^\infty \frac{\varphi(t)dt}{(t-z)^2}
				\end{equation}
				is absolutely convergent for all $z$ belonging to the real axis we obtain that $f(z)$ will be holomorphic in the upper and in the lower half plane.
				We denote $f_+:=f|_{\{\im z> 0\}}$ and $f_-:=f|_{\{\im z< 0\}}$. 
				Next consider the case when $z$ belongs to the real axis.
				We may write the principal value integral of Eqn.(\ref{pvvarphi0})		
				\begin{equation}
				\int_{-\infty}^\infty \frac{\varphi(t)dt}{t-z}=\lim_{\stackrel{N\to N}{\epsilon\to 0}}\left(
				\int_{-N}^{t-\epsilon} \frac{\varphi(t)dt}{t-z}
				+\int_{t+\epsilon}^N \frac{\varphi(t)dt}{t-z}
				\right)
				\end{equation}
				in particular
				\begin{equation}
				\int_{-\infty}^\infty \frac{dt}{t-z}=\lim_{\stackrel{N\to N}{\epsilon\to 0}}
				\frac{-\epsilon(N-t)}{(-N-t)\epsilon}
				\end{equation}
				which implies
				\begin{equation}
				\int_{-\infty}^\infty \frac{\varphi(t)dt}{t-z}=\lim_{\epsilon\to 0}\left(
				\int_{-\infty}^{t-\epsilon} \frac{(\varphi(t)-\varphi(\infty))dt}{t-z}
				+\int_{t+\epsilon}^\infty \frac{(\varphi(t)-\varphi(\infty))dt}{t-z}
				\right)
				\end{equation}
				Setting $A(t,z):=\frac{(\varphi(t)-\varphi(\infty))}{t-z}$ we may for a finite $c>0$ write
				\begin{equation}
				\left(\int_{-\infty}^{t-\epsilon} Adt +\int_{t+\epsilon}^{\infty} Adt\right)=
				\left(\int_{-\infty}^{t-c} Adt +\int_{t+c}^{\infty} Adt\right)+
				\int_{t-c}^{t-\epsilon} Adt +\int_{t+\epsilon}^{t+c} Adt
				\end{equation}
				The integrals in the parenthesis of the right hand side are independent of $\epsilon$, 
				and the domain of integration of the last two integrals is bounded thus repeating (for the last two integrals) the arguments of Proposition \ref{pvexistsprop} and Proposition \ref{gakovpreprop} we can verify 
				existence of the principal value integral of Eqn.(\ref{pvvarphi0}) for the case when $z$ belongs to the real axis and that the Sokhotsky-Plemelj formulas remain valid according to
				\begin{equation}\label{hobbobo}
				f_\pm(z)=\pm\frac{1}{2}\varphi(z)+\frac{1}{2\pi i}\int_{-\infty}^\infty\frac{\varphi(t)dt}{t-z} ,\quad \re z\in (-\infty,\infty),\quad \im z=0
				\end{equation}
				where $f_+$ ($f_-$) are the limits of $f(\zeta)$ as $\zeta$ tends to $z$ while remaining in the upper (lower) half plane.
				Now perform the change of variable $z=-\frac{1}{\zeta}$,
				$f^*(\zeta):=f(z)=f(-1/\zeta)$, $\varphi^*(\sigma):=\varphi(-1/\sigma),$
				where $\sigma=1/t$ traverses the real axis in positive direction whenever  
				$z=t$ does.
				We have
				\begin{equation}\label{hoboanvandes}
				f^*(\zeta)=\frac{\zeta}{2\pi i}\int_{-\infty}^\infty\frac{\varphi^*(\sigma)d\sigma}{\sigma(\sigma -\zeta)}=
				\frac{1}{2\pi i}\int_{-\infty}^\infty 
				\frac{\varphi^*(\sigma)d\sigma}{\sigma -\zeta}
				-\frac{1}{2\pi i}\int_{-\infty}^\infty 
				\frac{\varphi^*(\sigma)d\sigma}{\sigma}
				\end{equation}
				where the last integral is a constant.
				Suppose $\varphi^*(\sigma)$ satisfies the H\"older condition at $\sigma=0$ 
				in the sense that there exists a constant $C>0$ and $\beta\in (0,1)$ such that
				\begin{equation}
				\abs{\varphi^*(\sigma_2)-\varphi^*(\sigma_1)}\leq C\abs{\sigma_2-\sigma_1}^\alpha 
				\end{equation}
				so that
				\begin{equation}
				\abs{\varphi(t_2)-\varphi(t_1)}\leq C\abs{\frac{1}{t_2}-\frac{1}{t_2}}^\alpha 
				\end{equation}
				which means that $\varphi$ satisfies the condition of Eqn.(\ref{thusassumeekv}). If $z$ tends to infinity along a path in the upper (lower)
				half plane then $\zeta\to 0$ remaining in the upper (lower) half plane.
				Thus applying Eqn.(\ref{hobbobo}) with the choice of $+$ to the first integral in the right hand side of
				Eqn.(ref{hoboanvandes}) we get
				\begin{equation}
				f_+(\infty)=f^*_+(0)=\frac{\varphi^*(0)}{2}+\frac{1}{2\pi i}\int_{-\infty}^\infty\frac{\varphi^*(\sigma)d\sigma}{\sigma}-
			\frac{1}{2\pi i}\int_{-\infty}^\infty\frac{\varphi^*(\sigma)d\sigma}{\sigma}=\frac{\varphi(\infty)}{2}
		\end{equation}
	And the analogous calculation using the $-$ sign in Eqn.(\ref{hobbobo}) renders
$f_-(\infty)=-\frac{1}{2}\varphi(\infty),$
which verifies Eqn.(\ref{hobbobo})also for the (boundary) point at infinity
\\
\\
Let us also recall the Weierstrass preparation theorem and the Weierstrass division theorem. These involve normalized pseudopolynomials, also called Weierstrass polynomials.
\begin{definition}
	Let $n>1$, $q\in\Z_+,$ and set $z':=(z_1,\ldots,z_{n-1})$.
	A {\em Weierstrass polynomial}\index{Weierstrass polynomial} in $z_n$ of degree $q$ is a holomorphic function $W$ in a neighborhood of the origin
	in $\Cn$ of the form
	\begin{equation}
	W(z',z_n)=z_n^q +\sum_{j=1}^k a_j(z')z_n^{q-j}
	\end{equation}
\end{definition}

For reference we shall also give the definition of the discriminant.
Consider the symmetric polynomial (the subscript $V$ stands for Vandermonde)
\begin{equation}
p_V(w_1,\ldots,w_s):=\Pi_{i<j} (w_i-w_j)^2
\end{equation}
It is a known result (see e.g.\ van der Waerden \cite{vdw66}, parag 63) 
that for a given 
symmetric (i.e.\ it is invariant with respect to interchanging two indeterminants $w_i$ and $w_j$, $i\neq j$)
polynomial with integer coefficients there exists a unique symmetric polynomial
$Q_V(v_1,\ldots,v_s)$ with integer coefficients such that
$p_V(w_1,\ldots,w_s)=Q_V(\sigma_1(w_1,\ldots,w_s),\ldots,\sigma_s(w_1,\ldots,w_s))$ where
the $\sigma_i$ are the elementary symmetric polynomials
\begin{equation}
\begin{array}{ll}
\sigma_1(w_1,\ldots,w_s) & =\sum_{j=1}^s w_j\\
\sigma_2(w_1,\ldots,w_s) & =w_1\left(\sum_{j=2}^s w_j\right) +w_2\left(\sum_{j=3}^s w_j\right)+\cdots +w_s\\
&\vdots\\
\sigma_s(w_1,\ldots,w_s) & =w_1\cdots w_s
\end{array}
\end{equation}

\begin{definition}\index{Discriminant of a pseudopolynomial}
	The {\em discriminant}\index{Discriminant} of a pseudopolynomial of the form $P=w^s+h_1(z)w^{s-1}+\cdots+h_s(z)$
	in $\mathscr{O}(\Omega)[w]$ is defined as the holomorphic function on $\Omega$ given by
	\begin{equation}
	\Delta_P:=Q_V(-h_1,h_2,\ldots,(-1)^sh_s)
	\end{equation}
\end{definition}
Let $t_j$, $j=1,\ldots,s$ be the zeros of the polynomial
$w\mapsto P(w,z)$. Then $\Delta_P=0$ if and only if there exists a pair $j\neq k$ such that
$t_j=t_k.$ Let us look closer at this claim.
The elementary symmetric polynomials are the functions $\sigma_j(z)$ appearing in the expansion
\begin{equation}
\Pi_{i=1}^n (w-z_i)=\zeta^n-\sigma_1(z)w^{n-1}+\cdots +(-1)^n\sigma_n(z)
\end{equation}
Now consider the functions
\begin{equation}
s_k(z):=\sum_{j=1}^n z_j^2,\quad k=1,\ldots,n
\end{equation}
\begin{lemma}\label{taylleabor}
Each $\sigma_j(z),$ $j=1,\ldots, n$ can be written as a polynomial in the $s_1\ldots,s_n.$
\end{lemma}
\begin{proof}
For fixed value of $z=(z_1,\ldots,z_n)$	
set \begin{equation}
\psi(w):=\Pi_{i=1}^n (1-\zeta z_j)=1-\sigma_1\zeta +\sigma_2w^2+\cdots +(-1)^n\sigma_n w^n
\end{equation}
We have $\psi(0)=1$ and $\log \psi(\zeta)$ is holomorphic in a neighborhood of $\zeta=0.$ Then
\begin{equation}
\frac{\psi'(w)}{\psi(w)}=\sum_{i=1}^n \frac{-z}{1-\zeta z_i}=-s_1-s_2w-\cdots -s_n w^{n-1}-\cdots
\end{equation}
where the sum in the right hand side converges in a neighborhood of $\zeta=0$, which gives
$-\psi'(w)=\psi(w)\sum_{k=1}^\infty s_k w^{k-1}$ Comparing powers of $w$ gives
$\sigma_1=s_1,$ $\sigma_2=-\sigma_1s_1+s_2,$ $3\sigma_3=\sigma_2s_1-\sigma_1s_2 +s_3$ and so on.
For $j=1,\ldots,n$the $j$:th equation will express $\sigma_j$ as a polynomial in the $s_1,\ldots,s_j$ and $\sigma_1,\ldots,\sigma_{j-1},$
where the coefficients of these polynomials are independent of $z.$This implies that each 
$\sigma_j$ is a polynomial
in the $s_k$ with coefficients independent of $z.$ This completes the proof.
\end{proof}

\begin{theorem}[Weierstrass preparation theorem]\index{Weierstrass preparation theorem}
	Let $F(z,w)$ be analytic on $\{\abs{z-z_0}<r,\abs{w-w_0}<\rho\}$ for some point $(z_0,w_0)\in \C^{n-1}\times \C$ and let $\rho>0.$
	Suppose
	\begin{equation}\label{weierstrassekv344}
	F(z_0,w_0)=0,\quad F(z_0,w)\not\equiv 0
	\end{equation}
	Then there exists $r'>0, \rho>\rho'>0$ such that 
	on $\{\abs{z-z_0}<r',\abs{w-w_0}<\rho'\}$ the function $F(z,w)$ takes the form
	\begin{equation}
	F(z,w)=G(z,w)\sum_{j=0}^k A_j(z)w^{k-j}
	\end{equation}
	for holomorphic $A_j(z),$ $j=0,\ldots,k$ (the existence of $k\geq 1$
	follows from Eqn.(\ref{weierstrassekv344}) as $F(z_0,w)$ has a zero of order $k$ at $w_0$).
\end{theorem}
\begin{proof}
	Set $\gamma':=\{\abs{w-w_0}=\rho'\}$ with $\rho'$ chosen such that $F(z_0,w)$ has no zeros on 
	$\{\abs{w-w_0}\leq \rho'\}$ except $w_0$. Let $0<\mu :=\min_{\gamma'}\abs{F(z_0,w)}$ and let $r>r'>0$
	be sufficiently small such that
	\begin{equation}
	\abs{F(z,w)-F(z_0,w)}<\mu ,\quad z\in c':=\{\abs{z-z_0}=r'\} 
	\end{equation}
	By Rouch\'e's theorem $F(z_0,w)$ and $F(z,w)=F(z_0,w)+(F(z,w)-F(z_0,w))$ have, for a fixed
	$z\in \{\abs{z-z_0}\leq r'\},$
	the same
	number of zeros in $\{\abs{w-w_0}<\rho'\},$ which we know to be $k$, say $w_1(z),\ldots,w_k(z).$
	Without loss of generality we assume $(z_0,w_0)=(0,0).$	
	Set
	\begin{multline}
	P(z,w):=\Pi_{j=1}^k(w-w_j(z))=w^k -w^{k-1}(\sum_{j=1}^k w_j(z))
	+\\
	w^{k-2}(\sum_{j=2}^k w_{j-1}(z)w_j(z))-\cdots +(-1)^k w_1(z)\cdots w_k(z)=\\
	w^k+A_{k-1}(z)w^{k-1}+A_{k-2}(z)w^{k-2}+\cdots +A_0(z)
	\end{multline}
	By Lemma \ref{taylleabor} 
	each $A_j(z)$ is expressed as a polynomial in the $s_m(z)$.
	If the $s_j(z),$ $j=1,\ldots,k$ are analytic in a neighborhood of $z_0$ then so are the
	$A_l(z),$ $l=0,\ldots,k-1.$
	By the Residue theorem we have (since $F(z,w)$ has zeros at the points $w_1(z),\ldots,w_k(z)$ and no poles)
	\begin{equation}\label{weierstrassekv347}
	\frac{1}{2\pi i}\int_{\gamma'}w^m\frac{\partial_w F(z,w)dw}{F(z,w)}=\sum_{j=1}^k w_j^m(z)=s_m(z),\quad z\in \{\abs{z-z_0}\leq r'\}
	\end{equation}
	$m=1,\ldots,k.$
	Now $F(z,w)$ and $\partial_w F(z,w)$ are analytic on $\{\abs{z-z_0}\leq r'\}$
	for each $w\in \gamma'$ and since
	$\abs{F(z,w)-F(z_0,w)}<\mu\leq F(z_0,w)$ on $\{\abs{z-z_0}\leq r'\}\times \gamma'$
	we have that
	$F(z,w)$ does not vanish on $\{\abs{z-z_0}\leq r'\}\times \gamma'$. 
	Set
	\begin{equation}
	\Phi_m(z,w):=w^m\frac{\partial_w F(z,w)dw}{F(z,w)}
	\end{equation}
	so that Eqn.(\ref{weierstrassekv347}) is
	\begin{equation}\label{weierstrassekv347}
	s_m(z)=\frac{1}{2\pi i}\int_{\gamma'}\Phi_m(z,w),\quad z\in \{\abs{z-z_0}\leq r'\}
	\end{equation}
	$m=1,\ldots,k.$
	Note that  
	$\Phi_m(z,w)$ is continuous on the closed set
	$\{\abs{z-z_0}\leq r'\}\times\gamma'$ which yields uniform boundedness 
	on $\{\abs{z-z_0}\leq r'\}$,
	of the family $\{\Phi_m(z,w):w\in \gamma'\}$.
	\begin{lemma}\label{i1719}
		Let $\Omega\subset\C$ be a domain in the $z$-plane and
		let $L$ be a rectifiable curve in the $w$-plane.
		Let $\Phi(z,w)$ be a function that for each $w\in L$ is analytic in $\Omega$ and 
		for each $z\in \Omega$ is
		continuous with respect to $w$ on $L$. Suppose further that
		the family $\{\Phi(z,w):w\in L\}$ is uniformly bounded in $\Omega.$
		Then 
		\begin{equation}
		s(z)=\int_L\Phi(z,w)dw
		\end{equation}
		is analytic on $\Omega$. 
	\end{lemma}
	\begin{proof}
		Let $L$ be given by the equation $w=\lambda(t),$ $t\in [\alpha,\beta].$ Let $\phi^{(n)}:=\{t_0^{(n)},\ldots,t_n^{(n)}\},$
		$\alpha=t_0^{(0)}<\cdots <t_n^{(n)}=\beta,$
		$n\in \Z_+$ be a sequence of partitions of $[\alpha ,\beta]$, such that
		the norm of $\phi^{(n)}$ approaches zero as $n\to \infty.$ Set
		$w^{(n)}:=\lambda(t_k^{(n)}),$ $k\in \Z_+$ and set
		\begin{equation}
		f_n(z):=\sum_{k=1}^n \Phi(z,w_k^{(n)})(w_k^{(n)} -w_{k-1}^{(n)})
		\end{equation}
		For any compact $K\subset\Omega$ we have
		\begin{equation}
		\abs{f_n(z)}\leq C_K\sum_{k=1}^n \abs{w_k^{(n)} -w_{k-1}^{(n)}}\leq C_K \mbox{length}(L),\quad z\in K
		\end{equation}
		where $C_k:=\sup_{z\in K,w\in L}\abs{\Phi(z,w)}.$
		Thus $\{f_n\}_{n\in \Z_+}$ is uniformly bounded on $\omega$.
		By Montel's theorem it follows that $\{f_n(z)\}_{n\in \Z_+}$ is compact on $\Omega.$ Since $L$
		is rectifiable and $\Phi(z,w)$ is continuous with respect to $w$ on $L$ for each $z\in \omega$ we have that $f_n(z)$ converges, as $n\to \infty$, on $\Omega$
		the function $s(z)=\int_L\Phi(z,w)dw$. By Theorem \ref{vitaliporterthm} (Vitali's theorem) we may assume
		that $f_n(z)\to s(z)$ uniformly on $\Omega$ as $n\to \infty.$ A uniformly convergent sequence of holomorphic functions on a domain has 
		a limit function $s(z)$ which is holomorphic (see e.g.\ Theorem \ref{hormkonvthm}). This completes the proof
		of Lemma \ref{i1719}.
	\end{proof}
	By Lemma \ref{i1719}, with $\Omega$ replaced by $\{\abs{z-z_0}\leq r'\}$ and $L$ replaced by $\gamma'$, we have that 
	$s_m(z,w)$ is analytic on $\{\abs{z-z_0}\leq r'\},$ and the same is true for $A_j(z),$ $m=1,\ldots ,k$. Hence $P(z,w)$ is a Weierstrass polynomial of order $k$ with respect to $w$ and vanishes exactly on the same set as $F.$ 
	Set
	\begin{equation}
	G(z,w):=\frac{F(z,w)}{P(z,w)}
	\end{equation}
	This is defined and holomorphic in $\{\abs{z-z_0}< r'\}\times \{\abs{w-w_0}< \rho'\}$ outside the zero set of $F$ and $P$. Moreover for fixed $z$,
	$G(z,w)$ has only removable singularites in the disc $\abs{w-w_0}<\rho'$
	so that $G$ can be extended to a function in $\{\abs{z-z_0}< r'\}\times \{\abs{w-w_0}< \rho'\}$ that is analytic with respect to $w$ for each fixed $z$ as well as in the complement of the zero set.
	Furthermore, $G$ is bounded away from $0$ on $\{\abs{z-z_0}< r'\}\times \{\abs{w-w_0}= \rho'\}$. 
This together with the maximum principle implies that $F/P$ is bounded. 
	Writing
	\begin{equation}
	G(z,w)=\frac{1}{2\pi i}\int_{\abs{zeta}=\rho'}\frac{G(z,\zeta)d\zeta}{\zeta -w}
	\end{equation}
	we have that $G(z,w)$ is holomorphic with respect to $z$ as well.
This completes the proof.
\end{proof}	

\begin{theorem}[Weirerstrass's division theorem]\label{weierdivsion}
	Let $f$ be holomorphic in a neighborhood of the origin in $\Cn$ and let $W$ be an arbitrary Weierstrass polynomial
	in $z_n$ of degree $q.$ Then
	$f$ has a unique representation near $0$ of the form $f=QW+R$ where $Q$ is holomorphic and $R$ is a pseudopolynomial in $z_n$
	of degree $<q.$
\end{theorem}
\begin{proof}
	If $f$ and $W$ are holomorphic in $\{\abs{z_j}<r_j,j=1,\ldots,n\}$ for some $r\in \R_+^n$
	let $\delta$ be such that $\delta_j<r_j$, $j=1,\ldots,n,$ be such that $W(z)\neq 0$ for $z'\in \{\abs{z_j}<\delta_j,
	j=1,\ldots,n-1\}$ and $\abs{z_n}=\delta_n.$ Then one can verify that the function
	\begin{equation}
	Q(z):=\frac{1}{2\pi i}\int_{\abs{z_n}=\delta_n\}} \frac{f(z',\zeta)}{W(z',\zeta)}\frac{d\zeta}{\zeta-z_n}
	\end{equation}
	for $z\in \{\abs{z_j}<\delta_j, j=1,\ldots,n\},$
	is holomorphic and that the same holds true for $R(z):=f(z)-Q(z)W(z).$
	$R(z)$ is realized to be a pseudopolynomial in $z_n$ of degree $<q$ since by the formula
	\begin{equation}
	f(z)-Q(z)W(z)=\frac{1}{2\pi i}\int \frac{f(z',\zeta)}{W(z',\zeta)}\frac{W(z',\zeta)-W(z',z_n)}{\zeta-z_n}d\zeta
	\end{equation} 
	for $z\in \{\abs{z_j}<\delta_j, j=1,\ldots,n\},$ shows that it is a pseudopolynomial of degree strictly less than the degree of $W.$ If $f=Q_j W+R_j,$ $j=1,2.$ then $(Q_1-Q_2)W\equiv (R_1-R_2),$ hence uniqueness follows from uniqueness of holomorphic functions
	on open subsets of $\C.$
\end{proof}

\end{appendix}
\scriptsize

\bibliographystyle{amsplain}

\small
\printindex

\end{document}